\documentclass[reqno]{amsart}
\pdfoutput=1 

\usepackage[utf8]{inputenc}
\usepackage{Wave_Maps}
\usepackage{Wave_Maps_Paraproducts}

\usepackage{hyperref}

\title[Invariant Gibbs measures for $(1+1)$-dimensional wave maps into Lie groups]{Invariant Gibbs measures for $(1+1)$-dimensional wave maps \\ into Lie groups} 
\author[Bjoern Bringmann]{Bjoern Bringmann \vspace{-3ex}}
\address{\nopagebreak Bjoern Bringmann, Department of Mathematics, Princeton University, Princeton, NJ 08544}
\email{bringmann@princeton.edu}

\begin{document}

\maketitle

\vspace{-2ex}
\begin{abstract}
We discuss the $(1+1)$-dimensional wave maps equation with values in a compact Lie group. The corresponding Gibbs measure is given by a Brownian motion on the Lie group, which plays a central role in stochastic geometry. Our main theorem is the almost sure global well-posedness and invariance of the Gibbs measure for the wave maps equation. It is the first result of this kind for any geometric wave equation. 

\noindent Our argument relies on a novel finite-dimensional approximation of the wave maps equation which involves the so-called Killing renormalization. The main part of this article then addresses the global convergence of our approximation and the almost invariance of the Gibbs measure under the corresponding flow.
The proof of global convergence requires a  carefully crafted Ansatz which includes modulated linear waves, modulated bilinear waves, and mixed modulated objects. 
The interactions between the different objects in our Ansatz are analyzed using an intricate combination of analytic, geometric, and probabilistic ingredients. In particular, geometric aspects of the wave maps equation are utilized via orthogonality, which has previously been used in the deterministic theory of wave maps at critical regularity. The proof of almost invariance of the Gibbs measure under our approximation relies on conservative structures, which are a new framework for the approximation of Hamiltonian equations, and delicate estimates of the energy increment.
\end{abstract}
\vspace{-1.25ex}

\tableofcontents

\section{Introduction}\label{section:introduction}

\begin{figure}
    \centering
    \captionsetup[subfigure]{labelformat=simple}
\begin{subfigure}[t]{0.45\textwidth}    
    \centering
    \includegraphics[width=0.8\textwidth]{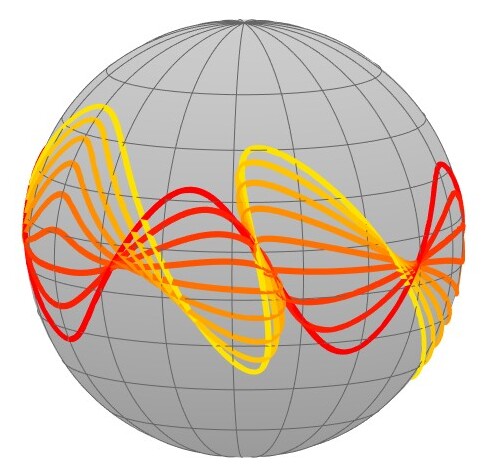}
    \caption{\small{Wave maps}}
    \label{fig:WM}
\end{subfigure}
\hspace{4ex}
\begin{subfigure}[t]{0.45\textwidth}
    \centering
    \includegraphics[width=0.8\textwidth]{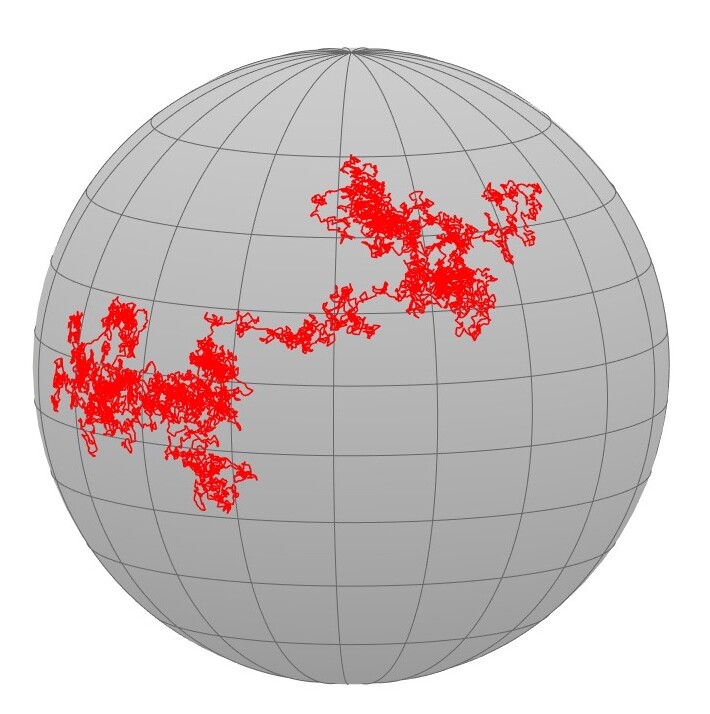}
    \caption{\small{Brownian motion}}
    \label{fig:BM}
\end{subfigure}
 \captionsetup{width=0.95\linewidth}
\caption{\small{In  $(a)$, we illustrate a wave map $\phi$ into $\mathbb{S}^2$ with smooth initial data, which is shown at nine different times. 
    The initial position of $\phi$, which is colored in red, can be written as $\phi(x)=(\sin(\theta(x)) \cos(x),\sin(\theta(x)) \sin(x),\cos(\theta(x))$, where  $\theta(x)=\frac{\pi}{2}+\frac{\pi}{8}cos(4x)$. The initial velocity $\partial_t \phi$ is given by the projection of $2 cos(4x) (1,0,1)$ onto the tangent space $T_{\phi(x)}(\mathbb{S}^2)$. 
    In $(b)$, we illustrate a sample path of a Brownian motion on $\mathbb{S}^2$.}}
\label{fig:sphere}
\end{figure}

The wave maps equation is the epitome of a geometric wave equation. In order to formulate the wave maps equation, we consider the Minkowski space $\R_t\times \R_x^d$ with signature $(- \, + \, \hdots + )$ and consider a compact Riemannian manifold $(\scrM,h)$. The wave maps equation for a map $\phi\colon \R_t \times \R_x^d \rightarrow \scrM$ is then given by
\begin{equation}\label{intro:eq-WM-covariant}
(\phi^\ast \nabla)^\mu \partial_\mu \phi =0. 
\end{equation}
In \eqref{intro:eq-WM-covariant}, $\nabla$ is the Levi-Civita connection on $(\scrM,h)$, $\phi^\ast$ is the pull-back map on $T\scrM$ induced by $\phi$, and Greek indices are raised with respect to the Minkowski metric. In local coordinates, the wave maps equation \eqref{intro:eq-WM-covariant} can be written as 
\begin{equation}\label{intro:eq-WM-M}
\partial_\mu \partial^\mu \phi^k = \Gamma^k_{ij}(\phi) \partial_\mu \phi^i \partial^\mu \phi^j,
\end{equation}
where $(\Gamma^k_{ij})$ are the Christoffel symbols of the Riemannian manifold $(\scrM,h)$. 
The initial data for the wave maps equation consists of the initial position $\phi_0 \colon \R_x^d\rightarrow \scrM$ and the initial velocity $\phi_1 \in \phi_0^\ast T\scrM$. The wave maps equation \eqref{intro:eq-WM-M} is the Euler-Lagrange equation of the Lagrangian
\begin{equation}\label{intro:eq-Lagrangian}
\mathscr{L}(\phi) := \int\displaylimits_{\R\times \R^d}   h_{\phi(t,x)} \big( \partial_\mu \phi(t,x), \partial^\mu \phi (t,x) \big)\, \dt \dx,
\end{equation}
and can therefore be seen as the natural generalization of the scalar-valued linear wave equation to the manifold-valued setting. 
Our main focus lies on $(1+1)$-dimensional wave maps, which are of particular significance in differential geometry and general relativity. In the $(1+1)$-dimensional setting, a wave map $\phi\colon \R_t \times \R_x \rightarrow \scrM$ describes the evolution of a string in the manifold $\scrM$. In the classical literature \cite{G80,GV82,KT98,LS81,MNT10,P76}, the initial position $\phi_0\colon \R_x \rightarrow \scrM$ and initial velocity $\phi_1^\ast \in \phi_0^\ast T \scrM$ are often taken as deterministic and smooth. For an illustration of wave maps with smooth initial data, we refer the reader to Figure \ref{fig:WM}. In this article, our main goal is to understand the probabilistic aspects of wave maps. Due to this, we study the evolution of a random string under the wave maps equation,  which requires us to go well beyond the classical smooth setting.

The most natural model for a random string on the manifold $\scrM$ is a Brownian motion (see Figure \ref{fig:BM}), which can be defined as the diffusion process induced by the Laplace-Beltrami operator. Brownian motions on manifolds were first studied by It\^{o} \cite{I50}, who considered Brownian motions on Lie groups, and are a central object at the interface of differential geometry and probability theory. In addition to their intrinsic significance, Brownian motions on manifolds also have important applications to heat kernel estimates and index theorems (see e.g. \cite{H02}). In the context of $(1+1)$-dimensional wave maps, their particular significance stems from the relationship with the Gibbs measure of the wave maps equation. At a formal level, the Gibbs measure at inverse temperature $\beta>0$ is given by 
\begin{equation}\label{intro:eq-Gibbs}
``\mathrm{d}\mu_\beta(\phi_0,\phi_1)= \mathcal{Z}^{-1} 
\exp\bigg( - \frac{\beta}{2} \int_{\R} \Big( 
\big| \partial_x \phi_0 \big|_{h(\phi_0)}^2 
+ \big| \phi_1 \big|^2_{h(\phi_0)} \Big) \mathrm{d}x \bigg) \mathrm{d}\phi_0 \mathrm{d}\phi_1",
\end{equation}
where $\phi_0$ and $\phi_1$ correspond to the initial data of \eqref{intro:eq-WM-M}. The links between the Gibbs measure $\mu_\beta$ and Brownian motions on $\scrM$ were first explored by Andersson and Driver \cite{AD99}, who showed\footnote{In \cite{AD99}, the authors only consider the initial position $\phi_0$ and do not consider the initial velocity $\phi_1$. While the results of \cite{AD99} therefore do not directly apply to wave maps, our main theorem only concerns the special case of Lie groups, in which the rigorous construction of the Gibbs measure is more elementary (see Definition \ref{lifting:def-Gibbs}).} that the first marginal of \eqref{intro:eq-Gibbs} is absolutely continuous with respect to the law of a (rescaled) Brownian motion on $\scrM$. 
Equipped with both the Gibbs measure and the wave maps equation, we arrive at the following two questions: 
\begin{enumerate}[label=(\Roman*)]
\item If the initial data is drawn from the Gibbs measure \eqref{intro:eq-Gibbs}, does the wave maps equation \eqref{intro:eq-WM-M} have a local solution? 
\item If so, is the solution of the wave maps equation \eqref{intro:eq-WM-M} global and is the Gibbs measure invariant under the global dynamics? 
\end{enumerate}
The first question has been answered in the affirmative in earlier work of the first author, L\"{u}hrmann, and Staffilani \cite{BLS21}. At least in the special case in which $\scrM$ is a compact Lie group with a bi-invariant Riemannian metric, such as the special orthogonal group $SO(n)$, we now answer the second question. 

\begin{theorem}[Informal version]\label{intro:thm-formal}
Let $\scrM$ be a compact Lie group with a bi-invariant Riemannian metric and let $\beta>0$ be the inverse temperature. Then, the wave maps equation \eqref{intro:eq-WM-M} is almost-surely globally well-posed for initial data drawn from the Gibbs measure $\mu_\beta$. Furthermore, the Gibbs measure $\mu_\beta$ is invariant under the wave maps equation. 
\end{theorem}

While we postpone a detailed discussion of Theorem \ref{intro:thm-formal} until Sections \ref{section:introduction-result} and \ref{section:overview}, we make three remarks. 
\begin{enumerate}[label=(\roman*)]
    \item In recent years, the invariance of Gibbs measures has been shown for several \emph{scalar} wave and Schr\"{o}dinger equations \cite{B94,B96,B23,BDNY22,DNY19,OOT21}. In contrast, our main theorem concerns a \emph{geometric} wave equation, and it is the first result of this kind. 
    \item In \cite{BLS21}, the first author, L\"{u}hrmann, and Staffilani proved the almost-sure local well-posedness for initial data drawn from the Gibbs measure (in the case of a general compact Riemannian manifold $\scrM$). In Theorem \ref{intro:thm-formal}, we strengthen local well-posedness to global well-posedness and invariance, which requires several novel ingredients. This passage from local to global well-posedness and invariance has already proven itself to be challenging even in the context of geometric stochastic  \emph{parabolic} equations. For examples of this, we refer the reader to the open problem for the geometric stochastic heat equation in \cite[Conjecture 4.5]{BGHZ22} and the recent breakthrough for the stochastic Yang-Mills equation in \cite{CS23}. 
    \item While Theorem \ref{intro:thm-formal} only covers compact Lie groups, most of our argument also applies to general compact Riemannian manifolds. However, there are certain limitations that currently prohibit us from treating the general case, and we refer to Remark
    \ref{intro:rem-general-case} for a more detailed discussion. 
\end{enumerate}

\subsection{Motion of deterministic and random strings}\label{section:intro-strings}

In this article, we are interested in the evolution of random strings in a Riemannian manifold, which can be described using two natural and prominent models: The $(1+1)$-dimensional wave maps equation with random initial data (as in Theorem \ref{intro:thm-formal}) and the one-dimensional geometric stochastic heat equation (as in \cite{BGHZ22,H16}). In the following, we focus entirely on these two models and postpone a broader discussion until Subsection \ref{section:introduction-general}. 

\subsubsection{Wave maps in $(1+1)$-dimensions}\label{section:intro-1d-wave-maps}

In the physical literature, $(1+d)$-dimensional wave maps are commonly referred to as $\sigma$-models, which were first introduced in \cite{GL60}. Since then, $\sigma$-models have played a fundamental role in particle physics and we refer the reader to \cite{S62,W88,Z89} and the references therein.  As previously mentioned, the $(1+1)$-dimensional wave maps equation is particularly relevant in differential geometry and general relativity. 
In differential geometry, its significance stems from the local equivalence of the following three models: 
\begin{enumerate}[label=(\roman*)]
    \item The Gauss-Codazzi equation for surfaces in $\mathbb{R}^3$ with constant negative Gaussian curvature, 
    \item the $(1+1)$-dimensional sine-Gordon equation,
    \item and the $(1+1)$-dimensional wave maps equation with target manifold $\scrM=\mathbb{S}^2$. 
\end{enumerate}
For detailed discussions of the relationships between the three models, we refer the reader to \cite{BS13,SS96,TU00,TU04}. As a result of this equivalence, $(1+1)$-dimensional wave maps into $\mathbb{S}^2$ correspond to solutions of one of the oldest problems in differential geometry, i.e., constructions of surfaces with constant negative Gaussian curvature. 
In general relativity, the central objects of interest are space-times, which are $(1+3)$-dimensional Lorentzian manifolds. A special class of spacetimes is given by the Gowdy vacuum spacetimes, which obey certain symmetry conditions. In the setting of Gowdy vacuum spacetimes, the Einstein vacuum equations can be reduced to the $(1+1)$-dimensional wave maps equation \cite{CB99,N07,R04}. This reduction has been used successfully to study the asymptotic behavior of Gowdy vacuum spacetimes.  \\

The deterministic theory of the $(1+1)$-dimensional wave maps equation has been studied in \cite{G80,GV82,KT98,LS81,MNT10,P76,Z99} and is rather well-understood. Local well-posedness in $H^s$-spaces was first shown for $s\geq 2$ in \cite{GV82,LS81}, $s\geq 1$ in \cite{Z99}, $s>\frac{3}{4}$ in \cite{KT98}, and finally for $s>\frac{1}{2}$ in \cite{MNT10}, which covers the entire scaling-subcritical regime. Due to the conserved energy of \eqref{intro:eq-WM-M}, which controls the $\dot{H}^1$-norm of the solution, the local theory directly implies the global well-posedness of \eqref{intro:eq-WM-M} in the energy space. Using the I-method, this has been extended to global well-posedness in $H^{s}$ for all $s>\frac{3}{4}$ in \cite{KT98}. In comparison, the deterministic well-posedness theory of the $(1+d)$-dimensional wave maps equation in dimension $d\geq 2$ is much more involved, and we further discuss this in Subsection \ref{intro:sec-high-dimensional-wm} below.
At and below the regularity $s=\frac{1}{2}$, deterministic well-posedness of \eqref{intro:eq-WM-M} breaks down entirely. In \cite{T00}, Tao showed that \eqref{intro:eq-WM-M} is ill-posed in both the critical Sobolev space $\dot{H}^{\frac{1}{2}}$ and the critical Besov space $\dot{B}^{\frac{1}{2}}_{2,1}$. The ill-posedness persists\footnote{In \cite{T00}, Tao showed norm-inflation in $\dot{H}^{\frac{1}{2}}$ and $\dot{B}^{\frac{1}{2}}_{2,1}$, which is a strong form of ill-posedness. In contrast, \cite{BLS21} only showed the unboundedness of the first Picard-iterate in $\C^s$, which is a weaker form of ill-posedness. Nevertheless, it rules out any proof of well-posedness in $\C^s$ using contraction-mapping arguments.}
even under stronger integrability conditions and, as shown in \cite[Theorem 1.2]{BLS21}, also occurs in the H\"{o}lder spaces $\C^{s}$ for $0<s\leq \frac{1}{2}$.

Compared to the deterministic theory of \eqref{intro:eq-WM-M}, the probabilistic theory of \eqref{intro:eq-WM-M} is much less understood. The most natural probabilistic questions for \eqref{intro:eq-WM-M} concerns the local and global well-posedness of \eqref{intro:eq-WM-M} on the support of the Gibbs measure and the invariance of the Gibbs measure. Since the Gibbs measure corresponds to Brownian motion, it lives at regularities $s<\frac{1}{2}$, and at such regularities the wave maps equation is both ill-posed and has no known conservation laws. Due to the absence of nonlinear smoothing in \eqref{intro:eq-WM-M}, 
 bridging the gap between the regularity of Brownian motion and the regularity of the deterministic well-posedness theory turns out to be a formidable challenge. As mentioned above, the author, L\"{u}hrmann, and Staffilani \cite{BLS21} previously obtained the local well-posedness of \eqref{intro:eq-WM-M} with Brownian initial data. As it requires more details, we postpone a further discussion of the ideas in \cite{BLS21} until Section \ref{section:overview}. In the case $\mathcal{M}=\mathbb{S}^d$, Brze\'{z}niak and Jendrej \cite{BJ22} introduced a lattice discretization of \eqref{intro:eq-WM-M} and examined the evolution of Brownian initial data under the discretized equation. 
 Using compactness arguments, Brze\'{z}niak and Jendrej showed that a subsequence of the solutions converges in law as the lattice spacing tends to zero. Furthermore, they proved that the law of the corresponding limit is invariant under time-like translations. This statement is quite different from our main result (Theorem \ref{intro:thm-rigorous-phi}),  which entails a strong rather than weak solution theory.

\subsubsection{Geometric stochastic heat equations in one dimension}

The parabolic counterpart of the wave maps equation is the harmonic map heat flow, which is a central object in differential geometry and partial differential equations \cite{ES64}. In the following discussion, we focus on the one-dimensional, periodic setting. For a map $\phi \colon [0,\infty)\times \mathbb{T}\rightarrow \scrM$, the harmonic map heat flow can be written in local coordinates as
\begin{equation}\label{intro:eq-harmonic-map-heat}
\partial_t \phi^k = \partial_x^2 \phi^k + \Gamma^k_{ij}(\phi) \partial_x \phi^i \partial_x \phi^j. 
\end{equation}
It is the $L^2$-gradient flow associated with the energy
\begin{equation}\label{intro:eq-energy}
E(\phi) = \frac{1}{2} \int_{\bT} h_{\phi(x)}\big( \partial_x \phi(x), \partial_x \phi(x) \big) \dx, 
\end{equation}
and can therefore be seen as the generalization of the scalar-valued, linear heat equation to the manifold-valued setting. For deterministic and smooth initial data $\phi_0 \colon \bT \rightarrow \scrM$, the local and global well-posedness of \eqref{intro:eq-harmonic-map-heat} is classical \cite{ES64,O85}. Similar as for the wave maps equation, it is possible to study \eqref{intro:eq-harmonic-map-heat} with random initial data, such as a Brownian loop on the manifold $\scrM$. However, since the harmonic map heat flow \eqref{intro:eq-harmonic-map-heat} is smoothing, many of the properties of Brownian loops are destroyed under the corresponding flow. To remedy this, it is natural to add a stochastic forcing term in \eqref{intro:eq-harmonic-map-heat}. The resulting model is called the geometric stochastic heat equation \cite{BGHZ22,H16} and can be written in local coordinates as\footnote{In \cite{BGHZ22}, the geometric stochastic heat equation \eqref{intro:eq-geometric-stochastic-heat} also includes a vector-field $h^k(\phi)$, but we omit it here.}
\begin{equation}\label{intro:eq-geometric-stochastic-heat}
\partial_t \phi^k = \partial_x^2 \phi^k + \Gamma^k_{ij}(\phi) \partial_x \phi^i \partial_x \phi^j + \sigma^k_\ell(\phi) \xi^\ell. 
\end{equation}
The coefficients $(\sigma^k_\ell)$ are related to the Riemannian metric (see \cite[(1.3)]{BGHZ22}) and the distributions $(\xi^\ell)$ are space-time white noises. The geometric stochastic heat equation \eqref{intro:eq-geometric-stochastic-heat} can be interpreted as the Langevin equation corresponding to the energy \eqref{intro:eq-energy} and is the natural generalization of the scalar-valued linear heat equation with additive space-time white noise to the manifold-valued setting. \\

Since the space-time white noises $(\xi^\ell)$ have parabolic regularity $-\frac{3}{2}-$, the parabolic regularity of the solution $\phi$ is at most $\frac{1}{2}-$. Thus, neither  $\Gamma^k_{ij}(\phi) \partial_x \phi^i \partial_x \phi^j$ nor $\sigma^k_\ell(\phi) \xi^\ell$ can be defined using only information on the parabolic regularity of $\phi$, and the local theory of \eqref{intro:eq-geometric-stochastic-heat} therefore has to rely on the random structure of $\phi$. The random structure of $\phi$ can be described using regularity structures \cite{BCCH21,BHZ19,CH16,H14}, which form a general framework for the local well-posedness of singular, stochastic partial differential equations. Instead of a single solution of \eqref{intro:eq-geometric-stochastic-heat}, however, regularity structures yield a family of solutions of renormalized versions of \eqref{intro:eq-geometric-stochastic-heat}, which is indexed by a $54$-dimensional renormalization group. In \cite[Theorem 1.6]{BGHZ22}, it is shown that a unique solution from this family can be selected by imposing natural conditions on the solution theory, such as equivariance under the action of the diffeomorphism group, an It\^{o} isometry, and a minimality condition (cf. \cite[Section 1]{BGHZ22}).  \\

While \cite{BGHZ22} settles the local well-posedness of \eqref{intro:eq-geometric-stochastic-heat}, the global well-posedness and invariance of Brownian loops for a renormalized version of \eqref{intro:eq-geometric-stochastic-heat} are still open and are discussed in \cite[Section 4.3]{BGHZ22}. One of the main difficulties, which is shared by Theorem \ref{intro:thm-formal}, is that a proof of invariance likely requires a finite-dimensional approximation of \eqref{intro:eq-geometric-stochastic-heat}, such as a lattice discretization. However, the BPHZ-theorem from \cite{CH16} is currently only available in the continuous setting, and can therefore not be applied to a lattice discretization of \eqref{intro:eq-geometric-stochastic-heat}. In simplified settings, such as when the Riemannian manifold $\scrM$ is given by a Lie group $\frkG$, it is likely possible to circumvent this problem. The reason is that the general results of \cite{CH16} can then likely be obtained by hand, which would be practically infeasible for general Riemannian manifolds $\scrM$. A similar approach has been used for a certain coupled KPZ equation in \cite{FH17}, and can likely be extended to geometric stochastic heat equations taking values in a Lie group $\frkG$. 

\subsection{A general class of deterministic and random models}\label{section:introduction-general}
We previously discussed the $(1+1)$-dimensional wave maps equation and one-dimensional geometric stochastic heat equation as individual models. We now take a different perspective and view them as part of a more general class of models. In order to avoid technicalities, our discussion of this general class will be purely formal. 
We let $\introstate$ be a given state space and let $\introcotangent$ be its cotangent bundle. Furthermore, we let  $\introfunctional \colon \introstate \rightarrow \R \medcup \{ \infty\}$ be a functional and $\introhamiltonian \colon \introcotangent  \rightarrow \R \medcup \{ \infty\}$ be a Hamiltonian. In many settings, the Hamiltonian is essentially of the form
\begin{equation*}
\introhamiltonian(\phi,\pi) = \introfunctional(\phi) + \tfrac{1}{2} \langle \pi, \pi \rangle_\phi, 
\end{equation*}
where $\langle \cdot, \cdot \rangle_\phi$ is an inner product on the fiber $\introcotangentphi$ of $\introcotangent$. For this reason, we phrase all examples below in terms of the functional $\introfunctional$. Equipped with $\introfunctional$ and $\introhamiltonian$, we can write the corresponding Langevin equation as 
\begin{equation}\label{intro:eq-Langevin} 
\dot{\phi} = - (\nabla \introfunctional)(\phi) + \zeta,
\end{equation}
where $\zeta$ denotes a stochastic forcing term\footnote{As can be seen from \eqref{intro:eq-energy}, the precise form of \eqref{intro:eq-Langevin} can be rather complicated. The reason is that both $\nabla$ and $\zeta$ depend on the metric structure of $\introstate$, which can itself be rather complicated.}. The corresponding Hamiltonian equation can be written as 
\begin{equation}\label{intro:eq-Hamiltonian}
    \dot{\phi} = (\nabla_\pi \introhamiltonian )(\phi,\pi), 
\hspace{12ex} 
\dot{\pi} = - (\nabla_\phi \introhamiltonian)(\phi,\pi). 
\end{equation}
While both the Langevin equation \eqref{intro:eq-Langevin} and Hamiltonian equation \eqref{intro:eq-Hamiltonian} can be formulated in this general setting, their analysis depends heavily on the choice of the state space $\introstate$ and the functional $\introfunctional$. For the wave maps equation and geometric stochastic heat equation in dimension $d\geq 1$, the state space consists of maps $\phi \colon \R^{d} \rightarrow \scrM$ and the functional $\introfunctional$ is given by
\begin{equation}\label{intro:eq-functional-sigma}
\introfunctional(\phi) = \int_{\R^d} \langle \partial_j \phi, \partial^j \phi \rangle_{h(\phi)} \, \dx. 
\end{equation}
In addition to the wave maps and geometric stochastic heat equations, however, many other important models\footnote{While the Langevin and Hamiltonian equations are often studied in the period setting, where $\R^d$ is replaced by $\mathbb{T}^d$, we do not emphasize this difference in our literature overview.} can be written in the form \eqref{intro:eq-Langevin} or \eqref{intro:eq-Hamiltonian}. For $\Phi^{4}_{d}$-models, the state space consists of functions $\phi \colon \R^d \rightarrow \R$  and the corresponding $\Phi^{4}_d$-functional is given by 
\begin{equation}\label{intro:eq-functional-Phi}
\introfunctional(\phi) = \int_{\R^d}  \bigg( \frac{|\nabla \phi|^2}{2} + \frac{|\phi|^{4}}{4} \bigg) \dx. 
\end{equation}
The corresponding Langevin and (real-valued) Hamiltonian equations are then given by the cubic stochastic heat equation and cubic wave equation, respectively. In Yang-Mills theory, the state space consists of connection one-forms $A\colon \R^d \rightarrow \frkg^d$, where $\frkg$ is a Lie algebra, and the Yang-Mills functional is given by
\begin{equation}\label{intro:eq-functional-YM}
\introfunctional(A) =  \int_{\R^d}  \big\langle (F_A)_{ij}, (F_A)^{ij} \big\rangle_{\frkg} \dx. 
\end{equation}
In \eqref{intro:eq-functional-YM}, $F_A$ denotes the curvature tensor corresponding to $A$. The corresponding Langevin and Hamiltonian equation are then given by the stochastic Yang-Mills flow and hyperbolic Yang-Mills equations, respectively. For detailed discussions of the Langevin and Hamiltonian equations for the three different models, we refer the reader to 
\cite{BGHZ22,C22,GH21} and \cite{BDNY22,KTV14,OT19,Tao06}, respectively. In the rest of this section, we focus our attention on the following four instances:
\begin{enumerate}[label=(\roman*)]
    \item Deterministic wave maps, 
    \item scalar-valued random dispersive equations,  
    \item stochastic Yang-Mills equations, 
    \item and random geometric wave equations.  
\end{enumerate}

\subsubsection{Deterministic aspects of wave maps in $(1+d)$-dimensions}\label{intro:sec-high-dimensional-wm}

Since the case $d=1$ was previously discussed in Subsection \ref{section:intro-1d-wave-maps}, we now focus on the case $d\geq 2$. The local well-posedness of the wave maps equation \eqref{intro:eq-WM-M} in the Sobolev spaces $(H^s\times H^{s-1})(\mathbb{R}^d)$ has now been established at all scaling sub-critical regularities $s>\frac{d}{2}$. At the higher regularities $s>\frac{d}{2}+1$ or $s>\frac{d}{2}+\frac{1}{2}$, the local well-posedness of \eqref{intro:eq-WM-M} can be obtained using Sobolev embedding or the $L_t^2L_x^\infty$-Strichartz estimate, respectively. Both arguments make no use of the geometric structure of \eqref{intro:eq-WM-M}, and also apply to wave equations with general quadratic derivative-nonlinearities. For the optimal range of regularities $s>\frac{d}{2}$, local well-posedness of \eqref{intro:eq-WM-M} was obtained by Klainerman and Machedon \cite{KM93,KM95,KM97} and Klainerman and Selberg \cite{KS97,KS02}. Their arguments heavily rely upon the algebraic structure of the null-form $\partial_\mu \phi^i \partial^\mu \phi^j$, which weakens the interactions between parallel waves, and the argument cannot be extended to wave equations with general quadratic derivative-nonlinearities \cite{L96}. At regularities $s\leq \frac{d}{2}$, the wave maps equation \eqref{intro:eq-WM-M} is expected to be ill-posed, and we refer the reader to \cite{DAG05,T00} and the references therein.\\

The global existence of solutions of \eqref{intro:eq-WM-M} in dimension $d\geq 2$ is a challenging problem even for small, smooth initial data. 
In the case $\scrM=\mathbb{S}^n$, it was solved in the seminal works of Tao, who first treated the high-dimensional case $d\geq 5$ in \cite{T01I} and then all dimensions $d\geq 2$ in \cite{T01II}. One of the main difficulties in \cite{T01I} is that, at the critical regularity $s=\frac{d}{2}$, the nonlinearity of the wave maps equation \eqref{intro:eq-WM-M} is non-perturbative. To overcome this difficulty, Tao relied on the gauge-freedom in \eqref{intro:eq-WM-M}, i.e., the freedom to work in any frame of the tangent bundle $T\mathbb{S}^n$. Using a carefully crafted microlocal gauge-transformation \cite[(9)]{T01I}, Tao then transformed the wave maps equation \eqref{intro:eq-WM-M} into a wave equation with a perturbative nonlinearity. One of the key ingredients in the argument from \cite{T01I} is the approximate orthogonality of the microlocal gauge transformation. This orthogonality is intimately tied to the geometric nature of the wave maps equation \eqref{intro:eq-WM-M} and has no counterpart in wave equations with general null-form nonlinearities. Further results on the global existence of solutions of \eqref{intro:eq-WM-M} with smooth, small initial data are available in the seminal articles of Tataru \cite{Tat98,Tat01}, who introduced part of the functional framework used in \cite{T01II}, and in \cite{KR01,K03,K04,NSU03,SS02,Tat05}, which consider more general target manifolds.\\ 

For large, smooth initial data, solutions of the wave maps equation \eqref{intro:eq-WM-M} can display diverse and intricate global dynamics. The global dynamics are best understood in the energy-critical case $d=2$ and are influenced by the geometry of the target manifold. Global existence of smooth solutions of \eqref{intro:eq-WM-M} is expected at all energies below the energy of any non-trivial harmonic map into the target manifold, and this has been obtained in the independent works \cite{KS12,ST10_1,ST10_2,TaoH3,TaoH4,TaoH5,TaoH6,TaoH7}. In particular, since there are no non-trivial harmonic maps from $\mathbb{R}^2$ into the hyperbolic spaces $\mathbb{H}^{m}$, this implies the global regularity of wave maps into hyperbolic spaces. At energies above the energy of the first non-trivial harmonic map into the target manifold, solutions may blow up in finite time \cite{KST08,RR12,RS10}. 
In the energy-supercritical case $d\geq 3$, much less is known about the long-time behavior of solutions with large, smooth initial data. At this time, most available results concern the existence and stability of self-similar blowup \cite{CDG17,D11,DW23,DW22,S88}.\\

Many aspects of the deterministic theory of wave maps, such as null-form estimates and the ideas behind Tao's microlocal gauge transformation, have already been used in the probabilistic theory of \cite{BLS21} and will continue to be used in this article. In contrast to \cite{BLS21}, however, this article also relies on the orthogonality properties we discussed in the context of \cite{T01II}, and therefore relies more heavily on the geometric aspects of the wave maps equation than \cite{BLS21}. 

\subsubsection{Probabilistic aspects of scalar dispersive equations}
The Hamiltonian evolution equations corresponding to the periodic real and complex-valued $\Phi^{4}_d$-functional are cubic wave and Schr\"{o}dinger equations, which are given by
\begin{align}
-\partial_t^2 \phi + \Delta \phi &= \phi^3  \hspace{10ex} (t,x) \in \mathbb{R} \times \mathbb{T}^d,  \label{intro:eq-NLW} \\ 
i \partial_t \phi + \Delta \phi &= |\phi|^{2} \phi \hspace{7.35ex} (t,x) \in \mathbb{R}\times \mathbb{T}^d.  \label{intro:eq-NLS}
\end{align}
The invariance of the Gibbs measure for \eqref{intro:eq-NLW} and \eqref{intro:eq-NLS} was first shown in dimension $d=1$ and $d=2$  by Friedlander \cite{F85}, Zhidkov \cite{Z94}, and Bourgain \cite{B94,B96,B99}. More recently, \cite{GKO18} introduced a para-controlled approach to wave equations, \cite{DNY19} introduced random averaging operators, and \cite{DNY22} introduced the method of random tensors, which can all be used to study random scalar dispersive equations in more singular settings than in \cite{B96}. In \cite{BDNY22}, the para-controlled approach and random tensor estimates were then used to prove the invariance of the Gibbs measure for \eqref{intro:eq-NLW} in dimension $d=3$. As of the time of writing of this article, the invariance of the Gibbs measure for the \eqref{intro:eq-NLS} in dimension $d=3$ is still open (cf. \cite[Section 9]{DNY22}). Despite the remaining open problems, our understanding of the probabilistic theory of \eqref{intro:eq-NLW} and \eqref{intro:eq-NLS} has improved significantly since the beginning of this field, and is certainly much more advanced than for geometric wave equations such as \eqref{intro:eq-WM-M}. 

While the results obtained for \eqref{intro:eq-NLW} and \eqref{intro:eq-NLS} are important inspirations for us, the wave maps equation is more closely related to scalar-valued wave equations with quadratic derivative-nonlinearities, such as 
\begin{equation}\label{intro:eq-DNLW}
-\partial_t^2 \phi + \Delta \phi = |\nabla \phi|^2 
\hspace{10ex} (t,x) \in \mathbb{R}\times\mathbb{R}^3. 
\end{equation}
In \cite{B21}, the author proved the probabilistic well-posedness of \eqref{intro:eq-DNLW} at regularities $s\geq 1.984$, which barely beats the Lorentz-critical regularity $s_L=2$ of \eqref{intro:eq-DNLW}. The main difficulty lies in the absence of nonlinear smoothing for \eqref{intro:eq-DNLW}, which is caused by low$\times$high-interactions. In \cite{B21}, the problematic low$\times$high-interactions are absorbed into so-called adapted linear evolutions, which capture the evolution of the high-frequency random initial data via a linearization of \eqref{intro:eq-DNLW} around the low-frequency components. The method of \cite{B21} later has been used in \cite{CGI22,KLS20} and partially motivates the modulated linear waves below. However, while the absence of nonlinear smoothing is a common feature of both \eqref{intro:eq-WM-M} and \eqref{intro:eq-DNLW}, the setting of this article is much more difficult than the setting of \cite{B21}. The reason is that the high$\times$high$\rightarrow$low-interactions in \eqref{intro:eq-WM-M} pose severe challenges, whereas the high$\times$high$\rightarrow$low-interactions in \eqref{intro:eq-DNLW} are harmless. 

\subsubsection{Stochastic Yang-Mills equations} 

We let $\frkG$ be a Lie group and $\frkg$ be the corresponding Lie algebra. For the notation used in the following discussion, we refer the reader to \cite{BC23,CCHS20,CCHS22}. The stochastic Yang-Mills equation for $A\colon [0,\infty) \times \bT^d \rightarrow \frkg^d$, which is the Langevin equation corresponding to \eqref{intro:eq-functional-YM}, is given by
\begin{equation}\label{intro:eq-SYM}
\partial_t A = - D_A^\ast F_A + \xi \qquad \qquad (t,x) \in (0,\infty)\times \bT^d, 
\end{equation}
where $D_A$ is the covariant derivative, $D_A^\ast$ is its adjoint, and $\xi$ is a $\frkg^d$-valued space-time white noise. For any smooth $g\colon \bT^d\rightarrow \frkG$, the corresponding gauge transformation is defined by  $A \mapsto A^g= g A g^{-1} - (\mathrm{d}g)  g^{-1}$. 
Formally, \eqref{intro:eq-SYM} is invariant under this gauge transformation, and can therefore not be parabolic. The lack of parabolicity can be addressed using the DeTurck trick \cite{DeT83}, which allows us to transform \eqref{intro:eq-SYM} into 
\begin{equation}\label{intro:eq-SYM-heat}
\partial_t A = - D_A^\ast F_A - D_A D_A^\ast A+ \xi \qquad \qquad (t,x) \in (0,\infty)\times \bT^d. 
\end{equation}
Since the principal term in $- D_A^\ast F_A - D_A D_A^\ast A$ is given by $\Delta A$, \eqref{intro:eq-SYM-heat} is commonly referred to as the stochastic Yang-Mills heat equation. \\ 

In dimension $d=2,3$, the local well-posedness of \eqref{intro:eq-SYM-heat} was obtained using regularity structures in \cite{CCHS20,CCHS22}. In addition to local well-posedness, \cite{CCHS20,CCHS22} also proved the gauge-covariance of \eqref{intro:eq-SYM-heat}, constructed a state space for solutions of \eqref{intro:eq-SYM-heat}, and constructed a Markov process corresponding to the projection of solutions to \eqref{intro:eq-SYM-heat} onto the space of gauge orbits. In dimension $d=2$, the local well-posedness and gauge-invariance have also been revisited using para-controlled calculus \cite{BC23}. In a recent seminal article \cite{CS23}, Chevyrev and Shen obtained the global well-posedness and invariance for \eqref{intro:eq-SYM} in dimension $d=2$. In addition to earlier methods from \cite{CCHS20}, the argument in \cite{CS23} relies on a lattice-approximation of \eqref{intro:eq-SYM}, properties of the two-dimensional Yang-Mills measure, and Bourgain's globalization argument \cite{B94}. In contrast, global well-posedness and invariance\footnote{In fact, even the construction of the Yang-Mills measure is still open in dimension $d=3$.} are still open in dimension $d=3$ and likely present a significant challenge. 

We emphasize that there is a striking parallel between the developments for the stochastic Yang-Mills equation and the wave maps equation with random initial data: In both settings, going from local well-posedness to global well-posedness and invariance poses a formidable challenge, and it goes far beyond a direct application of Bourgain's globalization argument \cite{B94}. 

\subsubsection{Random geometric wave equations}
The analysis of random geometric wave equations is still in its infancy and, at this moment, the only related works are the previously mentioned articles \cite{BLS21,BJ22} and \cite{B22EWM,BR23,KLS20}. In \cite{B22EWM}, the author studied the wave maps equation \eqref{intro:eq-WM-M} in dimension $d=3$, but restricted the model to equivariant maps on exterior domains. For each of the topological solitons $Q_n$ of this restricted model (see \cite[(1.7)]{B22EWM}), the author then proved the existence and invariance of a Gibbs measure which is supported on the homotopy class of $Q_n$. This is possible even in dimension $d=3$ since the exterior equivariant wave maps equation can be written as a one-dimensional wave equation with a sine-nonlinearity rather than a derivative-nonlinearity.

The hyperbolic Yang-Mills equation corresponding to \eqref{intro:eq-functional-YM} has not yet been studied from a probabilistic perspective, but related models have been studied in \cite{BR23,KLS20}. In \cite{KLS20}, the authors studied the energy-critical Maxwell-Klein-Gordon equation with random initial data below the energy space. One of the major achievements of \cite{KLS20} is to combine the probabilistic method from \cite{B21} with the deterministic functional framework of \cite{KST15} and induction procedure of \cite{KL15}. While \cite{KLS20} considers one of the most intricate and natural geometric wave equations, however, the random initial data in \cite{KLS20} stems from a Wiener-randomization and is not built in a gauge-covariant fashion. 
In \cite[(1.3)]{BR23}, the author and Rodnianski introduced a new model for gauge-covariant wave equations. While the deterministic aspects of this model are much simpler than for the energy-critical Maxwell-Klein-Gordon equation in \cite{KLS20}, the model includes a gauge-covariant, singular stochastic forcing. Using ideas from the probabilistic theory of scalar wave equations \cite{B23,DNY19,DNY22}, \cite{BR23} then proved the probabilistic global well-posedness of this gauge-covariant wave equation.

\subsection{Main result}\label{section:introduction-result}
We let $\frkG$ be a compact Lie group\footnote{Using $\frkG$ to denote a Lie group is rather uncommon and $G$ is used more widely. In this article, however, $G$ will be reserved for Gaussian random variables.} which is equipped with a bi-invariant Riemannian metric. Due to the Peter-Weyl theorem (see e.g. \cite[Corollary IV.4.22]{K02}), we may assume that $\frkG$ is a closed subgroup of $\textup{Gl}(n,\mathbb{C})$. For example, $\frkG$ may be taken as the special orthogonal group $\textup{SO}(n)$, the unitary group $\textup{U}(n)$, or the special unitary group $\textup{SU}(n)$. The Lie algebra corresponding to the Lie group $\frkG$ is denoted by $\frkg$ and is equipped with the induced Lie bracket 
$[\cdot,\cdot]\colon \frkg \times \frkg \rightarrow \frkg$, the induced inner-product $\langle \cdot, \cdot \rangle \colon \frkg \times \frkg \rightarrow \R$, and the induced norm 
$\| \cdot \|\colon \frkg \rightarrow [0,\infty)$. \\

In the case $\scrM=\frkG$,  the Lagrangian of the wave maps equation from \eqref{intro:eq-Lagrangian} can be written as 
\begin{equation}\label{intro:eq-Lagrangian-G}
\mathcal{L}(\phi) = \int_{\R^{1+1}} \Big( - \big\| \phi(t,x)^{-1} \partial_t \phi(t,x)  \big\|_{\frkg}^2 
+ \big\| \phi(t,x)^{-1} \partial_x \phi(t,x) \big\|_\frkg^2 \Big)  \dx \dt. 
\end{equation}
Since the wave maps equation is the Euler-Lagrange equation for \eqref{intro:eq-Lagrangian-G}, it then takes the simple form \begin{equation}\label{intro:eq-WM-G}
\partial_t \big( \phi^{-1} \partial_t \phi \big) = \partial_x \big( \phi^{-1} \partial_x \phi \big). 
\end{equation} 

Instead of working with the $\frkG$-valued wave map $\phi\colon \R^{1+1}\rightarrow \frkG$, we primarily work with the $\frkg$-valued maps $A,B\colon \R^{1+1}\rightarrow \frkg$, which are defined as 
\begin{equation}\label{intro:eq-A-B}
A:= \phi^{-1} \partial_t \phi \qquad \text{and} \qquad 
B:= \phi^{-1}  \partial_x \phi. 
\end{equation}
Due to our assumption that $\frkG$ is a closed subgroup of $\textup{Gl}(n,\mathbb{C})$, the products in \eqref{intro:eq-A-B} can be interpreted as matrix products. The reason for introducing $A$ and $B$ is that the Lie algebra $\frkg$ is a linear space and several analytical tools are therefore more easily used for $A$ and $B$ than for $\phi$. 
As shown in Lemma \ref{lifting:lem-equivalence}, the wave maps equation \eqref{intro:eq-WM-G} is equivalent to the system
\begin{equation}\label{intro:eq-system-A-B}
\partial_t A = \partial_x B \qquad 
\text{and} \qquad \partial_t B = \partial_x A - [ A , B]. 
\end{equation}
The Hamiltonian of \eqref{intro:eq-system-A-B} is given by
\begin{equation}
H(A,B) = \int_{\R} \dx \big( \| A(x) \|_{\frkg}^2 + \| B(x) \|_{\frkg}^2 \big)
\end{equation}
and the corresponding Gibbs measure at inverse temperature $\beta>0$ is formally given by
\begin{equation}\label{intro:eq-Gibbs-AB-formal}
``\mathrm{d}\mu_\beta(A,B) = \mathcal{Z}_\beta^{-1} \exp \Big( - \beta H(A,B) \Big)  \mathrm{d}A\mathrm{d}B\,".
\end{equation}
Since $A$ and $B$ take values in the linear space $\frkg$ and the energy is quadratic, the formal definition \eqref{intro:eq-Gibbs-AB-formal} can easily be made rigorous. Indeed, we can define
\begin{equation}\label{intro:eq-Gibbs-AB-rigorous}
\mu_\beta = \operatorname{Law} \Big( \big( \beta^{-\frac{1}{2}} W_0, \beta^{-\frac{1}{2}} W_1 \big) \Big),
\end{equation}
where $W_0,W_1\colon \R \rightarrow \frkg$ are two independent, $\frkg$-valued white noises (see Definition \ref{prelim:def-white-noise-g-valued}).  
Since the regularity of white noise is just below $-\frac{1}{2}$, the concept of a solution of \eqref{intro:eq-system-A-B} with initial data $(\beta^{-\frac{1}{2}}W_0,\beta^{-\frac{1}{2}} W_1)$ is difficult to define directly, and we instead use a limiting procedure. To this end, we first let $N\in \dyadic$ and define $(A_{\leq N},B_{\leq N})$ as the solution of the initial value problem
\begin{equation}\label{intro:eq-system-A-B-N}
\begin{cases}
\begin{aligned}
\partial_t A_{\leq N} &= \partial_x B_{\leq N},  &
\qquad \partial_t B_{\leq N} &= \partial_x A_{\leq N} - \big[ A_{\leq N}, B_{\leq N} \big], \\ 
A_{\leq N}(0) &= \beta^{-\frac{1}{2}} P_{\leq N}^x W_0, &
B_{\leq N}(0) &= \beta^{-\frac{1}{2}} P_{\leq N}^x W_1,
\end{aligned}
\end{cases}
\end{equation}
where $P_{\leq N}^x$ is the Littlewood-Paley operator from \eqref{prelim:eq-P-N} below. In other words, $(A_{\leq N},B_{\leq N})$ is the solution of the wave maps equation with frequency-truncated initial data. The local convergence result of \cite{BLS21} implies that, for any $R\geq 1$ and any sufficiently small-time $0<\tau\ll 1$, the sequence $(A_{\leq N},B_{\leq N})$ converges on the space-time domain $[-\tau,\tau]\times [-R,R]$ as $N$ tends to infinity. In the main theorem of this article, we now obtain the global convergence of $(A_{\leq N},B_{\leq N})$ and the invariance of the Gibbs measure under the limiting dynamics. 

\begin{theorem}[Global well-posedness and invariance]\label{intro:thm-rigorous-A-B}
Let $\beta\in (0,\infty)$, let $s<\frac{1}{2}$, and, for all $N\in \dyadic$, let $(A_{\leq N},B_{\leq N})\colon \R^{1+1} \rightarrow \frkg^2$ be the unique global solution of \eqref{intro:eq-system-A-B-N}. Then, the limit 
\begin{equation*}
(A,B) = \lim_{N\rightarrow \infty} (A_{\leq N},B_{\leq N}) 
\end{equation*}
almost surely exists in $C_t^0 \C_x^{s-1}([-T,T]\times [-R,R]\rightarrow \frkg)^2$ for all $T,R\geq 1$. Furthermore, for each time $t\in \R$, the law of $(A(t),B(t))$ is given by the Gibbs measure $\mu_\beta$ from \eqref{intro:eq-Gibbs-AB-rigorous}.
\end{theorem}

While we postpone an overview of the proof of Theorem \ref{intro:thm-rigorous-A-B} to Section \ref{section:overview}, we briefly describe the significance of both the theorem and our argument. As previously mentioned, Theorem \ref{intro:thm-rigorous-A-B} is the first result on the invariance of Gibbs measures under the evolution of a geometric wave equation. Its proof relies on techniques for geometric wave equations with deterministic data at critical regularity \cite{T01I} and techniques for scalar-valued wave equations with random initial data \cite{B21,DNY19,GKO18}. Most importantly, this article combines analytic, geometric, and probabilistic techniques in a single argument, and understanding the interplay of all three aspects constitutes a monumental task. In our argument, this interplay manifests itself in the following forms: 
\begin{enumerate}[label=(\roman*)]
\item Our Ansatz in \eqref{intro:eq-ansatz-U}-\eqref{intro:eq-ansatz-V} consists of modulated linear waves, modulated bilinear waves, and mixed modulated objects, which exhibit different levels of analytic and probabilistic structure. 
\item The absence of nonlinear smoothing and skew-symmetry of the Lie bracket both present themselves in the modulated linear waves (Subsection \ref{section:overview-ansatz}), which combine analytic and geometric aspects. 
\item The Killing map on the Lie algebra $\frkg$ and the covariance function of white noise and its anti-derivative enter into the Killing-renormalization (Subsection \ref{section:overview-killing}), which therefore combines geometric and probabilistic ingredients. 
\item High$\times$high$\rightarrow$low-interactions and the Jacobi identity on the Lie algebra $\frkg$ affect the Jacobi errors (Subsection \ref{section:overview-jacobi}), which therefore involve analytic and geometric aspects.
\item The proof of invariance (Subsection \ref{section:overview-almost}) relies on a geometric identity involving the Killing map and Lie bracket, and therefore uses geometric ingredients to prove a probabilistic statement.
\end{enumerate}

\revision{Our proof of Theorem \ref{intro:thm-rigorous-A-B} yields detailed information on the random structure of $(A_{\leq N},B_{\leq N})$, which is captured through our Ansatz (see e.g. \eqref{ansatz:eq-UN-rigorous-decomposition}-\eqref{ansatz:eq-VN-rigorous-decomposition} and Definitions \ref{ansatz:def-modulated-linear}, \ref{ansatz:def-modulated-bilinear}, \ref{ansatz:def-mixed}, and \ref{ansatz:def-modulated-linear-reversed}). Due to this,} it is not too difficult to transfer our results from the $\frkg$-valued derivatives back to the original, $\frkG$-valued map (see Section \ref{section:lifting}). The corresponding Gibbs measure is denoted by $\muGb$ and formally given by 
\begin{equation}\label{intro:eq-Gibbs-phi-G-formal}
``\mathrm{d}\muGb\big(\phi_0,\phi_1\big) = \mathcal{Z}_\beta^{-1} \exp \bigg( - \frac{\beta}{2} 
\int_{\R} \dx \Big( \big\| \phi_0^{-1} \partial_x \phi_0 \big\|_{\frkg}^2 + \big\| \phi_0^{-1} \phi_1 \big\|_{\frkg}^2 \Big) \bigg) 
\mathrm{d}\phi_0 \mathrm{d}\phi_1". 
\end{equation}
The rigorous definition of $\muGb$, which is more involved than the rigorous definition of $\mu_\beta$ in \eqref{intro:eq-Gibbs-AB-rigorous}, is postponed until Subsection \ref{section:lifting-Gibbs}. From the rigorous definition, it follows that $\phi_0(0)$ is distributed according to the Haar measure on $\frkG$, $\phi_0(0)^{-1} \phi_0(x)$ is a $\frkG$-valued Brownian motion at inverse temperature $\beta>0$, and $\phi_0(x)^{-1} \phi_1(x)$ is a $\frkg$-valued white noise at inverse temperature $\beta>0$. 

To introduce the mollified initial data, we let $W_0,W_1\colon \R \rightarrow \frkg$ be two independent, $\frkg$-valued white noises. Furthermore, we let $g$ be a $\frkG$-valued random variable that is independent of $(W_0,W_1)$ and whose law is given by the Haar measure on $\frkG$. For each $N\in \dyadic$, we then define $(\phi_{\leq N,0},\phi_{\leq N,1}) \colon \R \rightarrow T\frkG$ by 
\begin{equation}\label{intro:eq-phi-initial-mollified}
\phi_{\leq N,0}(0)=g, \quad \partial_x \phi_{\leq N,0}(x) = \beta^{-\frac{1}{2}} \phi_{\leq N,0}(x) P_{\leq N} W_0, \quad \text{and} \quad 
\phi_{\leq N,1}(x) = \beta^{-\frac{1}{2}} \phi_{\leq N,0}(x) P_{\leq N} W_1. 
\end{equation}
Equipped with \eqref{intro:eq-phi-initial-mollified}, we can now state our main theorem at the level of $\frkG$-valued wave maps.

\begin{theorem}[Global well-posedness and invariance for $\frkG$-valued maps]\label{intro:thm-rigorous-phi}
Let $\beta \in (0,\infty)$, let $s<\frac{1}{2}$, and, for all $N\in \dyadic$, let $\phi_{\leq N}\colon \R^{1+1}\rightarrow \frkG$ be the unique global solution of the wave maps equation \eqref{intro:eq-WM-G} with initial data as in \eqref{intro:eq-phi-initial-mollified}. Then, the limit
\begin{equation*}
\big(\phi, \partial_t \phi \big) = \lim_{N\rightarrow \infty} \big( \phi_{\leq N}, \partial_t \phi_{\leq N} \big)
\end{equation*}
almost surely exists in $(C_t^0 \C_x^{s}\times C_t^0 \C_x^{s-1})([-T,T]\times [-R,R]\rightarrow T\frkG)$ for all $T,R\geq 1$. Furthermore, for each $t\in \R$, the law of $(\phi(t),\partial_t \phi(t))$ is given by the Gibbs measure $\muGb$ from Definition \ref{lifting:def-Gibbs}.
\end{theorem}

As mentioned above, Theorem \ref{intro:thm-rigorous-phi} follows without difficulties from the proof of Theorem \ref{intro:thm-rigorous-A-B}. The only additional ingredients are lifts (Definition \ref{lifting:def-lifts}), which allow us to represent the $\frkG$-valued map $\phi_{\leq N}$ in terms of the $\frkg$-valued derivatives $A_{\leq N}$ and $B_{\leq N}$. 

\begin{remark}[Further properties of the limits]\label{intro:rem-properties}
Using a modification of our arguments, it is likely possible to show that the limits in Theorem \ref{intro:thm-rigorous-A-B} and Theorem \ref{intro:thm-rigorous-phi} do not depend on the choice of the symbol $\rho$, which appears in the definition of $P_{\leq N}$, and that the limits satisfy the flow property. For more detailed discussions of both properties, we refer the reader to Remark \ref{lipschitz:rem-kernel-rho} and Remark \ref{main:rem-flow}, respectively. 
\end{remark}

\begin{remark}[Complete integrability]
The $(1+1)$-dimensional wave maps equation is known to be completely integrable \cite{TU04}. However, our argument only relies on analytic, geometric, and probabilistic techniques, and does not directly rely on integrable methods. Nevertheless, it is an interesting (but likely very difficult) question whether Theorem \ref{intro:thm-rigorous-A-B} and Theorem \ref{intro:thm-rigorous-phi} can also be obtained using integrable techniques. 
\end{remark}

\textbf{Acknowledgements:}
First and foremost, the author thanks Igor Rodnianski for numerous insightful comments and helpful discussions during the work on this project. The author also thanks Sky Cao, Ilya Chevyrev, Jonas L\"{u}hrmann, Tadahiro Oh, Katharina Schratz, Gigliola Staffilani, and Nikolay Tzvetkov for interesting conversations. \revision{Finally, the author thanks the anonymous referees for their detailed and helpful comments.} The author was partially supported by the National Science Foundation under Grant No. DMS-1926686 \revision{and DMS-2453126}.

\section{Overview of the argument}\label{section:overview}

During the following overview of our argument, we focus on Theorem \ref{intro:thm-rigorous-A-B}, i.e., the $\frkg$-valued system \eqref{intro:eq-system-A-B}. For notational purposes, we introduce the temperature $\coup>0$, which is defined as 
\begin{equation}\label{intro:eq-coup}
\coup = \frac{1}{8\beta}.
\end{equation}
The factor of $\frac{1}{8}$ in \eqref{intro:eq-coup} is included for notational convenience since it leads to the simplest pre-factors in \eqref{intro:eq-system-U-V-N-data} below. For simplicity, we also assume that the $\frkg$-valued white noise $(W_0,W_1)$ is $2\pi$-periodic, in which case it can be written as 
\begin{equation}\label{intro:eq-W0-W1}
W_0 = \sum_{n\in \Z} G_{0,n} e^{inx} \qquad 
\text{and} \qquad 
W_1 = \sum_{n\in \Z} G_{1,n} e^{inx}. 
\end{equation}
Here, $(G_{0,n})_{n\in \Z}$ and $(G_{1,n})_{n\in \Z}$ are independent, standard, $\frkg$-valued Gaussian sequences (see  Definition \ref{prelim:def-standard-Gaussian-sequence}). 
Due to finite speed of propagation, the periodicity assumption is  not essential, and will later be removed.  
We split the rest of this overview of our argument into five parts, which address the following aspects: 
\begin{enumerate}[label=(\Roman*)]
\item Finite-dimensional approximation, 
\item Ansatz, 
\item Killing-renormalization and orthogonality, 
\item Jacobi errors and Bourgain-Bulut argument, 
\item and almost invariance and conservative structures.  
\end{enumerate}

\begin{remark}[The temperature $\coup$]
In many articles on random dispersive equations, the temperature $\coup$ is not too important and therefore fixed as a simple value (see, however, \cite{OOT21,OST22}). In this article, the main reason to keep the temperature $\coup$ as a parameter is to later use scaling (see e.g. Lemma \ref{ansatz:lem-scaling-symmetry}
and the proof of Proposition \ref{main:prop-refined-gwp}). Via scaling, we can avoid using cut-offs to time-scales $\tau\ll 1$ in our 
local theory, and instead only use cut-offs to time-scales $\sim 1$, which is technically convenient.
\end{remark}

\subsection{Finite-dimensional approximation}\label{section:overview-discretization}
We previously introduced \eqref{intro:eq-system-A-B-N}, which is a natural approximation of the initial-value problem for \eqref{intro:eq-system-A-B} with white noise initial data. In \cite{BLS21}, it has already been shown that \eqref{intro:eq-system-A-B-N} is locally well-posed. More precisely, it has been show that the limit of $(A_{\leq N},B_{\leq N})$ exists with high-probability on the space-time domain $[-\tau,\tau]\times [-R,R]$ as long as $\tau$ is sufficiently small (depending on $R$). The issue with the approximation \eqref{intro:eq-system-A-B-N} is that it is ill-suited for a proof of global well-posedness and invariance of white noise for \eqref{intro:eq-system-A-B}. This is because the initial data in \eqref{intro:eq-system-A-B-N} is not invariant under the flow, which makes it difficult to use Bourgain's globalization argument \cite{B94,B96}. In order to prove global well-posedness and invariance, we are therefore looking for a new finite-dimensional approximation, which should be different from but still have the same limit as \eqref{intro:eq-system-A-B-N}. In order to introduce and discuss this new approximation of \eqref{intro:eq-system-A-B}, we first need to introduce null-coordinates and null-variables. \\

We define the null-coordinates $(u,v)\in \mathbb{R}^{1+1}$ by $u:= x-t$ and $v:=x+t$. We also define the corresponding null-variables by $U:= \frac{1}{4} (A-B)$ and $V:=\frac{1}{4}(A+B)$. In the null-coordinates and null-variables, the system \eqref{intro:eq-system-A-B-N} takes the form 
\begin{equation}\label{intro:eq-system-U-V}
\partial_v U = \partial_u V = \big[ U,V\big]. 
\end{equation}
We emphasize that \eqref{intro:eq-system-U-V} only contains cross-interactions between $U$ and $V$ and no self-interactions of $U$ or $V$, which is due to the null structure of the wave maps equation. In null-coordinates, the initial time-slice $\{(t,x)\in \R^{1+1}\colon t=0\}$ corresponds to the diagonal $\{(u,v)\in \R^{1+1}\colon u=v\}$. The initial condition $(A(0),B(0))=(W_0,W_1)$ can then be written as 
\begin{equation}\label{intro:eq-initial-U-V}
U\big|_{u=v} = \hcoup W^+ \qquad \text{and} \qquad V\big|_{u=v} = \hcoup W^-, 
\end{equation}
where $W^+,W^-\colon \R \rightarrow \frkg$ are defined as
\begin{equation}\label{intro:eq-Wpm}
W^\pm = \frac{\sqrt{8}}{4} \big(W_0 \mp W_1\big) = \sum_{n\in \Z} \frac{G_{0,n}\mp G_{1,n}}{\sqrt{2}} e^{inx} =: \sum_{n\in \Z} G_n^\pm e^{inx}. 
\end{equation}  
Due to the rotation invariance of Gaussians, $(G_n^+)_{n\in \Z}$ and $(G_n^{-})_{n\in \Z}$ are independent, standard,  $\frkg$-valued Gaussian sequences. 
Equipped with the null-coordinates and null-variables, 
we can write the wave maps equation with frequency-truncated initial data \eqref{intro:eq-system-A-B-N} as 
\begin{equation}\label{intro:eq-system-U-V-N-data}
\begin{cases}
\partial_v U_{\leq N} = \partial_u V_{\leq N} = \big[ U_{\leq N}, V_{\leq N} \big], \\ 
U_{\leq N}\big|_{u=v}= \hcoup P_{\leq N}W^+, \quad V_{\leq N}\big|_{u=v} = \hcoup P_{\leq N}W^-.
\end{cases}
\end{equation}

We now revive our search for a new finite-dimensional approximation of \eqref{intro:eq-system-A-B} or, equivalently, of \eqref{intro:eq-system-U-V}. Our new approximation should, as best as possible, preserve the following two features of \eqref{intro:eq-system-A-B} and \eqref{intro:eq-system-U-V}: 
\begin{enumerate}[label=(\roman*)]
\item \label{intro:item-null}   The null structure, which is needed for well-posedness. 
\item \label{intro:item-Hamiltonian} The Hamiltonian structure, which is needed for the invariance of white noise. 
\end{enumerate}
It is rather difficult to 
simultaneously preserve \ref{intro:item-null} and \ref{intro:item-Hamiltonian} under approximations of \eqref{intro:eq-system-U-V}. The reason is that the null-structure relies on the equal treatment of time and space, whereas the Hamiltonian structure treats time and space differently. If we only needed to preserve either \ref{intro:item-null} or \ref{intro:item-Hamiltonian}, but not both, then life would be much easier. For example, \ref{intro:item-null} is preserved by the approximation 
\begin{equation}\label{intro:eq-wrong-approximation-null}
\partial_v U^{(N)} = \partial_u V^{(N)} = P_{\leq N}^{u,v} \Big[ P_{\leq N}^{u,v} U^{(N)} , P_{\leq N}^{u,v} V^{(N)} \Big]. 
\end{equation}
\revision{In \eqref{intro:eq-wrong-approximation-null}, $P_{\leq N}^{u,v}$ is the Littlewood-Paley projection in both the $u$ and $v$-variable, which is defined in \eqref{prelim:eq-P-N-uv} below.} However, while the local well-posedness of \eqref{intro:eq-wrong-approximation-null} can be shown using the arguments in \cite{BLS21}, \eqref{intro:eq-wrong-approximation-null} completely breaks the Hamiltonian structure of \eqref{intro:eq-system-A-B}. In fact, due to the $P_{\leq N}^{u,v}$-operators, the resulting system for $A^{(N)}$ and $B^{(N)}$ can no longer even be written as a Banach-space valued differential equation in time. Alternatively, \ref{intro:item-Hamiltonian} could be preserved\footnote{This is not obvious since \eqref{intro:eq-wrong-approximation-Hamiltonian} does not stem from an approximation of the Hamiltonian, but rather an approximation of the nonlinear structure of the state space. For more details, we refer the reader to Section \ref{section:conservative}.} by the approximation 
\begin{equation}\label{intro:eq-wrong-approximation-Hamiltonian}
\partial_v U^{(N)} = \partial_u V^{(N)} = P_{\leq N}^x \Big[ P_{\leq N}^x U^{(N)} , P_{\leq N}^x V^{(N)} \Big]. 
\end{equation}
The issues with \eqref{intro:eq-wrong-approximation-Hamiltonian} are more difficult to explain than for \eqref{intro:eq-wrong-approximation-null}. Loosely speaking, the $P_{\leq N}^x$-operators introduce shifts in the $x$-variables. In interactions between $\frkg$-valued white noises $W^\pm$ and their integrals $\Int W^\pm$, the shifts then lead to problematic probabilistic resonances which are not present in \eqref{intro:eq-system-U-V-N-data}. \revision{Here, we use the term ``problematic probabilistic resonances" to describe interactions involving Gaussian random variables which produce relatively large expectations or, more generally, interactions whose projections to lower-order Gaussian chaoses are relatively large. For a more detailed discussion of the problematic probabilistic resonances, we refer to Section \ref{section:overview-killing} below and Section \ref{section:ansatz-remainder}, Section \ref{section:killing-tensor}, and Section \ref{section:jacobi} in the main body of this article. The problematic probabilistic resonances} not only makes it difficult to establish well-posedness of \eqref{intro:eq-wrong-approximation-Hamiltonian}, but likely also prevents the convergence of the solutions of \eqref{intro:eq-system-U-V-N-data} and solutions of \eqref{intro:eq-wrong-approximation-Hamiltonian} to the same limit.  \\ 

Instead of \eqref{intro:eq-wrong-approximation-null} or \eqref{intro:eq-wrong-approximation-Hamiltonian}, we work with a new approximation of \eqref{intro:eq-system-U-V} which is given by 
\begin{equation}\label{intro:eq-system-U-V-Kil}
\partial_v U^{(N)} = \partial_u V^{(N)} = P_{\leq N}^x \Big[ P_{\leq N}^x U^{(N)} , P_{\leq N}^x V^{(N)} \Big] - \coup \Renorm[N] (U^{(N)}+V^{(N)}). 
\end{equation}
The goal behind the counterterm in \eqref{intro:eq-system-U-V-Kil} is to strike a delicate balance: On the one hand, we want that it removes the most problematic probabilistic resonances in \eqref{intro:eq-wrong-approximation-Hamiltonian} and thereby restores the well-posedness. On the other hand, while the counterterm breaks the conservation of the energy, we still want to be able to control the energy increment and thereby prove the almost invariance of white noise under \eqref{intro:eq-system-U-V-Kil}. 
To achieve this goal, we choose $\Renorm[N]$ as the so-called Killing-renormalization, which is natural from both probabilistic and geometric perspectives (see Definition \ref{ansatz:def-Killing}): In the $x$-variable, $\Renorm[N]$ acts as a convolution against the covariance function of white noise and its antiderivative. On the Lie algebra, it acts as the Killing map $\Kil\colon \frkg \rightarrow \frkg$, which is a central object from the theory of Lie algebras.

The daunting task ahead of us is to show that, with our choice of $\Renorm[N]$, \eqref{intro:eq-system-U-V-Kil} is indeed well-posed, converges to the same limit as \eqref{intro:eq-system-U-V-N-data}, and almost preserves $\frkg$-valued white noise \revision{(see Proposition \ref{main:prop-almost-invariance})}.  This presents many formidable challenges and requires several additional ideas, which will be discussed throughout the rest of this overview. 

\begin{remark}[Killing map]\label{intro:rem-killing}
The Killing map $\Kil\colon \frkg \rightarrow \frkg$, which enters into the definition of our Killing-renormalization $\Renorm[\Nscript][]$, is related to the Ricci curvature on $\frkG$ (see e.g. \cite[Proposition 21.19]{GQ20}). To be more precise, if $\textup{Ric}\colon \frkg \times \frkg \rightarrow \mathbb{R}$ is the restriction of the Ricci curvature to $\frkg=T_e \frkG$, then it holds for all $X,Y\in \frkg$ that
\begin{equation*}
\textup{Ric}\big( X , Y \big) = \big\langle \Kil \big( X \big) , Y \big \rangle. 
\end{equation*}
We also note that the Killing map not only appears in the renormalization of the finite-dimensional approximation \eqref{intro:eq-system-U-V-Kil}, but also appears in the renormalization of stochastic objects in the stochastic Yang-Mills heat equation \cite{BC23,CCHS20}.
\end{remark}

\begin{remark}[Lattice-approximation of wave maps]
While our approximation \eqref{intro:eq-system-U-V-Kil} of \eqref{intro:eq-system-U-V} relies on frequency-truncations and counterterms, one can also formulate lattice-approximations of \eqref{intro:eq-system-U-V}. Different lattice-approximations of the wave maps equation \eqref{intro:eq-WM-M} have previously been studied in several numerical works \cite{B09,B15,CV16,KW14,MS98} and, as mentioned above, in the probabilistic work \cite{BJ22}. However, it seems difficult to compare solutions of \eqref{intro:eq-system-A-B-N}, i.e., the wave maps equation with mollified initial data, and solutions of any lattice-approximation of the wave maps equation. Furthermore, it seems difficult to implement the analytic techniques of \cite{BLS21}, which heavily use frequency-space arguments, on a lattice. 
\end{remark}

\begin{remark}[Finite-dimensional approximation of scalar-valued dispersive equations]
In most works on invariant measures of scalar-valued dispersive equations, the choice of suitable finite-dimensional approximations is rather straightforward. However, even in the scalar-valued setting, this choice can be difficult if the well-posedness theory relies on gauge transformations (see e.g. \cite{DTV15,NORBS12}). 
\end{remark}

\begin{remark}[Discretization of singular parabolic SPDEs]
In \eqref{intro:eq-system-U-V-N-data} and \eqref{intro:eq-system-U-V-Kil}, we introduced two different approximations of the wave maps equation with white noise initial data. As part of the proof of Theorem \ref{intro:thm-rigorous-A-B}, we will show that both approximations converge to the same limit. The analysis of different approximations and their limits also plays an important role in singular parabolic stochastic partial differential equations, and we refer the reader to 
\cite{CM18,EH19,FH17,HM12,HM18} and the references therein. However, there is an important difference between the setting of this article and the parabolic setting, which is the absence of nonlinear smoothing (see Subsection \ref{section:overview-ansatz}). 
\end{remark}

\begin{remark}[General compact Riemannian manifolds]\label{intro:rem-general-case}
As mentioned below Theorem \ref{intro:thm-formal}, there are obstructions that prevent us from treating general compact Riemannian manifolds $\scrM$, in which case the wave maps equation is given by  
\begin{equation}\label{intro:eq-general-case}
\partial_\mu \partial^\mu \phi^k = \Gamma^k_{ij}(\phi) \partial_\mu \phi^i \partial^\mu \phi^j.
\end{equation}
In order to prove the global well-posedness and invariance for \eqref{intro:eq-general-case}, we would need to introduce a renormalized, finite-dimensional approximation of \eqref{intro:eq-general-case}. At least for Riemannian manifolds with non-constant scalar curvature, however, such finite-dimensional approximations are not even completely understood in the construction of the Gibbs measure or for the geometric stochastic heat equation. For a detailed discussion of this, we refer the reader to \cite[Section 4.3]{BGHZ22}. 
\end{remark}

\subsection{Ansatz}\label{section:overview-ansatz} 
The unknowns in our argument are the so-called modulation operators and nonlinear remainders, which are denoted by 
\begin{equation}\label{intro:eq-unknowns}
\SNintro[K][+][k], \SNintro[M][-][m] \colon \R^{1+1}_{u,v} \rightarrow \End(\frkg) \qquad \text{and} \qquad \UNintro[][\fs],\VNintro[][\fs] \colon \R^{1+1}_{u,v} \rightarrow \frkg. 
\end{equation}
In \eqref{intro:eq-unknowns}, $K,M\in \dyadic$ are dyadic scales, $k\in \Z_K$, and $m\in \Z_M$, where we set $\Z_N = \big\{ n\in \Z \colon \frac{N}{2}<|n|\leq N\big\}$.  
Using \eqref{intro:eq-unknowns}, our Ansatz for the solution of \eqref{intro:eq-system-U-V-Kil} can then be written as 
\begin{align}
\UNintro[][] &=   \UNintro[][+] +  \UNintro[][+-] +  \UNintro[][-] +  \UNintro[][+\fs] +  \UNintro[][\fs-] +\UNintro[][\fs],  \label{intro:eq-ansatz-U}\\
\VNintro[][] &= \VNintro[][-] +  \VNintro[][+-] +  \VNintro[][+] +  \VNintro[][\fs-] +  \VNintro[][+\fs] +\VNintro[][\fs]. \label{intro:eq-ansatz-V}
\end{align}
The first five terms in both \eqref{intro:eq-ansatz-U} and \eqref{intro:eq-ansatz-V} 
exhibit different levels of random structure and are all determined by the unknowns in 
\eqref{intro:eq-unknowns}. The last terms $\UNintro[][\fs]$ and $\VNintro[][\fs]$ in \eqref{intro:eq-ansatz-U} and \eqref{intro:eq-ansatz-V} previously appeared in \eqref{intro:eq-unknowns}. In our argument, we neither use nor prove any information on the structure of $\UNintro[][\fs]$ and $\VNintro[][\fs]$, and simply treat them as nonlinear, smooth remainders. 
Our discussion of \eqref{intro:eq-ansatz-U} and \eqref{intro:eq-ansatz-V} in this introduction will be rather brief and we refer to Section \ref{section:ansatz} for more details. \revision{Throughout this brief discussion, we use the symbol $\simeq$ rather informally, and it only suggests that the two terms on each side behave somewhat similarly.} In the following, we focus our attention on the two terms $\UNintro[][+]$ and $\UNintro[][+-]$, \revision{but briefly discuss the remaining random terms at the end of the subsection}.   

\subsubsection{The modulated linear wave $\UNintro[][+]$:}\label{section:overview-UNintro}
The modulated linear wave $\UNintro[][+]$ is defined as\footnote{For technical reasons,
we later restrict the dyadic sum in \eqref{intro:eq-Ansatz-Up-Sp} to dyadic scales $K\gg 1$, see e.g. \eqref{prelim:eq-sumlarge} and Definition \ref{ansatz:def-modulated-linear}.}
\begin{equation}\label{intro:eq-Ansatz-Up-Sp}
\UNintro[][+](u,v)= \sum_{K\in \dyadic} \UNintro[K][+](u,v) \qquad \text{and} \qquad \UNintro[K][+](u,v) = \hcoup \sum_{k\in \Z_K} \SNintro[K][+][k](u,v) G_k^+ \, e^{iku}. 
\end{equation}
Its purpose is to capture the evolution of the high-frequency initial data $W_K^+=\sum_{k\in \Z_K} G_k^+ e^{iku}$ on a low-frequency background. More 
precisely, $\UNintro[K][+]$ is meant to be an approximate solution of 
\begin{equation}\label{intro:eq-ansatz-Up}
\begin{cases}
    \begin{aligned}
    \partial_v \UNintro[K][+] &\simeq  P_{\leq N}^x \Big[ P_{\leq N}^x \UNintro[K][+], \big\{ \textup{low-frequency terms}\big\} \Big], \\
    \UNintro[K][+] \Big|_{u=v} &= \hcoup W_K^+, 
    \end{aligned}
    \end{cases}
\end{equation}
where the low-frequency terms in \eqref{intro:eq-ansatz-Up} depend on the terms in our Ansatz for $\VNintro[][]$ \revision{(see 
Definition \ref{ansatz:def-lo}, Definition \ref{ansatz:def-shhl}, \eqref{ansatz:eq-motivation-Up-again}, and Definition \ref{ansatz:def-modulation-equations})}.
Due to \eqref{intro:eq-Ansatz-Up-Sp}, the initial value problem \eqref{intro:eq-ansatz-Up} then suggests an initial value problem for the modulation operator $\SNintro[K][+][k] \colon \R^{1+1}\rightarrow \End(\frkg)$, which is given by 
\begin{equation}\label{intro:eq-ansatz-modulation}
\partial_v \SNintro[K][+][k] = \rho_{\leq N}^2(k)\ad\Big( \big\{ \textup{low-frequency terms}\big\} \Big) 
\circ \SNintro[K][+][k], \qquad \SNintro[K][+][k]\big|_{u=v} = \Id_\frkg.
\end{equation}
In \eqref{intro:eq-ansatz-modulation}, $\rho_{\leq N}(k)$ is the symbol of $P_{\leq N}^x$ and $\ad(\cdot)$ is the adjoint map on $\frkg$. The following three aspects of \eqref{intro:eq-Ansatz-Up-Sp}, \eqref{intro:eq-ansatz-Up}, and \eqref{intro:eq-ansatz-modulation} are most important for our argument: 
\begin{enumerate}[label=(\roman*),leftmargin=5ex]
    \item (Absence of nonlinear smoothing)
    The nonlinear part of the evolution of $\UNintro[K][+]$ can be written as the $v$-integral of the right-hand side of \eqref{intro:eq-ansatz-Up}. Thus, while $\UNintro[K][+]$ exhibits nonlinear smoothing in the $v$-variable, there is no nonlinear smoothing in the $u$-variable. This absence of nonlinear smoothing is the reason for the introduction of the modulated linear wave $\UNintro[K][+]$. In terms of the modulation operators $\SNintro[K][+][k]$, the absence of nonlinear smoothing stems from the fact that the difference between $\SNintro[K][+][k]$ and $\Id_\frkg$ cannot be bounded by an inverse power of $K$. 
    \item\label{intro:item-dependence} (Dependence on $k\in \Z_K$) In the regime $K\sim N$, the $\rho_{\leq N}^2(k)$-factor in \eqref{intro:eq-ansatz-modulation} necessitates the dependence of the $\SNintro[K][+][k]$-operators on $k\in \Z_K$. This dependence on $k\in \Z_K$ is dangerous since, in expressions such as \eqref{intro:eq-Ansatz-Up-Sp}, one has to prevent $(\SNintro[K][+][k])_{k\in \Z_K}$ from destroying\footnote{To see the potential issue, consider the following example: Let $(a_k)_{k\in \Z_K}$ be a deterministic sequence and let $(g_k)_{k\in \Z_K}$ be a sequence of independent, standard, real-valued Gaussians. Then, the sum $\sum_{k\in \Z_K} s_k g_k$ exhibits square-root cancellation if $s_k=a_k$, but may exhibit no cancellation if $s_k=a_k g_k$.} the randomness in $(G_k^+)_{k\in \Z_K}$. To address this, we prove a new chaos estimate with dependent coefficients (see Section \ref{section:chaos}), which relies on $\ell^1$-bounds for the discrete derivative of $\SNintro[K][+][k]$ in $k\in \Z_K$. 
    \item\label{intro:item-orthogonality} (Orthogonality) For any $A\in \frkg$, the adjoint map $\ad(A)\colon \frkg\rightarrow \frkg$ is skew-symmetric. As a consequence of the skew-symmetry and \eqref{intro:eq-ansatz-modulation}, we obtain that 
    \begin{equation*}
        \SNintro[K][+][k](u,v) \colon \frkg \rightarrow \frkg
    \end{equation*}
    is an orthogonal transformation\footnote{Due to additional frequency-truncations, the orthogonality actually only holds for the pure modulation operators, see Definition \ref{ansatz:def-pure} and Proposition \ref{modulation:prop-properties}.} for all $K\in \dyadic$, $k\in \Z_k$, and $u,v\in \R$. The orthogonality of $\SNintro[K][+][k]$ is an important aspect of our argument, as it is used to show that the Killing-renormalization cancels a resonant interaction that involves modulated linear waves. A similar orthogonality property has been used in the deterministic theory of wave maps  \cite{KT98,T01I,T01II}, but has not been used in the probabilistic theory of \cite{BLS21}. \revision{Orthogonality has previously also been used in the context of random nonlinear Schrödinger equations in \cite{B97,DNY21} but, unlike in our setting, the orthogonality used in \cite{B97,DNY21} is not pointwise in space-time.}
\end{enumerate}

In the following, we also need the modulated linear wave $\VNintro[][-]$, which is the direct counterpart of $\UNintro[][+]$. It consists of the dyadic components $\VNintro[M][-]$, which are given by 
\begin{equation}
\VNintro[M][-](u,v) = \hcoup \sum_{m\in \Z_M} \SNintro[M][-][m](u,v) G_m^- e^{imv}. 
\end{equation}

\begin{remark}[On absence of nonlinear smoothing]
The absence of nonlinear smoothing has been previously encountered in the context of random wave equations in \cite{B21,BLS21,BR23,CGI22,KLS20}. However, the aspects discussed in \ref{intro:item-dependence} and \ref{intro:item-orthogonality}, i.e., the dependence on $k\in \Z_K$ and orthogonality, have not previously been addressed in this context.
\end{remark}

\begin{remark}[On $k$-dependence]
The $k$-dependence of $\SNintro[K][+][k]$ was not present in \cite{BLS21}, since \cite{BLS21} only included a frequency-truncation in the initial data and not in the nonlinearity. Without this dependence, one could simply write $\UNintro[K][+]$ as 
\begin{equation*}
\UNintro[K][+] = \hcoup \SNintro[K][+] \Big( \sum_{k\in \Z_K} G_k^+ e^{iku} \Big),
\end{equation*}
and then the randomness in $(G_k^+)_{k\in \Z_K}$ easily yields square-root cancellation. In other works on random dispersive and wave equations \cite{B21,DNY19}, the Ansatz guarantees that the coefficients $(\SNintro[K][+][k])_{k\in \Z_K}$ and Gaussians $(G_k^+)_{k\in \Z_K}$ in expressions such as $\sum_{k\in \Z_K}\SNintro[K][+][k] G_k^+ e^{iku}$ are probabilistically independent, and then square-root cancellation can again easily be obtained. However, it would be difficult to enforce this probabilistic independence in our setting,  since the modulation operators $\SNintro[K][+][k]$ need to absorb certain high$\times$high$\rightarrow$low-interactions (see e.g. Definition \ref{ansatz:def-shhl} and Definition \ref{ansatz:def-modulation-equations}).
\end{remark}

\subsubsection{The modulated bilinear wave $\UNintro[][+-]$:} 
The modulated bilinear wave $\UNintro[][+-]$ is defined as
\begin{equation}\label{overview:eq-modulated-bilinear}
\UNintro[][+-] = \sum_{\substack{K \simeq_\delta M}} \UNintro[K,M][+-], \quad \text{where} \quad \UNintro[K,M][+-] =  \coup \sum_{\substack{k \in \Z_{K}}} \sum_{m \in \Z_M}
\bigg[ \SNintro[K][+][k] G_k^+ \, e^{iku},  \SNintro[M][-][m] G_m^- \frac{e^{imv}}{im} \bigg]_{\leq N}, 
\end{equation}
where $[\cdot,\cdot]_{\leq N}$ is the frequency-truncated Lie-bracket from \eqref{prelim:eq-Lie-bracket-truncated}. 
We note that $\UNintro[K,M][+-]$ is primarily supported on $u$-frequencies $\sim K$ and $v$-frequencies $\sim M$. 
The condition $K\simeq_\delta M$, which is defined in \eqref{prelim:eq-simeq-delta} below, then implies that $\UNintro[K,M][+-]$ is supported on high frequencies in both the $u$ and $v$-variable. 
Due to our definition of $\UNintro[K][+]$ and $\VNintro[M][-]$, it follows that $\UNintro[K,M][+-]$ satisfies
\begin{equation*}
\partial_v \UNintro[K,M][+-] \simeq \big[ \UNintro[K][+], \VNintro[M][-] \big]_{\leq N}.
\end{equation*}
Thus, $\UNintro[K,M][+-]$ captures the interactions of the two modulated linear waves $\UNintro[K][+]$ and $\VNintro[M][-]$, which travel in opposite directions. 
 The most difficult interactions between $\UNintro[][+-]$ and terms in $\VNintro[][]$, i.e., terms in the Ansatz \eqref{intro:eq-ansatz-V}, are high$\times$high$\rightarrow$low-interactions in the $v$-variable. Such interactions are responsible for both the Killing-renormalization (Subsection \ref{section:overview-killing}) and the Jacobi errors (Subsection \ref{section:overview-jacobi}), which are two of the main aspects of this article.

 \subsubsection{\protect{The terms $\UNintro[][-]$, $\UNintro[][+\fs]$ and $\UNintro[][\fs-]$:}}\label{section:overview-other-random-terms}
 \revision{We now briefly discuss the terms $\UNintro[][-]$, $\UNintro[][+\fs]$ and $\UNintro[][\fs-]$, but refer to Section \ref{section:ansatz-heuristic} and Section \ref{section:ansatz-rigorous} for more detailed discussions. From \eqref{intro:eq-Ansatz-Up-Sp}, one sees that the modulated linear wave $\UNintro[][+]$ only has high frequencies in the $u$-variable and that its random structure is inherited from $W^+$. The reversed modulated linear wave $\UNintro[][-]$ has similar properties, but the role of the $u$ and $v$-variables is reversed. Indeed, $\UNintro[][-]$ only has high frequencies in the $v$-variable and its random structure is inherited from $W^-$.} 

\revision{
 From \eqref{overview:eq-modulated-bilinear}, one sees that the modulated bilinear wave $\UNintro[][+-]$ has high frequencies in both variables. Moreover, its random structure in the $u$ and $v$-variables is inherited from $W^+$ and $W^-$, respectively. Like $\UNintro[][+-]$, the mixed modulated objects $\UNintro[][+\fs]$ and $\UNintro[][\fs-]$ also have high frequencies in both variables. However, they only exhibit random structure in one variable, and their behavior in the remaining variable is mostly dictated by the nonlinear, smooth remainders $\UNintro[][\fs]$ and $\VNintro[][\fs]$. 
 }

\subsection{Killing-renormalization and orthogonality}\label{section:overview-killing}

We now describe the motivation behind the Killing-renormalization in \eqref{intro:eq-system-U-V-Kil}. For expository purposes, we first describe it at the level of Picard iterates, and only afterward describe it at the level of our Ansatz. In the latter context, we also describe the significance of the orthogonality of $\SNintro[K][+][k]$ and $\SNintro[M][-][m]$, which was previously discussed in Subsection \ref{section:overview-ansatz}.

\subsubsection{The Killing-renormalization at the level of Picard-iterates.}
In the Picard iteration for \eqref{intro:eq-system-U-V-Kil}, one encounters the nonlinear term 
\begin{equation}\label{intro:eq-Picard}
\coup^{\frac{3}{2}} \Big[ \big[ W^+(u), (\Int^v_{u\rightarrow v} W^-)(u,v) \big]_{\leq N}, W^-(v) \Big]_{\leq N},
\end{equation}
where $\Int^v_{u\rightarrow v}$ denotes the integral in the $v$-variable from $u$ to $v$ (see Definition \ref{prelim:def-integral}). In order to further examine this expression, we let $\widecheck{\rho}_{\leq N}$ be the kernel of the Littlewood-Paley operator $P_{\leq N}$. Using the kernel $\widecheck{\rho}_{\leq N}$, the interaction \eqref{intro:eq-Picard} can be written as 
\begin{equation}\label{intro:eq-Picard-expanded}
\coup^{\frac{3}{2}} P_{\leq N}^x \int_{\R} \dy \big( \widecheck{\rho}_{\leq N} \ast \widecheck{\rho}_{\leq N} \big)(y) \Big[ 
\big[ (P_{\leq N} W^+)(u-y), (\Int^v_{u\rightarrow v} P_{\leq N} W^-)(u-y,v-y)\big], (P_{\leq N}W^-)(v) \Big].
\end{equation}
In \eqref{intro:eq-Picard-expanded}, there is the potential for a probabilistic resonance \revision{(see Sections \ref{section:ansatz-remainder}, \ref{section:killing-tensor}, and \ref{section:jacobi})} between 
$\Int^v_{u\rightarrow v} P_{\leq N} W^-$ and $P_{\leq N}W^-$. In fact, it holds\footnote{The formula \eqref{intro:eq-Picard-Killing} is an approximation since the contribution of a boundary term of the $\Int^v_{u\rightarrow v}$-integral has been omitted.} for all $E\in \frkg$ that 
\begin{equation}\label{intro:eq-Picard-Killing}
\mathbb{E}  \Big[ 
\big[ E , (\Int^v_{u\rightarrow v} P_{\leq N} W^-)(u-y,v-y)\big], (P_{\leq N}W^-)(v) \Big] \simeq \Cf^{(\Nscript)}(y) \Kil(E),
\end{equation}
where $\Cf^{(\Nscript)}$ is the covariance function of frequency-truncated, real-valued white noise and its antiderivative and $\Kil\colon \frkg\rightarrow \frkg$ is the Killing map (Definition \ref{prelim:def-Killing} and Definition \ref{ansatz:def-Killing}). We note that, since $\Cf^{(\Nscript)}(0)=0$, the contribution of the resonant term \eqref{intro:eq-Picard-Killing} to \eqref{intro:eq-Picard} only occurs due to the $P_{\leq N}^x$-operators and therefore does not appear in the wave maps equation with frequency-truncated initial data \eqref{intro:eq-system-U-V-N-data}. Since the contribution of the resonant part \eqref{intro:eq-Picard-Killing} to \eqref{intro:eq-Picard} is non-negligible, one may be afraid that the solutions of \eqref{intro:eq-system-U-V-N-data} and \eqref{intro:eq-system-U-V-Kil} do not converge to the same limit. However, we are rescued by the Killing-renormalization, which cancels the contribution of the resonant part \eqref{intro:eq-Picard-Killing} to \eqref{intro:eq-Picard}. It is defined as 
\begin{equation}
\Renorm[N] \big( W^+ \big)(u) = P_{\leq N}^x \int_{\R} \dy \big( \widecheck{\rho}_{\leq N} \ast \widecheck{\rho}_{\leq N} \big)(y) \Cf^{(\Nscript)}(y) \Kil \big( P_{\leq N}^x W^+ \big)(u-y) 
\end{equation}
and therefore the contribution of the resonant term \eqref{intro:eq-Picard-Killing} does not occur in the renormalized interaction
\begin{equation}\label{intro:eq-Picard-Killing-renormalized}
\coup^{\frac{3}{2}}  \Big[ \big[ W^+(u), (\Int^v_{u\rightarrow v} W^-)(u,v) \big]_{\leq N}, W^-(v) \Big]_{\leq N} - \coup  \Renorm[N] \big( \hcoup W^+\big)(u).
\end{equation}

\subsubsection{The Killing-renormalization at the level of our Ansatz.}
While our previous discussion explains how the Killing-renormalization acts as a counter-term in \eqref{intro:eq-Picard}, this alone is insufficient for our well-posedness theory. Instead of \eqref{intro:eq-Picard-Killing-renormalized}, we need to control 
\begin{equation}\label{intro:eq-interaction-Killing}
\begin{aligned}
&\Big[ \UNintro[][+-], \VNintro[][-] \Big]_{\leq N} - \coup \Renorm[N] \UNintro[][+]  \\
\simeq\, &  \Big[ \big[ \UNintro[][+], \Int^v_{u\rightarrow v} \VNintro[][-] \big]_{\leq N}, \VNintro[][-] \Big]_{\leq N} - \coup \Renorm[N] \UNintro[][+].
\end{aligned}
\end{equation}
As we discussed in Subsection \ref{section:overview-ansatz}, the wave maps equation lacks nonlinear smoothing, and therefore the difference between $\VNintro[][-]$ and $\hcoup W^-$ is no smoother than $\hcoup W^-$ itself. However, 
despite the lack of smoothing, 
the Killing-renormalization still eliminates the resonant term in  \eqref{intro:eq-interaction-Killing}. The reason is the orthogonality of $\SNintro[M][-][m]\colon \frkg\rightarrow \frkg$, which guarantees that the probabilistic resonances between $\hcoup P_{\leq N}^x W^-$ and $\hcoup \Int^v_{u\rightarrow v} P_{\leq N}^x W^-$ and between 
$P_{\leq N}^x\VNintro[][-]$ and $\Int^v_{u\rightarrow v} P_{\leq N}^x \VNintro[][-]$ are identical. As a result, the Killing-renormalization not only cancels the resonant term in \eqref{intro:eq-Picard-Killing-renormalized} but also in \eqref{intro:eq-interaction-Killing}.

\subsection{Jacobi errors and a Bourgain-Bulut argument}\label{section:overview-jacobi}

In this subsection, we first introduce the so-called Jacobi errors and then discuss how the Jacobi errors are controlled using a Bourgain-Bulut argument \cite{BB14}. The starting point of our discussion is a sum of two interactions, 
which is given by\footnote{In terms of our Ansatz from \eqref{intro:eq-ansatz-U}-\eqref{intro:eq-ansatz-V}, the interactions in \eqref{intro:eq-jacobi-interactions} stem from the Lie bracket of $\UNintro[][+-]$ and $\VNintro[][\fs]$ and the Lie bracket of $\UNintro[][+\fs]$ and $\VNintro[][-]$.}
\begin{equation}\label{intro:eq-jacobi-interactions}
\Big[ \Big[ \UNintro[K][+] , \Int^v_{u\rightarrow v} \VNintro[M][-] \Big]_{\leq N}, \VNintro[][\fs]\Big]_{\leq N} + \Big[ \Big[ \UNintro[K][+], \Int^v_{u\rightarrow v} \VNintro[][\fs] \Big]_{\leq N}, \VNintro[M][-]\Big]_{\leq N}, 
\end{equation}
where $K,M\in \dyadic$. In the following discussion, we are primarily concerned with the case $K,M \sim N$ and high$\times$high$\rightarrow$low-interactions in the $v$-variable between the $\VNintro[M][-]$ and $\VNintro[][\fs]$-terms.
Since the nonlinear remainder $\VNintro[][\fs]$ then appears at high frequencies, one may think that \eqref{intro:eq-jacobi-interactions} can be absorbed back into the nonlinear remainders. However, similar as in \cite[Section 1.3.2]{BLS21}, one quickly realizes that this is not possible using direct estimates. 
In \cite{BLS21}, the interactions in \eqref{intro:eq-jacobi-interactions} were instead absorbed back into the modulated linear waves $\UNintro[K][+]$ and $\VNintro[K][-]$. In the setting of this article, however, there is a new difficulty which was not present in \cite{BLS21}: If we absorb \eqref{intro:eq-jacobi-interactions} into $\UNintro[K][+]$, we then still have to preserve the orthogonality of the modulation operators $(\SNintro[K][+][k])_{k\in \Z_K}$ in \eqref{intro:eq-ansatz-Up}. Since the orthogonality relies on the skew-symmetry of the Lie bracket, we would therefore like to approximate the interactions in \eqref{intro:eq-jacobi-interactions} by a term of the form\footnote{The reason for the $(P_{\leq N}^x)^j$-operator in \eqref{intro:eq-jacobi-desired} is not obvious, but will become clear in Section \ref{section:jacobi}.}
\begin{equation}\label{intro:eq-jacobi-desired}
\Big[ (P_{\leq N}^x)^j \UNintro[K][+], \big\{ \, \textup{low-frequency terms} \big\}   \Big]_{\leq N},
\end{equation}
where $j\in \mathbb{N}$. In order to approximate \eqref{intro:eq-jacobi-interactions} using a term of the form \eqref{intro:eq-jacobi-desired}, it would be sufficient to control the
 so-called Jacobi errors\footnote{In Definition \ref{ansatz:def-jacobi-errors}, we give a precise definition of all Jacobi errors, which slightly differs from \eqref{intro:eq-jacobi-errors}.}, which are morally given by 
\begin{equation}\label{intro:eq-jacobi-errors}
\begin{aligned}
\JcbNErr[K,M][+]  &= \Big[ \Big[ \UNintro[K][+] , \Int^v_{u\rightarrow v} \VNintro[M][-] \Big]_{\leq N}, \VNintro[][\fs]\Big]_{\leq N}  \\
&+ \Big[ \Big[ \UNintro[K][+], \Int^v_{u\rightarrow v} \VNintro[][\fs] \Big]_{\leq N}, \VNintro[M][-]\Big]_{\leq N} \\ 
&- \Big[ (P_{\leq N}^x)^2 \UNintro[K][+], \Big[ \Int^v_{u\rightarrow v} \VNintro[M][-], \VNintro[][\fs] \Big]_{\leq N} \Big]_{\leq N}.
\end{aligned}
\end{equation}
The reasons why one may attempt to control the Jacobi errors in \eqref{intro:eq-jacobi-errors} are two-fold: First, as we are concerned with high$\times$high$\rightarrow$low-interactions in the $v$-variable, integration by parts allows us to move $\Int^v_{u\rightarrow v}$-operators between the $\VNintro[M][-]$ and $\VNintro[][\fs]$-terms. Second, the Lie bracket of $\frkg$ satisfies the Jacobi identity 
\begin{equation}\label{intro:eq-jacobi-identity}
\big[ \big[ A, B \big], C \big] 
+ \big[ \big[ B, C \big], A \big] 
+ \big[ \big[ C, A \big], B \big] = 0  
\end{equation}
for all $A,B,C\in \frkg$, which suggests the presence of a cancellation between the three terms in \eqref{intro:eq-jacobi-errors}. Unfortunately, our attempts to control the Jacobi errors \eqref{intro:eq-jacobi-errors} in this fashion were unsuccessful. The reason is that the Littlewood-Paley operators $P_{\leq N}^x$ in the frequency-truncated Lie bracket $[\cdot, \cdot]_{\leq N}$ break the Jacobi identity and, at least in the regime $K,M \sim N$, this cannot be repaired. To summarize, we are unable to use direct estimates to absorb  \eqref{intro:eq-jacobi-interactions}  in either the modulated linear waves or in the nonlinear remainders. \\ 

In order to control  \eqref{intro:eq-jacobi-interactions}, we instead use a Bourgain-Bulut argument \cite{BB14}, which is non-perturbative and relies on (almost) invariance. To this end, we first approximate the $v$-integrals in \eqref{intro:eq-jacobi-interactions} by $x$-integrals. This is possible since 
$\VNintro[M][-]$ and $\VNintro[][\fs]$ enter at much higher $v$-frequencies than $u$-frequencies (see Lemma \ref{jacobi:lem-replacement}). Furthermore, it is also necessary since the almost invariance of the Gibbs measure only yields information about the behavior at a fixed time, but no information about the joint behavior at different times. Then, we insert $\Pbd$-operators, which project the solution to frequency scales near the frequency boundary. This is possible since the Jacobi identity is only broken due to the $P_{\leq N}^x$-operators and can therefore be restored at frequency scales much lower than $N$. Finally, we recombine all terms from our Ansatz \eqref{intro:eq-ansatz-V} and thereby restore the full solution $\VNintro[][]$. This is necessary since the almost invariance of the Gibbs measure only yields information on the full solution $\VNintro[][]$, but yields no information on the individual terms in our Ansatz from \eqref{intro:eq-ansatz-V}. All in all, we then arrive at an expression of the form\footnote{While \eqref{intro:eq-jacobi-Bourgain-Bulut} is morally correct, the precise formulation includes additional frequency-projections, see e.g. Definition \ref{jacobi:def-chhl} and Proposition \ref{jacobi:prop-main}.}
\begin{equation}\label{intro:eq-jacobi-Bourgain-Bulut}
\Big[  \Big[ \UNintro[K][+], 
\Int^x_{0\rightarrow x} \Pbd \VNintro[][] \Big]_{\leq N}, \Pbd \VNintro[][] \Big]_{\leq N} - \coup \Renorm[N][]\UNintro[K][+],
\end{equation}
 By using the almost invariance of the Gibbs measure under the finite-dimensional approximation \eqref{intro:eq-system-U-V-Kil} and by taking into account the Killing-renormalization, we can then control the high$\times$high$\rightarrow$low-interaction between $\Int^x_{0\rightarrow x}\Pbd \VNintro[][]$ and $\Pbd \VNintro[][]$. In particular,  we obtain estimates which are sufficient to absorb \eqref{intro:eq-jacobi-Bourgain-Bulut} into the smooth remainders. 

\begin{remark}
The term \eqref{intro:eq-interaction-Killing} has been treated in two different ways: In Subsection \ref{section:overview-killing}, it was considered individually. In contrast, in Subsection \ref{section:overview-jacobi}, it has been included in \eqref{intro:eq-jacobi-Bourgain-Bulut}. The reason for this is that our argument will combine a short-time well-posedness result on time-scales $\sim N^{-\epsilon}$, where $\epsilon>0$ is a small parameter, and a local well-posedness result on time-scales $\sim 1$, see Section \ref{section:main}.  
\end{remark}

\begin{remark}\label{intro:rem-Bourgain-Bulut}
While we refer to our argument to control \eqref{intro:eq-jacobi-Bourgain-Bulut} as a Bourgain-Bulut argument, we believe that there is an important conceptual difference between our setting and the setting in \cite[Section 4]{BB14}. In \cite[Section 4]{BB14}, the local well-posedness relies on invariance and no local well-posedness results using only contraction-mapping arguments are obtained. In contrast, our argument proves the local well-posedness of the wave maps equation \eqref{intro:eq-system-U-V-N-data} with frequency-truncated initial data using a contraction-mapping argument. The (almost) invariance of the Gibbs measure is only used for the local well-posedness of the wave maps equation \eqref{intro:eq-system-U-V-Kil} with a frequency-truncated nonlinearity. However, the only reason for introducing \eqref{intro:eq-system-U-V-Kil} is to prove global well-posedness and invariance (Theorem \ref{intro:thm-rigorous-A-B}), and should therefore not be considered as part of our local theory.\\
A further difference between our setting and \cite{BB14} is that the Gibbs measures of \eqref{intro:eq-system-U-V-Kil} are given by white noise and therefore do not depend on $N\in \dyadic$. In contrast, the Gibbs measures of the frequency-truncated nonlinear wave equation in \cite[(1.2)]{BB14} depend on $N\in \dyadic$.
\end{remark}

\begin{remark}
The Bourgain-Bulut argument from \cite{BB14} was also used in recent work of Sun, Tzvetkov, and Xu. 
While the limiting solution $u$ from \cite[(4.3)]{STX22} is estimated using perturbative methods, the estimates of the sequence of solutions $u_N$ from \cite[(4.1)]{STX22} rely on a Bourgain-Bulut argument. 
\end{remark}

\subsection{Almost invariance and conservative structures}\label{section:overview-almost}

As part of the proof of Theorem \ref{intro:thm-rigorous-A-B}, we need to prove the almost invariance of the Gibbs measure under the finite-dimensional approximation \eqref{intro:eq-system-U-V-Kil}. For this problem, it is convenient to transform the finite-dimensional approximation \eqref{intro:eq-system-U-V-Kil} into Cartesian coordinates, which leads to the initial value problem
\begin{equation}\label{intro:eq-system-A-B-Kil}
\begin{cases}
\begin{aligned}
\partial_t \A[N] &= \partial_x \B[N], \\ 
\partial_t \B[N] &= \partial_x \A[N] - P_{\leq N}^x \Big[ P_{\leq N}^x \A[N], P_{\leq N}^x \B[N] \Big] + 2\coup \Renorm[N] \A[N], \\
\A[N](0) &= \sqrt{8\coup} W_0, \qquad \B[N](0) =  \sqrt{8\coup} W_1.   
\end{aligned}
\end{cases}
\end{equation}
We let $(\Aflow[\Nscript],\Bflow[\Nscript])$ be the flow of \eqref{intro:eq-system-A-B-Kil}, let $\mu^{(\coup)}$ be the 
law of $(\sqrt{8\coup} W_0, \sqrt{8\coup} W_1)$, and let 
\begin{equation*}
\mu^{(\Nscript,\coup)}_{t}:= \big( \Aflow[\Nscript](t),\Bflow[\Nscript](t)\big)_\# \mu^{(\coup)}
\end{equation*}
be the push-forward of $\mu^{(\coup)}$ under the map  $(\Aflow[\Nscript](t),\Bflow[\Nscript](t))$.
For all $T\geq 1$, we then want to prove that 
\begin{equation}\label{intro:eq-TV-bound}
\sup_{t\in [0,T]}\Big\| \mu^{(\Nscript,\coup)}_{t} - \mu^{(\coup)} \Big\|_{\TV} \lesssim_T N^{-\epsilon},     
\end{equation}
where $\epsilon>0$ is a small parameter. In order to prove \eqref{intro:eq-TV-bound}, we use an explicit formula for the Radon-Nikodym derivative of $\mu^{(\Nscript,\coup)}_{t}$ with respect to $\mu^{(\coup)}$, which is given by
\begin{equation}\label{intro:eq-explicit-Radon-Nikodym}
\frac{\mathrm{d}\mu^{(\Nscript,\coup)}_{t}}{\mathrm{d}\mu^{(\coup)}} = 
\exp \bigg( \frac{1}{4} \int_{-t}^{0} \int_{\bT} \, \ds \dx \, 
\big \langle \Renorm[N][] \A[N](s,x), \B[N](s,x) \big \rangle_\frkg \bigg),
\end{equation}
where $(\A[N],\B[N])$ are as in \eqref{intro:eq-system-A-B-Kil}. 
A similar explicit formula has been used by Debussche and Tsutsumi in \cite{DT21} to prove the quasi-invariance of Gaussian measures under the flow of nonlinear Schr\"{o}dinger equations. We now split the rest of this subsection into two parts: In the first part, we discuss the proof of the explicit formula \eqref{intro:eq-explicit-Radon-Nikodym}. In the second part, we discuss our estimates of the right-hand side in \eqref{intro:eq-explicit-Radon-Nikodym}. 

\subsubsection{Conservative structures}

In order to prove the explicit formula for the Radon-Nikodym derivative \eqref{intro:eq-explicit-Radon-Nikodym}, 
we need to understand how much of the Hamiltonian structure of the wave maps equation \eqref{intro:eq-system-A-B} is still present in the finite-dimensional approximation \eqref{intro:eq-system-A-B-Kil}. In the scalar setting \cite{B94,B96}, this is often a rather simple question since the finite-dimensional approximation is often the Hamiltonian equation corresponding to a frequency-truncated Hamiltonian. In our setting, however, this question is more complicated. The reason is that our finite-dimensional approximation \eqref{intro:eq-system-A-B-Kil} is not obtained from a frequency-truncated Hamiltonian, but rather from a frequency-truncated nonlinear structure on our state space. \\

In the following discussion, we use the inverse temperature $\beta>0$ rather than the temperature $\coup$ from \eqref{intro:eq-coup}. We also consider an abstract Hamiltonian evolution equation on a state space $\introstate$, which can be formally written as 
\begin{equation}\label{intro:eq-abstract-Hamiltonian}
\dot{\phi} = \Jc_\phi \mathrm{d}H_\phi. 
\end{equation}
In \eqref{intro:eq-abstract-Hamiltonian}, $H\colon \introstate \rightarrow \mathbb{R}$ is a Hamiltonian, $\mathrm{d}H$ is the exterior derivative of $H$, and $\Jc\colon T^\ast \introstate\rightarrow T\introstate$ is a bundle homomorphism. We also introduce an inverse temperature $\beta>0$, a reference measure\footnote{In Section \ref{section:conservative}, we primarily work with (volume) forms rather than measures, but we ignore this difference in this overview.} $\nu$ on $\introstate$, and a Gibbs measure $\mu_\beta$ on $\introstate$, which is given by
\begin{equation}\label{intro:eq-abstract-Gibbs}
\mu_\beta = \mathcal{Z}_\beta^{-1} \exp\big(-\beta H\big)\, \nu. 
\end{equation}
For the explicit expressions of $\introstate$, $H$, $\Jc$, $\nu$, and $\mu_\beta$ for the wave maps equation, we refer to Section \ref{section:conservative-wave-maps}. To obtain a finite-dimensional approximation of \eqref{intro:eq-abstract-Hamiltonian}, we let $\introstate^{(\Nscript)}$, $\Jc^{(\Nscript)}$, and $\nu^{(\Nscript)}$ be finite-dimensional approximations of $\introstate$, $\Jc$, and $\nu$, respectively. We let $H^{(\Nscript)}$ be the restriction of $H$ to $\introstate^{(\Nscript)}\subseteq \introstate$, which does not involve any finite-dimensional approximations, and let $\mu_\beta^{(\Nscript)}$ be the Gibbs measure\footnote{In this discussion, we work with a frequency-truncated Gibbs measure, whereas the explicit formula \eqref{intro:eq-explicit-Radon-Nikodym} concerns the full Gibbs measure. Since the initial data at frequencies larger than $N$ evolves linearly under  \eqref{intro:eq-system-A-B-Kil}, this can easily be remedied.} given by
\begin{equation}\label{intro:eq-abstract-Gibbs-N}
\mu^{(\Nscript)}_\beta = \mathcal{Z}_\beta^{-1} \exp\big(-\beta H^{(\Nscript)} \big)\, \nu^{(\Nscript)}. 
\end{equation}
Furthermore, we let $K^{(\Nscript)}$ be a vector field on $\introstate^{(\Nscript)}$, which will be used to capture the Killing-renormalization in \eqref{intro:eq-system-A-B-Kil}. Then, we can introduce the finite-dimensional Hamiltonian evolution equation
\begin{equation}\label{intro:eq-abstract-Hamiltonian-N}
\dot{\phi}^{(\Nscript)} = \Jc^{(\Nscript)}_{\phi^{(\Nscriptscript)}} \mathrm{d}H^{(\Nscript)}_{\phi^{(\Nscriptscript)}} + K^{(\Nscript)}. 
\end{equation}
We emphasize that \eqref{intro:eq-abstract-Hamiltonian-N} involves a finite-dimensional approximation $\Jc^{(\Nscript)}$ of the bundle homomorphism $\Jc$, which is often not needed in scalar settings. 
The explicit formula for the Radon-Nikodym derivative \eqref{intro:eq-explicit-Radon-Nikodym} then takes the form
\begin{equation}\label{intro:eq-abstract-explicit-Radon-Nikodym}
\frac{\mathrm{d}\big( \Phi^{(\Nscript)}(t)_\# \mu_\beta^{(\Nscript)}\big)}{\mathrm{d}\mu_\beta^{(\Nscript)}} = \exp \bigg(  \beta \int_{-t}^{0} \ds \, K^{(\Nscript)}_{\phi^{(\Nscriptscript)}(s)}(H^{(\Nscript)}) \bigg), 
\end{equation}
where $\Phi^{(\Nscript)}$ is the flow of \eqref{intro:eq-abstract-Hamiltonian-N}. 
Clearly, the explicit formula \eqref{intro:eq-abstract-explicit-Radon-Nikodym} cannot be obtained without making assumptions on the state space $\introstate^{(\Nscript)}$, bundle homomorphism $\Jc^{(\Nscript)}$,  measure $\nu^{(\Nscript)}$, Hamiltonian $H^{(\Nscript)}$, and vector field $K^{(\Nscript)}$. This is the motivation behind conservative structures, which consist of tuples $(\introstate^{(\Nscript)},\Jc^{(\Nscript)},\nu^{(\Nscript)})$. Conservative structures form a convenient framework for the finite-dimensional approximation of Hamiltonian evolution equations since they are more flexible than symplectic structures (Lemma \ref{structure:lem-almost-symplectic}). In order to obtain \eqref{intro:eq-abstract-explicit-Radon-Nikodym} using only mild assumptions on the Hamiltonian $H^{(\Nscript)}$ and vector-field $K^{(\Nscript)}$, we impose two conditions on conservative structures: The first condition requires that $\Jc^{(\Nscript)}$ is skew-symmetric and the second condition requires that $\Jc^{(\Nscript)}$ preserves the measure $\nu^{(\Nscript)}$, see Definition \ref{structure:def-conservative} and Lemma \ref{structure:lem-conservative-Riemannian}. For a more detailed discussion of conservative structures, we refer the reader to Subsection \ref{section:conservative-abstract}, which can be read independently of the rest of this article. 

\subsubsection{Energy increment}
In order to prove the almost invariance of the Gibbs measures, we want to show that, on an event with high probability, the Radon-Nikodym derivative in \eqref{intro:eq-explicit-Radon-Nikodym} is close to one. This requires us to show that the energy-increment, i.e., the argument of the exponential in \eqref{intro:eq-explicit-Radon-Nikodym}, is small. Using our null-coordinates and null-variables, the energy-increment can be written as 
\begin{equation}\label{intro:eq-energy-increment}
\frac{1}{8} \int_{\Domain_t} \du \dv  \, \big\langle \UNintro[][](u,v) , \Renorm[N]  \VNintro[][](u,v) \big\rangle_\frkg,
\end{equation}
where $\Domain_t=\big\{ (u,v)\in \R^{1+1}\colon (-u+v,u+v)\in [-2t,0]\times [-2\pi,2\pi]\big\}$. In order to control \eqref{intro:eq-energy-increment}, we insert our Ansatz from \eqref{intro:eq-ansatz-U} and \eqref{intro:eq-ansatz-V}, and then separately consider the terms 
\begin{equation}\label{intro:eq-energy-increment-decomposed}
\frac{1}{8} \int_{\Domain_t} \du \dv  \, \big\langle \UNintro[][\ast_1](u,v) , \Renorm[N]  \VNintro[][\ast_2](u,v) \big\rangle_\frkg, 
\end{equation}
where $\ast_1,\ast_2\in \{ +, - , +-, +\fs, \fs-, \fs \}$. For most interactions, \eqref{intro:eq-energy-increment-decomposed} can be controlled using similar arguments as in our well-posedness theory, but some interactions require additional, more delicate ingredients. In our argument for the $(+-)$$\times$$(+-)$-interaction, it is crucial that the covariance function of white noise and its antiderivative (Definition \ref{ansatz:def-Killing}) is odd. In our argument for the $(+)$$\times$$(+-)$ and $(+)$$\times$$(+)$-interactions, we use the symmetry of the Killing map $\Kil$, which is crucial for the Lie-algebraic identities in Lemma \ref{prelim:lem-Killing-vs-Bracket}. To be more precise, we use that 
\begin{equation*}
\big\langle E_a , \Kil \big[ E^a, X \big] \big \rangle_\frkg =0 \qquad \text{for all } X \in \frkg,
\end{equation*}
where $(E_a)$ is an orthonormal basis of $\frkg$.

\begin{remark}[Bootstrap argument]
Since we rely on a Bourgain-Bulut argument (Subsection \ref{section:overview-jacobi}), our proofs of local well-posedness and almost invariance of the Gibbs measure are intertwined. Since the corresponding bootstrap argument is rather technical, we postpone its discussion until Section \ref{section:main}. 
\end{remark}

\subsection{Structure of the article}

We briefly describe the structure of the rest of this article, which is also illustrated in Figure \ref{figure:structure}. 
In Section \ref{section:preliminaries}, we cover preliminaries, which include definitions and basic results from differential geometry, partial differential equations, and probability theory. All the material in this section is standard, and experts can likely skim or entirely skip much of this section. 

In Section \ref{section:ansatz}, which is the heart of this article, we rigorously introduce our Ansatz for \eqref{intro:eq-system-U-V-Kil}. Furthermore, we collect and classify all error terms, which sets the stage for most later sections. During the reading of this article, this section should be consulted repeatedly. 
In Section \ref{section:conservative}, we introduce conservative structures, whose main application is the explicit formula for the Radon-Nikodym derivative \eqref{intro:eq-explicit-Radon-Nikodym}. In Section \ref{section:chaos}, which is rather short, we state and prove a general chaos estimate with dependent coefficients. This chaos estimate will be one of the main ingredients in all of our probabilistic estimates. 

In Section \ref{section:Killing}, Section \ref{section:modulated-mixed}, and  Section \ref{section:modulation}, the main goal is to control the modulated objects, such as $\UNintro[][+]$ and $\UNintro[][-]$, and their interactions. Since the results of the three sections will be used repeatedly in the rest of this article, all three sections should be read in detail. 
In Sections \ref{section:hhl}-\ref{section:energy-increment},  we control all error terms in the remainder equations (Definition \ref{ansatz:def-remainder-equations}) and control the energy increment (Definition \ref{increment:def-energy-increment}). Since there are almost no dependencies between the sections, they can be read independently. The most interesting parts of this block are Section \ref{section:jacobi}, in which we control the Jacobi errors, and Section \ref{section:energy-increment}, in which we control the energy increment. In Section \ref{section:Lipschitz}, we state Lipschitz-variants of the estimates from Sections \ref{section:Killing}-\ref{section:renormalization}. Since the Lipschitz-variants can be obtained using minor modifications of earlier arguments, Section \ref{section:Lipschitz} contains no new ideas. 

In Section \ref{section:lifting}, we first discuss lifts of $\frkg$-valued maps to $\frkG$-valued maps and then discuss $\frkG$-valued Brownian motion, which serves as a preparation for the proof of Theorem \ref{intro:thm-rigorous-phi}.
In Section \ref{section:main}, we then prove both main theorems, i.e., Theorem \ref{intro:thm-rigorous-A-B} and Theorem \ref{intro:thm-rigorous-phi}. The arguments in this last section combine results from almost all earlier sections.

In Appendix \ref{section:appendix-auxiliary}, we record several auxiliary analytical estimates, which are more technical than the analytical estimates in the preliminaries. In Appendix \ref{section:classical}, we recall a classical well-posedness and stability result (Proposition \ref{classical:prop-main}), which will be used for a soft argument. In Appendix \ref{section:appendix-symbolic}, we also include a symbolic index, which may serve as a useful reference.

\tikzset{overviewrect/.style={rectangle, rounded corners, minimum width=2.5cm, minimum
height=1cm,text centered,align=center, draw=black, fill=white!10},
arrow/.style={thick,->,>=stealth}}

\tikzstyle{innerWhite} = [semithick, white,line width=1.4pt, shorten >= 4.5pt]

\begin{figure}[t]
\begin{center}
\scalebox{0.8}{
\begin{tikzpicture}[node distance=2cm,scale=1, every node/.style={scale=1}]

\node (ansatz)[overviewrect]  at (0,0) 
{Section \begin{NoHyper}\ref{section:ansatz}\end{NoHyper}};

\node (chaos)[overviewrect]  at (4.25,-1.5) 
{Section \begin{NoHyper}\ref{section:chaos}\end{NoHyper}};

\node (conservative)[overviewrect]  at (-4.25,-1.5) 
{Section \begin{NoHyper}\ref{section:conservative}\end{NoHyper}};

\begin{scope}[on background layer]
\draw [fill=blue!5] (2,-8) rectangle (-2,-3); 
\end{scope}
\node (killing)[overviewrect]  at (0,-4) 
{Section \begin{NoHyper}\ref{section:Killing}\end{NoHyper}};
\node (mixed)[overviewrect]  at (0,-5.5) 
{Section \begin{NoHyper}\ref{section:modulated-mixed}\end{NoHyper}};
\node (modulation)[overviewrect]  at (0,-7) 
{Section \begin{NoHyper}\ref{section:modulation}\end{NoHyper}};

\begin{scope}[on background layer]
\draw [fill=blue!5] (5.5,-12.5) rectangle (-5.5,-9); 
\end{scope}

\node (hhl)[overviewrect]  at (-3.5,-10) 
{Section \begin{NoHyper}\ref{section:hhl}\end{NoHyper}};

\node (jacobi)[overviewrect]  at (0,-10) 
{Section \begin{NoHyper}\ref{section:jacobi}\end{NoHyper}};

\node (perturbative)[overviewrect]  at (3.5,-10) 
{Section \begin{NoHyper}\ref{section:null}\end{NoHyper}};

\node (structural)[overviewrect]  at (-3.5,-11.5) 
{Section \begin{NoHyper}\ref{section:structural}\end{NoHyper}};

\node (renormalization)[overviewrect]  at (0,-11.5) 
{Section \begin{NoHyper}\ref{section:renormalization}\end{NoHyper}};

\node (energy)[overviewrect]  at (3.5,-11.5) 
{Section \begin{NoHyper}\ref{section:energy-increment}\end{NoHyper}};


\node (lipschitz)[overviewrect]  at (-3.5,-14) 
{Section \begin{NoHyper}\ref{section:Lipschitz}\end{NoHyper}};

\node (lifting)[overviewrect]  at (3.5,-14) 
{Section \begin{NoHyper}\ref{section:lifting}\end{NoHyper}};

\node (main)[overviewrect]  at (0,-15.5) 
{Section \begin{NoHyper}\ref{section:main}\end{NoHyper}};


\draw[arrow] (-1.25,-0.25) -- (-4.25,-0.25)--(conservative);
\draw[arrow] (ansatz) -- (0,-3);

\draw[arrow] (chaos) -- (4.25,-4) -- (2,-4);

\draw[arrow] (killing)--(mixed);
\draw[arrow] (mixed)--(modulation);

\draw[arrow] (0,-8)--(0,-9);

\draw[arrow] (-1.25,0.25) -- (-7,0.25) -- (-7,-15.75)-- (-1.25,-15.75);
\draw[arrow] (conservative) -- (-6.5,-1.5) -- (-6.5,-15.25) -- (-1.25,-15.25);

\draw[arrow] (2,-7) -- (7,-7) -- (7,-15.75) -- (1.25,-15.75);
\draw[arrow] (0,-12.5)--(main);

\draw[arrow] (-2,-7) -- (-6,-7) -- (-6,-14) -- (-4.75,-14);
\draw[arrow] (-3.5,-12.5) -- (-3.5,-13.5);

\draw[arrow] (lifting) -- (0.5,-14)-- (0.5,-15);
\draw[arrow] (lipschitz) -- (-0.5,-14)-- (-0.5,-15);

\end{tikzpicture}
}
\end{center}
\caption{\small{This figure illustrates the relationship between the different sections of this article.}}\label{figure:structure}
\end{figure} 

\section{Preliminaries}\label{section:preliminaries}

In this section, we make preparations for our main argument. In Subsection \ref{section:parameters}, we introduce parameters and general notation which will be used throughout the article. In Subsection \ref{section:function-spaces}, we introduce our function spaces and recall several para-product, commutator, and integral estimates. In Subsection \ref{section:Lie} and \ref{section:probability-theory}, we recall basic definitions and results from Lie algebra and probability theory, respectively.  

\subsection{Parameters and notation}\label{section:parameters}

In this subsection, we introduce our parameters, Littlewood-Paley operators, and cut-off functions.

\subsubsection{Parameters}
We choose \revision{six} parameters $\delta_0,\delta_1,\delta_2,\delta_3,\delta_4,\delta_5 \in (0,1)$ such that
\begin{equation}\label{prelim:eq-parameter-sizes}
\delta_0 \ll 1 \qquad \text{and} \qquad \delta_{j+1} \leq (\delta_{j})^{10} \textup{ for all } 0\leq j\leq 4. 
\end{equation}
We then define our main parameters $\delta\in (0,1)$, $r\in (1/2,1)$, $s\in (0,1/2)$, and $\eta \in (0,1)$ as 
\begin{equation}\label{prelim:eq-parameter-regularities}
\delta:=\delta_0, \qquad r:= 1/2+\delta_1, \qquad s:=1/2-\delta_2, \qquad \text{and} \qquad \eta:=\delta_3. 
\end{equation}
Furthermore, we define parameters $\vartheta \in (0,1)$ and $\varsigma \in (0,1)$ as 
\begin{equation}\label{prelim:eq-parameter-regularities-less}
\vartheta := \delta_4 \qquad \text{and} \qquad \varsigma = \delta_5,
\end{equation}
but they will be used less frequently than the main parameters from \eqref{prelim:eq-parameter-regularities}.
For notational convenience, we also define
\begin{equation*}
 \delta^\prime := \delta +  2 \delta_2 \qquad \text{and} \qquad r^\prime := r - 2\delta_2.
\end{equation*}
Thus, $\deltap$ is slightly bigger than $\delta$ and $r^\prime$ is slightly smaller than $r$. 
Furthermore, we define
\begin{equation}\label{prelim:eq-scrr}
\scrr := 1 - 10\delta. 
\end{equation}
This parameter is much larger than $r$, but serves a similar purpose. While $r$ will be used for high-regularity terms in the remainder equations (see Definition \ref{ansatz:def-remainder-equations}), $\scrr$ will be used for high-regularity terms in the modulation equations (see Definition \ref{ansatz:def-modulation-equations}).
In the following, we often write
\begin{equation}\label{prelim:eq-dependence-delta}
C=C(\delta_\ast) \qquad \text{and} \qquad c=c(\delta_\ast) 
\end{equation}
for sufficiently large and sufficiently small constants depending on $\delta_0,\delta_1,\delta_2$, $\delta_3$, $\delta_4$, and $\delta_5$, respectively. 

\begin{remark}
The two equivalent sets of parameters $(\delta,r,s,\eta,\vartheta,\varsigma)$ and $(\delta_0,\delta_1,\delta_2,\delta_3,\delta_4,\delta_5)$ serve two different purposes. The first set of parameters $(\delta,r,s,\eta,\vartheta,\varsigma)$ will be used to describe terms in our Ansatz and state our main estimates (see e.g. Section \ref{section:ansatz}). Since we use six different letters, it should be relatively easy to distinguish between the six different parameters. The second set of parameters $(\delta_0,\delta_1,\delta_2,\delta_3,\delta_4,\delta_5)$ will be used whenever we need to estimate complicated expressions involving our parameters, since \eqref{prelim:eq-parameter-sizes} makes it easy to remember their relative sizes. 
\end{remark}

\subsubsection{Dyadic scales}

We denote the set of dyadic integers by 
$\dyadic = \{ 2^n \colon n\in \mathbb{N}_0 \}$. We also introduce a large dyadic integer $\Nlarge\in \dyadic$ and a set of large dyadic integers $\dyadiclarge \subseteq \dyadic$ as 
\begin{equation}\label{prelim:eq-dyadiclarge}
\Nlarge := \inf \Big\{ N \in \dyadic \colon  N > 2^{100 \delta_4^{-1}} \Big\} 
\qquad \text{and} \qquad 
\dyadiclarge := \Big\{ N \in \dyadic \colon N \geq \Nlarge \Big\}.  
\end{equation}
The objects in \eqref{prelim:eq-dyadiclarge} are used to avoid technical problems with the zero frequency, see e.g. Lemma \ref{prelim:lem-Duhamel-integral} below. For every $M\in \dyadic$, the set of integers $\Z_M \subseteq \Z$ is defined as
\begin{alignat*}{3}
 \Z_M &:= \Big\{ m \in \Z \colon |m|\leq 1 \Big\} 
 \qquad &&\text{if } M=1, \\ 
 \Z_M &:= \Big\{ m \in \Z \colon M/2 < |m|\leq M \Big\} 
  \qquad  &&\text{if } M\geq 2.
\end{alignat*}
For any $K,M\in \dyadic$, we also define
\begin{align}
K \simeq_\delta M \qquad  &\Longleftrightarrow \qquad  \min(K,M) \geq \max(K,M)^{1-\delta}, \label{prelim:eq-simeq-delta} \\
K \ll_\delta M \qquad  &\Longleftrightarrow \qquad  K < M^{1-\delta}, \\ 
K \lesssim_\delta M \qquad  &\Longleftrightarrow \qquad  K \leq  M^{1/(1-\delta)}. 
\end{align}
Given a positive integer $k\geq 1$ and a map $f\colon \big( \dyadic \big)^k \rightarrow \R$, we often use the shorthand notation
\begin{equation}\label{prelim:eq-sum}
\sum_{M_1,\hdots,M_k} f\big( M_1, \hdots, M_k \big) := \sum_{M_1,\hdots,M_k \in \dyadic} f\big( M_1, \hdots, M_k \big).
\end{equation}
In the following, we often restrict the dyadic scales in \eqref{prelim:eq-sum} from $\dyadic$ to $\dyadiclarge$, i.e., we restrict to dyadic integers larger than or equal to $\Nlarge$. For notational convenience, we therefore write
\begin{equation}\label{prelim:eq-sumlarge}
\sumlarge_{M_1,\hdots,M_k} f\big( M_1, \hdots, M_k \big) := \sum_{M_1,\hdots,M_k \in \dyadiclarge} f\big( M_1, \hdots, M_k \big).
\end{equation}

\subsubsection{Harmonic analysis}
We define the one-dimensional torus $\bT$ as $\bT:=\R/ (2\pi \Z)$ and identify its fundamental domain with $(-\pi,\pi]$. 
For any $R\geq 1$, we also define $\bT_R:= \R / (2\pi R \Z)$. 
We let $\rho \colon  \R_\xi \rightarrow [0,1]$ be a smooth even cut-off function satisfying 
\begin{equation}\label{prelim:eq-rho}
\rho|_{[-7/8,7/8]}=1 \quad \text{and} \quad \rho|_{\R\backslash [-1,1]}=0. 
\end{equation}
Furthermore, we define
\begin{equation}\label{prelim:eq-rho-leqN}
\rho_{\leq N}(\xi) := \rho(\xi/N). 
\end{equation}
It is convenient to also define $(\rho_{N})_{N\in \dyadic}$ and $(\rho_{<N})_{N\in \dyadic}$ by 
\begin{alignat}{3}\label{prelim:eq-rhoN}
\rho_1 &:= \rho_{\leq 1},  &
\rho_N &:= \rho_{\leq N} - \rho_{\leq N/2} \quad \text{for all} ~ N \geq 2, \\ 
\rho_{<1} &:= 0,&  \qquad \text{and} \qquad 
\rho_{<N} &:= \rho_{\leq N/2}  \hspace{9.5ex}\text{for all} ~ N \geq 2. 
\end{alignat}

For any $y\in \R$, we define the shift operator $\Theta_y \colon C^\infty_b(\R) \rightarrow C^\infty_b(\R)$ by 
\begin{equation}\label{prelim:eq-shift}
\Theta_y f(x) := f(x-y). 
\end{equation}
Similarly, we define the shift operators $\Theta_y^x,\Theta_y^u,\Theta_y^v\colon C^\infty_b(\R^{1+1}_{u,v}) \rightarrow C^\infty_b(\R^{1+1}_{u,v})$ by 
\begin{equation*}
(\Theta_y^x f)(u,v):= f(u-y,v-y), \quad (\Theta_y^u f)(u,v):= f(u-y,v), \quad \text{and} \quad (\Theta_y^v f)(u,v):= f(u,v-y). 
\end{equation*}
For any $f\in C_b^\infty(\R)$ and $N \in \dyadic$, we define the Littlewood-Paley operators 
\begin{equation}\label{prelim:eq-P-N}
P_N f(x) := \int_{\R}\dy \, \widecheck{\rho}_N(y) (\Theta_y f)(x). 
\end{equation}
The operators $P_{\leq N}$ and $P_{<N}$ are defined similarly as in \eqref{prelim:eq-P-N} but with $\widecheck{\rho}_N$ replaced by $\widecheck{\rho}_{\leq N}$
or $\widecheck{\rho}_{<N}$, respectively.  
For any $w\in \{x,u,v\}$ and any $f\in C_b^\infty(\R_{u,v}^{1+1})$, we also define
\begin{equation*}
P_N^w f(u,v) := \int_{\R} \dy \, \widecheck{\rho}_N(y) (\Theta_y^w f)(u,v). 
\end{equation*}
For notational convenience, we further define
\begin{equation}\label{prelim:eq-P-N-uv}
P^{u,v}_{\leq N}  := P^u_{\leq N} P^v_{\leq N} , \qquad 
P^{u,v}_{N} := P^{u,v}_{\leq N} - P^{u,v}_{<N}, 
\qquad \text{and} \qquad P^{u,v}_{>N} = 1 - P^{u,v}_{\leq N}. 
\end{equation}
We define the Littlewood-Paley operators $(\Pbd)_{N\in \dyadic}$ by
\begin{equation}\label{prelim:eq-pbd}
\Pbd := \sum_{\substack{L\in \dyadic\colon \\ N^{1-2\delta_1}\leq L \leq N}} P_L^x. 
\end{equation}
The $\Pbd$-operators project to frequencies near the frequency-boundary of \eqref{intro:eq-system-A-B-N} at scale $\sim N$, and will be important in our treatment of the Jacobi errors (see Section \ref{section:jacobi}).
In addition to the Littlewood-Paley operators $(P_{N})_{N\in \dyadic}$ from \eqref{prelim:eq-P-N}, we sometimes also need sharp frequency-projection operators. For any $R\geq 1$, any $N\in \dyadic$, and any $f\colon \bT_R \rightarrow \bC$, we define 
\begin{equation}\label{prelim:eq-Psharp}
\Psharp_{R;N} f = \int_{\bT_R} \dy \, \widecheck{\rho}_N^{\,\,\sharp}(y) \Theta^x_y f(x),
\end{equation}
where 
\begin{equation*}
\widecheck{\rho}_N^{\,\,\sharp}(y) = 
\begin{cases}
\frac{1}{2\pi R} \sum_{\ell \in \Z} \mathbf{1} \Big\{ \big| \tfrac{\ell}{R} \big| \leq 1 \Big\} e^{i \frac{\ell}{R} y} 
\qquad \hspace{6.25ex} \text{if} \quad  N=1, \\ 
\frac{1}{2\pi R} \sum_{\ell \in \Z} \mathbf{1} \Big\{ \tfrac{N}{2} < \big| \tfrac{\ell}{R} \big| \leq N \Big\} e^{i \frac{\ell}{R} y} 
\qquad \text{if} \quad N\geq 2.
\end{cases}
\end{equation*}
Finally, for any $N \in \dyadic$ and $\Phi,\Psi\colon \R \rightarrow \frkg$, we write
\begin{equation}\label{prelim:eq-Lie-bracket-truncated}
\big[ \Phi, \Psi \big]_{\leq N} = P_{\leq N}^x \big[ P_{\leq N}^x \Phi, P_{\leq N}^x \Psi \big],
\end{equation}
where $[\cdot,\cdot] \colon \frkg \rightarrow \frkg$ is the Lie-bracket on $\frkg$. Thus, $[\cdot,\cdot]_{\leq N}$ is a frequency-truncated Lie bracket, which will be used repeatedly throughout this article. As it will be important later, we already emphasize that $[\cdot,\cdot]_{\leq N}$ does not satisfy the Jacobi identity for Lie brackets. 

\subsubsection{Time cut-off function and lattice partition}

In our local theory, we will make use of a cut-off function in time. To this end, we introduce the following sets of
functions. 

\begin{definition}[Cut-off functions]\label{prelim:def-cut-off}
We let $\chifixed\colon \R \rightarrow [0,1]$ be a fixed smooth function satisfying 
\begin{equation*}
\chifixed \, \big|_{[-2,2]}= 1 \qquad \text{and} \qquad \chifixed \, \big|_{\R \backslash [-3,3]} =0.
\end{equation*}
Then, we define 
\begin{equation}\label{prelim:eq-cuttilde}
\Cuttilde := \Big\{ \chifixed\big(\kappa t - t_0\big) \colon \kappa \in 2^{-\mathbb{N}_0},\, t_0 \in \R \Big\} \medcup \big\{ 1 \big\} 
\end{equation}
and define 
\begin{equation}\label{prelim:eq-cut}
\Cut := \Big\{ \chifixed( \cdot - \tau) \zeta(\cdot) \colon \tau \in [-1,1], \, \zeta \in \Cuttilde \Big\}.
\end{equation}
\end{definition}

\begin{remark}
We remark that the cut-off functions in $\Cut$ are localized on the unit timescale $\sim 1$ and the cut-off functions in $\Cuttilde$ are localized on time-scales $\gtrsim 1$. Thus, neither $\Cut$ nor $\Cuttilde$ contains cut-off functions supported on small timescales, which will only be accessed via a scaling argument (see Lemma \ref{ansatz:lem-scaling-symmetry}). 
While the cut-off functions from Definition \ref{prelim:def-cut-off} are far from central to this article, they are related to different technical aspects. First, the translation by $\tau \in [-1,1]$ in \eqref{prelim:eq-cut} is necessary to account for the time-shift in Proposition \ref{structure:prop-Gibbs} below. Second, the reason for using a larger class of functions than $\{ \chifixed(\cdot-\tau)\colon \tau \in [-1,1]\}$ is that this would be insufficient for our scaling argument, see e.g. Lemma \ref{ansatz:lem-scaling-symmetry} and Section \ref{section:main}.
\end{remark}

In the following, we also introduce a partition of unity $(\psix)_{x_0\in \Lambda}$, which will be used to decompose $\frkg$-valued white noise (see Lemma \ref{prelim:lem-white-noise-representation}).

\begin{definition}[Lattice partition]\label{prelim:def-lattice-partition}
We define a lattice $\Lambda\subseteq \R$ and finite lattice $\LambdaR \subseteq \R$, where $R\geq 1$, as 
\begin{equation*}
\Lambda := \Big\{ \tfrac{2\pi}{32} \ell \colon \ell \in \Z\Big \} 
\qquad \text{and} \qquad 
\LambdaR := \Big\{  \tfrac{2\pi}{32} \ell \colon \ell \in \Z, \,-16 R < \ell \leq 16 R \Big\}.
\end{equation*}
We also define $\psi \colon \R \rightarrow [0,\infty)$ by 
\begin{equation}\label{prelim:eq-psi}
\psi(x) := \Big( \sum_{x_0 \in \Lambda} \widetilde{\psi}(x-x_0)^2 \Big)^{-\frac{1}{2}} \widetilde{\psi}(x), \qquad \text{where} \quad \widetilde{\psi}(x) := 
\begin{cases}
\exp \Big( - \frac{1}{(\frac{\pi}{8})^2-|x|^2} \Big) \qquad \text{if } |x|<\tfrac{\pi}{8},\\ 
0 \hspace{22ex} \text{if } |x|\geq \tfrac{\pi}{8}.
\end{cases}
\end{equation}
Furthermore, we define the families $(\psix)_{x_0\in \Lambda}$ and $(\psiRx)_{x_0 \in \LambdaRR}$ by 
\begin{equation*}
\psix(x):= \psi(x-x_0) \qquad \text{and} \qquad 
\psiRx (x) := \sum_{\substack{\ell \in \Z }} \psix \big( x + 2\pi R \ell \big).
\end{equation*}
\end{definition}
As a consequence of Definition \ref{prelim:def-lattice-partition}, we have that
\begin{equation}\label{prelim:eq-psi-properties}
\operatorname{supp}\big( \psix \big) \subseteq \big[ x_0 - \tfrac{\pi}{8}, x_0 + \tfrac{\pi}{8} \big]
\qquad \text{and} \qquad 
\sum_{x_0 \in \Lambda} \psix(x)^2 = \sum_{x_0 \in \LambdaRR} \psiRx(x)^2 =1.
\end{equation}
In the following, it will often be convenient to work with a frequency-truncated version of $(\psiRx)_{x_0\in \LambdaRR}$, which is introduced in the next definition.

\begin{definition}[Frequency-truncation of lattice partition]\label{prelim:def-lattice-partition-truncated}
Let $R\geq 1$ and let $\LambdaR$ and $(\psiRx)_{x_0\in \LambdaRR}$ be as in Definition \ref{prelim:def-lattice-partition}. Then, for each $x_0 \in \LambdaR$ and $L\in \dyadic$, we define 
\begin{equation*}
\psiRxL := P_{\leq L^{\vartheta}} \psiRx. 
\end{equation*}
\end{definition}
We note that, due to the smoothness of $\psiRx$, it holds that 
$\| \psiRx - \psiRxL \|_{L_x^\infty(\R)} \lesssim L^{-100}$. Thus, while working with the frequency-truncated version $\psiRxL$ instead of $\psiRx$ can be convenient, it is never necessary. In order to obtain estimates for terms involving either $\chi\in \Cut$ or $(\psix)_{x\in \Lambda}$, it is convenient to introduce a Schwartz-type norm. For any smooth $\varphi\colon \R \rightarrow \R$, it is defined as 
\begin{equation}\label{prelim:eq-Schwartz}
\big\| \varphi \big\|_{\Schwartz} := \sup_{K\in \dyadic} K^{10} \big\| \langle z \rangle^{10} (P_K \varphi)(z) \big\|_{L^\infty}.
\end{equation}

\subsection{Function spaces, para-products, commutators, and integrals}\label{section:function-spaces}

We first introduce the functional framework which will be used throughout this article. Their definition relies on the Littlewood-Paley operators from Subsection \ref{section:parameters}.

\begin{definition}[H\"{o}lder and product-norms]\label{prelim:def-hoelder} 
For any $\gamma \in \R$ and any smooth, bounded function $f\colon \R \rightarrow \bC$, we define the Hölder-norm
\begin{equation}\label{prelim:eq-Hoelder}
\| f \|_{\C^\gamma_x(\R)} := \sup_{N\geq 1} N^\gamma \| P_N f(x) \|_{L^\infty_x(\R)}. 
\end{equation}
For any $\gamma_1,\gamma_2 \in \R$ and any smooth, bounded function $f\colon \R^{1+1}_{u,v} \rightarrow \bC$, we also define 
\begin{equation}\label{prelim:eq-Hoelder-product}
\| f \|_{\C_u^{\gamma_1} \C_v^{\gamma_2}(\R^{1+1})} 
:= \sup_{N_1, N_2 \geq 1} N_1^{\gamma_1} N_2^{\gamma_2}
\| P_{N_1}^u P_{N_2}^v f(u,v) \|_{L^\infty_{u,v}(\R^{1+1})}. 
\end{equation}
\end{definition}

In the following, we will also work with the $\C_x^{\gamma}$-norm on the torus $\bT_R$ or on compact intervals $J\subseteq \R$. For any smooth $f\colon \bT_R \rightarrow \bC$, we simply write
\begin{equation}\label{prelim:eq-Hoelder-periodic}
\| f \|_{\C_x^{\gamma}(\bT_R)} := \sup_{N\geq 1} N^\gamma \| P_N f \|_{L^\infty_x(\bT_R)},
\end{equation}
which also agrees with the $\C_x^{\gamma}(\R)$-norm of the $2\pi R$-periodic extension of $f$ to $\R$. For any smooth $f\colon J \rightarrow \bC$, we define
\begin{equation}\label{prelim:eq-Hoelder-local}
\| f \|_{\C_x^{\gamma}(J)} := \inf \Big\{ \big\| \widetilde{f} \, \big\|_{\C_x^{\gamma}(\R)} \colon   
\widetilde{f}|_J = f \Big\}.
\end{equation}

In the following, we not only need the norms from \eqref{prelim:eq-Hoelder}-\eqref{prelim:eq-Hoelder-local}, 
but also need the corresponding function spaces. For technical reasons, we would like the function spaces to be separable
(see Remark \ref{main:rem-null-lwp}). We therefore make the following definition, but emphasize that the precise details are not too important for this article. 

\begin{definition}[H\"{o}lder and product-spaces]\label{prelim:def-hoelder-spaces} 
For any $\gamma \in \R \backslash \{0\}$ and $R\geq 1$, we define $\C_x^\gamma(\bT_R)$ as the closure of $C^\infty(\bT_R)$ 
with respect to the norm in \eqref{prelim:eq-Hoelder-periodic}. For any compact interval $J\subseteq \R$, we define $\C_x^\gamma(J)$
as the closure of $C^\infty_c(\R)$ with respect to the norm in \eqref{prelim:eq-Hoelder-local}. Furthermore, we define 
$\C_x^\gamma(\R)$ as the closure of
\begin{equation*}
\textup{span} \Big( \Big( \bigcup_{R\in \mathbb{N}} C^\infty(\bT_R) \Big) \medcup C^\infty_c(\R) \Big)
\end{equation*}
with respect to the norm in \eqref{prelim:eq-Hoelder}. Finally, for all $\gamma_1,\gamma_2 \in \R \backslash \{0\}$, we define
$\Cprod{\gamma_1}{\gamma_2}(\R^{1+1})$ as the closure of 
\begin{equation*}
\textup{span} \Big( \Big\{ f_1(u) f_2(v) \colon f_1 \in \C^{\gamma_1}(\R), \, f_2 \in \C^{\gamma_2}(\R) \Big\} \Big)
\end{equation*}
with respect to the norm in \eqref{prelim:eq-Hoelder-product}.
\end{definition}

\subsubsection{Para-product estimates}\label{section:para-product}

We now introduce our para-product operators and state basic para-product estimates. Since we will use several different para-products, it is convenient to introduce a general framework, which is based on frequency-scale relations. 

\begin{definition}[Frequency-scale relations and para-products]\label{prelim:def-para}
A frequency-scale relation $\Rcal$ is a binary relation on $2^{\mathbb{N}_0}\times 2^{\mathbb{N}_0}$, i.e., pairs of frequency scales. For any 
$f,g \in C^\infty_b(\R)$, we define the corresponding para-product as 
\begin{equation*}
f \pararcal g := \sum_{\substack{ (M,N)\in \Rcal}} P_M f \cdot P_N g. 
\end{equation*}

If $\Rcal_1$ and $\Rcal_2$ are two frequency-scale relations and $f,g\in C^\infty_b(\R^{1+1}_{u,v})$, we define 
\begin{equation*}
f \parauvrcal g := 
\sum_{\substack{ (M_1,N_1)\in \Rcal_1}} 
\sum_{\substack{ (M_2,N_2)\in \Rcal_2}}
P_{M_1}^u P_{M_2}^v f \cdot P^u_{N_1} P^v_{N_2}  g. 
\end{equation*}
\end{definition}

For example, if $\Rcal=\Rcal_1=\Rcal_2=\ll$, then it holds that
\begin{align}
  f \parall g &:= 
  \sum_{\substack{M,N\colon \\ M \ll N}} P_M f \cdot P_N g, \notag \\ 
  f \parauvll g &:= 
  \sum_{\substack{M_1,N_1\colon \\ M_1 \ll N_1 }} 
  \sum_{\substack{M_2,N_2\colon \\ M_2 \ll N_2}}
P_{M_1}^u P_{M_2}^v f \cdot P^u_{N_1} P^v_{N_2}  g. \label{prelim:eq-para-uvll}
\end{align}
In \eqref{prelim:eq-para-uvll}, we added the subscripts $u$ and $v$ in $\ll_u$ and $\ll_v$, respectively. This further emphasizes which binary relation applies to which variable, and therefore serves expository purposes.

\begin{definition}[$\Rcal$-admissible regularities]\label{prelim:def-admissible}
Let $\Rcal$ be a frequency-scale relation and let $\alpha,\beta,\gamma \in \R\backslash \{ 0 \}$ be regularity parameters. Then, the tuple $(\alpha,\beta;\gamma)$ is called $\Rcal$-admissible if 
\begin{equation*}
 \sup_{K} \sum_{\substack{(M,N)\in \Rcal}} \bigg[ \Big(
1\big\{ K \sim M \gg N \big\} + 1 \big\{ K \sim N \gg M \big\}
+ 1 \big\{ M \sim N \gtrsim K \big\} \Big) K^\gamma M^{-\alpha} N^{-\beta}\bigg] \lesssim 1. 
\end{equation*}
\end{definition}

We now give sufficient conditions for admissible tuples for several important frequency-scale relations. 

\begin{lemma}[Admissible regularities]\label{prelim:lem-admissible}
Let $\Rcal$ be a frequency-scale relation, let $\alpha,\beta,\gamma\in \R\backslash \{0\}$, and define $\alpha_-:=\min(\alpha,0)$ and $\beta_-:= \min(\beta,0)$. Then, under either of the following conditions, the tuple $(\alpha,\beta;\gamma)$ is $\Rcal$-admissible:
\begin{align}
\Rcal&=\ll \qquad \text{and} \qquad \gamma\leq \alpha_- + \beta, \\
\Rcal&=\sim, \qquad \gamma \leq \alpha+\beta, \qquad \text{and} \qquad \alpha+\beta>0, \\ 
\Rcal&= \gg \qquad \text{and} \qquad \gamma \leq \alpha+\beta_-, \\
\Rcal&=\hspace{-0.5ex}\times, \qquad \gamma \leq \min(\alpha,\beta), \qquad \text{and} \qquad \alpha+\beta>0. 
\end{align}
\end{lemma}

\begin{remark} In Lemma \ref{prelim:lem-admissible}, $\Rcal=\times$ refers to the complete binary relation, i.e., the relation which connects all frequency-scales. The notation $\times$ is chosen since the operator $\paratimes$ coincides with the usual product. 
\end{remark}

\begin{proof}
The proof follows from a direct calculation (similar as in \cite[Proposition 2.5]{BLS21}) and we therefore omit the standard details. 
\end{proof}

We now state a general para-product estimate which will be used repeatedly throughout this article. 

\begin{lemma}[Para-product estimate]\label{prelim:lem-paraproduct}
Let $\Rcal$ be a frequency-scale relation, let $\alpha,\beta,\gamma \in \R \backslash \{ 0 \}$ be regularity parameters, and assume that $(\alpha,\beta;\gamma)$ is $\Rcal$-admissible. Then, it holds that 
\begin{equation*}
\big\| f \pararcal g \big\|_{\C^\gamma(\R)} 
\lesssim \| f \|_{\C^{\alpha}(\R)} \| g \|_{\C^{\beta}(\R)}. 
\end{equation*}
Furthermore, let $\Rcal_1$ and $\Rcal_2$ be two frequency-scale relations, let $\alpha_1,\alpha_2,\beta_1,\beta_2,\gamma_1,\gamma_2 \in \R \backslash \{ 0 \}$ be regularity parameters, and assume that $(\alpha_j,\beta_j;\gamma_j)$ is $\Rcal_j$-admissible for $j=1,2$. Then, 
\begin{equation*}
    \big\| f \parauvrcal g \big\|_{\Cprod{\gamma_1}{\gamma_2}}
\lesssim \| f \|_{\Cprod{\alpha_1}{\alpha_2}} \| g \|_{\Cprod{\beta_1}{\beta_2}}. 
\end{equation*}
\end{lemma}

\begin{proof}
The desired estimate follows directly from Definition \ref{prelim:def-hoelder} and Definition \ref{prelim:def-admissible}, and we omit the details (see also \cite[Proposition 2.5]{BLS21}). 
\end{proof}

For convenience, we also record the following corollary of Lemma \ref{prelim:lem-paraproduct}. 

\begin{corollary}[Product estimate]\label{prelim:cor-product}
 Let $\alpha_1,\alpha_2,\beta_1,\beta_2,\gamma_1,\gamma_2 \in \R \backslash \{ 0 \}$ be regularity parameters and assume that, for $j=1,2$, 
 \begin{equation*}
\gamma_j \leq \min(\alpha_j, \beta_j ) \qquad \text{and} \qquad 
\alpha_j + \beta_j >0. 
 \end{equation*}
 Then, it holds that 
\begin{equation}\label{prelim:eq-product}
\big\| f g \big\|_{\Cprod{\gamma_1}{\gamma_2}} 
\lesssim \big\| f \big\|_{\Cprod{\alpha_1}{\alpha_2}} \big\| g \big\|_{\Cprod{\beta_1}{\beta_2}}. 
\end{equation}
\end{corollary}

Finally, we record a trilinear (rather than bilinear) para-product estimate, which will be useful in Section \ref{section:modulation}.

\begin{lemma}[Trilinear para-product estimates]\label{prelim:lem-para-product-trilinear}
Let $\alpha,\beta,\gamma,\zeta\in (-1,1)\backslash\{0\}$ satisfy
\begin{equation*}
\alpha+\beta+\gamma > 0, \quad \beta+\gamma<0, \quad \text{and} \quad \zeta \in (-s,s). 
\end{equation*}
For all $f,g,h\colon \R^{1+1}_{u,v} \rightarrow \bC$, it then holds that
\begin{equation*}
\Big\| \big( f \Para[v][ll] g\big) \Para[v][sim] h - f \big( g \Para[v][sim] h \big) \Big\|_{\Cprod{\zeta}{\alpha+\beta+\gamma}} \lesssim \big\| f \big\|_{\Cprod{\zeta}{\alpha}} \big\| g \big\|_{\Cprod{s}{\beta}} \big\| h \big\|_{\Cprod{s}{\gamma}}. 
\end{equation*}
\end{lemma}
This lemma is a minor modification of \cite[Lemma 2.4]{GIP15}
and we omit the proof.

\subsubsection{Integral and trace estimates}\label{section:integral-trace}

We now introduce and estimate integral and trace operators, which both appear in the Duhamel integral for the linear wave equation. 

\begin{definition}[Integral operator]\label{prelim:def-integral}
For all $f\in C^\infty_b(\R)$, we define 
\begin{equation*}
    \Int [f]  (x):= \int_0^x \dy f(y). 
\end{equation*}
In addition, for all $f\in C^\infty_b(\R^{1+1})$, we define
\begin{equation*}
\Int^u[f](u,v) := \int_0^u \du^\prime f(u^\prime,v) \quad \text{and} \quad 
\Int^v[f](u,v) := \int_0^v \dv^\prime  f(u,v^\prime).
\end{equation*}
Finally, for all $a,b\colon \R^{1+1}\rightarrow \R$ and $f\in C^\infty_b(\R^{1+1})$, we also define 
\begin{align*}
\Int^u_{a\rightarrow b}[f](u,v)       := \Int^u[f](b,v) - \Int^u[f](a,v) 
\qquad \text{and} \qquad 
\Int^v_{a\rightarrow b}[f](u,v) := \Int^v[f](u,b) - \Int^v[f](u,a).
\end{align*}
\end{definition}

In the next lemma, we state our estimates for the integral operators from Definition \ref{prelim:def-integral}.

\begin{lemma}[Integral estimate]\label{prelim:lem-integral}
Let $\gamma\in (-1,\infty)\backslash \{0\}$, let $\varphi\in C^\infty_b(\R)$, and let $x_0\in \R$. 
Then, it holds for all $f\in C^\infty_b(\R)$ that 
\begin{equation}\label{prelim:eq-integral}
\big\| \Int [P_{>1}^x f] \big\|_{\C_x^{\gamma+1}}  + \big\|P_{>1}^x \Int [ f] \big\|_{\C_x^{\gamma+1}} \lesssim  \| f \|_{\C_x^\gamma}
\qquad \text{and} \qquad 
\big\| \Int [\varphi(\cdot - x_0) f] \big\|_{\C_x^{\gamma+1}} \lesssim \| \varphi \|_{\Schwartz} \| f \|_{\C_x^\gamma},
\end{equation}
where $\|\cdot\|_{\Schwartz}$ is as in \eqref{prelim:eq-Schwartz}. 
Similarly, let $\gamma_1,\gamma_2 \in (-1,\infty)\backslash \{0\}$ and 
$u_0,v_0\in \R$. Then, it holds for all $f\in C^\infty_b(\R^{1+1})$ that
\begin{alignat*}{3}
\big\| \Int^u[ P^u_{>1} f ] \big\|_{\Cprod{\gamma_1+1}{\gamma_2}}
+ \big\|  P^u_{>1} \Int^u[ f ] \big\|_{\Cprod{\gamma_1+1}{\gamma_2}}
&\lesssim \big\|  f \big\|_{\Cprod{\gamma_1}{\gamma_2}}, 
&\qquad 
\big\| \Int^u[ \varphi(u-u_0) f ] \big\|_{\Cprod{\gamma_1+1}{\gamma_2}}
&\lesssim \| \varphi \|_{\Schwartz} \big\|  f \big\|_{\Cprod{\gamma_1}{\gamma_2}}, \\
\big\| \Int^v[ P^v_{>1} f ] \big\|_{\Cprod{\gamma_1}{\gamma_2+1}}
+ \big\| P^v_{>1} \Int^v[  f ] \big\|_{\Cprod{\gamma_1}{\gamma_2+1}}
&\lesssim \big\|  f \big\|_{\Cprod{\gamma_1}{\gamma_2}}, 
&\qquad 
\big\| \Int^v[ \varphi(v-v_0) f ] \big\|_{\Cprod{\gamma_1}{\gamma_2+1}}
&\lesssim \| \varphi \|_{\Schwartz} \big\|  f \big\|_{\Cprod{\gamma_1}{\gamma_2}}. 
\end{alignat*}
\end{lemma}

\begin{proof}
Since all estimates in Lemma \ref{prelim:lem-integral} are standard, we only sketch the argument. It suffices to treat \eqref{prelim:eq-integral}, since the estimates in the $u$ and $v$-variables are similar. For the $\Int P^x_{>1}$-estimate, the estimate follows by writing $P^x_{>1}$ as $\partial_x \partial_x^{-1} P^x_{>1}$, where $\partial_x^{-1}P^x_{>1}$ is the Fourier multiplier with symbol $(i\xi)^{-1} \rho_{>1}(\xi)$. A similar argument works for $P^x_{>1} \Int$, since the argument then still has to enter at frequencies bounded away from zero. For the $\Int \varphi$-estimate, we refer to 
 \cite[Lemma A.10]{GIP15}. 
\end{proof}

We now turn our attention from integral to trace operators.

\begin{definition}[Trace operator]\label{prelim:def-trace} For any $f\in C^\infty_b(\R^{1+1}_{u,v})$, we define the traces $\Tru f \in C^\infty_b(\R_u)$ and $\Trv f \in C_b^\infty(\R_v)$ by 
\begin{equation*}
\Tru f(u):= f(u,u) \quad \text{and} \quad \Trv f(v):= f(v,v). 
\end{equation*}
\end{definition}

\begin{lemma}[Trace estimate]\label{prelim:lem-trace}
Let $\gamma,\gamma_1,\gamma_2 \in \R \backslash \{0\}$ and $f\in C^\infty_b(\R^{1+1}_{u,v})$. Then, the estimates
\begin{equation*}
\big\| \Tru f \big\|_{C_u^\gamma} \lesssim \| f \|_{\Cprod{\gamma_1}{\gamma_2}} 
\quad \text{and} \quad 
\big\| \Trv f \big\|_{C_v^\gamma} \lesssim \| f \|_{\Cprod{\gamma_1}{\gamma_2}} 
\end{equation*}
hold under either of the following two conditions:
\begin{itemize}
    \item[(i)] (General case) $\gamma_1+\gamma_2>0$ and $\gamma\leq \min(\gamma_1,\gamma_2)$.
    \item[(ii)] (Non-resonant case) $\gamma\leq \min(\gamma_1,\gamma_2,\gamma_1+\gamma_2)$ and $P_M^u P_N^v f=0$ for all $M\sim N$. 
\end{itemize}
\end{lemma}
This lemma can be found in \cite[Lemma 2.21]{BLS21}.

\begin{lemma}[Duhamel integral estimate]\label{prelim:lem-Duhamel-integral}
Let $\gamma_1,\gamma_2\in (-1,1)\backslash\{0\}$ and let $\varphi\in C^\infty_b(\R)$.  If $\gamma_1\leq \gamma_2+1$, then it holds for all $f\in C^\infty_b(\R^{1+1})$ that 
\begin{align*}
\big\| \Int^v_{u\rightarrow v}\big[ P^v_{>1} f \big] \big\|_{\Cprod{\gamma_1}{\gamma_2+1}} 
+ \big\| P^v_{>1} \Int^v_{u\rightarrow v}\big[  f \big] \big\|_{\Cprod{\gamma_1}{\gamma_2+1}} 
&\lesssim  \big\| f \big\|_{\Cprod{\gamma_1}{\gamma_2}}, \\ 
\big\| \Int^v_{u\rightarrow v}\big[ \varphi(v-u) f \big] \big\|_{\Cprod{\gamma_1}{\gamma_2+1}} &\lesssim \big\| \varphi \big\|_{\Schwartz} \big\| f \big\|_{\Cprod{\gamma_1}{\gamma_2}}. 
\end{align*}
Alternatively, if $\gamma_2 \leq \gamma_1+1$, then it also holds that
\begin{align*}
\big\| \Int^u_{v\rightarrow u}\big[ P^u_{>1} f \big] \big\|_{\Cprod{\gamma_1+1}{\gamma_2}} 
+ \big\| P^u_{>1}  \Int^u_{v\rightarrow u}\big[ f \big] \big\|_{\Cprod{\gamma_1+1}{\gamma_2}} 
&\lesssim  \big\| f \big\|_{\Cprod{\gamma_1}{\gamma_2}}, \\ 
\big\| \Int^u_{v\rightarrow u}\big[ \varphi(v-u) f \big] \big\|_{\Cprod{\gamma_1+1}{\gamma_2}} 
&\lesssim \big\| \varphi \big\|_{\Schwartz} \big\| f \big\|_{\Cprod{\gamma_1}{\gamma_2}}. 
\end{align*}
\end{lemma}

By writing $\Int^v_{u\rightarrow v}$ as $\Int^v- \Tr^u \Int^v$, this estimate  follows from 
 Lemma \ref{prelim:lem-integral} and Lemma \ref{prelim:lem-trace}. 

\subsubsection{Commutator estimates}\label{section:commutator}

We now consider commutators and commutator estimates. 

\begin{lemma}[Basic commutator estimate]\label{prelim:lem-commutator-basic}
Let $f,g\in L_x^\infty(\R)$, let $M,N\in \dyadic$, and let $y\in \R$. Then, it holds that 
\begin{equation*}
\big\| \tau_y \big( P_M f \, P_N g \big) - P_M f \tau_y P_N g \big\|_{L_x^\infty} \lesssim M |y| \, \big\| P_M f \big\|_{L_x^\infty} \big\| P_N g \big\|_{L^\infty_x},
\end{equation*}
where $\tau_y$ is the shift-operator from \eqref{prelim:eq-shift}. Similarly, if  $K\in \dyadic$, it holds that
\begin{equation*}
\big\| P_K \big( P_M f \, P_N g \big) - P_M f \, P_K P_N g \big\|_{L_x^\infty} \lesssim M N^{-1} \, \big\| P_M f \big\|_{L_x^\infty} \big\| P_N g \big\|_{L^\infty_x}. 
\end{equation*}
\end{lemma}
This follows from standard commutator estimates, see e.g. \cite[Section 2.10]{BCD11}. We now list several consequences of the basic commutator estimates, which are closer to our functional framework.

\begin{lemma}[Commutator estimate for $P_{\leq N}^x$]\label{prelim:commutator-PNx}
Let $\alpha_1,\alpha_2,\beta_1,\beta_2,\gamma_1,\gamma_2 \in \R \backslash \{0\}$ be regularity parameters and let $\Rcal_1$ and $\Rcal_2$ be two frequency-scale relations. For $j=1,2$, assume that $(\alpha_j,\beta_j;\gamma_j)$ is $\Rcal_j$-admissible. Furthermore, let $N\in \dyadic$. Then, it holds that 
\begin{equation*}
 \big\| P_{\leq N}^x  \big( f \parauvrcal g \big)
 - \big( P_{\leq N}^x f \big) \parauvrcal g \big\|_{\Cprod{\gamma_1}{\gamma_2}} 
 \lesssim N^{-1} \big\| f \big\|_{\Cprod{\alpha_1}{\alpha_2}} 
 \Big( \, \big\| g \big\|_{\Cprod{\beta_1+1}{\beta_2}} 
 + \big\| g \big\|_{\Cprod{\beta_1}{\beta_2+1}} \Big). 
\end{equation*}
\end{lemma}

\begin{proof}
Since paraproduct operators commute with translations, it holds that 
\begin{align*}
&P_{\leq N}^x  \big( f \parauvrcal g \big)
 - \big( P_{\leq N}^x f \big) \parauvrcal g \\
 =& \, \int_\R \dy \,  \widecheck{\rho}_{\leq N}(y) \, 
 \big( \Theta^x_y f \parauvrcal \Theta^x_y g \big) 
 - \int_\R \dy\,  \widecheck{\rho}_{\leq N}(y) \, 
 \big( \Theta^x_y f \parauvrcal  g \big) \\
 =&  \, \int_\R \dy \,  \widecheck{\rho}_{\leq N}(y) \, 
 \Big( \Theta^x_y f \parauvrcal \big( \Theta^x_y g - g \big) \Big). 
\end{align*}
Using the triangle inequality and the paraproduct estimate (Lemma \ref{prelim:lem-paraproduct}), it follows that
\begin{align*}
&\big\| P_{\leq N}^x  \big( f \parauvrcal g \big)
 - \big( P_{\leq N}^x f \big) \parauvrcal g \big\|_{\Cprod{\gamma_1}{\gamma_2}} \\
 \leq& \, \int_\R \dy \,  \big|\widecheck{\rho}_{\leq N}(y)\big| \, \Big\| 
 \Theta^x_y f \parauvrcal \big( \Theta^x_y g - g \big)  \Big\|_{\Cprod{\gamma_1}{\gamma_2}} \\ 
 \lesssim& \, \int_\R \dy \,  \big|\widecheck{\rho}_{\leq N}(y)\big| \, \big\|  \Theta^x_y f \big\|_{\Cprod{\alpha_1}{\alpha_2}} 
 \big\| \Theta^x_y g - g \big\|_{\Cprod{\beta_1}{\beta_2}}. 
\end{align*}
After using
\begin{equation*}
\big\|  \Theta^x_y f \big\|_{\Cprod{\alpha_1}{\alpha_2}}  \lesssim 
\big\|  f \big\|_{\Cprod{\alpha_1}{\alpha_2}} 
\qquad \text{and} \qquad 
\big\| \Theta^x_y g - g \big\|_{\Cprod{\beta_1}{\beta_2}} \lesssim |y|  \Big( \, \big\| g \big\|_{\Cprod{\beta_1+1}{\beta_2}} 
 + \big\| g \big\|_{\Cprod{\beta_1}{\beta_2+1}} \Big),
\end{equation*}
this easily yields the desired estimate. 
\end{proof}

Our next lemma is not concerned with a commutator estimate and instead states that certain operators commute exactly. 

\begin{lemma}\label{prelim:lem-commutativity}
For all $y\in \R$, $\Theta^x_y$ commutes with $\Int^v_{u\rightarrow v}$ and $\Int^u_{v\rightarrow u}$. Furthermore, $P_{\leq N}^x$ commutes with 
$\Int^v_{u\rightarrow v}$ and $\Int^u_{v\rightarrow u}$ for all $N\in \dyadic$.
\end{lemma}

\begin{proof}
Due to symmetry in the $u$ and $v$-variables and the definition of $P_{\leq N}^x$, it suffices to prove that $\Theta^x_y$ commutes with $\Int^v_{u\rightarrow v}$. For any smooth $f\colon \R_{u,v}^{1+1}\rightarrow \C$, it holds that
\begin{align*}
&\hspace{3ex} \big(\Theta^x_y \Int^v_{u\rightarrow v} f \big)(u,v) =  \big(  \Int^v_{u\rightarrow v} f \big)(u-y,v-y) =  \int\displaylimits_{u-y}^{v-y} \dz f(u-y,z) 
=\int\displaylimits_{u}^{v} \dz^\prime f(u-y,z^\prime-y)  \\ 
&=  \int\displaylimits_{u}^{v} \dz^\prime (\Theta^x_y f)(u,z^\prime)  
= \big( \Int^v_{u\rightarrow v} \Theta^x_y f\big)(u,v),
\end{align*}        
which yields the desired claim. 
\end{proof}

\subsection{Lie algebras}\label{section:Lie}

We first recall that our Lie group is denoted by $\frkG$ and the corresponding Lie algebra is denoted by $\frkg$. We also recall that, since $\frkG$ is equipped with a bi-invariant Riemannian metric, $\frkg$ comes with a natural inner product. Since the Riemannian metric on $\frkG$ is bi-invariant, the inner product and Lie bracket on $\frkg$ are compatible, i.e., the following identity is satisfied.
\begin{lemma}[Cyclicity] \label{appendix:lem-cyclic}
For all $X,Y,Z\in \frkg$, it holds that 
\begin{equation*}
\big \langle [X,Y], Z \big \rangle = \big \langle [Y,Z], X \big \rangle.  
\end{equation*}
\end{lemma}

For a proof of Lemma \ref{appendix:lem-cyclic}, we refer the reader to \cite[Proposition 21.8]{GQ20}. We now define the space of endomorphisms on $\frkg$ as 
\begin{equation}\label{prelim:eq-endomorphism}
\operatorname{End}(\frkg):= \big\{ \Phi \colon \frkg \rightarrow \frkg \big| \, \Phi \textup{ is linear} \big\}. 
\end{equation}
While many of the linear transformations on $\frkg$ in this article are Lie algebra homomorphisms, i.e., preserve the Lie bracket, this is not required in \eqref{prelim:eq-endomorphism}. 
For any $X\in \frkg$, we define the adjoint map $\Ad(X) \in \operatorname{End}(\frkg)$ by 
\begin{equation}\label{prelim:eq-adjoint-map}
\Ad(X)\colon \frkg \rightarrow \frkg, \, Y \mapsto [X,Y]. 
\end{equation}
Using the adjoint map, the Jacobi identity for the Lie algebra $\frkg$ can be written as 
\begin{equation}\label{prelim:eq-jacobi}
\Ad (X)  \Ad (Y) - \Ad(Y) \Ad(X) = \Ad ( [ X,Y] ).   
\end{equation} 
We recall that, due to Lemma \ref{appendix:lem-cyclic}, $\Ad(X)$ is a skew-symmetric operator on $\frkg$. More precisely, it holds that 
\begin{equation}\label{prelim:eq-Hermitian-adjoint}
\Ad(X)^\ast = - \Ad(X),
\end{equation}
where $\Ad(X)^\ast$ is the Hermitian \revision{conjugate} of the linear map $\Ad(X)$. For any $d\geq 1$, we further define 
\begin{equation*}
\frkg^{\otimes d} := \frkg \otimes \hdots \otimes \frkg, 
\end{equation*}
i.e., $\frkg^{\otimes d}$ is the $d$-fold tensor product of $\frkg$. For any permutation $\sigma \colon \{1,\hdots, d\} \rightarrow \{1,\hdots,d\}$, we also define the shuffle operator $\shuffle{\sigma}\colon \frkg^{\otimes d} \rightarrow \frkg^{\otimes d}$ by 
\begin{equation}\label{prelim:eq-shuffle} 
\shuffle{\sigma}\big( X_1 \otimes \hdots \otimes X_d \big) 
= X_{\sigma(1)} \otimes \hdots \otimes X_{\sigma(d)} \qquad \textup{for all } X_1,\hdots,X_d \in \frkg. 
\end{equation}

In order to work in coordinates of $\frkg$, we now further fix an orthonormal basis
\begin{equation}\label{prelim:eq-ONB}
\ONB = (E_a)_{a=1}^{\dim \frkg} \subseteq \frkg
\end{equation}
\revision{In the following, we will sometimes write $E^a:= \delta^{ab} E_b$, where the repeated index $b$ is summed over.} Equipped with the orthonormal basis from \eqref{prelim:eq-ONB}, we can now define the Casimir, which is a classical object from the study of Lie algebras.

\begin{definition}[Casimir]\label{prelim:def-casimir}
We define the Casmir of the Lie algebra $\frkg$ as 
\begin{equation*}
\Cas := E^a \otimes E_a \in \frkg \otimes \frkg,
\end{equation*}
\revision{where the repeated index $a$ is summed over.}
\end{definition}

Due to the next lemma, Definition \ref{prelim:def-casimir} does not depend on our choice of the orthonormal basis. 

\begin{lemma}[Coordinate-invariance of the Casimir]\label{prelim:lem-casimir}
Let $S\colon \frkg \rightarrow \frkg$ be a linear transformation. Then, it holds that 
\begin{equation}\label{prelim:eq-casimir}
S E^a \otimes S E_a - \Cas = \big( S S^\ast - \Id_\frkg \big)^{ab} \big( E_a \otimes E_b \big). 
\end{equation}
In particular, if $S\colon \frkg \rightarrow \frkg$ is orthogonal, then it holds that
$S E^a \otimes S E_a = \Cas$.
\end{lemma}

\begin{proof}
We write $(S^{ab})_{a,b=1,\hdots,\dim\frkg}$ for the coordinate matrix of $S$ with respect to the orthonormal basis $(E_a)_{a=1}^{\dim\frkg}$.
Since $(E_a)_{a=1}^{\dim\frkg}$ is orthonormal, the coordinate matrix of the Hermitian \revision{conjugate} $S^\ast$ is the transpose of the coordinate matrix of $S$. It then follows that 
\begin{equation*}
S E^a \otimes S E_a 
=\delta_{ab} \big( S E^a \otimes S E^b \big)
= \delta_{ab} S^{ca} S^{db} \big( E_c \otimes E_d \big) 
= \delta_{ab} S^{ca} (S^\ast)^{bd} \big( E_c \otimes E_d \big)
= (S S^\ast )^{cd} \big( E_c \otimes E_d \big).
\end{equation*}
By relabeling the summation indices $(c,d)$ as $(a,b)$ and using Definition \ref{prelim:def-casimir}, we obtain \eqref{prelim:eq-casimir}.
\end{proof}

\begin{definition}[Killing map]\label{prelim:def-Killing}
We define a linear map $\Kil \colon \frkg \rightarrow \frkg$ by
\begin{equation}\label{prelim:eq-Killing}
\Kil(X) := \big[ \big[ X , E^a \big], E_a \big]. 
\end{equation}
\end{definition}

The Killing map is the unique linear map corresponding to the Killing form, which is a central bilinear form in the study of Lie algebras (see e.g. \cite[Section 21]{GQ20}). The link between the Killing map and Killing form is part of the following lemma.

\begin{lemma}[Killing map and Killing form]\label{prelim:lem-Killing-form}
For all $X,Y\in \frkg$, it holds that 
\begin{equation}\label{prelim:eq-Killing-symmetry}
\big\langle \Kil X , Y \big\rangle = \big\langle X, \Kil Y \big\rangle. 
\end{equation}
Furthermore, it holds that 
\begin{equation}\label{prelim:eq-Killing-form}
\big\langle \Kil X , Y \big\rangle = \operatorname{tr}\big( \Ad(X) \Ad(Y) \big). 
\end{equation}
\end{lemma}
In the literature, the Killing form is usually defined as the right-hand side of \eqref{prelim:eq-Killing-form}. 

\begin{proof}
Due to the cyclicity of the trace, \eqref{prelim:eq-Killing-symmetry} easily follows from \eqref{prelim:eq-Killing-form}. As a result, it suffices to prove \eqref{prelim:eq-Killing-symmetry}. Using Lemma \ref{appendix:lem-cyclic}, the symmetry of the inner product, and the skew-symmetry of the Lie bracket, it holds that
\begin{align*}
\big\langle \Kil X , Y \big\rangle  
= \delta^{ab} \big\langle \big[ \big[ X, E_a \big], E_b \big], Y \big\rangle 
= \delta^{ab} \big\langle \big[ E_b, Y \big], \big[ X, E_a \big] \big\rangle 
= - \delta^{ab} \big\langle \Ad(X)E_a, \Ad(Y)E_b \big\rangle. 
\end{align*}
Since $(E_a)_{a=1}^{\dim\frkg}$ is an orthonormal basis, it holds that
\begin{align*}
- \delta^{ab} \big\langle \Ad(X)E_a, \Ad(Y)E_b \big\rangle
&= - \delta^{ab} \delta^{cd} \big\langle E_c, \Ad(X)E_a \big\rangle \,  
\big\langle E_d , \Ad(Y) E_b \big\rangle \\ 
&= - \delta^{ab} \delta^{cd} \Ad(X)_{ca} \Ad(Y)_{db}. 
\end{align*}
Since $\Ad(Y)$ is skew-symmetric, it follows that
\begin{equation*}
- \delta^{ab} \delta^{cd} \Ad(X)_{ca} \Ad(Y)_{db}
= -\Ad(X)_{ca} \Ad(Y)^{ca} = \Ad(X)_{ca} \Ad(Y)^{ac} =  \operatorname{tr}\big( \Ad(X) \Ad(Y) \big), 
\end{equation*}
which agrees with the right-hand side of \eqref{prelim:eq-Killing-form}.
\end{proof}

\begin{lemma}[Lie brackets involving the Killing map]\label{prelim:lem-Killing-vs-Bracket}
It holds that 
\begin{equation}\label{prelim:eq-Killing-vs-Bracket-e1}
\big[ \Kil E_a , E^a \big] =0.
\end{equation}
Furthermore, we have for all $X\in \frkg$ that
\begin{equation}\label{prelim:eq-Killing-vs-Bracket-e2}
\big\langle E_a , \Kil \big[ E^a  , X \big] \big\rangle =0. 
\end{equation}
\end{lemma}

\begin{proof} The first identity follows directly from the symmetry of the Killing map (Lemma \ref{prelim:lem-Killing-form}). To prove the second identity \eqref{prelim:eq-Killing-vs-Bracket-e2}, we combine Lemma \ref{appendix:lem-cyclic} and the symmetry of the Killing map, which yield
\begin{equation*}
\big\langle E_a , \Kil \big[ E_a , X \big] \big\rangle = 
\big\langle \Kil E_a , \big[ E_a , X \big] \big\rangle 
= \big\langle \big[ \Kil E_a , E^a \big], X \big \rangle. 
\end{equation*}
Thus, \eqref{prelim:eq-Killing-vs-Bracket-e2} can be deduced from \eqref{prelim:eq-Killing-vs-Bracket-e1}. 
\end{proof}
In this article, the Killing map and Casimir play similar roles, which is due to the following lemma.

\begin{lemma}[Link between Killing map and Casimir]\label{prelim:lem-Killing-Casimir}
Let $\Tc\colon \frkg \otimes \frkg \otimes \frkg \rightarrow \frkg$ be the unique linear map satisfying 
\begin{equation*}
\Tc \big( X \otimes Y \otimes Z \big) = \big[ \big[ X, Y \big] , Z \big]
\end{equation*}
for all $X,Y,Z\in \frkg$. Then, it holds that 
\begin{equation}\label{prelim:eq-Killing-Casimir}
\Tc \big( X \otimes \Cas \big) = \Kil X. 
\end{equation}
\end{lemma}

\begin{proof} Using Definition \ref{prelim:def-casimir} and Definition \ref{prelim:def-Killing}, it holds that
\begin{equation*}
\Tc \big( X \otimes \Cas \big) = \Tc \big( X \otimes E_a \otimes E^a \big)
= \big[ \big[ X , E_a \big], E^a \big] = \Kil X. 
\end{equation*}
This yields the desired identity \eqref{prelim:eq-Killing-Casimir}.  
\end{proof}

Since we make extensive use of Fourier analysis, we often work with complex-valued functions. Due to this, we cannot always work with the real Lie algebra $\frkg$ and also have to introduce its complexification, which is defined by
\begin{equation*}
\frkg_{\mathbb{C}}:= \frkg \otimes \mathbb{C}.
\end{equation*}
The complexification $\frkg_{\mathbb{C}}$ can be identified with $\frkg \oplus \frkg$, whose elements are written as 
\begin{equation*}
G = \big( \Re G , \Im G \big)
\end{equation*}
and which is equipped with the addition and Lie bracket
\begin{align}
G + H &:= \big(\Re G + \Re H, \Im G + \Im H\big), \\
\big[ G , H \big] &:= \Big( \big[ \Re G, \Re H\big] - \big[ \Im G, \Im H\big],
\big[ \Re G, \Im H\big] + \big[ \Im G, \Re H \big] \Big). \label{prelim:eq-gC-Lie-bracket}
\end{align}

\subsection{Probability theory}\label{section:probability-theory}

We recall several definitions and estimates from probability theory. 

\subsubsection{Total variation distance and Wasserstein metric}

We first recall the definition of the total variation distance of probability measures.

\begin{definition}\label{prelim:def-tv-distance}
Let $(\Omega,\Ec)$ be a measurable space and let $\mu$ and $\nu$ be two probability measures on $(\Omega,\Ec)$. Then, the total variation distance of $\mu$ and $\nu$ is defined as 
\begin{equation*}
\big\| \mu - \nu \big\|_{\TV} 
= 2 \sup_{A\in \Ec} \big| \mu(A) - \nu(A)\big|. 
\end{equation*}
\end{definition}

In Section \ref{section:main}, it will sometimes be convenient to work with a weaker metric than the $\TV$-distance. We therefore also introduce the following Wasserstein metric. 

\begin{definition}[Wasserstein metric]\label{prelim:def-Wasserstein}
Let $R\geq 1$ and let $\mu^{(\Rscript)}$ and $\nu^{(\Rscript)}$ be two probability measures on $(\C_x^{s-1}\times \C_x^{s-1})(\bT_R \rightarrow \frkg^2)$, which is equipped with the corresponding Borel $\sigma$-algebra. Then, we define their Wasserstein distance as 
\begin{align*}
&\,\,  \WassersteinR \big( \mu^{(\Rscript)}, \nu^{(\Rscript)}\big) \\ 
:=&\,  \inf_{\gamma^{(\Rscriptscript)}}
\int \min \Big( \Big\| \big(W^{(\Rscript)}_0,W^{(\Rscript)}_1\big)- \big(\widetilde{W}^{(\Rscript)}_0,\widetilde{W}^{(\Rscript)}_1\big) \Big\|_{\C^{s-1}\times \C^{s-1}}, \, 1 \Big) \, 
\mathrm{d}\gamma^{(\Rscript)}\big(W^{(\Rscript)}_0,W^{(\Rscript)}_1,\widetilde{W}^{(\Rscript)}_0,\widetilde{W}^{(\Rscript)}_1\big),
\end{align*}
where the infimum is taken over all couplings of $\mu^{(\Rscript)}$ and $\nu^{(\Rscript)}$. 
\end{definition}

In the following lemma, we list the basic properties of the Wasserstein metric $\WassersteinR$.

\begin{lemma}[Basic properties of Wasserstein metric]\label{prelim:lem-Wasserstein-properties}
Let $R\geq 1$ and let $\mu^{(\Rscript)}$ and $\nu^{(\Rscript)}$ be two probability measures on $(\C_x^{s-1}\times \C_x^{s-1})(\bT_R \rightarrow \frkg^2)$. Then, we have the following properties: 
\begin{enumerate}[label=(\roman*)]
    \item\label{prelim:item-Wasserstein-TV} It holds that $\WassersteinR(\mu^{(\Rscript)},\nu^{(\Rscript)})\leq \| \mu^{(\Rscript)} - \nu^{(\Rscript)}\|_{\TV}$.
    \item\label{prelim:item-Wasserstein-coupling} For any coupling $\gamma^{(\Rscript)}$ of $\mu^{(\Rscript)}$ and $\nu^{(\Rscript)}$ and any $\epsilon \in (0,1)$, it holds that
    \begin{equation*}
    \WassersteinR \big( \mu^{(\Rscript)}, \nu^{(\Rscript)}\big) \leq \epsilon + 
    \gamma^{(\Rscript)}\Big(\Big\| \big(W^{(\Rscript)}_0,W^{(\Rscript)}_1\big)- \big(\widetilde{W}^{(\Rscript)}_0,\widetilde{W}^{(\Rscript)}_1\big) \Big\|_{\C^{s-1}\times \C^{s-1}} \geq \epsilon \Big). 
    \end{equation*}
    \item\label{prelim:item-Wasserstein-Lipschitz} 
    For any $R^\prime \geq 1$ and any Lipschitz function $F\colon (\C_x^{s-1}\times \C_x^{s-1})(\bT_R\rightarrow \frkg^2) \rightarrow (\C_x^{s-1}\times \C_x^{s-1})(\bT_{R^\prime} \rightarrow \frkg^2)$, it holds that 
    \begin{equation*}
    \WassersteinRp \big( F_\sharp \mu^{(\Rscript)}, F_\sharp \nu^{(\Rscript)}\big)
    \leq \operatorname{Lip}(F) \, \WassersteinR \big( \mu^{(\Rscript)},  \nu^{(\Rscript)}\big),
    \end{equation*}
    where $F_\sharp$ denotes the push-forward under $F$.
\end{enumerate}
\end{lemma}

\begin{proof}
The estimate in \ref{prelim:item-Wasserstein-TV} is standard, see e.g. \cite[Theorem 6.15 and Case 6.16]{V09}. The estimate in \ref{prelim:item-Wasserstein-coupling} follows from the trivial estimate 
$\min \big( z, 1 \big) \leq \epsilon + \mathbf{1}\{ z\geq \epsilon \}$ for all $z\in [0,\infty)$. The estimate in \ref{prelim:item-Wasserstein-Lipschitz} holds since, for any coupling $\gamma^{(\Rscript)}$ of $\mu^{(\Rscript)}$ and $\nu^{(\Rscript)}$, $(F\otimes F)_\sharp \gamma^{(\Rscript)}$ is a coupling of $F_\sharp \mu^{(\Rscript)}$ and $F_\sharp \nu^{(\Rscript)}$.
\end{proof}

\subsubsection{White noise}

We now first recall the definition of white noise and then recall a Fourier-based representation of white noise.

\begin{definition}[Real-valued white noise]\label{prelim:def-white-noise-real-valued}
A random, tempered distribution $w^{(\coup)}\colon \R\rightarrow \R$ is called a real-valued white noise at temperature $\coup>0$ if it satisfies the following two properties: 
\begin{enumerate}[label=(\roman*)]
    \item\label{prelim:item-white-noise-1} For all $k\geq 1$ and all Schwartz functions $\varphi_1,\hdots,\varphi_k\in S(\R)$, 
    \begin{equation*}
    \bigg( \int_{\R} \dx \, w^{(\coup)}(x) \varphi_j(x) \bigg)_{j=1}^k
    \end{equation*}
    is a mean-zero, Gaussian random vector in $\R^k$. 
    \item\label{prelim:item-white-noise-2} For all Schwartz functions $\varphi,\psi\in S(\R)$, it holds that 
    \begin{equation*}
    \E \bigg[ \bigg(  \int_{\R} \dx \, w^{(\coup)}(x) \varphi(x) \bigg) 
    \bigg(  \int_{\R} \dy \, w^{(\coup)}(y) \psi(y) \bigg) \bigg] = \coup \int_{\R} \dx \varphi(x) \psi(x). 
    \end{equation*}
\end{enumerate}
For any $R\geq 1$, the definition of real-valued, $2\pi R$-periodic white noise $w^{(\Rscript,\coup)}\colon \bT_R \rightarrow \R$ at temperature $\coup>0$ is similar, but with the domain $\R$ replaced with $\bT_R$. 
\end{definition}

In the next definition, we generalize real-valued to $\frkg$-valued white noise. 

\begin{definition}[$\frkg$-valued white noise]\label{prelim:def-white-noise-g-valued} 
A random, tempered distribution $W^{(\coup)}\colon \R \rightarrow \frkg$ is called a $\frkg$-valued white noise at temperature $\coup>0$ if 
\begin{equation}\label{prelim:eq-white-noise-g-valued-coordinates}
\Big( \langle W^{(\coup)}(x), E^a \rangle_\frkg \Big)_{a=1}^{\dim \frkg}
\end{equation}
is a collection of independent, real-valued white noises at temperature $\coup>0$, where $(E_a)_{a=1}^{\dim \frkg}$ is our orthonormal basis of $\frkg$. Similarly, a random distribution $W^{(\Rscript,\coup)}\colon \bT_R \rightarrow \frkg$ is called a $\frkg$-valued, $2\pi R$-periodic white noise at temperature $\coup>0$ if 
\begin{equation}\label{prelim:eq-white-noise-g-valued-coordinates-periodic}
\Big( \langle W^{(\Rscript,\coup)}(x), E^a \rangle_\frkg \Big)_{a=1}^{\dim \frkg}
\end{equation}
is a collection of independent, real-valued, $2\pi R$-periodic white noises at temperature $\coup>0$. 
\end{definition}

While Definition \ref{prelim:def-white-noise-real-valued} and  Definition \ref{prelim:def-white-noise-g-valued} are natural from a conceptual perspective, it is difficult to directly work with them. Instead of the definition, we will mostly work with Fourier-based representations of white noise. In order to introduce the Fourier-based representations, we first make the following definition. 

\begin{definition}[Standard Gaussian sequences]\label{prelim:def-standard-Gaussian-sequence}
A random sequence $(g_n)_{n\in \Z}\subseteq \mathbb{C}$ is called a standard, $\mathbb{C}$-valued Gaussian sequence if it satisfies the following three properties: 
\begin{enumerate}[label=(\roman*)]
\item (Independence) For all $m,n\in \Z$ satisfying $m\neq \pm n$, $g_m$ and $g_n$ are independent. 
\item (Distribution) $g_0$ is a standard real-valued Gaussian and, for all $m\in \Z\backslash \{0\}$, $g_m$ is a standard complex-valued Gaussian. 
\item (Real-valuedness) For all $m\in \Z$, it holds that $\overline{g_m}=g_{-m}$. 
\end{enumerate}
A random sequence $(G_n)_{n\in \Z}\subseteq \frkg_\mathbb{C}$ is called a standard, $\frkg_{\mathbb{C}}$-valued Gaussian sequence if
\begin{equation*}
\Big( \langle G_m , E^a \rangle \Big)_{\substack{a=1,\hdots,\dim \frkg\\ m\in \Z }}
\end{equation*}
is a collection of independent, standard, $\mathbb{C}$-valued Gaussian sequences. 
\end{definition}

Using Definition \ref{prelim:def-standard-Gaussian-sequence}, we now first obtain a representation of $2\pi$-periodic white noise.

\begin{lemma}[Representation of $2\pi$-periodic white noise]\label{prelim:lem-white-noise-representation-2pi}
Let $\coup>0$ and let $(G_m)_{m\in \Z}$ be a standard, $\frkg_{\mathbb{C}}$-valued Gaussian sequence. Then,
\begin{equation*}
W^{(\coup)}=\hcoup \sum_{m\in \Z} G_m e^{imx}
\end{equation*}
is a $\frkg$-valued, $2\pi$-periodic white noise at temperature $\coup$. 
\end{lemma}

\begin{proof}
This follows directly from Definition \ref{prelim:def-white-noise-g-valued} and Plancherell's formula. 
\end{proof}

While Lemma \ref{prelim:lem-white-noise-representation-2pi} contains a convenient representation of $2\pi$-periodic white noise, the direct analog of Lemma \ref{prelim:lem-white-noise-representation-2pi} for $2\pi R$-periodic white noise is ill-suited for our purposes. The reason is that integration of Fourier modes such as $e^{\frac{i}{R} x}$ leads to a loss of $\sim R$, which is problematic. Instead, we rely on a different representation, which is built out of $2\pi$-periodic white noises and the lattice partition from Definition \ref{prelim:def-lattice-partition}.

\begin{lemma}[Representation of white noise]\label{prelim:lem-white-noise-representation}
Let $R\geq 1$ and let $(W^{(\coup)}_{x_0})_{x_0\in \Lambda}$ be a family of independent, $2\pi$-periodic, $\frkg$-valued white noises at temperature $\coup>0$. Then, we have the following properties: 
\begin{enumerate}[label=(\roman*)]
    \item\label{prelim:item-representation-1} $\sum_{x_0\in \Lambda} \psix W^{(\coup)}_{x_0}$ is a $\frkg$-valued white noise at temperature $\coup$. 
    \item\label{prelim:item-representation-2} $\sum_{x_0 \in \LambdaRR} \psiRx W^{(\coup)}_{x_0}$ is a $2\pi R$-periodic, $\frkg$-valued white noise at temperature $\coup$. 
    \item\label{prelim:item-representation-3} For all smooth $\varphi\colon \R \rightarrow \bC$ whose support is contained in $[-\tfrac{\pi}{2} R, \tfrac{\pi}{2}R]$, it holds that
    \begin{equation*}
    \varphi(x) \sum_{x_0\in \Lambda} \psix W^{(\coup)}_{x_0} 
    = \varphi(x) \sum_{x_0 \in \LambdaRR} \psiRx W^{(\coup)}_{x_0}. 
    \end{equation*}
\end{enumerate}
\end{lemma}

\begin{remark}
The statements in \ref{prelim:item-representation-1} and \ref{prelim:item-representation-2} not only yield representations of $2\pi R$-periodic and non-periodic, $\frkg$-valued white noises, but in fact yield a coupling between them, which will be convenient in the proof of Proposition \ref{main:prop-refined-gwp}.
\end{remark}

\begin{proof}[Proof of Lemma \ref{prelim:lem-white-noise-representation}:]
Since the properties in \ref{prelim:item-representation-1} and \ref{prelim:item-representation-2} follow directly from \eqref{prelim:eq-psi-properties} and a direct calculation of the covariance function, we omit the details. The property in \ref{prelim:item-representation-3} follows from $\varphi \psix=0$ for all $x_0 \in \Lambda \backslash \LambdaR$ and $\varphi \psix = \varphi \psiRx$ for all $x_0\in \LambdaR$.
\end{proof}

We now bound the difference between our representations of $2\pi R$-periodic and non-periodic white noise.

\begin{lemma}\label{prelim:lem-white-noise-difference}
Let $R\geq 1$, let $\Nd \in \dyadic$, let $\coup>0$, and let $(W^{(\coup)}_{x_0})_{x_0\in \Lambda}$ be as in Lemma \ref{prelim:lem-white-noise-representation}. Also, let 
\begin{equation*}
Z:= \sum_{x_0\in \Lambda} \psix W^{(\coup)}_{x_0}- \sum_{x_0 \in \LambdaRR} \psiRx W^{(\coup)}_{x_0}. 
\end{equation*}
and let $\varphi\colon \R \rightarrow [0,1]$ be any smooth cut-off function satisfying $\varphi|_{[-\frac{5}{4},\frac{5}{4}]}=1$
and \revision{$\varphi|_{\R\backslash [-\frac{11}{8},\frac{11}{8}]}=0$}. Then, it holds for all $p\geq 1$ that 
\begin{equation}\label{prelim:eq-white-noise-difference}
\E \Big[ \Big\| \varphi \big( \tfrac{x}{R}\big) P_{\leq \Nd} Z \Big\|_{L_x^\infty(\R)}^p\Big]^{\frac{1}{p}}
\lesssim \sqrt{p} \, (R\Nd)^{-100}.
\end{equation}
\end{lemma}

\begin{proof}
We first let $\widetilde{\varphi}\colon \R \rightarrow [0,1]$ be any smooth cut-off function satisfying $\widetilde{\varphi}|_{[-\frac{3}{2},\frac{3}{2}]}=1$
and \revision{$\widetilde{\varphi}|_{\R\backslash[-\frac{\pi}{2},\frac{\pi}{2}]}=0$}. 
Due to Lemma \ref{prelim:lem-white-noise-representation}.\ref{prelim:item-representation-3}, it holds that
\begin{equation*}
Z = \revision{\big( 1- \widetilde{\varphi} \big(  \tfrac{x}{R} \big) \big)} Z 
\end{equation*}
and therefore
\begin{equation*}
\varphi \big( \tfrac{x}{R}\big) P_{\leq \Nd} Z 
= \varphi \big( \tfrac{x}{R}\big) P_{\leq \Nd} \bigg( \Big( 1- \widetilde{\varphi} \big(  \tfrac{x}{R} \big) \Big) \Big)  Z \bigg).
\end{equation*}
The desired estimate \eqref{prelim:eq-white-noise-difference} now follows from standard mismatch estimates (see e.g. \cite[Lemma 2.3]{Bri20}) and standard moment estimates (see e.g. \cite[Lemma 2.29]{BLS21}).
\end{proof}

\subsubsection{Gibbs measures}
Equipped with the $\frkg$-valued white noise from Definition \ref{prelim:def-white-noise-g-valued}, we can now define the Gibbs measures $\muR$ and $\mu^{(\coup)}$. We will be more detailed than in the introduction and also define the corresponding sample spaces and $\sigma$-algebras. 

\begin{definition}[\protect{Sample spaces and $\sigma$-algebras}]\label{prelim:def-sample-space}
For any $R\geq 1$, we define the sample space
$\Omegas$ as  $\C_x^{s-1}(\bT_R \rightarrow \frkg^2)$ and equip it with the corresponding 
 Borel $\sigma$-Algebra $\mathcal{B}(\Omegas)$. We further define $\Omegasline$ as the space of all maps which belong to $\C_x^{s-1}(J\rightarrow \frkg^2)$ for all compact intervals $J\subseteq \R$ \revision{and for which the norm}
\begin{equation}\label{prelim:eq-sample-space-norm}
\big\| (A,B) \big\| := \sum_{R\in \dyadic} R^{-10} \Big( \big\| \chi\big( \tfrac{x}{R} \big) A \big\|_{\C_x^{s-1}(\R\rightarrow \frkg)}
+ \big\| \chi\big( \tfrac{x}{R} \big) B \big\|_{\C_x^{s-1}(\R\rightarrow \frkg)} \Big)
\end{equation}
\revision{is finite. It is equipped with the norm from \eqref{prelim:eq-sample-space-norm} and the corresponding Borel $\sigma$-algebra by $\mathcal{B}(\Omegasline)$.}
\end{definition}

\begin{definition}[Gibbs measure]\label{prelim:def-Gibbs}
For any $R \geq 1$ and $\coup>0$, the Gibbs measure $\muR \colon \mathcal{B}(\Omegas) \rightarrow [0,1]$ is defined as 
\begin{equation*}
\muR = \operatorname{Law} \Big( \big(  W_0^{(\Rscript,\coup)}, W_1^{(\Rscript,\coup)} \big) \Big), 
\end{equation*}
where $W_0^{(\Rscript,\coup)}, W_1^{(\Rscript,\coup)}\colon \bT_R \rightarrow \frkg$ are independent, $2\pi R$-periodic, $\frkg$-valued white noises at temperature $8\coup$. Similarly, the Gibbs measure $\mu^{(\coup)}\colon \mathcal{B}(\Omegasline) \rightarrow [0,1]$ is defined as 
\begin{equation*}
\mu^{(\coup)} = \operatorname{Law} \Big( \big(  W_0^{(\coup)}, W_1^{(\coup)} \big) \Big), 
\end{equation*}
where $W_0^{(\coup)}, W_1^{(\coup)}\colon \R \rightarrow \frkg$ are independent,  $\frkg$-valued white noises at temperature $8\coup$. 
\end{definition}

\subsubsection{Gaussian chaos}\label{section:preliminaries-chaos}
We now recall the definition and properties of Gaussian chaos. To this end, we let 
 $(G_{u_0,k}^+)_{k\in \Z}$ and $(G_{v_0,m}^-)_{m\in \Z}$, where $u_0,v_0\in \Lambda$, be independent, standard, $\frkg_\bC$-valued Gaussian sequences. The standard Gaussian sequences can be used to define independent, $2\pi $-periodic, $\frkg$-valued white noises 
 $W^{(\coup),+}_{u_0},W^{(\coup),-}_{v_0}\colon \bT \rightarrow \frkg$ at temperature $\coup$ via 
\begin{equation*}
W^{(\coup),+}_{u_0} = \hcoup \sum_{k\in \Z} G_{u_0,k}^+ e^{ikx} \qquad \text{and} \qquad 
W^{(\coup),-}_{v_0} = \hcoup \sum_{m\in \Z} G_{v_0,m}^- e^{imx}.
\end{equation*}
Using the orthonormal basis from \eqref{prelim:eq-ONB}, we can decompose 
\begin{equation}\label{prelim:eq-Gm-decomposition}
G_{u_0,k}^+ = G_{u_0,k}^{+,a} E_a \qquad \text{and} \qquad G_{v_0,m}^- = G_{v_0,m}^{-,b} E_b,
\end{equation}
where $(G_{u_0,k}^{+,a})_{k\in \Z}$ and $(G_{v_0,m}^{-,b})_{m\in \Z}$ are standard, $\mathbb{C}$-valued Gaussian sequences. For any degree $d\geq 1$, we define the space of Gaussian chaoses $\mathcal{G}\mathcal{C}_{\leq d}$ of degree less than or equal to $d$ as the closure of  
\begin{equation*}
\begin{aligned}
\operatorname{span} \bigg( \Big\{& 
G_{u_1,k_1}^{+,a_1} \hdots G_{u_{d_p},k_{d_p}}^{+,a_{d_p}} G_{v_1,m_1}^{-,b_1} \hdots G_{v_{d_m},m_{d_m}}^{-,b_{d_m}} \Big|  d_p, d_m \in \mathbb{N}_0, \, d_p + d_m \leq d,  u_1,\hdots,u_{d_p} \in \Lambda, \,   \\
& v_1,\hdots,v_{d_m}\in \Lambda, \, 1\leq a_1,\hdots, a_{d_p}, b_1,\hdots, b_{d_m} \leq \dim\frkg, \, k_1,\hdots,k_{d_p}, m_1,\hdots, m_{d_m}\in \Z \Big\} \bigg). 
\end{aligned}
\end{equation*}
For any homogeneous polynomial $q(G^+,G^-)$ of degree $d$, we define the Wick-ordered version $\mbox{$\lcol q(G^+,G^-) \rcol$}$ as the projection of $q(G^+,G^-) $ onto the orthogonal complement of $\mathcal{G}\mathcal{C}_{\leq d-1}$. In this article, we will mostly use that 
\begin{align*}
\lcol G_{u_1,k_1}^{+,a_1} G_{u_2,k_2}^{+,a_2} \rcol &=  
G_{u_1,k_1}^{+,a_1} G_{u_2,k_2}^{+,a_2} - \delta^{a_1 a_2} \delta_{u_1=u_2} \delta_{k_1+k_2=0}, \\ 
\lcol G_{v_1,m_1}^{-,b_1} G_{v_2,m_2}^{-,b_2} \rcol &= 
 G_{v_1,m_1}^{-,b_1} G_{v_2,m_2}^{-,b_2} -\delta^{b_1b_2} \delta_{v_1=v_2} \delta_{m_1+m_2=0}, \\ 
 \lcol G_{u_0,k}^{+,a} G_{v_0,m}^{-,b} \rcol &= G_{u_0,k}^{+,a} G_{v_0,m}^{-,b}. 
\end{align*}

In Section \ref{section:Killing}, we will need the following two properties of (Wick-ordered) polynomials in Gaussians. For their proofs, we refer to the monograph \cite{N06}. 

\begin{lemma}[$L^2_\omega$-orthogonality]\label{prelim:lem-orthogonality}
For any $d_p,d_m \in \mathbb{N}$, the family of Wick-ordered Gaussian chaoses 
\begin{equation*}
\begin{aligned}
\Big\{ \biglcol \, 
&G_{u_1,k_1}^{+,a_1} \hdots G_{u_{d_p},k_{d_p}}^{+,a_{d_p}} G_{v_1,m_1}^{-,b_1} \hdots G_{v_{d_m},m_{d_m}}^{-,b_{d_m}} \bigrcol \Big|  \,  
 u_1,\hdots,u_{d_p},v_1,\hdots, v_{d_m}\in \Lambda, \\ 
&\, 1\leq a_1,\hdots, a_{d_p}, b_1,\hdots, b_{d_m} \leq \dim\frkg, 
k_1,\hdots,k_{d_p}, m_1,\hdots, m_{d_m}\in \Z \Big\} \bigg)
\end{aligned}
\end{equation*}
is, up to permutations, pairwise orthogonal. 
\end{lemma}

\begin{lemma}[Gaussian hypercontractivity]\label{prelim:lem-hypercontractivity}
Let $d\in \mathbb{N}$ and let $q(G^+,G^-)$ be a polynomial of degree less than or equal to $d$. Then, it holds for all $p\geq 1$ that
\begin{equation*}
\E \Big[ \big|q(G^+,G^-) \big|^p \Big]^{\frac{1}{p}} \lesssim p^{\frac{d}{2}} 
\E \Big[ \big|q(G^+,G^-) \big|^2 \Big]^{\frac{1}{2}}.  
\end{equation*}
\end{lemma}

\subsubsection{Maxima of random variables} We now discuss estimates of the maxima of random variables. While we will only need to bound maxima of Gaussian chaoses, it is convenient to work with more general random variables. 
\begin{definition}\label{prelim:def-subgaussian}
Let $0<\gamma<\infty$. For any random variable $X$, we define 
\begin{equation*}
\big\| X\big\|_{\Psi_\gamma} := \sup_{p\geq 1} p^{-\gamma} \E \big[ |X|^p \big]^{1/p}. 
\end{equation*}
\end{definition}
Using Definition \ref{prelim:def-subgaussian}, we can now state the following lemma.

\begin{lemma}[Control of maxima]\label{prelim:lem-maxima}
Let $0<\gamma<\infty$, let $J\geq 1$, and let $(X_j)_{j=1}^J$ be a finite sequence of random variables. Then, it holds that 
\begin{equation*}
\big\| \max_{j=1,\hdots,J} |X_j| \big\|_{\Psi_\gamma} 
\leq e \log(2+J)^\gamma \max_{j=1,\hdots,J} \big\| X_j \big\|_{\Psi_\gamma}. 
\end{equation*}
\end{lemma}
For a proof of this lemma, we refer to \cite{V18} or \cite[Lemma 4.48]{B23}. 
Finally, we present a corollary of Lemma \ref{prelim:lem-maxima} which concerns the suprema of continuous random processes. In the literature on dispersive equations, such estimates are known as meshing arguments. 

\begin{corollary}[Meshing argument]\label{prelim:cor-meshing}
Let $F\colon \bT^2 \rightarrow \bC$ be a random process, let $\epsilon,\gamma>0$, and let $C,M\geq 1$. Then, it holds that 
\begin{equation}\label{prelim:eq-meshing}
\big\| \sup_{u,v\in \bT^2} |F(u,v)| \big\|_{\Psi_\gamma}
\lesssim_{\epsilon,\gamma,C} M^\epsilon \sup_{u,v\in \bT^2} \big\| F(u,v) \big\|_{\Psi_\gamma} + M^{-C} \big\| \Lip(F) \big\|_{\Psi_\gamma},
\end{equation}
where $\Lip(F)$ is the Lipschitz constant of $F$. 
\end{corollary}

\begin{remark}
In all applications of Corollary \ref{prelim:cor-meshing}, the $\epsilon$-loss in $M$ is acceptable and the second term in \eqref{prelim:eq-meshing} is completely negligible. Thus, Corollary \ref{prelim:cor-meshing} allows us to essentially replace moments of suprema by suprema of moments. 
\end{remark}

\begin{proof}[Proof of Corollary \ref{prelim:cor-meshing}:]
Let $h\in (0,1)$ remain to be chosen and let $\bT^{(h)}:= \bT \medcap (h \Z)$ be the periodic lattice with lattice spacing $h$. Using the triangle inequality, it follows that
\begin{equation}\label{prelim:eq-meshing-p1}
\sup_{u,v \in \bT} |F(u,v)| \leq \sup_{u,v\in \bT^{(h)}} |F(u,v)| + 2h \Lip(F). 
\end{equation}
Using \eqref{prelim:eq-meshing-p1} and  Lemma \ref{prelim:lem-maxima}, it follows that
\begin{align*}
\big\| \sup_{u,v \in \bT} |F(u,v)| \big\|_{\Psi_\gamma} 
&\lesssim \big\| \sup_{u,v\in \bT^{(h)}} |F(u,v)| \big\|_{\Psi_\gamma}  
+ h \big\| \Lip(F) \big\|_{\Psi_\gamma}  \\
&\lesssim \log(h^{-1})^\gamma  \sup_{u,v\in \bT^{(h)}} \big\|  F(u,v) \big\|_{\Psi_\gamma}  
+ h \big\| \Lip(F) \big\|_{\Psi_\gamma} \\
&\lesssim \log(h^{-1})^\gamma  \sup_{u,v\in \bT} \big\|  F(u,v) \big\|_{\Psi_\gamma}  
+ h \big\| \Lip(F) \big\|_{\Psi_\gamma}. 
\end{align*}
The desired estimate now follows by choosing $h:=M^{-C}$. 
\end{proof}

\section{Ansatz}\label{section:ansatz}

In this section, which is at the heart of this article, we derive our Ansatz for solutions of the discretized wave maps equation \eqref{intro:eq-system-U-V-Kil}. In Subsection \ref{section:ansatz-discretized-wave-maps},
we re-introduce the setting of this article, which was previously discussed in the introduction and overview of the argument (Subsection \ref{section:introduction-result} and Section \ref{section:overview}). In contrast to our earlier discussion, however, the definitions of all objects (such as the Killing-renormalization) will be made precise. 
In Subsection \ref{section:ansatz-heuristic} and Subsection \ref{section:ansatz-rigorous}, we then introduce our Ansatz for \eqref{intro:eq-system-U-V-Kil}. The first of the two subsections (Subsection \ref{section:ansatz-heuristic}) is heuristic and only used as a motivation, but hopefully makes the rest of this section more accessible. The rigorous definitions of the terms in our Ansatz are then included in the second of the two subsections (Subsection \ref{section:ansatz-rigorous}). 

In Subsection \ref{section:ansatz-modulation}, Subsection \ref{section:ansatz-remainder}, and Subsection \ref{section:ansatz-contraction}, we then make numerous definitions related to the evolution equations for the terms in our Ansatz. In Subsection \ref{section:ansatz-modulation}, we take a closer look at the modulation equations, which are the evolution equations for the pure modulation operators $\pSN[K][+][k]$ and $\pSN[M][-][m]$. In Subsection \ref{section:ansatz-remainder}, we take a closer look at the remainder equations, which are the evolution equations for the nonlinear remainders $\UN[][\fs]$ and $\VN[][\fs]$. In particular, we define all the errors terms which will be encountered in the rest of this article, such as the high$\times$high$\rightarrow$low-errors, Jacobi errors, structural errors, renormalization errors, and perturbative interactions. Finally, in Subsection \ref{section:ansatz-contraction}, we formalize the assumptions in our contraction-mapping arguments. 

\subsection{Setting}\label{section:ansatz-discretized-wave-maps} 
The discretized wave maps equation \eqref{intro:eq-system-U-V-Kil} was briefly introduced and discussed in Section \ref{section:overview}. We now repeat parts of our earlier discussion but also provide additional details and rigorous definitions. We recall from \eqref{intro:eq-A-B} that the $(1+1)$-dimensional wave maps equation can be written as 
\begin{equation}\label{ansatz:eq-wave-map}
\begin{cases}
\begin{aligned}
\partial_t A &= \partial_x B, \\ 
\partial_t B &= \partial_x B - \big[ A , B \big]. 
\end{aligned}
\end{cases}
\end{equation}
As our main goal is to prove Theorem \ref{intro:thm-rigorous-A-B}, we want to study \eqref{ansatz:eq-wave-map} with initial data given by a pair
\begin{equation}\label{ansatz:eq-initial-data}
\Big(  W_0^{(\coup)}, W_1^{(\coup)} \Big)
\end{equation}
of independent, $\frkg$-valued white noises at temperature $8\coup$ (see Definition \ref{prelim:def-white-noise-g-valued}). As previously mentioned in Section \ref{section:overview}, the prefactor $8$ in the temperature is convenient once we work in null coordinates. 

\subsubsection{Killing-renormalization}
In order to state our discretization of \eqref{ansatz:eq-wave-map}, we first need to introduce the covariance functions and  Killing-renormalizations. 

\begin{definition}[Killing-renormalization]\label{ansatz:def-Killing} 
For all $N,\Nd\in \dyadic$, we make the following definitions: 
\begin{enumerate}[label=(\roman*)]
\item\label{ansatz:item-covariance} (Covariance functions) For all $M\in \dyadic$ satisfying $M\geq 2$, we define  $\Cf_M, \Cf^{(\Nscript)}_M,\Cf^{(\Nscript,\Nd)}_M\colon \R \rightarrow \R$ by 
\begin{align*}
\Cf_M(y) 
&:= \int_{\R} \dxi \,  \mathbf{1} \Big\{ \tfrac{M}{2} \leq |\xi|\leq M \Big\} 
\frac{e^{-i\xi y}}{i\xi}  
= - \int_{\R} \dxi \, \mathbf{1} \Big\{ \tfrac{M}{2} \leq |\xi|\leq M \Big\} 
\frac{\sin(\xi y)}{\xi}, \\ 
\Cf_M^{(N)}(y) &:=   
 - \int_{\R} \dxi \,   \rho_{\leq N}^2(\xi) \mathbf{1} \Big\{ \tfrac{M}{2} \leq |\xi|\leq M \Big\} 
\frac{\sin(\xi y)}{\xi}, \\ 
\Cf^{(N,\Nd)}_M(y) &:=   
- \int_{\R} \dxi \,   \rho_{\leq N}^2(\xi) \rho_{\leq \Nd}^2(\xi) \mathbf{1} \Big\{ \tfrac{M}{2} \leq |\xi|\leq M \Big\} 
\frac{\sin(\xi y)}{\xi}.
\end{align*}
For $M=1$, the functions $\Cf_1, \Cf^{(\Nscript)}_1,\Cf^{(\Nscript,\Nd)}_1\colon \R \rightarrow \R$ are defined similarly, 
but with $\mathbf{1}\big\{ \tfrac{M}{2} \leq |\xi| \leq M \big\}$ replaced by $\mathbf{1}\big\{ |\xi| \leq 1 \big\}$. Furthermore, for all $K\in \dyadic$, we define $\Cf_{<K}, \Cf^{(\Nscript)}_{<K},\Cf^{(\Nscript,\Nd)}_{<K}\colon \R \rightarrow \R$ by
\begin{equation*}
\Cf_{<K} = \sum_{M<K} \Cf_M, \quad  \Cf^{(\Nscript)}_{<K} = \sum_{M<K} \Cf^{(\Nscript)}_M, \quad 
\text{and} \quad \Cf^{(\Nscript,\Nd)}_{<K} = \sum_{M<K} \Cf^{(\Nscript,\Nd)}_M. 
\end{equation*}
Finally, we define $\Cf^{(N)},\Cf^{(N,\Nd)}\colon \R \rightarrow \R$ by 
\begin{equation*}
\Cf^{(N)} = \sum_{M} \Cf^{(N)}_M \qquad \text{and} \qquad 
\Cf^{(N,\Nd)}  = \sum_{M} \Cf^{(N,\Nd)}_M. 
\end{equation*}
\item (Killing-renormalization) For all $M\in \dyadic$, we define the operators $\Renorm[N][M]$ and $\Renorm[N,\Nd][M]$ by
\begin{align*}
\Renorm[N][M] F &:=  P_{\leq N}^x \int_{\R} \dy \, 
\big(\widecheck{\rho}_{\leq N} \ast \widecheck{\rho}_{\leq N}\big)(y)  \Cf^{(\Nscript)}_M(y)  \Theta^x_y P_{\leq N}^x \Kil F, \\ 
\Renorm[N,\Nd][M] F &:=   P_{\leq N}^x \int_{\R} \dy \, 
\big(\widecheck{\rho}_{\leq N} \ast \widecheck{\rho}_{\leq N}\big)(y) \Cf^{(\Nscript,\Ndscript)}_M(y)  \Theta^x_y P_{\leq N}^x \Kil F
\end{align*}
for all $F\colon \R \rightarrow \frkg$. \revision{Here, $\Kil\colon \frkg\rightarrow \frkg$ is the Killing map from Definition \ref{prelim:def-Killing}.} The operators $\Renorm[N][]$ and $\Renorm[N,\Nd][]$ are defined similarly, but with the corresponding covariance functions from \ref{ansatz:item-covariance}.
\end{enumerate}
\end{definition}

In the discretized wave maps equation (see Definition \ref{ansatz:def-discretized-wave-maps} and Definition \ref{ansatz:def-discretized-wave-maps-periodic-truncated}), we only make use of $\Renorm[N][]$ and $\Renorm[N,\Nd][]$, but not $\Renorm[N][M]$ or $\Renorm[N,\Nd][M]$. However, the functions and operators with additional frequency-localization, such as $\Renorm[N][M]$ or $\Renorm[N,\Nd][M]$, will play an important role in decompositions and estimates of the discretized wave maps equation. 
We now record elementary properties of the covariance function and Killing-renormalization. 

\begin{lemma}[Properties of the covariance function and Killing-renormalization]\label{ansatz:lem-renormalization}
Let $M,N,\Nd\in \dyadic$. Then, it holds that 
\begin{equation}\label{ansatz:eq-renormalization-1}
\big| \Cf^{(N)}_M(y) \big| + \big| \Cf^{(N,\Nd)}_M(y) \big|\lesssim \min \big( 1, M |y| \big).  
\end{equation}
Furthermore, it holds for all $\alpha_1,\alpha_2\in \R\backslash \{0\}$, $K,L \in \dyadic$, and $F\colon \R^{1+1} \rightarrow \frkg$ that
\begin{align}
\big\| \Renorm[N][M] F \big\|_{\Cprod{\alpha_1}{\alpha_2}}
+ \big\| \Renorm[N,\Nd][M] F \big\|_{\Cprod{\alpha_1}{\alpha_2}}
&\lesssim \frac{M}{N} \big\| F \big\|_{\Cprod{\alpha_1}{\alpha_2}}, 
\label{ansatz:eq-renormalization-2} \\
\big\| \Renorm[N][M] P_K^u P_L^v F \big\|_{\Cprod{\alpha_1}{\alpha_2}}
+ \big\| \Renorm[N,\Nd][M] P_K^u P_L^v F \big\|_{\Cprod{\alpha_1}{\alpha_2}}
&\lesssim \frac{M}{N} \frac{\max(K,L)}{N} \big\|  F \big\|_{\Cprod{\alpha_1}{\alpha_2}}.
\label{ansatz:eq-renormalization-3}
\end{align}
Furthermore, it holds that 
\begin{align}
\bigg| \Cf_M(y) - \Big( -  \sum_{m\in \Z_M} \frac{\sin(my)}{m} \Big) \bigg| 
&\lesssim |y|, \label{ansatz:eq-renormalization-4} \\ 
\bigg| \Cf_M^{(\Nscript)}(y) - \Big( - \sum_{m\in \Z_M} \rho_{\leq N}^2(m) 
\frac{\sin(my)}{m} \Big) \bigg| &\lesssim |y|, \label{ansatz:eq-renormalization-5}\\ 
\bigg| \Cf_M^{(\Nscript,\Nd)}(y) - \Big( - \sum_{m\in \Z_M} \rho_{\leq N}^2(m) \rho_{\leq \Nd}^2(m)
\frac{\sin(my)}{m} \Big) \bigg| &\lesssim |y|. \label{ansatz:eq-renormalization-6}
\end{align}
\end{lemma}

\begin{remark}
Due to \eqref{ansatz:eq-renormalization-3}, the renormalization only accounts for boundary effects in frequency space. For functions $F\colon \R^{1+1}_{u,v}\rightarrow \frkg$ which are supported on $u$ and $v$-frequencies $\lesssim N^{1-\theta}$ for some $\theta\in (0,1)$, i.e., on frequencies far away from the frequency boundary, applying $\Renorm[N]$ leads to a gain of $\sim N^{-\theta}$. 
\end{remark}

\begin{remark}
While the definitions of $\Cf_M(y)$, $\Cf_M^{(\Nscript)}$, and $\Cf_M^{(\Nscript,\Nd)}$ are based on integrals, we will mostly work with the sums in \eqref{ansatz:eq-renormalization-4}, \eqref{ansatz:eq-renormalization-5}, and \eqref{ansatz:eq-renormalization-6}. The reason for still using integrals in Definition \ref{ansatz:def-Killing} is to obtain the exact scaling symmetry in Lemma \ref{ansatz:lem-scaling-symmetry}, which would otherwise only be satisfied approximately.
\end{remark}

\begin{proof}
We prove \eqref{ansatz:eq-renormalization-1}-\eqref{ansatz:eq-renormalization-3} only for $\Cf^{(N)}_M$ and $\Renorm[N][M]$, since 
the arguments for $\Cf^{(N,\Nd)}_M$ and $\Renorm[N,\Nd][M]$ are similar.
Using the definition of $\Cf^{(N)}_M$ and the trivial bound $|\sin(\zeta)|\leq \min(1,|\zeta|)$, it holds that 
\begin{equation*}
\Big| \Cf^{(N)}_M(y) \Big|
\leq \int_{\R} \dxi \, \mathbf{1}\big\{ |\xi| \sim M \big\}
\frac{|\sin(\xi y)|}{|\xi|}
\leq \int_{\R} \dxi \, \mathbf{1}\big\{ |\xi| \sim M \big\}
\min\Big( \frac{1}{|\xi|}, |y| \Big) \lesssim  \min \big( 1, M |y|\big). 
\end{equation*}
This completes the proof of \eqref{ansatz:eq-renormalization-1}.  Using \eqref{ansatz:eq-renormalization-1}, it follows that 
\begin{equation}\label{ansatz:eq-renormalization-p1}
\Big\| \big( \widecheck{\rho}_{\leq N}
\ast \widecheck{\rho}_{\leq N}\big)(y) \Cf^{(N)}_M(y) 
\Big\|_{L^1_y} \lesssim \Big\| \langle N |y| \rangle^{-10} M |y| \Big\|_{L_y^1} \lesssim MN^{-1}. 
\end{equation}
Due to the triangle inequality, we then obtain \eqref{ansatz:eq-renormalization-2} from \eqref{ansatz:eq-renormalization-p1}.  Since $\Cf^{(N)}_M(y)$ is odd and $\widecheck{\rho}_{\leq N}(y)$ is even, it holds that
\begin{equation*}
\int_{\R} \dy \, \big( \widecheck{\rho}_{\leq N}
\ast \widecheck{\rho}_{\leq N}\big)(y) \Cf^{(N)}_M(y)
=0.
\end{equation*}
As a result, we obtain that 
\begin{equation*}
\begin{aligned}
    \Renorm[N][M] P_K^u P_L^v F 
    &= P_{\leq N}^x \int_{\R} \dy \,  \big(\widecheck{\rho}_{\leq N} \ast \widecheck{\rho}_{\leq N}\big)(y)  \Cf^{(\Nscript)}_M(y) 
    \Theta^x_y P_{\leq N}^x P_K^u P_L^v F  \\ 
    &= P_{\leq N}^x \int_{\R} \dy \,  \big(\widecheck{\rho}_{\leq N} \ast \widecheck{\rho}_{\leq N}\big)(y)  \Cf^{(\Nscript)}_M(y) 
    \big( \Theta^x_y -1 \big) P_{\leq N}^x P_K^u P_L^v F. 
    \end{aligned}
\end{equation*}
Using the triangle inequality, it follows that 
\begin{align*}
\Big\| \Renorm[N][M] P_K^u P_L^v F \Big\|_{\Cprod{\alpha_1}{\alpha_2}}
&\lesssim \int_{\R} \dy \, \big| \big(\widecheck{\rho}_{\leq N} \ast \widecheck{\rho}_{\leq N}\big)(y)  \Cf^{(\Nscript)}_M(y) \big| 
\Big\| (\Theta^x_y-1)  P_K^u P_L^v F\Big\|_{\Cprod{\alpha_1}{\alpha_2}} \\
&\lesssim \bigg(  \int_{\R} \dy \,  \big| \big(\widecheck{\rho}_{\leq N} \ast \widecheck{\rho}_{\leq N}\big)(y)  \Cf^{(\Nscript)}_M(y)\big|\,\big|y\big| \bigg) \max(K,L) \big\| F \big\|_{\Cprod{\alpha_1}{\alpha_2}}.
\end{align*}
Since
\begin{equation*}
\int_{\R} \dy \,  \big| \big(\widecheck{\rho}_{\leq N} \ast \widecheck{\rho}_{\leq N}\big)(y)  \Cf^{(\Nscript)}_M(y)\big|\, \big|y\big| \lesssim \int_{\R} \dy \, \langle N y \rangle^{-10} M |y|^2 \lesssim MN^{-2}, 
\end{equation*}
this yields \eqref{ansatz:eq-renormalization-3}. Finally, it remains to prove the estimates \eqref{ansatz:eq-renormalization-4}, \eqref{ansatz:eq-renormalization-5}, and \eqref{ansatz:eq-renormalization-6}. Since the three estimates can easily be derived by treating the sums as Riemann sums for the integrals in $\Cf_M$, $\Cf_M^{(\Nscript)}$ and $\Cf_M^{(\Nscript,\Nd)}$, we omit the details. 
\end{proof}

\subsubsection{Discretized wave maps equation}
Equipped with Definition \ref{ansatz:def-Killing}, we can now introduce the discretized wave maps equation. 

\begin{definition}[Discretized wave maps equation]\label{ansatz:def-discretized-wave-maps} 
Let $N\in \dyadic$, let $\coup>0$, and let $A^{(\Nscript,\coup)},B^{(\Nscript,\coup)}\colon \R^{1+1} \rightarrow \frkg$. We call $(A^{(\Nscript,\coup)},B^{(\Nscript,\coup)})$ a solution of the discretized wave maps equation with
initial data  $(W_0^{(\coup)},W_1^{(\coup)})$ if 
\begin{equation}\label{ansatz:eq-discretized-wave-maps}
\begin{cases}
\begin{aligned}
\partial_t A^{(\Nscript,\coup)} &= \partial_x B^{(\Nscript,\coup)}, \\ 
\partial_t B^{(\Nscript,\coup)} &= \partial_x A^{(\Nscript,\coup)} -  \Big[  A^{(\Nscript,\coup)},  B^{(\Nscript,\coup)} \Big]_{\leq N} + 2\coup \Renorm[N] A^{(\Nscript,\coup)}, \\
A^{(\Nscript,\coup)}(0) &= W_0^{(\coup)}, \qquad B^{(\Nscript,\coup)}(0) =  W_1^{(\coup)}. 
\end{aligned}
\end{cases}
\end{equation}
Here, $[ \cdot, \cdot]_{\leq N}$ is the frequency-truncated Lie bracket from \eqref{prelim:eq-Lie-bracket-truncated}. 
\end{definition}

From a technical perspective, it is difficult to directly work with \eqref{ansatz:eq-discretized-wave-maps}. 
The reason is that we later want to choose $(W_0^{(\coup)},W_1^{(\coup)})\colon \R \rightarrow \frkg^2$ as $\frkg$-valued white noise at temperature $8\coup$, which not only has low regularity but is also unbounded in space (even in low-regularity spaces). 
In order to make rigorous statements, we need to introduce several cut-offs, which will later be removed by taking limits. First, we let $R\geq 1$ and let $(W_0^{(\Rscript,\coup)},W_1^{(\Rscript,\coup)})\colon \bT_R \rightarrow \frkg^2$, which will later be taken as a pair of $2\pi R$-periodic, $\frkg$-valued white noises. Second, we let $\Nd\in \dyadic$ be a dyadic scale satisfying $\Nd\lesssim N$ and let $P_{\leq \Nd}^x$ be the corresponding Littlewood-Paley operator from \eqref{prelim:eq-P-N}. Instead of $(W_0^{(\coup)},W_1^{(\coup)})$, we then consider
\begin{equation}\label{ansatz:eq-WRNd}
\Big( P_{\leq \Nd}^x W_0^{(\Rscript,\coup)}, P_{\leq \Nd}^x W_1^{(\Rscript,\coup)}\Big),
\end{equation}
which has been periodized and frequency-truncated. Furthermore, let $C_b^\infty(\R)$ be the space of all smooth functions from $\R$ to $\R$ which are bounded and have bounded derivatives. We now let $\chi \in C_b^\infty(\R)$ and refer to $\chi$ as a cut-off function in time. At least for now, we impose no further conditions on $\chi$, but see \mbox{Definition \ref{prelim:def-cut-off}} and Hypothesis \ref{hypothesis:pre}.

\begin{notation}[Dependence on parameters]\label{ansatz:notation}
To simplify the notation, we now indicate the dependence of objects on $(N, \Nd, R,\coup) \in \dyadic \times \dyadic \times [1,\infty) \times (0,\infty)$ and $\chi \in C^\infty_b(\R)$ by 
\begin{equation}\label{ansatz:eq-Ncs}
(\mathscr{N}) =  (N, \Nd, R, \coup ,\chi ). 
\end{equation}
Furthermore, we also write 
\begin{alignat*}{3}
\Cf^{(\Ncs)} &= \Cf^{(\Nscript,\Nd)}, \qquad & \qquad \Cf^{(\Ncs)}_M &= \Cf^{(\Nscript,\Nd)}_M, \\ 
\Renorm[\Ncs][] &= \Renorm[\Nscript,\Nd][], \qquad & \qquad \Renorm[\Ncs][M] &= \Renorm[\Nscript,\Nd][M].
\end{alignat*}
\end{notation}

Equipped with our new notation, we now introduce the discretized wave maps equation with periodized and frequency-truncated initial data. 

\begin{definition}\label{ansatz:def-discretized-wave-maps-periodic-truncated}
Let $N,\Nd\in \dyadic$, let $R\geq 1$, and let $\coup>0$. 
We call $(\A[\Ncs],\B[\Ncs])\colon \R \times \bT_R \rightarrow \frkg$ a solution of the discretized wave maps equation with cut-off function $\chi$ and initial data $(P_{\leq \Nd}^x W_0^{(\Rscript,\coup)},P_{\leq \Nd}^x W_1^{(\Rscript,\coup)})$ if 
\begin{equation}\label{ansatz:eq-discretized-wave-maps-truncated}
\begin{cases}
\begin{aligned}
\partial_t \A[\Ncs] &= \partial_x \B[\Ncs], \\ 
\partial_t \B[\Ncs] &= \partial_x \A[\Ncs] - \chi(t) \big[  \A[\Ncs],  \B[\Ncs] \big]_{\leq N} +2 \chi(t)^2 \coup  \Renorm[\Nc]\A[\Ncs], \\
\A[\Ncs](0) &=   P_{\leq \Nd}^x W_0^{(\Rscript,\coup)}, \qquad \B[\Ncs](0) =   P_{\leq \Nd}^x W_1^{(\Rscript,\coup)}. 
\end{aligned}
\end{cases}
\end{equation}
\end{definition}

\begin{remark}
The reason for introducing the additional frequency-truncation in the initial data in \revision{\eqref{ansatz:eq-discretized-wave-maps-truncated}}, even though the nonlinearity in \revision{\eqref{ansatz:eq-discretized-wave-maps-truncated}} has already been frequency-truncated, is that this makes it easier for us to compare the solutions of \revision{\eqref{ansatz:eq-discretized-wave-maps-truncated}} with solutions of \eqref{intro:eq-system-A-B-N}, which is important for the proof of Theorem~\ref{intro:thm-rigorous-A-B}. 
\end{remark}

\subsubsection{Scaling symmetry}
In the next lemma, we study the scaling symmetry of the discretized wave maps equation. To this end, we first introduce the scaling operator $\Scaling_\kappa$, where $\kappa \in (0,\infty)$. For any $f\colon \R\rightarrow \frkg$ and $g\colon \R^{1+1}\rightarrow \frkg$, it is given by  
\begin{equation}\label{ansatz:eq-scaling-transform}
\big( \Scaling_\kappa f \big) (x) = \kappa f(\kappa x) 
\qquad \text{and} \qquad \big( \Scaling_\kappa g\big)(t,x) = \kappa g(\kappa t,\kappa x).
\end{equation}
\begin{lemma}[Scaling symmetry]\label{ansatz:lem-scaling-symmetry}
Let $(N,\Nd,R,\coup)\in \dyadic \times \dyadic \times [1,\infty) \times (0,\infty)$, let $\chi\in C^\infty_b(\R)$, and let $(\A[\Ncs],\B[\Ncs])$ be the solution of \eqref{ansatz:eq-discretized-wave-maps-truncated} with initial data $(W_0^{(\Rscript,\coup)},W_1^{(\Rscript,\coup)})$. Let 
\begin{equation*}
(N^\prime,\Nd^\prime,R^\prime,\coup^\prime)\in \dyadic \times \dyadic \times [1,\infty) \times (0,\infty),
\end{equation*} 
let $\chi^\prime \in C^\infty_b(\R)$, and let $\kappa \in 2^{\Z}$ be such that 
\begin{equation*}
N^\prime = \kappa N, \qquad \Nd^\prime =\kappa \Nd, \qquad 
R^\prime = \frac{R}{\kappa}, \qquad \coup^\prime = \kappa \coup, \qquad \text{and} \qquad
\chi^\prime(t)=\chi(\kappa t). 
\end{equation*}
Then, $\big( \Scaling_\kappa \A[\Ncs], \Scaling_\kappa \B[\Ncs] \big)$
solves \eqref{ansatz:eq-discretized-wave-maps-truncated} but with the parameters $(N,\Nd,R,\coup,\chi)$ replaced by $(N^\prime,\Nd^\prime,R^\prime,\coup^\prime,\chi^\prime)$ and the initial data $(W_0^{(\Rscript,\coup)},W_1^{(\Rscript,\coup)})$ replaced by
\begin{equation}\label{ansatz:eq-scaled-white-noise}
\Big( \Scaling_\kappa W_0^{(\Rscript,\coup)}, \Scaling_\kappa W_1^{(\Rscript,\coup)}\Big). 
\end{equation}
Furthermore, if $(W_0^{(\Rscript,\coup)},W_1^{(\Rscript,\coup)})$ is $2\pi R$-periodic $\frkg$-valued white noise at temperature $8\coup$, then \eqref{ansatz:eq-scaled-white-noise} is $2\pi R^\prime$-periodic, $\frkg$-valued white noise at temperature $8\coup^\prime$. 
\end{lemma}

\begin{remark}
In Section \ref{section:main}, we will use Lemma \ref{ansatz:lem-scaling-symmetry} to reduce local well-posedness for a general temperature $\coup>0$ to the low-temperature case, i.e., $0<\coup \ll 1$. In our setting, this is a much more convenient 
approach to create smallness than other alternatives, such as inserting a cut-off function of the form $\chi(\coup^{-1}t)$ into \eqref{ansatz:eq-discretized-wave-maps-truncated}. Indeed, due to the low regularity of white noise, it would be technically difficult to include $\chi(\coup^{-1}t)$-factors in our Ansatz below. 
\end{remark}

\begin{proof} Similar as in \eqref{ansatz:eq-Ncs}, we write
\begin{equation*}
(\mathscr{N}^\prime) =  (N^\prime, \Nd^\prime, R^\prime, \coup^\prime ,\chi^\prime). 
\end{equation*}
From a direct calculation, we obtain for all $f,g\colon \R^{1+1} \rightarrow \R$ that 
\begin{equation}\label{ansatz:eq-scaling-p1}
\Scaling_\kappa P_{\leq N}^x f = P^x_{\leq N^\prime} \Scaling_\kappa f 
\qquad \text{and} \qquad \kappa \Scaling_\kappa (fg) = \Scaling_\kappa (f) \, \Scaling_\kappa(g). 
\end{equation}
Furthermore, it follows from Definition \ref{ansatz:def-Killing} that, for all $z\in \R$, 
\begin{equation}\label{ansatz:eq-scaling-p2}
\Cf^{(\Ncs^\prime)}\big( \tfrac{z}{\kappa} \big) = \Cf^{(\Ncs)}\big( z \big). 
\end{equation}
Using \eqref{ansatz:eq-scaling-p2}, it then follows for all $F\colon \R^{1+1}\rightarrow \frkg$ that
\begin{equation}\label{ansatz:eq-scaling-p3}
\Scaling_\kappa \Renorm[\Ncs][] F =  \Renorm[\Ncs^\prime][] \Scaling_\kappa F. 
\end{equation}
Equipped with \eqref{ansatz:eq-scaling-p1} and \eqref{ansatz:eq-scaling-p3}, it is then easy to see that $\big( \Scaling_\kappa \A[\Ncs], \Scaling_\kappa \B[\Ncs] \big)$ satisfies the discretized wave maps equation \eqref{ansatz:eq-discretized-wave-maps-truncated}, and we omit the details. Furthermore, it follows from \eqref{ansatz:eq-scaling-p1} that the initial data of $\big( \Scaling_\kappa \A[\Ncs], \Scaling_\kappa \B[\Ncs] \big)$ is given by 
\begin{equation*}
\big( \Scaling_\kappa P_{\leq \Nd}^x W_0^{(\Rscript,\coup)}, \Scaling_\kappa  P_{\leq \Nd}^x  W_1^{(\Rscript,\coup)}\big)
= \big( P_{\leq \Nd^\prime}^x \Scaling_\kappa  W_0^{(\Rscript,\coup)}, P_{\leq \Nd^\prime}^x 
 \Scaling_\kappa   W_1^{(\Rscript,\coup)}\big),
\end{equation*}
which is as desired. Finally, the claim regarding the distribution of \eqref{ansatz:eq-scaled-white-noise} follows directly from our definition of white noise, i.e., Definition \ref{prelim:def-white-noise-g-valued}. 
\end{proof}

\subsubsection{Null-coordinates}
As previously discussed in the introduction (Subsection \ref{section:overview}), most of our analysis of the discretized wave maps equation is performed in null-coordinates and null-variables. For any Cartesian coordinates $(t,x)\in \R^{1+1}$, the corresponding null-coordinates $(u,v)\in \R^{1+1}$ are defined as 
\begin{equation}\label{ansatz:eq-null-coordinates}
u:= x-t \qquad \text{and} \qquad v:= x+t. 
\end{equation}
For $\A[\Ncs][]$ and  $\B[\Ncs][]$ as in Definition \ref{ansatz:def-discretized-wave-maps-periodic-truncated}, we define the corresponding null-variables as 
\begin{equation}\label{ansatz:eq-null-unknowns}
    \U[\Ncs][] := \frac{\A[\Ncs][] - \B[\Ncs][]}{4} \qquad \text{and} \qquad \V[\Ncs][] := \frac{\A[\Ncs][]+\B[\Ncs][]}{4}. 
\end{equation}
Furthermore, we define the left and right-moving waves $W^{(\Rscript,\coup),\pm}\colon \bT_R \rightarrow \frkg$ as 
\begin{equation}\label{ansatz:eq-null-data}
W^{(\Rscript,\coup),\pm} :=  \frac{W_0^{(\Rscript,\coup)} \mp W_1^{(\Rscript,\coup)}}{4}. 
\end{equation}
In terms of the null coordinates and null unknowns, \eqref{ansatz:eq-discretized-wave-maps-truncated} can then be written as  
\begin{equation}\label{ansatz:eq-wave-maps} 
\begin{cases}
\partial_v \U[\Ncs][]= \partial_u \V[\Ncs][] = 
 \chinull \Big[  \U[\Ncs][],  \V[\Ncs][]\Big]_{\leq N} - \chinullsquare \coup \Renorm[\Ncs] \big( \U[\Ncs][] + \V[\Ncs][] \big), \\
\U[\Ncs][]\big|_{u=v}= P_{\leq \Nd}^x W^{(\Rscript,\coup),+}, \, \V[\Ncs][]\big|_{u=v}=   P_{\leq \Nd}^x W^{(\Rscript,\coup),-}. 
\end{cases}
\end{equation}

\begin{remark}[\protect{The cut-off function $\chi$}]\label{ansatz:rem-cut-off-null}
In \eqref{ansatz:eq-discretized-wave-maps-truncated}, the cut-off function $\chi\colon \R^{1+1} \rightarrow \R, (t,x) \mapsto \chi(t)$ depends only on the time-variable $t$. In \eqref{ansatz:eq-wave-maps}, we make a slight abuse of notation, and also write $\chi$ for the corresponding function  $\R^{1+1} \rightarrow \R, (u,v) \mapsto \chi(\frac{v-u}{2})$ in null-coordinates. 
\end{remark}

\begin{remark}
We note that if $W_0^{(\Rscript,\coup)},W_1^{(\Rscript,\coup)}\colon \bT_R\rightarrow \frkg$ are independent white noises at temperature $8\coup$, then $W^{(\Rscript,\coup),+},W^{(\Rscript,\coup),-}\colon \bT_R \rightarrow \frkg$ are independent white noises at temperature $\coup$. Furthermore, we note that the prefactors  of the nonlinearity and the renormalization term in \eqref{ansatz:eq-wave-maps} are equal to one. This is the reason for denoting the temperature of $W_0^{(\Rscript,\coup)}$ and $W_1^{(\Rscript,\coup)}$ by $8\coup$ and our choices of the prefactors in  \eqref{ansatz:eq-discretized-wave-maps-truncated} and  \eqref{ansatz:eq-null-unknowns}. 
\end{remark}

\subsubsection{Initial data}\label{section:ansatz-initial-data}
In order to prove our main theorem (Theorem \ref{intro:thm-rigorous-A-B}), we want to choose the initial data $(W_0^{(\Rscript,\coup)},W_1^{(\Rscript,\coup)})$ in Definition \ref{ansatz:def-discretized-wave-maps-periodic-truncated} as a pair of independent, $2\pi R$-periodic, $\frkg$-valued white noises at temperature $8\coup$. As a result, the right and left-moving waves $W^{(\Rscript,\coup),+},W^{(\Rscript,\coup),-}\colon \bT_R \rightarrow \frkg$ from \eqref{ansatz:eq-null-data} should be independent, $2\pi R$-periodic, $\frkg$-valued white noises at temperature $\coup$.
Due to Lemma \ref{prelim:lem-white-noise-representation-2pi} and Lemma \ref{prelim:lem-white-noise-representation}, we can therefore represent their frequency-truncations as 
\begin{equation}\label{ansatz:eq-Wpm}
\begin{aligned}
 P_{\leq \Nd}^x W^{(\Rscript,\coup),+}(x) &= \hcoup  P_{\leq \Nd}^x
\sum_{u_0 \in \LambdaRR} \sum_{K\in \dyadic} \sum_{k\in \Z_K} \psiRu(x) G_{u_0,k}^+ e^{ikx}, \\ 
  P_{\leq \Nd}^x W^{(\Rscript,\coup),-}(x) &= \hcoup  P_{\leq \Nd}^x
\sum_{v_0 \in \LambdaRR}  \sum_{M\in \dyadic} \sum_{m\in \Z_M} \psiRv(x) G_{v_0,m}^- e^{imx}, 
\end{aligned}
\end{equation}
where $(G^+_{u_0,k})_{k\in \Z}$ and $(G^-_{v_0,m})_{m\in \Z}$ are independent $\frkg$-valued Gaussian sequences (as in Definition \ref{prelim:def-standard-Gaussian-sequence}). However, in order to iterate our local theory (see Section \ref{section:main}), it will be necessary to introduce slightly more general initial data than in \eqref{ansatz:eq-Wpm}. 

\begin{definition}[Initial data]\label{ansatz:def-initial-data}
Let $N,\Nd\in \Dyadiclarge$, let  $R\geq 1$, let  $\coup>0$, and let 
\begin{equation*}
\big( \SNin[K][+] \big)_{K\in \dyadic}, \big( \SNin[M][-] \big)_{M\in \dyadic}\colon \bT_R \rightarrow \End(\frkg), 
\qquad \text{and} \qquad 
Z^{(\Ncs),+},Z^{(\Ncs),-}\colon \bT_R \rightarrow \frkg.
\end{equation*}
For all $K,M\in \dyadic$, we assume that the initial modulation operators satisfy the following two conditions:
\begin{alignat*}{3}
P^x_{\gg K^{1-\delta}} \SNin[K][+] &= P^x_{\gg M^{1-\delta}} \SNin[M][-] =0 &\qquad\qquad&\textup{for all } K,M\in \dyadic,  \\
\SNin[K][+] &= \SNin[M][-] = \Id_\frkg  &\qquad&\textup{for all } K,M\in \dyadic \textup{ satisfying } K,M >N^{1-\delta}.
\end{alignat*}
are satisfied. Then, we introduce 
\begin{equation}\label{ansatz:eq-Wpm-tilde}
\begin{aligned}
\widetilde{W}^{(\Ncs),+ }(x) 
&=   \sum_{\substack{K\in \dyadic}}   \SNin[K][+]  P^x_{\leq \Nd} P^{\sharp}_{R;K} W^{(\Rscript,\coup),+} 
+ Z^{(\Ncs),+}, \\
\widetilde{W}^{(\Ncs),-}(x) 
&=   \sum_{\substack{M\in \dyadic}}   \SNin[M][-] P^x_{\leq \Nd} P^{\sharp}_{R;M} W^{(\Rscript,\coup),-} 
+ Z^{(\Ncs),-},
\end{aligned}
\end{equation}
where $P^{\sharp}_{R;L}$ is as in \eqref{prelim:eq-Psharp}. 
\end{definition}

\begin{remark}\label{ansatz:rem-initial-data} We make the following remarks regarding Definition \ref{ansatz:def-initial-data}.
\begin{enumerate}[label=(\roman*)]
\item The expressions in \eqref{ansatz:eq-Wpm-tilde} involve $W^{(\Rscript,\coup),+}$ and $W^{(\Rscript,\coup),-}$, but make no explicit use of the representations of white noise from \eqref{ansatz:eq-Wpm}. This will be important in our globalization argument, see e.g. the proof of Lemma \ref{main:lem-gwp-small}.
\item In the case $\SNin[K][+]= \SNin[M][-]=\Id_\frkg$ for all $K,M\in \dyadic$ and $Z^{(\Ncs),+}=Z^{(\Ncs),-}=0$, the initial data in \eqref{ansatz:eq-Wpm-tilde} agrees with the initial data in \eqref{ansatz:eq-Wpm}.
\item In our local well-posedness results, we will later assume that $\SNin[K][+]$ and $\SNin[M][-]$ are close to the identity $\Id_\frkg$ and that $Z^{(\Ncs),+}$ and $Z^{(\Ncs),-}$ are small, smooth remainders (see Definition \ref{main:def-perturbations-null}).
\item In Proposition \ref{main:prop-null-lwp}, we only compare the solution $(U_{\leq \Nd},V_{\leq \Nd})$ of \eqref{intro:eq-system-U-V-N-data} and the solution $(U^{(\Nscript)},V^{(\Nscript)})$ of \eqref{intro:eq-system-U-V-Kil} in the regime $\Nd \leq N^{1-\delta}$. This is the reason for the different treatment of the frequencies $K,M\leq N^{1-\delta}$ and $K,M>N^{1-\delta}$ in Definition \ref{ansatz:def-initial-data}.
\item Due to Definition \ref{ansatz:def-initial-data}, the linear transformations $\SNin[K][+](u,v)\colon \frkg \rightarrow \frkg$ and $\SNin[M][-](u,v)\colon \frkg\rightarrow \frkg$ are orthogonal for $K,M>N^{1-\delta}$, but not necessarily for $K,M\leq N^{1-\delta}$. As a result, our pure modulation operators $(\pSN[K][+][k])_{K\in \Dyadiclarge,k\in \Z_K}$ and $(\pSN[M][-][m])_{M\in \Dyadiclarge,m\in\Z_M}$ from Definition \ref{ansatz:def-pure} below will also only be orthogonal for $K,M>N^{1-\delta}$. Since the orthogonality is only used to address probabilistic resonances stemming from $P_{\leq N}^x$, this will be sufficient for our argument.  
\end{enumerate} 
\end{remark}

Equipped with the initial data from Definition \ref{ansatz:def-initial-data}, we arrive at the initial value problem
\begin{equation*}
\begin{cases}
\partial_v \U[\Ncs][]= \partial_u \V[\Ncs][] = 
 \chinull \Big[  \U[\Ncs][],  \V[\Ncs][]\Big]_{\leq N} - \chinullsquare \coup \Renorm[\Ncs] \big( \U[\Ncs][] + \V[\Ncs][] \big), \\
\U[\Ncs][]\big|_{u=v}= \widetilde{W}^{(\Ncs),+}, \, \V[\Ncs][]\big|_{u=v}=   \widetilde{W}^{(\Ncs),-}, 
\end{cases}
\end{equation*}
which will be our main interest for the rest of this section. \revision{This is the same initial value problem as in \eqref{ansatz:eq-wave-maps}, except that the initial data $(P_{\leq \Nd}^x W^{(\Rscript,\coup),+},P_{\leq \Nd}^x W^{(\Rscript,\coup),-})$ has been replaced with the more general initial data $(\widetilde{W}^{(\Ncs),+},\widetilde{W}^{(\Ncs),-})$ from Definition \ref{ansatz:def-initial-data}. In particular, it has been obtained through (mostly technical) modifications of our finite-dimensional approximation \eqref{intro:eq-system-U-V-Kil}.}

\subsection{Heuristic description of the Ansatz}\label{section:ansatz-heuristic}
As already discussed in the introduction, the wave maps equation \eqref{ansatz:eq-wave-maps} cannot be solved directly via a contraction mapping argument. Instead, we use an Ansatz for the solution of \eqref{ansatz:eq-wave-maps} which captures the random structure of the solution. To be more precise, we write 
\begin{equation}\label{ansatz:eq-UN-decomposition}
\begin{aligned}
\UN &= \UN[][+] +  \UN[][+-] +  \UN[][-] +  \UN[][+\fs] +  \UN[][\fs-] +\UN[][\fs]. 
\end{aligned}
\end{equation}
\revision{The regularities of the terms in \eqref{ansatz:eq-UN-decomposition} are illustrated in Figure \ref{figure:ansatz-regularities}.} For the first five terms in \eqref{ansatz:eq-UN-decomposition}, we often rely on dyadic decompositions given by 
\begin{alignat*}{5}
\UN[][+] &= \Sumlarge_{\substack{ K \leq \Nd }} \UN[K][+], \qquad & 
\UN[][+-]&=\Sumlarge_{\substack{K_u,K_v \leq \Nd \colon \\ K_u \simeq_\delta K_v}} \UN[K_u,K_v][+-], \qquad & 
\UN[][-] &= \Sumlarge_{\substack{K \leq \Nd}} \UN[K][-] \\ 
\UN[][+\fs] &= \Sumlarge_{\substack{ K \leq \Nd}}  \UN[K][+\fs], \qquad &
\UN[][\fs-] &= \Sumlarge_{\substack{ K \leq \Nd}}\UN[K][\fs-].  
\end{alignat*}    
In the above dyadic decompositions, $\sum\limits^{\hspace{0.75pt}\scalebox{0.6}{$\bullet$}}$ is as in \eqref{prelim:eq-sumlarge}.
Similarly, we write 
\begin{equation}\label{ansatz:eq-VN-decomposition}
\begin{aligned}
\VN &= \VN[][-] +  \VN[][+-] +  \VN[][+] +  \VN[][\fs-] +  \VN[][+\fs] +\VN[][\fs] 
\end{aligned}
\end{equation}
and 
\begin{alignat*}{5}
\VN[][-] &= \Sumlarge_{\substack{ M \leq \Nd}} \VN[M][-], \qquad & 
\VN[][+-]&=\Sumlarge_{\substack{M_u,M_v \leq \Nd\colon  \\ M_u \simeq_\delta M_v}} \VN[M_u,M_v][+-], \qquad &
\VN[][+] &= \Sumlarge_{\substack{ M \leq \Nd}} \VN[M][+] , \\ 
\VN[][\fs-] &=  \Sumlarge_{\substack{ M \leq \Nd}} \VN[M][\fs-], \qquad &
\VN[][+\fs] &= \Sumlarge_{\substack{ M \leq \Nd}} \VN[M][+\fs]. 
\end{alignat*}

\begin{remark} 
Due to \eqref{prelim:eq-sumlarge}, all dyadic sums above are taken only over dyadic scales in $\Dyadiclarge$, i.e., only over large dyadic scales. This is for technical convenience, since it allows us to directly use bounds on the Duhamel integral (see Lemma \ref{prelim:lem-Duhamel-integral}) without relying on cut-off functions. The corresponding contributions for dyadic scales in $\dyadic \backslash \Dyadiclarge $, i.e., low dyadic scales, will be absorbed into the nonlinear remainders $\UN[][\fs]$ and $\VN[][\fs]$, which will always appear with cut-off functions.
\end{remark}

The nonlinear remainders $\UN[][\fs]$ and $\VN[][\fs]$ in \eqref{ansatz:eq-UN-decomposition} and \eqref{ansatz:eq-VN-decomposition} will be treated as arbitrary elements in $\Cprod{r-1}{r}$ and $\Cprod{r}{r-1}$, respectively. In particular, we do not require any knowledge of the random structure of $\UN[][\fs]$ or $\VN[][\fs]$. In contrast, the other terms in \eqref{ansatz:eq-UN-decomposition} and \eqref{ansatz:eq-VN-decomposition} all exhibit a useful random structure, which we now discuss at a heuristic level.  Throughout our heuristic discussion, we focus on the terms in \eqref{ansatz:eq-UN-decomposition} since, after reversing the roles of the $u$ and $v$-variables, the terms in \eqref{ansatz:eq-VN-decomposition} are similar. 
We further omit factors involving the cut-off function $\chinull$ and pretend that 
$\widetilde{W}^{(\Ncs),+}$ and $\widetilde{W}^{(\Ncs),-}$ are replaced with frequency-truncated $2\pi $-periodic, $\frkg$-valued white noises, i.e., can be written as 
\begin{equation*}
P_{\leq \Nd}^x W^{(\coup),\pm} = \hcoup \sum_{\ell \in \Z} \rhoND(\ell) G_\ell^{\pm} e^{i\ell x}.
\end{equation*}
For a rigorous treatment of our Ansatz,  we refer to Subsection \ref{section:ansatz-rigorous} below. 

\begin{figure}
\begin{tabular}{
!{\vrule width 1pt}>{\centering\arraybackslash}P{1.5cm}
!{\vrule width 1pt}>{\centering\arraybackslash}P{\regcolwidth}
!{\vrule width 1pt}>{\centering\arraybackslash}P{\regcolwidth}
!{\vrule width 1pt}>{\centering\arraybackslash}P{\regcolwidth}
!{\vrule width 1pt}>{\centering\arraybackslash}P{\regcolwidth}
!{\vrule width 1pt}>{\centering\arraybackslash}P{\regcolwidth}
!{\vrule width 1pt}>{\centering\arraybackslash}P{\regcolwidth}
!{\vrule width 1pt}} 
\noalign{\hrule height 1pt} 
Object &  $\UN[][+]$ & $\UN[][+-]$ & $\UN[][-]$   & $\UN[][+\fs]$ & $\UN[][\fs-]$ & $\UN[][\fs]$   
\\[4pt] \noalign{\hrule height 1pt} \rule{0pt}{14pt}
Space & $\Cprod{s-1}{s}$ & $\Cprod{s-1}{s}$ 
& $\Cprod{-\frac{1}{2}+\eta}{s}$ 
& $\Cprod{s-1}{r}$
& $\Cprod{r-1}{s}$
& $\Cprod{r-1}{r}$
\\[3pt] \noalign{\hrule height 1pt} 
\end{tabular}
\caption{\small{\revision{In this figure, we list product norms in which the terms in \eqref{ansatz:eq-UN-decomposition} can be controlled. The corresponding estimates will be proven in Lemmas \ref{modulation:lem-linear}, \ref{modulation:lem-bilinear}, \ref{modulation:lem-mixed}, and \ref{modulation:lem-linear-reversed} below. We emphasize, however, that the product norms listed above only provide rudimentary information on the terms in \eqref{ansatz:eq-UN-decomposition}, and that our argument provides much more detailed information on their properties.}}}
\label{figure:ansatz-regularities}
\end{figure}

\begin{enumerate}[leftmargin=5ex,label=(\roman*)]
    \item\label{ansatz:item-Up} 
    At a heuristic level, the modulated linear wave $\UN[K][+]$ will be chosen as a solution of the initial value problem\footnote{The precise form of the initial value problem is rather complicated, see e.g. Definition \ref{ansatz:def-modulated-linear} and Definition \ref{ansatz:def-modulation-equations}. \revision{The low-frequency terms will be made precise in Definition \ref{ansatz:def-lo} and Definition \ref{ansatz:def-shhl}.}}
    \begin{equation}\label{ansatz:eq-motivation-Up}
    \begin{cases}
    \begin{aligned}
    \partial_v \UN[K][+] &\simeq  P_{\leq N}^x \Big[ P_{\leq N}^x \UN[K][+], \big\{ \textup{low-frequency terms}\big\} \Big], \\
    \UN[K][+] \Big|_{u=v} &= \hcoup \sum_{k \in \Z_K} \rhoND(k) G_k^{+} e^{ik x}.  
    \end{aligned}
    \end{cases}
    \end{equation}
    The reason for including the term on the right-hand side of \eqref{ansatz:eq-motivation-Up} is the absence of nonlinear smoothing in the $u$-variable, which prevents us from simply using Bourgain's trick (see \cite{B21,BLS21} for related discussions). In order to analyze the interactions between $\UN[K][+]$ and terms from our Ansatz for $\VN[][]$, it is crucial that $\UN[K][+]$ exhibits a random structure. In fact, since \eqref{ansatz:eq-motivation-Up} will only be approximately satisfied, we will be able to write $\UN[K][+]$ as 
    \begin{equation}\label{ansatz:eq-motivation-Up-Sp}
    \UN[K][+] = \hcoup \sum_{k\in \Z_K} \rhoND(k) \SN[K][+][k](u,v) G_k^+ \, e^{iku}.
    \end{equation}
    The $(\SN[K][+][k])_{K\in \Dyadiclarge, k\in \Z_K}$ are called modulation operators and are one of the most central objects in this article. The expression for $\UN[K][+]$ in \eqref{ansatz:eq-motivation-Up-Sp} should be compared with the initial data in \eqref{ansatz:eq-motivation-Up}, in which $\SN[K][+][k]$ is replaced by the identity $\Id_\frkg$. We already emphasize here that the operators $(\SN[K][+][k])_{k\in \Z_K}$ and Gaussians $(G_k^+)_{k\in \Z_K}$ will be probabilistically dependent, which is addressed in Section \ref{section:chaos}. 
\end{enumerate}
The modulated linear wave $\UN[][+]$ addresses the absence of nonlinear smoothing in the $u$-variable. If the absence of nonlinear smoothing were our only problem, it would be sufficient to make the Ansatz $\UN[][]=\UN[][+] +\UN[][\fs]$, where $\UN[][\fs]$ is a smooth nonlinear remainder. In addition to the absence of nonlinear smoothing, however, we also need to address high$\times$high$\rightarrow$low-interactions between $\UN[][]$ and $\VN[][]$ in either the $u$ or $v$-variable. Due to this, we not only need detailed information on the behavior of $\UN[][]$ in the $u$-variable but also in the $v$-variable. This makes it necessary to also introduce the following four structured terms. 
\begin{enumerate}[leftmargin=5ex,label=(\roman*)]
   \setcounter{enumi}{1}
    \item\label{ansatz:item-Upm}  The modulated bilinear wave $\UN[K,M][+-]$ captures interactions between $\UN[K][+]$ and $\VN[M][-]$ when $K$ and $M$ are comparable, i.e., when $K\simeq_\delta M$. At a heuristic level, we choose $\UN[K,M][+-]$ as a solution of 
    \begin{equation}\label{ansatz:eq-motivation-Upm}
    \partial_v \UN[K,M][+-] \simeq  P_{\leq N}^x \Big[ P_{\leq N}^x \UN[K][+], P_{\leq N}^x \VN[M][-] \Big]. 
    \end{equation}
    Thus, $\UN[K,M][+-]$ behaves like $\UN[K][-]$ in the $u$-variable and like the $v$-integral of $\VN[M][-]$ in the $v$-variable. 
    \item\label{ansatz:item-Um}  The reversed modulated linear wave $\UN[M][-]$ captures interactions between low-frequency terms and the modulated linear wave $\VN[M][-]$. At a heuristic level, we choose $\UN[M][-]$ as a solution of 
    \begin{equation*}
    \partial_v \UN[M][-] \simeq  P_{\leq N}^x \Big[ \big\{ \textup{low-frequency terms} \big\}, P_{\leq N}^x \VN[M][-] \Big]. 
    \end{equation*}
    Thus, $\UN[M][-]$ will behave like a low-frequency function in the $u$-variable and like the $v$-integral of $\VN[M][-]$ in the $v$-variable. 
\end{enumerate}
Throughout this article, we refer to the terms from \ref{ansatz:item-Up}-\ref{ansatz:item-Um} as modulated objects. In both the $u$ and $v$-variable, the modulated objects either behave like low-frequency functions, the modulated linear waves, or integrals of the modulated linear waves. For the proof of well-posedness, however, it is insufficient to only isolate the modulated objects. In addition, we also have to isolate so-called mixed modulated objects whose behavior in either the $u$ or $v$-variable can also be dictated by the nonlinear remainders $\UN[][\fs]$ and $\VN[][\fs]$. 
Whereas the modulated objects are naturally placed\footnote{At this point it may seem counter-intuitive that $\UN[K][+]$ is placed only in $\Cprod{s-1}{s}$ and not in $\Cprod{s-1}{\infty}$. After all, the second term in the Lie bracket of \eqref{ansatz:eq-motivation-Up} is claimed to live at low frequencies. The reason why we still only obtain $\Cprod{s-1}{s}$-bounds is that our definition of low frequencies only requires the frequencies to be smaller than $\sim K^{1-}$, which is still quite high.}
in $\Cprod{s-1}{s}$, the mixed modulated objects will naturally be placed in $\Cprod{s-1}{r}$ or $\Cprod{r-1}{s}$. 
\begin{enumerate}[leftmargin=5ex,label=(\roman*)]
   \setcounter{enumi}{3} 
   \item \label{ansatz:item-Ups} The mixed modulated object $\UN[K][+\fs]$ captures interactions between $\UN[K][+]$ and $P^v_{\geq K^{1-\deltap}} \VN[][\fs]$, i.e., the nonlinear remainder in $\VN[][]$ at high $v$-frequencies. At a heuristic level, we choose $\UN[K][+\fs]$ as a solution\footnote{While the right-hand side of \eqref{ansatz:eq-motivation-Ups} is morally correct, it is not quite the right expression. It is missing both a $\VN[][+\fs]$-term in the second argument and a $\Para[v][ll]$-operator. The reasons for both of these changes are rather technical, and we omit the terms here to simplify the expression. For the rigorous definition, we refer the reader to Definition \ref{ansatz:def-mixed}.}  of 
    \begin{equation}\label{ansatz:eq-motivation-Ups}
    \partial_v \UN[K][+\fs] \simeq  P_{\leq N}^x \Big[ P_{\leq N}^x \UN[K][+], P_{\leq N}^x P^v_{\geq K^{1-\deltap}} \VN[][\fs] \Big]. 
    \end{equation}
    Thus, $\UN[K][+\fs]$ behaves like $\UN[K][+]$ in the $u$-variable and like the $v$-integral of $P^v_{\geq K^{1-\deltap}} \VN[][\fs]$ in the $v$-variable. In particular, $\UN[K][+\fs]$ can naturally be placed in $\Cprod{s-1}{r}$. 
    \item\label{ansatz:item-Usm}  The mixed modulated object \revision{$\UN[M][\fs-]$} captures interactions between $P^u_{\geq M^{1-\deltap}}\UN[][\fs]$, i.e., the nonlinear remainder in $\UN[][]$ at high $u$-frequencies, and $\VN[M][-]$. At a heuristic level, we choose $\UN[M][\fs-]$ as a solution\footnote{The same reservations as for \eqref{ansatz:eq-motivation-Ups} also apply to \eqref{ansatz:eq-motivation-Usm}.}  of 
    \begin{equation}\label{ansatz:eq-motivation-Usm}
    \partial_v \UN[M][\fs-] \simeq  P_{\leq N}^x \Big[ P_{\leq N}^x P^u_{\geq M^{1-\deltap}}\UN[][\fs],  P_{\leq N}^x \VN[M][-] \Big]. 
    \end{equation}
    Thus, $\UN[M][\fs-]$ behaves like $P^u_{\geq M^{1-\deltap}}\UN[][\fs]$ in the $u$-variable and like the $v$-integral of $\VN[M][-]$ in the $v$-variable. In particular, $\UN[M][\fs-]$ can naturally be placed in $\Cprod{r-1}{s}$. 
\end{enumerate}

\begin{remark}
We emphasize two differences between our Ansatz in \eqref{ansatz:eq-UN-decomposition}  and the Ansatz in the earlier work \cite[Section 3]{BLS21}. First, the modulation operators $\SN[K][+][k]$ in \eqref{ansatz:eq-motivation-Up-Sp} depend on the frequency variable $k\in \Z_K$, whereas the counterpart $A^+_K$ in \cite[Section 3]{BLS21} did not depend on the frequency variable. The reason for this dependence lies in the frequency-truncation of the wave maps equation (see Definition \ref{ansatz:def-modulation-equations}), and this dependence makes several new insights of this article necessary (see e.g. Section \ref{section:chaos}). Second, the Ansatz in \eqref{ansatz:eq-UN-decomposition} includes $\UN[][+\fs]$ and $\UN[][\fs-]$, whose counterparts in \cite{BLS21} were simply absorbed into the counterparts of $\UN[][+]$ and $\UN[][-]$, respectively. This is necessary since, especially for terms with frequency support near the frequency boundary, we require a more detailed understanding of high$\times$high$\rightarrow$low-interactions than in \cite{BLS21}. 
\end{remark}

\subsection{Rigorous definition of the Ansatz}\label{section:ansatz-rigorous} 
We now turn to a rigorous discussion of our Ansatz and therefore focus on the initial value problem 
\begin{equation}\label{ansatz:eq-wave-maps-tilde}
\begin{cases}
\partial_v \U[\Ncs][]= \partial_u \V[\Ncs][] = 
 \chinull \Big[  \U[\Ncs][],  \V[\Ncs][]\Big]_{\leq N} -  \coup \Renorm[\Ncs] \big( \U[\Ncs][] + \V[\Ncs][] \big), \\
\U[\Ncs][]\big|_{u=v}= \widetilde{W}^{(\Ncs),+}, \, \V[\Ncs][]\big|_{u=v}=   \widetilde{W}^{(\Ncs),-}, 
\end{cases}
\end{equation}
where the initial data $(\widetilde{W}^{(\Ncs),+},\widetilde{W}^{(\Ncs),-})$ is as in Definition \ref{ansatz:def-initial-data}.
To this end, we first recall from \eqref{ansatz:eq-UN-decomposition} and \eqref{ansatz:eq-VN-decomposition} that our Ansatz is given by
\begin{equation}\label{ansatz:eq-UN-rigorous-decomposition}
\begin{aligned}
\UN &= \UN[][+] + \UN[][+-] + \UN[][-] + \UN[][+\fs] + \UN[][\fs-] +\UN[][\fs]  
\end{aligned}
\end{equation}
and 
\begin{equation}\label{ansatz:eq-VN-rigorous-decomposition}
\begin{aligned}
\VN &= \VN[][-] + \VN[][+-] + \VN[][+] + \revision{\VN[][\fs-] + \VN[][+\fs]} + \VN[][\fs]. 
\end{aligned}
\end{equation}
As previously indicated in \eqref{ansatz:eq-motivation-Up-Sp}, our Ansatz also involves the modulation operators 
$(\SN[K][+][k])_{K\in \Dyadiclarge,k\in \Z_K}$ and $(\SN[M][-][m])_{M\in \Dyadiclarge, m \in \Z_M}$. For technical reasons, it is convenient to view $\SN[][+]$ and $\SN[][-]$ as frequency-truncations of so-called pure modulation operators. 

\begin{definition}[Modulation operators and pure modulation operators]\label{ansatz:def-pure}
We introduce the unknowns 
\begin{align}
\pSN[][+]=\big( \pSN[K][+][k] \big)_{K\in \Dyadiclarge, k \in \Z_K} &\colon \R_{u,v}^{1+1} \rightarrow \End(\frkg), \label{ansatz:eq-pure-Sp}\\ 
\pSN[][-]=\big( \pSN[M][-][m] \big)_{M\in \Dyadiclarge, m \in \Z_M} &\colon \R_{u,v}^{1+1} \rightarrow \End(\frkg),\label{ansatz:eq-pure-Sm}
\end{align}
which are called pure modulation operators. Then, the modulation operators 
\begin{align}
\SN[][+]=\big( \SN[K][+][k] \big)_{K\in \Dyadiclarge, k \in \Z_K} &\colon \R_{u,v}^{1+1} \rightarrow \End(\frkg), \label{ansatz:eq-Sp}\\ 
\SN[][-]=\big( \SN[M][-][m] \big)_{M\in \Dyadiclarge, m \in \Z_M} &\colon \R_{u,v}^{1+1} \rightarrow \End(\frkg),\label{ansatz:eq-Sm}
\end{align}
are defined in terms of \eqref{ansatz:eq-pure-Sp} and \eqref{ansatz:eq-pure-Sm} as 
\begin{equation}\label{ansatz:eq-pure-to-truncated}
\SN[K][+][k] = P^{u,v}_{\leq K^{1-\delta+\vartheta}} \pSN[K][+][k] 
\qquad \text{and} \qquad 
\SN[M][-][m] = P^{u,v}_{\leq M^{1-\delta+\vartheta}} \pSN[M][-][m].
\end{equation}
\end{definition}

\begin{remark}
For most of this article, we will only work with the modulation operators from \eqref{ansatz:eq-Sp} and \eqref{ansatz:eq-Sm}. The reason for still introducing the pure modulation operators from \eqref{ansatz:eq-pure-Sp} and \eqref{ansatz:eq-pure-Sm} is as follows: Since the frequency support properties of $\SN[][+]$ and $\SN[][-]$ are important in many of our arguments, we want to know their exact frequency support. At the same time, we want the modulation equations in Definition \ref{ansatz:def-modulation-equations} below to be as simple as possible and therefore prefer to avoid inserting several Littlewood-Paley operators into them. To achieve both, the simplest approach is to work with the two sets of unknowns from \eqref{ansatz:eq-pure-Sp}-\eqref{ansatz:eq-pure-Sm} and \eqref{ansatz:eq-Sp}-\eqref{ansatz:eq-Sm} and then link them via the identity \eqref{ansatz:eq-pure-to-truncated}.
\end{remark}

In total, \eqref{ansatz:eq-UN-rigorous-decomposition}, \eqref{ansatz:eq-VN-rigorous-decomposition}, \eqref{ansatz:eq-pure-Sp}, and \eqref{ansatz:eq-pure-Sm} leave us with fourteen different types of unknowns. The main goal of this subsection is to give explicit definitions of some of the terms in \eqref{ansatz:eq-UN-rigorous-decomposition} and \eqref{ansatz:eq-VN-rigorous-decomposition}, which reduces the number of unknowns. To give the reader a sense of direction, we already state the following fact. 

\begin{fact} Once Definition \ref{ansatz:def-modulated-linear}, Definition \ref{ansatz:def-modulated-bilinear}, Definition \ref{ansatz:def-mixed}, and Definition \ref{ansatz:def-modulated-linear-reversed} below have been made, the terms  
\begin{equation}\label{ansatz:eq-motivation-dependent-unknowns}
\begin{aligned}
&\UN[][+], \UN[][+-],  \UN[][-],  \UN[][+\fs], \UN[][\fs-], \\  
&\VN[][-], \VN[][+-], \VN[][+], \revision{\VN[][\fs-]}, \quad \textup{and} \quad   \revision{\VN[][+\fs]}
\end{aligned}
\end{equation}
are determined uniquely by 
\begin{equation*}
\SN[][+],  \SN[][-], \UN[][\fs], \quad \text{and} \quad \VN[][\fs].  
\end{equation*}
In particular, the definitions of \eqref{ansatz:eq-motivation-dependent-unknowns} are made in terms of the modulation operators $\SN[][\pm]$ rather than the pure modulation operators $\pSN[][\pm]$.
\end{fact}

In the initial stages of our argument (before Section \ref{section:modulation}), we will treat $\pSN[][+],  \pSN[][-], \UN[][\fs]$, and $\VN[][\fs]$ as unknowns. In Section \ref{section:modulation}, we will then solve the modulation equations (see Definition \ref{ansatz:def-modulation-equations}), which are evolution equations for $\pSN[][+]$ and $\pSN[][-]$. In later stages of our argument (after Section \ref{section:modulation}), we will then view  $\pSN[][+]$ and $\pSN[][-]$ (and therefore also $\SN[][+]$ and $\SN[][-]$) as functions of $\UN[][\fs]$ and $\VN[][\fs]$. For expository purposes, it is helpful to make the following definition.

\begin{definition}[Unknowns]\label{ansatz:def-unknowns}
We introduce the following two sets of unknowns:
\begin{enumerate}[label=(\roman*)]
\item (Pre-modulation unknowns) The first set of unknowns includes the pure modulation operators $\pSN[][+]$ and $\pSN[][-]$ from \eqref{ansatz:eq-pure-Sp} and \eqref{ansatz:eq-pure-Sm} and the nonlinear remainders $\UN[][\fs]$ and $\VN[][\fs]$. 
\item (Post-modulation unknowns) The second set of unknowns includes only the nonlinear remainders $\UN[][\fs]$ and $\VN[][\fs]$. 
\end{enumerate}
\end{definition}

We now turn towards the definitions of the terms in \eqref{ansatz:eq-motivation-dependent-unknowns}. To this end, we recall that the frequency-truncated lattice partition $(\psiRx)_{x_0\in \LambdaRR}$ was introduced in Definition \ref{prelim:def-lattice-partition-truncated}.

\begin{definition}[Modulated linear waves]\label{ansatz:def-modulated-linear}
For all $K,M,N,\Nd \in \Dyadiclarge$, we define
\begin{align}
\UN[K][+](u,v) &= \hcoup \sum_{u_0 \in \LambdaRR}  \psiRuK(u) \sum_{k\in \Z_K} \rhoND(k) \SN[K][+][k](u,v) G_{u_0,k}^+ \, e^{iku}, \label{ansatz:eq-Up}\\
\VN[M][-](u,v) &= \hcoup \sum_{v_0 \in \LambdaRR} \psiRvM(v) \sum_{m\in \Z_M} \rhoND(m) \SN[M][-][m](u,v) G_{v_0,m}^- \, e^{imv}. 
\label{ansatz:eq-Vm}
\end{align}
In particular, $\UN[][+]$ depends only on $\SN[][+]$ and  $\VN[][-]$ depends only on $\SN[][-]$.  
\end{definition}

\begin{remark}
We emphasize that the modulation operators $\SN[K][+][k]$ in \eqref{ansatz:eq-Up} depend on $k\in \Z_K$ but do  not depend on $u_0 \in \LambdaR$. The reason is that the Littlewood-Paley operators $P_{\leq N}^x$ do not at all commute with multiplication by $e^{iku}$, but do approximately commute with multiplication by $\psiRuK$.
\end{remark}

\begin{remark}
To avoid confusion, we note that $\UN[K][+]$ from \eqref{ansatz:eq-Up} does not yet satisfy \eqref{ansatz:eq-motivation-Up}, which was described as the motivation for $\UN[K][+]$. The initial value problem \eqref{ansatz:eq-motivation-Up} will only be satisfied due to the initial value problem for $\pSN[K][+]$, which is the subject of Definition \ref{ansatz:def-modulation-equations} below.  
\end{remark}

To simplify the notation below, we now introduce integrated versions of $\UN[K][+]$ and $\VN[M][-]$.

\begin{definition}[Integrated modulated linear waves]\label{ansatz:def-integrated-modulated-waves}
For all $K,M,N,\Nd \in \Dyadiclarge$ satisfying $K,M\leq \Nd$, we define
\begin{align}
\IUN[K][+](u,v) &=\hcoup \sum_{u_0 \in \LambdaRR} \psiRuK \sum_{k\in \Z_K} \rhoND(k) \SN[K][+][k](u,v) G_{u_0,k}^+ \, \frac{e^{iku}}{ik}, \\
\IVN[M][-](u,v) &= \hcoup \sum_{v_0 \in \LambdaRR} \psiRvM \sum_{m\in \Z_M} \rhoND(m) \SN[M][-][m](u,v) G_{v_0,m}^- \, \frac{e^{imv}}{im}. 
\end{align}
\end{definition}

Before we can introduce the modulated bilinear waves, we introduce frequency-truncated versions of the cut-off function $\chinull$.

\begin{definition}[Frequency-truncated cut-off functions]\label{ansatz:def-cutoff-frequency-truncated}
Let $\chi \in C^\infty_b(\R)$ and let $K,L,M\in \dyadic$. Then, we define
\begin{equation*}
\chi_L(t) := \big( P_{\leq L^{\vartheta}} \chi \big)(t) \qquad \text{and} \qquad \chi_{K,M}(t):= \chi_{\max(K,M)}(t). 
\end{equation*}
Similar as in Remark \ref{ansatz:rem-cut-off-null}, we also write $\chinull[L]$ and $\chinull[K,M]$ for the corresponding functions 
\begin{equation*}
    (u,v)\mapsto \chi_L(\frac{v-u}{2}) \qquad \text{and} \qquad (u,v) \mapsto \chi_{K,M}(\frac{v-u}{2}). 
\end{equation*} 
\end{definition}

\begin{remark}[Commutativity of $\chinull$ and $P_{\leq N}^x$]\label{ansatz:rem-commutativity-chipm}
Since the functions $\chinull$, $\chinull[L] $, and $\chinull[K,M]$ only depend on the time-variable $t$ or, equivalently, only the difference $v-u$, all of them commute with the Littlewood-Paley operators $P_{\leq N}^x$. This will be used repeatedly throughout this article. 
\end{remark}

\begin{remark}
Due to the cut-off functions in \eqref{ansatz:eq-wave-maps}, the rigorous definitions of our modulated and mixed modulated objects should include factors of $\chinull$. In order to exactly preserve the frequency-support properties of our objects, however, we will often replace $\chinull$ with $\chinull[L]$ or $\chinull[K,M]$. This is of course only for technical convenience, since the estimate $\| \chinull - \chinull[L] \|_{L^\infty}\lesssim L^{-100}$ will always allow us to interchange $\chinull$ and $\chinull[L]$.
\end{remark}

Equipped with Definition \ref{ansatz:def-integrated-modulated-waves} and Definition \ref{ansatz:def-cutoff-frequency-truncated}, we can now define the bilinear modulated linear waves. 

\begin{definition}[Modulated bilinear waves]\label{ansatz:def-modulated-bilinear} 
Let $K,M,N,\Nd\in \Dyadiclarge$ satisfy $K,M\leq \Nd$. Then, the modulated bilinear waves as  $\UN[K,M][+-]$ and $ \UN[M][-]$ are defined as 
\begin{align*}
\UN[K,M][+-] 
&:=  \chinull[K,M] \Big[ \UN[K][+], \IVN[M][-] \Big]_{\leq N}  \\
&=  \coup \, \chinull[K,M]  \hspace{-0.5ex} \sum_{u_0,v_0\in \LambdaRR} \sum_{k \in \Z_{K}} \sum_{m\in \Z_M} \rhoND(k) \rhoND(m)
\bigg[ \psiRuK \SN[K][+][k] G_{u_0,k}^+ \, e^{iku},    \psiRvM \SN[M][-][m] G_{v_0,m}^- \frac{e^{imv}}{im} \bigg]_{\leq N}, \\
\VN[K,M][+-] 
&:= \Big[ \IUN[K][+], \VN[M][-] \Big]_{\leq N}  \\ 
&= \coup \, \chinull[K,M]  \hspace{-0.5ex} \sum_{u_0,v_0\in \LambdaRR} \sum_{k \in \Z_{K}} \sum_{m\in \Z_M} \rhoND(k) \rhoND(m)
\bigg[ \psiRuK \SN[K][+][k] G_{u_0,k}^+ \, \frac{e^{iku}}{ik}, \psiRvM \SN[M][-][m] G_{v_0,m}^- e^{imv} \bigg]_{\leq N}. 
\end{align*}
In particular, $\UN[K,M][+-]$ and $\VN[K,M][+-]$ are determined by $\SN[K][+]$ and $\SN[M][-]$. 
\end{definition}

In \eqref{ansatz:eq-UN-rigorous-decomposition}, we first list the modulated objects $\UN[][+]$, $\UN[][+-]$, and $\UN[][-]$ and then list the mixed modulated objects $\UN[][+\fs]$ and $\UN[][\fs-]$. For our rigorous definitions, however, it is more convenient to first define the mixed modulated objects $\UN[][+\fs]$ and $\UN[][\fs-]$ and only then define the last modulated object $\UN[][-]$. The reason is that, as we will see in Definition \ref{ansatz:def-modulated-linear-reversed}, the definition of $\UN[][-]$ involves $\UN[][+\fs]$ and $\UN[][\fs-]$.

\begin{definition}[Mixed modulated objects]\label{ansatz:def-mixed}
Let $N,\Nd\in \Dyadiclarge$. Then, we define $(\UN[M][\fs-])_M$ and $(\VN[K][+\fs])_K$, where $K,M\in \Dyadiclarge$ satisfy $K,M\leq \Nd$, as the unique solution of the systems
\begin{align}
\UN[M][\fs-] &=  \chinull[M]
\bigg[ P^{u}_{\geq M^{1-\deltap}} \Big( P^v_{<M^{1-\deltap}} \UN[][\fs]
+ \Sumlarge_{\substack{K< M^{1-\delta}}} \UN[K][\fs-] \Big) \Para[u][gg] \IVN[M][-] \bigg]_{\leq N}, \label{ansatz:eq-rigorous-Usm} \\ 
\VN[K][+\fs]&= \chinull[K] 
\bigg[ \IUN[K][+] \Para[v][ll] P^v_{\geq K^{1-\deltap}} \Big( P^u_{<K^{1-\deltap}} \VN[][\fs] + \Sumlarge_{\substack{M< K^{1-\delta}}} \VN[M][+\fs] \Big) \bigg]_{\leq N}.\label{ansatz:eq-rigorous-Vps}
\end{align}
Furthermore, we define $(\UN[K][+\fs])_{K}$ and $(\VN[M][\fs-])_{M}$ explicitly as 
\begin{align}
\UN[K][+\fs] &=  \chinull[K]  \bigg[ \UN[K][+] \Para[v][ll] \Int^v_{u\rightarrow v} P^v_{\geq K^{1-\deltap}} \Big( P^u_{<K^{1-\deltap}} \VN[][\fs] + \Sumlarge_{\substack{ M< K^{1-\delta}}}  \VN[M][+\fs] \Big) \bigg]_{\leq N}, \label{ansatz:eq-rigorous-Ups} \\ 
\VN[M][\fs-] &= \chinull[M] \bigg[ \Int^u_{v\rightarrow u} P^{u}_{\geq M^{1-\deltap}} \Big( P^v_{<M^{1-\deltap}} \UN[][\fs]
+ \Sumlarge_{\substack{ K< M^{1-\delta}}} \UN[K][\fs-] \Big) \Para[u][gg]  \VN[M][-] \bigg]_{\leq N}. \label{ansatz:eq-rigorous-Vsm}
\end{align}
In particular, $\UN[][+\fs],\UN[][\fs-],\VN[][\fs-]$, and $\VN[][+\fs]$ are determined by $\SN[][+]$, $\SN[][-]$, $\UN[][\fs]$, and $\VN[][\fs]$. 
\end{definition}

\begin{remark}\label{ansatz:remark-recursive-mixed}
For $M=\Nlarge$, there is no $K\in \Dyadiclarge$ satisfying $K<M^{1-\delta}$, and thus the dyadic sum in \eqref{ansatz:eq-rigorous-Usm} is empty. As a result, \eqref{ansatz:eq-rigorous-Usm} yields an explicit definition of $\UN[\Nlarge][\fs-]$. For all $M\in \Dyadiclarge$ satisfying $M>\Nlarge$, \eqref{ansatz:eq-rigorous-Usm} defines $\UN[M][\fs-]$ as a function of $\UN[][\fs]$, $(\UN[K][\fs-])_{K\in \Dyadiclarge:\, K<M^{1-\delta}}$, and $\SN[M][-]$. In particular, \eqref{ansatz:eq-rigorous-Usm} can be used to define $(\UN[M][\fs-])_{M\in \Dyadiclarge}$ recursively as a function of $\UN[][\fs]$ and $(\SN[M][-])_{M\in \Dyadiclarge}$. Similarly, \eqref{ansatz:eq-rigorous-Vps} can be used to define  $(\VN[K][+\fs])_{K\in \Dyadiclarge}$ recursively. Once \eqref{ansatz:eq-rigorous-Usm} and \eqref{ansatz:eq-rigorous-Vps} have been solved, \eqref{ansatz:eq-rigorous-Ups} and \eqref{ansatz:eq-rigorous-Vsm} yield explicit definitions of $(\UN[K][+\fs])_{K\in \Dyadiclarge}$ and $(\VN[M][\fs-])_{M\in \Dyadiclarge}$, respectively. 
We also remark that while the recursive structure can be used to define the mixed modulated objects, it will not be essential for our quantitative estimates (see Section \ref{section:modulated-mixed}). 
\end{remark}

\begin{remark}
Since the first argument already contains $P^u_{\geq M^{1-\deltap}}$, 
the $\Para[u][gg]$-operator in \eqref{ansatz:eq-rigorous-Usm} may at first appear superfluous. However, since $\IVN[M][-]$ has $u$-frequencies $\lesssim M^{1-\delta+\vartheta}$, which includes the range between $M^{1-\deltap}$ and $M^{1-\delta+\vartheta}$, the  $\Para[u][gg]$-operator is necessary to avoid high$\times$high$\rightarrow$low-interactions in the $u$-variable. 
\end{remark}

\begin{remark}
We note that the rigorous definition of $\UN[K][+\fs]$ from \eqref{ansatz:eq-rigorous-Ups} differs slightly from the heuristic description in \eqref{ansatz:eq-motivation-Ups}. The most important difference is that the right-hand side of \eqref{ansatz:eq-rigorous-Ups} contains terms of the form $P^v_{\geq K^{1-\deltap}} \VN[M][+\fs]$, which were omitted in \eqref{ansatz:eq-motivation-Ups}. These terms have to be included for the same reason as the $\VN[][\fs]$-terms since, at least when $M$ is much smaller than $K$, $P^v_{\geq K^{1-\deltap}} \VN[M][+\fs]$ essentially behaves like $P^v_{\geq K^{1-\deltap}} P_M^u \VN[][\fs]$. Similarly, \eqref{ansatz:eq-rigorous-Usm} contains more terms than \eqref{ansatz:eq-motivation-Usm}. 
\end{remark}

Before we can define the last remaining terms in our Ansatz, i.e., $\UN[][-]$ and $\VN[][+]$, we need to introduce additional notation. In the next definition, we introduce $\LON[K][-]$ and $\LON[M][+]$, which combine several low-frequency terms.

\begin{definition}[Low-frequency terms]\label{ansatz:def-lo} 
Let $K,M,N,\Nd\in \Dyadiclarge$ satisfy $K,M\leq \Nd$. Then, we define
\begin{equation}\label{ansatz:eq-LON-p}
\begin{aligned}
\LON[K][-]
&:=  \Sumlarge_{\substack{ M <K^{1-\delta}}} \VN[M][-] 
+ \Sumlarge_{\substack{ M_u,M_v < K^{1-\delta}, \\M_u \simeq_\delta M_v}}\VN[M_u,M_v][+-] 
+ \Sumlarge_{\substack{ M <K^{1-\delta}}} \VN[M][+] \\ 
&+ \Sumlarge_{\substack{ M <K^{1-\delta}}} \Big( P^u_{<K^{1-\deltap}} \VN[M][\fs -] + P^v_{<K^{1-\deltap}} \VN[M][+\fs] \Big) 
+ P^{u,v}_{<K^{1-\deltap}} \VN[][\fs]. 
\end{aligned}
\end{equation}
Similarly, we define
\begin{equation}\label{ansatz:eq-LON-m}
\begin{aligned}
\LON[M][+] &:= 
\Sumlarge_{\substack{ K< M^{1-\delta}}} \UN[K][+] 
+ \Sumlarge_{\substack{ K_u,K_v< M^{1-\delta}, \\ K_u \simeq_\delta K_v}} \UN[K_u,K_v][+-] 
+ \Sumlarge_{\substack{K< M^{1-\delta}}} \UN[K][-] \\
&+ \Sumlarge_{\substack{ K< M^{1-\delta}}} \Big( P^v_{<M^{1-\deltap}} \UN[K][+\fs] + P^u_{<M^{1-\deltap}} \UN[K][\fs-] \Big) + P^{u,v}_{<M^{1-\deltap}} \UN[][\fs]. 
\end{aligned}
\end{equation}
\end{definition}

\begin{remark}[A technical comment regarding the low-frequency terms]\label{ansatz:rem-delta-1}
We emphasize that the definition of $\LON[K][-]$ in \eqref{ansatz:eq-LON-p} contains two different frequency-thresholds. The modulated linear wave $\VN[M][-]$ is included when $M< K^{1-\delta}$, i.e., when the $v$-frequency is bounded by $\lesssim K^{1-\delta}$.  In contrast, the smooth remainder $\VN[][\fs]$ only enters as $P^{u,v}_{<K^{1-\deltap}} \VN[][\fs]$, which restricts to  $v$-frequencies $\lesssim K^{1-\deltap}$. Since $\deltap > \delta$, we therefore include fewer $v$-frequencies of $\VN[][\fs]$ than of $\VN[][-]$. 
This is for technical reasons related to high$\times$high$\rightarrow$low-interactions which will be discussed in Remark \ref{ansatz:rem-delta-2}. 
\end{remark}

In the next definition, we define simplified high$\times$high$\rightarrow$low-interactions. Since a detailed discussion of the motivation behind this definition requires more background, it is postponed until Section \ref{section:ansatz-remainder}.  To simplify the notation in the simplified high$\times$high$\rightarrow$low-interactions, we first introduce the combined smooth remainders.

\begin{definition}[Combined smooth remainders]\label{ansatz:def-combined}
For all $K,M,N,\Nd \in \Dyadiclarge$ satisfying $K,M\leq \Nd$, we define 
\begin{align}
\UN[< M^{1-\delta}][\fsc] &:= P^v_{<M^{1-\deltap}} \UN[][\fs] 
+ \Sumlarge_{ L < M^{1-\delta}}\UN[L][\fs-], \label{ansatz:eq-UN-sast} \\
\VN[< K^{1-\delta}][\fcs] &:= P^u_{<K^{1-\deltap}} \VN[][\fs] + \Sumlarge_{\substack{ L < K^{1-\delta}}} \VN[L][+\fs]. \label{ansatz:eq-VN-asts}
\end{align}
\end{definition}

The superscript ``$\fsc"$ in \eqref{ansatz:eq-UN-sast} indicates that $\UN[][\fsc]$ consists of both $\UN[][\fs]$ and $\UN[][\fs-]$-terms, i.e., behaves like a smooth remainder in the $u$-variable but like a combination of smooth remainders and modulated objects in the $v$-variable. The superscript ``$\fcs$" has a similar motivation, but we flipped the order of $\fs$ and $\fc$ since  \eqref{ansatz:eq-VN-asts} contains $\VN[][+\fs]$-terms. 
We note that $\UN[< M^{1-\delta}][\fsc]$ is supported on low $v$-frequencies, i.e., on $v$-frequencies bounded by either $\lesssim M^{1-\deltap}$ (for the $\UN[][\fs]$-term) or $\lesssim M^{1-\delta}$ (for the $\UN[][+\fs]$-terms), but there is no restriction on the $u$-frequencies. 

\begin{definition}[Simplified high$\times$high$\rightarrow$low-term]\label{ansatz:def-shhl}
Let $K,M,N,\Nd \in \Dyadiclarge$. Then, we define
\begin{equation*}
\begin{aligned}
\SHHLN[M][u] &:=  \chinull[M]
\Sumlarge_{\substack{ K \leq \Nd\colon \\  K \lesssim_\delta M}} P^{u,v}_{<M^\delta}
\Big[ P_{\leq N}^x  P^u_{\geq M^{1-\deltap}} \UN[< M^{1-\delta}][\fsc], P_{\leq N}^x \IUN[K][+] \Big], \\
\SHHLN[K][v] &:= - \chinull[K] 
\Sumlarge_{\substack{ M \leq \Nd\colon \\  M \lesssim_\delta K}} P^{u,v}_{<K^\delta}
\Big[ P_{\leq N}^x P^v_{\geq K^{1-\deltap}}  \VN[< K^{1-\delta}][\fcs], 
P_{\leq N}^x \IVN[M][-] \Big]. 
\end{aligned}
\end{equation*}
\end{definition}

Equipped with Definition \ref{ansatz:def-integrated-modulated-waves}, Definition \ref{ansatz:def-lo}, and Definition \ref{ansatz:def-shhl}, we can now introduce $\UN[][-]$ and $\VN[][+]$. 

\begin{definition}[Reversed modulated linear waves]\label{ansatz:def-modulated-linear-reversed}
Let $N,\Nd \in \Dyadiclarge$. Then, the reversed modulated linear waves $(\UN[M][-])_{M\in \Dyadiclarge}$ and $(\VN[K][+])_{K\in \Dyadiclarge}$ are defined as the unique solution to the system of equations
\begin{align}
\UN[M][-] 
&=  \chinull[M]\Big[ \LON[M][+], \IVN[M][-] \Big]_{\leq N}
+  \chinull[M]\Big[ \SHHLN[M][u], (P_{\leq N}^x)^2  \IVN[M][-] \Big]_{\leq N}, \label{ansatz:eq-Um-rigorous} \\
\VN[K][+] &=  \chinull[K] \Big[  \IUN[K][+] , \LON[K][-] \Big]_{\leq N}
+ \chinull[K] \Big[   (P_{\leq N}^x)^2  \IUN[K][+] , \SHHLN[K][v] \Big]_{\leq N}. \label{ansatz:eq-Vm-rigorous} 
\end{align}
In particular, $\UN[][-]$ and $\VN[][+]$ are functions of 
\begin{equation}\label{ansatz:eq-Um-Vm-dependence}
\SN[][+], \SN[][-], \UN[][+\fs],\UN[][\fs-], \UN[][\fs], \VN[][\fs-], \VN[][+\fs], \quad \text{and} \quad \VN[][\fs]. 
\end{equation}
\end{definition}

\begin{remark}\label{ansatz:rem-PNx-SHHLN}
The additional $(P_{\leq N})^2$-operators in the $\SHHLN$-terms stem from iterated, frequency-truncated Lie brackets. For more details, see e.g. the proof of Proposition \ref{jacobi:prop-main}. 
\end{remark}

\begin{remark}\label{ansatz:remark-recursive}
We note that $\LON[M][+]$ and $\SHHLN[M][u]$ only depend on $(\UN[L][-])_{L \in \Dyadiclarge: \, L < M^{1-\delta}}$ and terms in \eqref{ansatz:eq-Um-Vm-dependence}. Thus, \eqref{ansatz:eq-Um-rigorous} yields an explicit definition of $\UN[\Nlarge][-]$ and a recursive definition of $\UN[M][-]$ for all $M> \Nlarge$.  Similarly, \eqref{ansatz:eq-Vm-rigorous} yields an explicit definition of $\VN[\Nlarge][+]$ and a recursive definition of $\VN[K][+]$ for all $K> \Nlarge$. 
\end{remark}

\begin{remark}\label{ansatz:rem-PNX-LON}
Due to Definition \ref{ansatz:def-lo} and Definition \ref{ansatz:def-shhl}, $\LON[K][-]$ and $\SHHLN[K][v]$ are supported on frequencies $\lesssim K^{1-\delta}$ and $\lesssim K^\delta$ \revision{in both the $u$ and $v$-variables}, respectively. In particular, since $K\lesssim N$, it always holds that
\begin{equation}\label{ansatz:eq-PNX-LON}
P_{\leq N}^x \LON[K][-] = \LON[K][-] \qquad \text{and} \qquad
P_{\leq N}^x \SHHLN[K][v] = \SHHLN[K][v]. 
\end{equation}
Thus, while the $P_{\leq N}^x$-operators are important in many terms, they affect neither $\LON[K][-]$ nor $\SHHLN[K][v]$.
\end{remark}

As a result of our definitions, all the modulated and mixed modulated objects obey certain frequency support properties. Since the frequency support properties will be used repeatedly in our argument, we collect them in the following lemma. 
 \begin{lemma}[Frequency support properties]\label{ansatz:lem-frequency-support}
  Let $K,M,N,\Nd\in \Dyadiclarge$ satisfy $K,M\leq \Nd$. Then, our modulation operators, modulated objects, and mixed modulated objects satisfy the following frequency support properties: 
  \begin{enumerate}[label=(\roman*)]
  \item\label{ansatz:item-freq-1} $\SN[K][+][k]$ is supported on $u$ and $v$-frequencies $\lesssim K^{1-\delta+\vartheta}$. 
  \item\label{ansatz:item-freq-2} $\UN[K][+]$ is supported on $u$-frequencies $\sim K$ and $v$-frequencies $\lesssim K^{1-\delta+\vartheta}$. 
  \item\label{ansatz:item-freq-3} $\IUN[K][+]$ is supported on $u$-frequencies $\sim K$ and $v$-frequencies $\lesssim K^{1-\delta+\vartheta}$. 
  \item\label{ansatz:item-freq-4} $\UN[K,M][+-]$ is supported on $u$-frequencies $\lesssim \max(K,M^{1-\delta+\vartheta})$ 
  and $v$-frequencies $\lesssim \max(K^{1-\delta+\vartheta},M)$.
 \item\label{ansatz:item-freq-5} $\UN[M][-]$ is supported on $u$-frequencies $\lesssim M^{1-\delta}$ and $v$-frequencies $\sim M$.
  \item\label{ansatz:item-freq-6} $\UN[K][+\fs]$ is supported on $u$-frequencies $\sim K$ and $v$-frequencies $\gtrsim K^{1-\deltap}$. 
  \item\label{ansatz:item-freq-7} $\UN[M][\fs-]$ is supported on $u$-frequencies $\gtrsim M^{1-\deltap}$ and $v$-frequencies $\sim M$. 
  \end{enumerate}
  The frequency support properties of $\SN[M][-][m]$ and the modulated and mixed modulated objects in \eqref{ansatz:eq-VN-decomposition} are similar but with reversed roles of the $u$ and $v$-variables. 
  \end{lemma}

\begin{remark}
On first reading, the reader may pretend that $\UN[K,M][+-]$ is supported on $u$-frequencies $\sim K$ (and $v$-frequencies $\sim M$), but this is not completely true. The reason is that $\VN[M][-]$ appears in the definition of $\UN[K,M][+-]$ and $\VN[M][-]$ is supported on $u$-frequencies $\lesssim M^{1-\delta+\vartheta}$. Since $\UN[K,M][+-]$ is defined in the regime $K\simeq_\delta M$, 
one can estimate $M^{1-\delta+\vartheta}$ from above by $M^\vartheta K$, but this is slightly bigger than $K$. 
\end{remark}

\begin{proof}
The frequency support properties directly follow from our definitions. Indeed, the property in \ref{ansatz:item-freq-1} follows from Definition \ref{ansatz:def-pure}. The properties in \ref{ansatz:item-freq-2} and \ref{ansatz:item-freq-3} follow from Definition \ref{ansatz:def-modulated-linear}, Definition \ref{ansatz:def-integrated-modulated-waves}, and \ref{ansatz:item-freq-1}. The property in \ref{ansatz:item-freq-4} follows from Definition \ref{ansatz:def-modulated-bilinear}, \ref{ansatz:item-freq-2}, and \ref{ansatz:item-freq-3}. The properties in \ref{ansatz:item-freq-6} and \ref{ansatz:item-freq-7} follow from Definition \ref{ansatz:def-mixed}, \ref{ansatz:item-freq-2}, and \ref{ansatz:item-freq-3}. Finally, the property in \ref{ansatz:item-freq-5} follows from Definition \ref{ansatz:def-lo}, Definition \ref{ansatz:def-shhl}, Definition \ref{ansatz:def-modulated-linear-reversed} and \ref{ansatz:item-freq-2}-\ref{ansatz:item-freq-6}.
\end{proof}

Equipped with Definition \ref{ansatz:def-modulated-linear}, Definition \ref{ansatz:def-modulated-bilinear}, Definition \ref{ansatz:def-mixed}, and Definition \ref{ansatz:def-modulated-linear-reversed}, we have now defined all terms in  \eqref{ansatz:eq-motivation-dependent-unknowns}. It thus remains to determine the evolution equations for the pre-modulation unknowns $\pSN[][+]$, $\pSN[][-]$, $\UN[][\fs]$, and $\VN[][\fs]$. To this end, we recall that our heuristic discussion of $\UN[K][+]$ involved the initial value problem \eqref{ansatz:eq-motivation-Up}. While the low-frequency terms in \eqref{ansatz:eq-motivation-Up} were left unspecified and the initial data did not match Definition \ref{ansatz:def-initial-data}, we can now be more precise. Equipped with Definition \ref{ansatz:def-lo} and \revision{Definition \ref{ansatz:def-shhl}}, we can reformulate \eqref{ansatz:eq-motivation-Up} as 
\begin{equation}\label{ansatz:eq-motivation-Up-again}
    \begin{cases}
    \begin{aligned}
    \partial_v \UN[K][+] &\simeq   \chinull \Big[  \UN[K][+], \LON[K][-] \Big]_{\leq N} + \chinull \Big[ (P_{\leq N}^x)^2  \UN[K][+], \SHHLN[K][v] \Big]_{\leq N}, \\
    \UN[K][+] \Big|_{u=v} &= \revision{\hcoup} \sum_{u_0\in \LambdaRR}  \sum_{k\in \Z_K}  \psiRuK(u) \rhoND(k) \SNin[K][+](u) G^+_{u_0,k} e^{iku}.
    \end{aligned}
    \end{cases}
\end{equation}
Since $\UN[K][+]$ is defined in terms of $\SN[K][+]$ and $\SN[K][+]$ is defined in terms of $\pSN[K][+]$ (see Definition \ref{ansatz:def-modulated-linear} and Definition \ref{ansatz:def-pure}), this suggests an initial value problem for $\pSN[K][+]$.

\begin{definition}[Modulation equations]\label{ansatz:def-modulation-equations}
Let $N,\Nd \in \Dyadiclarge$ and let 
\begin{equation*}
\big( \pSN[K][+][k] \big)_{K\in \Dyadiclarge, k \in \Z_K} \colon 
 \R_{u,v}^{1+1} \rightarrow \End(\frkg), 
\qquad \text{and} \qquad \big( \pSN[M][-][m] \big)_{M\in \Dyadiclarge, m \in \Z_M} \colon 
 \R_{u,v}^{1+1} \rightarrow \End(\frkg). 
\end{equation*}
Then, we say that $\pSN[][+]$ and $\pSN[][-]$ solve the modulation equations if the evolution equations
\begin{align}
\partial_v \pSN[K][+][k]  
&=  \chinull[K] \rho_{\leq N}^2(k)  \Big[ \pSN[K][+][k], \LON[K][-] \Big] 
+ \chinull[K]  \rho_{\leq N}^4(k)   \Big[ \pSN[K][+][k], \SHHLN[K][v] \Big], \label{ansatz:eq-modulation-Sp} \\ 
\partial_u \pSN[M][-][m]  
&=  \chinull[M]\rho_{\leq N}^2(m)   \Big[  \LON[M][+] , \pSN[M][-][m] \Big] 
+  \chinull[M]\rho_{\leq N}^4(m)  \Big[ \SHHLN[M][u], \pSN[M][-][m] \Big] \label{ansatz:eq-modulation-Sm}
\end{align}
and the initial conditions
\begin{equation*}
\pSN[K][+][k](u,v) \Big|_{v=u} = \SNin[K][+](u) \qquad \text{and} \qquad 
\pSN[M][-][m](u,v) \Big|_{u=v} = \SNin[M][-](v)
\end{equation*}
are satisfied. 
\end{definition}

\begin{remark}
Since $\UN[K][+](u,v)= \sum_{u_0 \in \LambdaRR} \sum_{k\in \Z_K} \psiRuK(u) \rhoND(k) \SN[K][+][k](u,v)G_{u_0,k}^+ e^{iku}$, we have replaced the actions of all $P_{\leq N}^x$-operators on $\UN[K][+]$ in  \eqref{ansatz:eq-motivation-Up-again} with $\rho_{\leq N}(k)$ in \eqref{ansatz:eq-modulation-Sp}. Without this replacement, \eqref{ansatz:eq-modulation-Sp} would not be an ordinary differential equation in the $v$-variable. 
We also remark that the dependence of $\pSN[K][+][k]$ and $\SN[K][+][k]$ on $k\in \Z_K$ is only due to the $\rho_{\leq N}^2(k)$ and $\rho_{\leq N}^4(k)$-factors in \eqref{ansatz:eq-modulation-Sp}. In particular, this dependence only occurs when $K\sim N$, i.e., when $K$ is near the frequency boundary. 
\end{remark}

\begin{remark}\label{ansatz:rem-initial-data-modulation-equations}
\revision{For now, we only impose the conditions from Definition \ref{ansatz:def-initial-data} on the initial data of the modulation equations. In fact, at least for now, the reader may simply think of $\SNin[K][+]$ and $\SNin[M][-]$ as being equal to $\Id_\frkg$. Later on, however, $\SNin[K][+]$ and $\SNin[M][-]$ will just be close to $\Id_\frkg$ in the $\C_x^s$-norm (see Definition \ref{main:def-perturbations-null}). The added flexibility in the initial data will be important in proving that our two different approximations
\eqref{intro:eq-system-U-V-N-data} and \eqref{intro:eq-system-U-V-Kil} converge to the same limit (see Lemma \ref{main:lem-abstract} and Remark \ref{main:rem-abstract-motivation}).}
\end{remark}

The modulation equations \eqref{ansatz:eq-modulation-Sp}-\eqref{ansatz:eq-modulation-Sm} are a system of first-order ordinary differential equations for the modulation operators. However, as will be described in Subsection \ref{section:ansatz-modulation} below, the regularity of $\LON[K][-]$ and $\LON[M][+]$ prohibits us from solving \eqref{ansatz:eq-modulation-Sp}-\eqref{ansatz:eq-modulation-Sm} using classical ODE-methods. Instead, our argument relies on the para-controlled approach to rough ODEs of \cite{GIP15}. \\

To complete this subsection, it only remains to state the evolution equations for $\UN[][\fs]$ and $\VN[][\fs]$. We recall that $\UN[][]$ and $\VN[][]$ are solutions of the discretized wave maps equation \eqref{ansatz:eq-wave-maps-tilde} and that, after solving the modulation equations from Definition \ref{ansatz:def-modulation-equations}, the terms 
\begin{equation*}
\begin{aligned}
&\UN[][+], \UN[][+-],  \UN[][-],  \UN[][+\fs], \UN[][\fs-], \\  
&\VN[][-], \VN[][+-], \VN[][+], \revision{\VN[][\fs-]}, \quad \textup{and} \quad   \revision{\VN[][+\fs]}
\end{aligned}
\end{equation*}
are functions of $\UN[][\fs]$ and $\VN[][\fs]$. As a result, the evolution equations for $\UN[]$ and $\VN[]$ are forced upon us by \eqref{ansatz:eq-wave-maps-tilde}, our Ansatz, and the definitions of this subsection.

\begin{definition}[Remainder equations]\label{ansatz:def-remainder-equations}
Let $N,\Nd\in \Dyadiclarge$ and let $\UN[][\fs],\VN[][\fs]\colon \R_{u,v}^{1+1}\rightarrow \frkg$. Then, we call $\UN[]$ and $\VN[]$ solutions of the remainder equations if the evolution equations
\begin{equation}\label{ansatz:eq-remainder-equation-Us}
\begin{aligned}
\partial_v \UN[][\fs] 
&= \chinull \big[ \UN, \VN \big]_{\leq N} -  \chinullsquare \coup  \Renorm[\Nc]  \big( \UN + \VN \big) \\
&- \partial_v \big( \UN[][+] + \UN[][+-] + \UN[][-] + \UN[][+\fs]+\UN[][\fs-] \big)
\end{aligned}
\end{equation}
and 
\begin{equation}\label{ansatz:eq-remainder-equation-Vs}
\begin{aligned}
\partial_u \VN[][\fs] 
&= \chinull \big[ \UN, \VN \big]_{\leq N} - \chinullsquare \coup \Renorm[\Nc]  \big( \UN + \VN \big) \\
&- \partial_u \big( \VN[][-] + \VN[][+-] + \VN[][+] + \VN[][\fs-] + \VN[][+\fs] \big)
\end{aligned}
\end{equation}
and the initial conditions 
\begin{align}
\UN[][\fs]  \big|_{u=v} &= \widetilde{W}^{(\Ncs),+} - \big( \UN[][+] + \UN[][+-] + \UN[][-] + \UN[][+\fs]+\UN[][\fs-] \big)\big|_{u=v}, \label{ansatz:eq-initial-data-remainder-e1} \\ 
\VN[][\fs] \big|_{u=v} &= \widetilde{W}^{(\Ncs),-} -  \big( \VN[][-] + \VN[][+-] + \VN[][+] + \VN[][\fs-] + \VN[][+\fs] \big)\big|_{u=v} \label{ansatz:eq-initial-data-remainder-e2}
\end{align}
are satisfied, where $\widetilde{W}^{(\Ncs),+}$ and $\widetilde{W}^{(\Ncs),-}$ are as in Definition \ref{ansatz:def-initial-data}.  
\end{definition}

\begin{remark}\label{ansatz:rem-initial-data-remainder-equations}
\revision{Since it depends on both the initial data from Definition \ref{ansatz:def-initial-data} and the terms in our Ansatz, the initial data in \eqref{ansatz:eq-initial-data-remainder-e1}-\eqref{ansatz:eq-initial-data-remainder-e2} is rather complicated. It will only be estimated at the end of our argument, see 
\eqref{main:eq-null-qb-1}-\eqref{main:eq-null-qb-4primeprime}. Similarly as in Remark \ref{ansatz:rem-initial-data-modulation-equations}, the remainders  $Z^{(\Ncs),\pm}$ appearing in \eqref{ansatz:eq-initial-data-remainder-e1}-\eqref{ansatz:eq-initial-data-remainder-e2} through 
$\widetilde{W}^{(\Ncs),\pm}$ will be needed to prove that our two approximations
\eqref{intro:eq-system-U-V-N-data} and \eqref{intro:eq-system-U-V-Kil} converge to the same limit.}
\end{remark}

The remainder equations \eqref{ansatz:eq-remainder-equation-Us} and \eqref{ansatz:eq-remainder-equation-Vs} are concise and, at least from a conceptual standpoint, rather pleasant. For estimates of $\UN[][\fs]$ and $\VN[][\fs]$, however, the form of the right-hand sides in \eqref{ansatz:eq-remainder-equation-Us} and \eqref{ansatz:eq-remainder-equation-Vs} is not helpful. This is because the Ansatz has not been inserted into the null-form $[ \UN, \VN]_{\leq N}$ and the derivatives of the structured terms have not been grouped with the interactions which are meant to be canceled. In Subsection \ref{section:ansatz-remainder}, we convert \eqref{ansatz:eq-remainder-equation-Us} and \eqref{ansatz:eq-remainder-equation-Vs} into a form that is more helpful for our estimates.

\subsection{A closer look at the modulation equations}\label{section:ansatz-modulation}

We recall from Definition \ref{ansatz:def-modulation-equations} that the modulation equation for $\pSN[K][+][k]$ is given by
\begin{equation}\label{ansatz:eq-Sp-motivation}
\partial_v \pSN[K][+][k] =  \chinull[K] \rho_{\leq N}^2(k) \Big[ \pSN[K][+][k], \LON[K][-] \Big] 
+ \chinull[K]  \rho_{\leq N}^4(k) \Big[ \pSN[K][+][k], \SHHLN[K][v] \Big].
\end{equation}
From an expository perspective, it is convenient to rewrite \eqref{ansatz:eq-Sp-motivation} using a notation that emphasizes that $\pSN[K][+][k]$ is a linear transformation. We therefore write \eqref{ansatz:eq-Sp-motivation} as 
\begin{equation}\label{ansatz:eq-Sp-matrix-form}
\partial_v \pSN[K][+][k] = - \chinull[K] \ad \Big(  \rho_{\leq N}^2(k) \LON[K][-] +  
\rho_{\leq N}^4(k) \SHHLN[K][v] \Big) \circ \pSN[K][+][k], 
\end{equation}
where $\ad$ is the adjoint map from \eqref{prelim:eq-adjoint-map}. Due to Definition \ref{ansatz:def-lo},  \eqref{ansatz:eq-Sp-matrix-form} includes the interaction 
\begin{equation}\label{ansatz:eq-Sp-motivation-para}
- \chinull[K] \rho_{\leq N}^2(k) \ad \big( P_{\leq N}^x \VN[M][-] \big) \circ \pSN[K][-][k]
\end{equation}
for all $M\in \Dyadiclarge$ satisfying\footnote{Since our goal is to control $\pSN[K][-][k]$ uniformly in $K\in \Dyadiclarge$, the constraint $M<K^{1-\delta}$ does not yield much useful information on \eqref{ansatz:eq-Sp-motivation-para}. It will only be used to limit the effect of the Littlewood-Paley operator $P_{\leq N}^x$.} the constraint $M<K^{1-\delta}$. Since $\VN[M][-]$ has $v$-regularity $-1/2-$, we therefore expect that $\pSN[K][+][k]$ has at most $v$-regularity $1/2-$. As a result, the sum of the $v$-regularities of $\pSN[K][+][k]$ and $\VN[M][-]$ is negative, and therefore the high$\times$high-interaction in the $v$-variable in \eqref{ansatz:eq-Sp-motivation-para} cannot be defined directly. This problem is addressed using a para-controlled formulation, which is the subject of the next definition. 

\begin{definition}[Para-controlled modulation equations]\label{ansatz:def-modulation-para}
Let $N,\Nd \in \Dyadiclarge$. Let 
\begin{equation}\label{ansatz:eq-X-Y}
\begin{aligned}
\big( \pXN[K][+][k] \big)_{K\in \Dyadiclarge, k \in \Z_K}, 
\big( \pYN[K][+][k] \big)_{K\in \Dyadiclarge, k \in \Z_K}
&\colon \mathbb{R}_{u,v}^{1+1} \rightarrow \frkg, \\ 
\big( \pXN[M][-][m] \big)_{M\in \Dyadiclarge, m \in \Z_M}, 
\big( \pYN[M][-][m] \big)_{M\in \Dyadiclarge, m \in \Z_M}
&\colon \mathbb{R}_{u,v}^{1+1} \rightarrow \frkg.
\end{aligned}
\end{equation}
Furthermore, define
\begin{equation}\label{ansatz:eq-para-to-original}
\pSN[K][+][k] := \pXN[K][+][k] + \pYN[K][+][k]
\qquad \text{and} \qquad 
\pSN[M][-][m] := \pXN[M][-][m] + \pYN[M][-][m].
\end{equation}
Then, we say that $\pXN[][+][]$, $\pYN[][+][]$, $\pXN[][-][]$, and $\pYN[][-][]$ solve the para-controlled modulation equations if the following integral equations are satisfied: 
\begin{enumerate}[label=(\roman*)]
    \item For all $K\in \Dyadiclarge$ and $k\in \Z_K$, $\pXN[K][+][k]$ satisfies
    \begin{align}
    \pXN[K][+][k] 
    = -  \chinull[K] \rho_{\leq N}^2(k) \, \Int^{v}_{u\rightarrow v} \Big( \ad \big( \LON[K][-] \big) \Big) \Para[v][gg] \pSN[K][+][k].  \label{ansatz:eq-X-p-1}
    \end{align}
    \item For all $K\in \Dyadiclarge$ and $k\in \Z_K$, $\pYN[K][+][k]$ satisfies 
    \begin{align}
     &\pYN[K][+][k] - \SNin[K][+] \notag \\ 
     =& \,   -  \rho_{\leq N}^2(k) \, 
     \bigg(   
     \Int^{v}_{u\rightarrow v} \Big( \chinull[K]  \Big( \ad \big( \LON[K][-] \big) \Para[v][gg] \pSN[K][+][k] \Big) \Big) 
     - \chinull[K] \Big( 
     \Int^{v}_{u\rightarrow v} \Big( \ad \big( \LON[K][-] \big) \Big) \Para[v][gg] \pSN[K][+][k] \Big)
     \bigg) \label{ansatz:eq-Y-p-1} \\
     +&\,  
     \rho_{\leq N}^4(k) \Int^v_{u\rightarrow v} \bigg(  \chinull[K] \bigg( \ad \big( \LON[K][-] \big) \Para[v][sim] \bigg( \chinull[K] \Big( \Int^v_{u\rightarrow v} \big( \ad\big( \LON[K][-] \big) \big) \Para[v][gg] \pSN[K][+][k] \Big) \bigg) \bigg) \bigg)  
     \label{ansatz:eq-Y-p-2} \\ 
     -&\,  \rho_{\leq N}^2(k) \, 
     \Int^v_{u\rightarrow v}  \bigg(  \chinull[K] \Big( \ad \big( \LON[K][-] \big)
     \Para[v][sim] \pYN[K][+][k] \Big)  \bigg)  \label{ansatz:eq-Y-p-3} \\
     -&\, \rho_{\leq N}^2(k) \, 
     \Int^v_{u\rightarrow v}  \bigg(  \chinull[K] \Big( \ad \big( \LON[K][-] \big)
     \Para[v][ll] \pSN[K][+][k] \Big)  \bigg)  \label{ansatz:eq-Y-p-4} \\ 
     -&\,  \rho_{\leq N}^4(k) \, 
     \Int^v_{u\rightarrow v}  \bigg(  \chinull[K] \ad \big( \SHHLN[K][v] \big)
      \pSN[K][+][k] \bigg) . \label{ansatz:eq-Y-p-5} 
    \end{align}
    \item For all $M\in \Dyadiclarge$ and $m\in \Z_M$, $\pXN[M][-][m]$ satisfies 
    \begin{equation}\label{ansatz:eq-X-m-1}
    \pXN[M][-][m] = \chinull[M]\rho_{\leq N}^2(m) \Int^u_{v\rightarrow u} \Big( \ad\big( \LON[M][+] \big) \Big) \Para[u][gg] \pSN[M][-][m]. 
    \end{equation}
    \item For all $M\in \Dyadiclarge$ and $m\in \Z_M$, $\pYN[M][-][m]$ satisfies 
    \begin{align}
    &\pYN[M][-][M] - \SNin[M][+] \notag \\ 
     =& \,   \rho_{\leq N}^2(m) 
     \bigg( 
     \Int^u_{v\rightarrow u} \Big(  \chinull[M]\Big( \ad\big( \LON[M][+] \big) \Para[u][gg] \pSN[M][-][m] \Big) \Big)
     - \chinull[M]\Big( \Int^u_{v\rightarrow u} \Big( \ad \big( \LON[M][+] \big) \Big) \Para[u][gg] \pSN[M][-][m] \Big) \bigg) \label{ansatz:eq-Y-m-1} \\ 
     +& \,  \rho_{\leq N}^4(m) \Int^u_{v\rightarrow u} \bigg( \chinull[M]\bigg( \ad \big( \LON[M][+] \big) \Para[u][sim] \bigg( \chinull[M]\Big( \Int^u_{v\rightarrow u} \big( \ad\big( \LON[M][+] \big) \big) \Para[u][gg] \pSN[M][-][m] \Big) \bigg) \bigg) \bigg) \label{ansatz:eq-Y-m-2}\\ 
     +& \, \rho_{\leq N}^2(m) \Int^u_{v\rightarrow u} \bigg( \chinull[M] \Big( \ad \big( \LON[M][+] \big) \Para[u][sim] \pYN[M][-][m] \Big) \bigg) \label{ansatz:eq-Y-m-3}\\
      +& \, \rho_{\leq N}^2(m) \Int^u_{v\rightarrow u} \bigg( \chinull[M]\Big( \ad \big( \LON[M][+] \big) \Para[u][ll] \pSN[M][-][m] \Big)  \bigg)\label{ansatz:eq-Y-m-4}\\
       +& \, \rho_{\leq N}^4(m) \Int^u_{v\rightarrow u} \bigg( \chinull[M] \ad \big( \SHHLN[M][u] \big)  \pSN[M][-][m]\bigg). \label{ansatz:eq-Y-m-5}
    \end{align}
\end{enumerate}

\end{definition}

The next lemma shows that solving the para-controlled modulation equations also yields a solution of the original modulation equations, which should come as no surprise. 

\begin{lemma}[Para-controlled modulation equations]\label{modulation:lem-para-controlled-modulation-equations}
Let $\pXN[][+][]$, $\pYN[][+][]$, $\pXN[][-][]$, and $\pYN[][-][]$  be as in \eqref{ansatz:eq-X-Y} and let $\pSN[][+][]$ and $\pSN[][-][]$ be as in  \eqref{ansatz:eq-para-to-original}. If $\pXN[][+][]$, $\pYN[][+][]$, $\pXN[][-][]$, and $\pYN[][-][]$ solve the para-controlled modulation equations (from Definition \ref{ansatz:def-modulation-para}), then $\pSN[][+][]$ and $\pSN[][-][]$ solve the modulation equations (from Definition \ref{ansatz:def-modulation-equations}). 
\end{lemma}

\begin{proof}
We only show that $\pSN[K][+][k]$ satisfies the modulation equation \eqref{ansatz:eq-Sp-matrix-form}, since the argument for $\pSN[M][-][m]$ is similar. We first note that
\begin{align*}
\eqref{ansatz:eq-X-p-1} + \eqref{ansatz:eq-Y-p-1} 
&= -   \rho_{\leq N}^2(k) \,   
     \Int^{v}_{u\rightarrow v} \bigg( \chinull[K]  \Big( \ad \big( \LON[K][-] \big) \Para[v][gg] \pSN[K][+][k] \Big) \bigg). 
\end{align*}

By using the integral equation for $\pXN[K][+][k]$, i.e., \eqref{ansatz:eq-X-p-1}, we also obtain that 
\begin{equation*}
\eqref{ansatz:eq-Y-p-2} = 
-  \rho_{\leq N}^2(k) 
     \Int^v_{u\rightarrow v} \bigg( \chinull[K] \Big( \ad \big( \LON[K][-] \big)
     \Para[v][sim] \pXN[K][+][k] \Big) \bigg). 
\end{equation*}
Due to the representation of $\pSN[K][+][k]$ from \eqref{ansatz:eq-para-to-original}, it then follows that
\begin{align*}
\eqref{ansatz:eq-Y-p-2} + \eqref{ansatz:eq-Y-p-3} 
&= -  \rho_{\leq N}^2(k) 
     \Int^v_{u\rightarrow v} \bigg( \chinull[K] \Big( \ad \big( \LON[K][-] \big)
     \Para[v][sim] \big( \pXN[K][+][k] + \pYN[K][+][k] \big) \Big) \bigg) \\
&= -  \rho_{\leq N}^2(k) 
     \Int^v_{u\rightarrow v} \bigg( \chinull[K]  \Big( \ad \big( \LON[K][-] \big)
     \Para[v][sim]  \pSN[K][+][k]  \Big) \bigg).
\end{align*}
By adding the integral equations for $\pXN[K][+][k]$ and $\pYN[K][+][k]$, we therefore obtain that
\begin{align*}
&\hspace{3ex}\pSN[K][+][k] - \SNin[K][+] \\ 
&= \Big( \eqref{ansatz:eq-X-p-1} + \eqref{ansatz:eq-Y-p-1}  \Big) 
+ \Big( \eqref{ansatz:eq-Y-p-2} + \eqref{ansatz:eq-Y-p-3} \Big) 
+ \eqref{ansatz:eq-Y-p-4} + \eqref{ansatz:eq-Y-p-5} \\ 
&= -  \rho_{\leq N}^2(k) 
     \Int^v_{u\rightarrow v} \bigg( \chinull[K]  \Big( 
     \ad \big( \LON[K][-] \big) \Para[v][gg]  \pSN[K][+][k] 
     + \ad \big( \LON[K][-] \big) \Para[v][sim]  \pSN[K][+][k] 
     + \ad \big( \LON[K][-] \big) \Para[v][ll]  \pSN[K][+][k]  \Big)  \bigg) \\
    &-  \rho_{\leq N}^4(k) \, 
     \Int^v_{u\rightarrow v} \Big( \chinull[K]  \ad \big( \SHHLN[K][v] \big)
      \pSN[K][+][k] \Big) \\ 
&= -  \rho_{\leq N}^2(k) 
     \Int^v_{u\rightarrow v} \Big(  \chinull[K] 
     \ad \big(\LON[K][-] \big)   \pSN[K][+][k] \Big) 
     -  \rho_{\leq N}^4(k) \, 
     \Int^v_{u\rightarrow v} \Big( \chinull[K] \ad \big( \SHHLN[K][v] \big)
      \pSN[K][+][k] \Big). 
\end{align*}
Thus, $\pSN[K][+][k]$ satisfies the integral form of the modulation equation \eqref{ansatz:eq-Sp-matrix-form}. 
\end{proof}

The well-posedness of the para-controlled modulation equations is one of the main steps in our argument. After first proving numerous estimates, it will be shown in Section \ref{section:modulation}. 

\subsection{A closer look at the remainder equations}\label{section:ansatz-remainder}

\begin{figure}
\begin{tabular}{
!{\vrule width 1pt}>{\centering\arraybackslash}P{1.5cm}
!{\vrule width 1pt}>{\centering\arraybackslash}P{\bigcolwidth}
!{\vrule width 1pt}>{\centering\arraybackslash}P{\bigcolwidth}
!{\vrule width 1pt}>{\centering\arraybackslash}P{\bigcolwidth}
!{\vrule width 1pt}>{\centering\arraybackslash}P{\bigcolwidth}
!{\vrule width 1pt}>{\centering\arraybackslash}P{\bigcolwidth}
!{\vrule width 1pt}>{\centering\arraybackslash}P{\bigcolwidth}
!{\vrule width 1pt}} 
\noalign{\hrule height 1pt} & & & & & & 
 \\[-5.8ex]
 \, \, \, \, \, \,   $V$ \vspace{-1ex} \newline \, \, \, \,  $U$  \hspace{-0.5ex}\vspace{-0.5ex} 
&  $(-)$ & $(+-)$ & $(+)$   & $(\fs-)$ & $(+\fs)$ & $(\fs)$  
\\[4pt] \noalign{\hrule height 1pt} \rule{0pt}{14pt}
$(+)$ 
& \cellcolor{Yellow!30} 
& \cellcolor{Red!30} $\HHLN[][u]$
& \cellcolor{Yellow!30} 
& \cellcolor{Red!30} $\HHLN[][u]$
& \cellcolor{Red!30} $ \HHLN[][v]$
& \cellcolor{Red!30} $ \HHLN[][v]$
\\[4pt] \noalign{\hrule height 1pt} \rule{0pt}{14pt}
$(+-)$
& \cellcolor{Red!30} $\HHLN[][v]$
& \cellcolor{Green!30} 
& \cellcolor{Green!30} 
& \cellcolor{Green!30} 
& \cellcolor{Red!30} $ \HHLN[][v]$
& \cellcolor{Red!30} $ \HHLN[][v]$
\\[4pt] \noalign{\hrule height 1pt} \rule{0pt}{14pt}
$(-)$ 
& \cellcolor{Yellow!30}
& \cellcolor{Green!30} 
& \cellcolor{Green!30} 
& \cellcolor{Green!30} 
& \cellcolor{Green!30} 
& \cellcolor{Green!30}  
\\[4pt] \noalign{\hrule height 1pt} \rule{0pt}{14pt}
$(+\fs)$ 
& \cellcolor{Red!30} $\HHLN[][v]$
& \cellcolor{Green!30} 
& \cellcolor{Green!30} 
& \cellcolor{Green!30} 
& \cellcolor{Green!30} 
& \cellcolor{Green!30}  
\\[4pt] \noalign{\hrule height 1pt} \rule{0pt}{14pt}
$(\fs-)$ 
& \cellcolor{Red!30} $\HHLN[][u]$
& \cellcolor{Red!30} $\HHLN[][u]$
& \cellcolor{Green!30} 
& \cellcolor{Green!30} 
& \cellcolor{Green!30} 
& \cellcolor{Green!30}  
\\[4pt] \noalign{\hrule height 1pt} \rule{0pt}{14pt}
$(\fs)$ 
& \cellcolor{Red!30} $\HHLN[][u]$
& \cellcolor{Red!30} $\HHLN[][u]$
& \cellcolor{Green!30} 
& \cellcolor{Green!30} 
& \cellcolor{Green!30} 
& \cellcolor{Green!30}  
\\[3pt] \noalign{\hrule height 1pt} 
\end{tabular}
\caption{\small{This figure displays the thirty-six different interactions in the nonlinearity \eqref{ansatz:eq-nonlinearity-motivation}. The green cells correspond to interactions that are entirely perturbative, i.e., which can be controlled in $\Cprod{r-1}{r-1}$. The yellow cells correspond to interactions that, at least in certain frequency regimes, are non-perturbative. These interactions will be absorbed by one of the structured terms in \eqref{ansatz:eq-UN-decomposition} or \eqref{ansatz:eq-VN-decomposition}. The red cells correspond to interactions which, possibly in addition to other non-perturbative interactions, enter into the high$\times$high$\rightarrow$low-matrices (see Definition \ref{ansatz:def-hhl} below). The high$\times$high$\rightarrow$low-matrices are some of the most complicated objects in our analysis and will be emphasized throughout this subsection.}}
\label{figure:ansatz-cases}
\end{figure}

In order to control the nonlinear remainders $\UN[][\fs]$ and $\VN[][\fs]$, we need to obtain more detailed information on the right-hand sides of \eqref{ansatz:eq-remainder-equation-Us} and \eqref{ansatz:eq-remainder-equation-Vs}. To this end, we need to insert our Ansatz from  \eqref{ansatz:eq-UN-decomposition} and \eqref{ansatz:eq-VN-decomposition} into the nonlinearity
\begin{equation}\label{ansatz:eq-nonlinearity-motivation}
 \chinull  \big[ \UN, \VN \big]_{\leq N} -  \chinullsquare \coup \Renorm[\Nc]  \big( \UN + \VN \big). 
\end{equation}
In total, this leads to thirty-six\footnote{In fact, the number of interactions is even larger once we take into account different frequency regimes.} different interactions, which creates organizational challenges. After inserting our Ansatz, we then need to identify the terms which will be canceled by either  
\begin{equation}\label{ansatz:eq-derivative-motivation}
\begin{aligned}
&\partial_v \big( \UN[][+] + \UN[][+-] + \UN[][-] + \UN[][+\fs]+\UN[][\fs-] \big) \\
\text{or} \qquad &
\partial_u \big( \VN[][-] + \VN[][+-] + \VN[][+] + \VN[][\fs-] + \VN[][\fs+] \big), 
\end{aligned}
\end{equation}
i.e., the terms which will be absorbed by the structured terms in our Ansatz. In order to give the reader a sense of direction, we first state our final decomposition and only then give precise definitions of all objects.

\begin{proposition}[Decomposition]\label{ansatz:prop-decomposition}
Let $\UN[][]$ and $\VN[][]$ be as in \eqref{ansatz:eq-UN-decomposition} and \eqref{ansatz:eq-VN-decomposition}, respectively. Furthermore, let 
 $(\sigma_u,\sigma_v)\in \{ (1,0),(0,1) \}$. Then, it holds that
\begin{equation}\label{ansatz:eq-decomposition}
\begin{aligned}
& \chinull \big[ \UN[][], \VN[][] \big]_{\leq N} - \chinullsquare \coup \Renorm[\Ncs] \big( \UN[][] + \VN[][] \big)  \\
-& \, \sigma_u \, \partial_v \big( \UN[][+] + \UN[][+-] + \UN[][-] + \UN[][+\fs]+\UN[][\fs-] \big)  \\
-& \, \sigma_v \, \partial_u \big( \VN[][-] + \VN[][+-] + \VN[][+] + \VN[][\fs-] + \VN[][\fs+] \big)  \\ 
=&\, \chinull \HHLNErr + \chinull \JcbNErr + \sigma_u \SEN[][u][] + \sigma_v \SEN[][v][]
+ \chinull \PIN[] + \chinull \REN .
\end{aligned}
\end{equation}
Here, $\HHLNErr[]$, $\JcbNErr$, $\SEN[][u][]$, $\SEN[][v][]$,  $\PIN[]$, and $\REN$ are as in Definition \ref{ansatz:def-hhlerr}, \mbox{Definition \ref{ansatz:def-jacobi-errors}}, Definition \ref{ansatz:def-structural-error}, Definition \ref{ansatz:def-perturbative}, and Definition \ref{ansatz:def-renormalization-error}.  
\end{proposition}

\begin{remark}
The purpose of the parameters $(\sigma_u,\sigma_v)\in \{ (1,0),(0,1) \}$ is to treat \eqref{ansatz:eq-remainder-equation-Us} and \eqref{ansatz:eq-remainder-equation-Vs} simultaneously. In contrast to the other error terms,  $\SEN[][u][]$ and $\SEN[][v][]$ do not come with a pre-factor of $\chinull $. The reason is that  $\SEN[][u][]$ and $\SEN[][v][]$ contain single factors of $\chinull -\chinull[K]$ or $\chinull -\chinull[M]$, in which we cannot factor out $\chinull$.
\end{remark}

With a lot of work, we will be able to show that the terms on the right-hand side of \eqref{ansatz:eq-decomposition} are all perturbative, i.e., can be controlled in $\Cprod{r-1}{r-1}$. The purpose of Proposition \ref{ansatz:prop-decomposition} is therefore to group the derivatives of the structured terms from \eqref{ansatz:eq-derivative-motivation} with the terms they are meant to cancel and organize the remaining error terms. 
We now define all objects from Proposition \ref{ansatz:prop-decomposition} and then provide its proof.  Throughout the rest of this subsection, we encourage the reader to consult with Figure \ref{figure:ansatz-cases} regarding the organization of the terms. \\

We first focus our attention on high$\times$high$\rightarrow$low-interactions. For expository purposes, we do not immediately state our rigorous definitions and instead first provide a heuristic motivation. To give a specific example, we consider the $(+-)$$\times$$(\fs)$-interaction, which is a dyadic sum of the terms
\begin{equation}\label{ansatz:eq-motivation-hhl-1}
 \Big[ \UN[K,M][+-], \VN[ ][\fs] \Big]_{\leq N} =  \chinull[K,M]  \Big[ \Big[ \UN[K][+] , \IVN[M][-] \Big]_{\leq N}, \VN[ ][\fs] \Big]_{\leq N}. 
\end{equation}
The most problematic interaction in \eqref{ansatz:eq-motivation-hhl-1} contains a high$\times$low-interaction in the $u$-variable and a high$\times$high-interaction in the $v$-variable. It can be written as 
\begin{equation}\label{ansatz:eq-motivation-hhl-2}
  \chinull[K,M] \Big[ \Big[ \UN[K][+] , \IVN[M][-] \Big]_{\leq N} \Para[v][sim] P^u_{<K^{1-\deltap}} \VN[ ][\fs] \Big]_{\leq N}. 
\end{equation}
The term \eqref{ansatz:eq-motivation-hhl-2} is non-perturbative and requires a rather careful analysis. In Definition \ref{ansatz:def-hhl} below, we will introduce high$\times$high$\rightarrow$low-operators such\footnote{To be precise, we never actually define $\HHLN[K,M][\revision{v},(\fs)\times (-)]$. The reason is that the $(+-)$$\times$$(\fs)$-interaction is immediately grouped together with the $(+-)$$\times$$(+\fs)$-interaction, and we then only need $\HHLN[K,M][\revision{v},(\fcs)\times (-)]$ from \eqref{ansatz:eq-hhl-m-2}.} as $\HHLN[K,M][\revision{v},(\fs)\times (-)]$, which capture the essence of \eqref{ansatz:eq-motivation-hhl-2}. Then, we can write
\begin{align}
&\,    \chinull[K,M] \Big[ \Big[ \UN[K][+] , \IVN[M][-] \Big]_{\leq N} \Para[v][sim] P^u_{<K^{1-\deltap}} \VN[ ][\fs] \Big]_{\leq N}  \notag \\ 
=& \,   \chinull[K,M]  \Big[ \Big[ \UN[K][+] , \IVN[M][-] \Big]_{\leq N} \Para[v][sim] P^u_{<K^{1-\deltap}} \VN[ ][\fs] \Big]_{\leq N}
-   P_{\leq N}^x \HHLN[K,M][\revision{v},(\fs)\times (-)] P_{\leq N}^x \UN[K][+] 
\label{ansatz:eq-motivation-hhl-3} \\ 
+& \,   P_{\leq N}^x \HHLN[K,M][\revision{v},(\fs)\times (-)] P_{\leq N}^x \UN[K][+]. 
\label{ansatz:eq-motivation-hhl-4} 
\end{align}
The first term \eqref{ansatz:eq-motivation-hhl-3} is always perturbative
and will be called a high$\times$high$\rightarrow$low-error (see Definition \ref{ansatz:def-hhlerr} below). The second term \eqref{ansatz:eq-motivation-hhl-4} is more delicate. In order to later use certain symmetry properties (see Section \ref{section:jacobi}), \eqref{ansatz:eq-motivation-hhl-4} is never considered individually and is first grouped together with other high$\times$high$\rightarrow$low-interactions. Then, we write 
\begin{equation}\label{ansatz:eq-motivation-hhl-5}
 P_{\leq N}^x \HHLN[K][v] P_{\leq N}^x \UN[K][+] =     \sum_{\ast} \Sumlarge_{M\leq \Nd} 
P_{\leq N}^x \HHLN[K,M][\revision{v},\ast] P_{\leq N}^x \UN[K][+],
\end{equation}
where $\ast$ ranges over four different labels (see Definition \ref{ansatz:def-hhl}). 
The right-hand side of \eqref{ansatz:eq-motivation-hhl-5} will then be split further into 
\begin{equation}\label{ansatz:eq-motivation-hhl-6}
\Big[ (P_{\leq N}^x)^2 \UN[K][+] , \SHHLN[K][+] \Big]_{\leq N} 
\qquad \text{and} \qquad 
\JcbNErr[K][+],
\end{equation}
where $\SHHLN$ is as in Definition \ref{ansatz:def-shhl} above and the Jacobi error $\JcbNErr$ is as in Definition \ref{ansatz:def-jacobi-errors} below. As can be seen from the modulation equations (Definition \ref{ansatz:def-modulation-equations}), the $\SHHLN$-term in \eqref{ansatz:eq-motivation-hhl-6} is absorbed into the modulated linear wave $\UN[K][+]$. The Jacobi error $\JcbNErr[K][+]$ is absorbed into the nonlinear remainder but, as we will see in Section \ref{section:jacobi}, this requires considerable effort. 

After this heuristic discussion, we now make three precise definitions related to high$\times$high$\rightarrow$low-interactions.

\begin{definition}[High$\times$high$\rightarrow$low-operators]\label{ansatz:def-hhl}
Let $K,M,N,\Nd \in \Dyadiclarge$ and let $y\in \R$. Then, we define
\begin{equation}\label{ansatz:eq-hhl-m}
\HHLN[K,M,y][v,(-)\times (-)],
\HHLN[K,M,y][v,(\fcs)\times (-)],
\HHLN[K,M,y][v,(-)\times (\fcs)],
\HHLN[K,M,y][v,\textup{kil}]
\colon  \R^{1+1} \rightarrow \End(\frkg) 
\end{equation} 
by 
\begin{align}
\HHLN[K,M,y][v,(-)\times (-)](u,v) &:= \chinull[K,M]  P^{u,v}_{<K^{\delta}} \bigg( 
\ad \Big( P_{\leq N}^x \VN[M][-] \Big) 
\ad \Big( \Theta^x_y P_{\leq N}^x \IVN[M][-] \Big)  
\bigg), \label{ansatz:eq-hhl-m-1}\allowdisplaybreaks[4]\\
\HHLN[K,M,y][v,(\fcs)\times (-)](u,v) &:= \chinull[K,M] P^{u,v}_{<K^{\delta}} \bigg( 
\ad \Big( P_{\leq N}^x  P^v_{\geq K^{1-\deltap}} \VN[< K^{1-\delta}][\fcs] \Big) 
\ad \Big( \Theta^x_y P_{\leq N}^x \IVN[M][-] \Big)
\bigg), \label{ansatz:eq-hhl-m-2}\allowdisplaybreaks[4]\\
\HHLN[K,M,y][v,(-)\times (\fcs)](u,v) &:=  \chinull[K] P^{u,v}_{<K^{\delta}} \bigg( 
\ad \Big( P_{\leq N}^x \VN[M][-] \Big) 
\ad \Big( \Theta^x_y P_{\leq N}^x  \Int^v_{u\rightarrow v} P^v_{\geq K^{1-\deltap}}   \VN[< K^{1-\delta}][\fcs] \Big) 
\bigg), \label{ansatz:eq-hhl-m-3}\allowdisplaybreaks[4]\\
\HHLN[M,y][v,\textup{kil}](u,v) &:=  \coup \chinull[K]  \Cf^{(\Ncs)}_M(y)  \Kil. \label{ansatz:eq-hhl-m-4}
\end{align}
Similarly, we define 
\begin{equation}\label{ansatz:eq-hhl-p}
\HHLN[M,K,y][u,(+)\times(+)],
\HHLN[M,K,y][u,(\fsc)\times(+)],
\HHLN[M,K,y][u,(+)\times(\fsc)],
\HHLN[K,y][u,\textup{kil}]
\colon  \R^{1+1} \rightarrow \End(\frkg) 
\end{equation} 
by 
\begin{align}
\HHLN[M,K,y][u,(+)\times(+)](u,v)  &:= \chinull[M,K] P^{u,v}_{<M^{\delta}} \bigg( 
\ad\Big( P_{\leq N}^x \UN[K][+] \Big)
\ad\Big( \Theta^x_y P_{\leq N}^x \IUN[K][+] \Big) 
\bigg),\label{ansatz:eq-hhl-p-1}\allowdisplaybreaks[4]\\
\HHLN[M,K,y][u,(\fsc)\times(+)](u,v) &:= \chinull[M,K] P^{u,v}_{<M^{\delta}} \bigg( 
\ad\Big( P_{\leq N}^x  P^u_{\geq M^{1-\deltap}} \UN[< M^{1-\delta}][\fsc] \Big)
\ad\Big( \Theta^x_y P_{\leq N}^x \IUN[K][+] \Big) 
\bigg),\label{ansatz:eq-hhl-p-2}\allowdisplaybreaks[4]\\
\HHLN[M,K,y][u,(+)\times(\fsc)](u,v)  &:=\chinull[M]   P^{u,v}_{<M^{\delta}} \bigg( 
\ad\Big( P_{\leq N}^x \UN[K][+] \Big) 
\ad\Big( \Theta^x_y P_{\leq N}^x \Int^u_{v\rightarrow u}  P^u_{\geq M^{1-\deltap}} \UN[< M^{1-\delta}][\fsc] \Big) 
\bigg),\label{ansatz:eq-hhl-p-3}\allowdisplaybreaks[4]\\
\HHLN[K,y][u,\textup{kil}](u,v)  &:=  \coup \chinull[M] \Cf^{(\Ncs)}_M(y) \Kil . \label{ansatz:eq-hhl-p-4} 
\end{align}
We also define the combined high$\times$high$\rightarrow$low-matrices 
\begin{equation}\label{ansatz:eq-hhl-c1}
\HHLN[K,y][v], \HHLN[M,y][u] \colon \R^{1+1} \rightarrow \End(\frkg)
\end{equation}
by 
\begin{align*}
 \hspace{1ex} \HHLN[K,y][v] 
 &:= 
\Sumlarge_{\substack{M  \leq \Nd\colon \\  M \simeq_\delta K}} 
\Big( \HHLN[K,M,y][v,(-)\times (-)] +  \HHLN[K,M,y][v,\textup{kil}] \Big) \\
&+ 
\Sumlarge_{\substack{ M \leq \Nd\colon\\ M \lesssim_\delta K}} \Big( 
\HHLN[K,M,y][v,(\fcs)\times (-)] + 
\HHLN[K,M,y][v,(-)\times (\fcs)] \Big)
\end{align*}
and 
\begin{align*}
 \hspace{1ex} \HHLN[M,y][u] 
 &:= 
\Sumlarge_{\substack{ K \leq \Nd\colon  \\ K \simeq_\delta M}} \Big( 
 \HHLN[M,K,y][u,(+)\times(+)] + 
 \HHLN[K,y][u,\textup{kil}] \Big) \\
 &+ \Sumlarge_{\substack{ K \leq \Nd\colon \\ K \lesssim_\delta M}} \Big( \HHLN[M,K,y][u,(\fsc)\times(+)] + 
\HHLN[M,K,y][u,(+)\times(\fsc)] \Big).
\end{align*}
Finally, for each $w\in \{u,v\}$ and $L\in \Dyadiclarge$, we define the operators corresponding to one of the matrices in \eqref{ansatz:eq-hhl-c1} by
\begin{equation}\label{ansatz:eq-hhl-c4}
\begin{aligned}
\HHLN[L][w] F (u,v) &= \int_{\R} \dy \, (\widecheck{\rho}_{\leq N} \ast \widecheck{\rho}_{\leq N})(y) \, \HHLN[L,y][w](u,v) \Theta^x_y F(u,v).  
\end{aligned}
\end{equation}
For each of the matrices in \eqref{ansatz:eq-hhl-m} and \eqref{ansatz:eq-hhl-p}, the corresponding operators are defined similarly as in \eqref{ansatz:eq-hhl-c4}. 
\end{definition}

Equipped with the high$\times$high$\rightarrow$low-operators, we now introduce the high$\times$high$\rightarrow$low-errors, which capture errors terms such as \eqref{ansatz:eq-motivation-hhl-3}.

\begin{definition}[High$\times$high$\rightarrow$low error terms]\label{ansatz:def-hhlerr} Let $N,\Nd\in \Dyadiclarge$. Then, we define 
\begin{align}
 \HHLNErr[][(+-)\times (-)]  
:=&\,   \Sumlarge_{\substack{ K_u,K_v,M \leq \Nd\colon\\   K_u \simeq_\delta K_v, \\ 
\max(K_u,K_v)\geq M^{1-\delta}}} 
\bigg( \Big[ \UN[K_u,K_v][+-] \Para[v][sim] \VN[M][-] \Big]_{\leq N} 
\label{ansatz:eq-hhlerr-d1}  \\  
&\hspace{2ex}-  \mathbf{1}\big\{ K_v = M \big\}  P_{\leq N}^x \HHLN[K_u,M][v,(-)\times (-)] P_{\leq N}^x \UN[K_u][+] \bigg), \notag 
\allowdisplaybreaks[4] \\
\HHLNErr[][(+-)\times (\fcs)]  
:=&\,  \Sumlarge_{\substack{K,M \leq \Nd \colon \\  K \simeq_\delta M}} 
\bigg( \Big[ \UN[K,M][+-] \Para[v][sim] \VN[<K^{1-\delta}][\fcs] \Big]_{\leq N} 
\label{ansatz:eq-hhlerr-d2}\\
&\hspace{2ex}-  P_{\leq N}^x \HHLN[K,M][v,(\fcs)\times (-)] P_{\leq N}^x \UN[K][+] \bigg), \notag \allowdisplaybreaks[4] \\ 
\HHLNErr[][(+)\times (\fcs)]  
:=&\,  \Sumlarge_{\substack{K \leq \Nd}}  
\bigg( \Big[ \UN[K][+] \Para[v][sim] P^v_{\geq K^{1-\deltap}} \VN[<K^{1-\delta}][\fcs] \Big]_{\leq N} 
\label{ansatz:eq-hhlerr-d3} \\
&\hspace{2ex}-  \Sumlarge_{\substack{M \leq \Nd \colon \\  M <K^{1-\delta}}} P_{\leq N}^x \HHLN[K,M][v,(\fcs)\times (-)] P_{\leq N}^x \UN[K][+] \bigg), \notag \allowdisplaybreaks[4] \\ 
\HHLNErr[][(+\fs)\times (-)]  
:=&\,  \Sumlarge_{\substack{K,M \leq \Nd \colon\\ M \lesssim_\delta K}}  
\bigg( \Big[ \UN[K][+\fs] \Para[v][sim] \VN[M][-] \Big]_{\leq N} 
\label{ansatz:eq-hhlerr-d4} \\
&\hspace{2ex}-   P_{\leq N}^x \HHLN[K,M][v,(-)\times (\fcs)] P_{\leq N}^x \UN[K][+] \bigg). \notag 
\end{align}
Similarly, we define 
\begin{align}
\HHLNErr[][(+)\times (+-)] 
:=&\, \Sumlarge_{\substack{K,M_u,M_v\leq \Nd\colon \\ M_u \simeq_\delta M_v, \\ \max(M_u,M_v)\geq K^{1-\delta}}} 
\bigg( \Big[ \UN[K][+] \Para[u][sim] \VN[M_u,M_v][+-] \Big]_{\leq N}
\label{ansatz:eq-hhlerr-d5}  \\
&\hspace{2ex} -  \mathbf{1} \big\{ K= M_u \big\} 
P_{\leq N}^x \HHLN[M_v,K][u,(+)\times (+)] P_{\leq N}^x \VN[M_v][-] \bigg), \notag \allowdisplaybreaks[4] \\ 
\HHLNErr[][(\fsc)\times (+-)] 
:=& \, \Sumlarge_{\substack{K,M \leq \Nd \colon \\ K \simeq_\delta M}} 
\bigg( \Big[ \UN[<M^{1-\delta}][\fsc] \Para[u][sim] \VN[K,M][+-] 
\Big]_{\leq N} \label{ansatz:eq-hhlerr-d6} \\
&\hspace{2ex} -  P_{\leq N}^x \HHLN[M,K][u,(\fsc)\times (+)] P_{\leq N}^x \VN[M][-] \bigg), \notag \allowdisplaybreaks[4] \\ 
\HHLNErr[][(\fsc)\times (-)] 
:=& \, \Sumlarge_{\substack{M \leq \Nd }} \bigg( 
\Big[ P^u_{\geq M^{1-\deltap}} \UN[<M^{1-\delta}][\fsc] \Para[u][sim] \VN[M][-] \Big]_{\leq N} \\ 
&\hspace{2ex} -  \Sumlarge_{\substack{K \leq \Nd \colon \\ K<M^{1-\delta}}} P_{\leq N}^x 
\HHLN[M,K][u,(\fsc)\times (+)] P_{\leq N}^x \VN[M][-] \bigg), \allowdisplaybreaks[4] \notag \\ 
\HHLNErr[][(+)\times (\fs-)] 
=:&\,  \Sumlarge_{\substack{K,M \leq \Nd \colon \\ K \lesssim_\delta M}} 
\bigg( \Big[ \UN[K][+] \Para[u][sim] \VN[M][\fs-] \Big]_{\leq N}\\
&\hspace{2ex} -  P_{\leq N}^x \HHLN[M,K][u,(+)\times (\fsc)] P_{\leq N}^x \VN[M][-] \bigg). 
\end{align}
Finally, we define the combined high$\times$high$\rightarrow$low-error by
\begin{align*}
&\, \, \HHLNErr[][] \\ 
:=&\,   \HHLNErr[][(+-)\times (-)] + \HHLNErr[][(+-)\times (\fcs)]  + \HHLNErr[][(+)\times (\fcs)] +   \HHLNErr[][(+\fs)\times (-)]  \\
+&\, \HHLNErr[][(+)\times (+-)] + \HHLNErr[][(\fsc)\times (+-)]
+ \HHLNErr[][(\fsc)\times (-)] + \HHLNErr[][(+)\times (\fs-)] . 
\end{align*}
\end{definition}

\begin{remark}
The definition of $\HHLNErr[][(+)\times (\fcs)]$ is currently hard to motivate, since it is based on the behavior of $\UN[K][+]$ in the $v$-variable.   The behavior of $\UN[K][+]$ in the $v$-variable, however, is hidden in  the modulation operators $(\SN[K][+][k])_{k\in \Z_K}$. Ultimately, the $\HHLN[K,M][\revision{v},(\fcs)\times (-)]$-terms will appear through the  para-controlled structure of  $(\SN[K][+][k])_{k\in \Z_K}$, which was briefly discussed in Subsection \ref{section:ansatz-modulation} and will be discussed in more detail in Section \ref{section:modulation}. 
\end{remark}

In our last definition related to high$\times$high$\rightarrow$low-interactions, we now define the so-called Jacobi errors, which were previously mentioned in the decomposition \eqref{ansatz:eq-motivation-hhl-6} of \eqref{ansatz:eq-motivation-hhl-5}.  We chose the terminology ``Jacobi error" since, as we will see in Section \ref{section:jacobi}, they primarily result from the failure of the Jacobi identity for the frequency-truncated Lie bracket $[\cdot,\cdot]_{\leq N}$.

\begin{definition}[Jacobi errors]\label{ansatz:def-jacobi-errors}
Let $K,M,N,\Nd\in \Dyadiclarge$ satisfy $K,M\leq \Nd$. Then, we define 
\begin{align*}
 \JcbNErr[K][+]&:=   P_{\leq N}^x \HHLN[K][v] P_{\leq N}^x \UN[K][+] 
 -   \big[  (P_{\leq N}^x)^2  \UN[K][+], \SHHLN[K][v] \big]_{\leq N}, \\  
 \JcbNErr[M][-]&:= P_{\leq N}^x \HHLN[M][+] P_{\leq N}^x \VN[M][-] 
 -    \big[ \SHHLN[M][u], (P_{\leq N}^x)^2   \VN[M][-] \big]_{\leq N}. 
\end{align*}
Furthermore, we define
\begin{alignat*}{3}
\JcbNErr[][\pm] := \Sumlarge_{L \leq \Nd} \JcbNErr[L][\pm]  \qquad 
\text{and} \qquad \JcbNErr[][] &:= \JcbNErr[][+]  + \JcbNErr[][-]. 
\end{alignat*}
\end{definition}

We now turn our attention away from high$\times$high$\rightarrow$low-interactions and discuss a different kind of error term. 
As discussed above (see e.g. \eqref{ansatz:eq-motivation-Up} and \eqref{ansatz:eq-motivation-Up-again}), the modulated linear wave $\UN[K][+]$ approximately solves 
\begin{equation}\label{ansatz:eq-motivation-Up-again-again}
    \partial_v \UN[K][+] \simeq \chinull   \Big[ \UN[K][+], \LON[K][-] \Big]_{\leq N} +  \chinull  \Big[ (P_{\leq N}^x)^2   \UN[K][+], \SHHLN[K][v] \Big]_{\leq N}.
\end{equation}
We recall that the reason for only solving \eqref{ansatz:eq-motivation-Up-again-again} approximately rather than exactly is that it allows us to write $\UN[K][+]$ in the simple form from Definition \ref{ansatz:def-modulated-linear}, which will be helpful in many of our estimates. The error term resulting from \eqref{ansatz:eq-motivation-Up-again-again} will be treated perturbatively and will soon be called a structural error (see Definition \ref{ansatz:def-structural-error} below). Similar error terms also result from the definitions of other structured terms in our Ansatz, which satisfy
\begin{flalign}
\partial_v \UN[K,M][+-] &\simeq  \chinull \Big[ \UN[K][+], \VN[M][-] \Big]_{\leq N} , 
\allowdisplaybreaks[4] \label{ansatz:eq-motivation-pm} \\[1ex]
\partial_v \UN[M][-] &\simeq   \chinull \Big[ \LON[M][+],  \VN[M][-] \Big]_{\leq N}
+  \chinull  \Big[  \SHHLN[M][u] , (P_{\leq N}^x)^2  \VN[M][-] \Big]_{\leq N}, \allowdisplaybreaks[4] \label{ansatz:eq-motivation-m} \\[1ex] 
 \partial_v \UN[K][+\fs] &\simeq \chinull  \Big[ \UN[K][+] \Para[v][ll]
    P^v_{\geq K^{1-\deltap}} \VN[<K^{1-\delta}][\fcs] \Big]_{\leq N} ,
  \allowdisplaybreaks[4] \label{ansatz:eq-motivation-ps} \\[1ex]
\partial_v \UN[M][\fs-] &\simeq  \chinull  \Big[ P^u_{\geq M^{1-\deltap}} \UN[<M^{1-\delta}][\fsc] \Para[u][gg] \VN[M][-] \Big]_{\leq N}. \label{ansatz:eq-motivation-sm}
\end{flalign}

\begin{definition}[Structural errors]\label{ansatz:def-structural-error} 
Let $N,\Nd \in \Dyadiclarge$. Then, we define 
\begin{equation*}
\SEN[][u][] := \sum_{\substack{\ast \in \{ +, +-, \\ -,+\fs,\fs-\}}} 
\SEN[][u][\ast],
\end{equation*}
where the individual structural errors are defined as follows: First, each individual structural error is defined as a dyadic sum, i.e.,
\begin{alignat*}{3}
\SEN[][u][+]  :=&\, 
\Sumlarge_{K\leq \Nd} \SEN[K][u][+], \qquad 
\SEN[][u][+-] :=&&\, 
\Sumlarge_{\substack{K,M \leq \Nd\colon \\ K \simeq_\delta M}} \SEN[K,M][u][+-], \\
\SEN[][u][-] :=&\, 
\Sumlarge_{M\leq \Nd} \SEN[M][u][-], \qquad 
\SEN[][u][+\fs] :=&&\, \Sumlarge_{K\leq \Nd} \SEN[K][u][+\fs],\\
\SEN[][u][\fs-] :=&\, \Sumlarge_{M\leq \Nd} \SEN[M][u][\fs-].
\end{alignat*}
Second, the dyadic components are defined as
\begin{flalign}
\SEN[K][u][+] &:= \chinull  \Big[ \UN[K][+], \LON[K][-] \Big]_{\leq N} +   \chinull  \Big[ (P_{\leq N}^x)^2 \UN[K][+], \SHHLN[K][v] \Big]_{\leq N} -    \partial_v \UN[K][+], 
\label{ansatz:eq-serr-p}\allowdisplaybreaks[4] \\[2ex]
\SEN[K,M][u][+-] &:= \chinull  \Big[ \UN[K][+], \VN[M][-] \Big]_{\leq N} -    \partial_v \UN[K,M][+-], 
\label{ansatz:eq-serr-pm}\allowdisplaybreaks[4] \\[2ex]
\SEN[M][u][-] &:=  \chinull  \Big[ \LON[M][+], \VN[M][-] \Big]_{\leq N} 
+  \chinull  \Big[  \SHHLN[M][u] ,  (P_{\leq N}^x)^2  \VN[M][-] \Big]_{\leq N} -   \partial_v \UN[M][-],   \label{ansatz:eq-serr-m}\allowdisplaybreaks[4] \\[2ex] 
\SEN[K][u][+\fs] &:= \chinull  \Big[ \UN[K][+] \Para[v][ll]
    P^v_{\geq K^{1-\deltap}} \VN[<K^{1-\delta}][\fcs] \Big]_{\leq N} -   \partial_v \UN[K][+\fs],
    \label{ansatz:eq-serr-ps}\allowdisplaybreaks[4] \\[2ex]
\SEN[M][u][\fs-] &:= \chinull \Big[ P^u_{\geq M^{1-\deltap}} \UN[<M^{1-\delta}][\fsc] \Para[u][gg] \VN[M][-] \Big]_{\leq N} -   \partial_v \UN[M][\fs-].
    \label{ansatz:eq-serr-fm}
\end{flalign}
We remark that for each $\ast \in \{ +, +-, -, +\fs, -\fs\}$, the structural error $\SEN[][u][\ast]$ corresponds to the structured term $\UN[][\ast]$. In other words, the super-scripts of $\SEN[][u][\ast]$ and $\partial_v \UN[][\ast]$ in 
\eqref{ansatz:eq-serr-p}-\eqref{ansatz:eq-serr-fm} always match. 
The structural error $\SEN[][v][]$ is defined similarly, and the corresponding dyadic components are given by 
\begin{flalign*}
\SEN[M][v][-] &:=  \chinull  \Big[ \LON[M][+], \VN[M][-] \Big]_{\leq N} 
+\chinull  \Big[  \SHHLN[M][u] , (P_{\leq N}^x)^2  \VN[M][-] \Big]_{\leq N} -   \partial_u \VN[M][-], \allowdisplaybreaks[4] \\[2ex] 
\SEN[K,M][v][+-] &:=\chinull  \Big[ \UN[K][+], \VN[M][-] \Big]_{\leq N} -   \partial_u \VN[K,M][+-],   \allowdisplaybreaks[4] \\[2ex]
\SEN[K][v][+] &:= \chinull  \Big[ \UN[K][+], \LON[K][-] \Big]_{\leq N} +\chinull   \Big[ (P_{\leq N}^x)^2  \UN[K][+], \SHHLN[K][v] \Big]_{\leq N} -      \partial_u \VN[K][+], \allowdisplaybreaks[4] \\[2ex]
\SEN[M][v][\fs-] &:= \chinull  \Big[ P^u_{\geq M^{1-\deltap}} \UN[<M^{1-\delta}][\fsc] \Para[u][gg] \VN[M][-] \Big]_{\leq N} -   \partial_u \VN[M][\fs-], \allowdisplaybreaks[4] \\[2ex]
\SEN[K][v][+\fs] &:= \chinull  \Big[ \UN[K][+] \Para[v][ll]
    P^v_{\geq K^{1-\deltap}} \VN[<K^{1-\delta}][\fcs] \Big]_{\leq N} -   \partial_u \VN[K][+\fs]. 
\end{flalign*}
\end{definition}

\begin{remark}
In the structural error $\SEN[K,M][u][+-]$, we include the full $\chinull$-factor, whereas Definition \ref{ansatz:def-modulated-bilinear} only contains a $\chinull[K,M]$-factor. As part of the proof of Lemma \ref{structural:lem-pm}, we will see that terms involving the difference $\chinull-\chinull[K,M]$ are negligible. For the other $\chinull$-factors in \eqref{ansatz:eq-serr-p}-\eqref{ansatz:eq-serr-fm}, the situation is similar.
\end{remark}

We now turn to perturbative interactions, which are the subject of the following definitions. In the first definition, we define the red, yellow, and green index sets, which are used to enumerate the corresponding cells in Figure \ref{figure:ansatz-cases}.

\begin{definition}[Red, yellow, and green index sets]
We define 
\begin{equation*}
\Red, \Yellow, \Green \subseteq \big\{+, -, +-, +\fs, \fs-, \fs\big\}^2
\end{equation*}
as the index sets corresponding to the red, yellow, and green cells in Figure \ref{figure:ansatz-cases}, respectively. More precisely, we define 
\begin{align*}
\Red &:= \Big\{ 
\big( +, +\fs \big), 
\big( + , \fs \big), 
\big( +-, - \big), 
\big( +-, +\fs \big), 
\big( +-, \fs \big), 
\big( +\fs, - \big) \Big\} \\
&\hspace{4ex} \bigcup \Big\{ 
\big( +, +- \big), 
\big( +, \fs- \big), 
\big( \fs-, - \big), 
\big( \fs-, +- \big),
\big( \fs, - \big), 
\big( \fs, +- \big) 
\Big\},  \\ 
\Yellow &:= \Big\{ 
\big( +, - \big),
\big( +, + \big), 
\big( -, - \big),
\Big\}, 
\\ 
\Green &:= \Big\{ 
\big( +-, +- \big), 
\big( +-, + \big), 
\big( +- , \fs - \big),
\big( -, +- \big), 
\big( +\fs, +- \big) \Big\} 
 \\
&\hspace{4ex} \bigcup
\Big\{ 
(\ast_1, \ast_2) \colon \ast_1=-, +\fs, \fs-, \fs \quad \text{and} \quad 
\ast_2= +, \fs-, +\fs, \fs \Big\}. 
\end{align*}
\end{definition}

\begin{remark}
In the next definition, we define the perturbative interactions associated with each of the thirty-six cells of Figure \ref{figure:ansatz-cases}. On first reading, we encourage the reader to not take this definition too seriously. After all, as the name suggests, all interactions included here will ultimately turn out to be perturbative. It may therefore make sense to only consult this definition when reading Section \ref{section:null}, in which all the perturbative interactions are controlled. 
\end{remark}

\begin{definition}[Perturbative interactions]\label{ansatz:def-perturbative} 
Let $N,\Nd \in \Dyadiclarge$. Then, we define
\begin{equation*}
\PIN[] := \PIN[\red]+ \PIN[\yellow] + \PIN[\green] 
\end{equation*}
where $\PIN[\red]$, $\PIN[\yellow]$, and $\PIN[\green]$ capture perturbative interactions associated with the red, yellow, and green cells in Figure \ref{figure:ansatz-cases}, respectively. More precisely, we define 
\begin{align*}
\PIN[\red] &:= \sum_{(\ast_1,\ast_2)\in \Red} \PIN[(\,\ast_1)\times (\,\ast_2)],\\
\PIN[\yellow]&:= \sum_{(\ast_1,\ast_2)\in \Yellow} \PIN[(\,\ast_1)\times (\,\ast_2)], \\ 
\PIN[\green]&:= \sum_{(\ast_1,\ast_2)\in \Green} \PIN[(\,\ast_1)\times (\,\ast_2)].
\end{align*}
 For the  red cells in Figure \ref{figure:ansatz-cases}, the individual perturbative interactions are defined as follows: For the $(+)$$\times$$(+-)$ and $(+-)$$\times$$(-)$-interaction, we define
\begin{align}
\PIN[(+)\times (+-)] 
&:=  \Sumlarge_{\substack{K,M_u,M_v \leq \Nd \colon \\ \max(M_u,M_v) \geq K^{1-\delta}}} 
\Big[ \UN[K][+] \Para[u][nsim] \VN[M_u,M_v][+-] \Big]_{\leq N}, 
\allowdisplaybreaks[4]\\
\PIN[(+-)\times (-)] 
&:= 
\Sumlarge_{\substack{K_u,K_v,M \leq \Nd \colon \\ \max(K_u,K_v) \geq M^{1-\delta} }}
 \Big[ \UN[K_u,K_v][+-] \Para[v][nsim] \VN[M][-] \Big]_{\leq N}. 
 \end{align}
For the $(+)$$\times$$(\fs-)$ and $(+\fs)$$\times$$(-)$-interaction, we define
\begin{align}
\PIN[(+)\times (\fs-)] &:= 
\Sumlarge_{\substack{K,M\leq \Nd \colon \\ M < K^{1-\delta}}} 
\Big[ \UN[K][+], P^u_{\geq K^{1-\deltap}} \VN[M][\fs -] \Big]_{\leq N}
\allowdisplaybreaks[4]\\
&+ \Sumlarge_{\substack{K,M\leq \Nd \colon \\ M \geq K^{1-\delta}}} 
\Big[ \UN[K][+] \Para[u][nsim] \VN[M][\fs -] \Big]_{\leq N}, 
\allowdisplaybreaks[4]\\
\PIN[(+\fs)\times (-)] &:= 
\Sumlarge_{\substack{K,M\leq \Nd\colon \\ K < M^{1-\delta}}}
\Big[ P^v_{\geq M^{1-\deltap}} \UN[K][+\fs], \VN[M][-] \Big]_{\leq N}
\allowdisplaybreaks[4]\\
&+ \Sumlarge_{\substack{K,M\leq \Nd\colon \\ K \geq M^{1-\delta}}}
\Big[  \UN[K][+\fs] \Para[v][nsim] \VN[M][-] \Big]_{\leq N}.
\end{align}
For the $(+)$$\times$$(+\fs)$ and the $(\fs-)$$\times$$(-)$-interaction, we define 
\begin{align}
\PIN[(+)\times (+\fs)] &:= 
\Sumlarge_{\substack{K,L \leq \Nd \colon \\ L \geq K^{1-\delta}}}
\Big( \Big[ \UN[K][+], \VN[L][+\fs] \Big]_{\leq N} 
- \chinull \coup \mathbf{1} \big\{ L = K \big\} \,  \Renorm[\Ncs][K] P^v_{\geq K^{1-\deltap}} 
\VN[<K^{1-\delta}][\fcs] 
\Big)
\allowdisplaybreaks[4]\\
&+ \Sumlarge_{\substack{K,L \leq \Nd \colon \\ L < K^{1-\delta}}}
 \Big[ \UN[K][+] \Para[v][gg] P^v_{>K^{1-\deltap}} \VN[L][+\fs] \Big]_{\leq N},
\allowdisplaybreaks[4]\\
\PIN[(\fs-)\times (-)] &:= 
\Sumlarge_{\substack{L,M \leq \Nd \colon \\ L \geq M^{1-\delta}}}
\Big( \Big[ \UN[L][\fs-], \VN[M][-] \Big]_{\leq N} 
-  \chinull \coup \mathbf{1}\big\{ L = M \big\} \, 
\Renorm[\Ncs][M] P^u_{\geq M^{1-\deltap}} \UN[<M^{1-\delta}][\fsc]
\Big) \\ 
&+ \Sumlarge_{\substack{L,M\leq \Nd \colon \\  L < M^{1-\delta}}} 
\Big[ P^u_{>M^{1-\deltap}} \UN[L][\fs-] \Para[u][ll] \VN[M][-] \Big]_{\leq N}. 
\end{align}
For the $(+)$$\times$$(\fs)$ and $(\fs)$$\times$$(-)$-interaction, we define
\begin{align}
\PIN[(+)\times (\fs)] &:= 
\Sumlarge_{K\leq \Nd} \Big[ \UN[K][+], P^u_{\geq K^{1-\deltap}} \VN[][\fs] \Big]_{\leq N}
\allowdisplaybreaks[4]\\ 
&+ \Sumlarge_{K\leq \Nd} \Big[ \UN[K][+] \Para[v][gg]   P^v_{\geq K^{1-\deltap}} P^u_{< K^{1-\deltap}} \VN[][\fs] \Big]_{\leq N}, 
\allowdisplaybreaks[4]\\ 
\PIN[(\fs)\times (-)] &:= 
\Sumlarge_{\substack{M \leq \Nd}} 
\Big[ P^u_{\geq M^{1-\deltap}} \UN[][\fs], \VN[M][-] \Big]_{\leq N} 
\allowdisplaybreaks[4]\\ 
&+ \Sumlarge_{M\leq \Nd} \Big[ P^u_{\geq M^{1-\deltap}} P^v_{<M^{1-\deltap}} \UN[][\fs]
\Para[u][ll] \VN[M][-] \Big]_{\leq N}. 
\end{align}
For the $(+-)$$\times$$(+\fs)$ and $(\fs-)$$\times$$(+-)$-interaction, we define
\begin{align}
\PIN[(+-)\times (+\fs)] &:= 
\Sumlarge_{\substack{K,M,L \leq \Nd \colon \\ K \simeq_\delta M, \\ L\geq K^{1-\delta}}}
\Big[ \UN[K,M][+-], \VN[L][+\fs] \Big]_{\leq N}  \\ 
&+ \Sumlarge_{\substack{K,M,L \leq \Nd \colon \\ K \simeq_\delta M, \\ L< K^{1-\delta}}}
\Big[ \UN[K,M][+-] \Para[v][nsim] \VN[L][+\fs] \Big]_{\leq N}, 
\allowdisplaybreaks[4]\\
\PIN[(\fs-)\times (+-)] &:= 
\Sumlarge_{\substack{K,L,M \leq \Nd \colon \\ K \simeq_\delta M, \\ 
L \geq M^{1-\delta}}} \Big[ \UN[L][\fs-], \VN[K,M][+-] \Big]_{\leq N} \\
&+  \Sumlarge_{\substack{K,L,M \leq \Nd \colon \\ K \simeq_\delta M, \\ 
L < M^{1-\delta}}} 
\Big[ \UN[L][\fs-] \Para[u][nsim] \VN[K,M][+-] \Big]_{\leq N}. 
\end{align} 
For the $(+-)$$\times$$(\fs)$ and $(\fs)$$\times$$(+-)$-interaction, we define
\begin{align} 
\PIN[(+-)\times (\fs)] &:= 
\Sumlarge_{\substack{K,M\leq \Nd\colon \\ K \simeq_\delta M}} 
\Big[ \UN[K,M][+-],  P^u_{\geq K^{1-\deltap}} \VN[][\fs] \Big]_{\leq N}
\allowdisplaybreaks[4]\\
&+ \Sumlarge_{\substack{K,M\leq \Nd\colon \\ K \simeq_\delta M}} 
\Big[ \UN[K,M][+-] \Para[v][nsim] P^u_{< K^{1-\deltap}} \VN[][\fs] \Big]_{\leq N}, 
\allowdisplaybreaks[4]\\
\PIN[(\fs)\times (+-)] &:= 
\Sumlarge_{\substack{K,M \leq  \Nd \colon \\ K \simeq_\delta M}}
\Big[ P^v_{\geq M^{1-\deltap}}  \UN[][\fs], 
\VN[K,M][+-] \Big]_{\leq N} 
\allowdisplaybreaks[4]\\
 &+ \Sumlarge_{\substack{K,M \leq  \Nd \colon \\ K \simeq_\delta M}}
\Big[ P^v_{< M^{1-\deltap}}  \UN[][\fs] \Para[u][nsim]
\VN[K,M][+-] \Big]_{\leq N}. 
\end{align}
Furthermore, for the yellow cells in Figure \ref{figure:ansatz-cases}, the perturbative interactions are defined as  
\begin{align} 
\PIN[(+)\times (-)] &:= 0, 
\allowdisplaybreaks[4]\\
\PIN[(+)\times (+)] &:= 
\Sumlarge_{\substack{K,M \leq \Nd \colon \\ M \geq K^{1-\delta}}}
 \Big[ \UN[K][+], \VN[M][+] \Big]_{\leq N}, 
\allowdisplaybreaks[4]\\
\PIN[(-)\times (-)] &:= 
\Sumlarge_{\substack{K,M \leq \Nd \colon \\ K \geq M^{1-\delta}}}
\Big[ \UN[K][-], \VN[M][-] \Big]_{\leq N}. 
\end{align}
Finally, for any of the green cells in Figure \ref{figure:ansatz-cases}, i.e., for any $(\ast_1,\ast_2) \in \Green$, all interactions are included in the perturbative interactions. That is, we define
\begin{equation*}
\PIN[(\,\ast_1)\times (\,\ast_2)] :=\Big[ \UN[][\ast_1], \VN[][\ast_2] \Big]_{\leq N}.
\end{equation*}
\end{definition}

Our next definition concerns the renormalization error $\REN$. The renormalization error collects all terms coming from the renormalization $\Renorm[\Ncs](\UN[][]+\VN[][])$ which are not needed to cancel resonances, i.e., which are perturbative. 

\begin{definition}[Renormalization error]\label{ansatz:def-renormalization-error} 
Let $N,\Nd\in \Dyadiclarge$. Then, we define the renormalization errors 
\begin{align}
\RenNErr[][u,+] &:=
- \coup  \Sumlarge_{\substack{K,M \leq \Nd \colon \\ K \simeq_\delta M}} \big( \chinull - \chinull[K] \big) \Renorm[\Ncs][M] \UN[K][+]  
-  \coup \chinull \Sumlarge_{\substack{K,M \leq \Nd \colon \\ K \not \simeq_\delta M}} \Renorm[\Ncs][M] \UN[K][+]  \\
&\hspace{3.25ex} - \coup \chinull \Renorm[\Ncs][<\Nlarge] \UN[][+], 
\allowdisplaybreaks[4] \notag \\ 
\RenNErr[][u,\fs-] &:=  
- \coup \chinull \Sumlarge_{\substack{L,M \leq \Nd \colon \\ L<M^{1-\delta}}}  \hspace{-1ex}
\Renorm[\Ncs][M] P^u_{\leq M^{1-\deltap}} \UN[L][\fs-] 
- \coup \chinull  \Sumlarge_{\substack{L,M \leq \Nd \colon \\ L \geq M^{1-\delta}}} \hspace{-1ex}
\Renorm[\Ncs][M] \UN[L][\fs-]   \\
&\hspace{3.25ex}-\coup \chinull \Renorm[\Ncs][<\Nlarge] \UN[][\fs-], \notag  
\allowdisplaybreaks[4] \\ 
\RenNErr[][u,\fs ] &:=  -  \coup \chinull \Sumlarge_{M\leq \Nd} \Renorm[\Ncs][M] 
\Big( P^v_{\geq M^{1-\deltap}} + P^u_{< M^{1-\deltap}} P^v_{< M^{1-\deltap}} \Big) \UN[][\fs]
-  \coup \chinull \Renorm[\Ncs][<\Nlarge] \UN[][\fs].
\end{align}
For all remaining $\ast \in \{+-,-,+\fs\}$, we define
\begin{equation}\label{ansatz:eq-renormalization-error-u-remaining}
\RenNErr[][u,\ast] := - \coup \chinull  \Renorm[\Ncs][] \UN[][\ast]. 
\end{equation}
Similarly, we define 
\begin{align}
\RenNErr[][v,-] &:= 
- \coup  \Sumlarge_{\substack{K,M \leq \Nd \colon \\ K \simeq_\delta M}} \big( \chinull - \chinull[M]\big)\Renorm[\Ncs][K] \VN[M][-]
-  \coup \chinull \Sumlarge_{\substack{K,M \leq \Nd \colon \\ K \not \simeq_\delta M}} \Renorm[\Ncs][K] \VN[M][-]\\
&\hspace{3.25ex}
- \coup \chinull \Renorm[\Ncs][<\Nlarge] \VN[][-]
\allowdisplaybreaks[4] \notag \\ 
\RenNErr[][v,+\fs] &:= 
-  \coup \chinull \Sumlarge_{\substack{L,K \leq \Nd \colon \\ L<K^{1-\delta}}} \hspace{-1ex}
\Renorm[\Ncs][K] P^v_{\leq K^{1-\deltap}} \VN[L][+\fs] 
-  \coup \chinull \Sumlarge_{\substack{L,K \leq \Nd \colon \\ L \geq K^{1-\delta}}}  \hspace{-1ex}
\Renorm[\Ncs][K] \VN[L][+\fs]\\
&\hspace{3.25ex}- \coup \chinull \Renorm[\Ncs][<\Nlarge] \VN[][+\fs], 
\allowdisplaybreaks[4] \notag \\ 
\RenNErr[][v,\fs] &:= 
 - \coup \chinull  \Sumlarge_{K\leq \Nd} \Renorm[\Ncs][K] 
\Big( P^u_{\geq K^{1-\deltap}} + P^u_{< K^{1-\deltap}} P^v_{< K^{1-\deltap}} \Big) \VN[][\fs]
- \coup \chinull \Renorm[\Ncs][<\Nlarge] \VN[][\fs]
\end{align}
For all remaining $\ast \in \{+-,+,\fs-\}$, we define
\begin{equation*}
\RenNErr[][v,\ast] := -  \coup \chinull \Renorm[\Ncs][] \VN[][\ast]. 
\end{equation*}
Finally, we define
\begin{align*}
\RenNErr[][u] &:= \sum_{\ast} \RenNErr[][u,\ast], \hspace{14ex}
\RenNErr[][v] := \sum_{\ast} \RenNErr[][v,\ast], \\
\RenNErr[][] &:= \RenNErr[][u] + \RenNErr[][v]. 
\end{align*}
\end{definition}

\begin{remark}
As is clear from Definition \ref{ansatz:def-renormalization-error}, many sub-terms of $\Renorm[\Ncs] \UN[][]$ are perturbative and therefore not necessary for well-posedness. The only necessary terms are of the form
\begin{equation*}
\coup \Renorm[\Ncs][M] P^u_{K} P^v_{L} \UN,
\end{equation*}
where $M$ and $K$ are comparable to $N$ and $L$ is comparable to one. This frequency regime is rather restrictive and cannot occur in the structured terms $\UN[][+-]$, $\UN[][-]$, and $\UN[][+\fs]$, which explains our choice in \eqref{ansatz:eq-renormalization-error-u-remaining}. 
\end{remark}

For expository purposes, we now state and prove a lemma which is a first step towards Proposition \ref{ansatz:prop-decomposition}. In this lemma, we only consider red cells in Figure \ref{figure:ansatz-cases} which involve the high$\times$high$\rightarrow$low-matrices $\HHLN[][v]$. 

\begin{lemma}[Red cells]\label{ansatz:lem-red-interactions}
Let $N,\Nd\in \Dyadiclarge$. Then, we have the decomposition 
\begin{align}
&\, 
\Big[ \UN[][+-], \VN[][-] \Big]_{\leq N} + \Big[ \UN[][+-], \VN[][+\fs] \Big]_{\leq N}    + \Big[ \UN[][+-], \VN[][\fs] \Big]_{\leq N} 
  \notag  \\
+&\,  \Big[ \UN[][+], \VN[][+\fs] \Big]_{\leq N} + \Big[ \UN[][+], \VN[][\fs] \Big]_{\leq N}   + \Big[ \UN[][+\fs], \VN[][-] \Big]_{\leq N}
\notag \\
-&\,  \chinull \coup  \Big( \Renorm[\Ncs][] \UN[][+]  + \Renorm[\Ncs] \VN[][+\fs] + \Renorm[\Ncs] \VN[][\fs] \Big) \notag \allowdisplaybreaks[4] \\
=&\, \Sumlarge_{M\leq \Nd}  \Big[  
\Sumlarge_{\substack{K_u,K_v < M^{1-\delta}\colon \\ K_u \simeq_\delta K_v }} 
\UN[K_u,K_v][+-] + \Sumlarge_{K< M^{1-\delta}} P^v_{<M^{1-\deltap}} \UN[K][+\fs]  ,   \VN[M][-] \Big]_{\leq N} 
\label{ansatz:eq-red-interaction-q1} \allowdisplaybreaks[4]  \\
+& \Sumlarge_{K \leq \Nd} \Big[ \UN[K][+] , P^v_{<K^{1-\deltap}} \VN[<K^{1-\delta}][\fcs]  \Big]_{\leq N}
\label{ansatz:eq-red-interaction-q2} \allowdisplaybreaks[4]\\ 
+& \Sumlarge_{K\leq \Nd} \Big[ \UN[K][+] \Para[v][ll] P^v_{\geq K^{1-\deltap}} \VN[<K^{1-\delta}][\fcs] \Big]_{\leq N} 
\label{ansatz:eq-red-interaction-q3} \allowdisplaybreaks[4]\\ 
+& \,  \Sumlarge_{K\leq \Nd}  \Big[  (P_{\leq N}^x)^2  \UN[K][+], \SHHLN[K][v] \Big]_{\leq N} 
\label{ansatz:eq-red-interaction-q4}  \allowdisplaybreaks[4] \\ 
+& \JcbNErr[][+] \label{ansatz:eq-red-interaction-q5} \allowdisplaybreaks[4] \\ 
+& \, \HHLNErr[][(+-)\times (-)] + \HHLNErr[][(+-)\times (\fcs)]
+ \HHLNErr[][(+)\times (\fcs)] + \HHLNErr[][(+\fs)\times (-)]
\label{ansatz:eq-red-interaction-q6}  \allowdisplaybreaks[4] \\ 
+&\, \PIN[(+-)\times (-)] + \PIN[(+-)\times (+\fs)] 
+ \PIN[(+-)\times (\fs)]   \label{ansatz:eq-red-interaction-q7} \allowdisplaybreaks[4] \\ 
+& \PIN[(+)\times (+\fs)]+  \PIN[(+)\times (\fs)] + \PIN[(+\fs)\times (-)] 
\label{ansatz:eq-red-interaction-q8} \allowdisplaybreaks[4] \\ 
+& \, \RenNErr[][u,+] +   \RenNErr[][v,+\fs] +  \RenNErr[][v,\fs].  \label{ansatz:eq-red-interaction-q9}
\end{align}
\end{lemma}

\begin{remark} 
We briefly discuss the nature of the terms in \eqref{ansatz:eq-red-interaction-q1}-\eqref{ansatz:eq-red-interaction-q8}.
The terms in 
\eqref{ansatz:eq-red-interaction-q1}-\eqref{ansatz:eq-red-interaction-q4} are non-perturbative, but will all be absorbed by structured terms in our Ansatz. 
The term \eqref{ansatz:eq-red-interaction-q1} is included in $ [ \LON[M][+], \VN[M][-]]_{\leq N}$ and therefore absorbed by either $\UN[M][-]$ or $\VN[M][-]$. The term \eqref{ansatz:eq-red-interaction-q2} is included in $ [  \UN[K][+], \LON[K][-] ]_{\leq N}$ and therefore absorbed by either $\UN[K][+]$ or $\VN[K][+]$. The term \eqref{ansatz:eq-red-interaction-q3} is absorbed by either $\UN[K][+\fs]$ or $\VN[K][+\fs]$. Finally, due to Definition \ref{ansatz:def-modulation-equations}, the term \eqref{ansatz:eq-red-interaction-q4} is also absorbed by either $\UN[K][+]$ or $\VN[K][+]$. The remaining terms in \eqref{ansatz:eq-red-interaction-q5}-\eqref{ansatz:eq-red-interaction-q9} will later shown to be perturbative.
\end{remark}

\begin{remark}
The proof of Lemma \ref{ansatz:lem-red-interactions} is rather technical and involves several terms. However, in all steps of this proof, there are two simple goals: The first goal is to isolate any non-perturbative interaction, such as the low$\times$high-interaction in \eqref{ansatz:eq-red-interaction-q1} or the high$\times$low-interaction in \eqref{ansatz:eq-red-interaction-q2}. The second goal is to isolate the problematic high$\times$high-interactions, which contribute to \eqref{ansatz:eq-red-interaction-q4}, \eqref{ansatz:eq-red-interaction-q5}, and \eqref{ansatz:eq-red-interaction-q6}. 
\end{remark}

\begin{proof}[Proof of Lemma \ref{ansatz:lem-red-interactions}]
We recall that, due to Definition \ref{ansatz:def-jacobi-errors}, it holds that
\begin{equation}\label{ansatz:eq-red-p1}
\eqref{ansatz:eq-red-interaction-q4} + \eqref{ansatz:eq-red-interaction-q5} 
=  \Sumlarge_{K\leq \Nd} P_{\leq N}^x \HHLN[K][v] P_{\leq N}^x  \UN[K][+]. 
\end{equation}
Due to Definition \ref{ansatz:def-hhl}, the $\HHLN[K][v]$-operator can be further decomposed as 
\begin{equation}
\begin{aligned}
\HHLN[K][v] =&\,  
\Sumlarge_{\substack{M \leq \Nd \colon \\ M \simeq_\delta K}} \HHLN[K,M][v,(-)\times (-)] 
+ \Sumlarge_{\substack{M \leq \Nd \colon \\ M \simeq_\delta K}} \HHLN[K,M][v,(\fcs)\times (-)] 
+ \Sumlarge_{\substack{M \leq \Nd \colon \\ M < K^{1-\delta}}} \HHLN[K,M][v,(\fcs)\times (-)] \\
+&\, \Sumlarge_{\substack{M\leq \Nd \colon \\ M \lesssim_\delta K}}
\HHLN[K,M][v,(-)\times (+\fs)]
+ \Sumlarge_{\substack{M \leq \Nd \colon \\ M \simeq_\delta K}}
\HHLN[M][v,\textup{kil}]. 
\label{ansatz:eq-red-p1'}
\end{aligned}
\end{equation}
We now separately consider the contributions of the $(+-)$$\times$$(-)$, $(+-)$$\times$$(+\fs)$, $(+-)$$\times$$(\fs)$, $(+)$$\times$$(+\fs)$,  $(+)$$\times$$(\fs)$, and $(+\fs)$$\times$$(-)$-interactions from the statement of the proposition. \\ 

\emph{Contribution of the $(+-)$$\times$$(-)$-interaction:} 
We first decompose
\begin{align}
\Big[ \UN[][+-], \VN[][-] \Big]_{\leq N} 
&= \Sumlarge_{\substack{K_u,K_v,M\leq \Nd \colon \\ K_u \simeq_\delta K_v }}
\Big[ \UN[K_u,K_v][+-], \VN[M][-] \Big]_{\leq N} \notag \\ 
&=  \Sumlarge_{\substack{K_u,K_v,M\leq \Nd \colon \\ K_u \simeq_\delta K_v, \\ \max(K_u,K_v)< M^{1-\delta} }}
\Big[ \UN[K_u,K_v][+-], \VN[M][-] \Big]_{\leq N}  \allowdisplaybreaks[4] \label{ansatz:eq-red-p2} \\ 
&+ \Sumlarge_{\substack{K_u,K_v,M\leq \Nd \colon \\ K_u \simeq_\delta K_v, \\ \max(K_u,K_v) \geq M^{1-\delta} }}
\Big[ \UN[K_u,K_v][+-] \Para[v][nsim] \VN[M][-] \Big]_{\leq N} \label{ansatz:eq-red-p3} \allowdisplaybreaks[4] \\
&+  \Sumlarge_{\substack{K_u,K_v,M\leq \Nd \colon \\ K_u \simeq_\delta K_v, \\ \max(K_u,K_v) \geq M^{1-\delta} }}
\Big[ \UN[K_u,K_v][+-] \Para[v][sim] \VN[M][-] \Big]_{\leq N}. \label{ansatz:eq-red-p4} 
\end{align}
The first term \eqref{ansatz:eq-red-p2} is non-perturbative and contained in \eqref{ansatz:eq-red-interaction-q1}. The second term  \eqref{ansatz:eq-red-p3} is included in the perturbative interaction $\PIN[(+-)\times (-)]$.  Using Definition \ref{ansatz:def-hhl} and Definition \ref{ansatz:def-hhlerr}, the third term can be written as  
\begin{align}
\eqref{ansatz:eq-red-p4}
&= \HHLNErr[][(+-)\times (-)] \label{ansatz:eq-red-p4'}\\
&+  \Sumlarge_{\substack{K,M\leq \Nd\colon \\ K \simeq_\delta M}} P_{\leq N}^x \HHLN[K,M][v,(-)\times (-)] P_{\leq N}^x \UN[K][+]. \label{ansatz:eq-red-p4''}
\end{align}
The first term  \eqref{ansatz:eq-red-p4'} is included in \eqref{ansatz:eq-red-interaction-q6}. The second term \eqref{ansatz:eq-red-p4''} corresponds to the contribution of the first summand in \eqref{ansatz:eq-red-p1'}, and is therefore included in \eqref{ansatz:eq-red-p1}. \\

\emph{Contribution of the $(+-)$$\times$$(+\fs)$ and $(+-)$$\times$$(\fs)$-interaction:} 
We first decompose 
\begin{align}
&\hspace{3ex}\Big[ \UN[][+-], \VN[][+\fs] \Big]_{\leq N} 
+ \Big[ \UN[][+-], \VN[\fs] \Big]_{\leq N} \notag \\ 
&= \Sumlarge_{\substack{K,M,L \leq \Nd \colon \\ K \simeq_\delta M}} \Big[ \UN[K,M][+-], \VN[L][+\fs] \Big]_{\leq N} 
+ \Sumlarge_{\substack{K,M \leq \Nd \colon \\ K \simeq_\delta M}} \Big[ \UN[K,M][+-], \VN[][\fs] \Big]_{\leq N}
\notag \\
&= \Sumlarge_{\substack{K,M,L \leq \Nd \colon \\ K \simeq_\delta M, \\ L\geq K^{1-\delta}}} \Big[ \UN[K,M][+-], \VN[L][+\fs] \Big]_{\leq N} 
+ \Sumlarge_{\substack{K,M \leq \Nd \colon \\ K \simeq_\delta M}} \Big[ \UN[K,M][+-], P^u_{\geq K^{1-\deltap}} \VN[][\fs] \Big]_{\leq N} \label{ansatz:eq-red-p7} \allowdisplaybreaks[4] \\ 
&+ \Sumlarge_{\substack{K,M,L \leq \Nd \colon \\ K \simeq_\delta M, \\ L< K^{1-\delta}}} \Big[ \UN[K,M][+-] \Para[v][nsim] \VN[L][+\fs] \Big]_{\leq N} 
+ \Sumlarge_{\substack{K,M \leq \Nd \colon \\ K \simeq_\delta M}} \Big[ \UN[K,M][+-] 
\Para[v][nsim] P^u_{< K^{1-\deltap}} \VN[][\fs] \Big]_{\leq N} \label{ansatz:eq-red-p8} \allowdisplaybreaks[4] \\ 
&+ \Sumlarge_{\substack{K,M,L \leq \Nd \colon \\ K \simeq_\delta M, \\ L< K^{1-\delta}}} \Big[ \UN[K,M][+-] \Para[v][sim] \VN[L][+\fs] \Big]_{\leq N} 
+ \Sumlarge_{\substack{K,M \leq \Nd \colon \\ K \simeq_\delta M}} \Big[ \UN[K,M][+-] 
\Para[v][sim] P^u_{< K^{1-\deltap}} \VN[][\fs] \Big]_{\leq N} \label{ansatz:eq-red-p9}. 
\end{align}
The  terms in both \eqref{ansatz:eq-red-p7} and \eqref{ansatz:eq-red-p8} are perturbative and included  in  
 $\PIN[(+-)\times (+\fs)]$ and  $\PIN[(+-)\times (\fs)]$. From the definition of $\VN[<K^{1-\delta}][\fcs]$ from \eqref{ansatz:eq-VN-asts}, it follows that 
 \begin{align}
 &\Sumlarge_{\substack{K,M,L \leq \Nd \colon \\ K \simeq_\delta M, \\ L< K^{1-\delta}}} \Big[ \UN[K,M][+-] \Para[v][sim] \VN[L][+\fs] \Big]_{\leq N} 
+ \Sumlarge_{\substack{K,M \leq \Nd \colon \\ K \simeq_\delta M}} \Big[ \UN[K,M][+-] 
\Para[v][sim] P^u_{< K^{1-\deltap}} \VN[][\fs] \Big]_{\leq N} \notag \\
=& \, \Sumlarge_{\substack{K,M \leq \Nd \colon \\ K \simeq_\delta M}} \Big[ \UN[K,M][+-] 
\Para[v][sim]  \VN[<K^{1-\delta}][\fcs] \Big]_{\leq N}. \label{ansatz:eq-red-p9'}
 \end{align}
 Using Definition \ref{ansatz:def-hhl} and Definition \ref{ansatz:def-hhlerr}, we then further decompose
 \begin{align}
\eqref{ansatz:eq-red-p9'}
&= \HHLNErr[][(+-)\times (\fcs)]  \label{ansatz:eq-red-p9''} \\ 
&+   \Sumlarge_{\substack{K,M \leq \Nd\colon \\ K \simeq_\delta M}} 
P_{\leq N}^x \HHLN[K,M][v,(\fcs)\times (-)] P_{\leq N}^x \UN[K][+]. \label{ansatz:eq-red-p9'''}
 \end{align}
 The first term  \eqref{ansatz:eq-red-p9''} is included in \eqref{ansatz:eq-red-interaction-q6}. The second term \eqref{ansatz:eq-red-p9'''} corresponds to the contribution of the second in summand in \eqref{ansatz:eq-red-p1'}, and is therefore included in \eqref{ansatz:eq-red-p1}. \\

\emph{Contribution of the $(+)$$\times$$(+\fs)$ and $(+)$$\times$$(\fs)$-interaction:} 
We first decompose
\begin{align}
&\hspace{3ex}\Big[ \UN[][+] , \VN[][+\fs] \Big]_{\leq N}
+ \Big[ \UN[][+] , \VN[][\fs] \Big]_{\leq N} \notag  \\
=& \Sumlarge_{\substack{K,L \leq \Nd}}
\Big[ \UN[K][+], \VN[L][+\fs] \Big]_{\leq N}
+ \Sumlarge_{K\leq \Nd} \Big[ \UN[K][+], \VN[][\fs] \Big]_{\leq N} \notag  \allowdisplaybreaks[4]\\
=& \Sumlarge_{\substack{K,L \leq \Nd \colon \\ L \geq K^{1-\delta}}}
\Big[ \UN[K][+], \VN[L][+\fs] \Big]_{\leq N} \label{ansatz:eq-red-q0} \\
+& \Sumlarge_{K\leq \Nd} \Big[ \UN[K][+],  P^u_{\geq K^{1-\deltap}} \VN[][\fs] \Big]_{\leq N}   \label{ansatz:eq-red-q1} \\ 
+& \Sumlarge_{\substack{K,L \leq \Nd \colon \\ L<  K^{1-\delta}}}
\Big[ \UN[K][+], P^v_{<K^{1-\deltap}} \VN[L][+\fs] \Big]_{\leq N}
+ \Sumlarge_{K\leq \Nd} \Big[ \UN[K][+],  P^u_{< K^{1-\deltap}} P^v_{< K^{1-\deltap}} \VN[][\fs] \Big]_{\leq N} \label{ansatz:eq-red-q2} \\
+& \Sumlarge_{\substack{K,L \leq \Nd \colon \\ L<  K^{1-\delta}}}
\Big[ \UN[K][+], P^v_{\geq K^{1-\deltap}} \VN[L][+\fs] \Big]_{\leq N}
+ \Sumlarge_{K\leq \Nd} \Big[ \UN[K][+],  P^u_{< K^{1-\deltap}} P^v_{\geq  K^{1-\deltap}} \VN[][\fs] \Big]_{\leq N} \label{ansatz:eq-red-q3}. 
\end{align}
The term in \eqref{ansatz:eq-red-q0} requires a renormalization which is included in $ \chinull \coup \Renorm[\Ncs] ( \VN[][+\fs]+\VN[][\fs])$ and will be isolated in \eqref{ansatz:eq-red-p13} below. Once the renormalization has been taken into account, \eqref{ansatz:eq-red-q0} is perturbative and included in $\PIN[(+)\times (+\fs)]$. 
The term in \eqref{ansatz:eq-red-q1} is perturbative and included in  $\PIN[(+)\times (\fs)]$. The terms in \eqref{ansatz:eq-red-q2} are non-perturbative and included in \eqref{ansatz:eq-red-interaction-q2}. Using the definition of $\VN[<K^{1-\delta}][\fcs]$ from \eqref{ansatz:eq-VN-asts}, we first write the last term \eqref{ansatz:eq-red-q3} as 
\begin{equation*}
\eqref{ansatz:eq-red-q3} 
= \Sumlarge_{K \leq \Nd} \Big[ \UN[K][+] , P^v_{\geq K^{1-\deltap}} \VN[<K^{1-\delta}][\fcs] \Big]_{\leq N}. 
\end{equation*}
 Using our para-product operators, we then further decompose
\begin{align}
\Sumlarge_{K\leq \Nd} \Big[ \UN[K][+], P^v_{\geq K^{1-\deltap}} \VN[<K^{1-\delta}][\fcs] \Big]_{\leq N} 
=& \Sumlarge_{K \leq \Nd} \Big[ \UN[K][+] \Para[v][gg] P^v_{\geq K^{1-\deltap}} \VN[<K^{1-\delta}][\fcs] \Big]_{\leq N} \label{ansatz:eq-red-q4}\\
+& \Sumlarge_{K \leq \Nd} \Big[ \UN[K][+] \Para[v][ll] P^v_{\geq K^{1-\deltap}} \VN[<K^{1-\delta}][\fcs] \Big]_{\leq N}\label{ansatz:eq-red-q5}\\
+& \Sumlarge_{K \leq \Nd} \Big[ \UN[K][+] \Para[v][sim] P^v_{\geq K^{1-\deltap}} \VN[<K^{1-\delta}][\fcs] \Big]_{\leq N}.  \label{ansatz:eq-red-q6}
\end{align}
The first term \eqref{ansatz:eq-red-q4} is perturbative. After reinserting the definition of $\VN[<K^{1-\delta}][\fcs]$, its contribution is distributed over $\PIN[(+)\times (+\fs)]$ and $\PIN[(+)\times (\fs)]$. The second term \eqref{ansatz:eq-red-q5} is non-perturbative and included in \eqref{ansatz:eq-red-interaction-q3}. Using Definition \ref{ansatz:def-hhl} and Definition \ref{ansatz:def-hhlerr}, the last term \eqref{ansatz:eq-red-q6}  can be written as 
\begin{align}
\eqref{ansatz:eq-red-q6} 
&=  \HHLNErr[][(+)\times (\fcs)] \label{ansatz:eq-red-q7} \\ 
&+  \Sumlarge_{\substack{K,M \leq \Nd \colon \\ M <K^{1-\delta}}} P_{\leq N}^x 
\HHLN[K,M][v,(\fcs)\times (-)] P_{\leq N}^x \UN[K][+]. 
\label{ansatz:eq-red-q8}
\end{align}
 The first term  \eqref{ansatz:eq-red-q7} is included in \eqref{ansatz:eq-red-interaction-q6}. The second term \eqref{ansatz:eq-red-q8} corresponds to the contribution of the third summand in \eqref{ansatz:eq-red-p1'}, and is therefore included in \eqref{ansatz:eq-red-p1}. \\

\emph{Contribution of the $(+\fs)$$\times$$(-)$-interaction:} 
We first decompose 
\begin{align}
\Big[ \UN[][+\fs], \VN[][-] \Big]_{\leq N}
&= \Sumlarge_{K,M\leq \Nd} \Big[ \UN[K][+\fs], \VN[M][-] \Big]_{\leq N} 
\notag \\ 
&= \Sumlarge_{\substack{K,M\leq \Nd \colon \\ K <M^{1-\delta} }}
\Big[ P^v_{<M^{1-\deltap}} \UN[K][+\fs] , \VN[M][-] \Big]_{\leq N}
\label{ansatz:eq-red-q9} \allowdisplaybreaks[4]\\ 
&+ \Sumlarge_{\substack{K,M\leq \Nd \colon \\ K <M^{1-\delta} }}
\Big[ P^v_{\geq M^{1-\deltap}} \UN[K][+\fs] , \VN[M][-] \Big]_{\leq N}
\label{ansatz:eq-red-q10} \allowdisplaybreaks[4]\\
&+ \Sumlarge_{\substack{K,M\leq \Nd \colon \\ K \geq M^{1-\delta} }}
\Big[  \UN[K][+\fs] , \VN[M][-] \Big]_{\leq N}.
\label{ansatz:eq-red-q11} 
\end{align}
The first term \eqref{ansatz:eq-red-q9} is non-perturbative but included in \eqref{ansatz:eq-red-interaction-q1}. The second term \eqref{ansatz:eq-red-q10} is perturbative and included in $\PIN[(+\fs)\times (-)]$. For the third term \eqref{ansatz:eq-red-q11}, we further decompose 
\begin{align}
\eqref{ansatz:eq-red-q11} &= 
\Sumlarge_{\substack{K,M\leq \Nd \colon \\ K \geq M^{1-\delta} }}
\Big[  \UN[K][+\fs] \Para[v][nsim] \VN[M][-] \Big]_{\leq N} 
\label{ansatz:eq-red-q12} \allowdisplaybreaks[4]\\
&+ 
\Sumlarge_{\substack{K,M\leq \Nd \colon \\ K \geq M^{1-\delta} }}
\Big[  \UN[K][+\fs] \Para[v][sim] \VN[M][-] \Big]_{\leq N}. 
\label{ansatz:eq-red-q13}
\end{align}
 The first summand \eqref{ansatz:eq-red-q12} is also perturbative and included in $\PIN[(+\fs)\times (-)]$. Using Definition \ref{ansatz:def-hhl} and Definition \ref{ansatz:def-hhlerr}, the second summand \eqref{ansatz:eq-red-q13} can be written as 
 \begin{align}
\eqref{ansatz:eq-red-q13} 
&= \HHLNErr[][(+\fs)\times (-)] \label{ansatz:eq-red-q14} \\ 
&+  \Sumlarge_{\substack{K,M\leq \Nd\colon \\ M \lesssim_\delta K}} 
P_{\leq N}^x \HHLN[K,M][(-)\times (\fcs)] P_{\leq N}^x \UN[K][+]
\label{ansatz:eq-red-q15}. 
 \end{align}
  The first term  \eqref{ansatz:eq-red-q4} is included in \eqref{ansatz:eq-red-interaction-q6}. The second term \eqref{ansatz:eq-red-q15} corresponds to the contribution of the fourth summand in \eqref{ansatz:eq-red-p1'}, and is therefore included in \eqref{ansatz:eq-red-p1}.\\

\emph{Contribution of the renormalization terms:} We first consider $-  \chinull \coup \Renorm[\Ncs][] \UN[][+]$, which is decomposed as 
\begin{equation}\label{ansatz:eq-red-p12}
\begin{aligned}
&\, -   \chinull \coup  \Renorm[\Ncs][] \UN[][+] \\
=&\,  
-   \coup  \Sumlarge_{\substack{K,M \leq \Nd \colon  \\ K \simeq_\delta M}} \chinull[K] 
\Renorm[\Ncs][M] \UN[K][+] 
-   \coup  \Sumlarge_{\substack{K,M \leq \Nd \colon  \\ K \simeq_\delta M}} \big( \chinull - \chinull[K] \big)
\Renorm[\Ncs][M] \UN[K][+] \\
&\, -   \chinull \coup  \Sumlarge_{\substack{K,M \leq \Nd \colon \\ K \not \simeq_\delta M}} 
\Renorm[\Ncs][M] \UN[K][+]
-\chinull \coup \Renorm[\Ncs][<\Nlarge] \UN[][+].
\end{aligned}
\end{equation}
Due to Definition \ref{ansatz:def-Killing} and Definition \ref{ansatz:def-hhl}, the first summand in \eqref{ansatz:eq-red-p12} corresponds to the contribution of the fifth summand in \eqref{ansatz:eq-red-p1'} and is therefore included in \eqref{ansatz:eq-red-p1}. The second, third, and fourth summand in \eqref{ansatz:eq-red-p12} agree with the renormalization error $\RenNErr[][u,+]$, which is included in \eqref{ansatz:eq-red-interaction-q9}. \\

We now consider $-  \chinull \coup  \Renorm[\Ncs][] \VN[][+\fs]$ and $- \chinull \coup  \Renorm[\Ncs][] \VN[][\fs]$. To this end, we decompose
\begin{align}
&-  \chinull \coup  \Renorm[\Ncs][] \VN[][+\fs] -  \chinull \coup  \Renorm[\Ncs][] \VN[][+\fs]  \notag \\
=& \, -  \chinull \coup    \Sumlarge_{K\leq \Nd} \Renorm[\Ncs][K] P^v_{\geq K^{1-\deltap}} 
\bigg( \Sumlarge_{\substack{L \leq \Nd\colon \\ L < K^{1-\delta}}} \VN[L][+\fs]
 + P^u_{<K^{1-\deltap}} \VN[][\fs] \bigg) 
 \label{ansatz:eq-red-p13} \allowdisplaybreaks[4] \\
 -&\,  \chinull \coup  \Sumlarge_{\substack{K,L \leq \Nd\colon \\ L < K^{1-\delta}}} \Renorm[\Ncs][K] P^v_{\leq K^{1-\deltap}} \VN[L][+\fs] 
 -  \chinull \coup  \Sumlarge_{\substack{K,L \leq \Nd\colon \\ L \geq K^{1-\delta}}} \Renorm[\Ncs][K] \VN[L][+\fs]
  \label{ansatz:eq-red-p14} \allowdisplaybreaks[4] \\
  -&  \chinull \coup  \Sumlarge_{K\leq \Nd} \Renorm[\Ncs][K] 
  \Big( P^u_{\geq K^{1-\deltap}} 
  + P^u_{< K^{1-\deltap}} P^v_{< K^{1-\deltap}} \Big)
  \VN[][\fs]  \label{ansatz:eq-red-p15}  \\ 
  -& \chinull \Renorm[\Ncs][<\Nlarge] \big( \VN[][+\fs] + \VN[][\fs] \big) \label{ansatz:eq-red-p16}.
\end{align}
After using the definition of $\VN[<K^{1-\delta}][\fcs]$ from \eqref{ansatz:eq-VN-asts}, the first term \eqref{ansatz:eq-red-p13} is exactly the renormalization of \eqref{ansatz:eq-red-q0} and therefore enters into $\PIN[(+)\times (+\fs)]$. The sum of the second, third, and fourth term \eqref{ansatz:eq-red-p14}, \eqref{ansatz:eq-red-p15}, and \eqref{ansatz:eq-red-p16} coincides with the renormalization errors  $\RenNErr[][v,+\fs]+ \RenNErr[][v,\fs]$, and is therefore included in  \eqref{ansatz:eq-red-interaction-q9}. \\

In total, we have now identified all terms in \eqref{ansatz:eq-red-interaction-q1}-\eqref{ansatz:eq-red-interaction-q9}, and therefore completed the proof. 
\end{proof}

Equipped with all of our previous preparations, we can now provide a proof of Proposition \ref{ansatz:prop-decomposition}. 

\begin{proof}[Proof of Proposition \ref{ansatz:prop-decomposition}]
This follows from inserting our Ansatz from \eqref{ansatz:eq-UN-decomposition} and \eqref{ansatz:eq-VN-decomposition} into the nonlinearity, organizing the terms as illustrated in Figure \ref{figure:ansatz-cases}, and using the definitions of this subsection. To be more precise, the contributions to the null-form are organized as follows: First, we decompose
\begin{align*}
\chinull \Big[ \UN[], \VN[] \Big]_{\leq N} 
=& \, \chinull \sum_{(\ast_1,\ast_2)\in \Red} 
\Big[\UN[][\ast_1], \VN[][\ast_2] \Big]_{\leq N}
+ \chinull \sum_{(\ast_1,\ast_2)\in \Yellow} 
\Big[\UN[][\ast_1], \VN[][\ast_2] \Big]_{\leq N} \\ 
+&\, \chinull  \sum_{(\ast_1,\ast_2)\in \Green} 
\Big[\UN[][\ast_1], \VN[][\ast_2] \Big]_{\leq N}.
\end{align*}
We now separately discuss the contributions of the red, yellow, and green cells. \\

\emph{Contributions of red cells:}  
Due to symmetry in the $u$ and $v$-variables, it suffices to treat the red cells from Figure \ref{figure:ansatz-cases} corresponding to 
\begin{equation*}
\Big\{ 
\big( +, +\fs \big), 
\big( + , \fs \big), 
\big( +-, - \big), 
\big( +-, +\fs \big), 
\big( +-, \fs \big), 
\big( +\fs, - \big) \Big\}. 
\end{equation*}
Their combined contribution was already organized in 
Lemma \ref{ansatz:lem-red-interactions}
and, combined with the derivatives of the structured terms, all contributions 
from \eqref{ansatz:eq-red-interaction-q1}-\eqref{ansatz:eq-red-interaction-q9} are as in \eqref{ansatz:eq-decomposition}. \\ 

\emph{Contributions of yellow cells:} 
Due to symmetry in the $u$ and $v$-variables, it suffices to treat the $(+)$$\times$$(-)$, and  $(+)$$\times$$(+)$-interactions. \\

For the $(+)$$\times$$(-)$-interaction, we decompose
\begin{equation}\label{ansatz:eq-decomposition-p0}
\begin{aligned}
\chinull  \Big[ \UN[][+], \VN[][-]\Big]_{\leq N}
&= \chinull  \Sumlarge_{\substack{K,M\leq \Nd\colon \\  K <M^{1-\delta} }}
\Big[ \UN[K][+] , \VN[M][-] \Big]_{\leq N} 
+\chinull  \Sumlarge_{\substack{K,M\leq \Nd\colon \\  M <K^{1-\delta} }} 
\Big[ \UN[K][+] , \VN[M][-] \Big]_{\leq N} \\
&+ \chinull \Sumlarge_{\substack{K,M\leq \Nd\colon \\  K \simeq_\delta M }}
\Big[ \UN[K][+] , \VN[M][-] \Big]_{\leq N} .
\end{aligned}
\end{equation}
Depending on $(\sigma_u,\sigma_v)$, the three summands in \eqref{ansatz:eq-decomposition-p0} are included in either
\begin{equation*}
\SEN[][u][+] + \SEN[][u][-] + \SEN[][u][+-]  
\qquad \text{or} \qquad 
\SEN[][v][+] + \SEN[][v][-] + \SEN[][v][+-]. 
\end{equation*}

For the $(+)$$\times$$(+)$-interaction, we decompose
\begin{equation}\label{ansatz:eq-decomposition-p1}
\chinull  \big[ \UN[][+], \VN[][+] \big]_{\leq N}
= \chinull  \Sumlarge_{\substack{K,M\leq \Nd \colon \\ M < K^{1-\delta} }}
 \big[ \UN[K][+], \VN[M][+] \big]_{\leq N}
  + \chinull  \Sumlarge_{\substack{K,M\leq \Nd \colon \\ M \geq K^{1-\delta} }}
 \big[ \UN[K][+], \VN[M][+] \big]_{\leq N}. 
\end{equation}
Depending on $(\sigma_u,\sigma_v)$, the first summand in \eqref{ansatz:eq-decomposition-p1} enters into either $\SEN[][u][+]$ or $\SEN[][v][+]$. The second summand in \eqref{ansatz:eq-decomposition-p1} is included in the perturbative interaction $\chinull \PIN[(+)\times (+)]$.  \\ 

\emph{Contributions of green cells:} For each green cell, i.e., each $(\ast_1,\ast_2)\in \Green$, the entire contribution
\begin{equation*}
   \Big[\UN[][\ast_1], \VN[][\ast_2] \Big]_{\leq N} 
\end{equation*}
is placed in the perturbative interaction $\chinull  \PIN[(\,\ast_1)\times (\,\ast_2)]$. 
\end{proof}

\begin{remark}\label{ansatz:rem-delta-2}
In Remark \ref{ansatz:rem-delta-1} above, we already mentioned that the choice of the parameters $\delta<\deltap$ impacts the high$\times$high$\rightarrow$low-interactions. If we had instead chosen $\delta>\deltap$, the problematic high$\times$high$\rightarrow$low-interactions would no longer occur in the $(+)$$\times$$(+\fs)$, $(+)$$\times$$(\fs)$, $(\fs-)$$\times$$(-)$, or $(\fs)$$\times$$(-)$-interaction, but would now occur in the $(+)$$\times$$(-)$-interaction. Since it is responsible for the modulated bilinear waves, the $(+)$$\times$$(-)$-interaction is complicated enough as it is, and this motivated our choice $\delta<\deltap$. 
\end{remark}

\subsection{Contraction-mapping hypothesis}\label{section:ansatz-contraction} 
The purpose of many sections of this article is to prove estimates which can be used to close our contraction-mapping arguments for the modulation and remainder equations (see Section \ref{section:modulation} and Section \ref{section:main}). For expository purposes, it is convenient to collect all assumptions on the unknowns which will be used in the contraction-mapping arguments in different hypotheses. In this way, we can then easily reference the different assumptions throughout this article. \\

We first record all probabilistic estimates which will be used in our argument. 
The statements of our probabilistic estimates involve our $\Wuv[\alpha][\beta]$ and $\Wfuv[\alpha][\beta]$-spaces, which are spaces for 
$\Cprod{\alpha}{\beta}$-valued sequences. Since our definitions of $\Wuv[\alpha][\beta]$ and $\Wfuv[\alpha][\beta]$ require additional notation, we postpone them until Definition \ref{chaos:def-wc} and Definition \ref{killing:def-FW}, respectively.

\begin{hypothesis}[Probabilistic hypothesis]\label{hypothesis:probabilistic} 
Let $\Ac \geq 1$ be a parameter, let $R\geq 1$,  and let $(G_{u_0,k}^+)_{k\in \Z}$ and $(G_{v_0,m}^-)_{m\in \Z}$ be the standard, $\frkg_\bC$-valued Gaussian sequences from Subsection \ref{section:ansatz-initial-data}, i.e., from our definition of the initial data. We then assume that the following estimates are satisfied: 
\begin{enumerate}[leftmargin=4ex,label=(\roman*)]
\item \label{ansatz:item-hypothesis-crude} (Crude estimate) For all $u_0,v_0\in \LambdaR$ and $k,m\in \Z$, it holds that
\begin{equation*}
\big\|  G^+_{u_0,k} \big\|_{\frkg}\leq \Ac \langle k \rangle^{\frac{\eta}{2}} \qquad \text{and} \qquad 
\big\| G^-_{v_0,m} \big\|_{\frkg} \leq \Ac \langle m \rangle^{\frac{\eta}{2}}. 
\end{equation*}
\item \label{ansatz:item-hypothesis-sharp} (Sharp frequency-truncations) 
Let $\alpha,\beta \in (-1,1)\backslash \{0\}$, let $K,L,M\in \dyadic$, and let $W^{(\Rscript),\pm}:= W^{(\Rscript,1),\pm}$ be as in \eqref{ansatz:eq-Wpm}. For all $f\in \C_x^s(\R)$ satisfying $\| f \|_{\C_x^s}\leq 1$, it then holds that 
\begin{align*}
\big\| \Psharp_{R;K} W^{(\Rscript),+} \big\|_{\C_x^\alpha} &\leq  \Ac K^{\alpha+\frac{1}{2}+\frac{\eta}{2}}, \\ 
\big\| P_L \Psharp_{R;K} W^{(\Rscript),+} -  \sum_{x_0\in \LambdaRR} \sum_{k\in \Z_K} \rho_L(k) \psiRxK 
G_{x_0,k}^+ e^{ikx} \big\|_{\C_x^\alpha} &\leq  \Ac  K^{\alpha+\frac{1}{2}+\frac{\eta}{2}} \max(K,L)^{-\frac{1}{2}}, \\ 
\big\| \Com \big( \Psharp_{R;L}, (P_{\ll K} f) \big) \Psharp_{R;K} W^{(\Rscript),+} \big\|_{\C_x^\alpha}
&\leq \Ac K^{\alpha+1-s+\frac{\eta}{2}} \max(K,L)^{-\frac{1}{2}}, 
\end{align*}
where $\Com$ is the commutator. Similarly, it holds that 
\begin{align*}
\big\| \Psharp_{R;M} W^{(\Rscript),-} \big\|_{\C_x^\beta} &\leq  \Ac M^{\beta+\frac{1}{2}+\frac{\eta}{2}}, \\ 
\big\| P_L \Psharp_{R;M} W^{(\Rscript),-} -  \sum_{x_0\in \LambdaRR} \sum_{m\in \Z_M} \rho_L(m) \psiRxM 
G_{x_0,m}^- e^{imx} \big\|_{\C_x^\beta} &\leq  \Ac  M^{\beta+\frac{1}{2}+\frac{\eta}{2}} \max(M,L)^{-\frac{1}{2}}, \\ 
\big\| \Com \big( \Psharp_{R;L}, (P_{\ll M} f) \big) \Psharp_{R;M} W^{(\Rscript),-} \big\|_{\C_x^\beta}
&\leq \Ac M^{\beta+1-s+\frac{\eta}{2}} \max(M,L)^{-\frac{1}{2}}.
\end{align*}
\item \label{ansatz:item-hypothesis-linear} (Linear estimates) 
Let $\alpha,\beta \in (-1,1)\backslash \{0\}$ and let $u_0,v_0\in \LambdaR$. Furthermore, 
let $(\Smod[L][+])_{L\in \Dyadiclarge} \subseteq \Wfuv[s][\beta]$, let $(\Smod[L][-])_{L\in \Dyadiclarge} \subseteq \Wfuv[\alpha][s]$, 
and let $K,M,\Nd\in \Dyadiclarge$. Then, it holds that 
\begin{align*}
 \Big\|
\sum_{k\in \Z_K} \rhoND(k) \Smod[K][+][k] G_{u_0,k}^+ \, e^{iku} 
\Big\|_{\Cprod{\alpha}{\beta}} 
&\leq \Ac  K^{\alpha+\frac{1}{2}+\frac{\eta}{2}} \big\| \Smod[K][+] \big\|_{\Wuv[s][\beta]}, \\ 
 \Big\|
\sum_{k\in \Z_K} \rhoND(k) \Smod[K][+][k] G_{u_0,k}^+ \, \frac{e^{iku}}{ik} 
\Big\|_{\Cprod{\alpha}{\beta}}  
&\leq \Ac K^{\alpha-\frac{1}{2}+\frac{\eta}{2}} \big\| \Smod[K][+] \big\|_{\Wuv[s][\beta]} , \\ 
 \Big\| 
\sum_{m\in \Z_M} \rhoND(m) \Smod[M][-][m] G_{v_0,m}^- \, e^{imv} \Big\|_{\Cprod{\alpha}{\beta}} 
&\leq \Ac M^{\beta+\frac{1}{2}+\frac{\eta}{2}} \big\| \Smod[M][-] \big\|_{\Wuv[\alpha][s]} , \\
 \Big\| 
\sum_{m\in \Z_M} \rhoND(m) \Smod[M][-][m] G_{v_0,m}^- \, \frac{e^{imv}}{im} \Big\|_{\Cprod{\alpha}{\beta}} 
&\leq \Ac  M^{\beta-\frac{1}{2}+\frac{\eta}{2}} \big\| \Smod[M][-] \big\|_{\Wuv[\alpha][s]}. 
\end{align*}
\item  \label{ansatz:item-hypothesis-tensor-product} (Main tensor product)  Let $\gamma \in (-s,0)$
and let $u_0,u_1,v_0,v_1\in \LambdaR$. Furthermore, let $(\Smod[L][+])_{L\in \Dyadiclarge}$, $(\Smodtil[L][+])_{L\in \Dyadiclarge}\subseteq \Wfuv[s][s]$, let  $(\Smod[L][-])_{L\in \Dyadiclarge},(\Smodtil[L][-])_{L\in \Dyadiclarge}\subseteq \Wfuv[s][s]$, 
let $K,M,L,\Nd\in \Dyadiclarge$, and let $y,z\in \R$. 
Then, it holds that 
\begin{align*}
&\bigg\| 
\Theta^x_{y} \Big( \sum_{\ell \in \Z_L} \rhoND(\ell) \Smod[L][+][\ell] G_{u_0,\ell}^+ \frac{e^{i\ell u}}{i\ell} \Big) \otimes 
\Theta^x_{z} \Big( \sum_{k\in \Z_K} \rhoND(k) \Smodtil[K][+][k] G_{u_1,k}^+ \, e^{iku} \Big)  \\
&\hspace{2ex} - \mathbf{1}\big\{ K=L\big\} \mathbf{1} \big\{ u_0=u_1 \big\}  
\sum_{k\in \Z_K}  \delta^{ab} \rhoNDsquare(k) \Big( \big( \Theta^x_y \Smod[K][+][k] E_a \big) \otimes \big( \Theta^x_z \Smodtil[K][+][k] E_b \big) \Big) \frac{e^{ik(y-z)}}{(-ik)} 
\bigg\|_{\Cprod{\gamma}{s}}  \\ 
\leq&\,  \Ac^2  \max\big( K, L \big)^{\gamma+\frac{1}{2}+\frac{\eta}{2}} L^{-\frac{1}{2}} \big\| \Smod[L][+] \big\|_{\Wuv[s][s]} \big\| \Smodtil[K][+] \big\|_{\Wuv[s][s]}
\end{align*}
and
\begin{align*}
&\bigg\| 
\Theta^x_{y} \Big( \sum_{\ell \in \Z_L} \rhoND(\ell) \Smod[L][-][\ell] G_{v_0,\ell}^- \frac{e^{i\ell v}}{i\ell} \Big) \otimes 
\Theta^x_{z} \Big( \sum_{m\in \Z_M} \rhoND(m) \Smodtil[M][-][m] G_{v_1,m}^- \, e^{imv} \Big)  \\ 
&\hspace{2ex} - \mathbf{1}\big\{ M=L\big\} \mathbf{1} \big\{ v_0=v_1 \big\}  
\sum_{m\in \Z_M}  \delta^{ab} \rhoNDsquare(m) \Big( \big( \Theta^x_y \Smod[M][-][m] E_a \big) \otimes \big( \Theta^x_z \Smodtil[M][-][m] E_b \big) \Big) \frac{e^{im(y-z)}}{(-im)} \bigg\|_{\Cprod{s}{\gamma}} \notag  \\ 
\leq&\,  \Ac^2  \max\big( L, M \big)^{\gamma+\frac{1}{2}+\frac{\eta}{2}} L^{-\frac{1}{2}} 
\big\| \Smod[L][-] \big\|_{\Wuv[s][s]} \big\| \Smodtil[M][-] \big\|_{\Wuv[s][s]}. 
\end{align*}
\item \label{ansatz:item-hypothesis-trace} (Trace estimate)  Let $\gamma \in (-s,0)$ and let $u_0,v_0\in \LambdaR$. 
Furthermore, let  $(\Smod[L][+])_{L\in \Dyadiclarge}\subseteq \Wfuv[s][s]$, let  $(\Smod[L][-])_{L\in \Dyadiclarge} \subseteq \Wfuv[s][s]$, and let $K,M,\Nd\in \Dyadiclarge$. Finally, let $y,z\in \R$. Then, it holds that 
\begin{align*}
&\, \Big\| 
\Theta^x_{y} \Big( \sum_{k\in \Z_K} \rhoND(k) \Smod[K][+][k](x-t,x+t) G_{u_0,k}^+ \, e^{ik(x-t)} \Big) \\
&\hspace{2ex}\otimes 
\Theta^x_{z} \Big( \sum_{m\in \Z_M} \rhoND(m) \Smod[M][-][m](x-t,x+t) G_{v_0,m}^- \frac{e^{im(x+t)}}{im} \Big) 
\Big\|_{C_t^0 \C_x^{\gamma}} \notag \\
\leq&\, \Ac^2  \max(K,M)^{\gamma+\frac{1}{2}+\frac{\eta}{2}} M^{-\frac{1}{2}} \big\| \Smod[K][+] \big\|_{\Wuv[s][s]}  \big\| \Smod[M][-] \big\|_{\Wuv[s][s]} 
\end{align*}
and 
\begin{align*}
&\, \Big\| 
\Theta^x_{y} \Big( \sum_{k\in \Z_K} \rhoND(k) \Smod[K][+][k](x-t,x+t) G_{u_0,k}^+ \, \frac{e^{ik(x-t)}}{ik} \Big) \\
&\hspace{2ex} \otimes 
\Theta^x_{z} \Big( \sum_{m\in \Z_M} \rhoND(m) \Smod[M][-][m](x-t,x+t) G_{v_0,m}^- e^{im(x+t)} \Big) 
\Big\|_{C_t^0 \C_x^{\gamma}} \notag \\
\leq&\, \Ac^2  \max(K,M)^{\gamma+\frac{1}{2}+\frac{\eta}{2}} K^{-\frac{1}{2}} \big\| \Smod[K][+] \big\|_{\Wuv[s][s]}  \big\| \Smod[M][-] \big\|_{\Wuv[s][s]}.  
\end{align*} 
\item \label{ansatz:item-hypothesis-for-energy-increment} (Energy increment) 
   Let $\gamma \in (-s,0)$ and let $u_0,v_0\in \LambdaR$. 
Furthermore, let  $(\Smod[L][+])_{L\in \Dyadiclarge}\subseteq \Wfuv[s][s]$, let  $(\Smod[L][-])_{L\in \Dyadiclarge} \subseteq \Wfuv[s][s]$, and let $K,M,\Nd\in \Dyadiclarge$. Finally, let $\iota \in \{ \pm 1 \}$  
and let $c,y,z\in \R$. Then, it holds that
\begin{align*}
&\bigg\| \bigg( 
\Theta^x_{y} \Big( \sum_{k\in \Z_K} \rhoND(k) \Smod[K][+][k] G_{u_0,k}^+ \, e^{iku} \Big)  \otimes 
\Theta^x_{z} \Big( \sum_{m\in \Z_M} \rhoND(m) \Smod[M][-][m]G_{v_0,m}^- \frac{e^{imv}}{im} \Big) 
\bigg)(u,\iota u+c)
\bigg\|_{\C_u^{\gamma}}  \notag \\
\leq\, & \Ac^2  \max(K,M)^{\gamma+\frac{1}{2}+\frac{\eta}{2}} M^{-\frac{1}{2}} \big\| \Smod[K][+] \big\|_{\Wuv[s][s]}  \big\| \Smod[M][-] \big\|_{\Wuv[s][s]}
\end{align*}
and 
\begin{align*}
&\bigg\| \bigg( 
\Theta^x_{y} \Big( \sum_{k\in \Z_K} \rhoND(k) \Smod[K][+][k] G_{u_0,k}^+ \, \frac{e^{iku}}{ik} \Big)  \otimes 
\Theta^x_{z} \Big( \sum_{m\in \Z_M} \rhoND(m) \Smod[M][-][m]G_{v_0,m}^- e^{imv} \Big) 
\bigg)(\iota v+c,v )
\bigg\|_{\C_v^{\gamma}} \\
\leq\, & \Ac^2  \max(K,M)^{\gamma+\frac{1}{2}+\frac{\eta}{2}} K^{-\frac{1}{2}} \big\| \Smod[K][+] \big\|_{\Wuv[s][s]}  \big\| \Smod[M][-] \big\|_{\Wuv[s][s]}. 
\end{align*}
\end{enumerate}
\end{hypothesis}

\begin{remark}
Due to the uniformity of Hypothesis \ref{hypothesis:probabilistic} in $u_0,u_1,v_0,v_1\in \LambdaR$, we will later be required to take $\Ac \geq R^{\eta}$, see e.g. Proposition \ref{killing:prop-probabilistic-hypothesis}. The reason for referring to the 
estimate in \ref{ansatz:item-hypothesis-for-energy-increment} of Hypothesis \ref{hypothesis:probabilistic} as ``energy increment" is that it will only be needed in the proof of Lemma \ref{increment:lem-increment-object-I}, which is used to control the energy increment from Definition \ref{increment:def-energy-increment}. 
\end{remark}

The next hypothesis is the so-called pre-modulation hypothesis, which will be used in the contraction-mapping argument for the modulation equations. In other words, it will be used before and while solving the modulation equations from Definition \ref{ansatz:def-modulation-equations}.

\begin{hypothesis}[Pre-modulation hypothesis]\label{hypothesis:pre} 
Let $N,\Nd\in \Dyadiclarge$, let $R\geq 1$, let $\coup \in (0,1)$, and let $\Ac,\Bc\geq 1$. Then, the corresponding pre-modulation hypothesis  consists of the following assumptions:
\begin{enumerate}[label=(\roman*)]
    \item (Probabilistic estimates) The probabilistic estimates from Hypothesis \ref{hypothesis:probabilistic} are satisfied. 
    \item\label{ansatz:item-pre-truncation} (Frequency-truncation parameters) It holds that $\Nd \leq 4 N$. 
    \item (Cut-off function) The function $\chi=\chi(t)$ from Definition \ref{ansatz:def-discretized-wave-maps-periodic-truncated} is an element of $\Cut$, which was introduced in Definition \ref{prelim:def-cut-off}.
    \item (Small coupling) It holds that
    \begin{equation*}
    \Dc  := \hcoup \Ac \Bc \leq c,
    \end{equation*}
    where $c=c(\delta_\ast)>0$ is a sufficiently small constant. 
    \item\label{ansatz:item-pre-modulation} (Pure modulation operators) 
    For all $K,M\in \Dyadiclarge$, it holds that $\pSN[K][+],\pSN[M][-]\in \Wuv[s][s]$  and
    \begin{equation*}
      \big\| \pSN[K][+][k] \big\|_{\Wuv[s][s]},   \big\| \pSN[M][-][m] \big\|_{\Wuv[s][s]} \leq \Bc, 
    \end{equation*}
    where $\Wuv[s][s]$-norm and space are as in Definition \ref{chaos:def-wc} below. 
    \item (Nonlinear remainders) The nonlinear remainders $\UN[][\fs],\VN[][\fs]\colon \R^{1+1}_{u,v}\rightarrow \frkg$ satisfy
    \begin{equation*}
    \big\| \UN[][\fs] \big\|_{\Cprod{r-1}{r}}, 
    \big\| \VN[][\fs] \big\|_{\Cprod{r}{r-1}} 
    \leq \hcoup \Ac \Bc=\Dc. 
    \end{equation*}
\end{enumerate}
\end{hypothesis}

\begin{remark}[Frequency-truncation parameters]\label{ansatz:rem-frequency-truncation}
Despite the assumption $\Nd\leq 4N$ in Hypothesis \ref{hypothesis:pre}.\ref{ansatz:item-pre-truncation}, our results can later be applied to all $\Nd \in \Dyadiclarge$. The reason is that the component of the initial data at frequencies bigger than $4N$ evolves linearly under the finite-dimensional approximation \eqref{ansatz:eq-wave-maps-tilde}.

We also recall that definitions of our modulated objects, such as Definition \ref{ansatz:def-modulated-linear}, Definition \ref{ansatz:def-modulated-bilinear}, Definition \ref{ansatz:def-mixed}, and Definition \ref{ansatz:def-modulated-linear-reversed}, only include dyadic scales $K,M\leq \Nd$. Together with \ref{hypothesis:pre}.\ref{ansatz:item-pre-truncation}, it therefore always holds that $K,M\lesssim N$, which will be used repeatedly.  
\end{remark}

\begin{remark}[Role of the parameters]\label{ansatz:rem-parameters}
The parameter $\coup$ in Hypothesis \ref{hypothesis:pre} is related via \eqref{ansatz:eq-Wpm} to the size of the initial data.  The parameter $\Ac$ determines the probability of the event on which Hypothesis \ref{hypothesis:probabilistic} is satisfied, see Proposition \ref{killing:prop-probabilistic-hypothesis} below. The parameter $\Bc$ serves as an upper bound on the size of the pure modulation parameters. Finally, the parameter $\Dc$ serves as an upper bound on the size of the nonlinear remainders. Throughout this article, we will always take $\Bc \lesssim 1$ and $\Dc \ll 1$. 
\end{remark}

\begin{remark}[\protect{The parameter $\Bc$}]\label{ansatz:rem-parameters-B}
Due to Definition \ref{ansatz:def-pure}, the pre-modulation hypothesis (Hypothesis \ref{hypothesis:pre}) directly implies for all $K,M\in \Dyadiclarge$ that 
\begin{equation}\label{ansatz:eq-parameter-bc}
 \big\| \SN[K][+][k] \big\|_{\Wuv[s][s]},   \big\| \SN[M][-][m] \big\|_{\Wuv[s][s]}
\leq \| \widecheck{\rho} \,\|_{L^1(\R)} \Bc \lesssim \Bc, 
\end{equation}
where $\rho$ is as in \eqref{prelim:eq-rho}. As will be shown in Lemma \ref{modulation:lem-linear}, \eqref{ansatz:eq-parameter-bc} implies that the modulated linear waves  satisfy the estimates 
\begin{equation*}
\big\| \UN[][+] \big\|_{\Cprod{s-1}{s}} \lesssim \hcoup \Ac \Bc \qquad \text{and} \qquad
\big\| \VN[][-] \big\|_{\Cprod{s}{s-1}} \lesssim \hcoup \Ac \Bc. 
\end{equation*}
Thus, the bounds on the modulation operators and nonlinear remainders from Hypothesis \ref{hypothesis:pre} are chosen such that $\UN[][+]$, $\VN[][-]$, $\UN[][\fs]$, and $\VN[][\fs]$ have the same size (with respect to the corresponding norms). 
\end{remark}

Once the modulation equations have been solved, Proposition \ref{modulation:prop-main} below yields additional information on the modulation operators. Since this additional information is needed to solve the remainder equations, we record some of the additional information in a second hypothesis. This second hypothesis is called the post-modulation hypothesis.

\begin{hypothesis}[Post-modulation hypothesis]\label{hypothesis:post} Let $N,\Nd\in \Dyadiclarge$, let $R\geq 1$, let 
$\coup \in (0,1)$, and let $\Ac,\Bc\geq 1$. Then, the \revision{post-modulation} hypothesis consists of the following assumptions:
\begin{enumerate}[label=(\roman*)]
\item The pre-modulation hypothesis (Hypothesis \ref{hypothesis:pre}) is satisfied. 
\item\label{ansatz:item-post-2} The \revision{pure modulation operators $(\pSN[K][+][k]\big)_{K\in \Dyadiclarge,k\in \Z_K}$ and $\big(\pSN[M][-][m]\big)_{M\in \Dyadiclarge,m\in \Z_M}$} satisfy the modulation equations from Definition \ref{ansatz:def-modulation-equations}. 
\item\label{ansatz:item-post-orthogonality} The modulation operators satisfy the orthogonality estimate
\begin{equation*}
\sup_{K\in \Dyadiclarge} K^{100} \Big\| \SN[K][+][k] \big( \SN[K][+][k] \big)^\ast - \Id_\frkg \Big\|_{\Wuv[s][s]}, 
\sup_{M\in \Dyadiclarge} M^{100} \Big\| \SN[M][-][m] \big( \SN[M][-][m] \big)^\ast - \Id_\frkg \Big\|_{\Wuv[s][s]}
\leq C \Bc^2 \Dc, 
\end{equation*}
where $C=C(\delta_\ast)$ is a large constant. 
\end{enumerate}
\end{hypothesis}
We now introduce one last definition, which concerns gains. The gains are purely notational devices which allow us to combine the statements of several different estimates.  

\begin{definition}[Gains]\label{ansatz:def-Gains}
Let $K,K_u,K_v,M,M_u,M_v,N,\Nd\in \Dyadiclarge$. Then, we define
\begin{alignat*}{3}
&\Gain\big( \UN[K][+] \big)= \Gain\big(  \UN[K][-] \big) 
=&& \,\Gain\big(  \UN[K][+\fs] \big) =  \Gain\big( \UN[K][\fs-] \big) = K^{-\eta}, \\
&\Gain\big( \VN[M][-] \big)= \Gain\big( \VN[M][+] \big) 
= && \,\Gain\big( \VN[M][\fs-] \big) = \Gain\big( \VN[M][+\fs]\big) = M^{-\eta}, \\
&\Gain \big( \UN[K_u,K_v][+-] \big) = (K_u K_v)^{-\eta}, \quad
&&\, \Gain \big( \VN[M_u,M_v][+-] \big) = (M_u M_v)^{-\eta}, \\
&\Gain \big( \UN[][\fs] \big)=1, &&\, \Gain \big( \VN[][\fs] \big) =1. 
\end{alignat*}
\end{definition}

\section{Conservative structures}\label{section:conservative}
The goal of this section is to understand the push-forward of white noise under the discretized wave maps equation. 
We therefore first recall from Definition \ref{ansatz:def-discretized-wave-maps-periodic-truncated} that the discretized wave maps equation with cut-off function $\chi \in C^\infty_b(\R)$ and initial time $t_0\in \R$ is given by 
\begin{equation}\label{structure:eq-discretized-wave-maps}
\begin{cases}
\begin{aligned}
\partial_t A^{(\Nscript,\Rscript,\coup)} &= \partial_x B^{(\Nscript,\Rscript,\coup)}, \\ 
\partial_t B^{(\Nscript,\Rscript,\coup)} &= \partial_x A^{(\Nscript,\Rscript,\coup)} - \chi(t)  \Big[ A^{(\Nscript,\Rscript,\coup)},  B^{(\Nscript,\Rscript,\coup)} \Big]_{\leq N} + 2 \chi(t)^2 \coup \Renorm[N] A^{(\Nscript,\Rscript,\coup)}, \\
A^{(\Nscript,\Rscript,\coup)}(t_0) &= W_0^{(\Rscript,\coup)}, \qquad B^{(\Nscript,\Rscript,\coup)}(t_0) = W_1^{(\Rscript,\coup)},  
\end{aligned}
\end{cases}
\end{equation}
where $W_0^{(\Rscript,\coup)},W_1^{(\Rscript,\coup)}\colon \bT_R \rightarrow \frkg$. In this section, the dependence of the solution on the initial time $t_0\in \R$ and cut-off function $\chi$ will sometimes be important. Due to this, we emphasize this dependence by writing
\begin{equation}\label{structure:eq-dependence}
A^{(\Nscript,\Rscript,\coup)}(t,x;t_0,\chi) \qquad \text{and} \qquad  B^{(\Nscript,\Rscript,\coup)}(t,x;t_0,\chi).
\end{equation}
For any fixed $N\in \dyadic$, it is easy to see that \eqref{structure:eq-discretized-wave-maps} is globally well-posed (see e.g. Lemma \ref{structure:lem-global-flow}). Thus, \eqref{structure:eq-discretized-wave-maps} induces a global flow, which we denote by $(\Aflow[\Nscript,\Rscript,\coup], \Bflow[\Nscript,\Rscript,\coup])(t;t_0,\chi)$.

We recall from Definition \ref{prelim:def-Gibbs} that the Gibbs measure $\muR$ is defined as the law of a pair of independent, $2\pi R$-periodic, $\frkg$-valued white noises at temperature $8\coup$. 
In the next definition, we introduce the Gibbs measure $\muNR_{\chi,t}$, which captures the evolution of $\muR$ under \eqref{structure:eq-discretized-wave-maps}. 

\begin{definition}[\protect{Time-evolved Gibbs measure}]\label{structure:def-Gibbs-time-evolved}
For any $N\in \dyadic$, $R\geq 1$, $\chi \in C^\infty_b(\R)$, and $t \in \R$, we define the time-evolved Gibbs measure $\muNR_{\chi,t}$ as the push-forward of $\muR$  under the discretized wave maps equation \eqref{structure:eq-discretized-wave-maps}, i.e., 
\begin{equation*}
\muNR_{\chi,t}= \big(\Aflow[\Nscript,\Rscript,\coup], \Bflow[\Nscript,\Rscript,\coup]\big)(t;0,\chi)_\# \muR. 
\end{equation*}
\end{definition}

\begin{remark}
In Definition \ref{structure:def-Gibbs-time-evolved}, we indexed the Gibbs measures $\muNR_{\chi,t}$ and $\muR$ using $\coup$. In Section \ref{section:conservative-abstract} below, we index abstract Gibbs measures using the inverse temperature $\beta$ (see e.g. Definition \ref{structure:def-Gibbs-form} and Definition \ref{structure:def-measures}). Since $\coup$ and $\beta$ are always related by the identity $\coup = (8\beta)^{-1}$, this should not cause any confusion. 
\end{remark}

Equipped with Definition \ref{structure:def-Gibbs-time-evolved}, we can now state the main proposition of this section. 

\begin{proposition}[Push-forward of the Gibbs measure under the discretized wave maps equation]\label{structure:prop-Gibbs} 
Let $N\in \dyadic$, let $R\geq 1$, let $\coup>0$, and let $\chi\in C^\infty_b(\R)$. Then, it holds for all $\tau \in \R$ that 
\begin{equation}\label{structure:eq-Gibbs-formula}
\begin{aligned}
&\, \frac{\mathrm{d}\muNR_{\chi,\tau} }{ \mathrm{d}\muR} 
\Big( W_0^{(\Rscript,\coup)},W_1^{(\Rscript,\coup)}\Big) \\
=&\,  \exp \bigg( \, 
\frac{1}{4} \int_{0}^\tau \int_{\bT_R} \ds \dx \, \chi(s)^2 
\big \langle \Renorm[N] A^{(\Nscript,\Rscript,\coup)}(s,x;\tau,\chi), 
B^{(\Nscript,\Rscript,\coup)}(s,x;\tau,\chi) \big \rangle
\bigg), 
\end{aligned}
\end{equation}
where $A^{(\Nscript,\Rscript,\coup)}$ and $B^{(\Nscript,\Rscript,\coup)}$ are as in \eqref{structure:eq-dependence}.
\end{proposition}

\begin{remark}
We note that, while the flow $(\Aflow[\Nscript,\Rscript,\coup],\Bflow[\Nscript,\Rscript,\coup])$ is based on the frequency-truncated evolution equation, the initial measure $\muR$ does not contain any frequency-truncation. We also note that, by using the change of variables $s^\prime:= s-\tau$, the integral in \eqref{structure:eq-Gibbs-formula} can be written as 
\begin{equation*}
\frac{1}{4} \int_{-\tau}^0 \int_{\bT_R} \ds^\prime \dx \, \big( \Theta^t_{-\tau} \chi\big)(s)^2 \big \langle \Renorm[N] A^{\Nscript,\Rscript,\coup}(s^\prime,x;0,\Theta^t_{-\tau} \chi), 
B^{(\Nscript,\Rscript,\coup)}(s^\prime,x;0,\Theta^t_{-\tau} \chi) \big \rangle,
\end{equation*}
where $\Theta^t_{-\tau} \chi(s)= \chi(s+\tau)$. 
However, the formulation in \eqref{structure:eq-Gibbs-formula} will be more convenient in the proof of Proposition \ref{main:prop-almost-invariance} below. 
\end{remark}

The proof of Proposition \ref{structure:prop-Gibbs} is postponed until the end of Subsection \ref{section:conservative-wave-maps}. A similar explicit formula as in \eqref{structure:eq-Gibbs-formula} was recently used by Debussche and Tsutsumi \cite{DT21} in the context of quasi-invariant Gaussian measures for nonlinear Schr\"{o}dinger equations, and has since found further applications in \cite{FS22}. Despite the similarities between the formulas in Proposition \ref{structure:prop-Gibbs} and \cite{DT21}, our arguments differ significantly from \cite{DT21}. 
The reason is that the nonlinearity in the discretized wave maps equation is induced by a nonlinear structure on our function space, whereas the nonlinearity in \cite{DT21} is induced by a non-quadratic Hamiltonian. In order to prove Proposition \ref{structure:prop-Gibbs}, we introduce so-called conservative structures, which are inspired by symplectic structures but more flexible. This additional flexibility allows us to discretize the conservative structure of the original wave maps equation (see Subsection \ref{section:conservative-wave-maps}), which ultimately leads to a proof of Proposition \ref{structure:prop-Gibbs}. We remark that our treatment in Subsection \ref{section:conservative-abstract} and Subsection \ref{section:conservative-wave-maps} is more detailed than is required for a proof of Proposition \ref{structure:prop-Gibbs},  but we believe our arguments here to be of independent interest.

\subsection{Abstract conservative structures}\label{section:conservative-abstract}
In this subsection, we introduce abstract conservative structures. In introducing conservative structures, our goal is to create a general framework for proving the invariance of Gibbs measures under (discretized and/or geometric) partial differential equations. Based on this goal, our treatment of conservative structures goes beyond the necessary properties needed in the proof of Proposition \ref{structure:prop-Gibbs} and is sprinkled with additional examples, motivations, and properties. In particular, we include a detailed discussion of the relationship between symplectic structures and conservative structures. \\

 Before stating the precise definition of conservative structures (Definition \ref{structure:def-conservative}), we recall basic definitions and facts from differential geometry. We let $M$ be a smooth, $d$-dimensional manifold. We denote the space of smooth vector fields by $\Gamma(TM)$, the space of smooth covector fields (or one-forms) by $\Gamma(T^\ast M)$, and the space of smooth $k$-forms by $\Omega^k(M)$.  For any smooth vector field $X \in \Gamma(TM)$, smooth $k$-form $\omega \in \Omega^k(M)$, and smooth $l$-form $\eta\in \Omega^l(M)$, we define 
\begin{enumerate}[label=(\roman*)]
    \item the Lie derivative $\Lie_X$, 
    \item the exterior derivative $\mathrm{d}\omega$, 
    \item the interior product $X \intprod \omega$, 
    \item and the wedge product $\omega \wedge \eta$
\end{enumerate}
as in \cite[Section 9 and 14]{Lee13}. We also recall the identities 
\begin{align}
\mathrm{d} \big( \omega \wedge \eta \big)&= \mathrm{d}\omega \wedge \eta + (-1)^k \omega \wedge \mathrm{d}\eta, \label{prelim:eq-diff-geo-d} \\
X \intprod \big( \omega \wedge \eta \big) &= \big( X \intprod \omega \big) \wedge \eta + (-1)^k \omega \wedge \big( X \intprod \eta \big). \label{prelim:eq-diff-geo-intprod} 
\end{align}

We also let $\boldsymbol{\nu}$ be a $d$-form, which we refer to as a volume form. Throughout this subsection, all volume forms will be denoted by bold symbols. The purpose of this is to distinguish volume forms from measures, which will appear in Definition \ref{structure:def-measures} below. While we do not assume this here, most of the volume forms in this section will be nowhere vanishing, i.e., orientation forms. Furthermore, we consider a smooth, time-dependent bundle homomorphism $J\colon \R \times T^\ast M \rightarrow TM$. That is, for each $t\in \R$ and $p\in M$, we have a linear operator
 \begin{equation*}
\Jt_p \colon T_p^\ast M \rightarrow T_p M,
 \end{equation*}
 which depends smoothly on $t$ and $p$. Equivalently, we may view $J$ as a contravariant $2$-tensor field, which can be written in coordinates as
 \begin{equation*}
    \Jt = \Jt^{ij} \frac{\partial}{\partial x^i} \otimes \frac{\partial}{\partial x^j}. 
 \end{equation*}
 For any $t\in \R$ and any smooth differential form $\varphi \in \Gamma(T^\ast M)$, we define a vector field $\Jt\varphi \in \Gamma(TM)$ by $(\Jt\varphi)_p := \Jt_p \varphi_p$. In coordinates, it holds that 
$(\Jt\varphi)^i = \Jt^{ij} \varphi_j$. 

\begin{remark}
The time dependence of the bundle homomorphism $J$ is needed to account for the cut-off function $\chi$ in \eqref{structure:eq-discretized-wave-maps}, but it is not the main aspect of this section.
\end{remark}
 
 \begin{definition}[Conservative structure]\label{structure:def-conservative}
 A triple $(M,\boldsymbol{\nu},J)$, which consists of a smooth  manifold $M$, a volume form $\boldsymbol{\nu}$, and a smooth, time-dependent bundle homomorphism $J\colon \R \times T^\ast M \rightarrow TM $, is called a conservative structure if the following two conditions are satisfied: 
 \begin{enumerate}[label=(\alph*)]
     \item (Skew-symmetry) \label{structure:item-conservative-a} The bundle homomorphism $J$ is skew-symmetric. In other words, for all $t\in \R$, $p\in M$, and  $\varphi,\psi \in T_p^\ast M$, it holds that 
     $\varphi( \Jt_p \psi)= - \psi ( \Jt_p \varphi)$. 
     \item (Volume-preservation) \label{structure:item-conservative-b} For any $t\in \R$ and $H\in C^\infty(M)$, the vector field $\Jt \del H$ preserves the volume form $\boldsymbol{\nu}$, i.e., 
$\Lie_{\Jt \del H} \big( \boldsymbol{\nu} \big) = 0$.
 \end{enumerate}
 \end{definition}
Far from just  a single vector field $X\in \Gamma(TM)$ satisfying $\Lie_X (\boldsymbol{\nu})=0$, the second property in Definition \ref{structure:def-conservative} guarantees the existence of a  family of such vector fields. \\

We now start our discussion of the connection between symplectic and conservative structures.

\begin{definition}[Almost symplectic and symplectic structure]
Let $M$ be a smooth manifold. Then, we make the following definitions: 
\begin{enumerate}[label=(\roman*)] 
\item A smooth $2$-form $\omega \in \Omega^2(M)$ is called an \emph{almost symplectic structure} if it is non-singular. This means that,  for all  $p\in M$ and $u \in T_p M$, there exists a $v \in T_p M$ satisfying 
\begin{equation*}
w_p ( u, v ) \neq 0. 
\end{equation*}
\item A smooth $2$-form $\omega \in \Omega^2(M)$ is called a \emph{symplectic structure} if it is an almost symplectic structure and closed, i.e., satisfies $\mathrm{d}\omega =0$. 
\end{enumerate}
\end{definition}

In the next lemma, we characterize which almost symplectic structures induce conservative structures. 

\begin{lemma}[On almost symplectic and conservative structures]\label{structure:lem-almost-symplectic}  
Let $M$ be a smooth $2n$-dimensional manifold, let $\omega$ be an almost symplectic structure on $M$, and let $\boldsymbol{\nu}:= \omega^n$ be the corresponding volume form. Finally, let $J\colon T^\ast M\rightarrow TM$ be the unique time-independent bundle homomorphism satisfying 
\begin{equation}\label{structure:eq-almost-symplectic-J}
\omega_p \big( J_p \varphi, v \big) = \varphi(v)
\end{equation}
for all $p\in M$, $\varphi \in T_p^\ast M$, and $v\in T_p M$. 
Then, 
$(M,\boldsymbol{\nu},J)$ is a conservative structure if and only if $\mathrm{d}\omega^{n-1}=0$. 
In particular, $(M,\boldsymbol{\nu},J)$ is a conservative structure if $\omega$ is symplectic.  
\end{lemma} 

\begin{remark}
In Lemma \ref{structure:lem-almost-symplectic}, we only consider time-independent almost symplectic structures.
The reason is that Definition \ref{structure:def-conservative} only allows for time-independent volume forms. 
\end{remark}

\begin{remark}
Due to Lemma \ref{structure:lem-almost-symplectic}, the reader may wonder why, instead of working with symplectic structures, we introduced conservative structures. We found that symplectic structures $\omega$ (or almost symplectic structures $\omega$ satisfying $\mathrm{d}\omega^{n-1}=0$) are more difficult to preserve under discretization than conservative structures. As our main argument deals with the discretized wave maps equation from Definition \ref{ansatz:def-discretized-wave-maps}, it is therefore convenient for us to work with conservative structures.  
\end{remark}

\begin{proof}[Proof of Lemma \ref{structure:lem-almost-symplectic}]
Due to the non-singularity and smoothness of $\omega$, it is clear that \eqref{structure:eq-almost-symplectic-J} uniquely determines a smooth bundle homomorphism $J\colon T^\ast M\rightarrow TM$. Furthermore, since $\omega$ is skew-symmetric,  $J$ is also skew-symmetric.  In order to match the notation in the symplectic geometry literature, we now write $\nabla_\omega H := J \mathrm{d} H$. \\

It  remains to prove that the condition 
\begin{equation}\label{structure:eq-almost-symplectic-p1}
\Lie_{\nabla_\omega H} \omega^n = 0 \qquad \textup{for all} \, H \in C^\infty(M)
\end{equation}
 is equivalent to $\mathrm{d}\omega^{n-1}=0$. Using Cartan's magic formula, it holds that 
\begin{equation}\label{structure:eq-almost-symplectic-p2}
\Lie_{\nabla_\omega H} \omega^{n} = \mathrm{d} \big( \nabla_\omega H \intprod \omega^n \big) + \nabla_\omega H \intprod \mathrm{d}\omega^n. 
\end{equation}
Since $\omega^{n}$ is a $2n$-form, it holds that $\mathrm{d}\omega^n =0$. Due to the definition of $\nabla_\omega H$, 
it holds that $\nabla_\omega H \intprod \omega = \mathrm{d}H$. 
Together with \eqref{prelim:eq-diff-geo-intprod}, this yields $\nabla_\omega H \intprod \omega^n = n \, \mathrm{d}H \wedge \omega^{n-1}$. Using \eqref{prelim:eq-diff-geo-d}, we then obtain that
\begin{equation}\label{structure:eq-almost-symplectic-p3}
\begin{aligned}
\mathrm{d} \big( \nabla_\omega H \intprod \omega^n \big) + \nabla_\omega H \intprod \mathrm{d}\omega^n 
&= n \mathrm{d}\big( \mathrm{d} H \wedge \omega^{n-1} \big)  \\ 
&= n \, \mathrm{d}^2 H \wedge \omega^{n-1} - n \, \mathrm{d}H \wedge \mathrm{d}\omega^{n-1}\\
&= - n \, \mathrm{d}H \wedge \mathrm{d}\omega^{n-1}. 
\end{aligned}
\end{equation}
By combining \eqref{structure:eq-almost-symplectic-p2} and \eqref{structure:eq-almost-symplectic-p3}, it follows that \eqref{structure:eq-almost-symplectic-p1} is equivalent to 
\begin{equation}\label{structure:eq-almost-symplectic-p4}
\mathrm{d}H \wedge \mathrm{d}\omega^{n-1} = 0 \qquad \textup{for all} \, H \in C^\infty(M). 
\end{equation}
Since $\mathrm{d}\omega^{n-1}$ is a $(2n-1)$-form, it is easy to see that \eqref{structure:eq-almost-symplectic-p4} is equivalent to $\mathrm{d}\omega^{n-1}=0$. 
\end{proof}

In contrast to symplectic structures, conservative structures can be highly degenerate, which is illustrated by the following example.

\begin{example}[Trivial conservative structure]\label{structure:example-trivial}   
Let $M$ be any  manifold and let $\boldsymbol{\nu}$ be a volume form on $M$. Furthermore, let $J\colon T^\ast M \rightarrow TM$ be the trivial bundle homomorphism, i.e., let  $J_p \varphi =0$
    for all $p\in M$ and $\varphi \in T^\ast_p M$. Then, as is clear from the definitions, $(M,\boldsymbol{\nu},J)$ is a conservative structure.
\end{example}

After discussing the relationship between symplectic and conservative structures, we now discuss basic operations on conservative structures. 

\begin{lemma}[Operations on conservative structures]\label{structure:lem-operations-conservative}
~ 
\begin{enumerate}[label=(\roman*)]
    \item\label{structure:item-products} (Products) Let $(M_1,\boldsymbol{\nu}_1,J_1)$ and $(M_2,\boldsymbol{\nu}_2,J_2)$ be two conservative structures. Then, 
    \begin{equation*}
        (M_1 \times M_2, \boldsymbol{\nu}_1 \wedge \boldsymbol{\nu}_2 , J_1 \oplus J_2) 
    \end{equation*}is a conservative structure. Here, $J_1 \oplus J_2$ is defined by \begin{equation*}
    \big( J_1 \oplus J_2 \big) \big( \varphi_1 \otimes \varphi_2 \big) 
    = J_1 \varphi_1 \otimes \varphi_2 + \varphi_1 \otimes J_2 \varphi_2
    \end{equation*}
    for all $\varphi_1 \in \Gamma(T^\ast M_1)$ and $\varphi_2 \in \Gamma(T^\ast M_2)$. 
    \item\label{structure:item-linear-combination} (Linear combination of bundle homomorphisms) Let $(M,\boldsymbol{\nu},J_1)$ and $(M,\boldsymbol{\nu},J_2)$ be two conservative structures with the same smooth manifold $M$ and the same volume form $\boldsymbol{\nu}$. For any smooth $\theta_1,\theta_2\colon \R \rightarrow \R$, we  define 
    \begin{equation*}
     \Jt:= \theta_1(t) J_1(t)+ \theta_2(t) J_2(t).
    \end{equation*}
    Then, $(M,\boldsymbol{\nu},J)$ is a conservative structure. 
    \item\label{structure:item-multiplication} (Multiplication) Let $(M,\boldsymbol{\nu},J)$ be a conservative structure. Furthermore, let $f\in C^\infty(\R \times M)$ satisfy $\operatorname{Im}(J)f=\{0 \}$, i.e., 
        $(\Jt_p\varphi)f(t)=0$
    for all $t\in \R$, $p \in M$, and $\varphi \in T^\ast_p M$. Then, $(M,\boldsymbol{\nu},fJ)$ is a conservative structure.
\end{enumerate}
\end{lemma}

\begin{proof}
We separately treat the three operations \ref{structure:item-products}, \ref{structure:item-linear-combination}, and \ref{structure:item-multiplication}. \\

\emph{Proof of \ref{structure:item-products}:} 
The skew-symmetry of $J_1 \oplus J_2$ is clear. Thus, it remains to prove for all $t\in \R$ and all $H\in C^\infty(M_1 \times M_2)$ that 
\begin{equation*}
\Lie_{(J_1(t) \oplus J_2(t)) \mathrm{d}H}\big( \boldsymbol{\nu}_1 \wedge \boldsymbol{\nu}_2\big) = 0.
\end{equation*}
Using the linearity of the Lie derivative, we obtain that
\begin{equation}\label{structure:eq-products-p1} 
\begin{aligned}
&\Lie_{(J_1(t) \oplus J_2(t)) \mathrm{d}H}\big( \boldsymbol{\nu}_1 \wedge \boldsymbol{\nu}_2\big)  \\ 
=& \Lie_{J_1(t) \mathrm{d}H} \boldsymbol{\nu}_1 \wedge \boldsymbol{\nu}_2  
+ \boldsymbol{\nu}_1 \wedge \Lie_{J_2(t) \mathrm{d}H} \boldsymbol{\nu}_2. 
\end{aligned}
\end{equation}
The first and second summand in \eqref{structure:eq-products-p1} vanish since $(M_1,\boldsymbol{\nu}_1,J_1)$ and $(M_2,\boldsymbol{\nu}_2,J_2)$ are conservative structures, respectively. \\

\emph{Proof of \ref{structure:item-linear-combination}:} Clearly, $J$ satisfies the skew-symmetry condition. Due to the linearity of the Lie derivatives in the vector field, $J$ also preserves the volume form $\boldsymbol{\nu}$. \\

\emph{Proof of \ref{structure:item-multiplication}:} Clearly, $fJ$ satisfies the skew-symmetry condition. Thus, it remains to verify that 
\begin{equation*}
\Lie_{f(t) \Jt \del H}\big( \boldsymbol{\nu} \big)=0
\end{equation*}
for all $t\in \R$ and $H\in C^\infty(M)$. We first prove that 
\begin{equation}\label{structure:eq-operations-p1}
\Lie_{f(t) \Jt \del H}\big( \boldsymbol{\nu} \big)= f(t) \Lie_{\Jt \del H}\big( \boldsymbol{\nu} \big) + \big( \Jt \del H\big)(f(t))) \, \boldsymbol{\nu}. 
\end{equation}
Using Cartan's magic formula, we obtain\footnote{The first term on the right-hand side vanishes since $\boldsymbol{\nu}$ is a top-order form, but this is not relevant here.} 
\begin{equation}\label{structure:eq-operations-p2}
\Lie_{f(t) \Jt \mathrm{d}H} \boldsymbol{\nu} = \big( f(t) \Jt \mathrm{d}H \big) \intprod \mathrm{d} \boldsymbol{\nu} + \mathrm{d} \big( ( f(t) \Jt \mathrm{d} H )\intprod \boldsymbol{\nu} \big). 
\end{equation}
Using the product rule for the exterior derivative \eqref{prelim:eq-diff-geo-d}, it holds that
\begin{equation}\label{structure:eq-operations-p3}
\mathrm{d} \big( ( f(t) \Jt \mathrm{d} H )\intprod \boldsymbol{\nu} \big) 
= f(t) \mathrm{d} \big( (  \Jt \mathrm{d} H )\intprod \boldsymbol{\nu} \big) 
+ \mathrm{d}f(t) \wedge \big( (\Jt \mathrm{d}H) \intprod \boldsymbol{\nu} \big). 
\end{equation}
After combining \eqref{structure:eq-operations-p2} and \eqref{structure:eq-operations-p3}, the proof of \eqref{structure:eq-operations-p1} is reduced to the identity 
\begin{equation}\label{structure:eq-operations-p4}
\mathrm{d}f(t) \wedge \big( (\Jt \mathrm{d}H) \intprod \boldsymbol{\nu} \big)
= (\Jt \mathrm{d}H)(f(t))\, \boldsymbol{\nu}. 
\end{equation}
Since $\mathrm{d}f(t) \wedge \boldsymbol{\nu}$ is a $(d+1)$-form, it holds that $\mathrm{d}f(t) \wedge \boldsymbol{\nu}=0$. From \eqref{prelim:eq-diff-geo-intprod}, it follows that 
\begin{equation*}
0 =  \Jt \mathrm{d}H \intprod \big( \mathrm{d} f(t) \wedge \boldsymbol{\nu} \big) = - \mathrm{d}f(t) \wedge \big( (\Jt\mathrm{d}H) \intprod \boldsymbol{\nu} \big) + \big( (\Jt\mathrm{d}H) \intprod \mathrm{d}f(t)) \big) \boldsymbol{\nu}. 
\end{equation*}
Since $(\Jt\mathrm{d}H) \intprod \mathrm{d}f(t)= (\Jt\mathrm{d}H)(f(t))$, this proves \eqref{structure:eq-operations-p4} and hence completes the proof of \eqref{structure:eq-operations-p1}. It now remains to show that the right-hand side of \eqref{structure:eq-operations-p1} vanishes. The first summand in \eqref{structure:eq-operations-p1} vanishes since $J$ preserves the volume form $\boldsymbol{\nu}$. The second summand in \eqref{structure:eq-operations-p1} vanishes due to our assumption on $f$. 
\end{proof}

We now turn to dynamics on conservative structures. To this end, let $H\in C^\infty(M)$ be a smooth function. Then, we obtain a time-dependent vector field $\Jt \del H$, which generates a smooth local flow. In order to obtain a smooth global flow, it is convenient to require a growth condition, which is stated in the following definition.

 \begin{definition}[Coercive functions]\label{structure:def-coercive}
 A Hamiltonian $H\in C^\infty(M)$ is called coercive if $H$ is nonnegative and if, for all $\theta \in \R$, the sublevel sets
$\big\{ p \in M \colon H(p) \leq \theta \big\}$
 are compact. 
 \end{definition}
 
 Given a conservative structure $(M,\boldsymbol{\nu},J)$ and a coercive Hamiltonian $H$, we arrive at the differential equation 
 \begin{equation}\label{structure:eq-JgradH}
\dot{\gamma}(t) = (J(t) \del H)_{\gamma(t)}, \qquad \gamma(t_0)=\gamma_0 \in M. 
 \end{equation}
 Unfortunately, \eqref{structure:eq-JgradH} does not capture the discretized wave maps equation (Definition \ref{ansatz:def-discretized-wave-maps}). In order to fit the discretized wave maps equation into our framework, we will introduce an additional smooth, time-dependent vector field $K$, which will be as in the following definition.
 
\begin{definition}[Conservative perturbation]\label{structure:def-perturbation}
A smooth, time-dependent vector field $K$ is called a conservative perturbation with respect to a conservative structure $(M,\boldsymbol{\nu},J)$ and a coercive Hamiltonian $H$ if the following two assumptions are satisfied: 
\begin{enumerate}[label=(\alph*)]
\item (Growth) For all $t\in \R$ and $p\in M$, it holds that 
$| \Kt_p(H)| \lesssim \big(1+H(p) \big)$.
\item (Volume-preservation) $K$ preserves the volume form $\boldsymbol{\nu}$, i.e., $\Lie_{\Kt} ( \boldsymbol{\nu} )=0$ for all $t\in \R$. 
\end{enumerate}
\end{definition}

Equipped with a conservative structure $(M,\boldsymbol{\nu},J)$, a coercive Hamiltonian $H$, and a conservative perturbation $K$, we can now state the main evolution equation of this subsection. It is given by 
\begin{equation}\label{structure:eq-JgradH-K}
\dot{\gamma}(t) = (\Jt \del H)_{\gamma(t)} + \Kt_{\gamma(t)}, \qquad \gamma(t_0)=\gamma_0 \in M. 
\end{equation}
 
 In the following lemma, we prove that \eqref{structure:eq-JgradH-K} is globally well-posed. Furthermore, we obtain an energy identity and show that the flow preserves the volume form. 
 
 \begin{lemma}[Global flow, energy identity, and volume-preservation]\label{structure:lem-global-flow}
 Let $(M,\boldsymbol{\nu},J)$ be a conservative structure, let $H$ be a coercive Hamiltonian, and let $K$ be a conservative perturbation. Then, we have the following properties: 
 \begin{enumerate}[label=(\roman*)]
     \item (Global flow) The ordinary differential equation \eqref{structure:eq-JgradH-K} is globally well-posed and induces a unique, smooth global flow $\Gamma \colon \R \times \R \times M \rightarrow M, (t,t_0,\gamma_0)\mapsto \Gamma(t;t_0) \gamma_0$. 
     \item (Energy identity) For all $t\in \R$, $t_0\in \R$, and $\gamma_0 \in M$, it holds that
     \begin{equation*}
     H\big( \Gamma(t;t_0) \gamma_0 \big) = H\big( \gamma_0 \big) + \int_{t_0}^t \mathrm{d}s \, \Ks_{\Gamma(s;t_0)\gamma_0}(H).
     \end{equation*}
     In particular, if $K\equiv 0$, the energy is conserved. 
     \item (Volume-preservation) The flow $\Gamma$ is volume-preserving, i.e., for all $t,t_0\in \R$, it holds that 
     \begin{equation*}
        \Gamma(t;t_0)^\ast \boldsymbol{\nu} = \boldsymbol{\nu}.
     \end{equation*}
 \end{enumerate}
 \end{lemma}
 
 \begin{proof} We first prove the existence, uniqueness, and smoothness of the global flow $\Gamma$. To this end, 
let $\gamma_0 \in M$, let $I$ be the maximal time interval of existence for \eqref{structure:eq-JgradH-K}, and let $\gamma\colon I \rightarrow M$ be the corresponding maximal solution. Using the skew-symmetry of $J$, we obtain for all $t\in I$ that
\begin{equation}\label{structure:eq-gwp-p1}
\frac{\mathrm{d}}{\mathrm{d}t} H(\gamma(t)) = 
\del H_{\gamma(t)} \big( \dot{\gamma}(t) \big) = \del H_{\gamma(t)} \big( \Jt_{\gamma(t)} \del H_{\gamma(t)} + \Kt_{\gamma(t)} \big) = \del H_{\gamma(t)} \big( \Kt_{\gamma(t)} \big)=\Kt_{\gamma(t)}(H). 
\end{equation}
Since $K$ is a conservative perturbation, it follows that 
\begin{equation*}
\Big|\Kt_{\gamma(t)}(H) \Big| \leq C \Big( 1+ H(\gamma(t))\Big). 
\end{equation*}
Using Gronwall's inequality, this implies that 
\begin{equation}\label{structure:eq-gwp-p2}
H(\gamma(t)) \leq \exp\big( C |t|\big) \big( 1+ H(\gamma_0) \big).
\end{equation}
Due to Definition \ref{structure:def-coercive}, the sublevel sets of $H$ are compact, and thus \eqref{structure:eq-gwp-p2} implies the global well-posedness of  \eqref{structure:eq-JgradH-K}. Using classical ODE theory (see e.g. \cite[Section 9]{Lee13}), this also implies the existence,  uniqueness, and smoothness of the global flow $\Gamma \colon \R \times \R \times M \rightarrow M$. 

The desired energy identity follows directly from \eqref{structure:eq-gwp-p1}. Finally, due to our assumptions on $J$ and $K$, the time-dependent vector field $J \del H+K$ preserves the volume form $\boldsymbol{\nu}$. As a result, the corresponding flow $\Gamma$ satisfies 
$\Gamma(t;t_0)^\ast \boldsymbol{\nu} = \boldsymbol{\nu}$
for all $t,t_0\in \R$. 
 \end{proof}
 
In Definition \ref{structure:def-conservative-Riemannian} below, the place of the volume form $\boldsymbol{\nu}$ will be taken by the Riemannian volume form. While the Riemannian volume form is a central object in finite-dimensional Riemannian geometry, it does not have an infinite-dimensional analogue. In order to work in the infinite-dimensional setting, we introduce the Gibbs form, which has an infinite-dimensional analogue and corresponds to a weighted version of $\boldsymbol{\nu}$. 

\begin{definition}[Gibbs form]\label{structure:def-Gibbs-form} 
Let $(M,\boldsymbol{\nu},J)$ be a conservative structure and let $H$ be a coercive Hamiltonian. Furthermore, let $\beta>0$ and $\Zc>0$ be two parameters. Then, we define the Gibbs form $\Mu_\beta$ by 
\begin{equation*}
\Mu_\beta= \Zc^{-1} \exp\big( - \beta H \big) \boldsymbol{\nu}. 
\end{equation*}
\end{definition}

\begin{remark}
At this point, we have not imposed any integrability conditions on $H$, but the Gibbs form $\Mu_\beta$ is still well-defined as a volume form. However, once we consider the Gibbs measure $\mu_\beta$ (rather than the Gibbs form $\Mu_\beta$), an integrability condition will be imposed. \end{remark}

To avoid confusion, we note that, if $(M,\boldsymbol{\nu},J)$ is a conservative structure, then $(M,\Mu_\beta,J)$ is typically not a  conservative structure. The reason is that the second condition in Definition \ref{structure:def-conservative} will typically be violated. We also note that, as will be clear from Lemma \ref{structure:lem-Gibbs-form}, the time-dependent vector field $K$ in our setting preserves the volume form $\boldsymbol{\nu}$ but does not preserve the Gibbs form $\Mu_\beta$. 

\begin{lemma}[Pull-back of the Gibbs form]\label{structure:lem-Gibbs-form}
Let $(M,\boldsymbol{\nu},J)$ be a conservative structure, let $H$ be a coercive Hamiltonian, and let $K$ be a conservative perturbation. Furthermore, let $\Gamma \colon \R \times \R \times M \rightarrow M$ be the global flow from Lemma \ref{structure:lem-global-flow} and let $\Mu_\beta$ be the Gibbs form from Definition \ref{structure:def-Gibbs-form}. Then, it holds that 
\begin{equation*}
\big( \Gamma(t;t_0)^\ast \Mu_\beta \big)_p = \exp \bigg( - \beta  \int_{t_0}^t \mathrm{d}s \, \Ks_{\Gamma(s;t_0)p}(H) \bigg) \, (\Mu_\beta)_p
\end{equation*}
for all $p\in M$ and $t\in \R$.
\end{lemma}

\begin{proof}
For all $t\in \R$ and $p\in M$, it follows from the definition of the Lie derivative that 
\begin{equation}\label{structure:eq-transport-p1}
\frac{\mathrm{d}}{\mathrm{d}t} \Big( \Gamma(t;t_0)^\ast \Mu_\beta\Big)_p 
= \Big( \Gamma(t;t_0)^\ast \Lie_{\Jt \del H+ \Kt }( \Mu_\beta) \Big)_p.
\end{equation}
Using the product rule for Lie derivatives, using that $(M,\boldsymbol{\nu},J)$ is a conservative structure, and using that $K$ is a conservative perturbation, it follows that
\begin{align}
\Lie_{\Jt \del H + \Kt} ( \Mu_\beta ) 
&= \Zc^{-1} \big( \Jt \del H + \Kt \big)\big( e^{-\beta H} \big) \boldsymbol{\nu} 
+ \Zc^{-1} e^{-\beta H} \Lie_{\Jt \del H + \Kt }\big(\boldsymbol{\nu}\big) \notag \\
&= - \beta \Zc^{-1} \Kt\big(H\big) e^{-\beta H} \boldsymbol{\nu} \notag \\
&= - \beta \Kt(H) \Mu_\beta \label{structure:eq-transport-p2}.
\end{align}
By combining \eqref{structure:eq-transport-p1} and \eqref{structure:eq-transport-p2}, it follows that 
\begin{equation}\label{structure:eq-transport-p3}
\frac{\mathrm{d}}{\mathrm{d}t} \big( \Gamma(t;t_0)^\ast \Mu_\beta \big)_p = - \beta \Big( \Gamma(t;t_0)^\ast \big( \Kt\big(H\big) \, \Mu_\beta \big) \Big)_p = -\beta \Kt_{\Gamma(t;t_0)p}(H) \big( \Gamma(t;t_0)^\ast \Mu_\beta\big)_p. 
\end{equation}
By solving the (time-dependent) ordinary differential equation \eqref{structure:eq-transport-p3}, we obtain the desired conclusion.
\end{proof}

In the framework of abstract conservative structures, it is natural to work with differential forms. Once we move to the discretized wave maps equation, however, it is ultimately more natural (or at least common) to work with probability measures instead of differential forms. We therefore now convert Lemma \ref{structure:lem-Gibbs-form} into a more measure-theoretic formulation. To this end, we first make the following definition. 

\begin{definition}\label{structure:def-measures} 
Let $M$ be a smooth manifold and let $\boldsymbol{\nu}$ be a (positively oriented) orientation form on $M$. Then, we make the following definitions. 
\begin{enumerate}[label=(\roman*)]
    \item Let $\Borel(M)$ be the Borel $\sigma$-algebra of $M$. We define the volume measure $\Vol\colon \Borel(M) \rightarrow [0,\infty]$ by 
    \begin{equation*}
        \int_M f(p) \mathrm{d}\Vol(p)= \int_M f \boldsymbol{\nu} 
    \end{equation*}
    for all $f\in C^\infty_c(M)$, where the right-hand side denotes the integral of a volume form. 
    \item (Quantitative coercivity) A coercive Hamiltonian $H\colon M \rightarrow [0,\infty)$ is called quantitatively coercive if 
    \begin{equation*}
    \exp\big( - \beta H \big) \in L^1(M,\Vol)
    \end{equation*}
    for all $\beta>0$. 
    \item (Gibbs measure) For any $\beta>0$ and quantitatively coercive Hamiltonian $H\colon \R \rightarrow [0,\infty)$, we define the Gibbs measure $\mu_\beta\colon \Borel(M) \rightarrow [0,1]$ by 
    \begin{equation*}
    \int_M f(p) \, \mathrm{d}\mu_\beta(p) =\int_M f \Mu_\beta
    \end{equation*}
    for all $f\in C^\infty_c(M)$, where $\Mu_\beta$ is as in Definition \ref{structure:def-Gibbs-form} and $\Zc=\Zc_\beta>0$ is chosen such that $\mu_\beta(M)=1$. 
\end{enumerate}
\end{definition}

Equipped with Definition \ref{structure:def-measures}, we now prove the following corollary of Lemma \ref{structure:lem-Gibbs-form}. 

\begin{corollary}[Push-forward of Gibbs measure]\label{structure:cor-Gibbs-measure}
Let $M$ be an oriented smooth manifold, let $\boldsymbol{\nu}$ be a positively oriented volume form on $M$, let $J\colon \R\times  T^\ast M\rightarrow TM$ be a smooth, time-dependent bundle homomorphism, and assume that $(M,\boldsymbol{\nu},J)$ is a conservative structure. Furthermore, let $H\colon M\rightarrow [0,\infty)$ be a quantitatively coercive Hamiltonian and let $K$ be a conservative perturbation. Finally, let $\Gamma\colon \R \times \R \times  M\rightarrow M$ be the corresponding global flow, let $\beta>0$, and let $\mu_\beta$ be the Gibbs measure from Definition \ref{structure:def-measures}. Then, it holds for all $t,t_0\in \R$ that 
\begin{equation}\label{structure:eq-cor-Gibbs-measure}
\mathrm{d}\Big( \Gamma(t;t_0)_\# \mu_\beta\Big)
= \exp \bigg( \beta \int_{t_0}^{t} \ds (K(s)(H)\circ \Gamma(s;t)) \bigg) \mathrm{d}\mu_\beta,
\end{equation}
where $\Gamma(t;t_0)_\#$ is the push-forward. 
\end{corollary}

\begin{proof}
For any $f\in C^\infty_c(M)$, we first convert the integral against the measure $\Gamma(t;t_0)_\# \mu_\beta$ into an integral against the Gibbs form $\Mu_\beta$. It holds that 
\begin{equation}\label{structure:eq-measure-p1}
\int_M f(p) \mathrm{d} \big( \Gamma(t;t_0)_\# \mu_\beta \big)(p) 
= \int_M f \big( \Gamma(t;t_0) p \big) \mathrm{d} \mu_\beta (p)
= \int_M (f\circ \Gamma(t;t_0)) \Mu_\beta. 
\end{equation} 
In \eqref{structure:eq-measure-p1}, the last term is the integral of a volume form over $M$. 
Using that $\Gamma(t;t_0)^{-1} = \Gamma(t_0;t)$ and the coordinate-invariance of integrals of differential forms, it follows that 
\begin{equation}\label{structure:eq-measure-p2} 
\begin{aligned}
 \int_M (f\circ \Gamma(t;t_0)) \Mu_\beta 
 &= \int_M (f\circ \Gamma(t;t_0)) \Gamma(t;t_0)^\ast \Gamma(t_0;t)^\ast \Mu_\beta \allowdisplaybreaks[3] \\ 
 &= \int_M \Gamma(t;t_0)^\ast \Big( f \Gamma(t_0;t)^\ast \Mu_\beta \Big) \allowdisplaybreaks[3] \\ 
 &= \int_M f \Gamma(t_0;t)^\ast \Mu_\beta. 
 \end{aligned}
\end{equation}
By first using Lemma \ref{structure:lem-Gibbs-form} and then converting back from differential forms into measures, we obtain that 
\begin{equation}
\begin{aligned}\label{structure:eq-measure-p3}
    \int_M f \,\Gamma(t_0;t)^\ast \Mu_\beta 
    &= \int_M f \exp\bigg( - \beta \int_t^{t_0} \ds \big( K(s)(H) \circ \Gamma(s;t)\big) \bigg) \Mu_\beta \allowdisplaybreaks[3] \\
    &= \int_M f \exp\bigg(  \beta \int_{t_0}^{t} \ds \big( K(s)(H) \circ \Gamma(s;t)\big) \bigg) \Mu_\beta \allowdisplaybreaks[3] \\
    &= \int_M f(p) \exp\bigg( \beta \int_{t_0}^{t} \ds \big( K(s)(H) \circ \Gamma(s;t)\big)(p) \bigg) \mathrm{d}\mu_\beta(p). 
\end{aligned}
\end{equation}
By combining \eqref{structure:eq-measure-p1}, \eqref{structure:eq-measure-p2}, and \eqref{structure:eq-measure-p3}, we ultimately arrive at the desired identity \eqref{structure:eq-cor-Gibbs-measure}.
\end{proof}

As briefly mentioned above, we will now restrict ourselves to a Riemannian manifold $(M,g)$ and the corresponding Riemannian volume form $\boldsymbol{\nu}$. While the main results of this section (Lemma \ref{structure:lem-global-flow} and Lemma \ref{structure:lem-Gibbs-form}) can be obtained without the Riemannian structure, the Riemannian structure simplifies the verification of the conditions in Definition \ref{structure:def-conservative} and Definition \ref{structure:def-perturbation}. Furthermore, as will be discussed in Section \ref{section:conservative-wave-maps}, a simple Riemannian structure is available in the setting of wave maps.  

\begin{definition}[Conservative Riemannian structure]\label{structure:def-conservative-Riemannian}
A tuple $(M,g,\boldsymbol{\nu},J)$ is called a conservative Riemannian structure if $M$ is a smooth manifold, $g$ is a Riemannian metric on $M$, $\boldsymbol{\nu}$ is the corresponding Riemannian volume form, and $(M,\boldsymbol{\nu},J)$ is a conservative structure.
\end{definition}

In the following lemma, we provide a criterion for conservative Riemannian structures.

\begin{lemma}[Coordinate-dependent criteria for conservative Riemannian structures]
\label{structure:lem-conservative-Riemannian}
Let $M$ be a smooth  manifold, let $g$ be a Riemannian metric on $M$, and let $\boldsymbol{\nu}$ be the corresponding Riemannian volume form. Let $J\colon \R \times  T^\ast M\rightarrow TM$ be a smooth, time-dependent bundle homomorphism. Furthermore, assume that, for every $p\in M$, there exist local coordinates $(x^j)_{j=1}^d$ in a neighborhood of $p$ such that the following two properties are satisfied:
\begin{enumerate}[label=(\alph*)]
    \item (Skew-symmetry) \label{structure:item-riemannian-conservative-a} For all $t\in \R$ and $1\leq i,j\leq d$, it holds that $\Jt^{ij} = - \Jt^{ji}$. 
    \item (Volume-preservation) \label{structure:item-riemannian-conservative-b} For all $t\in \R$ and $1\leq j \leq d$, it holds that
    \begin{equation*}
    \frac{1}{\sqrt{\det(g)}} \partial_i \big( \sqrt{\det(g)} \Jt^{ij} \big)=0. 
    \end{equation*}
\end{enumerate}
Then, $(M,g,\boldsymbol{\nu},J)$ is a conservative Riemannian structure.
\end{lemma}

\begin{proof}
We only need to prove that the conditions \ref{structure:item-conservative-a} and \ref{structure:item-conservative-b} in Definition \ref{structure:def-conservative} are satisfied. Clearly, \ref{structure:item-conservative-a} in Definition \ref{structure:def-conservative} follows from our Assumption \ref{structure:item-riemannian-conservative-a}. In order to prove condition \ref{structure:item-conservative-b} from Definition \ref{structure:def-conservative}, we first note that
\begin{equation*}
\Lie_{\Jt\del H}\big( \boldsymbol{\nu} \big) = \operatorname{div} \big( \Jt \del H\big) \boldsymbol{\nu}.
\end{equation*}
In local coordinates, it holds that 
\begin{align}
\operatorname{div} \big( \Jt \del H \big) 
&= \frac{1}{\sqrt{\det(g)}} \partial_i \big( \sqrt{\det(g)} \Jt^{ij} \partial_j H \big) \notag \\
&= \frac{1}{\sqrt{\det(g)}} \partial_i \big( \sqrt{\det(g)} \Jt^{ij}  \big) \partial_j H 
+  \frac{1}{\sqrt{\det(g)}}  \sqrt{\det(g)} \Jt^{ij}   \partial_i \partial_j H. 
\label{structure:eq-conservative-Riemannian-p1}
\end{align}
The first summand in \eqref{structure:eq-conservative-Riemannian-p1} vanishes due to our Assumption \ref{structure:item-riemannian-conservative-b}. The second summand in \eqref{structure:eq-conservative-Riemannian-p1} vanishes due to the skew-symmmetry of $J$, i.e., our Assumption \ref{structure:item-riemannian-conservative-a}, and the symmetry of $\partial_i \partial_j H$. 
\end{proof}

\subsection{Conservative structure of wave maps}\label{section:conservative-wave-maps}
In this subsection, we study the conservative structure of the $2\pi R$-periodic, $(1+1)$-dimensional wave maps equation. We define the state space $\State_R$ and its smooth subspace $\State_R^\infty\subseteq \State_R$ by 
\begin{align*}
\State_R &:= L^2(\bT_R \rightarrow \frkg) \times L^2(\bT_R \rightarrow \frkg), \\
\State_R^\infty &:= \bigcap\displaylimits_{k=0}^\infty \Big( H^k(\bT_R \rightarrow \frkg) \times H^k(\bT_R \rightarrow \frkg) \Big). 
\end{align*}

The state space $\State_R$ is equipped with the $L^2$-inner product 
\begin{equation}\label{structure:eq-inner-product}
\left\langle \begin{pmatrix} A_1 \\ B_1 \end{pmatrix} , \begin{pmatrix} A_2 \\ B_2 \end{pmatrix} \right \rangle 
:= \int_{\bT_R } \dx \Big( \big \langle  A_1(x) , A_2(x) \big \rangle_\frkg +
 \big \langle  B_1(x) , B_2(x) \big \rangle_\frkg \Big),
\end{equation}
which turns $\State_R$ into a Hilbert space. As a result, we now identify
\begin{equation*}
T_{(A,B)} \State_R \simeq T^\ast_{(A,B)} \State_R \simeq \State_R
\end{equation*}
for all $(A,B)\in \State_R$. We define the Hamiltonian by 
\begin{equation}\label{structure:eq-Hamiltonian}
H(A,B) := \frac{1}{2} \int_{\bT_R } \dx \Big( \| A(x) \|_{\frkg}^2 + \| B(x) \|_{\frkg}^2 \Big). 
\end{equation}
The gradient of $H$ with respect to the $L^2$-inner product is given by 
\begin{equation}\label{structure:eq-Hamiltonian-gradient}
\big( \nabla_{L^2} H \big)_{(A,B)} = \begin{pmatrix} A \\ B \end{pmatrix}. 
\end{equation}

Finally, we define $\J$-operator, which, together with the Hamiltonian $H$ from above, will yield the wave maps equations \eqref{ansatz:eq-wave-map}. Since the $\J$-operator involves both derivatives and products, however, we need to work on $\State_R^\infty$ instead of $\State_R$.  
For any $t\in \R$ and $(A,B) \in \State_R^\infty$, we define the linear operator $\J(t)_{(A,B)} \colon \State_R^\infty \rightarrow \State_R^\infty$ by 
\begin{equation}\label{structure:eq-J}
\J(t)_{(A,B)} = 
\begin{pmatrix} 0 & \partial_x \\ \partial_x & 0 \end{pmatrix}
- \frac{1}{2} \chi(t) \begin{pmatrix} [A,\cdot\,] & [\,\cdot,B] \\ [\,\cdot,B] & [A,\cdot\,] \end{pmatrix}. 
\end{equation}
More precisely, for any $(C,D)\in \State_R^\infty$, we define 
\begin{equation*}
\J(t)_{(A,B)} \begin{pmatrix} C \\ D \end{pmatrix} 
= \begin{pmatrix} \partial_x D \\ \partial_x C \end{pmatrix} 
- \frac{1}{2} \chi(t) \begin{pmatrix} 
[A,C]+ [D,B] \\ 
[C,B] + [A,D] 
\end{pmatrix} . 
\end{equation*}
In particular, it holds that
\begin{equation*}
    \big( \J(t) \nabla_{L^2} H \big)_{(A,B)} = \begin{pmatrix} \partial_x B \\ \partial_x A - \chi(t) [A,B] \end{pmatrix},
\end{equation*}
which coincides with the right-hand side of the wave maps equation \eqref{ansatz:eq-wave-map}. 
If $\boldsymbol{\nu}$ formally represents the Lebesgue volume form on $\State_R$, then the tuple
$(\State_R,\boldsymbol{\nu}, \J)$ formally yields a conservative structure. However, since the infinite-dimensional Lebesgue volume form $\boldsymbol{\nu}$ cannot be defined rigorously, this aspect can only be discussed rigorously in the discretized setting. \\

We now turn to the discretized wave maps equation from \eqref{structure:eq-discretized-wave-maps}. 
While we make no assumption on the frequency-support of $\A[\Nscript,\Rscript,\coup],\B[\Nscript,\Rscript,\coup]\colon \R\times \bT_R \rightarrow \frkg$, we note that the discretized wave maps equation \eqref{ansatz:eq-discretized-wave-maps} is linear at frequencies $\gg N$. To be more precise, we write
\begin{align*}
\A[\Nscript,\Rscript,\coup] &= P_{\lesssim N}^x \A[\Nscript,\Rscript,\coup] + P_{\gg N}^x \A[\Nscript,\Rscript,\coup] =: \A[\Nscript,\Rscript,\coup][\lesssim N] + \A[\Nscript,\Rscript,\coup][\gg N], \\
\B[\Nscript,\Rscript,\coup] &= P_{\lesssim N}^x \B[\Nscript,\Rscript,\coup] + P_{\gg N}^x \B[\Nscript,\Rscript,\coup] =: \B[\Nscript,\Rscript,\coup][\lesssim N] + \B[\Nscript,\Rscript,\coup][\gg N]. 
\end{align*}
Then, the two evolution equations in \eqref{ansatz:eq-discretized-wave-maps} are equivalent to the four evolution equations 
\begin{align}
\partial_t \A[\Nscript,\Rscript,\coup][\lesssim N] &= \partial_x \B[\Nscript,\Rscript,\coup][\lesssim N],\label{structure:eq-dDWM-hl-1} \\
\partial_t \B[\Nscript,\Rscript,\coup][\lesssim N] &= \partial_x \A[\Nscript,\Rscript,\coup][\lesssim N] -\chi(t)  P_{\leq N}^x \Big[ P_{\leq N}^x \A[\Nscript,\Rscript,\coup][\lesssim N], P_{\leq N}^x \B[\Nscript,\Rscript,\coup][\lesssim N] \Big] + 2 \chi(t)^2 \coup \Renorm[N] \A[\Nscript,\Rscript,\coup][\lesssim N], \label{structure:eq-dDWM-hl-2} \\
\partial_t \A[\Nscript,\Rscript,\coup][\gg N] &= \partial_x \B[\Nscript,\Rscript,\coup][\gg N], \label{structure:eq-dDWM-hl-3} \\
\partial_t \B[\Nscript,\Rscript,\coup][\gg N] &= \partial_x \A[\Nscript,\Rscript,\coup][\gg N]. \label{structure:eq-dDWM-hl-4}
\end{align}
Thus, the discretized wave maps equation \eqref{ansatz:eq-discretized-wave-maps} can be decoupled into two finite-dimensional, nonlinear evolution equations \eqref{structure:eq-dDWM-hl-1}-\eqref{structure:eq-dDWM-hl-2} and two infinite-dimensional, linear evolution equations \eqref{structure:eq-dDWM-hl-3}-\eqref{structure:eq-dDWM-hl-4}. For all practical purposes, the discretized wave maps equation \eqref{ansatz:eq-discretized-wave-maps} can therefore be viewed as being finite-dimensional, and we now discuss its conservative structure. 
We define a finite-dimensional subspace $\StateN_R \subseteq \State_R$ by 
\begin{equation}\label{structure:eq-state-N}
\StateN_R:= \Big\{ (A,B) \in \State_R \colon P_{\gg N}^x A = P_{\gg N}^x B =0 \Big\}. 
\end{equation}
The finite-dimensional subspace $\StateN_R$ is also equipped with the inner product from \eqref{structure:eq-inner-product}. 
Furthermore, we define $\nuN$ as the canonical Lebesgue form on $\StateN_R$ (with any choice of orientation). 
We define a frequency-truncation approximation of $\J$ from \eqref{structure:eq-J} by 
\begin{equation}\label{structure:eq-J-N}
(\JNt)_{(A,B)} = \begin{pmatrix} 0 & \partial_x \\ \partial_x & 0 \end{pmatrix}
- \frac{1}{2} \chi(t) \begin{pmatrix} 
P_{\leq N}^x \big[ P_{\leq N}^x  A, P_{\leq N}^x \big( \, \cdot \big) \big] 
& P_{\leq N}^x\big[ P_{\leq N}^x  \big( \, \cdot \big),P_{\leq N}^x B\big] \\[1ex]
P_{\leq N}^x\big[ P_{\leq N}^x \big( \, \cdot \big),P_{\leq N}^x B\big] 
& P_{\leq N}^x\big[ P_{\leq N}^x  A,P_{\leq N}^x\big( \, \cdot \big)\big] \end{pmatrix}.
\end{equation}
With a slight abuse of notation, we also denote by $H$ the restriction of the Hamiltonian from  \eqref{structure:eq-Hamiltonian} to $\StateN_R$.
Finally, in order to capture the Killing-renormalization, we define a vector field $\KNt$ on $\StateN_R$ by 
\begin{equation}\label{structure:eq-KN}
\big( \KNt\big)_{(A,B)} := 
\begin{pmatrix} 0 \\ 2 \chi(t)^2 \coup \Renorm[N][] A \end{pmatrix}. 
\end{equation}
Equipped with $H$, $\JN$, and $\KN$, we arrive at the evolution equation
\begin{equation*}
\partial_t \begin{pmatrix} \A[\Nscript,\Rscript,\coup][\lesssim N] \\[0.5ex] \B[\Nscript,\Rscript,\coup][\lesssim N] \end{pmatrix} = \Big( \JNt \nabla_{L^2} H \Big)_{(\A[\Nscriptscript,\Rscriptscript,\coup][\lesssim N],\B[\Nscriptscript,\Rscriptscript,\coup][\lesssim N])} + \big( \KNt\big)_{(\A[\Nscriptscript,\Rscriptscript,\coup][\lesssim N],\B[\Nscriptscript,\Rscriptscript,\coup][\lesssim N])}, 
\end{equation*}
which is equivalent to the (nonlinear part of the) discretized wave maps equation, i.e., \eqref{structure:eq-dDWM-hl-1}-\eqref{structure:eq-dDWM-hl-2}. 
In the next proposition, we show that \eqref{structure:eq-dDWM-hl-1}-\eqref{structure:eq-dDWM-hl-2} comes from a conservative Riemannian structure, coercive Hamiltonian, and conservative perturbation.

\begin{proposition}[Conservative structure of discretized wave maps] \label{structure:prop-conservative-discrete}
Let $N\in \dyadic$, let $R\geq 1$, and let $\StateN_R$, $\nuN$, $\JN$, $H$, and $\KN$ be as above. Then, we have the following properties:
\begin{enumerate}[label=(\roman*)]
    \item\label{structure:item-conservative-wm-1} $\big(\StateN_R,\langle \cdot, \cdot\rangle_{L^2}, \nuN, \JN\big)$ is a conservative Riemannian structure.  
    \item\label{structure:item-conservative-wm-2} $H\colon \StateN_R \rightarrow \R$ is a coercive Hamiltonian. 
    \item\label{structure:item-conservative-wm-3} $\KN$ is a conservative perturbation with respect to $H$. 
\end{enumerate}
\end{proposition}

\begin{proof}
We prove the three properties \ref{structure:item-conservative-wm-1}, \ref{structure:item-conservative-wm-2}, and \ref{structure:item-conservative-wm-3} separately. Throughout the proof, we omit the time-dependence of $\JN$ and $\KN$ from our notation. \\

\emph{Proof of \ref{structure:item-conservative-wm-1}:} By the definition of the inner product and the Lebesgue form $\nuN$, it only remains to prove that $(\StateN_R,\nuN,\JN)$ is a conservative structure. To this end, we decompose
\begin{align*}
\JN_{(A,B)} &= \begin{pmatrix} 0 & \partial_x \\ \partial_x & 0 \end{pmatrix}
- \frac{1}{2} \chi(t) \begin{pmatrix} 
P_{\leq N}^x \big[ P_{\leq N}^x  A, P_{\leq N}^x \big( \, \cdot \big) \big] 
& P_{\leq N}^x\big[ P_{\leq N}^x  \big( \, \cdot \big),P_{\leq N}^x B\big]  \\[1ex]
P_{\leq N}^x\big[ P_{\leq N}^x \big( \, \cdot \big),P_{\leq N}^x B\big] 
& P_{\leq N}^x\big[ P_{\leq N}^x  A,P_{\leq N}^x\big( \, \cdot \big)\big]  \end{pmatrix} \\
&=: \JNI_{(A,B)} - \frac{1}{2} \chi(t) \JNII_{(A,B)}. 
\end{align*}
Using Lemma \ref{structure:lem-operations-conservative}, it suffices to prove that $(\StateN_R,\nuN,\JNI)$ and $(\StateN_R,\nuN,\JNII)$ are conservative structures. For $(\StateN_R,\nuN,\JNI)$, this is rather easy. Indeed, the skew-symmetry is clearly satisfied and, since $\JNI_{(A,B)}$ does not actually depend on $(A,B)$, the divergence-condition (see Lemma \ref{structure:lem-conservative-Riemannian}) is clearly satisfied as well. \\

It remains to prove that $(\StateN_R,\nuN,\JNII)$ is a conservative structure. We first prove the skew-symmetry of $\JNII$. Using the adjoint map from $\eqref{prelim:eq-adjoint-map}$ and the skew-symmetry of the Lie bracket, we write
\begin{equation*}
\JNII_{(A,B)} = 
\begin{pmatrix} 
P_{\leq N}^x \circ \Ad(P_{\leq N}^x A) \circ P_{\leq N}^x 
& -P_{\leq N}^x \circ \Ad(P_{\leq N}^x B) \circ P_{\leq N}^x  \\[1ex]
-P_{\leq N}^x \circ \Ad(P_{\leq N}^x B) \circ P_{\leq N}^x 
& P_{\leq N}^x \circ \Ad(P_{\leq N}^x A) \circ P_{\leq N}^x  \end{pmatrix}. 
\end{equation*}
Due to the self-adjointness of $P_{\leq N}^x$ and the skew-symmetry of the adjoint map (as in \eqref{prelim:eq-Hermitian-adjoint}), we see that 
\begin{equation*}
P_{\leq N}^x \circ \Ad(P_{\leq N}^x A) \circ P_{\leq N}^x \qquad \text{and} \qquad 
P_{\leq N}^x \circ \Ad(P_{\leq N}^x B) \circ P_{\leq N}^x 
\end{equation*}
are skew-symmetric, which implies the skew-symmetry of $\JNII$. Thus, it remains to prove that $\JNII$ is volume-preserving. Due to Lemma \ref{structure:lem-conservative-Riemannian}, it suffices to prove that, for all $(C,D)\in \StateN_R$, the vector field 
\begin{equation}\label{structure:eq-conservative-wm-p1}
\begin{aligned}
\begin{pmatrix} A \\ B \end{pmatrix} \mapsto 
\JNII_{(A,B)} \begin{pmatrix} C \\ D \end{pmatrix} 
= \begin{pmatrix} 
P_{\leq N}^x \big[ P_{\leq N}^x  A, P_{\leq N}^x C\big]    \\[1ex]
P_{\leq N}^x\big[ P_{\leq N}^x  C ,P_{\leq N}^x B\big]   \end{pmatrix}
+ 
\begin{pmatrix} 
 P_{\leq N}^x\big[ P_{\leq N}^x   D,P_{\leq N}^x B\big]  \\[1ex]
 P_{\leq N}^x\big[ P_{\leq N}^x  A,P_{\leq N}^x D\big] \end{pmatrix}.
 \end{aligned}
\end{equation} 
is divergence-free. This can be done via a direct calculation in coordinates, but we present a more conceptual argument. For this, we treat the two summands in \eqref{structure:eq-conservative-wm-p1} separately. For the first summand, we note that 
\begin{equation*}
\begin{pmatrix} 
P_{\leq N}^x \big[ P_{\leq N}^x  A, P_{\leq N}^x C\big]    \\[1ex]
P_{\leq N}^x\big[ P_{\leq N}^x  C ,P_{\leq N}^x B\big]   \end{pmatrix} 
= \begin{pmatrix} 
P_{\leq N}^x \big[ P_{\leq N}^x  A, P_{\leq N}^x C\big]    \\[1ex]
-P_{\leq N}^x\big[ P_{\leq N}^x  B ,P_{\leq N}^x C\big]   \end{pmatrix}. 
\end{equation*}
Thus, when taking the divergence, the contributions of the $A$ and $B$-components cancel each other. For the second summand in \eqref{structure:eq-conservative-wm-p1}, we only note the off-diagonal structure, i.e., the output in the $A$-component only depends on the $B$ and the output in the $B$-component only depends on $A$. This guarantees that each summand in the definition of the divergence is zero. \\ 

\emph{Proof of \ref{structure:item-conservative-wm-2}:} Since $H(A,B)= \tfrac{1}{2}\cdot \big( \| A \|_{L^2}^2 + \| B \|_{L^2}^2 \big)$,  the sublevel sets of $H$ are balls in $\StateN_R$. Since $\StateN_R$ is finite-dimensional, this implies the compactness of the sublevel sets. \\

\emph{Proof of \ref{structure:item-conservative-wm-3}:} We first prove the growth estimate on $\KN$. To this end, we note that 
\begin{align*}
\big( \KN H \big)_{(A,B)}
= \big\langle \KN_{(A,B)}, \nabla H_{(A,B)} \big\rangle 
=  2 \chi(t)^2 \coup  \, \left\langle \begin{pmatrix} 0 \\ \Renorm[N] A  \end{pmatrix}, \begin{pmatrix} A \\ B \end{pmatrix} \right\rangle. 
\end{align*}
As  a result, it follows that 
\begin{equation*}
\Big| \big( \KN H \big)_{(A,B)} \Big| \leq 2 \|\chi\|_{L^\infty}^2 \coup \big\| \Renorm[N]\big\|_{L^2 \rightarrow L^2} \big\| A \big\|_{L^2} \big\| B \big\|_{L^2} \leq 2\|\chi\|_{L^\infty}^2 \coup  \big\| \Renorm[N]\big\|_{L^2 \rightarrow L^2} H(A,B). 
\end{equation*}
Due to the boundedness\footnote{At this point, it is irrelevant whether the boundedness of $\Renorm[N]$ is uniform in $N$.} of $\Renorm[N]$ on $L^2$, we obtain the desired estimate. Thus, it remains to show that $\KN$ preserves the Lebesgue measure $\nuN$, i.e., it remains to show that $\KN$ is divergence-free. However, since $\KN$ induces a shear flow, this is\footnote{The output of $\KN$ in the $A$-component is zero and the output of $\KN$ in the $B$-component only depends on $A$.} clear.  
\end{proof}

The main result of this section, i.e., Proposition \ref{structure:prop-Gibbs}, is now a direct consequence of our previous considerations. 

\begin{proof}[Proof of Proposition \ref{structure:prop-Gibbs}]
Using Corollary \ref{structure:cor-Gibbs-measure} and Proposition \ref{structure:prop-conservative-discrete}, 
we obtain that 
\begin{equation*}
\begin{aligned}
&\, \frac{\mathrm{d}\muNR_{\chi,\tau} }{ \mathrm{d}\muR} 
\Big( W_0^{(\Rscript,\coup)},W_1^{(\Rscript,\coup)}\Big) \\
=&\,  \exp \bigg( \, 
2\beta \coup \int_{0}^\tau \int_{\bT_R} \ds \dx \, \chi(s)^2  
\big \langle \Renorm[N] A^{(\Nscript,\Rscript,\coup)}(s,x;\tau, \chi), 
B^{(\Nscript,\Rscript,\coup)}(s,x;\tau, \chi) \big \rangle
\bigg), 
\end{aligned}
\end{equation*}
Due to \eqref{intro:eq-coup}, it holds that $2\beta\coup=\frac{1}{4}$, which then yields the desired formula \eqref{structure:eq-Gibbs-formula}.
\end{proof}

\section{Abstract chaos estimate with dependent coefficients}\label{section:chaos}

In this section, we obtain an abstract chaos estimate with dependent coefficients, which will be one of the main ingredients in our stochastic estimates of Section \ref{section:Killing}.
To state the estimate, we first introduce convenient notation. For any $d_u,d_v\geq 1$, $m_u \in \Z^{d_u}$, and $m_v \in \Z^{d_v}$, we write 
\begin{equation*}
m_{u,\scrs} := \sum_{j=1}^{d_u} (m_u)_j \quad \text{and} \quad 
m_{v,\scrs} := \sum_{j=1}^{d_v} (m_v)_j. 
\end{equation*}
The subscript ``$\scrs$" stands for ``sum". 

\begin{proposition}[Chaos estimate with dependent coefficients]\label{chaos:prop-chaos}
Let $d_u,d_v \geq 1$, let $1\leq a_1,\hdots,a_{d_u}\leq \dim\frkg$, and let $1\leq b_1,\hdots,b_{d_v}\leq \dim \frkg$. Let $M\in \dyadic$ and let $\cs \colon \Z^{d_u}\times \Z^{d_v}\rightarrow \bC$ be a deterministic sequence supported on frequencies $m_u\in \Z^{d_u}$ and $m_v \in \Z^{d_v}$ satisfying $|m_u|,|m_v|\lesssim M$. Let $\alpha_1,\alpha_2,\beta_1,\beta_2,\gamma_1,\gamma_2\in (-1,1) \backslash \{0\}$ be regularity parameters, let $\Rcal_1$ and $\Rcal_2$ be frequency-scale relations, and assume that, for $j=1,2$, $(\alpha_j,\beta_j;\gamma_j)$ be $\Rcal_j$-admissible. Finally, let $\Wuv[\alpha_1][\alpha_2]$ be as in Definition \ref{chaos:def-wc} below. Then, it holds for all $\epsilon>0$ and $p\geq 1$ that 
\begin{equation}\label{chaos:eq-chaos}
\begin{aligned}
&\E \bigg[ \sup_{\| f \|_{\Wuv[\alpha_1][\alpha_2][]}\leq 1} \bigg\| \sum_{\substack{m_u\in \Z^{d_u} \\  m_v \in \Z^{d_v} }}  \big( P^{u,v}_{\leq M} f_{m_u m_v} \big) (u,v) \parauvrcal  \Big( \cs_{m_u m_v} \biglcol \prod_{j=1}^{d_u} g^{+,a_j}_{m_{u,j}} \prod_{j=1}^{d_v} g^{-,b_j}_{m_{v,j}} \, \bigrcol \\
&\hspace{1ex}\times \exp\big( i m_{u,\scrs} u \big)   \exp\big( i m_{v,\scrs} v \big) \Big) \bigg\|_{\Cprod{\gamma_1}{\gamma_2}}^p  \bigg]^{1/p} \\
\lesssim_{d_u,d_v,\epsilon}& \, \,  M^\epsilon  p^{(d_u+d_v)/2}\max_{L_u,L_v} \bigg[ L_u^{\beta_1} L_v^{\beta_2}  
\Big\|
\cs_{m_u m_v} \mathbf{1} \Big\{ \big|  m_{u,\scrs} \big| \sim L_u \Big\}  \mathbf{1} \Big\{ \big|  m_{v,\scrs} \big| \sim L_v \Big\} \Big\|_{\ell^2} \bigg].  
\end{aligned}
\end{equation}
\end{proposition}

\begin{remark}
We make the following remarks regarding Proposition \ref{chaos:prop-chaos}. 
\begin{enumerate}[label=(\roman*)]
    \item The chaos estimate with dependent coefficients will be one of the main ingredients in our analysis of modulated linear waves and modulated stochastic objects (see e.g. Section \ref{section:Killing}). 
    \item The gain from randomness in \eqref{chaos:eq-chaos} lies in the \revision{use of the $\ell^2$-norm rather than the $\ell^1$-norm of $\cs_{m_u m_v}$ on the right-hand side.}
    \item The most basic form of the chaos estimate concerns expressions of the form
    \begin{equation*}
        \Big| \sum_{m\in \Z} f_m g_m \Big|. 
    \end{equation*}
Here, $(f_m)_{m\in \Z}$ is a general random sequence and $(g_m)_{m\in \Z}$ is a sequence of standard Gaussians. The most important aspect of our general chaos estimate is that it does not require the probabilistic independence of $(f_m)_{m\in \Z}$ and $(g_m)_{m\in \Z}$. Instead, it requires $\ell^1$-control over discrete derivatives of  $(f_m)_{m\in \Z}$, \revision{which is encoded in the $\Wc$-norm (and $\Wuv[\alpha_1][\alpha_2][]$-norms) from Definition \ref{chaos:def-wc} below.}   
\item The idea behind our chaos estimate is similar as in the construction of the Wiener integral. In this construction, the integral of a smooth function $f\colon [0,1]\rightarrow \R$ against the derivative of a Brownian motion $B\colon [0,1]\rightarrow \R$ is defined as 
\begin{equation*}
\int_{0}^1 f(x) \mathrm{d}B(x) = f(1) B(1) - \int_0^1 f^\prime(x) B(x) \dx. 
\end{equation*}
\end{enumerate}
\end{remark}

The rest of this section is devoted to the proof of Proposition \ref{chaos:prop-chaos}. We first introduce discrete derivatives and shift operators. 

\begin{definition}[Discrete derivatives and shift operators]\label{chaos:def-derivative}
Let $d \geq 1$ be a positive integer and let $1\leq j \leq d$. We define the shift operator $\Theta_j$ and discrete derivative $\Der_j$ by
\begin{equation*}
\big( \Theta_j f \big)_m := f_{m+e_j} \qquad \text{and} \qquad \big( \Der_j f \big)_m := f_{m+e_j} - f_m
\end{equation*}
for all $f\colon \Z^d \rightarrow \bC$, where $e_j$ is the $j$th canonical basis vector. For any subset $\Jc=\{j_1,\hdots,j_k \} \subseteq \{ 1, \hdots,d\}$, we also define
\begin{equation}\label{chaos:eq-DJ-ThetaJ}
\Theta_\Jc f := \Theta_{j_1} \Theta_{j_2} \hdots \Theta_{j_k} f \qquad \text{and} \qquad
\Der_\Jc f:= \Der_{j_1} \Der_{j_2} \hdots \Der_{j_k} f. 
\end{equation}
\end{definition}
The discrete derivatives and shift operators satisfy a product rule and commute. To be more precise, let $f,h\colon \Z^d \rightarrow \bC$ and $1\leq i,j \leq d$. Then, it holds that 
\begin{equation}\label{chaos:eq-product-rule}
\Der_j \big( f h\big) = \Der_j f \cdot \Theta_j h + f \cdot \Der_j h
\end{equation}
and
\begin{equation}\label{chaos:eq-commutativity}
\Der_i \Der_j f = \Der_j \Der_i f, \quad \Theta_i \Der_j f = \Der_j \Theta_i f, \quad \text{and} \quad
\Theta_i \Theta_j f = \Theta_j \Theta_i f. 
\end{equation}
In particular, \eqref{chaos:eq-commutativity} implies that the order of the operators in \eqref{chaos:eq-DJ-ThetaJ} is irrelevant. For any $1\leq j \leq d$, the adjoint of $\Der_j$ on $\ell^2(\Z^d)$ is given by
\begin{equation}\label{chaos:eq-Dj-adjoint}
\big( \Der_j^\ast f \big)_m = - (f_m - f_{m-e_j} ). 
\end{equation}
The discrete derivatives and shift operators can also be extended from scalar-valued sequences $f=(f_m)_{m\in \Z^d}$ to function-valued sequences $f(u,v)=(f_m(u,v))_{m\in \Z^d}$. 

\begin{definition}[\protect{$\Wc$ and $\Wuv[\alpha_1][\alpha_2][]$}-norms and spaces]\label{chaos:def-wc}
Let $d\geq 1$ and $f\colon \Z^d \rightarrow \bC$. Then, we define 
\begin{equation}\label{chaos:eq-wc}
\big\| f \big\|_{\Wc} = \big\| f \big\|_{\Wc_m(\Z^d)} := \max_{\Jc \subseteq \{ 1, \hdots, d \}} \sum_{m_\Jc \in \Z^\Jc} \sup_{m_{\Jc^c} \in \Z^{\Jc^c}} \big| (\Der_\Jc f)_m \big|. 
\end{equation}
The corresponding space $\Wc$ is defined as the closure of 
\begin{equation}\label{chaos:eq-wc-space}
\Big\{  f\colon \Z^d   \rightarrow \bC \, \Big| \, f_m = f_{-m} \textup{ for all } m\in \Z^d, \, f_{m} = 0 \textup{ for all but finitely many } m \in \Z^d\Big\}
\end{equation}
with respect to the norm in \eqref{chaos:eq-wc}. 
For all $\alpha_1,\alpha_2\in \R \backslash \{0\}$ and $f\colon \Z^d \times \R^{1+1} \rightarrow \bC$, we also define 
\begin{equation}\label{chaos:eq-wc-uv}
\big\| f \big\|_{\Wuv[\alpha_1][\alpha_2][]} = \big\| f \big\|_{\Wuv[\alpha_1][\alpha_2][m](\Z^d)} :=  \max_{\Jc \subseteq \{ 1, \hdots, d \}} \sum_{m_\Jc \in \Z^\Jc} \sup_{m_{\Jc^c} \in \Z^{\Jc^c}} \big\| (\Der_\Jc f)_m(u,v) \big\|_{\Cprod{\alpha_1}{\alpha_2}}. 
\end{equation}
The corresponding space $\Wuv[\alpha_1][\alpha_2]$ is defined as the closure of 
\begin{equation}\label{chaos:eq-wc-uv-space}
\begin{aligned}
\Big\{& f\colon \Z^d \times \R^{1+1}   \rightarrow \bC \, \Big|\, f_m \equiv f_{-m} \textup{ and }  f_m \in \Cprod{\alpha_1}{\alpha_2} \textup{ for all } m \in \Z^d,  \\ 
&f_{m} \equiv 0 \textup{ for all but finitely many } m \in \Z^d\Big\}
\end{aligned}
\end{equation}
with respect to the norm in \eqref{chaos:eq-wc-uv}.
\end{definition}

\begin{remark}[Symmetry]
The symmetry condition $f_m = f_{-m}$ in \eqref{chaos:eq-wc-space} and \eqref{chaos:eq-wc-uv-space}  
is not needed for the proof of Proposition \ref{chaos:prop-chaos}. However, it will be needed in the proof of Proposition \ref{killing:prop-resonant} below, and we found it convenient to directly include it in Definition \ref{chaos:def-wc}.
\end{remark}

\begin{remark}\label{chaos:rem-wc}
\revision{While Definition \ref{chaos:def-wc} is stated only for $\bC$-valued functions, it can be extended to the vector-valued case. Indeed, let $\mathsf{V}$ be a finite-dimensional vector space over $\R$ or $\bC$. If $f\colon \Z^d\rightarrow \mathsf{V}$ or $f\colon \Z^d \times \R^{1+1} \rightarrow \mathsf{V}$, then the $\Wc$ and $\Wuv[\alpha_1][\alpha_2][]$-norms can be defined exactly as in \eqref{chaos:eq-wc} or \eqref{chaos:eq-wc-uv}, except that all absolute values (including those in Definition \ref{prelim:def-hoelder-spaces}) are replaced with the norm on $\mathsf{V}$. In particular, by choosing $\mathsf{V}=\End(\frkg)$, we can consider the $\Wuv[\alpha_1][\alpha_2]$-norms of the modulation operators from Definition \ref{ansatz:def-pure}.}
\end{remark}

As can be seen directly from Definition \ref{chaos:def-wc}, the $\Wc$ and $\Wuv[\alpha_1][\alpha_2][]$-norms consist of mixed $\ell^1\ell^\infty$-norms. The $\ell^1$-norm is taken in all variables involved in the discrete derivatives and the $\ell^\infty$-norm is taken in all remaining variables. In the following lemma, we record elementary properties of the $\Wc$ and $\Wuv[\alpha_1][\alpha_2][]$-spaces.

\begin{lemma}[\protect{Properties of $\Wc$ and $\Wuv[\alpha_1][\alpha_2][]$}]\label{chaos:lem-wc}
The $\Wc$-norm satisfies the following estimates:
\begin{enumerate}[label=(\roman*)]
\item\label{chaos:item-wc-boundedness} (Boundedness of discrete derivatives and shifts) For all $d\geq 1$, $f\colon \Z^d \rightarrow \bC$, and $\Jc \subseteq \{1,\hdots, d\}$, it holds that
\begin{equation*}
\big\| \Theta_\Jc f \big\|_{\Wc} + \big\| \Der_\Jc f \big\|_{\Wc} \lesssim \big\| f \big\|_{\Wc}. 
\end{equation*}
\item\label{chaos:item-wc-tensor} (Tensor products) Let $d\geq 1$ and, for each $1\leq j \leq d$, let $f^{(j)}\colon \Z \rightarrow \bC$. Then, it holds that
\begin{equation*}
\Big\| \,   \prod_{j=1}^d f^{(j)}_{m_j} \,  \Big\|_{\Wc(\Z^d)} \lesssim  \prod_{j=1}^d \big\| f^{(j)}  \big\|_{\Wc(\Z)}.
\end{equation*}
\item\label{chaos:item-wc-product} (Products) Let $d\geq 1$ and $f,h \colon \Z^d \rightarrow \bC$. Then,
\begin{equation*}
\big\| f h \big\|_{\Wc(\Z^d)} \lesssim \big\| f \big\|_{\Wc(\Z^d)}\big\| h \big\|_{\Wc(\Z^d)}. 
\end{equation*}
\item\label{chaos:item-wc-indicator} (Indicator) Let $a,b \in \R$. Then, it holds that 
\begin{equation*}
\big\| \mathbf{1} \big\{ a\leq |m| \leq b\big\} \big\|_{\Wc(\Z)} \lesssim 1. 
\end{equation*}
\item \label{chaos:item-wc-cutoff} (Cut-off) Let $M$ be a dyadic integer and let $\rho_{\leq M}$ be as in \eqref{prelim:eq-rho-leqN}. Then, it holds that
\begin{equation*}
    \big\| \rho_{\leq M}(m) \big\|_{\Wc(\Z)} \lesssim 1.
\end{equation*}
\end{enumerate}
Similar estimates also hold if $\Wc$ is replaced by $\Wuv[\alpha_1][\alpha_2][]$.
\end{lemma}

\begin{proof}
Since the proofs of all five estimates \ref{chaos:item-wc-boundedness}-\ref{chaos:item-wc-cutoff} are elementary, we omit the details.
\end{proof}

Equipped with Definition \ref{chaos:def-wc} and Lemma \ref{chaos:lem-wc}, we now prove a chaos estimate for scalar-valued sequences. As we will see below, this is the main step in the proof of Proposition \ref{chaos:prop-chaos}. 

\begin{lemma}[Scalar-valued chaos estimate with dependent coefficients] \label{chaos:lem-scalar}
Let $\mathcal{A}$ be an index set and for all indices $a\in \mathcal{A}$, let $(g_m^a)_{m\in \Z^d}$ be an independent copy 
of a standard, $\mathbb{C}$-valued Gaussian sequence (as in Definition \ref{prelim:def-standard-Gaussian-sequence}). 
Let $d\geq 1$, let $a_1,\hdots,a_d \in \mathcal{A}$, let $M$ be a dyadic integer, and let $\cs \colon \Z^d \rightarrow \mathbb{C}$ be a deterministic sequence supported on frequencies satisfying $|m_j|\leq M$ for all $1\leq j \leq d$. 
 For all $\epsilon>0$ and $p\geq 1$, it then holds that 
\begin{equation*}
\E \Big[ \sup_{\| f \|_{\Wc(\Z^d)} \leq 1} \Big| \sum_{m\in \Z^d} f_m \cs_m \lcol  \prod_{j=1}^d g_{m_j}^{a_j}  \rcol \Big|^p \Big]^{1/p}
\lesssim_\epsilon M^\epsilon p^{d/2} \big\| \cs_m \big\|_{\ell^2_m}.
\end{equation*}
\end{lemma}

The main idea behind Lemma \ref{chaos:lem-scalar} is summation by parts. While the Gaussian chaos $(\cs_m \lcol \prod_{j=1}^d g_{m_j}^{a_j} \rcol)_{m\in\Z^d}$ heavily oscillates in $m\in \Z^d$, the random sequences $f=(f_m)_{m\in \Z^d}$ vary slowly in $m\in \Z^d$, which leads to a gain.

\begin{proof}[Proof of Lemma \ref{chaos:lem-scalar}]
Due to the support condition on the deterministic sequence $(\cs_m)_{m\in \Z^d}$ and Lemma \ref{chaos:lem-wc}, we can always replace $f \in \Wc$ with
\begin{equation*}
\Big( \prod_{j=1}^d \mathbf{1} \big\{ |m_j|\leq M \big\} \Big) f_m. 
\end{equation*}
As a result, we now simply assume that all random sequences $f=(f_m)_{m\in \Z^d}$ satisfy the same support condition as $(\cs_m)_{m\in \Z^d}$. We now define a Gaussian chaos $\mathcal{G}\colon \Z^d \rightarrow \bC$ by 
\begin{equation}\label{chaos:eq-scalar-p1}
\mathcal{G}_m := \sum_{\substack{n\in \Z^d \colon \\  n \leq m}} \cs_{n} \lcol \prod_{j=1}^d g_{n_j}^{a_j} \rcol. 
\end{equation}
In \eqref{chaos:eq-scalar-p1}, the inequality $n\leq m$ is componentwise, i.e., $n_j \leq m_j$ for all $1\leq j \leq d$. From \eqref{chaos:eq-Dj-adjoint} and \eqref{chaos:eq-scalar-p1}, it follows for all $m\in \Z^d$ that 
\begin{equation}\label{chaos:eq-scalar-p2}
\big( \Der_1^\ast \Der_2^\ast \hdots \Der_d^\ast G \big)_m = (-1)^d \cs_m  \lcol \prod_{j=1}^d g_{m_j}^{a_j} \rcol
\end{equation}
Using summation by parts, it follows that
\begin{equation}\label{chaos:eq-scalar-p3}
\begin{aligned}
&\Big| \sum_{m\in \Z^d} f_m \cs_m \lcol \prod_{j=1}^d g_{m_j}^{a_j} \rcol \Big|
= \Big| \sum_{m\in \Z^d} f_m \big( \Der_1^\ast \Der_2^\ast \hdots \Der_d^\ast G \big)_m   \Big| 
= \Big| \sum_{m\in \Z^d} \big( \Der_1 \Der_2 \hdots \Der_d f\big)_m \mathcal{G}_m \Big| \\
\leq&\,  \sum_{m\in \Z^d}\big|  \big( \Der_1 \Der_2 \hdots \Der_d f\big)_m \big| \times \sup_{m\in \Z^d} |\mathcal{G}_m| 
\leq \| f \|_{\Wc(\Z^d)} \sup_{m\in \Z^d} |\mathcal{G}_m|. 
\end{aligned}
\end{equation}
As a result, it follows that 
\begin{equation}\label{chaos:eq-scalar-p4}
\E \Big[ \sup_{\| f \|_{\Wc(\Z^d)} \leq 1} \Big| \sum_{m\in \Z^d} f_m \cs_m \lcol  \prod_{j=1}^d g_{m_j}^{a_j}  \rcol \Big|^p \Big]^{1/p} \leq 
\E \Big[ \sup_{m\in \Z^d} |\mathcal{G}_m|^p \Big]^{1/p}. 
\end{equation}
Due to Gaussian hypercontractivity (Lemma \ref{prelim:lem-hypercontractivity}), Lemma \ref{prelim:lem-maxima},  and the fact that $(\mathcal{G}_m)_{m\in \Z^d}$ is constant outside of the support of $(\cs_m)_{m\in \Z^d}$, it then suffices to prove that
\begin{equation*}
\sup_{m\in \Z^d} \E \Big[ |\mathcal{G}_m|^2 \Big] \lesssim \big\| \cs_m \big\|_{\ell^2_m}^2.
\end{equation*}
This follows directly from the definition of $\mathcal{G}$ and  orthogonality (as in \ref{prelim:lem-orthogonality}).
\end{proof}

Equipped with Lemma \ref{chaos:lem-scalar}, we are now ready to prove the main result of this section.

\begin{proof}[Proof of Proposition \ref{chaos:prop-chaos}] 
Using the definition of the para-product operator (Definition \ref{prelim:def-para}), we can decompose the argument from \eqref{chaos:eq-chaos} as 
\begin{equation}\label{chaos:eq-decomposition}
\begin{aligned}
&\sum_{\substack{m_u\in \Z^{d_u} \\  m_v \in \Z^{d_v} }}  P^{u,v}_{\leq M} f_{m_u m_v}(u,v) \parauvrcal  \Big( \cs_{m_u m_v} \biglcol \prod_{j=1}^{d_u} g^{+,a_j}_{m_{u,j}} \prod_{j=1}^{d_v} g^{-,b_j}_{m_{v,j}} \, \bigrcol \,  \exp\big( i m_{u,\scrs} u \big)   \exp\big( i m_{v,\scrs} v \big) \Big)\\
=& \sum_{\substack{K_u,K_v,L_u,L_v \lesssim M\colon \\ 
(K_u,L_u) \in \Rcal_1, \\ (K_v,L_v) \in \Rcal_2}}
\Chaos_{K_u,K_v,L_u,L_v}[f](u,v),
\end{aligned}
\end{equation}
where 
\begin{align*}
\Chaos_{K_u,K_v,L_u,L_v}[f](u,v) :=&  \sum_{\substack{m_u \in \Z^{d_u} \\ m_v\in \Z^{d_v}}} \hspace{-1.5ex} \bigg( \big( P_{K_u}^u P_{K_v}^v P^{u,v}_{\leq M} f_{m_u m_v}\big)(u,v) 
\cs_{m_u m_v}^{(L_u,L_v)} \lcol \prod_{j=1}^{d_u} g_{m_{u,j}}^{+,a_j} \prod_{j=1}^{d_v} g_{m_{v,j}}^{-,b_j} \rcol 
e^{im_{u,\scrs} u } e^{i m_{v,\scrs} v } \bigg), \\ 
\cs_{m_u m_v}^{(L_u,L_v)} :=& \cs_{m_u m_v} \mathbf{1} \Big\{ \big|  m_{u,\scrs} \big| \sim L_u \Big\}  \mathbf{1} \Big\{ \big|  m_{v,\scrs} \big| \sim L_v \Big\}.
\end{align*}
We emphasize that the restriction $K_u,K_v \lesssim M$ in \eqref{chaos:eq-decomposition} is due to the $P^{u,v}_{\leq M}$-operator acting on $f$ in \eqref{chaos:eq-chaos} and the restriction $L_u,L_v \lesssim M$ in \eqref{chaos:eq-decomposition} is due to our assumption on the support of $\cs_{m_u m_v}$. \revision{By definition, the $\Cprod{\gamma_1}{\gamma_2}$-norm of \eqref{chaos:eq-decomposition} is given by  
\begin{equation}\label{chaos:eq-decomposition-norm}
\sup_{N_u, N_v \geq 1} N_u^{\gamma_1} N_v^{\gamma_2} 
\bigg\| P_{N_u}^u P_{N_v}^v \sum_{\substack{K_u,K_v,L_u,L_v \lesssim M\colon \\ 
(K_u,L_u) \in \Rcal_1, \\ (K_v,L_v) \in \Rcal_2}}
\Chaos_{K_u,K_v,L_u,L_v}[f](u,v) \bigg\|_{L^\infty_{u,v}}.
\end{equation}
From frequency-support considerations, one sees that 
$P_{N_u}^u P_{N_v}^v \Chaos_{K_u,K_v,L_u,L_v}[f]$ is only non-zero if, for both $w=u$ and $w=v$, it holds that $N_w \sim K_w \gg L_w$, $N_w \sim L_w \gg K_w$, or $K_w \sim L_w \gtrsim N_w$. We then obtain 
\begin{align}
\eqref{chaos:eq-decomposition-norm} &\lesssim \sup_{N_u, N_v \geq 1} N_u^{\gamma_1} N_v^{\gamma_2} \hspace{-3ex} \sum_{\substack{K_u,K_v,L_u,L_v \lesssim M\colon \\ 
(K_u,L_u) \in \Rcal_1, \\ (K_v,L_v) \in \Rcal_2}}    
\Big( \big( 
1\big\{ N_u \sim K_u \gg L_u \big\} + 1 \big\{ N_u \sim L_u \gg K_u \big\}
+ 1 \big\{ K_u \sim L_u \gtrsim N_u \big\}\big) \notag \\ 
&\quad \times \big( 
1\big\{ N_v \sim K_v \gg L_v \big\} + 1 \big\{ N_v \sim L_v \gg K_v \big\}
+ 1 \big\{ K_v \sim L_v \gtrsim N_v \big\}\big) 
\big\|  \Chaos_{K_u,K_v,L_u,L_v}[f](u,v) \big\|_{L^\infty_{u,v}} \Big) \notag  \\ 
&\lesssim   \sup_{\substack{K_u,K_v, \\ L_u,L_v \lesssim M}} 
K_u^{\alpha_1} K_v^{\alpha_2} L_u^{\beta_1} L_v^{\beta_2} 
\Big\| \Chaos_{K_u,K_v,L_u,L_v}[f](u,v)\Big\|_{L^\infty_{u,v}} \notag \\ 
&\lesssim    \sum_{\substack{K_u,K_v, \\ L_u,L_v \lesssim M}} 
K_u^{\alpha_1} K_v^{\alpha_2} L_u^{\beta_1} L_v^{\beta_2} 
\Big\| \Chaos_{K_u,K_v,L_u,L_v}[f](u,v)\Big\|_{L^\infty_{u,v}}. \label{chaos:eq-chaos-q0} 
\end{align}
In the second inequality above, we also used that $(\alpha_j,\beta_j;\gamma_j)$ is $\Rcal_j$-admissible for $j=1,2$.}
To simplify the notation, we now set $d:= d_u + d_v$. Using Gaussian hypercontractivity (Lemma \ref{prelim:lem-hypercontractivity}) and crude estimates, one can easily prove that
\begin{equation}\label{chaos:eq-chaos-q1}
K_u^{\alpha_1} K_v^{\alpha_2}  \sup_{p\geq 1} p^{-d/2}
\E \Big[ \Lip \Big( \sup_{\| f \|_{\Wuv[\alpha_1][\alpha_2][]} \leq 1} \big| \Chaos_{K_u,K_v,L_u,L_v}[f](u,v)\big| \Big)^p \Big]^{1/p} 
\lesssim M^{10d} \big\|  \cs_{m_u m_v}^{(L_u,L_v)} \big\|_{\ell^2},
\end{equation}
where $\Lip$ refers to the Lipschitz constant in the $u$ and $v$-variables. By combining \eqref{chaos:eq-chaos-q1} and our meshing argument (Corollary \ref{prelim:cor-meshing}), it follows that 
\begin{align}
&K_u^{\alpha_1} K_v^{\alpha_2} L_u^{\beta_1} L_v^{\beta_2} 
\sup_{p\geq 1} p^{-d/2}
\E \Big[  \sup_{\| f \|_{\Wuv[\alpha_1][\alpha_2][]}\leq 1}
\Big\| \Chaos_{K_u,K_v,L_u,L_v}[f](u,v)\Big\|_{L^\infty_{u,v}}^p \Big]^{1/p} \notag \\
\lesssim&\, M^{\epsilon/4} K_u^{\alpha_1} K_v^{\alpha_2} L_u^{\beta_1} L_v^{\beta_2}  \sup_{u,v\in \R} \sup_{p\geq 1} p^{-d/2} 
\E \Big[  \sup_{\| f \|_{\Wuv[\alpha_1][\alpha_2][]} \leq 1}
\Big| \Chaos_{K_u,K_v,L_u,L_v}[f](u,v)\Big|^p \Big]^{1/p} \label{chaos:eq-chaos-q2}  \\
+&\, M^{-100d} L_u^{\beta_1} L_v^{\beta_2}\big\|  \cs_{m_u m_v}^{(L_u,L_v)} \big\|_{\ell^2}. \label{chaos:eq-chaos-q3}
\end{align}
The second term \eqref{chaos:eq-chaos-q3} clearly yields an acceptable contribution, and it therefore remains to estimate \eqref{chaos:eq-chaos-q2}. 
Using our scalar-valued chaos estimate (Lemma \ref{chaos:lem-scalar}), it holds that 
\begin{align*}
\eqref{chaos:eq-chaos-q2} 
&\lesssim M^{\epsilon/2} 
 \bigg(  K_u^{\alpha_1} K_v^{\alpha_2} \sup_{u,v\in \R}  \sup_{\| f \|_{\Wuv[\alpha_1][\alpha_2][]} \leq 1} 
\big\| P_{K_u}^u P_{K_v}^v P_{\leq M}^{u,v} f_{m_u m_v}(u,v) \big\|_{\Wc(\Z^d)}  \bigg)  L_u^{\beta_1} L_v^{\beta_2} \big\|  \cs_{m_u m_v}^{(L_u,L_v)} \big\|_{\ell^2} \\ 
&\lesssim M^{\epsilon/2} L_u^{\beta_1} L_v^{\beta_2}\big\|  \cs_{m_u m_v}^{(L_u,L_v)} \big\|_{\ell^2}. 
\end{align*}
Since we have an additional $M^{\epsilon/2}$-factor at our disposal, we can perform the sum in \eqref{chaos:eq-chaos-q0}, and then obtain the desired estimate. 
\end{proof}

\section{Chaos estimates for modulated objects}\label{section:Killing}

In Subsection \ref{section:ansatz-contraction}, we stated our probabilistic hypothesis (Hypothesis \ref{hypothesis:probabilistic}), which collects several probabilistic estimates. The main goal of this section is to verify that the probabilistic hypothesis is satisfied with high probability and to examine its implications. To this end, we make the following definition, \revision{which involves the parameter $\Ac$ (see Hypothesis \ref{hypothesis:probabilistic} and Remark \ref{ansatz:rem-parameters}).}

\begin{definition}\label{killing:def-probabilistic-hypothesis} 
Let $R\geq 1$ and let $\Ac \geq 1$. 
Then, we define the event $\PHA\subseteq \Omega$ as the event on which all estimates stated in our probabilistic hypothesis (Hypothesis \ref{hypothesis:probabilistic}) are satisfied. 
\end{definition}

In the following proposition, we control the probability of $\PHA$.  

\begin{proposition}[Probabilistic hypothesis]\label{killing:prop-probabilistic-hypothesis} 
Let $c=c(\delta_\ast)$ be a sufficiently small constant. For all $R\geq 1$ and $\Ac \geq 1$, it then holds that
\begin{equation}\label{killing:eq-probabilistic-hypothesis}
\mathbb{P} \big( \Omega \backslash \PHA \big) \leq  \revision{c^{-1}} \exp\Big( - c R^{-2\eta} \Ac^2 \Big).  
\end{equation}
\end{proposition}

Before the proof of Proposition \ref{killing:prop-probabilistic-hypothesis}, which is postponed until Subsection \ref{section:killing-proof}, we first need to prove several chaos estimates involving modulated linear waves and their tensor products\footnote{In this section, tensor products always refer to tensor products in the Lie algebra $\frkg$. In particular, the tensor products should not be confused with the random tensors of \cite{DNY22}.}. Our reasons for working with tensor products are two-fold: 
First, the tensor product $\frkg\otimes \frkg$ is a canonical object, and, in particular, its construction does not rely on a choice of coordinates on $\frkg$. This implies that certain objects in $\frkg\otimes\frkg$ obey natural transformation laws (see e.g. Lemma \ref{prelim:lem-casimir}).  Second, the tensor product can be used to control several other expressions, such as Lie brackets or iterated Lie brackets. This perspective is motivated by recent works on the stochastic Yang-Mills heat equation, see e.g. \cite{BC23,CCHS22}.
As part of our analysis of tensor products (see Subsection \ref{section:killing-resonant}), we will also naturally encounter the Casimir from Definition \ref{prelim:def-casimir}, which will then lead to the Killing-renormalization.\\ 

In order to state our estimates, we still need to define the $\Wfuv[\alpha][\beta]$-spaces,
which previously appeared in our probabilistic hypothesis (Hypothesis \ref{hypothesis:probabilistic}). The $\Wfuv[\alpha][\beta]$-spaces are subspaces of the $\Wuv[\alpha][\beta]$-spaces from Definition \ref{chaos:def-wc}, but also require an additional frequency-support condition. 

\begin{definition}[\protect{$\Wfuv[\alpha][\beta]$-spaces}]\label{killing:def-FW}
Let $\alpha,\beta \in (-1,1)\backslash\{0\}$. Then, we define $\Wfuv[\alpha][\beta]$
as the subspace of $\Wuv[\alpha][\beta]=\Wuv[\alpha][\beta](\Z)$ consisting of all $(S_\ell)_{\ell \in \Z}$ which satisfy the following frequency-support condition: 
For all $L\in \dyadic$ and all $\ell \in \Z_L$, it holds that $P^{u,v}_{\gg L^{1-\delta+\vartheta}} S_\ell =0$.
\end{definition}

\subsection{Modulated linear waves}\label{section:killing-linear}
The goal of this short subsection is to prove first-order chaos estimates. Their primary application will be regularity estimates for the modulated linear waves $\UN[K][+]$ and $\VN[M][-]$ and the integrated modulated linear waves $\IUN[K][+]$ and $\IVN[M][-]$.

\begin{lemma}\label{killing:lem-regularity-modulated-linear}
For all $\alpha,\beta \in (-1,1)\backslash\{0\}$, all $\epsilon>0$, and all $p\geq 1$, it holds that 
\begin{align}
\E \bigg[ \sup_{u_0 \in \LambdaRR} \sup_{\Nd \in \Dyadiclarge} \sup_{K \in \Dyadiclarge} 
\sup_{\substack{\Smod[K][+]\in \Wfuv[s][\beta]\colon \\ \| \Smod[K][+] \|_{\Wuv[s][\beta][]} \leq 1}}
\Big( K^{-\alpha-\frac{1}{2}-\epsilon} \Big\|
\sum_{k\in \Z_K} \rhoND(k) \Smod[K][+][k] G_{u_0,k}^+ \, e^{iku} 
\Big\|_{\Cprod{\alpha}{\beta}} \Big)^p \bigg]^{1/p} &\lesssim \sqrt{p} R^\epsilon,
\label{killing:eq-modulated-linear-U}\\ 
\E \bigg[ \sup_{v_0 \in \LambdaRR} \sup_{\Nd \in \Dyadiclarge} \sup_{M \in \Dyadiclarge}
\sup_{\substack{\Smod[M][-]\in \Wfuv[\alpha][s]\colon \\ \| \Smod[M][-] \|_{\Wuv[\alpha][s][]} \leq 1}} 
 \Big( M^{-\beta-\frac{1}{2}-\epsilon} \Big\| 
\sum_{m\in \Z_M} \rhoND(m) \Smod[M][-][m] G_{v_0,m}^- \, e^{imv} \Big\|_{\Cprod{\alpha}{\beta}} \Big)^p \bigg]^{1/p} &\lesssim \sqrt{p} R^\epsilon. 
\label{killing:eq-modulated-linear-V}
\end{align}
Similarly, it holds that 
\begin{align}
\E \bigg[ \sup_{u_0 \in \LambdaRR} \sup_{\Nd \in \Dyadiclarge} \sup_{K \in \Dyadiclarge} 
\sup_{\substack{\Smod[K][+]\in \Wfuv[s][\beta]\colon \\ \| \Smod[K][+] \|_{\Wuv[s][\beta][]} \leq 1}}
\Big( K^{-\alpha+\frac{1}{2}-\epsilon}\Big\|
\sum_{k\in \Z_K} \rhoND(k) \Smod[K][+][k] G_{u_0,k}^+ \, \frac{e^{iku}}{ik} 
\Big\|_{\Cprod{\alpha}{\beta}}  \Big)^p \bigg]^{1/p} &\lesssim \sqrt{p} R^\epsilon, 
\label{killing:eq-integrated-modulated-linear-U} \\ 
\E \bigg[  \sup_{v_0 \in \LambdaRR} \sup_{\Nd \in \Dyadiclarge} \sup_{M \in \Dyadiclarge}
\sup_{\substack{\Smod[M][-]\in \Wfuv[\alpha][s]\colon \\ \| \Smod[M][-] \|_{\Wuv[\alpha][s][]} \leq 1}} 
\Big( M^{-\beta+\frac{1}{2}-\epsilon}  \Big\| 
\sum_{m\in \Z_M}  \rhoND(m) \Smod[M][-][m] G_{v_0,m}^- \, \frac{e^{imv}}{im} \Big\|_{\Cprod{\alpha}{\beta}} \Big)^p \bigg]^{1/p} &\lesssim \sqrt{p} R^\epsilon.
\label{killing:eq-integrated-modulated-linear-V}
\end{align}
\end{lemma}

\begin{remark}
In light of Definition \ref{ansatz:def-modulated-linear}, \eqref{killing:eq-modulated-linear-U} directly yields a regularity estimate for the modulated linear waves $\UN[K][+]$. The reason for stating our estimate as in \eqref{killing:eq-modulated-linear-U}, and not simply as an estimate for $\UN[K][+]$, is that \eqref{killing:eq-modulated-linear-U} is more versatile. In addition to the regularity estimate for $\UN[K][+]$, it will also be used for commutator estimates (see e.g. Lemma \ref{modulation:lem-PNX-modulated} and Lemma \ref{modulation:lem-integration}). 
\end{remark}

\begin{proof}
We only prove \eqref{killing:eq-modulated-linear-U}, since the proofs of \eqref{killing:eq-modulated-linear-V}, \eqref{killing:eq-integrated-modulated-linear-U}, and \eqref{killing:eq-integrated-modulated-linear-V} are similar. Due to Lemma \ref{prelim:lem-maxima}, we may reduce to a single $u_0\in \LambdaR$. Due to Lemma \ref{chaos:lem-wc}, it holds that 
\begin{equation*}
\Big\| \rhoND(k) \Smod[K][+][k] \Big\|_{\Wuv[s][s][k]} 
\lesssim \big\| \rhoND(k) \big\|_{\Wc_k} \Big\| \Smod[K][+][k] \Big\|_{\Wuv[s][s][k]}
\lesssim \Big\| \Smod[K][+][k] \Big\|_{\Wuv[s][s][k]}.
\end{equation*}

Furthermore, due to Definition \ref{killing:def-FW} and $\Smod[K][+][]\in \Wfuv[s][\beta]$, the expression
\begin{equation*}
 \sum_{k\in \Z_K} \rhoND(k) \Smod[K][+][k](u,v)G_{u_0,k}^+ \, e^{iku}
\end{equation*}
only contains low$\times$high-interactions in the $u$-variable and only contains $v$-frequencies at scales $\lesssim K$. Using our chaos estimate (Proposition \ref{chaos:prop-chaos}), it then follows that 
\begin{align*}
    \E \bigg[  \sup_{\Nd\in \Dyadiclarge} \sup_{\substack{\Smod[K][+]\in \Wfuv[s][\beta]\colon \\ \| \Smod[K][+] \|_{\Wuv[s][\beta][]} \leq 1}}
 \Big\| \sum_{k\in \Z_K}\rhoND(k) \Smod[K][+][k](u,v)G_{u_0,k}^+ \, e^{iku}\Big\|_{\Cprod{\alpha}{\beta}}^p \bigg]^{1/p}
 \lesssim 
 \sqrt{p} K^{\frac{\epsilon}{2}} K^{\alpha} \big\| \mathbf{1} \big\{ k \in \Z_K \big\} \big\|_{\ell_k^2}  
 \lesssim \sqrt{p} K^{\alpha+\frac{1}{2}+\frac{\epsilon}{2}}.
\end{align*}
By summing over $K\in \Dyadiclarge$ and using the $K^{-\frac{\epsilon}{2}}$-factor which is still at our disposal, this yields the desired estimate in \eqref{killing:eq-modulated-linear-U}. 
\end{proof}

\subsection{Tensor products of modulated linear waves}\label{section:killing-tensor}

In this subsection, we study the tensor products of translated
modulated linear waves and translated integrated modulated linear waves.

\begin{proposition}[The main tensor product]\label{killing:prop-main-tensor}
Let $K,M,\Nd\in \Dyadiclarge$ and let $\gamma \in (-s,0)$. For all $\epsilon>0$ and $p\geq 1$, it then holds that  
\begin{align}
&\E \bigg[ \sup_{u_0,u_1 \in \LambdaRR} \sup_{\substack{y,z\in \R}} \hspace{-1ex} 
\sup_{\substack{\Smod[L][+],\Smodtil[K][+]\in \Wfuv[s][s]\colon \\ \| \Smod[L][+] \|_{\Wuv[s][s][]},\| \Smodtil[K][+] \|_{\Wuv[s][s][]}\leq 1}}\hspace{-0.5ex}
\bigg\| 
\Theta^x_{y} \Big( \sum_{\ell \in \Z_L} \hspace{-1ex}\rhoND(\ell) \Smod[L][+][\ell] G_{u_0,\ell}^+ \frac{e^{i\ell u}}{i \ell} \Big) \otimes 
\Theta^x_{z} \Big( \sum_{k\in \Z_K} \hspace{-0.75ex} \rhoND(k) \Smodtil[K][+][k] G_{u_1,k}^+ \, e^{iku} \Big) \notag \\ 
&\hspace{8ex} - \mathbf{1}\big\{ K=L\big\} \mathbf{1} \big\{ u_0=u_1 \big\}  
\sum_{k\in \Z_K}  \delta^{ab} \rhoNDsquare(k) \Big( \big( \Theta^x_y \Smod[K][+][k] E_a \big) \otimes \big( \Theta^x_z \Smodtil[K][+][k] E_b \big) \Big) \frac{e^{ik(y-z)}}{(-ik)} 
\bigg\|_{\Cprod{\gamma}{s}}^p \bigg]^{1/p}  \notag \\
&\lesssim p R^\epsilon  \max\big( K, L \big)^{\gamma+\frac{1}{2}+\epsilon} L^{-\frac{1}{2}}. \label{killing:eq-main-tensor}
\end{align}
\end{proposition}

\begin{remark}\label{killing:rem-resonant}
In \eqref{killing:eq-main-tensor}, we subtract the term 
\begin{equation}\label{killing:eq-resonant}
\mathbf{1}\big\{ K=L\big\} \mathbf{1} \big\{ u_0=u_1 \big\}  
\sum_{k\in \Z_K} \delta^{ab} \rhoNDsquare(k) \Big( \big( \Theta^x_y \Smod[K][+][k] E_a \big) \otimes \big( \Theta^x_z \Smodtil[K][+][k] E_b \big) \Big) \frac{e^{ik(y-z)}}{(-ik)},
\end{equation}
which coincides with the resonant part of the second-order Gaussian chaos in \eqref{killing:eq-main-tensor}. 
The resonant term \eqref{killing:eq-resonant} will be examined in Proposition \ref{killing:prop-resonant} below. In particular, it will be linked to the Casimir $\Cas \in \frkg \otimes \frkg$ from Definition \ref{prelim:def-casimir} and the covariance function $\Cf$ from Definition \ref{ansatz:def-Killing}.  
\end{remark}

\begin{proof}[Proof of Proposition \ref{killing:prop-main-tensor}]
Using Lemma \ref{prelim:lem-maxima}, we may fix $u_0,u_1\in \LambdaR$. By using the  elementary identity 
$\Theta^x_y f \cdot \Theta^x_z g = \Theta^x_z \big( \Theta^x_{y-z} f \cdot g \big)$ and the translation-invariance of the $\Cprod{\gamma}{s}$-norm, we may further reduce to $z=0$.
 Using the decomposition of $(G^+_{u_0,m})_{m\in \Z}$ and $(G^+_{u_1,k})_{k\in \Z}$ into coordinates from \eqref{prelim:eq-Gm-decomposition}, we now decompose the second-order Gaussian chaos as 
\begin{align}
&\Big( \Theta^x_{y}  \sum_{\ell \in \Z_L} \rhoND(\ell) \Smod[L][+][\ell] G_{u_0,\ell}^+ \frac{e^{i\ell u}}{i \ell} \Big) \otimes  \Big( \sum_{k\in \Z_K} \rhoND(k) \Smodtil[K][+][k] G_{u_1,k}^+ \, e^{iku} \Big) \notag \\
=&  \sum_{\substack{ k \in \Z_{K}  \\ \ell \in \Z_{L} }} \rhoND(k)  \rhoND(\ell)
\Big( \Theta^x_y \Smod[L][+][\ell] G_{u_0,\ell}^+ \otimes \Smodtil[K][+][k] G_{u_1,k}^+ \Big) \frac{e^{-i\ell y}}{i \ell} e^{i (k+\ell) u}  \notag  \\ 
=&  \sum_{\substack{ k \in \Z_{K} \\ \ell \in \Z_{L}  }} \rhoND(k)  \rhoND(\ell)
\Big( \Theta^x_y \Smod[L][+][\ell]  E_{a} \otimes \Smodtil[K][+][k] E_{b} \Big) G_{u_0,\ell}^{+,a} G_{u_1,k}^{+,b} \frac{e^{-i\ell y}}{i \ell} e^{i(k+\ell )u}  \label{killing:eq-main-tensor-p1}. 
\end{align}
We then insert the chaos decomposition 
\begin{equation*}
G_{u_0,\ell}^{+,a} G_{u_1,k}^{+,b} = \lcol G_{u_0,\ell}^{+,a} G_{u_1,k}^{+,b} \rcol \hspace{-0.3ex}+ \mathbf{1}\big\{ u_0 = u_1 \big\} \delta^{a b} \delta_{k+\ell=0}  
\end{equation*}
from Subsection \ref{section:preliminaries-chaos} into \eqref{killing:eq-main-tensor-p1}, which yields the decomposition
\begin{align}
\eqref{killing:eq-main-tensor-p1} 
&= \sum_{\substack{ k \in \Z_{K} \\ \ell \in \Z_{L}}} \rhoND(k)  \rhoND(\ell)
 \Big( \Theta^x_y\Smod[L][+][\ell]  E_{a} \otimes \Smodtil[K][+][k] E_{b} \Big) \lcol G_{u_0,\ell}^{+,a} G_{u_1,k}^{+,b}\rcol \frac{e^{-i\ell y}}{i \ell} e^{i(k+\ell )u}  \label{killing:eq-main-tensor-p2}\\ 
&+  \mathbf{1}\big\{ u_0 = u_1 \big\} \sum_{\substack{k \in \Z_{K} \\ \ell \in \Z_{L}}} \rhoND(k)  \rhoND(\ell)
 \Big( \Theta^x_y\Smod[L][+][\ell]  E_{a} \otimes \Smodtil[K][+][k] E_{b} \Big)  \delta^{a b} \delta_{k+\ell=0}  \frac{e^{-i\ell y}}{i \ell} e^{i(k+\ell )u}. \label{killing:eq-main-tensor-p3}
\end{align}
Using the symmetry condition $\Smod[L][+][\ell]=\Smod[L][+][-\ell]$, which holds for all $\ell \in \Z_L$, the resonant term \eqref{killing:eq-main-tensor-p3} can be written as 
\begin{equation*}
\eqref{killing:eq-main-tensor-p3}
= \mathbf{1}\big\{ K=L\big\} \mathbf{1} \big\{ u_0=u_1 \big\}  
\sum_{k\in \Z_K} \delta^{ab} \rhoNDsquare(k) \Big( \big( \Theta^x_y \Smod[K][+][k] E_a \big) \otimes \big(  \Smodtil[K][+][k] E_b \big) \Big) \frac{e^{iky}}{(-ik)}. 
\end{equation*}
As a result, \eqref{killing:eq-main-tensor-p3} is canceled by the second summand in \eqref{killing:eq-main-tensor}, and it therefore only remains to prove that
\begin{equation}\label{killing:eq-main-tensor-p4}
\E \bigg[ \sup_{y\in \R} \, 
\sup_{\substack{\Smod[L][+]\in \Wfuv[s][s]\colon \\\| \Smod[L][+] \|_{\Wuv[s][s][]}\leq 1 }}
\sup_{\substack{\Smodtil[K][+]\in \Wfuv[s][s]\colon \\\| \Smodtil[K][+] \|_{\Wuv[s][s][]}\leq 1 }}\hspace{-0.5ex} 
\big\| \, \eqref{killing:eq-main-tensor-p2} \big\|_{\Cprod{\gamma}{s}}^p \bigg]^{1/p} \lesssim p \max(L,K)^{\gamma+\frac{1}{2}+\epsilon} L^{-\frac{1}{2}}. 
\end{equation} 
To this end, we first use Lemma \ref{chaos:lem-wc}, which yields that
\begin{equation*}
\max_{a,b=1,\hdots,\dim\frkg} \, \sup_{y\in \R}  \Big\| \rhoND(k) \rhoND(\ell) \Big( \big( \Theta^x_y\Smod[L][+][\ell]  E_{a} \big) \otimes \big( \Smodtil[K][+][k] E_{b}  \big) \Big)\Big\|_{\Wuv[s][s][](\Z^2)} \lesssim  
\big\| \Smod[L][+][\ell] \big\|_{\Wuv[s][s][](\Z)} \big\| \Smodtil[K][+][k] \big\|_{\Wuv[s][s][](\Z)}. 
\end{equation*}
Using our abstract chaos estimate (Proposition \ref{chaos:prop-chaos}), it then follows that
\begin{equation}\label{killing:eq-main-tensor-p5}
\begin{aligned}
&\E \bigg[ \sup_{y\in \R} \, 
\sup_{\substack{\Smod[L][+]\in \Wfuv[s][s]\colon \\\| \Smod[L][+] \|_{\Wuv[s][s][]}\leq 1 }}
\sup_{\substack{\Smodtil[K][+]\in \Wfuv[s][s]\colon \\\| \Smodtil[K][+] \|_{\Wuv[s][s][]}\leq 1 }}\hspace{-0.5ex} 
\big\| \, \eqref{killing:eq-main-tensor-p2} \big\|_{\Cprod{\gamma}{s}}^p \bigg]^{1/p} \\
\lesssim& \, \revision{p} \max(L,K)^{\epsilon} \max_{M \in \dyadic} M^{\gamma} \bigg\| \frac{1}{i\ell} \mathbf{1} \big\{ \,  k \in \Z_{K},  \ell \in \Z_{L} \big\} \, \mathbf{1} \big\{ |k+\ell|\sim M \big\} \bigg\|_{\ell^2_{k} \ell^2_{\ell}}. 
\end{aligned}
\end{equation}
From a direct calculation, it follows that 
\begin{equation}\label{killing:eq-main-tensor-p6}
\begin{aligned}
\bigg\| \frac{1}{i\ell} \mathbf{1} \big\{ \,  k \in \Z_{K},  \ell \in \Z_{L} \big\} \, \mathbf{1} \big\{ |k+\ell|\sim M \big\} \bigg\|_{\ell^2_{k} \ell^2_{\ell}}
\lesssim  L^{-1} L^{\frac{1}{2}} \min(K,M)^{\frac{1}{2}} \lesssim K^{\frac{1}{2}+\gamma} M^{-\gamma} L^{-\frac{1}{2}}. 
\end{aligned}
\end{equation}
\revision{In the last inequality in \eqref{killing:eq-main-tensor-p6}, we used that $\gamma \in (-s,0)\subseteq (-\frac{1}{2},0)$ and bounded the minimum of $K$ and $M$ using a product. By inserting \eqref{killing:eq-main-tensor-p6} into \eqref{killing:eq-main-tensor-p5}}, we then obtain the desired estimate \eqref{killing:eq-main-tensor-p4}. 
\end{proof} 

In the rest of this subsection, we record three different variants of Proposition \ref{killing:prop-main-tensor}. The first variant is phrased in terms of Cartesian rather than null coordinates and will be used to control traces of the modulated bilinear waves $\UN[K,M][+-]$ and $\VN[K,M][+-]$. 

\begin{lemma}[Tensor product estimate in Cartesian coordinates]\label{killing:lem-Cartesian}
Let $K,M,\Nd\in \Dyadiclarge$ and let $\gamma\in (-s,0)$. For all $\epsilon>0$ and $p\geq 1$, it then holds that  
\begin{equation}\label{killing:eq-Cartesian}
\begin{aligned}
&\E \bigg[  \sup_{u_0,v_0\in \LambdaRR} \sup_{\substack{y,z\in \R }}  
\sup_{\substack{\Smod[K][+]\in \Wfuv[s][s][] \colon \\ \| \Smod[K][+] \|_{\Wuv[s][s][]}\leq 1}}
\sup_{\substack{\Smod[M][-]\in \Wfuv[s][s][] \colon \\ \| \Smod[M][-] \|_{\Wuv[s][s][]}\leq 1}}
\Big\| 
\Theta^x_{y} \Big( \sum_{k\in \Z_K} \rhoND(k) \Smod[K][+][k](x-t,x+t) G_{u_0,k}^+ \, e^{ik(x-t)} \Big)  \\
&\hspace{5ex}\otimes 
\Theta^x_{z} \Big( \sum_{m\in \Z_M} \rhoND(m) \Smod[M][-][m](x-t,x+t) G_{v_0,m}^- \frac{e^{im(x+t)}}{im} \Big) 
\Big\|_{C_t^0 \C_x^{\gamma}}^p \bigg]^{\frac{1}{p}}  \\
\lesssim&\, p R^\epsilon \max(K,M)^{\gamma+\frac{1}{2}+\epsilon} M^{-\frac{1}{2}}. 
\end{aligned}
\end{equation}
\end{lemma}

\begin{proof}
Since  the argument in \eqref{killing:eq-Cartesian} contains tensor products of $G_{u_0,k}^+$ and $G_{v_0,m}^-$, which are probabilistically independent, there is no resonant part. Thus, \eqref{killing:eq-Cartesian} can be shown using the same estimate as for the non-resonant part \eqref{killing:eq-main-tensor-p2} in the proof of Proposition \ref{killing:prop-main-tensor} and we thus omit the details. 
\end{proof}

The second variant of Proposition \ref{killing:prop-main-tensor} is rather technical. 
At this stage of the argument, it is difficult to motivate the expression in \eqref{killing:eq-for-energy-increment} below, but it appears naturally in our estimates of the energy increment (see Lemma \ref{increment:lem-increment-object-I}). 

\begin{lemma}[Tensor product estimate for energy increment]\label{killing:lem-for-energy-increment}
Let $K,M,\Nd \in \Dyadiclarge$ and let $\gamma \in (-s,0)$. Furthermore, let $\iota \in \{\pm 1\}$ and let $c\in \R$. 
For all $\epsilon>0$ and $p\geq 1$, it then holds that  
\begin{equation}\label{killing:eq-for-energy-increment}
\begin{aligned}
&\E \bigg[ \sup_{u_0,v_0\in \LambdaRR}  \sup_{\substack{y,z\in \R }}  
\sup_{\substack{\Smod[K][+]\in \Wfuv[s][s][] \colon \\ \| \Smod[K][+] \|_{\Wuv[s][s][]}\leq 1}}
\sup_{\substack{\Smod[M][-]\in \Wfuv[s][s][] \colon \\ \| \Smod[M][-] \|_{\Wuv[s][s][]}\leq 1}}
\bigg\| \bigg( 
\Theta^x_{y} \Big( \sum_{k\in \Z_K} \rhoND(k) \Smod[K][+][k] G_{u_0,k}^+ \, e^{iku} \Big)  \\
&\hspace{33ex}\otimes 
\Theta^x_{z} \Big( \sum_{m\in \Z_M} \rhoND(m) \Smod[M][-][m]G_{v_0,m}^- \frac{e^{imv}}{im} \Big) 
\bigg)(u,\iota u+c)
\bigg\|_{\C_u^{\gamma}}^p \bigg]^{\frac{1}{p}}   \\
\lesssim&\, p R^\epsilon \max(K,M)^{\gamma+\frac{1}{2}+\epsilon} M^{-\frac{1}{2}}. 
\end{aligned}
\end{equation}
\end{lemma}

\begin{proof}
Similar as \eqref{killing:eq-Cartesian} from Lemma \ref{killing:lem-Cartesian}, \eqref{killing:eq-for-energy-increment} can be shown using the argument for the non-resonant part \eqref{killing:eq-main-tensor-p2} in the proof of Proposition \ref{killing:prop-main-tensor}. 
\end{proof}

The third variant of Proposition \ref{killing:prop-main-tensor} is much simpler than the original formulation, since
it does not involve modulation operators. The only new aspect is that it involves the Cartesian integral operator $\Int^x_{0\rightarrow x}$. In Section \ref{section:main}, it will be used to control
the probability of the Bourgain-Bulut event (see Corollary \ref{main:cor-bourgain-bulut-probability}).

\begin{lemma}[Tensor product estimate involving $\Int^x_{0\rightarrow x}$]\label{killing:lem-bourgain-bulut} 
Let $N\in \Dyadiclarge$, let $\Pbd$ be as in \eqref{prelim:eq-pbd}, and let $\CNbd$ be as in Definition \ref{jacobi:def-cf} below. 
Let $R\geq 1$, let $\coup>0$, and let $W^{(\Rscript,\coup)},W^{(\Rscript,\coup)}_0,W^{(\Rscript,\coup)}_1 \colon \bT_R \rightarrow \frkg$ be independent, $2\pi R$-periodic, $\frkg$-valued white noises at temperature $8\coup$. Then, it holds for all $p\geq 1$ that 
\begin{align}
&\E \bigg[ \sup_{y\in \R} \, \Big( \langle N y \rangle^{-2} 
 \Big\| P_{\leq N^{\deltap}}^x \Big( 
\Pbd W^{(\Rscript,\coup)} \otimes  (\Theta^x_y-1) P^x_{>1} \Int^x_{0\rightarrow x} \Pbd W^{(\Rscript,\coup)}  - 8\coup \CNbd(y) \Cas \Big) \Big\|_{L_x^\infty(\bT_R)} \Big)^p \bigg]^{\frac{1}{p}} \notag \\
\lesssim&\, p  R^\epsilon \coup N^{-\frac{1}{2}+\frac{\delta}{2}+2\delta_1}
\label{killing:eq-bourgain-bulut-e1}
\end{align}
and 
\begin{equation}\label{killing:eq-bourgain-bulut-e2}
\begin{aligned}
&\,\E \bigg[ \sup_{y\in \R} \, \Big( \langle N y \rangle^{-2}
 \Big\| P_{\leq N^{\deltap}}^x \Big( 
\Pbd W^{(\Rscript,\coup)}_0 \otimes (\Theta^x_y-1) P^x_{>1} \Int^x_{0\rightarrow x} \Pbd W^{(\Rscript,\coup)}_1 \Big) \Big\|_{L_x^\infty(\bT_R)}\Big)^p \bigg]^{\frac{1}{p}} \\
\lesssim&\, p R^\epsilon \coup N^{-\frac{1}{2}+\frac{\delta}{2}+2\delta_1}.
\end{aligned}
\end{equation}
\end{lemma}

\begin{proof}
Since \eqref{killing:eq-bourgain-bulut-e1} and \eqref{killing:eq-bourgain-bulut-e2} follow from the representation of white noise in Lemma \ref{prelim:lem-white-noise-representation} and minor modifications of the proof of Proposition \ref{killing:prop-main-tensor}, we only sketch the argument. As in the proof of Proposition \ref{killing:prop-main-tensor}, we first perform a decomposition into non-resonant and resonant parts. The non-resonant parts in \eqref{killing:eq-bourgain-bulut-e1} and \eqref{killing:eq-bourgain-bulut-e2} can be controlled by first using the $P_{\leq N^{\deltap}}^x$-operator to move from $L_x^\infty$ to $\C_x^{-s+\eta}$, which costs a factor of $N^{\deltap (s-\eta)}\sim N^{\frac{\delta}{2}+\mathcal{O}(\delta_2)}$, and then arguing as in the proof of \eqref{killing:eq-main-tensor-p4}. Due to the independence of $W^{(\Rscript,\coup)}_0$ and $W^{(\Rscript,\coup)}_1$, there is no resonant part in 
\eqref{killing:eq-bourgain-bulut-e2}, and therefore it remains to treat the resonant part in \eqref{killing:eq-bourgain-bulut-e1}. Up to simple error terms, which stem from the $\psix$-functions in Lemma \ref{prelim:lem-white-noise-representation} and integration by parts, the resonant part in \eqref{killing:eq-bourgain-bulut-e1} is given by 
\begin{equation}\label{killing:eq-bourgain-bulut-p}
\begin{aligned}
&\, 8 \coup \, P_{\leq N^{\deltap}}^x 
\bigg( \sum_{n\in \Z} \sum_{\substack{K,L\in \dyadic\colon \\ N^{1-2\delta_1}\leq K,L \leq N}} \hspace{-4ex}
\rho_{K}(n) \rho_L(n) e^{-inx} \big( \Theta^x_y - 1 \big) \Big( \frac{e^{inx}-1}{in} \Big) - \CNbd(y) \bigg) \Cas \\
=&\, 8 \coup \, 
\bigg( \sum_{n\in \Z} \sum_{\substack{K,L\in \dyadic\colon \\ N^{1-2\delta_1}\leq K,L \leq N}} \hspace{-4ex}
\rho_{K}(n) \rho_L(n) \frac{e^{-iny}}{in} - \CNbd(y) \bigg) \Cas.
\end{aligned}
\end{equation}
Using a minor variant of \eqref{ansatz:eq-renormalization-4}-\eqref{ansatz:eq-renormalization-6} from Lemma \ref{ansatz:lem-renormalization}
and Definition \ref{jacobi:def-cf}, it follows that
\begin{equation*}
\big| \eqref{killing:eq-bourgain-bulut-p}\big| \lesssim \coup |y| \lesssim \coup N^{-1} \langle N y \rangle,
\end{equation*}
which is acceptable. 
\end{proof}

\subsection{Proof of Proposition \ref{killing:prop-probabilistic-hypothesis}}\label{section:killing-proof} 
In order to prove Proposition \ref{killing:prop-probabilistic-hypothesis}, we need to show that the estimates from Hypothesis \ref{hypothesis:probabilistic} are satisfied on an event satisfying the estimate \eqref{killing:eq-probabilistic-hypothesis}. 
This follows from our earlier estimates in Subsection \ref{section:killing-linear} and Subsection \ref{section:killing-tensor}
as well as simpler estimates from Section \ref{section:preliminaries} and Section \ref{section:appendix-auxiliary}.
In the following, we reference the results corresponding to Hypothesis \ref{hypothesis:probabilistic}.\ref{ansatz:item-hypothesis-crude}-\ref{ansatz:item-hypothesis-for-energy-increment} separately. 
\begin{itemize}
\item[\ref{ansatz:item-hypothesis-crude}:] This follows from Lemma \ref{prelim:lem-maxima}.
\item[\ref{ansatz:item-hypothesis-sharp}:] This follows from Lemma \ref{prelim:lem-commutator-sharp}.
\item[\ref{ansatz:item-hypothesis-linear}:] This follows from Lemma \ref{killing:lem-regularity-modulated-linear}.
\item[\ref{ansatz:item-hypothesis-tensor-product}:] This follows from Proposition \ref{killing:prop-main-tensor}.
\item[\ref{ansatz:item-hypothesis-trace}:] This follows from Lemma \ref{killing:lem-Cartesian}.
\item[\ref{ansatz:item-hypothesis-for-energy-increment}:] This follows from Lemma \ref{killing:lem-for-energy-increment}. 
\end{itemize}

\subsection{The resonant term}\label{section:killing-resonant}

In this subsection, we analyze the resonant term from Proposition \ref{killing:prop-main-tensor} and Remark \ref{killing:rem-resonant}. As part of this analysis, we encounter the Casimir from Definition \ref{prelim:def-casimir}, which is a central object in this article.

\begin{proposition}[Resonant tensor]\label{killing:prop-resonant}
For all $K,\Nd\in \Dyadiclarge$ and $\Smod[K][+],\Smodtil[K][+]\in \Wfuv[s][s]$, it holds that
\begin{equation}\label{killing:eq-resonant-e1}
\bigg\| \sum_{k\in \Z_K} \delta^{ab} \rhoNDsquare(k) \Big( \big( \Theta^x_y \Smod[K][+][k] E_a \big) \otimes \big( \Theta^x_z \Smodtil[K][+][k] E_b \big) \Big)
\frac{e^{ik(y-z)}}{(-ik)} \bigg\|_{\Cprod{s}{s}} 
\lesssim K |y-z| \big\| \Smod[K][+] \big\|_{\Wuv[s][s]} \big\| \Smodtil[K][+] \big\|_{\Wuv[s][s]}.
\end{equation}
Furthermore, it also holds that 
\begin{equation}\label{killing:eq-resonant-e2}
\begin{aligned}
&\,\bigg\| \sum_{k\in \Z_K}  \delta^{ab}  \rhoNDsquare(k) \Big( \big( \Theta^x_y \Smod[K][+][k] E_a \big) \otimes \big( \Theta^x_z \Smod[K][+][k] E_b \big) \Big)
\frac{e^{ik(y-z)}}{(-ik)} - \Cf^{(\Nd)}_K(y-z) \Cas \bigg\|_{\Cprod{s}{s}} \\
\lesssim&\, \Big\| \Smod[K][+][k]  \big( \Smod[K][+][k] \big)^\ast - \Id_\frkg \Big\|_{\Wuv[s][s]} + |y-z| + 
K^{1-\delta+\vartheta} |y-z| \big\| \Smod[K][+] \big\|_{\Wuv[s][s]}^2 
\end{aligned} 
\end{equation}
\end{proposition}

\begin{remark}[Almost orthogonality]
As will be shown in Section \ref{section:modulation}, and as was previously mentioned in Subsection \ref{section:overview-ansatz} and Subsection \ref{section:ansatz}, the modulation operators stemming from our modulation equations (Definition \ref{ansatz:def-modulation-equations}) are almost orthogonal. As a result, the first summand in \eqref{killing:eq-resonant-e2} will turn out to be small.
\end{remark}

\begin{proof} 
Due to Definition \ref{chaos:def-wc}, we have for all $k\in \Z_K$ that $\Smod[K][+][k]=\Smod[K][+][-k]$ and $\Smodtil[K][+][k]=\Smodtil[K][+][-k]$. We can therefore replace the complex exponential $(-ik)^{-1} e^{ik(y-z)}$ in \eqref{killing:eq-resonant-e1} with 
$-k^{-1}\sin( k(y-z))$. Together with the estimate $|\sin(\xi)|\leq |\xi|$,  it follows that
\begin{align*}
\textup{LHS of } \eqref{killing:eq-resonant-e1}
&\lesssim 
\sum_{k\in \Z_K} \Big\| \delta^{ab} \Big( \big( \Theta^x_y \Smod[K][+][k] E_a \big) \otimes \big( \Theta^x_z \Smod[K][+][k] E_b \big) \Big) \Big\|_{\Cprod{s}{s}} \bigg| \frac{\sin(k(y-z))}{k}\bigg| \\ 
&\lesssim \sum_{k\in \Z_K} \big\| \Smod[K][+][k] \big\|_{\Cprod{s}{s}} \big\| \Smodtil[K][+][k] \big\|_{\Cprod{s}{s}} |y-z| \\
&\lesssim K |y-z| \sup_{k\in \Z_K}\big\| \Smod[K][+][k] \big\|_{\Cprod{s}{s}} \sup_{k\in \Z_K} \big\| \Smodtil[K][+][k] \big\|_{\Cprod{s}{s}}.
\end{align*}
Due to the definition of the $\Wuv[s][s]$-norm, this directly implies the first estimate \eqref{killing:eq-resonant-e1}. We now turn to the second estimate \eqref{killing:eq-resonant-e2}. To this end, we first use Lemma \ref{ansatz:lem-renormalization}, which implies
\begin{equation*}
\bigg| \Cf^{(\Nd)}_K(y-z) - \bigg( - \sum_{k\in \Z_K} \rhoNDsquare(k) \frac{\sin(k(y-z))}{k} \bigg) \bigg| \lesssim |y-z|. 
\end{equation*}

Using Definition \ref{prelim:def-casimir} and Lemma \ref{prelim:lem-casimir}, we then write 
\begin{align}
 &-\sum_{k\in \Z_K}\rhoNDsquare(k)   \delta^{ab} \Big( \big( \Theta^x_y \Smod[K][+][k] E_a \big) \otimes \big( \Theta^x_z \Smod[K][+][k] E_b \big) \Big)
\frac{\sin(k(y-z))}{k}
+  \sum_{k\in \Z_K} \rhoNDsquare(k) \frac{\sin(k(y-z))}{k} \Cas \notag \\
=&\, - \sum_{k\in \Z_K} \rhoNDsquare(k)  \delta^{ab}  \Big( \big( \Theta^x_y \Smod[K][+][k] E_a \big) \otimes \big( \Theta^x_z \Smod[K][+][k] E_b \big)  - E_a \otimes E_b \Big)\frac{\sin(k(y-z))}{k} \notag \\ 
=&\, - \sum_{k\in \Z_K} \rhoNDsquare(k) \Big( \big(\Theta^x_y \Smod[K][+][k]\big) \big(\Theta^x_z \Smod[K][+][k] \big)^\ast - \Id_\frkg\Big)^{\hspace{-0.2ex}ab}  
\big(E_a \otimes E_b) \, \frac{\sin(k(y-z))}{k}.\label{killing:eq-resonant-p0}
\end{align}
Together with the estimate $|\sin(\xi)|\leq 1$ and Definition \ref{chaos:def-wc}, it follows that 
\begin{align*}
\big\| \eqref{killing:eq-resonant-p0} \big\|_{\Cprod{s}{s}}
&\lesssim K^{-1} \sum_{k\in \Z_K}  \Big\| \big(\Theta^x_y \Smod[K][+][k]\big) \big(\Theta^x_z \Smod[K][+][k] \big)^\ast - \Id_\frkg \Big\|_{\Cprod{s}{s}} \\ 
&\lesssim \sup_{k\in \Z_K} \Big\| \big(\Theta^x_y \Smod[K][+][k]\big) \big(\Theta^x_z \Smod[K][+][k] \big)^\ast - \Id_\frkg \Big\|_{\Cprod{s}{s}} \\
&\lesssim \Big\|  \big(\Theta^x_y \Smod[K][+][k]\big) \big(\Theta^x_z \Smod[K][+][k] \big)^\ast - \Id_\frkg \Big\|_{\Wuv[s][s]}. 
\end{align*}
Using the translation-invariance of the $\Wuv[s][s]$-norm, we can further estimate
\begin{align}
 \Big\|  \big(\Theta^x_y \Smod[K][+][k]\big) \big(\Theta^x_z \Smod[K][+][k] \big)^\ast - \Id_\frkg \Big\|_{\Wuv[s][s]} 
=&\, \Big\|  \big(\Theta^x_{y-z} \Smod[K][+][k]\big) \big(\Smod[K][+][k] \big)^\ast - \Id_\frkg \Big\|_{\Wuv[s][s]} \notag \\ 
=&\, \Big\|  \big( \Smod[K][+][k]\big) \big(\Smod[K][+][k] \big)^\ast - \Id_\frkg \Big\|_{\Wuv[s][s]}
\hspace{-1ex}+ \Big\|   \big(\Theta^x_{y-z} \Smod[K][+][k]- \Smod[K][+][k]\big) \big(\Smod[K][+][k] \big)^\ast  \Big\|_{\Wuv[s][s]}. \label{killing:eq-resonant-p1}
\end{align}
The first summand in \eqref{killing:eq-resonant-p1} is as on the right-hand side of \eqref{killing:eq-resonant-e2}. 
To control the second summand in \eqref{killing:eq-resonant-p1}, we first recall from Definition \ref{killing:def-FW} that $\Smod[K][+][k]$ is supported on frequencies $\lesssim K^{1-\delta+\vartheta}$. Using Lemma \ref{chaos:lem-wc}, the second summand can then be estimated by 
\begin{equation}\label{killing:eq-resonant-p2}
\Big\|   \big(\Theta^x_{y-z} \Smod[K][+][k]- \Smod[K][+][k]\big) \big(\Smod[K][+][k] \big)^\ast  \Big\|_{\Wuv[s][s]}
\lesssim \Big\|   \Theta^x_{y-z} \Smod[K][+][k]- \Smod[K][+][k] \Big\|_{\Wuv[s][s]} \Big\| \Smod[K][+][k]  \Big\|_{\Wuv[s][s]}
\lesssim K^{1-\delta+\vartheta} |y-z| \Big\| \Smod[K][+][k]  \Big\|_{\Wuv[s][s]}^2,
\end{equation}
which completes the proof of \eqref{killing:eq-resonant-e2}.
\end{proof}

\subsection{Applications}\label{section:killing-consequences} 
In this subsection, we obtain several estimates for interactions of modulated objects. In order to make it easier to use these estimates in later sections, we will not state the estimates as moment bounds (as in Lemma \ref{killing:lem-regularity-modulated-linear} or Proposition \ref{killing:prop-main-tensor}) but rather as consequences of the post-modulation hypothesis (Hypothesis \ref{hypothesis:post}), which is a pointwise statement. 

From a notational perspective, it is convenient to make the following definition. 

\begin{definition}[Wick-ordered quadratic tensor product]\label{killing:def-Wick} 
Let $K,L,N,\Nd \in \Dyadiclarge$ and let $y,z \in \R$. Then, we define
\begin{equation}\label{killing:eq-wick-d1}
\begin{aligned}
&\biglcol \,  \Theta^x_{y} \IUN[L][+] \otimes \Theta^x_{z} \UN[K][+] \bigrcol \\
\overset{\textup{def}}{=}&  \,  \Theta^x_{y} \IUN[L][+] \otimes \Theta^x_{z} \UN[K][+] - \coup \mathbf{1}\big\{ K = L\big\} \Cf^{(\Ndscript)}_{K}(y-z) \Cas
\end{aligned}
\end{equation}
and 
\begin{equation}\label{killing:eq-wick-d2}
\begin{aligned}
&\biglcol \,  \Theta^x_{y} P_{\leq N}^x \IUN[L][+] \otimes \Theta^x_{z} P_{\leq N}^x \UN[K][+] \bigrcol \\
\overset{\textup{def}}{=}&  \,  \Theta^x_{y} P_{\leq N}^x \IUN[L][+] \otimes \Theta^x_{z} P_{\leq N}^x \UN[K][+] - \coup \mathbf{1}\big\{ K = L\big\} \Cf^{(\Ncs)}_{K}(y-z) \Cas.
\end{aligned}
\end{equation}
\end{definition}
The $\coup$-factors in \eqref{killing:eq-wick-d1} and \eqref{killing:eq-wick-d2} are due to the $\hcoup$-factors in both Definition \ref{ansatz:def-modulated-linear} and Definition \ref{ansatz:def-integrated-modulated-waves}. 
We note that the only two differences between \eqref{killing:eq-wick-d1} and \eqref{killing:eq-wick-d2} are the two $P_{\leq N}^x$-operators and that \eqref{killing:eq-wick-d2} contains $\Cf^{(\Ncs)}$ rather than $\Cf^{(\Ndscript)}$. 

Equipped with Definition \ref{killing:def-Wick}, we can now state and prove the following lemma.

\begin{lemma}[Tensor product of modulated linear waves and their integrals]\label{killing:lem-tensor-modulated-linear}
Let the post-modulation hypothesis (Hypothesis \ref{hypothesis:post}) be satisfied, let $K,L\in \Dyadiclarge$, and let $\gamma \in (-s,0)$. Furthermore, let $h\in (0,\infty)$ and let $y,z\in \R$ satisfy $|y|,|z|\leq h$. Then, it holds that
\begin{align}
 \Big\| \hspace{-0.5ex} \biglcol \,   \Theta^x_{y} \IUN[L][+] \hspace{-0.5ex} \otimes \Theta^x_{z}  \UN[K][+] \bigrcol  \hspace{-0.5ex}\Big\|_{\Cprod{r-1}{s}} 
&\lesssim \hspace{-0.5ex} \Big( \hspace{-0.5ex}\max\big( K, L \big)^{r-\frac{1}{2}+\eta} L^{-\frac{1}{2}} \hspace{-0.25ex} +\hspace{-0.5ex} \mathbf{1} \big\{ K = L \big\} N^{-\delta+\vartheta} N h \Big) \Dc^2, \label{killing:eq-tensor-truncated-1}\\ 
 \Big\| \hspace{-0.5ex} \biglcol \,   \Theta^x_{y} P_{\leq N}^x \IUN[L][+] \hspace{-0.5ex} \otimes \Theta^x_{z} P_{\leq N}^x \UN[K][+] \bigrcol  \hspace{-0.5ex}\Big\|_{\Cprod{r-1}{s}} 
&\lesssim \hspace{-0.5ex} \Big( \hspace{-0.5ex}\max\big( K, L \big)^{r-\frac{1}{2}+\eta} L^{-\frac{1}{2}}  \hspace{-0.25ex} +\hspace{-0.5ex} \mathbf{1} \big\{ K = L \big\} N^{-\delta+\vartheta} \langle N h \rangle \Big) \Dc^2. \label{killing:eq-tensor-truncated-2}
\end{align}
\end{lemma}

\begin{proof} 
We first prove \eqref{killing:eq-tensor-truncated-1}.
To simplify the notation, we write Definition \ref{ansatz:def-modulated-linear} and Definition \ref{ansatz:def-integrated-modulated-waves} as 
\begin{equation}\label{killing:eq-tensor-truncated-p2}
\IUN[L][+] = \sum_{u_0\in \LambdaRR} \psiRuL \IUN[u_0,L][+] 
\qquad \text{and} \qquad 
\UN[K][+] = \sum_{u_1 \in \LambdaRR} \psiRuKone \UN[u_1,K][+],
\end{equation}
where 
\begin{align}\label{killing:eq-tensor-truncated-p3}
\IUN[u_0,L][+]= \hcoup \sum_{\ell \in \Z_L} \rhoND(\ell)  \SN[L][+][\ell] G^+_{u_0,\ell} \frac{e^{i\ell u}}{i\ell}
\qquad \text{and} \qquad 
\UN[u_1,K][+] = \hcoup \sum_{k\in \Z_K} \rhoND(k) \SN[K][+][k] G^+_{u_1,k} e^{iku}.
\end{align}
Using \eqref{killing:eq-tensor-truncated-p2} and using the properties of $(\psiRu)_{u_0\in \LambdaRR}$ from
\eqref{prelim:eq-psi-properties}, we then obtain that 
\begin{align}
 &\hspace{2ex} \biglcol \,   \Theta^x_{y} \IUN[L][+] \otimes \Theta^x_{z}  \UN[K][+] \bigrcol \notag \\
 &= \sum_{u_0,u_1 \in \LambdaRR} \bigg(  
\big( \Theta^x_{y} \psiRuL \big) \big( \Theta^x_{z} \psiRuKone \big)   
\big( \Theta^x_{y} \IUN[u_0,L][+] \otimes \Theta^x_{z}  \UN[u_1,K][+]\big) \notag \\ 
&\hspace{15ex}- \coup \psiRu \psiRuone \mathbf{1}\big\{ u_0 = u_1 \big\} \mathbf{1}\big\{ K=L\big\} \Cf^{(\Ndscript)}_K(y-z) \Cas   \bigg) \notag \\ 
&= \sum_{u_0,u_1 \in \LambdaRR} 
\bigg( 
\big( \Theta^x_{y} \psiRuL \big) \big( \Theta^x_{z} \psiRuKone \big) \label{killing:eq-tensor-truncated-p4} \\ 
&\hspace{8ex} 
\times \Big( \Theta^x_{y} \IUN[u_0,L][+] \otimes \Theta^x_{z}  \UN[u_1,K][+] 
- \coup  \mathbf{1}\big\{ u_0 = u_1 \big\} \mathbf{1}\big\{ K=L\big\} \Cf^{(\Ndscript)}_K(y-z) \Cas \Big) \bigg) \notag \\ 
&+ \coup  \mathbf{1}\big\{ K=L\big\} \Cf^{(\Ndscript)}_K(y-z) 
\sum_{u_0,u_1 \in \LambdaRR} \bigg( \mathbf{1}\big\{ u_0 = u_1 \big\}   
\Big( \big( \Theta^x_{y} \psiRuL \big) \big( \Theta^x_{z} \psiRuKone \big) 
- \psiRu \psiRuone \Big) \bigg)  \Cas . \label{killing:eq-tensor-truncated-p5}
\end{align}
We now estimate \eqref{killing:eq-tensor-truncated-p4} and \eqref{killing:eq-tensor-truncated-p5} separately. 
Using Lemma \ref{prelim:lem-psi-sum}, we first obtain that
\begin{equation}\label{killing:eq-tensor-truncated-q4}
\begin{aligned}
&\, \big\| \eqref{killing:eq-tensor-truncated-p4} \big\|_{\Cprod{r-1}{s}} \\ 
\lesssim&\,    \sup_{u_0,u_1\in \LambdaRR} \Big\| \Theta^x_{y} \IUN[u_0,L][+] \otimes \Theta^x_{z}  \UN[u_1,K][+] - \coup  \mathbf{1}\big\{ u_0 = u_1 \big\} \mathbf{1}\big\{ K=L\big\} \Cf^{(\Ndscript)}_K(y-z) \Cas \Big\|_{\Cprod{r-1}{s}}.  
\end{aligned}
\end{equation}
We now estimate \eqref{killing:eq-tensor-truncated-q4} in two different ways. In the regime $K<N^{1-\delta+\vartheta}$, we make use of Lemma \ref{ansatz:lem-renormalization}, Hypothesis \ref{hypothesis:probabilistic}.\ref{ansatz:item-hypothesis-tensor-product}, the estimate \eqref{killing:eq-resonant-e1} from Proposition \ref{killing:prop-resonant}, and Hypothesis \ref{hypothesis:pre}.\ref{ansatz:item-pre-modulation}, which yield that 
\begin{equation}\label{killing:eq-tensor-truncated-q5}
\begin{aligned}
\eqref{killing:eq-tensor-truncated-q4} 
\lesssim&\, \Big\| \Theta^x_{y} \IUN[u_0,L][+] \otimes \Theta^x_{z}  \UN[u_1,K][+]  \Big\|_{\Cprod{r-1}{s}}
+ \coup \mathbf{1} \big\{ K= L \big\} \, \Big| \Cf^{(\Ndscript)}_K(y-z) \Big| \allowdisplaybreaks[3] \\ 
\lesssim&\, \coup \, \Ac^2 \Big( \max(K,L)^{r-\frac{1}{2}+\eta} L^{-\frac{1}{2}} + \mathbf{1} \big\{ K=L \big\} K \big| y-z\big| \Big) \big\| \SN[L][+] \big\|_{\Wuv[s][s]} \big\| \SN[K][+] \big\|_{\Wuv[s][s]}  \\
&\hspace{0ex}+ \coup \mathbf{1} \big\{ K=L \big\} K \big| y-z\big|\\ 
\lesssim&\, \Big( \max(K,L)^{r-\frac{1}{2}+\eta} L^{-\frac{1}{2}} + \mathbf{1} \big\{ K=L \big\} N^{1-\delta+\vartheta} h \Big) \Dc^2.
\end{aligned}
\end{equation}
In the regime $K\geq N^{1-\delta+\vartheta}$, we make use of Hypothesis \ref{hypothesis:probabilistic}.\ref{ansatz:item-hypothesis-tensor-product}, the estimate \eqref{killing:eq-resonant-e2} from Proposition \ref{killing:prop-resonant}, and Hypothesis \ref{hypothesis:post}.\ref{ansatz:item-post-orthogonality}, which yield that 
\begin{equation}\label{killing:eq-tensor-truncated-q6}
\begin{aligned}
\eqref{killing:eq-tensor-truncated-q4} 
\lesssim&\, \coup \Ac^2 \Big( \max(K,L)^{r-\frac{1}{2}+\eta} L^{-\frac{1}{2}} + \mathbf{1} \big\{ K=L \big\} K^{1-\delta+\vartheta} \big| y-z\big| \Big) \big\| \SN[L][+] \big\|_{\Wuv[s][s]} \big\| \SN[K][+] \big\|_{\Wuv[s][s]} \\ 
&+ \coup \Ac^2 \mathbf{1}\big\{ K=L \big\} \Big\| \SN[K][+][k] \big( \SN[K][+][k] \big)^\ast - \Id_\frkg \Big\|_{\Wuv[s][s]} 
+\coup \Ac^2 \mathbf{1}\big\{ K=L \big\} |y-z| \allowdisplaybreaks[3] \\
\lesssim&\, \Big( \max(K,L)^{r-\frac{1}{2}+\eta} L^{-\frac{1}{2}} + \mathbf{1} \big\{ K=L \big\} N^{1-\delta+\vartheta}  h \Big) \Dc^2 + \mathbf{1}\big\{ K=L \big\} K^{-100} \Dc^2 
+ \mathbf{1}\big\{ K=L \big\} h \Dc^2\allowdisplaybreaks[3]  \\ 
\lesssim&\, \Big( \max(K,L)^{r-\frac{1}{2}+\eta} L^{-\frac{1}{2}} + \mathbf{1} \big\{ K=L \big\} N^{1-\delta+\vartheta} h \Big) \Dc^2 .
\end{aligned}
\end{equation}
By combining \eqref{killing:eq-tensor-truncated-q4}, \eqref{killing:eq-tensor-truncated-q5}, and \eqref{killing:eq-tensor-truncated-q6}, it follows that \eqref{killing:eq-tensor-truncated-p4} makes an acceptable contribution to \eqref{killing:eq-tensor-truncated-1}. In order to control \eqref{killing:eq-tensor-truncated-p5}, we first estimate 
\begin{align*}
 \mathbf{1}\big\{ K=L \big\} \sum_{u_0,u_1\in \LambdaRR} \mathbf{1}\big\{ u_0 = u_1 \big\} 
\Big| \big( \Theta^x_{y} \psiRuL \big) \big( \Theta^x_{z} \psiRuKone \big) 
- \psiRu \psiRuone \Big| 
\lesssim \mathbf{1}\big\{ K=L \big\} \big( h + K^{-100} \big). 
\end{align*} 
Together with Lemma \ref{ansatz:lem-renormalization}, it then follows that
\begin{align*}
\big| \eqref{killing:eq-tensor-truncated-p5}\big| 
&\lesssim  \coup \mathbf{1}\big\{ K=L\big\} \min \big( 1, K h \big) \big( h +K^{-100} \big) \lesssim \coup  \mathbf{1}\big\{ K=L\big\} h, 
\end{align*}
which is an acceptable contribution to \eqref{killing:eq-tensor-truncated-1}. This completes the proof of \eqref{killing:eq-tensor-truncated-1} and it only remains to prove \eqref{killing:eq-tensor-truncated-2}. Due to Definition \ref{ansatz:def-Killing}, it holds that 
\begin{equation*}
\Cf^{(\Ncs)}_K(y-z) = \int_{\R^2} \dy^\prime \dz^\prime \, \widecheck{\rho}_{\leq N}(y^\prime) \widecheck{\rho}_{\leq N}(z^\prime)  \Cf^{(\Ndscript)}_K(y+y^\prime-z-z^\prime).
\end{equation*}
Together with the definition of $P_{\leq N}^x$ and Definition \ref{killing:def-Wick}, we then obtain that
\begin{equation}\label{killing:eq-tensor-truncated-p6}
\begin{aligned}
&\biglcol \, \Theta^x_{y} P_{\leq N}^x \IUN[L][+] \otimes \Theta^x_{z} P_{\leq N}^x \UN[K][+]  \bigrcol  \\
=& \, \int_{\R^2} \dy^\prime \dz^\prime \, \Big( \widecheck{\rho}_{\leq N}(y^\prime) \widecheck{\rho}_{\leq N}(z^\prime) \biglcol \,  
\Theta^x_{y+y^\prime} \IUN[L][+] \otimes \Theta^x_{z+z^\prime}  \UN[K][+] \bigrcol \Big) . 
\end{aligned}
\end{equation}
Using \eqref{killing:eq-tensor-truncated-p6}, the second estimate \eqref{killing:eq-tensor-truncated-2} can then easily be derived from the first estimate \eqref{killing:eq-tensor-truncated-1}. 
\end{proof}

\begin{proposition}[Killing-renormalized iterated Lie-bracket]\label{killing:prop-quadratic} 
Let the post-modulation hypothesis (Hypothesis \ref{hypothesis:post}) be satisfied, let $K,L\in \Dyadiclarge$, and let $\gamma\in (-s,0)$. For all $y\in \R$, it then holds that
\begin{equation}\label{killing:eq-product}
\begin{aligned}
&\sup_{\substack{E \in \frkg \colon \\ \| E \|_{\frkg} \leq 1}}
\Big\| \Big[ P_{\leq N}^x \UN[K][+], \Big[  \Theta^x_y P_{\leq N}^x\IUN[L][+] , E \Big] \Big] -  \coup \mathbf{1}\big\{ K = L \big\} \Cf^{(\Ncs)}_{K}(y) \Kil(E)  \Big\|_{\Cprod{\gamma}{s}} \\
\lesssim& \,  \Big( \max(K,L)^{\gamma+\frac{1}{2}+\eta} L^{-\frac{1}{2}} + \mathbf{1}\big\{ K=L \big\} N^{-\delta+\vartheta} \langle N y \rangle \Big) \Dc^2. 
\end{aligned}
\end{equation}
\end{proposition}

\begin{proof}[Proof of Proposition \ref{killing:prop-quadratic}:] 
In order to prove \eqref{killing:eq-product}, we only have to reduce it to an estimate of the Wick-ordered tensor product from Definition \ref{killing:def-Wick}. To this end, we write 
\begin{align}
 \Big[ P_{\leq N}^x \UN[K][+], \Big[  \Theta^x_y P_{\leq N}^x\IUN[L][+] , E \Big] \Big]
=& -  \Big[ \big(   \Id_\frkg \otimes \ad(E) \big) 
\Big(   P_{\leq N}^x \UN[K][+] \otimes   \Theta^x_y P_{\leq N}^x\IUN[L][+]  \Big) \Big]_\frkg \notag \\
=& -  \Big[ \big(   \Id_\frkg \otimes \ad(E) \big) 
\Big(   \biglcol \, P_{\leq N}^x \UN[K][+] \otimes   \Theta^x_y P_{\leq N}^x\IUN[L][+]  \bigrcol \Big) \Big]_\frkg  \label{killing:eq-quadratic-p1} \\
-& \coup \, \mathbf{1} \big\{ K = L \big\} 
\Cf^{(\Ncs)}_{K}(y) 
\Big[  \big(  \Id_\frkg \otimes \ad(E)  \big) \Cas \Big]_\frkg.  \label{killing:eq-quadratic-p2}
\end{align}
Due to Lemma \ref{killing:lem-tensor-modulated-linear}, the contribution of \eqref{killing:eq-quadratic-p1} is acceptable. Using Definition \ref{prelim:def-casimir} and Definition \ref{prelim:def-Killing}, it holds that
\begin{equation*}
-\coup \Big[  \big(  \Id_\frkg \otimes \ad(E)  \big) \Cas \Big]_\frkg 
= -\coup \Big[ E_a \otimes \ad(E) E^a \Big]_\frkg= - \coup \Big[ E_a, \Big[ E ,E^a \Big] \Big] = \coup \Kil(E). 
\end{equation*}
Thus, \eqref{killing:eq-quadratic-p2} is exactly canceled by the Killing-renormalization.
\end{proof}

In the following corollary, we show that the quadratic tensor estimate implies several estimates for higher-order tensors.

\begin{corollary}[Higher-order tensors]\label{killing:cor-tensors}
Let the post-modulation hypothesis (Hypothesis \ref{hypothesis:post}) be satisfied and let $h \in (0,\infty)$. Then, we have the following estimates: 
\begin{enumerate}[label=(\roman*)]
\item (Cubic tensor) For all $K,L,M\in \Dyadiclarge$, it holds that
\begin{align}
& \sup_{\substack{y_1,y_2,y_3\in \R \colon \\ |y_1|,|y_2|,|y_3|\leq h}} 
 \Big\| \Big( \biglcol \, \Theta^x_{y_1} P_{\leq N}^x \UN[K][+]  \otimes \Theta^x_{y_2} P_{\leq N}^x \IUN[L][+]   \bigrcol \Big)  \otimes  \Theta^x_{y_3} \VN[M][-] \Big\|_{\Cprod{r-1}{r-1}}  \notag \\
\lesssim&\,  \Big( \max\big( K, L \big)^{r-s} L^{-\frac{1}{2}} + \mathbf{1} \big\{ K = L \big\} N^{-\delta+\vartheta} \langle N h \rangle \Big) M^{r-s} \Dc^3.  \label{killing:eq-tensors-cubic-e1} 
\end{align}
\item (Quartic tensor) For all $K_u,K_v,M_u,M_v\in \Dyadiclarge$ satisfying $K_u \simeq_\delta K_v$ and $M_u\simeq_\delta M_v$, it holds that
\begin{align}
& \sup_{\substack{y_1,\hdots, y_4\in \R \colon \\ |y_1|,\hdots, |y_4|\leq h}}
\Big\| \Theta^x_{y_1} P_{\leq N}^x \UN[K_u][+]  
\otimes \Theta^x_{y_2} P_{\leq N}^x  \IVN[K_v][-]
\otimes \Theta^x_{y_3} P_{\leq N}^x  \IUN[M_u][+] 
\otimes \Theta^x_{y_4} P_{\leq N}^x \VN[M_v][-] \notag \\
&\hspace{3ex}- \coup^2 \mathbf{1} \big\{ K_u = M_u \big\}  \mathbf{1} \big\{ K_v = M_v \big\} 
\Cf^{(\Ncs)}_{K_u}(y_3-y_1) \Cf^{(\Ncs)}_{K_v}(y_2-y_4) \,   \shuffle{(2,3)} \big( 
\Cas \otimes \Cas \big) \Big\|_{\Cprod{r-1}{r-1}} \notag \\
&\lesssim \,  \max\big( K_u,K_v,M_u,M_v)^{-(1-4\delta)\delta} \langle N h \rangle^2 \Dc^4, \label{killing:eq-tensors-quartic-e1}
\end{align}
where $\shuffle{(2,3)}$ is the shuffle operator from \eqref{prelim:eq-shuffle}. 
\end{enumerate}
\end{corollary}

\begin{proof}[Proof of Corollary \ref{killing:cor-tensors}] We prove the estimates \eqref{killing:eq-tensors-cubic-e1} and  \eqref{killing:eq-tensors-quartic-e1} separately. \\

\emph{Proof of \eqref{killing:eq-tensors-cubic-e1}:}
Using our product estimate (Lemma \ref{prelim:cor-product}), it holds that 
\begin{align*}
&\, \Big\| \Big( \biglcol \, \Theta^x_{y_1} P_{\leq N}^x \UN[K][+]  \otimes \Theta^x_{y_2} P_{\leq N}^x \IUN[L][+]   \bigrcol \Big)  \otimes  \Theta^x_{y_3} \VN[M][-] \Big\|_{\Cprod{r-1}{r-1}} \\
\lesssim&\, 
\Big\| \biglcol \, \Theta^x_{y_1} P_{\leq N}^x \UN[K][+]  \otimes \Theta^x_{y_2} P_{\leq N}^x \IUN[L][+]   \bigrcol 
\Big\|_{\Cprod{r-1}{s}} 
\Big\| \Theta^x_{y_3} \VN[M][-] \Big\|_{\Cprod{s}{r-1}}. 
\end{align*}
The first factor can be controlled using Lemma \ref{killing:lem-tensor-modulated-linear}. Using Hypothesis \ref{hypothesis:probabilistic}.\ref{ansatz:item-hypothesis-linear}, Hypothesis \ref{hypothesis:pre}.\ref{ansatz:item-pre-modulation}, and Lemma \ref{prelim:lem-psi-sum}, the second factor can be estimated by 
$\| \Theta^x_{y_3} \VN[M][-] \|_{\Cprod{s}{r-1}} \lesssim M^{r-s} \Dc$.
For a more general estimate of $\VN[M][-]$, we also refer to Lemma \ref{modulation:lem-linear} below. \\ 

\emph{Proof of \eqref{killing:eq-tensors-quartic-e1}:} 
For notational purposes, we first write
\begin{equation}\label{killing:eq-tensors-quartic-p0}
\begin{aligned}
&\Theta^x_{y_1} P_{\leq N}^x \UN[K_u][+]  
\otimes \Theta^x_{y_2} P_{\leq N}^x  \IVN[K_v][-]
\otimes \Theta^x_{y_3} P_{\leq N}^x  \IUN[M_u][+] 
\otimes \Theta^x_{y_4} P_{\leq N}^x \VN[M_v][-] \\
&- \coup^2 \mathbf{1} \big\{ K_u = M_u \big\}  \mathbf{1} \big\{ K_v = M_v \big\} 
\Cf^{(\Ncs)}_{K_u}(y_3-y_1) \Cf^{(\Ncs)}_{K_v}(y_2-y_4) \,   \shuffle{(2,3)} \big( 
\Cas \otimes \Cas \big)  \\
=&  \shuffle{(2,3)} \bigg( \Theta^x_{y_1} P_{\leq N}^x \UN[K_u][+]  
\otimes \Theta^x_{y_3} P_{\leq N}^x  \IUN[M_u][+] 
\otimes \Theta^x_{y_2} P_{\leq N}^x  \IVN[K_v][-]
\otimes \Theta^x_{y_4} P_{\leq N}^x \VN[M_v][-] \\
&- \coup^2 \mathbf{1} \big\{ K_u = M_u \big\}  \mathbf{1} \big\{ K_v = M_v \big\} 
\Cf^{(\Ncs)}_{K_u}(y_3-y_1) \Cf^{(\Ncs)}_{K_v}(y_2-y_4) \big( \Cas \otimes \Cas \big) \bigg). 
\end{aligned}
\end{equation}
We now estimate the argument of the shuffle operator in \eqref{killing:eq-tensors-quartic-p0}. 
Using Definition \ref{killing:def-Wick}, we write 
\begin{align}
&\Theta^x_{y_1} P_{\leq N}^x \UN[K_u][+]  
\otimes \Theta^x_{y_3} P_{\leq N}^x  \IUN[M_u][+] 
\otimes \Theta^x_{y_2} P_{\leq N}^x  \IVN[K_v][-]
\otimes \Theta^x_{y_4} P_{\leq N}^x \VN[M_v][-]  \notag \\
&-  \coup^2 \mathbf{1} \big\{ K_u = M_u \big\}  \mathbf{1} \big\{ K_v = M_v \big\} 
\Cf^{(\Ncs)}_{K_u}(y_3-y_1) \Cf^{(\Ncs)}_{K_v}(y_2-y_4) \big( \Cas \otimes \Cas \big)  \notag \\ 
=&  \Big( \biglcol \, 
\Theta^x_{y_1} P_{\leq N}^x \UN[K_u][+]  
\otimes \Theta^x_{y_3} P_{\leq N}^x  \IUN[M_u][+] \bigrcol \Big)
\otimes \Big( \biglcol \, 
 \Theta^x_{y_2} P_{\leq N}^x  \IVN[K_v][-]
\otimes \Theta^x_{y_4} P_{\leq N}^x \VN[M_v][-]  \bigrcol  \Big)  \label{killing:eq-tensors-quartic-p1} \\
+&  \coup \mathbf{1} \big\{ K_v = M_v \big\} \Cf^{(\Ncs)}_{K_v}(y_2-y_4)   \bigg( \Big(\biglcol \, 
\Theta^x_{y_1} P_{\leq N}^x \UN[K_u][+]  
\otimes \Theta^x_{y_3} P_{\leq N}^x  \IUN[M_u][+] \bigrcol  \Big) 
\otimes \Cas \bigg)  \label{killing:eq-tensors-quartic-p2} \\
+& \coup \mathbf{1} \big\{ K_u = M_u \big\} \Cf^{(\Ncs)}_{K_u}(y_3-y_1) 
\bigg( \Cas \otimes \Big( \biglcol \, 
\Theta^x_{y_2} P_{\leq N}^x  \IVN[K_v][-]
\otimes \Theta^x_{y_4} P_{\leq N}^x \VN[M_v][-] \bigrcol  \Big) \bigg) \label{killing:eq-tensors-quartic-p3}. 
\end{align}
By symmetry in the $u$ and $v$-variables, it suffices to estimate \eqref{killing:eq-tensors-quartic-p1} and \eqref{killing:eq-tensors-quartic-p2}. For \eqref{killing:eq-tensors-quartic-p1}, we first use the product estimate (Corollary \ref{prelim:cor-product})
and then use our quadratic tensor estimate (Lemma \ref{killing:lem-tensor-modulated-linear}), which yield 
\begin{align}
\big\| \eqref{killing:eq-tensors-quartic-p1} \big\|_{\Cprod{r-1}{r-1}} 
&\lesssim \Big\| \biglcol \, 
\Theta^x_{y_1} P_{\leq N}^x \UN[K_u][+]  
\otimes \Theta^x_{y_3} P_{\leq N}^x  \IUN[M_u][+] \bigrcol  \Big\|_{\Cprod{r-1}{s}} \notag \\
&\times \Big\|  \biglcol \, 
 \Theta^x_{y_2} P_{\leq N}^x  \IVN[K_v][-]
\otimes \Theta^x_{y_4} P_{\leq N}^x \VN[M_v][-]   \bigrcol  \Big\|_{\Cprod{s}{r-1}} \notag \\ 
&\lesssim \bigg( \max(K_u,M_u)^{r-s} M_u^{-\frac{1}{2}} + \mathbf{1} \big\{ K_u = M_u \big\} N^{-\delta+\vartheta} \langle N h \rangle \bigg) \label{killing:eq-tensors-quartic-q1} \\
&\times \bigg( \max(K_v,M_v)^{r-s} K_v^{-\frac{1}{2}} + \mathbf{1} \big\{ K_v = M_v \big\} N^{-\delta+\vartheta} \langle N h \rangle \bigg) \times \Dc^4. \notag 
\end{align}
In order to simplify \eqref{killing:eq-tensors-quartic-q1}, we recall that $r-s=\delta_1+\delta_2$ and $\vartheta=\delta_4$. Furthermore, we use that 
\begin{alignat}{5}
\max(K_u,M_u) &\lesssim K_u M_u, \qquad& M_u^{-\frac{1}{2}} &\lesssim M_u^{-\delta}, \qquad& N^{-\delta+\vartheta} \lesssim M_u^{-\delta+\vartheta}, \label{killing:eq-tensors-quartic-q2} \\ 
\max(K_v,M_v) &\lesssim K_v M_v, \qquad& K_v^{-\frac{1}{2}} &\lesssim K_v^{-\delta}, \qquad& N^{-\delta+\vartheta} \lesssim K_v^{-\delta+\vartheta}.\label{killing:eq-tensors-quartic-q3}
\end{alignat}
As a result, we obtain that
\begin{equation*}
\eqref{killing:eq-tensors-quartic-q1} \lesssim \big( K_u K_v M_u M_v\big)^{2\delta_1} \big( K_v M_u \big)^{-\delta} \langle N h \rangle^2 \Dc^4. 
\end{equation*}
Since $K_u \simeq_\delta K_v$ and $M_u \simeq_\delta M_v$, this yields an acceptable contribution to \eqref{killing:eq-tensors-quartic-e1}.\\

For \eqref{killing:eq-tensors-quartic-p2}, we first use the bound from Lemma \ref{ansatz:lem-renormalization}, which yields 
\begin{align}
\big\| \eqref{killing:eq-tensors-quartic-p2} \big\|_{\Cprod{r-1}{r-1}}
\lesssim& \, \big\| \eqref{killing:eq-tensors-quartic-p2} \big\|_{\Cprod{r-1}{s}}  \notag \\
\lesssim& \, \lambda \mathbf{1} \big\{ K_v = M_v \big\} \langle N h \rangle 
\Big\|  \biglcol \, 
\Theta^x_{y_1} P_{\leq N}^x \UN[K_u][+]  
\otimes \Theta^x_{y_3} P_{\leq N}^x  \IUN[M_u][+]   \bigrcol \Big\|_{\Cprod{r-1}{s}}. \label{killing:eq-tensors-quartic-p5}
\end{align}
Using Lemma \ref{killing:lem-tensor-modulated-linear} and \eqref{killing:eq-tensors-quartic-q2}, we then further obtain that
\begin{equation}\label{killing:eq-tensors-quartic-p6}
\eqref{killing:eq-tensors-quartic-p5}\lesssim \lambda \Dc^2 \mathbf{1} \big\{ K_v = M_v \big\}  (K_u M_u)^{2\delta_1} M_u^{-\delta} \langle Nh \rangle^2. 
\end{equation}
Due to Hypothesis \ref{hypothesis:post}, it holds that $\coup \Dc^2\lesssim \Dc^4$. Furthermore, by combining $K_u \simeq_\delta K_v$, $M_u \simeq_\delta M_v$, and $K_v=M_v$, we obtain the lower bound
\begin{equation*}
M_u \gtrsim \max(K_u,K_v,M_u,M_v)^{(1-\delta)^2}.
\end{equation*}
Thus, \eqref{killing:eq-tensors-quartic-p6} yields an acceptable contribution to \eqref{killing:eq-tensors-quartic-e1}.
\end{proof}

In the last proposition of this section, we consider the Lie bracket of two modulated bilinear waves, which is a quartic interaction. In contrast to the quartic interaction from Corollary \ref{killing:cor-tensors}, however, this interaction exhibits a cancellation and therefore does not require a renormalization.

\begin{proposition}[The quartic term]\label{killing:prop-quartic}
Let the post-modulation hypothesis (Hypothesis \ref{hypothesis:post}) be satisfied and let 
 $K_u,K_v,M_u,M_v\in \Dyadiclarge$ satisfy  $K_u \simeq_\delta K_v$ and $M_u \simeq_\delta M_v$. Then, it holds that 
\begin{equation}\label{killing:eq-quartic}
\Big\| \Big[ \UN[K_u,K_v][+-],  \VN[M_u,M_v][+-] \Big]_{\leq N} \Big\|_{\Cprod{r-1}{r-1}} \lesssim (K_u K_v M_u M_v)^{-\eta} \Dc^4. 
\end{equation}
\end{proposition}

\begin{remark}
The most dangerous aspect of \eqref{killing:eq-quartic} is the potential for two simultaneous probabilistic resonances, i.e., probabilistic resonances in both the $G^+_{u_0,k}$ and $G^-_{v_0,m}$-variables hidden in \eqref{killing:eq-quartic}. Fortunately, as we will see in the proof below, this contribution vanishes due to the properties of the Lie bracket. 
\end{remark}
\begin{proof}
Using the definition of the modulated bilinear waves (Definition \ref{ansatz:def-modulated-bilinear}) and the integral representation of $P_{\leq N}^x$, we first write 
\begin{equation}\label{killing:eq-quartic-p0}
\begin{aligned}
& \Big[  \UN[K_u,K_v][+-],  \VN[M_u,M_v][+-] \Big]_{\leq N} \\
=&P_{\leq N}^x \int_{\R^2} \dy \dz  \bigg( 
\big( \widecheck{\rho}_{\leq N} \ast \widecheck{\rho}_{\leq N} \big)(y) 
\big( \widecheck{\rho}_{\leq N} \ast \widecheck{\rho}_{\leq N} \big)(z) \\ 
&\times \Big[ \Big[ \Theta^x_y P_{\leq N}^x \UN[K_u][+], 
\Theta^x_y P_{\leq N}^x \IVN[K_v][-] \Big]_\frkg, 
\Big[ \Theta^x_z P_{\leq N}^x \IUN[M_u][+], 
\Theta^x_z P_{\leq N}^x \VN[M_v][-] \Big]_\frkg \Big]_\frkg  \bigg). 
\end{aligned}
\end{equation}
In order to analyze \eqref{killing:eq-quartic-p0}, we want to utilize our quartic tensor estimate from 
 Corollary \ref{killing:cor-tensors}. To this end, we define the iterated Lie bracket
\begin{equation}\label{killing:eq-quartic-p1}
\revision{\mathcal{L}} \colon \frkg^{\otimes 4} \rightarrow \frkg, \, \,  
A \otimes B \otimes C \otimes D 
\mapsto \big[ \big[ A , B \big], \big[ C, D \big] \big]. 
\end{equation}  
Using the shuffle operator from \eqref{prelim:eq-shuffle} and \eqref{killing:eq-quartic-p1}, we can write
\begin{align}
& \Big[ \Big[ \Theta^x_y P_{\leq N}^x \UN[K_u][+], 
\Theta^x_y P_{\leq N}^x \IVN[K_v][-] \Big]_\frkg, 
\Big[ \Theta^x_z P_{\leq N}^x \IUN[M_u][+], 
\Theta^x_z P_{\leq N}^x \VN[M_v][-] \Big]_\frkg \Big]_\frkg \notag \\
=& \, \revision{\mathcal{L}} \bigg(
 \Theta^x_y P_{\leq N}^x \UN[K_u][+]
 \otimes \Theta^x_y P_{\leq N}^x \IVN[K_v][-] 
\otimes  \Theta^x_z P_{\leq N}^x \IUN[M_u][+]
\otimes \Theta^x_z P_{\leq N}^x \VN[M_v][-] \bigg) \notag \\
=& \, \revision{\mathcal{L}} \bigg( 
 \Theta^x_y P_{\leq N}^x \UN[K_u][+]
 \otimes \Theta^x_y P_{\leq N}^x \IVN[K_v][-] 
\otimes  \Theta^x_z P_{\leq N}^x \IUN[M_u][+]
\otimes \Theta^x_z P_{\leq N}^x \VN[M_v][-] \label{killing:eq-quartic-p2} \\ 
&\hspace{5ex} - \coup^2 \mathbf{1} \big\{ K_u = M_u \big\} \mathbf{1} \big\{ K_v = M_v \big\} \Cf^{(\Ncs)}_{K_u}(z-y)  \Cf^{(\Ncs)}_{K_v}(y-z) \shuffle{(2,3)}\big( \Cas \otimes \Cas \big) \bigg) \notag \\  
+&  \, \coup^2 \mathbf{1} \big\{ K_u = M_u \big\} \mathbf{1} \big\{ K_v = M_v \big\} \Cf^{(\Ncs)}_{K_u}(z-y) \Cf^{(\Ncs)}_{K_v}(y-z)  \revision{\mathcal{L}} \Big(  \shuffle{(2,3)}\big( \Cas \otimes \Cas \big) \Big). \label{killing:eq-quartic-p3}
\end{align}
The contribution of \eqref{killing:eq-quartic-p2} can be controlled directly using Corollary \ref{killing:cor-tensors}. To treat \eqref{killing:eq-quartic-p3}, we use the definition of $\revision{\mathcal{L}}$ and $\Cas$, which yield that 
\begin{align}
\revision{\mathcal{L}} \Big(  \shuffle{(2,3)}\big( \Cas \otimes \Cas \big) \Big)
= \revision{\mathcal{L}} \Big( E_a \otimes E_b \otimes E^a \otimes E^b \Big) 
= \big[ \big[ E_a , E_b \big]_\frkg, \big[ E^a , E^b \big]_\frkg \big]_\frkg. \label{killing:eq-quartic-p6}
\end{align}
Due to the skew-symmetry of the Lie bracket, \eqref{killing:eq-quartic-p6} is identically zero, which completes the proof of \eqref{killing:eq-quartic}. We emphasize that the fact that \eqref{killing:eq-quartic-p6} is zero is crucial for the final estimate \eqref{killing:eq-quartic}. Indeed, at least in the regime $K_u,K_v,M_u,M_v\sim N$, neither the covariance functions in \eqref{killing:eq-quartic-p3} nor the integrals in \eqref{killing:eq-quartic-p0} would yield the desired decay in \eqref{killing:eq-quartic}.
\end{proof}

\section{Modulated and mixed modulated objects}\label{section:modulated-mixed}

In this section we obtain basic estimates of modulated objects, mixed modulated objects, and their interactions. Throughout this section, we only rely on the pre-modulation hypothesis (Hypothesis \ref{hypothesis:pre}) and never use the post-modulation hypothesis (Hypothesis \ref{hypothesis:post}). The reason is that we intend to use the estimates of this section in the contraction-mapping argument for the modulation equations. In Subsection \ref{section:modulation-structure}, we then improve certain estimates of this section by replacing using the post-modulation rather than pre-modulation hypothesis, i.e., making use of orthogonality (Proposition \ref{modulation:prop-properties}). 

\subsection{Modulated linear and bilinear waves}\label{section:modulated-sub-linear}

In the first subsection, we control the modulated linear and bilinear waves. The main ingredients in all the following estimates are the probabilistic bounds from Hypothesis \ref{hypothesis:pre}.

\begin{lemma}[Regularity of modulated linear  waves]\label{modulation:lem-linear}
Let the pre-modulation hypothesis (Hypothesis \ref{hypothesis:pre}) be satisfied.  Furthermore, let $K,M\in \Dyadiclarge$ and let $\gamma \in (-1,1)\backslash \{0\}$. Then, it holds that 
\begin{alignat}{3}
\big\| \UN[K][+] \big\|_{\Cprod{\gamma}{s}} &\lesssim \Dc   K^{\gamma+\frac{1}{2}+\eta}, 
& \qquad \qquad \big\| \IUN[K][+] \big\|_{\Cprod{\gamma}{s}} &\lesssim  \Dc   K^{\gamma-\frac{1}{2}+\eta}, 
\label{modulation:eq-lb-1}\\ 
\big\| \VN[M][-]\big\|_{\Cprod{s}{\gamma}}
&\lesssim \Dc   M^{\gamma+\frac{1}{2}+\eta}, & \qquad  \qquad 
\big\| \IVN[M][-]\big\|_{\Cprod{s}{\gamma}}
&\lesssim \Dc   M^{\gamma-\frac{1}{2}+\eta}. \label{modulation:eq-lb-2}
\end{alignat}
\end{lemma}

\begin{proof}
We only prove the estimate for $\UN[K][+]$, since the remaining estimates can be obtained similarly. Using Definition \ref{ansatz:def-modulated-linear}, it holds that
\begin{equation}\label{modulation:eq-lb-p0}
\UN[K][+] = \hcoup \sum_{u_0\in \LambdaRR} \sum_{k\in \Z_K} \psiRuK  \rhoND(k) \SN[K][+][k] G_{u_0,k}^+ \, e^{iku}. 
\end{equation}
Using Lemma \ref{prelim:lem-psi-sum}, we obtain that 
\begin{equation}\label{modulation:eq-lb-p1}
\Big\| \UN[K][+] \Big\|_{\Cprod{\gamma}{s}} 
\lesssim \hcoup \sup_{u_0\in \LambdaRR} \Big\| \sum_{k\in \Z_K} \rhoND(k)  \SN[K][+][k] G_{u_0,k}^+ \, e^{iku} \Big\|_{\Cprod{\gamma}{s}}
\end{equation}
Using  Hypothesis \ref{hypothesis:probabilistic}.\ref{ansatz:item-hypothesis-linear} and the bound on the modulation operators from Hypothesis \ref{hypothesis:pre}, it follows that
\begin{equation*}
\eqref{modulation:eq-lb-p1}
\lesssim \hcoup \Ac K^{\gamma+1/2+\eta} \Big\| \SN[K][+][k] \Big\|_{\Wuv[s][s]} \lesssim \hcoup \Ac \Bc K^{\gamma+1/2+\eta} = \Dc K^{\gamma+1/2+\eta}. \qedhere
\end{equation*}
\end{proof}

\begin{lemma}[Regularity of modulated bilinear waves]\label{modulation:lem-bilinear}
Assume that the pre-modulation hypothesis (Hypothesis \ref{hypothesis:pre}) is satisfied. Furthermore, let $K,M\in \Dyadiclarge$ satisfy $K\simeq_\delta M$. Then, it holds that  
\begin{align}
\big\| \UN[K,M][+-] \big\|_{\Cprod{s-1}{s}} 
+ \big\| \UN[K,M][+-] \big\|_{\Cprod{\scrr-1}{\eta}}
&\lesssim \Dc^2 (KM)^{-\eta}, 
\label{modulation:eq-lb-3} \\ 
\big\| \VN[K,M][+-]\big\|_{\Cprod{s}{s-1}}
+ \big\| \VN[K,M][+-]\big\|_{\Cprod{\eta}{\scrr-1}}
&\lesssim \Dc^2 (KM)^{-\eta}. \label{modulation:eq-lb-4}
\end{align}
 In addition, it holds for all  $\alpha,\beta\in (-s,s)$ \revision{that}
\begin{align}
\big\| \UN[K,M][+-] \big\|_{\Cprod{\alpha}{\beta}} 
&\lesssim \Dc^2 K^{\alpha+1/2+\eta} M^{\beta-1/2+\eta} \label{modulation:eq-bilinear-e1}, \\
\big\| \VN[K,M][+-] \big\|_{\Cprod{\alpha}{\beta}} 
&\lesssim  \Dc^2  K^{\alpha-1/2+\eta} M^{\beta+1/2+\eta} \label{modulation:eq-bilinear-e2}. 
\end{align}
\end{lemma} 

\begin{proof}
 By symmetry, it suffices to prove \eqref{modulation:eq-lb-3} and \eqref{modulation:eq-bilinear-e1}.  In order to prove \eqref{modulation:eq-lb-3}, we decompose
\begin{align*}
\UN[K,M][+-] &= \chinull[K,M] \Big[ \UN[K][+], \IVN[M][-] \Big]_{\leq N} \\
&= \chinull[K,M] \Big[ \UN[K][+] \Para[u][nsim] \IVN[M][-] \Big]_{\leq N}
 + \chinull[K,M] \Big[ \UN[K][+] \Para[u][sim] \IVN[M][-] \Big]_{\leq N}. 
\end{align*}
Using our para-product estimate (Lemma \ref{prelim:lem-paraproduct}), \eqref{modulation:eq-lb-1}, and \eqref{modulation:eq-lb-2}, the non-resonant part can be estimated by 
\begin{equation*}
\Big\| \Big[ \UN[K][+] \Para[u][nsim] \IVN[M][-] \Big]_{\leq N} 
\Big\|_{\Cprod{s-1}{s}} 
\lesssim \big\| \UN[K][+] \big\|_{\Cprod{s-1}{s}} 
\big\| \IVN[M][-] \big\|_{\Cprod{s}{s}} 
\lesssim \Dc^2 (KM)^{s-1/2+\eta}.
\end{equation*}
Since $s-1/2+\eta=-\delta_2+\delta_3$ and $\delta_3\ll \delta_2$, this is acceptable.  Similarly, it holds that
\begin{equation*}
\Big\| \Big[ \UN[K][+] \Para[u][nsim] \IVN[M][-] \Big]_{\leq N} 
\Big\|_{\Cprod{\scrr-1}{\eta}} 
\lesssim \big\| \UN[K][+] \big\|_{\Cprod{\scrr-1}{\eta}} 
\big\| \IVN[M][-] \big\|_{\Cprod{s}{\eta}} 
\lesssim \Dc^2 K^{\scrr-s} M^{\eta-s}.
\end{equation*}
Since $K\simeq_\delta M$ and since  
\begin{equation*}
(1-\delta)^{-1}(\scrr-s)+(1-\delta)(\eta-s) 
= (1-\delta)^{-1} (\tfrac{1}{2}-10\delta+\delta_2) + (1-\delta) (-\tfrac{1}{2}+\delta_2+\delta_3) = -9\delta + \mathcal{O}(\delta^2),
\end{equation*}
this is acceptable. 
In order to treat the resonant product, we first recall that $\UN[K][+]$ is supported on $u$-frequencies $\sim K$ and $\IVN[M][-]$ is supported on $u$-frequencies $\lesssim M^{1-\delta+\vartheta}$. Thus, \revision{for} $\UN[K][+] \Para[u][sim] \IVN[M][-]$ to be non-zero, it needs to hold that $K\lesssim M^{1-\delta+\vartheta}$. Using our para-product estimate (Lemma \ref{prelim:lem-paraproduct}), \eqref{modulation:eq-lb-1}, and \eqref{modulation:eq-lb-2}, it then follows that 
\begin{align*}
&\Big\| \Big[ \UN[K][+] \Para[u][sim] \IVN[M][-] \Big]_{\leq N} 
\Big\|_{\Cprod{s-1}{s}} 
+ \Big\| \Big[ \UN[K][+] \Para[u][sim] \IVN[M][-] \Big]_{\leq N} 
\Big\|_{\Cprod{\scrr-1}{\eta}} \\
\lesssim \, & \Big\| \Big[ \UN[K][+] \Para[u][sim] \IVN[M][-] \Big]_{\leq N} 
\Big\|_{\Cprod{\eta}{s}}  \\ 
\lesssim \, & \big\| \UN[K][+] \big\|_{\Cprod{-s+\eta}{s}}
\big\| \IVN[M][-] \big\|_{\Cprod{s}{s}} \\ 
\lesssim&\, \Dc^2  K^{-s+\frac{1}{2}+2\eta} M^{s-\frac{1}{2}+\eta}  
\lesssim \Dc^2 M^{-\delta (\frac{1}{2}-s)+3\eta+\vartheta}. 
\end{align*}
Since $-\delta (\frac{1}{2}-s)+3\eta +\vartheta =-\delta\delta_2+3\delta_3 + \delta_4$ and $\delta_3,\delta_4\ll \delta\delta_2$, this is acceptable. 
It remains to prove \eqref{modulation:eq-bilinear-e1}.  Using product estimates (Lemma \ref{prelim:lem-paraproduct}) and Lemma \ref{modulation:lem-linear}, it holds that
\begin{align*}
 &\big\| \UN[K,M][+-] \big\|_{\Cprod{\alpha}{\beta}} 
 \lesssim \Big\| \Big[ \UN[K][+] , \IVN[M][-] \Big]_{\leq N} \Big\|_{\Cprod{\alpha}{\beta}} \\
 \lesssim\,& \Big\| \UN[K][+] \Big\|_{\Cprod{\alpha}{s}} 
 \Big\| \IVN[M][-] \Big\|_{\Cprod{s}{\beta}} 
\lesssim \Dc^2 K^{\alpha+1/2+\eta} M^{\beta-1/2+\eta}. \qedhere
\end{align*}
\end{proof}

In our previous estimates (Lemma \ref{modulation:lem-linear} and Lemma \ref{modulation:lem-bilinear}), we obtained regularity estimates for the modulated linear and bilinear waves. We now turn our attention to commutator-type estimates.

\begin{lemma}\label{modulation:lem-PNX-modulated}
Assume that the pre-modulation hypothesis (Hypothesis \ref{hypothesis:pre}) is satisfied. Furthermore, let $K,L,M\in \Dyadiclarge$, let $j\in \mathbb{N}\backslash \{0\}$, and let $\gamma \in (-1,1)\backslash \{0\}$. Then, it holds that
\begin{align}
\Big\| (P_{\leq L}^x)^j \UN[K][+] 
- \hcoup \sum_{u_0\in \LambdaRR}  \sum_{k\in \Z_K} \psiRuK \rho_{\leq L}^{\, j}(k) \rhoND(k) \SN[K][+][k] G_{u_0,k}^+ \, e^{iku} \Big\|_{\Cprod{\gamma}{s}} &\lesssim \Dc K^{\gamma+1/2+\eta} K^{-\delta+\vartheta}, 
\label{modulation:lem-PNX-modulated-e1} \\ 
\Big\| (P_{\leq L}^x)^j \VN[M][-] 
-  \hcoup \hspace{-1ex} \sum_{v_0\in \LambdaRR} \sum_{m\in \Z_M} \psiRvM \rho_{\leq L}^{\, j}(m) \rhoND(m)  \SN[M][-][m] G_{v_0,m}^- \, e^{imv} \Big\|_{\Cprod{s}{\gamma}} &\lesssim \Dc M^{\gamma+1/2+\eta} M^{-\delta+\vartheta}. 
\label{modulation:lem-PNX-modulated-e2}
\end{align}
\end{lemma}

In comparison to \eqref{modulation:eq-lb-1} and \eqref{modulation:eq-lb-2}, the right-hand sides of \eqref{modulation:lem-PNX-modulated-e1} and \eqref{modulation:lem-PNX-modulated-e2} gain a factor of $K^{-\delta+\vartheta}$ and $M^{-\delta+\vartheta}$, respectively. 

\begin{proof}
Due to symmetry in the $u$ and $v$-variables, it suffices to prove \eqref{modulation:lem-PNX-modulated-e1}. Due to Lemma \ref{ansatz:lem-frequency-support}, we only have to treat the case $K\sim L$.
Due to pre-modulation hypothesis and Lemma \ref{prelim:lem-psi-sum} (see also the proof of Lemma \ref{modulation:lem-linear}), it holds that 
\begin{align*}
&\Big\| (P_{\leq L}^x)^j \UN[K][+] - \hcoup \sum_{u_0\in \LambdaRR}  \sum_{k\in \Z_K} \psiRuK  (P_{\leq L}^x)^j \Big( \rhoND(k) \SN[K][+][k] G_{u_0,k}^+ \, e^{iku} \Big) \Big\|_{\Cprod{\gamma}{s}} \\
\lesssim&\, \hcoup L^{-1} \sup_{u_0\in \LambdaRR} \Big\| \sum_{k\in \Z_K} \rhoND(k) \SN[K][+][k] G_{u_0,k}^+ e^{iku} \Big\|_{\Cprod{\gamma}{s}} 
\lesssim \Dc K^{\gamma+\frac{1}{2}+\eta-1},
\end{align*}
which is acceptable. Due to Lemma \ref{prelim:lem-psi-sum}, we may now fix $u_0 \in \LambdaR$.
To simplify the notation, we let $\widecheck{\rho}_{\leq L}^{\,\,\ast j}$ be the inverse Fourier-transform of $\rho_{\leq L}^{\, j}$. We now write
\begin{align*}
&\hcoup (P_{\leq L}^x)^j \sum_{k\in \Z_K} \rhoND(k) \SN[K][+][k] G_{u_0,k}^+ e^{iku}
- \hcoup \sum_{k\in \Z_K} \rho_{\leq L}^{\, j}(k) \rhoND(k) \SN[K][+][k] G_{u_0,k}^+ \, e^{iku} \\
=& \hcoup \int_{\R} \dy \, \widecheck{\rho}_{\leq L}^{\,\,\ast j}(y) \sum_{k\in \Z_K} \rhoND(k) \Big( \Theta^x_y \SN[K][+][k] - \SN[K][+][k] \Big) G_{u_0,k}^+ \, e^{ik(u+y)}. 
\end{align*}
Using Hypothesis \ref{hypothesis:probabilistic} and that  $\SN[K][+][k]$ is supported on frequencies $\lesssim K^{1-\delta+\vartheta}$, it then follows that
\begin{align*}
&\Big\| \hcoup (P_{\leq L}^x)^j \sum_{k\in \Z_K} \rhoND(k) \SN[K][+][k] G_{u_0,k}^+ e^{iku}
- \sum_{k\in \Z_K} \rho_{\leq L}^{\, j}(k) \rhoND(k) \SN[K][+][k] G_{u_0,k}^+ \, e^{iku} \Big\|_{\Cprod{\gamma}{s}} \\
\lesssim&\, \hcoup \Ac  K^{\gamma+1/2+\eta}
\bigg( \int_{\R}  \dy \, \big| \widecheck{\rho}_{\leq L}^{\,\,\ast j}(y) \big| \, \Big\| \rhoND(k) \Big( \Theta^x_y \SN[K][+][k] - \SN[K][+][k] \Big) \Big\|_{\Wuv[s][s][]} \bigg) \\
\lesssim&\, \hcoup \Ac \Bc K^{\gamma+1/2+\eta} K^{1-\delta+\vartheta} \Big( \int_{\R}  \dy \, \big| \widecheck{\rho}_{\leq L}^{\,\,\ast j}(y) \big| \, \big| y \big| \Big)  
\lesssim \Dc K^{\gamma+1/2+\eta} K^{1-\delta+\vartheta} K^{-1},
\end{align*}
which is acceptable. 
\end{proof}

In the next lemma, we control the difference of $\Int^u_{v\rightarrow u} \UN[K][+]$ and  $\IUN[K][+]$. 

\begin{lemma}[Integration error]\label{modulation:lem-integration}
Assume that the pre-modulation hypothesis (Hypothesis \ref{hypothesis:pre}) is satisfied and let $K,M\in \Dyadiclarge$.  Then, it holds that 
\begin{align}
\Big\| \partial_u \IUN[K][+] - \UN[K][+] \Big\|_{\Cprod{r-1}{r}}
&\lesssim \Dc K^{-1/2}, \label{killing:eq-integration-e1}\\
\Big\| \partial_v \IVN[M][-] - \VN[M][-] \Big\|_{\Cprod{r}{r-1}}
&\lesssim \Dc M^{-1/2}. \label{killing:eq-integration-e2}
\end{align}
Furthermore, it holds that 
\begin{align}
\Big\|  \IUN[K][+] - \Int^u_{v\rightarrow u} \UN[K][+] \Big\|_{\Cprod{2s}{s}}
&\lesssim \Dc K^{-\eta}, \label{killing:eq-integration-e3}\\
\Big\| \IVN[M][-] - \Int^v_{u\rightarrow v} \VN[M][-] \Big\|_{\Cprod{s}{2s}}
&\lesssim  \Dc M^{-\eta}. \label{killing:eq-integration-e4}
\end{align}
\end{lemma}

\begin{proof}
Due to symmetry, it suffices to prove  \eqref{killing:eq-integration-e1} and \eqref{killing:eq-integration-e3}. 
Using the definitions of $\UN[K][-]$ and $\IUN[K][-]$, we obtain that
\begin{align}
&\, \partial_u \IUN[K][+] - \UN[K][+] \notag \\
=&\, \hcoup \partial_u \sum_{u_0\in \LambdaRR} \sum_{k\in \Z_K} \psiRuK \rhoND(k) \SN[K][+][k] G_{u_0,k}^+ \frac{e^{iku}}{ik} 
- \hcoup \sum_{u_0\in \LambdaRR} \sum_{k\in \Z_K} \psiRuK \rhoND(k)  \SN[K][+][k] G_{u_0,k}^+ e^{iku} \notag \\ 
=&\, \hcoup  \sum_{u_0\in \LambdaRR} \sum_{k\in \Z_K} \rhoND(k)  \partial_u \Big( \psiRuK \SN[K][+][k] G_{u_0,k}^+ \Big) \frac{e^{iku}}{ik}. 
\label{killing:eq-integration-p1}
\end{align}
Using Lemma \ref{prelim:lem-psi-sum}, it then holds that 
\begin{equation}\label{killing:eq-integration-q1}
\begin{aligned}
\big\| \eqref{killing:eq-integration-p1} \big\|_{\Cprod{r-1}{r}} 
\lesssim&\, \revision{\hcoup} \sup_{u_0\in \LambdaRR} \Big\| \sum_{k\in \Z_K} \rhoND(k)  \SN[K][+][k] G_{u_0,k}^+ \revision{\frac{e^{iku}}{ik}}  \Big\|_{\Cprod{r-1}{r}} \\
+&\, 
\revision{\hcoup} \sup_{u_0\in \LambdaRR} \Big\| \sum_{k\in \Z_K} \rhoND(k)  \big( \partial_u \SN[K][+][k] \big) G_{u_0,k}^+ \revision{\frac{e^{iku}}{ik}}  \Big\|_{\Cprod{r-1}{r}}.
\end{aligned}
\end{equation}
Using the pre-modulation hypothesis (Hypothesis \ref{hypothesis:pre}) and Lemma \ref{ansatz:lem-frequency-support}, which implies that $\SN[K][+][k]$ is supported on frequencies $\lesssim K^{1-\delta+\vartheta}$, it follows that 
\begin{align*}
\eqref{killing:eq-integration-q1} 
\lesssim \hcoup \Ac K^{r-\frac{3}{2}+\eta} \Big( \Big\|   \SN[K][+][k] \Big\|_{\Wuv[\eta][r][k]} + \Big\| \partial_u  \SN[K][+][k] \Big\|_{\Wuv[\eta][r][k]} \Big)  
\lesssim \Dc K^{r-\frac{3}{2}+\eta}  \big( K^{1-\delta+\vartheta} \big)^{1+\eta-s} \big( K^{1-\delta+\vartheta} \big)^{r-s}.
\end{align*}
Since
\begin{align*}
&r-\frac{3}{2}+\eta + (1-\delta+\vartheta)(1+\eta-s)+(1-\delta+\vartheta)(r-s) \\  =\,&  \delta_1 - 1 + \delta_3 + (1-\delta+\delta_4)(\tfrac{1}{2}+\delta_2+\delta_3)+(1-\delta+\delta_4)(\delta_1+\delta_2) = -\frac{1}{2}-\frac{\delta}{2} + \mathcal{O}(\delta_1), 
\end{align*}
this completes the proof of \eqref{killing:eq-integration-e1}.
 It now only remains to prove \eqref{killing:eq-integration-e3}. Due to the definition of the $\Cprod{\gamma_1}{\gamma_2}$-norms, it holds that 
\begin{align*}
&\Big\| \,  \IUN[K][+] - \Int^u_{v\rightarrow u} \UN[K][+] \Big\|_{\Cprod{2s}{s}} \\ 
\lesssim \, &  
\Big\|  \, \IUN[K][+] -  \Int^u_{v\rightarrow u} \UN[K][+] \Big\|_{\Cprod{s}{s}}
 + 
\Big\| \, \partial_u \Big( \IUN[K][+] -  \Int^u_{v\rightarrow u} \UN[K][+] \Big) \Big\|_{\Cprod{2s-1}{s}}. 
\end{align*}
For the first term, we use our Duhamel integral estimate (Lemma \ref{prelim:lem-Duhamel-integral}) and Lemma \ref{modulation:lem-linear}, which yield that 
\begin{equation*}
\Big\|  \, \IUN[K][+] -  \Int^u_{v\rightarrow u} \UN[K][+] \Big\|_{\Cprod{s}{s}} 
\lesssim \Big\| \, \IUN[K][+] \Big\|_{\Cprod{s}{s}}
+ \Big\|  \, \UN[K][+] \Big\|_{\Cprod{s-1}{s}} \lesssim \Dc K^{s-1/2+\eta}. 
\end{equation*}
For the second term, we use \eqref{killing:eq-integration-e1}, which yields that 
\begin{equation*}
\Big\| \, \partial_u \Big( \IUN[K][+] -  \Int^u_{v\rightarrow u} \UN[K][+] \Big) \Big\|_{\Cprod{2s-1}{s}}
\lesssim 
\Big\| \, \partial_u \IUN[K][+] - \UN[K][+] \Big\|_{\Cprod{2s-1}{s}}
\lesssim \Dc K^{2s-r-1/2}.
\end{equation*}
Since $s-1/2+\eta=-\delta_2+\delta_3$ and $2s-r-1/2=-\delta_1-2\delta_2$, both contributions are acceptable, and thus it completes the proof of \eqref{killing:eq-integration-e3}. 
\end{proof}

In the next lemma, we control an integration error in the modulated bilinear wave. 

\begin{lemma}[Bilinear integration error]\label{modulation:lem-integration-bilinear} 
Assume that the pre-modulation hypothesis (Hypothesis \ref{hypothesis:pre}) is satisfied. Furthermore, let  $K,M\in \Dyadiclarge$ satisfy $K\simeq_\delta M$. Then, it holds that
\begin{equation}\label{modulation:eq-integration-bilinear}
\Big\| \Int^v_{u\rightarrow v} \VN[K,M][+-] - \chinull[K,M] \Big[ \IUN[K][+],  \IVN[M][-] \Big]_{\leq N} \Big\|_{\Cprod{s}{r}} 
\lesssim \Dc^2 (KM)^{-\eta}. 
\end{equation}
\end{lemma}

\begin{proof}
We first estimate 
\begin{align}
&\Big\| \Int^v_{u\rightarrow v} \VN[K,M][+-] -  \chinull[K,M] \Big[ \IUN[K][+],  \IVN[M][-] \Big]_{\leq N} \Big\|_{\Cprod{s}{r}} 
\notag \\ 
\lesssim\, &  \Big\| \Int^v_{u\rightarrow v} \VN[K,M][+-] -  \chinull[K,M] \Big[ \IUN[K][+],  \IVN[M][-] \Big]_{\leq N} \Big\|_{\Cprod{s}{s}} 
\notag \\ 
+\, &  \Big\|  \partial_v \Big( \Int^v_{u\rightarrow v} \VN[K,M][+-] -  \chinull[K,M] \Big[ \IUN[K][+],  \IVN[M][-] \Big]_{\leq N} \Big) \Big\|_{\Cprod{s}{r-1}} 
\notag \\ 
\lesssim\, & \Big\| \Int^v_{u\rightarrow v} \VN[K,M][+-] \Big\|_{\Cprod{s}{s}} 
+ \Big\| \Big[ \IUN[K][+],  \IVN[M][-] \Big]_{\leq N} \Big\|_{\Cprod{s}{s}} 
\label{modulation:eq-integration-bilinear-p1}\\
+\, & \Big\| \chinull[K,M] \Big[ \IUN[K][+], \VN[M][-] \Big]_{\leq N} 
- \partial_v \Big(  \chinull[K,M] \Big[ \IUN[K][+],  \IVN[M][-] \Big]_{\leq N}  \Big)
\Big\|_{\Cprod{s}{r-1}}. \label{modulation:eq-integration-bilinear-p2}
\end{align}
We now estimate the terms in \eqref{modulation:eq-integration-bilinear-p1} and \eqref{modulation:eq-integration-bilinear-p2} separately. Using our integral estimate (Lemma \ref{prelim:lem-Duhamel-integral}), \revision{Definition \ref{ansatz:def-modulated-bilinear}}, and\footnote{To be more precise, we use the proof of Lemma \ref{modulation:lem-bilinear}, in which we controlled $\big[\IUN[K][+],\VN[M][-] \big]_{\leq N}$.} Lemma \ref{modulation:lem-bilinear}, the first summand in \eqref{modulation:eq-integration-bilinear-p1} can be bounded by
\begin{equation*}
  \Big\| \Int^v_{u\rightarrow v} \VN[K,M][+-] \Big\|_{\Cprod{s}{s}}  
  \lesssim \Big\|  \Big[ \IUN[K][+], \VN[M][-] \Big]_{\leq N}  \Big\|_{\Cprod{s}{s-1}} \lesssim \Dc^2 (KM)^{-\eta},  
\end{equation*}
which is acceptable. Using our paraproduct estimate (Lemma \ref{prelim:lem-paraproduct}) and Lemma \ref{modulation:lem-linear}, the second summand in \eqref{modulation:eq-integration-bilinear-p1} can be bounded by
\begin{equation}\label{modulation:eq-integration-bilinear-q1}
\Big\| \Big[ \IUN[K][+],  \IVN[M][-] \Big]_{\leq N} \Big\|_{\Cprod{s}{s}} 
\lesssim \Big\| \IUN[K][+] \Big\|_{\Cprod{s}{s}} 
\Big\| \IVN[M][-] \Big\|_{\Cprod{s}{s}} \lesssim \Dc^2  (KM)^{-\eta},
\end{equation}
which is acceptable. It now remains to control  \eqref{modulation:eq-integration-bilinear-p2}. 
To this end, we first estimate
\begin{align}
\eqref{modulation:eq-integration-bilinear-p2}
&\leq \Big\| \big( \partial_v \chinull[K,M] \big) \Big[ \IUN[K][+], \IVN[M][-]\Big]_{\leq N} \Big\|_{\Cprod{s}{r-1}} 
\label{modulation:eq-integration-bilinear-q2} \\ 
&+ \Big\|  \chinull[K,M] \Big[ \IUN[K][+], \VN[M][-] - \partial_v \IVN[M][-] \Big]_{\leq N} \Big\|_{\Cprod{s}{r-1}} 
\label{modulation:eq-integration-bilinear-p3} \\ 
&+ \Big\| \chinull[K,M]  \Big[ \partial_v \IUN[K][+] \Para[v][nsim]  \IVN[M][-] \Big]_{\leq N} 
\Big\|_{\Cprod{s}{r-1}} \label{modulation:eq-integration-bilinear-p4} 
\\
&+ \Big\| \chinull[K,M]  \Big[ \partial_v \IUN[K][+] \Para[v][sim]  \IVN[M][-] \Big]_{\leq N} 
\Big\|_{\Cprod{s}{r-1}} \label{modulation:eq-integration-bilinear-p5}. 
\end{align}
The first term \eqref{modulation:eq-integration-bilinear-q2} can be estimated as in \eqref{modulation:eq-integration-bilinear-q1}. 
Using Lemma \ref{modulation:lem-linear} and Lemma \ref{modulation:lem-integration}, it holds that
\begin{equation*}
\eqref{modulation:eq-integration-bilinear-p3} 
\lesssim \Big\| \IUN[K][+] \Big\|_{\Cprod{s}{s}} 
\Big\| \VN[M][-] - \partial_v \IVN[M][-] \Big\|_{\Cprod{s}{r-1}} 
\lesssim \Dc^2  K^{-\eta} M^{-1/2},
\end{equation*}
which is acceptable. Using our para-product estimate (Lemma \ref{prelim:lem-paraproduct}) and Lemma \ref{modulation:lem-linear}, it holds that
\begin{align*}
\eqref{modulation:eq-integration-bilinear-p4}
\lesssim 
\Big\|  \partial_v \IUN[K][+]  \Big\|_{\Cprod{s}{r-1}} 
\Big\|   \IVN[M][-] \Big\|_{\Cprod{s}{\eta}} 
\lesssim \Dc^2  K^{r-s} M^{\eta-s}.
\end{align*}
Since $K\simeq_\delta M$ and $(r-s)+(\eta-s)=-\frac{1}{2}+\mathcal{O}(\delta_1)$, this is more than acceptable. It now remains to treat \eqref{modulation:eq-integration-bilinear-p5}. Since $\IUN[K][+]$ is supported on $v$-frequencies $\lesssim K^{1-\delta+\vartheta}$ and $\VN[M][-]$ is supported on $v$-frequencies $\sim M$,  \eqref{modulation:eq-integration-bilinear-p5} is only non-trivial when $M\lesssim K^{1-\delta+\vartheta}$. Using our paraproduct estimate (Lemma \ref{prelim:lem-paraproduct}) and Lemma \ref{modulation:lem-linear}, it then follows that
\begin{equation*}
    \eqref{modulation:eq-integration-bilinear-p5}
\lesssim 
\Big\|  \partial_v \IUN[K][+]  \Big\|_{\Cprod{s}{s-1}} 
\Big\|   \IVN[M][-] \Big\|_{\Cprod{s}{1-s+\eta}} 
\lesssim \Dc^2 K^{s-1/2+\eta} M^{1/2-s+2\eta} 
\lesssim  \Dc^2 K^{-\delta (1/2-s)+3\eta+\vartheta}.
\end{equation*}
Since
\begin{equation*}
-\delta \big( \tfrac{1}{2}-s \big)+3\eta +\vartheta = -\delta \delta_2 + 3 \delta_3 + \delta_4 \qquad \text{and} \qquad \delta_3,\delta_4  \ll \delta \delta_2, 
\end{equation*}
this is acceptable.
\end{proof}

\subsection{Mixed modulated objects}\label{section:modulated-sub-mixed}

In this subsection, we control the mixed modulated linear objects $\UN[][+\fs]$, $\UN[][\fs-]$, $\VN[][\fs-]$, and $\VN[][+\fs]$, which were previously defined in Definition \ref{ansatz:def-mixed}. The basic estimates for the mixed modulated objects are contained in Lemma \ref{modulation:lem-mixed} below, which is proven using a kickback argument. 

\begin{lemma}[Control of mixed modulated objects]\label{modulation:lem-mixed}
Let the pre-modulation hypothesis (Hypothesis \ref{hypothesis:pre}) be satisfied.  Then, we have the following estimates:
\begin{align}
\sup_{M\in \Dyadiclarge} M^\eta \big\| \UN[M][\fs-] \big\|_{\Cprod{r-1}{s}} &\lesssim \Dc^2, 
\label{mixed:eq-control-1}  \\
\sup_{K\in \Dyadiclarge} K^\eta \big\| \VN[K][+\fs] \big\|_{\Cprod{s}{r-1}} &\lesssim \Dc^2, 
\label{mixed:eq-control-2} \\
\sup_{K\in \Dyadiclarge} K^\eta \big\| \UN[K][+\fs] \big\|_{\Cprod{s-1}{r}} &\lesssim \Dc^2, 
\label{mixed:eq-control-3} \\ 
\sup_{M\in \Dyadiclarge} M^\eta \big\| \VN[M][\fs-] \big\|_{\Cprod{r}{s-1}} &\lesssim \Dc^2. 
\label{mixed:eq-control-4} 
\end{align}
\end{lemma}

For the reader's convenience, we recall from Definition \ref{ansatz:def-mixed} that
\begin{align}
\UN[M][\fs-] &= \chinull[M] 
\bigg[ P^{u}_{\geq M^{1-\deltap}} \Big( P^v_{<M^{1-\deltap}} \UN[][\fs]
+ \Sumlarge_{K< M^{1-\delta}} \UN[K][\fs-] \Big) \Para[u][gg] \IVN[M][-] \bigg]_{\leq N}, \label{mixed:eq-def-1}
\\ 
\VN[K][+\fs]&= \chinull[K]
\bigg[ \IUN[K][+] \Para[v][ll] P^v_{\geq K^{1-\deltap}} \Big( P^u_{<K^{1-\deltap}} \VN[][\fs] + \Sumlarge_{M< K^{1-\delta}} \VN[M][+\fs] \Big) \bigg]_{\leq N}, \label{mixed:eq-def-2} \\ 
\UN[K][+\fs] &= \chinull[K] \bigg[ \UN[K][+] \Para[v][ll]  \Int^v_{u\rightarrow v} P^v_{\geq K^{1-\deltap}} \Big( P^u_{<K^{1-\deltap}} \VN[][\fs] + \Sumlarge_{M< K^{1-\delta}} \VN[M][+\fs] \Big) \bigg]_{\leq N},  \label{mixed:eq-def-3}  \\ 
\VN[M][\fs-] &= \chinull[M] \bigg[ \Int^u_{v\rightarrow u} P^{u}_{\geq M^{1-\deltap}} \Big( P^v_{<M^{1-\deltap}} \UN[][\fs]
+ \Sumlarge_{K< M^{1-\delta}} \UN[K][\fs-] \Big) \Para[u][gg] \VN[M][-] \bigg]_{\leq N}.  \label{mixed:eq-def-4}
\end{align}

\begin{proof}[Proof of Lemma \ref{modulation:lem-mixed}:]
We first prove \eqref{mixed:eq-control-1}, i.e., we first prove uniform estimates for the mixed modulated objects $(\UN[M][\fs-])_{M\in \Dyadiclarge}$. Using a product estimate (Lemma \ref{prelim:lem-paraproduct}) and Lemma \ref{modulation:lem-linear}, we obtain that 
\begin{align}
\big\| \UN[M][\fs-] \big\|_{\Cprod{r-1}{s}} 
&\lesssim \Big\| P^v_{<M^{1-\deltap}} \UN[][\fs]
+ \Sumlarge_{K< M^{1-\delta}} \UN[K][\fs-] \Big\|_{\Cprod{r-1}{s}} 
\big\| \IVN[M][-] \big\|_{\Cprod{s}{s}}  \notag \\
&\lesssim \Dc M^{-\eta} \Big( \big\| \UN[][\fs] \big\|_{\Cprod{r-1}{r}} 
+ \Sumlarge_{\substack{ K < M^{1-\delta}}} \big\| \UN[K][\fs-] \big\|_{\Cprod{r-1}{s}} \Big) \notag \\
&\lesssim \Dc M^{-\eta} \Big( \Dc+ \sup_{\substack{K\in \Dyadiclarge \colon \\ K < M^{1-\delta}}} K^\eta 
\big\| \UN[K][\fs-] \big\|_{\Cprod{r-1}{s}} \Big). \label{mixed:eq-control-p1}
\end{align}
By multiplying \eqref{mixed:eq-control-p1} with $M^\eta$, it follows that 
\begin{equation}\label{mixed:eq-control-p2} 
 M^\eta \big\| \UN[M][\fs-] \big\|_{\Cprod{r-1}{s}} 
\lesssim \Dc^2 +  \Dc \sup_{\substack{K\in \Dyadiclarge \colon \\ K < M^{1-\delta}}} K^\eta 
\big\| \UN[K][\fs-] \big\|_{\Cprod{r-1}{s}}. 
\end{equation}

Using Hypothesis \ref{hypothesis:post}, which yields that $\Dc=\hcoup \Ac \Bc \ll 1$, this implies the desired estimate \eqref{mixed:eq-control-1}. Due to symmetry in the $u$ and $v$-variables, the proof of \eqref{mixed:eq-control-2} is similar and we omit the details.

We now turn to the proof of \eqref{mixed:eq-control-3}. Due to the Lemma \ref{ansatz:lem-frequency-support}, the right-hand side of \eqref{mixed:eq-def-3} only contains high$\times$low-interactions in the $u$-variable. Furthermore, it contains only low$\times$high-interactions in the $v$-variable. From the para-product estimate (Lemma \ref{prelim:lem-paraproduct}) and the Duhamel integral estimate  (Lemma \ref{prelim:lem-Duhamel-integral}), it follows that 
\begin{align}
\big\| \UN[K][+\fs] \big\|_{\Cprod{s-1}{r}}  
&\lesssim \big\|  \UN[K][+] \big\|_{\Cprod{s-1}{s}} 
\Big\| \Int^v_{u\rightarrow v} P^v_{\geq K^{1-\deltap}} \Big( P^u_{<K^{1-\deltap}} \VN[][\fs] + \Sumlarge_{M< K^{1-\delta}} \VN[M][+\fs] \Big) \Big\|_{\Cprod{s}{r}} \\ 
&\lesssim  \big\|  \UN[K][+] \big\|_{\Cprod{s-1}{s}}  
\Big( \big\| \VN[][\fs]\big\|_{\Cprod{s}{r-1}} + \Sumlarge_{M<K^{1-\delta}} \big\| \VN[M][+\fs]\big\|_{\Cprod{s}{r-1}} \Big).
\label{mixed:eq-control-q1}
\end{align}
Using Lemma \ref{modulation:lem-linear} and \eqref{mixed:eq-control-2}, it then follows that 
\begin{align*}
\eqref{mixed:eq-control-q1}
\lesssim \Dc K^{-\eta} \Big( \Dc + \Dc \Sumlarge_{M<K^{1-\delta}} M^{-\eta} \Big) \lesssim \Dc K^{-\eta}. 
\end{align*}
Due to symmetry in the $u$ and $v$-variables, the proof of \eqref{mixed:eq-control-4} is similar. 
\end{proof}

As a direct consequence of Lemma \ref{modulation:lem-mixed}, we also obtain control of $\UN[][\fsc]$ and $\VN[][\fcs]$ from \eqref{ansatz:eq-UN-sast} and \eqref{ansatz:eq-VN-asts}. 

\begin{corollary}\label{modulation:cor-control-combined}
Let the pre-modulation hypothesis (Hypothesis \ref{hypothesis:pre}) be satisfied. Then, it holds that 
\begin{align*}
\sup_{M \in \Dyadiclarge} \big\| \UN[<M^{1-\delta}][\fsc] \big\|_{\Cprod{r-1}{s}} \lesssim \Dc
\qquad \text{and} \qquad 
\sup_{K\in \Dyadiclarge} \big\| \VN[<K^{1-\delta}][\fcs] \big\|_{\Cprod{s}{r-1}} \lesssim \Dc. 
\end{align*}
\end{corollary}

\begin{proof}
Due to symmetry in the $u$ and $v$-variables, it suffices to prove the estimate of $\UN[<M^{1-\delta}][\fsc]$. Using the definition of $\UN[<M^{1-\delta}][\fsc]$ from \eqref{ansatz:eq-UN-sast} and Lemma \ref{modulation:lem-mixed}, it holds that
\begin{equation*}
\big\| \UN[<M^{1-\delta}][\fsc] \big\|_{\Cprod{r-1}{s}} 
\lesssim \big\| \UN[][\fs] \big\|_{\Cprod{r-1}{s}} 
+ \Sumlarge_{\substack{L\leq \Nd \colon \\ L <M^{1-\delta}}} \big\| \UN[L][\fs-] \big\|_{\Cprod{r-1}{s}} \lesssim \Dc + \Sumlarge_{\substack{L\leq \Nd \colon \\ L <M^{1-\delta}}} L^{-\eta} \Dc^2 \lesssim \Dc. \qedhere
\end{equation*}
\end{proof}

Equipped with Lemma \ref{modulation:lem-mixed} and Corollary \ref{modulation:cor-control-combined}, it is relatively easy to also bound the simplified high$\times$high$\rightarrow$low-interactions from Definition \ref{ansatz:def-shhl}. 

\begin{lemma}[Control of simplified high$\times$high$\rightarrow$low-interactions]\label{modulation:lem-shhl}
Let the pre-modulation hypothesis (Hypothesis \ref{hypothesis:pre}) be satisfied. Then, it holds that
\begin{align}
\sup_{K\in \Dyadiclarge} K^\eta \big\| \SHHLN[K][v] \big\|_{\Cprod{s}{\eta}} &\lesssim \Dc^2, 
\label{mixed:eq-shhl-1} \\
\sup_{M\in \Dyadiclarge} M^\eta \big\| \SHHLN[M][u] \big\|_{\Cprod{\eta}{s}} &\lesssim \Dc^2. 
\label{mixed:eq-shhl-2}
\end{align}
\end{lemma}

\begin{proof}[Proof of Lemma \ref{modulation:lem-shhl}:]
We only prove \eqref{mixed:eq-shhl-1}, since the proof of \eqref{mixed:eq-shhl-2} is similar. Due to its definition (and the frequency-support properties of  $\IVN[M][-]$), $\SHHLN[K][v]$ only contains high$\times$high-interactions in the $v$-variable. Using the high$\times$high-paraproduct estimate (Lemma \ref{prelim:lem-paraproduct}), Lemma \ref{modulation:lem-linear}, and Corollary \ref{modulation:cor-control-combined}, it then follows that  
\begin{align*}
&\big\| \SHHLN[K][v] \big\|_{\Cprod{s}{\eta}} 
\lesssim \Sumlarge_{\substack{M \leq \Nd \colon \\ M \simeq_\delta K}} 
\big\| \VN[<K^{1-\delta}][\fcs] \big\|_{\Cprod{s}{r-1}} 
\big\| \IVN[M][-] \big\|_{\Cprod{s}{1-r^\prime}}  \\ 
\lesssim\,&   \Sumlarge_{\substack{M \leq \Nd \colon \\ M \simeq_\delta K}} M^{1-r^\prime-s}  \Dc^2 
\lesssim K^{(1-\delta)(1-r^\prime-s)} \Dc^2. 
\end{align*}
Since $(1-\delta)(1-r^\prime-s)=(1-\delta)(-\delta_1+3\delta_2)=-\delta_1+\mathcal{O}(\delta\delta_1)$, this is acceptable.
\end{proof}

\subsection{Reversed modulated linear waves}\label{section:modulated-sub-reversed} 
In this subsection, 
we turn to the reversed modulated linear waves $\UN[M][-]$ and $\VN[K][+]$. 

\begin{lemma}[Reversed modulated linear waves]\label{modulation:lem-linear-reversed}
Let the pre-modulation hypothesis (Hypothesis \ref{hypothesis:pre}) be satisfied. Then, it holds that
\begin{align}
\sup_{M\in \Dyadiclarge} M^{\delta (1/2-s)-10\eta} \Big\| \, \UN[M][-] \Big\|_{\Cprod{-1/2+\eta}{s}} &\lesssim \Dc^2, \label{modulation:eq-lf-1} \\ 
\sup_{K\in \Dyadiclarge} K^{\delta (1/2-s)-10\eta}\Big\| \, \VN[K][+] \Big\|_{\Cprod{s}{-1/2+\eta}} &\lesssim \Dc^2. \label{modulation:eq-lf-2}
\end{align}
\end{lemma}

\begin{proof}
Due to symmetry in the $u$ and $v$-variables, it suffices to prove \eqref{modulation:eq-lf-1}. In order to prove \eqref{modulation:eq-lf-1}, it is sufficient to prove that 
\begin{align}
\Big\| \, \UN[M][-] \Big\|_{\Cprod{-1/2+\eta}{s}} &\lesssim  M^{-\delta (1/2-s)+10\eta} 
 \Big( \Dc^2 + \Dc \sup_{\substack{ K \in \Dyadiclarge \colon \\ K<M^{1-\delta}}} K^{r+\eta-2s} \big\| \UN[K][-] \big\|_{\Cprod{s-1}{s}} \Big)
\label{modulation:eq-lf-p9}  
\end{align}
To this end, we first use the definitions of  $\UN[M][-]$ and $\LON[M][+]$, which yield that
\begin{align}
\UN[M][-] &=  \chinull[M] \Big[ \LON[M][+], \IVN[M][-] \Big]_{\leq N}  +  \chinull[M] \Big[ \SHHLN[M][u], (P_{\leq N}^x)^2 \IVN[M][-] \Big]_{\leq N} \notag  \\ 
&=  \chinull[M] \Sumlarge_{K<M^{1-\delta}} 
 \Big[ \UN[K][+], \IVN[M][-] \Big]_{\leq N} 
 \label{modulation:eq-lf-p1} \allowdisplaybreaks[4] \\ 
&+  \chinull[M] \Sumlarge_{K<M^{1-\delta}} 
 \Big[ \UN[K][-], \IVN[M][-] \Big]_{\leq N} 
 \label{modulation:eq-lf-p2} \allowdisplaybreaks[4] \\ 
 &+  \chinull[M] \Sumlarge_{K_u,K_v<M^{1-\delta}} 
 \Big[ \UN[K_u,K_v][+-], \IVN[M][-] \Big]_{\leq N} 
 \label{modulation:eq-lf-p3}  \allowdisplaybreaks[4] \\ 
 &+  \chinull[M]\Sumlarge_{K<M^{1-\delta}} 
 \Big[ P^v_{<M^{1-\deltap}} \UN[K][+\fs], \IVN[M][-] \Big]_{\leq N} 
 \label{modulation:eq-lf-p4}  \allowdisplaybreaks[4] \\ 
 &+  \chinull[M] \Sumlarge_{K<M^{1-\delta}} 
 \Big[ P^u_{<M^{1-\deltap}} \UN[K][\fs-], \IVN[M][-] \Big]_{\leq N} 
 \label{modulation:eq-lf-p5}  \allowdisplaybreaks[4] \\ 
 &+  \chinull[M] \Big[ P^{u,v}_{<M^{1-\deltap}}  \UN[][\fs], \IVN[M][-] \Big]_{\leq N} 
 \label{modulation:eq-lf-p6}  \allowdisplaybreaks[4] \\ 
 &+   \chinull[M] \Big[ \SHHLN[M][u], (P_{\leq N}^x)^2 \IVN[M][-] \Big]_{\leq N} \label{modulation:eq-lf-p7}. 
\end{align}
We now separately estimate \eqref{modulation:eq-lf-p1} in the $\Cprod{-1/2+\eta}{s}$-norm and \eqref{modulation:eq-lf-p2}-\eqref{modulation:eq-lf-p7} in the stronger $\Cprod{r-1}{s}$-norm. \\

\emph{Estimate of \eqref{modulation:eq-lf-p1}:} Using product estimates, using that $\IVN[M][-]$ is supported on $u$-frequencies $\lesssim M^{1-\delta+\vartheta}$, and Lemma \ref{modulation:lem-linear}, it holds that
\begin{align*}
\Big\| \Big[ \UN[K][+] , \IVN[M][-] \Big]_{\leq N} 
\Big\|_{\Cprod{-1/2+\eta}{s}}
&\lesssim \big\| \UN[K][+] \big\|_{\Cprod{-1/2+\eta}{s}} 
\big\| \IVN[M][-] \big\|_{\Cprod{1/2+\eta}{s}} \\ 
&\lesssim K^{2\eta} \big(M^{1-\delta+\vartheta}\big)^{\frac{1}{2}+\eta-s} M^{s-\frac{1}{2}+\eta} \Dc^2 
\lesssim K^{2\eta} M^{-\delta (\frac{1}{2}-s)+2\eta+\vartheta} \Dc^2. 
\end{align*}
Since $K<M^{1-\delta}\lesssim M$, this is acceptable. \\ 

\emph{Estimate of \eqref{modulation:eq-lf-p2}:} We recall that $\UN[K][-]$ is supported on $u$-frequencies $\lesssim K$ and $v$-frequencies $\sim K$. In particular, since $K<M^{1-\delta}$, \eqref{modulation:eq-lf-p2} contains only low$\times$high-interactions in the $v$-variable. Using our para-product estimate (Lemma \ref{prelim:lem-paraproduct}) and Lemma \ref{modulation:lem-linear}, it follows that 
\begin{align*}
 \Big\| \Big[ \UN[K][-], \IVN[M][-] \Big]_{\leq N}  \Big\|_{\Cprod{r-1}{s}}
 &\lesssim \big\| \UN[K][-] \big\|_{\Cprod{r-1}{\eta}}
\big\| \IVN[M][-] \big\|_{\Cprod{s}{s}}  \\ 
&\lesssim K^{r-s} K^{\eta-s} M^{s-1/2+\eta} 
\Dc  \big\| \UN[K][-] \big\|_{\Cprod{s-1}{s}}. 
\end{align*}
Since $s-\frac{1}{2}+\eta=-\delta_2+\delta_3$, the exponent of $M$ is better than the exponent in the desired estimate \eqref{modulation:eq-lf-p9}, and thus this yields an acceptable contribution.  \\

\emph{Estimate of \eqref{modulation:eq-lf-p3}, \eqref{modulation:eq-lf-p4}, and \eqref{modulation:eq-lf-p5}:} 
To treat the estimates of \eqref{modulation:eq-lf-p3}, \eqref{modulation:eq-lf-p4}, and \eqref{modulation:eq-lf-p5} simultaneously, we let $K,K_u,K_v < M^{1-\delta}$ and let 
\begin{equation}\label{modulation:eq-lf-p8}
\UN[][\ast] \in \Big\{ \UN[K_u,K_v][+-], P^v_{<M^{1-\deltap}} \UN[K][+\fs], P^u_{<M^{1-\deltap}} \UN[K][\fs-] \Big\}. 
\end{equation}
We note that all terms in \eqref{modulation:eq-lf-p8} have $v$-frequencies $\lesssim M^{1-\deltap}\ll M$. Using the low$\times$high-paraproduct estimate (Lemma \ref{prelim:lem-paraproduct}) and Lemma \ref{modulation:lem-linear}, it follows that 
\begin{align*}
\Big\| \Big[ \UN[][\ast], \IVN[M][-] \Big]_{\leq N} \Big\|_{\Cprod{r-1}{s}} 
\lesssim \big\| \UN[][\ast] \big\|_{\Cprod{r-1}{\eta}}
\big\| \IVN[M][-]\big\|_{\Cprod{s}{s}} 
\lesssim M^{s-1/2+\eta} \Dc \big\| \UN[][\ast] \big\|_{\Cprod{r-1}{\eta}}. 
\end{align*}
Using Lemma \ref{ansatz:lem-frequency-support}, Lemma \ref{modulation:lem-bilinear}, and Lemma \ref{modulation:lem-mixed}, it follows that
\begin{alignat*}{3}
\big\| \UN[K_u,K_v][+-] \big\|_{\Cprod{r-1}{\eta}}
&\lesssim (K_u K_v)^{-\eta} \Dc^2, \hspace{10ex} &
\big\| P^v_{<M^{1-\deltap}} \UN[K][+\fs] \big\|_{\Cprod{r-1}{\eta}}
&\lesssim K^{r-s} K^{(1-\deltap)(\eta-s)} \Dc^2, \\
\big\| P^u_{<M^{1-\deltap}} \UN[K][\fs-] \big\|_{\Cprod{r-1}{\eta}}
&\lesssim K^{\eta-s} \Dc^2.
\end{alignat*}
As a result, the total contribution is given by $M^{s-1/2+\eta} \Dc^3$, which is acceptable. \\

\emph{Estimate of \eqref{modulation:eq-lf-p6}:} 
Using our paraproduct estimate (Lemma \ref{prelim:lem-paraproduct}) and Lemma \ref{modulation:lem-linear}, it holds that 
\begin{equation*}
\Big\| 
\Big[ P^{u,v}_{<M^{1-\deltap}}  \UN[][\fs], \IVN[M][-] \Big]_{\leq N}  
\Big\|_{\Cprod{r-1}{s}} 
\lesssim \big\| \UN[][\fs] \big\|_{\Cprod{r-1}{s}} 
\big\|  \IVN[M][-] \big\|_{\Cprod{s}{s}}
\lesssim M^{s-1/2+\eta} \Dc^2,
\end{equation*}
which is acceptable. \\

\emph{Estimate of \eqref{modulation:eq-lf-p7}:} Using our paraproduct estimate (Lemma \ref{prelim:lem-paraproduct}), Lemma \ref{modulation:lem-linear}, and Lemma \ref{modulation:lem-shhl}, it holds that 
\begin{equation*}
\Big\|  \Big[ \SHHLN[M][u], (P_{\leq N}^x)^2 \IVN[M][-] \Big]_{\leq N} \Big\|_{\Cprod{r-1}{s}} 
\lesssim \big\| \SHHLN[M][u] \big\|_{\Cprod{r-1}{s}} 
\big\| \IVN[M][-] \big\|_{\Cprod{s}{s}} \lesssim M^{s-1/2+\eta} \Dc^3, 
\end{equation*}
which is acceptable.  
\end{proof}

As a result of the estimates of this section, we now control all modulated and mixed modulated objects in $\UN[][]$ and $\VN[][]$ in $\Cprod{s-1}{s}$ or $\Cprod{s}{s-1}$, respectively. As a direct consequence, we obtain the following corollary.

\begin{corollary}[\protect{Control of $\LON[ ][\pm]$}]\label{modulation:cor-LON}
Let the pre-modulation hypothesis (Hypothesis \ref{hypothesis:pre}) be satisfied and let $K,M\in \Dyadiclarge$. Then, it holds that 
\begin{equation*}
\big\| \LON[K][-] \big\|_{\Cprod{s}{s-1}} 
+ \big\| \LON[M][+] \big\|_{\Cprod{s-1}{s}} \lesssim \Dc. 
\end{equation*}
\end{corollary}

\begin{proof}
This follows directly from the definitions of $\LON[K][+]$ and $\LON[M][-]$, Lemma \ref{modulation:lem-linear}, Lemma \ref{modulation:lem-bilinear}, Lemma \ref{modulation:lem-mixed}, and Lemma \ref{modulation:lem-linear-reversed}.
\end{proof}

\subsection{Estimates in Cartesian coordinates}\label{section:modulated-sub-Cartesian} 

In the last subsection, we prove regularity estimates of the modulated objects, mixed modulated objects, and nonlinear remainder in Cartesian rather than in null-coordinates. The estimates in Cartesian coordinates will not be needed to treat the modulation equations (Section \ref{section:modulation}), but will be needed to treat the Jacobi errors (Section \ref{section:jacobi}) and perform our globalization argument (Section \ref{section:main}). Our reason for placing the regularity estimates in Cartesian coordinates in this section is that the arguments are similar as in Subsections \ref{section:modulated-sub-linear}-\ref{section:modulated-sub-reversed}.

\begin{lemma}[Estimates in Cartesian coordinates]\label{modulation:lem-Cartesian}
Assume that the post-modulation hypothesis (Hypothesis \ref{hypothesis:post}) is satisfied and let  $K,M\in \Dyadiclarge$. Then, it holds that
\begin{alignat}{12} 
K^\eta \big\| \UN[K][+] \big\|_{\Ctx{s-1}} \lesssim \, && \Dc,&& \quad 
(KM)^\eta \big\| \UN[K,M][+-] \big\|_{\Ctx{\scrr-1}} \lesssim \, && \Dc^2,&& \quad 
M^\eta \big\| \UN[M][-] \big\|_{\Ctx{\scrr-1}} \lesssim \, && \Dc^2,&& \label{modulation:eq-cartesian-e3} \\ 
K^\eta \big\| \UN[K][+\fs] \big\|_{\Ctx{\scrr-1}} \lesssim \, && \Dc^2,&& \quad
M^\eta \big\| \UN[M][\fs-] \big\|_{\Ctx{r-1}} \lesssim \, && \Dc^2,&& \quad  
\big\| \UN[][\fs] \big\|_{\Ctx{r-1}} \lesssim  && \hspace{-2ex} \Dc. && \label{modulation:eq-cartesian-e4} 
\end{alignat}
Furthermore, it also holds that 
\begin{alignat}{12} 
M^\eta \big\| \VN[M][-] \big\|_{\Ctx{s-1}} \lesssim \, && \Dc,&& \quad 
(KM)^\eta \big\| \VN[K,M][+-] \big\|_{\Ctx{\scrr-1}} \lesssim \, && \Dc^2,&& \quad 
K^\eta \big\| \VN[K][+] \big\|_{\Ctx{\scrr-1}} \lesssim \, && \Dc^2,&& \label{modulation:eq-cartesian-e1} \\ 
M^\eta \big\| \VN[M][\fs-] \big\|_{\Ctx{\scrr-1}} \lesssim \, && \Dc^2,&& \quad
K^\eta \big\| \VN[K][+\fs] \big\|_{\Ctx{r-1}} \lesssim \, && \Dc^2,&& \quad  
\big\| \VN[][\fs] \big\|_{\Ctx{r-1}} \lesssim  && \hspace{-2ex} \Dc. && \label{modulation:eq-cartesian-e2} 
\end{alignat}
\end{lemma} 

\begin{remark}
We emphasize that, as already observed in \cite{BLS21}, the bound on $\VN[K,M][+-]$ relies crucially on the probabilistic independence of $W^+$ and $W^-$. 
\end{remark}

\begin{proof}
Due to symmetry in the $u$ and $v$-variables, it suffices to prove the estimates in \eqref{modulation:eq-cartesian-e1} and \eqref{modulation:eq-cartesian-e2}. 
One of the main ingredients in the proof is our trace estimate (Lemma \ref{prelim:lem-trace}), which contains a general and non-resonant case.\\ 

\emph{Estimate of $\VN[M][-]$:}
Due to Lemma \ref{ansatz:lem-frequency-support}, we see that
$\VN[M][-]$ is supported on $u$-frequencies $\lesssim M^{1-\delta+\vartheta}$ and $v$-frequencies $\sim M$. Using the non-resonant trace estimate and Lemma \ref{modulation:lem-linear}, it follows that 
\begin{equation*}
\big\| \VN[M][-] \big\|_{\Ctx{s-1}} \lesssim  \big\| \VN[M][-] \big\|_{\Cprod{s}{s-1}} \lesssim M^{-\eta} \Dc.
\end{equation*}

\emph{Estimate of $\VN[K,M][+-]$:} For the $(+-)$-term, it follows from Definition \ref{ansatz:def-modulated-bilinear} that
\begin{align*}
&\VN[K,M][+-] \\ 
=&\, \coup \chinull[K,M]  P_{\leq N}^x \sum_{u_0,v_0\in \LambdaRR} \int_{\R^2} \dy \dz \, \bigg( \widecheck{\rho}_{\leq N}(y) \widecheck{\rho}_{\leq N}(z)  \psiRuK(x-y-t) \psiRvM(x-y+t) \\ 
\times& \Big[ \Theta^x_y  \Big( \hspace{-0.2ex} \sum_{k\in \Z_K} \rhoND(k) \SN[K][+][k](x-t,x+t) G_k^+ \frac{e^{ik(x-t)}}{ik} \Big), 
\Theta^x_z \Big( \hspace{-0.2ex} \sum_{m\in \Z_M} \rhoND(m) \SN[M][-][m](x-t,x+t) G_m^- e^{im(x+t)}  \Big) \Big] \bigg). 
\end{align*}
Using Lemma \ref{prelim:lem-psi-sum}, Hypothesis \ref{hypothesis:probabilistic}.\ref{ansatz:item-hypothesis-trace}, and the bound on the modulation operators from Hypothesis \ref{hypothesis:pre}, it then follows that
\begin{equation}\label{modulation:eq-cartesian-p1}
\big\| \VN[K,M][+-] \big\|_{C_t^0 \C_x^{\scrr-1}} 
\lesssim \coup \Ac^2 \Bc^2 \max(K,M)^{\scrr-\frac{1}{2}+\eta} K^{-\frac{1}{2}}
= \Dc^2 \max(K,M)^{\scrr-\frac{1}{2}+\eta} K^{-\frac{1}{2}}. 
\end{equation}
Since $\scrr=1-10\delta$ and $K\gtrsim M^{1-\delta}$, it holds that 
\begin{equation}\label{modulation:eq-cartesian-p2} 
\max(K,M)^{\scrr-\frac{1}{2}+\eta} K^{-\frac{1}{2}} \lesssim \max(K,M)^{-10\delta+\frac{\delta}{2}+\eta}
\lesssim (KM)^{-\eta}, 
\end{equation}
which is acceptable. \\

\emph{Estimate of $\VN[K][+]$:} 
Due to Lemma \ref{ansatz:lem-frequency-support}, we see that $\VN[K][+]$ is supported on $u$-frequencies $\sim K$ and $v$-frequencies $\lesssim K^{1-\delta}$. Using the non-resonant trace estimate and Lemma \ref{modulation:lem-linear-reversed}, it then follows that  
\begin{align*}
\big\| \VN[K][+] \big\|_{\Ctx{\scrr-1}} 
&\lesssim \big\| \VN[K][+] \big\|_{\Cprod{\eta}{\scrr-1}}  
\lesssim K^{\eta-s}  K^{\scrr-1-(s-1)} \big\| \VN[K][+] \big\|_{\Cprod{s}{s-1}} \lesssim K^{-\eta} \Dc^2. 
\end{align*}

\emph{Estimate of $\VN[M][\fs-]$:} 
Since $\VN[M][\fs-]$ is supported on $u$-frequencies $\gtrsim M^{1-\deltap}$, we first decompose 
\begin{equation*}
\big\| \VN[M][\fs-] \big\|_{\Ctx{\scrr-1}} \lesssim \Sumlarge_{\substack{ L \gtrsim M^{1-\deltap}}} 
\big\| P_L^u \VN[M][\fs-] \big\|_{\Ctx{\scrr-1}}. 
\end{equation*}
In the case $L\sim M$, we use the general trace estimate and Lemma \ref{modulation:lem-mixed}, which yield that
\begin{equation*}
\big\| P_L^u \VN[M][\fs-] \big\|_{\Ctx{\scrr-1}} \lesssim \big\| P_L^u \VN[M][\fs-] \big\|_{\Cprod{\eta}{\eta}} 
\lesssim M^{\eta-r} M^{\eta-(s-1)} \big\| \VN[M][\fs-] \big\|_{\Cprod{r}{s-1}} \lesssim M^{\eta-r+\eta-(s-1)} \Dc^2. 
\end{equation*}
Since $\eta-r+\eta-(s-1)=-\delta_1+\delta_2+2\delta_3\leq -\delta_3$, this is acceptable. In the case $M^{1-\deltap} \lesssim L\not\sim M$, we use the non-resonant trace estimate and Lemma \ref{modulation:lem-mixed}, which yield that
\begin{align*}
   \big\| P_L^u \VN[M][\fs-] \big\|_{\Ctx{\scrr-1}} \lesssim \big\| P_L^u \VN[M][\fs-] \big\|_{\Cprod{\eta}{\scrr-1}} 
   \lesssim L^{\eta-r} M^{\scrr-1-(s-1)} \big\| P_L^u \VN[M][\fs-] \big\|_{\Cprod{r}{s-1}} 
   \lesssim M^{(1-\deltap)(\eta-r)+\scrr-s} \Dc^2. 
\end{align*}
Since $\scrr=1-10\delta$, it holds that 
\begin{equation*}
(1-\deltap)(\eta-r)+\scrr-s 
= (1-\delta) \big( - \tfrac{1}{2} \big) + \tfrac{1}{2} - 10\delta + \mathcal{O}(\delta_1) = - \big( 10-\tfrac{1}{2}\big) \delta  + \mathcal{O}(\delta_1), 
\end{equation*}
which is acceptable. \\ 

\emph{Estimate of $\VN[K][+\fs]$ and $\VN[][\fs]$:} Using the general case of the trace estimate and Lemma \ref{modulation:lem-mixed}, we obtain that 
\begin{align*}
\big\| \VN[K][+\fs] \big\|_{\Ctx{r-1}} 
&\lesssim  \big\| \VN[K][+\fs] \big\|_{\Cprod{s}{r-1}} 
\lesssim  K^{-\eta} \Dc^2, \\ 
\big\| \VN[][\fs] \big\|_{\Ctx{r-1}} 
&\lesssim \big\| \VN[][\fs] \big\|_{\Cprod{r}{r-1}} 
\lesssim \Dc. \qedhere
\end{align*}
\end{proof}

In our treatment of the Jacobi errors (Section \ref{section:jacobi}), we also need the following refinement of Lemma \ref{modulation:lem-Cartesian}. In this refinement, $\Pbd$ is as in \eqref{prelim:eq-pbd} above, i.e., it is given by 
\begin{equation*}
\Pbd := \sum_{\substack{ N^{1-2\delta_1}\leq L \leq N}} P_L^x. 
\end{equation*}

\begin{lemma}\label{modulation:lem-Cartesian-frequency-boundary}
Assume that the post-modulation hypothesis (Hypothesis \ref{hypothesis:post}) is satisfied. Furthermore, let $K\in \Dyadiclarge$ satisfy $N^{1-2\delta_1}\lesssim K\lesssim N$. Then, it holds that 
\begin{equation}\label{modulation:eq-Cartesian-frequency-boundary}
\Big\| \Pbd \Big( P^v_{>K^{1-\deltap}} \VN[\leq K^{1-\delta}][\fcs] - \VN[][+\fs] - \VN[][\fs] \Big) \Big\|_{C_t^0 \C_x^{\scrr-1}} 
\lesssim N^{-\delta} \Dc. 
\end{equation}
\end{lemma}

\begin{proof}
Since $K\lesssim N$, Lemma \ref{ansatz:lem-frequency-support} implies that $\VN[\leq K^{1-\delta}][\fcs] $ is supported on $u$-frequencies $\lesssim K^{1-\delta} \ll N^{1-2\delta_1}$. Thus, it holds that 
\begin{equation*}
\Pbd  P^v_{\leq K^{1-\deltap}} \VN[\leq K^{1-\delta}][\fcs]=0,  
\end{equation*}
and we can therefore replace $P^v_{>K^{1-\deltap}} \VN[\leq K^{1-\delta}][\fcs]$ in \eqref{modulation:eq-Cartesian-frequency-boundary} by $\VN[\leq K^{1-\delta}][\fcs]$. Using the definition of $\VN[\leq K^{1-\delta}][\fcs]$, it then follows that 
\begin{align*}
&\hspace{1ex}\Big\| \Pbd \Big( P^v_{>K^{1-\deltap}} \VN[\leq K^{1-\delta}][\fcs] - \VN[][+\fs] - \VN[][\fs] \Big) \Big\|_{C_t^0 \C_x^{\scrr-1}}  \\ 
\lesssim&\,  \Sumlarge_{\substack{L\leq \Nd \colon \\ L>K^{1-\delta}}} \big\| \Pbd \VN[L][+\fs] \big\|_{C_t^0 \C_x^{\scrr-1}}  
+ \big\| \Pbd P^u_{>K^{1-\deltap}} \VN[][\fs] \big\|_{C_t^0 \C_x^{\scrr-1}}. 
\end{align*}
Since $\VN[L][+\fs]$ is supported on $u$-frequencies of size $\sim L \gtrsim K^{1-\delta} \gg K^{1-\deltap}$, it therefore suffices to prove for all $\VN[][\ast]\in \Cprod{s}{r-1}$ that
\begin{equation}\label{modulation:eq-Cartesian-frequency-boundary-p1}
\big\|  \Pbd P^u_{>K^{1-\deltap}} \VN[][\ast] \big\|_{C_t^0 \C_x^{\scrr-1}} \lesssim N^{-\delta} \big\| \VN[][\ast] \big\|_{\Cprod{s}{r-1}}. 
\end{equation}
In order to prove \eqref{modulation:eq-Cartesian-frequency-boundary-p1}, we use a dyadic decomposition of $\VN[][\ast]$ into $P_{L_1}^u P_{L_2}^v \VN[][\ast]$, where $L_1 \gtrsim K^{1-\delta}$. If $L_1\nsim L_2$, we use the non-resonant trace estimate (Lemma \ref{prelim:lem-trace}), which yields that
\begin{align*}
\big\|  \Pbd P_{L_1}^u P_{L_2}^v \VN[][\ast] \VN[][\ast] \big\|_{C_t^0 \C_x^{\scrr-1}}
&\lesssim N^{\scrr-1-(r-1)+\eta} \big\| P_{L_1}^u P_{L_2}^v \VN[][\ast] \big\|_{C_t^0 \C_x^{r-1-\eta}} \\ 
&\lesssim N^{\scrr-r+\eta} \big\|  P_{L_1}^u P_{L_2}^v \VN[][\ast] \big\|_{\Cprod{\eta}{r-1-\eta}} \\ 
&\lesssim N^{\scrr-r+\eta} L_1^{\eta-s} L_2^{-\eta} \big\| \VN[][\ast] \big\|_{\Cprod{s}{r-1}}. 
\end{align*}
Since $L_1 \gtrsim K^{1-\delta}$ and $K\geq N^{1-2\delta_1}$, it holds that
\begin{equation*}
N^{\scrr-r+\eta} L_1^{\eta-s} L_2^{-\eta} \lesssim N^{\scrr-r+\eta} N^{(1-\delta)(1-2\delta_1) (2\eta-s)} (L_1 L_2)^{-\eta}.
\end{equation*}
Since 
\begin{equation*}
\scrr - r +\eta + (1-\delta) (1-2\delta_1) (2\eta-s) = 1- 10\delta - \frac{1}{2} + (1-\delta) \big( - \tfrac{1}{2} \big) + \mathcal{O}(\delta_1) = - 10\delta + \tfrac{\delta}{2} + \mathcal{O}(\delta_1), 
\end{equation*}
this is acceptable. If $L_1 \sim L_2$, we use the resonant trace estimate (Lemma \ref{prelim:lem-trace}), which yields that 
\begin{align*}
\big\|  \Pbd P_{L_1}^u P_{L_2}^v \VN[][\ast]  \big\|_{C_t^0 \C_x^{\scrr-1}} 
&\lesssim \big(N^{1-2\delta_1}\big)^{-10\delta} \big\|  \Pbd P_{L_1}^u P_{L_2}^v \VN[][\ast]  \big\|_{C_t^0 \C_x^{0}} \\
&\lesssim N^{-10(1-2\delta_1)\delta} \big\|  P_{L_1}^u P_{L_2}^v \VN[][\ast] \big\|_{\Cprod{s-\eta}{r-1-\eta}} \\ 
&\lesssim N^{-10(1-2\delta_1)\delta} (L_1 L_2)^{-\eta} \big\| \VN[][\ast]\big\|_{\Cprod{s}{r-1}}, 
\end{align*}
which is also acceptable. 
\end{proof}

In the last lemma of this section, we control certain high$\times$high-interaction involving Cartesian coordinates and Cartesian integrals.
\begin{lemma}\label{modulation:lem-Cartesian-high-high}
Assume that the post-modulation hypothesis (Hypothesis \ref{hypothesis:post}) is satisfied. Furthermore, let 
$L_1,L_2\in \dyadic$ satisfy $L_1 \sim L_2$ and $L_1,L_2 \geq 2$. Then, it holds that 
\begin{equation}\label{modulation:eq-Cartesian-high-high}
\begin{aligned}
& \, \Big\| P_{L_1}^x \UN[][] \otimes \Int^x_{0\rightarrow x} P^x_{L_2} \UN[][] \Big\|_{C_t^0 \C_x^{2s-1}} 
+ \Big\| P_{L_1}^x \UN[][] \otimes \Int^x_{0\rightarrow x} P^x_{L_2} \VN[][] \Big\|_{C_t^0 \C_x^{2s-1}} \\
+&\, \Big\| P_{L_1}^x \VN[][] \otimes \Int^x_{0\rightarrow x} P^x_{L_2} \UN[][] \Big\|_{C_t^0 \C_x^{2s-1}} 
+  \Big\| P_{L_1}^x \VN[][] \otimes \Int^x_{0\rightarrow x} P^x_{L_2} \VN[][] \Big\|_{C_t^0 \C_x^{2s-1}} 
\lesssim (L_1 L_2)^{-\eta} \Dc^2.
\end{aligned}
\end{equation}
\end{lemma}

\begin{proof}
We only prove the estimate 
\begin{equation}\label{modulation:eq-Cartesian-hh-p1}
 \Big\| P_{L_1}^x \UN[][] \otimes \Int^x_{0\rightarrow x} P^x_{L_2} \UN[][] \Big\|_{C_t^0 \C_x^{2s-1}}  \lesssim (L_1 L_2)^{-\eta} \Dc^2,
\end{equation}
since the remaining estimates are similar. To this end, we first insert our Ansatz from \eqref{ansatz:eq-UN-rigorous-decomposition}. Using the triangle inequality, it then follows that
\begin{equation}\label{modulation:eq-Cartesian-hh-p2}
 \Big\| P_{L_1}^x \UN[][] \otimes \Int^x_{0\rightarrow x} P^x_{L_2} \UN[][] \Big\|_{C_t^0 \C_x^{2s-1}} \leq 
 \sum_{\ast_1,\ast_2} \Big\|   P_{L_1}^x \UN[][\ast_1] \otimes \Int^x_{0\rightarrow x} P^x_{L_2} \UN[][\ast_2] \Big\|_{C_t^0 \C_x^{2s-1}},
\end{equation}
where $\ast_1,\ast_2 \in \{ +, +-, -, +\fs, \fs-, \fs \}$. In the case $(\ast_1,\ast_2)\neq (+,+)$, Lemma \ref{modulation:lem-Cartesian} implies that $\UN[][\ast_1]$ or $\UN[][\ast_2]$ can be replaced in $C_t^0 \C_x^{r-1}$. Since $r+s>1$, the corresponding term in \eqref{modulation:eq-Cartesian-hh-p2} can be controlled using Lemma \ref{prelim:lem-paraproduct} and Lemma \ref{prelim:lem-integral}. 
It therefore remains to treat the case $(\ast_1,\ast_2)=(+,+)$. Using Definition \ref{ansatz:def-modulated-linear} and Lemma \ref{ansatz:lem-frequency-support}, it suffices to prove for all $K_1,K_2 \in \Dyadiclarge$ satisfying $K_1 \sim K_2 \sim L_1 \sim L_2$ that 
\begin{equation}\label{modulation:eq-Cartesian-hh-p3}
\Big\| P_{L_1}^x \UN[K_1][+] \otimes \Int^x_{0\rightarrow x} P^x_{L_2} \UN[K_2][+] \Big\|_{C_t^0 \C_x^{2s-1}}
\lesssim (L_1 L_2)^{-\eta} \Dc^2.
\end{equation}
Using Lemma \ref{modulation:lem-PNX-modulated} and using a similar argument as in the proof of Lemma \ref{modulation:lem-integration}, we obtain that 
\begin{align}
\Big\|  P_{L_1}^x \UN[K_1][+]  - \hcoup \sum_{u_1 \in \LambdaRR} \sum_{k_1 \in \Z_{K_1}} 
\psiRuKoneone \rho_{L_1}(k_1) \rhoND(k_1) \SN[K_1][+][k_1] G^+_{u_1,k_1} e^{ik_1 u} \Big\|_{C_t^0 \C_x^{r-1}} 
&\lesssim K_1^{-\eta} \Dc, \label{modulation:eq-Cartesian-hh-p4} \\ 
\Big\| \Int^x_{0\rightarrow x} P_{L_2}^x \UN[K_2][+]  - \hcoup \sum_{u_2 \in \LambdaRR} \sum_{k_2 \in \Z_{K_2}} 
\psiRuKtwotwo \rho_{L_2}(k_2) \rhoND(k_2) \SN[K_2][+][k_2] G^+_{u_2,k_2} \frac{e^{ik_2 u}}{ik_2} \Big\|_{C_t^0 \C_x^{r-1}} 
&\lesssim K_2^{-\eta} \Dc. \label{modulation:eq-Cartesian-hh-p5}
\end{align}
The contribution of the error terms coming from the approximations in \eqref{modulation:eq-Cartesian-hh-p4} and \eqref{modulation:eq-Cartesian-hh-p5} to \eqref{modulation:eq-Cartesian-hh-p3} can be controlled using Lemma \ref{prelim:lem-paraproduct} and Lemma \ref{prelim:lem-integral}. As a result, it then suffices to prove\footnote{The proof of the analogue of \eqref{modulation:eq-Cartesian-hh-p6} for the $P_{L_1}^x \UN[][] \otimes \Int^x_{0\rightarrow x} P^x_{L_2} \VN[][]$ and $P_{L_1}^x \VN[][] \otimes \Int^x_{0\rightarrow x} P^x_{L_2} \UN[][]$ in Lemma \ref{modulation:lem-Cartesian-high-high} differs slightly from the argument below. Instead of the trace estimate (Lemma \ref{prelim:lem-trace}), Hypothesis \ref{hypothesis:probabilistic}.\ref{ansatz:item-hypothesis-tensor-product}, and Proposition \ref{killing:prop-resonant}, we then directly use Hypothesis \ref{hypothesis:probabilistic}.\ref{ansatz:item-hypothesis-trace}.} that 
\begin{equation}\label{modulation:eq-Cartesian-hh-p6}
\begin{aligned}
&\, \coup \bigg\| \sum_{u_1,u_2\in \LambdaRR} \sum_{\substack{k_1 \in \Z_{K_1}\\ k_2 \in \Z_{K_2}}}
\Big( \prod_{j=1}^2 \psiRuKj \rho_{L_j}(k_j) \rhoND(k_j) \Big) 
\Big(  \SN[K_1][+][k_1] G^+_{u_1,k_1} \otimes \SN[K_2][+][k_2] G^+_{u_2,k_2} \Big) e^{ik_1 u}  \frac{e^{ik_2 u}}{ik_2}
\bigg\|_{C_t^0 \C_x^{2s-1}} \\ 
\lesssim&\, (K_1 K_2)^{-\eta} \Dc^2. 
\end{aligned}
\end{equation}
By using the trace estimate (Lemma \ref{prelim:lem-trace}) and Lemma \ref{prelim:lem-psi-sum}, we first obtain that
\begin{align}
&\,\textup{LHS of \eqref{modulation:eq-Cartesian-hh-p6}} \notag \\ 
\lesssim&\,  \coup \Big\| \sum_{u_1,u_2\in \LambdaRR} \sum_{\substack{k_1 \in \Z_{K_1}\\ k_2 \in \Z_{K_2}}}
\Big( \prod_{j=1}^2 \psiRuKj \rho_{L_j}(k_j) \rhoND(k_j) \Big) 
\Big(  \SN[K_1][+][k_1] G^+_{u_1,k_1} \otimes \SN[K_2][+][k_2] G^+_{u_2,k_2} \Big) e^{ik_1 u}  \frac{e^{ik_2 u}}{ik_2}
\Big\|_{\Cprod{s}{2s-1}} \notag  \\ 
\lesssim&\, \coup \sup_{u_1,u_2 \in \LambdaRR} 
\Big\|  \sum_{\substack{k_1 \in \Z_{K_1}\\ k_2 \in \Z_{K_2}}}
\Big( \prod_{j=1}^2 \rho_{L_j}(k_j) \rhoND(k_j) \Big) 
\Big(  \SN[K_1][+][k_1] G^+_{u_1,k_1} \otimes \SN[K_2][+][k_2] G^+_{u_2,k_2} \Big) e^{ik_1 u}  \frac{e^{ik_2 u}}{ik_2}
\Big\|_{\Cprod{s}{2s-1}}. \label{modulation:eq-Cartesian-hh-p7}
\end{align}
We now note that the argument on the right-hand side of \eqref{modulation:eq-Cartesian-hh-p6} contains no $\Theta^x_y$-operators. 
Due to \eqref{killing:eq-resonant-e1} from Proposition \ref{killing:prop-resonant}, its resonant part therefore must equal zero. 
Together with Hypothesis \ref{hypothesis:probabilistic}.\ref{ansatz:item-hypothesis-tensor-product}, we then obtain that
\begin{equation}\label{modulation:eq-Cartesian-hh-p7p}
\eqref{modulation:eq-Cartesian-hh-p7} \lesssim \coup \Ac^2 \max(K_1,K_2)^{2s-1-\frac{1}{2}+\eta} K_2^{-\frac{1}{2}} 
\Big\| \rho_{L_1}(k_1) \SN[K_1][+][k_1] \Big\|_{\Wuv[s][s]}
\Big\| \rho_{L_2}(k_2) \SN[K_1][+][k_2] \Big\|_{\Wuv[s][s]}.
\end{equation}
By inserting our choice of parameters from \eqref{prelim:eq-parameter-regularities}, using $K_1 \sim K_2$, and using Lemma \ref{chaos:lem-wc}, it follows that
\begin{equation*}
\eqref{modulation:eq-Cartesian-hh-p7p} \lesssim (K_1 K_2)^{-\delta_2+\frac{\delta_3}{2}}  \Dc^2 ,
\end{equation*}
which is acceptable. 
\end{proof}

\begin{remark}[Estimates with reversed roles of $t$ and $x$]\label{modulation:rem-reversed-roles}
We first note that the estimate from Lemma \ref{killing:lem-Cartesian} is still satisfied after replacing $C_t^0 \C_x^\gamma$ by $C_x^0 \C_t^\gamma$, i.e., after reversing the roles of $t$ and $x$. After including the $C_x^0 \C_t^\gamma$-estimate in 
Hypothesis \ref{hypothesis:probabilistic}.\ref{ansatz:item-hypothesis-trace}, the same estimates as in \eqref{modulation:eq-cartesian-e1}, \eqref{modulation:eq-cartesian-e2}, and \eqref{modulation:eq-Cartesian-high-high} can then be proven 
with $C_t^0 \C_x^\gamma$ replaced by $C_x^0 \C_t^\gamma$, $\Int^x_{0\rightarrow x}$ replaced by $\Int^t_{0\rightarrow t}$, and 
$P_{L_1}^x$ and $P_{L_2}^x$ replaced by $P_{L_1}^t$ and $P_{L_2}^t$. To be more precise, it can be proven that 
\begin{alignat*}{12} 
M^\eta \big\| \VN[M][-] \big\|_{\Cxt{s-1}} \lesssim \, && \Dc,&& \quad 
(KM)^\eta \big\| \VN[K,M][+-] \big\|_{\Cxt{\scrr-1}} \lesssim \, && \Dc^2,&& \quad 
K^\eta \big\| \VN[K][+] \big\|_{\Cxt{\scrr-1}} \lesssim \, && \Dc^2,&&  \\ 
M^\eta \big\| \VN[M][\fs-] \big\|_{\Cxt{\scrr-1}} \lesssim \, && \Dc^2,&& \quad
K^\eta \big\| \VN[K][+\fs] \big\|_{\Cxt{r-1}} \lesssim \, && \Dc^2,&& \quad  
\big\| \VN[][\fs] \big\|_{\Cxt{r-1}} \lesssim  && \hspace{-2ex} \Dc, && 
\end{alignat*}
and 
\begin{align*}
& \, \Big\| P_{L_1}^t \UN[][] \otimes \Int^t_{0\rightarrow t} P^t_{L_2} \UN[][] \Big\|_{\Cxt{2s-1}} 
+ \Big\| P_{L_1}^t \UN[][] \otimes \Int^t_{0\rightarrow t} P^t_{L_2} \VN[][] \Big\|_{\Cxt{2s-1}} \\
+&\, \Big\| P_{L_1}^t \VN[][] \otimes \Int^t_{0\rightarrow t} P^t_{L_2} \UN[][] \Big\|_{\Cxt{2s-1}} 
+  \Big\| P_{L_1}^t \VN[][] \otimes \Int^t_{0\rightarrow t} P^t_{L_2} \VN[][] \Big\|_{\Cxt{2s-1}} 
\lesssim (L_1 L_2)^{-\eta} \Dc^2.
\end{align*}
\end{remark}
\section{Modulation equations}\label{section:modulation}

We now analyze the modulation equations, which were previously introduced in Definition \ref{ansatz:def-modulation-equations}. Using the adjoint map from \eqref{prelim:eq-adjoint-map}, the modulation equations can be written as  
\begin{equation}\label{modulation:eq-motivation}
\begin{cases}
\begin{aligned}
\partial_v \pSN[K][+][k]
&= -  \chinull[K] \Ad \Big(  \rho_{\leq N}^2(k) \LON[K][-]  + 
\rho_{\leq N}^4(k) \SHHLN[K][v] \Big) \circ  \pSN[K][+][k], \\
\partial_u \pSN[M][-][m]
&= \chinull[M]  \Ad \Big(  \rho_{\leq N}^2(m) \LON[M][+]  + 
\rho_{\leq N}^4(m) \SHHLN[M][u] \Big) \circ  \pSN[M][-][M], \\
\pSN[K][+][k]\big|_{v=u} &= \SNin[K][+], \qquad \pSN[M][-][m]\big|_{u=v} = \SNin[M][-].
\end{aligned}
\end{cases}
\end{equation}
 Here, $\LON$ and $\SHHLN$ are as in Definition \ref{ansatz:def-lo} and  Definition \ref{ansatz:def-shhl}, respectively. 
 In the main result of this section (Proposition \ref{modulation:prop-main}), we obtain the well-posedness of the modulation equations and several properties of the solutions. In order to state Proposition \ref{modulation:prop-main}, we first make the following definition. 

\begin{definition}[\protect{The parameter $\Bcin$}]\label{modulation:def-bin}
Let $(\SNin[K][+])_{K\in \dyadic}$ and $(\SNin[M][-])_{M\in \dyadic}$ be initial modulation operators as in Definition \ref{ansatz:def-initial-data}. Then, we define
\begin{equation*}
\Bcin := \sup_{K\in \dyadic} \big\| \SNin[K][+] \big\|_{\C_x^s} + \sup_{M\in \dyadic}\big\| \SNin[M][-] \big\|_{\C_x^s}. 
\end{equation*}
\end{definition}

For most of this article, the reader may pretend that the parameters $\Bc$ and $\Bcin$ from Hypothesis \ref{hypothesis:pre} and Definition \ref{modulation:def-bin} are related via the identity $\Bc=4\Bcin$. The exception is the proof of Proposition \ref{main:prop-null-lwp}, where the quotient of $\Bcin$ and $\Bc$ is used to absorb a constant $C=C(\delta_\ast)$. 

\begin{proposition}[Well-posedness of modulation equations]\label{modulation:prop-main} 
Let $\Nd,N\in \Dyadiclarge$, let $\Ac,\Bc\geq 1$, let $\coup>0$, and let $\Dc= \hcoup \Ac \Bc$.   
In addition, let $(\SNin[K][+])_{K\in \dyadic}$ and $(\SNin[M][-])_{M\in \dyadic}$ be initial modulation operators as in Definition \ref{ansatz:def-initial-data} and let  $\UN[][\fs],\VN[][\fs]\colon \R^{1+1} \rightarrow \frkg$ be nonlinear remainders.
Assume that the probabilistic hypothesis (Hypothesis \ref{hypothesis:probabilistic}) and the conditions 
\begin{equation}\label{modulation:eq-well-posedness-assumption}
\Bcin \leq \tfrac{1}{4} \Bc, \qquad \Dc \leq c, \qquad \text{and} \qquad \big\| \UN[][\fs] \big\|_{\Cprod{r-1}{r}}, \big\| \VN[][\fs] \big\|_{\Cprod{r}{r-1}} \leq \Dc, 
\end{equation}
are satisfied, where $\Bcin$ is as in Definition \ref{modulation:def-bin} and $c=c(\delta_\ast)$ is as in Hypothesis \ref{hypothesis:pre}.
Then, there exist unique pure modulation operators \revision{$( \pSN[K][+][k])_{K\in \Dyadiclarge,k\in \Z_K}\colon \R^{1+1}\rightarrow \End(\frkg)$} and \revision{$( \pSN[M][-][m] )_{M\in \Dyadiclarge,m\in \Z_M}\colon \R^{1+1}\rightarrow \End(\frkg)$} which solve the modulation equations \eqref{modulation:eq-motivation} and which, for all $K,M\in \Dyadiclarge$, satisfy the estimates
\begin{equation}\label{modulation:eq-main-bounds}
\big\| \pSN[K][+][k] \big\|_{\Wuv[s][s][k]} ,  \big\| \pSN[M][-][m]  \big\|_{\Wuv[s][s][m]} \leq 2 \Bcin. 
\end{equation}
Furthermore, we have the following properties:  
\begin{enumerate}[leftmargin=5.5ex,label=(\roman*)]
\item \label{modulation:item-difference} (Frequency-truncation) 
For all $K,M\in \Dyadiclarge$, it holds that 
\begin{equation*}
K^{100} \Big\| \SN[K][+][k] - \pSN[K][+][k] \Big\|_{\Wuv[s][s][k]}, M^{100} \Big\| \SN[M][-][m] - \pSN[M][-][m] \Big\|_{\Wuv[s][s][m]} \lesssim    \Dc \Bcin,
\end{equation*}
\revision{where $\SN[K][+][k]$ and $\SN[M][-][m]$ are as in \eqref{ansatz:eq-pure-to-truncated}.}
\item \label{modulation:item-distance-initial} (Distance to initial data) For all $K,M\in \Dyadiclarge$, it holds that 
\begin{equation*}
\Big\| \SN[K][+][k] - \SNin[K][+] \Big\|_{\Wuv[s][s][k]},
\Big\| \SN[M][-][m] - \SNin[M][-] \Big\|_{\Wuv[s][s][m]} \lesssim  \Dc \Bcin.
\end{equation*}
\item \label{modulation:item-orthogonality} 
(Orthogonality) For all $K,M\in \Dyadiclarge$ satisfying $K,M>N^{1-\delta}$, it holds that 
\begin{equation*}
K^{100} \Big\| \big(\SN[K][+][k]\big)^\ast \SN[K][+][k] - \Id_\frkg \Big\|_{\Wuv[s][s][k]} 
,M^{100} \Big\| \big(\SN[M][-][m]\big)^\ast \SN[M][-][m] - \Id_\frkg \Big\|_{\Wuv[s][s][m]}\lesssim   \Dc \Bcin^2 .
\end{equation*}
\item \label{modulation:item-structure}
(Para-controlled structure) 
The modulation operators can be written as 
\begin{align*}
\SN[K][+][k] 
&= - \chinull[K] \rho_{\leq N}^2(k) \, \Big( \Int^v_{u\rightarrow v} \big( \ad\big(\LON[K][-]\big) \big) \Para[v][gg] \SN[K][+][k] \Big) + \YN[K][+][k] , \\ 
\SN[M][-][m] 
&=  \chinull[M] \rho_{\leq N}^2(m)\, \Big( 
\Int^u_{v\rightarrow u} \big( \ad \big( \LON[M][+] \big) \big) \Para[u][gg] \SN[M][-][m] \Big) +  \YN[M][-][m], 
\end{align*}
where, for all $K,M\in \dyadiclarge$, the nonlinear remainders satisfy 
\begin{equation*}
\Big\| \YN[K][+][k] \Big\|_{\Wuv[s][\scrr][k]} , \Big\| \YN[M][-][m] \Big\|_{\Wuv[\scrr][s][m]} \lesssim \Bcin. 
\end{equation*}
\end{enumerate}
\end{proposition}

\begin{remark}
In earlier work \cite[Section 8]{BLS21}, the first author, L\"{u}hrmann, and Staffilani proved the well-posedness of a similar system of modulation equations. Since the pure modulation operators $(\pSN[K][+][k])_{k\in \Z_K}$ and $(\pSN[M][-][m])_{m\in \Z_M}$ are morally supported on frequencies $\lesssim K^{1-\delta}$ and $\lesssim M^{1-\delta}$, which are much smaller than $N$, our well-posedness argument is closely related to the argument from \cite{BLS21}. In comparison with \cite[Section 8]{BLS21}, the main novelties of Proposition \ref{modulation:prop-main} lie in the $\Wuv[\alpha][\beta]$-norms and orthogonality.
\end{remark}

As already discussed in Section \ref{section:ansatz-modulation}, the modulation equations \eqref{modulation:eq-motivation} cannot be solving using classical methods for ordinary differential equations. The reason is that the $u$-regularity of $\LON[M][+]$ and the $v$-regularity of $\LON[K][-]$ are below $-\frac{1}{2}$. Instead of classical arguments, we rely on the para-controlled approach to rough ordinary differential equations from \cite{GIP15}. We therefore make the Ansatz 
\begin{equation*}
\pSN[K][+][k] := \pXN[K][+][k] + \pYN[K][+][k]
\qquad \text{and} \qquad 
\pSN[M][-][m] := \pXN[M][-][m] + \pYN[M][-][m],
\end{equation*}
where $\pXN[K][+][k]$ and $\pXN[M][-][m]$ are para-controlled components and 
$\pYN[K][+][k]$ and $\pYN[M][-][m]$ are nonlinear remainders. In order to obtain solutions $\pSN[K][+][k]$ and $\pSN[M][-][m]$ of the modulation equations, we would like $\pXN[K][+][k]$, $\pXN[M][-][m]$, $\pYN[K][+][k]$, and $\pYN[M][-][m]$ to be solutions of the para-controlled modulation equations from Definition \ref{ansatz:def-modulation-para}.
This requires control of the terms in \eqref{ansatz:eq-X-p-1}-\eqref{ansatz:eq-Y-m-5} and the corresponding estimates are the subject of the next two subsections. 

\subsection{Resonant term}\label{section:modulation-resonant}

In this subsection, we control the resonant products 
\begin{equation}\label{modulation:eq-resonant-motivation-0}
\Ad \big( \LON[K][-] \big) \Para[v][sim] \Ad \big( P^v_{>1} \Int^{v}_{u\rightarrow v} \LON[K][-] \big)
\qquad \text{and} \qquad 
\Ad \big( \LON[M][+] \big) \Para[u][sim] \Ad \big( P^{u}_{>1} \Int^{u}_{v\rightarrow u} \LON[M][+] \big). 
\end{equation}
The $P^v_{>1}$ and $P^u_{>1}$-operators in \eqref{modulation:eq-resonant-motivation-0} are inserted for technical convenience and allow us to repeatedly use Lemma \ref{prelim:lem-Duhamel-integral}.
The estimates of \eqref{modulation:eq-resonant-motivation-0} 
are the most important ingredients in our estimates of 
\eqref{ansatz:eq-Y-p-2} and \eqref{ansatz:eq-Y-m-2}, i.e., 
the double Duhamel-terms in the para-controlled modulation equations. 

\begin{proposition}[Resonant product]\label{modulation:prop-resonant}
Let the pre-modulation hypothesis (Hypothesis \ref{hypothesis:pre}) be satisfied. Then, it holds that 
\begin{align}
\sup_{K\in \Dyadiclarge} \Big\| \Ad \big( \LON[K][-] \big) \Para[v][sim] \Ad \big( P^v_{>1} \Int^{v}_{u\rightarrow v} \LON[K][-] \big) \Big\|_{\Cprod{s}{\scrr-1}} &\lesssim \Dc^2, \label{modulation:eq-resonant} \\
\sup_{M\in \Dyadiclarge} \Big\| \Ad \big( \LON[M][+] \big) \Para[u][sim] \Ad \big( P^u_{>1} \Int^{u}_{v\rightarrow u} \LON[M][+] \big) \Big\|_{\Cprod{\scrr-1}{s}} 
&\lesssim \Dc^2. 
\end{align}
\end{proposition}

Due to symmetry in the $u$ and $v$-variables, we focus our attention on \eqref{modulation:eq-resonant}. 
After inserting the definition of $\LON[K][-]$, Proposition \ref{modulation:prop-resonant} requires estimates of
\begin{equation}\label{modulation:eq-resonant-motivation}
\Ad \Big(  \VN[][\ast_1] \Big) \Para[v][sim] 
\Ad \Big( P^v_{>1} \Int^v_{u\rightarrow v}  \VN[][\ast_2] \Big),
\end{equation}
where, up to additional frequency-truncations, $\VN[][\ast_1]$ and $\VN[][\ast_2]$ are terms from our Ansatz. The proof of Proposition \ref{modulation:prop-resonant} is split over three different lemmas.

\begin{lemma}\label{modulation:lem-resonant-m}
Let the pre-modulation hypothesis (Hypothesis \ref{hypothesis:pre}) be satisfied and let
$K,L,L_u,L_v,M\in \Dyadiclarge$ satisfy $L,L_u,L_v,M \lesssim N^{1-\delta}$. Furthermore, let 
\begin{equation}\label{modulation:eq-resonant-m}
\VN[][\ast] \in \Big\{ \VN[L][-], \VN[L_u,L_v][+-], \VN[L][+], P_{<K^{1-\deltap}}^{u} \VN[L][\fs-] \Big\}.
\end{equation}
Then, it holds that
\begin{equation}\label{modulation:eq-resonant-m-estimate}
\Big\| \Ad \big(  \VN[M][-] \big) \Para[v][sim] \Ad \big(  P^v_{>1} \Int^v_{u\rightarrow v} \VN[][\ast] \big) \Big\|_{\Cprod{s}{\scrr-1}} \lesssim M^{-\eta} \Gain \big( \VN[][\ast] \big) \Dc^2. 
\end{equation}
\end{lemma}

\begin{proof}
We treat the four cases in \eqref{modulation:eq-resonant-m} separately. Since the $\VN[M][-]$-term in \eqref{modulation:eq-resonant-m-estimate} is supported on $v$-frequencies $\sim M \geq \Nlarge$ and we consider high$\times$high-interactions, the $P^v_{>1} \Int^v_{u\rightarrow v}$-operator in \eqref{modulation:eq-resonant-m-estimate} can always be replaced by $\Int^v_{u\rightarrow v}$, which will be done in the first and second cases below.  \\

\emph{The $(-)$$\times$$(-)$-interaction:}
We estimate 
\begin{align}
&\, \Big\|  \Ad \big(  \VN[M][-] \big)
\Para[v][sim]  \Ad \big( \Int^v_{u\rightarrow v} \VN[L][-] \big) \Big\|_{\Cprod{s}{\scrr-1}} \notag \\
\lesssim&\, \Big\|  \Ad \big(  \VN[M][-] \big)
\Para[v][sim]  \Ad \big( \IVN[L][-] \big) \Big\|_{\Cprod{s}{\scrr-1}} \label{modulation:eq-resonant-m-p1} \\
+&\,  \Big\|  \Ad \big(  \VN[M][-] \big)
\Para[v][sim]  \Ad \Big(  \big( \Int^v_{u\rightarrow v} \VN[L][-] - \IVN[L][-] \big) \Big) \Big\|_{\Cprod{s}{\scrr-1}} 
\label{modulation:eq-resonant-m-p2}. 
\end{align}
We first estimate \eqref{modulation:eq-resonant-m-p1}. In this estimate, we can assume that $L\sim M$, since otherwise \eqref{modulation:eq-resonant-m-p1} equals zero. Using Lemma \ref{prelim:lem-insertion-parasim}, it then follows that 
\begin{align}
\eqref{modulation:eq-resonant-m-p1} 
&\lesssim 
\Big\| \Ad \big( \VN[M][-] \big)
\,  \Ad \big( \IVN[L][-] \big) \Big\|_{\Cprod{s}{\scrr-1}} \label{modulation:eq-resonant-m-p3} \\
&+\,  M^{\scrr-2s} \big\|  \VN[M][-] \big\|_{\Cprod{s}{s-1}} \big\| \IVN[L][-]\big\|_{\Cprod{s}{s}}.  \label{modulation:eq-resonant-m-p4}
\end{align}
Since $\scrr-2s<-\delta$, the contribution of \eqref{modulation:eq-resonant-m-p4} can be controlled using Lemma \ref{modulation:lem-linear}, and it therefore remains to control \eqref{modulation:eq-resonant-m-p3}. 
Using\footnote{While Lemma \ref{killing:lem-tensor-modulated-linear} requires the post-modulation hypothesis, it is clear from the proof that the estimate \eqref{killing:eq-tensor-truncated-1} for $y=z=0$ only requires the pre-modulation hypothesis.} Lemma \ref{killing:lem-tensor-modulated-linear} and that \eqref{modulation:eq-resonant-m-p3} contains no $\Theta^x_y$-operators, it follows that  
\begin{align*}
\Big\| \Ad \big( \VN[M][-] \big)
\,  \Ad \big(  \IVN[L][-] \big) \Big\|_{\Cprod{s}{\scrr-1}}
\lesssim  \max(M,L)^{\scrr-1+\frac{1}{2}+\eta} L^{-\frac{1}{2}} \Dc^2. 
\end{align*}
Since $\scrr=1-10\delta$ and $L\sim M$, this 
yields an acceptable contribution. It now remains to control \eqref{modulation:eq-resonant-m-p2}. Using our paraproduct estimate (Lemma \ref{prelim:lem-paraproduct}), Lemma \ref{modulation:lem-linear}, and Lemma \ref{modulation:lem-integration}, it holds that
\begin{align*}
 &\Big\|  \Ad \big(  \VN[M][-] \big)
\Para[v][sim]  \Ad \Big(  \big( \Int^v_{u\rightarrow v} \VN[L][-] - \IVN[L][-] \big) \Big) \Big\|_{\Cprod{s}{\scrr-1}} \\
\lesssim\, & \big\| \VN[M][-] \big\|_{\Cprod{s}{s-1}} 
\big\| \Int^v_{u\rightarrow v} \VN[L][-] - \IVN[L][-]  \big\|_{\Cprod{s}{r}} 
\lesssim (ML)^{-\eta} \Dc^2, 
\end{align*}
which is acceptable. \\

\emph{The $(-)$$\times$$(+-)$-interaction:} 
In this step we control the stronger $\Cprod{s}{\eta}$-norm. First, we estimate 
\begin{align}
&\Big\|  \Ad \big(  \VN[M][-] \big)
\Para[v][sim]  \Ad \big( \Int^v_{u\rightarrow v} \VN[L_u,L_v][+-] \big) \Big\|_{\Cprod{s}{\eta}} \notag \\
\lesssim& \, \Big\|  \Ad \big(  \VN[M][-] \big)
\Para[v][sim]  \Ad \Big( \Big[ \IUN[L_u][+], \IVN[L_v][-] \Big]_{\leq N} \Big)  \Big\|_{\Cprod{s}{\eta}} 
\label{modulation:eq-resonant-m-p5} \\ 
+& \, \Big\|  \Ad \big(  \VN[M][-] \big)
\Para[v][sim]  \Ad \Big(  \Big( \Int^v_{u\rightarrow v} \VN[L_u,L_v][+-] - \Big[ \IUN[L_u][+], \IVN[L_v][-] \Big]_{\leq N} \Big) \Big)  \Big\|_{\Cprod{s}{\eta}}. 
\label{modulation:eq-resonant-m-p6}
\end{align}
We now estimate \eqref{modulation:eq-resonant-m-p5} and \eqref{modulation:eq-resonant-m-p6} separately and start with \eqref{modulation:eq-resonant-m-p5}. 
Due to Lemma \ref{ansatz:lem-frequency-support}, we can assume that $\max(L_u,L_v)\gtrsim M$, since 
otherwise \eqref{modulation:eq-resonant-m-p5} equals zero. Using Lemma \ref{prelim:lem-paraproduct} and Lemma \ref{modulation:lem-linear}, we obtain
\begin{align*}
\eqref{modulation:eq-resonant-m-p5} 
\lesssim\, & \big\| \VN[M][-] \big\|_{\Cprod{s}{-s+\eta}} 
\Big\| \Big[ \IUN[L_u][+], \IVN[L_v][-]\Big]_{\leq N} \Big\|_{\Cprod{s}{s}} \allowdisplaybreaks[3] \\
\lesssim \, & \big\| \VN[M][-] \big\|_{\Cprod{s}{-s+\eta}}  
\big\| \IUN[L_u][+] \big\|_{\Cprod{s}{s}} 
\big\| \IVN[L_v][-] \big\|_{\Cprod{s}{s}} \\
\lesssim \, & M^{-s+1/2+2\eta} L_u^{\revision{s-1/2}+\eta} L_v^{\revision{s-1/2}+\eta} \Dc^3.
\end{align*}
Since $L_u \simeq_\delta L_v$ and $\max(L_u,L_v)\gtrsim M$, it follows that
\begin{align*}
&\, M^{-s+1/2+2\eta} L_u^{s-1/2+\eta} L_v^{s-1/2+\eta}
\lesssim (M L_u^{-1} L_v^{-1})^{\delta_2} \max(L_u,L_v,M)^{4\delta_3}  \\
\lesssim&\,  \min(L_u,L_v)^{-\delta_2} \max(L_u,L_v,M)^{4\delta_3} 
\lesssim \max(L_u,L_v,M)^{-\delta_2+\delta \delta_2 + 4\delta_3}, 
\end{align*}
which yields an acceptable contribution. It remains to control \eqref{modulation:eq-resonant-m-p6}. Using Lemma \ref{prelim:lem-paraproduct} and Lemma \ref{modulation:lem-integration-bilinear}, it holds that
\begin{align*}
\eqref{modulation:eq-resonant-m-p6} \lesssim
\big\| \VN[M][-] \big\|_{\Cprod{s}{s-1}} 
\Big\|  \Int^v_{u\rightarrow v} \VN[L_u,L_v][+-] - \Big[ \IUN[L_u][+], \IVN[L_v][-] \Big]_{\leq N} \Big\|_{\Cprod{s}{r}} 
\lesssim (ML_uL_v)^{-\eta} \Dc^3,
\end{align*}
which is acceptable. \\ 

\emph{The $(-)$$\times$$(+)$-interaction:} 
In this case we control the stronger $\Cprod{s}{\eta}$-norm. 
Using Lemma \ref{ansatz:lem-frequency-support}, it follows that the resonant product is only non-zero when $L^{1-\delta}\gtrsim M$.
Using our paraproduct estimate (Lemma \ref{prelim:lem-paraproduct}), Duhamel integral estimate (Lemma \ref{prelim:lem-Duhamel-integral}), Lemma \ref{modulation:lem-linear}, and Lemma \ref{modulation:lem-linear-reversed}, it follows that
\begin{align*}
    \Big\|  \Ad \big(  \VN[M][-] \big)
\Para[v][sim]  \Ad \big( P^v_{>1} \Int^v_{u\rightarrow v} \VN[L][+] \big) \Big\|_{\Cprod{s}{\eta}} 
\lesssim&\, \big\| \VN[M][-] \big\|_{\Cprod{s}{-1/2+\eta}} 
\big\| P^v_{>1} \Int^v_{u\rightarrow v} \VN[L][+] \big\|_{\Cprod{s}{1/2+\eta}} \\ \lesssim&\, \big\| \VN[M][-] \big\|_{\Cprod{s}{-
1/2+\eta}} 
\big\|  \VN[L][+] \big\|_{\Cprod{s}{-1/2+\eta}} \\ 
\lesssim& \, M^{2\eta} L^{-\delta (1/2-s)+10\eta} \Dc^3. 
\end{align*}
Since $L^{1-\delta}\gtrsim M$, this yields an acceptable contribution. \\

\emph{The $(-)$$\times$$(\fs-)$-interaction:} In this case we control the stronger $\Cprod{s}{\eta}$-norm. 
Using Lemma \ref{ansatz:lem-frequency-support}, it follows that the resonant product is only non-zero when $L\sim M$. Using our para-product estimate (Lemma \ref{prelim:lem-paraproduct}), Lemma \ref{prelim:lem-Duhamel-integral}, and Lemma \ref{modulation:lem-linear}, it follows that 
\begin{align*}
&\Big\|  \Ad \big(  \VN[M][-] \big)
\Para[v][sim]  \Ad \big( P^v_{>1} \Int^v_{u\rightarrow v} P_{<K^{1-\deltap}}^{u} \VN[L][\fs-] \big) \Big\|_{\Cprod{s}{\eta}} \\
\lesssim\, & \big\| \VN[M][-] \big\|_{\Cprod{s}{\revision{-s+\eta}}}
\big\| P^v_{>1} \Int^v_{u\rightarrow v} P_{<K^{1-\deltap}}^{u} \VN[L][\fs-] \big\|_{\Cprod{s}{s}}
\lesssim M^{1-2s+\eta} \Dc \big\| \VN[L][\fs-] \big\|_{\Cprod{s}{s-1}}.
\end{align*}
Since $\VN[L][\fs-]$ is supported on $u$-frequencies $\gtrsim L^{1-\deltap}$, it follows from Lemma \ref{modulation:lem-mixed} that
\begin{equation*}
 \big\| \VN[L][\fs-] \big\|_{\Cprod{s}{s-1}} 
 \lesssim L^{(1-\deltap)(s-r)} \big\|  \VN[L][\fs-] \big\|_{\Cprod{r}{s-1}}  
 \lesssim L^{(1-\deltap)(s-r)} \Dc. 
\end{equation*}
Since $L\sim M$ and $(1-\deltap)(s-r) = -\delta_1 + \mathcal{O}(\delta \delta_1)$, we obtain an acceptable contribution. 
\end{proof}

\begin{lemma}\label{modulation:lem-two-high-u}  
Let the pre-modulation hypothesis (Hypothesis \ref{hypothesis:pre}) be satisfied. 
Furthermore, let $K,\allowbreak L,\allowbreak L_u,\allowbreak L_v,\allowbreak M,\allowbreak M_u,\allowbreak M_v \in \Dyadiclarge$
and let 
\begin{equation*}
\VN[][\ast_1] \in \Big\{ \VN[M_u,M_v][+-], \VN[M][+], P_{<K^{1-\deltap}}^{u}\VN[M][\fs-] \Big\}
\qquad \text{and} \qquad 
\VN[][\ast_2] \in \Big\{ \VN[L_u,L_v][+-], \VN[L][+], P_{<K^{1-\deltap}}^{u}\VN[L][\fs-] \Big\}. 
\end{equation*}
Then, it holds that 
\begin{equation}\label{modulation:eq-two-high-u}
    \Big\| \Ad \big(  \VN[][\ast_1] \big) \Para[v][sim] \Ad \big(  P^v_{>1} \Int^v_{u\rightarrow v} \VN[][\ast_2] \big) \Big\|_{\Cprod{s}{\eta}} \lesssim \Gain(\VN[][\ast_1]) \Gain \big( \VN[][\ast_2] \big) \Dc^4. 
\end{equation}
\end{lemma}

\begin{proof}[Proof of Lemma \ref{modulation:lem-two-high-u}]
Using our paraproduct estimates (Lemma \ref{prelim:lem-paraproduct}) and integral estimates (Lemma \ref{prelim:lem-Duhamel-integral}), we obtain that 
\begin{align*}
&\Big\| \Ad \big(  \VN[][\ast_1] \big) \Para[v][sim] \Ad \big( P^v_{>1} \Int^v_{u\rightarrow v} \VN[][\ast_2] \big) \Big\|_{\Cprod{s}{\eta}} \allowdisplaybreaks[3]\\
\lesssim\, &  \Big\| \Ad \big(  \VN[][\ast_1] \big) \DPara[ll][sim] \Ad \big( P^v_{>1} \Int^v_{u\rightarrow v} \VN[][\ast_2] \big) \Big\|_{\Cprod{s}{\eta}} 
+  \Big\| \Ad \big(  \VN[][\ast_1] \big) \DPara[sim][sim] \Ad \big(  P^v_{>1} \Int^v_{u\rightarrow v} \VN[][\ast_2] \big) \Big\|_{\Cprod{s}{\eta}} \\ 
&+\,  \Big\| \Ad \big(  \VN[][\ast_1] \big) \DPara[gg][sim] \Ad \big(  P^v_{>1} \Int^v_{u\rightarrow v} \VN[][\ast_2] \big) \Big\|_{\Cprod{s}{\eta}} \allowdisplaybreaks[4]\\
\lesssim\, & \big\| \VN[][\ast_1] \big\|_{\Cprod{\eta}{-s+\eta}} 
\big\| P^v_{>1} \Int^v_{u\rightarrow v} \VN[][\ast_2] \big\|_{\Cprod{\revision{s}}{s}} 
+ \big\| \VN[][\ast_1] \big\|_{\Cprod{s}{s-1}} 
\big\| P^v_{>1} \Int^v_{u\rightarrow v} \VN[][\ast_2] \big\|_{\Cprod{\eta}{1-s+\eta}} \allowdisplaybreaks[3]\\
\lesssim\, &  \big\| \VN[][\ast_1] \big\|_{\Cprod{\eta}{-s+\eta}} 
\big\|  \VN[][\ast_2] \big\|_{\Cprod{\revision{s}}{s-1}} 
+ \big\| \VN[][\ast_1] \big\|_{\Cprod{s}{s-1}} 
\big\|  \VN[][\ast_2] \big\|_{\Cprod{\eta}{-s+\eta}}. 
\end{align*}
Since the options for $\VN[][\ast_1]$ and $\VN[][\ast_2]$ consist of the same types of terms, it then only remains to prove that 
\begin{equation*}
 \big\| \VN[][\ast_1] \big\|_{\Cprod{s}{s-1}} + \big\| \VN[][\ast_1] \big\|_{\Cprod{\eta}{-s+\eta}}   \lesssim \Gain(\VN[][\ast_1]) \Dc^2. 
\end{equation*}
For $\VN[M_u,M_v][+-]$, this follows directly from Lemma \ref{modulation:lem-bilinear}. For $\VN[M][+]$, this follows directly from Lemma \ref{modulation:lem-linear-reversed} and that $\VN[M][+]$ is supported on $u$-frequencies $\sim M$ and $v$-frequencies $\lesssim M$. For $\VN[M][\fs-]$, this follows directly from Lemma \ref{modulation:lem-mixed} and that $\VN[M][\fs-]$ is supported on $u$-frequencies $\gtrsim M^{1-\deltap}$ and $v$-frequencies $\sim M$. 
\end{proof}

\begin{lemma}[Interactions involving $(+\fs)$ or $(\fs)$ in the modulation equations]
\label{modulation:lem-ps-or-s}
Let the pre-modulation hypothesis (Hypothesis \ref{hypothesis:pre}) be satisfied and let $K,L,M,M_u,M_v\in \Dyadiclarge$. Furthermore, let 
\begin{align*}
&\VN[][\ast_1] \in \Big\{ \VN[M][-], \VN[M_u,M_v][+-], \VN[M][+], P_{<K^{1-\deltap}}^{u}\VN[M][\fs-], P_{<K^{1-\deltap}}^{v} \VN[M][+\fs], P_{<K^{1-\deltap}}^{u,v} \VN[][\fs] \Big\} \\ 
\text{and} \quad 
&\VN[][\ast_2] \in \Big\{ P_{<K^{1-\deltap}}^{v} \VN[L][+\fs], P_{<K^{1-\deltap}}^{u,v}\VN[][\fs] \Big\}.
\end{align*}
Then, it holds that 
\begin{equation*}
\begin{aligned}
\Big\| \Ad \big(    \VN[][\ast_1] \big) \Para[v][sim]
\Ad \big(  P^v_{>1} \Int^{v}_{u\rightarrow v}  \VN[][\ast_2] \big)
\Big\|_{\Cprod{s}{\scrr-1}} 
\lesssim  \Gain \big( \VN[][\ast_1] \big) \Gain \big( \VN[][\ast_2] \big) \Dc^3.
\end{aligned}
\end{equation*}
\end{lemma}

\begin{proof}
Using the paraproduct and integral estimates (Lemma \ref{prelim:lem-paraproduct} and Lemma \ref{prelim:lem-Duhamel-integral}), we obtain that 
\begin{align*}
 \Big\| \Ad \big(  \VN[][\ast_1] \big) \Para[v][sim]
\Ad \big( P^v_{>1} \Int^{v}_{u\rightarrow v}  \VN[][\ast_2] \big)
\Big\|_{\Cprod{s}{\scrr-1}}  
\lesssim&\, \big\| \VN[][\ast_1] \big\|_{\Cprod{s}{s-1}} 
\big\| P^v_{>1} \Int^{v}_{u\rightarrow v} \VN[][\ast_2] \big\|_{\Cprod{s}{r}} \\
\lesssim&\, \big\| \VN[][\ast_1] \big\|_{\Cprod{s}{s-1}}  \big\|  \VN[][\ast_2] \big\|_{\Cprod{s}{r-1}}.
\end{align*}
Due to Lemma \ref{modulation:lem-linear}, Lemma \ref{modulation:lem-bilinear}, Lemma \ref{modulation:lem-mixed}, and Lemma \ref{modulation:lem-linear-reversed}, this yields the desired estimates. 
\end{proof}

\begin{figure}
\scalebox{0.9
}{
\begin{tabular}{
!{\vrule width 1pt}>{\centering\arraybackslash}P{1.5cm}
!{\vrule width 1pt}>{\centering\arraybackslash}P{\smallcolwidth}
!{\vrule width 1pt}>{\centering\arraybackslash}P{\smallcolwidth}
!{\vrule width 1pt}>{\centering\arraybackslash}P{\smallcolwidth}
!{\vrule width 1pt}>{\centering\arraybackslash}P{\smallcolwidth}
!{\vrule width 1pt}>{\centering\arraybackslash}P{\smallcolwidth}
!{\vrule width 1pt}>{\centering\arraybackslash}P{\smallcolwidth}
!{\vrule width 1pt}} 
\noalign{\hrule height 1pt} & & & & & & 
 \\[-5.8ex]
 \begin{tabular}{ll}
 & \hspace{-2.25ex}$ V^{\ast_2}$ \\[-2.5ex]
 $V^{\ast_1} \hspace{1ex}$ &  
 \end{tabular}
&  $(-)$ & $(+-)$ & $(+)$   & $(\fs-)$ & $(+\fs)$ & $(\fs)$  
\\[4pt] \noalign{\hrule height 1pt} \rule{0pt}{14pt}
$(-)$ 
& \cellcolor{Green!30} \ref{modulation:lem-resonant-m}
& \cellcolor{Green!30} \ref{modulation:lem-resonant-m}
& \cellcolor{Green!30} \ref{modulation:lem-resonant-m}
& \cellcolor{Green!30} \ref{modulation:lem-resonant-m}
& \cellcolor{magenta!30}  \ref{modulation:lem-ps-or-s}
& \cellcolor{magenta!30}  \ref{modulation:lem-ps-or-s}  
\\[4pt] \noalign{\hrule height 1pt} \rule{0pt}{14pt}
$(+-)$
& \cellcolor{Gray!80} 
& \cellcolor{Orange!30} \ref{modulation:lem-two-high-u}
& \cellcolor{Orange!30} \ref{modulation:lem-two-high-u}
& \cellcolor{Orange!30} \ref{modulation:lem-two-high-u}
& \cellcolor{magenta!30}  \ref{modulation:lem-ps-or-s}  
& \cellcolor{magenta!30}  \ref{modulation:lem-ps-or-s}  
\\[4pt] \noalign{\hrule height 1pt} \rule{0pt}{14pt}
$(+)$ 
& \cellcolor{Gray!80} 
& \cellcolor{Gray!80} 
& \cellcolor{Orange!30} \ref{modulation:lem-two-high-u}
& \cellcolor{Orange!30} \ref{modulation:lem-two-high-u}
& \cellcolor{magenta!30}  \ref{modulation:lem-ps-or-s}  
& \cellcolor{magenta!30}  \ref{modulation:lem-ps-or-s}  
\\[4pt] \noalign{\hrule height 1pt} \rule{0pt}{14pt}
$(\fs-)$ 
& \cellcolor{Gray!80} 
& \cellcolor{Gray!80} 
& \cellcolor{Gray!80} 
& \cellcolor{Orange!30} \ref{modulation:lem-two-high-u}
& \cellcolor{magenta!30}  \ref{modulation:lem-ps-or-s}  
& \cellcolor{magenta!30}  \ref{modulation:lem-ps-or-s}  
\\[4pt] \noalign{\hrule height 1pt} \rule{0pt}{14pt}
$(+\fs)$ 
& \cellcolor{Gray!80} 
& \cellcolor{Gray!80} 
& \cellcolor{Gray!80} 
& \cellcolor{Gray!80} 
& \cellcolor{magenta!30}  \ref{modulation:lem-ps-or-s}  
& \cellcolor{magenta!30}  \ref{modulation:lem-ps-or-s}  
\\[4pt] \noalign{\hrule height 1pt} \rule{0pt}{14pt}
$(\fs)$ 
& \cellcolor{Gray!80} 
& \cellcolor{Gray!80} 
& \cellcolor{Gray!80} 
& \cellcolor{Gray!80} 
& \cellcolor{Gray!80} 
& \cellcolor{magenta!30}  \ref{modulation:lem-ps-or-s}   
\\[3pt] \noalign{\hrule height 1pt} 
\end{tabular}
}
\caption{\small{
In this figure we display the different cases in the proof of Proposition \ref{modulation:prop-resonant}. Due to Lemma \ref{prelim:lem-integration-by-parts}, the estimates of \eqref{modulation:eq-resonant-motivation} are morally symmetric in $\VN[][\ast_1]$ and $\VN[][\ast_2]$, and it therefore suffices to treat the cases on or above the diagonal. In each list, we list the number of the lemma which controls the corresponding interaction. For visual purposes, we also colored cells containing the same lemma in the same color. 
}}
\label{figure:resonant-cases}
\end{figure}

Equipped with Lemma \ref{modulation:lem-resonant-m}, Lemma \ref{modulation:lem-two-high-u}, and Lemma \ref{modulation:lem-ps-or-s}, we can now prove Proposition \ref{modulation:prop-resonant}. 

\begin{proof}[Proof of Proposition \ref{modulation:prop-resonant}:] 
Due to symmetry in the $u$ and $v$-variables, it suffices to prove \eqref{modulation:eq-resonant}. After inserting the definition of $\LON[K][-]$, it suffices to control 
\begin{equation*}
\Ad \big(  \VN[][\ast_1] \big) \Para[v][sim] 
\Ad \big( P^v_{>1} \Int^v_{u\rightarrow v}  \VN[][\ast_2] \big),
\end{equation*}
where 
\begin{align}
\VN[][\ast_1] 
&\in \Big\{ \VN[M][-], \VN[M_u,M_v][+-], \VN[M][+], P^u_{<K^{1-\deltap}} \VN[M][\fs-], P^v_{<K^{1-\deltap}}\VN[M][+\fs], P^{u,v}_{<K^{1-\deltap}}
\VN[][\fs] \Big\}, \label{modulation:eq-resonant-p1} \\
\VN[][\ast_1] 
&\in \Big\{ \VN[L][-], \VN[L_u,L_v][+-], \VN[L][+], P^u_{<K^{1-\deltap}} \VN[L][\fs-], P^v_{<K^{1-\deltap}}\VN[L][+\fs], P^{u,v}_{<K^{1-\deltap}}
\VN[][\fs] \Big\}, \label{modulation:eq-resonant-p2}
\end{align}
and $L,L_u,L_v,M,M_u,M_v\lesssim K^{1-\delta} \lesssim N^{1-\delta}$.  
The different cases resulting from \eqref{modulation:eq-resonant-p1} and \eqref{modulation:eq-resonant-p2} are illustrated in Figure \ref{figure:resonant-cases}. 
Using Lemma \ref{prelim:lem-integration-by-parts}, it is straightforward to show that
\begin{equation*}
\begin{aligned}
&\bigg\| \Ad \big(  \VN[][\ast_1] \big) \Para[v][sim] 
\Ad \big( P^v_{>1} \Int^v_{u\rightarrow v}  \VN[][\ast_2] \big)
+ \Big(  \Ad \big(  \VN[][\ast_2] \big) \Para[v][sim] 
\Ad \big( P^v_{>1} \Int^v_{u\rightarrow v}  \VN[][\ast_1] \big) \Big)^\ast \bigg\|_{\Cprod{s}{\scrr-1}} \\
\lesssim&\, \Gain\big( \VN[][\ast_1]\big) 
 \Gain\big( \VN[][\ast_2]\big) \Dc^3,
\end{aligned}
\end{equation*}
where $(\hdots)^\ast$ denotes the Hermitian \revision{conjugate}. 
Thus, it suffices to treat the cells in Figure \ref{figure:resonant-cases} on or above the diagonal. For each of these cells, the contribution is estimated by using the lemma listed in the cell. 
\end{proof}

\subsection{Nonlinear estimates}\label{section:modulation-nonlinear} 

In Section \ref{section:modulated-mixed} and Subsection \ref{section:modulation-resonant}, we estimated the modulated objects, mixed modulated objects, and the resonant term from \eqref{modulation:eq-resonant-motivation-0}. Equipped with these estimates, we can now control all nonlinear terms in the para-controlled modulation equations.
In the first lemma, we control the right-hand side of \eqref{ansatz:eq-X-p-1}. 

\begin{lemma}\label{modulation:lem-X-estimate}
Let the pre-modulation hypothesis (Hypothesis \ref{hypothesis:pre}) be satisfied and let $K\in \Dyadiclarge$.  Then, it holds that
\begin{equation}\label{modulation:eq-X-estimate}
\Big\| \rho_{\leq N}^2(k) \Big( \Int^v_{u\rightarrow v} \Ad\big( \LON[K][-] \big) \Big) \Para[v][gg] \pSN[K][+][k] \Big\|_{\Wuv[s][s][k]} 
\lesssim \Dc\,  \big\| \pSN[K][+] \big\|_{\Wuv[s][s][k]}.
\end{equation}
\end{lemma}

\begin{proof}
Due to the $\Para[v][gg]$-operator, we may replace $\Int^v_{u\rightarrow v}$-operator in \eqref{modulation:eq-X-estimate} by $P^v_{>1} \Int^v_{u\rightarrow v}$. 
Using the para-product estimate (Lemma \ref{prelim:lem-paraproduct}), Duhamel integral estimate (Lemma \ref{prelim:lem-Duhamel-integral}), and Lemma \ref{chaos:lem-wc}, it holds that
\begin{align*}
\Big\| \rho_{\leq N}^2(k) \Big( P^v_{>1} \Int^v_{u\rightarrow v} \Ad\big( \LON[K][-] \big) \Big) \Para[v][gg] \pSN[K][+][k] \Big\|_{\Wuv[s][s][k]} 
\lesssim& \, 
\Big\|  P^v_{>1} \Int^v_{u\rightarrow v} \LON[K][-]   \Big\|_{\Cprod{s}{s}}
\Big\| \rho_{\leq N}^2(k) \pSN[K][+][k] \Big\|_{\Wuv[s][s][k]}  \\
\lesssim& \, \big\| \LON[K][-]  \big\|_{\Cprod{s}{s-1}}
\big\| \rho_{\leq N}^2(k) \big\|_{\Wc_k} \big\| \pSN[K][+][k] \big\|_{\Wuv[s][s][k]}. 
\end{align*}
Using Lemma \ref{chaos:lem-wc} and Corollary \ref{modulation:cor-LON}, it holds that $\| \rho_{\leq N}^2(k) \|_{\Wc_k}\lesssim 1$ and $\| \LON[K][-]  \|_{\Cprod{s}{s-1}} \lesssim \Dc$, and thus we obtain an acceptable contribution. 
\end{proof}

In the next lemma, we control \eqref{ansatz:eq-Y-p-1}, i.e., the integral commutator term.

\begin{lemma}\label{modulation:lem-Y-integral-commutator}
Let the pre-modulation hypothesis (Hypothesis \ref{hypothesis:pre}) be satisfied and let $K\in \Dyadiclarge$. Then, it holds that
\begin{align*}
&\bigg\|  \rho_{\leq N}^2(k) \, 
     \Int^{v}_{u\rightarrow v}  \Big( \chinull[K] \Big( \Ad \big( \LON[K][-] \big) \Para[v][gg] \pSN[K][+][k] \Big) \Big)
     - \chinull[K]   \rho_{\leq N}^2(k) 
     \Big( \Big( \Int^{v}_{u\rightarrow v} \Ad \big( \LON[K][-] \big) \Big) \Para[v][gg] \pSN[K][+][k]
     \Big) \Big)  \bigg\|_{\Wuv[s][\scrr][k]} \\
     &\lesssim \Dc \, \big\| \pSN[K][+] \big\|_{\Wuv[s][s][k]}.
\end{align*}
\end{lemma}

\begin{proof} 
We first decompose
\begin{align}
&\rho_{\leq N}^2(k) \, 
     \Int^{v}_{u\rightarrow v}  \Big( \chinull[K] \Big( \Ad \big( \LON[K][-] \big) \Para[v][gg] \pSN[K][+][k] \Big) \Big)
     - \chinull[K]   \rho_{\leq N}^2(k) 
     \Big( \Big( \Int^{v}_{u\rightarrow v} \Ad \big( \LON[K][-] \big) \Big) \Para[v][gg] \pSN[K][+][k]
     \Big) \Big) \notag  \\ 
=&\,\rho_{\leq N}^2(k) \, 
     \bigg( \Int^{v}_{u\rightarrow v}  \Big( \chinull[K] \Big( \Ad \big( \LON[K][-] \big) \Para[v][gg] \pSN[K][+][k] \Big) \Big) 
      - \chinull[K]  \Int^{v}_{u\rightarrow v}  \Big(  \Ad \big( \LON[K][-] \big) \Para[v][gg] \pSN[K][+][k] \Big) \bigg) 
      \label{modulation:eq-Y-integral-commutator-q1} \\ 
+&\, \chinull[K] \rho_{\leq N}^2(k) 
\bigg( \Int^{v}_{u\rightarrow v}  \Big(  \Ad \big( \LON[K][-] \big) \Para[v][gg] \pSN[K][+][k] \Big) 
- \Big( \Int^{v}_{u\rightarrow v}   \Ad \big( \LON[K][-] \big) \Big) \Para[v][gg] \pSN[K][+][k] \Big) \bigg).
  \label{modulation:eq-Y-integral-commutator-q2} 
\end{align}
We first treat \eqref{modulation:eq-Y-integral-commutator-q1}. Using Lemma \ref{prelim:lem-commutator-integral} and Lemma \ref{chaos:lem-wc}, 
it holds that 
\begin{align}
\big\| \eqref{modulation:eq-Y-integral-commutator-q1} \big\|_{\Wuv[s][\scrr][k]}
&\lesssim \Big\| \rho_{\leq N}^2(k) \Big( \ad \big( \LON[K][-] \big) \Para[v][gg] \pSN[K][+][k] \Big) \Big\|_{\Wuv[s][s-1][k]} \notag \\
&\lesssim \big\| \rho_{\leq N}^2(k) \big\|_{\Wc_k} 
\Big\| \ad \big( \LON[K][-] \big) \Para[v][gg] \pSN[K][+][k] \Big\|_{\Wuv[s][s-1][k]}. \label{modulation:eq-Y-integral-commutator-q3} 
\end{align}
The first factor in \eqref{modulation:eq-Y-integral-commutator-q3} is controlled due to Lemma \ref{chaos:lem-wc}. Using Lemma \ref{prelim:lem-paraproduct} and Corollary \ref{modulation:cor-control-combined}, it then follows that
\begin{equation*}
    \Big\| \ad \big( \LON[K][-] \big) \Para[v][gg] \pSN[K][+][k] \Big\|_{\Wuv[s][s-1][k]}
    \lesssim \big\| \LON[K][-] \big\|_{\Cprod{s}{s-1}} \big\| \pSN[K][+][k] \big\|_{\Wuv[s][s][k]}
    \lesssim \Dc \big\| \pSN[K][+][k] \big\|_{\Wuv[s][s][k]}. 
\end{equation*}
Thus, the contribution of \eqref{modulation:eq-Y-integral-commutator-q1} is acceptable, and it remains to treat \eqref{modulation:eq-Y-integral-commutator-q2}. Using Lemma \ref{prelim:lem-commutator-integral} and using Lemma \ref{chaos:lem-wc}, we obtain that
\begin{align*}
\Big\|  \eqref{modulation:eq-Y-integral-commutator-q2} \Big\|_{\Wuv[s][\scrr][k]} 
     \lesssim 
     \Big\| \Ad \big( \LON[K][-] \big) \Big\|_{\Cprod{s}{s-1}} \Big\| \rho_{\leq N}^2(k) \pSN[K][+][k] \Big\|_{\Wuv[s][s][k]}
     \lesssim  \big\|  \LON[K][-]  \big\|_{\Cprod{s}{s-1}}
     \big\|  \rho_{\leq N}^2(k) \big\|_{\Wc_k}
     \big\| \pSN[K][+][k] \big\|_{\Wuv[s][s][k]}. 
\end{align*}
Using Lemma \ref{chaos:lem-wc} and Corollary \ref{modulation:cor-LON}, it holds that $\| \rho_{\leq N}^2(k) \|_{\Wc_k}\lesssim 1$ and $\| \LON[K][-]  \|_{\Cprod{s}{s-1}} \lesssim \Dc$, and thus this yields an acceptable contribution. 
\end{proof}

We now estimate the term in the integral equation for $\pYN[K][+][k]$ which came from the double Duhamel trick, i.e., the term in \eqref{ansatz:eq-Y-p-2}. 

\begin{lemma}\label{modulation:lem-Y-double-duhamel}
Let the pre-modulation hypothesis (Hypothesis \ref{hypothesis:pre}) be satisfied and let $K\in \Dyadiclarge$. Then, it holds that
\begin{equation*}
\bigg\| \rho_{\leq N}^4(k) \Ad \big( \LON[K][-] \big) \Para[v][sim] \Big( \chinull[K] \Big( \Int^v_{u\rightarrow v} \big( \Ad\big( \LON[K][-] \big) \big) \Para[v][gg] \pSN[K][+][k] \Big)  \Big)
     \bigg\|_{\Wuv[s][\scrr-1][k]} \lesssim \Dc^2 \, \big\| \pSN[K][+] \big\|_{\Wuv[s][s][k]}. 
\end{equation*}
\end{lemma}

\begin{proof}
We split the argument into three steps. \\ 

\emph{First step:} In the first step, we prove that 
\begin{equation}\label{modulation:eq-double-Duhamel-q1}
\begin{aligned}
\bigg\|& \rho_{\leq N}^4(k) \Ad \big( \LON[K][-] \big) \Para[v][sim] \Big( \chinull[K] \Big( \Int^v_{u\rightarrow v} \big( \Ad\big( \LON[K][-] \big) \big) \Para[v][gg] \pSN[K][+][k] \Big)  \Big)  \\
&- 
\chinull[K] \rho_{\leq N}^4(k) \Ad \big( \LON[K][-] \big) \Para[v][sim] \Big( \Int^v_{u\rightarrow v} \big( \Ad\big( \LON[K][-] \big) \big) \Para[v][gg] \pSN[K][+][k] \Big)  
     \bigg\|_{\Wuv[s][\scrr-1][k]}  
     \lesssim \Dc^2 \big\| \pSN[K][+] \big\|_{\Wuv[s][s][k]}. 
\end{aligned}
\end{equation}
Using Lemma \ref{prelim:lem-chi-commutator}, Lemma \ref{chaos:lem-wc}, and Corollary \ref{modulation:cor-control-combined}, it holds that
\begin{align*}
\textup{LHS of } \eqref{modulation:eq-double-Duhamel-q1}
&\lesssim \big\| \rho_{\leq N}^4(k) \big\|_{\Wc_k} \big\| \chinull[K] \big\|_{\Cprod{s}{s}} 
\big\| \ad \big( \LON[K][-] \big) \big\|_{\Cprod{s}{s-1}} 
\Big\| \Int^v_{u\rightarrow v} \big( \Ad\big( \LON[K][-] \big) \big) \Para[v][gg] \pSN[K][+][k] \Big\|_{\Wuv[s][s][k]} \\
&\lesssim \Dc \Big\| \Int^v_{u\rightarrow v} \big( \Ad\big( \LON[K][-] \big) \big) \Para[v][gg] \pSN[K][+][k] \Big\|_{\Wuv[s][s][k]}.
\end{align*}
Using the definition of $\Para[v][gg]$, our para-product estimate (Lemma \ref{prelim:lem-paraproduct}), and Duhamel integral estimate (Lemma \ref{prelim:lem-Duhamel-integral}), it then holds that 
\begin{align*}
 \Big\| \Int^v_{u\rightarrow v} \big( \Ad\big( \LON[K][-] \big) \big) \Para[v][gg] \pSN[K][+][k] \Big\|_{\Wuv[s][s][k]}
 &\lesssim \big\| P^v_{>1} \Int^v_{u\rightarrow v}  \LON[K][-] \big\|_{\Cprod{s}{s}}
 \big\| \pSN[K][+][k] \big\|_{\Wuv[s][s][k]}  \\
 &\lesssim \big\|  \LON[K][-] \big\|_{\Cprod{s}{s-1}}
 \big\| \pSN[K][+][k] \big\|_{\Wuv[s][s][k]}.
\end{align*}
Together with Corollary \ref{modulation:cor-control-combined}, this completes the proof of \eqref{modulation:eq-double-Duhamel-q1}. \\

\emph{Second step:} In the second step, we prove that 
\begin{equation}\label{modulation:eq-double-Duhamel-p0}
\begin{aligned}
\bigg\| &\rho_{\leq N}^4(k) \Ad \big( \LON[K][-] \big) \Para[v][sim] \Big( \Int^v_{u\rightarrow v} \big( \Ad\big( \LON[K][-] \big) \big) \Para[v][gg] \pSN[K][+][k] \Big)  \\
&- \rho_{\leq N}^4(k) 
\Big( \Ad \big( \LON[K][-] \big) \Para[v][sim] P^v_{>1}  \Int^v_{u\rightarrow v} \big( \Ad\big( \LON[K][-] \big) \big) \Big) \pSN[K][+][k]
     \bigg\|_{\Wuv[s][\scrr-1][k]} \lesssim \Dc^2 \big\| \pSN[K][+] \big\|_{\Wuv[s][s][k]}. 
\end{aligned}
\end{equation}
Using the definition of $\Para[v][gg]$, we may replace the first $\Int^v_{u\rightarrow v}$-operator in \eqref{modulation:eq-double-Duhamel-p0} by $P^v_{>1} \Int^v_{u\rightarrow v}$. 
Using the trilinear para-product estimate from Lemma \ref{prelim:lem-para-product-trilinear}, the integral estimate (Lemma \ref{prelim:lem-Duhamel-integral}),  and Lemma \ref{chaos:lem-wc}, we then have that 
\begin{align}
& \bigg\| \rho_{\leq N}^4(k) \Ad \big( \LON[K][-] \big) \Para[v][sim] \Big( P^v_{>1} \Int^v_{u\rightarrow v} \big( \Ad\big( \LON[K][-] \big) \big) \Para[v][gg] \pSN[K][+][k] \Big)  \notag \\
&\hspace{3ex} - \rho_{\leq N}^4(k) 
\Big( \Ad \big( \LON[K][-] \big) \Para[v][sim]  P^v_{>1} \Int^v_{u\rightarrow v} \big( \Ad\big( \LON[K][-] \big) \big) \Big) \pSN[K][+][k]
     \bigg\|_{\Wuv[s][\scrr-1][k]}  \notag \\ 
    \lesssim& \, \Big\|  \Ad \big( \LON[K][-] \big) \Big\|_{\Cprod{s}{s-1}} 
    \Big\| P^v_{>1} \Int^v_{u\rightarrow v} \Ad \big( \LON[K][-] \big) \Big\|_{\Cprod{s}{s}}
    \Big\| \rho_{\leq N}^4(k) \pSN[K][+][k] \Big\|_{\Wuv[s][s][k]} \notag \\
    \lesssim& \, \big\|  \LON[K][-]  \big\|_{\Cprod{s}{s-1}}^2 
    \big\|  \rho_{\leq N}^4(k) \big\|_{\Wc_k}
     \big\| \pSN[K][+][k] \big\|_{\Wuv[s][s][k]}. \label{modulation:eq-double-Duhamel-p1}
\end{align}
Using Lemma \ref{chaos:lem-wc} and Corollary \ref{modulation:cor-LON}, we obtain that 
$ \|  \rho_{\leq N}^4(k)\|_{\Wc_k}\lesssim 1$ and  $\| \LON[K][-]  \|_{\Cprod{s}{s-1}} \lesssim \Dc$. 
 This completes the proof of \eqref{modulation:eq-double-Duhamel-p0} and hence the second step.\\

\emph{Third step:} In the third step, we prove that 
\begin{equation}\label{modulation:eq-double-Duhamel-p2}
\bigg\|  \rho_{\leq N}^4(k) 
\Big( \Ad \big( \LON[K][-] \big) \Para[v][sim]  P^v_{>1} \Int^v_{u\rightarrow v} \big( \Ad\big( \LON[K][-] \big) \big) \Big) \pSN[K][+][k]
     \bigg\|_{\Wuv[s][\scrr-1][k]} \lesssim \Dc^2 \, \big\| \pSN[K][+] \big\|_{\Wuv[s][s][k]}. 
\end{equation}
Using a product estimate (Lemma \ref{prelim:lem-paraproduct}), which requires that  $\scrr-1+s>0$, and Lemma \ref{chaos:lem-wc}, we obtain that
\begin{align}
&\bigg\|  \rho_{\leq N}^4(k) 
\Big( \Ad \big( \LON[K][-] \big) \Para[v][sim]  P^v_{>1} \Int^v_{u\rightarrow v} \big( \Ad\big( \LON[K][-] \big) \big)  \Big) \pSN[K][+][k]
     \bigg\|_{\Wuv[s][\scrr-1][k]}  \notag \\
\lesssim&\, \Big\|  
 \Ad \big( \LON[K][-] \big) \Para[v][sim]  P^v_{>1} \Int^v_{u\rightarrow v} \big( \Ad\big( \LON[K][-] \big) \big)    \Big\|_{\Cprod{s}{\scrr-1}} \Big\| \rho_{\leq N}^4(k) \pSN[K][+][k]
     \Big\|_{\Wuv[s][s][k]} \notag \\
     \lesssim&\,  \Big\|  
 \Ad \big( \LON[K][-] \big) \Para[v][sim]  P^v_{>1} \Int^v_{u\rightarrow v} \big( \Ad\big( \LON[K][-] \big) \big)   \Big\|_{\Cprod{s}{\scrr-1}}  \big\|  \rho_{\leq N}^4(k) \big\|_{\Wc_k}
     \big\| \pSN[K][+][k] \big\|_{\Wuv[s][s][k]}. \label{modulation:eq-double-Duhamel-p3}
\end{align}
Using Lemma \ref{chaos:lem-wc} and Proposition \ref{modulation:prop-resonant}, we obtain that 
\begin{equation*}
\big\|  \rho_{\leq N}^4(k) \big\|_{\Wc_k} \lesssim 1 \qquad \text{and} \qquad \Big\|  
 \Ad \big( \LON[K][-] \big) \Para[v][sim]  P^v_{>1} \Int^v_{u\rightarrow v} \big( \Ad\big( \LON[K][-] \big)\big)    \Big\|_{\Cprod{s}{\scrr-1}} \lesssim \Dc^2. 
\end{equation*}
This completes the proof of \eqref{modulation:eq-double-Duhamel-p2} and hence the proof of this lemma. 
\end{proof}

In the next lemma, we estimate the remaining terms in the para-controlled modulation equation for $\pYN[K][+][k]$, i.e., \eqref{ansatz:eq-Y-p-3}, \eqref{ansatz:eq-Y-p-4}, and \eqref{ansatz:eq-Y-p-5}.

\begin{lemma}\label{modulation:lem-Y-remaining} 
Let the pre-modulation hypothesis (Hypothesis \ref{hypothesis:pre}) be satisfied, let $K\in \Dyadiclarge$, and let $(\pYN[K][+][k])_{k\in \Z_K}\colon \R^{1+1}\rightarrow \frkg$. Then, we have the following estimates: 
\begin{align}
\Big\| \rho_{\leq N}^2(k) \Ad \big( \LON[K][-] \big) \Para[v][sim] \pYN[K][+][k] \Big\|_{\Wuv[s][\scrr-1][k]} &\lesssim \Dc \, \big\| \pYN[K][+][k] \big\|_{\Wuv[s][\scrr][k]}
\label{modulation:eq-Y-remaining-e1}, \\
\Big\| \rho_{\leq N}^2(k) \Ad \big( \LON[K][-] \big) \Para[v][ll] \pSN[K][+][k] \Big\|_{\Wuv[s][\scrr-1][k]} &\lesssim \Dc 
\, \big\| \pSN[K][+] \big\|_{\Wuv[s][s][k]}
\label{modulation:eq-Y-remaining-e2}, \\
\Big\| \rho_{\leq N}^4(k) \Ad \big( \SHHLN[K][v] \big) \, \pSN[K][+][k] \Big\|_{\Wuv[s][\scrr-1][k]} &\lesssim \Dc^2 
\, \big\| \pSN[K][+] \big\|_{\Wuv[s][s][k]}
\label{modulation:eq-Y-remaining-e3}. 
\end{align}
\end{lemma}

\begin{proof}
We prove \eqref{modulation:eq-Y-remaining-e1}, \eqref{modulation:eq-Y-remaining-e2}, and \eqref{modulation:eq-Y-remaining-e3} separately. \\

\emph{Proof of \eqref{modulation:eq-Y-remaining-e1}:} Using the para-product estimate (Lemma \ref{prelim:lem-paraproduct}), which requires that $s-1+\scrr>0$, and Lemma \ref{chaos:lem-wc}, we obtain that
\begin{align}
\Big\| \rho_{\leq N}^2(k) \Ad \big( \LON[K][-] \big) \Para[v][sim] \pYN[K][+][k] \Big\|_{\Wuv[s][\scrr-1][k]} 
&\lesssim \Big\|  \Ad \big( \LON[K][-] \big) \Big\|_{\Cprod{s}{s-1}} 
\Big\| \rho_{\leq N}^2(k) \pYN[K][+][k] \Big\|_{\Wuv[s][\scrr][k]}  \notag \\
&\lesssim \big\| \LON[K][-] \big\|_{\Cprod{s}{s-1}} 
\big\| \rho_{\leq N}^2(k) \big\|_{\Wc_k} \big\| \pYN[K][+][k] \big\|_{\Wuv[s][\scrr][k]}. 
\label{modulation:eq-Y-remaining-p1}
\end{align}
Using Lemma \ref{chaos:lem-wc} and Corollary \ref{modulation:cor-LON}, we have that 
\begin{equation}\label{modulation:eq-Y-remaining-p2}
     \big\| \rho_{\leq N}^2(k) \big\|_{\Wc_k}  \lesssim 1 \qquad \text{and} \qquad \big\| \LON[K][-] \big\|_{\Cprod{s}{s-1}}  \lesssim \Dc. 
\end{equation}
Thus, the contribution of \eqref{modulation:eq-Y-remaining-p1} is acceptable. \\

\emph{Proof of \eqref{modulation:eq-Y-remaining-e2}:} By using our para-product estimate (Lemma \ref{prelim:lem-paraproduct}), which requires that $\scrr<2s$, and using Lemma \ref{chaos:lem-wc}, we obtain that
\begin{align}
\Big\| \rho_{\leq N}^2(k) \Ad \big( \LON[K][-] \big) \Para[v][ll] \pSN[K][+][k] \Big\|_{\Wuv[s][\scrr-1][k]} 
&\lesssim \Big\|  \Ad \big( \LON[K][-] \big) \Big\|_{\Cprod{s}{s-1}} 
\Big\| \rho_{\leq N}^2(k) \pSN[K][+][k] \Big\|_{\Wuv[s][s][k]}  \notag \\
&\lesssim \big\| \LON[K][-] \big\|_{\Cprod{s}{s-1}} 
\big\| \rho_{\leq N}^2(k) \big\|_{\Wc_k} \big\| \pSN[K][+][k] \big\|_{\Wuv[s][s][k]}. 
\label{modulation:eq-Y-remaining-p3}
\end{align}
Due to \eqref{modulation:eq-Y-remaining-p2}, this yields an acceptable contribution. \\

\emph{Proof of \eqref{modulation:eq-Y-remaining-e3}:} Using the product estimate (Lemma \ref{prelim:lem-paraproduct}), which requires that $\scrr-1+s>0$, and using Lemma \ref{chaos:lem-wc}, we obtain that
\begin{align}
\Big\| \rho_{\leq N}^4(k) \Ad \big( \SHHLN[K][v] \big) \, \pSN[K][+][k] \Big\|_{\Wuv[s][\scrr-1][k]} 
&\lesssim \Big\|  \Ad \big( \SHHLN[K][v] \big)  \Big\|_{\Cprod{s}{\scrr-1}} 
\Big\| \rho_{\leq N}^4(k) \pSN[K][+][k] \Big\|_{\Wuv[s][s][k]} \notag  \\
&\lesssim \big\|   \SHHLN[K][v] \big\|_{\Cprod{s}{\scrr-1}} 
\big\| \rho_{\leq N}^4(k) \big\|_{\Wc_k} \big\|  \pSN[K][+][k] \big\|_{\Wuv[s][s][k]}. \label{modulation:eq-Y-remaining-p4}
\end{align}
Due to Lemma \ref{chaos:lem-wc} and Lemma \ref{modulation:lem-shhl}, we have that 
\begin{equation*}
\big\| \rho_{\leq N}^4(k) \big\|_{\Wc_k} \lesssim 1 \qquad \text{and} \qquad \big\|   \SHHLN[K][v] \big\|_{\Cprod{s}{\scrr-1}}  \lesssim \Dc^2, 
\end{equation*}
and thus \eqref{modulation:eq-Y-remaining-p4} yields an acceptable contribution. 
\end{proof}

\subsection{Well-posedness of the para-controlled modulation equations}\label{section:modulation-wellposedness}
In this subsection, we use our previous estimates (from Subsection \ref{section:modulation-resonant} and Subsection \ref{section:modulation-nonlinear}) to prove the well-posedness of the para-controlled modulation equations. 

\begin{proposition}[Well-posedness of the para-controlled modulation equations]\label{modulation:prop-para-well-posedness}
Let the same assumptions as in Proposition \ref{modulation:prop-main} be satisfied.  
Then, there exist unique
\begin{equation}\label{modulation:eq-para-operators}
\begin{aligned}
\big( \pXN[K][+][k] \big)_{K\in \Dyadiclarge,k\in \Z_K},
\big( \pYN[K][+][k] \big)_{K\in \Dyadiclarge,k\in \Z_K} 
&\colon \R^{1+1} \rightarrow \End(\frkg), \\ 
\big( \pXN[M][-][m] \big)_{M\in \Dyadiclarge,m\in \Z_M},
\big( \pYN[M][-][m] \big)_{M\in \Dyadiclarge,m\in \Z_M}
&\colon \R^{1+1}\rightarrow \End(\frkg) 
\end{aligned}
\end{equation}
which satisfy the para-controlled modulation equations (Definition \ref{ansatz:def-modulation-equations}) and the estimates 
\begin{equation}\label{modulation:eq-para-bounds}
\begin{aligned}
\sup_{K\in \Dyadiclarge } \Big\| \pXN[K][+][k]\Big\|_{\Wuv[s][s][k]},
\sup_{M\in \Dyadiclarge } \Big\| \pXN[M][-][m] \Big\|_{\Wuv[s][s][m]}
&\lesssim \Dc \Bcin, \\ 
\sup_{K\in \Dyadiclarge } \Big\| \pYN[K][+][k] - \SNin[K][+] \Big\|_{\Wuv[s][\scrr][k]},
\sup_{M\in \Dyadiclarge } \Big\| \pYN[M][-][m]- \SNin[M][-] \Big\|_{\Wuv[\scrr][s][m]} 
&\lesssim \Dc \Bcin.
\end{aligned}
\end{equation}
\end{proposition}

In Subsection \ref{section:modulation-resonant} and Subsection \ref{section:modulation-nonlinear}, we already obtained all estimates which are needed in the proof of Proposition \ref{modulation:prop-para-well-posedness}. Thus, it only remains to implement a contraction-mapping argument. 

\begin{proof}[Proof of Proposition \ref{modulation:prop-para-well-posedness}:]
We first define the set 
\begin{equation}\label{modulation:eq-para-p1}
\begin{aligned}
\BallN := \Big\{ &\big( \pXN[][+], \pXN[][-], \pYN[][+], \pYN[][-] \big) \colon  
\, \max \Big( \sup_{\substack{K\in \Dyadiclarge}} \big\| \pXN[K][+][k] \big\|_{\Wuv[s][s][k]}, 
\sup_{\substack{M\in \Dyadiclarge}}\big\| \pXN[M][-][m] \big\|_{\Wuv[s][s][m]}, \\
&\sup_{\substack{K\in \Dyadiclarge}}\big\| \pYN[K][+][k] \big\|_{\Wuv[s][\scrr][k]}, 
\sup_{\substack{M\in \Dyadiclarge}} \big\| \pYN[M][-][m] \big\|_{\Wuv[\scrr][s][m]}
\Big) 
\leq 2 \Bcin
\Big\}. 
\end{aligned}
\end{equation}
Throughout this proof, we view $\pSN[][+]$ and $\pSN[][-]$ as functions of $\pXN[][+]$, $\pXN[][-]$, $\pYN[][+]$, and $\pYN[][-]$, which are given by 
\begin{equation*}
\pSN[K][+][k] = \pXN[K][+][k] + \pYN[K][+][k] \qquad 
\text{and} \qquad 
\pSN[M][-][m] = \pXN[M][-][m] + \pYN[M][-][m].
\end{equation*}
As a result, $(\pXN[][+],\pXN[][-],\pYN[][+],\pYN[][-])\in \BallN$ directly implies the estimates
\begin{equation*}
\sup_{\substack{K\in \Dyadiclarge}} \big\| \pSN[K][+][k] \big\|_{\Wuv[s][s][k]},
\sup_{\substack{M\in \Dyadiclarge}} \big\| \pSN[M][-][m] \big\|_{\Wuv[s][s][m]} \leq 4 \Bcin \leq \Bc. 
\end{equation*}
In particular, the bound in Hypothesis \ref{hypothesis:pre}.\ref{ansatz:item-pre-modulation} is satisfied. 
We now define $\XNcal[][\pm]$ and $\YNcal[][\pm]$ as the maps which encode the para-controlled modulation equations 
\eqref{ansatz:eq-X-p-1}-\eqref{ansatz:eq-Y-m-5}. More precisely, we define
\begin{align*}
\XNcal[K][+][k](\pXN[][+],\pXN[][-],\pYN[][+],\pYN[][-]) &= \eqref{ansatz:eq-X-p-1}, \\
\XNcal[M][-][m](\pXN[][+],\pXN[][-],\pYN[][+],\pYN[][-]) &= \eqref{ansatz:eq-X-m-1}, \\
\YNcal[K][+][k](\pXN[][+],\pXN[][-],\pYN[][+],\pYN[][-]) &= 
\SNin[K][+] + \eqref{ansatz:eq-Y-p-1} +\eqref{ansatz:eq-Y-p-2}+\eqref{ansatz:eq-Y-p-3}+\eqref{ansatz:eq-Y-p-4}+\eqref{ansatz:eq-Y-p-5}, \\
\YNcal[M][-][m](\pXN[][+],\pXN[][-],\pYN[][+],\pYN[][-]) &= 
\SNin[M][-] + \eqref{ansatz:eq-Y-m-1} +\eqref{ansatz:eq-Y-m-2}+\eqref{ansatz:eq-Y-m-3}+\eqref{ansatz:eq-Y-m-4}+\eqref{ansatz:eq-Y-m-5}.
\end{align*}
To simplify the notation, we also write\footnote{The operator $\ZNcal$ should not be confused with the partition function $\Zc_\beta$ in the definitions of our Gibbs measures.} 
\begin{align*}
\ZN &= \big( \pXN[][+], \pXN[][-], \pYN[][+], \pYN[][-] \big), \\ 
\big\| \ZN \big\| &=  \max \Big( \sup_{\substack{K\in \Dyadiclarge}} \big\| \pXN[K][+][k] \big\|_{\Wuv[s][s][k]}, 
\sup_{\substack{M\in \Dyadiclarge}} \big\| \pXN[M][-][m] \big\|_{\Wuv[s][s][m]}, 
\sup_{\substack{K\in \Dyadiclarge}} \big\| \pYN[K][+][k] \big\|_{\Wuv[s][\scrr][k]}, 
\sup_{\substack{M\in \Dyadiclarge}} \big\| \pYN[M][-][m] \big\|_{\Wuv[\scrr][s][m]}
\Big), \\ 
\ZNcal &= \big( \XNcal[][+], \XNcal[][-], \YNcal[][+], \YNcal[][-] \big). 
\end{align*}
To solve the para-controlled modulation equations, we first have to show that $\ZNcal$ maps $\BallN$ into $\BallN$ and then have to show that $\ZNcal$ is a contraction, which implies that $\ZNcal$ has a unique fixed point in $\BallN$. For expository purposes, we split the rest of this proof into two steps. \\ 

\emph{Step 1: $\ZNcal$ maps $\BallN$ into $\BallN$.} Due to the symmetry of our estimates in the $u$ and $v$-variables, it suffices to control $\XNcal[][+][]$ and $\YNcal[][+][]$. Using Lemma \ref{modulation:lem-X-estimate}, we obtain that for all $K\in \Dyadiclarge$ that 
\begin{equation}\label{modulation:eq-para-X}
\big\| \XNcal[K][+][k](\ZN) \big\|_{\Wuv[s][s][k]} \leq \big\| \eqref{ansatz:eq-X-p-1} \big\|_{\Wuv[s][s][k]} 
\lesssim \Dc \Bcin. 
\end{equation}
Since $\Dc$ is small, this easily yields the desired upper bound by $2\Bcin$. Using the definition of $\YNcal[K][+][k]$ and the triangle inequality, we also obtain for all $K\in \Dyadiclarge$ that
\begin{equation*}
\big\| \YNcal[K][+][k](\ZN) - \SNin[K][+] \big\|_{\Wuv[s][\scrr][k]} 
\leq 
\big\| \eqref{ansatz:eq-Y-p-1} \big\|_{\Wuv[s][\scrr][k]} 
+ \big\| \eqref{ansatz:eq-Y-p-2} \big\|_{\Wuv[s][\scrr][k]} 
+ \big\| \eqref{ansatz:eq-Y-p-3} \big\|_{\Wuv[s][\scrr][k]} 
+ \big\| \eqref{ansatz:eq-Y-p-4} \big\|_{\Wuv[s][\scrr][k]} 
+ \big\| \eqref{ansatz:eq-Y-p-5} \big\|_{\Wuv[s][\scrr][k]}.
\end{equation*}
We now estimate \eqref{ansatz:eq-Y-p-1} using Lemma \ref{modulation:lem-Y-integral-commutator}, \eqref{ansatz:eq-Y-p-2} using Lemma \ref{modulation:lem-Y-double-duhamel}, and \eqref{ansatz:eq-Y-p-3}, \eqref{ansatz:eq-Y-p-4}, and \eqref{ansatz:eq-Y-p-5} using Lemma \ref{modulation:lem-Y-remaining}. In total, it then follows that 
\begin{equation}\label{modulation:eq-para-Y}
\big\| \YNcal[K][+][k](\ZN) - \SNin[K][+] \big\|_{\Wuv[s][\scrr][k]}   \lesssim \Dc \Bcin.
\end{equation}
Since $\SNin[K][+]=\SNin[K][+](u)$ neither depends on $v\in \R$ nor $k\in \Z_K$, our assumptions imply that
\begin{equation*}
\big\| \SNin[K][+] \big\|_{\Wuv[s][\scrr][k]} = \big\| \SNin[K][+] \big\|_{\C_u^{s}} \leq \Bcin
\end{equation*}
for all $K\in \Dyadiclarge$. 
Since $\Dc$ is small, this  yields the desired upper bound of $\| \YNcal[K][+][k](\ZN) \big\|_{\Wuv[s][\scrr][k]}$  by $2\Bcin$. Thus, the maps $\XNcal[][+]$ and $\YNcal[][+]$ satisfy the desired upper bounds, and hence $\ZNcal$ maps $\BallN$ into $\BallN$. From
\eqref{modulation:eq-para-X} and \eqref{modulation:eq-para-Y} it further follows that, once the para-controlled modulation equations have been solved, the solution satisfies \eqref{modulation:eq-para-bounds}.\\ 

\emph{Step 2: $\ZNcal$ is a contraction.} In order to prove that $\ZNcal$ is a contraction on $\BallN$, it suffices to prove the estimate
\begin{equation}\label{modulation:eq-para-p2}  
\big\| \ZNcal \big( \ZN \big) - \ZNcal \big( \ZNtil \big) \big\| \lesssim \Dc\,  \big\| \ZN - \ZNtil \big\|
\end{equation}
for all $\ZN,\ZNtil\in \BallN$. This follows from a standard (but tedious) modification of the estimates from the first step.
We omit the details here, but refer to Section \ref{section:Lipschitz} for related modifications. 
\end{proof}

We now record a consequence of the estimates of Subsection \ref{section:modulation-nonlinear} and the proof of Proposition \ref{modulation:prop-para-well-posedness}.

\begin{corollary}[Moral frequency support]\label{modulation:cor-moral-frequency}
Let the assumptions of Proposition \ref{modulation:prop-para-well-posedness} be satisfied, let $K\in \Dyadiclarge$, and let $L_u,L_v\in \dyadic$.
Then, we have the refined estimates 
\begin{equation}\label{modulation:eq-refined-bounds}
\Big( \tfrac{\max( L_u,L_v)}{K^{1-\delta}} \Big)^{100\vartheta^{-1}}  
\Big( \Big\| P^u_{L_u} P^v_{L_v} \pXN[K][+][k] \Big\|_{\Wuv[s][s][k]} + 
\Big\| P^u_{L_u} P^v_{L_v} \big( \pYN[K][+][k] - \SNin[K][+] \big) \Big\|_{\Wuv[s][\scrr][k]} \Big)
\lesssim     \Dc \Bcin.
\end{equation}
Similar estimates also hold for $\XN[M][-][m]$ and $\YN[M][-][m]$, but with $K$ replaced by $M$. 
\end{corollary}

Loosely speaking, Corollary \ref{modulation:cor-moral-frequency} states that the operators $( \pXN[K][+][k] \big)_{k\in \Z_K}$
and $(\pYN[K][+][k] \big)_{k\in \Z_K}$ are morally supported on frequencies $\lesssim K^{1-\delta}$. Due to Definition \ref{ansatz:def-lo}, Definition \ref{ansatz:def-shhl}, and Definition \ref{ansatz:def-modulation-para}, the para-controlled modulation equations for $( \pXN[K][+][k] \big)_{k\in \Z_K}$
and $(\pYN[K][+][k] \big)_{k\in \Z_K}$ are ordinary differential equations whose coefficients are supported on frequencies $\lesssim K^{1-\delta}$, and this is therefore to be expected. 

\begin{proof}[Proof of Corollary \ref{modulation:cor-moral-frequency}:]
The refined estimates \eqref{modulation:eq-refined-bounds} can be proven using a modification of the functional framework used in the proof of Proposition \ref{modulation:prop-para-well-posedness}. Since the modifications are purely technical and do not require any additional ideas, we only sketch them. For a more detailed argument involving a similar functional framework, we also refer the reader to \cite{B21}. 

We first introduce a modification of our $\Cprod{\alpha}{\beta}$-norm from Definition \ref{prelim:def-hoelder}. For any $\alpha,\beta\in (-1,1)\backslash \{0\}$, $L\in \dyadic$, and $f\colon \R^{1+1}\rightarrow \bC$, we define
\begin{equation}\label{modulation:eq-moral-frequency-e1}
\big\| f \big\|_{(\Cprod{\alpha}{\beta})_{\lesssim L}} := \sup_{L_u,L_v\in \dyadic} L_u^\alpha L_v^\beta \max \bigg( \frac{L_u}{L}, \frac{L_v}{L}, 1 \bigg)^{100\vartheta^{-1}} \big\| P^u_{L_u} P^v_{L_v} f \big\|_{L^\infty}.
\end{equation}
Under the same conditions as in Lemma \ref{prelim:lem-paraproduct}, it holds for all $L\in \dyadic$ and $f,g\colon \R^{1+1}\rightarrow \bC$ that 
\begin{equation}\label{modulation:eq-moral-frequency-p1}
\big\| f \parauvrcal g \big\|_{(\Cprod{\gamma_1}{\gamma_2})_{\lesssim L}}
\lesssim \| f \|_{(\Cprod{\alpha_1}{\alpha_2})_{\lesssim L}} \| g \|_{(\Cprod{\beta_1}{\beta_2})_{\lesssim L}}. 
\end{equation}
The reason is that 
$P^u_{K_u} P^v_{K_v} ( P^u_{L_u} P^v_{L_v} f \cdot P^u_{M_u} P^v_{M_v} g)$
is only non-zero if $K_u \lesssim \max(L_u,M_u)$ and $K_v\lesssim \max(L_v,M_v)$. Thus, the maxima 
from \eqref{modulation:eq-moral-frequency-e1} which occur on the right-hand side of \eqref{modulation:eq-moral-frequency-p1} can always be used to control the corresponding factor on the left-hand side of \eqref{modulation:eq-moral-frequency-p1}. We also note that, for all $L\in \dyadic$ and $f\colon \R^{1+1}\rightarrow \bC$, it holds that
\begin{equation}\label{modulation:eq-moral-frequency-p2}
\| P^{u,v}_{\leq L} f \|_{(\Cprod{\alpha_1}{\alpha_2})_{\lesssim L}} \lesssim \| f \|_{\Cprod{\alpha_1}{\alpha_2}}. 
\end{equation}

We also introduce a variant of the $\Wuv[\alpha][\beta](\Z)$-norm from Definition \ref{chaos:def-wc}. For any $L\in \dyadic$ and $f\colon \Z \times \R^{1+1}\rightarrow \bC$, it is defined as 
\begin{equation}\label{modulation:eq-moral-frequency-e2} 
\big\| f_\ell(u,v) \big\|_{\Wuv[\alpha][\beta][\lesssim L](\Z)} 
:= \sup_{\ell \in \Z} \big\| f_\ell(u,v) \big\|_{(\Cprod{\alpha}{\beta})_{\lesssim L}} 
+ \sum_{\ell \in \Z} \big\| \big( \Der f\big)_\ell (u,v) \big\|_{(\Cprod{\alpha}{\beta})_{\lesssim L}},
\end{equation}
where $\Der$ is the discrete derivative from Definition \ref{chaos:def-derivative}.

In order to prove the refined estimates \eqref{modulation:eq-refined-bounds}, we need to estimate $\pXN[K][+][k]$ and $\pYN[K][+][k]$ in $\Wuv[s][s][\lesssim K^{1-\delta}]$ and $\Wuv[s][\scrr][\lesssim K^{1-\delta}]$ instead of $\Wuv[s][s][]$ and $\Wuv[s][\scrr][]$, respectively. Since the only ingredients in the proof of Proposition \ref{modulation:prop-para-well-posedness} were the estimates from  Lemma \ref{modulation:lem-X-estimate}, Lemma \ref{modulation:lem-Y-integral-commutator}, Lemma \ref{modulation:lem-Y-double-duhamel}, and Lemma \ref{modulation:lem-Y-remaining}, we then only need to extend the latter estimates from $\Wuv[\alpha][\beta][]$ to $\Wuv[\alpha][\beta][\lesssim K^{1-\delta}]$-norms. The reason this is possible is that, by definition, $\LON[K][-]$ and $\SHHLN[K][v]$ are supported on frequencies $\lesssim K^{1-\delta}$ in both variables. Thus, this extension only requires us to replace each application of Lemma \ref{prelim:lem-paraproduct} with applications of \eqref{modulation:eq-moral-frequency-p1} and \eqref{modulation:eq-moral-frequency-p2}.
\end{proof}

\subsection{Geometric properties and proof of Proposition \ref{modulation:prop-main}}
\label{section:modulation-properties}

In the previous subsections, we dealt with analytical aspects of the modulation equations. In this subsection, we now address geometric aspects of the modulation equations. 

\begin{proposition}[Properties of modulation operators]\label{modulation:prop-properties} 
Let $K,M\in \Dyadiclarge$, let $k\in \Z_K$, let $m\in \Z_M$, and let $u,v\in \R$. Then, the solutions of the modulation equations \eqref{modulation:eq-motivation} satisfy the following properties: 
\begin{enumerate}[label=(\roman*),leftmargin=6ex]
\item\label{modulation:item-properties-1} 
(Orthogonality) The linear transformations
\begin{equation}\label{modulation:eq-multiplied}
 \pSN[K][+][k](u,v) \big( \SNin[K][+](u) \big)^{-1} \colon \frkg \rightarrow \frkg 
 \qquad \text{and} \qquad
 \pSN[M][-][m](u,v) \big( \SNin[M][-](v) \big)^{-1} \colon \frkg \rightarrow \frkg
\end{equation}
are orthogonal. 
\item\label{modulation:item-properties-2}
(Automorphism)  The linear transformations in \eqref{modulation:eq-multiplied} are Lie algebra automorphisms. 
\item\label{modulation:item-properties-3} ($k$ and $m$-dependence) 
If $K\leq \tfrac{1}{2} N$, then $\pSN[K][+][k]$ is constant in $k\in \Z_{K}$. Similarly, if $M\leq \tfrac{1}{2} N$, then $\pSN[M][-][m]$ is constant in $m\in \Z_M$.
\end{enumerate}
\end{proposition}

\begin{remark} 
The automorphism property is not used in any part of this article, and it is only recorded here to underline the geometric nature of the modulation operators.  
\end{remark}
\begin{proof}
It suffices to prove the properties for $\pSN[K][+][k]$, since the arguments for $\pSN[M][-][m]$ are similar. To simplify the notation, we set 
\begin{align*}
\Qk(u,v) &:= - \chinull[K](\tfrac{v-u}{2}) \rho_{\leq N}^2(k) \LON[K][-](u,v) - \chinull[K](\tfrac{v-u}{2}) \rho_{\leq N}^4(k) \SHHLN[K][v](u,v), \\ 
\pSNhat[K][+][k](u,v) &:= \pSN[K][+][k](u,v) \big( \SNin[K][+](u) \big)^{-1}. 
\end{align*}
Then, the initial value problems for $\pSN[K][+][k]$ and $\pSNhat[K][+][k]$ can be written as 
\begin{alignat}{3}
\partial_v \pSN[K][+][k](u,v) 
&= \Ad \big( \Qk(u,v) \big)  \pSN[K][+][k](u,v),  &\qquad \pSN[K][+][k](u,v)\big|_{v=u}&= \SNin[K][+](u), 
\label{modulation:eq-Sp-rewritten}\\
\partial_v \pSNhat[K][+][k](u,v) 
&= \Ad \big( \Qk(u,v) \big)  \pSNhat[K][+][k](u,v),  &\qquad \pSNhat[K][+][k](u,v)\big|_{v=u}&= \Id_\frkg. 
\label{modulation:eq-Sp-hat-rewritten}
\end{alignat}
We now treat \ref{modulation:item-properties-1}, \ref{modulation:item-properties-2}, and \ref{modulation:item-properties-3} separately.\\

\emph{Proof of \ref{modulation:item-properties-1}:} Due to \eqref{prelim:eq-Hermitian-adjoint}, 
\begin{equation*}
\Ad \big( \Qk(u,v) \big) \colon \frkg \rightarrow \frkg
\end{equation*}
is skew-symmetric for all $u,v\in \R$. Using the evolution equation in \eqref{modulation:eq-Sp-hat-rewritten}, it follows that 
\begin{align*}
\partial_v \Big( \big( \pSNhat[K][+][k] \big)^\ast \pSNhat[K][+][k] \Big) 
&= \big(\pSNhat[K][+][k]\big)^\ast \partial_v \pSNhat[K][+][k] + \partial_v \big( \pSNhat[K][+][k] \big)^\ast \pSNhat[K][+][k] \\
&= \big(\pSNhat[K][+][k]\big)^\ast \Big( \Ad \big( \Qk \big) + \Ad \big( \Qk \big)^\ast \Big)  \pSNhat[K][+][k] \\
&=0. 
\end{align*}
Together with $\pSNhat[K][+][k](u,u)=\Id_\frkg$, this implies the desired orthogonality. \\ 

\emph{Proof of \ref{modulation:item-properties-2}:} 
Due to orthogonality, $\pSNhat[K][+][k](u,v)$ is clearly invertible.  Thus, it only remains to prove that $\pSNhat[K][+][k](u,v)$ preserves the Lie bracket. That is, it remains to prove for all $A,B\in \frkg$ that 
\begin{equation}\label{modulation:eq-properties-p1}
\Big[ \pSNhat[K][+][k](u,v) A, \pSNhat[K][+][k](u,v) B \Big] 
= \pSNhat[K][+][k](u,v) \Big[ A, B \Big]. 
\end{equation}
In order to prove \eqref{modulation:eq-properties-p1}, we define
\begin{equation*}
\pSNhat[K][+][k,A,B](u,v) := \Big[ \pSNhat[K][+][k](u,v) A, \pSNhat[K][+][k](u,v) B \Big] 
- \pSNhat[K][+][k](u,v) \Big[ A, B \Big]. 
\end{equation*}
Using the evolution equation in \eqref{modulation:eq-Sp-hat-rewritten} and the Jacobi identity on the Lie algebra $\frkg$, one easily obtains that
\begin{equation*}
\partial_v \pSNhat[K][+][k,A,B](u,v) = \Ad \Big( \Qk(u,v) \Big) \pSNhat[K][+][k,A,B](u,v). 
\end{equation*}
Since the initial condition in \eqref{modulation:eq-Sp-hat-rewritten} guarantees that $\pSNhat[K][+][k,A,B](u,u)=0$, this yields $\pSNhat[K][+][k,A,B](u,v)=0$ and therefore \eqref{modulation:eq-properties-p1}.  \\

\emph{Proof of \ref{modulation:item-properties-3}:} Due to \eqref{prelim:eq-rho} and \eqref{prelim:eq-rho-leqN}, 
the condition $K\leq \frac{1}{2} N$ implies that $\rho_{\leq N}(k)$ is constant on $\Z_K$. Since the solution of the modulation equations is unique, it then follows that $\pSN[K][+][k]$ is also constant on $\Z_K$.
\end{proof}

Equipped with both Proposition \ref{modulation:prop-para-well-posedness} (from Subsection \ref{section:modulation-wellposedness}) and Proposition \ref{modulation:prop-properties}, we can now prove the main result of this section. 

\begin{proof}[Proof of Proposition \ref{modulation:prop-main}:] 
Using Proposition \ref{modulation:prop-para-well-posedness}, we obtain unique $\pXN[][+]$, $\pXN[][-]$, $\pYN[][+]$, and $\pYN[][-]$ which satisfy the para-controlled modulation equations (Definition \ref{ansatz:def-modulation-para}) and satisfy the bounds from \eqref{modulation:eq-para-bounds}. We now set 
\begin{equation}\label{modulation:eq-main-p1}
\pSN[K][+][k] := \pXN[K][+][k] + \pYN[K][+][k] \qquad \text{and} \qquad 
\pSN[M][-][m] := \pXN[M][-][m] + \pYN[M][-][m].
\end{equation}
Due to Proposition \ref{modulation:lem-para-controlled-modulation-equations}, the pure modulation operators from \eqref{modulation:eq-main-p1} solve the modulation equations from Definition \ref{ansatz:def-modulation-equations}. Furthermore, due to \eqref{modulation:eq-para-bounds}, the bounds in \eqref{modulation:eq-main-bounds} are satisfied. 
Thus, this already addresses the claims regarding the existence of $\SN[][+]$ and $\SN[][-]$ and the uniform bounds in Proposition \ref{modulation:prop-main}. The uniqueness statement follows from soft arguments since, for any fixed $\Nd,N\in \Dyadiclarge$, the modulation equations \eqref{modulation:eq-motivation} form a system of ordinary differential equations with smooth coefficients. It therefore only remains to verify the properties in \ref{modulation:item-difference}, \ref{modulation:item-distance-initial}, \ref{modulation:item-orthogonality}, and  \ref{modulation:item-structure}. Since the arguments for $\SN[M][-][m]$ are similar, we focus our attention on $\SN[K][+][k]$.  

To obtain \ref{modulation:item-difference}, we use Definition \ref{ansatz:def-initial-data} and Definition \ref{ansatz:def-pure}, which implies that
\begin{equation*}
\SN[K][+][k] - \pSN[K][+][k] =  - P^{u,v}_{>K^{1-\delta+\vartheta}} \big( \pSN[K][+][k] \big) 
= -  P^{u,v}_{>K^{1-\delta+\vartheta}} \big( \pSN[K][+][k] -\SNin[K][+] \big).
\end{equation*}
Thus, the estimate in \ref{modulation:item-difference} follows directly from Corollary \ref{modulation:cor-moral-frequency}. 
To obtain \ref{modulation:item-distance-initial}, we split
\begin{equation*}
\SN[K][+][k] - \SNin[K][+] = \big( \SN[K][+][k] - \pSN[K][+][k] \big) + \big( \pSN[K][+][k] -  \SNin[K][+]\big).
\end{equation*}
The first and second summand can be estimated using \ref{modulation:item-difference} and  using Proposition \ref{modulation:prop-para-well-posedness}, respectively.
To \mbox{obtain \ref{modulation:item-orthogonality}}, we use \mbox{Definition \ref{ansatz:def-initial-data}}, which implies that $\SNin[K][+]$ is orthogonal for all $K>N^{1-\delta}$. Due to Proposition \ref{modulation:prop-properties}, it then follows that $\pSN[K][+][k]$ is orthogonal for all $K>N^{1-\delta}$ and $k\in \Z_K$. As a result, we can write 
\begin{equation*}
\big( \SN[K][+][k] \big)^\ast \SN[K][+][k] - \Id_\frkg 
= \big( \SN[K][+][k] \big)^\ast \SN[K][+][k] - \big( \pSN[K][+][k] \big)^\ast \pSN[K][+][k].
\end{equation*}
The desired estimate in \ref{modulation:item-orthogonality} can now be obtained from the estimate in \ref{modulation:item-difference} and Lemma \ref{chaos:lem-wc}. To obtain \ref{modulation:item-structure}, we first use Proposition \ref{modulation:prop-para-well-posedness} and write
\begin{align}
\SN[K][+][k] 
&= \pSN[K][+][k] + \big( \SN[K][+][k] - \pSN[K][+][k] \big) \notag \\ 
&=  -  \chinull[K] \rho_{\leq N}^2(k) \, \Int^{v}_{u\rightarrow v} \Big( \ad \big( \LON[K][-] \big) \Big) \Para[v][gg] \pSN[K][+][k] 
+ \pYN[K][+][k] + \big( \SN[K][+][k] - \pSN[K][+][k] \big) \notag \\ 
&= -  \chinull[K] \rho_{\leq N}^2(k) \, \Int^{v}_{u\rightarrow v} \Big( \ad \big( \LON[K][-] \big) \Big) \Para[v][gg] \SN[K][+][k] \notag \\ 
&+ \pYN[K][+][k] + \big( \SN[K][+][k] - \pSN[K][+][k] \big)
- \chinull[K] \rho_{\leq N}^2(k) \, \Int^{v}_{u\rightarrow v} \Big( \ad \big( \LON[K][-] \big) \Big) \Para[v][gg] \big( \pSN[K][+][k] - \SN[K][+][k] \big). \label{modulation:eq-main-p2}
\end{align}
The three terms in \eqref{modulation:eq-main-p2} can easily be controlled in $\Wuv[s][\scrr]$ using Proposition \ref{modulation:prop-para-well-posedness} and Corollary \ref{modulation:cor-moral-frequency}, and we therefore obtain \ref{modulation:item-structure}. 
\end{proof}

\subsection{Structure of modulated linear waves}\label{section:modulation-structure} 

In Subsection \ref{section:modulation-wellposedness}, we proved the well-posedness of the para-controlled modulation equations. In particular, we established that the modulation operators $\SN[][\pm]$ exhibit a para-controlled structure. In this subsection, we examine the consequences of this para-controlled structure of $\SN[][\pm]$ for the modulated linear waves $\UN[][+]$ and $\VN[][-]$.

\begin{lemma}\label{modulation:lem-structure}
Let the post-modulation hypothesis (Hypothesis \ref{hypothesis:post}) be satisfied and let $\gamma \in (-1,1)\backslash\{0\}$. For all $K\in \Dyadiclarge$, it then holds that 
\begin{align}
\Big\| 
\UN[K][+] - \chinull[K] \Big[ \UN[K][+] \Para[v][ll] \IVN[<K^{1-\delta}][-] \Big]_{\leq N} \Big\|_{\Cprod{\gamma}{r}} 
&\lesssim K^{\gamma+1/2+\eta} \Dc, 
\label{modulation:eq-structure-e1} \\
\Big\| 
\UN[K][+] 
- \chinull[K]  \Big[ \UN[K][+] \Para[v][ll] \Big(  \IVN[<K^{1-\delta}][-]
+ \Int^v_{u\rightarrow v} P^v_{<K^{1-\deltap}} \VN[<K^{1-\delta}][\fcs]  \Big) \Big]_{\leq N}
\Big\|_{\Cprod{\gamma}{1/2+\delta}} 
&\lesssim K^{\gamma+1/2+\eta} \Dc. 
\label{modulation:eq-structure-e2}
\end{align}
Similarly, for all $M\in \Dyadiclarge$, it holds that 
\begin{align}
\Big\| \VN[M][-] 
- \chinull[M] \Big[ \IUN[<M^{1-\delta}][+] \Para[u][gg] \VN[M][-] \Big]_{\leq N}
\Big\|_{\Cprod{r}{\gamma}} &\lesssim M^{\gamma+1/2+\eta} \Dc, \\ 
\Big\| \VN[M][-] 
- \chinull[M] \Big[ \Big( \IUN[<M^{1-\delta}][+]
+ \Int^u_{v\rightarrow u} P^u_{<M^{1-\delta}} \UN[<M^{1-\delta}][\fsc] \Big) \Para[u][gg] \VN[M][-] \Big]_{\leq N}
\Big\|_{\Cprod{1/2+\delta}{\gamma}} &\lesssim M^{\gamma+1/2+\eta} \Dc. 
\end{align}
\end{lemma}

\begin{remark}
Since the $v$-regularity in \eqref{modulation:eq-structure-e1} is lower than the $v$-regularity in \eqref{modulation:eq-structure-e2}, \eqref{modulation:eq-structure-e2} is a more accurate description of $\UN[K][+]$. In some situations, however, the additional $v$-regularity is not necessary, and then the simpler expansion from \eqref{modulation:eq-structure-e1} will be more convenient. 
\end{remark}

The following argument relies heavily on the para-controlled structure of the modulation operators $(\SN[K][+][k])_{k\in \Z}$. 

\begin{proof} 
Due to the symmetry of the estimates in the $u$ and $v$-variables, it suffices to prove \eqref{modulation:eq-structure-e1} and \eqref{modulation:eq-structure-e2}. 
For expository purposes, we separate the proof into four different steps. \\

\emph{Step 1:} In the first step, we prove that 
\begin{align}\label{modulation:eq-structure-p1} 
\Big\| \UN[K][+] 
-  \chinull[K] \Big[  \UN[K][+]  \Para[v][ll] 
\Int^v_{u\rightarrow v} \LON[K][-] \Big]_{\leq N}
\Big\|_{\Cprod{\gamma}{1/2+\delta}}
\lesssim K^{\gamma+1/2+\eta} \Dc. 
\end{align}
Using Proposition \ref{modulation:prop-main}, we obtain that 
\begin{equation}\label{modulation:eq-structure-p2}
\SN[K][+][k] = - \chinull[K] \rho_{\leq N}^2(k) \Big( \Int^v_{u\rightarrow v} \big( \Ad \big( \LON[K][-] \big) \big) \Para[v][gg] \SN[K][+][k] \Big) + \YN[K][+][k],
\end{equation}
where the remainder satisfies 
\begin{equation}\label{modulation:eq-structure-p3}
\big\| \, \YN[K][+][] \big\|_{\Wuv[s][\scrr](\Z)} \lesssim \Bc. 
\end{equation}
Using Definition \ref{ansatz:def-modulated-linear}, it then follows that 
\begin{align}
\UN[K][+] 
&= \hcoup \sum_{u_0\in \LambdaRR} \sum_{k\in \Z_K} \psiRuK  \rhoND(k) \SN[K][+][k] G_{u_0,k}^+ \, e^{iku}  \notag \\ 
&= \chinull[K] \hcoup \Big[ \sum_{u_0\in \LambdaRR} \sum_{k\in \Z_K} \psiRuK  \rho_{\leq N}^2(k)  \rhoND(k) \SN[K][+][k] G_{u_0,k}^+ e^{iku} \Para[v][ll] \Int^v_{u\rightarrow v}\LON[K][-]  \Big] \label{modulation:eq-structure-p4} \\ 
&+ \hcoup \sum_{u_0\in \LambdaRR} \sum_{k\in \Z_K} \psiRuK  \rhoND(k) \YN[K][+][k] G_{u_0,k}^+ \, e^{iku}. \label{modulation:eq-structure-p4prime} 
\end{align}
The term in \eqref{modulation:eq-structure-p4} can then be further decomposed as 
\begin{align}
&\chinull[K]  \hcoup \Big[ \sum_{u_0\in \LambdaRR} \sum_{k\in \Z_K} \psiRuK  \rho_{\leq N}^2(k)  \rhoND(k) \SN[K][+][k] G_{u_0,k}^+ e^{iku} \Para[v][ll] \Int^v_{u\rightarrow v}\LON[K][-]  \Big] \notag\\ 
=\,& \chinull[K]  \Big[ \Big( \hcoup \hspace{-1.25ex}\sum_{u_0\in \LambdaRR} \sum_{k\in \Z_K} \psiRuK  \rho_{\leq N}^2(k)  \rhoND(k) \SN[K][+][k] G_{u_0,k}^+ e^{iku} - (P_{\leq N}^x)^2 \UN[K][+] \Big) \Para[v][ll] P^v_{>1}\Int^v_{u\rightarrow v}\LON[K][-]  \Big] 
\label{modulation:eq-structure-p5} \\
+\,& \chinull[K]  \Big( \Big[ (P_{\leq N}^x)^2 \UN[K][+]  \Para[v][ll]P^v_{>1} \Int^v_{u\rightarrow v}\LON[K][-]  \Big] 
- \Big[  \UN[K][+]  \Para[v][ll] P^v_{>1} \Int^v_{u\rightarrow v}\LON[K][-]  \Big]_{\leq N} \Big)
\label{modulation:eq-structure-p6} \\
+\,& \chinull[K]   \Big[  \UN[K][+]  \Para[v][ll] \Int^v_{u\rightarrow v}\LON[K][-]  \Big]_{\leq N}. 
\end{align}
In \eqref{modulation:eq-structure-p5} and \eqref{modulation:eq-structure-p6}, we also used the definition of $\Para[v][ll]$ to insert additional $P^v_{>1}$-operators. 
In order to prove our desired estimate \eqref{modulation:eq-structure-p1}, it thus suffices to prove the three estimates 
\begin{align}
\hcoup \Big\| \sum_{u_0\in \LambdaRR} \sum_{k\in \Z_K} \psiRuK  \rhoND(k) \YN[K][+][k] G_{u_0,k}^+ \, e^{iku} \Big\|_{\Cprod{\gamma}{\scrr}} &\lesssim K^{\gamma+1/2+\eta} \Dc, 
\label{modulation:eq-structure-p7} \\ 
\big\| \eqref{modulation:eq-structure-p5}  \big\|_{\Cprod{\gamma}{1/2+\delta}} &\lesssim K^{\gamma+1/2+\eta} \Dc^2,
\label{modulation:eq-structure-p8} \\ 
\big\| \eqref{modulation:eq-structure-p6}  \big\|_{\Cprod{\gamma}{1/2+\delta}} &\lesssim K^{\gamma+1/2+\eta} \Dc^2. 
\label{modulation:eq-structure-p9}
\end{align}
The first estimate \eqref{modulation:eq-structure-p7} follows directly from \eqref{modulation:eq-structure-p3} and  \ref{ansatz:item-hypothesis-linear} in Hypothesis \ref{hypothesis:probabilistic}. To prove the second estimate \eqref{modulation:eq-structure-p8}, we first note that \eqref{modulation:eq-structure-p5} only contains high$\times$low-interactions in the $u$-variable. Using our paraproduct estimate (Lemma \ref{prelim:lem-paraproduct}), Lemma \ref{prelim:lem-Duhamel-integral}, and Lemma \ref{modulation:lem-PNX-modulated}, it then follows that
\begin{align*}
\big\| \eqref{modulation:eq-structure-p5}  \big\|_{\Cprod{\gamma}{1/2+\delta}}
&\lesssim 
\Big\| \hcoup \sum_{u_0\in \LambdaRR} \sum_{k\in \Z_K} \psiRuK  \rho_{\leq N}^2(k)  \rhoND(k) \SN[K][+][k] G_{u_0,k}^+ e^{iku} - (P_{\leq N}^x)^2 \UN[K][+] \Big\|_{\Cprod{\gamma}{s}} \\ 
&\hspace{3ex} \times 
\Big\| P^v_{>1} \Int^v_{u\rightarrow v} \LON[K][-] \Big\|_{\Cprod{s}{1/2+\delta}} \\ 
&\lesssim K^{\gamma\revision{+1/2}+\eta} K^{-\delta+\vartheta} \Dc
\Big\|  \LON[K][-] \Big\|_{\Cprod{s}{-1/2+\delta}} . 
\end{align*}
Due to Lemma \ref{ansatz:lem-frequency-support}, $\LON[K][-]$ is supported on $v$-frequencies $\lesssim K^{1-\delta}$. Together with Corollary \ref{modulation:cor-LON}, we obtain
\begin{align*}
K^{-\delta+\vartheta} 
\Big\|  \LON[K][-] \Big\|_{\Cprod{s}{-1/2+\delta}} 
\lesssim K^{-\delta+\vartheta} K^{(1-\delta)(\delta+\frac{1}{2}-s)} 
\Big\| \LON[K][-] \Big\|_{\Cprod{s}{s-1}} 
\lesssim K^{-\delta^2} K^{\frac{1}{2}-s+\vartheta} \Dc.  
\end{align*}
Since $\frac{1}{2}-s+\vartheta=\delta_2+\delta_4\ll \delta^2$, this completes the proof of \eqref{modulation:eq-structure-p8}. 
It remains to prove the third estimate \eqref{modulation:eq-structure-p9}. We first note that, as discussed in Remark \ref{ansatz:rem-PNX-LON}, $P^v_{>1}\Int^v_{u\rightarrow v} \LON[K][-] = P_{\leq N}^x P^v_{>1}\Int^v_{u\rightarrow v} \LON[K][-]$. Just as \eqref{modulation:eq-structure-p5}, \eqref{modulation:eq-structure-p6} only contains high$\times$low-interactions in the $u$-variable. Using Lemma \ref{prelim:commutator-PNx}, it then follows that
\begin{align*}
\Big\| \, \eqref{modulation:eq-structure-p6} \Big\|_{\Cprod{\gamma}{1/2+\delta}}
&\lesssim N^{-1} \Big\| \UN[K][+] \Big\|_{\Cprod{\gamma}{s}} 
\Big( \Big\| P^v_{>1} \Int^v_{u\rightarrow v} \LON[K][-] \Big\|_{\Cprod{s+1}{1/2+\delta}}
+ \Big\| P^v_{>1} \Int^v_{u\rightarrow v} \LON[K][-] \Big\|_{\Cprod{s+1}{3/2+\delta}} \Big) \\ 
&\lesssim K^{\gamma+1/2+\eta} N^{-1} \Dc \Big( \Big\| P^v_{>1} \Int^v_{u\rightarrow v} \LON[K][-] \Big\|_{\Cprod{s+1}{1/2+\delta}}
+ \Big\| P^v_{>1} \Int^v_{u\rightarrow v} \LON[K][-] \Big\|_{\Cprod{s+1}{3/2+\delta}} \Big). 
\end{align*}
Using that $\LON[K][-]$ is supported on frequencies $\lesssim K^{1-\delta}$, Lemma \ref{prelim:lem-Duhamel-integral}, and Corollary \ref{modulation:cor-LON}, it follows that 
\begin{align*}
    &N^{-1} \Big( \Big\| P^v_{>1} \Int^v_{u\rightarrow v} \LON[K][-] \Big\|_{\Cprod{s+1}{1/2+\delta}}
+ \Big\| P^v_{>1} \Int^v_{u\rightarrow v} \LON[K][-] \Big\|_{\Cprod{s+1}{3/2+\delta}} \Big) \\ 
\lesssim&\,    N^{-1} K^{1-\delta} \Big\| P^v_{>1} \Int^v_{u\rightarrow v} \LON[K][-] \Big\|_{\Cprod{s}{1/2+\delta}} 
\lesssim N^{-1}   K^{1-\delta} \Big\| \LON[K][-] \Big\|_{\Cprod{s}{-1/2+\delta}} 
\lesssim N^{-1} K^{1-\delta} K^{(1-\delta)(\delta+1/2-s)} \Dc. 
\end{align*}
Since $K\lesssim N$ and $\delta^2\gg 1/2-s$, this yields an acceptable contribution to \eqref{modulation:eq-structure-p9}. This completes our proofs of \eqref{modulation:eq-structure-p7}, \eqref{modulation:eq-structure-p8}, and \eqref{modulation:eq-structure-p9}, and therefore our proof of \eqref{modulation:eq-structure-p1}. \\ 

\emph{Step 2:} In the second step, we prove that 
\begin{align}
\Big\|  \Big[  \UN[K][+] \Para[v][ll] \Big( \Int^v_{u\rightarrow v} \LON[K][-] 
-  \IVN[<K^{1-\delta}][-] 
-  \Int^v_{u\rightarrow v} P^v_{<K^{1-\deltap}} 
\VN[<K^{1-\deltap}][\fcs] \Big) \Big]_{\leq N} \Big\|_{\Cprod{\gamma}{1/2+\delta}} 
\lesssim K^{\gamma+1/2+\eta} \Dc^2. \label{modulation:eq-structure-p10}
\end{align}
Due to the definition of $\Para[v][ll]$, we may replace $\Para[v][ll]$ in \eqref{modulation:eq-structure-p10} with $\Para[v][ll]P^v_{>1}$.
Due to Lemma \ref{ansatz:lem-frequency-support}, the argument in \eqref{modulation:eq-structure-p10} only contains high$\times$low-interactions in the $u$-variable. Due to  our paraproduct estimate (Lemma \ref{prelim:lem-paraproduct}), it suffices to prove that
\begin{align}
\Big\| P^v_{>1} \Big( \Int^v_{u\rightarrow v} \LON[K][-] 
-  \IVN[<K^{1-\delta}][-] 
-  \Int^v_{u\rightarrow v} P^v_{<K^{1-\deltap}} 
\VN[<K^{1-\deltap}][\fcs] \Big) \Big\|_{\Cprod{\eta}{1/2+\delta}} 
\lesssim \Dc. 
\end{align}
Using the definition of $\LON[K][-]$ and Lemma \ref{prelim:lem-Duhamel-integral}, it then remains to prove the four estimates 
\begin{align}
\Sumlarge_{\substack{M\leq \Nd \colon \\ M<K^{1-\delta}}} 
\Big\| \Int^v_{u\rightarrow v} \VN[M][-] - \IVN[M][-] \Big\|_{\Cprod{\eta}{1/2+\delta}} &\lesssim \Dc,  \label{modulation:eq-structure-p11}  \\ 
\Sumlarge_{\substack{M_u,M_v\leq \Nd \colon \\ M_u \simeq_\delta M_v, \\ M_u,M_v<K^{1-\delta}}} 
\Big\| \VN[M_u,M_v][+-] \Big\|_{\Cprod{\eta}{-1/2+\delta}} &\lesssim \Dc^2, 
\label{modulation:eq-structure-p12}  \\ 
\Sumlarge_{\substack{M\leq \Nd \colon \\ M<K^{1-\delta}}} 
\Big\| \VN[M][+] \Big\|_{\Cprod{\eta}{-1/2+\delta}} &\lesssim \Dc^2, 
\label{modulation:eq-structure-p13}  \\ 
\Sumlarge_{\substack{M\leq \Nd \colon \\ M<K^{1-\delta}}} 
\Big\| \VN[M][\fs-] \Big\|_{\Cprod{\eta}{-1/2+\delta}} &\lesssim \Dc^2. 
\label{modulation:eq-structure-p14}  
\end{align}
The first estimate \eqref{modulation:eq-structure-p11} follows directly from Lemma \ref{modulation:lem-integration}. The second estimate \eqref{modulation:eq-structure-p12} follows directly from Lemma \ref{modulation:lem-bilinear}. The third estimate \eqref{modulation:eq-structure-p13} follows directly from the Lemma \ref{ansatz:lem-frequency-support} and Lemma \ref{modulation:lem-linear-reversed}. Finally, the fourth estimate \eqref{modulation:eq-structure-p14} follows directly from the Lemma \ref{ansatz:lem-frequency-support} and Lemma \ref{modulation:lem-mixed}. This completes the proof of \eqref{modulation:eq-structure-p11}-\eqref{modulation:eq-structure-p14} and hence the proof of \eqref{modulation:eq-structure-p10}. \\

\emph{Step 3:} In the third step, we prove that
\begin{equation}\label{modulation:eq-structure-p15}
\Big\| \Big[ \UN[K][+] \Para[v][ll] 
\Int^v_{u\rightarrow v} P^v_{<K^{1-\deltap}} \VN[<K^{1-\delta}][\fcs] \Big]_{\leq N} \Big\|_{\Cprod{\gamma}{r}}
\lesssim K^{\gamma+1/2+\eta} \Dc^2. 
\end{equation}
As before, we may replace $\Para[v][ll]$ in \eqref{modulation:eq-structure-p15} by $\Para[v][ll]P^v_{>1}$. 
Due to the Lemma \ref{ansatz:lem-frequency-support}, the argument in \eqref{modulation:eq-structure-p15} contains only high$\times$low-interactions in the $u$-variable. Using our paraproduct estimate (Lemma \ref{prelim:lem-paraproduct}) and Corollary \ref{modulation:cor-control-combined}, it then follows that 
\begin{align*}
&\Big\| \Big[ \UN[K][+] \Para[v][ll] 
P^v_{>1} \Int^v_{u\rightarrow v} P^v_{<K^{1-\deltap}} \VN[<K^{1-\delta}][\fcs] \Big]_{\leq N} \Big\|_{\Cprod{\gamma}{r}} \\
\lesssim\,& \big\| \UN[K][+] \big\|_{\Cprod{\gamma}{s}} 
\Big\| P^v_{>1}\Int^v_{u\rightarrow v} P^v_{<K^{1-\deltap}} \VN[<K^{1-\delta}][\fcs]  \Big\|_{\Cprod{s}{r}} 
\lesssim \big\| \UN[K][+] \big\|_{\Cprod{\gamma}{s}} 
\big\| \VN[<K^{1-\delta}][\fcs]  \Big\|_{\Cprod{s}{r-1}}
\lesssim  K^{\gamma+1/2+\eta} \Dc^2. 
\end{align*}

\emph{Step 4:} In the fourth and final step, we finish up the proof. The first estimate \eqref{modulation:eq-structure-e1} follows from \eqref{modulation:eq-structure-p1}, \eqref{modulation:eq-structure-p10}, and the triangle inequality. The second estimate \eqref{modulation:eq-structure-e2} then follows \eqref{modulation:eq-structure-e1}, \eqref{modulation:eq-structure-p15}, and the triangle inequality.
\end{proof}

As a first application of Lemma \ref{modulation:lem-structure}, we now control resonant parts of the $(+)$$\times$$(-)$-interaction. 

\begin{lemma}\label{modulation:lem-pm-para}
Let the post-modulation hypothesis (Hypothesis \ref{hypothesis:post}) be satisfied and let $K,M \in \Dyadiclarge$ satisfy $K\simeq_\delta M$. Then, it holds that 
\begin{equation}\label{modulation:eq-pm-para}
\Big\| \Big[ \UN[K][+] \Para[v][gtrsim] \VN[M][-] \Big]_{\leq N} 
\Big\|_{\Cprod{r-1}{r-1}} 
\lesssim M^{-\delta +2 (r-s)} \Dc^2.
\end{equation}
\end{lemma}

\begin{proof}
Using Lemma \ref{modulation:lem-structure}, it holds that
\begin{align}\label{modulation:eq-pm-para-p0}
\UN[K][+] 
= \chinull[K] \Big[ \UN[K][+] \Para[v][ll] \IVN[<K^{1-\delta}][-] \Big]_{\leq N}
+ \chinull[K] \Big[ \UN[K][+] \Para[v][ll] \Int^v_{u\rightarrow v} P^v_{<K^{1-\deltap}} \VN[<K^{1-\delta}][\fcs] \Big]_{\leq N} 
+ \UN[K][\#], 
\end{align}
where the remainder $\UN[K][\#]$ satisfies
\begin{equation}\label{modulation:eq-pm-para-p1}
\big\| \UN[K][\#] \big\|_{\Cprod{r-1}{1/2+\delta}} \lesssim K^{r-s} \Dc. 
\end{equation}
Using Lemma \ref{prelim:lem-chi-commutator}, the $\chinull[K]$-factors in \eqref{modulation:eq-pm-para-p0} can easily be commuted with the $\Para[v][gtrsim]$-operator in \eqref{modulation:eq-pm-para}. Thus, it remains to estimate 
\begin{align}
&\, \chinull[K] \Big[ \Big[ \UN[K][+] \Para[v][ll] \IVN[<K^{1-\delta}][-] \Big]_{\leq N} \Para[v][gtrsim] \VN[M][-] \Big]_{\leq N} \label{modulation:eq-pm-para-p2} \\ 
 +&\, \chinull[K] \Big[ \Big[ \UN[K][+] \Para[v][ll] \Int^v_{u\rightarrow v} P^v_{<K^{1-\deltap}} \VN[<K^{1-\delta}][\fcs] \Big]_{\leq N}
  \Para[v][gtrsim] \VN[M][-] \Big]_{\leq N}
  \label{modulation:eq-pm-para-p3} \\
  +&\,  \Big[ \UN[K][\#] \Para[v][gtrsim] \VN[M][-] \Big]_{\leq N}
 \label{modulation:eq-pm-para-p4}. 
\end{align}
We now treat \eqref{modulation:eq-pm-para-p2}, \eqref{modulation:eq-pm-para-p3}, and \eqref{modulation:eq-pm-para-p4} separately. We start with our estimate of \eqref{modulation:eq-pm-para-p2}, which is the most difficult term. To estimate \eqref{modulation:eq-pm-para-p2}, we further decompose
\begin{align}
\eqref{modulation:eq-pm-para-p2}
&=  \chinull[K] \Big[ \Big[ \UN[K][+] \Para[v][ll] \IVN[<K^{1-\delta}][-] \Big]_{\leq N} \Para[v][gg] \VN[M][-] \Big]_{\leq N} 
\label{modulation:eq-pm-para-p5} \\ 
&+ \chinull[K] \Big[ \Big[ \UN[K][+] \Para[v][ll] \IVN[<K^{1-\delta}][-] \Big]_{\leq N} \Para[v][sim] \VN[M][-] \Big]_{\leq N}. 
\label{modulation:eq-pm-para-p6} 
\end{align}
Using our para-product estimate (Lemma \ref{prelim:lem-paraproduct}) and Lemma \ref{modulation:lem-linear}, we obtain that
\begin{align*}
\big\| \eqref{modulation:eq-pm-para-p5} \big\|_{\Cprod{r-1}{r-1}}
&\lesssim 
\Big\|  \UN[K][+] \Para[v][ll] \IVN[<K^{1-\delta}][-] 
\Big\|_{\Cprod{r-1}{s}} 
\Big\| \VN[M][-] \Big\|_{\Cprod{s}{r-1-s+\eta}} \\ 
&\lesssim \Big\| \UN[K][+] \Big\|_{\Cprod{r-1}{s}} 
\Big\| \IVN[<K^{1-\delta}][-]  \Big\|_{\Cprod{s}{s}}
\Big\| \VN[M][-] \Big\|_{\Cprod{s}{r-1-s+\eta}} \\
&\lesssim K^{r-s} M^{-1/2+r-s+2\eta} \Dc^3. 
\end{align*}
Since $K\simeq_\delta M$, this is acceptable. By inserting the integral representation of $P_{\leq N}^x$, it follows that
\begin{align*}
&\Big\| \eqref{modulation:eq-pm-para-p6} \big\|_{\Cprod{r-1}{r-1}} \\ 
\lesssim\, \, \, & \int_{\R^3}  \dy_1 \dy_2 \dy_3 \Big( \prod_{j=1}^3 N \langle N y_j \rangle^{-10} \Big) 
\Big\| \Big( \Theta^x_{y_1} \UN[K][+] \Para[v][ll] \Theta^x_{y_2} \IVN[<K^{1-\delta}][-] \Big) \Para[v][sim] \Theta^x_{y_3} \VN[M][-] \Big\|_{\Cprod{r-1}{r-1}}. 
\end{align*}
Using Lemma \ref{prelim:lem-para-product-trilinear}, it holds that 
\begin{align}
&\Big\| \Big( \Theta^x_{y_1} \UN[K][+] \Para[v][ll] \Theta^x_{y_2} \IVN[<K^{1-\delta}][-] \Big) \Para[v][sim] \Theta^x_{y_3} \VN[M][-] \Big\|_{\Cprod{r-1}{r-1}} \notag \\
\lesssim&\, \Big\| \Theta^x_{y_1} \UN[K][+] \Big\|_{\Cprod{r-1}{s}}
\Big\| \Theta^x_{y_2} \IVN[<K^{1-\delta}][-]  \Para[v][sim] \Theta^x_{y_3} \VN[M][-] \Big\|_{\Cprod{s}{r-1}} \label{modulation:eq-pm-para-p7} \\
 +&\,  \Big\| \Theta^x_{y_1} \UN[K][+] \Big\|_{\Cprod{r-1}{s}} 
\Big\| \Theta^x_{y_2} \IVN[<K^{1-\delta}][-] \Big\|_{\Cprod{s}{s}} 
\Big\| \Theta^x_{y_3} \VN[M][-] \Big\|_{\Cprod{s}{-2s+\eta}}
\label{modulation:eq-pm-para-p8}. 
\end{align}
Using Lemma \ref{modulation:lem-linear}, it is easy to see that \eqref{modulation:eq-pm-para-p8} yields an acceptable contribution, and we therefore turn our attention to  \eqref{modulation:eq-pm-para-p7}. 
Using Lemma \ref{modulation:lem-linear} and Lemma \ref{prelim:lem-insertion-parasim}, it holds that 
\begin{align}
\eqref{modulation:eq-pm-para-p7}
&\lesssim K^{r-s} \Dc
\Sumlarge_{\substack{L<K^{1-\delta}\colon \\ L\sim M}} 
\Big\| \Theta^x_{y_2} \IVN[L][-]  \Para[v][sim] \Theta^x_{y_3} \VN[M][-] \Big\|_{\Cprod{s}{r-1}} \notag \\ 
&\lesssim K^{r-s} \Dc \Sumlarge_{\substack{L<K^{1-\delta}\colon \\ L\sim M}}
\Big( \Big\| \Theta^x_{y_2} \IVN[L][-]  \otimes \Theta^x_{y_3} \VN[M][-] \Big\|_{\Cprod{s}{r-1}} + M^{r-2s} \Dc^2 \Big). \label{modulation:eq-pm-para-p9}
\end{align}
Using Lemma \ref{killing:lem-tensor-modulated-linear} and $L<K^{1-\delta}\leq M$, which rules out a probabilistic resonance between $\IVN[L][-]$ and $\VN[M][-]$, it follows that 
\begin{equation*}
\eqref{modulation:eq-pm-para-p9} \lesssim 
K^{r-s} \Big( M^{r-1+\eta} + M^{r-2s} \Big) \Dc^3 \lesssim K^{r-s} M^{r-2s} \Dc^3,
\end{equation*}
which is acceptable. 
This completes our estimate of \eqref{modulation:eq-pm-para-p2} and we now turn to \eqref{modulation:eq-pm-para-p3}. In fact, we do not need to estimate  \eqref{modulation:eq-pm-para-p3} since it is identically zero. The reason is that $ P^v_{<K^{1-\deltap}} \VN[<K^{1-\delta}][\fcs]$ is supported on $v$-frequencies $\lesssim K^{1-\deltap}$, $\VN[M][-]$ is supported on $v$-frequencies $\sim M \gtrsim K^{1-\delta}$, and $\deltap>\delta$. It now only remains to control \eqref{modulation:eq-pm-para-p4}. Using our paraproduct estimate (Lemma \ref{prelim:lem-paraproduct}), Lemma \ref{modulation:lem-linear}, and \eqref{modulation:eq-pm-para-p1}, it holds that 
\begin{equation*}
\Big\| \Big[ \UN[K][\#] \Para[v][gtrsim] \VN[M][-] \Big]_{\leq N} 
\Big\|_{\Cprod{r-1}{r-1}}
\lesssim \Big\|  \UN[K][\#]  \Big\|_{\Cprod{r-1}{1/2+\delta}}
\Big\| \VN[M][-] \Big\|_{\Cprod{s}{-1/2-\delta+\eta}}
\lesssim K^{r-s} M^{-\delta+2\eta} \Dc^2. 
\end{equation*}
Since $K\simeq_\delta M$, this is acceptable. 
\end{proof}

The next lemma is related to Proposition \ref{modulation:prop-resonant}, but assumes stronger conditions on the terms involved in the interaction and their frequency-support properties. Due to this, this next lemma yields better decay in the highest frequency scale. 

\begin{lemma}\label{modulation:lem-ivm-lo}
Assume that the post-modulation hypothesis (Hypothesis \ref{hypothesis:post}) is satisfied. Furthermore,  let $K,M\in \Dyadiclarge$ satisfy
\begin{equation}\label{modulation:eq-sharp-condition}
M\geq K^{1-\delta}. 
\end{equation}
Then, it holds that
\begin{equation}\label{modulation:eq-ivm-lo}
\sup_{y,z\in \R} \Big\| \Theta^x_y \IVN[M][-] \otimes \Theta^x_z \LON[K][-] \Big\|_{\Cprod{s}{r-1}} \lesssim M^{-1/2+\delta} \Dc^2. 
\end{equation}
\end{lemma}

\begin{remark}
The precise form of the condition $M\geq K^{1-\delta}$ in \eqref{modulation:eq-sharp-condition} is crucial. Since $\LON[K][-]$ contains $\VN[L][-]$-terms for all $L<K^{1-\delta}$, \eqref{modulation:eq-ivm-lo} then only contains interactions between $\IVN[M][-]$ and  $\VN[L][-]$ when $L<M$, but not when $L=M$. 
\end{remark}

\begin{proof}
We recall from Definition \ref{ansatz:def-lo} that 
\begin{equation}\label{modulation:eq-ivm-lo-p1}
\begin{aligned}
\LON[K][-] 
&= \Sumlarge_{L<K^{1-\delta}} \VN[L][-] 
+ \Sumlarge_{L_u,L_v <K^{1-\delta}} \VN[L_u,L_v][+-] 
+ \Sumlarge_{L<K^{1-\delta}} \VN[L][+] \\ 
&+ \Sumlarge_{L<K^{1-\delta}} P^u_{<K^{1-\deltap}} \VN[L][\fs-] 
+ \Sumlarge_{L<K^{1-\delta}} P^v_{<K^{1-\deltap}} \VN[L][+\fs] 
+ P^{u,v}_{<K^{1-\deltap}} \VN[][\fs]. 
\end{aligned}
\end{equation}
We now distinguish four different cases. \\ 

\emph{Case 1: Estimate for $\VN[L][-]$.} 
Since $L<K^{1-\delta}\leq M$, it follows from Lemma \ref{killing:lem-tensor-modulated-linear}
\begin{equation}\label{modulation:eq-ivm-vm}
\sup_{y,z\in \R} \Big\| \Theta^x_y \IVN[M][-] \otimes \Theta^x_z \VN[L][-] \Big\|_{\Cprod{s}{r-1}}
\lesssim \coup \Ac^2 \Bc^2 M^{r-1/2+\eta} M^{-1/2} = M^{-1/2+\delta_1+\delta_3} \Dc^2,
\end{equation}
which is acceptable. \\ 

\emph{Case 2: Estimate for $\VN[L_u,L_v][+-]$.} Since $L_v <K^{1-\delta}\leq M$, it follows from \eqref{modulation:eq-ivm-vm} and Lemma \ref{modulation:lem-linear} that
\begin{align*}
&\sup_{y,z\in \R} \Big\| \Theta^x_y \IVN[M][-] \otimes \Theta^x_z \VN[L_u,L_v][+-] \Big\|_{\Cprod{s}{r-1}} \\
\lesssim\, & \sup_{y,z_1,z_2\in \R} 
\Big\| \Theta^x_y \IVN[M][-] \otimes \Theta^x_{z_1} \IUN[L_u][+] \otimes \Theta^x_{z_2} \VN[L_v][-] 
\Big\|_{\Cprod{s}{r-1}} \\ 
\lesssim\, & \sup_{z_1\in \R} \Big\|\Theta^x_{z_1} \IUN[L_u][+]  \Big\|_{\Cprod{s}{s}} 
\sup_{y,z_2 \in \R} \Big\| \Theta^x_y \IVN[M][-] \otimes \Theta^x_{z_2} \VN[L_v][-] \Big\|_{\Cprod{s}{r-1}} \\
\lesssim \, & \Dc \times M^{r-1/2+\eta} M^{-1/2} \Dc^2 = M^{-1/2+\delta_1+\delta_3} \Dc^3,
\end{align*}
which is acceptable. \\ 

\emph{Case 3: Estimate for $\VN[L][\fs-]$.} From the definition of $\VN[L][\fs -]$ (Definition \ref{ansatz:def-mixed}), it follows that
\begin{align}
\VN[L][\fs -] &=  \chinull[L] \Big[ \Int^u_{v\rightarrow u} P^u_{\geq L^{1-\deltap}} \UN[<L^{1-\delta}][\fsc] 
\Para[u][gg] \VN[L][-] \Big]_{\leq N} \notag \\
&=  \chinull[L] \Big[ \Int^u_{v\rightarrow u} P^u_{\geq L^{1-\deltap}} \UN[<L^{1-\delta}][\fsc], \VN[L][-] \Big]_{\leq N}
-  \chinull[L] \Big[ \Int^u_{v\rightarrow u} P^u_{\geq L^{1-\deltap}} \UN[<L^{1-\delta}][\fsc] 
\Para[u][lesssim] \VN[L][-] \Big]_{\leq N}. 
\label{modulation:eq-ivm-lo-q1}
\end{align}
We now treat the first and second summand in \eqref{modulation:eq-ivm-lo-q1} separately.  To treat the first summand, we write 
\begin{align}
&\Theta^x_y \IVN[M][-] \otimes \Theta^x_z \Big[ \Int^u_{v\rightarrow u} P^u_{\geq L^{1-\deltap}} \UN[<L^{1-\delta}][\fsc], \VN[L][-] \Big]_{\leq N} \label{modulation:eq-ivm-lo-p2} \\
=&\, \int_{\R^3}  \dz_1 \dz_2 \dz_3 \bigg( \Big( \prod_{j=1}^3 \widecheck{\rho}_{\leq N}(z_j) \Big)   
\times \mathcal{T} \Big( \Theta^x_y \IVN[M][-] \otimes \Theta^x_{z+z_1+z_2} \VN[L][-] 
\otimes \Theta^x_{z+z_1+z_3}  \Int^u_{v\rightarrow u} P^u_{\geq L^{1-\deltap}} \UN[<L^{1-\delta}][\fsc] \Big) \bigg), \notag 
\end{align}
where $\mathcal{T}\colon \frkg^{\otimes 3} \rightarrow \frkg^{\otimes 2}$ is the unique linear map satisfying
\begin{equation*}
\mathcal{T}(A\otimes B \otimes C) = A \otimes [C,B]
\end{equation*}
for all $A,B,C\in \frkg$. Using our product estimate (Corollary \ref{prelim:cor-product}), it follows that
\begin{align}
\big\| \eqref{modulation:eq-ivm-lo-p2} \big\|_{\Cprod{s}{r-1}} 
\lesssim&\,  \sup_{\substack{z_1,z_2,z_3\in \R}} 
\Big\| \Theta^x_y \IVN[M][-] \otimes \Theta^x_{z+z_1+z_2} \VN[L][-] 
\otimes \Theta^x_{z+z_1+z_3}  \Int^u_{v\rightarrow u} P^u_{\geq L^{1-\deltap}} \UN[<L^{1-\delta}][\fsc]  \Big\|_{\Cprod{s}{r-1}} \notag \\
\lesssim&\, \sup_{z_1,z_2\in \R}\Big\| \Theta^x_y \IVN[M][-] \otimes \Theta^x_{z+z_1+z_2} \VN[L][-] \Big\|_{\Cprod{s}{r-1}} 
\times    \Big\|  \Int^u_{v\rightarrow u} P^u_{\geq L^{1-\deltap}} \UN[<L^{1-\delta}][\fsc]  \Big\|_{\Cprod{s}{s}}. \label{modulation:eq-ivm-lo-p3}
\end{align}
Using Lemma \ref{prelim:lem-Duhamel-integral} and Corollary \ref{modulation:cor-control-combined}, it holds that
\begin{equation}\label{modulation:eq-ivm-lo-p4}
 \Big\|  \Int^u_{v\rightarrow u} P^u_{\geq L^{1-\deltap}} \UN[<L^{1-\delta}][\fsc]  \Big\|_{\Cprod{s}{s}}
 \lesssim  \Big\| P^u_{\geq L^{1-\deltap}} \UN[<L^{1-\delta}][\fsc]  \Big\|_{\Cprod{s}{s-1}} 
 \lesssim L^{-(1-\deltap)(r-s)} \Dc. 
\end{equation}
Using Lemma \ref{killing:lem-tensor-modulated-linear} and $L<K^{1-\delta}\leq M$, it also holds that 
\begin{equation}\label{modulation:eq-ivm-lo-p5}
\begin{aligned}
\Big\| \Theta^x_y \IVN[M][-] \otimes \Theta^x_{z+z_1+z_2}  \VN[L][-] \Big\|_{\Cprod{s}{r-1}} 
\lesssim  \max(M,L)^{r-1+\frac{1}{2}+\eta} M^{-\frac{1}{2}} \Dc^2 
\lesssim M^{-\frac{1}{2}+\delta_1+\delta_2} \Dc^2. 
\end{aligned}
\end{equation}
By inserting \eqref{modulation:eq-ivm-lo-p4} and \eqref{modulation:eq-ivm-lo-p5} into \eqref{modulation:eq-ivm-lo-p3}, we obtain that the contribution of the first summand in \eqref{modulation:eq-ivm-lo-q1} to \eqref{modulation:eq-ivm-lo} is acceptable. It therefore remains to treat the second summand in \eqref{modulation:eq-ivm-lo-q1}. Using our para-product estimate (Lemma \ref{prelim:lem-paraproduct}), Duhamel integral estimate (Lemma \ref{prelim:lem-Duhamel-integral}), Lemma \ref{modulation:lem-linear}, and Corollary \ref{modulation:cor-control-combined}, we first obtain that 
\begin{align}
&\Big\| \Big[ \Int^u_{v\rightarrow u} P^u_{\geq L^{1-\deltap}} \UN[<L^{1-\delta}][\fsc] 
\Para[u][lesssim] \VN[L][-] \Big]_{\leq N} \Big\|_{\Cprod{s}{\eta}} 
\lesssim \Big\| \Int^u_{v\rightarrow u} P^u_{\geq L^{1-\deltap}} \UN[<L^{1-\delta}][\fsc]  \Big\|_{\Cprod{\eta}{\eta}} 
\big\| \VN[L][-] \big\|_{\Cprod{s}{\eta}} \notag \\ 
\lesssim&\, \big\|  P^u_{\geq L^{1-\deltap}} \UN[<L^{1-\delta}][\fsc]  \big\|_{\Cprod{-1+\eta}{\eta}} 
\big\| \VN[L][-] \big\|_{\Cprod{s}{\eta}} 
\lesssim L^{(1-\deltap)(-1+\eta-(r-1))} L^{\frac{1}{2}+\eta} \Dc^2. \label{modulation:eq-ivm-lo-p5prime} 
\end{align}
We now note that 
\begin{equation*}
(1-\deltap) \big(-1+\eta-(r-1)\big) + \tfrac{1}{2} + \eta = (1-\delta) \big( -\tfrac{1}{2}- \delta_1 \big) + \tfrac{1}{2} + \mathcal{O}(\delta_2) \leq \tfrac{\delta}{2}.
\end{equation*}
Together with $L\lesssim M$, we can therefore bound the $L$-factor in \eqref{modulation:eq-ivm-lo-p5prime} by $M^{\frac{\delta}{2}}$. We then obtain that 
\begin{align*}
&\, \Big\| \Theta^x_y \IVN[M][-] \otimes \Theta^x_z \Big[ \Int^u_{v\rightarrow u} P^u_{\geq L^{1-\deltap}} \UN[<L^{1-\delta}][\fsc] \Para[u][lesssim] \VN[L][-] \Big]_{\leq N} \Big\|_{\Cprod{s}{r-1}} \\ 
\lesssim&\, \Big\|  \IVN[M][-]  \Big\|_{\Cprod{s}{\eta}} \Big\| \Big[ \Int^u_{v\rightarrow u} P^u_{\geq L^{1-\deltap}} \UN[<L^{1-\delta}][\fsc] \Para[u][lesssim] \VN[L][-] \Big]_{\leq N}  \Big\|_{\Cprod{s}{\eta}} \\ 
\lesssim&\,  M^{-\frac{1}{2}+2\eta}M^{\frac{\delta}{2}} \Dc^3 \lesssim M^{-\frac{1}{2}+\delta} \Dc^3. 
\end{align*}
Thus, the contribution of the second summand in \eqref{modulation:eq-ivm-lo-q1} to \eqref{modulation:eq-ivm-lo} is acceptable.\\

\emph{Case 4: $ \VN[L][+], \VN[L][+\fs]$, and $\VN[][\fs]$.} To treat all three subcases simultaneously, we let 
\begin{equation*}
\VN[][\ast] \in \Big\{ \VN[L][+], P^v_{<K^{1-\deltap}} \VN[L][+\fs] , P^{u,v}_{<K^{1-\deltap}}\VN[][\fs] \Big\}. 
\end{equation*}
Using Hypothesis \ref{hypothesis:post}, Lemma \ref{modulation:lem-mixed}, and Lemma \ref{modulation:lem-linear-reversed}, we obtain that
\begin{equation*}
\big\| \VN[][\ast] \big\|_{\Cprod{s}{s-1}} \lesssim \Dc. 
\end{equation*}
We now claim that $\VN[][\ast]$ is supported on $v$-frequencies $\ll M$. Indeed, due to Lemma \ref{ansatz:lem-frequency-support} and the conditions $L<K^{1-\delta}\leq M$, $\VN[L][+]$ is supported on $v$-frequencies $\lesssim L^{1-\delta+\vartheta} \lesssim M^{1-\delta+\vartheta}$.
Furthermore, $P^v_{<K^{1-\deltap}} \VN[L][+\fs]$ and $P^{u,v}_{<K^{1-\deltap}}\VN[][\fs]$ are supported on $v$-frequencies $\lesssim K^{1-\deltap}$. Due to \eqref{modulation:eq-sharp-condition} and $\deltap=\delta+2\delta_2$, this is also much smaller than $M$. Since $\VN[][\ast]$ is supported on $v$-frequencies much smaller than $M$, the high$\times$low-estimate (Lemma \ref{prelim:lem-paraproduct}) yields that
\begin{align*}
\sup_{y,z\in \R} \Big\| \Theta^x_y \IVN[M][-] \otimes \Theta^x_z \VN[][\ast] \Big\|_{\Cprod{s}{r-1}} 
\lesssim\,  \big\| \IVN[M][-] \big\|_{\Cprod{s}{\eta}} \big\| \VN[][\ast] \big\|_{\Cprod{s}{r-1}} 
\lesssim  M^{-1/2+2\eta} M^{r-s} \Dc^2, 
\end{align*}
which is more than acceptable. 
\end{proof}

\subsection{Initial data of modulated linear waves}

In the last part of this section, we show that the initial data of modulated linear waves is close to (the low-regularity part of) the initial data from Definition \ref{ansatz:def-initial-data}.

\begin{lemma}\label{modulation:lem-initial-modulated-linear}
Let the post-modulation hypothesis (Hypothesis \ref{hypothesis:post}) be satisfied. Then, it holds that 
\begin{equation}\label{modulation:eq-initial-modulated-linear}
\begin{aligned}
 \Big\|  \sum_{K\in \dyadic} \SNin[K][+](u) \big( P_{\leq \Nd} P^\sharp_{R;K} W^{(\Rscript,\coup),+}\big)(u)  - \UN[][+] (u,u) \Big\|_{\C_u^{r-1}} 
\lesssim \hcoup \Ac \Bcin.
\end{aligned}
\end{equation}
\end{lemma}

\begin{remark}
Since $\Bcin \lesssim \Bc$, the right-hand side of \eqref{modulation:eq-initial-modulated-linear} can be bounded from above by $\Dc$. The reason for not using this estimate is that, in the proof of Proposition \ref{main:prop-null-lwp}, the quotient of $\Bcin$ and $\Bc$ will be used to absorb a constant $C=C(\delta_\ast)$. 
\end{remark}

\begin{proof}
Using  Definition \ref{ansatz:def-modulated-linear}, we write the argument in \eqref{modulation:eq-initial-modulated-linear} as  
\begin{align}
 &\sum_{K\in \dyadic} \SNin[K][+](u) \big( P_{\leq \Nd} P^\sharp_{R;K} W^{(\Rscript,\coup),+}\big)(u)  - \UN[][+] (u,u) \notag \\
=& \,  \sum_{\substack{K\colon \\[0.1ex] K < \Nlarge}}  \SNin[K][+](u) \big( P_{\leq \Nd} P^\sharp_{R;K} W^{(\Rscript,\coup),+}\big)(u)
\label{ansatz:eq-modulated-linear-p1} \\ 
+&  \Sumlarge_{K} \SNin[K][+](u) \Big( \big( P_{\leq \Nd} P^\sharp_{R;K} W^{(\Rscript,\coup),+}\big)(u) 
- \hcoup \sum_{u_0 \in \LambdaRR} \sum_{k\in \Z_K} \psiRuK(u) \rhoND(k) G_{u_0,k}^+ e^{iku} \Big) 
\label{ansatz:eq-modulated-linear-p2} \\ 
+& \hcoup \Sumlarge_K  \sum_{u_0 \in \LambdaRR} \sum_{k\in \Z_K}  \rhoND(k) \psiRuK \Big( \SNin[K][+](u) - \SN[K][+][k](u,u) \Big)  G_{u_0,k}^+ e^{iku}. 
\label{ansatz:eq-modulated-linear-p3}
\end{align}
Using a product estimate (Lemma \ref{prelim:lem-paraproduct}) and Hypothesis \ref{hypothesis:probabilistic}.\ref{ansatz:item-hypothesis-sharp}, it holds that
\begin{align*}
 \big\| \eqref{ansatz:eq-modulated-linear-p1} \big\|_{\C_u^{r-1}} 
\lesssim \sum_{\substack{K\colon \\[0.1ex] K < \Nlarge}} \big\|\SNin[K][+] \big\|_{\C_u^s} 
\big\| P^\sharp_{R;K} W^{(\Rscript,\coup),+}\big\|_{\C_u^{r-1}}  
\lesssim  \hcoup \Ac \Bcin  \Big( \sum_{\substack{K\colon \\[0.1ex] K < \Nlarge}}  K^{r-\frac{1}{2}+\eta} \Big) 
\lesssim \hcoup \Ac\Bcin , 
\end{align*}
which is acceptable. Using a product estimate (Lemma \ref{prelim:lem-paraproduct}) and Hypothesis \ref{hypothesis:probabilistic}.\ref{ansatz:item-hypothesis-sharp}, it also holds that
\begin{align*}
&\,  \big\| \eqref{ansatz:eq-modulated-linear-p2} \big\|_{\C_u^{r-1}} \\
\lesssim&\, \Sumlarge_{K} \Big\| \SNin[K][+] \Big\|_{\C_u^s} \Big\|  \big( P_{\leq \Nd} P^\sharp_{R;K} W^{(\Rscript,\coup),+}\big)(u) 
- \hcoup \sum_{u_0 \in \LambdaRR} \sum_{k\in \Z_K} \psiRuK(u) \rhoND(k) G_{u_0,k}^+ e^{iku} \Big\|_{\C_u^{r-1}} \\
\lesssim&\, \hcoup \Ac \Bcin \Big( \Sumlarge_K K^{r-1+\eta} \Big) 
\lesssim  \hcoup \Ac \Bcin, 
\end{align*}
which is also acceptable. 
In order to treat the term \eqref{ansatz:eq-modulated-linear-p3}, we first recall the initial condition $\pSN[K][+][k](u,u)=\SNin[K][+](u)$.
Using Hypothesis \ref{hypothesis:probabilistic}.\ref{ansatz:item-hypothesis-crude}, Proposition \ref{modulation:prop-main}, and Lemma \ref{prelim:lem-psi-sum}, it then follows that 
\begin{align*}
&\, \big\| \eqref{ansatz:eq-modulated-linear-p3} \big\|_{\C_u^{r-1}} 
\lesssim \hcoup \Ac \Sumlarge_K \sum_{k\in \Z_K} \big\|  \pSN[K][+][k](u,u) - \SN[K][+][k](u,u) \big\|_{L_u^\infty} \langle k \rangle^\eta \\ 
\lesssim&\,  \hcoup \Ac \Big( \Sumlarge_K  K  K^{-100} K^\eta \Big) \sup_{K\in \Dyadiclarge} K^{100} \big\| \pSN[K][+][k] - \SN[K][+][k] \big\|_{\Wuv[s][s]} 
\lesssim \hcoup \Ac \Bcin.
\end{align*}
This completes our estimates of the three terms \eqref{ansatz:eq-modulated-linear-p1}, \eqref{ansatz:eq-modulated-linear-p2}, and \eqref{ansatz:eq-modulated-linear-p3}, and hence completes the proof of \eqref{modulation:eq-initial-modulated-linear}.
\end{proof}

\section{\texorpdfstring{\protect{High$\times$high$\rightarrow$low-errors}}{High high to low-errors}}\label{section:hhl}

In this section, we control the high$\times$high$\rightarrow$low-errors from Definition \ref{ansatz:def-hhlerr}, which are the first of five errors listed in Proposition \ref{ansatz:prop-decomposition}. The main estimate of this section is the subject of the following proposition. 

\begin{proposition}[\protect{High$\times$high$\rightarrow$low-errors}]\label{hhl:prop-main}
Let the post-modulation hypothesis (Hypothesis \ref{hypothesis:post}) be satisfied. Then, it holds that
\begin{equation*}
\Big\| \HHLNErr[][] \Big\|_{\Cprod{r-1}{r-1}} \lesssim \Dc^2. 
\end{equation*}
\end{proposition}

The proof of Proposition \ref{hhl:prop-main} relies on combinations of several para-product estimates. One of these para-product estimates, which will be used repeatedly in the proof of Proposition \ref{hhl:prop-main}, is isolated in our next lemma. 

\begin{lemma}\label{hhl:lem-paraproduct}
Let $K,L\in \dyadic$ satisfy $L\gtrsim K^\delta$ and let $(\beta,\gamma)\in \{ (r,s), (s,r) \}$. Then, it holds that 
\begin{equation}\label{hhl:eq-paraproduct}
\begin{aligned}
&\Big\| \Big( P_{\lesssim K}^u f 
\Para[v][ll] P^v_{\geq L} g \Big) \Para[v][sim] h 
- P_{\lesssim K}^u f  \, P^{u,v}_{<K^\delta} \Big( P^v_{\geq L} g \Para[v][sim] h \Big) 
\Big\|_{\Cprod{r-1}{r-1}} \\ 
\lesssim&\, (KL)^{-\eta} 
\big\| f \big\|_{\Cprod{s-1}{s}} 
\big\| g \big\|_{\Cprod{s}{\beta}}
\big\| h \big\|_{\Cprod{s}{\gamma-1}}. 
\end{aligned}
\end{equation}
\end{lemma}

We emphasize that $P_{\lesssim K}^u$ acts in the $u$-variable, whereas $ P^v_{\geq L}$ acts in the $v$-variable. As we will see in the proof of Lemma \ref{hhl:lem-paraproduct}, the condition $L\gtrsim K^\delta$ allows us to trade $v$-regularity of $g$ for $u$-regularity of $f$. 

\begin{proof}[Proof of Lemma \ref{hhl:lem-paraproduct}:]
In order to prove \eqref{hhl:eq-paraproduct}, we first decompose
\begin{align}
& \Big( P_{\lesssim K}^u f 
\Para[v][ll] P^v_{\geq L} g \Big) \Para[v][sim] h 
- P_{\lesssim K}^u f  \, \,  P^{u,v}_{<K^\delta} \Big( P^v_{\geq L} g \Para[v][sim] h \Big)  \notag \\
=\,& \Big( P_{\lesssim K}^u f 
\Para[v][ll] P^v_{\geq L} g \Big) \Para[v][sim] h 
- P_{\lesssim K}^u f  \, \,  \Big( P^v_{\geq L} g \Para[v][sim] h \Big)\label{hhl:eq-para-p1} \\ 
+\,& P_{\lesssim K}^u f  \Para[u][nsim] \Big( \big(1-P^{u,v}_{<K^\delta}\big)  \Big( P^v_{\geq L} g \Para[v][sim] h \Big) \Big) \label{hhl:eq-para-p3} \\
+\, & P_{\lesssim K}^u f  \Para[u][sim] \Big(  \big(1-P^{u,v}_{<K^\delta}\big)  \Big( P^v_{\geq L} g \Para[v][sim] h \Big) \Big)\label{hhl:eq-para-p4}. 
\end{align}
We now estimate \eqref{hhl:eq-para-p1}, \eqref{hhl:eq-para-p3}, and \eqref{hhl:eq-para-p4} separately. \\

\emph{Estimate of \eqref{hhl:eq-para-p1}:} 
 Using Lemma \ref{prelim:lem-para-product-trilinear}, the contribution of \eqref{hhl:eq-para-p1} can be estimated by
\begin{align*}
\big\| \eqref{hhl:eq-para-p1} \big\|_{\Cprod{r-1}{r-1}}
&\lesssim 
\big\| P_{\lesssim K}^u f \big\|_{\Cprod{r-1}{s}}
\big\| P^v_{\geq L} g \big\|_{\Cprod{s}{1-2s+\eta}} 
\big\| h \big\|_{\Cprod{s}{s-1}} \\ 
&\lesssim K^{r-s} L^{1-2s+\eta-s} 
\big\| f \big\|_{\Cprod{s-1}{s}}
\big\| g \big\|_{\Cprod{s}{s}} 
\big\| h \big\|_{\Cprod{s}{s-1}}. 
\end{align*}
Since $\beta,\gamma \geq s$ and $L\gtrsim K^\delta$, this yields an acceptable contribution. \\

\emph{Estimate of \eqref{hhl:eq-para-p3}:} 
Using our non-resonant estimate (Lemma \ref{prelim:lem-paraproduct}), we obtain that
\begin{align*}
\big\| \eqref{hhl:eq-para-p3}\big\|_{\Cprod{r-1}{r-1}}
&\lesssim 
\big\| P_{\lesssim K}^u f \big\|_{\Cprod{r-1}{s}}
\big\| \big(1-P^{u,v}_{<K^\delta}\big)  \big( P^v_{\geq L} g \Para[v][sim] h \big) \big\|_{\Cprod{\eta}{r-1}} \\ 
&\lesssim K^{r-s} \max \big( K^{\delta (\eta-s)}, K^{\delta (r-1-\eta)} \big) 
\big\| f \big\|_{\Cprod{s-1}{s}}
\big\| P^v_{\geq L} g \Para[v][sim] h  \big\|_{\Cprod{s}{\eta}} \\  
&\lesssim K^{r-s} \max \big( K^{\delta (\eta-s)}, K^{\delta (r-1-\eta)} \big)  L^{-\eta} 
\big\| f \big\|_{\Cprod{s-1}{s}} 
\big\| g \big\|_{\Cprod{s}{\beta}}
\big\| h \big\|_{\Cprod{s}{\gamma-1}}. 
\end{align*}
In the last line, we used that $(\beta-\eta)+(\gamma-1)>\eta$. Since $r-1-\eta,\eta-s\leq -\frac{1}{2}+2\delta_1$, it follows that 
\begin{equation*}
K^{r-s} \max \big( K^{\delta (\eta-s)}, K^{\delta (r-1-\eta)} \big) \lesssim K^{-\frac{\delta}{2}+4\delta_1} \lesssim K^{-\eta}. 
\end{equation*}
Thus, the contribution of \eqref{hhl:eq-para-p3} is acceptable.\\

\emph{Estimate of \eqref{hhl:eq-para-p4}:} 
Using our resonant estimate (Lemma \ref{prelim:lem-paraproduct}), we obtain that
\begin{align*}
\big\| \eqref{hhl:eq-para-p4} \big\|_{\Cprod{r-1}{r-1}} 
&\lesssim 
\big\| P_{\lesssim K}^u f \big\|_{\Cprod{-s+\eta}{s}}
\big\| P^v_{\geq L} g \Para[v][sim] h  \big\|_{\Cprod{s}{r-1}}\\ 
&\lesssim K^{1-2s+\eta} L^{1-\beta-\gamma+\eta}
\big\| f \big\|_{\Cprod{s-1}{s}} 
\big\| g \big\|_{\Cprod{s}{\beta}}
\big\| h \big\|_{\Cprod{s}{\gamma-1}}. 
\end{align*}
Since $1-\beta-\gamma=1-s-r$ and $L \gtrsim K^\delta$, 
we obtain that 
\begin{equation*}
K^{1-2s+\eta} L^{1-\beta-\gamma+\eta}
\lesssim K^{1-2s+\eta+\delta (1-s-r+2\eta)} L^{-\eta}. 
\end{equation*}
Using our parameter conditions, it holds that 
\begin{equation*}
1-2s+\eta+\delta (1-s-r+2\eta) 
= 2\delta_2 + \delta_3 + \delta \big( -\delta_1 + \delta_2 + 2\delta_3 \big) = - \delta \delta_1 + \mathcal{O}(\delta_2),
\end{equation*}
and thus the contribution of \eqref{hhl:eq-para-p4} is acceptable.
\end{proof}

\begin{proof}[Proof of Proposition \ref{hhl:prop-main}:]
Due to Definition \ref{ansatz:def-hhlerr} and the symmetry in the $u$ and $v$-variables, it suffices to prove the four estimates
\begin{align}
\Big\| \HHLNErr[][(+-)\times (-)] \Big\|_{\Cprod{r-1}{r-1}}
&\lesssim \Dc^3 \label{hhl:eq-main-p-a}, \\ 
\Big\| \HHLNErr[][(+-)\times (\fcs)] \Big\|_{\Cprod{r-1}{r-1}}
&\lesssim \Dc^3 \label{hhl:eq-main-p-b}, \\
\Big\| \HHLNErr[][(+)\times (\fcs)] \Big\|_{\Cprod{r-1}{r-1}}
&\lesssim \Dc^2 \label{hhl:eq-main-p-c}, \\
\Big\| \HHLNErr[][(+\fs)\times (-)] \Big\|_{\Cprod{r-1}{r-1}}
&\lesssim \Dc^3 \label{hhl:eq-main-p-d}. 
\end{align}
We now address \eqref{hhl:eq-main-p-a}, \eqref{hhl:eq-main-p-b}, \eqref{hhl:eq-main-p-c}, and \eqref{hhl:eq-main-p-d} separately. \\ 

\emph{Proof of \eqref{hhl:eq-main-p-a}:} We first recall from Definition \ref{ansatz:def-hhlerr} that
\begin{align}
 &\HHLNErr[][(+-)\times (-)]  \label{ansatz:eq-main-p-a-0} \\ 
 =\, & \Sumlarge_{\substack{K,L,M\leq \Nd \colon \\ K \simeq_\delta L, \\ 
 \max(K,L)\geq M^{1-\delta} }} 
 \bigg( \Big[ \UN[K,L][+-] \Para[v][sim] \VN[M][-] \Big]_{\leq N} 
 - \mathbf{1} \big\{ L= M \big\} P_{\leq N}^x \HHLN[K,M][v,(-)\times (-)] P_{\leq N}^x \UN[K][+] \bigg). \notag 
 \end{align}
 We now estimate the dyadic components in \eqref{ansatz:eq-main-p-a-0} separately. To this end, we further decompose 
 \begin{align}
&\Big[ \UN[K,L][+-] \Para[v][sim] \VN[M][-] \Big]_{\leq N} 
 - \mathbf{1} \big\{ L= M \big\} P_{\leq N}^x \HHLN[K,M][v,(-)\times (-)] P_{\leq N}^x \UN[K][+] \notag \\ 
=&\, \Big[ \UN[K,L][+-] , \VN[M][-] \Big]_{\leq N} 
 - \mathbf{1} \big\{ L= M \big\} P_{\leq N}^x \HHLN[K,M][v,(-)\times (-)] P_{\leq N}^x \UN[K][+] 
 \label{ansatz:eq-main-p-a-1}  \\ 
-& \, \Big[ \UN[K,L][+-] \Para[v][nsim] \VN[M][-] \Big]_{\leq N} \label{ansatz:eq-main-p-a-2}. 
 \end{align}
 We first treat \eqref{ansatz:eq-main-p-a-1}, which is the more difficult term. Using the definition of $\UN[K,L][+-]$, the integral representation of $P_{\leq N}^x$, and Remark \ref{ansatz:rem-commutativity-chipm}, it holds that
 \begin{align}
&\Big[ \UN[K,L][+-] , \VN[M][-] \Big]_{\leq N}   \notag \\ 
=\, & \chinull[K,M] P_{\leq N}^x \Big[ (P_{\leq N}^x)^2 \Big[ P_{\leq N}^x \UN[K][+], P_{\leq N}^x \IVN[L][-] \Big], 
P_{\leq N}^x \VN[M][-] \Big] \notag \\ 
=\, & \chinull[K,M] P_{\leq N}^x \int_{\R} \dy \big( \widecheck{\rho}_{\leq N} \ast \widecheck{\rho}_{\leq N}\big)(y) 
\Ad \big( P_{\leq N}^x \VN[M][-] \big) \Ad \big( \Theta^x_y P_{\leq N}^x \IVN[L][-] \big) 
\Theta^x_y P_{\leq N}^x \UN[K][+]. 
 \end{align}
Together with the definition of $\HHLN[K,M][v,(-)\times (-)]$, it then follows that
\begin{align*}
\eqref{ansatz:eq-main-p-a-1} 
=\,& \chinull[K,M] P_{\leq N}^x \int_{\R} \dy  \bigg( \big( \widecheck{\rho}_{\leq N} \ast \widecheck{\rho}_{\leq N}\big)(y)\\ 
&\times \big( 1- \mathbf{1} \big\{ L=M \big\} P^{u,v}_{\leq K^{\delta}} \big) \Big( 
\Ad \big( P_{\leq N}^x \VN[M][-] \big) \Ad \big( \Theta^x_y P_{\leq N}^x \IVN[L][-] \big) \Big)
\Theta^x_y P_{\leq N}^x \UN[K][+] \bigg). 
\end{align*}
Since neither the covariance function $\Cf^{(\Ncs)}_M(y)$ nor the Casimir $\Cas$ depends on the null-coordinates $u,v\in \R$, it follows from Definition \ref{killing:def-Wick} that 
\begin{align*}
&\big( 1- \mathbf{1} \big\{ L=M \big\} P^{u,v}_{\leq K^{\delta}} \big) \Big( 
\Ad \big( P_{\leq N}^x \VN[M][-] \big) \Ad \big( \Theta^x_y P_{\leq N}^x \IVN[L][-] \big) \Big) \\
=\, & \big( 1- \mathbf{1} \big\{ L=M \big\} P^{u,v}_{\leq K^{\delta}} \big) \Big( \biglcol 
\Ad \big( P_{\leq N}^x \VN[M][-] \big) \Ad \big( \Theta^x_y P_{\leq N}^x \IVN[L][-] \big) \bigrcol \Big). 
\end{align*}
Using Lemma \ref{killing:lem-tensor-modulated-linear} and Lemma \ref{modulation:lem-linear}, it then follows that
\begin{align*}
\Big\| \eqref{ansatz:eq-main-p-a-1} \Big\|_{\Cprod{r-1}{r-1}} 
\lesssim\,&
\int_{\R} \dy \,  \bigg(
\Big| \big( \widecheck{\rho}_{\leq N} \ast \widecheck{\rho}_{\leq N}\big)(y) \Big| \\
&\times \Big\| \biglcol 
\Ad \big( P_{\leq N}^x \VN[M][-] \big) \Ad \big( \Theta^x_y P_{\leq N}^x \IVN[L][-] \big) \bigrcol 
\Big\|_{\Cprod{s}{r-1}} \Big\| \Theta^x_y \UN[K][+] \Big\|_{\Cprod{r-1}{s}} \bigg) \\
\lesssim \, & \int_{\R} \dy \, N\langle Ny \rangle^{-10} 
\Big( \max(M,L)^{r-1/2+\eta} L^{-1/2} + N^{-\delta+\vartheta} \langle Ny \rangle \Big) K^{r-s} \Dc^3 \\ 
\lesssim\, & \Big( \max(M,L)^{r-1/2+\eta} L^{-1/2} + N^{-\delta+\vartheta} \Big) K^{r-s} \Dc^3. 
\end{align*}
Since $K\simeq_\delta L$ and $\max(K,L)\geq M^{1-\delta}$, this easily yields an acceptable contribution. It now only remains to estimate \eqref{ansatz:eq-main-p-a-2}. Using our non-resonant estimate (Lemma \ref{prelim:lem-paraproduct}), Lemma \ref{modulation:lem-linear}, and Lemma \ref{modulation:lem-bilinear}, it holds that
\begin{equation*}
\Big\| \eqref{ansatz:eq-main-p-a-2} \Big\|_{\Cprod{r-1}{r-1}} 
\lesssim \big\| \UN[K,L][+-] \big\|_{\Cprod{r-1}{\eta}} \big\| \VN[M][-] \big\|_{\Cprod{s}{r-1}} 
\lesssim K^{r-1/2+\eta} L^{-1/2+2\eta} M^{r-1/2+\eta} \Dc^3. 
\end{equation*}
Since $K\simeq_\delta L$ and $\max(K,L)\geq M^{1-\delta}$, this clearly yields an acceptable contribution. \\

\emph{Proof of \eqref{hhl:eq-main-p-b}:}
We first recall from Definition \ref{ansatz:def-hhlerr} that 
\begin{equation}\label{hhl:eq-main-p-b-0}
\begin{aligned}
    &\HHLNErr[][(+-)\times (\fcs)] \\  
    =& \, \Sumlarge_{\substack{K,M \leq \Nd \colon \\  K \simeq_\delta M}} 
\bigg( \Big[ \UN[K,M][+-] \Para[v][sim] \VN[<K^{1-\delta}][\fcs] \Big]_{\leq N} - P_{\leq N}^x \HHLN[K,M][v,(\fcs)\times (-)] P_{\leq N}^x \UN[K][+] \bigg). 
\end{aligned}
\end{equation}
We now estimate the dyadic summands in \eqref{hhl:eq-main-p-b-0} individually. To this end, we first decompose
\begin{align}
&\, \Big[ \UN[K,M][+-] \Para[v][sim] \VN[<K^{1-\delta}][\fcs] \Big]_{\leq N}
-  P_{\leq N}^x \HHLN[K,M][\revision{v},(\fcs)\times (-)] P_{\leq N}^x \UN[K][+] \notag \\
=&\, \Big[ \Big( \chinull[K,M] \Big[ \UN[K][+], \IVN[M][-] \Big]_{\leq N} \Big) \Para[v][sim]  \VN[<K^{1-\delta}][\fcs] \Big] 
- \chinull[K,M] \Big[   \Big[ \UN[K][+], \IVN[M][-] \Big]_{\leq N}  \Para[v][sim]  \VN[<K^{1-\delta}][\fcs] \Big] \label{hhl:eq-main-q-b-1} \\ 
+&\, \chinull[K,M] \Big[   \Big[ \UN[K][+], \IVN[M][-] \Big]_{\leq N}  \Para[v][sim]  \VN[<K^{1-\delta}][\fcs] \Big]
-  P_{\leq N}^x \HHLN[K,M][\revision{v},(\fcs)\times (-)] P_{\leq N}^x \UN[K][+]. \label{hhl:eq-main-q-b-2}
\end{align}
The term \eqref{hhl:eq-main-q-b-1} can be easily controlled using Lemma \ref{prelim:lem-chi-commutator} and it remains to treat \eqref{hhl:eq-main-q-b-2}. 
We then write 
\begin{align}
\eqref{hhl:eq-main-q-b-2} =\, & \chinull[K,M]  P_{\leq N}^x \Big[ (P_{\leq N}^x)^2 
\Big[ P_{\leq N}^x \UN[K][+], P_{\leq N}^x \IVN[M][-] \Big] \Para[v][sim] P_{\leq N}^x \VN[<K^{1-\delta}][\fcs] \Big] \notag \\ 
 & -  P_{\leq N}^x \HHLN[K,M][v,(\fcs)\times (-)] P_{\leq N}^x \UN[K][+] \notag \\
 =\, &  \chinull[K,M]  P_{\leq N}^x \Big[ (P_{\leq N}^x)^2 
\Big[ P_{\leq N}^x \UN[K][+] \Para[v][gtrsim] P_{\leq N}^x \IVN[M][-] \Big] \Para[v][sim] P_{\leq N}^x \VN[<K^{1-\delta}][\fcs] \Big] \label{hhl:eq-main-p-b-1} \\ 
+\, &  \bigg(   \chinull[K,M] P_{\leq N}^x \Big[ (P_{\leq N}^x)^2 
\Big[ P_{\leq N}^x \UN[K][+] \Para[v][ll] P_{\leq N}^x \IVN[M][-] \Big] \Para[v][sim] P_{\leq N}^x \VN[<K^{1-\delta}][\fcs] \Big]  \label{hhl:eq-main-p-b-2}  \\ 
 & - P_{\leq N}^x \HHLN[K,M][v,(\fcs)\times (-)] P_{\leq N}^x \UN[K][+] \bigg). \notag 
\end{align}
We now estimate \eqref{hhl:eq-main-p-b-1} and \eqref{hhl:eq-main-p-b-2} separately. For \eqref{hhl:eq-main-p-b-1}, we obtain from our paraproduct estimates (Lemma \ref{prelim:lem-paraproduct}), Lemma \ref{modulation:lem-linear}, and Corollary \ref{modulation:cor-control-combined} that 
\begin{align*}
\Big\| \eqref{hhl:eq-main-p-b-1}\Big\|_{\Cprod{r-1}{r-1}} 
&\lesssim 
\Big\| P_{\leq N}^x \UN[K][+] \Para[v][gtrsim] P_{\leq N}^x \IVN[M][-] \Big\|_{\Cprod{r-1}{s}} 
\Big\|  P_{\leq N}^x \VN[<K^{1-\delta}][\fcs] \Big\|_{\Cprod{s}{r-1}} \\ 
&\lesssim \Big\| \UN[K][+] \Big\|_{\Cprod{r-1}{s}} 
\Big\| \IVN[M][-] \Big\|_{\Cprod{s}{\eta}} 
\Big\|  \VN[<K^{1-\delta}][\fcs] \Big\|_{\Cprod{s}{r-1}} \\ 
&\lesssim  K^{r-s} M^{\eta-s} \Dc^3.
\end{align*}
Since $K\simeq_\delta M$, this is acceptable. 
It now remains to control \eqref{hhl:eq-main-p-b-2}. By inserting the integral representation of $(P_{\leq N}^x)^2$, we obtain that
\begin{align*}
\eqref{hhl:eq-main-p-b-2} 
=&   P_{\leq N}^x \int_{\R} \dy \big( \widecheck{\rho}_{\leq N} \ast \widecheck{\rho}_{\leq N} \big)(y) 
\bigg( \chinull[K,M] \Big[ \Big[ \Theta^x_y P_{\leq N}^x \UN[K][+] \Para[v][ll] \Theta^x_y P_{\leq N}^x \IVN[M][-] \Big] \Para[v][sim] P_{\leq N}^x \VN[<K^{1-\delta}][\fcs] \Big]  \\
&- \HHLN[K,M,y][v,(\fcs)\times (-)] \Theta^x_y P_{\leq N}^x \UN[K][+] \bigg).
\end{align*}
Due to Lemma \ref{ansatz:lem-frequency-support} and $K\simeq_\delta M$, the term $\Theta^x_y P_{\leq N}^x \IVN[M][-]$ is supported on $v$-frequencies $\sim M \gtrsim K^{1-\delta}$. Due to the $\Para[v][sim]$-operator and $\deltap>\delta$, we may therefore replace 
 $\VN[<K^{1-\delta}][\fcs]$ with $P^v_{>K^{1-\deltap}} \VN[<K^{1-\delta}][\fcs]$. Using the definition of $\HHLN[K,M,y][\revision{v},(\fcs)\times (-)]$ and Lemma \ref{hhl:lem-paraproduct}, it then follows that 
\begin{equation*}
\big\| \eqref{hhl:eq-main-p-b-2} \big\|_{\Cprod{r-1}{r-1}}
\lesssim (KM)^{-\eta} \big\| \UN[K][+]\big\|_{\Cprod{s-1}{s}}
\big\| \IVN[M][-]\big\|_{\Cprod{s}{s}} 
\big\| \VN[<K^{1-\delta}][\fcs] \big\|_{\Cprod{s}{r-1}}
\lesssim (KM)^{-\eta} \Dc^3, 
\end{equation*}
which is an acceptable contribution. \\ 

\emph{Proof of \eqref{hhl:eq-main-p-c}:}
We first recall from Definition \ref{ansatz:def-hhlerr} that
\begin{equation}\label{hhl:eq-main-p-c-0}
\begin{aligned}
&\HHLNErr[][(+)\times (\fcs)] \\   
=&\,   \Sumlarge_{\substack{K \leq \Nd}}  
\bigg( \Big[ \UN[K][+] \Para[v][sim] P^v_{\geq K^{1-\deltap}} \VN[<K^{1-\delta}][\fcs] \Big]_{\leq N} - \Sumlarge_{\substack{M \leq \Nd \colon \\  M <K^{1-\delta}}} P_{\leq N}^x \HHLN[K,M][v,(\fcs)\times (-)] P_{\leq N}^x \UN[K][+] \bigg). 
\end{aligned}
\end{equation}
We now control the dyadic summands in \eqref{hhl:eq-main-p-c-0} individually. To this end, we first use Lemma \ref{modulation:lem-structure}, which yields the decomposition 
\begin{equation}\label{hhl:eq-main-p-c-1}
\UN[K][+] =  \Sumlarge_{\substack{M\leq \Nd \colon \\ M<K^{1-\delta}}} \chinull[K,M]
\Big[ \UN[K][+] \Para[v][ll] \IVN[M][-] \Big]_{\leq N} 
+ \UN[K][\#].
\end{equation}
In our derivation of \eqref{hhl:eq-main-p-c-1} from Lemma \ref{modulation:lem-structure}, we also used that $\chinull[K,M]=\chinull[K]$ for all $M<K^{1-\delta}$. The remainder term $\UN[K][\#]$ in \eqref{hhl:eq-main-p-c-1} satisfies
\begin{equation*}
\Big\| \UN[K][\#] \Big\|_{\Cprod{r-1}{r}} \lesssim K^{r-s} \Dc. 
\end{equation*}
Equipped with \eqref{hhl:eq-main-p-c-1}, we then write the dyadic summands from \eqref{hhl:eq-main-p-c-0} as 
\begin{align}
&\Big[ \UN[K][+] \Para[v][sim] P^v_{\geq K^{1-\deltap}} \VN[<K^{1-\delta}][\fcs] \Big]_{\leq N} 
- \Sumlarge_{\substack{M \leq \Nd \colon \\  M <K^{1-\delta}}} P_{\leq N}^x \HHLN[K,M][v,(\fcs)\times (-)] P_{\leq N}^x \UN[K][+] 
\notag \\
=&\, \Sumlarge_{\substack{M\leq \Nd \colon \\ M<K^{1-\delta}}} \bigg( 
\Big[ \Big( \chinull[K,M] \Big[ \UN[K][+] \Para[v][ll] \IVN[M][-] \Big]_{\leq N} \Big) \Para[v][sim]  P^v_{\geq K^{1-\deltap}} \VN[<K^{1-\delta}][\fcs] \Big]_{\leq N} \label{hhl:eq-main-p-c-2}  \\
&\hspace{10ex}-  P_{\leq N}^x \HHLN[K,M][v,(\fcs)\times (-)] P_{\leq N}^x \UN[K][+]  \bigg) \notag \\
+&\, \Big[ \UN[K][\#] \Para[v][sim] P^v_{\geq K^{1-\deltap}} \VN[<K^{1-\delta}][\fcs] \Big]_{\leq N}. \label{hhl:eq-main-p-c-3}  
\end{align}
We now note that, due to the $\Para[v][sim]$-operator, only the terms in \eqref{hhl:eq-main-p-c-2} satisfying $M\gtrsim K^{1-\deltap}$ yield a non-zero contribution. As a result, \eqref{hhl:eq-main-p-c-2} can be estimated exactly as in the proof of \eqref{hhl:eq-main-p-b}. It then remains to control \eqref{hhl:eq-main-p-c-3}, which can be estimated by 
\begin{align*}
\Big\| \Big[ \UN[K][\#] \Para[v][sim] P^v_{\geq K^{1-\deltap}} \VN[<K^{1-\delta}][\fcs] \Big]_{\leq N} \Big\|_{\Cprod{r-1}{r-1}} 
\lesssim& \,  \big\| \UN[K][\#]  \big\|_{\Cprod{r-1}{r}}
\big\| P^v_{\geq K^{1-\deltap}} \VN[<K^{1-\delta}][\fcs]  \big\|_{\Cprod{s}{-r+\eta}} \\ 
\lesssim& \,   \big( K^{1-\deltap} \big)^{-r+\eta-(r-1)} 
\big\| \UN[K][\#]  \big\|_{\Cprod{r-1}{r}}
\big\| P^v_{\geq K^{1-\deltap}} \VN[<K^{1-\delta}][\fcs]  \big\|_{\Cprod{s}{r-1}}  \\
\lesssim&\, K^{r-s} K^{(1-\deltap)(1-2r+\eta)} \Dc^2. 
\end{align*}
Since 
\begin{equation*}
r-s+ (1-\delta)(1-2r+\eta) 
= \delta_1 - \delta_2 + (1-\delta-2\delta_2) (-2\delta_1+\delta_3) = -\delta_1 + \mathcal{O}(\delta \delta_1), 
\end{equation*}
this yields an acceptable contribution. \\ 

\emph{Proof of \eqref{hhl:eq-main-p-d}:}  
We first recall from Definition \ref{ansatz:def-hhl} that
\begin{equation}\label{hhl:eq-main-p-d-00}
\begin{aligned}
&\HHLNErr[][(+\fs)\times (-)] \\   
=&\,  \Sumlarge_{\substack{K,M \leq \Nd \colon\\ M \lesssim_\delta K}}  
\bigg( \Big[ \UN[K][+\fs] \Para[v][sim] \VN[M][-] \Big]_{\leq N} 
-  P_{\leq N}^x \HHLN[K,M][v,(-)\times (\fcs)] P_{\leq N}^x \UN[K][+] \bigg). 
\end{aligned}
\end{equation}
We also recall from Definition \ref{ansatz:def-mixed} that
\begin{equation}\label{hhl:eq-main-p-d-0}
\UN[K][+\fs] 
=  \chinull[K] \Big[  \UN[K][+] \Para[v][ll] 
 \Int^v_{u\rightarrow v} P^v_{\geq K^{1-\deltap}} \VN[<K^{1-\delta}][\fcs] \Big]_{\leq N}. 
\end{equation}
Due to Lemma \ref{prelim:lem-chi-commutator}, the $\chinull[K]$-factor can easily be commuted with the $\Para[v][sim]$-operator in \eqref{hhl:eq-main-p-d-00}. Then, it remains to consider the dyadic summands 
\begin{align} 
&\chinull[K] \Big[ \Big[  \UN[K][+] \Para[v][ll] 
 \Int^v_{u\rightarrow v} P^v_{\geq K^{1-\deltap}} \VN[<K^{1-\delta}][\fcs] \Big]_{\leq N}
 \Para[v][sim] \VN[M][-] \Big]_{\leq N}
- P_{\leq N}^x \HHLN[K,M][v,(-)\times (\fcs)] P_{\leq N}^x \UN[K][+] \label{hhl:eq-main-p-d-1} \\ 
=& P_{\leq N}^x \int_{\R} \dy \, 
(\widecheck{\rho}_{\leq N}\ast \widecheck{\rho}_{\leq N})(y) 
\bigg( \chinull[K]\Big[ \Big[ \Theta^x_y P_{\leq N}^x \UN[K][+] 
\Para[v][ll] \Theta^x_y P_{\leq N}^x \Int^v_{u\rightarrow v} P^v_{\geq K^{1-\deltap}} \VN[<K^{1-\delta}][\fcs] \Big] 
\Para[v][sim] P_{\leq N}^x \VN[M][-] \Big] \notag \\ 
&- \HHLN[K,M,y][v,(-)\times (\fcs)] \Theta^x_y P_{\leq N}^x \UN[K][+]  \bigg).  \notag 
\end{align}
Using the definition of $\HHLN[K,M,y][v,(-)\times (\fcs)]$ and Lemma \ref{hhl:lem-paraproduct}, it then follows that
\begin{align*}
\big\| \eqref{hhl:eq-main-p-d-1} \big\|_{\Cprod{r-1}{r-1}}
\lesssim K^{-\eta} 
\big\| \UN[K][+] \big\|_{\Cprod{s}{s}}
\big\| \Int^v_{u\rightarrow v} P^v_{\geq K^{1-\deltap}} \VN[<K^{1-\delta}][\fcs] \big\|_{\Cprod{s}{r}}
\big\| \VN[M][-] \big\|_{\Cprod{s}{s-1}}. 
\end{align*}
Using Lemma \ref{prelim:lem-Duhamel-integral}, Lemma \ref{modulation:lem-linear}, and Corollary \ref{modulation:cor-control-combined}, it follows that 
\begin{equation*}
K^{-\eta} 
\big\| \UN[K][+] \big\|_{\Cprod{s}{s}}
\big\| \Int^v_{u\rightarrow v} P^v_{\geq K^{1-\deltap}} \VN[<K^{1-\delta}][\fcs] \big\|_{\Cprod{s}{r}}
\big\| \VN[M][-] \big\|_{\Cprod{s}{s-1}}\lesssim K^{-\eta} \Dc^3,
\end{equation*}
which is acceptable. 
\end{proof}

\section{Jacobi errors}\label{section:jacobi}

In this subsection, we examine the Jacobi error $\JcbNErr$ from Definition \ref{ansatz:def-jacobi-errors}. Out of all terms in \mbox{Proposition \ref{ansatz:prop-decomposition}},  the Jacobi error requires the most delicate argument. As previously discussed in Subsection \ref{section:overview-jacobi}, the Jacobi error cannot be controlled using only perturbative arguments and will instead be controlled in parts using a Bourgain-Bulut argument. 

In order to make our bounds on the Jacobi error $\JcbNErr$ precise, we need to make three additional definitions. In the first definition, we introduce variants of the covariance function and Killing-renormalization from Definition \ref{ansatz:def-Killing}.

\begin{definition}[Covariance function and Killing-renormalization near the frequency-boundary]\label{jacobi:def-cf} 
For all $\Nd,N \in \Dyadiclarge$ and $M\in \dyadic$ satisfying $M\geq 2$, we define
\begin{align*}
\CNbd_M(y) &:= - \int_{\R} \dxi \, \Big( \rho_{\leq N}(\xi) - \rho_{<N^{1-2\delta_1}}(\xi) \Big)^2
\mathbf{1} \big\{ \tfrac{M}{2} \leq |\xi|  \leq M \big\} \frac{\sin(\xi y)}{\xi},  \\ 
\CNcsbd_M(y) &:= - \int_{\R} \dxi \, \Big( \rho_{\leq N}(\xi) - \rho_{<N^{1-2\delta_1}}(\xi) \Big)^2 \rho_{\leq \Nd}(\xi)^2
\mathbf{1} \big\{ \tfrac{M}{2} \leq |\xi|  \leq M \big\} \frac{\sin(\xi y)}{\xi}.
\end{align*}
In addition, we define $\CNbd_1=\CNcsbd_1=0$, 
\begin{equation*}
\CNbd(y) := \sum_{\substack{M\in \dyadic }} \CNbd_M(y),
\qquad \text{and} \qquad 
\CNcsbd(y) := \sum_{\substack{M\in \dyadic }} \CNcsbd_M(y). 
\end{equation*}
For any $F\colon \R_x \rightarrow \frkg$, we also define 
\begin{align*}
\Renormbd[\Nscript] F &:= P_{\leq N}^x \int_{\R} \dy\,  (\widecheck{\rho}_{\leq N}\ast \widecheck{\rho}_{\leq N} )(y) \CNbd(y) \Theta^x_y P_{\leq N}^x \Kil F, \\ 
\Renormbd[\Ncs] F &:= P_{\leq N}^x \int_{\R} \dy\,  (\widecheck{\rho}_{\leq N}\ast \widecheck{\rho}_{\leq N} )(y) \CNcsbd(y) \Theta^x_y P_{\leq N}^x \Kil F. 
\end{align*}
\end{definition}

As previously mentioned in Subsection \ref{section:overview-jacobi}, the Bourgain-Bulut argument requires us to use Cartesian-coordinates. For this reason, we will need Cartesian versions of the high$\times$high$\rightarrow$low-matrices and operators from Definition \ref{ansatz:def-hhl}.

\begin{definition}[\protect{Cartesian high$\times$high$\rightarrow$low-matrices and operators}]\label{jacobi:def-chhl}
Let $N\in \Dyadiclarge$. For all $A,B\colon \R_x \rightarrow \frkg$ and $y\in \R$,  we define the Cartesian high$\times$high$\rightarrow$low-matrix $\CHHLN_y (A,B)\colon \R_x \rightarrow \End(\frkg)$ by 
\begin{equation}\label{jacobi:eq-chhl-y}
\CHHLN_y(A,B) := P^x_{\leq N^{\deltap}} 
\Big( \ad\big( \Pbd A \big) \ad \big( (\Theta^x_y -1) P^x_{>1}\Int^x_{0\rightarrow x} \Pbd B\big) \Big).
\end{equation}
For any $F\colon \R_x \rightarrow \frkg$, we also define
\begin{equation}\label{jacobi:eq-chhl}
\CHHLN(A,B) F :=  \int_{\R} \dy\,  (\widecheck{\rho}_{\leq N}\ast \widecheck{\rho}_{\leq N} )(y) \CHHLN_y(A,B) \Theta^x_y F.
\end{equation}
\end{definition}

\begin{remark}\label{jacobi:rem-acting-spatial-variable}
In the proof of Proposition \ref{main:prop-null-lwp}, it is crucial that \eqref{jacobi:eq-chhl-y} and \eqref{jacobi:eq-chhl} only act on the spatial coordinates. In particular, it is crucial that neither \eqref{jacobi:eq-chhl-y} nor \eqref{jacobi:eq-chhl} contains Littlewood-Paley operators in the time-coordinate $t$ or the null-coordinates $u$ or $v$. 

A technical difference between Definition \ref{ansatz:def-hhl} and Definition \ref{jacobi:def-chhl} is that Definition \ref{ansatz:def-hhl} includes $\chi$-factors, whereas Definition \ref{jacobi:def-chhl} does not include $\chi$-factors. This choice has been made to emphasize that $\CHHLN$ only acts on the spatial variables, but does lead to inconvenient $\chi$-factors in Proposition \ref{jacobi:prop-shift}.
\end{remark}

In the last definition needed for Proposition \ref{jacobi:prop-main}, we introduce modified Jacobi errors.

\begin{definition}[Modified Jacobi errors]\label{jacobi:def-modified-jacobi} 
We define the modified Jacobi-error $\JcbNErr[][+,\dagger]$ by
\begin{equation}\label{jacobi:eq-modified-jacobi-plus}
\JcbNErr[][+,\dagger] = \JcbNErr[][+] - \chi 
\Big( P_{\leq N}^x \CHHLN \big( \VN[][], \VN[][] \big) P_{\leq N}^x - \coup \Renormbd[\Ncs] \Big) \Big( \hspace{-2ex} \Sumlarge_{\substack{\hspace{1.5ex} K \geq N^{1-2\delta_1}}} \hspace{-2ex}  \UN[K][+]\Big).
\end{equation}
Similarly, we define
\begin{equation}\label{jacobi:eq-modified-jacobi-minus}
\JcbNErr[][-,\dagger] = \JcbNErr[][-] - \chi 
\Big( P_{\leq N}^x \CHHLN \big( \UN[][], \UN[][] \big) P_{\leq N}^x - \coup \Renormbd[\Ncs] \Big) \Big( \hspace{-2ex} \Sumlarge_{\substack{\hspace{1.5ex}M \geq N^{1-2\delta_1}}} \hspace{-2ex}  \VN[M][-] \Big).
\end{equation}
Finally, we define
\begin{equation}\label{jacobi:eq-modified-jacobi}
\JcbNErr[][\dagger] := \JcbNErr[][+,\dagger]  + \JcbNErr[][-,\dagger] . 
\end{equation}
\end{definition}

\begin{remark}[Motivation behind modified Jacobi errors]
Due to Definition \ref{jacobi:def-modified-jacobi}, we can write 
\begin{align}
\JcbNErr[][+]  
&:= \JcbNErr[][+,\dagger]  \label{jacobi:eq-modified-jacobi-rem-e1} \\ 
&+ \chi 
\Big( P_{\leq N}^x \CHHLN \big( \VN[][], \VN[][] \big) P_{\leq N}^x - \coup \Renormbd[\Ncs] \Big) \Big( \hspace{-2ex} \Sumlarge_{\substack{\hspace{1.5ex} K \geq N^{1-2\delta_1}}} \hspace{-2ex}  \UN[K][+] \Big).  \label{jacobi:eq-modified-jacobi-rem-e2}
\end{align}
The first term \eqref{jacobi:eq-modified-jacobi-rem-e1} will be controlled using perturbative arguments. In contrast, the second term \eqref{jacobi:eq-modified-jacobi-rem-e2} will be controlled using a Bourgain-Bulut argument. For this, it is crucial that 
\begin{equation*}
P_{\leq N}^x \CHHLN \big( \VN[][], \VN[][] \big) P_{\leq N}^x - \coup \Renormbd[\Ncs] 
\end{equation*}
only depends on the full solution $\VN[][]$, but not any of the individual terms from our Ansatz in \eqref{ansatz:eq-VN-rigorous-decomposition}.
\end{remark}

Equipped with Definition \ref{jacobi:def-cf}, Definition \ref{jacobi:def-chhl}, and Definition \ref{jacobi:def-modified-jacobi}, we can now state the main result of this section. 

\begin{proposition}[Jacobi errors]\label{jacobi:prop-main} 
Let the post-modulation hypothesis (Hypothesis \ref{hypothesis:post}) be satisfied. Then, it holds that 
\begin{align}
\Big\| \JcbNErr[][] \Big\|_{\Cprod{r-1}{r-1}} &\lesssim \bigg( 1+ \frac{\Nd}{N} \,  N^{\delta_2} \bigg) \Dc^3,
\label{jacobi:eq-main-estimate-1} \\
\Big\| \JcbNErr[][\dagger] \Big\|_{\Cprod{r-1}{r-1}} &\lesssim  \Dc^3. 
\label{jacobi:eq-main-estimate-2} 
\end{align}
Furthermore, for all $U,V\colon \R^{1+1}\rightarrow \frkg$, it holds that 
\begin{equation}\label{jacobi:eq-bourgain-bulut-1}
\begin{aligned}
&\, \Sumlarge_{\substack{\hspace{1.5ex} K \geq N^{1-2\delta_1}}} \Big\| \, \chi 
\Big( P_{\leq N}^x \CHHLN \big( V, V \big) P_{\leq N}^x - \coup \Renormbd[\Ncs] \Big) \big(   \UN[K][+] \big) \Big\|_{\Cprod{r-1}{r-1}} \\
\lesssim&\, N^{\delta_1+5\delta_3} \Dc \sup_{y\in \R} \, \langle Ny \rangle^{-10} \Big\|  \chi \Big( \CHHLN_y\big( V, V \big) - \coup \CNcsbd(y) \Kil \Big) \Big\|_{L_t^\infty L_x^\infty}
\end{aligned}
\end{equation}
and 
\begin{equation}\label{jacobi:eq-bourgain-bulut-2}
\begin{aligned}
&\, \Sumlarge_{\substack{\hspace{1.5ex}M \geq N^{1-2\delta_1}}} \Big\| \, \chi
\Big( P_{\leq N}^x \CHHLN \big( U, U \big) P_{\leq N}^x - \coup \Renormbd[\Ncs] \Big) \big(  \VN[M][-] \big) \Big\|_{\Cprod{r-1}{r-1}} \\ 
\lesssim&\, N^{\delta_1+5\delta_3} \Dc \sup_{y\in \R} \, \langle Ny \rangle^{-10} \Big\|  \chi \Big( \CHHLN_y\big( U, U \big) - \coup \CNcsbd(y) \Kil \Big) \Big\|_{L_t^\infty L_x^\infty}. 
\end{aligned}
\end{equation}
\end{proposition}

We first comment on the nature of \eqref{jacobi:eq-main-estimate-1} and \eqref{jacobi:eq-main-estimate-2}. Due to the $N^{\delta_2}$-loss, the first estimate \eqref{jacobi:eq-main-estimate-1} generally cannot be used to close our contraction-mapping argument for the remainder equations. The only exceptions are when either $\Nd$ is much smaller than $N$ or when $N^{\delta_2}\Dc^2 \ll 1$, i.e., the regime in which $\Dc$ is small depending on $N$. By scaling, the latter can be used to prove that our Ansatz from \eqref{ansatz:eq-UN-decomposition} and \eqref{ansatz:eq-VN-decomposition} holds on the timescale $\sim N^{-\delta_2}$, i.e., on a timescale which depends on $N$. This is useful since this timescale is still larger than the size of the energy increment  (see Section \ref{section:energy-increment}), which allows us to prove the almost invariance of the Gibbs measure.

In contrast to the right-hand side of \eqref{jacobi:eq-main-estimate-1}, the right-hand side of \eqref{jacobi:eq-main-estimate-2} is always under control. In order to control the original Jacobi error $\JcbNErr[][]$, it therefore remains to control the right-hand sides of
\eqref{jacobi:eq-bourgain-bulut-1} and \eqref{jacobi:eq-bourgain-bulut-2} with $U=\UN[][]$ and $V=\VN[][]$.
While this cannot be done using our Ansatz from \eqref{ansatz:eq-UN-rigorous-decomposition}-\eqref{ansatz:eq-VN-rigorous-decomposition} and earlier estimates, it can be done using the almost invariance of the Gibbs measure (see Proposition \ref{main:prop-almost-invariance} and Corollary \ref{main:cor-bourgain-bulut-probability}). Thus, the almost invariance of the Gibbs measure is used even in the proof of  local well-posedness of the discretized wave maps equation \eqref{ansatz:eq-wave-maps}, which is reminiscent of a Bourgain-Bulut argument \cite{BB14}. \\

We now comment on the proof of Proposition \ref{jacobi:prop-main}. Due to Definition \ref{ansatz:def-hhl} and Definition \ref{ansatz:def-jacobi-errors}, we need to understand the sum $\HHLN[K,y][v]+\ad\big( \SHHLN[K][y]\big)$ of the high$\times$high$\rightarrow$low-matrix $\HHLN[K,y][v]$ and the adjoint map of the simplified high$\times$high$\rightarrow$low-interaction $\SHHLN[K][y]$, which is done in three stages. In Subsection 
\ref{section:jacobi-with-shift}, we control $\HHLN[K,\bfzero][v]+\ad\big( \SHHLN[K][y]\big)$. The key ingredient in our estimate is the Jacobi identity on the \mbox{Lie algebra $\frkg$}, which leads to a crucial cancellation. In Subsection \ref{section:jacobi-with-shift}, we obtain estimates of the difference $\HHLN[K,y][v]-\HHLN[K,\bfzero][v]$ and refined estimates of the difference 
$\HHLN[K,y][v]-\HHLN[K,\bfzero][v]- \chinull[K] P^{u,v}_{<K^\delta}\CHHLNc_{K,y}$, where $\CHHLNc_{K,y}$ is the combined Cartesian high$\times$high$\rightarrow$low-matrix from \eqref{jacobi:eq-chhl-combined}. The latter estimate requires us to switch from null-coordinates back into Cartesian coordinates, which turns out to be rather delicate. In Subsection \ref{section:jacobi-cartesian}, we then examine the difference between $P^{u,v}_{<K^\delta}\CHHLNc_{K,y}$ and the matrix $\CHHLN_y ( \VN[][], \VN[][])- \coup \CNcsbd(y) \Kil$ from Proposition \ref{jacobi:prop-main}. This requires a detailed understanding of the Cartesian high$\times$high$\rightarrow$low-matrices corresponding to all terms from our Ansatz \eqref{ansatz:eq-VN-rigorous-decomposition}.

\subsection{\protect{High$\times$high$\rightarrow$low-matrices without shift}}\label{section:jacobi-without-shift}

We first analyze the high$\times$high$\rightarrow$low-matrices $\HHLN[K,\bfzero][v]$, i.e., the high$\times$high$\rightarrow$low-matrices without a shift. In the following proposition, we show that  $\HHLN[K,\bfzero][v]$ can be approximated by the adjoint map $\Ad\big( \SHHLN[K][v] \big)$, where $\SHHLN[K][v]$ is as in Definition \ref{ansatz:def-shhl}.

\begin{proposition}\label{jacobi:prop-HHL-without-shift} 
Let the post-modulation hypothesis (Hypothesis \ref{hypothesis:post}) be satisfied. Then, it holds that
\begin{equation}\label{jacobi:eq-HHL-without-shift}
\Big\| \HHLN[K,\bfzero][v] + \Ad \big( \SHHLN[K][v] \big) \Big\|_{\Cprod{\eta}{r-1}} 
\lesssim \big( K^{-\frac{1}{2}+\delta} + N^{-\delta+\vartheta} \big) \Dc^2. 
\end{equation}
\end{proposition}

We first prove the following lemma, which will later be used to simplify $\HHLN[K,\bfzero][v]$. 

\begin{lemma}\label{jacobi:lem-integration-parts}
Assume that the post-modulation hypothesis (Hypothesis \ref{hypothesis:post}) is satisfied and let 
$K,M \in \Dyadiclarge$. 
Then, it holds that
\begin{align}
&\bigg\| P^{u,v}_{<K^\delta} \Big( P_{\leq N}^x \VN[M][-] \otimes  P_{\leq N}^x \Int^v_{u\rightarrow v} P^v_{\geq K^{1-\deltap}} \VN[<K^{1-\delta}][\fcs]   + P_{\leq N}^x \IVN[M][-] \otimes  P_{\leq N}^x  P^v_{\geq K^{1-\deltap}} \VN[<K^{1-\delta}][\fcs]  \Big)
\bigg\|_{\Cprod{\eta}{r-1}} \notag \\
\lesssim \, &  K^{-\frac{1}{2}+\delta} M^{-\eta} \Dc^2. \label{jacobi:eq-low-e1}
\end{align}
\end{lemma}

\begin{proof}
We separate the proof of \eqref{jacobi:eq-low-e1} into two steps. \\ 

\emph{Step 1: Integration by parts.} In this step, we prove the estimate 
\begin{equation}
\begin{aligned}\label{jacobi:eq-low-p6}
&\bigg\| P^{u,v}_{<K^{1-\delta}} \Big( P_{\leq N}^x \VN[M][-] \otimes  P_{\leq N}^x \Int^v_{u\rightarrow v} P^v_{\geq K^{1-\deltap}} \VN[<K^{1-\delta}][\fcs]  \\
&\hspace{3ex}+  \Int^v_{u\rightarrow v} P_{\leq N}^x \VN[M][-] \otimes P_{\leq N}^x P^v_{\geq K^{1-\deltap}} \VN[<K^{1-\delta}][\fcs]\Big) \bigg\|_{\Cprod{\eta}{r-1}}  \\ 
\lesssim&\,  K^{-\frac{1}{2}+\delta} M^{-\eta} \Dc^2. 
\end{aligned}
\end{equation}
We first use Lemma \ref{prelim:lem-commutativity}, i.e., that $ P_{\leq N}^x$ and  $\Int^v_{u\rightarrow v}$ commute, which allows us to replace $P_{\leq N}^x \Int^v_{u\rightarrow v}$ in the first tensor product by  $ \Int^v_{u\rightarrow v} P_{\leq N}^x$. 
Then, using Lemma \ref{prelim:lem-integration-by-parts}, Lemma \ref{prelim:lem-Duhamel-integral}, Lemma \ref{modulation:lem-linear}, and Corollary \ref{modulation:cor-control-combined}, it follows that 
\begin{align*}
\textup{LHS of } \eqref{jacobi:eq-low-p6} 
&\lesssim  \big\| \Int^v_{u\rightarrow v} \VN[M][-] \big\|_{\Cprod{\eta}{s}} 
\big\| \Int^v_{u\rightarrow v} P^v_{\geq K^{1-\deltap}} \VN[<K^{1-\delta}][\fcs] \big\|_{\Cprod{\eta}{r-s+\eta}} \\
&\lesssim \big\| \VN[M][-] \big\|_{\Cprod{\eta}{s-1}} 
\big\| P^v_{\geq K^{1-\deltap}} \VN[<K^{1-\delta}][\fcs] \big\|_{\Cprod{\eta}{r-s+\eta-1}} 
\lesssim M^{-\eta} K^{(1-\deltap)(-s+\eta)} \Dc^2.
\end{align*}
Since $(1-\deltap)(-s+\eta)=-\frac{1}{2}+\frac{\delta}{2}+\mathcal{O}(\delta_2)$, this is acceptable. \\ 

\emph{Step 2: Replacing $\Int^v_{u\rightarrow v}  \VN[M][-]$ by $\IVN[M][-]$.} In this step, we prove the estimate 
\begin{equation}\label{jacobi:eq-low-p7}
\begin{aligned}
&\bigg\| P^{u,v}_{<K^\delta} \Big( \big( \Int^v_{u\rightarrow v} P_{\leq N}^x \VN[M][-]- P_{\leq N}^x \IVN[M][-] \big) \otimes P_{\leq N}^x P^v_{\geq K^{1-\deltap}} \VN[<K^{1-\delta}][\fcs] \Big) 
\bigg\|_{\Cprod{\eta}{r-1}} 
\lesssim  K^{-\frac{1}{2}+\delta} M^{-\eta}\Dc^2. 
\end{aligned}
\end{equation}
In order to prove \eqref{jacobi:eq-low-p7}, we first make use of the $P^{u,v}_{<K^\delta}$-operator, which yields the estimate
\begin{equation}\label{jacobi:eq-low-q1}
\textup{LHS of } \eqref{jacobi:eq-low-p7} 
\lesssim K^{\delta (r-1+2s-\eta)} \Big\| \big( \Int^v_{u\rightarrow v} P_{\leq N}^x \VN[M][-]- P_{\leq N}^x \IVN[M][-] \big) \otimes P_{\leq N}^x P^v_{\geq K^{1-\deltap}} \VN[<K^{1-\delta}][\fcs] \Big\|_{\Cprod{\eta}{-2s+\eta}}.
\end{equation}
Using our para-product estimate (Lemma \ref{prelim:lem-paraproduct}), Lemma \ref{prelim:lem-commutativity}, Lemma \ref{modulation:lem-integration}, and Corollary \ref{modulation:cor-control-combined}, it then follows that 
\begin{align*}
\eqref{jacobi:eq-low-q1}
&\lesssim   K^{\delta (r-1+2s-\eta)}  \Big\| \Int^v_{u\rightarrow v} P_{\leq N}^x \VN[M][-]- P_{\leq N}^x \IVN[M][-] \Big\|_{\Cprod{\eta}{2s}} \big\| P^v_{\geq K^{1-\deltap}} \VN[<K^{1-\delta}][\fcs] \big\|_{\Cprod{\eta}{-2s+\eta}} \\ 
&\lesssim  K^{\delta (r-1+2s-\eta)}  M^{-\eta} K^{(1-\deltap)(-2s+\eta-(r-1))} \Dc^2. 
\end{align*}
Since
\begin{equation*}
  \delta (r-1+2s-\eta)  + (1-\deltap)(-2s+\eta-(r-1)) 
  = \delta \big( \tfrac{1}{2} + \delta_1 \big) + (1-\delta) \big( - \tfrac{1}{2} - \delta_1 \big) + \mathcal{O}(\delta_2) = - \frac{1}{2} + \mathcal{O}(\delta_1), 
\end{equation*}
this is acceptable. 
By combining \eqref{jacobi:eq-low-p6} and  \eqref{jacobi:eq-low-p7}, we then obtain \eqref{jacobi:eq-low-e1}. 
\end{proof}

Equipped with Lemma \ref{jacobi:lem-integration-parts}, we now prove the main result of this subsection. 

\begin{proof}[Proof of Proposition \ref{jacobi:prop-HHL-without-shift}:] 
Using Definition \ref{ansatz:def-shhl}, using Definition \ref{ansatz:def-hhl},  and using that $\Cf_{M}^{(\Ncs)}(0)=0$, which eliminates the Killing-renormalization, it holds that 
\begin{align}
&\HHLN[K,\bfzero][v] + \Ad\big( \SHHLN[K][v] \big) \notag \\
=\, & 
\Sumlarge_{\substack{M \simeq_\delta K}} 
\chinull[K,M] P^{u,v}_{<K^{\delta}}  \bigg( \Ad \big( P_{\leq N}^x \VN[M][-] \big) \Ad \big(  P_{\leq N}^x \IVN[M][-] \big) \bigg) \label{jacobi:eq-p8} \\ 
+\, & 
\Sumlarge_{\substack{ M \lesssim_\delta K}} 
\bigg( \chinull[K,M] P^{u,v}_{<K^{\delta}} \Big( 
\Ad \big( P_{\leq N}^x P^v_{\geq K^{1-\deltap}} \VN[<K^{1-\delta}][\fcs] \big) 
\Ad \big(  P_{\leq N}^x \IVN[M][-]  \big) \Big) \label{jacobi:eq-p9} \\ 
&\hspace{10ex}+  \chinull[K] P^{u,v}_{<K^{\delta}} \Big( \Ad \big( P_{\leq N}^x \VN[M][-]  \big) 
\Ad \big(  P_{\leq N}^x \Int^v_{u\rightarrow v} P^v_{\geq K^{1-\deltap}} \VN[<K^{1-\delta}][\fcs] \big) \Big) \notag \\ 
&\hspace{10ex}- \chinull[K] P^{u,v}_{<K^{\delta}} \Big( \Ad \big( \big[ P^x_{\leq N} P^v_{\geq K^{1-\deltap}} \VN[<K^{1-\delta}][\fcs] , P^x_{\leq N} 
\IVN[M][-] \big] \big) \Big)\bigg). \notag 
\end{align}
We now estimate \eqref{jacobi:eq-p8} and \eqref{jacobi:eq-p9} separately. \\ 

\emph{Step 1: Estimate of \eqref{jacobi:eq-p8}.} 
We first write
\begin{equation}\label{jacobi:eq-q2}
\begin{aligned}
\Ad \big( P_{\leq N}^x \VN[M][-] \big) \Ad \big( P_{\leq N}^x \IVN[M][-] \big) 
= \mathcal{T} \Big(  P_{\leq N}^x \VN[M][-] \otimes   P_{\leq N}^x \IVN[M][-]  \Big),
\end{aligned}
\end{equation}
where $\mathcal{T}\colon \frkg^{\otimes 2}\rightarrow \operatorname{End}(\frkg)$ is the unique linear map satisfying
$\mathcal{T}(A\otimes B) = \Ad(A) \Ad(B)$
for all $A,B\in \frkg$. Using Definition \ref{killing:def-Wick} and $\Cf^{(\Ncs)}_M(0)=0$, we can further write
\begin{equation}\label{jacobi:eq-q2prime}
P_{\leq N}^x \VN[M][-] \otimes   P_{\leq N}^x \IVN[M][-]  
= \biglcol \, P_{\leq N}^x \VN[M][-] \otimes   P_{\leq N}^x \IVN[M][-]  \bigrcol. 
\end{equation}
By combining \eqref{jacobi:eq-q2}, \eqref{jacobi:eq-q2prime}, and Lemma \ref{killing:lem-tensor-modulated-linear}, we then obtain 
\begin{align}
&\Big\| 
\Ad \big( P_{\leq N}^x \VN[M][-] \big) \Ad \big( P_{\leq N}^x \IVN[M][-] \big) \Big\|_{\Cprod{\eta}{r-1}} 
\lesssim \, \Big( M^{r-1+\frac{1}{2}+\eta} M^{-\frac{1}{2}} + N^{-\delta+\vartheta} \Big)  \Dc^2. \label{jacobi:eq-p10} 
\end{align}
Using \eqref{prelim:eq-parameter-regularities}, \eqref{prelim:eq-parameter-regularities-less}, and $M\gtrsim K^{1-\delta}$, it follows that
\begin{equation*}
\eqref{jacobi:eq-p10} \lesssim \big(  K^{-(1-\delta) (-\frac{1}{2}+\delta_1+\delta_3)}  + N^{-\delta+\vartheta} \big) \Dc^2 
\lesssim \big( K^{-\frac{1}{2}+\delta} + N^{-\delta+\vartheta} \big) \Dc^2, 
\end{equation*}
which is an acceptable contribution to \eqref{jacobi:eq-HHL-without-shift}. \\ 

\emph{Step 2: Estimate of \eqref{jacobi:eq-p9}.}
Due to the smoothness of $\chinull$, Lemma \ref{modulation:lem-linear}, and Corollary \ref{modulation:cor-control-combined}, it holds that 
\begin{align*}
&\Big\| (\chinull[K,M]-\chinull[K]) P^{u,v}_{<K^{\delta}} \Big( 
\Ad \big( P_{\leq N}^x P^v_{\geq K^{1-\deltap}} \VN[<K^{1-\delta}][\fcs] \big) 
\Ad \big(  P_{\leq N}^x \IVN[M][-]  \big) \Big) \Big\|_{\Cprod{\eta}{r-1}} \\ 
\lesssim&\,  K^{-100} \big\|  \VN[<K^{1-\delta}][\fcs] \big\|_{\Cprod{\eta}{r-1}} 
\big\| \IVN[M][-] \big\|_{\Cprod{\eta}{s}} \lesssim K^{-100} M^{-\eta} \Dc^2,
\end{align*}
which is a more than acceptable contribution to \eqref{jacobi:eq-HHL-without-shift}. We may therefore replace the $\chinull[K,M]$-factor in the first summand of \eqref{jacobi:eq-p9} with a $\chinull[K]$-factor. Since all three summands in \eqref{jacobi:eq-p9} then share a $\chinull[K]$-factor, it can be factored out, and it then suffices to treat
\begin{equation}\label{jacobi:eq-q}
\begin{aligned}
P^{u,v}_{<K^{\delta}} &\Sumlarge_{\substack{M \lesssim_\delta K}} 
\bigg(  
\Ad \big( P_{\leq N}^x P^v_{\geq K^{1-\deltap}} \VN[<K^{1-\delta}][\fcs] \big) \Ad \big(  P_{\leq N}^x \IVN[M][-]  \big)  \\
&\hspace{10ex}  + \Ad \big( P_{\leq N}^x \VN[M][-]  \big) 
\Ad \big(  P_{\leq N}^x \Int^v_{u\rightarrow v} P^v_{\geq K^{1-\deltap}} \VN[<K^{1-\delta}][\fcs] \big)  \\ 
&\hspace{10ex}-  \Ad \Big( \big[ P^x_{\leq N} P^v_{\geq K^{1-\deltap}} \VN[<K^{1-\delta}][\fcs] , P^x_{\leq N} 
\IVN[M][-] \big] \Big) \bigg). 
\end{aligned}
\end{equation}
Using the Jacobi identity for the Lie algebra $\frkg$, which has been stated in terms of the adjoint map in \eqref{prelim:eq-jacobi}, it holds that 
\begin{equation}\label{jacobi:eq-p11} 
\begin{aligned}
&- \Ad \Big( \big[ P^x_{\leq N} P^v_{\geq K^{1-\deltap}} \VN[<K^{1-\delta}][\fcs] , P^x_{\leq N} 
\IVN[M][-] \big] \Big) \\ 
= \, & \Ad \Big( \big[ P^x_{\leq N} 
\IVN[M][-],  P^x_{\leq N} P^v_{\geq K^{1-\deltap}} \VN[<K^{1-\delta}][\fcs] \big] \Big) \\ 
=\, & 
\Ad \Big( P_{\leq N}^x \IVN[M][-] \Big) 
\Ad \Big( P^x_{\leq N} P^v_{\geq K^{1-\deltap}} \VN[<K^{1-\delta}][\fcs]\Big) \\ 
-\, &
\Ad \Big( P^x_{\leq N} P^v_{\geq K^{1-\deltap}} \VN[<K^{1-\delta}][\fcs]\Big) 
\Ad \Big( P_{\leq N}^x \IVN[M][-] \Big).  
\end{aligned}
\end{equation}
Since the last term in \eqref{jacobi:eq-p11} cancels the first summand in \eqref{jacobi:eq-q},  the $\Cprod{\eta}{r-1}$-norm of the dyadic summands in \eqref{jacobi:eq-q} is bounded by 
\begin{equation}\label{jacobi:eq-p12} 
\begin{aligned}
& \bigg\| \Ad \big( P_{\leq N}^x \VN[M][-]  \big) 
\Ad \big(  P_{\leq N}^x \Int^v_{u\rightarrow v} P^v_{\geq K^{1-\deltap}} \VN[<K^{1-\delta}][\fcs] \big) \\ &\hspace{2ex} + \Ad \big( P_{\leq N}^x \IVN[M][-] \big) 
\Ad \big( P^x_{\leq N} P^v_{\geq K^{1-\deltap}} \VN[<K^{1-\delta}][\fcs]\big) \bigg\|_{\Cprod{\eta}{r-1}}. 
\end{aligned}
\end{equation}
Using Lemma \ref{jacobi:lem-integration-parts}, it follows that
\begin{equation*}
\eqref{jacobi:eq-p12} \lesssim  K^{-\frac{1}{2}+\delta} M^{-\eta} \Dc^2,
\end{equation*}
which yields an acceptable contribution to \eqref{jacobi:eq-HHL-without-shift}. 
\end{proof} 

\subsection{\protect{High$\times$high$\rightarrow$low-matrices with shift}}\label{section:jacobi-with-shift}

We now turn to the difference between $\HHLN[K,y][v]$ and $\HHLN[K,\bfzero][v]$. For notational purposes, it is convenient to introduce the combined Cartesian high$\times$high$\rightarrow$low-operators 
\begin{equation}\label{jacobi:eq-chhl-combined}
\begin{aligned}
 \CHHLNc_{K,y}
 &:= \CHHLN\big( \VN[][-], \VN[][-] \big) 
 + \CHHLN \big( P^v_{\geq K^{1-\deltap}} \VN[\leq K^{1-\delta}][\fcs], \VN[][-] \big) \\
 &+ \CHHLN \big( P^v_{\geq K^{1-\deltap}} \VN[\leq K^{1-\delta}][\fcs], \VN[][-] \big) 
 - \coup \CNcsbd(y) \Kil. 
\end{aligned}
\end{equation}

\begin{proposition}\label{jacobi:prop-shift}
Let the post-modulation hypothesis (Hypothesis \ref{hypothesis:post}) be satisfied and let $K\in  \Dyadiclarge$. Then, it holds that
\begin{equation}\label{jacobi:eq-shift-general}
 \Big\| \HHLN[K,y][v] - \HHLN[K,\bfzero][v]\Big\|_{\Cprod{\eta}{r-1}} \lesssim 
 \bigg( K^{-\frac{1}{2}+\delta} + N^{-\delta+\vartheta} + \tfrac{\Nd}{N}\,  \Nd^{-\delta_1+2\delta_3} \bigg) \langle Ny \rangle \Dc^2. 
\end{equation}
Furthermore, if $K\geq N^{1-2\delta_1}$, then it also holds that 
\begin{equation}\label{jacobi:eq-shift-boundary}
\begin{aligned}
&\Big\| \HHLN[K,y][v] - \HHLN[K,\bfzero][v]
- \chinull[K] P^{u,v}_{<K^\delta} \CHHLNc_{K,y} \Big\|_{\Cprod{\eta}{r-1}} 
\lesssim \Big( N^{-\frac{1}{2}+\delta} + N^{-2\delta_1 + 4\delta_2} \Big) \langle N y \rangle \Dc^2.
\end{aligned}
\end{equation}
\end{proposition}

\begin{remark} In the regime $\Nd \sim N$, the right-hand sides of \eqref{jacobi:eq-shift-general} and \eqref{jacobi:eq-shift-boundary} contain terms which, up to insignificant factors, are of the form $N^{-\delta_1}$ and $N^{-2\delta_1}$, respectively. While the two terms only differ by a factor of $N^{-\delta_1}$, this difference is significant. The reason is that, in order to control the Jacobi error, we need to absorb a loss of the form $K^{\delta_1+\mathcal{O}(\delta_3)}$. If $K\sim \Nd \sim N$, this cannot be done using $N^{-\delta_1}$.
\end{remark}

The proof of Proposition \ref{jacobi:prop-shift} occupies the remainder of this subsection. We first prove \eqref{jacobi:eq-shift-general}, which is contained in Lemma \ref{jacobi:lem-general} and whose proof is rather simple. We then prove \eqref{jacobi:eq-shift-boundary}, which is more difficult and whose proof is distributed over several lemmas.

\subsubsection{\protect{General estimate}} 

We first prove \eqref{jacobi:eq-shift-general}, which is a general estimate for all $K\in \Dyadiclarge$.

\begin{lemma}\label{jacobi:lem-general}
Let the post-modulation hypothesis (Hypothesis \ref{hypothesis:post}) be satisfied
and let $K\in \Dyadiclarge$. Then, it holds that 
\begin{equation}\label{jacobi:eq-general}
\Big\| \HHLN[K,y][v] - \HHLN[K,\bfzero][v]\Big\|_{\Cprod{\eta}{r-1}} \lesssim \bigg( K^{-\frac{1}{2}+\delta} + N^{-\delta+\vartheta} + \tfrac{\Nd}{N} \Nd^{-\delta_1+2\delta_3} \bigg) \langle Ny \rangle  \Dc^2. 
\end{equation}
\end{lemma}

\begin{proof}
We first recall from Definition \ref{ansatz:def-hhl} that, for all $z\in \R$, 
\begin{align}
 \HHLN[K,z][v]  &:= 
\Sumlarge_{\substack{M \leq \Nd\colon \\ M \simeq_\delta K}} 
\Big( \HHLN[K,M,z][v,(-)\times (-)] +  \HHLN[K,M,z][v,\textup{kil}] \Big) \label{jacobi:eq-general-p1}  \\
&+ 
\Sumlarge_{\substack{M \leq \Nd\colon \\ M \lesssim_\delta K}} \Big( 
\HHLN[K,M,z][v,(\fcs)\times (-)] + 
\HHLN[K,M,z][v,(-)\times (\fcs)] \Big). \label{jacobi:eq-general-p2} 
\end{align}
We now control the contributions of \eqref{jacobi:eq-general-p1} and \eqref{jacobi:eq-general-p2} to \eqref{jacobi:eq-general} separately. \\

\emph{Contribution of \eqref{jacobi:eq-general-p1}:} Using Definition \ref{ansatz:def-hhl}, it follows that
\begin{align}
&\Big(  \HHLN[K,M,y][v,(-)\times (-)] +  \HHLN[K,M,y][v,\textup{kil}] \Big) - \Big(  \HHLN[K,M,\bfzero][v,(-)\times (-)] +  \HHLN[K,M,\bfzero][v,\textup{kil}] \Big) \notag \\ 
=&\,\chinull[K,M] P^{u,v}_{<K^\delta}   \bigg( 
\Ad \Big( P_{\leq N}^x \VN[M][-] \Big) 
\Ad \Big( \Theta^x_y P_{\leq N}^x \IVN[M][-] \Big)  - \coup \Cf^{(\Ncs)}_M(y) \Kil \bigg) \label{jacobi:eq-general-p3} \\ 
-&\,\chinull[K,M] P^{u,v}_{<K^\delta}   \bigg( 
\Ad \Big( P_{\leq N}^x \VN[M][-] \Big) 
\Ad \Big(  P_{\leq N}^x \IVN[M][-] \Big)  - \coup \Cf^{(\Ncs)}_M(0) \Kil \bigg) \label{jacobi:eq-general-p4} \\ 
+&\, \coup \big( \chinull[K,M] - \chinull[K] \big) \big( \Cf^{(\Ncs)}_M(y)-\Cf^{(\Ncs)}_M(0)\big) \label{jacobi:eq-general-q}. 
\end{align}
Using Proposition \ref{killing:prop-quadratic} (and symmetry of our estimates in the $u$ and $v$-variables), we obtain that
\begin{align*}
\Big\| \eqref{jacobi:eq-general-p3} \Big\|_{\Cprod{\eta}{r-1}} 
+ \Big\| \eqref{jacobi:eq-general-p4} \Big\|_{\Cprod{\eta}{r-1}} 
\lesssim \Big( M^{r-\frac{1}{2}+\eta} M^{-\frac{1}{2}} + N^{-\delta+\vartheta} \langle N y \rangle \Big) \Dc^2.
\end{align*}
Since $M\simeq_\delta K$, it follows that
\begin{equation*}
M^{r-\frac{1}{2}+\eta} M^{-\frac{1}{2}} + N^{-\delta+\vartheta} \langle N y \rangle \lesssim K^{-\frac{1}{2}+\delta} + N^{-\delta+\vartheta} \langle N y\rangle ,
\end{equation*}
which yields an acceptable contribution to \eqref{jacobi:eq-general}. Using the smoothness of $\chinull$ and Lemma \ref{ansatz:lem-renormalization}, we also have that 
\begin{equation*}
\big\| \eqref{jacobi:eq-general-q}\big\|_{\Cprod{\eta}{r-1}} 
\lesssim K^{-100} \Dc^2,
\end{equation*}
which yields a more than acceptable contribution to \eqref{jacobi:eq-general}.\\

\emph{Contribution of \eqref{jacobi:eq-general-p2}:} Since the argument for the $(-)$$\times$$(\fcs)$-interaction is similar, we only control the $(\fcs)$$\times$$(-)$-interaction.  
To this end, we first recall that 
\begin{align}
&\hspace{1ex} \HHLN[K,M,y][v,(\fcs)\times (-)] - \HHLN[K,M,\bfzero][v,(\fcs)\times (-)]   \notag \\ 
=&\, \chinull[K,M] P^{u,v}_{<K^{\delta}} \bigg( 
\Ad \Big( P_{\leq N}^x  P^v_{\geq K^{1-\deltap}} \VN[< K^{1-\delta}][\fcs] \Big) 
\Ad \Big( \big( \Theta^x_y -1 \big) P_{\leq N}^x \IVN[M][-] \Big) \bigg) \label{jacobi:eq-general-p5}. 
\end{align}
Using the high$\times$high-paraproduct estimate (Lemma \ref{prelim:lem-paraproduct}) and Corollary \ref{modulation:cor-control-combined}, it follows that
\begin{align*}
\Big\| \eqref{jacobi:eq-general-p5} \Big\|_{\Cprod{\eta}{r-1}} 
&\lesssim \Big\| \VN[< K^{1-\delta}][\fcs] \Big\|_{\Cprod{s}{r-1}} \Big\| \big( \Theta^x_y -1 \big)  \IVN[M][-] \Big\|_{\Cprod{s}{1-r+\eta}} \lesssim \Dc  \Big\| \big( \Theta^x_y -1 \big)  \IVN[M][-] \Big\|_{\Cprod{s}{1-r+\eta}}. 
\end{align*}
Due to Lemma \ref{ansatz:lem-frequency-support}, 
$\IVN[M][-]$ is supported on $u$ and $v$-frequencies $\lesssim M$. Using Lemma \ref{modulation:lem-linear}, we obtain that
\begin{equation*}
    \Dc  \Big\| \big( \Theta^x_y -1 \big)  \IVN[M][-] \Big\|_{\Cprod{s}{1-r+\eta}} 
    \lesssim M |y| \Dc \Big\| \IVN[M][-] \Big\|_{\Cprod{s}{1-r+\eta}} 
    \lesssim M^{\frac{3}{2}-r+2\eta} |y| \Dc^2. 
\end{equation*}
Since $M\lesssim \Nd$ and $\frac{3}{2}-r+2\eta=1-\delta_1+2\delta_3$, we obtain that
\begin{equation*}
 M^{\frac{3}{2}-r+2\eta} |y| \Dc^2
 \lesssim \Nd^{1-\delta_1+2\delta_3} N^{-1} \langle N y \rangle \Dc^2 = \tfrac{\Nd}{N} \, \Nd^{-\delta_1+2\delta_3} \langle Ny \rangle \Dc^2, 
\end{equation*}
which is an acceptable contribution to \eqref{jacobi:eq-general}.
\end{proof}

\subsubsection{\protect{Refined estimate for $N^{1-2\delta_1} \leq K\lesssim  N$}} 
In the following lemma, we show how to replace the $P_{\leq N}^x$ and $\Int^v_{u\rightarrow v}$-operators in the high$\times$high$\rightarrow$low-matrices. 

\begin{lemma}[Replacing $P_{\leq N}^x$ and $\Int^v_{u\rightarrow v}$]\label{jacobi:lem-replacement} 
Assume that the post-modulation hypothesis (Hypothesis \ref{hypothesis:post}) is satisfied. Furthermore, let $K,M\in \Dyadiclarge$ satisfy $K\geq N^{1-2\delta_1}$ and let 
\begin{equation*}
\VN[][\ast_1],\VN[][\ast_2] \in \Big\{ \VN[M][-], P^v_{\geq K^{1-\deltap}} \VN[<K^{1-\delta}][\fcs] \Big\},
\end{equation*}
Then, it holds that 
\begin{equation}\label{jacobi:eq-replacement}
\begin{aligned}
&\Big\| P^{u,v}_{<K^\delta} \Big( 
P_{\leq N}^x \VN[][\ast_1] \otimes (\Theta^x_y -1) \Int^v_{u\rightarrow v} P_{\leq N}^x \VN[][\ast_2] \Big) 
- P^{u,v}_{<K^\delta} \Big( 
\CHHLN_y \big( \VN[][\ast_1] , \VN[][\ast_2] \big) \Big) 
\Big\|_{\Cprod{\eta}{r-1}} \\ 
\lesssim &\, \Big( N^{-\frac{1}{2}+\delta} +  N^{-2\delta_1 + 4 \delta_2}\, N|y|  \Big) \Dc^2.
\end{aligned}
\end{equation}
\end{lemma}

\begin{proof} In order to prove \eqref{jacobi:eq-replacement}, it suffices to prove the two separate estimates
\begin{equation}\label{jacobi:eq-replacement-e1}
\begin{aligned}
&\Big\| P^{u,v}_{<K^\delta} \Big( 
P_{\leq N}^x \VN[][\ast_1] \otimes (\Theta^x_y -1) \Int^v_{u\rightarrow v} P_{\leq N}^x \VN[][\ast_2] \Big) \\
&\hspace{3ex} 
- P^{u,v}_{<K^\delta} \Big( 
\Pbd \VN[][\ast_1] \otimes (\Theta^x_y -1) \Int^v_{u\rightarrow v} \Pbd \VN[][\ast_2] \Big) 
\Big\|_{\Cprod{\eta}{r-1}} 
\lesssim  N^{-2\delta_1 + 4 \delta_2}\, N|y| \Dc^2 
\end{aligned}
\end{equation}
and 
\begin{equation}\label{jacobi:eq-replacement-e2} 
\begin{aligned}
 &\,\Big\| P^{u,v}_{<K^\delta} \Big( 
\Pbd \VN[][\ast_1] \otimes (\Theta^x_y -1) \Int^v_{u\rightarrow v} \Pbd \VN[][\ast_2] \Big) -  P^{u,v}_{<K^\delta} \Big( 
\CHHLN_y \big( \VN[][\ast_1] , \VN[][\ast_2] \big) \Big) \Big\|_{\Cprod{\eta}{r-1}} \\
\lesssim&\, N^{-\frac{1}{2}+\delta} \Dc^2. 
\end{aligned}
\end{equation}
\emph{Proof of \eqref{jacobi:eq-replacement-e1}:}  Using the definition of $\Pbd$, we obtain that 
\begin{align}
&P^{u,v}_{<K^\delta} \Big( 
P_{\leq N}^x \VN[][\ast_1] \otimes (\Theta^x_y -1) \Int^v_{u\rightarrow v} P_{\leq N}^x \VN[][\ast_2] \Big) \notag \\
&\hspace{3ex} 
- P^{u,v}_{<K^\delta} \Big( 
\Pbd \VN[][\ast_1] \otimes (\Theta^x_y -1) \Int^v_{u\rightarrow v} \Pbd \VN[][\ast_2] \Big) \notag \\ 
=&\, \sum_{\substack{L_1,L_2 \in \dyadic\colon \\ \min(L_1,L_2)<N^{1-2\delta_1}}} 
P^{u,v}_{<K^\delta} \Big( 
P_{L_1}^x \VN[][\ast_1] \otimes (\Theta^x_y -1) \Int^v_{u\rightarrow v} P_{L_2}^x \VN[][\ast_2] \Big). \label{jacobi:eq-replace-PNX-p1}
\end{align}
Due to the $P^{u,v}_{<K^\delta}$-operator outside the tensor products in \eqref{jacobi:eq-replace-PNX-p1}, the tensor products are localized to frequencies $\lesssim K^\delta \lesssim N^\delta$ in the $x$-variable. As a result, the condition $\min(L_1,L_2)<N^{1-2\delta_1}$ implies that $\max(L_1,L_2)\lesssim N^{1-2\delta_1}$. Using the high$\times$high-paraproduct estimate (Lemma \ref{prelim:lem-paraproduct}), the commutativity of  $\Theta^x_y$ and $\Int^v_{u\rightarrow v}$, and Lemma \ref{prelim:lem-Duhamel-integral}, we now estimate
\begin{align}
&\Big\| P^{u,v}_{<K^\delta} \Big( 
P_{L_1}^x \VN[][\ast_1] \otimes (\Theta^x_y -1) \Int^v_{u\rightarrow v} P_{L_2}^x \VN[][\ast_2] \Big) 
\Big\|_{\Cprod{\eta}{r-1}} \notag \\
\lesssim&\, \Big\| P_{L_1}^x \VN[][\ast_1] \Big\|_{\Cprod{\eta}{-s+\eta}} 
\Big\| \Int^v_{u\rightarrow v} (\Theta^x_y -1)  P_{L_2}^x \VN[][\ast_2]\Big\|_{\Cprod{\eta}{s}}  \notag \\
\lesssim&\, N^{1-2s+\eta} \Big\| \VN[][\ast_1] \Big\|_{\Cprod{s}{s-1}}  
\Big\| (\Theta^x_y -1)  P_{L_2}^x \VN[][\ast_2]\Big\|_{\Cprod{\eta}{s-1}}.\label{jacobi:eq-replace-PNX-p2}
\end{align}
The first factor in \eqref{jacobi:eq-replace-PNX-p2} is given by $N^{2\delta_2+\delta_3}$ and the second factor in \eqref{jacobi:eq-replace-PNX-p2} is controlled by $\Dc$. Using the fundamental theorem of calculus and using that $\max(L_1,L_2)\lesssim N^{1-2\delta_1}$, the third factor in \eqref{jacobi:eq-replace-PNX-p2} can be estimated by
\begin{equation*}
\Big\| (\Theta^x_y -1)  P_{L_2}^x \VN[][\ast_2]\Big\|_{\Cprod{\eta}{s-1}}
\lesssim L_2 |y| \Big\|  \VN[][\ast_2]\Big\|_{\Cprod{s}{s-1}}
\lesssim N^{1-2\delta_1} |y| \Dc. 
\end{equation*}
By inserting this back into \eqref{jacobi:eq-replace-PNX-p2}, we obtain the desired estimate \eqref{jacobi:eq-replacement-e1}. \\ 

\emph{Proof of \eqref{jacobi:eq-replacement-e2}:} 
Using Definition \ref{jacobi:def-chhl}, $P^{u,v}_{<K^\delta}=P^{u,v}_{<K^\delta} P^x_{\leq N^{\deltap}}$,  and our product estimate (Corollary \ref{prelim:cor-product}), it holds that 
\begin{align*}
 \textup{LHS of } \eqref{jacobi:eq-replacement-e2}\leq &\, \Big\| \Pbd \VN[][\ast_1] \otimes ( \Theta^x_y-1)(\Int^v_{u\rightarrow v} - P^x_{>1} \Int^x_{0\rightarrow x} ) \Pbd \VN[][\ast_2] \Big\|_{\Cprod{\eta}{r-1}} \\
 \lesssim&\, 
 \Big\| \Pbd \VN[][\ast_1] \Big\|_{\Cprod{\eta}{r-1}} 
 \Big\|  (\Int^v_{u\rightarrow v} - P^x_{>1} \Int^x_{0\rightarrow x} ) \Pbd \VN[][\ast_2] \Big\|_{\Cprod{\eta}{1-r^\prime}}.
\end{align*}
By using either Lemma \ref{modulation:lem-linear} or Corollary \ref{modulation:cor-control-combined}, it follows that
\begin{equation*}
    \Big\| \Pbd \VN[][\ast_1] \Big\|_{\Cprod{\eta}{r-1}} 
    \lesssim N^{r-s} \Big\| \Pbd \VN[][\ast_1] \Big\|_{\Cprod{s}{s-1}} 
    \lesssim N^{\delta_1+\delta_2} \Dc. 
\end{equation*}
We now intend to use Lemma \ref{prelim:lem-Cartesian-integral}, for which we need to verify that the $v$-frequency of $\Pbd \VN[][\ast_2]$ is much larger than its $u$-frequency. For $\VN[][\ast_2]=\VN[M][-]$, this is clear since $\VN[M][-]$ is supported on $u$-frequencies $\lesssim M^{1-\delta}$ and $v$-frequencies $\sim M$. For $\VN[][\ast_2]=P^v_{\geq K^{1-\deltap}} \VN[<K^{1-\delta}][\fcs]$, we first note that $\VN[][\ast_2]$ is supported on $u$-frequencies $\lesssim K^{1-\delta}\lesssim N^{1-\delta}$. Due to the definition of $\Pbd$, this implies that $\Pbd P^v_{\geq K^{1-\deltap}} \VN[<K^{1-\delta}][\fcs]$ can only be non-zero if $ \VN[<K^{1-\delta}][\fcs]$ enters at $v$-frequencies $\gtrsim N^{1-2\delta_1}$. Since $N^{1-2\delta_1}\gg N^{1-\delta}$, this implies that $\Pbd P^v_{\geq K^{1-\deltap}} \VN[<K^{1-\delta}][\fcs]$ is also supported on $v$-frequencies much larger than its $u$-frequencies. Thus, the conditions in Lemma \ref{prelim:lem-Cartesian-integral} are satisfied, and we then obtain that
 \begin{equation*}
     \Big\| (\Int^v_{u\rightarrow v} - P^x_{>1} \Int^x_{0\rightarrow x} ) \Pbd \VN[][\ast_2] \Big\|_{\Cprod{\eta}{1-r^\prime}} 
     \lesssim \Big\| \Pbd \VN[][\ast_2] \Big\|_{\Cprod{s}{-1+2\delta_1}}. 
 \end{equation*}
 Since $\Pbd \VN[][\ast_2]$ is supported on $v$-frequencies $\gtrsim N^{1-2\delta_1}$, it then follows that
 \begin{equation*}
    \Big\| \Pbd \VN[][\ast_2] \Big\|_{\Cprod{s}{-1+2\delta_1}}
    \lesssim \big( N^{1-2\delta_1} \big)^{-1+2\delta_1-(s-1)}
    \big\| \VN[][\ast_2] \big\|_{\Cprod{s}{s-1}} 
    \lesssim N^{(1-2\delta_1)(-\scalebox{0.8}{$\tfrac{1}{2}$}+2\delta_1+\delta_2)} \Dc. 
 \end{equation*}
 Since 
 \begin{equation*}
    \delta_1+\delta_2 
    + (1-2\delta_1)\Big(-\tfrac{1}{2}+2\delta_1+\delta_2\Big) 
    = -\frac{1}{2}+ \mathcal{O}(\delta_1),
 \end{equation*}
 this yields an acceptable contribution.
\end{proof}

In the next lemma, we examine the high$\times$high$\rightarrow$low-matrices for the $(-)$$\times$$(-)$-interaction. 

\begin{lemma}[\protect{The $(-)$$\times$$(-)$-interaction}]\label{jacobi:lemma-mm}
Assume that the post-modulation hypothesis (Hypothesis \ref{hypothesis:post}) is satisfied. Furthermore, let $K\in \Dyadiclarge$ satisfy $K\geq N^{1-2\delta_1}$. Then, it holds that
\begin{equation}\label{jacobi:eq-mm}
\begin{aligned}
&\Big\| \Sumlarge_{\substack{M\leq \Nd \colon \\  M\simeq_\delta K}} 
\Big( \HHLN[K,M,y][v,(-)\times (-)] - \HHLN[K,M,\bfzero][v,(-)\times (-)] \Big) - \chinull[K] P^{u,v}_{<K^\delta} 
\Big(\CHHLN_y \big( \VN[][-], \VN[][-] \big) \Big) \Big\|_{\Cprod{\eta}{r-1}} \\ 
\lesssim&\, \Big( N^{-\frac{1}{2}+\delta} + N^{-2\delta_1+4\delta_2} \, N |y| \Big) \Dc^2. 
\end{aligned}
\end{equation}
\end{lemma}

\begin{proof}
Since $K\geq N^{1-2\delta_1}$, the conditions $M\leq \Nd$ and $M\simeq_\delta K$ are equivalent to $K^{1-\delta}\leq M\leq \Nd$. As a result, it follows that 
\begin{equation}\label{jacobi:eq-mm-p1} 
\begin{aligned}
&\Sumlarge_{\substack{M\leq \Nd \colon \\  M\simeq_\delta K}} 
\Big( \HHLN[K,M,y][v,(-)\times (-)] - \HHLN[K,M,\bfzero][v,(-)\times (-)] \Big) \\ 
=&\, 
\Sumlarge_{\substack{M\colon \\ K^{1-\delta}\leq M\leq \Nd}} \chinull[K,M] P^{u,v}_{<K^\delta} 
\Big( \Ad \big( P_{\leq N}^x \VN[M][-] \big)
\Ad \big( (\Theta^x_y -1) P_{\leq N}^x \IVN[M][-] \big) \Big). 
\end{aligned}
\end{equation}
We now make several approximations of \eqref{jacobi:eq-mm-p1}, which will eventually lead us to the desired estimate \eqref{jacobi:eq-mm}. First, using the smoothness of $\chinull$, we can easily approximate the $\chinull[K,M]$-factor in \eqref{jacobi:eq-mm-p1} using $\chinull[K]$, which can then be pulled out of the dyadic sum. 
Next, using Lemma \ref{killing:lem-tensor-modulated-linear}, we obtain that
\begin{align}
&\Big\| \Sumlarge_{\substack{M\colon \\ K^{1-\delta}\leq M\leq \Nd}} P^{u,v}_{<K^\delta} 
\Big( \Ad \big( P_{\leq N}^x \VN[M][-] \big)
\Ad \big( (\Theta^x_y -1) P_{\leq N}^x \IVN[M][-] \big) \Big) \notag \\
&\hspace{3ex} - \Sumlarge_{\substack{L,M \colon  \\ K^{1-\delta}\leq L,M\leq \Nd}} P^{u,v}_{<K^\delta} 
\Big( \Ad \big( P_{\leq N}^x \VN[L][-] \big)
\Ad \big( (\Theta^x_y -1) P_{\leq N}^x \IVN[M][-] \big) \Big) 
\Big\|_{\Cprod{\eta}{r-1}} \allowdisplaybreaks[3] \label{jacobi:eq-mm-p2} \\ 
\lesssim&\, \Sumlarge_{\substack{L,M \colon  \\ K^{1-\delta}\leq L,M\leq \Nd,\\ L\neq M}} 
\max(L,M)^{r-\frac{1}{2}+\eta} M^{-\frac{1}{2}} \Dc^2 
\lesssim N^{\delta_1+\delta_3} K^{-\frac{1-\delta}{2}} \Dc^2. \notag 
\end{align}
Since $K\gtrsim N^{1-2\delta_1}$, the error term is acceptable. Second, using the high$\times$high-paraproduct estimate and Lemma \ref{modulation:lem-integration}, it follows that 
\begin{align}
&\Big\| \Sumlarge_{\substack{L,M \colon  \\ K^{1-\delta}\leq L,M\leq \Nd}} P^{u,v}_{<K^\delta} 
\Big( \Ad \big( P_{\leq N}^x \VN[L][-] \big)
\Ad \big( (\Theta^x_y -1) P_{\leq N}^x \IVN[M][-] \big) \Big) 
\Big\|_{\Cprod{\eta}{r-1}}  
\notag \\
&\hspace{3ex} - \Sumlarge_{\substack{L,M \colon  \\ K^{1-\delta}\leq L,M\leq \Nd}} P^{u,v}_{<K^\delta} 
\Big( \Ad \big( P_{\leq N}^x \VN[L][-] \big)
\Ad \big( (\Theta^x_y -1) P_{\leq N}^x \Int^v_{u\rightarrow v} \VN[M][-] \big) \Big) 
\Big\|_{\Cprod{\eta}{r-1}} 
\allowdisplaybreaks[3] \label{jacobi:eq-mm-p3}\\ 
\lesssim&\, \Sumlarge_{\substack{L,M \colon  \\ K^{1-\delta}\leq L,M\leq \Nd}} 
\Big\| \VN[L][-] \Big\|_{\Cprod{s}{-2s+\eta}} 
\Big\| \IVN[M][-] - \Int^v_{u\rightarrow v} \VN[M][-] \Big\|_{\Cprod{s}{2s}} 
\notag \\ 
\lesssim&\, \Sumlarge_{\substack{L,M \colon  \\ K^{1-\delta}\leq L,M\leq \Nd}} 
L^{-2s+\frac{1}{2}+2\eta} M^{-\eta} \Dc^2 
\lesssim K^{(1-\delta) \big(-\tfrac{1}{2}+2\delta_2 +2 \delta_3\big)} \Dc^2. 
\notag 
\end{align}
Since $K\gtrsim N^{1-2\delta_1}$, this is acceptable. Since $P_{\leq N}^x$ and $\Int^v_{u\rightarrow v}$ commute, we can also replace $P_{\leq N}^x \Int^v_{u\rightarrow v} \VN[M][-]$ by $ \Int^v_{u\rightarrow v} P_{\leq N}^x \VN[M][-]$. Third, using Lemma \ref{jacobi:lem-replacement}, we can replace $P_{\leq N}^x$ by $\Pbd$ and $\Int^v_{u\rightarrow v}$ by $P^x_{>1} \Int^x_{0\rightarrow x}$, i.e., we can estimate 
\begin{equation}\label{jacobi:eq-mm-p4}
\begin{aligned}
&\Big\| \Sumlarge_{\substack{L,M \colon  \\ K^{1-\delta}\leq L,M\leq \Nd}} P^{u,v}_{<K^\delta} 
\Big( \Ad \big( P_{\leq N}^x \VN[L][-] \big)
\Ad \big( (\Theta^x_y -1) \Int^v_{u\rightarrow v} P_{\leq N}^x  \VN[M][-] \big) \Big) \\
&\hspace{3ex} - \Sumlarge_{\substack{L,M \colon  \\ K^{1-\delta}\leq L,M\leq \Nd}} P^{u,v}_{<K^\delta} 
\Big( \Ad \big( \Pbd \VN[L][-] \big)
\Ad \big( (\Theta^x_y -1)  P^x_{>1} \Int^x_{0\rightarrow x} \Pbd \VN[M][-] \big) \Big) 
\Big\|_{\Cprod{\eta}{r-1}} \\ 
\lesssim&\, \big( N^{-\frac{1}{2}+\delta}+ N^{-2\delta_1+4\delta_2}\, N|y| \big) \Dc^2, 
\end{aligned}
\end{equation}
which is acceptable.  Due to Lemma \ref{ansatz:lem-frequency-support}
and the definition of $\Pbd$, it holds that
 $\Pbd \VN[L][-]=0$ and $\Pbd \VN[M][-]=0$ for all $L,M<K^{1-\delta}$. 
 Together with Definition \ref{jacobi:def-chhl}, we then obtain
 
\begin{equation}\label{jacobi:eq-mm-p6}
\begin{aligned}
  &\Sumlarge_{\substack{L,M \colon  \\ K^{1-\delta}\leq L,M\leq \Nd}} P^{u,v}_{<K^\delta} 
\Big( \Ad \big( \Pbd \VN[L][-] \big)
\Ad \big( (\Theta^x_y -1)  P^x_{>1} \Int^x_{0\rightarrow x} \Pbd \VN[M][-] \big) \Big) \allowdisplaybreaks[4] \\
=&\,   \Sumlarge_{\substack{L,M \colon \\ L,M\leq \Nd}} P^{u,v}_{<K^\delta} 
\Big( \Ad \big( \Pbd \VN[L][-] \big)
\Ad \big( (\Theta^x_y -1) P^x_{>1} \Int^x_{0\rightarrow x}  \Pbd \VN[M][-] \big) \Big) \allowdisplaybreaks[4] \\
=&\, P^{u,v}_{<K^\delta} 
\Big( \Ad \big( \Pbd \VN[][-] \big)  \Ad \big( (\Theta^x_y -1) P^x_{>1} \Int^x_{0\rightarrow x} \Pbd \VN[][-] \big) \Big) \allowdisplaybreaks[4] \\
=&\, P^{u,v}_{<K^\delta}  \Big( \CHHLN_y \big( \VN[][-], \VN[][-]\big) \Big).
\end{aligned}
\end{equation}
By combining the three estimates \eqref{jacobi:eq-mm-p2}, \eqref{jacobi:eq-mm-p3}, and \eqref{jacobi:eq-mm-p4} and the identity \eqref{jacobi:eq-mm-p6}, we obtain the desired estimate \eqref{jacobi:eq-mm}. 
\end{proof}

\begin{lemma}[\protect{The $(\fsc)$$\times$$(-)$ and $(-)$$\times$$(\fcs)$-interactions}]\label{jacobi:lem-fcsm}
Assume that the post-modulation hypothesis (Hypothesis \ref{hypothesis:post}) is satisfied. Furthermore, let $K\in \Dyadiclarge$ satisfy $K\geq N^{1-2\delta_1}$. Then, it holds that
\begin{align}
&\Big\| \Sumlarge_{\substack{M\leq \Nd \colon \\ M \lesssim_\delta K}} \Big( \HHLN[K,M,y][v,(\fcs)\times (-)] - \HHLN[K,M,\bfzero][v,(\fcs)\times (-)] \Big) 
- \chinull[K] P^{u,v}_{<K^\delta} \Big( \CHHLN_y \big( P^v_{\geq K^{1-\deltap}} \VN[\leq K^{1-\delta}][\fcs], \VN[][-] \big) \Big)
\Big\|_{\Cprod{\eta}{r-1}} \notag \\ 
\lesssim&\,  \big( N^{-\frac{1}{2}+\delta}+ N^{-2\delta_1+4\delta_2}\, N|y| \big) \Dc^2 \label{jacobi:eq-fcsm-e1}
\end{align}
and 
\begin{align}
&\Big\| \Sumlarge_{\substack{M\leq \Nd \colon \\ M \lesssim_\delta K}} \Big( \HHLN[K,M,y][v,(-)\times (\fcs)] - \HHLN[K,M,\bfzero][v,(-)\times (\fcs)] \Big) - \chinull[K] P^{u,v}_{<K^\delta} \Big( \CHHLN_y \big( \VN[][-], P^v_{\geq K^{1-\deltap}} \VN[\leq K^{1-\delta}][\fcs] \big) \Big) \Big\|_{\Cprod{\eta}{r-1}} \notag \\ 
\lesssim&\,  \big( N^{-\frac{1}{2}+\delta}+ N^{-2\delta_1+4\delta_2}\, N|y| \big) \Dc^2. 
\label{jacobi:eq-fcsm-e2} 
\end{align}
\end{lemma}

\begin{proof}
We first prove \eqref{jacobi:eq-fcsm-e1}. From Definition \ref{ansatz:def-hhl}, we obtain that
\begin{equation}\label{jacobi:eq-fcsm-q}
\begin{aligned}
&\Sumlarge_{\substack{M\leq \Nd \colon \\ M \lesssim_\delta K}} \Big( \HHLN[K,M,y][v,(\fcs)\times (-)] - \HHLN[K,M,\bfzero][v,(\fcs)\times (-)] \Big) \\
=&\, \Sumlarge_{\substack{M\leq \Nd \colon \\ M \lesssim_\delta K}}  \chinull[K,M]  P^{u,v}_{<K^\delta} 
\Big( \Ad \big( P_{\leq N}^x P^v_{\geq K^{1-\deltap}} \VN[\leq K^{1-\delta}][\fcs] \big)
\Ad \big( (\Theta^x_y -1) P_{\leq N}^x \IVN[M][-] \big) \Big). 
\end{aligned}
\end{equation}
Due to the smoothness of $\chinull$, we can easily approximate the $\chinull[K,M]$-factor in \eqref{jacobi:eq-fcsm-q} by $\chinull[K]$, which can then be pulled out of the dyadic sum. 
Since $K\geq N^{1-2\delta_1}$, the condition $M\leq \Nd$ already implies that $M\lesssim_\delta K$. Thus, we can replace the sum in \eqref{jacobi:eq-fcsm-q} by the sum over all $M\leq \Nd$. 
We now want to replace $\IVN[M][-]$ in \eqref{jacobi:eq-fcsm-q} by $\Int^v_{u\rightarrow v}\VN[M][-]$. We therefore use the high$\times$high-paraproduct estimate (Lemma \ref{prelim:lem-paraproduct}), Lemma \ref{modulation:lem-linear}, and Corollary \ref{modulation:cor-control-combined}, which yield 
\begin{equation}\label{jacobi:eq-fcsm-p1}
\begin{aligned}
&\Big\| \Sumlarge_{M\leq \Nd} P^{u,v}_{<K^\delta} 
\Big( \Ad \big( P_{\leq N}^x P^v_{\geq K^{1-\deltap}} \VN[\leq K^{1-\delta}][\fcs] \big)
\Ad \big( (\Theta^x_y -1) P_{\leq N}^x \IVN[M][-] \big) \Big) 
\\
&\hspace{3ex} - \Sumlarge_{M\leq \Nd} P^{u,v}_{<K^\delta} 
\Big( \Ad \big( P_{\leq N}^x P^v_{\geq K^{1-\deltap}} \VN[\leq K^{1-\delta}][\fcs] \big)
\Ad \big( (\Theta^x_y -1) P_{\leq N}^x \Int^v_{u\rightarrow v} \VN[M][-] \big) \Big) 
\Big\|_{\Cprod{\eta}{r-1}} \\ 
\lesssim&\, \Sumlarge_{M\leq \Nd} 
\Big\| P^v_{\geq K^{1-\deltap}} \VN[\leq K^{1-\delta}][\fcs] \Big\|_{\Cprod{s}{-2s+\eta}} 
\Big\| \IVN[M][-] - \Int^v_{u\rightarrow v} \VN[M][-] \Big\|_{\Cprod{s}{2s}} \\ 
\lesssim&\, \Sumlarge_{M\leq \Nd} 
K^{(1-\deltap)(-2s+\eta-(r-1))} M^{-\eta} \Dc^2 
\lesssim K^{(1-\delta-2\delta_2) \big(-\tfrac{1}{2}-\delta_1+2\delta_2 + \delta_3\big)} \Dc^2. 
\end{aligned}
\end{equation}
Since $K\gtrsim N^{1-2\delta_1}$, this is acceptable. We now replace $P_{\leq N}^x$ by $\Pbd$ and $\Int^v_{u\rightarrow v}$ by $P^x_{>1} \Int^x_{0\rightarrow x}$. Indeed, using Lemma \ref{jacobi:lem-replacement}, we obtain that 
\begin{align*}
&\Big\| \Sumlarge_{M\leq \Nd}  P^{u,v}_{<K^\delta} 
\Big( \Ad \big( P_{\leq N}^x P^v_{\geq K^{1-\deltap}} \VN[\leq K^{1-\delta}][\fcs] \big)
\Ad \big( (\Theta^x_y -1) \Int^v_{u\rightarrow v} P_{\leq N}^x  \VN[M][-] \big) \Big) \\
&\hspace{3ex} - \Sumlarge_{M\leq \Nd}  P^{u,v}_{<K^\delta} 
\Big( \Ad \big( \Pbd P^v_{\geq K^{1-\deltap}} \VN[\leq K^{1-\delta}][\fcs] \big)
\Ad \big( (\Theta^x_y -1)  P^x_{>1} \Int^x_{0\rightarrow x} \Pbd \VN[M][-] \big) \Big) 
\Big\|_{\Cprod{\eta}{r-1}} \\ 
\lesssim&\, \big( N^{-\frac{1}{2}+\delta}+ N^{-2\delta_1+4\delta_2}\, N|y| \big) \Dc^2, 
\end{align*}
which is acceptable. Since the sum of $\VN[M][-]$ over all $M\in \Dyadiclarge$ satisfying $M\leq \Nd$ is equal to $\VN[][-]$, this completes the proof of \eqref{jacobi:eq-fcsm-e1}. The second estimate is similar, and we omit the details. 
\end{proof}

\begin{lemma}[Killing-term]\label{jacobi:lem-killing}
Let $K,N,\Nd\in\Dyadiclarge$ satisfy $K\geq N^{1-2\delta_1}$. Then, it holds that
\begin{equation*}
\Big\| \Sumlarge_{\substack{M\leq \Nd\colon \\ M \simeq_\delta K}} \Big( \HHLN[K,M,y][v,\textup{kil}]- \HHLN[K,M,\bfzero][v,\textup{kil}] \Big) - \big( - \chinull[K] \coup \CNcsbd(y) \Kil \big) \Big\|_{\Cprod{\eta}{r-1}} \lesssim N^{-2\delta_1} N|y| \Dc^2. 
\end{equation*}
\end{lemma}

\begin{proof}
Since $K\geq N^{1-2\delta_1}$ and $M\leq \Nd$, the condition $M\simeq_\delta K$ is equivalent to $M\geq K^{1-\delta}$. Together with the identity $\Cf^{(\Ncs)}_M(0)=0$, it then follows that
\begin{align}
&\, \Sumlarge_{\substack{M\leq \Nd\colon \\ M \simeq_\delta K}} \hspace{-1ex}\Big( \HHLN[K,M,y][v,\textup{kil}]- \HHLN[K,M,\bfzero][v,\textup{kil}] \Big) - \big( - \chinull[K] \coup \CNcsbd(y) \Kil \big)  \notag \\ 
=& \, \chinull[K] \coup \bigg( - \Sumlarge_{\substack{M \leq \Nd\colon \\ M \geq K^{1-\delta}}} \hspace{-1ex} \Cf^{(\Ncs)}_M(y) 
+ \sum_{\substack{M \leq \Nd}} \hspace{-1ex} \CNcsbd_M(y) \bigg)  \notag \\ 
=&\, \chinull[K] \coup \bigg( \Sumlarge_{\substack{M \leq \Nd\colon \\ M < K^{1-\delta}}} \hspace{-1ex} \Cf^{(\Ncs)}_M(y) 
+ \Sumlarge_{\substack{M \leq \Nd }} \hspace{-1ex} \big( \CNcsbd_M -\Cf^{(\Ncs)}_M \big) (y) 
+ \sum_{\substack{M<\Nlarge}} \CNcsbd_M (y) \bigg). 
\label{jacobi:eq-killing-p1}
\end{align}
Since $\chinull$ is smooth and $\coup\lesssim \Dc^2$, it remains to control the dyadic sums in \eqref{jacobi:eq-killing-p1}. Using Lemma \ref{ansatz:lem-renormalization}, the first dyadic sum in \eqref{jacobi:eq-killing-p1} can be estimated by
\begin{equation*}
 \bigg|  \Sumlarge_{\substack{M \leq \Nd\colon \\ M < K^{1-\delta}}} \hspace{-1ex} \Cf^{(\Ncs)}_M(y)  \bigg|
 \lesssim   \Sumlarge_{\substack{M \leq \Nd\colon \\ M < K^{1-\delta}}} \hspace{-1ex} M |y| \lesssim K^{1-\delta} |y| \lesssim  N^{1-\delta} |y|,
\end{equation*}
which is acceptable. To treat the second and third sum in \eqref{jacobi:eq-killing-p1}, we first note that $\CNcsbd_M -\Cf^{(\Ncs)}_M =0$ for all $M\gg N^{1-2\delta_1}$. Using a minor variant of Lemma \ref{ansatz:lem-renormalization}, it therefore follows that
\begin{equation*}
\bigg| \Sumlarge_{\substack{M \leq \Nd }} \hspace{-1ex} \big( \CNcsbd_M -\Cf^{(\Ncs)}_M \big) (y) 
\bigg| 
+ \bigg| \sum_{\substack{M<\Nlarge}} \CNcsbd_M (y) \bigg| 
\lesssim \Sumlarge_{\substack{M \leq \Nd\colon \\  M\lesssim N^{1-2\delta_1} }} \hspace{-2ex} M |y|
+ \sum_{\substack{M<\Nlarge}} M |y| 
\lesssim N^{1-2\delta_1} |y|, 
\end{equation*}
which is acceptable. 
\end{proof}

\subsubsection{Proof of Proposition \ref{jacobi:prop-shift}} 

We now have all necessary estimates to analyze the difference of $\HHLN[K,y][v]$ and $\HHLN[K,\bfzero][v]$. 

\begin{proof}[Proof of Proposition \ref{jacobi:prop-shift}] 
The first estimate \eqref{jacobi:eq-shift-general} is already contained in Lemma \ref{jacobi:lem-general}. It thus remains to prove that, for all $K\geq N^{1-2\delta_1}$, the second estimate \eqref{jacobi:eq-shift-boundary} is satisfied. Using Definition \ref{ansatz:def-hhl} and \eqref{jacobi:eq-chhl-combined}, it follows that 
\begin{align}
&\HHLN[K,y][v] - \HHLN[K,\bfzero][v]
- \chinull[K] P^{u,v}_{<K^\delta} \CHHLNc_{K,y} \notag \\ 
=&\, \Sumlarge_{\substack{M\leq \Nd \colon \\  M\simeq_\delta K}} 
\Big( \HHLN[K,M,y][v,(-)\times (-)] - \HHLN[K,M,\bfzero][v,(-)\times (-)] \Big) - \chinull[K] P^{u,v}_{<K^\delta} 
\Big(\CHHLN_y \big( \VN[][-], \VN[][-] \big) \Big) 
\label{jacobi:eq-shift-p1} \\
+&\, \Sumlarge_{\substack{M\leq \Nd \colon \\ M \lesssim_\delta K}} \Big( \HHLN[K,M,y][v,(\fcs)\times (-)] - \HHLN[K,M,\bfzero][v,(\fcs)\times (-)] \Big) 
- \chinull[K] P^{u,v}_{<K^\delta} \Big( \CHHLN_y \big( P^v_{\geq K^{1-\deltap}} \VN[\leq K^{1-\delta}][\fcs], \VN[][-] \big) \Big) \label{jacobi:eq-shift-p2} \\
+&\, \Sumlarge_{\substack{M\leq \Nd \colon \\ M \lesssim_\delta K}} \Big( \HHLN[K,M,y][v,(-)\times (\fcs)] - \HHLN[K,M,\bfzero][v,(-)\times (\fcs)] \Big) - \chinull[K] P^{u,v}_{<K^\delta} \Big( \CHHLN_y \big( \VN[][-], P^v_{\geq K^{1-\deltap}} \VN[\leq K^{1-\delta}][\fcs] \big) \Big) 
\label{jacobi:eq-shift-p3} \\
+&\, \Sumlarge_{\substack{M\leq \Nd\colon \\ M \simeq_\delta K}} \Big( \HHLN[K,M,y][v,\textup{kil}]- \HHLN[K,M,\bfzero][v,\textup{kil}] \Big) - \big( - \chinull[K] \coup \CNcsbd(y) \Kil \big). 
\label{jacobi:eq-shift-p4}
\end{align}
The first term \eqref{jacobi:eq-shift-p1} has been estimated in Lemma \ref{jacobi:lemma-mm}, the second and third term \eqref{jacobi:eq-shift-p2} and \eqref{jacobi:eq-shift-p3} have been estimated in Lemma \ref{jacobi:lem-fcsm}, and the fourth term \eqref{jacobi:eq-shift-p4} has been estimated in \eqref{jacobi:lem-killing}. 
\end{proof}

\subsection{\protect{The Cartesian high$\times$high$\rightarrow$low-matrices}}\label{section:jacobi-cartesian}

In Proposition \ref{jacobi:prop-shift}, we showed that the difference 
$\HHLN[K,y][v] - \HHLN[K,\bfzero][v]$ can be approximated using $\chinull[K] P^{u,v}_{<K^\delta} \CHHLNc_{K,y}$. 
In order to prove \eqref{jacobi:eq-main-estimate-2}, we still need to show that $ P^{u,v}_{<K^\delta} \CHHLNc_{K,y}$
can be replaced by $\CHHLN_y(\VN[][],\VN[][])-\coup \CNcsbd(y) \Kil$, which is the subject of the next proposition.

\begin{proposition}[\protect{The Cartesian high$\times$high$\rightarrow$low-matrices}]\label{jacobi:prop-cartesian}
Let the post-modulation hypothesis (Hypothesis \ref{hypothesis:post}) be satisfied and let $K\in \Dyadiclarge$ satisfy $K\geq N^{1-2\delta_1}$.
Then, it holds that 
\begin{equation}\label{jacobi:eq-cartesian}
\begin{aligned}
&\bigg\| \bigg( P^{u,v}_{<K^{\delta}} \Big( \CHHLNc_{K,y} \Big)  - \Big( \CHHLN_y \big( \VN[][], \VN[][] \big) - \coup \CNcsbd(y) \Kil \Big) \bigg) \Theta^x_y P_{\leq N}^x \UN[K][+] \bigg\|_{\Cprod{r-1}{r-1}} \\
\lesssim&\,  N^{-\delta_1+\delta \delta_1} \Dc^3. 
\end{aligned}
\end{equation}
\end{proposition}

In order to prove of Proposition \ref{jacobi:prop-cartesian}, we first need to prove the following three lemmas. In the first lemma, we show that the $P^{u,v}_{<K^\delta}$-operator in  \eqref{jacobi:eq-cartesian} can be removed. 

\begin{lemma}[Removing $P^{u,v}_{< K^\delta}$]\label{jacobi:lem-Puv-Px}
Let the post-modulation hypothesis (Hypothesis \ref{hypothesis:post}) be satisfied. Furthermore, let
$K,M\in \Dyadiclarge$ satisfy $K\geq N^{1-2\delta_1}$ and let
\begin{equation*}
\VN[][\ast_1], \VN[][\ast_2] \in \Big\{ \VN[][-], P^v_{\geq K^{1-\deltap}} \VN[<K^{1-\delta}][\fcs] \Big\}. 
\end{equation*}
For all $y\in \R$, it then holds that 
\begin{align}
\Big\|  \CHHLN_y \big( \VN[][\ast_1], \VN[][\ast_2] \big) 
\Big\|_{\Cprod{s}{2r-1-\eta}} &\lesssim N^{2\delta_1+2\delta_3} \Dc^2, 
\label{jacobi:eq-Puv-Px-p1} \\ 
 \bigg\|  P^{u,v}_{\geq K^{\delta}} \Big( \CHHLN_y \big( \VN[][\ast_1], \VN[][\ast_2] \big) \Big) 
\Theta^x_y P_{\leq N}^x \UN[K][+] \bigg\|_{\Cprod{r-1}{r-1}} &\lesssim N^{-\delta_1}\Dc^3. 
\label{jacobi:eq-Puv-Px}
\end{align}
\end{lemma}

\begin{proof}
For expository purposes, we split the argument into two steps. \\

\emph{First step:} In the first step, we prove \eqref{jacobi:eq-Puv-Px-p1}. 
To this end, we recall from the proof of Lemma \ref{jacobi:lem-replacement} that, on the frequency-support of 
$\Pbd \VN[][\ast_1]$ and $\Pbd \VN[][\ast_2]$, the $v$-frequencies are much larger than the $u$-frequencies. Due to the $P^x_{\leq N^{\deltap}}$-operator,  it follows that $ \CHHLN_y \big( \VN[][\ast_1], \VN[][\ast_2] \big) $ can only contain high$\times$high$\rightarrow$low-interactions in the $v$-variable. As a result, we obtain the identity
\begin{align*}
 \CHHLN_y \big( \VN[][\ast_1], \VN[][\ast_2] \big) 
=  P^x_{<N^{\deltap}}\Big( 
\ad\big( \Pbd \VN[][\ast_1] \big) \Para[v][sim] \ad\big( (\Theta^x_y -1) P^x_{>1} \Int^x_{0\rightarrow x} \Pbd \VN[][\ast_2] \big) \Big).
\end{align*}
Using the high$\times$high-estimate (Lemma \ref{prelim:lem-paraproduct}) and Lemma \ref{prelim:lem-Cartesian-integral}, we then obtain that
\begin{align*}
&\, \Big\|  \CHHLN_y \big( \VN[][\ast_1], \VN[][\ast_2] \big) 
\Big\|_{\Cprod{s}{2r-1-\eta}}\\ 
\lesssim&\,  \Big\| \ad\big( \Pbd \VN[][\ast_1] \big) \Para[v][sim] \ad\big( (\Theta^x_y -1) P^x_{>1} \Int^x_{0\rightarrow x} \Pbd \VN[][\ast_2] \big) 
\Big\|_{\Cprod{s}{2r-1-\eta}} \\
\lesssim&\, \big\|  \Pbd \VN[][\ast_1] \big\|_{\Cprod{s}{r-1}} \big\| P^x_{>1} \Int^x_{0\rightarrow x} \Pbd \VN[][\ast_2] \big\|_{\Cprod{s}{r}} 
\lesssim \big\| \VN[][\ast_1] \big\|_{\Cprod{s}{r-1}} \big\| \VN[][\ast_2] \big\|_{\Cprod{s}{r-1}}. 
\end{align*}
Using Lemma \ref{modulation:lem-linear} and Corollary \ref{modulation:cor-control-combined}, it holds that
\begin{equation*}
\big\| \VN[][\ast_1] \big\|_{\Cprod{s}{r-1}}, \big\| \VN[][\ast_2] \big\|_{\Cprod{s}{r-1}} 
\lesssim N^{r-\frac{1}{2}+\eta} \Dc = N^{\delta_1+\delta_2} \Dc, 
\end{equation*}
which completes the proof of \eqref{jacobi:eq-Puv-Px-p1}.\\

\emph{Second step:} In the second step, we prove the estimate \eqref{jacobi:eq-Puv-Px}. To this end, we decompose 
\begin{align}
  &\, \big( P^{u,v}_{\geq K^\delta} \CHHLN_y \big( \VN[][\ast_1], \VN[][\ast_2] \big) \big) \Theta^x_y P_{\leq N}^x \UN[K][+] \notag \\
  =&\, \big( P^{u,v}_{\geq K^\delta} \CHHLN_y \big( \VN[][\ast_1], \VN[][\ast_2] \big) \big) \Para[u][nsim] \Theta^x_y P_{\leq N}^x \UN[K][+] \label{jacobi:eq-Puv-Px-p2}\\
  +&\, \big( P^{u,v}_{\geq K^\delta} \CHHLN_y \big( \VN[][\ast_1], \VN[][\ast_2] \big) \big) \Para[u][sim] \Theta^x_y P_{\leq N}^x \UN[K][+]. \label{jacobi:eq-Puv-Px-p3}
\end{align}
We now estimate \eqref{jacobi:eq-Puv-Px-p2} and \eqref{jacobi:eq-Puv-Px-p3} separately. Using Lemma \ref{prelim:lem-paraproduct}, we estimate the non-resonant term by 
\begin{align}
&\, \big\| \eqref{jacobi:eq-Puv-Px-p2} \big\|_{\Cprod{r-1}{r-1}} \notag\\ 
\lesssim&\, \big\| P^{u,v}_{\geq K^\delta} \CHHLN_y \big( \VN[][\ast_1], \VN[][\ast_2] \big) \big\|_{\Cprod{\eta}{r-1}} 
\big\|  \UN[K][+] \big\|_{\Cprod{r-1}{s}}  \notag \\ 
\lesssim&\, \max \Big( K^{\delta (\eta-s)}, K^{\delta (r-1-(2r-1-\eta))} \Big)
\big\|  \CHHLN_y \big( \VN[][\ast_1], \VN[][\ast_2] \big)\big\|_{\Cprod{s}{2r-1-\eta}} \big\| \UN[K][+] \big\|_{\Cprod{r-1}{s}}. \label{jacobi:eq-Puv-Px-p4} 
\end{align}
Using \eqref{jacobi:eq-Puv-Px-p1} and Lemma \ref{modulation:lem-linear}, we obtain that
\begin{equation}\label{jacobi:eq-Puv-Px-p5}
\eqref{jacobi:eq-Puv-Px-p4} \lesssim \max \Big( K^{\delta (\eta-s)}, K^{\delta (r-1-(2r-1-\eta))} \Big) N^{2\delta_1+2\delta_3} K^{r-\frac{1}{2}+\eta} \Dc^3. 
\end{equation}
Using the definition of our parameters from \eqref{prelim:eq-parameter-regularities} and $K\geq N^{1-2\delta_1}$, we obtain that
\begin{equation*}
\eqref{jacobi:eq-Puv-Px-p5} \lesssim N^{-\frac{\delta}{2}+\mathcal{O}(\delta_1)}, 
\end{equation*}
which is more than acceptable. It remains to treat the resonant-term \eqref{jacobi:eq-Puv-Px-p3}. Due to the high$\times$high-product $\Para[u][sim]$,
$ \CHHLN_y( \VN[][\ast_1], \VN[][\ast_2])$ only enters into \eqref{jacobi:eq-Puv-Px-p3} at $u$-frequencies $\gtrsim K$. Due to the $P^x_{\leq N^{\deltap}}$-operator in the definition of $ \CHHLN_y( \VN[][\ast_1], \VN[][\ast_2])$ and $K\geq N^{1-2\delta_1}$, this implies that $ \CHHLN_y( \VN[][\ast_1], \VN[][\ast_2])$ then also only enters into \eqref{jacobi:eq-Puv-Px-p3} at $v$-frequencies $\gtrsim K$. Using Lemma \ref{prelim:lem-paraproduct}, we then obtain that
\begin{align}
\big\| \eqref{jacobi:eq-Puv-Px-p3} \big\|_{\Cprod{r-1}{r-1}}
&\lesssim \big\|  P^{v}_{\gtrsim K}  \CHHLN_y\big( \VN[][\ast_1], \VN[][\ast_2]\big) \big\|_{\Cprod{s}{r-1}} 
\big\| \UN[K][+] \big\|_{\Cprod{r-1}{s}} \notag  \\
&\lesssim K^{r-1-(2r-1-\eta)} \big\|  \CHHLN_y\big( \VN[][\ast_1], \VN[][\ast_2]\big) \big\|_{\Cprod{s}{2r-1-\eta}} \big\| \UN[K][+] \big\|_{\Cprod{r-1}{s}}. \label{jacobi:eq-Puv-Px-p6}
\end{align}
Using definition of our parameters from \eqref{prelim:eq-parameter-regularities}, $K\geq N^{1-2\delta_1}$, 
Lemma \ref{modulation:lem-linear}, and \eqref{jacobi:eq-Puv-Px-p1}, we obtain that
\begin{equation*}
\eqref{jacobi:eq-Puv-Px-p6} \lesssim N^{-\frac{1}{2}+\mathcal{O}(\delta_1)} \Dc^3,
\end{equation*}
which is more than acceptable.
\end{proof}

In the next lemma, we prove a general estimate which relies on the $L_t^\infty L_x^\infty$-norm.

\begin{lemma}\label{jacobi:lem-F-Up}
Let the post-modulation hypothesis (Hypothesis \ref{hypothesis:post}) be satisfied and let $K\in \Dyadiclarge$ satisfy $N^{1-2\delta_1}\leq K \lesssim N$. Furthermore, let $y\in \R$ and let $F\colon \R^{1+1}  \rightarrow \End(\frkg)$. Then, it holds that
\begin{equation}\label{jacobi:eq-F-Up}
\begin{aligned}
 \Big\|  \big( P_{\leq N^{\deltap}}^x F  \big) \, \Theta^x_y P_{\leq N}^x \UN[K][+] \Big\|_{\Cprod{r-1}{r-1}} 
\lesssim \Dc N^{\delta_1+4\delta_3} \big\| P_{\leq N^{\deltap}}^x F \big\|_{L_t^\infty L_x^\infty}.
\end{aligned}
\end{equation}
\end{lemma}

\begin{proof}
We first use a dyadic decomposition of $F$, i.e., we decompose
\begin{equation}\label{jacobi:eq-F-Up-p1}
\begin{aligned}
\big( P_{\leq N^{\deltap}}^x F \big) \, \Theta^x_y P_{\leq N}^x \UN[K][+] 
= \sum_{L,M} 
\big( P^u_L P^v_M P_{\leq N^{\deltap}}^x F \big) \, \Theta^x_y P_{\leq N}^x \UN[K][+]. 
\end{aligned}
\end{equation}
In order to bound \eqref{jacobi:eq-F-Up-p1}, we distinguish two cases. In the case $L\ll K$, \eqref{jacobi:eq-F-Up-p1} only contains low$\times$high-interactions in the $u$-variable. Using Lemma \ref{prelim:lem-paraproduct} and Lemma \ref{modulation:lem-linear}, 
it then follows that 
\begin{equation}\label{jacobi:eq-F-Up-p2}
\begin{aligned}
&\, \Big\| \big( P^u_L P^v_M P_{\leq N^{\deltap}}^x F \big) \,  \Theta^x_y P_{\leq N}^x \UN[K][+] \Big\|_{\Cprod{r-1}{r-1}}  
\lesssim  \big\| P^u_L P^v_M P_{\leq N^{\deltap}}^x F \big\|_{\Cprod{\eta}{r-1}} \big\| \UN[K][+] \big\|_{\Cprod{r-1}{s}} \\ 
\lesssim & \, \Dc K^{r-\frac{1}{2}+\eta} L^\eta M^{r-1} \big\| P_{\leq N^{\deltap}}^x F  \big\|_{L^\infty_u L^\infty_v} 
\lesssim \Dc K^{\delta_1+4\delta_3} (LM)^{-\eta} \big\| P_{\leq N^{\deltap}}^x F  \big\|_{L^\infty_u L^\infty_v} 
\end{aligned}
\end{equation}
In the last inequality, we used $L\ll K$ and $r-1\leq -\eta$. It therefore remains to treat the case $L\gtrsim K$. Since $K\geq N^{1-2\delta_1}$,
$P^u_L P^v_M P_{\leq N^{\deltap}}^x F $ can then only be non-zero in the case $L\sim M$. Using our para-product estimate (Lemma \ref{prelim:lem-paraproduct}) and Lemma \ref{modulation:lem-linear}, it then follows that 
\begin{equation}\label{jacobi:eq-F-Up-p3}
\begin{aligned}
&\, \Big\| \big( P^u_L P^v_M P_{\leq N^{\deltap}}^x F \big) \,  \Theta^x_y P_{\leq N}^x \UN[K][+] \Big\|_{\Cprod{r-1}{r-1}}  
\lesssim \big\| P^u_L P^v_M P_{\leq N^{\deltap}}^x F \big\|_{\Cprod{1-r-2\eta}{r-1}} \big\| \UN[K][+] \big\|_{\Cprod{r-1+3\eta}{s}} \\
\lesssim&\, \Dc K^{r-\frac{1}{2}+4\eta} L^{1-r-2\eta} M^{r-1} \big\| P_{\leq N^{\deltap}}^x F  \big\|_{L^\infty_u L^\infty_v} 
\lesssim K^{\delta_1+4\delta_3} (LM)^{-\eta}  \big\| P_{\leq N^{\deltap}}^x F  \big\|_{L^\infty_u L^\infty_v} .
\end{aligned}
\end{equation}
By combining \eqref{jacobi:eq-F-Up-p1}, \eqref{jacobi:eq-F-Up-p2}, and \eqref{jacobi:eq-F-Up-p3}, it follows that 
\begin{align*}
 \Big\|  \big( P_{\leq N^{\deltap}}^x F  \big) \, \Theta^x_y P_{\leq N}^x \UN[K][+] \Big\|_{\Cprod{r-1}{r-1}}  
 \lesssim \Dc K^{\delta_1+4\delta_3} \big\| P_{\leq N^{\deltap}}^x F  \big\|_{L^\infty_u L^\infty_v}.
\end{align*}
Since $K\lesssim N$ and $\big\| P_{\leq N^{\deltap}}^x F  \big\|_{L^\infty_u L^\infty_v}=\big\| P_{\leq N^{\deltap}}^x F  \big\|_{L^\infty_t L^\infty_x}$, this implies \eqref{jacobi:eq-F-Up}.
\end{proof}

In the last lemma of this subsection, we bound the $L_t^\infty L_x^\infty$-norm of $\CHHLN_y(\VN[][\ast_1],\VN[][\ast_2])$ for several choices of $\VN[][\ast_1]$ and $\VN[][\ast_2]$.

\begin{lemma}[\protect{$L_t^\infty L_x^\infty$-bound of $\CHHLN_y$}]\label{jacobi:lem-remaining-interactions}
Assume that the post-modulation hypothesis (Hypothesis \ref{hypothesis:post}) is satisfied. Furthermore, let $\ast_1,\ast_2\in \{-,+-,+,\fs-,+\fs,\fs\}$ satisfy
\begin{equation}\label{jacobi:eq-remaining-interactions-a}
\big( \ast_1,\ast_2 \big) \neq \big( -, - \big), \big( -, +\fs\big), \big( -, \fs\big), \big( +\fs, -\big), \big( \fs, -\big). 
\end{equation}
For all $y\in \R$, it then holds that 
\begin{equation}\label{jacobi:eq-remaining-interactions}
\Big\| \CHHLN_y\big(\VN[][\ast_1],\VN[][\ast_2]\big)  \Big\|_{L_t^\infty L_x^\infty} 
\lesssim N^{-2(1-4\delta_1)\delta_1} \Dc^2.
\end{equation}
Furthermore, for all $K\in \Dyadiclarge$ satisfying $K\geq N^{1-2\delta_1}$, it holds that
\begin{align}
&\, \Big\|  \CHHLN_y \Big( \VN[][-] , P^v_{\geq K^{1-\deltap}} \VN[<K^{1-\delta}][\fcs] - \VN[][+\fs] - \VN[][\fs] \Big) \Big\|_{L_t^\infty L_x^\infty} \notag \\
+&\, \Big\|  \CHHLN_y \Big(P^v_{\geq K^{1-\deltap}} \VN[<K^{1-\delta}][\fcs] - \VN[][+\fs] - \VN[][\fs],  \VN[][-] \Big) \Big\|_{L_t^\infty L_x^\infty} 
\hspace{1ex}\lesssim N^{-\frac{1}{2}+20\delta} \Dc^2. \label{jacobi:eq-remaining-interactions-mcs}
\end{align}
\end{lemma}

\begin{proof}
We first prove the three separate estimates
\begin{align}
\textup{LHS of } \eqref{jacobi:eq-remaining-interactions}
&\lesssim N^{-\frac{1}{2}+20\delta} \big\| \VN[][\ast_1] \big\|_{C_t^0 \C_x^{s-1}} \big\| \VN[][\ast_2] \big\|_{C_t^0 \C_x^{\scrr-1}}. \label{jacobi:eq-remaining-interactions-p1} \\ 
\textup{LHS of } \eqref{jacobi:eq-remaining-interactions}
&\lesssim N^{-\frac{1}{2}+20\delta} \big\| \VN[][\ast_1] \big\|_{C_t^0 \C_x^{\scrr-1}} \big\| \VN[][\ast_2] \big\|_{C_t^0 \C_x^{s-1}}, \label{jacobi:eq-remaining-interactions-p2}\\ 
\textup{LHS of } \eqref{jacobi:eq-remaining-interactions}
&\lesssim N^{-2 (1-4\delta_1) \delta_1} \big\| \VN[][\ast_1] \big\|_{C_t^0 \C_x^{r-1}} \big\| \VN[][\ast_2] \big\|_{C_t^0 \C_x^{r-1}}. \label{jacobi:eq-remaining-interactions-p3}  
\end{align}
To this end, we first note that the left-hand side of \eqref{jacobi:eq-remaining-interactions} only contains high$\times$high-interactions. Using paraproduct estimates and Lemma \ref{prelim:lem-integral}, it therefore follows that 
\begin{align*}
\textup{LHS of } \eqref{jacobi:eq-remaining-interactions}
\lesssim&\, \Big\| \Pbd \VN[][\ast_1] \Big\|_{C_t^0 \C_x^{-\frac{1}{2}+\eta}} \Big\|  \big( \Theta^x_y - 1 \big) P^x_{>1} \Int^x_{0\rightarrow x} \Pbd \VN[][\ast_2] \Big\|_{C_t^0 \C_x^{\frac{1}{2}+\eta}} \\ 
\lesssim&\, \Big\| \Pbd \VN[][\ast_1] \Big\|_{C_t^0 \C_x^{-\frac{1}{2}+\eta}}  \Big\| \Pbd \VN[][\ast_2] \Big\|_{C_t^0 \C_x^{-\frac{1}{2}+\eta}}. 
\end{align*}
In order to obtain \eqref{jacobi:eq-remaining-interactions-p1}, we now bound
\begin{align*}
&\Big\| \Pbd \VN[][\ast_1] \Big\|_{C_t^0 \C_x^{-\frac{1}{2}+\eta}}  \Big\| \Pbd \VN[][\ast_2] \Big\|_{C_t^0 \C_x^{-\frac{1}{2}+\eta}} \\ 
\lesssim&\,  N^{-\frac{1}{2}+\eta-(s-1)} \big( N^{1-2\delta_1}\big)^{-\frac{1}{2}+\eta-(\scrr-1)}  \big\| \VN[][\ast_1] \big\|_{C_t^0 \C_x^{s-1}} \big\| \VN[][\ast_1] \big\|_{C_t^0 \C_x^{\scrr-1}}.
\end{align*}
Since $\scrr=1-10\delta$, this yields an acceptable contribution. The second estimate \eqref{jacobi:eq-remaining-interactions-p2} follows by using the same argument with reversed roles of $\VN[][\ast_1]$ and $\VN[][\ast_2]$. For the third inequality \eqref{jacobi:eq-remaining-interactions-p3}, we bound
\begin{align*}
&\Big\| \Pbd \VN[][\ast_1] \Big\|_{C_t^0 \C_x^{-\frac{1}{2}+\eta}}  \Big\| \Pbd \VN[][\ast_2] \Big\|_{C_t^0 \C_x^{-\frac{1}{2}+\eta}} \\ 
\lesssim&\, \big( N^{1-2\delta_1} \big)^{-\frac{1}{2}+\eta-(r-1)} \big( N^{1-2\delta_1} \big)^{-\frac{1}{2}+\eta-(r-1)} 
\Big\| \Pbd \VN[][\ast_1] \Big\|_{C_t^0 \C_x^{r-1}}  \Big\| \Pbd \VN[][\ast_2] \Big\|_{C_t^0 \C_x^{r-1}}.
\end{align*}
Since
\begin{equation*}
2 (1-2\delta_1) \Big( - \frac{1}{2}+\eta - (r-1) \Big) = - 2 (1-2\delta_1) \delta_1 + \mathcal{O}(\delta_3), 
\end{equation*}
this is acceptable. The first desired inequality \eqref{jacobi:eq-remaining-interactions} now follows from Lemma \ref{modulation:lem-Cartesian}, our Assumption \eqref{jacobi:eq-remaining-interactions-a}, and the three estimates \eqref{jacobi:eq-remaining-interactions-p1}, \eqref{jacobi:eq-remaining-interactions-p2}, and \eqref{jacobi:eq-remaining-interactions-p3}. The second desired inequality \eqref{jacobi:eq-remaining-interactions-mcs} follows from Lemma \ref{modulation:lem-Cartesian-frequency-boundary}, \eqref{jacobi:eq-remaining-interactions-p1}, and \eqref{jacobi:eq-remaining-interactions-p2}.
\end{proof}

Equipped with Lemma \ref{jacobi:lem-Puv-Px}, Lemma \ref{jacobi:lem-F-Up}, and Lemma \ref{jacobi:lem-remaining-interactions}, we can now prove Proposition \ref{jacobi:prop-cartesian}.

\begin{proof}[Proof of Proposition \ref{jacobi:prop-cartesian}]
We decompose the argument in \eqref{jacobi:eq-cartesian} as 
\begin{align}
&\bigg( P^{u,v}_{<K^{\delta}} \Big( \CHHLNc_{K,y} \Big)  - \Big( \CHHLN_y \big( \VN[][], \VN[][] \big) - \coup \CNcsbd(y) \Kil \Big) \bigg) \Theta^x_y P_{\leq N}^x \UN[K][+] \notag \\
=&\, - P^{u,v}_{\geq K^\delta} \Big( \CHHLNc_{K,y} \Big) \Theta^x_y P_{\leq N}^x \UN[K][+] \label{jacobi:eq-cartesian-p1}\\
+&\, \Big( \CHHLNc_{K,y} - \Big( \CHHLN_y \big( \VN[][], \VN[][] \big) - \coup \CNcsbd(y) \Kil \Big)  \Big)\Theta^x_y P_{\leq N}^x \UN[K][+].
\label{jacobi:eq-cartesian-p2}
\end{align}
The first term \eqref{jacobi:eq-cartesian-p1} has already been estimated in Lemma \ref{jacobi:lem-Puv-Px}. To estimate the second term \eqref{jacobi:eq-cartesian-p2}, we first use Lemma \ref{jacobi:lem-F-Up}, which yields that
\begin{equation}\label{jacobi:eq-cartesian-p3}
\big\| \eqref{jacobi:eq-cartesian-p2} \big\|_{\Cprod{r-1}{r-1}} 
\lesssim N^{\delta_1+4\delta_3} \Dc \, 
\Big\| \CHHLNc_{K,y} - \Big( \CHHLN_y \big( \VN[][], \VN[][] \big) - \coup \CNcsbd(y) \Kil \Big) \Big\|_{L_t^\infty L_x^\infty}.
\end{equation}
Using the definition of $\CHHLNc_{K,y}$ from \eqref{jacobi:eq-chhl-combined}, it then follows that
\begin{equation}\label{jacobi:eq-cartesian-p4}
\begin{aligned}
\eqref{jacobi:eq-cartesian-p3} 
&\lesssim N^{\delta_1+4\delta_3} \Dc \, 
\Big\|\, \CHHLN\big( \VN[][-], \VN[][-] \big) 
 + \CHHLN \big( P^v_{\geq K^{1-\deltap}} \VN[\leq K^{1-\delta}][\fcs], \VN[][-] \big) \\ 
 &\hspace{14ex}+ \CHHLN \big( P^v_{\geq K^{1-\deltap}} \VN[\leq K^{1-\delta}][\fcs], \VN[][-] \big) 
 - \CHHLN_y \big( \VN[][], \VN[][] \big) \Big\|_{L_t^\infty L_x^\infty}.
\end{aligned}
\end{equation}
By inserting the Ansatz for $\VN[][]$ from \eqref{ansatz:eq-VN-rigorous-decomposition} and using Lemma \ref{jacobi:lem-remaining-interactions}, it then follows that
\begin{equation*}
\eqref{jacobi:eq-cartesian-p4} 
\lesssim N^{\delta_1+4\delta_3} \Big( N^{-2 (1-4\delta_1)\delta_1} + N^{-\frac{1}{2}+20\delta}\Big) \Dc^3
\lesssim N^{-\delta_1+8\delta_1^2 + 4\delta_3} \Dc^3 \lesssim N^{-\delta_1+\delta \delta_1} \Dc^3, 
\end{equation*}
which is acceptable.
\end{proof}

\subsection{Proof of Proposition \ref{jacobi:prop-main}}

Equipped with Proposition \ref{jacobi:prop-HHL-without-shift} and Proposition \ref{jacobi:prop-shift}, we can now prove the main result of this section. 
\begin{proof}[Proof of Proposition \ref{jacobi:prop-main}:] 
We first prove \eqref{jacobi:eq-bourgain-bulut-1} and \eqref{jacobi:eq-bourgain-bulut-2}. Due to symmetry in the $u$ and $v$-variables, it suffices to prove \eqref{jacobi:eq-bourgain-bulut-1}. Using Definition \ref{jacobi:def-cf} and Definition \ref{jacobi:def-chhl}, we obtain that 
\begin{equation}\label{jacobi:eq-main-p0}
\begin{aligned}
& \, \chi\,  
\Big( P_{\leq N}^x \CHHLN \big( V, V \big) P_{\leq N}^x - \coup \Renormbd[\Ncs] \Big)  \UN[K][+]  \\
=&\, P_{\leq N}^x \int_{\R} \dy \, \big( \widecheck{\rho}_{\leq N} \ast \widecheck{\rho}_{\leq N} \big)(y) 
\chi  \big( \CHHLN_y \big( V , V \big) - \coup \CNcsbd(y) \Kil \big) \Theta^x_y P_{\leq N}^x \UN[K][+]. 
\end{aligned}
\end{equation}
Due to the $P_{\leq N^{\deltap}}^x$-operator in Definition \ref{jacobi:def-chhl} and the fact that the cut-off function $\chi$ only depends on the time-variable $t$, it follows that $\chi  \big( \CHHLN_y \big( V , V \big) - \coup \CNcsbd(y) \Kil \big)$ is supported on $x$-frequencies $\lesssim N^{\deltap}$. Using Lemma \ref{jacobi:lem-F-Up}, it then follows that
\begin{align*}
\big\| \eqref{jacobi:eq-main-p0}\big\|_{\Cprod{r-1}{r-1}}
&\lesssim \Dc N^{\delta_1+4\delta_3} 
\int_{\R} \dy \, \big|  \big( \widecheck{\rho}_{\leq N} \ast \widecheck{\rho}_{\leq N} \big)(y)  \big| 
\big\| \chi  \big( \CHHLN_y \big( V , V \big) - \coup \CNcsbd(y) \Kil \big) \big\|_{L_t^\infty L_x^\infty}  \\
&\lesssim \Dc N^{\delta_1+4\delta_3} \sup_{y\in \R} \, \langle N y \rangle^{-10} \big\| \chi  \big( \CHHLN_y \big( V , V \big) - \coup \CNcsbd(y) \Kil \big) \big\|_{L_t^\infty L_x^\infty}.
\end{align*}
After summing over $N^{1-2\delta_1}\leq K \lesssim N$, this yields \eqref{jacobi:eq-bourgain-bulut-1}.
It remains to prove \eqref{jacobi:eq-main-estimate-1} and \eqref{jacobi:eq-main-estimate-2}, which are more difficult. We first recall from Definition \ref{ansatz:def-jacobi-errors} and Definition \ref{jacobi:def-modified-jacobi} that 
\begin{equation*}
\JcbNErr[][] = \JcbNErr[][+] + \JcbNErr[][-]\qquad \text{and} \qquad 
\JcbNErr[][\dagger] = \JcbNErr[][+,\dagger] + \JcbNErr[][-,\dagger]. 
\end{equation*}
Due to the symmetry in the $u$ and $v$-variables, we only examine $\JcbNErr[][+]$ and $\JcbNErr[][+,\dagger]$. 
To this end, we first recall from Definition \ref{ansatz:def-jacobi-errors} that 
\begin{equation}\label{jacobi:eq-main-p1} 
\JcbNErr[K][+] = P_{\leq N}^x \HHLN[K][v] P_{\leq N}^x \UN[K][+] -  \Big[ (P_{\leq N}^x)^2 \UN[K][+], \SHHLN[K][v] \Big]_{\leq N}. 
\end{equation}
From Definition \ref{ansatz:def-hhl}, it follows that
\begin{equation}\label{jacobi:eq-main-p2} 
\begin{aligned}
P_{\leq N}^x \HHLN[K][v] P_{\leq N}^x \UN[K][+] 
= P_{\leq N}^x \int_{\R} \dy \big(\widecheck{\rho}_{\leq N} \ast \widecheck{\rho}_{\leq N}\big)(y)\, \HHLN[K,y][v] \Theta^x_y P_{\leq N}^x \UN[K][+]. 
\end{aligned}
\end{equation}
From the definition of $\SHHLN[K][v]$, it directly follows that $P_{\leq N}^x \SHHLN[K][v]=\SHHLN[K][v]$. Using the integral representation of $P_{\leq N}^x$ and the adjoint map from \eqref{prelim:eq-adjoint-map}, it then follows that 
\begin{equation}\label{jacobi:eq-main-p3} 
\begin{aligned}
- \Big[ (P_{\leq N}^x)^2 \UN[K][+], \SHHLN[K][v] \Big]_{\leq N} 
=&\,  - P_{\leq N}^x \Big[ (P_{\leq N}^x)^2 P_{\leq N}^x \UN[K][+] ,  \SHHLN[K][v] \Big] \\ 
=&\, P_{\leq N}^x \int_{\R} \dy\,  (\widecheck{\rho}_{\leq N} \ast \widecheck{\rho}_{\leq N})(y)  
\Ad\big( \SHHLN[K][v] \big) \Theta^x_y P_{\leq N}^x \UN[K][+]. 
\end{aligned}
\end{equation}
By combining \eqref{jacobi:eq-main-p1}, \eqref{jacobi:eq-main-p2}, and \eqref{jacobi:eq-main-p3}, it follows that
\begin{align}
\JcbNErr[K][+] 
&= P_{\leq N}^x 
\int_{\R} \dy \, \big(\widecheck{\rho}_{\leq N} \ast \widecheck{\rho}_{\leq N}\big)(y)\, 
\Big( \HHLN[K,y][v] + \Ad\big( \SHHLN[K][v] \big) \Big) 
\Theta^x_y P_{\leq N}^x \UN[K][+] \notag \\ 
&= P_{\leq N}^x 
\int_{\R} \dy \, \big(\widecheck{\rho}_{\leq N} \ast \widecheck{\rho}_{\leq N}\big)(y)\, 
\Big( \HHLN[K,\bfzero][v] + \Ad\big( \SHHLN[K][v] \big) \Big) 
\Theta^x_y P_{\leq N}^x \UN[K][+] 
\label{jacobi:eq-main-p4}\\
&+P_{\leq N}^x 
\int_{\R} \dy \, \big(\widecheck{\rho}_{\leq N} \ast \widecheck{\rho}_{\leq N}\big)(y)\, 
\Big( \HHLN[K,y][v] -  \HHLN[K,\bfzero][v]  \Big) 
\Theta^x_y P_{\leq N}^x \UN[K][+]. 
\label{jacobi:eq-main-p5}
\end{align} 
We first control \eqref{jacobi:eq-main-p4}. From the definitions of $\HHLN[K,\bfzero][v]$ and $\SHHLN[K][v]$, it follows that both are supported on $u$-frequencies $\lesssim K^\delta$. Using the low$\times$high-estimate (Lemma \ref{prelim:lem-paraproduct}), Lemma \ref{modulation:lem-linear}, and Proposition \ref{jacobi:prop-HHL-without-shift}, we obtain that 
\begin{align*}
\big\| \eqref{jacobi:eq-main-p4} \big\|_{\Cprod{r-1}{r-1}}
&\lesssim \int_{\R} \dy \, \big| \big( \widecheck{\rho}_{\leq N} \ast \widecheck{\rho}_{\leq N}\big)(y)\big| \, 
\big\| \HHLN[K,\bfzero][v] + \Ad\big( \SHHLN[K][v] \big) \big\|_{\Cprod{\eta}{r-1}} 
\big\| \Theta^x_y P_{\leq N}^x \UN[K][+] \big\|_{\Cprod{r-1}{s}} \\
&\lesssim \int_{\R} \dy \, \big| \big( \widecheck{\rho}_{\leq N} \ast \widecheck{\rho}_{\leq N}\big)(y)\big| \, 
\big( K^{-\frac{1}{2}+\delta} + N^{-\delta+\vartheta} \big) K^{r-\frac{1}{2}+\eta} \Dc^3 \\ 
&\lesssim \big( K^{-\frac{1}{2}+\delta+\delta_1+\delta_3} + N^{-\delta+\delta_4} K^{\delta_1+\delta_3} \big) \Dc^3. 
\end{align*}
Since $K\lesssim N$, this yields an acceptable contribution to both \eqref{jacobi:eq-main-estimate-1} and \eqref{jacobi:eq-main-estimate-2}. In order to complete the proof of \eqref{jacobi:eq-main-estimate-1} and \eqref{jacobi:eq-main-estimate-2}, it then suffices to prove the three estimates 
\begin{align}
\big\| \eqref{jacobi:eq-main-p5} \big\|_{\Cprod{r-1}{r-1}} 
&\lesssim \Big( K^{-\delta+2\delta_1} + \frac{\Nd}{N} N^{3\delta_3} \Big) \Dc^3, \label{jacobi:eq-main-p6} \\ 
\mathbf{1} \big\{ K < N^{1-2\delta_1} \big\} \big\| \eqref{jacobi:eq-main-p5} \big\|_{\Cprod{r-1}{r-1}} 
&\lesssim \Big( K^{-\delta+2\delta_1}  + N^{-\delta_2} \Big) \Dc^3, \label{jacobi:eq-main-p7} 
\end{align}
and 
\begin{equation}\label{jacobi:eq-main-p8}
\begin{aligned}
&\, \mathbf{1} \big\{ K \geq N^{1-2\delta_1} \big\} \Big\| \eqref{jacobi:eq-main-p5} - \chi \Big( P_{\leq N}^x \CHHLN\big( \VN[][], \VN[][] \big) P_{\leq N}^x - \coup \Renormbd[\Ncs] \Big) \UN[K][+] \Big\|_{\Cprod{r-1}{r-1}} \\
\lesssim&\,   N^{-\delta_1+2\delta \delta_1} \Dc^3,
\end{aligned} 
\end{equation}
which are shown separately. \\ 

\emph{Proof of \eqref{jacobi:eq-main-p6}:}
From the definitions of $\HHLN[K,y][v]$ and $\HHLN[K,\bfzero][v]$, it follows that both of them are supported on $u$-frequencies $\lesssim K^\delta$. Using the low$\times$high-estimate (Lemma \ref{prelim:lem-paraproduct}), Lemma \ref{modulation:lem-linear}, and 
\eqref{jacobi:eq-shift-general} from \mbox{Proposition \ref{jacobi:prop-shift}}, it then follows that 
\begin{equation}\label{jacobi:eq-main-p9}
\begin{aligned}
\big\| \eqref{jacobi:eq-main-p5} \big\|_{\Cprod{r-1}{r-1}} 
&\lesssim \int_{\R} \dy \, \big| \big(\widecheck{\rho}_{\leq N} \ast \widecheck{\rho}_{\leq N}\big)(y)\big| \, 
\big\|  \HHLN[K,y][v] -  \HHLN[K,\bfzero][v]  \big\|_{\Cprod{\eta}{r-1}}
\big\| \Theta^x_y P_{\leq N}^x \UN[K][+] \big\|_{\Cprod{r-1}{s}}  \\ 
&\lesssim \int_{\R} \, \big| \big(\widecheck{\rho}_{\leq N} \ast \widecheck{\rho}_{\leq N}\big)(y)\big| \, \langle N y \rangle 
\Big( K^{-\frac{1}{2}+\delta} + N^{-\delta+\vartheta} + \frac{\Nd}{N} \Nd^{-\delta_1+2\delta_3} \Big) K^{r-\frac{1}{2}+\eta} \Dc^3 \\
&\lesssim \Big( K^{-\frac{1}{2}+\delta + \delta_1 + \delta_3} + N^{-\delta+\delta_4} K^{\delta_1+\delta_3} 
+  \frac{\Nd}{N} \Nd^{-\delta_1+2\delta_3} K^{\delta_1+\delta_3} \Big) \Dc^3. 
\end{aligned}
\end{equation}
Using only $K\lesssim N$, the first two summands in \eqref{jacobi:eq-main-p9} can both be bounded by $K^{-\delta+2\delta_1}$, and it follows that
\begin{equation}\label{jacobi:eq-main-p10}
\big\| \eqref{jacobi:eq-main-p5} \big\|_{\Cprod{r-1}{r-1}}  \lesssim 
 \Big( K^{-\delta+2\delta_1} +  \frac{\Nd}{N} \Nd^{-\delta_1+2\delta_3} K^{\delta_1+\delta_3} \Big) \Dc^3.
\end{equation}
In order to obtain \eqref{jacobi:eq-main-p6}, we now only use $K\lesssim \Nd \lesssim N$, which yields
\begin{equation*}
 \frac{\Nd}{N} \Nd^{-\delta_1+2\delta_3} K^{\delta_1+\delta_3}  \lesssim  \frac{\Nd}{N} \Nd^{3\delta_3} \lesssim \frac{\Nd}{N} N^{3\delta_3}.
\end{equation*}
\emph{Proof of \eqref{jacobi:eq-main-p7}:} The estimate \eqref{jacobi:eq-main-p7} can also be derived from \eqref{jacobi:eq-main-p10}. Indeed, if $K\lesssim N^{1-2\delta_1}$, then \eqref{jacobi:eq-main-p10} implies 
\begin{equation*}
\big\| \eqref{jacobi:eq-main-p5} \big\|_{\Cprod{r-1}{r-1}}  \lesssim 
 \Big( K^{-\delta+2\delta_1} +  \frac{\Nd}{N} \Nd^{-\delta_1+2\delta_3} N^{(1-2\delta_1)(\delta_1+\delta_3)} \Big) \Dc^3.
\end{equation*}
Since $\Nd\lesssim N$, it holds that 
\begin{equation*}
 \frac{\Nd}{N} \Nd^{-\delta_1+2\delta_3} N^{(1-2\delta_1)(\delta_1+\delta_3)} \lesssim N^{-\delta_1+2\delta_3}N^{(1-2\delta_1)(\delta_1+\delta_3)} \lesssim N^{-2\delta_1^2+2\delta_3}.
\end{equation*}
Since $\delta_1^2 \gg \delta_2$, this yields \eqref{jacobi:eq-main-p7}. \\

\emph{Proof of \eqref{jacobi:eq-main-p8}:} Using Definition \ref{jacobi:def-cf}, using Definition \ref{jacobi:def-chhl}, and inserting the expression for \eqref{jacobi:eq-main-p5}, the argument on the left-hand side of \eqref{jacobi:eq-main-p8} can be written as 
\begin{align}
&\eqref{jacobi:eq-main-p5}
- \chi \Big( P_{\leq N}^x \CHHLN\big( \VN[][], \VN[][] \big) P_{\leq N}^x - \coup \Renormbd[\Ncs] \Big) \UN[K][+] \notag \\
=&\, P_{\leq N}^x 
\int_{\R} \dy \, \bigg(  \big(\widecheck{\rho}_{\leq N} \ast \widecheck{\rho}_{\leq N}\big)(y)\, 
\Big( \HHLN[K,y][v] -  \HHLN[K,\bfzero][v]  - \chinull \big( \CHHLN_y\big( \VN[][],\VN[][] \big) - \coup \CNcsbd(y) \Kil \big) \Big) \notag \\
&\times \Theta^x_y P_{\leq N}^x \UN[K][+]   \bigg). \notag 
\end{align}
In order to estimate this term, we decompose the sum of high$\times$high$\rightarrow$low-matrices in the integrand as 
\begin{align}
 &\HHLN[K,y][v] -  \HHLN[K,\bfzero][v]  - \chinull \Big( \CHHLN_y\big( \VN[][],\VN[][] \big) - \coup \CNcsbd(y) \Kil \big) \Big)\notag \\
=& \,  \HHLN[K,y][v] -  \HHLN[K,\bfzero][v]  - \chinull[K]  P^{u,v}_{\leq K^\delta} \CHHLNc_{K,y} \label{jacobi:eq-main-p12} \\
+&\,  \big( \chinull[K] - \chinull \big) P^{u,v}_{\leq K^\delta} \CHHLNc_{K,y} \label{jacobi:eq-main-p13}  \\ 
+&\, \chinull \Big( P^{u,v}_{\leq K^\delta} \CHHLNc_{K,y} - \big(\CHHLN_y\big( \VN[][],\VN[][] \big) - \coup \CNcsbd(y) \Kil \big) \Big). 
\label{jacobi:eq-main-p14} 
\end{align}
We recall that $\HHLN[K,y][v]$, $\HHLN[K,\bfzero][v]$, and $P^{u,v}_{\leq K^\delta} \CHHLNc_{K,y}$ are supported on $u$-frequencies $\lesssim K^{\delta}$.
The contribution of \eqref{jacobi:eq-main-p12} to \eqref{jacobi:eq-main-p8} can then be controlled using the low$\times$high-estimate (Lemma \ref{prelim:lem-paraproduct}), Lemma \ref{modulation:lem-linear}, and Proposition \ref{jacobi:prop-shift}. The contribution of \eqref{jacobi:eq-main-p13} to \eqref{jacobi:eq-main-p8} can easily be controlled using the \mbox{smoothness of $\chinull$}, Lemma \ref{modulation:lem-linear}, and \eqref{jacobi:eq-Puv-Px-p1} from Lemma \ref{jacobi:lem-Puv-Px}. Finally, the contribution of \eqref{jacobi:eq-main-p14} to \eqref{jacobi:eq-main-p8} has been controlled in Proposition \ref{jacobi:prop-cartesian}.
\end{proof}

\section{Perturbative interactions}\label{section:null}

In this section we control all perturbative interactions from Definition \ref{ansatz:def-perturbative}. That is, we control $\PIN[(\ast_1)\times (\ast_2)]$ for all 
\begin{equation*}
\ast_1 \in \big\{ +, +-, - , +\fs, \fs-, \fs \big\} 
\qquad \text{and} \qquad
\ast_2 \in \big\{-, +-, +, \fs-, +\fs, \fs \big\}. 
\end{equation*}
In total, this corresponds to thirty-six different cases. Fortunately, the number of cases can be reduced rather quickly. As will be shown in Lemma \ref{null:lem-other}, all cases in which
\begin{equation}\label{null:eq-motivation-other}
\ast_1 \in \big\{  - , +\fs, \fs-, \fs \big\} 
\qquad \text{and} \qquad
\ast_2 \in \big\{ +, \fs-, +\fs, \fs \big\}
\end{equation}
can be controlled rather easily, which eliminates sixteen cases. We then utilize the  symmetry of our estimates in the $u$ and $v$-variables, which allows us to reduce to the diagonal and upper-diagonal cells in Figure \ref{figure:ansatz-cases}. Finally, we recall from Definition \ref{ansatz:def-perturbative} that $\PIN[(+)\times (-)]=0$, which eliminates the case $(\ast_1,\ast_2)=(+,-)$. In the end, this leaves us with the ten different cases that are illustrated in Figure \ref{figure:null-cases}. 
The main estimate of this section, which relies on this case analysis, is stated in the following proposition.

\begin{figure}
\begin{tabular}{
!{\vrule width 1pt}>{\centering\arraybackslash}P{1.5cm}
!{\vrule width 1pt}>{\centering\arraybackslash}P{\percolwidth}
!{\vrule width 1pt}>{\centering\arraybackslash}P{\percolwidth}
!{\vrule width 1pt}>{\centering\arraybackslash}P{\percolwidth}
!{\vrule width 1pt}>{\centering\arraybackslash}P{\percolwidth}
!{\vrule width 1pt}>{\centering\arraybackslash}P{\percolwidth}
!{\vrule width 1pt}>{\centering\arraybackslash}P{\percolwidth}
!{\vrule width 1pt}} 
\noalign{\hrule height 1pt} & & & & & & 
 \\[-5.8ex]
 \, \, \, \, \, \,   $V$ \vspace{-1ex} \newline \, \, \, \,  $U$  \hspace{-0.5ex}\vspace{-0.5ex} 
&  $(-)$ & $(+-)$ & $(+)$   & $(\fs-)$ & $(+\fs)$ & $(\fs)$  
\\[4pt] \noalign{\hrule height 1pt} \rule{0pt}{14pt}
$(+)$ 
& \cellcolor{Gray!80} 
& \cellcolor{Red!30} Lem. \ref{null:lem-p-pm} 
& \cellcolor{Yellow!30} Lem. \ref{null:lem-p-p} 
& \cellcolor{Red!30} Lem. \ref{null:lem-p-s}
& \cellcolor{Yellow!30} Lem. \ref{null:lem-p-p} 
& \cellcolor{Yellow!30}  Lem. \ref{null:lem-p-s}
\\[4pt] \noalign{\hrule height 1pt} \rule{0pt}{14pt}
$(+-)$ 
& \cellcolor{Gray!80} 
& \cellcolor{Green!30} Prop. \ref{killing:prop-quartic}
& \cellcolor{Green!30} Lem. \ref{null:lem-pm-p}
& \cellcolor{Green!30} Lem. \ref{null:lem-pm-sm}
& \cellcolor{Red!30}  Lem. \ref{null:lem-pm-s}
& \cellcolor{Red!30} Lem. \ref{null:lem-pm-s}
\\[3pt] \noalign{\hrule height 1pt} 
\end{tabular}
\caption{\small{
In this figure, we illustrate the ten different main cases in our estimates of the perturbative interactions. Our definition implies $\PIN[(+)\times (-)]=0$, and the corresponding cell has been colored in gray. Due to the symmetry in the $u$ and $v$-variables, the cell corresponding to the $(+-)$$\times$$(-)$-interaction has also been colored in gray. The color of the remaining cells agree with Figure \ref{figure:ansatz-cases} from Subsection \ref{section:ansatz-remainder}. 
}}
\label{figure:null-cases}
\end{figure}

\begin{proposition}[Perturbative interactions]\label{null:prop-perturbative}
Let the post-modulation hypothesis (Hypothesis \ref{hypothesis:post}) be satisfied. Then, it holds that
\begin{equation*}
\big\| \PIN \big\|_{\Cprod{r-1}{r-1}} \lesssim \Dc^2. 
\end{equation*}
\end{proposition}
The proof of Proposition \ref{null:prop-perturbative}, which consists only of references to other lemmas and propositions, is presented in Subsection \ref{section:null-perturbative} below.

\subsection{\texorpdfstring{\protect{Contributions of $\UN[][+]$}}{}} 
In this subsection, we state and prove all estimates which involve $\UN[][+]$ and therefore correspond to the first row in  Figure \ref{figure:null-cases}.

\begin{lemma}[The $(+)$$\times$$(+-)$-interaction]
\label{null:lem-p-pm}
Let the post-modulation hypothesis (Hypothesis \ref{hypothesis:post}) be satisfied and let $K,M_u,M_v\in \Dyadiclarge$ satisfy $\max(M_u,M_v)\geq K^{1-\delta}$ and $M_u \simeq_\delta M_v$. Furthermore, let $y,z\in \R$. Then, it holds that
\begin{equation*}
\big\| \Theta^x_y \UN[K][+] \Para[u][nsim] \Theta^x_z 
\VN[M_u,M_v][+-] \big\|_{\Cprod{r-1}{r-1}} 
\lesssim (KM_u M_v)^{-\eta} \Dc^3. 
\end{equation*}
\end{lemma}

\begin{proof}
Using our para-product estimate (Lemma \ref{prelim:lem-paraproduct}), Lemma \ref{modulation:lem-linear}, and Lemma \ref{modulation:lem-bilinear}, it holds that 
\begin{align*}
\big\| \Theta^x_y \UN[K][+] \Para[u][nsim] \Theta^x_z 
\VN[M_u,M_v][+-] \big\|_{\Cprod{r-1}{r-1}} 
&\lesssim \big\| \UN[K][+] \big\|_{\Cprod{r-1}{s}} 
\big\| \VN[M_u,M_v][+-] \big\|_{\Cprod{\eta}{r-1}} \\ 
&\lesssim K^{r-s} M_u^{\eta-s} M_v^{r-s} \Dc^3. 
\end{align*}
Since $\max(M_u,M_v)\geq K^{1-\delta}$ and $M_u \simeq_\delta M_v$, this is more than acceptable. 
\end{proof}

\begin{lemma}[The $(+)$$\times$$(+)$ and $(+)$$\times$$ (+\fs)$-interactions]\label{null:lem-p-p}
Assume that the post-modulation hypothesis (Hypothesis \ref{hypothesis:post}) is satisfied  and let $K,L\in \Dyadiclarge$. If $L\geq K^{1-\delta}$, then it holds that
\begin{align}
\Big\| \Big[ \UN[K][+] , \VN[L][+] \Big]_{\leq N} 
\Big\|_{\Cprod{r-1}{r-1}} 
&\lesssim (KL)^{-\eta} \Dc^3, 
\label{null:eq-p-p-e1} \\ 
\Big\| \Big[ \UN[K][+] , \VN[L][+\fs] \Big]_{\leq N}
- \chinull \coup \mathbf{1} \big\{ K=L \big\} \Renorm[\Ncs][K] P^v_{\geq K^{1-\deltap}} \VN[<K^{1-\delta}][\fcs] 
\Big\|_{\Cprod{r-1}{r-1}} 
&\lesssim (KL)^{-\eta} \Dc^3. 
\label{null:eq-p-p-e2} 
\end{align}
Alternatively, if $L<K^{1-\delta}$ and $y,z\in \R$, then
\begin{equation}\label{null:eq-p-p-e3}
\Big\| \Theta^x_y \UN[K][+] \Para[v][gg] \Theta^x_z P^v_{>K^{1-\deltap}} \VN[L][+\fs] \Big\|_{\Cprod{r-1}{r-1}}
\lesssim (KL)^{-\eta} \Dc^3. 
\end{equation}
\end{lemma}

\begin{proof}
We prove \eqref{null:eq-p-p-e1}, \eqref{null:eq-p-p-e2}, and \eqref{null:eq-p-p-e3} separately. \\

\emph{Step 1: Proof of \eqref{null:eq-p-p-e1}.}
From Remark \ref{ansatz:rem-commutativity-chipm} and Definition \ref{ansatz:def-modulated-linear-reversed}, we obtain that 
\begin{align}
&\, \Big[ \UN[K][+] , \VN[L][+] \Big]_{\leq N} \notag \\ 
=&\, \chinull[L] \Big[ \UN[K][+] , \Big[ \IUN[L][+], \LON[L][-] \Big]_{\leq N} \Big]_{\leq N}
+ \chinull[L] \Big[ \UN[K][+] , \Big[ (P_{\leq N}^x)^2 \IUN[L][+], \SHHLN[L][v] \Big]_{\leq N} \Big]_{\leq N}. \label{null:eq-p-p-decomposition}
\end{align}
We first treat the contribution of the $\LON[L][-]$-term from \eqref{null:eq-p-p-decomposition}. Using the integral representation of the $P_{\leq N}^x$-operator, it holds that
\begin{align}
 &\Big[ \UN[K][+] , \Big[ \IUN[L][+], \LON[L][-] \Big]_{\leq N} \Big]_{\leq N} \notag \\
 =&\, P_{\leq N}^x \int_{\R} \dy \, (\widecheck{\rho}_{\leq N} \ast \widecheck{\rho}_{\leq N})(y) \Big[ P_{\leq N}^x \UN[K][+], \Big[ \Theta^x_y P_{\leq N}^x \IUN[L][+], \Theta^x_y P_{\leq N}^x \LON[L][-]\Big]\Big] \notag \\
 =&  P_{\leq N}^x \int_{\R} \dy \, (\widecheck{\rho}_{\leq N} \ast \widecheck{\rho}_{\leq N})(y) \bigg( \Big[ P_{\leq N}^x \UN[K][+], \Big[ \Theta^x_y P_{\leq N}^x \IUN[L][+], \Theta^x_y P_{\leq N}^x \LON[L][-]\Big]\Big] \label{null:eq-p-p-nonresonant} \\
 &\hspace{29ex}- \coup \mathbf{1}\big\{ K=L \big\} \Cf^{(\Ncs)}_K(y) \Kil \big( \Theta^x_y P_{\leq N}^x \LON[L][-] \big) \bigg) \notag \\ 
 +& \coup \mathbf{1}\big\{ K=L \big\} P_{\leq N}^x \int_{\R} \dy \, (\widecheck{\rho}_{\leq N} \ast \widecheck{\rho}_{\leq N})(y) \Cf^{(\Ncs)}_K(y) \Kil \big( \Theta^x_y P_{\leq N}^x \LON[L][-] \big). \label{null:eq-p-p-resonant}
\end{align}
Using Proposition \ref{killing:prop-quadratic} and Corollary \ref{modulation:cor-LON}, the non-resonant term \eqref{null:eq-p-p-nonresonant} can be controlled by 
\begin{equation*}
\big\| \eqref{null:eq-p-p-nonresonant} \big\|_{\Cprod{r-1}{r-1}} 
\lesssim \Big( \max(K,L)^{r-\frac{1}{2}+\eta} L^{-\frac{1}{2}} + N^{-\delta+\vartheta} \Big) L^{r-s} \Dc^3. 
\end{equation*}
Since $L\gtrsim K^{1-\delta}$, this contribution is acceptable. Using Definition \ref{ansatz:def-Killing}, the resonant term \eqref{null:eq-p-p-resonant} can be written as $\coup \mathbf{1}\{ K=L\} \Renorm[N][K] \LON[L][-]$.
Due to Lemma \ref{ansatz:lem-frequency-support} and Definition \ref{ansatz:def-lo}, 
$\LON[L][-]$ is supported on frequencies $\lesssim L^{1-\delta}\lesssim N^{1-\delta}$. Using Lemma \ref{ansatz:lem-renormalization}, we then obtain that 
\begin{align*}
\coup \Big\|  \eqref{null:eq-p-p-resonant} \Big\|_{\Cprod{r-1}{r-1}} 
\lesssim \coup N^{r-s} N^{-\delta} 
 \Big\| \LON[L][-] \Big\|_{\Cprod{s}{s-1}} 
\lesssim \, \coup \Dc N^{r-s-\delta}. 
\end{align*}
Since $r-s=\delta_1+\delta_2\ll \delta$ and $\coup \leq \coup \Ac^2 \Bc^2\leq \Dc^2$, this is acceptable.  \\ 

The estimate for the $\SHHLN[L][v]$-term in \eqref{null:eq-p-p-decomposition} is similar. The only difference is that, due to the additional $(P_{\leq N}^x)^2$-operator, the covariance function $\Cf^{(\Ncs)}_K(y)$ in \eqref{null:eq-p-p-nonresonant}-\eqref{null:eq-p-p-resonant} has to be replaced by 
\begin{equation*}
\widebar{\Cf}^{(\Ncs)}_K(y) := \int_{\R} \dz (\widecheck{\rho}_{\leq N} \ast \widecheck{\rho}_{\leq N})(z) \Cf^{(\Ncs)}_K(y+z). 
\end{equation*}
Just like $\Cf^{(\Ncs)}_K$, however, $\widebar{\Cf}^{(\Nscript)}_K$ is odd and satisfies \eqref{ansatz:eq-renormalization-1}, and therefore this does not affect any of our estimates. \\

\emph{Step 2: Proof of \eqref{null:eq-p-p-e2}.} 
Using Lemma \ref{ansatz:lem-renormalization} and Corollary \ref{modulation:cor-control-combined}, we first estimate  
\begin{align*}
&\Big\| \big( \chinull[L] - \chinull \big)  \coup \mathbf{1} \big\{ K=L \big\} \Renorm[\Ncs][K] P^v_{\geq K^{1-\deltap}} \VN[<K^{1-\delta}][\fcs] 
\Big\|_{\Cprod{r-1}{r-1}} \\ 
\lesssim&\,  \coup \mathbf{1}\big\{ K=L \big\} L^{-100} \big\| \VN[<K^{1-\delta}][\fcs] \big\|_{\Cprod{s}{r-1}} 
\lesssim (KL)^{-\eta} \coup \Dc.
\end{align*}
Since $\coup \lesssim \Dc^2$, this is an acceptable contribution to \eqref{null:eq-p-p-e2}.
As a result, we can replace the $\chinull$-factor in \eqref{null:eq-p-p-e2} by $\chinull[L]$. 
We then recall from Definition \ref{ansatz:def-mixed} that
\begin{align}
\VN[L][+\fs] &= \chinull[L] \Big[ \IUN[L][+] \Para[v][ll] P^v_{\geq L^{1-\deltap}} \VN[<L^{1-\delta}][\fcs] \Big]_{\leq N} \notag \\
&=  \chinull[L] \Big[ \IUN[L][+] , P^v_{\geq L^{1-\deltap}} \VN[<L^{1-\delta}][\fcs] \Big]_{\leq N}
-  \chinull[L] \Big[ \IUN[L][+] \Para[v][gtrsim] P^v_{\geq L^{1-\deltap}} \VN[<L^{1-\delta}][\fcs] \Big]_{\leq N}.
\label{null:eq-p-p-q}
\end{align}
We first control the contribution of the first summand in \eqref{null:eq-p-p-q}. 
Using Definition \ref{ansatz:def-Killing}, Proposition \ref{killing:prop-quadratic} and Corollary \ref{modulation:cor-control-combined}, we obtain that  
\begin{equation}\label{null:eq-p-p-qq}
\begin{aligned}
&\Big\| \chinull[L]  \Big[ \UN[K][+],  \Big[ \IUN[L][+], P^v_{\geq L^{1-\deltap}} \VN[<L^{1-\delta}][\fcs] \Big]_{\leq N} \Big]_{\leq N} - \chinull[L] \coup \mathbf{1} \big\{ K=L \big\} \Renorm[\Ncs][K] P^v_{\geq K^{1-\deltap}} \VN[<K^{1-\delta}][\fcs] 
\Big\|_{\Cprod{r-1}{r-1}} \\
\lesssim&\, \int_{\R} \dy \bigg( \big| \big( \widecheck{\rho}_{\leq N}\ast \widecheck{\rho}_{\leq N} \big) \big|(y) 
\, \, \Big\| \Theta^x_y P^v_{\geq L^{1-\deltap}} \VN[<L^{1-\delta}][\fcs] \Big\|_{\Cprod{s}{r-1}} \\ 
&\hspace{2ex} \times  \sup_{\substack{E\in \frkg\colon \\ \| E \|_{\frkg}\leq 1}} 
\Big\|  \Big[ P_{\leq N}^x \UN[K][+],  \Big[ \Theta^x_y P_{\leq N}^x \IUN[L][+], E \Big] \Big]
- \coup \mathbf{1} \big\{ K=L \big\}  \Cf^{(\Ncs)}_K(y) \Kil(E) \Big\|_{\Cprod{r-1}{s}} \bigg)\\
\lesssim&\,  \big( \max(K,L)^{r-\frac{1}{2}+\eta} L^{-\frac{1}{2}} + N^{-\delta+\vartheta} \big) \Dc^3.
\end{aligned}
\end{equation}
Since $L\gtrsim K^{1-\delta}$, this is acceptable.  It therefore remains to treat the contribution of the second summand in \eqref{null:eq-p-p-q}. To this end, we first use the para-product estimate from Lemma \ref{prelim:lem-paraproduct}, which yields 
\begin{align}
&\Big\| \chinull[L] \Big[ \UN[K][+],  \Big[ \IUN[L][+] \Para[v][gtrsim] P^v_{\geq L^{1-\deltap}} \VN[<L^{1-\delta}][\fcs] \Big]_{\leq N} \Big]_{\leq N} \Big\|_{\Cprod{r-1}{r-1}} \notag \\ 
\lesssim&\, \big\| \UN[K][+] \big\|_{\Cprod{r-1}{s}} 
\Big\|  \Big[ \IUN[L][+] \Para[v][gtrsim] P^v_{\geq L^{1-\deltap}} \VN[<L^{1-\delta}][\fcs] \Big]_{\leq N} \Big\|_{\Cprod{1-r^\prime}{r-1}} 
\notag \\
\lesssim&\, \big\| \UN[K][+] \big\|_{\Cprod{r-1}{s}} 
\big\| \IUN[L][+]\big\|_{\Cprod{1-r^\prime}{s}} 
\big\| P^v_{\geq L^{1-\deltap}} \VN[<L^{1-\delta}][\fcs] \big\|_{\Cprod{1-r^\prime}{-s+\eta}}.
\label{null:eq-p-p-q2}
\end{align}
Using Lemma \ref{modulation:lem-linear} and Corollary \ref{modulation:cor-control-combined}, it then follows that
\begin{align}\label{null:eq-p-p-q3}
\eqref{null:eq-p-p-q2} 
\lesssim K^{r-s} L^{1-r^\prime-s} L^{(1-\deltap)(-s+\eta-(r-1))} \Dc^3. 
\end{align} 
Due to our choice of parameters, it holds that 
\begin{equation*}
1-r^\prime-s + (1-\deltap) \big(-s+\eta-(r-1)\big) = - \delta_1 + (1-\delta) ( - \delta_1 ) + \mathcal{O}(\delta_2) = - 2 \delta_1 + \mathcal{O}(\delta \delta_1). 
\end{equation*}
As a result, it follows that
\begin{equation*}
\eqref{null:eq-p-p-q3} \lesssim K^{\delta_1} L^{-2\delta_1} \big( KL \big)^{\mathcal{O}(\delta \delta_1)} \Dc^3.
\end{equation*}
Since $L\gtrsim K^{1-\delta}$, this is acceptable. \\

\emph{Step 3: Proof of \eqref{null:eq-p-p-e3}.} Using our para-product estimate (Lemma \ref{prelim:lem-paraproduct}), Lemma \ref{modulation:lem-linear}, and Lemma \ref{modulation:lem-mixed}, it holds that
\begin{align*}
\Big\| \Theta^x_y \UN[K][+] \Para[v][gg] \Theta^x_z P^v_{>K^{1-\deltap}} \VN[L][+\fs] \Big\|_{\Cprod{r-1}{r-1}} 
&\lesssim \big\|  \UN[K][+]  \big\|_{\Cprod{r-1}{s}}
\big\| P^v_{>K^{1-\deltap}} \VN[L][+\fs] \big\|_{\Cprod{1-r^\prime}{r-1-s}} \\
&\lesssim K^{r-s} L^{1-r^\prime-s} K^{-(1-\deltap) s} \Dc^3,
\end{align*}
which is more than acceptable. 
\end{proof}

\begin{lemma}[The $(+)$$\times$$(\fs-)$ and $(+)$$\times$$ (\fs)$-interactions]\label{null:lem-p-s}
Assume that the post-modulation hypothesis (Hypothesis \ref{hypothesis:post}) is satisfied and let $K,M\in \Dyadiclarge$. Furthermore, let $y,z\in \R$. Then, it holds that 
\begin{align}
\mathbf{1}\big\{ M < K^{1-\delta} \big\} 
 \Big\| \Theta^x_y \UN[K][+] 
\otimes \Theta^x_z P^u_{\geq K^{1-\deltap}} \VN[M][\fs-] 
\Big\|_{\Cprod{r-1}{r-1}} 
&\lesssim (MK)^{-\eta} \Dc^2, 
\label{null:eq-p-s-e1} \\ 
 \Big\| \Theta^x_y \UN[K][+] 
\otimes \Theta^x_z P^u_{\geq K^{1-\deltap}} \VN[][\fs] 
\Big\|_{\Cprod{r-1}{r-1}} 
&\lesssim K^{-\eta} \Dc^2, \label{null:eq-p-s-e2} \\ 
\mathbf{1}\big\{ M \geq K^{1-\delta} \big\} 
 \Big\| \Theta^x_y \UN[K][+] 
\Para[u][nsim] \Theta^x_z  \VN[M][\fs-] 
\Big\|_{\Cprod{r-1}{r-1}} 
&\lesssim (MK)^{-\eta} \Dc^2, \label{null:eq-p-s-e3} \\ 
 \Big\| \Theta^x_y \UN[K][+] 
\Para[v][gg] \Theta^x_z P^v_{\geq K^{1-\deltap}} P^u_{< K^{1-\deltap}} \VN[][\fs] 
\Big\|_{\Cprod{r-1}{r-1}} 
&\lesssim K^{-\eta} \Dc^2. 
\label{null:eq-p-s-e4}
\end{align}
\end{lemma}

The main ingredients in the proof of Lemma \ref{null:lem-p-s} are  para-product estimates (Lemma \ref{prelim:lem-paraproduct}). Compared with the proof of Lemma \ref{null:lem-p-p}, the argument relies much less on the properties of the modulated stochastic objects.

\begin{proof}[Proof of Lemma \ref{null:lem-p-s}]
We separate the proof into three main steps. \\

\emph{Step 1: Proof of \eqref{null:eq-p-s-e1} and \eqref{null:eq-p-s-e2}.} In order to treat \eqref{null:eq-p-s-e1} and \eqref{null:eq-p-s-e2}, we let $M<K^{1-\delta}$ and 
\begin{equation}\label{null:eq-p-s-p1}
\VN[][\ast] \in \Big\{ \VN[M][\fs-], \VN[][\fs] \Big\}. 
\end{equation}
In either case, it holds for all $0<\gamma\leq r$ that
\begin{equation}\label{null:eq-p-s-p2} 
\begin{aligned}
\big\| P^u_{>K^{1-\deltap}} \VN[][\ast] \big\|_{\Cprod{\gamma}{r-1}} 
&\lesssim K^{(1-\deltap)(\gamma-r)} \big\| P^u_{>K^{1-\deltap}} \VN[][\ast] \big\|_{\Cprod{r}{r-1}} \\ 
&\lesssim K^{(1-\deltap) (\gamma-r)} K^{(1-\delta)(r-s)} \Dc.
\end{aligned}
\end{equation}
We then decompose
\begin{equation}\label{null:eq-p-s-p3}
\begin{aligned}
&\Theta^x_y \UN[K][+] 
\otimes \Theta^x_z P^u_{\geq K^{1-\deltap}} \VN[][\ast]  \\
=\, & \Theta^x_y \UN[K][+] 
\Para[u][sim] \Theta^x_z P^u_{\geq K^{1-\deltap}} \VN[][\ast]  
+  \Theta^x_y \UN[K][+] 
\Para[u][nsim] \Theta^x_z P^u_{\geq K^{1-\deltap}} \VN[][\ast]. 
\end{aligned}
\end{equation}
For the resonant term in \eqref{null:eq-p-s-p3}, we use the high$\times$high-estimate (Lemma \ref{prelim:lem-paraproduct}), Lemma \ref{modulation:lem-linear}, and \eqref{null:eq-p-s-p2}, which yields that 
\begin{align*}
\Big\| \Theta^x_y \UN[K][+] 
\Para[u][sim] \Theta^x_z P^u_{\geq K^{1-\deltap}} \VN[][\ast] 
\Big\|_{\Cprod{r-1}{r-1}}  
\lesssim  \big\| \UN[K][+] \big\|_{\Cprod{-r^\prime}{s}} 
\big\| \VN[][\ast] \big\|_{\Cprod{r}{r-1}} 
\lesssim K^{1-r^\prime-s} K^{(1-\delta)(r-s)} \Dc^2. 
\end{align*}
Since $1-r^\prime-s+(1-\delta)(r-s)=-\delta \delta_1 + \mathcal{O}(\delta_2)$, this is acceptable. 
For the non-resonant term in \eqref{null:eq-p-s-p3}, we also use the para-product estimate (Lemma \ref{prelim:lem-paraproduct}), Lemma \ref{modulation:lem-linear}, and \eqref{null:eq-p-s-p2}, which yields that 
\begin{align*}
\Big\| \Theta^x_y \UN[K][+] 
\Para[u][nsim] \Theta^x_z P^u_{\geq K^{1-\deltap}} \VN[][\ast] 
\Big\|_{\Cprod{r-1}{r-1}} 
\lesssim\, & \big\| \UN[K][+] \big\|_{\Cprod{r-1}{s}} 
\big\| \VN[][\ast] \big\|_{\Cprod{\eta}{r-1}} \\
\lesssim\, & K^{r-s} K^{(1-\deltap) (\eta-r)} K^{(1-\delta)(r-s)} \Dc^2.
\end{align*}
Since 
\begin{equation*}
r-s+(1-\deltap) (\eta-r) + (1-\delta)(r-s) = - \frac{1}{2} + \mathcal{O}(\delta), 
\end{equation*}
this is more than acceptable. \\

\emph{Step 2: Proof of \eqref{null:eq-p-s-e3}.} Using the non-resonant estimate (Lemma \ref{prelim:lem-paraproduct}), Lemma \ref{modulation:lem-linear}, and Lemma \ref{modulation:lem-mixed}, we have that 
\begin{align*}
\Big\| \Theta^x_y \UN[K][+] 
\Para[u][nsim] \Theta^x_z  \VN[M][\fs-] 
\Big\|_{\Cprod{r-1}{r-1}} 
&\lesssim \, \big\| \UN[K][+] \big\|_{\Cprod{r-1}{s}} 
\big\| \VN[M][\fs-] \big\|_{\Cprod{\eta}{r-1}} \\
&\lesssim \, K^{r-s} M^{(1-\delta)(\eta-r)} M^{r-s} \Dc^2. 
\end{align*}
Since $M\geq K^{1-\deltap}$, this easily yields an acceptable contribution. \\ 

\emph{Step 3: Proof of \eqref{null:eq-p-s-e4}.} Using the high$\times$low-estimate (Lemma \ref{prelim:lem-paraproduct}) and Lemma \ref{modulation:lem-linear}, we obtain that 
\begin{align*}
&\Big\| \Theta^x_y \UN[K][+] 
\Para[v][gg] \Theta^x_z P^v_{\geq K^{1-\deltap}} P^u_{< K^{1-\deltap}} \VN[][\fs] 
\Big\|_{\Cprod{r-1}{r-1}} \\ 
\lesssim\, &\big\| \UN[K][+] \big\|_{\Cprod{r-1}{s}} 
\big\| P^v_{\geq K^{1-\deltap}} P^u_{< K^{1-\deltap}} \VN[][\fs] 
\big\|_{\Cprod{r}{r-1-s}} 
\lesssim K^{r-s} K^{-(1-\deltap)s} \Dc^2,
\end{align*}
which is more than acceptable. 
\end{proof}

\subsection{\texorpdfstring{\protect{Contributions of $\UN[][+-]$}}{}} 

In this subsection, we state and prove all estimates which involve $\UN[][+-]$. In Figure \ref{figure:null-cases}, this corresponds to the second row.

\begin{lemma}[The $(+-)$$\times$$(+)$-interaction] \label{null:lem-pm-p}
Let the post-modulation hypothesis (Hypothesis \ref{hypothesis:post}) be satisfied and let $K,L,M \in \Dyadiclarge$ satisfy $L\simeq_\delta M$. Then, it holds that
\begin{equation*}
\Big\| P_{\leq N}^x \UN[L,M][+-] \otimes P_{\leq N}^x \VN[K][+] \Big\|_{\Cprod{r-1}{r-1}} 
\lesssim (KLM)^{-\eta} \Dc^3. 
\end{equation*}
\end{lemma}

\begin{proof}
Using the definition of $\UN[L,M][+-]$, the definition of $\VN[K][+]$, the integral representation of $P_{\leq N}^x$, and Remark \ref{ansatz:rem-commutativity-chipm}, it holds that
\begin{align}
&\Big\| P_{\leq N}^x \UN[L,M][+-] \otimes P_{\leq N}^x \VN[K][+] \Big\|_{\Cprod{r-1}{r-1}} \notag \\ 
\lesssim&\, \int_{\R^4} \dy_1 \hdots \dy_4 \bigg(
\Big( \prod_{j=1}^4 N \langle N y_j \rangle^{-10} \Big) \notag \\
&\times \Big( 
\Big\| \Theta^x_{y_1} \UN[L][+] \otimes 
\Theta^x_{y_2} \IVN[M][-] \otimes
\Theta^x_{y_3} \IUN[K][+] \otimes 
\Theta^x_{y_4} \LON[K][-]  \Big\|_{\Cprod{r-1}{r-1}} 
\label{null:eq-pm-p-1}\\ 
&+\Big\| \Theta^x_{y_1} \UN[L][+] \otimes 
\Theta^x_{y_2} \IVN[M][-] \otimes
\Theta^x_{y_3} \IUN[K][+] \otimes 
\Theta^x_{y_4} \SHHLN[K][v] \Big\|_{\Cprod{r-1}{r-1}} 
\Big) \bigg) \label{null:eq-pm-p-2}. 
\end{align}
We now estimate \eqref{null:eq-pm-p-1} and \eqref{null:eq-pm-p-2} separately. \\

\emph{Estimate of \eqref{null:eq-pm-p-1}:} 
In order to control \eqref{null:eq-pm-p-1}, we distinguish the two cases $K\neq L$ and $K=L$. If $K\neq L$, we first use our paraproduct estimates (Lemma \ref{prelim:lem-paraproduct}), which yield that
\begin{align}
\eqref{null:eq-pm-p-1}
&\lesssim 
\Big\| \Theta^x_{y_1} \UN[L][+] \otimes \Theta^x_{y_3} \IUN[K][+] \Big\|_{\Cprod{r-1}{s}}
\Big\| \Theta^x_{y_2} \IVN[M][-] \Big\|_{\Cprod{s}{\eta}} 
\Big\| \LON[K][-] \Big\|_{\Cprod{s}{\eta}}. 
\label{null:eq-pm-p-3}
\end{align}
We estimate the first factor in \eqref{null:eq-pm-p-1} using Lemma \ref{killing:lem-tensor-modulated-linear} and the second factor using Lemma \ref{modulation:lem-linear}. To control the third factor,  we use that $ \LON[K][-] $ is supported on frequencies $\lesssim K^{1-\delta}$ (see Definition \ref{ansatz:def-lo} and Lemma \ref{ansatz:lem-frequency-support}) and Corollary \ref{modulation:cor-LON}. In total, this yields
\begin{equation*}
\eqref{null:eq-pm-p-3} 
\lesssim \max(K,L)^{r-\frac{1}{2}+\eta} K^{-\frac{1}{2}} M^{-\frac{1}{2}+2\eta} K^{(1-\delta)(\eta-(s-1))} \Dc^3. 
\end{equation*}
We then estimate $\max(K,L)\leq KL$ and insert the parameter expressions from Section \ref{section:parameters}, which yields that 
\begin{equation*}
\max(K,L)^{r-\frac{1}{2}+\eta} K^{-\frac{1}{2}} M^{-\frac{1}{2}+2\eta} K^{(1-\delta)(\eta-(s-1))} 
\lesssim K^{-\frac{\delta}{2}} M^{-\frac{1}{2}} (KLM)^{\mathcal{O}(\delta_1)}. 
\end{equation*}
Since $L\simeq_\delta M$, this is acceptable. It now remains to treat the case $K=L$. Using our paraproduct estimate from Lemma \ref{prelim:lem-paraproduct}, we first estimate
\begin{equation}
\eqref{null:eq-pm-p-1} 
\lesssim 
\Big\| \Theta^x_{y_1} \UN[L][+] \Big\|_{\Cprod{r-1}{s}}
\Big\| \Theta^x_{y_3} \IUN[K][+] \Big\|_{\Cprod{1-r^\prime}{s}} 
\Big\| \Theta^x_{y_2} \IVN[M][-] \otimes \Theta^x_{y_4} \LON[K][-] 
\Big\|_{\Cprod{s}{r-1}}. 
\label{null:eq-pm-p-4} 
\end{equation}
The first and second factor in \eqref{null:eq-pm-p-4} are estimated using Lemma \ref{modulation:lem-linear}. We now note that, since $K=L$ and $L\simeq_\delta M$, it holds that $K\simeq_\delta M$ and hence $M\geq K^{1-\delta}$. Due to this, the third factor in \eqref{null:eq-pm-p-4} can be estimated using Lemma \ref{modulation:lem-ivm-lo}. In total, it follows that
\begin{equation*}
\eqref{null:eq-pm-p-4} \lesssim L^{r-\frac{1}{2}+\eta} K^{\frac{1}{2}-r^\prime+\eta} M^{-\frac{1}{2}+\delta} \Dc^3. 
\end{equation*}
Since $K=L\simeq_\delta M$, this is easily seen to be acceptable. \\ 

\emph{Estimate of \eqref{null:eq-pm-p-2}:} 
We first use our paraproduct estimate (Lemma \ref{prelim:lem-paraproduct}), which yields that
\begin{align}
\eqref{null:eq-pm-p-2} 
&\lesssim \Big\|\Theta^x_{y_1} \UN[L][+] \Big\|_{\Cprod{r-1}{s}}
\Big\| \Theta^x_{y_2} \IVN[M][-] \Big\|_{\Cprod{s}{\eta}}
\Big\| \Theta^x_{y_3} \IUN[K][+] \Big\|_{\Cprod{s}{s}} 
\Big\| \Theta^x_{y_4} \SHHLN[K][v] \Big\|_{\Cprod{s}{\eta}}. 
\label{null:eq-pm-p-x}
\end{align}
We estimate the first three factors in \eqref{null:eq-pm-p-x} using Lemma \ref{modulation:lem-linear} and the last factor using Lemma \ref{modulation:lem-shhl}, which yields that
\begin{equation*}
\eqref{null:eq-pm-p-x} 
\lesssim L^{r-s} M^{-\frac{1}{2}+2\eta} K^{s-\frac{1}{2}+\eta} \Dc^3 = L^{\delta_1+\delta_2} M^{-\frac{1}{2}+\delta_3} K^{-\delta_2+\delta_3} \Dc^3. 
\end{equation*}
Since $L \simeq_\delta M$, the $M^{-\frac{1}{2}}$-factor easily makes up for the $L^{\delta_1+\delta_2}$-loss, and therefore this is acceptable. 
\end{proof}

\begin{lemma}[The $(+-)$$\times$$(\fs-)$-case]\label{null:lem-pm-sm} 
Let the post-modulation hypothesis (Hypothesis \ref{hypothesis:post}) be satisfied and let $K_u,K_v,M\in \Dyadiclarge$ satisfy $K_u \simeq_\delta K_v$. Furthermore, let $y,z\in \R$. Then,
\begin{equation}\label{null:eq-pm-sm}
 \Big\| \Theta^x_y \UN[K_u,K_v][+-] \otimes \Theta^x_z \VN[M][\fs-] \Big\|_{\Cprod{r-1}{r-1}} \lesssim (K_u K_v M)^{-\eta} \Dc^4. 
\end{equation}
\end{lemma}

\begin{proof}
In order to prove \eqref{null:eq-pm-sm}, we first decompose
\begin{equation}\label{null:eq-pm-sm-p1}
\begin{aligned}
&\Theta^x_y \UN[K_u,K_v][+-] \otimes \Theta^x_z \VN[M][\fs-]
=& \, \Theta^x_y \UN[K_u,K_v][+-] \Para[v][sim] \Theta^x_z \VN[M][\fs-]
+ \Theta^x_y \UN[K_u,K_v][+-] \Para[v][nsim] \Theta^x_z \VN[M][\fs-].
\end{aligned}
\end{equation}
We now estimate the resonant and non-resonant term in \eqref{null:eq-pm-sm-p1} separately. 
Due to Lemma \ref{ansatz:lem-frequency-support}, $\UN[K_u,K_v][+-]$ is supported on $v$-frequencies $\lesssim K_v$ and 
$\VN[M][\fs-]$ is supported on $v$-frequencies $\sim M$. Thus, the resonant term in \eqref{null:eq-pm-sm-p1} can only be non-zero if $K_v \gtrsim M$. Using our para-product estimate (Lemma \ref{prelim:lem-paraproduct}), Lemma \ref{modulation:lem-bilinear}, and Lemma \ref{modulation:lem-mixed}, the resonant term in \eqref{null:eq-pm-sm-p1} can be controlled by 
\begin{equation*}
\Big\| \Theta^x_y \UN[K_u,K_v][+-] \Para[v][sim] \Theta^x_z \VN[M][\fs-] \Big\|_{\Cprod{r-1}{r-1}} 
\lesssim \big\|  \UN[K_u,K_v][+-] \big\|_{\Cprod{r-1}{\eta}}
\big\|  \VN[M][\fs-] \big\|_{\Cprod{1-r^\prime}{\eta}}. 
\end{equation*}
Using Lemma \ref{modulation:lem-linear}, Lemma \ref{modulation:lem-bilinear}, $K_v \gtrsim M$, and $K_u \simeq_\delta K_v$, it follows that
\begin{align*}
\big\|  \UN[K_u,K_v][+-] \big\|_{\Cprod{r-1}{\eta}}
\big\|  \VN[M][\fs-] \big\|_{\Cprod{1-r^\prime}{\eta}}
&\lesssim K_u^{r-s} K_v^{\eta-s} 
M^{(1-\deltap)(1-r^\prime-r)} M^{\eta-(s-1)} \Dc^4\\ 
&\lesssim K_u^{\delta_2} K_v^{-\frac{1}{2}} M^{\frac{1}{2}-2\delta_2} (K_u K_v M)^{\mathcal{O}(\delta \delta_2)} \Dc^4. 
\end{align*}
Since $M\lesssim K_v$ and $K_u \simeq_\delta K_v$, it follows that
\begin{align*}
K_u^{\delta_2} K_v^{-\frac{1}{2}} M^{\frac{1}{2}-2\delta_2} (K_u K_v M)^{\mathcal{O}(\delta \delta_2)} 
&\lesssim K_u^{\delta_2} K_v^{-2\delta_2} (K_u K_v M)^{\mathcal{O}(\delta \delta_2)}\\ 
&\lesssim \max(K_u,K_v,M)^{-\delta_2} (K_u K_v M)^{\mathcal{O}(\delta \delta_2)},
\end{align*}
which is acceptable. 
For the non-resonant term in \eqref{null:eq-pm-sm-p1}, our para-product estimate (Lemma \ref{prelim:lem-paraproduct}) implies that
\begin{equation*}
\Big\| \Theta^x_y \UN[K_u,K_v][+-] \Para[v][nsim] \Theta^x_z \VN[M][\fs-] \Big\|_{\Cprod{r-1}{r-1}} 
\lesssim \big\| \UN[K_u,K_v][+-] \big\|_{\Cprod{r-1}{\eta}}
\big\| \VN[M][\fs-] \big\|_{\Cprod{1-r^\prime}{r-1}}. 
\end{equation*}
Using Lemma \ref{modulation:lem-bilinear}, it holds that \begin{align*}
\big\| \UN[K_u,K_v][+-] \big\|_{\Cprod{r-1}{\eta}}
&\lesssim K_u^{r-s} K_v^{\eta-s} \Dc^2 \lesssim (K_u K_v)^{-\eta}  \Dc^2 , \\
\big\| \VN[M][\fs-] \big\|_{\Cprod{1-r^\prime}{r-1}} 
&\lesssim M^{(1-\deltap)(1-r^\prime-r)} M^{r-s}  \Dc^2  \lesssim M^{-\eta}  \Dc^2 , 
\end{align*}
which are both acceptable. 
\end{proof}

\begin{lemma}[The $(+-)$$\times$$(+\fs)$ and $(+-)$$\times$$(\fs)$-interaction]\label{null:lem-pm-s}  
Assume that the post-modulation hypothesis (Hypothesis \ref{hypothesis:post}) is satisfied and let $K,M,L\in \Dyadiclarge$ satisfy $K \simeq_\delta M$. For all $y,z\in \R$, it then holds that 
\begin{align}
\mathbf{1} \big\{ L \geq K^{1-\delta} \big\} 
 \big\|
\Theta^x_y \UN[K,M][+-] 
\otimes \Theta^x_z \VN[L][+\fs] \big\|_{\Cprod{r-1}{r-1}} 
\lesssim& (KLM)^{-\eta} \Dc^3, \label{null:eq-pm-s-e1} \\ 
 \big\|
\Theta^x_y \UN[K,M][+-] 
\otimes \Theta^x_z P^{u}_{\geq K^{1-\deltap}} \VN[][\fs] \big\|_{\Cprod{r-1}{r-1}} 
\lesssim& (KM)^{-\eta} \Dc^3,  \label{null:eq-pm-s-e2} \\
\mathbf{1} \big\{ L < K^{1-\delta} \big\} 
 \big\|
\Theta^x_y \UN[K,M][+-] 
\Para[v][nsim] \Theta^x_z \VN[L][+\fs] \big\|_{\Cprod{r-1}{r-1}} 
\lesssim& (KLM)^{-\eta} \Dc^3, \label{null:eq-pm-s-e3} \\ 
 \big\|
\Theta^x_y \UN[K,M][+-] 
\Para[v][nsim] \Theta^x_z P^{u}_{< K^{1-\deltap}} \VN[][\fs] \big\|_{\Cprod{r-1}{r-1}} 
\lesssim& (KM)^{-\eta} \Dc^3. \label{null:eq-pm-s-e4}
\end{align}
\end{lemma}

\begin{proof} 
We separate the proof into two steps. \\

\emph{Step 1: Proof of \eqref{null:eq-pm-s-e1} and \eqref{null:eq-pm-s-e2}.} 
For any $L\geq K^{1-\delta}$ and 
$\VN[][\ast] \in \big\{ \VN[L][+\fs] , P^u_{\geq K^{1-\deltap}} \VN[][\fs] \big\}$, 
it holds that
\begin{equation}\label{null:eq-pm-s-p1}
\big\| \VN[][\ast] \big\|_{\Cprod{1-r^\prime}{r-1}} 
\lesssim K^{(1-\delta)(1-r^\prime-s)} \Gain\big( \VN[][\ast] \big)\Dc.
\end{equation}
Using the product estimate (Corollary \ref{prelim:cor-product}), Lemma \ref{modulation:lem-bilinear}, and \eqref{null:eq-pm-s-p1}, it then follows that 
\begin{equation*}
\begin{aligned}
\big\| \Theta^x_y \UN[K,M][+-] 
\otimes \Theta^x_z \VN[][\ast] \big\|_{\Cprod{r-1}{r-1}}  
&\lesssim \big\| \UN[K,M][+-] \big\|_{\Cprod{r-1}{1-r^\prime}}
\big\| \VN[][\ast] \big\|_{\Cprod{1-r^\prime}{r-1}} \\ 
&\lesssim K^{r-s} M^{1-r^\prime-s} K^{(1-\delta)(1-r^\prime-s)} \Gain\big( \VN[][\ast] \big) \Dc^{3}.
\end{aligned}
\end{equation*}
Since $K\simeq_\delta M$, this yields an acceptable contribution.\\

\emph{Step 2: Proof of \eqref{null:eq-pm-s-e3} and \eqref{null:eq-pm-s-e4}.} 
For any 
$\VN[][\ast] \in \big\{ \VN[L][+\fs] , P^{u}_{< K^{1-\deltap}} \VN[][\fs] \big\}$, it holds that
\begin{equation}\label{null:eq-pm-s-p2}
\big\| \VN[][\ast] \big\|_{\Cprod{1-r^\prime}{r-1}} 
\lesssim \Gain(\VN[][\ast]) \Dc. 
\end{equation}
Using the non-resonant para-product estimate (Lemma \ref{prelim:lem-paraproduct}), Lemma \ref{modulation:lem-bilinear}, and \eqref{null:eq-pm-s-p2}, we obtain that
\begin{align*}
    \big\|
\Theta^x_y \UN[K,M][+-] 
\Para[v][nsim] \Theta^x_z \VN[][\ast] \big\|_{\Cprod{r-1}{r-1}}
&\lesssim \big\| \UN[K,M][+-] \big\|_{\Cprod{r-1}{\eta}} 
\big\| \VN[][\ast] \big\|_{\Cprod{1-r^\prime}{r-1}} \\
&\lesssim K^{r-s} M^{\eta-s} 
\Gain(\VN[][\ast]) \Dc^{3}. 
\end{align*}
Together with $K\simeq_\delta M$, this yields the desired estimate. 
\end{proof}

\subsection{Other contributions}

In this subsection, we treat all cases listed in \eqref{null:eq-motivation-other}. As discussed at the beginning of this section, the cases are relatively easy and are not listed in Figure \ref{figure:null-cases}.

\begin{lemma}[Other contributions]\label{null:lem-other} 
Let the post-modulation hypothesis (Hypothesis \ref{hypothesis:post}) be satisfied and let $K,M \in \Dyadiclarge$. Furthermore, let $y,z\in \R$ and let
\begin{align*}
&\UN[][\ast_1] \in \Big\{ \UN[K][-], \UN[K][\fs-], \UN[K][+\fs], \UN[][\fs] \Big\} 
\quad \text{and} \quad \\ 
&\VN[][\ast_2] \in \Big\{ \VN[M][+], \VN[M][+\fs], \VN[M][\fs-], \VN[][\fs] \Big\}. 
\end{align*}
Then, it holds that 
\begin{equation}\label{null:eq-fs}
\begin{aligned}
\Big\| \Theta^x_y \UN[][\ast_1] \otimes \Theta^x_z \VN[][\ast_2] \Big\|_{\Cprod{r-1}{r-1}} 
\lesssim \Gain\big(  \UN[][\ast_1] \big) \Gain\big(  \VN[][\ast_2] \big). 
\end{aligned}
\end{equation}
\end{lemma}

\begin{proof}[Proof of Lemma \ref{null:lem-other}:]
We first use the product estimate (Corollary \ref{prelim:cor-product}), which yields
\begin{align*}
 \Big\| \Theta^x_y \UN[][\ast_1] \otimes \Theta^x_z \VN[][\ast_2] \Big\|_{\Cprod{r-1}{r-1}} 
&\lesssim  \Big\| \Theta^x_y \UN[][\ast_1] \Big\|_{\Cprod{r-1}{1-r^\prime}} 
\Big\|  \Theta^x_z \VN[][\ast_2] \Big\|_{\Cprod{r-1}{1-r^\prime}}  \\ 
&\lesssim \Big\|  \UN[][\ast_1] \Big\|_{\Cprod{r-1}{1-r^\prime}} 
 \Big\|   \VN[][\ast_2] \Big\|_{\Cprod{r-1}{1-r^\prime}}. 
\end{align*}
Using Lemma \ref{modulation:lem-mixed} and Lemma \ref{modulation:lem-linear-reversed}, we obtain\footnote{In cases in which $\ast_1\neq \fs$ or $\ast_2\neq \fs$, one could improve the power of $\Dc$, but this would not affect Proposition \ref{null:prop-perturbative}.} in all cases that
\begin{equation*}
\Big\|  \UN[][\ast_1] \Big\|_{\Cprod{r-1}{1-r^\prime}}  \lesssim \Gain\big(   \UN[][\ast_1] \big) \Dc
\quad \text{and} \quad \Big\|   \VN[][\ast_2] \Big\|_{\Cprod{r-1}{1-r^\prime}} \lesssim \Gain \big( \VN[][\ast_2] \big) \Dc,
\end{equation*}
which yields the desired estimate \eqref{null:eq-fs}. 
\end{proof}

\subsection{Proof of Proposition \ref{null:prop-perturbative}}\label{section:null-perturbative}

We are now ready to prove the main result of this section. 

\begin{proof}[Proof of Proposition \ref{null:prop-perturbative}:]
This follows directly from the definitions of the perturbative interactions (from Definition \ref{ansatz:def-perturbative}) and the estimates of this section. To be slightly more precise, we first decompose
\begin{equation*}
\PIN[] := \sum_{(\ast_1,\ast_2)} \PIN[(\,\ast_1)\times (\,\ast_2)]. 
\end{equation*}
If $\ast_1\in \{+,+-\}$ and $\ast_2\neq -$, $ \PIN[(\,\ast_1)\times (\,\ast_2)]$ can be controlled using the designated estimate from Figure \ref{figure:null-cases}. If $\ast_2 \in \{ -, +-\}$ and $\ast_1 \neq +$, the estimate follows by first using the symmetry in the $u$ and $v$-variable and then the designated estimate from Figure \ref{figure:null-cases}. Thus, it only remains to treat the cases in which $\ast_1 \neq +,+-$ and $\ast_2 \neq -,+-$, which can all be handled using Lemma \ref{null:lem-other}. 
\end{proof}

\section{Structural errors}\label{section:structural}

In this section we control all structural error terms, i.e., all error terms from Definition \ref{ansatz:def-structural-error}. While this has been previously discussed in Subsection \ref{section:ansatz-rigorous}, we  remind the reader of the origin of the structural errors. To this end, 
we briefly restrict ourselves to just one example and consider the modulated object $\UN[K,M][+-]$. The primary purpose of this object is to absorb the contribution of the  $(+)$$\times$$(-)$-interaction
\begin{equation}\label{structural:eq-motivation-0}
\chinull  \Big[ \UN[K][+], \VN[M][-] \Big]_{\leq N}. 
\end{equation}
In addition to absorbing \eqref{structural:eq-motivation-0}, however, we also want $\UN[K,M][+-]$ to have an expression which is as simple as possible. Instead of choosing $\UN[K,M][+-]$ as the $v$-integral of \eqref{structural:eq-motivation-0}, which would absorb \eqref{structural:eq-motivation-0} completely, 
we therefore chose the simpler expression
\begin{equation}\label{structural:eq-motivation-2}
\UN[K,M][+-]= \chinull[K,M] \Big[ \UN[K][+], \IVN[M][-] \Big]_{\leq N}.
\end{equation}
Due to this choice, the remainder equations then include the error term
\begin{equation*}
 \chinull \Big[ \UN[K][+], \VN[M][-] \Big]_{\leq N} - \partial_v \UN[K,M][+-],
\end{equation*}
which is the structural error term corresponding to $\UN[K,M][+-]$. Naturally, similar structural errors occur as a result of our definitions of other modulated and mixed modulated objects, and all the structural errors have already been introduced in Definition \ref{ansatz:def-structural-error}. 
The main estimate of this section shows that all structural errors are perturbative, i.e., in the weighted H\"{o}lder space $\WCprod{r-1}{r-1}$ from Definition \ref{prelim:def-weighted-hoelder}.

\begin{proposition}\label{structural:prop-main}
Let the post-modulation hypothesis (Hypothesis \ref{hypothesis:post}) be satisfied.  Then, it holds that
\begin{equation*}
\Big\|\,  \SEN[][u][] \Big\|_{\WCprod{r-1}{r-1}} 
+ \Big\|\, \SEN[][v][] \Big\|_{\WCprod{r-1}{r-1}} 
\lesssim \Dc^2. 
\end{equation*}
\end{proposition}

\begin{remark}
In Proposition \ref{ansatz:prop-decomposition}, the structural errors have no pre-factors of $\chinull$. In order to control the $\Int^v_{u\rightarrow v}$ and $\Int^u_{v\rightarrow u}$-integrals of the structural errors, we therefore need to control the 
$\WCprod{r-1}{r-1}$-norm from Definition \ref{prelim:def-weighted-hoelder} instead of the $\Cprod{r-1}{r-1}$-norm from Definition \ref{prelim:def-hoelder}.
\end{remark}

Due to the symmetry of our estimates in the $u$ and $v$-variables, we focus on the estimate of $\SEN[][u][]$. Due to Definition \ref{ansatz:def-structural-error}, this requires estimates of the five individual structural errors
\begin{equation*}
    \SEN[][u][+], \, \SEN[][u][+-], \, 
    \SEN[][u][-], \, \SEN[][u][+\fs], \quad  \text{and} \quad  
    \SEN[][u][\fs-]. 
\end{equation*}
The first lemma concerns $\SEN[][u][+]$, which is the structural error corresponding to the modulated linear wave $\UN[K][+]$. 

\begin{lemma}[\protect{Control of $\SEN[K][u][+]$}]\label{structural:lem-p}
Let the post-modulation hypothesis (Hypothesis \ref{hypothesis:post}) be satisfied and let $K\in \Dyadiclarge$. Then, it holds that
\begin{equation}\label{structural:eq-p}
\Big\| \SEN[K][u][+] \Big\|_{\WCprod{r-1}{r-1}} \lesssim K^{-\eta} \Dc^2. 
\end{equation}
\end{lemma}

The main ingredients in the following proof are the modulation equations (Definition \ref{ansatz:def-modulation-equations}) and our chaos estimate (Proposition \ref{chaos:prop-chaos}). 

\begin{proof}
Using the definition of $\SEN[K][u][+]$, we first write 
\begin{align}
\SEN[K][u][+]
&= \big( \chinull - \chinull[K] \big) \Big( \Big[ \UN[K][+], \LON[K][-] \Big]_{\leq N} 
+ \Big[  (P_{\leq N}^x)^2 \UN[K][+] , \SHHLN[K][v] \Big]_{\leq N} \Big) \label{structural:eq-p-q1} \\ 
&+  \chinull[K] \Big[ \UN[K][+], \LON[K][-] \Big]_{\leq N} 
+ \chinull[K] \Big[  (P_{\leq N}^x)^2 \UN[K][+] , \SHHLN[K][v] \Big]_{\leq N} - \partial_v \UN[K][+].  \label{structural:eq-p-q2} 
\end{align}
Using Lemma \ref{prelim:lem-weighted-hoelder-properties} and Lemma \ref{modulation:lem-linear}, Lemma \ref{modulation:lem-shhl}, and Corollary \ref{modulation:cor-LON}, we then easily obtain that
\begin{align*}
\big\| \eqref{structural:eq-p-q1} \big\|_{\WCprod{r-1}{r-1}} 
&\lesssim K^{-10} \Big( \Big\| \Big[ \UN[K][+], \LON[K][-] \Big]_{\leq N}  \Big\|_{\Cprod{r-1}{r-1}}  
+ \Big\|  \Big[  (P_{\leq N}^x)^2 \UN[K][+] , \SHHLN[K][v] \Big]_{\leq N} \Big\|_{\Cprod{r-1}{r-1}} \Big)\\ 
&\lesssim K^{-10} \big\| \UN[K][+] \big\|_{\Cprod{r-1}{s}} \Big( \big\| \LON[K][-] \big\|_{\Cprod{s}{r-1}} 
+ \big\| \SHHLN[K][v]\big\|_{\Cprod{s}{r-1}} \Big) \\ 
&\lesssim K^{-9} \Dc^2.
\end{align*}
Since this is an acceptable contribution to \eqref{structural:eq-p}, it remains to treat \eqref{structural:eq-p-q2}. 
Using the definition of $\UN[K][+]$ and Remark \ref{ansatz:rem-PNX-LON}, we then write 
\begin{align}
&\hspace{3ex} \eqref{structural:eq-p-q2} \notag \\ 
&= \chinull[K] \hcoup  \sum_{u_0\in \LambdaRR} \sum_{k\in \Z_K} 
P_{\leq N}^x \Big[ P_{\leq N}^x \Big( \psiRuK \rhoND(k) \SN[K][+][k] G^+_{u_0,k}  e^{iku} \Big), \LON[K][-] \Big] \notag \\
&+  \chinull[K] \hcoup  \sum_{u_0\in \LambdaRR} \sum_{k\in \Z_K} 
P_{\leq N}^x \Big[ (P_{\leq N}^x)^3 \Big( \psiRuK \rhoND(k)  \SN[K][+][k] G^+_{u_0,k}  e^{iku} \Big), \SHHLN[K][v] \Big] \notag \\
&- \hcoup  \sum_{u_0\in \LambdaRR} \sum_{k\in \Z_K}  \partial_v \Big( \psiRuK \rhoND(k)  \SN[K][+][k] G^+_{u_0,k}  e^{iku} \Big) \notag \\
&= \chinull[K] \hcoup  \sum_{u_0\in \LambdaRR} \sum_{k\in \Z_K} 
\Big( P_{\leq N}^x \Big[ P_{\leq N}^x \big( \psiRuK \rhoND(k)  \SN[K][+][k] G^+_{u_0,k} e^{iku}\big) , \LON[K][-] \Big]  
\label{structural:eq-p-q3} \\ 
&\hspace{23ex}-\psiRuK \rho_{\leq N}^2(k) \Big[ \SN[K][+][k] G^+_{u_0,k} e^{iku}, \LON[K][-] \Big] \Big) \notag\\ 
&+ \chinull[K] \hcoup  \sum_{u_0\in \LambdaRR} \sum_{k\in \Z_K} 
\Big( (P_{\leq N}^x) \Big[ (P_{\leq N}^x)^3 \big(  \psiRuK \rhoND(k)  \SN[K][+][k] G^+_{u_0,k} e^{iku} \big), \SHHLN[K][v] \Big] 
\label{structural:eq-p-q4} \\ 
&\hspace{23ex}- \psiRuK \rho_{\leq N}^4(k) \rhoND(k)  \Big[ \SN[K][+][k] G^+_{u_0,k} e^{iku}, \SHHLN[K][v] \Big] \Big)  
\notag \\ 
&+ \hcoup  \sum_{u_0\in \LambdaRR} \sum_{k\in \Z_K}  \Big( \psiRuK  \rhoND(k)  \Big( \Big[ \SN[K][+][k] G_{u_0,k}^+ e^{iku}, 
\rho_{\leq N}^2(k) \LON[K][-] + \rho_{\leq N}^4(k) \SHHLN[K][v] \Big]  \label{structural:eq-p-qq5} \\ 
&\hspace{34ex}- \partial_v \SN[K][+][k] G_{u_0,k}^+ e^{iku} \Big) \Big). \notag
\end{align}
The last term \eqref{structural:eq-p-qq5}  can easily be controlled using Lemma \ref{modulation:prop-main}.\ref{modulation:item-difference} and Lemma \ref{prelim:lem-psi-sum}. Since the estimates of \eqref{structural:eq-p-q3} and \eqref{structural:eq-p-q4} are similar, we only treat \eqref{structural:eq-p-q3}. Using the integral representation of $P_{\leq N}^x$, it holds that 
\begin{align}
\eqref{structural:eq-p-q3} 
=&\, \chinull[K] \hcoup \sum_{u_0\in \LambdaRR} \sum_{k\in \Z_K} \int_{\R^2} \dy_1\dy_2 \, \bigg( 
\widecheck{\rho}_{\leq N}(y_1) \widecheck{\rho}_{\leq N}(y_2) \rhoND(k) e^{ik (u+y_1+y_2)}  \notag \\
&\times 
\Big( \big(\Theta^x_{y_1+y_2} \psiRuK\big) \Big[ \big( \Theta^x_{y_1+y_2} \SN[K][+][k] \big) G^+_{u_0,k}, \Theta^x_{y_1} \LON[K][-] \Big] 
- \psiRuK \Big[ \SN[K][+][k] G_{u_0,k}^+, \LON[K][-] \Big] \Big) \bigg) \notag \allowdisplaybreaks[3]\\ 
=&\, \chinull[K] \hcoup \sum_{u_0\in \LambdaRR} \sum_{k\in \Z_K} \int_{\R^2} \dy_1\dy_2 \, \bigg( 
\widecheck{\rho}_{\leq N}(y_1) \widecheck{\rho}_{\leq N}(y_2) \rhoND(k) e^{ik (u+y_1+y_2)}  \notag \\
&\times 
 \Big( \Theta^x_{y_1+y_2} \psiRuK - \psiRuK \Big)  \Big[ \big( \Theta^x_{y_1+y_2} \SN[K][+][k] \big) G^+_{u_0,k}, \Theta^x_{y_1} \LON[K][-] \Big]  \bigg) 
 \allowdisplaybreaks[3 ]\label{structural:eq-p-q6}  \\ 
 +&\, \chinull[K] \hcoup \sum_{u_0\in \LambdaRR} \sum_{k\in \Z_K} \psiRuK \int_{\R^2} \dy_1\dy_2 \, \bigg( 
\widecheck{\rho}_{\leq N}(y_1) \widecheck{\rho}_{\leq N}(y_2)  \rhoND(k) e^{ik (u+y_1+y_2)}  \notag \\
&\times 
\Big(  \Big[ \big( \Theta^x_{y_1+y_2} \SN[K][+][k] \big) G^+_{u_0,k}, \Theta^x_{y_1} \LON[K][-] \Big] 
- \Big[ \SN[K][+][k] G_{u_0,k}^+, \LON[K][-] \Big] \Big) \bigg).
\allowdisplaybreaks[3] \label{structural:eq-p-q7}
\end{align}
Due to the smoothness and support properties of $(\psiRuK)_{u_0\in \LambdaRR}$, the estimate of \eqref{structural:eq-p-q6} is much easier than the estimate of \eqref{structural:eq-p-q7}, and we therefore focus on \eqref{structural:eq-p-q7}. By first using Lemma  \ref{prelim:lem-weighted-hoelder-properties} and then \mbox{Lemma \ref{prelim:lem-psi-sum}}, we obtain that
\begin{equation}\label{structural:eq-p-q8}
\begin{aligned}
\big\| \eqref{structural:eq-p-q7} \big\|_{\WCprod{r-1}{r-1}} 
\lesssim& \, \hcoup \sup_{u_0\in \LambdaRR} \bigg\| \sum_{k\in \Z_K}  \int_{\R^2} \dy_1\dy_2 \, \bigg( 
\widecheck{\rho}_{\leq N}(y_1) \widecheck{\rho}_{\leq N}(y_2) \rhoND(k) e^{ik (u+y_1+y_2)} \\
&\times 
\Big(  \Big[ \big( \Theta^x_{y_1+y_2} \SN[K][+][k] \big) G^+_{u_0,k}, \Theta^x_{y_1} \LON[K][-] \Big] 
- \Big[ \SN[K][+][k] G_{u_0,k}^+, \LON[K][-] \Big] \Big) \bigg) \bigg\|_{\Cprod{r-1}{r-1}}.
\end{aligned}
\end{equation}
In order to estimate \eqref{structural:eq-p-q8}, we first prove that 
\begin{equation}\label{structural:eq-p-p4}
\Big\| \Theta^x_{y_1+y_2} \SN[K][+][k] \otimes \Theta^x_{y_1} \LON[K][-] -  \SN[K][+][k] \otimes \LON[K][-] \Big\|_{\Wuv[s][r-1][k]} \lesssim K^{1-\delta+2\delta_1}\big( |y_1|+|y_2| \big) \Bc \Dc. 
\end{equation}
To obtain \eqref{structural:eq-p-p4}, we use the fundamental theorem of calculus and product estimates, which yield that
\begin{align}
&\big\|  \Theta^x_{y_1+y_2} \SN[K][+][k] \otimes \Theta^x_{y_1} \LON[K][-] -  \SN[K][+][k] \otimes \LON[K][-] \big\| _{\Wuv[s][r-1][k]} \notag \\
\lesssim\,& \max_{y\in \R} \big\|  \Theta^x_y (\partial_u + \partial_v) \SN[K][+][k] \otimes \Theta^x_{y_1} \LON[K][-] \big\| _{\Wuv[s][r-1][k]}
|y_1+y_2|
\notag \\ 
+\,& \max_{y\in \R} \big\|  \ \SN[K][+][k] \otimes \Theta^x_{y} (\partial_u+\partial_v) \LON[K][-] \big\| _{\Wuv[s][r-1][k]} |y_1|
\notag \\ 
\lesssim\, &  \big\|  (\partial_u + \partial_v) \SN[K][+][k]  \big\| _{\Wuv[s][s][k]} \big\|  \LON[K][-] \big\| _{\Cprod{s}{r-1}}
 |y_1+y_2| \label{structural:eq-p-p5} \\ 
+\, &  \big\| \SN[K][+][k]  \big\| _{\Wuv[s][s][k]}
\big\|  (\partial_u + \partial_v)  \LON[K][-] \big\|_{\Cprod{s}{r-1}}  
 |y_1| \label{structural:eq-p-p6}. 
\end{align}
Due to Lemma \ref{ansatz:lem-frequency-support}, $\SN[K][+][k]$ and $\LON[K][-]$ are supported on $u$ and $v$-frequencies $\lesssim K^{1-\delta+\vartheta}$. As a result, we obtain that
\begin{align*}
\eqref{structural:eq-p-p5} + \eqref{structural:eq-p-p6} 
&\lesssim K^{1-\delta+\vartheta} K^{(1-\delta+\vartheta)(r-s)} 
\big\| \SN[K][+][k]  \big\| _{\Wuv[s][s][k]}
\big\| \LON[K][-] \big\|_{\Cprod{s}{s-1}}   \big(  |y_1| +  |y_2| \big) \\
&\lesssim K^{1-\delta+2\delta_1} \big(  |y_1| +  |y_2| \big) \Bc \Dc. 
\end{align*}
This completes the proof of \eqref{structural:eq-p-p4}. By combining \eqref{structural:eq-p-p4} and \ref{ansatz:item-hypothesis-linear} from Hypothesis \ref{hypothesis:probabilistic}, it then follows that
\begin{align*}
\eqref{structural:eq-p-q8}
\lesssim\, & \hcoup \Ac K^{r-\frac{1}{2}+\eta}\int_{\R^2} \dy_1  \dy_2 \,\big| \widecheck{\rho}_{\leq N}(y_1) \widecheck{\rho}_{\leq N}(y_2) \big|
 \big\|  \Theta^x_{y_1+y_2} \SN[K][+][k] \otimes \Theta^x_{y_1} \LON[K][-] -  \SN[K][+][k] \otimes \LON[K][-] \big\| _{\Wuv[s][r-1][k]} \\
\lesssim\, & \hcoup \Ac \Bc \Dc K^{r-\frac{1}{2}+\eta} K^{1-\delta+2\delta_1}
\int_{\R^2} \dy_1  \dy_2 \, \big| \widecheck{\rho}_{\leq N}(y_1) \widecheck{\rho}_{\leq N}(y_2) \big| \big( \big|y_1\big| + \big| y_2\big| \big) \\ 
\lesssim\, & \Dc^2 K^{1-\delta+4\delta_1} N^{-1}. 
\end{align*}
Since $K\lesssim N$ and $\delta_1 \ll \delta$, this is acceptable. 
\end{proof}

\begin{lemma}[\protect{Control of $\SEN[K,M][u][+-]$}]\label{structural:lem-pm}
Let the post-modulation hypothesis (Hypothesis \ref{hypothesis:post}) be satisfied and let $K,M\in \Dyadiclarge$ satisfy $K\simeq_\delta M$. Then, it holds that
\begin{equation}\label{structural:eq-pm}
\Big\| \SEN[K,M][u][+-] \Big\|_{\WCprod{r-1}{r-1}} 
\lesssim (KM)^{-\eta} \Dc^2. 
\end{equation}
\end{lemma}

Out of all the estimates in this section, Lemma \ref{structural:lem-pm} is the most delicate. 
Fortunately, the main difficulties have already been addressed in Lemma \ref{modulation:lem-pm-para} above, which controls resonances in the $(+)$$\times$$(-)$-interaction. 
\begin{proof}
Using the definitions of $\SEN[K,M][u][+-]$ and $\UN[K,M][+-]$, it holds that
\begin{align}
\SEN[K,M][u][+-]
&= \chinull \Big[ \UN[K][+], \VN[M][-] \Big]_{\leq N}
- \partial_v \UN[K,M][+-]  \notag \\
&= \chinull \Big[ \UN[K][+], \VN[M][-] \Big]_{\leq N} - 
\partial_v \Big( \chinull[K,M] \Big[ \UN[K][+], \IVN[M][-] \Big]_{\leq N} \Big) \notag \\
&= \big( \chi - \chinull[K,M]\big) \Big[ \UN[K][+], \VN[M][-] \Big]_{\leq N} 
- \big( \partial_v \chinull[K,M] \big) \Big[ \UN[K][+], \IVN[M][-] \Big]_{\leq N} \label{structural:eq-pm-q1} \\
&+ \chinull[K,M] \Big[ \UN[K][+] \Para[v][ll] \Big( \VN[M][-]  - \partial_v \IVN[M][-] \Big) \Big]_{\leq N} \label{structural:eq-pm-p1}\\ 
&- \chinull[K,M] \Big[ \partial_v \UN[K][+] \Para[v][ll] \IVN[M][-] \Big]_{\leq M} \label{structural:eq-pm-q} \\
&+ \chinull[K,M] \Big[ \UN[K][+] \Para[v][gtrsim]  \VN[M][-]    \Big]_{\leq N} \label{structural:eq-pm-p2} \\ 
&- \chinull[K,M] \partial_v \Big[ \UN[K][+] \Para[v][gtrsim] \IVN[M][-] \Big]_{\leq N}
\label{structural:eq-pm-p3}. 
\end{align}
We now estimate \eqref{structural:eq-pm-q1}-\eqref{structural:eq-pm-p3} separately.
For \eqref{structural:eq-pm-q1}, we first obtain from Lemma \ref{prelim:lem-weighted-hoelder-properties} that
\begin{equation}\label{structural:eq-pm-q2}
\big\| \eqref{structural:eq-pm-q1} \big\|_{\WCprod{r-1}{r-1}} 
\lesssim (KM)^{-100} \Big\| \Big[ \UN[K][+], \VN[M][-] \Big]_{\leq N} \Big\|_{\Cprod{r-1}{r-1}} 
+ \Big\| \Big[ \UN[K][+], \IVN[M][-] \Big]_{\leq N} \Big\|_{\Cprod{r-1}{r-1}}.
\end{equation}
Using our product estimate (Corollary \ref{prelim:cor-product}) and Lemma \ref{modulation:lem-linear}, we then 
obtain that
\begin{equation*}
\eqref{structural:eq-pm-q2} 
\lesssim (KM)^{-100} \big\| \UN[K][+] \big\|_{\Cprod{r-1}{s}} \big\| \VN[M][-] \big\|_{\Cprod{s}{r-1}} 
+ \big\| \UN[K][+] \big\|_{\Cprod{r-1}{s}} \big\| \IVN[M][-] \big\|_{\Cprod{s}{r-1}} 
\lesssim (KM)^{-\eta} \Dc^2.
\end{equation*}

For \eqref{structural:eq-pm-p1}, we obtain from Lemma \ref{prelim:lem-weighted-hoelder-properties}, Corollary \ref{prelim:cor-product}, and Lemma \ref{modulation:lem-integration} that
\begin{equation*}
\big\| \eqref{structural:eq-pm-p1} \big\|_{\WCprod{r-1}{r-1}} 
\lesssim \big\| \UN[K][+] \big\|_{\Cprod{r-1}{s}} 
\big\|  \VN[M][-]  - \partial_v \IVN[M][-] \big\|_{\Cprod{s}{r-1}} 
\lesssim K^{r-s} M^{-1/2} \Dc^2.
\end{equation*}
Since $K\simeq_\delta M$, this yields an acceptable contribution. For \eqref{structural:eq-pm-q}, we obtain from Lemma \ref{prelim:lem-weighted-hoelder-properties}, Lemma \ref{prelim:lem-paraproduct}, and Lemma \ref{modulation:lem-linear} that
\begin{equation*}
\Big\| \eqref{structural:eq-pm-q} \Big\|_{\WCprod{r-1}{r-1}} 
\lesssim \Big\| \partial_v \UN[K][+] \Big\|_{\Cprod{r-1}{\eta}} 
\Big\| \IVN[M][-] \Big\|_{\Cprod{s}{r-1}} \lesssim K^{r-s} K^{1+\eta-s} M^{r-1-s} \Dc^2.
\end{equation*}
Since $K\simeq_\delta M$, this easily yields an acceptable contribution. 
For \eqref{structural:eq-pm-p2}, the estimate follows from Lemma \ref{prelim:lem-weighted-hoelder-properties} and Lemma \ref{modulation:lem-pm-para}. 
For \eqref{structural:eq-pm-p3}, we obtain from Lemma \ref{prelim:lem-weighted-hoelder-properties}, our para-product estimate (Lemma \ref{prelim:lem-paraproduct}), and Lemma \ref{modulation:lem-linear} that
\begin{equation*}
\begin{aligned}
\Big\| \eqref{structural:eq-pm-p3} \Big\|_{\WCprod{r-1}{r-1}} 
&\lesssim \Big\|  \Big[ \UN[K][+] \Para[v][gtrsim] \IVN[M][-] \Big]_{\leq N} \Big\|_{\Cprod{r-1}{r}} \\ 
&\lesssim \big\| \UN[K][+] \big\|_{\Cprod{r-1}{r}} 
\big\| \IVN[M][-] \big\|_{\Cprod{s}{\eta}} 
\lesssim (K^{r-s})^2 M^{\eta-s} \Dc^2.
\end{aligned}
\end{equation*}
Since $K\simeq_\delta M$, this is acceptable. 
\end{proof}

\begin{lemma}[\protect{Control of $\SEN[M][u][-]$}]\label{structural:lem-m}
Let the post-modulation hypothesis (Hypothesis \ref{hypothesis:post}) be satisfied and let $M\in \Dyadiclarge$. Then, it holds that
\begin{equation}\label{structural:eq-m}
\Big\| \SEN[M][u][-] \Big\|_{\WCprod{r-1}{r-1}} 
\lesssim M^{-1/2+\delta} \Dc^2.
\end{equation}
\end{lemma}
The main ingredients in the following proof are Lemma \ref{modulation:lem-integration} and our para-product estimates (Lemma \ref{prelim:lem-paraproduct}). 

\begin{proof} 
Using the definitions of $\SEN[M][u][-]$ and $\UN[M][-]$, it follows that
\begin{align}
\SEN[M][u][-]  
 &=  \chinull  \Big[ \LON[M][+], \VN[M][-] \Big]_{\leq N} 
+ \chinull  \Big[ \SHHLN[M][u],  (P_{\leq N}^x)^2 \VN[M][-] \Big]_{\leq N}  \notag \\
&- \partial_v  \Big( \chinull[M] \Big[ \LON[M][+], \IVN[M][-] \Big]_{\leq N} \Big) 
- \partial_v  \Big( \chinull[M] \Big[ \SHHLN[M][u],  (P_{\leq N}^x)^2 \IVN[M][-] \Big]_{\leq N} \Big)\notag \\
&= \big( \chinull - \chinull[M] \big) \Big( \Big[ \LON[M][+], \VN[M][-] \Big]_{\leq N} + \Big[ \SHHLN[M][u], (P_{\leq N}^x)^2 \VN[M][-] \Big]_{\leq N} \Big)  \label{structural:eq-m-q1} \\
&- \big( \partial_v \chinull[M] \big) \Big(  \Big[ \LON[M][+], \IVN[M][-] \Big]_{\leq N} 
+ \Big[ \SHHLN[M][u], (P_{\leq N}^x)^2 \IVN[M][-] \Big]_{\leq N}\Big) \label{structural:eq-m-q2} \\ 
&+ \chinull[M]  \Big[ \LON[M][+],  \big( \VN[M][-] - \partial_v \IVN[M][-] \big) \Big]_{\leq N} \label{structural:eq-m-p1} \\ 
&+ \chinull[M]  \Big[ \SHHLN[M][u], (P_{\leq N}^x)^2 \big( \VN[M][-] - \partial_v \IVN[M][-] \big) \Big]_{\leq N} \label{structural:eq-m-p2} \\ 
&+  \chinull[M]  \Big[ \partial_v \LON[M][+],  \IVN[M][-] \Big]_{\leq N} \label{structural:eq-m-p3} \\ 
&+ \chinull[M]   \Big[ \partial_v \SHHLN[M][u],  (P_{\leq N}^x)^2 \IVN[M][-] \Big]_{\leq N}  \label{structural:eq-m-p4}. 
\end{align}
In the following argument, we repeatedly use Lemma \ref{prelim:lem-weighted-hoelder-properties} to move from $\Cprod{r-1}{r-1}$ to $\WCprod{r-1}{r-1}$ without explicit reference. 
Using Lemma \ref{modulation:lem-linear}, Lemma \ref{modulation:lem-shhl}, and Corollary \ref{modulation:cor-LON}, it holds that
\begin{align*}
&\big\| \eqref{structural:eq-m-q1}\big\|_{\WCprod{r-1}{r-1}} + \big\| \eqref{structural:eq-m-q2}\big\|_{\WCprod{r-1}{r-1}} \\
\lesssim&\, \Big(  \big\| \LON[M][+] \big\|_{\Cprod{r-1}{s}} 
+ \big\| \SHHLN[M][u] \big\|_{\Cprod{r-1}{s}} \Big) 
\Big( M^{-100} \big\| \VN[M][-] \big\|_{\Cprod{s}{r-1}} + \big\| \IVN[M] \big\|_{\Cprod{s}{r-1}} \Big) \\ 
\lesssim& \, M^{r-\frac{3}{2}+\eta} \Dc^2, 
\end{align*}
which is more than acceptable. 
Using Lemma \ref{modulation:lem-integration}, Lemma  \ref{modulation:lem-shhl}, and Corollary \ref{modulation:cor-LON},  we have that
\begin{align*}
&\big\| \eqref{structural:eq-m-p1} \big\|_{\WCprod{r-1}{r-1}} 
+ \big\| \eqref{structural:eq-m-p2} \big\|_{\WCprod{r-1}{r-1}} \\ 
\lesssim\, & \Big( \big\| \LON[M][+] \big\|_{\Cprod{r-1}{s}} 
+ \big\| \SHHLN[M][u] \big\|_{\Cprod{r-1}{s}} \Big) 
\big\|  \VN[M][-] - \partial_v \IVN[M][-] \big\|_{\Cprod{r}{r-1}} 
\lesssim M^{r-s} M^{-1/2} \Dc^2,
\end{align*}
which is acceptable. Due to Lemma \ref{ansatz:lem-frequency-support} and the definitions of $\LON[M][+]$ and $\SHHLN[M][u]$, \eqref{structural:eq-m-p3} and \eqref{structural:eq-m-p4} only contain low$\times$high-interactions in the $v$-variable. As a result, our para-product estimate (Lemma \ref{prelim:lem-paraproduct}), Lemma \ref{modulation:lem-shhl}, and Corollary \ref{modulation:cor-LON} yield that
\begin{align*}
   \big\| \eqref{structural:eq-m-p3} \big\|_{\WCprod{r-1}{r-1}} 
+ \big\| \eqref{structural:eq-m-p4} \big\|_{\WCprod{r-1}{r-1}} 
\lesssim\, & \Big( \big\| \partial_v \LON[M][+] \big\|_{\Cprod{r-1}{\eta}} 
+ \big\| \partial_v \SHHLN[M][u] \big\|_{\Cprod{r-1}{\eta}} \Big) 
\big\|  \IVN[M][-] \big\|_{\Cprod{s}{r-1}}   \\ 
\lesssim\, & M^{r-s} M^{1+\eta-s} M^{r-1-s} \Dc^2 \lesssim M^{-1/2+\delta} \Dc^2,
\end{align*}
which is acceptable. 
\end{proof}

\begin{lemma}[\protect{Control of $\SEN[K][u][+\fs]$}]\label{structural:lem-ps}
Let the post-modulation hypothesis (Hypothesis \ref{hypothesis:post}) be satisfied and let $K\in \Dyadiclarge$. Then, it holds that
\begin{equation}\label{structural:eq-ps}
\Big\| \SEN[K][u][+\fs] \Big\|_{\Cprod{r-1}{r-1}} 
\lesssim K^{-1/2+\delta} \Dc^2. 
\end{equation}
\end{lemma}

The following proof is rather simple and only utilizes our paraproduct estimates (Lemma \ref{prelim:lem-paraproduct}). 

\begin{proof}
Using the definition of $\SEN[K][u][+\fs]$ and $\UN[K][+\fs]$, it holds that
\begin{align}
&\,\SEN[K][u][+\fs] \notag \\ 
=&\, \chinull \Big[ \UN[K][+] \Para[v][ll] P^v_{\geq K^{1-\deltap}} \VN[<K^{1-\delta}][\fcs] \Big]_{\leq N}
- \partial_v \Big( \chinull[K] \Big[ \UN[K][+] \Para[v][ll] \Int^v_{u \rightarrow v} P^v_{\geq K^{1-\deltap}} \VN[<K^{1-\delta}][\fcs] \Big]_{\leq N} \Big) \notag \\ 
=&\, \big( \chinull - \chinull[K] \big) \Big[ \UN[K][+] \Para[v][ll] P^v_{\geq K^{1-\deltap}} \VN[<K^{1-\delta}][\fcs] \Big]_{\leq N} 
- \big( \partial_v \chinull[K] \big) \Big[ \UN[K][+] \Para[v][ll] \Int^v_{u \rightarrow v} P^v_{\geq K^{1-\deltap}} \VN[<K^{1-\delta}][\fcs] \Big]_{\leq N} \label{structural:eq-ps-q1} \\
-&\, \chinull[K] \Big[ \partial_v \UN[K][+] \Para[v][ll] \Int^v_{u \rightarrow v} P^v_{\geq K^{1-\deltap}} \VN[<K^{1-\delta}][\fcs] \Big]_{\leq N} . 
\label{structural:eq-ps-p1}
\end{align}
The first term \eqref{structural:eq-ps-q1} can be controlled similarly as \eqref{structural:eq-pm-q1} in the proof of Lemma \ref{structural:lem-pm} and we omit the details. 
Using the low$\times$high-paraproduct estimate (Lemma \ref{prelim:lem-paraproduct}), Lemma \ref{prelim:lem-Duhamel-integral}, Lemma \ref{modulation:lem-linear}, and Corollary \ref{modulation:cor-control-combined}, we obtain that
\begin{align*}
\big\| \eqref{structural:eq-ps-p1} \big\|_{\WCprod{r-1}{r-1}} 
&\lesssim \Big\|  \partial_v \UN[K][+] \Big\|_{\Cprod{r-1}{s-1}} 
\Big\| \Int^v_{u \rightarrow v} P^v_{\geq K^{1-\deltap}} \VN[<K^{1-\delta}][\fcs] \Big\|_{\Cprod{s}{r-s+\eta}} \\ 
&\lesssim K^{r-s} K^{(1-\deltap)(-s+\eta)}
 \Big\|  \partial_v \UN[K][+] \Big\|_{\Cprod{s-1}{s}} 
\Big\|  \VN[<K^{1-\delta}][\fcs]  \Big\|_{\Cprod{s}{r-1}} \\
&\lesssim K^{r-s} K^{(1-\deltap)(-s+\eta)} \Dc^2, 
\end{align*}
which is acceptable. 
\end{proof}

\begin{lemma}[\protect{Control of $\SEN[M][u][\fs-]$}]\label{structural:lem-sm}
Let the post-modulation hypothesis (Hypothesis \ref{hypothesis:post}) be satisfied and let $M\in \Dyadiclarge$. Then, it holds that
\begin{equation}\label{structural:eq-sm}
\Big\| \SEN[M][u][\fs-] \Big\|_{\WCprod{r-1}{r-1}} 
\lesssim M^{-1/2+\delta} \Dc^2. 
\end{equation}
\end{lemma}

The main ingredients in the following proof are Lemma \ref{modulation:lem-integration} and our para-product estimates (Lemma \ref{prelim:lem-paraproduct}). 

\begin{proof}
Using the definitions of $\SEN[M][u][\fs-]$ and $\UN[M][\fs-]$, it holds that 
\begin{align}
&\SEN[M][u][\fs-] \notag \\ 
=&\,  \chinull \Big[ P^u_{\geq M^{1-\deltap}} \UN[<M^{1-\delta}][\fsc] \Para[u][gg] \VN[M][-] \Big]_{\leq N}
- \partial_v  \Big( \chinull[M] \Big[ P^u_{\geq M^{1-\deltap}} \UN[<M^{1-\delta}][\fsc] \Para[u][gg] \IVN[M][-] \Big]_{\leq N} \Big) \notag \\
=&\, \big( \chinull - \chinull[M] \big) \Big[ P^u_{\geq M^{1-\deltap}} \UN[<M^{1-\delta}][\fsc] \Para[u][gg] \VN[M][-] \Big]_{\leq N}
- \big( \partial_v \chinull[M] \big) \Big[ P^u_{\geq M^{1-\deltap}} \UN[<M^{1-\delta}][\fsc] \Para[u][gg] \IVN[M][-] \Big]_{\leq N} 
\label{structural:eq-sm-q1} \\ 
+&\,  \Big[ P^u_{\geq M^{1-\deltap}} \UN[<M^{1-\delta}][\fsc] \Para[u][gg] \Big( \VN[M][-] - \partial_v \IVN[M][-] \Big) \Big]_{\leq N} \label{structural:eq-sm-p1}\\
-&\, \Big[ P^u_{\geq M^{1-\deltap}} \partial_v \UN[<M^{1-\delta}][\fsc] \Para[u][gg] \IVN[M][-] \Big]_{\leq N}  \label{structural:eq-sm-p2}. 
\end{align}
The first term \eqref{structural:eq-sm-q1} can be controlled similarly as \eqref{structural:eq-pm-q1} in the proof of Lemma \ref{structural:lem-pm} and we omit the details. 
We now estimate \eqref{structural:eq-sm-p1} and \eqref{structural:eq-sm-p2} separately. Using our para-product estimate (Lemma \ref{prelim:lem-paraproduct}), Lemma \ref{modulation:lem-integration}, and Corollary \ref{modulation:cor-control-combined}, we have that
\begin{equation*}
\big\| \eqref{structural:eq-sm-p1} \big\|_{\WCprod{r-1}{r-1}}
\lesssim \big\|  \UN[<M^{1-\delta}][\fsc] \big\|_{\Cprod{r-1}{s}}
\big\|  \VN[M][-] - \partial_v \IVN[M][-]\big\|_{\Cprod{r}{r-1}} 
\lesssim M^{-1/2} \Dc^2,
\end{equation*}
which is acceptable. Due to Lemma \ref{ansatz:lem-frequency-support}, \eqref{structural:eq-sm-p2} only contains low$\times$high-interactions in the $v$-variable. Thus, our para-product estimate (Lemma \ref{prelim:lem-paraproduct}) yields that 
\begin{align*}
\big\| \eqref{structural:eq-sm-p2} \big\|_{\WCprod{r-1}{r-1}}
&\lesssim \big\| \partial_v \UN[<M^{1-\delta}][\fsc] \big\|_{\Cprod{r-1}{\eta}} \big\| \IVN[M][-] \big\|_{\Cprod{s}{r-1}}  \\
&\lesssim M^{1+\eta-s} M^{r-1-s} 
\big\| \UN[<M^{1-\delta}][\fsc] \big\|_{\Cprod{r-1}{s}} \big\| \IVN[M][-] \big\|_{\Cprod{s}{s}} \\
&\lesssim M^{1+\eta-s} M^{r-1-s} \Dc^2, 
\end{align*}
which is acceptable. 
\end{proof}

Equipped with our earlier lemmas, we can now prove the main proposition of this section. 

\begin{proof}[Proof of Proposition \ref{structural:prop-main}]
Due to the symmetry of our estimates in the $u$ and $v$-variables, it suffices to control $\SEN[][u][]$. Due to Definition \ref{ansatz:def-structural-error}, it then suffices to control the terms
\begin{equation}\label{structural:eq-main-p1}
\SEN[][u][+], \, \SEN[][u][+-], \, 
    \SEN[][u][-], \, \SEN[][u][+\fs], \quad  \text{and} \quad  
    \SEN[][u][\fs-],
\end{equation}
which can be done using Lemma \ref{structural:lem-p}, Lemma \ref{structural:lem-pm}, Lemma \ref{structural:lem-m}, Lemma \ref{structural:lem-ps}, and Lemma \ref{structural:lem-sm}. 
\end{proof}

\section{Renormalization errors} \label{section:renormalization}

In this section, we control all renormalization errors, i.e., all errors from Definition \ref{ansatz:def-renormalization-error}. Our main estimate is recorded in the following proposition. 

\begin{proposition}[Renormalization errors]\label{ren:prop-main}
Let the post-modulation hypothesis (Hypothesis \ref{hypothesis:post}) be satisfied. Then, it holds that 
\begin{equation}\label{ren:eq-main}
\Big\| \, \REN \Big\|_{\Cprod{r-1}{r-1}} \lesssim N^{-\eta} \coup \Dc \lesssim N^{-\eta} \Dc^3. 
\end{equation}
\end{proposition}

In comparison with estimates of the high$\times$high$\rightarrow$low-errors or Jacobi errors, the estimates of the renormalization errors are rather simple. In fact, the proof of Proposition \ref{ren:prop-main} only relies on the frequency-support properties of the modulated and mixed modulated objects and Lemma \ref{ansatz:lem-renormalization}. 

\begin{proof}[Proof of Proposition \ref{ren:prop-main}:]
Due to the symmetry of our estimates in the $u$ and $v$-variables, it suffices to estimate $\RenNErr[][u]$. Due to Definition \ref{ansatz:def-renormalization-error}, it then remains to bound the six terms
\begin{equation*}
\RenNErr[][u,+], \, \RenNErr[][u,+-], \, \RenNErr[][u,-], \, \RenNErr[][u,+\fs], \, \RenNErr[][u,\fs-], \quad \text{and} \quad \RenNErr[][u,\fs],
\end{equation*}
which are addressed in four separate steps. \\ 

\emph{Estimate of $\RenNErr[][u,+]$}: From Definition \ref{ansatz:def-renormalization-error}, it follows that
\begin{equation*}
\RenNErr[][u,+] = 
-  \coup \Sumlarge_{\substack{K,M\leq \Nd\colon \\ K  \simeq_\delta M}} \big( \chinull- \chinull[K] \big) \Renorm[\Ncs][M] \UN[K][+]
-  \coup \chinull \Sumlarge_{\substack{K,M\leq \Nd\colon \\ K \not \simeq_\delta M}} \Renorm[\Ncs][M] \UN[K][+]
- \coup \chinull \Renorm[\Ncs][<\Nlarge] \UN[][+]. 
\end{equation*}
Using Lemma \ref{ansatz:lem-renormalization}, Definition \ref{ansatz:def-cutoff-frequency-truncated}, Lemma \ref{modulation:lem-linear}, and that $\UN[K][+]$ is supported on $u$ and $v$-frequencies $\lesssim K$, we obtain that 
\begin{align*}
&\coup\,   \bigg\|  \Sumlarge_{\substack{K,M\leq \Nd\colon \\ K \simeq_\delta M}} \big( \chinull- \chinull[K] \big) \Renorm[\Ncs][M] \UN[K][+] 
\bigg\|_{\Cprod{r-1}{r-1}}  \\
\lesssim&\, \coup  \Sumlarge_{\substack{K,M\leq \Nd\colon \\ K \simeq_\delta M}}  K^{-10} \big\| \Renorm[\Ncs][M] \UN[K][+] \big\|_{\Cprod{r-1}{r-1}} 
\lesssim \coup  \Sumlarge_{\substack{K,M\leq \Nd\colon \\ K \simeq_\delta M}}  K^{-10} N^{-1} K K^{r-s} \big\| \UN[K][+] \big\|_{\Cprod{s-1}{s}} 
\lesssim \coup N^{-1} \Dc, 
\end{align*}
which is more than acceptable.
Using that $\UN[K][+]$ is supported on $u$ and $v$-frequencies $\lesssim K$, using \mbox{Lemma \ref{ansatz:lem-renormalization}}, and using Lemma \ref{modulation:lem-linear},  it follows that
\begin{align*}
\coup \bigg\| \Sumlarge_{\substack{K,M\leq \Nd\colon \\ K \not \simeq_\delta M}} \Renorm[\Ncs][M] \UN[K][+] \bigg\|_{\Cprod{r-1}{r-1}}
&\lesssim \coup \Sumlarge_{\substack{K,M\leq \Nd\colon \\ K \not \simeq_\delta M}} 
\frac{M}{N} \frac{K}{N} \, \big\| \UN[K][+] \big\|_{\Cprod{r-1}{r-1}} \\ 
&\lesssim  \coup \Dc \Sumlarge_{\substack{K,M\leq \Nd\colon \\ K \not \simeq_\delta M}} 
\frac{M}{N} \frac{K}{N} K^{r-s}. 
\end{align*}
 Furthermore, 
since $K,M\leq \Nd$, $\Nd\lesssim N$, and $K\not\simeq_\delta M$, it holds that $\min(K,M)\lesssim N^{1-\delta}$. As a result, 
\begin{equation*}
   \coup \Dc \Sumlarge_{\substack{K,M\leq \Nd\colon \\ K \not \simeq_\delta M}} 
\frac{M}{N} \frac{K}{N} K^{r-s} 
\lesssim\coup \Dc N^{-\delta+r-s} . 
\end{equation*}
which is acceptable. Using $\Nlarge \lesssim 1$ and using Lemma \ref{ansatz:lem-renormalization}, Lemma \ref{ansatz:lem-frequency-support}, and Lemma \ref{modulation:lem-linear}, it also holds that 
\begin{equation}\label{ren:eq-less-Nlarge}
\coup \Big\| \Renorm[\Ncs][<\Nlarge] \UN[][+] \Big\|_{\Cprod{r-1}{r-1}} 
\lesssim \coup  N^{-1} \big\| \UN[][+] \big\|_{\Cprod{r-1}{r-1}} 
\lesssim \coup \Dc N^{-1+r-s},
\end{equation}
which is also acceptable. \\ 

\emph{Estimate of $\RenNErr[][u,\fs-]$}: Due to Definition \ref{ansatz:def-renormalization-error}, it holds that
\begin{equation}\label{ren:eq-p3}
\begin{aligned}
&\RenNErr[][u,\fs-] \\ 
=& - \coup \chinull \Sumlarge_{\substack{L,M\leq \Nd \colon \\ L<M^{1-\delta}}}
\Renorm[\Ncs][M] P^u_{\leq M^{1-\deltap}} \UN[L][\fs-]
- \coup \chinull \Sumlarge_{\substack{L,M\leq \Nd \colon \\ L\geq M^{1-\delta}}}
\Renorm[\Ncs][M] \UN[L][\fs-]
- \coup \chinull \Renorm[\Ncs][<\Nlarge] \UN[][\fs-]. 
\end{aligned}
\end{equation}
We first control the first summand in \eqref{ren:eq-p3}. 
By combining Lemma \ref{ansatz:lem-frequency-support} and the constraints 
$L,M\leq \Nd$ and $L<M^{1-\delta}$, we obtain that $P^u_{\leq M^{1-\deltap}} \UN[L][\fs-]$ is supported on $u$ and $v$-frequencies $\lesssim N^{1-\deltap}$. Using Lemma \ref{ansatz:lem-renormalization} and Lemma \ref{modulation:lem-mixed}, it then follows that
\begin{equation*}
\coup \Big\| \Renorm[\Ncs][M] P^u_{\leq M^{1-\deltap}} \UN[L][\fs-] \Big\|_{\Cprod{r-1}{r-1}} 
\lesssim  \coup N^{-\deltap} \Big\| \UN[L][\fs-] \|_{\Cprod{r-1}{s}} 
\lesssim \coup \Dc^2 N^{-\deltap},
\end{equation*}
which yields an acceptable contribution. We now turn to the second summand in \eqref{ren:eq-p3}. From Lemma \ref{ansatz:lem-frequency-support},  it follows that $\UN[L][\fs-]$ is supported on $v$-frequencies $\sim L$. Using Lemma \ref{ansatz:lem-renormalization} and Lemma \ref{modulation:lem-mixed}, we then obtain
\begin{equation*}
\coup \Big\| \Renorm[\Ncs][M] \UN[L][\fs-] \Big\|_{\Cprod{r-1}{r-1}} 
\lesssim \coup MN^{-1} \Big\|  \UN[L][\fs-] \Big\|_{\Cprod{r-1}{r-1}}  
\lesssim \coup \Dc^2 L^{r-s-1} MN^{-1}. 
\end{equation*}
Due to the constraint $L\geq M^{1-\delta}$ in the second summand in \eqref{ren:eq-p3}, this yields an acceptable contribution. Finally, for the third summand in \eqref{ren:eq-p3}, it holds that 
\begin{equation*}
\coup \Big\|  \Renorm[\Ncs][<\Nlarge] \UN[][\fs-] \Big\|_{\Cprod{r-1}{r-1}}  
\lesssim \coup N^{-1} \Big\| \UN[][\fs-] \Big\|_{\Cprod{r-1}{r-1}} 
\lesssim \coup \Dc N^{-1},
\end{equation*}
which is more than acceptable.\\

\emph{Estimate of $\RenNErr[][u,\fs]$}: Due to Definition \ref{ansatz:def-renormalization-error}, it holds that
\begin{equation*}
\RenNErr[][u,\fs]:= - \coup \chinull \Sumlarge_{M\leq \Nd} 
\Renorm[\Ncs][M] \Big( P^v_{\geq M^{1-\deltap}} + P^u_{\leq M^{1-\deltap}} P^v_{\leq M^{1-\deltap}} \Big) \UN[][\fs]
- \coup \chinull \Renorm[\Ncs][<\Nlarge] \UN[][\fs].
\end{equation*}
Using Lemma \ref{ansatz:lem-renormalization}, it holds that 
\begin{align*}
&\coup \Big\| \Renorm[\Ncs][M] \Big( P^v_{\geq M^{1-\deltap}} + P^u_{\leq M^{1-\deltap}} P^v_{\leq M^{1-\deltap}} \Big) \UN[][\fs] \Big\|_{\Cprod{r-1}{r-1}} \\
\lesssim&\, \coup 
MN^{-1} M^{-(1-\deltap)} 
\big\|  P^v_{\geq M^{1-\deltap}} \UN[][\fs] \big\|_{\Cprod{r-1}{r}}
+ \coup N^{-\deltap} \big\|  P^u_{\leq M^{1-\deltap}} P^v_{\leq M^{1-\deltap}} \UN[][\fs] \big\|_{\Cprod{r-1}{r}} \\
\lesssim&\, \coup \Dc \big( M^{\deltap}N^{-1} + N^{-\deltap} \big),
\end{align*}
which is acceptable. Furthermore, we also obtain from Lemma \ref{ansatz:lem-renormalization} that
\begin{equation*}
\coup \Big\| \Renorm[\Ncs][<\Nlarge] \UN[][\fs] \Big\|_{\Cprod{r-1}{r-1}} 
\lesssim \coup N^{-1} \big\| \UN[][\fs]\big\|_{\Cprod{r-1}{r-1}} \lesssim \coup \Dc N^{-1},
\end{equation*}
which is more than acceptable. \\ 
 
\emph{Estimate of $\RenNErr[][u,+-]$, $\RenNErr[][u,-]$, and $\RenNErr[][u,+\fs]$:} 
From Definition \ref{ansatz:def-renormalization-error}, Lemma \ref{ansatz:lem-renormalization}, and the triangle inequality, it follows that
\begin{equation}\label{ren:eq-p1}
\begin{aligned}
&\Big\| \RenNErr[][u,+-] \Big\|_{\Cprod{r-1}{r-1}} 
+ \Big\| \RenNErr[][u,-] \Big\|_{\Cprod{r-1}{r-1}} 
+ \Big\| \RenNErr[][u,+\fs] \Big\|_{\Cprod{r-1}{r-1}} \\
\lesssim& \, \coup N^{-\eta} \bigg(
\Sumlarge_{\substack{K,M\leq \Nd \colon \\ K \simeq_\delta M}}
\Big\| \UN[K,M][+-] \Big\|_{\Cprod{r+\eta-1}{r+\eta-1}} 
+ \Sumlarge_{\substack{M\leq \Nd}}
\Big\|\UN[M][-] \Big\|_{\Cprod{r+\eta-1}{r+\eta-1}}  \\ 
+& \,  \Sumlarge_{\substack{K\leq \Nd}}
\Big\| \UN[K][+\fs] \Big\|_{\Cprod{r+\eta-1}{r+\eta-1}} \bigg). 
\end{aligned}
\end{equation}
For the $(+-)$-term, it follows from Lemma \ref{modulation:lem-bilinear} and $K\simeq_\delta M$ that 
\begin{equation*}
\Big\| \UN[K,M][+-] \Big\|_{\Cprod{r+\eta-1}{r+\eta-1}}  
\lesssim K^{r-1/2+\eta} M^{r-3/2+\eta} \Dc^2 \lesssim \max(K,M)^{-1+2\delta} \Dc^2,
\end{equation*}
which is acceptable. 
For the $(-)$ and $(+\fs)$-terms, it follows from Lemma \ref{ansatz:lem-frequency-support}, Lemma \ref{modulation:lem-mixed}, and Lemma \ref{modulation:lem-linear-reversed} that 
\begin{align*}
\Big\|\UN[M][-] \Big\|_{\Cprod{r+\eta-1}{r+\eta-1}}   
&\lesssim M^{(1-\delta) (r+\eta-s)} M^{r+\eta-1-s} \Big\|\UN[M][-] \Big\|_{\Cprod{s-1}{s}} \lesssim M^{-1+2\delta} \Dc^2, \\ 
\Big\| \UN[K][+\fs] \Big\|_{\Cprod{r+\eta-1}{r+\eta-1}} 
&\lesssim K^{r+\eta-1} K^{(1-\delta)(r+\eta-1-r)}
\Big\| \UN[K][+\fs] \Big\|_{\Cprod{s-1}{r}} \lesssim K^{-1+2\delta} \Dc^2,  
\end{align*}
which both lead to acceptable contributions.
\end{proof}

\section{Energy increment}\label{section:energy-increment}

The goal of this section is to control the energy increment for solutions of the finite-dimensional approximation of the wave maps equation \eqref{ansatz:eq-wave-maps}, which is used to control the Radon-Nikodym derivative in Proposition \ref{structure:prop-Gibbs}. In turn, the estimate of the Radon-Nikodym derivative can then be used to prove the almost invariance of the Gibbs measure (see Proposition \ref{main:prop-almost-invariance}). The ingredients used in this section are mostly similar as for our estimates of the error terms in the remainder equations (Sections \ref{section:hhl}-\ref{section:renormalization}). However, there are two important new ingredients, which are two key cancellations in the proof of Lemma \ref{increment:lem-increment-object-II}. 

Before we state the main result of this section (Proposition \ref{increment:prop-energy-increment}), we make the following definition. 

\begin{definition}[Energy increment]\label{increment:def-energy-increment} Let
$N\in \Dyadiclarge$, let $R\geq 1$, and let 
$U,V \colon \R \times \bT_R \rightarrow \frkg$. Furthermore, let $\zeta\in \Cuttilde$, where $\Cuttilde$ is as in Definition \ref{prelim:def-cut-off}, and let $\tau,t_0,t_1\in \R$. Then, we define the energy increment 
\begin{equation*}
\EIN \big(U,V;t_0,t_1\big) := \frac{1}{4} \int_{t_0}^{t_1} \dt \int_{\bT_R} \dx \,  \zeta(t)^2 \Big \langle U(t,x) , \Renorm[\Nscript][] V(t,x) \Big\rangle_\frkg. 
\end{equation*}
We also define
\begin{equation*}
\EIN\big(U,V;\tau\big) := \EIN\big(U,V;-\tau,0\big),
\end{equation*}
which matches the energy increment in Proposition \ref{structure:prop-Gibbs}.  
\end{definition}

For the rest of this section, it will be convenient to work with the kernel $\widecheck{\Gamma}^{(\Nscript)}$ and operator $\Gamma^{(\Nscript),x}$, which are defined as
\begin{equation}\label{increment:eq-Gamma}
\widecheck{\Gamma}^{(\Nscript)}(y) := \big( \widecheck{\rho}_{\leq N} \ast \widecheck{\rho}_{\leq N}\big)(y) 
\Cf^{(\Nscript)}(y) 
\qquad \text{and} \qquad \Gamma^{(\Nscript),x}  := \int_{\R} \dy \, \widecheck{\Gamma}^{(\Nscript)}(y) \Theta^x_y.
\end{equation}
Due to Definition \ref{ansatz:def-Killing} and \eqref{increment:eq-Gamma}, the Killing-renormalization $\Renorm[\Nscript]$ can then be written as
\begin{equation}\label{increment:eq-Gamma-Killing}
\Renorm[\Nscript] = P_{\leq N}^x \Gamma^{(\Nscript),x} P_{\leq N}^x \Kil. 
\end{equation}
In the next lemma, we obtain a representation of the energy increment in which  the Killing-renormalization $\Renorm[\Nscript]$ has been replaced using \eqref{increment:eq-Gamma-Killing} and the integration is performed in the $u$ and $v$-variables. 

\begin{lemma}[Representation of the energy-increment]\label{increment:lem-representation-increment}
Let $N\in \Dyadiclarge$, let $R\geq 1$, and let $\zeta\in \Cuttilde$. Furthermore, let $U,V \colon \R_{u,v}^{1+1} \rightarrow \frkg$, let $t_0,t_1\in \R$, and let  
\begin{equation*}
\Domain_{t_0,t_1}^{(\Rscript)} := \Big\{ (u,v) \in \R \times \R \colon 
-\pi R \leq \frac{u+v}{2} \leq \pi R, \, t_0 \leq \frac{v-u}{2} \leq t_1 \Big\}. 
\end{equation*}
Then, it holds that 
\begin{equation}\label{increment:eq-representation-increment}
\begin{aligned}
\EIN \big(U,V;t_0,t_1\big) 
&= \frac{1}{8} \int\displaylimits_{\Domain_{t_0,t_1}^{(\Rscriptscript)}} \du \dv \, \zeta\big( \tfrac{v-u}{2} \big)^2
\Big\langle P_{\leq N}^x U(u,v) , \Gamma^{(\Nscript),x} \Kil P_{\leq N}^x V(u,v) \Big\rangle_\frkg \\
&= -\frac{1}{8} \int\displaylimits_{\Domain_{t_0,t_1}^{(\Rscriptscript)}} \du \dv \, \zeta\big( \tfrac{v-u}{2} \big)^2
\Big\langle \Gamma^{(\Nscript),x} \Kil P_{\leq N}^x U(u,v) ,  P_{\leq N}^x V(u,v) \Big\rangle_\frkg. 
\end{aligned}
\end{equation}
\end{lemma}

\begin{proof}
This follows directly from \eqref{increment:eq-Gamma-Killing}, the skew-adjointness of $\Gamma^{(\Nscript),x}$, the symmetry of $\Kil$, and our change of variables from Cartesian into null coordinates. 
\end{proof}

In the main proposition of this subsection, we control the energy increment of the solution to the finite-dimensional approximation of the wave maps equation. 

\begin{proposition}[Energy increment]\label{increment:prop-energy-increment}
Let the post-modulation hypothesis (Hypothesis \ref{hypothesis:post}) be satisfied. Then, it holds that 
\begin{equation*}
\sup_{\zeta \in \Cuttilde} \sup_{\tau \in [-1,1]} \Big| \EIN\big(\UN[][],\VN[][];\tau\big) \Big| 
\lesssim R N^{-(1-10\delta)\delta_1} \Dc^2. 
\end{equation*}
\end{proposition}

\begin{figure}
\scalebox{0.925}{
\begin{tabular}{
!{\vrule width 1pt}>{\centering\arraybackslash}P{1.5cm}
!{\vrule width 1pt}>{\centering\arraybackslash}P{\colwidth}
!{\vrule width 1pt}>{\centering\arraybackslash}P{\colwidth}
!{\vrule width 1pt}>{\centering\arraybackslash}P{\colwidth}
!{\vrule width 1pt}>{\centering\arraybackslash}P{\colwidth}
!{\vrule width 1pt}>{\centering\arraybackslash}P{\colwidth}
!{\vrule width 1pt}>{\centering\arraybackslash}P{\colwidth}
!{\vrule width 1pt}} 
\noalign{\hrule height 1pt} & & & & & & 
 \\[-5.8ex]
 \, \, \, \, \, \,   $V$ \vspace{-1ex} \newline \, \, \, \,  $U$  \hspace{-0.5ex}\vspace{-0.5ex} 
&  $(-)$ & $(+-)$ & $(+)$   & $(\fs-)$ & $(+\fs)$ & $(\fs)$  
\\[4pt] \noalign{\hrule height 1pt} \rule{0pt}{14pt}
$(+)$ 
& \cellcolor{Blue!30}  \ref{increment:lem-increment-object-I}
& \cellcolor{Orange!30}  \ref{increment:lem-increment-object-II}
& \cellcolor{Orange!30}  \ref{increment:lem-increment-object-II}
& \cellcolor{magenta!30}  \ref{increment:lem-increment-fs}
& \cellcolor{magenta!30}  \ref{increment:lem-increment-fs}
& \cellcolor{magenta!30}  \ref{increment:lem-increment-fs}
\\[4pt] \noalign{\hrule height 1pt} \rule{0pt}{14pt}
$(+-)$ 
& \cellcolor{Gray!80} 
& \cellcolor{Orange!30}  \ref{increment:lem-increment-object-II}
& \cellcolor{Orange!30}  \ref{increment:lem-increment-object-II}
& \cellcolor{magenta!30}  \ref{increment:lem-increment-fs}
& \cellcolor{magenta!30}  \ref{increment:lem-increment-fs}
& \cellcolor{magenta!30}  \ref{increment:lem-increment-fs}
\\[4pt] \noalign{\hrule height 1pt} \rule{0pt}{14pt}
$(-)$ 
& \cellcolor{Gray!80}
& \cellcolor{Gray!80}
& \cellcolor{Green!30}  \ref{increment:lem-increment-remaining}
& \cellcolor{Green!30}  \ref{increment:lem-increment-remaining}
& \cellcolor{Green!30}  \ref{increment:lem-increment-remaining}
& \cellcolor{Green!30}  \ref{increment:lem-increment-remaining}
\\[4pt] \noalign{\hrule height 1pt} \rule{0pt}{14pt}
$(+\fs)$ 
& \cellcolor{Gray!80}
& \cellcolor{Gray!80}
& \cellcolor{Gray!80}
& \cellcolor{Green!30}  \ref{increment:lem-increment-remaining}
& \cellcolor{Green!30}  \ref{increment:lem-increment-remaining}
& \cellcolor{Green!30}  \ref{increment:lem-increment-remaining}
\\[4pt] \noalign{\hrule height 1pt} \rule{0pt}{14pt}
$(\fs-)$ 
& \cellcolor{Gray!80}
& \cellcolor{Gray!80}
& \cellcolor{Gray!80}
& \cellcolor{Gray!80}
& \cellcolor{Green!30}  \ref{increment:lem-increment-remaining}
& \cellcolor{Green!30}  \ref{increment:lem-increment-remaining}
\\[4pt] \noalign{\hrule height 1pt} \rule{0pt}{14pt}
$(\fs)$ 
& \cellcolor{Gray!80}
& \cellcolor{Gray!80}
& \cellcolor{Gray!80}
& \cellcolor{Gray!80}
& \cellcolor{Gray!80}
& \cellcolor{Green!30}  \ref{increment:lem-increment-remaining}
\\[3pt] \noalign{\hrule height 1pt} 
\end{tabular}
}
\caption{\small{In this figure, we give an overview of the cases in the estimate of the energy increment. The different cases are treated using either Lemma \ref{increment:lem-increment-object-I}, Lemma \ref{increment:lem-increment-object-II}, Lemma \ref{increment:lem-increment-fs}, or Lemma \ref{increment:lem-increment-remaining}. The corresponding cells are colored in blue, orange, magenta, and green, respectively. 
}} 
\label{figure:increment-cases}
\end{figure}

The proof of Proposition \ref{increment:prop-energy-increment} occupies the rest of this subsection.
The different cases in the proof of Proposition \ref{increment:prop-energy-increment}, which are split over four separate lemmas, are illustrated in Figure \ref{figure:increment-cases}. Before we treat specific interactions, however, we prove a general lemma which allows us to bound the energy increment by $\Cprod{\gamma_1}{\gamma_2}$-norms of the integrand.

\begin{lemma}\label{increment:lem-increment-integrand}
Let $R\geq 1$, let $t_0,t_1\in [-1,1]$, and let $\gamma_1,\gamma_2 \in (-1,1)\backslash\{0\}$ satisfy $\gamma_1+\gamma_2>-1$. For all $F\colon \R^{1+1}_{u,v}\rightarrow \frkg$, it then holds that
\begin{equation}\label{increment:eq-increment-integrand-F}
\bigg| \int\displaylimits_{\Domain_{t_0,t_1}^{(\Rscriptscript)}} \du \dv F(u,v) \bigg| \lesssim  R \big\| F \big\|_{\Cprod{\gamma_1}{\gamma_2}}. 
\end{equation}
In particular, for all  $U,V \colon \R_{u,v}^{1+1} \rightarrow \frkg$ and $N \in \Dyadiclarge$, it holds that 
\begin{equation}\label{increment:eq-increment-integrand}
\sup_{\zeta \in \Cuttilde} \Big| \EIN\big( U,V;t_0,t_1 \big) \Big| \lesssim  R
\Big\|  \Big\langle P_{\leq N}^x U , \Gamma^{(\Nscript),x} \Kil P_{\leq N}^x V \Big\rangle_\frkg \Big\|_{\Cprod{\gamma_1}{\gamma_2}}. 
\end{equation}
\end{lemma}

\begin{proof} Due to Lemma \ref{increment:lem-representation-increment}, it suffices to prove the first estimate \eqref{increment:eq-increment-integrand-F}. To this end, we let $\varphi\in C^\infty_c(\R)$ be such that $\{ \varphi(\cdot-x_0)\}_{x_0\in \Z}$ is a partition of unity. Then, we decompose 
\begin{equation}\label{increment:eq-increment-integrand-p1}
\int\displaylimits_{\Domain_{t_0,t_1}^{(\Rscriptscript)}} \du \dv F(u,v) 
= \sum_{u_0,v_0\in \Z} \int\displaylimits_{\Domain_{t_0,t_1}^{(\Rscriptscript)}} \du \dv \, \varphi(u-u_0) \varphi(v-v_0) F(u,v).
\end{equation}
Due to the definition of $\Domain_{t_0,t_1}^{(\Rscript)}$ and $t_0,t_1 \in [-1,1]$, we may restrict $u_0,v_0\in \Z$ in 
\eqref{increment:eq-increment-integrand-p1} to $|u_0|,|v_0|\lesssim R$ and $|u_0-v_0|\lesssim 1$. Together with Lemma \ref{prelim:lem-integral}
and Lemma \ref{prelim:lem-trace}, it then follows that
\begin{align*}
\big| \eqref{increment:eq-increment-integrand-p1} \big| 
&\lesssim \sum_{\substack{u_0,v_0\in \Z\colon \\ |u_0|,|v_0|\lesssim R, \\ |u_0-v_0|\lesssim 1}} 
\bigg| \int_{-\pi R -t_1}^{\pi R - t_0} \du \,  
\varphi(u-u_0) \Int^v \Big(  \varphi(v-v_0) F(u,v) \Big)\Big|_{v=\max(u-2 t_1,-u-2\pi R)}^{\min(u-2t_0,-u+2\pi R)}
\bigg| \\\
&\lesssim \sum_{\substack{u_0,v_0\in \Z\colon \\ |u_0|,|v_0|\lesssim R, \\ |u_0-v_0|\lesssim 1}} \big\| F \big\|_{\Cprod{\gamma_1}{\gamma_2}}
\lesssim R \big\| F \big\|_{\Cprod{\gamma_1}{\gamma_2}}. \qedhere
\end{align*}
\end{proof}

In the next lemma, which uses no structural information on the solution, we show that the energy increment is controlled if either $\UN[][]$ or $\VN[][]$ is projected to frequencies much smaller than $N$. 

\begin{lemma}[Eliminating low-frequency terms]\label{increment:lem-low-frequency}
Let $K_u,K_v,M_u,M_v \in \dyadic$ satisfy $\max(K_u,K_v) \lesssim N^{1-\delta}$ or $\max(M_u,M_v)\lesssim N^{1-\delta}$ and let $R\geq 1$. Then, it holds 
for all $\UN[][\ast]\in \Cprod{s-1}{s}$ and $\VN[][\ast]\in \Cprod{s}{s-1}$ that 
\begin{equation*}
\begin{aligned}
 &\sup_{\zeta\in \Cuttilde} \sup_{\tau \in [-1,1]} \Big| \EIN\big( P_{K_u}^u P_{K_v}^v \UN[][\ast], P^u_{M_u} P^v_{M_v}\VN[][\ast];\tau \big) \Big| \\  
\lesssim&\, R  N^{-\delta+10\delta_2} 
\big\| \UN[][\ast] \big\|_{\Cprod{s-1}{s}} \big\| \VN[][\ast] \big\|_{\Cprod{s}{s-1}}. 
\end{aligned}
\end{equation*}
\end{lemma}

\begin{proof}
By symmetry (as in Lemma \ref{increment:lem-representation-increment}), it suffices to treat the case $\max(M_u,M_v) \lesssim N^{1-\delta}$. Using Lemma \ref{increment:lem-increment-integrand}, it suffices to control the integrand in $\EIN$ in the $\Cprod{s-1}{-s+\eta}$-norm.  Using the product estimate from Corollary \ref{prelim:cor-product}, we obtain that
\begin{align*}
&\Big\| \Big\langle P_{\leq N}^x  P_{K_u}^u P_{K_v}^v \UN[][\ast], \Gamma^{(\Nscript),x} \Kil P_{\leq N}^x  P^u_{M_u} P^v_{M_v}\VN[][\ast] \Big\rangle_\frkg \,  \Big\|_{\Cprod{s-1}{-s+\eta}} \\
\lesssim& \, \big\| P_{\leq N}^x  P_{K_u}^u P_{K_v}^v \UN[][\ast] \big\|_{\Cprod{s-1}{s}} \big\| \Gamma^{(\Nscript),x} P^u_{M_u} P^v_{M_v}\VN[][\ast] \big\|_{\Cprod{1-s+\eta}{-s+\eta}}.  
\end{align*}
Using Lemma \ref{ansatz:lem-renormalization}, we obtain that
\begin{align*}
    \big\| \Gamma^{(\Nscript),x} P^u_{M_u} P^v_{M_v}\VN[][\ast] \big\|_{\Cprod{1-s+\eta}{-s+\eta}}
    &\lesssim M_u^{1-s+\eta-s} M_v^{-s+\eta-(s-1)}  \big\| \Gamma^{(\Nscript),x} P^u_{M_u} P^v_{M_v}\VN[][\ast] \big\|_{\Cprod{s}{s-1}}  \\
    &\lesssim (M_u M_v)^{1-2s+\eta} \max\big( M_u , M_v \big) N^{-1} \big\| \VN[][\ast] \big\|_{\Cprod{s}{s-1}}.
\end{align*}
Since $\max(M_u,M_v) \lesssim N^{1-\delta}$, the prefactor can be estimated by
\begin{equation*}
 (M_u M_v)^{1-2s+\eta} \max\big( M_u , M_v \big) N^{-1} 
 \lesssim N^{4\delta_2+2\delta_3} N^{1-\delta} N^{-1} \lesssim N^{-\delta+4\delta_2+2\delta_3},
\end{equation*}
which is acceptable. 
\end{proof}

We now start the detailed case analysis which is illustrated in Figure \ref{figure:increment-cases}. In our first lemma, we control the energy increment which comes from the interaction of $\UN[][+]$ and $\VN[][-]$. 

\begin{lemma}[Energy increment from modulated objects I]\label{increment:lem-increment-object-I}
Assume that the post-modulation hypothesis (Hypothesis \ref{hypothesis:post}) is satisfied. For all $K,M\in \Dyadiclarge$, it holds that
\begin{equation*}
\sup_{\zeta\in \Cuttilde} \sup_{\tau \in [-1,1]} \Big| \EIN \Big( \UN[K][+],\VN[M][-];\tau \Big) \Big| \lesssim R N^{-\delta+10\delta_2} \Dc^2. 
\end{equation*}
\end{lemma}

\begin{proof}
By symmetry, we may assume that $M\geq K$.  Due to Lemma \ref{increment:lem-low-frequency}, we may also assume that $M\gtrsim N^{1-\delta}$. Using Lemma \ref{increment:lem-representation-increment}, we first write 
\begin{align}
 \EIN \Big( \UN[K][+],\VN[M][-];\tau \Big) 
=\,& - \frac{1}{8} \int\displaylimits_{\Domain_{-\tau,0}^{(\Rscriptscript)}} \du \dv \, \zeta\big( \tfrac{v-u}{2} \big)^2 
\Big\langle \Renorm[N] \UN[K][+] ,  \VN[M][-] \Big\rangle_\frkg \notag \\
=\, - &\frac{1}{8} \int\displaylimits_{\Domain_{-\tau,0}^{(\Rscriptscript)}} \du \dv \, \zeta\big( \tfrac{v-u}{2} \big)^2
\Big\langle \Renorm[N] \UN[K][+] ,  \partial_v \IVN[M][-] \Big\rangle_\frkg 
\label{increment:eq-object-I-p1} \\ 
-\, & \frac{1}{8} \int\displaylimits_{\Domain_{-\tau,0}^{(\Rscriptscript)}} \du \dv \, \zeta\big( \tfrac{v-u}{2} \big)^2
\Big\langle \Renorm[N] \UN[K][+] , \VN[M][-]-   \partial_v \IVN[M][-]\Big\rangle_\frkg 
\label{increment:eq-object-I-p2}. 
\end{align}
Using integration by parts, we then write \eqref{increment:eq-object-I-p1} as 
\begin{align}
\eqref{increment:eq-object-I-p1} 
&=  \frac{1}{8} \int\displaylimits_{\Domain_{-\tau,0}^{(\Rscriptscript)}} \du \dv \, \zeta\big( \tfrac{v-u}{2} \big)^2
\Big\langle \partial_v \Renorm[N] \UN[K][+] ,   \IVN[M][-] \Big\rangle_\frkg 
\label{increment:eq-object-I-p3} \\ 
&+ \frac{1}{8}  \int\displaylimits_{\Domain_{-\tau,0}^{(\Rscriptscript)}} \du \dv \, \partial_v \Big( \zeta\big( \tfrac{v-u}{2} \big)^2 \Big)
\Big\langle \Renorm[N] \UN[K][+] ,   \IVN[M][-] \Big\rangle_\frkg  \label{increment:eq-object-I-q} \\ 
&- \frac{1}{8} \int_{-\pi R}^{\pi R+\tau} \du \, \zeta\big( \tfrac{v-u}{2} \big)^2 \Big\langle  \Renorm[N] \UN[K][+] ,   \IVN[M][-] \Big\rangle_\frkg\Big|_{v=\max(u-2\tau,-u-2\pi R)}^{\min(u,-u+2\pi R)}.
\label{increment:eq-object-I-p4}
\end{align}
It remains to estimate \eqref{increment:eq-object-I-p2}, \eqref{increment:eq-object-I-p3}, \eqref{increment:eq-object-I-q}, and \eqref{increment:eq-object-I-p4}, which are treated separately. \\ 

\emph{Estimate of \eqref{increment:eq-object-I-p2}:} 
Using our para-product estimate (Lemma \ref{prelim:lem-paraproduct}), Lemma \ref{ansatz:lem-renormalization}, and Lemma \ref{increment:lem-increment-integrand}, we obtain that
\begin{align*}
\big| \eqref{increment:eq-object-I-p2}\big| 
 &\lesssim R \, \Big\| \Big\langle \Renorm[N] \UN[K][+] , \VN[M][-]-   \partial_v \IVN[M][-]\Big\rangle_\frkg \,  \Big\|_{\Cprod{s-1}{r-1}}  \\
&\lesssim R\, \big\| \UN[K][+] \big\|_{\Cprod{s-1}{s}} 
\big\|\VN[M][-]-   \partial_v \IVN[M] \big\|_{\Cprod{r}{r-1}}. 
\end{align*}
Using Lemma \ref{modulation:lem-linear} and Lemma \ref{modulation:lem-integration}, we further estimate 
\begin{equation*}
\big\| \UN[K][+] \big\|_{\Cprod{s-1}{s}} 
\big\|\VN[M][-]-   \partial_v \IVN[M] \big\|_{\Cprod{r}{r-1}} 
\lesssim M^{-\frac{1}{2}} \Dc^2. 
\end{equation*}
Since $M\gtrsim N^{1-\delta}$, this is acceptable. \\ 

\emph{Estimate of \eqref{increment:eq-object-I-p3}:}
Using Lemma \ref{ansatz:lem-frequency-support} and $M\geq K$, it follows that $\UN[K][+]$ is supported on $v$-frequencies $\lesssim K^{1-\delta+\vartheta}\lesssim M^{1-\delta+\vartheta}$, and thus the integrand in \eqref{increment:eq-object-I-p3} only contains low$\times$high-interactions in the $v$-variable. Using our low$\times$high-estimate (Lemma \ref{prelim:lem-paraproduct}) and Lemma \ref{increment:lem-increment-integrand}, it then follows that 
\begin{align*}
\big| \eqref{increment:eq-object-I-p3}\big| 
\lesssim  R\, \Big\| \Big\langle \partial_v \Renorm[N] \UN[K][+] ,   \IVN[M][-] \Big\rangle_\frkg  \, \Big\|_{\Cprod{s-1}{r-1}}  
\lesssim R\, \big\| \partial_v \UN[K][+] \big\|_{\Cprod{s-1}{\eta}} \big\| \IVN[M][-] \big\|_{\Cprod{s}{r-1}}. 
\end{align*}
Using Lemma \ref{ansatz:lem-frequency-support}, Lemma \ref{modulation:lem-linear}, and $M\geq K$, it then follows that 
\begin{align*}
\big\| \partial_v \UN[K][+] \big\|_{\Cprod{r-1}{\eta}} \big\| \IVN[M][-] \big\|_{\Cprod{s}{r-1}}
&\lesssim K^{\eta+1-s} M^{r-1-s} \big\| \UN[K][+] \big\|_{\Cprod{s-1}{s}} \big\| \IVN[M][-] \big\|_{\Cprod{s}{s}} \lesssim M^{r-2s+\eta} \Dc^2. 
\end{align*}
Since $r-2s+\eta=-\frac{1}{2}+\mathcal{O}(\delta_1)$ and $M\gtrsim N^{1-\delta}$, this is acceptable. \\ 

\emph{Estimate of \eqref{increment:eq-object-I-q}:} Using Lemma \ref{increment:lem-increment-integrand}, Lemma \ref{ansatz:lem-renormalization}, and Lemma \ref{modulation:lem-linear}, it holds that 
\begin{align*}
\big| \eqref{increment:eq-object-I-q} \big| 
\lesssim R \Big\| \Big\langle  \Renorm[N] \UN[K][+] ,   \IVN[M][-] \Big\rangle_\frkg  \, \Big\|_{\Cprod{s-1}{r-1}}  
\lesssim R \big\| \UN[K][+] \big\|_{\Cprod{s-1}{s}} \big\| \IVN[M][-] \big\|_{\Cprod{s}{r-1}} 
\lesssim R M^{r-1-s} \Dc^2. 
\end{align*}
Since $M\gtrsim N^{1-\delta}$, this is acceptable. \\

\emph{Estimate of \eqref{increment:eq-object-I-p4}:} 
Using integral estimates (Lemma \ref{prelim:lem-integral}), we obtain that 
\begin{equation}\label{increment:eq-object-I-p5}
\big| \eqref{increment:eq-object-I-p4}\big| 
\lesssim R \, \bigg\| \zeta\big( \tfrac{v-u}{2} \big)^2 \Big\langle  \Renorm[N] \UN[K][+] ,   \IVN[M][-] \Big\rangle_\frkg\Big|_{v=\max(u-2\tau,-u-2\pi R)}^{\min(u,-u+2\pi R)} \bigg\|_{\C_u^{r-1}}.  
\end{equation}
\begin{equation*}
\eqref{increment:eq-object-I-p5} \lesssim \coup \Ac^2 \Bc^2 \max(K,M)^{r-\frac{1}{2}+\eta} M^{-\frac{1}{2}} \lesssim \Dc^2 M^{-1+r+\eta}. 
\end{equation*}
Since $-1+r+\eta=-\frac{1}{2}+\mathcal{O}(\delta_1)$ and $M\gtrsim N^{1-\delta}$, this is acceptable. 
\end{proof}

In the previous lemma, we treated the $(+)$$\times$$(-)$-interaction. In the following lemma, we treat all other interactions 
which only involve modulated objects.  

\begin{lemma}[Energy increment from modulated objects II]\label{increment:lem-increment-object-II}
Assume that the post-modulation hypothesis (Hypothesis \ref{hypothesis:post}) is satisfied. Let 
 $K,M,K_u,K_v,M_u,M_v \in \Dyadiclarge$ and let 
\begin{equation*}
\UN[][\ast] \in \Big\{ \UN[K][+], \UN[K_u,K_v][+-]\Big\} 
\qquad \text{and} \qquad 
\VN[][\ast] \in \Big\{ \VN[M_u,M_v][+-], \VN[M][+] \Big\}.
\end{equation*}
Then, it holds that
\begin{equation*}
\sup_{\zeta\in \Cuttilde}  \sup_{\tau \in [-1,1]} \Big| \EIN \Big( \UN[][\ast],\VN[][\ast];\tau \Big) \Big| \lesssim R N^{-\delta+10\delta_1} \Dc^3. 
\end{equation*}
\end{lemma}

\begin{remark}
The proof of Lemma \ref{increment:lem-increment-object-II} is rather delicate. In our estimate of the $(+-)$$\times$$(+-)$-interaction, it is crucial that the covariance function $\Cf^{(N)}$ from Definition \ref{ansatz:def-Killing} is odd. In our estimates of the $(+)$$\times$$(+-)$ and $(+)$$\times$$(+)$-interactions, we rely on a cancellation stemming from the Lie bracket and the Killing map (see Lemma \ref{prelim:lem-Killing-vs-Bracket}). 
\end{remark}

\begin{proof} In contrast to the proof of Lemma \ref{increment:lem-increment-object-I}, this argument controls all energy increments by estimating the integrands in $\Cprod{\gamma_1}{\gamma_2}$-norms. Due to Lemma \ref{increment:lem-low-frequency}, we may assume that 
\begin{equation}\label{increment:eq-increment-modulated2-p0}
K,M,\max(K_u,K_v),\max(M_u,M_v) \gtrsim N^{1-\delta}.
\end{equation}

\emph{Case 1: $(+-)\times(+-)$-interaction.} Due to the definition of $\UN[K_u,K_v][+-]$ and $\VN[M_u,M_v][+-]$, the integral representation of $(P_{\leq N}^x)^2$ and $\Gamma^{(\Nscript),x}$, and Remark \ref{ansatz:rem-commutativity-chipm}, it holds that
\begin{equation}\label{increment:eq-increment-modulated2-p1}
\begin{aligned}
&\Big\langle P_{\leq N}^x \UN[K_u,K_v][+-] , \Gamma^{(\Nscript),x} \Kil P_{\leq N}^x \VN[M_u,M_v][+-] \Big\rangle_\frkg \\
=& \chinull[K_u,K_v] \chinull[M_u,M_v] \int_{\R^3} \dw \dy \dz \bigg(  \big( \widecheck{\rho}_{\leq N} \ast \widecheck{\rho}_{\leq N}\big)(y)  \widecheck{\Gamma}^{(\Nscript)}(w) \big( \widecheck{\rho}_{\leq N} \ast \widecheck{\rho}_{\leq N}\big)(z-w )  \\
&\times \mathcal{T} \Big( \Theta^x_y P_{\leq N}^x \UN[K_u][+] \otimes \Theta^x_y P_{\leq N}^x \IVN[K_v][-] 
\otimes \Theta^x_z P_{\leq N}^x \IUN[M_u][+] \otimes \Theta^x_z P_{\leq N}^x \VN[M_v][-] \Big) \bigg),
\end{aligned}
\end{equation}
where
\begin{equation}\label{increment:eq-increment-modulated2-p2}
\mathcal{T} \colon \frkg^{\otimes 4} \rightarrow \R, \, \, A \otimes B \otimes C \otimes D \mapsto \big \langle \big[ A, B \big] , \Kil \big[ C, D \big] \big \rangle_\frkg. 
\end{equation}
Due to Corollary \ref{killing:cor-tensors} and Lemma \ref{increment:lem-increment-integrand}, the non-resonant part of \eqref{increment:eq-increment-modulated2-p1} yields an acceptable contribution. Thus, it remains to treat the resonant part of \eqref{increment:eq-increment-modulated2-p1}. Since multiplication by $\chinull[K_u,K_v] \chinull[M_u,M_v]$ is bounded in $\Cprod{\gamma_1}{\gamma_2}$-spaces, we may omit the $\chinull[K_u,K_v] \chinull[M_u,M_v]$-factor in \eqref{increment:eq-increment-modulated2-p1}. Then, the resonant part is given by
\begin{equation}\label{increment:eq-increment-modulated2-p3}
\begin{aligned}
&\coup \int_{\R^3} \dw \dy \dz \bigg(  \big( \widecheck{\rho}_{\leq N} \ast \widecheck{\rho}_{\leq N}\big)(y)  \widecheck{\Gamma}^{(\Nscript)}(w) \big( \widecheck{\rho}_{\leq N} \ast \widecheck{\rho}_{\leq N}\big)(z-w )  \Cf^{(\Ncs)}_{K_u}(y-z) \Cf^{(\Ncs)}_{K_v}(z-y) \bigg) \\
&\times \, 
\mathcal{T} \Big( E_a \otimes E_b \otimes E^a \otimes E^b \Big) . 
\end{aligned}
\end{equation}
In fact, we now show that the resonant part \eqref{increment:eq-increment-modulated2-p3} equals zero. 
Using the definition of $\mathcal{T}$ from \eqref{increment:eq-increment-modulated2-p2}, Lemma \ref{prelim:lem-Killing-form}, and the skew-adjointness of the adjoint map, we obtain that 
\begin{equation}\label{increment:eq-increment-modulated2-p4}
\begin{aligned}
\mathcal{T} \Big( E_a \otimes E_b \otimes E^a \otimes E^b \Big) 
&= \Big\langle \big[ E_a, E_b \big], \Kil \big[ E^a, E^b \big] \Big \rangle_\frkg \\ 
&= \Tr \Big( \Ad\big( \big[ E_a, E_b \big] \big)  \Ad\big( \big[ E^a, E^b \big] \big) \Big) \\ 
&= - \sum_{a,b=1}^{\dim \frkg} \Big\| \Ad\big( \big[ E_a, E_b \big] \big)  \Big\|_{\textup{HS}}^2,
\end{aligned}
\end{equation}
where $\| \cdot\|_{\textup{HS}}$ denotes the Hilbert-Schmidt norm. 
For a general Lie algebra $\frkg$, such as $\frkg= \mathfrak{s}\mathfrak{o}(n)$,  \eqref{increment:eq-increment-modulated2-p4} is clearly non-zero. In order to prove that \eqref{increment:eq-increment-modulated2-p3} is zero, we therefore have to show that the integral vanishes. To this end, we examine the integrand in \eqref{increment:eq-increment-modulated2-p3}, i.e., 
\begin{equation}\label{increment:eq-increment-modulated2-p4'}
\big( w,y,z \big) \in \R^3 \mapsto \big( \widecheck{\rho}_{\leq N} \ast \widecheck{\rho}_{\leq N}\big)(y)  \widecheck{\Gamma}^{(\Nscript)}(w) \big( \widecheck{\rho}_{\leq N} \ast \widecheck{\rho}_{\leq N}\big)(z-w )  \Cf^{(\Ncs)}_{K_u}(y-z) \Cf^{(\Ncs)}_{K_v}(z-y). 
\end{equation}
Since 
$\big( \widecheck{\rho}_{\leq N} \ast \widecheck{\rho}_{\leq N}\big)(y)$ and  $\big( \widecheck{\rho}_{\leq N} \ast \widecheck{\rho}_{\leq N}\big)(z-w )$
are even and 
$\widecheck{\Gamma}^{(\Nscript)}(w)$, $\Cf^{(\Ncs)}_{K_u}(y-z)$,  and $\Cf^{(\Ncs)}_{Kv}(z-y)$
are odd, the function in \eqref{increment:eq-increment-modulated2-p4'} is odd. Thus,  \eqref{increment:eq-increment-modulated2-p3} has an odd integrand and the integral therefore vanishes. \\

\emph{Case 2: $(+)\times(+-)$-interaction.} Due to the definitions of $\UN[K][+]$ and $\VN[M_u,M_v][+-]$, the integral representations of $(P_{\leq N}^x)^2$ and $\Gamma^{(\Nscript),x}$, and Remark \ref{ansatz:rem-commutativity-chipm}, it holds that
\begin{align*}
&\Big \langle P_{\leq N}^x \UN[K][+], \Gamma^{(\Nscript),x} \Kil P_{\leq N}^x \VN[M_u,M_v][+-] \Big\rangle_\frkg \\
=& \chinull[M_u,M_v] \int_{\R^2} \dw \dz \, \widecheck{\Gamma}^{(\Nscript)}(w) \big( \widecheck{\rho}_{\leq N} \ast \widecheck{\rho}_{\leq N}\big)(z-w) 
\mathcal{T} \Big( P_{\leq N}^x \UN[K][+] \otimes \Theta^x_z P_{\leq N}^x  \IUN[M_u][+] \otimes \Theta^x_z P_{\leq N}^x \VN[M_v][-] \Big),
\end{align*}
where 
\begin{equation}\label{increment:eq-increment-modulated2-p5}
\mathcal{T}\colon \frkg^{\otimes 3} \rightarrow \R, \, (A,B,C) \mapsto \Big \langle A, \Kil \big[ B, C \big] \Big\rangle_\frkg. 
\end{equation}
We now decompose 
\begin{align}
&\mathcal{T} \bigg( P_{\leq N}^x \UN[K][+] \otimes \Theta^x_z P_{\leq N}^x \IUN[M_u][+] \otimes \Theta^x_z P_{\leq N}^x \VN[M_v][-] \bigg) \notag \\ 
=& \,  
\mathcal{T}  \bigg( \Big( \biglcol \, P_{\leq N}^x \UN[K][+] \otimes \Theta^x_z P_{\leq N}^x\IUN[M_u][+]  \bigrcol \Big) \otimes \Theta^x_z P_{\leq N}^x \VN[M_v][-] \bigg) \label{increment:eq-increment-modulated2-p5-nonres} \\
+&\, \coup \mathbf{1} \big\{ K=M_u \big\} \Cf^{(\Ncs)}_{K}(z) \mathcal{T}  \Big( \Cas \otimes \Theta^x_z P_{\leq N}^x \VN[M_v][-] \Big) \label{increment:eq-increment-modulated2-p5-res}. 
\end{align}
The contribution of \eqref{increment:eq-increment-modulated2-p5-nonres} can easily be controlled using  Corollary \ref{killing:cor-tensors} and Lemma \ref{increment:lem-increment-integrand}. Thus, it remains to treat the contribution of the resonant part \eqref{increment:eq-increment-modulated2-p5-res}. 
 We show that, due to the properties of the Lie bracket, the contribution of the resonant part \eqref{increment:eq-increment-modulated2-p5-res} vanishes. To see this, it is sufficient to show that, for all $X\in \frkg$, 
\begin{equation}\label{increment:eq-increment-modulated2-p8}
\mathcal{T} \big( \Cas \otimes X \big) =0. 
\end{equation}
Indeed, using the definition of the Casimir, the definition of $\mathcal{T}$ from \eqref{increment:eq-increment-modulated2-p5},  and Lemma \ref{prelim:lem-Killing-vs-Bracket}, it holds that
\begin{equation*}
    \mathcal{T} \big( \Cas \otimes X \big) = \mathcal{T} \big(  E_a \otimes E^a \otimes X \big) = \Big \langle E_a, \Kil \big[ E^a, X \big] \Big\rangle_\frkg =0. 
\end{equation*} 

\emph{Case 3: $(+)\times (+)-interaction$}. 
We first recall from Definition \ref{ansatz:def-modulated-linear-reversed} that
\begin{equation}\label{increment:eq-increment-modulated2-q7m}
\VN[M][+] = \chinull[M] \Big[ \IUN[M][+], \LON[M][-]  \Big]_{\leq N}
+ \chinull[M] \Big[ (P_{\leq N}^x)^2 \IUN[M][+],  \SHHLN[M][v] \Big]_{\leq N}. 
\end{equation}
We only treat the contribution of the $\LON[M][-]$-term in \eqref{increment:eq-increment-modulated2-q7m}, since the argument for the $\SHHLN[M][v]$-term is similar. 
Using the integral representations of $(P_{\leq N}^x)^2$ and $\Gamma^{(\Nscript),x}$ and using Remark \ref{ansatz:rem-commutativity-chipm}, it then follows that
\begin{equation}\label{increment:eq-increment-modulated2-p7m}
\begin{aligned}
&\Big \langle P_{\leq N}^x \UN[K][+], \Gamma^{(\Nscript),x} \Kil P_{\leq N}^x \chinull[M]
 \Big[ \IUN[M][+], \LON[M][-]  \Big]_{\leq N} \Big\rangle_\frkg \\
=& \chinull[M] \int_{\R^2} \dw \dz \, \bigg(  \widecheck{\Gamma}(w) \big( \widecheck{\rho}_{\leq N} \ast \widecheck{\rho}_{\leq N}\big)(z-w) \\
&\times 
\mathcal{T} \Big( P_{\leq N}^x \UN[K][+] \otimes \Theta^x_z P_{\leq N}^x \IUN[M][+] \otimes \Theta^x_z P_{\leq N}^x \LON[M][-]  \Big) \bigg),
\end{aligned}
\end{equation}
where $\mathcal{T}$ is as in \eqref{increment:eq-increment-modulated2-p5}. Using Definition \ref{killing:def-Wick}, we decompose
\begin{align}
&\mathcal{T} \bigg( P_{\leq N}^x \UN[K][+] \otimes \Theta^x_z P_{\leq N}^x \IUN[M][+] \otimes \Theta^x_z P_{\leq N}^x \LON[M][-]  \bigg) \notag \\
=&\,  \mathcal{T} \bigg( \Big( \biglcol \, P_{\leq N}^x \UN[K][+] \otimes \Theta^x_z P_{\leq N}^x \IUN[M][+] \bigrcol \Big) \otimes \Theta^x_z P_{\leq N}^x   \LON[M][-] \bigg) \label{increment:eq-increment-modulated2-p6} \\ 
+&\, \coup \mathbf{1} \big\{ K=M \big\} \Cf^{(\Ncs)}_M(z) \mathcal{T} \Big( \Cas \otimes \Theta^x_z P_{\leq N}^x  \LON[M][-]  \Big). \label{increment:eq-increment-modulated2-p7}
\end{align}
Exactly as in Case 2, i.e., the $(+)\times (+-)$-interaction, we can show that the resonant term \eqref{increment:eq-increment-modulated2-p7} vanishes. 
It therefore remains to estimate the contribution of non-resonant part \eqref{increment:eq-increment-modulated2-p6}. Using Lemma \ref{killing:lem-tensor-modulated-linear} and Corollary \ref{modulation:cor-LON}, it holds that 
\begin{align*}
&\Big\| \Big( \biglcol \, P_{\leq N}^x \UN[K][+] \otimes \Theta^x_z P_{\leq N}^x \IUN[M][+] \bigrcol \Big) \otimes \Theta^x_z P_{\leq N}^x  \LON[M][-]
\Big\|_{\Cprod{r-1}{r-1}} \\
\lesssim& \, 
\Big\|  \biglcol \, P_{\leq N}^x \UN[K][+] \otimes \Theta^x_z P_{\leq N}^x \IUN[M][+] \bigrcol \Big\|_{\Cprod{r-1}{s}} 
\Big\|  \Theta^x_z P_{\leq N}^x   \LON[M][-]  \Big\|_{\Cprod{s}{r-1}} \\
\lesssim&\, \bigg( \max(K,M)^{r-s} M^{-1/2} + \mathbf{1}\big\{ K=M  \big\} N^{-\delta+\vartheta} \langle N |z| \rangle \bigg) M^{r-s} \Dc^3.
\end{align*}
Since $N^{1-\delta}\lesssim M \lesssim N$, the resulting contribution to \eqref{increment:eq-increment-modulated2-p7m} can easily be bounded by $N^{-\delta+10 \delta_1}$. \\

\emph{Case 4: $(+-)\times (+)-interaction$}. Due to the integral representations of $P_{\leq N}^x$ and $\Gamma^{(\Nscript),x}$ and  Remark \ref{ansatz:rem-commutativity-chipm}, it holds that
\begin{align*}
&\Big \langle P_{\leq N}^x \UN[K_u,K_v][+-], \Gamma^{(\Nscript),x} \Kil P_{\leq N}^x \VN[M][+] \Big\rangle_\frkg \\
=& \chinull[K_u,K_v] \chinull[M] \int_{\R^3} \dw \dy \dz \, \bigg( 
\big( \widecheck{\rho}_{\leq N} \ast \widecheck{\rho}_{\leq N}\big)(y) \widecheck{\Gamma}(w) \big( \widecheck{\rho}_{\leq N} \ast \widecheck{\rho}_{\leq N}\big)(z-w) \\
&\times 
\mathcal{T} \Big( \Theta^x_y P_{\leq N}^x \UN[K_u][+] \otimes \Theta^x_z P_{\leq N}^x \IUN[M][+] \otimes \Theta^x_y P_{\leq N}^x \IVN[K_v][-] \otimes   \Theta^x_z P_{\leq N}^x \big(  \LON[M][-] + \SHHLN[M][v] \big) \Big) \bigg),
\end{align*}
where 
\begin{equation*}
\mathcal{T}\colon \frkg^{\otimes 4}\rightarrow \R, \, A \otimes B \otimes C \otimes D \mapsto \Big \langle \big [ A,C \big], \Kil \big[ B,D \big] \Big \rangle 
\end{equation*}
Using Definition \ref{killing:def-Wick}, we then decompose
\begin{align}
 &\Theta^x_y P_{\leq N}^x \UN[K_u][+] \otimes \Theta^x_z P_{\leq N}^x \IUN[M][+] \otimes \Theta^x_y P_{\leq N}^x \IVN[K_v][-] \otimes   \Theta^x_z P_{\leq N}^x \big(  \LON[M][-] + \SHHLN[M][v] \big)  \notag \\ 
 =&  \Big( \biglcol  \, \Theta^x_y P_{\leq N}^x \UN[K_u][+] \otimes \Theta^x_z P_{\leq N}^x \IUN[M][+] \bigrcol \Big) \otimes \Theta^x_y P_{\leq N}^x \IVN[K_v][-]  
  \otimes   \Theta^x_z P_{\leq N}^x \big(  \LON[M][-] + \SHHLN[M][v] \big) \label{increment:eq-increment-modulated2-p9} \\ 
 +& \coup \mathbf{1} \big\{ K_u = M \big\} \Cf_{K_u}^{(\Ncs)}(z-y) \, 
 \Big( \Cas \otimes  \Theta^x_y P_{\leq N}^x \IVN[K_v][-] \otimes   \Theta^x_z P_{\leq N}^x \big(  \LON[M][-] + \SHHLN[M][v] \big) \Big).
 \label{increment:eq-increment-modulated2-p10}
\end{align}
We first estimate \eqref{increment:eq-increment-modulated2-p9}. Using Lemma \ref{killing:lem-tensor-modulated-linear}, Lemma \ref{modulation:lem-linear}, Lemma \ref{modulation:lem-shhl}, and Corollary \ref{modulation:cor-LON}, it holds that
\begin{align*}
\Big\| \eqref{increment:eq-increment-modulated2-p9} \Big\|_{\Cprod{-s}{-s}} 
&\lesssim \Big\| \biglcol  \, \Theta^x_y P_{\leq N}^x \UN[K_u][+] \otimes \Theta^x_z P_{\leq N}^x \IUN[M][+] \bigrcol \Big\|_{\Cprod{-s+\eta}{s}} \\ 
&\times \Big\| \Theta^x_y P_{\leq N}^x \IVN[K_v][-]   \Big\|_{\Cprod{s}{s}} 
\Big\|  \Theta^x_z P_{\leq N}^x \big(  \LON[M][-] + \SHHLN[M][v] \big) \Big\|_{\Cprod{s}{-s+\eta}} \\ 
&\lesssim \Big( \max(K_u,M)^{1-2s+2\eta} M^{-1/2} + N^{-\delta+\vartheta} \langle N(z-y) \rangle \Big)  M^{1-2s+\eta} \Dc^4. 
\end{align*}
Since $N^{1-\delta}\lesssim M \lesssim N$, this yields an acceptable contribution. It remains to estimate \eqref{increment:eq-increment-modulated2-p10}. 
To this end, we first use Lemma \ref{ansatz:lem-renormalization}, which yields that $|\Cf_{K_u}^{(\Ncs)}(z-y)|\lesssim 1$. 
We further observe that the conditions $K_u=M$ and $K_u \simeq_\delta K_v$ imply that $K_v \geq M^{1-\delta}$. Using Lemma \ref{modulation:lem-linear}, Lemma \ref{modulation:lem-shhl}, and Lemma \ref{modulation:lem-ivm-lo}, it then follows that
\begin{align*}
\Big\| \eqref{increment:eq-increment-modulated2-p10} \Big\|_{\Cprod{-s}{-s}} 
\lesssim\, & 
 \coup \Big\| \Theta^x_y P_{\leq N}^x \IVN[K_v][-] \otimes   \Theta^x_z P_{\leq N}^x\LON[M][-]  \Big\|_{\Cprod{-s}{-s}} \\
 + \, & \coup 
 \Big\|  \Theta^x_y P_{\leq N}^x \IVN[K_v][-]  \Big\|_{\Cprod{s}{\eta}}
 \Big\|  \Theta^x_z P_{\leq N}^x  \SHHLN[M][v] \Big\|_{\Cprod{s}{\eta}}\\
 \lesssim\, &  K_v^{-1/2+\delta} \lambda \Dc^2. 
\end{align*}
Since $\coup \lesssim \Dc^2$ and $M=K_u \simeq_\delta K_v$, this is easily seen to yield an acceptable contribution. 
\end{proof}

Equipped with Lemma \ref{increment:lem-increment-object-I} and Lemma \ref{increment:lem-increment-object-II}, it remains to control interactions that involve at least one mixed modulated object or smooth remainder. To this end, we first prove the following auxiliary lemma. 
 
\begin{lemma}\label{increment:lem-aux-ast} Assume that the post-modulation hypothesis (Hypothesis \ref{hypothesis:post}). Furthermore,  let $K,M\in \Dyadiclarge$ and let 
\begin{equation*}
\UN[][\ast] \in \big\{ \UN[K][+\fs], \UN[K][\fs-], \UN[][\fs] \big\} 
\qquad \text{and} \qquad
\VN[][\ast] \in \big\{ \VN[M][\fs-], \VN[M][+\fs], \VN[][\fs] \big\}. 
\end{equation*}
Then, it holds that
\begin{align}
\Big\| P^{u,v}_{\geq N^{1-\delta}} \UN[][\ast] \Big\|_{\Cprod{-s}{1/2}} 
&\lesssim N^{-(1-4\delta) \delta_1} \Dc, \label{increment:eq-aux-ast-e1} \\ 
\Big\| P^{u,v}_{\geq N^{1-\delta}} \VN[][\ast] \Big\|_{\Cprod{1/2}{-s}} 
&\lesssim N^{-(1-4\delta) \delta_1} \Dc. \label{increment:eq-aux-ast-e2} 
\end{align}
\end{lemma}

\begin{proof}
Due to symmetry in the $u$ and $v$-variables, it suffices to prove \eqref{increment:eq-aux-ast-e2}. In order to prove \eqref{increment:eq-aux-ast-e2}, we distinguish the three different options for $\VN[][\ast]$. \\ 

\emph{Estimate of $\VN[M][\fs-]$:} 
Since $\VN[M][\fs-]$ is supported on $u$-frequencies $\gtrsim M^{1-\deltap}$ and $v$-frequencies $\sim M$, the term $P^{u,v}_{\geq N^{1-\delta}}\VN[M][\fs-]$ is supported on $u$-frequencies $\gtrsim N^{(1-\deltap)(1-\delta)}$ and $v$-frequencies $\lesssim N$. It then follows that
\begin{equation}\label{increment:eq-aux-ast-p1}
\Big\| P^{u,v}_{\geq N^{1-\delta}}\VN[M][\fs-] \Big\|_{\Cprod{1/2}{-s}} \lesssim N^{-(1-\deltap)(1-\delta)(r-1/2)} N^{1-2s} 
\Big\| \VN[M][\fs-] \Big\|_{\Cprod{r}{s-1}}.
\end{equation}
To simplify the pre-factor, we note that 
\begin{equation*}
-(1-\deltap)(1-\delta)(r-1/2) + 1-2s = -(1-2\delta) \delta_1 + \mathcal{O}(\delta^2 \delta_1). 
\end{equation*}
Together with Lemma \ref{modulation:lem-mixed}, it follows that \eqref{increment:eq-aux-ast-p1} can be controlled by $N^{-(1-4\delta)\delta_1}\Dc^2$, which is acceptable. \\  

\emph{Estimate of $\VN[M][+\fs]$:} 
Since $\VN[M][+\fs]$ is supported on $u$-frequencies $\sim M$ and $v$-frequencies $\gtrsim M^{1-\deltap}$, the term $P^{u,v}_{\geq N^{1-\delta}}\VN[M][+\fs]$ is supported on $u$-frequencies $\lesssim N$ and  $v$-frequencies $\gtrsim N^{(1-\deltap)(1-\delta)}$. It then follows that
\begin{equation*}
\Big\| P^{u,v}_{\geq N^{1-\delta}}\VN[M][+\fs] 
\Big\|_{\Cprod{1/2}{-s}} 
\lesssim N^{1-2s}  N^{-(1-\deltap)(1-\delta)(r-1/2)} 
\Big\| \VN[M][+\fs] \Big\|_{\Cprod{s}{r-1}}.
\end{equation*}
Since the pre-factor is as in \eqref{increment:eq-aux-ast-p1} and the $\Cprod{r}{s-1}$-norm of $\VN[M][\fs-]$ can be controlled using Lemma \ref{modulation:lem-mixed}, this is acceptable.\\

\emph{Estimate of $\VN[M][\fs]$:} 
It clearly holds that
\begin{equation*}
  \Big\| P^{u,v}_{\geq N^{1-\delta}}\VN[][\fs] 
\Big\|_{\Cprod{1/2}{-s}}   
\lesssim N^{-(1-\delta) (r-1+s)} \Big\| \VN[][\fs] \Big\|_{\Cprod{r}{r-1}}. 
\end{equation*}
Using that $-(1-\delta)(r-1+s)=-(1-\delta)\delta_1 + \mathcal{O}(\delta_2)$ and using Hypothesis \ref{hypothesis:post} to control $\VN[][\fs]$, this yields an acceptable contribution. 
\end{proof}

In the next lemma, we treat the remaining cases in the first two rows of Figure \ref{figure:increment-cases}. 

\begin{lemma}\label{increment:lem-increment-fs} 
Assume that the post-modulation hypothesis (Hypothesis \ref{hypothesis:post}) is satisfied. Furthermore, let 
 $K,L,M \in \Dyadiclarge$ satisfy $K\simeq_\delta L$ and let 
\begin{equation}\label{increment:eq-increment-fs-cases}
\UN[][\ast] \in \Big\{ \UN[K][+], \UN[K,L][+-]\Big\} 
\qquad \text{and} \qquad 
\VN[][\ast] \in \Big\{ \VN[M][\fs-], \VN[M][+\fs], \VN[][\fs] \Big\}. 
\end{equation}
Then, it holds that
\begin{equation*}
\sup_{\zeta \in \Cuttilde} \sup_{\tau \in [-1,1]} \Big| \EIN \Big( \UN[][\ast],\VN[][\ast];\tau \Big) \Big| \lesssim R N^{-(1-8\delta)\delta_1} \Dc^3. 
\end{equation*}
\end{lemma}

\begin{proof}
We decompose 
\begin{equation}\label{increment:eq-increment-fs-p1}
\begin{aligned}
& \EIN \Big( \UN[][\ast],\VN[][\ast];\tau  \Big) \\ 
=\, & \EIN \Big( \UN[][\ast], P^{u,v}_{< N^{1-\delta}} \VN[][\ast];\tau \Big)
+  \EIN \Big( \UN[][\ast], P^{u,v}_{\geq N^{1-\delta}} \VN[][\ast];\tau \Big). 
\end{aligned}
\end{equation}
The first summand in \eqref{increment:eq-increment-fs-p1} can be easily controlled by using Lemma \ref{increment:lem-low-frequency} and our estimates for the terms in \eqref{increment:eq-increment-fs-cases}, i.e., Lemma \ref{modulation:lem-bilinear}, Lemma \ref{modulation:lem-mixed}, and Lemma \ref{modulation:lem-linear-reversed}. 
Thus, it remains to control the second summand in \eqref{increment:eq-increment-fs-p1}. Due to Lemma \ref{increment:lem-increment-integrand}, it suffices to prove that
\begin{equation}\label{increment:eq-increment-fs-p2}
\Big\| \Big\langle P_{\leq N}^x \UN[][\ast], \Gamma^{(\Nscript),x} \Kil P_{\leq N}^x P^{u,v}_{\geq N^{1-\delta}} \VN[][\ast] \Big\rangle_\frkg \,  \Big\|_{\Cprod{-s}{-s}} 
\lesssim N^{-(1-8\delta)\delta_1} \Dc^3. 
\end{equation} 
In order to prove \eqref{increment:eq-increment-fs-p2}, we first use our product estimate (Lemma \ref{prelim:lem-paraproduct}), which yields that
\begin{align*}
\Big\| \Big\langle P_{\leq N}^x \UN[][\ast], \Gamma^{(\Nscript),x} \Kil P_{\leq N}^x P^{u,v}_{\geq N^{1-\delta}} \VN[][\ast] \Big\rangle_\frkg \,  \Big\|_{\Cprod{-s}{-s}} 
\lesssim \big\| \UN[][\ast] \big\|_{\Cprod{-s}{1/2}} 
\big\| P^{u,v}_{\geq N^{1-\delta}}  \VN[][\ast] \big\|_{\Cprod{1/2}{-s}}. 
\end{align*}
For $\UN[][\ast]$ as in \eqref{increment:eq-increment-fs-cases}, it follows from Lemma \ref{modulation:lem-bilinear} and Lemma \ref{modulation:lem-linear-reversed} that
\begin{equation*}
\big\| \UN[][\ast] \big\|_{\Cprod{-s}{1/2}}  \lesssim 
N^{-s-(s-1)} N^{\frac{1}{2}-s} \Dc^2 = N^{3\delta_2} \Dc^2.   
\end{equation*}
Together with Lemma \ref{increment:lem-aux-ast}, it then follows that
\begin{equation*}
    \big\| \UN[][\ast] \big\|_{\Cprod{-s}{1/2}} 
\big\| P^{u,v}_{\geq N^{1-\delta}}  \VN[][\ast] \big\|_{\Cprod{1/2}{-s}} 
\lesssim N^{3\delta_2} N^{-(1-4\delta)\delta_1} \Dc^3, 
\end{equation*}
which is acceptable. 
\end{proof}

In the last lemma of this subsection, we treat all cases from Figure \ref{figure:increment-cases} which involve neither the modulated linear wave $\UN[][+]$ nor the modulated bilinear wave $\UN[][+-]$. 

\begin{lemma}\label{increment:lem-increment-remaining}
Assume that the post-modulation hypothesis (Hypothesis \ref{hypothesis:post}) is satisfied. Furthermore, let $K,M\in \Dyadiclarge$ and let
\begin{equation*}
\UN[][\ast] \in \big\{ \UN[K][-],\UN[K][+\fs], \UN[K][\fs-], \UN[][\fs] \big\} 
\quad \text{and} \quad
\VN[][\ast] \in \big\{ \VN[M][+], \VN[M][\fs-], \VN[M][+\fs], \VN[][\fs] \big\}. 
\end{equation*}
Then, it holds that 
\begin{equation*}
\sup_{\zeta \in \Cuttilde} \sup_{\tau\in [-1,1]} \Big| \EIN \Big( \UN[][\ast],\VN[][\ast];\tau \Big) \Big| \lesssim R N^{-(1-8\delta) \delta_1} \Dc^2. 
\end{equation*}
\end{lemma}

\begin{proof}
We first decompose 
\begin{align}
&\EIN \Big( \UN[][\ast],\VN[][\ast];\tau \Big) \notag  \\ 
=\,&  \bigg( \EIN \Big( \UN[][\ast],\VN[][\ast];\tau \Big) - 
\EIN \Big( P^{u,v}_{>N^{1-\delta}} \UN[][\ast],P^{u,v}_{>N^{1-\delta}} \VN[][\ast];\tau\Big) \bigg) 
\label{increment:eq-remaining-p1} \\ 
+\,& \EIN \Big( P^{u,v}_{>N^{1-\delta}} \UN[][\ast],P^{u,v}_{>N^{1-\delta}} \VN[][\ast];\tau \Big).
\label{increment:eq-remaining-p2}
\end{align}
The first summand \eqref{increment:eq-remaining-p1} can be controlled using Lemma \ref{increment:lem-low-frequency}, and it remains to control the second summand \eqref{increment:eq-remaining-p2}. To this end, we first use Lemma \ref{increment:lem-increment-integrand}, which yields that 
\begin{equation}\label{increment:eq-remaining-p2p}
\begin{aligned}
&\bigg| \EIN \Big( P^{u,v}_{>N^{1-\delta}} \UN[][\ast],P^{u,v}_{>N^{1-\delta}} \VN[][\ast];\tau \Big) \bigg| \\
\lesssim&\, R \, \Big\| 
\Big\langle P_{\leq N}^x P^{u,v}_{>N^{1-\delta}} \UN[][\ast], 
\Gamma^{(\Nscript),x} \Kil P_{\leq N}^x P^{u,v}_{>N^{1-\delta}} \VN[][\ast] \Big\rangle_\frkg  
\Big\|_{\Cprod{-s}{-s}}. 
\end{aligned}
\end{equation}
We now first consider the special case
\begin{equation}\label{increment:eq-remaining-p3}
\big( \UN[][\ast], \VN[][\ast] \big) = \big( \UN[K][-], \VN[M][+] \big). 
\end{equation}
By symmetry, we may assume that $K\leq M$. 
Due to Lemma \ref{ansatz:lem-frequency-support}, 
$\UN[K][-]$ is supported on $u$-frequencies $\lesssim K^{1-\delta+\vartheta}$ and $\VN[M][+]$ is supported on $u$-frequencies $\sim M$, and thus the argument in \eqref{increment:eq-remaining-p2p} only contains low$\times$high-interactions in the $u$-variable. Using our low$\times$high-estimate (Lemma \ref{prelim:lem-paraproduct}), it follows that 
\begin{equation}\label{increment:eq-remaining-p4}
\eqref{increment:eq-remaining-p2p} 
\lesssim R \big\|  P^{u,v}_{>N^{1-\delta}} \UN[K][-] \big\|_{\Cprod{-s}{1/2}} 
\big\| P^{u,v}_{>N^{1-\delta}} \VN[M][+]  \big\|_{\Cprod{\eta}{-s}}. 
\end{equation}
Due to the frequency-support properties of $\UN[K][-]$ and $\VN[M][+]$, the right-hand side of \eqref{increment:eq-remaining-p4} is non-trivial only when $K,M\gtrsim N^{1-\delta}$. Furthermore, using Lemma \ref{modulation:lem-linear-reversed}, 
\begin{align*}
\big\|  P^{u,v}_{>N^{1-\delta}} \UN[K][-] \big\|_{\Cprod{-s}{1/2}}  
&\lesssim K^{1-2s} K^{1/2-s} \big\| \UN[K][-] \big\|_{\Cprod{s}{s-1}} \lesssim N^{3\delta_2} \Dc^2, \\ 
\big\| P^{u,v}_{>N^{1-\delta}} \VN[M][+]  \big\|_{\Cprod{\eta}{-s}} 
&\lesssim M^{\eta-s} M^{1-2s} \big\| \VN[M][+]\big\|_{\Cprod{s}{s-1}} 
\lesssim N^{-\frac{1}{2}+2\delta} \Dc^2. 
\end{align*}
By inserting this back into \eqref{increment:eq-remaining-p4}, we obtain a more than acceptable contribution and therefore complete the case \eqref{increment:eq-remaining-p3}. Due to symmetry, it now remains to treat the cases in which
\begin{equation}\label{increment:eq-remaining-p5}
\VN[][\ast] \in \big\{  \VN[M][\fs-], \VN[M][+\fs], \VN[][\fs] \big\}. 
\end{equation}
Using our paraproduct estimate (Lemma \ref{prelim:lem-paraproduct}) and Lemma \ref{increment:lem-aux-ast}, it follows that
\begin{align*}
\eqref{increment:eq-remaining-p2p} 
\lesssim R \big\| P^{u,v}_{>N^{1-\delta}} \UN[][\ast] \big\|_{\Cprod{-s}{1/2}} 
\big\| P^{u,v}_{>N^{1-\delta}}  \VN[][\ast] \big\|_{\Cprod{1/2}{-s}} \lesssim R N^{2(1-2s)} N^{-(1-4\delta)\delta_1} \Dc^2.  
\end{align*}
Since $2(1-2s)=4\delta_2$, this is acceptable. 
\end{proof}

Equipped with all lemmas above, we can now prove the main proposition of this section. 

\begin{proof}[Proof of Proposition \ref{increment:prop-energy-increment}:] 
Using our Ansatz from \eqref{ansatz:eq-UN-decomposition} and \eqref{ansatz:eq-VN-decomposition}, we decompose 
\begin{equation}\label{increment:eq-energy-increment-p1}
 \Big| \EIN\big(\UN[][],\VN[][];\tau\big) \Big|  
 \leq \sum_{\ast_1,\ast_2}  \Big| \EIN\big(\UN[][\ast_1],\VN[][\ast_2];\tau\big) \Big|.
\end{equation}
In \eqref{increment:eq-energy-increment-p1}, the sum is taken over all $\ast_1,\ast_2\in \{+,-,+-,+\fs,\fs-,\fs\}$, which correspond to all cells in Figure \ref{figure:increment-cases}. Using the symmetry of our estimates in the $u$ and $v$-variables, it suffices to treat the cells in Figure \ref{figure:increment-cases} which are on or above the diagonal. For each of these cells, the contribution to \eqref{increment:eq-energy-increment-p1} can be controlled using the lemma listed in the cell, i.e., using Lemma \ref{increment:lem-increment-object-I}, Lemma \ref{increment:lem-increment-object-II}, Lemma \ref{increment:lem-increment-fs}, or Lemma \ref{increment:lem-increment-remaining}. 
\end{proof}
\section{Lipschitz-continuous dependence}\label{section:Lipschitz}

In Sections \ref{section:Killing}-\ref{section:renormalization}, we obtained several estimates of our modulated and mixed modulated
objects, the terms in the modulation equations, and the terms in the remainder equations. In Section \ref{section:main}, however, we not only need estimates of all of these terms, but also need their Lipschitz-continuous dependence on the modulation operators and nonlinear remainders. In this section, we therefore state Lipschitz-variants of our earlier estimates. While the statements themselves will be complete and precise, many of the proofs will only be sketched, since they consist of minor modifications of our arguments from Sections \ref{section:Killing}-\ref{section:renormalization}.

Throughout the rest of this section, we let 
\begin{align}
&\big( \Nd, \SNin[][\pm] , \SN[][\pm] , \UN[][\fs], \VN[][\fs] \big) \label{lipschitz:eq-original}\\ 
\qquad \text{and} \qquad &\big( \Ndtil, \SNtilin[][\pm] , \SNtil[][\pm] , \UNtil[][\fs], \VNtil[][\fs] \big)
\label{lipschitz:eq-tilde}
\end{align} 
be two collections of frequency-truncation parameters, initial modulation operators (as in Definition \ref{ansatz:def-initial-data}), modulation operators, and nonlinear remainders. In all results below, we will assume that either the pre-modulation hypothesis (Hypothesis \ref{hypothesis:pre}) or post-modulation hypothesis (Hypothesis \ref{hypothesis:post}) is satisfied for both \eqref{lipschitz:eq-original} and 
\eqref{lipschitz:eq-tilde}. Previously, all objects in our Ansatz, such as the modulated linear waves $\UN[K][+]$, were defined using \eqref{lipschitz:eq-original}. In the following, the corresponding objects for \eqref{lipschitz:eq-tilde} will be denoted using similar notation, but with an additional tilde. For example, we write\footnote{For typographical reasons, the $\mathscr{N}$-superscript in $\UNtil[K][+]$ is not changed to $\widetilde{\mathscr{N}}$, and instead kept in its original form.}
\begin{equation*}
\UNtil[K][+] = \hcoup \sum_{u_0 \in \LambdaRR} \sum_{k\in \Z_K} \psiRuK \rhoNDtil(k) \SNtil[K][+][k] G_{u_0,k}^+ e^{iku}.
\end{equation*}
In order to state our Lipschitz-estimates, it is convenient to make the following definitions.

\begin{definition}[Distances]\label{lipschitz:def-distances}
We first define a distance on dyadic integers by 
\begin{equation*}
\ddyadic \big( \Nd,\Ndtil \big) 
= 
\begin{cases}
\begin{tabular}{cl}
$0$ & if $\Nd = \Ndtil$, \\
$\min\big( \Nd,\Ndtil\big)^{-\varsigma}$  & if  $\Nd \neq \Ndtil$. 
\end{tabular}
\end{cases}
\end{equation*}
We then define the pre-modulation distance of \eqref{lipschitz:eq-original} and \eqref{lipschitz:eq-tilde} by 
\begin{equation}\label{lipschitz:eq-dpre}
\begin{aligned}
\dpre &:= \Dc\,  \mathsf{d}_{\dyadic}\big( \Nd, \Ndtil \big) 
+ \hcoup \Ac \sup_{K\in \Dyadiclarge} \big\| \pSN[K][+][k] - \pSNtil[K][+][k] \big\|_{\Wuv[s][s]}
+ \hcoup \Ac \sup_{M\in \Dyadiclarge} \big\| \pSN[M][-][m] - \pSNtil[M][-][m] \big\|_{\Wuv[s][s]}\\ 
&+ \big\| \UN[][\fs] - \UNtil[][\fs] \big\|_{\Cprod{r-1}{r}}
+ \big\| \VN[][\fs] - \VNtil[][\fs] \big\|_{\Cprod{r}{r-1}}. 
\end{aligned}
\end{equation}
Finally, we define the post-modulation distance of \eqref{lipschitz:eq-original} and \eqref{lipschitz:eq-tilde} by 
\begin{equation}\label{lipschitz:eq-dpost}
\begin{aligned}
\dpost &:= \Dc\,  \mathsf{d}_{\dyadic}\big( \Nd, \Ndtil \big) 
+ \hcoup \Ac \sup_{K\in \dyadic} \big\| \SNin[K][+] - \SNtilin[K][+] \big\|_{\C_u^s}
+ \hcoup \Ac \sup_{M\in \dyadic} \big\| \SNin[M][-] - \SNtilin[M][-] \big\|_{\C_v^s}\\ 
&+ \big\| \UN[][\fs] - \UNtil[][\fs] \big\|_{\Cprod{r-1}{r}}
+ \big\| \VN[][\fs] - \VNtil[][\fs] \big\|_{\Cprod{r}{r-1}}. 
\end{aligned}
\end{equation}
\end{definition}

In both \eqref{lipschitz:eq-dpre} and \eqref{lipschitz:eq-dpost}, the pre-factors are chosen such that the upper bounds on the nonlinear remainders and modulation operators differ by a factor of $\hcoup \Ac$, which is consistent with the bounds in Hypothesis \ref{hypothesis:pre}. In Lemma \ref{lipschitz:lem-distance} below, we will show that the distances $\dpre$ and $\dpost$ are related via the inequality $\dpre \lesssim \dpost$.

\begin{remark}\label{lipschitz:rem-kernel-rho}
In order to obtain the Lipschitz-continuous dependence on $\Nd$, we estimate terms which contain the difference $\rho_{\leq \Nd}-\rho_{\leq \Ndtil}$. Similar estimates should remain valid for $\rho_{\leq \Nd}-\widetilde{\rho}_{\leq \Nd}$, where $\widetilde{\rho}$ is any symbol satisfying the same properties as $\rho$ (see Section \ref{section:parameters}). Due to this, the limits from Theorem \ref{intro:thm-rigorous-A-B} and Theorem
\ref{intro:thm-rigorous-phi} should not depend on the precise choice of $\rho$, which was previously mentioned in Remark \ref{intro:rem-properties}. 
\end{remark}

\subsection{Modulated and mixed modulated objects}\label{section:lipschitz-modulated}
In this subsection, we obtain Lipschitz-variants of our estimates of the modulated and mixed modulated objects. 
In the following lemma, we record the Lipschitz-variants of Lemma \ref{modulation:lem-linear}, Lemma \ref{modulation:lem-bilinear}, Lemma \ref{modulation:lem-mixed}, Lemma \ref{modulation:lem-shhl}, Lemma \ref{modulation:lem-linear-reversed}, and Corollary \ref{modulation:cor-LON}.

\begin{lemma}[Lipschitz-estimates of modulated and mixed modulated objects]\label{lipschitz:lem-modulated-mixed}
Assume that the pre-modulation hypothesis (Hypothesis \ref{hypothesis:pre}) is satisfied for \eqref{lipschitz:eq-original} and \eqref{lipschitz:eq-tilde} and let $K,M\in \Dyadiclarge$. Then, we have the following Lipschitz estimates:
\begin{align}
\big\| \UN[K][+] - \UNtil[K][+] \big\|_{\Cprod{\gamma}{s}} 
&\lesssim K^{\gamma+\frac{1}{2}+\eta} \dpre,  \label{lipschitz:eq-modulated-p} \\ 
\big\| \IUN[K][+]- \IUNtil[K][+] \big\|_{\Cprod{\gamma}{s}} 
&\lesssim  K^{\gamma-\frac{1}{2}+\eta} \dpre, \label{lipschitz:eq-modulated-ip} \\ 
\big\| \UN[K,M][+-] - \UNtil[K,M][+-] \big\|_{\Cprod{s-1}{s}} 
&\lesssim  (KM)^{-\eta} \Dc \dpre, \label{lipschitz:eq-modulated-pm-1} \\ 
 \big\| \UN[K,M][+-] - \UNtil[K,M][+-]\big\|_{\Cprod{\scrr-1}{\eta}}
 &\lesssim  (KM)^{-\eta} \Dc \dpre, \label{lipschitz:eq-modulated-pm-2} \\ 
\big\| \UN[K,M][+-] - \UNtil[K,M][+-] \big\|_{\Cprod{\alpha}{\beta}} 
&\lesssim K^{\alpha+1/2+\eta} M^{\beta-1/2+\eta} \Dc \dpre, \label{lipschitz:eq-modulated-pm-3} \\
\big\| \UN[M][\fs-] - \UNtil[M][\fs-] \big\|_{\Cprod{r-1}{s}} 
&\lesssim M^{-\eta} \Dc \dpre, \label{lipschitz:eq-modulated-sm} \\
\big\| \UN[K][+\fs] - \UNtil[K][+\fs] \big\|_{\Cprod{s-1}{r}} 
&\lesssim K^{-\eta} \Dc \dpre , \label{lipschitz:eq-modulated-ps} \\
 \big\| \SHHLN[M][u] - \SHHLNtil[M][u]\big\|_{\Cprod{\eta}{s}} 
 &\lesssim M^{-\eta} \Dc \dpre, \label{lipschitz:eq-modulated-shhl} \\
 \Big\| \, \UN[M][-] - \UNtil[M][-] \Big\|_{\Cprod{-1/2+\eta}{s}} 
 &\lesssim  M^{-\delta (\frac{1}{2}-s)+10\eta} \Dc \dpre ,\label{lipschitz:eq-modulated-m} \\
\big\| \LON[M][+] - \LONtil[M][+]  \big\|_{\Cprod{s-1}{s}} &\lesssim \dpre. \label{lipschitz:eq-modulated-lo}
\end{align}
Similar estimates also hold for the corresponding $\VN[][\ast]$-objects.
\end{lemma}

\begin{remark}
There are also Lipschitz-variants of the estimates from Lemma \ref{modulation:lem-PNX-modulated}, Lemma \ref{modulation:lem-integration}, Lemma \ref{modulation:lem-integration-bilinear}, and Corollary \ref{modulation:cor-control-combined}, but we do not explicitly record them here.
\end{remark}

\begin{proof}
Since the estimates can be obtained using similar arguments as in the proofs of Lemma \ref{modulation:lem-linear}, \mbox{Lemma \ref{modulation:lem-bilinear}}, Lemma \ref{modulation:lem-mixed}, Lemma \ref{modulation:lem-shhl}, Lemma \ref{modulation:lem-linear-reversed}, and Corollary \ref{modulation:cor-LON}, we only sketch the proof. To obtain \eqref{lipschitz:eq-modulated-p}, we write
\begin{align}
\UN[K][+] - \UNtil[K][+] 
&= \hcoup \sum_{u_0 \in \LambdaRR} \sum_{k\in \Z_K} \psiRuK \Big( \rhoND(k)  \SN[K][+][k] - \rhoNDtil(k) \SNtil[K][+][k] \Big) G^+_{u_0,k} e^{iku} 
\label{lipschitz:eq-modulated-q1}.
\end{align}
Similar as in the proof of Lemma \ref{modulation:lem-linear}, we obtain from Lemma \ref{prelim:lem-psi-sum} and Hypothesis \ref{hypothesis:probabilistic}.\ref{ansatz:item-hypothesis-linear} that
\begin{equation}\label{lipschitz:eq-modulated-q2}
\big\| \eqref{lipschitz:eq-modulated-q1} \big\|_{\Cprod{\gamma}{s}} 
\lesssim K^{\gamma+\frac{1}{2}+\frac{\eta}{2}} \hcoup \Ac \Big\|  \rhoND(k)  \SN[K][+][k] - \rhoNDtil(k) \SNtil[K][+][k] \Big\|_{\Wuv[s][s]}.
\end{equation}
We now note that the difference of $\rhoND(k)$ and $\rhoNDtil(k)$ can only make a non-zero contribution in the case $\Nd\neq \Ndtil$ and $K\gtrsim \min(\Nd,\Ndtil)$. 
Using Lemma \ref{chaos:lem-wc}, we then obtain that
\begin{align*}
\eqref{lipschitz:eq-modulated-q2} 
&\lesssim K^{\gamma+\frac{1}{2}+\frac{\eta}{2}} \hcoup \Ac \,  
\Big( \mathbf{1}\big\{ \Nd \neq \Ndtil, \, K \gtrsim \min(\Nd,\Ndtil) \big\} \big\| \SN[K][+][k] \big\|_{\Wuv[s][s]}
+ \big\| \SN[K][+][k] - \SNtil[K][+][k] \big\|_{\Wuv[s][s]} \Big).  
\end{align*}
Since 
\begin{equation*}
K^{\frac{\eta}{2}} \mathbf{1}\big\{ \Nd \neq \Ndtil, \, K \gtrsim \min(\Nd,\Ndtil) \big\} \lesssim K^{\eta} \ddyadic \big( \Nd, \Ndtil \big), 
\end{equation*}
this yields an acceptable contribution. The estimates in \eqref{lipschitz:eq-modulated-ip}, \eqref{lipschitz:eq-modulated-pm-1}, \eqref{lipschitz:eq-modulated-pm-2}, and \eqref{lipschitz:eq-modulated-pm-3} can be obtained using similar modifications as in the proof of \eqref{lipschitz:eq-modulated-p}. 
In order to obtain \eqref{lipschitz:eq-modulated-sm}, we use similar arguments as in \eqref{mixed:eq-control-p1} and \eqref{mixed:eq-control-p2}, which yield that 
\begin{equation}\label{lipschitz:eq-modulated-q3}
\begin{aligned}
&\, \sup_{M\in \Dyadiclarge} M^\eta \big\| \UN[M][\fs-] - \UNtil[M][\fs-] \big\|_{\Cprod{r-1}{s}}\\
\lesssim&\,  \Dc\,  \big\| \UN[][\fs]-\UNtil[][\fs] \big\|_{\Cprod{r-1}{r}}
+ \Dc \sup_{M\in \Dyadiclarge} M^\eta \big\| \IVN[M][-] - \IVNtil[M][-] \big\|_{\Cprod{s}{s}}.
\end{aligned}
\end{equation}
Together with the analogue of \eqref{lipschitz:eq-modulated-ip} for the difference of $\IVN[M][-]$ and $\IVNtil[M][-]$, this implies \eqref{lipschitz:eq-modulated-sm}. The remaining estimates in 
\eqref{lipschitz:eq-modulated-ps}, 
\eqref{lipschitz:eq-modulated-shhl},
\eqref{lipschitz:eq-modulated-m}, 
and \eqref{lipschitz:eq-modulated-lo}
follow from similar modifications. 
\end{proof}

We also record Lipschitz-variants of the estimates from Lemma \ref{modulation:lem-Cartesian}.

\begin{lemma}[Lipschitz-estimates in Cartesian coordinates]\label{lipschitz:lem-modulated-cartesian}
Assume that the pre-modulation hypothesis (Hypothesis \ref{hypothesis:pre}) is satisfied for \eqref{lipschitz:eq-original} and \eqref{lipschitz:eq-tilde} and let $K,M\in \Dyadiclarge$. Then, it holds that:
\begin{alignat}{3}
K^\eta \, \big\|\UN[K][+] - \UNtil[K][+] \big\|_{\Ctx{s-1}} &\lesssim  \dpre ,  
\quad& \quad 
M^\eta \, \big\|\VN[M][-]  - \VNtil[M][-] \big\|_{\Ctx{s-1}} &\lesssim  \dpre,  
\label{lipschitz:eq-modulated-cartesian-e1} \\  
(KM)^\eta \, \big\|\UN[K,M][+-] - \UNtil[K,M][+-]  \big\|_{\Ctx{\scrr-1}} &\lesssim \Dc \dpre , 
\quad& \quad 
(KM)^\eta \, \big\|\VN[K,M][+-] - \VNtil[K,M][+-] \big\|_{\Ctx{\scrr-1}} &\lesssim \Dc \dpre, 
\label{lipschitz:eq-modulated-cartesian-e2} \\ 
M^\eta \, \big\|\UN[M][-] - \UNtil[M][-] \big\|_{\Ctx{\scrr-1}} &\lesssim  \Dc \dpre, 
\quad& \quad 
K^\eta \, \big\|\VN[K][+] - \VNtil[K][+] \big\|_{\Ctx{\scrr-1}} &\lesssim  \Dc \dpre, 
\label{lipschitz:eq-modulated-cartesian-e3} \\ 
K^\eta \, \big\|\UN[K][+\fs] - \UNtil[K][+\fs] \big\|_{\Ctx{\scrr-1}} &\lesssim \Dc \dpre,
\quad& \quad 
M^\eta \, \big\|\VN[M][\fs-] - \VNtil[M][\fs-] \big\|_{\Ctx{\scrr-1}} &\lesssim \Dc \dpre, 
\label{lipschitz:eq-modulated-cartesian-e4} \\ 
M^\eta \, \big\|\UN[M][\fs-] - \UNtil[M][\fs-] \big\|_{\Ctx{r-1}} &\lesssim \Dc \dpre, 
\quad& \quad 
K^\eta \, \big\|\VN[K][+\fs] - \VNtil[K][+\fs] \big\|_{\Ctx{r-1}} &\lesssim \Dc \dpre, 
\label{lipschitz:eq-modulated-cartesian-e5} \\ 
\, \big\|\UN[][\fs] - \UNtil[][\fs]\big\|_{\Ctx{r-1}} &\lesssim  \dpre,  
\quad& \quad 
\, \big\| \VN[][\fs] - \VNtil[][\fs]\big\|_{\Ctx{r-1}} &\lesssim  \dpre. 
\label{lipschitz:eq-modulated-cartesian-e6}
\end{alignat}
\end{lemma}

\begin{remark}\label{lipschitz:rem-cartesian}
Similar Lipschitz-variants can be derived for Lemma \ref{modulation:lem-Cartesian-frequency-boundary} and Lemma \ref{modulation:lem-Cartesian-high-high}, but we do not explicitly record them here.
\end{remark}

\begin{proof}
The estimates in \eqref{lipschitz:eq-modulated-cartesian-e1} and \eqref{lipschitz:eq-modulated-cartesian-e3}-\eqref{lipschitz:eq-modulated-cartesian-e6} can easily be derived by combining the argument from the proof of Lemma \ref{modulation:lem-Cartesian} with a Lipschitz-estimate from Lemma \ref{lipschitz:lem-modulated-mixed}, and we omit the details. The proof of \eqref{lipschitz:eq-modulated-cartesian-e2} does not rely on Lemma \ref{lipschitz:lem-modulated-mixed}, but can be obtained from a modification of the proof of \mbox{Lemma \ref{modulation:lem-Cartesian}}, which we now sketch. Due to symmetry, it suffices to control the difference of $\VN[K,M][+-]$ and $\VNtil[K,M][+-]$. Similar as in the proof of Lemma \ref{modulation:lem-Cartesian}, we obtain the identity
\begin{align}
&\, \VN[K,M][+-] - \VNtil[K,M][+-] \label{lipschitz:eq-cartesian-proof-1} \\  
=&\, \coup \chinull[K,M]  P_{\leq N}^x \sum_{u_0,v_0\in \LambdaRR} \int_{\R^2} \dy \dz \, \bigg( \widecheck{\rho}_{\leq N}(y) \widecheck{\rho}_{\leq N}(z)  \psiRuK(x-y-t) \psiRvM(x-y+t) \notag \\ 
\times& \bigg( \Big[ \Theta^x_y  \Big( \hspace{-0.2ex} \sum_{k\in \Z_K} \hspace{-1ex} \rhoND(k) \SN[K][+][k](x-t,x+t) G_k^+ \frac{e^{ik(x-t)}}{ik} \Big), 
\Theta^x_z \Big( \hspace{-0.2ex} \sum_{m\in \Z_M} \hspace{-1ex} \rhoND(m) \SN[M][-][m](x-t,x+t) G_m^- e^{im(x+t)}  \Big) \Big] \notag \\ 
-&\, \Big[ \Theta^x_y  \Big( \hspace{-0.2ex} \sum_{k\in \Z_K} \hspace{-1ex} \rhoNDtil(k) \SNtil[K][+][k](x-t,x+t) G_k^+ \frac{e^{ik(x-t)}}{ik} \Big), 
\Theta^x_z \Big( \hspace{-0.2ex} \sum_{m\in \Z_M} \hspace{-1ex} \rhoNDtil(m) \SNtil[M][-][m](x-t,x+t) G_m^- e^{im(x+t)}  \Big) \Big]\bigg)\bigg). \notag  
\end{align}
In order to control this difference, we use the decompositions 
\begin{align*}
\rhoND(k) \SN[K][+][k] &= \rhoNDtil(k) \SNtil[K][+][k] + \big( \rhoND(k) \SN[K][+][k] - \rhoNDtil(k) \SNtil[K][+][k] \big), \\ 
\rhoND(m) \SN[M][-][m] &= \rhoNDtil(m) \SNtil[M][-][m] + \big( \rhoND(m) \SN[M][-][m] - \rhoNDtil(m) \SNtil[M][-][m] \big).
\end{align*}
Due to the bilinearity of the Lie-bracket, we can then decompose \eqref{lipschitz:eq-cartesian-proof-1} into terms which contain the difference of $ \rhoND(k) \SN[K][+][k]$ and $\rhoNDtil(k) \SNtil[K][+][k]$ or the difference of  $\rhoND(m) \SN[M][-][m]$ and $\rhoNDtil(m) \SNtil[M][-][m]$.
 Using a similar argument as in \eqref{modulation:eq-cartesian-p1} and \eqref{modulation:eq-cartesian-p2}, we then obtain
\begin{align}
&\, \Big\|  \VN[K,M][+-] - \VNtil[K,M][+-] \Big\|_{\Ctx{\scrr-1}} \notag \\
\lesssim&\,  \coup \Ac^2 \Bc (KM)^{-2\eta} 
 \Big( \big\| \rhoND(k) \SN[K][+][k] - \rhoNDtil(k)  \SNtil[K][+][k] \big\|_{\Wuv[s][s]}
+ \big\| \rhoND(m) \SN[M][-][m] - \rhoNDtil(m) \SNtil[M][-][m] \big\|_{\Wuv[s][s]}  \Big).
\label{lipschitz:eq-cartesian-proof-2}
\end{align}
We now note that the contributions
$\rhoND(k) -\rhoNDtil(k)$ or $\rhoND(m)-\rhoNDtil(m)$ can only be non-zero in the case $\Nd\neq \Ndtil$ and $\max(K,M)\gtrsim \min(\Nd,\Ndtil)$. Together with Lemma \ref{chaos:lem-wc}, it then follows that 
\begin{align*}
\eqref{lipschitz:eq-cartesian-proof-2} 
&\lesssim \coup \Ac^2 \Bc^2 (KM)^{-2\eta} \mathbf{1} \Big\{  \Nd \neq \Ndtil, \, \max(K,M) \geq \min(\Nd,\Ndtil) \Big\} \\
&+ \coup \Ac^2 \Bc (KM)^{-2\eta} \Big( \big\|  \SN[K][+][k] -   \SNtil[K][+][k] \big\|_{\Wuv[s][s]}
+ \big\|  \SN[M][-][m] - \SNtil[M][-][m] \big\|_{\Wuv[s][s]}  \Big). 
\end{align*}
Due to the definition of $\dpre$, this yields \eqref{lipschitz:eq-modulated-cartesian-e2}.
\end{proof}

\subsection{Modulation operators}

In this subsection, we derive a Lipschitz-variant  of Proposition \ref{modulation:prop-main}, i.e., the Lipschitz-estimates for solutions of the modulation equations.

\begin{proposition}[Lipschitz-estimates for modulation operators]\label{lipschitz:prop-modulation}
Assume that the post-modulation hypothesis (Hypothesis \ref{hypothesis:post}) is satisfied for \eqref{lipschitz:eq-original} and \eqref{lipschitz:eq-tilde}. Then, it holds that 
\begin{equation}\label{lipschitz:eq-modulation}
\begin{aligned}
&\, \sup_{K\in \Dyadiclarge} \Big\| \pSN[K][+][k] - \pSNtil[K][+][k] \Big\|_{\Wuv[s][s]}
+ \sup_{M\in \Dyadiclarge} \Big\| \pSN[M][-][m] - \pSNtil[M][-][m] \Big\|_{\Wuv[s][s]} \\
\lesssim&\, \sup_{K\in \dyadic} \Big\| \SNin[K][+] - \SNtilin[K][+] \Big\|_{\Wuv[s][s]}
+ \sup_{M\in \dyadic} \Big\| \SNin[M][-] - \SNtilin[M][-] \Big\|_{\Wuv[s][s]} \\
+&\, \Bc \Big( \Dc \ddyadic(\Nd,\Ndtil) +  \big\| \UN[][\fs] - \UNtil[][\fs] \big\|_{\Cprod{r-1}{r}}
+ \big\| \VN[][\fs] - \VNtil[][\fs] \big\|_{\Cprod{r}{r-1}} \Big).
\end{aligned}
\end{equation}
Furthermore, we also have Lipschitz-variants of the estimates in \ref{modulation:item-difference}, \ref{modulation:item-orthogonality}, and \ref{modulation:item-structure} from Proposition \ref{modulation:prop-main}. For the $\SN[K][+][k]$-terms, the Lipschitz-variants are given by
\begin{align*}
&\, \sup_{K\in \Dyadiclarge} K^{100} \Big\|  \Big( \SN[K][+][k] - \pSN[K][+][k] \Big) - \Big( \SNtil[K][+][k]- \pSNtil[K][+][k] \Big) \Big\|_{\Wuv[s][s]}  \\
+&\, \sup_{\substack{K\in \Dyadiclarge\colon\\ K>N^{1-\delta}}} \Bc^{-1} K^{100} 
\Big\| \big(\SN[K][+][k]\big)^\ast \SN[K][+][k]  - \big(\SNtil[K][+][k]\big)^\ast \SNtil[K][+][k] \Big\|_{\Wuv[s][s]} 
+ \sup_{K\in \Dyadiclarge} \Big\| \YN[K][+][k] - \YNtil[K][+][k] \Big\|_{\Wuv[s][\scrr]} \\
\lesssim&\, \sup_{K\in \dyadic} \Big\| \SNin[K][+] - \SNtilin[K][+] \Big\|_{\Wuv[s][s]}
+ \sup_{M\in \dyadic} \Big\| \SNin[M][-] - \SNtilin[M][-] \Big\|_{\Wuv[s][s]} \\
+&\, \Bc \Big( \Dc \ddyadic(\Nd,\Ndtil) +  \big\| \UN[][\fs] - \UNtil[][\fs] \big\|_{\Cprod{r-1}{r}}
+ \big\| \VN[][\fs] - \VNtil[][\fs] \big\|_{\Cprod{r}{r-1}} \Big). 
\end{align*}
Similar estimates also hold for the  $\SN[M][-][m]$-terms.
\end{proposition}

\begin{remark}
Due to Hypothesis \ref{hypothesis:post}, the parameter $\Dc =\hcoup \Ac \Bc$ is small. As a result, \eqref{lipschitz:eq-modulation} is equivalent to the estimate
\begin{equation}\label{lipschitz:eq-modulation-rewritten}
\begin{aligned}
&\, \sup_{K\in \Dyadiclarge} \Big\| \pSN[K][+][k] - \pSNtil[K][+][k] \Big\|_{\Wuv[s][s]}
+ \sup_{M\in \Dyadiclarge} \Big\| \pSN[M][-][m] - \pSNtil[M][-][m] \Big\|_{\Wuv[s][s]} \\
\lesssim&\, \sup_{K\in \dyadic} \Big\| \SNin[K][+] - \SNtilin[K][+] \Big\|_{\Wuv[s][s]}
+ \sup_{M\in \dyadic} \Big\| \SNin[M][-] - \SNtilin[M][-] \Big\|_{\Wuv[s][s]} 
+ \Bc \dpre. 
\end{aligned}
\end{equation}
The $\SNin[][\pm]$-terms and $\Bc \dpre$-terms in \eqref{lipschitz:eq-modulation-rewritten} can be traced back to different terms from the proof of Proposition \ref{modulation:prop-para-well-posedness}. While the $\SNin[][\pm]$-terms in \eqref{lipschitz:eq-modulation-rewritten} come from the $\SNin[][\pm]$-terms in 
$\YNcal[][\pm]$, the $\Bc \dpre$-term in \eqref{lipschitz:eq-modulation-rewritten} comes from the nonlinear terms in $\XNcal[][\pm]$ and $\YNcal[][\pm]$.
\end{remark}

The Lipschitz-estimates from Proposition \ref{lipschitz:prop-modulation} can be obtained from minor modifications of the proof of Proposition \ref{modulation:prop-main}. We illustrate the necessary modifications using a simple example, which is a Lipschitz-variant of Lemma \ref{modulation:lem-X-estimate}. 

\begin{lemma}[Lipschitz-variant of Lemma \ref{modulation:lem-X-estimate}]\label{lipschitz:lem-X-estimate}
Let the pre-modulation hypothesis (Hypothesis \ref{hypothesis:pre}) be satisfied for \eqref{lipschitz:eq-original} and \eqref{lipschitz:eq-tilde} and let $K\in \Dyadiclarge$.  Then, it holds that
\begin{equation*}
\begin{aligned}
\Big\| \rho_{\leq N}^2(k) \Big( \Int^v_{u\rightarrow v} \Ad\big( \LON[K][-] \big) \Big) \Para[v][gg] \pSN[K][+][k]
-  \rho_{\leq N}^2(k) \Big( \Int^v_{u\rightarrow v} \Ad\big( \LONtil[K][-] \big) \Big) \Para[v][gg] \pSNtil[K][+][k] \Big\|_{\Wuv[s][s][k]} 
\lesssim \Bc \dpre. 
\end{aligned}
\end{equation*}
\end{lemma}

\begin{proof} 
The estimate follows from a modification of the proof of Lemma \ref{modulation:lem-X-estimate}. As in the proof of \mbox{Lemma \ref{modulation:lem-X-estimate}}, we first replace $\Int^v_{u\rightarrow v}$ by $P^v_{>1} \Int^v_{u\rightarrow v}$. Then, we decompose
\begin{align}
&\, \rho_{\leq N}^2(k) \Big( P^v_{>1}  \Int^v_{u\rightarrow v} \Ad\big( \LON[K][-] \big) \Big) \Para[v][gg] \pSN[K][+][k]
-  \rho_{\leq N}^2(k) \Big( P^v_{>1}  \Int^v_{u\rightarrow v} \Ad\big( \LONtil[K][-] \big) \Big) \Para[v][gg] \pSNtil[K][+][k]  \notag \\ 
=&\, \rho_{\leq N}^2(k) \Big( P^v_{>1}  \Int^v_{u\rightarrow v} \Ad\big( \LON[K][-] - \LONtil[K][-] \big) \Big) \Para[v][gg] \pSN[K][+][k] 
\label{lipschitz:eq-X-estimate-p1} \\ 
+&\, \rho_{\leq N}^2(k) \Big( P^v_{>1}  \Int^v_{u\rightarrow v} \Ad\big( \LONtil[K][-] \big) \Big) \Para[v][gg] \big( \pSN[K][+][k] - \pSNtil[K][+][k] \big). 
\label{lipschitz:eq-X-estimate-p2} 
\end{align}
Using the same estimates as in the proof of Lemma \ref{modulation:lem-X-estimate}, it then follows that 
\begin{equation}\label{lipschitz:eq-X-estimate-p3}
\begin{aligned}
&\, \big\| \eqref{lipschitz:eq-X-estimate-p1} \|_{\Wuv[s][s][k]}
+ \big\| \eqref{lipschitz:eq-X-estimate-p2} \|_{\Wuv[s][s][k]} \\ 
\lesssim&\,  \big\| \LONtil[K][-] - \LON[K][-] \big\|_{\Cprod{s}{s-1}} \big\| \pSN[K][+][k] \big\|_{\Wuv[s][s][k]}
+ \big\| \LONtil[K][-]  \big\|_{\Cprod{s}{s-1}} \big\| \pSN[K][+][k] - \pSNtil[K][+][k] \big\|_{\Wuv[s][s][k]}.
\end{aligned}
\end{equation}
Using Hypothesis \ref{hypothesis:pre}, Corollary \ref{modulation:cor-LON}, and \eqref{lipschitz:eq-modulated-lo} from Lemma \ref{lipschitz:lem-modulated-mixed}, it follows that 
\begin{equation*}
\eqref{lipschitz:eq-X-estimate-p3} \lesssim \dpre \Bc + \Dc \big\| \pSN[K][+][k] - \pSNtil[K][+][k] \big\|_{\Wuv[s][s][k]} 
\lesssim \Bc \dpre,
\end{equation*}
which yields \eqref{lipschitz:lem-X-estimate}.
\end{proof}

The necessary modifications in the statements and proofs of Proposition \ref{modulation:prop-resonant}, Lemma \ref{modulation:lem-Y-integral-commutator}, Lemma \ref{modulation:lem-Y-double-duhamel}, and Lemma \ref{modulation:lem-Y-remaining} are similar, and we omit the details. With these modifications, Proposition \ref{lipschitz:prop-modulation} then follows from similar arguments as in the proof of Proposition \ref{modulation:prop-main}, and we again omit the details. For a more systematic description of the necessary modifications, we also refer to Algorithm \ref{lipschitz:algorithm} below.

As an application of Proposition \ref{lipschitz:prop-modulation}, we obtain the following estimate for the distances from \mbox{Definition \ref{lipschitz:def-distances}}.

\begin{lemma}[\protect{Relationship of $\dpre$ and $\dpost$}]\label{lipschitz:lem-distance}
Let the post-modulation hypothesis \mbox{(Hypothesis \ref{hypothesis:post})} be satisfied for both \eqref{lipschitz:eq-original} and \eqref{lipschitz:eq-tilde}. Then, it holds that 
\begin{equation}\label{lipschitz:eq-distance}
\dpre \lesssim \dpost. 
\end{equation}
\end{lemma}

\begin{proof}
The estimate \eqref{lipschitz:eq-distance} can easily be obtained by applying Proposition \ref{lipschitz:prop-modulation} to the $\pSN[][\pm]$-terms in \eqref{lipschitz:eq-dpre} and using $\Dc = \hcoup \Ac \Bc \leq 1$.
\end{proof}

At the end of this subsection, we also record a Lipschitz-variant of Lemma \ref{modulation:lem-initial-modulated-linear}.

\begin{lemma}[Lipschitz-estimate for initial data of modulated linear wave]\label{lipschitz:lem-modulated-initial}
Assume that the post-modulation hypothesis (Hypothesis \ref{hypothesis:post}) is satisfied for \eqref{lipschitz:eq-original} and \eqref{lipschitz:eq-tilde}. Then, it holds that 
\begin{equation}\label{lipschitz:eq-modulated-initial}
\begin{aligned}
&\, \Big\|  \Big( \sum_{K\in \dyadic} \SNin[K][+](u) \big( P_{\leq \Nd} P^\sharp_{R;K} W^{(\Rscript,\coup),+}\big)(u)  - \UN[][+] (u,u) \Big)  \\
&-\hspace{1ex}  \Big( \sum_{K\in \dyadic} \SNtilin[K][+](u) \big( P_{\leq \Ndtil} P^\sharp_{R;K} W^{(\Rscript,\coup),+}\big)(u)  - \UNtil[][+] (u,u) \Big) 
 \Big\|_{\C_u^{r-1}} \\  
\lesssim&\,  
\hcoup \Ac \sup_{K\in \dyadic} \big\| \SNin[K][+] - \SNtilin[K][+] \big\|_{\C_u^s}
+ \hcoup \Ac \sup_{M\in \dyadic} \big\| \SNin[M][-] - \SNtilin[M][-] \big\|_{\C_v^s} 
+  \Dc\,  \mathsf{d}_{\dyadic}\big( \Nd, \Ndtil \big) \\ 
+&\,  \, \Dc \Big( \big\| \UN[][\fs] - \UNtil[][\fs] \big\|_{\Cprod{r-1}{r}}
+ \big\| \VN[][\fs] - \VNtil[][\fs] \big\|_{\Cprod{r}{r-1}}\Big). 
\end{aligned}
\end{equation}
\end{lemma}

\begin{remark}
We emphasize that the $\SNin[][\pm]$ and $\Nd$-terms in \eqref{lipschitz:eq-modulated-initial} are as in the definition of $\dpost$ from \eqref{lipschitz:eq-dpost}, but the $\UN[][\fs]$ and $\VN[][\fs]$-terms in \eqref{lipschitz:eq-modulated-initial} have an additional $\Dc$-factor. This directly follows from the $\hcoup \Ac$-factors from our bounds on $W^{(\Rscript,\coup),\pm}$ (see e.g. Hypothesis \ref{hypothesis:probabilistic}) and the $\Bc$-factor in \eqref{lipschitz:eq-modulation}.
\end{remark}

Since \eqref{lipschitz:eq-modulated-initial} follows from standard modifications of the proof of Lemma \ref{modulation:lem-initial-modulated-linear}, we omit the proof.

\subsection{Tensor products of modulated linear waves}
In this subsection, we state Lipschitz-variants of our estimates of tensor products of modulated linear waves. 
Due to the Killing-renormalization, the necessary modifications are slightly more involved than in Subsection \ref{section:lipschitz-modulated}.
In the following lemma, we record a Lipschitz-variant of Lemma \ref{killing:lem-tensor-modulated-linear}.

\begin{lemma}[Lipschitz-estimate for the quadratic tensor-product of modulated linear waves]\label{lipschitz:lem-quadratic-tensor}
Assume that the post-modulation hypothesis (Hypothesis \ref{hypothesis:post}) is satisfied for \eqref{lipschitz:eq-original} and \eqref{lipschitz:eq-tilde} and let $K,L\in \Dyadiclarge$. 
Furthermore, let $h\in (0,\infty)$ and let $y,z\in \R$ satisfy $|y|,|z|\leq h$. Then, we have the following Lipschitz-estimates: 
\begin{equation}\label{lipschitz:eq-quadratic-tensor}
\begin{aligned}
&\, \Big\| \biglcol \hspace{0.5ex} \Theta^x_y \IUN[L][+] \otimes \Theta^x_z \UN[K][+] \bigrcol 
- \biglcol \hspace{0.5ex} \Theta^x_y \IUNtil[L][+] \otimes \Theta^x_z \UNtil[K][+] \bigrcol \Big\|_{\Cprod{\gamma}{s}} \\ 
\lesssim&\, \Big( \max(K,L)^{r-\frac{1}{2}+\eta} L^{-\frac{1}{2}}
+ \mathbf{1} \big\{ K=L \big\} N^{-\delta+\vartheta+\varsigma} N h \Big) \Dc \dpost
\end{aligned}
\end{equation}
and 
\begin{equation}\label{lipschitz:eq-quadratic-tensor-truncated}
\begin{aligned}
&\, \Big\| \biglcol \hspace{0.5ex} \Theta^x_y P_{\leq N}^x \IUN[L][+] \otimes \Theta^x_z P_{\leq N}^x \UN[K][+] \bigrcol 
- \biglcol \hspace{0.5ex} \Theta^x_y P_{\leq N}^x \IUNtil[L][+] \otimes \Theta^x_z P_{\leq N}^x \UNtil[K][+] \bigrcol \Big\|_{\Cprod{\gamma}{s}} \\ 
\lesssim&\, \Big( \max(K,L)^{r-\frac{1}{2}+\eta} L^{-\frac{1}{2}}
+ \mathbf{1} \big\{ K=L \big\} N^{-\delta+\vartheta+\varsigma} \langle N h \rangle \Big) \Dc \dpost.
\end{aligned}
\end{equation}
\end{lemma}

\begin{remark}
We remark that Lemma \ref{lipschitz:lem-quadratic-tensor} contains a $N^{-\delta+\vartheta+\varsigma}$-factor, whereas Lemma \ref{killing:lem-tensor-modulated-linear} only contains a $N^{-\delta+\vartheta}$-factor. Since $\varsigma$ is much smaller than $\vartheta$, this does not affect any subsequent arguments.
\end{remark}

\begin{remark}
Just like Lemma \ref{killing:lem-tensor-modulated-linear} cannot be obtained from Lemma \ref{modulation:lem-linear}, the Lipschitz-estimate of the Wick-ordered tensor product $\biglcol \hspace{0.5ex} \Theta^x_y \IUN[L][+] \otimes \Theta^x_z \UN[K][+] \bigrcol$ from Lemma \ref{lipschitz:lem-quadratic-tensor} cannot be obtained from the Lipschitz-estimates of $\IUN[L][+]$ and $\UN[K][+]$ from Lemma \ref{lipschitz:lem-modulated-mixed}. 
\end{remark}

\begin{proof}
Using the same identity as in \eqref{killing:eq-tensor-truncated-p6} from the proof of Lemma \ref{killing:lem-tensor-modulated-linear}, the second estimate \eqref{lipschitz:eq-quadratic-tensor-truncated} can be derived from the first estimate \eqref{lipschitz:eq-quadratic-tensor}. We therefore focus on the proof of \eqref{lipschitz:eq-quadratic-tensor}, which is a variant of the proof of Lemma \ref{killing:lem-tensor-modulated-linear}.

We let $\IUN[u_0,L][+]$ and $\UN[u_1,K][+]$ be as in \eqref{killing:eq-tensor-truncated-p3} and let $\IUNtil[u_0,L][+]$ and $\UNtil[u_1,K][+]$ be defined similarly as in \eqref{killing:eq-tensor-truncated-p3}, but with $\rhoND$ and $\SN[][+]$ replaced by $\rhoNDtil$ and $\SNtil[][+]$. Then, we decompose
\begin{align}
&\biglcol \hspace{0.5ex} \Theta^x_y \IUN[L][+] \otimes \Theta^x_z \UN[K][+] \bigrcol 
- \biglcol \hspace{0.5ex} \Theta^x_y \IUNtil[L][+] \otimes \Theta^x_z \UNtil[K][+] \bigrcol \notag\\
=& \, \sum_{u_0,u_1\in \LambdaRR} \big(  \Theta^x_{y} \psiRuL \big) \big( \Theta^x_{z} \psiRuKone \big)  
\bigg(\Theta^x_y \IUN[u_0,L][+] \otimes \Theta^x_z \UN[u_1,K][+] - \Theta^x_y \IUNtil[u_0,L][+] \otimes \Theta^x_z \UNtil[u_1,K][+] 
\label{lipschitz:eq-quadratic-p1} \\
&\hspace{33ex}- \coup \mathbf{1} \big\{ K=L \big\} \mathbf{1} \big\{ u_0= u_1 \big\} 
\Big( \Cf^{(\Nd)}_{K}(y-z) - \Cf^{(\Ndtil)}_{K}(y-z) \Big) \Cas \bigg) \notag \\
+&\coup \mathbf{1} \big\{ K=L \big\} \sum_{u_0\in \LambdaRR} 
\Big( \big(  \Theta^x_{y} \psiRuL \big) \big( \Theta^x_{z} \psiRuK \big) 
- \psiRuL \psiRuK \Big)
\Big( \Cf^{(\Nd)}_{K}(y-z) - \Cf^{(\Ndtil)}_{K}(y-z) \Big) \Cas.
\label{lipschitz:eq-quadratic-p2}
\end{align}
Using the smoothness and decay of $\psiRuL$ and $\psiRuK$ and Lemma \ref{ansatz:lem-renormalization}, we easily obtain
\begin{align*}
\big| \eqref{lipschitz:eq-quadratic-p2}\big| 
&\lesssim \coup \mathbf{1} \big\{ K=L \big\} \big( |y|+ |z| \big) \Big( \big| \Cf^{(\Nd)}_{K}(y-z)\big| + \big| \Cf^{(\Ndtil)}_{K}(y-z) \big| \Big) \mathbf{1} \big\{ \Nd \neq \Ndtil \big\} \\
&\lesssim \coup h \mathbf{1} \big\{ K=L \big\}   \mathbf{1} \big\{ \Nd \neq \Ndtil \big\}. 
\end{align*}
Since $\Nd,\Ndtil \lesssim N$, this is acceptable. 
Thus, it remains to control \eqref{lipschitz:eq-quadratic-p1}, which is more difficult. Due to Lemma \ref{prelim:lem-psi-sum}, it suffices to treat the contribution from \eqref{lipschitz:eq-quadratic-p1} for fixed $u_0,u_1 \in \LambdaR$. In order to simplify the notation, we introduce the bilinear maps
\begin{align*}
\mathcal{T} \big( S^+, \widetilde{S}^+ \big) 
&:= \coup \Theta^x_y \Big( \sum_{\ell \in \Z_L} S_\ell^+ G_{u_0,\ell}^+ \frac{e^{i\ell u}}{i \ell} \Big) 
\otimes \Theta^x_z \Big( \sum_{k\in \Z_K} \widetilde{S}_k^+ G_{u_1,k}^+ e^{iku} \Big) \\ 
&\hspace{3ex}- \coup \mathbf{1} \big\{ K=L \big\} \mathbf{1} \big\{ u_0 =u_1 \big\} 
\sum_{k\in \Z_K} \delta^{ab} \big( \Theta^x_y S_k^+ E_a \otimes \Theta^x_z \widetilde{S}_k^+ E_b\big) \frac{e^{ik(y-z)}}{(-ik)}, \\ 
\mathcal{R} (S^+, \widetilde{S}^+ ) &:= \coup 
\sum_{k\in \Z_K} \delta^{ab} \big( \Theta^x_y S_k^+ E_a \otimes \Theta^x_z \widetilde{S}_k^+ E_b\big) \frac{e^{ik(y-z)}}{(-ik)}.
\end{align*}
Using $\mathcal{T}$ and $\mathcal{R}$, the summands in \eqref{lipschitz:eq-quadratic-p1} can be written as 
\begin{align}
&\mathcal{T} \big( \rhoND \SN[L][+], \rhoND \SN[K][+] \big) - 
\mathcal{T} \big( \rhoNDtil \SNtil[L][+], \rhoNDtil \SNtil[K][+] \big) \label{lipschitz:eq-quadratic-p3} \\ 
+& \mathbf{1} \big\{ K=L \big\} \mathbf{1} \big\{ u_0 =u_1 \big\}  \Big( \mathcal{R} \big( \rhoND \SN[L][+], \rhoND \SN[K][+] \big) - 
\mathcal{R} \big( \rhoNDtil \SNtil[L][+], \rhoNDtil \SNtil[K][+] \big) \label{lipschitz:eq-quadratic-p4} \\ 
&\hspace{24ex}- \coup \big( \Cf^{(\Nd)}_{K}(y-z) - \Cf^{(\Ndtil)}_{K}(y-z) \big) \Cas \Big). \notag
\end{align}
We first control \eqref{lipschitz:eq-quadratic-p3}. 
Since $\mathcal{T}$ is bilinear, it follows from Hypothesis \ref{hypothesis:probabilistic}.\ref{ansatz:item-hypothesis-tensor-product} that 
\begin{equation*}
\begin{aligned}
&\Big\| \mathcal{T} \big( \rhoND \SN[L][+], \rhoND \SN[K][+] \big) - 
\mathcal{T} \big( \rhoNDtil \SNtil[L][+], \rhoNDtil \SNtil[K][+] \big) \Big\|_{\Cprod{\gamma}{s}} \\ 
\lesssim&\,  \max(K,L)^{\gamma+\frac{1}{2}+\frac{\eta}{2}} L^{-\frac{1}{2}}  \coup \Ac^2 \Bc 
\big\| \rhoND(k) \SN[K][+][k] - \rhoNDtil(k)  \SNtil[K][+][k] \big\|_{\Wuv[s][s]}. 
\end{aligned}
\end{equation*}
By arguing similarly as for \eqref{lipschitz:eq-cartesian-proof-2} from the proof of Lemma \ref{lipschitz:lem-modulated-cartesian}, this yields an acceptable contribution. In order to control \eqref{lipschitz:eq-quadratic-p4}, we further decompose
\begin{align}
 &\, \mathcal{R} \big( \rhoND \SN[L][+], \rhoND \SN[K][+] \big) - \mathcal{R} \big( \rhoNDtil \SNtil[L][+], \rhoNDtil \SNtil[K][+] \big)
 - \coup \big( \Cf^{(\Nd)}_{K}(y-z) - \Cf^{(\Ndtil)}_{K}(y-z) \big) \Cas \notag \\ 
=& \, 
\coup \delta^{ab} \sum_{k\in \Z_K} \big( \rhoNDsquare(k) - \rhoNDsquaretil(k) \big) 
\Big( \big( \Theta^x_y \SN[K][+][k] E_a \big) \otimes \big( \Theta^x_z \SN[K][+][k] E_b \big) - E_a \otimes E_b \Big) \frac{e^{ik(y-z)}}{(-ik)}
\label{lipschitz:eq-quadratic-p5} \\ 
+&\,\coup \delta^{ab} \sum_{k\in \Z_K}  \rhoNDsquaretil(k) 
\Big( \big( \Theta^x_y \SN[K][+][k] E_a \big) \otimes \big( \Theta^x_z \SN[K][+][k] E_b \big)
- \big( \Theta^x_y \SNtil[K][+][k] E_a \big) \otimes \big( \Theta^x_z \SNtil[K][+][k] E_b \big) \Big) \frac{e^{ik(y-z)}}{(-ik)}
\label{lipschitz:eq-quadratic-p6} \\ 
-&\, \coup \bigg( \Cf^{(\Nd)}_{K}(y-z) - \Cf^{(\Ndtil)}_{K}(y-z) - 
\sum_{k\in \Z_K} \big( \rhoNDsquare(k) - \rhoNDsquaretil(k) \big) \frac{e^{ik(y-z)}}{(-ik)} \bigg) \Cas. 
\label{lipschitz:eq-quadratic-p7}
\end{align}
By arguing similarly as in our estimate of \eqref{killing:eq-resonant-p0}, it follows that
\begin{equation*}
\big\| \eqref{lipschitz:eq-quadratic-p5} \big\|_{\Cprod{s}{s}}
\lesssim \coup \Big\| \big( \Theta^x_y \SN[K][+][k] \big) \big( \Theta^x_z \SN[K][+][k] \big)^\ast - \Id_\frkg \Big\|_{\Wuv[s][s]}
\mathbf{1} \big\{ \Nd \neq \Ndtil,\, K \geq \min(\Nd,\Ndtil) \big\}. 
\end{equation*}
Using the estimates from \eqref{killing:eq-resonant-p1} and \eqref{killing:eq-resonant-p2} and Hypothesis \ref{hypothesis:post}.\ref{ansatz:item-post-orthogonality}, this easily yields an acceptable contribution. By arguing similarly as in our estimate of \eqref{killing:eq-resonant-p0}, it also follows that
\begin{equation*}
\big\| \eqref{lipschitz:eq-quadratic-p5} \big\|_{\Cprod{s}{s}}
\lesssim \coup \Big\| \big( \Theta^x_y \SN[K][+][k] \big) \big( \Theta^x_z \SN[K][+][k] \big)^\ast
- \big( \Theta^x_y \SNtil[K][+][k] \big) \big( \Theta^x_z \SNtil[K][+][k] \big)^\ast \Big\|_{\Wuv[s][s]}. 
\end{equation*}
Using the estimates from \eqref{killing:eq-resonant-p1} and \eqref{killing:eq-resonant-p2} and the Lipschitz-estimate from Proposition \ref{lipschitz:prop-modulation}, this also yields an acceptable contribution. Finally, using Lemma \ref{ansatz:lem-renormalization}, $\Nd,\Ndtil\lesssim N$, and $\coup \lesssim \Dc^2$, we obtain
\begin{equation*}
\big| \eqref{lipschitz:eq-quadratic-p7} \big| \lesssim \coup h \mathbf{1} \big\{ \Nd \neq \Ndtil \big\}
\lesssim \big( \Dc N^{1-\delta+\vartheta+\varsigma} h \big) \times \big( \Dc \max(\Nd,\Ndtil)^{-1+\delta-\vartheta+\varsigma} 
\mathbf{1}\big\{ \Nd \neq \Ndtil \big\} \big),
\end{equation*}
which is acceptable.
\end{proof}

In the following lemma, we record Lipschitz-variants of the estimates from Proposition \ref{killing:prop-quadratic}, Corollary \ref{killing:cor-tensors}, and Proposition \ref{killing:prop-quartic}.

\begin{lemma}\label{lipschitz:lem-higher-tensors}
Let the post-modulation hypothesis (Hypothesis  \ref{hypothesis:post}) be satisfied for \eqref{lipschitz:eq-original} and \eqref{lipschitz:eq-tilde}. Then, we have the following Lipschitz-estimates: 
\begin{enumerate}[label=(\roman*)]
\item Let $K,L\in \Dyadiclarge$, let $E\in \frkg$ satisfy $\| E\|_{\frkg}\leq 1$, and  let $y\in \R$. Then, 
\begin{equation}\label{lipschitz:eq-killing}
\begin{aligned}
&\Big\| \Big( \Big[ P_{\leq N}^x \UN[K][+], \Big[  \Theta^x_y P_{\leq N}^x\IUN[L][+] , E \Big] \Big] -  \coup \mathbf{1}\big\{ K = L \big\} \Cf^{(\Ncs)}_{K}(y) \Kil(E)  \Big) \\
&\hspace{2ex} - 
\Big( \Big[ P_{\leq N}^x \UNtil[K][+], \Big[  \Theta^x_y P_{\leq N}^x\IUNtil[L][+] , E \Big] \Big] -  \coup \mathbf{1}\big\{ K = L \big\} \Cftil^{(\Ncs)}_{K}(y) \Kil(E)  \Big) \Big\|_{\Cprod{\gamma}{s}} \\
\lesssim& \,  \Big( \max(K,L)^{\gamma+\frac{1}{2}+\eta} L^{-\frac{1}{2}} + \mathbf{1}\big\{ K=L \big\} N^{-\delta+\vartheta+\varsigma} \langle N y \rangle \Big) \Dc \dpost. 
\end{aligned}
\end{equation}
\item Let $K,L,M\in \Dyadiclarge$, let $h\in (0,\infty)$, and let $y_1,y_2,y_3\in \R$ satisfy $|y_1|,|y_2|,|y_3|\leq h$. Then, it holds that 
\begin{align*}
& \Big\| \Big( \biglcol \, \Theta^x_{y_1} P_{\leq N}^x \UN[K][+]  \otimes \Theta^x_{y_2} P_{\leq N}^x \IUN[L][+]   \bigrcol \Big)  \otimes  \Theta^x_{y_3} \VN[M][-]  \notag \\
&\hspace{2ex} - \Big( \biglcol \, \Theta^x_{y_1} P_{\leq N}^x \UNtil[K][+]  \otimes \Theta^x_{y_2} P_{\leq N}^x \IUNtil[L][+]   \bigrcol \Big)  \otimes  \Theta^x_{y_3} \VNtil[M][-]  \Big\|_{\Cprod{r-1}{r-1}} \\ 
\lesssim&\,  \Big( \max\big( K, L \big)^{r-s} L^{-\frac{1}{2}} + \mathbf{1} \big\{ K = L \big\} N^{-\delta+\vartheta+\varsigma} \langle N h \rangle \Big) M^{r-s} \Dc^2 \dpost. 
\end{align*}
\item
Let $K_u,K_v,M_u,M_v \in \Dyadiclarge$ satisfy $K_u \simeq_\delta K_v$ and $M_u\simeq_\delta M_v$.
Furthermore, let $h\in (0,\infty)$ and let $y_1,y_2,y_3,y_4\in \R$ satisfy $|y_1|,|y_2|,|y_3|,|y_4|\leq h$.
Then, it holds that 
\begin{align*}
& \bigg\| \bigg( \Theta^x_{y_1} P_{\leq N}^x \UN[K_u][+]  
\otimes \Theta^x_{y_2} P_{\leq N}^x  \IVN[K_v][-]
\otimes \Theta^x_{y_3} P_{\leq N}^x  \IUN[M_u][+] 
\otimes \Theta^x_{y_4} P_{\leq N}^x \VN[M_v][-] \notag \\
&\hspace{6ex}- \coup^2 \mathbf{1} \big\{ K_u = M_u \big\}  \mathbf{1} \big\{ K_v = M_v \big\} 
\Cf^{(\Ncs)}_{K_u}(y_3-y_1) \Cf^{(\Ncs)}_{K_v}(y_2-y_4) \,   \shuffle{(2,3)} \big( 
\Cas \otimes \Cas \big) \bigg) \\
&-\bigg( \Theta^x_{y_1} P_{\leq N}^x \UNtil[K_u][+]  
\otimes \Theta^x_{y_2} P_{\leq N}^x  \IVNtil[K_v][-]
\otimes \Theta^x_{y_3} P_{\leq N}^x  \IUNtil[M_u][+] 
\otimes \Theta^x_{y_4} P_{\leq N}^x \VNtil[M_v][-] \notag \\
&\hspace{6ex}- \coup^2 \mathbf{1} \big\{ K_u = M_u \big\}  \mathbf{1} \big\{ K_v = M_v \big\} 
\Cftil^{(\Ncs)}_{K_u}(y_3-y_1) \Cftil^{(\Ncs)}_{K_v}(y_2-y_4) \,   \shuffle{(2,3)} \big( 
\Cas \otimes \Cas \big) \bigg)
 \bigg\|_{\Cprod{r-1}{r-1}} \notag\\
&\lesssim \,  \max\big( K_u,K_v,M_u,M_v)^{-(1-4\delta)\delta} \langle N h \rangle^2 \Dc^3 \dpost.
\end{align*}
\item Let $K_u,K_v,M_u,M_v \in \Dyadiclarge$ satisfy $K_u \simeq_\delta K_v$ and $M_u\simeq_\delta M_v$.
Then, it holds that 
\begin{align*}
 \Big\| \Big[  \UN[K_u,K_v][+-],  \VN[M_u,M_v][+-] \Big]_{\leq N}
-  \Big[  \UNtil[K_u,K_v][+-],  \VNtil[M_u,M_v][+-] \Big]_{\leq N} \Big\|_{\Cprod{r-1}{r-1}} 
\lesssim  (K_u K_v M_u M_v)^{-\eta} \Dc^3 \dpost. 
\end{align*}
\end{enumerate}
\end{lemma}

We recall that the main ingredient in the proofs of 
Proposition \ref{killing:prop-quadratic}, Corollary \ref{killing:cor-tensors}, and Proposition \ref{killing:prop-quartic}
is the estimate of $ \biglcol \hspace{0.5ex} \Theta^x_y P_{\leq N}^x \IUN[L][+] \otimes \Theta^x_z P_{\leq N}^x \UN[K][+] \bigrcol $
from Lemma \ref{killing:lem-tensor-modulated-linear}. Since the Lipschitz-variant of the latter has been established in Lemma \ref{lipschitz:lem-quadratic-tensor}, the estimates in Lemma \ref{lipschitz:lem-higher-tensors} can now be obtained from minor modifications of the proofs of 
Proposition \ref{killing:prop-quadratic}, Corollary \ref{killing:cor-tensors}, and Proposition \ref{killing:prop-quartic}. Due to this, we omit the details. 

\subsection{Lipschitz-estimates for remainder equations}

In this subsection, we record Lipschitz-variants of the terms in the remainder equations, i.e., the terms in Proposition \ref{ansatz:prop-decomposition}. In the following lemma, we record Lipschitz-variants of Proposition \ref{hhl:prop-main}, \eqref{jacobi:eq-main-estimate-1} and \eqref{jacobi:eq-main-estimate-2} from Proposition \ref{jacobi:prop-main}, Proposition \ref{null:prop-perturbative}, \mbox{Proposition \ref{structural:prop-main}}, and Proposition \ref{ren:prop-main}.

\begin{proposition}\label{lipschitz:prop-remainder-no-jcb}
Let the post-modulation hypothesis (Hypothesis \ref{hypothesis:post}) be satisfied for \eqref{lipschitz:eq-original} and \eqref{lipschitz:eq-tilde}. Then, we have the following Lipschitz-estimates: 
\begin{align}
\Big\| \HHLNErr - \HHLNErrtil \Big\|_{\Cprod{r-1}{r-1}} &\lesssim \Dc \dpost, \label{lipschitz:eq-remainder-hhl} \\
\Big\| \JcbNErr - \JcbNErrtil \Big\|_{\Cprod{r-1}{r-1}} &\lesssim \Big( 1+ \tfrac{\max(\Nd,\Ndtil)}{N} N^{\delta_2} \Big) \Dc^2 \dpost, 
\label{lipschitz:eq-remainder-jcb} \\
\Big\| \JcbNErr[][\dagger] - \JcbNErrtildagger \Big\|_{\Cprod{r-1}{r-1}}
&\lesssim \Dc^2 \dpost, \label{lipschitz:eq-remainder-jcb-dagger}\\ 
\Big\| \PIN - \PINtil \Big\|_{\Cprod{r-1}{r-1}} &\lesssim \Dc \dpost, 
\label{lipschitz:eq-remainder-pin}\\ 
\Big\| \SEN[][u][]  - \SENtilu \Big\|_{\WCprod{r-1}{r-1}} &\lesssim \Dc \dpost, 
\label{lipschitz:eq-remainder-sen-u}\\ 
\Big\| \SEN[][v][]  - \SENtilv \Big\|_{\WCprod{r-1}{r-1}} &\lesssim \Dc \dpost,
\label{lipschitz:eq-remainder-sen-v} \\ 
\Big\| \RenNErr - \RenNErrtil \Big\|_{\Cprod{r-1}{r-1}} &\lesssim N^{-\eta} \Dc^2 \dpost. 
\label{lipschitz:eq-remainder-ren}
\end{align}
\end{proposition}

Proposition \ref{lipschitz:prop-remainder-no-jcb} follows from standard modifications of our arguments from Sections \ref{section:hhl}-\ref{section:renormalization} and, as a result, we do not provide all details. Instead, we describe the necessary modifications in two illustrative examples. In addition, we describe a general algorithm which can be used to perform the necessary modifications for all other estimates. 

The first example is a Lipschitz-variant of Lemma \ref{null:lem-p-pm}, which was used to control $\PIN$. 

\begin{lemma}[Lipschitz-variant of Lemma \ref{null:lem-p-pm}]\label{lipschitz:lem-example-1}
Let the post-modulation hypothesis (Hypothesis \ref{hypothesis:post}) be satisfied for \eqref{lipschitz:eq-original} and \eqref{lipschitz:eq-tilde} and let 
$K,M_u,M_v\in \Dyadiclarge$ satisfy $\max(M_u,M_v)\geq K^{1-\delta}$ and $M_u \simeq_\delta M_v$. Furthermore, 
let $y,z\in \R$. Then, it holds that
\begin{equation*}
\Big\| \Theta^x_y \UN[K][+] \Para[u][nsim] \Theta^x_z \VN[M_u,M_v][+-] 
- \Theta^x_y \UNtil[K][+] \Para[u][nsim] \Theta^x_z \VNtil[M_u,M_v][+-] \Big\|_{\Cprod{r-1}{r-1}}
\lesssim (KM_u M_v)^{-\eta} \Dc^2 \dpost. 
\end{equation*}
\end{lemma}

\begin{proof}
The argument is a minor modification of the proof of Lemma \ref{null:lem-p-pm}. We first decompose
\begin{align}
&\Theta^x_y \UN[K][+] \Para[u][nsim] \Theta^x_z \VN[M_u,M_v][+-] 
- \Theta^x_y \UNtil[K][+] \Para[u][nsim] \Theta^x_z \VNtil[M_u,M_v][+-] \notag \\ 
=&\, \Theta^x_y \Big( \UN[K][+] - \UNtil[K][+] \Big)  \Para[u][nsim] \Theta^x_z \VN[M_u,M_v][+-] \label{lipschitz:eq-example-1-a} \\
+&\, \Theta^x_y \UNtil[K][+] \Para[u][nsim]\Theta^x_z \Big( \VN[M_u,M_v][+-] - \VNtil[M_u,M_v][+-] \Big).  \label{lipschitz:eq-example-1-b}
\end{align}
Using our para-product estimate (Lemma \ref{prelim:lem-paraproduct}), Lemma \ref{modulation:lem-bilinear}, and the Lipschitz-estimate from \eqref{lipschitz:eq-modulated-p}, we obtain that
\begin{equation*}
\big\| \, \eqref{lipschitz:eq-example-1-a} \big\|_{\Cprod{r-1}{r-1}}
\lesssim \big\| \UN[K][+] - \UNtil[K][+] \big\|_{\Cprod{r-1}{s}} 
\big\| \VN[M_u,M_v][+-] \big\|_{\Cprod{\eta}{r-1}} 
\lesssim K^{r-s} \dpost \times M_u^{\eta-s} M_v^{r-s} \Dc^2.
\end{equation*}
Since $\max(M_u,M_v)\geq K^{1-\delta}$ and $M_u \simeq_\delta M_v$, this is acceptable. Similarly, using our 
para-product estimate (Lemma \ref{prelim:lem-paraproduct}), Lemma \ref{modulation:lem-linear}, and the Lipschitz-estimate from \eqref{lipschitz:eq-modulated-pm-3}, we obtain that 
\begin{equation*}
\big\| \, \eqref{lipschitz:eq-example-1-b} \big\|_{\Cprod{r-1}{r-1}}
\lesssim \big\| \UN[K][+]  \big\|_{\Cprod{r-1}{s}} 
\big\| \VN[M_u,M_v][+-] -  \VNtil[M_u,M_v][+-] \big\|_{\Cprod{\eta}{r-1}} 
\lesssim K^{r-s} \Dc \times M_u^{\eta-s} M_v^{r-s} \Dc \dpost,
\end{equation*}
which is also acceptable.
\end{proof}

The next example is a Lipschitz-variant of \eqref{null:eq-p-p-qq}, which was a crucial step in the proof of Lemma \ref{null:lem-p-p}.

\begin{lemma}[\protect{Lipschitz-variant of \eqref{null:eq-p-p-qq}}]\label{lipschitz:lem-example-2}
Let the post-modulation hypothesis (Hypothesis \ref{hypothesis:post}) be satisfied for \eqref{lipschitz:eq-original} and \eqref{lipschitz:eq-tilde}. Furthermore, let $K,L\in \Dyadiclarge$ satisfy $L\geq K^{1-\delta}$. Then, it holds that 
\begin{align}
&\,\bigg\|  \bigg(  \Big[ \UN[K][+],  \Big[ \IUN[L][+], P^v_{\geq L^{1-\deltap}} \VN[<L^{1-\delta}][\fcs] \Big]_{\leq N} \Big]_{\leq N} -  \coup \mathbf{1} \big\{ K=L \big\} \Renorm[\Ncs][K] P^v_{\geq K^{1-\deltap}} \VN[<K^{1-\delta}][\fcs] \bigg) \label{lipschitz:eq-example-2} \\
&\, \hspace{3ex}  - \bigg(  \Big[ \UNtil[K][+],  \Big[ \IUNtil[L][+], P^v_{\geq L^{1-\deltap}} \VNtil[<L^{1-\delta}][\fcs] \Big]_{\leq N} \Big]_{\leq N} -  \coup \mathbf{1} \big\{ K=L \big\} \Renormtil[\Ncs][K] P^v_{\geq K^{1-\deltap}} \VNtil[<K^{1-\delta}][\fcs] \bigg) 
\bigg\|_{\Cprod{r-1}{r-1}} \notag  \\ 
\lesssim&\, (KL)^{-\eta} \Dc^2 \dpost. \notag
\end{align}
\end{lemma}

\begin{proof}
We first decompose the term on the left-hand side of \eqref{lipschitz:eq-example-2} into 
\begin{equation}\label{lipschitz:eq-example-2-a}
\begin{aligned}
 \, &\Big[ \UN[K][+],  \Big[ \IUN[L][+], P^v_{\geq L^{1-\deltap}} \VN[<L^{1-\delta}][\fcs] \Big]_{\leq N} \Big]_{\leq N} -  \coup \mathbf{1} \big\{ K=L \big\} \Renorm[\Ncs][K] P^v_{\geq K^{1-\deltap}} \VN[<K^{1-\delta}][\fcs]  \\ 
 -\, \bigg(  &\Big[ \UNtil[K][+],  \Big[ \IUNtil[L][+], P^v_{\geq L^{1-\deltap}} \VN[<L^{1-\delta}][\fcs] \Big]_{\leq N} \Big]_{\leq N} -  \coup \mathbf{1} \big\{ K=L \big\} \Renormtil[\Ncs][K] P^v_{\geq K^{1-\deltap}} \VN[<K^{1-\delta}][\fcs] \bigg) 
\end{aligned}
\end{equation}
and 
\begin{equation}\label{lipschitz:eq-example-2-b}
\begin{aligned}
 &\, \Big[ \UNtil[K][+],  \Big[ \IUNtil[L][+], P^v_{\geq L^{1-\deltap}} \big( \VN[<L^{1-\delta}][\fcs]  -  \VNtil[<L^{1-\delta}][\fcs] \big) \Big]_{\leq N} \Big]_{\leq N}  \\
 -&\, \coup \mathbf{1} \big\{ K=L \big\} \Renormtil[\Ncs][K] P^v_{\geq K^{1-\deltap}} 
 \big( \VN[<K^{1-\delta}][\fcs] - \VNtil[<K^{1-\delta}][\fcs] \big). 
\end{aligned}
\end{equation}
Using a similar argument as in \eqref{null:eq-p-p-qq}, but using \eqref{lipschitz:eq-killing} from Lemma \ref{lipschitz:lem-higher-tensors} instead of Proposition \ref{killing:prop-quadratic}, the $\Cprod{r-1}{r-1}$-norm of \eqref{lipschitz:eq-example-2-a} can be estimated by
\begin{align*}
&\int_{\R} \dy \bigg( \big|\big( \widecheck{\rho}_{\leq N}\ast \widecheck{\rho}_{\leq N}\big)(y)\big| \, 
 \Big\| \Theta^x_y P^v_{\geq L^{1-\deltap}} \VN[<L^{1-\delta}][\fcs] \Big\|_{\Cprod{s}{r-1}} \allowdisplaybreaks[3] \\
&\times \sup_{\substack{E\in \frkg\colon \\ \| E \|_{\frkg}\leq 1}}
\Big\| \Big( \Big[ P_{\leq N}^x \UN[K][+],  \Big[ \Theta^x_y P_{\leq N}^x  \IUN[L][+], E \Big] \Big] -  \coup\mathbf{1} \big\{ K=L \big\} \Cf^{(\Ncs)}_K(y) \Kil(E) \Big) \\
&\hspace{10ex}- \Big( \Big[ P_{\leq N}^x \UNtil[K][+],  \Big[ \Theta^x_y P_{\leq N}^x \IUNtil[L][+], E \Big] \Big] -  \coup\mathbf{1} \big\{ K=L \big\} \Cftil^{(\Ncs)}_K(y) \Kil(E) \Big) \Big\|_{\Cprod{r-1}{s}} \bigg) \allowdisplaybreaks[3] \\
\lesssim&\,  \Dc \times \Big( \max(K,L)^{r-\frac{1}{2}+\eta} L^{-\frac{1}{2}} +  N^{-\delta+\vartheta+\varsigma} \Big) 
\Dc \dpost.
\end{align*}
Since $L\geq K^{1-\delta}$, this is acceptable. Using a similar argument as in \eqref{null:eq-p-p-qq}, but using Lemma \ref{lipschitz:lem-modulated-mixed} instead of Corollary \ref{modulation:cor-control-combined}, we can also estimate the $\Cprod{r-1}{r-1}$-norm of \eqref{lipschitz:eq-example-2-b} by 
\begin{align*}
& \int_{\R} \dy \bigg( \big|\big( \widecheck{\rho}_{\leq N}\ast \widecheck{\rho}_{\leq N}\big)(y)\big| \,\, 
\Big\| \Theta^x_y P^v_{\geq L^{1-\deltap}} \big( \VN[<L^{1-\delta}][\fcs] - \VNtil[<L^{1-\delta}][\fcs] \big)\Big\|_{\Cprod{s}{r-1}}  
\allowdisplaybreaks[3] \\
&\times \sup_{\substack{E\in \frkg\colon \\ \| E \|_{\frkg}\leq 1}}
\Big\|  \Big[ P_{\leq N}^x \UN[K][+],  \Big[ \Theta^x_y P_{\leq N}^x  \IUN[L][+], E \Big] \Big] -  \coup\mathbf{1} \big\{ K=L \big\} \Cf^{(\Ncs)}_K(y) \Kil(E) \Big\|_{\Cprod{r-1}{s}} \bigg)\\ 
\lesssim&\,  \dpost \times \Big( \max(K,L)^{r-\frac{1}{2}+\eta} L^{-\frac{1}{2}} + \mathbf{1}\big\{ K=L \big\} N^{-\delta+\vartheta+\varsigma} \Big) 
\Dc^2.
\end{align*}
Since $L\geq K^{1-\delta}$, this is also acceptable. 
\end{proof}

\begin{remark}
In the estimate of \eqref{lipschitz:eq-example-2-a}, it is important to use the Lipschitz-estimate of
\begin{equation}\label{lipschitz:eq-example-2-a-discussion}
\R \times \frkg \rightarrow \Cprod{r-1}{s}, \, \, (y,E) \mapsto \Big[ P_{\leq N}^x \UN[K][+],  \Big[ \Theta^x_y P_{\leq N}^x  \IUN[L][+], E \Big] \Big] -  \coup\mathbf{1} \big\{ K=L \big\} \Cf^{(\Ncs)}_K(y) \Kil(E) 
\end{equation}
from \eqref{lipschitz:eq-killing} in Lemma \ref{lipschitz:lem-higher-tensors} rather than the individual Lipschitz-estimates of $\UN[K][+]$ and $\IUN[L][+]$ from Lemma \ref{lipschitz:lem-modulated-mixed}. This can readily be seen from our earlier estimate \eqref{null:eq-p-p-qq}, which controlled \eqref{lipschitz:eq-example-2-a-discussion} using Proposition \ref{killing:prop-quadratic} rather than the individual estimates of 
$\UN[K][+]$ and $\IUN[L][+]$ from Lemma \ref{modulation:lem-linear}.
\end{remark}

After the two examples from Lemma \ref{lipschitz:lem-example-1} and Lemma \ref{lipschitz:lem-example-2}, we now give a general algorithm which can be used to modify all remaining estimates needed for the proof of Proposition \ref{lipschitz:prop-remainder-no-jcb}.

\begin{algorithm}[From a-priori to Lipschitz estimates]\label{lipschitz:algorithm}
We describe an algorithm which can be followed to convert our earlier a-priori estimates into Lipschitz estimates. To this end,
we first schematically describe the proof of an a-priori estimate.\\

\emph{A-priori estimate:} We were given a finite index $J$ and terms $(X_j)_{j=1}^J$ from our Ansatz in \eqref{ansatz:eq-UN-rigorous-decomposition} and \eqref{ansatz:eq-VN-rigorous-decomposition}. We were further given a norm $\| \cdot \|_{\mathsf{N}}$, an expression $\mathscr{E}(X_1,X_2,\hdots,X_J)$, and controlled 
\begin{equation*}
\big\| \mathscr{E}\big(X_1,X_2,\hdots,X_J\big)\big\|_{\mathsf{N}}
\end{equation*}
using the following argument. By using the definitions of our modulated and mixed modulated objects and grouping terms together, we wrote
\begin{equation}\label{lipschitz:algorithm-apriori-decomposition}
\mathscr{E}\big( X_1,X_2,\hdots,X_J\big)
= \sum_{i=1}^I \mathscr{M}_i \big( Y_{i,1}, Y_{i,2}, \hdots, Y_{i,J_i} \big).
\end{equation}
In \eqref{lipschitz:algorithm-apriori-decomposition}, $I$ is a finite index, $(J_i)_{i=1}^I$ are finite indices, and $(Y_{i,j})$ are objects which are determined by $X_1,\hdots,X_J$  and have previously been estimated. Furthermore, $(\mathscr{M}_i)_{i=1}^I$ are multi-linear forms. 
Using para-product, integral, trace, and commutator estimates as well as frequency-support considerations, we controlled 
\begin{equation}\label{lipschitz:algorithm-apriori-control-M}
\big\| \mathscr{M}_i \big( Y_{i,1}, Y_{i,2}, \hdots, Y_{i,J_i} \big) \big\|_{\mathsf{N}} 
\leq C_i \prod_{j=1}^{J_i} \big\| Y_{i,j} \big\|_{\mathsf{N}_{i,j}},
\end{equation}
where $C_i>0$ are constants and $\| \cdot \|_{\mathsf{N}_{i,j}}$ are certain norms. Using previous estimates, we then controlled each of the $Y_{i,j}$-terms by 
\begin{equation}\label{lipschitz:algorithm-apriori-control-Y}
\big\| Y_{i,j} \big\|_{\mathsf{N}_{i,j}} \leq C_{i,j},
\end{equation}
where $C_{i,j}>0$ are constants. In total, we therefore obtained the estimate
\begin{equation}\label{lipschitz:algorithm-apriori-final}
\big\| \mathscr{E}\big(X_1,X_2,\hdots,X_J\big)\big\|_{\mathsf{N}} \leq \sum_{i=1}^I C_i \prod_{j=1}^{J_i} C_{i,j}.
\end{equation}

\emph{Lipschitz estimate:} We now describe the necessary modifications to turn \eqref{lipschitz:algorithm-apriori-final} into a Lipschitz estimate. To this end, let $(\widetilde{X}_j)_{j=1}^J$ be terms from our Ansatz in \eqref{ansatz:eq-UN-rigorous-decomposition} and \eqref{ansatz:eq-VN-rigorous-decomposition}, but defined using \eqref{lipschitz:eq-tilde} instead of \eqref{lipschitz:eq-original}. We can then estimate 
\begin{equation*}
\big\|  \mathscr{E}\big(X_1,X_2,\hdots,X_J\big) -  \mathscr{E}\big(\widetilde{X}_1,\widetilde{X}_2,\hdots,\widetilde{X}_J\big) 
\big\|_{\mathsf{N}}
\end{equation*}
using the following argument. We first decompose $\mathscr{E}\big(X_1,X_2,\hdots,X_J\big)$ and  $\mathscr{E}\big(\widetilde{X}_1,\widetilde{X}_2,\hdots,\widetilde{X}_J\big)$ separately using the decomposition from \eqref{lipschitz:algorithm-apriori-decomposition}, which yields 
\begin{equation}\label{lipschitz:algorithm-lipschitz-decomposition}
\begin{aligned}
&\, \mathscr{E}\big(X_1,X_2,\hdots,X_J\big) -  \mathscr{E}\big(\widetilde{X}_1,\widetilde{X}_2,\hdots,\widetilde{X}_J\big) \\ 
=& \, \sum_{i=1}^{I} \Big( \mathscr{M}_i \big( Y_{i,1}, Y_{i,2}, \hdots, Y_{i,J_i} \big)
-  \mathscr{M}_i \big( \widetilde{Y}_{i,1}, \widetilde{Y}_{i,2}, \hdots, \widetilde{Y}_{i,J_i} \big) \Big). 
\end{aligned}
\end{equation}
We then use the multi-linearity of $\mathscr{M}_i$, which allows us to write
\begin{equation}\label{lipschitz:algorithm-lipschitz-multilinear}
\begin{aligned}
&\,  \mathscr{M}_i \big( Y_{i,1}, Y_{i,2}, \hdots, Y_{i,J_i} \big)
-  \mathscr{M}_i \big( \widetilde{Y}_{i,1}, \widetilde{Y}_{i,2}, \hdots, \widetilde{Y}_{i,J_i} \big) \\
=&\, \sum_{j=1}^J \mathscr{M}_i \big( \widetilde{Y}_{i,1}, \hdots, \widetilde{Y}_{i,j-1}, Y_{i,j}- \widetilde{Y}_{i,j}, 
Y_{i,j+1},\hdots, Y_{i,J_i} \big). 
\end{aligned}
\end{equation}
To control the summands in \eqref{lipschitz:algorithm-lipschitz-multilinear}, we then follow the
exact same arguments used in our proof of \eqref{lipschitz:algorithm-apriori-control-M}.
From this, we obtain 
\begin{equation}\label{lipschitz:algorithm-lipschitz-control-M}
\begin{aligned}
&\,\Big\|  \mathscr{M}_i \big( \widetilde{Y}_{i,1}, \hdots, \widetilde{Y}_{i,j-1}, Y_{i,j}- \widetilde{Y}_{i,j}, 
Y_{i,j+1},\hdots, Y_{i,J_i} \big)\Big\|_{\mathsf{N}} \\
\leq&\,  C_i \Big( \prod_{k=1}^{j-1} \big\| \widetilde{Y}_{i,k}\big\|_{\mathsf{N}_{i,k}} \Big) 
\times \big\| Y_{i,j}- \widetilde{Y}_{i,j}\big\|_{\mathsf{N}_{i,j}} 
\times \Big( \prod_{k=j+1}^{J_i} \big\| Y_{i,k} \big\|_{\mathsf{N}_{i,k}} \Big).
\end{aligned}
\end{equation}
Using the same bound as in \eqref{lipschitz:algorithm-apriori-control-Y}, we then control 
\begin{equation}\label{lipschitz:algorithm-lipschitz-control-Y}
 \prod_{k=1}^{j-1} \big\| \widetilde{Y}_{i,k}\big\|_{\mathsf{N}_{i,k}}\leq  \prod_{k=1}^{j-1} C_{i,k}
 \qquad \text{and} \qquad 
 \prod_{k=j+1}^{J_i} \big\| Y_{i,k} \big\|_{\mathsf{N}_{i,k}} \leq 
 \prod_{k=j+1}^{J_i} C_{i,k}. 
\end{equation}
In total, we therefore obtain the estimate
\begin{align}\label{lipschitz:algorithm-lipschitz-final}
\Big\| \mathscr{E}\big(X_1,X_2,\hdots,X_J\big) -  \mathscr{E}\big(\widetilde{X}_1,\widetilde{X}_2,\hdots,\widetilde{X}_J\big) \Big\|_{\mathsf{N}}
\leq \sum_{i=1}^I \sum_{j=1}^{J_i} C_i \, \Big( \prod_{\substack{1\leq k \leq J_i\colon \\ k\neq j}} \hspace{-2ex} C_{i,k} \Big) 
\big\| Y_{i,j} - \widetilde{Y}_{i,j} \big\|_{\mathsf{N}_{i,j}}.
\end{align}
In particular, we can derive Lipschitz estimates of $\mathscr{E}(X_1,X_2,\hdots,X_J)$ with respect to the $\|\cdot\|_{\mathsf{N}}$-norm from Lipschitz estimates of the $Y_{i,j}$-terms with respect to the $\|\cdot\|_{\mathsf{N}_{i,j}}$-norms.  
\end{algorithm}

\begin{remark}
In our proofs of Lemma \ref{lipschitz:lem-example-1} and Lemma \ref{lipschitz:lem-example-2}, we already follow Algorithm \ref{lipschitz:algorithm}. In both proofs, we only need one multi-linear form, i.e., it holds that $I=1$. In the proof of Lemma \ref{lipschitz:lem-example-1}, we neither insert the definition of a modulated or mixed modulated object nor group different terms together. As a result, both the original expression and multi-linear form can be written in terms of 
\begin{equation*}
X_1 = Y_{1,1} = \UN[K][+] \qquad \text{and} \qquad 
X_2 = Y_{1,2} = \VN[M_u,M_v][+-]. 
\end{equation*}
In the proof of Lemma \ref{lipschitz:lem-example-2}, however, we treat the tensor product of $P_{\leq N}^x \UN[K][+]$ and $\Theta^x_y P_{\leq N}^x \IUN[L][+]$ as a combined object. While the original expression can be written in terms of 
\begin{equation*}
X_1= \UN[K][+], \qquad X_2 = \IUN[L][+], \qquad \text{and} \qquad X_3 = \VN[<L^{1-\delta}][\fcs], 
\end{equation*}
we choose $Y_{1,1}$ as the map 
\begin{equation*}
\R \times \frkg \rightarrow \Cprod{r-1}{s}, \, \, (y,E) \mapsto \Big[ P_{\leq N}^x \UN[K][+],  \Big[ \Theta^x_y P_{\leq N}^x  \IUN[L][+], E \Big] \Big] -  \coup\mathbf{1} \big\{ K=L \big\} \Cf^{(\Ncs)}_K(y) \Kil(E) 
\end{equation*}
from \eqref{lipschitz:eq-example-2-a-discussion} and choose $Y_{1,2}=\VN[<L^{1-\delta}][\fcs]$.
\end{remark}

At the end of this section, we record a Lipschitz-variant of \eqref{jacobi:eq-bourgain-bulut-1} from Proposition \ref{jacobi:prop-main}.

\begin{lemma}[\protect{Lipschitz-variant of \eqref{jacobi:eq-bourgain-bulut-1}}]\label{lipschitz:lem-bourgain-bulut}
Let the post-modulation hypothesis (Hypothesis \ref{hypothesis:post}) be satisfied for  \eqref{lipschitz:eq-original} and \eqref{lipschitz:eq-tilde}. For all $V\colon \R^{1+1}\rightarrow \frkg$, it then holds that 
\begin{equation}\label{lipschitz:eq-bourgain-bulut}
\begin{aligned}
&\, \Sumlarge_{\substack{\hspace{1.5ex} K \geq N^{1-2\delta_1}}} \Big\| \, \chi 
\Big( P_{\leq N}^x \CHHLN \big( V, V \big) P_{\leq N}^x - \coup \Renormbd[\Ncs] \Big) \big(   \UN[K][+] - \UNtil[K][+]\big) \Big\|_{\Cprod{r-1}{r-1}} \\
\lesssim&\, N^{\delta_1+5\delta_3}  \sup_{y\in \R} \Big( \, \langle Ny \rangle^{-10} \Big\|  \chi \Big( \CHHLN_y\big( V, V \big) - \coup \CNcsbd(y) \Kil \Big) \Big\|_{L_t^\infty L_x^\infty} \Big) \,  \dpost.
\end{aligned}
\end{equation}
\end{lemma}

\begin{remark}
We emphasize that \eqref{lipschitz:eq-bourgain-bulut} only yields the Lipschitz-continuous dependence on $\UN[K][+]$, but makes no claim regarding the Lipschitz-continuous dependence on $V$. The reason is that the Cartesian high$\times$high$\rightarrow$low-interaction $\CHHLN_y\big( V, V \big) - \coup \CNcsbd(y) \Kil $ will not be controlled using perturbative arguments, but will instead be controlled using the almost invariance of the Gibbs measure (see Proposition \ref{main:prop-almost-invariance} and Corollary \ref{main:cor-bourgain-bulut-probability}).
\end{remark}

\begin{proof}[Proof of Lemma \ref{lipschitz:lem-bourgain-bulut}]
The desired inequality \eqref{lipschitz:eq-bourgain-bulut} can be obtained using a minor variant of the proof of \eqref{jacobi:eq-bourgain-bulut-1}. Indeed, by first using the identity \eqref{jacobi:eq-main-p0} and then using \eqref{jacobi:eq-F-Up-p2} and \eqref{jacobi:eq-F-Up-p3} from the proof of Lemma \ref{jacobi:lem-F-Up}, we obtain that
\begin{align*}
&\,  \Big\| \, \chi 
\Big( P_{\leq N}^x \CHHLN \big( V, V \big) P_{\leq N}^x - \coup \Renormbd[\Ncs] \Big) \big(   \UN[K][+] - \UNtil[K][+]\big) \Big\|_{\Cprod{r-1}{r-1}} \\
\lesssim&\, N^{\delta_1+4\delta_3} \sup_{y\in \R} \Big( \, \langle Ny \rangle^{-10} \Big\|  \chi \Big( \CHHLN_y\big( V, V \big) - \coup \CNcsbd(y) \Kil \Big) \Big\|_{L_t^\infty L_x^\infty} \Big) \\
&\times \Big( K^{-(r-\frac{1}{2}+\eta)} \big\| \UN[K][+] - \UNtil[K][+] \big\|_{\Cprod{r-1}{s}} 
+ K^{-(r-\frac{1}{2}+4\eta)}  \big\| \UN[K][+] - \UNtil[K][+] \big\|_{\Cprod{r-1+3\eta}{s}} \Big). 
\end{align*}
By using Lemma \ref{lipschitz:lem-modulated-mixed}, we then obtain the desired inequality \eqref{lipschitz:eq-bourgain-bulut}.
\end{proof}

\section{Lifting}\label{section:lifting}

In Sections \ref{section:ansatz}-\ref{section:energy-increment}, we only worked with the derivatives $A,B\colon \R^{1+1}\rightarrow \frkg$. 
In order to prove Theorem \ref{intro:thm-rigorous-phi}, we will later need to transfer our result from the derivatives $A,B\colon \R^{1+1}\rightarrow \frkg$ to the original map $\phi \colon \R^{1+1}\rightarrow \frkG$, which is the subject of this section. 
In Subsection \ref{section:lifting-lifts}, we focus on the smooth setting and examine the relationship between the derivatives $A$ and $B$ and the map $\phi$ using lifts. In Subsection \ref{section:lifting-enhanced}, we introduce the enhanced data spaces, which are low-regularity spaces, and obtain the Lipschitz continuity of our lifts with respect to the corresponding metric. 
In Subsection \ref{section:lifting-Gibbs}, we introduce the Gibbs measure $\muGb$, which is the Gibbs measure corresponding to the $\frkG$-valued map $\phi$. 

\subsection{Lifts}\label{section:lifting-lifts}

As mentioned above, we momentarily work in the smooth setting, and postpone a discussion of the low-regularity setting until Subsection \ref{section:lifting-enhanced}. In Section \ref{section:introduction}, we previously introduced two different formulations of the wave maps equation: In the first formulation, we let $(\phi_0,\phi_1) \colon \R \rightarrow  T\frkG$ be smooth maps satisfying $\phi_1(x) \in T_{\phi_0(x)} \frkG$ for all $x\in \R$, i.e., $\phi_0(x)^{-1} \phi_1(x)\in \frkg$ for all $x\in \R$. The wave maps equation with initial data $(\phi_0,\phi_1)$ for a map $\phi\colon \R^{1+1}\rightarrow \frkG$ can then be written as 
\begin{equation}\label{lifting:eq-WM-phi}
\begin{cases}
\begin{aligned}
\partial_t \big( \phi^{-1} \partial_t \phi \big) = \partial_x \big( \phi^{-1} \partial_x \phi\big), \\ 
\phi(0,\cdot) = \phi_0, \, \partial_t \phi(0,\cdot) = \phi_1. 
\end{aligned}
\end{cases}
\end{equation}
In the second formulation, which is at the level of derivatives, we let $A_0,B_0\colon \R \rightarrow \frkg$ be smooth maps. The wave maps equation with initial data $(A_0,B_0)$ for $A,B\colon \R^{1+1}\rightarrow \frkg$ can then be written as 
\begin{equation}\label{lifting:eq-WM-AB}
\begin{cases}
\begin{aligned}
\partial_t A &= \partial_x B, \\ 
\partial_t B &= \partial_x A - \big[ A , B \big], \\ 
A(0,\cdot) &= A_0,\,\, B(0,\cdot)=B_0.
\end{aligned}
\end{cases}
\end{equation}
In Section \ref{section:introduction}, we claimed that the solutions of the original formulation \eqref{lifting:eq-WM-phi} and the derivative formulation \eqref{lifting:eq-WM-AB} can be linked via the identities
\begin{equation}\label{lifting:eq-link}
A = \phi^{-1} \partial_t \phi \qquad \text{and} \qquad B= \phi^{-1} \partial_x \phi.
\end{equation}
In the following lemma, we now make this link between \eqref{lifting:eq-WM-phi} and \eqref{lifting:eq-WM-AB} precise. 

\begin{lemma}[Equivalence]\label{lifting:lem-equivalence}
Let $\phi \colon \R^{1+1}\rightarrow \frkG$ be a smooth solution of \eqref{lifting:eq-WM-phi} and let $A,B\colon \R^{1+1}\rightarrow \frkg$ be defined as in \eqref{lifting:eq-link}. Then, $A$ and $B$ are smooth solutions of \eqref{lifting:eq-WM-AB} with initial data given by $A_0 = \phi_0^{-1} \phi_1$ and $B_0=\phi_0^{-1} \partial_x \phi_0$. 

Alternatively, let $A,B\colon \R^{1+1}\rightarrow \frkg$ be smooth solutions of \eqref{lifting:eq-WM-AB} and let $g\in \frkG$. Then, there exist unique smooth maps $(\phi_0,\phi_1) \colon \R\rightarrow T \frkG$ satisfying 
\begin{equation*}
\phi_0(0)=g, \quad \partial_x \phi_0(x) = \phi_0(x) B_0(x), \quad \text{and} \quad  \phi_1(x) = \phi_0(x) A_0(x)
\end{equation*}
and a unique smooth map $\phi \colon \R^{1+1} \rightarrow \frkG$ satisfying $\phi(0,0)=g$ and \eqref{lifting:eq-link}. Furthermore, the map $\phi$ is a smooth solution of \eqref{lifting:eq-WM-phi}.
\end{lemma}

The result from Lemma \ref{lifting:lem-equivalence} is mostly contained in \cite[Proposition 2.1]{TU04}. For the reader's benefit, we still present a self-contained proof.

\begin{proof}[Proof of Lemma \ref{lifting:lem-equivalence}:]
First, let $\phi$ be a smooth solution of \eqref{lifting:eq-WM-phi} and let $A$ and $B$ be defined as in \eqref{lifting:eq-link}. The first equation in \eqref{lifting:eq-WM-AB} then follows directly from the evolution equation in \eqref{lifting:eq-WM-phi}. The second equation in \eqref{lifting:eq-WM-AB} can be derived from the identity $\partial_t \partial_x \phi = \partial_x \partial_t \phi$. More precisely, it holds that
\begin{align*}
\partial_t B - \partial_x A = \partial_t \big( \phi^{-1} \partial_x \phi \big) - \partial_x \big( \phi^{-1} \partial_x \phi \big) 
= - \phi^{-1} (\partial_t \phi) \phi^{-1} (\partial_x \phi) + \phi^{-1} (\partial_x \phi) \phi^{-1} (\partial_t \phi)  
= - \big[ A , B \big].
\end{align*}
Finally, the initial condition in \eqref{lifting:eq-WM-AB} follows directly from the initial condition in \eqref{lifting:eq-WM-phi} and the definitions of $A_0$ and $B_0$, which completes the proof of the first part of the lemma. \\ 

Second, let $g\in \frkG$, let $A_0,B_0\colon \R \rightarrow \frkG$ be smooth, and let $A,B \colon \R^{1+1} \rightarrow \frkG$ be smooth solutions of \eqref{lifting:eq-WM-AB}. The global existence and uniqueness of $\phi_0$ and $\phi_1$ follow directly from the Picard-Lindel\"{o}f theorem and the compactness of $\frkG$. The conditions on the map $\phi\colon \R^{1+1} \rightarrow \frkG$ can be written as 
\begin{equation}\label{lifting:eq-equivalence-p1}
\partial_t \phi = \phi A, \quad \partial_x \phi = \phi B, \quad \text{and} \quad \phi(0,0)=g. 
\end{equation}
The second equation in \eqref{lifting:eq-WM-AB} implies the compatibility of the two differential equations $\partial_t \phi = \phi A$ and $\partial_x \phi = \phi B$. Together with the compactness of $\frkG$, this implies the global existence and uniqueness of a solution to \eqref{lifting:eq-equivalence-p1}. From \eqref{lifting:eq-link} and the first equation in \eqref{lifting:eq-WM-AB}, it also follows that
\begin{equation*}
\partial_t \big( \phi^{-1} \partial_t \phi \big) = \partial_t A = \partial_x B = \partial_x \big( \phi^{-1} \partial_x \phi\big),
\end{equation*}
and thus $\phi$ solves the evolution equation in \eqref{lifting:eq-WM-phi}. Finally, the initial conditions in \eqref{lifting:eq-WM-phi} follow directly from the definitions of $\phi_0$, $\phi_1$, and $\phi$.
\end{proof}

While Lemma \ref{lifting:lem-equivalence} implies that $A$ and $B$ determine $\phi$ (modulo multiplication by an element of $\frkG$), this alone is insufficient for our purposes since we will need to transfer our estimates from $A$ and $B$ to $\phi$. In order to transfer our estimates, we need a more explicit representation\footnote{The representation from Lemma \ref{lifting:lem-representation} below has implicitly been used in the proof of Lemma \ref{lifting:lem-equivalence}, since it is typically used to prove that \eqref{lifting:eq-equivalence-p1} indeed has a global solution $\phi$.} of $\phi$ in terms of $A$ and $B$. To this end, we first make the following definition.

\begin{definition}[Lifts]\label{lifting:def-lifts}
Let $A\colon \R_t \rightarrow \frkG$ and let $B\colon \R_x \rightarrow \frkG$ be smooth. Then, we define $\ac=\ac[A]\colon \R_t \rightarrow \frkG$
and $\bc=\bc[B]\colon \R_x \rightarrow \frkG$ as the unique solutions of 
\begin{alignat}{3}
\partial_t \ac(t) &= \ac(t) A(t), &\qquad \ac(0) = \ec, \\ 
\partial_x \bc(x) &= \bc(x) B(x), &\qquad \bc(0) = \ec,
\end{alignat}
where $\ec$ is the unit element of $\frkG$. Furthermore, for all $t_0\in \R$, we also define $\ac[t_0,A]\colon \R_t \rightarrow \frkG$ as the unique solution of 
\begin{equation*}
\partial_t \ac[t_0,A](t) = \ac[t_0,A](t) A(t), \qquad \ac[t_0,A](t_0)=\ec. 
\end{equation*}
\end{definition}

Equipped with Definition \ref{lifting:def-lifts}, we now derive the desired representation of $\phi$ in terms of $A$ and $B$.

\begin{lemma}[Representation using lifts]\label{lifting:lem-representation}
Let $\phi \colon \R^{1+1} \rightarrow \frkG$ be a smooth solution of \eqref{lifting:eq-WM-phi}, let $g=\phi(0,0)\in \frkG$, and let $A,B\colon \R^{1+1}\rightarrow \frkg$ be as in \eqref{lifting:eq-link}. Then, it holds that
\begin{align}
\phi(t,x) &= g \cdot \ac[A(\cdot,0)](t)\cdot \bc[B(t,\cdot)](x), \label{lifting:eq-representation-1} \\
\partial_t \phi(t,x) &= g \cdot \ac[A(\cdot,0)](t)  \cdot \bc[B(t,\cdot)](x)\cdot A(t,x). \label{lifting:eq-representation-2}
\end{align}
\end{lemma}

\begin{proof}[Proof of Lemma \ref{lifting:lem-representation}:]
We first write the map $\phi$ as 
\begin{equation}\label{lifting:eq-representation-p1}
\phi(t,x) = \phi(0,0) \phi(0,0)^{-1} \phi(t,0) \phi(t,0)^{-1} \phi(t,x).     
\end{equation}
Using \eqref{lifting:eq-link}, we obtain that 
\begin{equation*}
\partial_t \big( \phi(0,0)^{-1} \phi(t,0) \big) = \phi(0,0)^{-1} \partial_t \phi(t,0) = \phi(0,0)^{-1} \phi(t,0) A(t,0). 
\end{equation*}
From Definition \ref{lifting:def-lifts}, it then follows that $\phi(0,0)^{-1} \phi(t,0)=\ac[A(\cdot,0)](t)$. Using a similar argument, it also follows that $\phi(t,0)^{-1} \phi(t,x)=\bc[B(t,\cdot)](x)$. By inserting both identities into \eqref{lifting:eq-representation-p1}, we obtain the desired first identity \eqref{lifting:eq-representation-1}. In order to obtain \eqref{lifting:eq-representation-2}, we first calculate
\begin{equation}\label{lifting:eq-representation-p2}
\partial_t \phi(t,x) = \phi(t,x) \phi(t,x)^{-1} \partial_t \phi(t,x) = \phi(t,x) A(t,x). 
\end{equation}
After inserting \eqref{lifting:eq-representation-1} into \eqref{lifting:eq-representation-p2}, we then obtain the desired second identity \eqref{lifting:eq-representation-2}.
\end{proof}

\subsection{Enhanced data space and properties of lifts}\label{section:lifting-enhanced}

We first introduce the enhanced data space. While elements of the enhanced data space can live at low regularities, we impose enough conditions to later define lifts on the enhanced data space (see Lemma \ref{lifting:lem-lift-properties}).

\begin{definition}[Enhanced data space]\label{main:def-enhanced-data}
For any $R\geq 1$ and smooth $A,\widetilde{A}\colon \bT_R \rightarrow \frkg$, we define the metric\footnote{The condition $K,L\geq 2$ in the dyadic sum is added to ensure the periodicity of the $x$-integrals. Of course, at the cost of powers of $R$, the contribution of $K,L\lesssim 1$ can always be bounded using only the $\C_x^{s-1}$-norms of $A$ and $\widetilde{A}$.}
\begin{equation}\label{lifting:eq-enhanced-data-A}
\begin{aligned}
&\,\big\| A \,; \widetilde{A}\,\big\|_{\Xi^{s-1}_x(\bT_R\rightarrow \frkg)} \\ 
:=&\, \big\| A - \widetilde{A} \,  \big\|_{\C_x^{s-1}(\bT_R\rightarrow \frkg)}
+ \sum_{\substack{K,L\geq 2 \colon \\ K \sim L}}\big\| P_K^x A \otimes \Int^x_{0\rightarrow x} P_L^x A - P_K^x \widetilde{A} \otimes \Int^x_{0\rightarrow x} P_L^x \widetilde{A} \, \big\|_{\C_x^{2s-1}(\bT_R\rightarrow \frkg \otimes \frkg)}.
\end{aligned}
\end{equation}
For any smooth $A,B,\widetilde{A},\widetilde{B}\colon \bT_R \rightarrow \frkg$, we also define the metric
\begin{equation}\label{lifting:eq-enhanced-data-AB}
\begin{aligned}
\big\| (A,B)\, ; (\widetilde{A},\widetilde{B}) \big\|_{\Xi^{s-1}_x(\bT_R\rightarrow \frkg^2)}
:=&\,  \big\| A \,; \widetilde{A}\,\big\|_{\Xi^{s-1}_x(\bT_R\rightarrow \frkg)} 
+ \big\| B \,; \widetilde{B}\,\big\|_{\Xi^{s-1}_x(\bT_R\rightarrow \frkg)} \\ 
+&\,\sum_{\substack{K,L\geq 2 \colon \\ K \sim L}}\big\| P_K^x A \otimes \Int^x_{0\rightarrow x} P_L^x B - P_K^x \widetilde{A} \otimes \Int^x_{0\rightarrow x} P_L^x \widetilde{B} \, \big\|_{\C_x^{2s-1}(\bT_R\rightarrow \frkg  \otimes \frkg)} \\ 
+&\,  \sum_{\substack{K,L\geq 2\colon \\ K \sim L}}\big\| P_K^x B \otimes \Int^x_{0\rightarrow x} P_L^x A - P_K^x \widetilde{B} \otimes \Int^x_{0\rightarrow x} P_L^x \widetilde{A} \, \big\|_{\C_x^{2s-1}(\bT_R\rightarrow \frkg  \otimes \frkg)}.
\end{aligned}
\end{equation}
The corresponding function spaces $\Xi^{s-1}_x(\bT_R\rightarrow \frkg)$ and $\Xi^{s-1}_x(\bT_R\rightarrow \frkg^2)$ are defined as the completions of  $C^\infty_x(\bT_R\rightarrow \frkg)$ and $C^\infty_x(\bT_R\rightarrow \frkg^2)$. 
\end{definition}

While Definition \ref{main:def-enhanced-data} captures the important aspects of our current functional framework, we also need the following related definitions.

\begin{definition}[Enhanced data spaces again]\label{lifting:def-enhance-data-again}
We make the following definitions.
\begin{enumerate}[label=(\roman*)]
    \item For smooth $A,B,\widetilde{A},\widetilde{B}\colon \R \rightarrow \frkg$, we define the $\Xi^{s-1}_x(\R\rightarrow \frkg)$ and
    $\Xi^{s-1}_x(\R\rightarrow \frkg^2)$-metrics as in \eqref{lifting:eq-enhanced-data-A} and \eqref{lifting:eq-enhanced-data-AB}, but with $\bT_R$ replaced with $\R$.
    \item Similar as in \eqref{prelim:eq-Hoelder-local}, for each compact interval $J\subseteq \R$, the $\Xi^{s-1}_x(J\rightarrow \frkg)$ and
    $\Xi^{s-1}_x(J\rightarrow \frkg^2)$-metrics are defined via restriction of the $\Xi^{s-1}_x(\R\rightarrow \frkg)$ and
    $\Xi^{s-1}_x(\R\rightarrow \frkg^2)$-metrics. 
    \item The function spaces $\Xi^{s-1}_t(\R\rightarrow \frkg)$ and
    $\Xi^{s-1}_t(\R\rightarrow \frkg^2)$, which are viewed as function spaces for time-dependent maps, are defined exactly as for maps depending on the spatial variable.
\end{enumerate}
\end{definition}
 
In Definition \ref{main:def-enhanced-data} and Definition \ref{lifting:def-enhance-data-again}, 
we introduced function spaces for $\frkg$-valued maps which will be used for the arguments of our lifts. In order to state the mapping properties of lifts, we will also need to introduce function spaces for the $\frkG$-valued maps. To this end, we first let $T\geq 1$ and let $C_t^0([0,T]\rightarrow \frkG)$ be the space of continuous maps from $[0,T]$ to $\frkG$. This space is equipped with the metric
\begin{equation*}
\big\| \phi_1 ; \phi_2 \big\|_{C_t^0([0,T]\rightarrow \frkG)} := \sup_{t\in [0,T]} d_{\frkG}\big( \phi_1(t),\phi_2(t) \big), 
\end{equation*}
where $d_\frkG$ is the intrinsic distance on the Riemannian manifold $\frkG$. In order to define additional function spaces, we now take an extrinsic perspective. 
Due to the Peter-Weyl theorem, we have previously assumed that $\frkG$ is a closed subgroup of $\textup{Gl}(n,\bC)$. Due to the canonical embedding of the group $\textup{Gl}(n,\bC)$ into the ambient vector space $\bC^{n\times n}$, we can therefore embed $\frkG$ into $\bC^{n\times n}$ and $T\frkG$ into $\bC^{n\times n}\times \bC^{n\times n}$. In particular, we may then define all H\"{o}lder spaces for $\frkG$-valued and $T\frkG$-valued maps via these embeddings.

\begin{lemma}\label{lifting:lem-lift-properties}
The lifts $\ac=\ac[\cdot]$ and $\bc=\bc[\cdot]$ satisfy the following properties: 
\begin{enumerate}[label=(\roman*)]
\item Let $A,\widetilde{A}\colon \R_t \rightarrow \frkg$ satisfy 
$\| A; 0\|_{\Xi_t^{s-1}(\R \rightarrow \frkg)}, \| \widetilde{A}; 0\|_{\Xi_t^{s-1}(\R \rightarrow \frkg)} \lesssim 1$. Then, it holds for all $t_0 \in \R$ that 
\begin{equation*}
\big\| \ac[t_0,A](t); \ac[t_0,\widetilde{A}](t) \big\|_{C_t^0 ([t_0-1,t_0+1]\rightarrow \frkG)} 
\lesssim \big\| A; \widetilde{A} \big\|_{\Xi_t^{s-1}(\R\rightarrow \frkg)}.
\end{equation*}
\item The map $B\mapsto \bc[B]$ is locally Lipschitz as a map from $\Xi^{s-1}_x([-R,R]\rightarrow \frkg)$ into $\C_x^{s}([-R,R]\rightarrow \frkG)$
for all $R\geq 1$. 
\item The map $(A,B)\mapsto (\bc[B],\bc[B]A )$ is locally Lipschitz as a map from $\Xi^{s-1}_x([-R,R]\rightarrow \frkg^2)$ into the space $(\C_x^s\times \C_x^{s-1})([-R,R]\rightarrow T\frkG)$ for all $R\geq 1$.
\end{enumerate}
\end{lemma}

\begin{remark} 
The statement of Lemma \ref{lifting:lem-lift-properties} is intended to match the statements in Section \ref{section:main} below. In particular, this is the reason for stating quantitative estimates for $\ac$, but only stating qualitative estimates for $\bc$.
\end{remark}

The claims in Lemma \ref{lifting:lem-lift-properties} can easily be obtained from the para-controlled approach to rough ordinary differential equations, see e.g. \cite[Section 3]{GIP15}, and we therefore omit the details.

\subsection{\protect{The Gibbs measure $\muGb$}}\label{section:lifting-Gibbs}

In this subsection, we give the precise definition of the Gibbs measure $\muGb$, which already appeared in the statement of Theorem \ref{intro:thm-rigorous-phi}. 
The Gibbs measure $\muGb$ is built out of the Haar measure on $\frkG$, $\frkG$-valued Brownian motion, and $\frkg$-valued white noise. While we have extensively worked with $\frkg$-valued white noise, we have yet to work with the Haar measure on $\frkG$ or $\frkG$-valued Brownian motion, and we now briefly recall their definitions and properties. 

\begin{remark}[Inverse temperature vs. temperature]
In contrast to Sections \ref{section:ansatz}-\ref{section:energy-increment}, it is now more convenient to work with the inverse temperature $\beta>0$. The corresponding temperature is then given by $\beta^{-1}$ which, due to \eqref{intro:eq-coup}, can also be written as $8\coup$.
\end{remark}

We first define the Haar measure on $\frkG$. To this end, we let 
$\mathcal{B}(\frkG)$ be the Borel $\sigma$-algebra on the Lie group $\frkG$. 

\begin{definition}[Haar measure]\label{lifting:def-Haar}
We define the Haar measure $\Haar$ as the unique probability measure on $(\frkG,\mathcal{B}(\frkG))$ which is invariant under left and right-multiplication, i.e., which satisfies $\Haar(g\cdot E)=\Haar(E\cdot g)=\Haar(E)$ for all $g\in \frkG$ and $E\in \mathcal{B}(\frkG)$. 
\end{definition}

For the existence and uniqueness of the Haar measure $\Haar$, we refer the reader to \cite[Chapter 1.4]{T14}. In the literature, one often considers Haar measures which are invariant under left-multiplication or right-multiplication, but not necessarily both. On compact Lie groups, however, invariance under left-multiplication and right-multiplication are equivalent. \\

We now turn to the next building block of the Gibbs measure $\muGb$, which is $\frkG$-valued Brownian motion.

\begin{definition}[Brownian motion on $\frkG$]\label{lifting:def-BM}
A stochastic process $\psi_\beta \colon \R \rightarrow \frkG$ is called a $\frkG$-valued Brownian motion at inverse temperature $\beta>0$ if it can be written as the solution of the Stratonovich stochastic differential equation
\begin{equation*}
\mathrm{d}\psi_\beta(x) =  \psi_\beta(x) \circ \mathrm{d}W_\beta(x), \qquad \psi_\beta(0)=\ec,
\end{equation*}
where $W_\beta\colon \R \rightarrow \frkg$ is a $\frkg$-valued white noise at inverse temperature $\beta>0$ and $\ec$ is the unit element of $\frkG$.
\end{definition}

Many further aspects of $\frkG$-valued Brownian motion, such as the connection to the heat kernel on $\frkG$, are beyond the scope of this article and we instead refer the reader to \cite{C04,D17,I50}. In the following, we primarily work with a characterization of $\frkG$-valued Brownian motion using the lifts from Definition \ref{lifting:def-lifts}, which is the subject of the next lemma. 

\begin{lemma}[Brownian motion on $\frkG$ and lifts]\label{lifting:lem-BM-lifting}
If $W_\beta\colon \R \rightarrow \frkg$ is a $\frkg$-valued white noise at inverse temperature $\beta>0$, then $\bc[W_\beta]\colon \R \rightarrow \frkG$ is a $\frkG$-valued Brownian motion at inverse temperature $\beta>0$. Conversely, if $\psi_\beta\colon \R \rightarrow \frkG$ is a $\frkG$-valued Brownian motion at inverse temperature $\beta>0$, then there exists a $\frkg$-valued white noise $W_\beta\colon \R \rightarrow \frkg$ at inverse temperature $\beta>0$ such that $\psi_\beta=\bc[W_\beta]$. 
\end{lemma}

\begin{proof}
Let $W_\beta\colon \R \rightarrow \frkg$ be a $\frkg$-valued white noise at inverse temperature $\beta>0$ and let $\psi_\beta\colon \R \rightarrow \frkG$ be the solution of the Stratonovich stochastic differential equation 
\begin{equation*}
\mathrm{d}\psi_\beta(x)= \psi_\beta(x) \circ \mathrm{d}W_\beta(x), \qquad \psi_\beta(0)=\ec. 
\end{equation*}
In order to prove the lemma, it suffices to prove that $\bc[W_\beta]=\psi_\beta$. Using Lemma \ref{lifting:lem-lift-properties}, we first obtain that 
\begin{equation}\label{lifting:eq-BM-lifting-p1}
\bc[W_\beta] = \lim_{N\rightarrow \infty} \bc[P_{\leq N} W_\beta], 
\end{equation}
where the convergence holds in $\C_x^s([-R,R]\rightarrow \frkG)$ for all $R \geq 1$. From the definition of the lift $\bc$, it follows that $\bc[P_{\leq N} W_\beta]$ is a solution of the ordinary differential equation 
\begin{equation*}
\partial_x \bc[P_{\leq N} W_\beta](x) = \bc[P_{\leq N} W_\beta](x) P_{\leq N}W_\beta(x), \qquad \bc[P_{\leq N} W_\beta](0)=\ec. 
\end{equation*}
Using the Wong-Zakai theorem (see e.g. \cite[Chapter 6]{IW89}), it then follows that 
\begin{equation}\label{lifting:eq-BM-lifting-p2}
\lim_{N\rightarrow \infty} \bc[P_{\leq N}W_\beta]= \psi_\beta, 
\end{equation}
where the convergence holds in $\C_x^s([-R,R]\rightarrow \frkG)$ for all $R \geq 1$. By combining \eqref{lifting:eq-BM-lifting-p1} and \eqref{lifting:eq-BM-lifting-p2}, we then obtain the desired identity $\bc[W_\beta]=\psi_\beta$, which completes the proof.
\end{proof}

Equipped with both the Haar measure (Definition \ref{lifting:def-Haar}) and $\frkG$-valued Brownian motion (Definition \ref{lifting:def-BM}), we can now give a precise definition of the Gibbs measure $\muGb$.

\begin{definition}[Gibbs measure]\label{lifting:def-Gibbs}
Let $\beta>0$ be an inverse temperature. Then, we define 
\begin{equation*}
\muGb = \operatorname{Law} \Big( \big( g \psi_\beta(x), g \psi_\beta(x) W_\beta(x) \big) \Big), 
\end{equation*}
where 
\begin{enumerate}[label=(\roman*)]
    \item $g$ is a $\frkG$-valued random variable whose law is given by the Haar measure $\Haar$, 
    \item $\psi_\beta$ is a $\frkG$-valued Brownian motion at inverse temperature $\beta>0$, 
    \item $W_\beta$ is a $\frkg$-valued white noise at inverse temperature $\beta>0$, 
    \item and $g$, $\psi_\beta$, and $W_\beta$ are independent. 
\end{enumerate}
\end{definition}

\begin{remark}
We note that, due to the independence of $\psi_\beta$ and $W_\beta$, the product $\psi_\beta(x) W_\beta(x)$ can easily be defined as a spatial distribution.
\end{remark}

\begin{lemma}[Representation of Gibbs measure]\label{lifting:lem-Gibbs-representation}
Let $\beta>0$ be an inverse temperature and let 
\begin{enumerate}[label=(\alph*)]
\item\label{lifting:item-rep-a} $g$ be a $\frkG$-valued random variable whose law is given by the Haar measure $\Haar$, 
\item\label{lifting:item-rep-b} $h$ be a $\frkG$-valued random variable, 
\item\label{lifting:item-rep-c} $W_{\beta,0},W_{\beta,1}\colon \R \rightarrow \frkg$ be $\frkg$-valued white noises at inverse temperature $\beta>0$.
\end{enumerate}
Furthermore, assume that 
\begin{enumerate}[label=(\alph*)]
\setcounter{enumi}{3}
\item\label{lifting:item-rep-d} $g$, $W_{\beta,0}$, and $W_{\beta,1}$ are independent 
\item\label{lifting:item-rep-e} and $h$ is determined by $W_{\beta,0}$ and $W_{\beta,1}$. 
\end{enumerate}
Then, it holds that
\begin{equation*}
\operatorname{Law} \Big( \big( g \cdot h \cdot \bc\big[ W_{\beta,0}\big], g \cdot h \cdot \bc\big[ W_{\beta,0}\big] \cdot W_{\beta,1} \big) \Big) = \muGb.
\end{equation*}
\end{lemma}

\begin{proof}
Due to Definition \ref{lifting:def-Gibbs}, it suffices to show that 
\begin{enumerate}[label=(\roman*)]
    \item\label{lifting:item-rep-1} the law of $g\cdot h$ is given by the Haar measure $\Haar$, 
    \item\label{lifting:item-rep-2} $\bc[W_{\beta,0}]$ is a $\frkG$-valued Brownian motion at inverse temperature $\beta>0$, 
    \item\label{lifting:item-rep-3} $W_{\beta,1}$ is a $\frkg$-valued white noise at inverse temperature $\beta>0$, 
    \item\label{lifting:item-rep-4} and $g\cdot h$, $W_{\beta,0}$, and $W_{\beta,1}$ are independent. 
\end{enumerate}
The properties in \ref{lifting:item-rep-2} and \ref{lifting:item-rep-3} follow directly from \ref{lifting:item-rep-c} and Lemma \ref{lifting:lem-BM-lifting}. Furthermore, due to \ref{lifting:item-rep-d}, $W_{\beta,0}$ and $W_{\beta,1}$ are independent. In order to obtain \ref{lifting:item-rep-1} and \ref{lifting:item-rep-4}, it therefore remains to prove that the law of $g\cdot h$ is given by the Haar measure and that $g\cdot h$ is independent of the pair $(W_{\beta,0},W_{\beta,1})$. To this end, it suffices to prove for all $E\in \mathcal{B}(\frkG)$ that
\begin{align}\label{lifting:eq-rep-p}
\mathbb{P} \big( g\cdot h \in E \big| \,W_{\beta,0},W_{\beta,1} \big) 
= \Haar\big( E \big) = \mathbb{P} \big( g\cdot h \in E \big).
\end{align}
Due to \ref{lifting:item-rep-e}, $h$ is measurable with respect to the $\sigma$-algebra generated by $W_{\beta,0}$ and $W_{\beta,1}$. Then, \eqref{lifting:eq-rep-p} follows directly from \ref{lifting:item-rep-a}, \ref{lifting:item-rep-d}, and the invariance of the Haar measure under right-multiplication.
\end{proof}

\section{Global well-posedness and invariance}\label{section:main}
In this section, we prove both Theorem \ref{intro:thm-rigorous-A-B} and Theorem \ref{intro:thm-rigorous-phi}, i.e., the global well-posedness and invariance of the Gibbs measure for the $\frkg$-valued derivatives $A,B\colon \R^{1+1}\rightarrow \frkg$ and for the $\frkG$-valued map $\phi\colon \R^{1+1} \rightarrow \frkG$. To this end, we first recall the most relevant formulations of the wave maps equation. In Cartesian coordinates, the wave maps equation with time cut-off function $\zeta \in C^\infty_b(\R)$ can be written as 
\begin{equation}\label{main:eq-WM}
\partial_t A^{(\zetascript)} = \partial_x B^{(\zetascript)} \qquad \text{and} \qquad \partial_t B^{(\zetascript)} =\partial_x A^{(\zetascript)} - \zeta(t) \big[ A^{(\zetascript)}, B^{(\zetascript)} \big].
\end{equation}
Our finite-dimensional approximation of \eqref{main:eq-WM}, which involves the Killing-renormalization, is given by
\begin{equation}\label{main:eq-WM-N}
\begin{aligned}
\partial_t A^{(\Nscript,\coup,\zetascript)} &= \partial_x B^{(\Nscript,\coup,\zetascript)}, \\ 
\partial_t B^{(\Nscript,\coup,\zetascript)} &= \partial_x A^{(\Nscript,\coup,\zetascript)} 
- \zeta(t) \big[ A^{(\Nscript,\coup,\zetascript)}, B^{(\Nscript,\coup,\zetascript)} \big]_{\leq N} + 2 \coup \zeta(t)^2  \Renorm[N] A^{(\Nscript,\coup,\zetascript)}.
\end{aligned}
\end{equation}
In null-coordinates, the counterpart of \eqref{main:eq-WM-N} is given by
\begin{equation}\label{main:eq-WM-null-N}
\partial_v U^{(\Nscript,\coup,\zetascript)} = \partial_u V^{(\Nscript,\coup,\zetascript)} 
= \zeta \big[ U^{(\Nscript,\coup,\zetascript)} , V^{(\Nscript,\coup,\zetascript)} \big] - \coup \zeta^2 \Renorm[N] \big( U^{(\Nscript,\coup,\zetascript)} + V^{(\Nscript,\coup,\zetascript)} \big).
\end{equation}
Whereas Section \ref{section:ansatz} and Sections \ref{section:Killing}-\ref{section:Lipschitz} primarily dealt with \eqref{main:eq-WM-null-N}, all three of the evolution equations \eqref{main:eq-WM}, \eqref{main:eq-WM-N}, and \eqref{main:eq-WM-null-N} will be important in this section. In Subsection \ref{section:main-lwp-remainder}, we prove the well-posedness of the remainder equations (Proposition \ref{main:prop-null-lwp}), 
which is the natural conclusion of our estimates from \mbox{Sections \ref{section:Killing}-\ref{section:Lipschitz}}. Due to our Ansatz from \eqref{ansatz:eq-UN-decomposition}-\eqref{ansatz:eq-VN-decomposition},  this yields detailed information on the random structure of the solution of \eqref{main:eq-WM-null-N}. In Subsection \ref{section:main-lwp}, we convert Proposition \ref{main:prop-null-lwp} into a well-posedness result for \eqref{main:eq-WM-N}, which is the subject of Proposition \ref{main:prop-lwp}. While Proposition \ref{main:prop-lwp} contains less detailed information on the random structure of the solution than Proposition \ref{main:prop-null-lwp}, it will be much easier to iterate in time, which is crucial for the rest of our argument. In Subsection \ref{section:main-invariance}, we prove the almost invariance of the Gibbs measure under the finite-dimensional approximation \eqref{main:eq-WM-N}. The main ingredients are the explicit formula for the Radon-Nikodym derivative from Proposition \ref{structure:prop-Gibbs} and our estimate of the energy-increment from \mbox{Section \ref{section:energy-increment}}. In Subsection \ref{section:main-gwp}, we use Bourgain's globalization argument to prove Theorem \ref{intro:thm-rigorous-A-B}, i.e., the global well-posedness and invariance of the Gibbs measure for \eqref{main:eq-WM}. Finally, in Subsection \ref{section:lifting-theorem}, we prove \mbox{Theorem \ref{intro:thm-rigorous-phi}}, i.e., the global well-posedness and invariance of the Gibbs measure for the $\frkG$-valued map $\phi$ from \eqref{intro:eq-WM-G}. \revision{For an overview of the main results of this section, we also refer the reader to Figure \ref{figure:structure-main-argument}.}\\

\tikzset{overviewrect/.style={rectangle, rounded corners, minimum width=2.5cm, minimum
height=1cm,text centered,align=center, draw=black, fill=white!10},
arrow/.style={thick,->,>=stealth}}
\tikzset{bigoverviewrect/.style={rectangle, rounded corners, minimum width=6cm, minimum
height=1cm,text centered,align=center, draw=black, fill=white!10},
arrow/.style={thick,->,>=stealth}}
\tikzset{medoverviewrect/.style={rectangle, rounded corners, minimum width=5cm, minimum
height=1cm,text centered,align=center, draw=black, fill=white!10},
arrow/.style={thick,->,>=stealth}}

\tikzstyle{innerWhite} = [semithick, white,line width=1.4pt, shorten >= 4.5pt]

\begin{figure}[t]
\begin{center}
\scalebox{0.8}{
\begin{tikzpicture}[node distance=2cm,scale=1, every node/.style={scale=1}]

\node (ansatz)[overviewrect]  at (0,-0.5) 
{Section \begin{NoHyper}\ref{section:ansatz}\end{NoHyper}};

\node (conservative)[overviewrect]  at (-4.25,-2) 
{Section \begin{NoHyper}\ref{section:conservative}\end{NoHyper}};

\node (body)[overviewrect]  at (4.25,-2) 
{Sections \begin{NoHyper}\ref{section:chaos}\,--\ref{section:lifting}\end{NoHyper}};

\begin{scope}[on background layer]
\draw [fill=blue!5] (8,-17.5) rectangle (-8,-3); 
\end{scope}
\node (LWPnullshort)[bigoverviewrect]  at (-4,-4.5) 
{\begin{NoHyper}Proposition \ref{main:prop-null-lwp}.\ref{main:item-null-lwp-short}\end{NoHyper}: \\ Short-time well-posedness \\ in null-coordinates};
\node (LWPnulllocal)[bigoverviewrect]  at (4,-4.5) 
{\begin{NoHyper}Proposition \ref{main:prop-null-lwp}.\ref{main:item-null-lwp-local}\end{NoHyper}: \\ Local well-posedness on $\BBAchi$  \\   in null-coordinates};
\node (LWPshort)[bigoverviewrect]  at (-4,-7) 
{\begin{NoHyper}Proposition \ref{main:prop-lwp}.\ref{main:item-lwp-short}\end{NoHyper}: \\ Short-time well-posedness};
\node (LWPlocal)[bigoverviewrect]  at (4,-7) 
{\begin{NoHyper}Proposition \ref{main:prop-lwp}.\ref{main:item-lwp-local}\end{NoHyper}: \\ Local well-posedness on $\BBAchi$    };

\node (almostinvariance)[bigoverviewrect]  at (-4,-9) 
{\begin{NoHyper}Proposition \ref{main:prop-almost-invariance}\end{NoHyper}: \\ Almost invariance};

\node (BBevent)[bigoverviewrect]  at (-4,-11) 
{\begin{NoHyper}Corollary \ref{main:cor-bourgain-bulut-probability}: \end{NoHyper} \\ Probability of $\BBAchi$};

\node (gwpsmall)[bigoverviewrect]  at (0,-13.5) 
{\begin{NoHyper}Lemma \ref{main:lem-gwp-small}: \end{NoHyper} \\ Almost global well-posedness and \\  invariance at low temperature};

\node (gwplarge)[bigoverviewrect]  at (0,-16) 
{\begin{NoHyper}Proposition \ref{main:prop-refined-gwp}: \end{NoHyper} \\ Refined global well-posedness and \\  invariance};


\draw[arrow] (-1.25,-0.5) -- (-4.25,-0.5)--(conservative);
\draw[arrow] (ansatz) -- (0,-3);
\draw[arrow] (body) -- (4.25,-3);
\draw[arrow] (conservative) -- (-9,-2) -- (-9,-9) -- (almostinvariance);

\draw[arrow] (LWPnullshort) -- (LWPshort); 
\draw[arrow] (LWPnulllocal) -- (LWPlocal); 
\draw[arrow] (LWPshort) -- (almostinvariance); 
\draw[arrow] (almostinvariance) -- (BBevent); 
\draw[arrow] (gwpsmall)--(gwplarge);

\draw[arrow] (BBevent) -- (-4,-13.5) -- (gwpsmall); 
\draw[arrow] (LWPlocal) -- (4,-13.5) -- (gwpsmall); 
\end{tikzpicture}
}
\end{center}
\caption{\small{\revision{This figure illustrates the main results of this section and their dependencies.}}}\label{figure:structure-main-argument}
\end{figure}

In Subsections \ref{section:main-lwp-remainder}-\ref{section:main-gwp} below, we state several lemmas and propositions which depend on each other. In many of them, we need to introduce large constants $C=C(\delta_\ast)$ and small constants $c=c(\delta_\ast)$. In order to organize the relative sizes of the constants, we introduce \revision{constants}
\begin{equation}\label{main:eq-constants}
\revision{(c_j)_{j=0}^{4} \qquad \text{and} \qquad (C_j)_{j=0}^{4}} 
\end{equation}
such that, for all $1\leq j \leq 4$, $C_j$ and $c_j$ are chosen as sufficiently large and sufficiently small depending on $(C_k)_{k=0}^{j-1}$, $(c_k)_{k=0}^{j-1}$, and the parameters from Section \ref{section:parameters}. 

\subsection{Well-posedness in null-coordinates}\label{section:main-lwp-remainder}

In this subsection, we prove the well-posedness of \eqref{main:eq-WM-null-N}. Due to our Ansatz from Section \ref{section:ansatz} and the well-posedness of the modulation equations (Proposition \ref{modulation:prop-main}), this boils down to the well-posedness of the remainder equations from Definition \ref{ansatz:def-remainder-equations}. 

\subsubsection{Bourgain-Bulut event and perturbations}

In the following definition, we introduce the Bourgain-Bulut event $\BBAzeta$. As the probability of the Bourgain-Bulut event $\BBAzeta$ is later controlled using the almost invariance of the Gibbs measure (Corollary \ref{main:cor-bourgain-bulut-probability}), it is best to describe the event in terms of Cartesian coordinates.

\begin{definition}[Event for Bourgain-Bulut estimate]\label{main:def-Bourgain-Bulut}
Let $N\geq 1$, let $R\geq 1$,  let $\coup>0$, and let $\zeta\in \C^\infty_b(\R)$.
For each $(W^{(\Rscript,\coup)}_{0},W^{(\Rscript,\coup)}_{1})\in (\C_x^{s-1}\times \C_x^{s-1})(\bT_R\rightarrow \frkg^2)$, let $(A^{(\Nscript,\Rscript,\coup,\zetascript)},B^{(\Nscript,\Rscript,\coup,\zetascript)})\colon \R \times \bT_R \rightarrow \frkg^2$ be the unique global solution of \eqref{main:eq-WM-N} with initial data $(W^{(\Rscript,\coup)}_{0},W^{(\Rscript,\coup)}_{1})$. Then, we define
\begin{align*}
&\BBAzeta  \\ 
:= \bigg\{& \big( W^{(\Rscript,\coup)}_{0},W^{(\Rscript,\coup)}_{1}\big)\in (\C_x^{s-1}\times \C_x^{s-1})(\bT_R\rightarrow \frkg^2) \colon \\ 
&\sup_{y\in \R} \, \langle Ny \rangle^{-4} \max\bigg( \Big\| 
\CHHLN_y \big( A^{(\Nscript,\Rscript,\coup,\zetascript)}, B^{(\Nscript,\Rscript,\coup,\zetascript)} \big)  \Big) 
\Big\|_{L_t^\infty L_x^\infty},\\ 
\hspace{1ex} &
 \Big\| 
\CHHLN_y \big( B^{(\Nscript,\Rscript,\coup,\zetascript)}, A^{(\Nscript,\Rscript,\coup,\zetascript)} \big)
 \Big\|_{L_t^\infty L_x^\infty} , \Big\| 
\CHHLN_y \big( A^{(\Nscript,\Rscript,\coup,\zetascript)}, A^{(\Nscript,\Rscript,\coup,\zetascript)} \big)
- 8 \coup \CNbd(y) \Kil \Big\|_{L_t^\infty L_x^\infty}, \\ 
\hspace{1ex} &\Big\| 
\CHHLN_y \big( B^{(\Nscript,\Rscript,\coup,\zetascript)}, B^{(\Nscript,\Rscript,\coup,\zetascript)} \big) 
- 8 \coup \CNbd(y) \Kil  \Big\|_{L_t^\infty L_x^\infty} \bigg)
\leq N^{-\frac{1}{2}+\delta} \coup \bigg\},
\end{align*}
where $\CNbd(y)$  is as in Definition \ref{jacobi:def-cf}, $\CHHLN_y$ is as in Definition \ref{jacobi:def-chhl}, and all $L_t^\infty L_x^\infty$-norms are taken over the region $[-4,4]\times \bT_R$.
\end{definition}

\begin{remark}
We emphasize that the definition of the Bourgain-Bulut event $\BBAzeta$ 
only relies on the solution $(A^{(\Nscript,\Rscript,\coup,\zetascript)},B^{(\Nscript,\Rscript,\coup,\zetascript)})$, which is defined in terms of
$(W^{(\Rscript,\coup)}_{0},W^{(\Rscript,\coup)}_{1})$. In particular, 
it relies neither on our Ansatz from \eqref{ansatz:eq-UN-rigorous-decomposition}-\eqref{ansatz:eq-VN-rigorous-decomposition} nor our representations of white noise from Lemma \ref{prelim:lem-white-noise-representation-2pi} and Lemma \ref{prelim:lem-white-noise-representation}.
\end{remark}

In the next definition, we introduce the class of perturbations of the initial data which will be considered in Proposition \ref{main:prop-null-lwp}.

\begin{definition}[Perturbations of initial modulation operators and nonlinear remainders]\label{main:def-perturbations-null} 
Let $N\in \dyadiclarge$, let $R \geq 1$, let $\coup>0$, let $\Ac\geq 1$,  and  let $\epsilon>0$. 
Then, we define 
$\Perpm(\Ac,\epsilon)$
as the set of all 
\begin{equation*}
\big(S^{(\Nscript,\Rscript,\coup),+}_K \big)_{K\in \dyadic},\big(S^{(\Nscript,\Rscript,\coup),-}_M \big)_{M\in \dyadic}\colon \bT_R \rightarrow \End(\frkg) \qquad \text{and} \qquad  
Z^{(\Nscript,\Rscript,\coup),+}, Z^{(\Nscript,\Rscript,\coup),-} \colon \bT_R \rightarrow \frkg
\end{equation*}
which satisfy the following conditions: 
\begin{enumerate}[label=(\roman*)]
    \item (Frequency support) For each $L\in \dyadic$, it holds that $P^v_{\gg L^{1-\delta}} S^{(\Nscript,\Rscript,\coup),\pm}_L=0$.
    \item (Frequency boundary) For all $L\in \dyadic$ satisfying $L>N^{1-\delta}$, it holds that $S^{(\Nscript,\Rscript,\coup),\pm}_L=\Id_\frkg$.
    \item (Distance to neutral element) It holds that
    \begin{equation*}
    \sup_{L\in \dyadic} \big\| S^{(\Nscript,\Rscript,\coup),\pm}_L  - \Id_\frkg \big\|_{\C_x^{s}}\leq \epsilon \qquad\text{and} \qquad 
    \big\| Z^{(\Nscript,\Rscript,\coup),\pm} \big\|_{\C_x^{r-1}} \leq \hcoup \Ac \epsilon.
    \end{equation*}
\end{enumerate}
\end{definition}

\begin{remark}
Similar as in Hypothesis \ref{hypothesis:pre} and Definition \ref{lipschitz:def-distances}, the estimates of the modulation operators and nonlinear remainders in Definition \ref{main:def-perturbations-null} differ by a factor of $\hcoup \Ac$.
\end{remark}

\subsubsection{Statement and proof of well-posedness in null-coordinates}

We now state and prove our well-posedness result for \eqref{main:eq-WM-null-N}. To this end, we \revision{recall that the parameters $\Ac$ and $\Bc$, which appear in the following definition, were previously motivated in Remark \ref{ansatz:rem-parameters} and Remark \ref{ansatz:rem-parameters-B}.} 

\begin{definition}[Well-posedness in null-coordinates]\label{main:def-null-lwp}
Let $N\in \Dyadiclarge$, let $R\geq 1$, let $\coup>0$,  let $\Ac,\Bc \geq 1$, let $\Dc=\hcoup \Ac \Bc$, 
and let $\zeta\in C^\infty_b(\R)$. Then, we define 
\begin{equation*}
\nullLWPDzeta \subseteq \Omega
\end{equation*}
as the event on which, under the assumptions \ref{main:item-null-lwp-A1}-\ref{main:item-null-lwp-A4}, the conclusions \ref{main:item-null-lwp-C1}-\ref{main:item-null-lwp-C4} hold. Here, the assumptions are defined 
as follows: 
\begin{enumerate}[label={[A\arabic*]}]
\item \label{main:item-null-lwp-A1} Let $\Nd \in \Dyadiclarge$ satisfy $\Nd \leq N^{1-\delta}$ and let $0<\epsilon \leq c_0$, \revision{where $c_0$ is as in \eqref{main:eq-constants}.}
\item \label{main:item-null-lwp-A2} Let the initial pure modulation operators and initial remainders
\begin{alignat}{3}
\big( S^{(\Nscript,\Rscript,\coup),\textup{in},+}_K \big)_{K\in \dyadic}, 
\big( S^{(\Nscript,\Rscript,\coup),\textup{in},-}_M \big)_{M\in \dyadic} &\colon \bT_R \rightarrow \End(\frkg), 
&\qquad 
Z^{(\Nscript,\Rscript,\coup),\pm}&\colon \bT_R \rightarrow \frkg, 
\label{main:eq-null-lwp-initial}\\ 
\big( \widetilde{S}^{(\Nscript,\Ndscript,\Rscript,\coup),\textup{in},+}_{K} \big)_{K\in \dyadic}, 
\big( \widetilde{S}^{(\Nscript,\Ndscript,\Rscript,\coup),\textup{in},-}_{M} \big)_{M\in \dyadic} &\colon \bT_R \rightarrow \End(\frkg), 
&\qquad 
\widetilde{Z}^{(\Nscript,\Ndscript,\Rscript,\coup),\pm} &\colon\bT_R \rightarrow \frkg, 
\label{main:eq-null-lwp-initial-tilde}
\end{alignat}
be such that 
\begin{equation*}
S^{(\Nscript,\Rscript,\coup),\textup{in},+}_K=S^{(\Nscript,\Rscript,\coup),\textup{in},-}_M=\Id_\frkg, 
\qquad Z^{(\Nscript,\Rscript,\coup),+}=Z^{(\Nscript,\Rscript,\coup),-}=0, 
\end{equation*}
and 
\begin{equation*}
\Big( \widetilde{S}^{(\Nscript,\Ndscript,\Rscript,\coup),\textup{in},+}, \widetilde{S}^{(\Nscript,\Ndscript,\Rscript,\coup),\textup{in},-}, 
\widetilde{Z}^{(\Nscript,\Ndscript,\Rscript,\coup),+}, \widetilde{Z}^{(\Nscript,\Ndscript,\Rscript,\coup),-} \Big) 
\in \Perpm(\Ac,\epsilon).
\end{equation*}
\item \label{main:item-null-lwp-A3} Let the initial data be given by 
\begin{align}
W^{(\Rscript,\coup),\pm} &= \hcoup \sum_{x_0\in \LambdaRR}\sum_{L\in \dyadic} \sum_{\ell \in \Z_L}\psiRx G_{x_0,\ell}^{\pm} e^{i\ell x},
\label{main:eq-null-lwp-initial-W} \\ 
\widetilde{W}^{(\Nscript,\Ndscript,\Rscript,\coup),\pm}&=   \sum_{L\in \dyadic}  \widetilde{S}^{(\Nscript,\Ndscript,\Rscript,\coup),\textup{in},\pm}_{L} P_{\leq \Nd}^x P^{\sharp}_{R;L} W^{(\Rscript,\coup),\pm} +\widetilde{Z}^{(\Nscript,\Ndscript,\Rscript,\coup),\pm}. 
\label{main:eq-null-lwp-initial-W-tilde}
\end{align}
\item \label{main:item-null-lwp-A4} Let $(U^{(\Nscript,\Rscript,\coup,\zetascript)},V^{(\Nscript,\Rscript,\coup,\zetascript)})$ and 
$(\widetilde{U}^{(\Nscript,\Ndscript,\Rscript,\coup,\zetascript)},\widetilde{V}^{(\Nscript,\Ndscript,\Rscript,\coup,\zetascript)})$ 
be the solutions of \eqref{main:eq-WM-null-N} with initial data from \eqref{main:eq-null-lwp-initial-W} and \eqref{main:eq-null-lwp-initial-W-tilde}, respectively.
\end{enumerate}
Furthermore, with $\epsilon^\prime:= C_0 (\epsilon+ \Bc \Nd^{-\varsigma})$ \revision{and $C_0$ as in \eqref{main:eq-constants}}, the conclusions \ref{main:item-null-lwp-C1}-\ref{main:item-null-lwp-C4} are defined as follows: 
\begin{enumerate}[label={[C\arabic*]}]
\item \label{main:item-null-lwp-C1} (Modulation and remainder equations) The modulation and remainder equations with initial data as in \eqref{main:eq-null-lwp-initial} or \eqref{main:eq-null-lwp-initial-tilde} have unique solutions
\begin{align}
&\Big( \pSNz[][+], \pSNz[][-], U^{(\Nscript,\Rscript,\coup,\zetascript),\fs}, V^{(\Nscript,\Rscript,\coup,\zetascript),\fs} \Big) 
\label{main:eq-null-lwp-S-rem-N} \\
\text{and} \qquad &\Big( \pStilz[][+], \pStilz[][-], \widetilde{U}^{(\Nscript,\Ndscript,\Rscript,\coup,\zetascript),\fs}, \widetilde{V}^{(\Nscript,\Ndscript,\Rscript,\coup,\zetascript),\fs}\Big)
\label{main:eq-null-lwp-S-rem-til} 
\end{align}
which satisfy the estimates
\begin{align*}
\sup_{K\in \Dyadiclarge} \big\| \pSNz[K][+] \big\|_{\Wuv[s][s]},  \sup_{K\in \Dyadiclarge} \big\| \pStilz[K][+] \big\|_{\Wuv[s][s]} 
&\leq \Bc, \\ 
\sup_{M\in \Dyadiclarge}\big\| \pSNz[M][-] \big\|_{\Wuv[s][s]}, \sup_{M\in \Dyadiclarge}\big\| \pStilz[M][-] \big\|_{\Wuv[s][s]} &\leq \Bc, \\ 
\big\|  U^{(\Nscript,\Rscript,\coup,\zetascript),\fs} \big\|_{\Cprod{r-1}{r}} 
, \big\| \widetilde{U}^{(\Nscript,\Ndscript,\Rscript,\coup,\zetascript),\fs} \big\|_{\Cprod{r-1}{r}}, 
\big\| V^{(\Nscript,\Rscript,\coup,\zetascript),\fs} \big\|_{\Cprod{r}{r-1}}, 
\big\| \widetilde{V}^{(\Nscript,\Ndscript,\Rscript,\coup,\zetascript),\fs} \big\|_{\Cprod{r}{r-1}} 
&\leq \Dc. 
\end{align*}
\item \label{main:item-null-lwp-C2} (Lipschitz estimates) The solutions \eqref{main:eq-null-lwp-S-rem-N} and \eqref{main:eq-null-lwp-S-rem-til} of the modulation and remainder equations satisfy the Lipschitz estimates 
\begin{align*}
\sup_{K\in \dyadiclarge} \big\| \pSNz[K][+]-  \pStilz[K][+] \big\|_{\Wuv[s][s]} &\leq \epsilon^\prime, \\  
\sup_{M\in \dyadiclarge}\big\| \pSNz[M][-] -  \pStilz[M][-] \big\|_{\Wuv[s][s]} &\leq  \epsilon^\prime, \\ 
\big\|  U^{(\Nscript,\Rscript,\coup,\zetascript),\fs} -  
\widetilde{U}^{(\Nscript,\Ndscript,\Rscript,\coup,\zetascript),\fs} \big\|_{\Cprod{r-1}{r}}
&\leq  \hcoup \Ac  \epsilon^\prime, \\ 
\big\| V^{(\Nscript,\Rscript,\coup,\zetascript),\fs} - 
\widetilde{V}^{(\Nscript,\Ndscript,\Rscript,\coup,\zetascript),\fs} \big\|_{\Cprod{r}{r-1}} 
 &\leq \hcoup \Ac \epsilon^\prime.
\end{align*}
\item \label{main:item-null-lwp-C3} (Parameter trick) The nonlinear remainders $U^{(\Nscript,\Rscript,\coup,\zetascript),\fs}$
and $V^{(\Nscript,\Rscript,\coup,\zetascript),\fs}$ satisfy the higher-regularity estimates 
\begin{equation*}
\big\|  U^{(\Nscript,\Rscript,\coup,\zetascript),\fs}\big\|_{\Cprod{r-1+\eta}{r+\eta}}, 
\big\|  V^{(\Nscript,\Rscript,\coup,\zetascript),\fs} \big\|_{\Cprod{r+\eta}{r-1+\eta}} 
\leq \Dc. 
\end{equation*}
\item \label{main:item-null-lwp-C4} (Structure of full solution) 
The solutions $(U^{(\Nscript,\Rscript,\coup,\zetascript)},V^{(\Nscript,\Rscript,\coup,\zetascript)})$ and $(\widetilde{U}^{(\Nscript,\Ndscript,\Rscript,\coup,\zetascript)},\widetilde{V}^{(\Nscript,\Ndscript,\Rscript,\coup,\zetascript)})$ from \ref{main:item-null-lwp-A4} can be decomposed as in our Ansatz \eqref{ansatz:eq-UN-rigorous-decomposition}-\eqref{ansatz:eq-VN-rigorous-decomposition}, where the modulated and mixed modulated objects are defined using the pure modulation operators and nonlinear remainders from \eqref{main:eq-null-lwp-S-rem-N} and \eqref{main:eq-null-lwp-S-rem-til}, respectively. 
\end{enumerate}
Finally, we define
\begin{equation}\label{main:eq-lwpd-null-cut}
\nullLWPDcut := \bigcap_{\chi \in \Cut} \nullLWPDchi.
\end{equation}
\end{definition}

\begin{remark}\label{main:rem-null-lwp}
\revision{To avoid confusion, let us also put Definition \ref{main:def-null-lwp} into slightly different words: We are in the event $\nullLWPDzeta$ if and only if, given any additional inputs satisfying the assumptions \ref{main:item-null-lwp-A1}-\ref{main:item-null-lwp-A4}, the corresponding solutions must satisfy the conclusions in \ref{main:item-null-lwp-C1}-\ref{main:item-null-lwp-C4}. In particular, one should think of the event $\nullLWPDzeta$ as imposing conditions on the variables $(G_{x_0,\ell}^{\pm})_{x_0\in\LambdaRR,\ell \in \Z}$ appearing in \eqref{main:eq-null-lwp-initial-W}.
Since the Assumption \ref{main:item-null-lwp-A2} and the intersection in \eqref{main:eq-lwpd-null-cut} involve uncountable sets, one may worry about the measurability of the event $\nullLWPDzeta$}. However, all the conditions can be rephrased via an approximation argument using only countable sets, and therefore the measurability is guaranteed.
\end{remark}

Equipped with Definition \ref{main:prop-null-lwp}, we can now state and prove our well-posedness result in null-coordinates. 

\begin{proposition}[Well-posedness in null-coordinates]\label{main:prop-null-lwp}
Let $N\in \dyadiclarge$, let $R\geq 1$, let $\coup>0$, let $\Ac,\Bc \geq C_0$, and let $\Dc=\hcoup \Ac \Bc$. Then, we have the following two estimates:
\begin{enumerate}[label=(\Roman*)]
\item \label{main:item-null-lwp-short} 
If $\Dc \leq c_0 N^{-\delta_2}$, then
\begin{equation*}
\mathbb{P} \Big(  \nullLWPDcut \Big) \geq 1- c_0^{-1} \exp\Big( - c_0 R^{-2\eta} \Ac^2 \Big). 
\end{equation*}
\item \label{main:item-null-lwp-local} If $\Dc \leq c_0$ and $\chi \in \Cut$, then it holds that  
\begin{equation}\label{main:eq-null-lwp-probability}
\mathbb{P} \Big(  \nullLWPDchi \Big) \geq \muR \Big( \BBAchi \Big) - c_0^{-1} \exp\Big( - c_0 R^{-2\eta} \Ac^2 \Big).
\end{equation}
\end{enumerate}
\end{proposition}

\begin{remark}
We recall that, due to the scaling-symmetry proven in Lemma \ref{ansatz:lem-scaling-symmetry}, the temperature $\coup$ can be interpreted as a timescale. Due to this, the conditions in \ref{main:item-lwp-short} and \ref{main:item-lwp-local} can be interpreted as restrictions to the short timescale $\sim N^{-2\delta_2}$ or the local timescale $\sim 1$, respectively. 
In Subsection \ref{section:main-invariance}, we will use \ref{main:item-lwp-short} to control the probability of $\BBAchi$, and therefore use \ref{main:item-lwp-short} to show that the lower-bound in \eqref{main:eq-null-lwp-probability} from \ref{main:item-lwp-local} is useful. 
\end{remark}

\begin{proof}[Proof of Proposition \ref{main:prop-null-lwp}:] 
All the important estimates and ideas needed for the proof of this proposition are included in Sections \ref{section:ansatz}-\ref{section:Lipschitz}. Nevertheless, due to the length of this article, combining them into a proof requires a certain amount of effort. In the first and second step of this proof, we separately analyze the solutions
$(\widetilde{U}^{(\Nscript,\Ndscript,\Rscript,\coup,\chiscript)},\widetilde{V}^{(\Nscript,\Ndscript,\Rscript,\coup,\chiscript)})$
and $(U^{(\Nscript,\Rscript,\coup,\chiscript)},V^{(\Nscript,\Rscript,\coup,\chiscript)})$, respectively. In particular, we obtain the bounds from \ref{main:item-null-lwp-C1} and the decompositions from \ref{main:item-null-lwp-C4} in Definition \ref{main:def-null-lwp}. In the third step, we consider the difference of $(\widetilde{U}^{(\Nscript,\Ndscript,\Rscript,\coup,\chiscript)},\widetilde{V}^{(\Nscript,\Ndscript,\Rscript,\coup,\chiscript)})$
and $(U^{(\Nscript,\Rscript,\coup,\chiscript)},V^{(\Nscript,\Rscript,\coup,\chiscript)})$ and obtain the Lipschitz-estimates from \ref{main:item-null-lwp-C2}. In the fourth step, we then describe a parameter trick, which allows us to prove the estimate in \ref{main:item-null-lwp-C3}. \\

\emph{Step 1: Solutions of the wave maps equation \eqref{main:eq-WM-null-N} with frequency-truncated nonlinearity and initial data.} 
In this step, we analyze the solution $(\widetilde{U}^{(\Nscript,\Ndscript,\Rscript,\coup,\chiscript)},\widetilde{V}^{(\Nscript,\Ndscript,\Rscript,\coup,\chiscript)})$ of \eqref{main:eq-WM-null-N} with the frequency-truncated initial data from \eqref{main:eq-null-lwp-initial-W-tilde}. For expository purposes, we split the argument into six sub-steps.\\

\emph{Step 1.(a): Notation.}
Let $\Nd\in \Dyadiclarge$ satisfy $\Nd \leq N^{1-\delta}$ and let $\chi \in \Cut$.
To simplify the notation, we write 
\begin{equation}\label{main:eq-null-pa-1}
\begin{aligned}
&\, \Big( S^{\fd,\diamond,+}, S^{\fd,\diamond,-} , U^{\fd,\fs} , V^{\fd,\fs} \Big) \\
:=&\, \Big(   
\widetilde{S}^{(\Nscript,\Ndscript,\Rscript,\coup,\chiscript),\diamond,+}, 
\widetilde{S}^{(\Nscript,\Ndscript,\Rscript,\coup,\chiscript),\diamond,-},  
\widetilde{U}^{(\Nscript,\Ndscript,\Rscript,\coup,\chiscript),\fs}, 
\widetilde{V}^{(\Nscript,\Ndscript,\Rscript,\coup,\chiscript),\fs}\Big). 
\end{aligned}
\end{equation}
In other words, we use the super-script ``$\fd$" to represent the objects from Definition \ref{main:def-null-lwp} corresponding to wave maps with frequency-truncated initial data (and nonlinearity). The modulated and mixed modulated objects $U^{\fd,\ast}$ and $V^{\fd,\ast}$, where 
$\ast \in \{+,-,+-,+\fs,\fs-\}$, are defined using \eqref{main:eq-null-pa-1} as in Definitions \ref{ansatz:def-pure}, \ref{ansatz:def-modulated-linear}, \ref{ansatz:def-modulated-bilinear}, \ref{ansatz:def-mixed}, and \ref{ansatz:def-modulated-linear-reversed}. 
We also define
\begin{equation}\label{main:eq-null-pa-2}
U^{\fd} = \sum_{\ast} U^{\fd,\ast} \qquad \text{and} \qquad V^{\fd} = \sum_{\ast} V^{\fd,\ast},
\end{equation}
where the sums are taken over all $\ast \in \{ +, - , +-, +\fs, \fs- , \fs \}$.
A similar notation as for $U^{\fd,\ast}$ and $V^{\fd,\ast}$ is used for the error terms from \mbox{Proposition \ref{ansatz:prop-decomposition}}.\\

\emph{Step 1.(b): Setting.}
For fixed nonlinear remainders $U^{\fd,\fs}$ and $V^{\fd,\fs}$, we consider the modulation equations 
\begin{equation}\label{main:eq-null-pb-2}
\begin{cases}
\begin{aligned}
\partial_v S^{\fd,\diamond,+}_{K,k}  
&=  \chinull[K] \rho_{\leq N}^2(k)   \Big[  S^{\fd,\diamond,+}_{K,k} , \LOd[K][-] \Big] + \chinull[K]  \rho_{\leq N}^4(k)  \Big[  S^{\fd,\diamond,+}_{K,k} , \SHHLd[K][v] \Big],   \\
\partial_u S^{\fd,\diamond,-}_{M,m}  
&=  \chinull[M]    \rho_{\leq N}^2(m)\Big[  \LOd[M][+] , S^{\fd,\diamond,-}_{M,m}  \Big] +  \chinull[M] \rho_{\leq N}^4(m)  \Big[ \SHHLd[M][u], S^{\fd,\diamond,-}_{M,m}  \Big], \\ 
S^{\fd,\diamond,+}_{K,k}  \Big|_{v=u} &= \widetilde{S}^{(\Nscript,\Ndscript,\Rscript,\coup),\textup{in},+}_{K}, 
\qquad S^{\fd,\diamond,-}_{M,m}  \Big|_{u=v} = \widetilde{S}^{(\Nscript,\Ndscript,\Rscript,\coup),\textup{in},-}_{M}
\end{aligned}
\end{cases}
\end{equation}
from Definition \ref{ansatz:def-modulation-equations}, where the initial pure modulation operators 
are as in Definition \ref{main:def-null-lwp}.\ref{main:item-null-lwp-A2}. We also consider the remainder equations 
\begin{equation}\label{main:eq-null-pb-3}
\begin{cases}
\begin{aligned}
\partial_v U^{\fd,\fs}  &= \chinull  \big[ U^{\fd}, V^{\fd} \big]_{\leq N} - \coup \chinullsquare \Renormd \big( U^{\fd} + V^{\fd} \big) \\ 
&\hspace{2.5ex}- \partial_v \big( U^{\fd,+} + U^{\fd,+-} + U^{\fd,-} + U^{\fd,+\fs} + U^{\fd,\fs-} \big) ,\\ 
\partial_u V^{\fd,\fs}  &= \chinull  \big[ U^{\fd}, V^{\fd} \big]_{\leq N} - \coup \chinullsquare \Renormd \big( U^{\fd} + V^{\fd} \big) \\ 
&\hspace{2.5ex}- \partial_u \big( V^{\fd,+} + V^{\fd,+-} + V^{\fd,-} + V^{\fd,+\fs} + V^{\fd,\fs-} \big) ,\\ 
U^{\fd,\fs} \big|_{v=u} 
&= \widetilde{W}^{(\Nscript,\Ndscript,\Rscript,\coup),+}(u)  
- \big( U^{\fd,+} + U^{\fd,+-} + U^{\fd,-} + U^{\fd,+\fs} + U^{\fd,\fs-} \big)\big|_{v=u}, \\  
V^{\fd,\fs} \big|_{u=v} 
&= \widetilde{W}^{(\Nscript,\Ndscript,\Rscript,\coup),-}(v)  
- \big( V^{\fd,-} + V^{\fd,+-} + V^{\fd,+} + V^{\fd,\fs-} + V^{\fd,+\fs} \big)\big|_{u=v}
\end{aligned}
\end{cases}
\end{equation}
from Definition \ref{ansatz:def-remainder-equations}. 
In keeping with \eqref{main:eq-null-pa-1}, we wrote $\Renormd$ for $\Renorm[\Nscript,\Ndscript]$ from Definition \ref{ansatz:def-Killing}.
For the rest of this step, we assume the smallness condition $\Dc \leq c_0$. We then show that, on the event $\PHA$, the modulation equations \eqref{main:eq-null-pb-2} and remainder equations \eqref{main:eq-null-pb-3} have unique solutions which satisfy the estimates in Definition \ref{main:def-null-lwp}.\ref{main:item-null-lwp-C1}. Furthermore, we show that $(U^\fd,V^\fd)$ from \eqref{main:eq-null-pa-2} solves the wave maps equation \eqref{main:eq-WM-null-N} with the initial data from \eqref{main:eq-null-lwp-initial-W-tilde}.\\

\emph{Step 1.(c): The maps $\Uscr^{\fd,\fs}$ and $\Vscr^{\fd,\fs}$.}
We first introduce the ball 
\begin{equation}\label{main:eq-null-pc-4}
\mathbb{B}^{\fd,\fs} := \Big\{ \big( U^{\fd,\fs}, V^{\fd,\fs} \big) \in (\Cprod{r-1}{r} \times \Cprod{r}{r-1})(\R^{1+1})\colon 
\big\| U^{\fd,\fs} \big\|_{\Cprod{r-1}{r}}, \big\| V^{\fd,\fs} \big\|_{\Cprod{r}{r-1}} 
\leq \Dc \Big\}.
\end{equation}
For any pair $(U^{\fd,\fs},V^{\fd,\fs})\in \mathbb{B}^{\fd,\fs}$, Proposition \ref{modulation:prop-main} implies the existence of a unique solution 
\begin{equation*}
\big( S^{\fd,\diamond,+}, S^{\fd,\diamond,-} \big)
= \big( S^{\fd,\diamond,+}, S^{\fd,\diamond,-} \big)\Big( \widetilde{S}^{(\Nscript,\Ndscript,\Rscript,\coup),\textup{in},\pm}, U^{\fd,\fs}, V^{\fd,\fs} \Big) 
\end{equation*}
of  the modulation equations \eqref{main:eq-null-pb-2}. It then only remains to solve the remainder equations \eqref{main:eq-null-pb-3}. To this end, we first recall from Proposition \ref{ansatz:prop-decomposition} that the right-hand sides of \eqref{main:eq-null-pb-3} can be written 
as
\begin{equation*}
\chinull \HHLERRd + \chinull \JcbERRd + \sigma_u \SEd{u} + \sigma_v \SEd{v} + \chinull \PId + \chinull \RenNErrd, 
\end{equation*}
where $(\sigma_u,\sigma_v)\in \{ (1,0),(0,1)\}$. We then define the map $\Uscr^{\fd,\fs}$ by 
\begin{equation}\label{main:eq-null-pc-2}
\begin{aligned}
\Uscr^{\fd,\fs}\big( U^{\fd,\fs}, V^{\fd,\fs} \big) 
&=  \widetilde{W}^{(\Nscript,\Ndscript,\Rscript,\coup),+}(u) 
- \big( U^{\fd,+} + U^{\fd,+-} + U^{\fd,-} + U^{\fd,+\fs} + U^{\fd,\fs-} \big)(u,u) \\ 
&+ \Int^v_{u\rightarrow v} \Big( \chinull \HHLERRd + \chinull \JcbERRd + \SEd{u} + \chinull \PId  + \chinull \RenNErrd \Big)(u,v).
\end{aligned}
\end{equation}
Similarly, we define the map $\Vscr^{\fd,\fs}$ by 
\begin{equation}\label{main:eq-null-pc-3}
\begin{aligned}
\Vscr^{\fd,\fs}\big( U^{\fd,\fs}, V^{\fd,\fs} \big) 
&= \widetilde{W}^{(\Nscript,\Ndscript,\Rscript,\coup),-}(v)  
- \big( V^{\fd,-} + V^{\fd,+-} + V^{\fd,+} + V^{\fd,\fs-} + V^{\fd,+\fs}  \big)(v,v) \\ 
&+ \Int^u_{v\rightarrow u} \Big( \chinull \HHLERRd + \chinull \JcbERRd + \SEd{v} + \chinull \PId + \chinull \RenNErrd \Big)(u,v).
\end{aligned}
\end{equation}
In order to solve the remainder equations \eqref{main:eq-null-pb-3}, it then suffices to show that $(\Uscr^{\fd,\fs},\Vscr^{\fd,\fs})$ is a contraction on the ball $\mathbb{B}^{\fd,\fs}$.\\

\emph{Step 1.(d): Self-mapping.} 
We first show that $(\Uscr^{\fd,\fs},\Vscr^{\fd,\fs})$ maps $\mathbb{B}^{\fd,\fs}$ into $\mathbb{B}^{\fd,\fs}$. Due to the symmetry of our estimates in the $u$ and $v$-variables, it suffices to treat $\Uscr^{\fd,\fs}$. To estimate the contribution of the initial data, we first use Definition \ref{ansatz:def-modulated-linear} and \eqref{main:eq-null-lwp-initial-W-tilde}, which imply that 
\begin{align} 
&\, \widetilde{W}^{(\Nscript,\Ndscript,\Rscript,\coup),+}(u) 
- \big( U^{\fd,+} + U^{\fd,+-} + U^{\fd,-} + U^{\fd,+\fs} + U^{\fd,\fs-} \big)(u,u) \notag \\ 
=&\,   \sum_{K\in \dyadic} 
\widetilde{S}^{(\Nscript,\Ndscript,\Rscript,\coup),\textup{in},+}_{K}(u) \big( P_{\leq \Nd} P^{\sharp}_{R;K} W^{(\Rscript,\coup),+}\big)(u) - U^{\fd,+}(u,u) \label{main:eq-null-qb-1} \\ 
+&\, \widetilde{Z}^{(\Nscript,\Ndscript,\Rscript,\coup),+}(u) \label{main:eq-null-qb-2} \\ 
-&\, \big( U^{\fd,+-} + U^{\fd,-} + U^{\fd,+\fs} + U^{\fd,\fs-} \big)(u,u). \label{main:eq-null-qb-3}
\end{align}
We estimate the two terms \eqref{main:eq-null-qb-1} and \eqref{main:eq-null-qb-2} using Assumption \ref{main:item-null-lwp-A2} and Lemma \ref{modulation:lem-initial-modulated-linear}, which yield that 
\begin{equation}\label{main:eq-null-qb-4}
\begin{aligned}
&\, \big\| \eqref{main:eq-null-qb-1} \big\|_{\C_u^{r-1}}
+ \big\| \eqref{main:eq-null-qb-2} \big\|_{\C_u^{r-1}} \\
\lesssim&\, \hcoup \Ac  \, 
\Big( \sup_{K\in \dyadic} \big\| \SNin[K][+] \big\|_{\C_u^s} + \sup_{M\in \dyadic}\big\| \SNin[M][-] \big\|_{\C_v^s} \Big) + \hcoup \Ac \epsilon
\lesssim  \hcoup \Ac \big(1+\epsilon\big). 
\end{aligned}
\end{equation}
Due to Assumption \ref{main:item-null-lwp-A1} and our assumption $\Bc\geq C_1$, we further obtain that
\begin{equation}\label{main:eq-null-qb-4prime}
    \hcoup \Ac \big(1+\epsilon\big) \lesssim \hcoup \Ac \lesssim C_0^{-1} \hcoup \Ac \Bc = C_0^{-1} \Dc. 
\end{equation}

The last term \eqref{main:eq-null-qb-3} is estimated using\footnote{In order to use Lemma \ref{modulation:lem-Cartesian}, we note that $\| U^{\fd,\ast}(u,u) \|_{\C_u^{r-1}} \leq \| U^{\fd,\ast}(x-t,x+t) \|_{C_t^0 \C_x^{r-1}}$.} Lemma \ref{modulation:lem-Cartesian}, which yields that
\begin{equation}\label{main:eq-null-qb-4primeprime}
\begin{aligned}
\big\| \eqref{main:eq-null-qb-3} \big\|_{\C_u^{r-1}} \lesssim \Dc^2. 
\end{aligned}
\end{equation}

We now estimate the contribution of the nonlinearity. To control the $\HHLERRd$, $\PId$, and $\RenNErrd$-terms, we combine Lemma \ref{prelim:lem-Duhamel-integral}, Proposition \ref{hhl:prop-main}, Proposition \ref{null:prop-perturbative}, and Proposition \ref{ren:prop-main}, which yield that 
\begin{equation}\label{main:eq-null-pd-1}
\begin{aligned}
&\, \Big\| \Int^v_{u\rightarrow v}\Big( \chinull \HHLERRd + \chinull \PId + \chinull \RenNErrd  \Big) \Big\|_{\Cprod{r-1}{r}}\\ 
\lesssim&\,   \big\| \HHLERRd \big\|_{\Cprod{r-1}{r-1}} + \big\|\, \PId \big\|_{\Cprod{r-1}{r-1}} 
+ \big\| \RenNErrd \big\|_{\Cprod{r-1}{r-1}}
\lesssim \Dc^2. 
\end{aligned}
\end{equation}
In order to control the $\SEd{u}$-term, we use Lemma \ref{prelim:lem-weighted-hoelder-properties} and Proposition \ref{structural:prop-main}, which yield that 
\begin{equation}\label{main:eq-null-pd-2}
\Big\| \Int^v_{u\rightarrow v} \Big( \SEd{u} \Big) \Big\|_{\Cprod{r-1}{r}}
\lesssim \big\| \SEd{u} \big\|_{\WCprod{r-1}{r-1}} \lesssim \Dc^2. 
\end{equation}
Finally, it remains to control the Jacobi error $\JcbERRd$. Using Lemma \ref{prelim:lem-Duhamel-integral} and \eqref{jacobi:eq-main-estimate-1} from Proposition \ref{jacobi:prop-main}, we obtain that 
\begin{equation}\label{main:eq-null-pd-3}
\Big\| \Int^v_{u\rightarrow v}\Big( \chinull \JcbERRd \Big) \Big\|_{\Cprod{r-1}{r}}
\lesssim  \big\| \JcbERRd\big\|_{\Cprod{r-1}{r-1}} 
\lesssim \Big( 1+ \frac{\Nd}{N} N^{\delta_2} \Big) \Dc^3. 
\end{equation}
Due to the condition $\Nd\leq N^{1-\delta}$ from Assumption \ref{main:item-null-lwp-A1}, the right-hand side of \eqref{main:eq-null-pd-4} can be further estimated by 
\begin{equation}\label{main:eq-null-pd-4}
\Big( 1+ \frac{\Nd}{N} N^{\delta_2} \Big) \Dc^3
\lesssim \Big( 1 + N^{-\delta} N^{\delta_2}\Big) \Dc^3 \lesssim \Dc^3. 
\end{equation}

By combining \eqref{main:eq-null-qb-4}-\eqref{main:eq-null-pd-4}, it follows that
\begin{equation}\label{main:eq-null-pd-5}
\big\| \Uscr^{\fd,\fs} \big( U^{\fd,\fs}, V^{\fd,\fs} \big) \big\|_{\Cprod{r-1}{r}} 
\lesssim C_0^{-1} \Dc + \Dc^2.
\end{equation}
Since $C_0=C_0(\delta_\ast)$ is sufficiently large, $c_0=c_0(\delta_\ast)$ is sufficiently small,  and $\Dc \leq c_0$,
\eqref{main:eq-null-pd-4} proves that $(\Uscr^{\fd,\fs},\Vscr^{\fd,\fs})$ maps $\mathbb{B}^{\fd,\fs}$ into $\mathbb{B}^{\fd,\fs}$. \\

\emph{Step 1.(e): Contraction.} It remains to prove that $(\Uscr^{\fd,\fs},\Vscr^{\fd,\fs})$ not only maps $\mathbb{B}^{\fd,\fs}$ into $\mathbb{B}^{\fd,\fs}$, but is also a contraction on $\mathbb{B}^{\fd,\fs}$. Due to the symmetry of our estimates in the $u$ and $v$-variables, it suffices to prove for all $(U^{\fd,\fs},V^{\fd,\fs}),(\widetilde{U}^{\fd,\fs},\widetilde{V}^{\fd,\fs})\in \mathbb{B}^{\fd,\fs}$ that
\begin{equation}\label{main:eq-null-pe-1}
\begin{aligned}
 \Big\| \Uscr^{\fd,\fs}\big(U^{\fd,\fs}, V^{\fd,\fs}\big) - \Uscr^{\fd,\fs} \big( \widetilde{U}^{\fd,\fs}, \widetilde{V}^{\fd,\fs} \big) \Big\|_{\Cprod{r-1}{r}} 
\lesssim \Dc  \Big( 
\big\| U^{\fd,\fs} - \widetilde{U}^{\fd,\fs}\big\|_{\Cprod{r-1}{r}}
+ \big\| V^{\fd,\fs} - \widetilde{V}^{\fd,\fs}\big\|_{\Cprod{r}{r-1}} \Big).
\end{aligned}
\end{equation}
This follows from a similar argument as for \eqref{main:eq-null-pd-5}, but with the estimates from Sections \ref{section:modulated-mixed}-\ref{section:renormalization} replaced by their Lipschitz-variants from Lemma \ref{lipschitz:lem-modulated-cartesian}, Lemma \ref{lipschitz:lem-modulated-initial}, and Proposition \ref{lipschitz:prop-remainder-no-jcb}.\\

\emph{Step 1.(f): Conclusion.} In Step 1.(d) and Step 1.(e), we proved that the map $(\Uscr^{\fd,\fs},\Vscr^{\fd,\fs})$ is a contraction on the ball $\mathbb{B}^{\fd,\fs}$ from \eqref{main:eq-null-pc-4}. Due to the contraction-mapping principle, it follows that there exists a unique solution $(U^{\fd,\fs},V^{\fd,\fs})$ of the remainder equations \eqref{main:eq-null-pb-3} in $\mathbb{B}^{\fd,\fs}$. Due to the definition of $\mathbb{B}^{\fd,\fs}$ and  \mbox{Proposition \ref{modulation:prop-main}},  it follows that 
\begin{equation*}
\sup_{K\in \Dyadiclarge} \big\| S^{\fd,\diamond,+}_K \big\|_{\Wuv[s][s]}\leq \Bc, \quad  \sup_{M\in \Dyadiclarge} \big\| S^{\fd,\diamond,-}_M \big\|_{\Wuv[s][s]} \leq \Bc, \quad  \big\| U^{\fd,\fs} \big\|_{\Cprod{r-1}{r}}\leq \Dc, \quad \text{and} \quad \big\| V^{\fd,\fs} \big\|_{\Cprod{r}{r-1}} \leq \Dc.
\end{equation*}
Thus, $(S^{\fd,\diamond,+},S^{\fd,\diamond,-},U^{\fd,\fs},V^{\fd,\fs})$ satisfy the desired bounds from \ref{main:item-null-lwp-C1}. As the last part of this step, we verify the part of \ref{main:item-null-lwp-C4} corresponding to $U^{\fd,\fs}$ and $V^{\fd,\fs}$.
Due to soft arguments, \eqref{main:eq-WM-null-N} has a unique global solution (see e.g. Lemma \ref{structure:lem-global-flow}),
and it therefore suffices to verify\footnote{While this is almost trivial for $U^{\fd}$ and $V^{\fd}$, it is less obvious for $\widebar{U}^{\fn}$ and $\widebar{V}^{\fn}$ from the second step.} that $U^{\fd}$ and $V^{\fd}$ from \eqref{main:eq-null-pa-2} solve the frequency-truncated wave maps equation \eqref{main:eq-WM-null-N}. Using \eqref{main:eq-null-pa-2} and  \eqref{main:eq-null-pb-3}, it holds that 
\begin{align*}
\partial_v U^{\fd} 
&= \partial_v \big( U^{\fd,+} + U^{\fd,+-} + U^{\fd,-} + U^{\fd,+\fs} + U^{\fd,\fs-} \big) + \partial_v U^{\fd,\fs} \\
&= \chinull  \big[ U^{\fd}, V^{\fd} \big]_{\leq N} - \coup \chinullsquare \Renormd \big( U^{\fd} + V^{\fd} \big).
\end{align*}
Together with a similar calculation for $V^{\fd}$, we obtain that $U^{\fd}$ and $V^{\fd}$ solve \eqref{main:eq-WM-null-N}.\\

\emph{Step 2: Solution of the wave maps equation \eqref{main:eq-WM-null-N} with a frequency-truncated nonlinearity.}
In the second step, we analyze the solution $(U^{(\Nscript,\Rscript,\coup,\chiscript)},V^{(\Nscript,\Rscript,\coup,\chiscript)})$
of the frequency-truncated wave maps equation \eqref{main:eq-WM-null-N} with the initial data from \eqref{main:eq-null-lwp-initial-W}, i.e., 
the initial data without frequency-truncations. For expository purposes, we split the argument into five sub-steps. \\

\emph{Step 2.(a): The case $\Dc\leq c_0 N^{-\delta_2}$.}
This case can be treated almost exactly as in the first step. The only difference is that the estimate from 
\eqref{main:eq-null-pd-4} is replaced by\footnote{We recall from Remark \ref{ansatz:rem-frequency-truncation} that we can always reduce to $\Nd\lesssim N$.}
\begin{equation}\label{main:eq-null-2-pa-1}
\Big( 1 + \frac{\Nd}{N}  N^{\delta_2} \Big) \Dc^3 \lesssim \big( 1+ N^{\delta_2}\big) \Dc^3 \lesssim \Dc^2,
\end{equation}
where we used the condition $\Dc\leq c_0 N^{-\delta_2}$. For the rest of the second step, we only assume $\Dc \leq c_0$, but restrict ourselves to the set 
\begin{equation}\label{main:eq-null-2-pa-2}
\Big\{ \big( W^{(\Rscript,\coup),+}, W^{(\Rscript,\coup),-} \big) \colon 
\big( W^{(\Rscript,\coup)}_0, W^{(\Rscript,\coup)}_1 \big) \in \BBAchi 
\Big\} 
\bigcap \PHA,
\end{equation}
where $(W^{(\Rscript,\coup)}_0, W^{(\Rscript,\coup)}_1)$ is as in \eqref{ansatz:eq-null-data} and $\BBAchi$ is as in Definition \ref{main:def-Bourgain-Bulut}. \\

\emph{Step 2.(b): Notation and setting.} We recall from Subsection \ref{section:overview-jacobi} that the remainder equations corresponding to $(U^{(\Nscript,R,\coup,\chiscript)},V^{(\Nscript,R,\coup,\chiscript)})$ cannot be solved using only a contraction-mapping argument, and instead also require a Bourgain-Bulut estimate. For this reason, we need to temporarily replace the remainder equations with an auxiliary problem, which requires a certain notational effort. 

For the rest of this proof, we use the super-script ``$\fn$" to represent objects which correspond to the frequency-truncated wave maps equation \eqref{main:eq-WM-null-N} with the initial data from \eqref{main:eq-null-lwp-initial-W}, i.e., the initial data without frequency-truncation. We let $U^{\fn}$ and $V^{\fn}$ be the unique global solutions of 
\begin{equation}\label{main:eq-null-2-pb-1}
\begin{cases}
\partial_v U^\fn=\partial_u V^{\fn} = \chinull \big[ U^{\fn}, V^{\fn} \big]_{\leq N} - \coup \chinullsquare \Renormn \big( U^{\fn} + V^{\fn} \big), \\
U^{\fn} \big|_{v=u}= W^{(\Rscript,\coup),+}, \quad V^{\fn} \big|_{u=v}= W^{(\Rscript,\coup),-}. 
\end{cases}
\end{equation}
In keeping with our notation, we wrote $\Renormn$ for $\Renorm[\Nscript]$. The existence and uniqueness of a global solution of \eqref{main:eq-null-2-pb-1} follows from soft arguments (see e.g. Lemma \ref{structure:lem-global-flow}), but does not directly yield any bounds on 
$U^{\fn}$ and $V^{\fn}$ which are uniform in $N$. \\

Since we soon work with auxiliary equations, we momentarily decorate all further objects with a bar. We therefore denote pure modulation operators and nonlinear remainders by 
\begin{equation}\label{main:eq-null-2-pb-2}
\begin{aligned}
 \Big( \widebar{S}^{\fn,\diamond,+}, \widebar{S}^{\fn,\diamond,-} , \widebar{U}^{\fn,\fs} , \widebar{V}^{\fn,\fs} \Big). 
\end{aligned}
\end{equation}
The corresponding modulated and mixed modulated objects $\widebar{U}^{\fn,\ast}$ and $\widebar{V}^{\fn,\ast}$, where 
$\ast \in \{+,-,+-,+\fs,\fs-\}$, are defined using \eqref{main:eq-null-2-pb-2} as in Definitions \ref{ansatz:def-pure}, \ref{ansatz:def-modulated-linear}, \ref{ansatz:def-modulated-bilinear}, \ref{ansatz:def-mixed}, and \ref{ansatz:def-modulated-linear-reversed}. 
We also define
\begin{equation}\label{main:eq-null-2-pb-3}
\widebar{U}^{\fn} = \sum_{\ast} \widebar{U}^{\fn,\ast} \qquad \text{and} \qquad \widebar{V}^{\fn} = \sum_{\ast} \widebar{V}^{\fn,\ast},
\end{equation}
where the sum is taken over all $\ast \in \{ +, - , +-, +\fs, \fs- , \fs \}$.
A similar notation as for $\widebar{U}^{\fn,\ast}$ and $\widebar{V}^{\fn,\ast}$ is used for the error terms from \mbox{Proposition \ref{ansatz:prop-decomposition}}. For fixed nonlinear remainders $\widebar{U}^{\fn,\fs}$ and $\widebar{V}^{\fn,\fs}$, we consider the modulation equations
\begin{equation}\label{main:eq-null-2-pb-4}
\begin{cases}
\begin{aligned}
\partial_v \widebar{S}^{\fn,\diamond,+}_{K,k}  
&=  \chinull[K] \rho_{\leq N}^2(k)   \Big[  \widebar{S}^{\fn,\diamond,+}_{K,k} , \LObarn[K][-] \Big] + \chinull[K]  \rho_{\leq N}^4(k)  \Big[  \widebar{S}^{\fn,\diamond,+}_{K,k} , \SHHLbarn[K][v] \Big],   \\
\partial_u \widebar{S}^{\fn,\diamond,-}_{M,m}  
&=  \chinull[M]    \rho_{\leq N}^2(m)\Big[  \LObarn[M][+] , \widebar{S}^{\fn,\diamond,-}_{M,m}  \Big] +  \chinull[M] \rho_{\leq N}^4(m)  \Big[ \SHHLbarn[M][u], \widebar{S}^{\fn,\diamond,-}_{M,m}  \Big], \\ 
\widebar{S}^{\fn,\diamond,+}_{K,k}  \Big|_{v=u} &= S^{(\Nscript,\Rscript,\coup),\textup{in},+}_{K}, 
\qquad \widebar{S}^{\fn,\diamond,-}_{M,m}  \Big|_{u=v} = S^{(\Nscript,\Rscript,\coup),\textup{in},-}_{M}.
\end{aligned}
\end{cases}
\end{equation}
To introduce the auxiliary equations for $(\widebar{U}^{\fn,\fs},\widebar{V}^{\fn,\fs})$, we recall from Definition \ref{jacobi:def-chhl} and Definition \ref{jacobi:def-modified-jacobi} that 
\begin{equation}\label{main:eq-null-2-pb-5}
\begin{aligned}
\JcbERRbarn &= \JcbERRbarndagger 
+  \Sumlarge_{\substack{\hspace{1.5ex} K \geq N^{1-2\delta_1}}} \chinull \Big( P_{\leq N}^x \CHHLN \big( \widebar{V}^{\fn}, \widebar{V}^{\fn} \big) P_{\leq N}^x - \coup \Renormnbd \Big) \big( \widebar{U}^{\fn,+}_K \big) \\
& + \Sumlarge_{\substack{\hspace{1.5ex}M \geq N^{1-2\delta_1}}} \chinull \Big( P_{\leq N}^x \CHHLN \big(\widebar{U}^{\fn}, \widebar{U}^{\fn}\big) P_{\leq N}^x - \coup \Renormnbd \Big) \big(  \widebar{V}^{\fn,-}_M \big). 
\end{aligned}
\end{equation}
In our auxiliary equations, we insert the decomposition of $\JcbERRbarn$ from \eqref{main:eq-null-2-pb-5}, but replace
$\CHHLN \big(\widebar{U}^{\fn}, \widebar{U}^{\fn}\big)$ and $\CHHLN \big(\widebar{V}^{\fn}, \widebar{V}^{\fn}\big)$ by 
$\CHHLN \big(U^{\fn}, U^{\fn}\big)$ and $\CHHLN \big(V^{\fn},V^{\fn}\big)$, respectively. To avoid confusion, we recall that 
$\widebar{U}^{\fn}$ and $\widebar{V}^{\fn}$ are defined using $(\widebar{U}^{\fn,\fs}, \widebar{V}^{\fn,\fs})$ 
as in \eqref{main:eq-null-2-pb-3}, whereas $U^{\fn}$ and $V^{\fn}$ are defined as the solution of \eqref{main:eq-null-2-pb-1}. 
In particular, $U^{\fn}$ and $V^{\fn}$ do not depend on $(\widebar{U}^{\fn,\fs}, \widebar{V}^{\fn,\fs})$. The auxiliary problem for the nonlinear remainders  $(\widebar{U}^{\fn,\fs}, \widebar{V}^{\fn,\fs})$ is therefore given by the evolution equation
\begin{equation}\label{main:eq-null-2-pb-6}
\begin{aligned}
\partial_v \widebar{U}^{\fn,\fs} &= \chinull \HHLERRbarn + \chinull \JcbERRbarndagger + \SEbarn{u} + \chinull \PIbarn  + \chinull \RenNErrbarn  \\ 
&+ \Sumlarge_{\substack{\hspace{1.5ex} K \geq N^{1-2\delta_1}}} \chinullsquare \Big( P_{\leq N}^x \CHHLN \big( V^{\fn}, V^{\fn} \big) P_{\leq N}^x - \coup \Renormnbd \Big) \big( \widebar{U}^{\fn,+}_K \big) \\ 
&+ \Sumlarge_{\substack{\hspace{1.5ex}M \geq N^{1-2\delta_1}}} \chinullsquare \Big( P_{\leq N}^x \CHHLN \big(U^{\fn}, U^{\fn}\big) P_{\leq N}^x - \coup \Renormnbd \Big) \big(  \widebar{V}^{\fn,-}_M \big),  
\end{aligned}
\end{equation}
the evolution equation 
\begin{equation}\label{main:eq-null-2-pb-7}
\begin{aligned}
\partial_u \widebar{V}^{\fn,\fs} &= \chinull \HHLERRbarn + \chinull \JcbERRbarndagger + \SEbarn{v} + \chinull \PIbarn  + \chinull \RenNErrbarn  \\ 
&+ \Sumlarge_{\substack{\hspace{1.5ex} K \geq N^{1-2\delta_1}}} \chinullsquare \Big( P_{\leq N}^x \CHHLN \big( V^{\fn}, V^{\fn} \big) P_{\leq N}^x - \coup \Renormnbd \Big) \big( \widebar{U}^{\fn,+}_K \big) \\ 
&+ \Sumlarge_{\substack{\hspace{1.5ex}M \geq N^{1-2\delta_1}}} \chinullsquare \Big( P_{\leq N}^x \CHHLN \big(U^{\fn}, U^{\fn}\big) P_{\leq N}^x - \coup \Renormnbd \Big) \big(  \widebar{V}^{\fn,-}_M \big), \\ 
\end{aligned}
\end{equation}
and the initial conditions 
\begin{align}
\widebar{U}^{\fn,\fs} \big|_{v=u} 
&= W^{(\Rscript,\coup),+}(u)  
- \big( \widebar{U}^{\fn,+} + \widebar{U}^{\fn,+-} + \widebar{U}^{\fn,-} + \widebar{U}^{\fn,+\fs} + \widebar{U}^{\fn,\fs-} \big)\big|_{v=u}, 
\label{main:eq-null-2-pb-8}\\  
 \widebar{V}^{\fn,\fs} \big|_{u=v} 
&= W^{(\Rscript,\coup),-}(v)  
- \big( \widebar{V}^{\fn,-} + \widebar{V}^{\fn,+-} + \widebar{V}^{\fn,+} + \widebar{V}^{\fn,\fs-} + \widebar{V}^{\fn,+\fs} \big)\big|_{u=v}.
\label{main:eq-null-2-pb-9}
\end{align}
In the next sub-step, we solve the auxiliary equations \eqref{main:eq-null-2-pb-6}-\eqref{main:eq-null-2-pb-9} using a contraction-mapping argument. \\ 

\emph{Step 2.(c): Contraction-mapping argument for the auxiliary equations.}
Similar as in Step 1.(c), we first introduce the ball 
\begin{equation}\label{main:eq-null-2-pc-1}
\widebar{\mathbb{B}}^{\fn,\fs} := \Big\{ \big( \widebar{U}^{\fn,\fs}, \widebar{V}^{\fn,\fs} \big) \in (\Cprod{r-1}{r} \times \Cprod{r}{r-1})(\R^{1+1})\colon 
\big\| \widebar{U}^{\fn,\fs} \big\|_{\Cprod{r-1}{r}}, \big\| \widebar{V}^{\fn,\fs} \big\|_{\Cprod{r}{r-1}} 
\leq \Dc \Big\}.
\end{equation}
For any pair $(\widebar{U}^{\fn,\fs},\widebar{V}^{\fn,\fs})\in \widebar{\mathbb{B}}^{\fn,\fs}$, Proposition \ref{modulation:prop-main} implies the existence of a unique solution 
\begin{equation*}
\big( \widebar{S}^{\fn,\diamond,+}, \widebar{S}^{\fn,\diamond,-} \big)
= \big( \widebar{S}^{\fn,\diamond,+}, \widebar{S}^{\fn,\diamond,-} \big)\Big( S^{(\Nscript,\Rscript,\coup),\textup{in},\pm}, \widebar{U}^{\fn,\fs}, \widebar{V}^{\fn,\fs} \Big) 
\end{equation*}
of  the modulation equations \eqref{main:eq-null-2-pb-4}. Similar as in Step 1.(c), we also define the 
 maps $\widebar{\Uscr}^{\fn,\fs}$ and $\widebar{\Vscr}^{\fn,\fs}$ by  
\begin{align*}
\widebar{\Uscr}^{\fn,\fs}\big(\widebar{U}^{\fn,\fs},\widebar{V}^{\fn,\fs}\big)
&= \textup{ RHS of } \eqref{main:eq-null-2-pb-8} + \Int^v_{u\rightarrow v} \Big( \textup{RHS of } \eqref{main:eq-null-2-pb-6} \Big), \\ 
\widebar{\Vscr}^{\fn,\fs}\big(\widebar{U}^{\fn,\fs},\widebar{V}^{\fn,\fs}\big)
&= \textup{ RHS of } \eqref{main:eq-null-2-pb-9} + \Int^u_{v\rightarrow u} \Big( \textup{RHS of } \eqref{main:eq-null-2-pb-7} \Big).
\end{align*}
In order to solve the auxiliary equations \eqref{main:eq-null-2-pb-6}-\eqref{main:eq-null-2-pb-9}, it then only remains to prove that 
$(\widebar{\Uscr}^{\fn,\fs},\widebar{\Vscr}^{\fn,\fs})$ is contraction on $\widebar{\mathbb{B}}^{\fn,\fs}$. While the argument is mostly similar as in Step 1.(c), there is a crucial difference in the treatment of the Jacobi errors. Instead of \eqref{main:eq-null-pd-3} and \eqref{main:eq-null-pd-4}, we use the three estimates\footnote{We emphasize that the estimates \eqref{main:eq-null-2-pc-3} and \eqref{main:eq-null-2-pc-4}, which rely on the Bourgain-Bulut estimate from Definition \ref{main:def-Bourgain-Bulut} and Corollary \ref{main:cor-bourgain-bulut-probability}, are much better than any of the available perturbative estimates of the left-hand side.}
\begin{align}
\Big\| \JcbERRbarndagger  \Big\|_{\Cprod{r-1}{r-1}} 
&\lesssim 
\label{main:eq-null-2-pc-2} \Dc^3, \\
\Sumlarge_{\substack{\hspace{1.5ex} K \geq N^{1-2\delta_1}}}  
\Big\| \chinull \Big( P_{\leq N}^x \CHHLN \big( V^{\fn}, V^{\fn} \big) P_{\leq N}^x - \coup \Renormnbd \Big) \big( \widebar{U}^{\fn,+}_K \big) \Big\|_{\Cprod{r-1}{r-1}} 
&\lesssim N^{-\frac{1}{2}+2\delta} \Dc^3
\label{main:eq-null-2-pc-3}, \\ 
\Sumlarge_{\substack{\hspace{1.5ex} M \geq N^{1-2\delta_1}}}   
\Big\| \chinull \Big( P_{\leq N}^x \CHHLN \big( U^{\fn}, U^{\fn} \big) P_{\leq N}^x - \coup \Renormnbd \Big) \big( \widebar{V}^{\fn,+}_M \big) \Big\|_{\Cprod{r-1}{r-1}} 
&\lesssim N^{-\frac{1}{2}+2\delta} \Dc^3
\label{main:eq-null-2-pc-4}. 
\end{align}
The first estimate \eqref{main:eq-null-2-pc-2} follows directly from \eqref{jacobi:eq-main-estimate-2} in Proposition \ref{jacobi:prop-main}.
To obtain the second estimate, we first use \eqref{jacobi:eq-bourgain-bulut-1} from Proposition \ref{jacobi:prop-main}, which yields that
\begin{equation}\label{main:eq-null-2-pc-5}
\begin{aligned}
&\Sumlarge_{\substack{\hspace{1.5ex} K \geq N^{1-2\delta_1}}}  
\Big\| \chinull \Big( P_{\leq N}^x \CHHLN \big( V^{\fn}, V^{\fn} \big) P_{\leq N}^x - \coup \Renormnbd \Big) \big( \widebar{U}^{\fn,+}_K \big) \Big\|_{\Cprod{r-1}{r-1}} \\
\lesssim&\, N^{\delta_1+5\delta_2} \Dc \sup_{y\in \R} \Big( \langle N y \rangle^{-10} 
\Big\| \, \chi \Big( \CHHLN_y \big( V^{\fn}, V^{\fn} \big) - \coup \CNbd(y) \Kil \Big) \Big\|_{L_t^\infty L_x^\infty} \Big).
\end{aligned}
\end{equation}
Since $\chi \in \Cut$, $\chi$ is supported on the time-interval $[-4,4]$. Using Definition \ref{main:def-Bourgain-Bulut} and the restriction to the event from \eqref{main:eq-null-2-pa-2}, it then follows that 
\begin{equation*}
\eqref{main:eq-null-2-pc-5}\lesssim N^{\delta_1+5\delta_2} N^{-\frac{1}{2}+\delta} \Dc^3 \lesssim N^{-\frac{1}{2}+2\delta} \Dc^3,
\end{equation*}
which implies \eqref{main:eq-null-2-pc-3}. The proof of the third  estimate \eqref{main:eq-null-2-pc-4} is similar. Using \eqref{main:eq-null-2-pc-2}, \eqref{main:eq-null-2-pc-3}, and \eqref{main:eq-null-2-pc-4}, as well as similar arguments for $\HHLERRbarn$, $\SEbarn{u}$, $\PIbarn$, and  $\RenNErrbarn$ as in Step 1.(d), it follows that $(\widebar{\Uscr}^{\fn,\fs},\widebar{\Vscr}^{\fn,\fs})$ maps $\widebar{\mathbb{B}}^{\fn,\fs}$ into $\widebar{\mathbb{B}}^{\fn,\fs}$. \\

In order to obtain that $(\widebar{\Uscr}^{\fn,\fs},\widebar{\Vscr}^{\fn,\fs})$ is a contraction on $\widebar{\mathbb{B}}^{\fn,\fs}$, we use a similar argument as in Step 1.(e). The only new ingredients are the Lipschitz-variants of \eqref{main:eq-null-2-pc-2}, \eqref{main:eq-null-2-pc-3}, and \eqref{main:eq-null-2-pc-4}, which can be obtained from \eqref{lipschitz:eq-remainder-jcb-dagger} in Proposition \ref{lipschitz:prop-remainder-no-jcb} and Lemma \ref{lipschitz:lem-bourgain-bulut}. All in all, the contraction-mapping principle then implies the existence of a unique solution of \eqref{main:eq-null-2-pb-6}-\eqref{main:eq-null-2-pb-9} in $\widebar{\mathbb{B}}^{\fn,\fs}$.\\

\emph{Step 2.(d): From the auxiliary equations to the wave maps and remainder equations.} 
In Step 2.(c), we solved the auxiliary equations \eqref{main:eq-null-2-pb-6}-\eqref{main:eq-null-2-pb-9} using a contraction-mapping argument. In this sub-step, we show that $(\widebar{U}^{\fn},\widebar{V}^{\fn})=(U^{\fn},V^{\fn})$, i.e., that $(\widebar{U}^{\fn},\widebar{V}^{\fn})$ defined as in \eqref{main:eq-null-2-pb-3} is a solution of the wave maps equation \eqref{main:eq-null-2-pb-1}. Furthermore, we show that $(\widebar{U}^{\fn,\fs},\widebar{V}^{\fn,\fs})$ solves the remainder equations corresponding to \eqref{main:eq-null-2-pb-1}. \\

From Proposition \ref{ansatz:prop-decomposition}, we first obtain that 
\begin{equation}\label{main:eq-null-2-pd-1}
\begin{aligned}
&\, \chinull \HHLERRbarn + \chinull \JcbERRbarndagger + \SEbarn{u} + \chinull \PIbarn  + \chinull \RenNErrbarn \\ 
=& \,\chinull  \big[ \widebar{U}^{\fn},\widebar{V}^{\fn}\big]_{\leq N}   
- \coup \chinullsquare \Renormn \big( \widebar{U}^{\fn} + \widebar{V}^{\fn} \big)
- \partial_v \big( \widebar{U}^{\fn,+} + \widebar{U}^{\fn,+-}+ \widebar{U}^{\fn,-} + \widebar{U}^{\fn,+\fs}+ \widebar{U}^{\fn,\fs- }\big) \\
+&\,   \chinull \big( \JcbERRbarndagger - \JcbERRbarn \big).
\end{aligned}
\end{equation}
Together with \eqref{main:eq-null-2-pb-3}, it then follows from \eqref{main:eq-null-2-pb-6} and  \eqref{main:eq-null-2-pd-1} that 
\begin{equation}\label{main:eq-null-2-pd-2}
\begin{aligned}
\partial_v \widebar{U}^{\fn} 
&= \partial_v \big( \widebar{U}^{\fn,+} + \widebar{U}^{\fn,+-}+ \widebar{U}^{\fn,-} + \widebar{U}^{\fn,+\fs}+ \widebar{U}^{\fn,\fs- }\big) 
+ \partial_v \big( \widebar{U}^{\fn,\fs}\big) \\ 
&= \chinull  \big[ \widebar{U}^{\fn},\widebar{V}^{\fn}\big]_{\leq N}   
- \coup \chinullsquare \Renormn \big( \widebar{U}^{\fn} + \widebar{V}^{\fn} \big) \\
&+ \chinull \big( \JcbERRbarndagger - \JcbERRbarn \big) \\ 
&+ \Sumlarge_{\substack{\hspace{1.5ex} K \geq N^{1-2\delta_1}}} \chinullsquare \Big( P_{\leq N}^x \CHHLN \big( V^{\fn}, V^{\fn} \big) P_{\leq N}^x - \coup \Renormnbd \Big) \big( \widebar{U}^{\fn,+}_K \big) \\ 
&+ \Sumlarge_{\substack{\hspace{1.5ex}M \geq N^{1-2\delta_1}}} \chinullsquare \Big( P_{\leq N}^x \CHHLN \big(U^{\fn}, U^{\fn}\big) P_{\leq N}^x - \coup \Renormnbd \Big) \big(  \widebar{V}^{\fn,-}_M \big).
\end{aligned}
\end{equation}
By inserting the formula from \eqref{main:eq-null-2-pb-5} for $ \JcbERRbarndagger - \JcbERRbarn$ into \eqref{main:eq-null-2-pd-2}, and by repeating a similar computation for $\widebar{V}^{\fn}$, we then obtain that
\begin{equation}\label{main:eq-null-2-pd-3}
\begin{aligned}
\partial_v \widebar{U}^{\fn} = \partial_u \widebar{V}^{\fn} 
&= \chinull  \big[ \widebar{U}^{\fn},\widebar{V}^{\fn}\big]_{\leq N}   
- \coup \chinullsquare \Renormn \big( \widebar{U}^{\fn} + \widebar{V}^{\fn} \big) \\
&+ \Sumlarge_{\substack{\hspace{1.5ex} K \geq N^{1-2\delta_1}}} \chinullsquare P_{\leq N}^x 
\Big( \CHHLN \big( V^{\fn}, V^{\fn} \big) - \CHHLN \big( \widebar{V}^{\fn}, \widebar{V}^{\fn} \big) \Big) 
P_{\leq N}^x  \big( \widebar{U}^{\fn,+}_K \big) \\ 
&+ \Sumlarge_{\substack{\hspace{1.5ex}M \geq N^{1-2\delta_1}}} \chinullsquare  P_{\leq N}^x 
\Big(\CHHLN \big(U^{\fn}, U^{\fn}\big)  - \CHHLN \big( \widebar{U}^{\fn}, \widebar{U}^{\fn} \big) \Big)
P_{\leq N}^x  \big(  \widebar{V}^{\fn,-}_M \big).
\end{aligned}
\end{equation}
Due to Definition \ref{jacobi:def-chhl} and Remark \ref{jacobi:rem-acting-spatial-variable}, \eqref{main:eq-null-2-pd-3} can be written in Cartesian coordinates as an ordinary differential equation for $(\widebar{U}^{\fn},\widebar{V}^{\fn})$ on the space of $2\pi R$-periodic functions with $x$-frequencies $\lesssim N$. Due to Picard-Lindel\"{o}f, it then follows that the solution of \eqref{main:eq-null-2-pd-3} is unique.
Since $(U^\fn,V^\fn)$ solves \eqref{main:eq-null-2-pb-1}, it clearly also solves \eqref{main:eq-null-2-pd-3}, and the uniqueness of the solution of \eqref{main:eq-null-2-pd-3} then implies that $(\widebar{U}^{\fn},\widebar{V}^{\fn})=(U^{\fn},V^{\fn})$.\\

Since $(\widebar{U}^{\fn},\widebar{V}^{\fn})=(U^\fn,V^\fn)$, it then follows from \eqref{main:eq-null-2-pb-5}, \eqref{main:eq-null-2-pb-6}, and \eqref{main:eq-null-2-pb-7} that
\begin{align}
\partial_v \widebar{U}^{\fn,\fs} &= \chinull \HHLERRbarn + \chinull \JcbERRbarn + \SEbarn{u} + \chinull \PIbarn  + \chinull \RenNErrbarn, 
\label{main:eq-null-2-pd-4}\\ 
\partial_u \widebar{V}^{\fn,\fs} &= \chinull \HHLERRbarn + \chinull \JcbERRbarn + \SEbarn{v} + \chinull \PIbarn  + \chinull \RenNErrbarn.
\label{main:eq-null-2-pd-5}
\end{align}
Due to Proposition \ref{ansatz:prop-decomposition}, it therefore follows that $(\widebar{U}^{\fn,\fs},\widebar{V}^{\fn,\fs})$ 
solves the remainder equations corresponding to the frequency-truncated wave maps equation \eqref{main:eq-null-2-pb-1}.\\

\emph{Step 2.(e): Notation revisited.} Due to Step 2.(d), we know that the solution $(\widebar{U}^{\fn,\fs},\widebar{V}^{\fn,\fs})$ of the auxiliary equations \eqref{main:eq-null-2-pb-6}-\eqref{main:eq-null-2-pb-9} actually solves the remainder equations corresponding to \eqref{main:eq-null-2-pb-1}. Since the bar in our notation for $(\widebar{U}^{\fn,\fs},\widebar{V}^{\fn,\fs})$ was only used to distinguish between the auxiliary and remainder equations, the bar is now obsolete, and we simply write
\begin{equation*}
\big( U^{\fn,\fs} V^{\fn,\fs}\big) := \big(\widebar{U}^{\fn,\fs},\widebar{V}^{\fn,\fs}\big). 
\end{equation*}
Similarly, we omit the bar from our notation for the modulated and mixed modulated objects in \eqref{main:eq-null-2-pb-3} and write
$\HHLERRn$, $\JcbERRn$, $\JcbERRndagger$, $\SEn{u}$, $\SEn{v}$, $\PIn$, and $\RenNErrn$ for the error terms in \eqref{main:eq-null-2-pb-5}-\eqref{main:eq-null-2-pb-7}.  \\

\emph{Step 3: Lipschitz estimate.} 
We now want to prove the Lipschitz estimates in \ref{main:item-null-lwp-C2}. Due to \mbox{Proposition \ref{lipschitz:prop-modulation}}, 
it suffices to prove the Lipschitz estimates for $(U^{\fd,\fs},V^{\fd,\fs})$ and $(U^{\fn,\fs},V^{\fn,\fs})$. 
To this end, we write the difference of $U^{\fd,\fs}$ and $U^{\fn,\fs}$ as
\begin{align}
&\,U^{\fd,\fs} - U^{\fn,\fs} \notag \\ 
=&\,  \big( \widetilde{W}^{(\Nscript,\Ndscript,\Rscript,\coup),+}(u) - U^{\fd,+}(u,u) \big) 
- \big( W^{(\Rscript,\coup),+}(u) - U^{\fn,+}(u,u) \big) 
\label{main:item-null-p3-1}\\
-&\,  \big( U^{\fd,+-}  + U^{\fd,-} +  U^{\fd,+\fs} + U^{\fd,\fs-}  - U^{\fn,+-}  -U^{\fn,-}  -U^{\fn,+\fs}  
- U^{\fn,\fs-}\big)(u,u) 
 \label{main:item-null-p3-2}\\
+&\,  \Int^v_{u\rightarrow v} \Big( \chinull
\big( \HHLERRd - \HHLERRn \big) + \SEd{u} - \SEn{u} + \chinull \big( \PId - \PIn \big) 
+ \chinull \big( \RenNErrd - \RenNErrn \big) \Big) 
\label{main:item-null-p3-3}\\
+&\, \Int^v_{u\rightarrow v} \Big( \chinull \big( \JcbERRd - \JcbERRn \big) \Big). 
\label{main:item-null-p3-4}
\end{align}
The initial data terms in \eqref{main:item-null-p3-1} can be estimated using 
\eqref{main:eq-null-lwp-initial-W}, \eqref{main:eq-null-lwp-initial-W-tilde}, and 
Lemma \ref{lipschitz:lem-modulated-initial}, which yield that
\begin{equation}\label{main:item-null-p3-5}
\begin{aligned}
\big\| \eqref{main:item-null-p3-1} \big\|_{\Cprod{r-1}{r}} 
&\lesssim 
\hcoup \Ac \sup_{L\in \dyadic} \big\| \widetilde{S}^{(\Nscript,\Ndscript,\Rscript,\coup),\textup{in},\pm}_L - \Id_\frkg \big\|_{\C^{s-1}_x}
+ \Dc \Nd^{-\varsigma}
+\big\| \widetilde{Z}^{(\Nscript,\Ndscript,\Rscript,\coup),\pm} \big\|_{\C_x^{r-1}} \\
&+  \Dc \Big( \big\| U^{\fd,\fs} - U^{\fn,\fs} \big\|_{\Cprod{r-1}{r}} + \big\| V^{\fd,\fs} - V^{\fn,\fs} \big\|_{\Cprod{r}{r-1}} \Big).
\end{aligned}
\end{equation}
In order to simplify the notation, we now use a similar notation as in Definition \ref{lipschitz:def-distances} and write
\begin{equation}\label{main:eq-dpost}
\begin{aligned}
\dpost &:=  \Dc \Nd^{-\varsigma} 
+ \hcoup \Ac  \sup_{K\in \dyadic} \big\| \widetilde{S}^{(\Nscript,\Ndscript,\Rscript,\coup),\textup{in},+}_{K} - \Id_\frkg \big\|_{\C_u^s}
+ \hcoup \Ac \sup_{M\in \dyadic} \big\| \widetilde{S}^{(\Nscript,\Ndscript,\Rscript,\coup),\textup{in},-}_{M} - \Id_\frkg \big\|_{\C_v^s} \\ 
&+ \big\| U^{\fd,\fs} - U^{\fn,\fs} \big\|_{\Cprod{r-1}{r}} 
+ \big\| V^{\fd,\fs} - V^{\fn,\fs} \big\|_{\Cprod{r}{r-1}}. 
\end{aligned}
\end{equation}
The terms in \eqref{main:item-null-p3-2} can be estimated using Lemma \ref{lipschitz:lem-modulated-cartesian} and the terms in \eqref{main:item-null-p3-3} can be estimated using Proposition \ref{lipschitz:prop-remainder-no-jcb}, which yield that
\begin{equation}\label{main:item-null-p3-6}
\big\| \eqref{main:item-null-p3-2} \big\|_{\Cprod{r-1}{r}}
+ \big\| \eqref{main:item-null-p3-3} \big\|_{\Cprod{r-1}{r}}
\lesssim \Dc \dpost. 
\end{equation}
In order to treat \eqref{main:item-null-p3-4}, we distinguish the cases in \ref{main:item-null-lwp-short}
and \ref{main:item-null-lwp-local} of the statement of this proposition. Under the condition $\Dc\leq c_0 N^{-\delta_2}$, 
we use Proposition \ref{lipschitz:prop-remainder-no-jcb} and argue as in \eqref{main:eq-null-2-pa-1}, which yields
\begin{equation}\label{main:item-null-p3-7}
\big\|  \eqref{main:item-null-p3-4} \big\|_{\Cprod{r-1}{r}} 
\lesssim  \Dc \dpost.
\end{equation}
If we only have the condition $\Dc\leq c_0$, then we restrict to the event in \eqref{main:eq-null-2-pa-2}. Since
$\Nd\leq N^{1-\delta}$, it holds that $U^{\fd,+}_K=V^{\fd,-}_M=0$ for all $K,M>N^{1-2\delta_1}$, and it therefore follows from Definition \ref{jacobi:def-modified-jacobi} that 
\begin{equation}\label{main:item-null-p3-8}
\JcbERRd = \JcbERRddagger. 
\end{equation}
As a result, we may decompose
\begin{align}
\JcbERRd- \JcbERRn 
&= \JcbERRddagger - \JcbERRndagger 
\label{main:item-null-p3-9} \\
& + \Sumlarge_{\substack{\hspace{1.5ex} K \geq N^{1-2\delta_1}}} \chinull \Big( P_{\leq N}^x \CHHLN \big( V^{\fn}, V^{\fn} \big) P_{\leq N}^x - \coup \Renormnbd \Big) \big( U^{\fn,+}_K \big) 
\label{main:item-null-p3-10} \\
& + \Sumlarge_{\substack{\hspace{1.5ex}M \geq N^{1-2\delta_1}}} \chinull \Big( P_{\leq N}^x \CHHLN \big(U^{\fn}, U^{\fn}\big) P_{\leq N}^x - \coup \Renormnbd \Big) \big(  V^{\fn,-}_M \big).
\label{main:item-null-p3-11}
\end{align}
By estimating \eqref{main:item-null-p3-9} using Proposition \ref{lipschitz:prop-remainder-no-jcb} and estimating \eqref{main:item-null-p3-10} and \eqref{main:item-null-p3-11} as before in \eqref{main:eq-null-2-pc-3} and \eqref{main:eq-null-2-pc-4}, we obtain that
\begin{equation}\label{main:item-null-p3-12} 
\big\|  \eqref{main:item-null-p3-3} \big\|_{\Cprod{r-1}{r}} 
\lesssim  \Dc^2 \dpost 
+ \Dc^3 N^{-\frac{1}{2}+2\delta}.
\end{equation}
Due to the assumption $\Nd \leq N^{1-\delta}$, the $N^{-\frac{1}{2}+2\delta}$-term in \eqref{main:item-null-p3-12} is negligible compared to the $\Nd^{-\varsigma}$-term coming from \eqref{main:eq-dpost}.
Together with a similar argument for $V^{\fd,\fs}-V^{\fn,\fs}$,
we obtain from our estimates in \eqref{main:item-null-p3-5}, \eqref{main:item-null-p3-6}, \eqref{main:item-null-p3-7}, and \eqref{main:item-null-p3-12} that
\begin{align*}
&\, \big\| U^{\fd,\fs} - U^{\fn,\fs} \big\|_{\Cprod{r-1}{r}} 
+ \big\| V^{\fd,\fs} - V^{\fn,\fs} \big\|_{\Cprod{r}{r-1}}  \\
\lesssim&\,
\hcoup \Ac \sup_{L\in \dyadic} \big\| \widetilde{S}^{(\Nscript,\Ndscript,\Rscript,\coup),\textup{in},\pm}_L - \Id_\frkg \big\|_{\C^{s-1}_x}
+ \Dc \Nd^{-\varsigma} +\big\| \widetilde{Z}^{(\Nscript,\Ndscript,\Rscript,\coup),\pm} \big\|_{\C_x^{r-1}} \\
+&\,  \Dc \Big( \big\| U^{\fd,\fs} - U^{\fn,\fs} \big\|_{\Cprod{r-1}{r}} + \big\| V^{\fd,\fs} - V^{\fn,\fs} \big\|_{\Cprod{r}{r-1}} \Big). 
\end{align*}
Using $\Dc \leq c_0$ and Assumption \ref{main:item-null-lwp-A2}, it then follows that 
\begin{align*}
  \big\| U^{\fd,\fs} - U^{\fn,\fs} \big\|_{\Cprod{r-1}{r}} 
+ \big\| V^{\fd,\fs} - V^{\fn,\fs} \big\|_{\Cprod{r}{r-1}} 
\lesssim   \hcoup \Ac \epsilon + \Dc \Nd^{-\varsigma} = \hcoup \Ac \big( \epsilon + \Bc \Nd^{-\varsigma}\big),
\end{align*}
which implies the Lipschitz estimates in \ref{main:item-null-lwp-C2}.\\

\emph{Step 4: Parameter trick.} 
In Definition \ref{main:def-null-lwp}.\ref{main:item-null-lwp-C3}, we claim that the nonlinear remainders $U^{\fn,\fs}$ 
and $V^{\fn,\fs}$ from the second step satisfy
\begin{equation}\label{main:eq-null-p4-1}
\big\| U^{\fn,\fs} \big\|_{\Cprod{r+\eta-1}{r+\eta}} , \big\| V^{\fn,\fs} \big\|_{\Cprod{r+\eta}{r+\eta-1}} \lesssim \Dc. 
\end{equation}
In other words, we claim that in our estimates of $U^{\fn,\fs}$ and $V^{\fn,\fs}$, the regularity parameter $r$ can be replaced by $r+\eta$. 
The reason is that, in our analysis of \eqref{main:eq-null-2-pb-1}, \eqref{main:eq-null-2-pb-4}, and \eqref{main:eq-null-2-pb-6}-\eqref{main:eq-null-2-pb-9}, we can replace our parameter
\begin{equation}\label{main:eq-null-p4-2}
\delta_1 \quad \text{by} \quad \delta_1+\delta_3,
\end{equation}
which can be seen as follows: The definition of the initial data in \eqref{main:eq-null-lwp-initial-W} only involves the Gaussians $(G^\pm_{x_0,\ell})_{x_0\in \LambdaRR,\ell \in \Z}$ and functions $(\psiRx)_{x_0\in \LambdaRR}$, and therefore depends on none of our parameters. The definitions of our modulation operators, modulated objects, and mixed modulation objects (Definitions \ref{ansatz:def-pure}, \ref{ansatz:def-modulated-linear}, \ref{ansatz:def-modulated-bilinear}, \ref{ansatz:def-mixed}, and \ref{ansatz:def-modulated-linear-reversed}) depend on the parameters $\delta$, $\deltap$, $\vartheta$, and $\Nlarge$, which are defined in terms of $\delta_0,\delta_2$, and $\delta_4$. Thus, none of the definitions of the terms in our Ansatz \eqref{ansatz:eq-UN-rigorous-decomposition}-\eqref{ansatz:eq-VN-rigorous-decomposition} depend on $\delta_1$, and our Ansatz is therefore not altered by the replacement \eqref{main:eq-null-p4-2}.
Furthermore, since the parameters 
\begin{equation*}
\big( \delta_0, \delta_1 + \delta_3 , \delta_2, \delta_3, \delta_4, \delta_5 \big)
\end{equation*}
still satisfy the conditions from \eqref{prelim:eq-parameter-sizes}, all of our estimates remain valid after the replacement \eqref{main:eq-null-p4-1}.

We emphasize that \eqref{main:eq-null-p4-1} does not hold for $U^{\fd,\fs}$ and $V^{\fd,\fs}$, i.e., the nonlinear remainders for the wave maps equation with frequency-truncated initial data. The reason is that the initial data in \eqref{main:eq-null-lwp-initial-W-tilde} depends on 
$\widetilde{Z}^{(\Nscript,\Ndscript,\Rscript,\coup),\pm}$, and our assumptions on $\widetilde{Z}^{(\Nscript,\Ndscript,\Rscript,\coup),\pm}$
from Definition \ref{main:def-null-lwp}.\ref{main:item-null-lwp-A2} depend on the regularity parameter $r$. 
\end{proof}

\subsection{Well-posedness}\label{section:main-lwp}

In this subsection, we obtain the well-posedness of \eqref{main:eq-WM-N}, i.e., the finite-dimensional approximation of the wave maps equation in Cartesian coordinates. In comparison to Proposition \ref{main:prop-null-lwp} above, the main result of this subsection (Proposition \ref{main:prop-lwp}) contains less detailed information on the random structure of the solution but can be iterated in time more easily. In Subsection  \ref{section:main-abstract}, we first prove an abstract estimate for our random structure. In Subsection \ref{section:main-lwp-statement}, we then state and prove Proposition \ref{main:prop-lwp}.

\subsubsection{Abstract estimate for our random structure}\label{section:main-abstract}
We now state and prove an abstract lemma regarding modulated objects.  

\begin{lemma}[Abstract estimate for our random structure]\label{main:lem-abstract} 
Let $R\geq 1$, let $\coup>0$, let $\Ac,\Bc \geq 1$, and let $\Dc := \hcoup \Ac \Bc$. 
Let the probabilistic hypothesis (Hypothesis \ref{hypothesis:probabilistic}) be satisfied 
and let $\big((G_{x_0,\ell}),W^{(\Rscript,\coup)}\big)$  be given by either $\big((G_{x_0,\ell}^+),W^{(\Rscript,\coup),+}\big)$ or 
$\big((G_{x_0,\ell}^-),W^{(\Rscript,\coup),-}\big)$. Then, the assumptions \ref{main:item-abstract-A1}-\ref{main:item-abstract-A4} below imply the conclusions \ref{main:item-abstract-C1}-\ref{main:item-abstract-C3} below. Here, the assumptions are defined as follows: 
\begin{enumerate}[label={[A\arabic*]}]
    \item \label{main:item-abstract-A1} \emph{Structure:} Let $\Nd \in \Dyadiclarge$, let 
    \begin{equation*}
    (\Sn[K])_{\substack{ K \lesssim \Nd}}, (\Sd[K])_{K\lesssim \Nd} \colon \bT_R\rightarrow \End(\frkg)
    \end{equation*}
    and let $\Zn,\Zd\colon \bT_R\rightarrow \frkg$. Furthermore, let $\Xd,\Xn\colon \bT_R\rightarrow \frkg$ and assume that  
    \begin{align*}
    P_{\leq \Nd} \Xn &=  
    P_{\leq \Nd} \sum_{\substack{K \lesssim \Nd}} \sum_{x_0\in \LambdaRR}\sum_{k\in \Z_K}  \Sn[K][](x) \psiRxK G_{x_0,k} e^{ikx}
    + P_{\leq \Nd} \Zn, \\ 
    \Xd &=  \sum_{\substack{K \lesssim \Nd}}  \sum_{x_0\in \LambdaRR}\sum_{k\in \Z_K}   \Sd[K](x) \psiRxK \rhoND(k) G_{x_0,k} e^{ikx} + \Zd, 
    \end{align*}
    \item \label{main:item-abstract-A2} \emph{Frequency-support:} For all $K\lesssim \Nd$, it holds that $\Sn[K][]=P_{\ll K} \Sn[K][]$ and $\Sd[K][]=P_{\ll K} \Sd[K][]$. 
    \item \label{main:item-abstract-A3}  \emph{Bounds:} Assume that 
    \begin{alignat*}{3}
    \sup_{\substack{K\in \dyadic\colon \\ K \lesssim \Nd}} \max \Big( 
    &\big\|  \Sn[K][] \big\|_{\C_x^s}, \big\|  \Sd[K][] \big\|_{\C_x^s}, 
    \big\|  ( \Sn[K][] )^{-1} \big\|_{\C_x^s}, &\big\|  ( \Sd[K][] )^{-1} \big\|_{\C_x^s} \Big)
    \leq 2, \\ 
    &\big\| \Zn \big\|_{\C_x^{r+\eta-1}} \leq \Dc, \qquad  \text{and}
    &\big\| \Zd \big\|_{\C_x^{r-1}} \leq \Dc. 
    \end{alignat*}
    \item \label{main:item-abstract-A4}  \emph{Difference:} Let $\epsilon>0$ and assume that
    \begin{equation*}
    \sup_{\substack{K\in \dyadic\colon \\ K \lesssim \Nd}} 
    \Big\|\Sn[K] - \Sd[K] \Big\|_{\C_x^{s}} \leq \epsilon 
    \qquad \text{and} \qquad \big\| \Zn - \Zd \big\|_{\C_x^{r-1}} \leq \hcoup \Ac \epsilon. 
    \end{equation*}
\end{enumerate}
Furthermore, the conclusions \ref{main:item-abstract-C1}-\ref{main:item-abstract-C3} are defined as follows:
\begin{enumerate}[label={[C\arabic*]}]
    \item \label{main:item-abstract-C1} \emph{Structure:} It holds that 
    \begin{equation*}
    \Xd = \sum_{\substack{K \lesssim \Nd }} \Sg[K] (x) \Psharp_{R;K} P_{\leq \Nd} \Xn + \Zg,
    \end{equation*}
    where $(\Psharp_{R;K})_{K\in \dyadic}$ are the sharp Fourier-cutoffs from  \eqref{prelim:eq-Psharp}
    and 
    $(\Sg[K])_{K\lesssim \Nd}\colon \bT_R\rightarrow \End(\frkg)$ and $Z\colon \bT_R\rightarrow \frkg$ satisfy the properties in \ref{main:item-abstract-C2} and \ref{main:item-abstract-C3} below. 
    \item \label{main:item-abstract-C2} \emph{Frequency-support:} For all $K\lesssim \Nd$, it holds that $P_{\gg K^{1-\delta}} \Sg[K]=0$.
    \item \label{main:item-abstract-C3}  \emph{Distance to neutral element:} For $\epsilon^\prime := C_0 \Bc R^\eta  (\epsilon+ \Bc \Nd^{-\varsigma})$, it holds that 
    \begin{equation*}
    \sup_{\substack{K\in \dyadic \colon \\ K\lesssim \Nd }} \Big\| \Sg[K] - \Id_\frkg \Big\|_{\C_x^s} \leq \epsilon^\prime 
    \qquad \text{and} \qquad \big\| \Zg \big\|_{C_x^{r-1}} \leq \hcoup \Ac \epsilon^\prime. 
    \end{equation*}
\end{enumerate}
\end{lemma}

\begin{remark}
The significance of Lemma \ref{main:lem-abstract} is that the structure in \ref{main:item-abstract-C1} is similar as the structure in Definition \ref{ansatz:def-initial-data}, i.e., as in the definition of our initial data. Due to this, we will be able to iterate our local theory (see Proposition \ref{main:prop-lwp} and the proof of Lemma \ref{main:lem-gwp-small}).
\end{remark}

\begin{remark}\label{main:rem-abstract-motivation}
As in the proof of Proposition \ref{main:prop-null-lwp}, 
the superscripts ``$\mathfrak{d}$" and ``$\mathfrak{n}$" in Lemma \ref{main:lem-abstract} stand for ``data" and ``nonlinearity", respectively. 
As Lemma \ref{main:lem-abstract} is rather technical, we also briefly state its claim at a more basic level. Ignoring all frequency projections and dyadic sums, we may write \ref{main:item-abstract-A1} schematically as 
\begin{equation*}
\Xn = \Sn[][] W^{(\Rscript,\coup)} + \Zn \qquad \text{and} \qquad \Xd = \Sd[][] W^{(\Rscript,\coup)} + \Zd. 
\end{equation*}
As long as $\Sn$ is invertible, this can be rewritten as 
\begin{equation*}
\Xd = \Sd (\Sn)^{-1} \Xn + \big( \Zd - \Sd (\Sn)^{-1} \Zn \big) =: S \Xn + Z. 
\end{equation*}
 Now, provided that $\| \Sn - \Sd \| \leq \epsilon$, $\| \Zn - \Zd \|\leq \hcoup \Ac \epsilon$, $\| \Zn \| \leq \Dc$, and $\| \Zd \| \leq \Dc$, we should then have $\| S -\Id \| \lesssim \epsilon$ and 
 \begin{equation*}
 \| Z \| \lesssim \| \Zd - \Zn \| + \| 1 - \Sd (\Sn)^{-1} \| \cdot  \| \Zn \| \lesssim \hcoup \Ac \epsilon + \epsilon \Dc \sim \hcoup \Ac \Bc \epsilon.
 \end{equation*}  
 Up to the additional $\Nd^{-\varsigma}$-term and $R^\eta$-factor, this is exactly the claim in \ref{main:item-abstract-C3}. 
\end{remark}

\begin{proof}
For expository purposes, we split the argument into three steps.\\ 

\emph{First step: Expressions for $\Xd$ and $\Xn$ in terms of $\Psharp_{R;K}$.}
We first write $\Xd$ and $P_{\leq \Nd} \Xn$ as 
\begin{align}
\Xd &= \sum_{K\lesssim \Nd} \Sd[K] P_{\leq \Nd} \Psharp_{R;K} W^{(\Rscript,\coup)} + \Zdtil, 
\label{main:eq-abstract-s1} \\ 
P_{\leq \Nd} \Xn &= P_{\leq \Nd} \sum_{K\lesssim \Nd} \Sn[K]  \Psharp_{R;K} W^{(\Rscript,\coup)} +P_{\leq \Nd} \Zntil, 
 \label{main:eq-abstract-s2} 
\end{align}
where 
\begin{align}
\Zdtil &= \Zd + \sum_{K\lesssim \Nd} \Sd[K] 
\Big(  \sum_{x_0\in \LambdaRR}\sum_{k\in \Z_K}   \psiRxK \rhoND(k) G_{x_0,k} e^{ikx}
- P_{\leq \Nd} \Psharp_{R;K} W^{(\Rscript,\coup)} \Big), 
\label{main:eq-abstract-s3}\\ 
\Zntil &= \Zn + \sum_{K\lesssim \Nd} \Sn[K] 
\Big(  \sum_{x_0\in \LambdaRR}\sum_{k\in \Z_K}   \psiRxK  G_{x_0,k} e^{ikx}
- \Psharp_{R;K} W^{(\Rscript,\coup)} \Big). 
\label{main:eq-abstract-s4} 
\end{align}
We recall that, due to Assumption \ref{main:item-abstract-A2}, $\Sd[K]$ and $\Sn[K]$ are supported on frequencies much smaller than $K$. 
Using the low$\times$high-estimate (Lemma \ref{prelim:lem-paraproduct}), Hypothesis \ref{hypothesis:probabilistic}.\ref{ansatz:item-hypothesis-sharp}, 
Assumption \ref{main:item-abstract-A3}, and Assumption \ref{main:item-abstract-A4}, it then easily follows that 
\begin{align}
&\big\| \Zdtil \big\|_{\C_x^{r-1}}\lesssim \big\| \Zd \big\|_{\C_x^{r-1}} + \hcoup \Ac \lesssim \Dc, 
\label{main:eq-abstract-s5} \\ 
&\big\| \Zntil \big\|_{\C_x^{r+\eta-1}}\lesssim \big\| \Zn \big\|_{\C_x^{r+\eta-1}} + \hcoup \Ac \lesssim \Dc, 
\label{main:eq-abstract-s6} 
\end{align}
and 
\begin{equation}\label{main:eq-abstract-s7}
\begin{aligned}
\big\| \Zntil - \Zdtil  \big\|_{\C_x^{r-1}}
&\lesssim 
\big\| \Zd - \Zntil \big\|_{\C_x^{r-1}} + \hcoup \Ac \sup_{K\in \dyadic}\big\| \Sn[K] - \Sd[K] \big\|_{\C_x^s} 
+  \hcoup \Ac \Nd^{-\varsigma} \\
&\lesssim \hcoup \Ac \big( \epsilon + \Nd^{-\varsigma} \big) 
\lesssim \hcoup \Ac \big( \epsilon + \Bc \Nd^{-\varsigma} \big). 
\end{aligned}
\end{equation}
Thus, $\Zdtil$ and $\Zntil$ effectively satisfy the same estimates as $\Zd$ and $\Zn$. \\ 

\emph{Second step: Expressions for $\Sg[K]$ and $\Zg$.}
To simplify the notation, we write $\Sg[K]=\Id_\frkg+\Eg[K]$, where $(\Eg[K])_{K\lesssim \Nd}\colon \bT_R\rightarrow \End(\frkg)$ remains to be chosen. In order to satisfy \ref{main:item-abstract-C1}, we then have to choose 
\begin{align}
\Zg &= \Xd -  \sum_{\substack{ K \lesssim \Nd }} \big( \Id_\frkg + \Eg[K](x) \big)  \Psharp_{R;K} P_{\leq \Nd} \Xn \notag \\ 
&= \big( \Xd - P_{\leq \Nd} \Xn \big)  
- \sum_{\substack{ K \lesssim \Nd}} \Eg[K](x) \Psharp_{R;K} P_{\leq \Nd} \Xn.
\label{main:eq-abstract-p1} 
\end{align}
By inserting the two decompositions from \eqref{main:eq-abstract-s1} and \eqref{main:eq-abstract-s2}, we can write the first summand in \eqref{main:eq-abstract-p1} as 
\begin{align}
\Xd - P_{\leq \Nd} \Xn  
&=  \sum_{\substack{K \lesssim \Nd}} \Sd[K] P_{\leq \Nd} \Psharp_{R;K}  W^{(\Rscript,\coup)} -  P_{\leq \Nd} \sum_{\substack{K\lesssim \Nd}} \Sn[K] \Psharp_{R;K}  W^{(\Rscript,\coup)} 
\label{main:eq-abstract-p2} \\
&+ \Zdtil - P_{\leq \Nd} \Zntil. \label{main:eq-abstract-p4}  
\end{align}
In order to simplify the notation, we now denote the commutator of two operators $A$ and $B$ by $\Com(A,B)$. Using this notation, we then write the operator  $P_{\leq \Nd}( \Sn[K] \, \cdot )$ in \eqref{main:eq-abstract-p2} as $\Sn[K] P_{\leq \Nd} + \Com(P_{\leq \Nd}, \Sn[K])$ and insert this expression into \eqref{main:eq-abstract-p2}. As a result, we then obtain the identity 
\begin{equation}
\begin{aligned}
& \,  \Xd - P_{\leq \Nd} \Xn  \notag  \\
=&\, 
 \sum_{\substack{K \lesssim \Nd}} \big( \Sd[K] - \Sn[K]\big) \Psharp_{R;K} P_{\leq \Nd} W^{(\Rscript,\coup)} 
+ \sum_{\substack{K \lesssim \Nd}}  \Com \big( \Sn[K], P_{\leq \Nd} \big) \Psharp_{R;K} W^{(\Rscript,\coup)} +  \Zdtil - P_{\leq \Nd} \Zntil. \label{main:eq-abstract-p5} 
\end{aligned}
\end{equation}
After rewriting the first summand in \eqref{main:eq-abstract-p1}, we now rewrite the second summand in \eqref{main:eq-abstract-p1}. Using a similar argument as for \eqref{main:eq-abstract-p2}, we obtain that 
\begin{align}
\sum_{\substack{ K \lesssim \Nd}} \Eg[K](x) \Psharp_{R;K} P_{\leq \Nd} \Xn 
&= 
 \sum_{\substack{ K,L \lesssim \Nd}} \Eg[K] \Sn[L] \Psharp_{R;K} \Psharp_{R;L} P_{\leq \Nd} W^{(\Rscript,\coup)} \label{main:eq-abstract-p6}  \\ 
&+  \sum_{\substack{ K,L \lesssim \Nd}} \Eg[K] \Com\big( \Psharp_{R;K}  P_{\leq \Nd} , \Sn[L] \big)   \Psharp_{R;L}  W^{(\Rscript,\coup)}  
+ \sum_{\substack{K \lesssim \Nd}} \Eg[K] \Psharp_{R;K} \Zntil. \notag 
\end{align}
Since $\Psharp_{R;K} \Psharp_{R;L} = \mathbf{1}\big\{ K=L \big\} \Psharp_{R;K}$, we can write the dyadic double sum in \eqref{main:eq-abstract-p6} as 
\begin{equation*}
\sum_{\substack{ K \lesssim \Nd}} \Eg[K] \Sn[K] \Psharp_{R;K} P_{\leq \Nd} W^{(\Rscript,\coup)}. 
\end{equation*}
All in all, inserting our expressions back into \eqref{main:eq-abstract-p1} yields the decomposition 
\begin{align}
\Zg &= \sum_{\substack{K \lesssim \Nd}} \big( \Sd[K] - \Sn[K] - \Eg[K] \Sn[K] \big) \Psharp_{R;K} P_{\leq \Nd} W^{(\Rscript,\coup)} \label{main:eq-abstract-q1}  \\
&+ \sum_{\substack{K \lesssim \Nd}}  \Com \big( \Sn[K], P_{\leq \Nd} \big) \Psharp_{R;K} W^{(\Rscript,\coup)} 
- \sum_{\substack{ K,L \lesssim \Nd}} \Eg[K] 
\Com\big( \Psharp_{R;K} P_{\leq \Nd} , \Sn[L] \big)   \Psharp_{R;L}  W^{(\Rscript,\coup)}   \label{main:eq-abstract-q2} \\ 
&+  \big( \Zdtil - P_{\leq \Nd} \Zntil \big)  
- \sum_{\substack{K \lesssim \Nd}} \Eg[K] \Psharp_{R;K} \Zntil.\label{main:eq-abstract-q3}  
\end{align}
In order to turn \eqref{main:eq-abstract-q1} into a remainder, and yet still satisfy the frequency-support condition from \ref{main:item-abstract-C2}, we now choose 
\begin{equation}\label{main:eq-abstract-E}
\Eg[K](x) := P_{\leq K^{1-\delta}} \Big( \big( \Sd[K](x)- \Sn[K](x)\big) \big( \Sn[K](x) \big)^{-1} \Big). 
\end{equation}
It now remains to prove the estimates in \ref{main:item-abstract-C3}, for which it suffices to prove the estimates\footnote{We note that the $\Bc$ and $R^{\eta}$-factors in $\epsilon^\prime=C_0 \Bc R^\eta  (\epsilon+\Bc \Nd^{-\varsigma})$ are only needed in the estimate of $Z$, and the additional $\Bc$-factor ultimately stems from the lack of $\Bc$-factors in our bounds on $S^{\fn}$ and $S^{\fd}$ in \ref{main:item-abstract-A3}. This choice, which will be justified using Proposition \ref{modulation:prop-main}, had to be made in order to obtain bounds on the inverses of $S^{\fn}$ and $S^{\fd}$.}
\begin{align*}
\big\| \Eg[K] \big\|_{\C_x^s} \lesssim \epsilon \qquad \text{and} \qquad 
\big\| \eqref{main:eq-abstract-q1} \big\|_{\C_x^{r-1}} 
+ \big\| \eqref{main:eq-abstract-q2} \big\|_{\C_x^{r-1}} 
+ \big\| \eqref{main:eq-abstract-q3} \big\|_{\C_x^{r-1}} 
\lesssim \hcoup \Ac \Bc R^\eta  \big( \epsilon + \Bc \Nd^{-\varsigma} \big).
\end{align*}

\emph{Third step: Estimates of $\Eg[K]$ and $\Zg$.}
We now separately estimate $\Eg[K]$ from \eqref{main:eq-abstract-E} and the terms in \eqref{main:eq-abstract-q1}, \eqref{main:eq-abstract-q2}, and \eqref{main:eq-abstract-q3}. \\ 

\emph{Estimate of $\Eg[K]$:} Using a product estimate, \ref{main:item-abstract-A3}, and \ref{main:item-abstract-A4}, it holds that 
\begin{equation}\label{main:eq-abstract-E-estimate}
\big\| \Eg[K] \big\|_{\C_x^s} \lesssim \big\| \Sd[K](x)- \Sn[K](x) \big\|_{\C_x^s} \big\| \big( \Sn[K](x) \big)^{-1} \big\|_{\C_x^s} \lesssim \epsilon. 
\end{equation}  

\emph{Estimate of \eqref{main:eq-abstract-q1}:} Due to our choice of $\Eg[K]$ in \eqref{main:eq-abstract-E}, it holds that 
\begin{equation*}
\Sd[K] - \Sn[K] - \Eg[K] \Sn[K]  
= P_{>K^{1-\delta}} \Big( \big( \Sd[K](x)- \Sn[K](x)\big) \big( \Sn[K](x) \big)^{-1} \Big) \Sn[K]. 
\end{equation*}
Using product estimates, \ref{main:item-abstract-A3}, and \ref{main:item-abstract-A4}, we then obtain for all $0<\alpha \leq s$ that 
\begin{align*}
 \Big\| \Sd[K] - \Sn[K] - \Eg[K] \Sn[K]   \Big\|_{\C_x^\alpha} 
\lesssim&\,  \Big\|  P_{>K^{1-\delta}} \Big( \big( \Sd[K]- \Sn[K]\big) \big( \Sn[K] \big)^{-1} \Big) \Big\|_{\C_x^\alpha} \big\| \Sn[K] \big\|_{\C_x^\alpha} \\ 
\lesssim&\,  K^{-(1-\delta)(s-\alpha)} \big\| \Sd[K] -\Sn[K] \big\|_{\C_x^s} 
\big\| \big( \Sn[K] \big)^{-1} \big\|_{\C_x^s} \big\| \Sn[K] \big\|_{\C_x^s} \\ 
\lesssim&\, \epsilon K^{-(1-\delta)(s-\alpha)}. 
\end{align*}
Since $\Eg[K],\Sd[K]$, and $\Sn[K]$ are supported on frequencies $\ll K$, \eqref{main:eq-abstract-q1} only contains low$\times$high-interactions. Using a low$\times$high-paraproduct estimate and Hypothesis \ref{hypothesis:probabilistic}.\ref{ansatz:item-hypothesis-sharp}, it then follows that
\begin{align*}
\big\| \eqref{main:eq-abstract-q1}  \big\|_{\C_x^{r-1}} 
&\lesssim \sum_{\substack{ K \lesssim \Nd}} 
\big\| \Sd[K] - \Sn[K] - \Eg[K] \Sn[K] \big\|_{\C_x^\eta} \big\| \Psharp_{R;K} W^{(\Rscript,\coup)} \big\|_{\C_x^{r-1}} \notag \\ 
&\lesssim  \hcoup \Ac   \epsilon
\sum_{\substack{K \lesssim \Nd}} 
K^{-(1-\delta)(s-\eta)} K^{r-\frac{1}{2}+\eta} 
\lesssim \hcoup \Ac   \epsilon.
\end{align*} 

\emph{Estimate of \eqref{main:eq-abstract-q2}:} We estimate the two summands in \eqref{main:eq-abstract-q2} separately. For the first summand in \eqref{main:eq-abstract-q2},  we use Lemma \ref{prelim:lem-commutator-basic}, Hypothesis \ref{hypothesis:probabilistic}.\ref{ansatz:item-hypothesis-sharp}, and Assumption \ref{main:item-abstract-A2}, which yields that 
\begin{align*}
&\big\| \Com \big( \Sn[K], P_{\leq \Nd} \big) \Psharp_{R;K} W^{(\Rscript,\coup)} \big\|_{\C_x^{r-1}} 
\lesssim   \Nd^{-1} \big\| \Sn[K] \big\|_{\C_x^1} \big\| \Psharp_{R;K} W^{(\Rscript,\coup)} \big\|_{\C_x^{r-1}} \\
\lesssim &\,  K^{1-s} \Nd^{-1} \big\| \Sn[K] \big\|_{\C_x^s} \big\| \Psharp_{R;K} W^{(\Rscript,\coup)} \big\|_{\C_x^{r-1}}
\lesssim  \hcoup \Ac K^{1-s} K^{r-\frac{1}{2}+\eta} \Nd^{-1}. 
\end{align*}
Since $K^{1-s} K^{r-\frac{1}{2}+\eta} \Nd^{-1} \lesssim \Nd^{-\frac{1}{2}+r-s+\eta} \lesssim \Nd^{-\frac{1}{2}+2\delta_1}$, this is acceptable. For the second summand in \eqref{main:eq-abstract-q2}, we use Lemma \ref{prelim:lem-paraproduct}, Lemma \ref{prelim:commutator-PNx}, Hypothesis \ref{hypothesis:probabilistic}.\ref{ansatz:item-hypothesis-sharp}, and \eqref{main:eq-abstract-E-estimate}.  We then obtain that
\begin{align*}
& \big\| \Eg[K] \Com \big( \Psharp_{R;K} P_{\leq \Nd}, \Sn[L] \big) \Psharp_{R;L} W^{(\Rscript,\coup)} \big\|_{\C_x^{r-1}} 
\lesssim \big\| \Eg[K] \big\|_{\C_x^s} \big\| \Com \big( \Psharp_{R;K} P_{\leq \Nd}, \Sn[L] \big) \Psharp_{R;L} W^{(\Rscript,\coup)} \big\|_{\C_x^{r-1}}   \\ 
\lesssim&\,  \hcoup \Ac \max(K,L)^{-\frac{1}{2}+\delta} \big\| \Eg[K] \big\|_{\C_x^s} 
\big\| \Sn[L] \big\|_{\C_x^s} \lesssim \hcoup \Ac  \epsilon  \max(K,L)^{-\frac{1}{2}+\delta},
\end{align*}
which is acceptable. \\ 

\emph{Estimate of \eqref{main:eq-abstract-q3}:} 
We treat the two summands in \eqref{main:eq-abstract-q3} separately. In order to control the first summand in \eqref{main:eq-abstract-q3}, we use \eqref{main:eq-abstract-s6} and \eqref{main:eq-abstract-s7}, which yield that
\begin{equation*}
\big\| \Zdtil - P_{\leq \Nd} \Zntil \big\|_{\C_x^{r-1}}
\lesssim \big\| \Zdtil - \Zntil \big\|_{\C_x^{r-1}} + \Nd^{-\eta} \big\| P_{>\Nd} \Zntil \big\|_{\C_x^{r+\eta-1}} \lesssim 
\hcoup \Ac   \big( \epsilon +  \Bc \Nd^{-\varsigma}\big). 
\end{equation*}
For the second summand in \eqref{main:eq-abstract-q3}, we use Lemma \ref{prelim:lem-PKsharp}, \eqref{main:eq-abstract-s6}, and  \eqref{main:eq-abstract-E}, which yield that
\begin{align}\label{main:eq-abstract-R-loss}
&\big\| \Eg[K] \Psharp_{R;K} \Zntil \big\|_{\C_x^{r-1}} 
\lesssim \big\| \Eg[K] \big\|_{\C_x^s} \big\| \Psharp_{R;K} \Zntil \big\|_{\C_x^{r-1}} 
\lesssim R^\eta K^{-\frac{\eta}{2}} \big\| \Eg[K] \big\|_{\C_x^s}  \big\| \Zntil \big\|_{\C_x^{r+\eta-1}} \lesssim \big( R^\eta K^{-\frac{\eta}{2}} \big)  \Dc \epsilon. 
\end{align}
Since $\Dc=\hcoup \Ac \Bc$, this is acceptable. 
We remark that \eqref{main:eq-abstract-R-loss} is the only estimate responsible for the $R^\eta$-factors in \ref{main:item-abstract-C3}.
\end{proof}

\subsubsection{Statement and proof of well-posedness}\label{section:main-lwp-statement}

In order to state Proposition \ref{main:prop-lwp}, we first define a variant of the set $\Perpm(W^{(\Rscript,\coup),+},W^{(\Rscript,\coup),-};\Ac,\epsilon)$ from Definition \ref{main:def-perturbations-null}.

\begin{definition}[Perturbations]\label{main:def-perturbations}
Let $N,\Nd\in \dyadiclarge$, let $R \geq 1$, let $\coup>0$, let $\Ac,\Bc\geq 1$, and let $\Dc=\hcoup \Ac \Bc$. Furthermore, let $(W_0^{(\Rscript,\coup)},W_1^{(\Rscript,\coup)})\colon \bT_R\rightarrow \frkg^2$ be initial data, let 
\begin{equation*}
W^{(\Rscript,\coup),\pm} = \tfrac{1}{4} \big( W_0^{(\Rscript,\coup)} \mp W_1^{(\Rscript,\coup)} \big)
\end{equation*}
and let $\epsilon>0$.
Then, we define $\Per \big( W^{(\Rscript,\coup)}_{0}, W^{(\Rscript,\coup)}_1;\Ac,\epsilon)$ as the set of all $(\widetilde{W}^{(\Nscript,\Ndscript,\Rscript,\coup)}_{0},\widetilde{W}^{(\Nscript,\Ndscript,\Rscript,\coup)}_{1})$ 
which can be written as 
\begin{align*}
\widetilde{W}^{(\Nscript,\Ndscript,\Rscript,\coup)}_{0} 
&= 2 \, \big( \widetilde{W}^{(\Nscript,\Ndscript,\Rscript,\coup),+} + \widetilde{W}^{(\Nscript,\Ndscript,\Rscript,\coup),-}\big), \\ 
\widetilde{W}^{(\Nscript,\Ndscript,\Rscript,\coup)}_{1}
&= 2 \, \big( - \widetilde{W}^{(\Nscript,\Ndscript,\Rscript,\coup),+} + \widetilde{W}^{(\Nscript,\Ndscript,\Rscript,\coup),-}\big),
\end{align*}
where 
\begin{equation*}
\widetilde{W}^{(\Nscript,\Ndscript,\Rscript,\coup),\pm}
= \sum_{L\in \dyadic} \widetilde{S}^{(\Nscript,\Rscript,\coup),\pm}_L P_{\leq \Nd}^x  \Psharp_{R;L} W^{(\Rscript,\coup),\pm} + \widetilde{Z}^{(\Nscript,\Rscript,\coup),\pm}
\end{equation*}
and
\begin{equation*}
\Big( \big( \widetilde{S}^{(\Nscript,\Rscript,\coup),+}_K \big)_{K\in \dyadic},
\big( \widetilde{S}^{(\Nscript,\Rscript,\coup),-}_M \big)_{M\in \dyadic}, 
\widetilde{Z}^{(\Nscript,\Rscript,\coup),+},
\widetilde{Z}^{(\Nscript,\Rscript,\coup),-} \Big) 
\in \Perpm (\Ac,\epsilon). 
\end{equation*}
\end{definition}

In the next definition, 
we introduce an event which captures local well-posedness in Cartesian coordinates
and which is a variant of the event from Definition \ref{main:def-null-lwp}.

\begin{definition}[Well-posedness]\label{main:def-lwp}
Let $N \in \dyadiclarge$, let $R\geq 1$, let $\coup>0$, let $\Ac,\Bc \geq 1$, let $\Dc=\hcoup \Ac \Bc$, and let $\zeta\in C^\infty_b(\R)$.  Then, we define the event
\begin{equation*}
\LWPDzeta \subseteq (\C_x^{s-1}\times \C_x^{s-1})(\bT_R \rightarrow \frkg^2) 
\end{equation*}
as the set of all $(W^{(\Rscript,\coup)}_0,W^{(\Rscript,\coup)}_1)$ for which, under the assumptions \ref{main:item-lwp-A1}-\ref{main:item-lwp-A3}, the conclusions \ref{main:item-lwp-C1}-\ref{main:item-lwp-C5} hold. Here, the assumptions \ref{main:item-lwp-A1}-\ref{main:item-lwp-A3} are defined as follows: 
\begin{enumerate}[leftmargin=7ex,label={[A\arabic*]}]
    \item\label{main:item-lwp-A1} Let $\Nd \in \dyadiclarge$ satisfy $\Nd \leq N^{1-\delta}$ and let $0<\epsilon \leq c_1$. \\[-2.5ex]
    \item\label{main:item-lwp-A2} Let $(\widetilde{W}^{(\Nscript,\Ndscript,\Rscript,\coup)}_{0},\widetilde{W}^{(\Nscript,\Ndscript,\Rscript,\coup)}_{1})\colon \bT_R \rightarrow \frkg^2$ satisfy
    \begin{equation*}
    \big( \widetilde{W}^{(\Nscript,\Ndscript,\Rscript,\coup)}_{0},\widetilde{W}^{(\Nscript,\Ndscript,\Rscript,\coup)}_{1} \big) 
    \in \Per \big( W^{(\Rscript,\coup)}_{0}, W^{(\Rscript,\coup)}_1;\Ac,\epsilon).
    \end{equation*}
    \item\label{main:item-lwp-A3} Let $(A^{(\Nscript,\Rscript,\coup,\zetascript)},B^{(\Nscript,\Rscript,\coup,\zetascript)})$
    and $(\widetilde{A}^{(\Nscript,\Ndscript,\Rscript,\coup,\zetascript)},\widetilde{B}^{(\Nscript,\Ndscript,\Rscript,\coup,\zetascript)})$ be 
    the solutions of \eqref{main:eq-WM-N} with initial data 
    $(W^{(\Rscript,\coup)}_{0}, W^{(\Rscript,\coup)}_1)$ and 
    $( \widetilde{W}^{(\Nscript,\Ndscript,\Rscript,\coup)}_{0},\widetilde{W}^{(\Nscript,\Ndscript,\Rscript,\coup)}_{1})$, 
    respectively. 
\end{enumerate}
Furthermore, with $\epsilon^\prime := C_1 \Bc R^\eta (\epsilon+ \Bc \Nd^{-\varsigma})$, the conclusions \ref{main:item-lwp-C1}-\ref{main:item-lwp-C5} are defined as follows:
\begin{enumerate}[leftmargin=7ex,label={[C\arabic*]}]
\item \label{main:item-lwp-C1} ($C_t^0 \Xi_x^{s-1}$-bound) It holds that
\begin{align*}
\big\| (\widetilde{A}^{(\Nscript,\Ndscript,\Rscript,\coup,\zetascript)},\widetilde{B}^{(\Nscript,\Ndscript,\Rscript,\coup,\zetascript)});(0,0)\big\|_{C_t^0 \Xi_x^{s-1}([-1,1] \times \bT_R \rightarrow \frkg^2)} &\leq C_1 \Dc, \\ 
\big\| (A^{(\Nscript,\Rscript,\coup,\zetascript)},B^{(\Nscript,\Rscript,\coup,\zetascript)});(0,0)\big\|_{C_t^0 \Xi_x^{s-1}([-1,1] \times \bT_R \rightarrow \frkg^2)}
&\leq C_1 \Dc.
\end{align*}
\item \label{main:item-lwp-C2} ($C_t^0 \Xi_x^{s-1}$-difference estimate) It holds that
\begin{equation*}
\big\| (\widetilde{A}^{(\Nscript,\Ndscript,\Rscript,\coup,\zetascript)},\widetilde{B}^{(\Nscript,\Ndscript,\Rscript,\coup,\zetascript)}); 
(A^{(\Nscript,\Rscript,\coup,\zetascript)},B^{(\Nscript,\Rscript,\coup,\zetascript)})
\big\|_{C_t^0 \Xi_x^{s-1}([-1,1] \times \bT_R \rightarrow \frkg^2)} \leq \hcoup \Ac  \epsilon^\prime.
\end{equation*}
\item \label{main:item-lwp-C3} (Difference estimate for lifts) It holds that
\begin{equation*}
\Big\| \ac\big[ \widetilde{A}^{(\Nscript,\Ndscript,\Rscript,\coup,\zetascript)}(\cdot,0)\big] ; 
\ac\big[ A^{(\Nscript,\Rscript,\coup,\zetascript)}(\cdot,0)\big] \Big\|_{C_t^0([-1,1]\rightarrow \frkG)} \leq \hcoup \Ac \epsilon^\prime. 
\end{equation*}
\item \label{main:item-lwp-C4} (Perturbation) For all $t \in [-1,1]$, it holds that 
\begin{equation*}
\big( \widetilde{A}^{(\Nscript,\Ndscript,\Rscript,\coup,\zetascript)}(t),\widetilde{B}^{(\Nscript,\Ndscript,\Rscript,\coup,\zetascript)}(t) \big)
\in 
\Per \big( A^{(\Nscript,\Rscript,\coup,\zetascript)}(t),B^{(\Nscript,\Rscript,\coup,\zetascript)}(t) ;\Ac,\epsilon^\prime).
\end{equation*}
\item \label{main:item-lwp-C5} (Energy increment) It holds that
\begin{equation*}
\sup_{\widetilde{\zeta}\in \Cuttilde} \sup_{t\in [-1,1]} \bigg| \int_{-t}^0 \int_{\bT_R} \ds \dx\,  \widetilde{\zeta}(s)^2\big \langle \Renorm[N] A^{(\Nscript,\Rscript,\coup,\zetascript)}(s,x), B^{(\Nscript,\Rscript,\coup,\zetascript)}(s,x) \big \rangle_\frkg \bigg| \leq C_1 R N^{-(1-10\delta)\delta_1} \Dc^2. 
\end{equation*}
\end{enumerate}
Furthermore, we define
\begin{align}
\LWPDcut &:= \bigcap_{\chi \in \Cut} \LWPDchi.\label{main:eq-lwpd-cut}
\end{align}
Finally, we define $\LWPD$ as the local well-posedness event for $\zeta=1$, i.e., 
\begin{align}
\LWPD &:= \LWPDone.  \label{main:eq-lwpd-one}
\end{align}
\end{definition}

\begin{remark}\label{main:rem-lwp-definition}
In the following, we make several remarks regarding Definition \ref{main:def-lwp}. 
\begin{enumerate}[label=(\alph*)]
\item\label{main:item-rem-cutoff} We emphasize that \ref{main:item-lwp-C1}-\ref{main:item-lwp-C5} only concern the behavior of the solution on $[-1,1]\times \bT_R$ and therefore only depend on the values of $\zeta\in C^\infty_b(\R)$ on $[-1,1]$. As a consequence, it therefore follows for all $\zeta \in C^\infty_b(\R)$ that $\LWPDzeta=\LWPDzetaprime$, where $\zeta^\prime:= \chi \zeta$ and $\chi$ is any smooth function satisfying $\chi(t)=1$ for all $t\in [-1,1]$. 
\item A technical difference between Definition \ref{main:def-null-lwp} and Definition \ref{main:def-lwp} lies in the choice of the $\epsilon^\prime$-parameter. Compared to the $\epsilon^\prime$-parameter from Definition \ref{main:def-null-lwp}, the $\epsilon^\prime$-parameter in Definition \ref{main:def-lwp} contains additional $\Bc$ and $R^{\eta}$-factors. This technical nuisance, which stems from 
 Lemma \ref{main:lem-abstract}, will fortunately not cause any problems below. 
\item Whereas our bound for the distance of $\ac[ \widetilde{A}^{(\Nscript,\Ndscript,\Rscript,\coup,\zetascript)}(\cdot,0)]$ and  
$\ac[ A^{(\Nscript,\Rscript,\coup,\zetascript)}(\cdot,0)]$ is explicitly contained in \ref{main:item-lwp-C3}, our bound for the distance of $\bc[ \widetilde{B}^{(\Nscript,\Ndscript,\Rscript,\coup,\zetascript)}(t,\cdot)]$ and  
$\bc[ B^{(\Nscript,\Rscript,\coup,\zetascript)}(t,\cdot)]$ is not explicitly contained in Definition \ref{main:def-lwp} and will later be derived from Lemma \ref{lifting:lem-lift-properties} and \ref{main:item-lwp-C2}. The technical reason for this is that, if we omitted \ref{main:item-lwp-C3}, then we would need to include local $\Xi^{s-1}_t$-bounds for the distance of $\widetilde{A}^{(\Nscript,\Ndscript,\Rscript,\coup,\zetascript)}$ and $A^{(\Nscript,\Rscript,\coup,\zetascript)}$. In the proof of Lemma \ref{main:lem-gwp-small}, we would then need to glue local $\Xi_t^{s-1}$-bounds together, which would be technically inconvenient.
\item In Definition \ref{main:def-null-lwp}, we only consider solutions of the finite-dimensional approximation of the wave maps equation \eqref{main:eq-WM-N}. In order to prove Theorem \ref{intro:thm-rigorous-A-B}, we also need to compare solutions of \eqref{main:eq-WM-N} with solutions of the original wave maps equation \eqref{main:eq-WM}. This will be done 
using classical estimates (Proposition \ref{classical:prop-main}), but is postponed until the last part of our argument (Proposition \ref{main:prop-refined-gwp}). 
\end{enumerate}
\end{remark}

Equipped with Definition \ref{main:def-lwp}, we can now state the main result of this subsection.

\begin{proposition}[Short-time and local well-posedness]\label{main:prop-lwp} 
Let $N \in \Dyadiclarge$, let $R \geq 1$, let $\coup>0$, let $\Ac,\Bc \geq C_1$, and let $\Dc=\hcoup \Ac \Bc$. Then, we have the following two estimates: 
\begin{enumerate}[label=(\Roman*)]
\item \label{main:item-lwp-short} 
If $\Dc \leq c_1 N^{-\delta_2}$, then
\begin{equation*}
\muR \Big(  \LWPDcut \Big) \geq 1- c_1^{-1} \exp\Big( - c_1 R^{-2\eta} \Ac^2 \Big). 
\end{equation*}
\item \label{main:item-lwp-local} If $\Dc \leq c_1$ and $\chi=\chifixed$ is as in Definition \ref{prelim:def-cut-off}, then
\begin{equation}\label{main:eq-lwp-probability}
\muR \Big(  \LWPD\Big) \geq \muR \Big( \BBAchi \Big) - c_1^{-1} \exp\Big( - c_1 R^{-2\eta} \Ac^2 \Big).
\end{equation}
\end{enumerate}
\end{proposition} 

The idea behind the proof of Proposition \ref{main:prop-lwp} is to first use the well-posedness in null-coordinates (Proposition \ref{main:prop-null-lwp}) and then derive the conclusions in Definition \ref{main:def-lwp} from the conclusions in Definition \ref{main:def-null-lwp}.

\begin{proof}[Proof of Proposition \ref{main:prop-lwp}:] 
For expository purposes, we split the argument into two steps.\\

\emph{First step: Using Proposition \ref{main:prop-null-lwp}.}
We define $W^{(\Rscript,\coup),+},W^{(\Rscript,\coup),-}\colon \bT_R \rightarrow \frkg$ as 
\begin{equation}\label{main:eq-lwp-proof-0}
W^{(\Rscript,\coup),\pm} := \tfrac{1}{4} \big( W_0^{(\Rscript,\coup)} \mp W_1^{(\Rscript,\coup)} \big),
\end{equation}
which are $2\pi R$-periodic, $\frkg$-valued white noises at temperature $\coup$. Using our representation from Lemma \ref{prelim:lem-white-noise-representation}, we may then write 
\begin{equation}\label{main:eq-lwp-proof-1}
W^{(\Rscript,\coup),\pm}(x) = \hcoup \sum_{x_0 \in \LambdaRR} \sum_{\ell \in \Z} \psiRx G_{x_0,\ell}^{\pm} e^{i \ell x}. 
\end{equation}
Similar as in \eqref{main:eq-lwp-proof-0}, we define $\widetilde{W}^{(\Nscript,\Ndscript,\Rscript,\coup),+},\widetilde{W}^{(\Nscript,\Ndscript,\Rscript,\coup),-}\colon \bT_R \rightarrow \frkg$ as 
\begin{equation*}
\widetilde{W}^{(\Nscript,\Ndscript,\Rscript,\coup),+}
:= \tfrac{1}{4} \big(  \widetilde{W}^{(\Nscript,\Ndscript,\Rscript,\coup)}_0 \mp  \widetilde{W}^{(\Nscript,\Ndscript,\Rscript,\coup)}_1 \big). 
\end{equation*}
Due to Assumption \ref{main:item-lwp-A2} in Definition \ref{main:def-lwp}, we may then write
\begin{equation}\label{main:eq-lwp-proof-2}
\widetilde{W}^{(\Nscript,\Ndscript,\Rscript,\coup),\pm}=   \sum_{L\in \dyadic}   \widetilde{S}^{(\Nscript,\Rscript,\coup),\textup{in},\pm}_L P_{\leq \Nd}^x \Psharp_{R;L} W^{(\Rscript,\coup),\pm} +\widetilde{Z}^{(\Nscript,\Rscript,\coup),\pm},
\end{equation}
where 
\begin{equation*}
\Big( \big(\widetilde{S}^{(\Nscript,\Rscript,\coup),\textup{in},+}_{K} \big)_{K\in \dyadic}, 
\big(\widetilde{S}^{(\Nscript,\Rscript,\coup),\textup{in},-}_{M} \big)_{M\in \dyadic}, 
\widetilde{Z}^{(\Nscript,\Rscript,\coup),+},
\widetilde{Z}^{(\Nscript,\Rscript,\coup),-} \Big) \in \Perpm \big( \Ac, \epsilon \big).
\end{equation*}
We now let $(U^{(\Nscript,\Rscript,\coup,\chiscript)},V^{(\Nscript,\Rscript,\coup,\chiscript)})$ 
and $(\widetilde{U}^{(\Nscript,\Ndscript,\Rscript,\coup,\chiscript)},\widetilde{V}^{(\Nscript,\Ndscript,\Rscript,\coup,\chiscript)})$
be the solutions of \eqref{main:eq-WM-null-N} with the initial data from
\eqref{main:eq-lwp-proof-1} and \eqref{main:eq-lwp-proof-2}, respectively. Then, we can write the solutions from Assumption \ref{main:item-lwp-A3} in Definition \ref{main:def-lwp} as 
\begin{align}
\Big( A^{(\Nscript,\Rscript,\coup,\chiscript)}, B^{(\Nscript,\Rscript,\coup,\chiscript)} \Big) 
&= 2\,  \Big( U^{(\Nscript,\Rscript,\coup,\chiscript)} + V^{(\Nscript,\Rscript,\coup,\chiscript)}, 
- U^{(\Nscript,\Rscript,\coup,\chiscript)} + V^{(\Nscript,\Rscript,\coup,\chiscript)} \Big), 
\label{main:eq-lwp-proof-3} \\
\Big( \widetilde{A}^{(\Nscript,\Ndscript,\Rscript,\coup,\chiscript)}, \widetilde{B}^{(\Nscript,\Ndscript,\Rscript,\coup,\chiscript)} \Big) 
&= 2 \, \Big( U^{(\Nscript,\Ndscript,\Rscript,\coup,\chiscript)} + V^{(\Nscript,\Ndscript,\Rscript,\coup,\chiscript)}, 
- U^{(\Nscript,\Ndscript,\Rscript,\coup,\chiscript)} + V^{(\Nscript,\Ndscript,\Rscript,\coup,\chiscript)} \Big). 
\label{main:eq-lwp-proof-4} 
\end{align}
We now note that, due to Remark \ref{main:rem-lwp-definition}.\ref{main:item-rem-cutoff}, it holds for $\chi=\widebar{\chi}$ as in \ref{main:item-lwp-local} that 
\begin{equation*}
\LWPD = \LWPDchi.
\end{equation*}
Thus, \eqref{main:eq-lwp-probability} can be converted into a form similar as in \eqref{main:eq-null-lwp-probability}, but with $\LWPDchi$ instead of $\nullLWPDchi$. Due to Proposition \ref{main:prop-null-lwp}, the events $\nullLWPDcut$ and $\nullLWPDchi$ satisfy the probability estimates in \ref{main:item-lwp-short} and \ref{main:item-lwp-local}. Thus, we may therefore utilize the pure modulation operators and nonlinear remainders 
\begin{align*}
&\Big( S^{(\Nscript,\Rscript,\coup,\chiscript),\diamond,+} ,S^{(\Nscript,\Rscript,\coup,\chiscript),\diamond,-}, U^{(\Nscript,\Rscript,\coup,\chiscript),\fs}, V^{(\Nscript,\Rscript,\coup,\chiscript),\fs} \Big) \\ 
\text{and} \qquad &\Big( \widetilde{S}^{(\Nscript,\Ndscript,\Rscript,\coup,\chiscript),\diamond,+}, \widetilde{S}^{(\Nscript,\Ndscript,\Rscript,\coup,\chiscript),\diamond,-}, \widetilde{U}^{(\Nscript,\Ndscript,\Rscript,\coup,\chiscript),\fs}, \widetilde{V}^{(\Nscript,\Ndscript,\Rscript,\coup,\chiscript),\fs}\Big)
\end{align*}
from Definition \ref{main:def-null-lwp} and may write our solutions  $(U^{(\Nscript,\Rscript,\coup,\chiscript)},V^{(\Nscript,\Rscript,\coup,\chiscript)})$ 
 and $(U^{(\Nscript,\Ndscript,\Rscript,\coup,\chiscript)},V^{(\Nscript,\Ndscript,\Rscript,\coup,\chiscript)})$ as
\begin{alignat*}{3}
U^{(\Nscript,\Rscript,\coup,\chiscript)} &= \sum_{\ast} U^{(\Nscript,\Rscript,\coup,\chiscript),\ast}, 
\quad&\quad 
V^{(\Nscript,\Rscript,\coup,\chiscript)} &= \sum_{\ast} V^{(\Nscript,\Rscript,\coup,\chiscript),\ast}, \\ 
U^{(\Nscript,\Ndscript,\Rscript,\coup,\chiscript)} &= \sum_{\ast} U^{(\Nscript,\Ndscript,\Rscript,\coup,\chiscript),\ast}, 
\quad&\quad 
V^{(\Nscript,\Ndscript,\Rscript,\coup,\chiscript)} &= \sum_{\ast} V^{(\Nscript,\Ndscript,\Rscript,\coup,\chiscript),\ast},
\end{alignat*}
where $\ast \in \{ +,+-,-,+\fs,\fs-,\fs\}$. 
It remains to show that 
the conclusions in Definition \ref{main:def-null-lwp} can be used to derive the conclusions in Definition \ref{main:def-lwp}. \\

\emph{Second step: The conclusions in Definition \ref{main:def-lwp}.} 
We now separately address the conclusions from \ref{main:item-lwp-C1}-\ref{main:item-lwp-C5} in Definition \ref{main:def-lwp}. \\

\emph{\ref{main:item-lwp-C1}: $C_t^0 \Xi_x^{s-1}$-bound.} This follows directly from Definition \ref{main:def-enhanced-data}, Lemma \ref{modulation:lem-Cartesian}, and Lemma \ref{modulation:lem-Cartesian-high-high}. In fact, due to the cut-off function in time, the estimate holds with $[-1,1]\times \bT_R$ replaced by $\mathbb{R}\times \bT_R$. \\ 

\emph{\ref{main:item-lwp-C2}: $C_t^0 \Xi_x^{s-1}$-difference estimate.} The argument is similar as for \ref{main:item-lwp-C1}. The only difference is that, instead of Lemma \ref{modulation:lem-Cartesian} and Lemma \ref{modulation:lem-Cartesian-high-high}, we use their Lipschitz-variants (from Lemma \ref{lipschitz:lem-modulated-cartesian} and  Remark \ref{lipschitz:rem-cartesian}), which yield that 
\begin{equation*}
\big\| (\widetilde{A}^{(\Nscript,\Ndscript,\Rscript,\coup,\zetascript)},\widetilde{B}^{(\Nscript,\Ndscript,\Rscript,\coup,\zetascript)}); 
(A^{(\Nscript,\Rscript,\coup,\zetascript)},B^{(\Nscript,\Rscript,\coup,\zetascript)})
\big\|_{C_t^0 \Xi_x^{s-1}([-1,1] \times \bT_R \rightarrow \frkg^2)} 
\lesssim C_0 \hcoup \Ac \big( \epsilon + \Bc \Nd^{-\varsigma} \big). 
\end{equation*}

\emph{\ref{main:item-lwp-C3}: Difference estimate for lifts.}
As discussed in Remark \ref{modulation:rem-reversed-roles}, similar estimates as in Lemma \ref{modulation:lem-Cartesian} and Lemma \ref{modulation:lem-Cartesian-high-high} are still satisfied after switching the $t$ and $x$-variables. Similar as in \ref{main:item-lwp-C2}, we then obtain that 
\begin{equation*}
\big\| (\widetilde{A}^{(\Nscript,\Ndscript,\Rscript,\coup,\chiscript)},\widetilde{B}^{(\Nscript,\Ndscript,\Rscript,\coup,\chiscript)}); 
(A^{(\Nscript,\Rscript,\coup,\chiscript)},B^{(\Nscript,\Rscript,\coup,\chiscript)})
\big\|_{C_x^0 \Xi_t^{s-1}( \bT_R \times \R \rightarrow \frkg^2)} \lesssim C_0 \hcoup \Ac \big( \epsilon + \Bc \Nd^{-\varsigma} \big). 
\end{equation*}
In particular, it holds that 
\begin{equation*}
\big\| \widetilde{A}^{(\Nscript,\Ndscript,\Rscript,\coup,\chiscript)}(\cdot,0); 
A^{(\Nscript,\Rscript,\coup,\chiscript)}(\cdot,0)
\big\|_{C_x^0 \Xi_t^{s-1}(  \R \rightarrow \frkg)} \lesssim C_0 \hcoup \Ac \big( \epsilon + \Bc \Nd^{-\varsigma} \big).
\end{equation*}
Together with Lemma \ref{lifting:lem-lift-properties}, this implies the desired estimates for the lifts. \\ 

\emph{\ref{main:item-lwp-C4}: Perturbation.} Similar as in the proof of Proposition \ref{main:prop-null-lwp}, we now replace the super-scripts 
$(N,R,\coup,\chi)$ and $(N,\Nd,R,\coup,\chi)$ by ``$\fn$" and ``$\fd$", respectively. Using our Ansatz \eqref{ansatz:eq-UN-rigorous-decomposition}, Definition \ref{ansatz:def-modulated-linear}, and Lemma \ref{ansatz:lem-frequency-support}, we can write 
\begin{equation}\label{main:eq-lwp-p1}
\begin{aligned}
P_{\leq \Nd}^x U^{\fn} (x-t,x+t) 
&= P_{\leq \Nd}^x \Sumlarge_{K\lesssim \Nd} \sum_{u_0 \in \LambdaRR} \sum_{k\in \Z_K} \rhoND(k) \psiRuK(x-t) S^{\fn,+}_{K,k}(x-t,x+t)  G^+_{u_0,K} e^{ik(x-t)} \\ 
&+ P_{\leq \Nd}^x \big( U^{\fn,+-} + U^{\fn,-} + U^{\fn,+\fs} + U^{\fn,\fs-} + U^{\fn,\fs} \big)(x-t,x+t) 
\end{aligned}
\end{equation}
and 
\begin{equation}\label{main:eq-lwp-p2}
\begin{aligned}
 U^{\fd} (x-t,x+t) 
&=  \Sumlarge_{K\lesssim \Nd} \sum_{u_0 \in \LambdaRR} \sum_{k\in \Z_K} \psiRuK(x-t) S^{\fd,+}_{K,k}(x-t,x+t)  G^+_{u_0,K} e^{ik(x-t)} \\ 
&+  \big( U^{\fd,+-} + U^{\fd,-} + U^{\fd,+\fs} + U^{\fd,\fs-} + U^{\fd,\fs} \big)(x-t,x+t).
\end{aligned}
\end{equation}
In order to apply Lemma \ref{main:lem-abstract} to \eqref{main:eq-lwp-p1} and \eqref{main:eq-lwp-p2}, we need to verify that the assumptions in Lemma \ref{main:lem-abstract} are satisfied. 
Due to the condition $\Nd\leq N^{1-\delta}$ from Assumption \ref{main:item-lwp-A1} and Proposition \ref{modulation:prop-properties}.\ref{modulation:item-properties-3}, the $S^{\fn,+}_{K,k}$ and $S^{\fd,+}_{K,k}$-operators appearing in \eqref{main:eq-lwp-p1} and \eqref{main:eq-lwp-p2} are constant in $k\in \Z_K$, and we therefore now omit the subscript $k$. 
Due to Proposition \ref{modulation:prop-main}.\ref{modulation:item-distance-initial}, we have for all $K\lesssim \Nd$ that 
\begin{equation}\label{main:eq-lwp-p3} 
\big\| S^{\fn,+}_{K} - \Id_\frkg \big\|_{\Cprod{s}{s}} \lesssim \Dc 
\qquad \text{and} \qquad \big\| S^{\fd,+}_{K} - \Id_\frkg \big\|_{\Cprod{s}{s}} 
\lesssim \Dc  + \epsilon. 
\end{equation}
Since $\Dc \leq c_1$ and $\epsilon\leq c_1$, it then follows from Lemma \ref{prelim:lem-trace} and \eqref{main:eq-lwp-p3} that  
\begin{equation}\label{main:eq-lwp-p4}
\big\| S^{\fn,+}_{K}\big\|_{C_t^0 \C_x^s}, 
\big\| (S^{\fn,+}_{K})^{-1} \big\|_{C_t^0 \C_x^s},
\big\| S^{\fd,+}_{K}\big\|_{C_t^0 \C_x^s}, 
\big\| (S^{\fd,+}_{K})^{-1} \big\|_{C_t^0 \C_x^s} \leq 2.
\end{equation}
Furthermore, due to (a minor variant of) Lemma \ref{modulation:lem-Cartesian} and Definition \ref{main:def-null-lwp}.\ref{main:item-null-lwp-C3}, it holds that
\begin{align}
\Big\| U^{\fn,+-} + U^{\fn,-} + U^{\fn,+\fs} + U^{\fn,\fs-} + U^{\fn,\fs} \Big\|_{C_t^0 \C_x^{r+\eta-1}}
&\lesssim \Dc, 
\label{main:eq-lwp-p5} \\ 
\Big\| U^{\fd,+-} + U^{\fd,-} + U^{\fd,+\fs} + U^{\fd,\fs-} + U^{\fd,\fs}  \Big\|_{C_t^0 \C_x^{r-1}} &\lesssim \Dc. 
\label{main:eq-lwp-p6}
\end{align}

Due to \eqref{main:eq-lwp-p4}, \eqref{main:eq-lwp-p5}, and \eqref{main:eq-lwp-p6}, we can 
apply\footnote{We remark that Lemma \ref{main:lem-abstract} contains sums over $\dyadic$, whereas \eqref{main:eq-lwp-p1} and \eqref{main:eq-lwp-p2} contain sums over $\Dyadiclarge$. This causes no problem, since terms corresponding to $K\in \dyadic\backslash\Dyadiclarge$ can always be absorbed into the $Z$-terms in Lemma \ref{main:lem-abstract}.} Lemma \ref{main:lem-abstract} to (translations of) $P_{\leq \Nd}^x U^{\fn} $ and $U^{\fd}$. From Lemma \ref{main:lem-abstract}, we then obtain that
\begin{equation*}
U^{\fd}(x-t,x+t) = \sum_K S^+_K(x-t,x+t) P_{\leq \Nd}^x \Psharp_{R;K} U^{\fn}(x-t,x+t) + Z^+(x-t,x+t), 
\end{equation*}
where the modulation operators and remainder satisfy
\begin{equation*}
\big\| S_K^+ - \Id_\frkg \big\|_{C_t^0 \C_x^s} \lesssim C_0^2 \Bc R^\eta\big( \epsilon + \Bc \Nd^{-\varsigma} \big) 
\qquad \text{and} \qquad 
\big\| Z^+ \big\|_{C_t^0 \C_x^{r-1}} \lesssim \hcoup \Ac \times C_0^2 \Bc R^\eta\big( \epsilon + \Bc \Nd^{-\varsigma} \big) .
\end{equation*}
Together with a similar argument for $V^{\fn}$ and $V^{\fd}$, it then follows for all $t\in \R$ that 
\begin{equation*}
\big( A^{\fd}(t),B^{\fd}(t) \big)
\in 
\Per \Big( A^{\fn}(t),B^{\fn}(t) ;\Ac, C_1 \Bc R^\eta\big( \epsilon + \Bc \Nd^{-\varsigma} \big)\Big),
\end{equation*}
which completes the proof of \ref{main:item-lwp-C4}. \\

\emph{\ref{main:item-lwp-C5}: Energy increment.} This follows directly from Proposition \ref{increment:prop-energy-increment}, i.e., our earlier estimate of the energy increment.
\end{proof}

\subsection{Almost invariance}\label{section:main-invariance}
In this subsection, we prove the almost invariance of the Gibbs measure under our finite-dimensional approximation of the wave maps equation \eqref{main:eq-WM-N}. To this end, we first recall from Definition \ref{prelim:def-Gibbs} and Definition \ref{structure:def-Gibbs-time-evolved} that the Gibbs measure and its time-evolution under 
\eqref{main:eq-WM-N} are denoted by 
$\muR$ and $\muNR_{\zeta,t}$, respectively. 

\begin{proposition}[Almost invariance]\label{main:prop-almost-invariance} 
Let $N\in \Dyadiclarge$, let $T,R\geq 1$,  let $\coup>0$, and let $\zeta \in \Cuttilde$. If  $\max(T,R,\lambda) \leq N^{\delta_2}$, then it holds that 
\begin{equation}\label{main:eq-almost-invariance-TV}
\sup_{t\in [0,T]} \Big\| \muNR_{\zeta,t} - \muR \Big\|_{\TV} \leq C_2 N^{-\frac{\delta_1}{2}}. 
\end{equation}
Furthermore, it holds for all events $E\subseteq \C_x^{s-1}(\bT_R \rightarrow \frkg^2)$ that 
\begin{equation}\label{main:eq-almost-invariance-event}
\sup_{t\in [0,T]} \muNR_{\zeta,t} ( E ) \leq C_2 \Big( \muR(E) + \exp \big( - c_2 N^{\delta_2}\big) \Big). 
\end{equation}
\end{proposition}

Before we turn to the proof of Proposition \ref{main:prop-almost-invariance}, we record the following scaling identity which involves the scaling transform from \eqref{ansatz:eq-scaling-transform}. 

\begin{lemma}[Scaling]\label{main:lem-scaling}
Let $(N,R,\coup),(N^\prime,R^\prime,\coup^\prime)\in \Dyadiclarge \times [1,\infty)\times (0,\infty)$, let $\zeta,\zeta^\prime \in C^\infty_b(\R)$ and assume that, for some $\kappa \in 2^{\Z}$, 
\begin{equation*}
N^\prime = \kappa N, \qquad R^\prime = \frac{R}{\kappa}, \quad 
\coup^\prime = \kappa \coup, \qquad \text{and} \qquad  \zeta^\prime(\cdot) = \zeta(\kappa \, \cdot).
\end{equation*}
Then, it holds for all $t \in \R$ that 
\begin{equation*}
\muNR_{\zeta,t} = ( \Scaling_{\kappa^{-1}})_\# \mupNR_{\zeta^\prime,\kappa^{-1}t}. 
\end{equation*}
\end{lemma}

\begin{proof}
This follows directly from Lemma \ref{ansatz:lem-scaling-symmetry}.
\end{proof}

\begin{proof}[Proof of Proposition \ref{main:prop-almost-invariance}:]
For expository purposes, we will soon split the argument into four steps. In the first three steps, we only treat the 
low-temperature case
\begin{equation}\label{main:eq-inv-low-temperature}
\coup \leq N^{-4\delta_2}. 
\end{equation}
In the last step, the condition \eqref{main:eq-inv-low-temperature} will be removed via a scaling argument. Throughout the four steps, 
we can always assume that $N$ is sufficiently large depending on $C_1$, $c_1$, and $(\delta_j)_{j=0}^5$, i.e., 
\begin{equation}\label{main:eq-large-N}
N \gg_{C_1,c_1,\delta_\ast} 1, 
\end{equation}
since otherwise \eqref{main:eq-almost-invariance-TV} and \eqref{main:eq-almost-invariance-event} are trivial. Furthermore, we use the following notation:
We let $C_1^\prime=C_1^\prime(C_1,\delta_\ast)$ and $C_1^{\prime\prime}=C_1^{\prime\prime}(C_1^\prime,\delta_\ast)$ be sufficiently large constants and let $c_1^\prime=c_1^\prime(c_1,\delta_\ast)$ and $c_1^{\prime\prime}=c_1^{\prime\prime}(c_1^\prime,\delta_\ast)$
be sufficiently small constants. Furthermore, we write $\Wn=(W_0^{(\Rscript,\coup)},W_1^{(\Rscript,\coup)})$ and let
\begin{equation*}
\Big( A^{(\Nscript,\Rscript,\coup,\zetascript)}(t,x;t_0,\Wn), B^{(\Nscript,\Rscript,\coup,\zetascript)}(t,x;t_0,\Wn) \Big)
\end{equation*}
be the solution of \eqref{main:eq-WM-N} with initial time $t_0$ and initial data $\Wn$. Finally, we let $\Phin_\zeta(t;t_0)$ be the corresponding flow. \\ 

\emph{First step: Upper bound at low temperature on the timescale $\sim 1$.}
Under the assumption \eqref{main:eq-inv-low-temperature}, we prove for all $R\geq 1$, all events $E\subseteq \Omegas$,  and all $\zeta \in \Cuttilde$ that
\begin{equation}\label{main:eq-inv-1-p1}
\sup_{\tau \in [0,1]} \muNR_{\zeta,\tau}(E) 
\leq \exp \Big( C_1^\prime R N^{-\frac{3\delta_1}{4}} \Big) \muR(E) + C_1^\prime \exp\Big( -c_1^\prime R^{-2\eta} N^{2\delta_2}\Big).
\end{equation}
In order to prove \eqref{main:eq-inv-1-p1}, we first let $\tau \in [0,1]$ and define 
\begin{align}
&\EInc \label{main:eq-inv-1-p2} \\[1ex]
:=&\,  \bigg\{\Wn \in \Omegas \colon 
\bigg| \int_{0}^{\tau} \ds \zeta(s)^2 \big \langle \Renorm[N] A^{(\Nscript,\Rscript,\coup,\zetascript)}(s,x;\tau,\Wn),
B^{(\Nscript,\Rscript,\coup,\zetascript)}(s,x;\tau,\Wn)\big \rangle \bigg| \leq C_1^\prime R N^{-\frac{3\delta_1}{4}} \bigg\}. \notag 
\end{align} 
Since 
\begin{equation*}
\big( A^{(\Nscript,\Rscript,\coup,\zetascript)}(s,x;\tau,\zeta,\Wn), B^{(\Nscript,\Rscript,\coup,\zetascript)}(s,x;\tau,\Wn) \big) 
= \Big( \Phin_\zeta (s;\tau) \Wn \Big)(x)
\end{equation*}
and $\Phin_\zeta(s;\tau)\Phin_\zeta(\tau;0)=\Phin_\zeta(s;0)$, we then obtain that
\begin{align}
&\, \Big\{ \Wn  \in \Omegas \colon \Phin_\zeta(\tau;0) \Wn \in \EInc \Big\} \label{main:eq-inv-1-q0} \\ 
=&\, \bigg\{ \Wn \in \Omegas \colon 
\bigg| \int_{0}^{\tau} \ds \zeta(s)^2 \big \langle \Renorm[N] A^{(\Nscript,\Rscript,\coup,\zetascript)}(s,x;0,\Wn),
B^{(\Nscript,\Rscript,\zetascript,\coup)}(s,x;0,\Wn)\big \rangle \bigg| \leq C_1^\prime R N^{-\frac{3\delta_1}{4}} \bigg\}. \notag 
\end{align}
We now choose parameters 
\begin{equation}\label{main:eq-inv-1-q1}
\Ac := \tfrac{c_1}{C_1} \coup^{-\frac{1}{2}} N^{-\delta_2} \qquad \text{and} \qquad \Bc:=C_1,
\end{equation}
where $c_1$ is as in Proposition \ref{main:prop-lwp}. Due to \eqref{main:eq-inv-low-temperature} and \eqref{main:eq-large-N}, it holds that $\Ac\geq C_1$, which is one of the assumptions in Proposition \ref{main:prop-lwp}. Due to \eqref{main:eq-inv-1-q1}, it then holds that
\begin{equation}\label{main:eq-inv-1-q2}
\Dc = \hcoup \Ac \Bc = \hcoup \times \tfrac{c_1}{C_1} \coup^{-\frac{1}{2}} N^{-\delta_2} \times C_1 = c_1 N^{-\delta_2}.
\end{equation}
Using Definition \ref{main:def-lwp}, Remark \ref{main:rem-lwp-definition}.\ref{main:item-rem-cutoff}, \eqref{main:eq-inv-1-q0}, and $\Dc \leq 1$, it then follows that  
\begin{equation}\label{main:eq-inv-1-p3}
\LWPDcut \subseteq \Phin_\zeta(\tau;0)^{-1} \EInc.
\end{equation}
Equipped with both \eqref{main:eq-inv-1-p2} and \eqref{main:eq-inv-1-p3}, we now estimate
\begin{equation}\label{main:eq-inv-1-p4}
\muNR_{\zeta,\tau}\big( E \big) 
\leq \muNR_{\zeta,\tau} \Big( E \medcap \EInc \Big) + \muNR_{\zeta,\tau} \Big( \Omegas \backslash \EInc \Big).
\end{equation}
To control the first summand in \eqref{main:eq-inv-1-p4}, we use Proposition \ref{structure:prop-Gibbs} and \eqref{main:eq-inv-1-p2}, which yield that
\begin{align*}
\muNR_{\zeta,\tau} \Big( E \medcap \EInc \Big) 
&= \int \mathbf{1} \Big\{ \Wn \in E \medcap \EInc \Big\}\,  \frac{\mathrm{d}\muNR_{\zeta,\tau}}{\mathrm{d}\muR}(\Wn) \,  \mathrm{d}\muR (\Wn) \\
&\leq  \int \mathbf{1} \Big\{ \Wn \in E \medcap \EInc \Big\} \exp\Big( C_1^\prime R N^{-\frac{3\delta_1}{4}} \Big) \mathrm{d}\muR (\Wn) \\
&= \exp\Big( C_1^\prime R N^{-\frac{3\delta_1}{4}} \Big) \muR\Big( E \medcap \EInc  \Big) \\ 
&\leq \exp\Big( C_1^\prime R N^{-\frac{3\delta_1}{4}} \Big) \muR\big( E  \big). 
\end{align*}
Using  \eqref{main:eq-inv-1-q2}, \eqref{main:eq-inv-1-p3}, and Proposition \ref{main:prop-lwp}.\ref{main:item-lwp-short}, the second summand in \eqref{main:eq-inv-1-p4} can be estimated by 
\begin{equation*}
\muNR_{\zeta,\tau} \Big( \Omegas \backslash \EInc \Big)
\leq \muNR \big( \Omegas \backslash \LWPDcut \big) \leq c_1^{-1} \exp \big( - c_1 R^{-2\eta} \Ac^2 \big).
\end{equation*}
After inserting \eqref{main:eq-inv-1-q1} and using the low-temperature condition \eqref{main:eq-inv-low-temperature}, this proves \eqref{main:eq-inv-1-p1}.\\

\emph{Second step: Upper bound at low temperature.}
Under the condition \eqref{main:eq-inv-low-temperature}, we now prove for all $T,R\geq 1$, all events $E\subseteq \Omegas$, and all $\zeta\in \Cuttilde$ that 
\begin{equation}\label{main:eq-inv-2-p1}
\begin{aligned}
\sup_{t\in [0,T]} \muNR_{\zeta,t}\big( E \big) 
&\leq \exp \Big( C_1^\prime (1+T) R N^{-\frac{3\delta_1}{4}} \Big) \muR\big( E \big) \\ 
&+ C_1^\prime (1+T) \exp \Big( C_1^\prime (1+T) R N^{-\frac{3\delta_1}{4}} \Big) \exp\big( -c_1^\prime R^{-2\eta} N^{2\delta_2}\big).
\end{aligned}
\end{equation}
For this, it suffices to prove for all $j\in \mathbb{N}$ that
\begin{equation}\label{main:eq-inv-2-p2}
\begin{aligned}
\sup_{t\in [0,j]} \muNR_{\zeta,t}\big( E \big) 
\leq \exp \Big( C_1^\prime j R N^{-\frac{3\delta_1}{4}} \Big) \muR\big( E \big)  
+ C_1^\prime j \exp \Big( C_1^\prime j R N^{-\frac{3\delta_1}{4}} \Big) \exp\big( -c_1^\prime R^{-2\eta} N^{2\delta_2}\big).
\end{aligned}
\end{equation}
In order to prove \eqref{main:eq-inv-2-p2}, we proceed via induction. The base case $j=1$ follows directly from \eqref{main:eq-inv-1-p1}, and it thus remains to perform the induction step. To this end, let $t\in [j,j+1]$. Then, it holds that 
\begin{equation}\label{main:eq-inv-2-p3}
\begin{aligned}
\muNR_{\zeta,t}\big( E \big) 
&= \muR \big( \Phin_{\zeta}(t;0) \Wn \in E \big) \\ 
&= \muR \big( \Phin_{\zeta}(t;1) \Phin_{\zeta}(1;0) \Wn \in E \big) \\
&= \muNR_{\zeta,1} \big( \Phin_{\zeta}(t;1) \Wn \in E \big). 
\end{aligned}
\end{equation}
Using the base case, i.e., \eqref{main:eq-inv-1-p1}, it then follows that
\begin{equation}\label{main:eq-inv-2-p3p}
\muNR_{\zeta,1} \big( \Phin_{\zeta}(t;1) \Wn \in E \big) 
\leq \exp\Big( C_1^\prime R N^{-\frac{3\delta_1}{4}} \Big) \muR\big(  \Phin_{\zeta}(t;1) \Wn \in E \big)
+ C_1^\prime \exp \big( - c_1^\prime R^{-2\eta} N^{2\delta_2} \big).
\end{equation}
As is clear from \eqref{main:eq-WM-N}, translation in time yields the identity
\begin{equation}\label{main:eq-inv-2-p4} 
\Phin_{\zeta}(t;1) = \Phin_{\Theta^t_{-1}\zeta} (t-1;0).
\end{equation}
Since the set $\Cuttilde$ is invariant under translation (see Definition \ref{prelim:def-cut-off}), it then follows from 
\eqref{main:eq-inv-2-p4} and the induction hypothesis that
\begin{equation}\label{main:eq-inv-2-p5}
\begin{aligned}
&\muR\big(  \Phin_{\zeta}(t;1) \Wn \in E \big)
= \muR \big( \Phin_{\Theta^t_{-1}\zeta} (t-1;0) \Wn \in E \big) 
= \muNR_{\Theta^t_{-1}\zeta,t-1}\big( E \big) \\
\leq&\, \exp \Big( C_1^\prime j R N^{-\frac{3\delta_1}{4}} \Big) \muR\big( E \big)  
+ C_1^\prime j \exp \Big( C_1^\prime j R N^{-\frac{3\delta_1}{4}} \Big) \exp\big( -c_1^\prime R^{-2\eta} N^{2\delta_2}\big). 
\end{aligned}
\end{equation}
By combining \eqref{main:eq-inv-2-p3}, \eqref{main:eq-inv-2-p3p}, and \eqref{main:eq-inv-2-p5}, we then obtain that
\begin{align*}
\muNR_{\zeta,t}\big( E \big)  
\leq \hspace{-0.25ex} \exp \Big( C_1^\prime (j+1) R N^{-\frac{3\delta_1}{4}} \Big) \muR\big( E \big)
+ \hspace{-0.25ex} C_1^\prime (j+1) \exp \Big( C_1^\prime (j+1) R N^{-\frac{3\delta_1}{4}} \Big) \hspace{-0.25ex} \exp\big( -c_1^\prime R^{-2\eta} N^{2\delta_2}\big),
\end{align*}
which completes the induction step.\\

\emph{Third step: The desired estimates at low temperature.}
Under the low-temperature condition \eqref{main:eq-inv-low-temperature}, we now prove for all $T,R \leq  N^{10\delta_2}$, all $E\subseteq \Omegas$, all $\zeta\in \Cuttilde$, and all $0\leq t \in [0,T]$ that 
\begin{equation}\label{main:eq-inv-3-p1}
\muNR_{\zeta,t}\big( E \big) \leq C_1^{\prime\prime} \Big( \muR\big(E\big) 
+ \exp\big( - c_1^{\prime\prime} N^{2(\delta_2-\eta)} \big) \Big)
\end{equation}
and 
\begin{equation}\label{main:eq-inv-3-p2}
\big\| \muNR_{\zeta,t} - \muR \big\|_{\TV} \leq C_1^{\prime\prime} N^{-\frac{5\delta_1}{8}}. 
\end{equation}
To obtain \eqref{main:eq-inv-3-p1} and \eqref{main:eq-inv-3-p2}, we first use \eqref{main:eq-inv-2-p1} and the conditions $ T,R \leq  N^{10\delta_2}$, which directly yield
\begin{equation}\label{main:eq-inv-3-p3}
\begin{aligned}
\muNR_{\zeta,t} \big( E \big) 
&\leq \exp \Big( C_1^\prime (1+N^{10\delta_2}) N^{10\delta_2} N^{-\frac{3\delta_1}{4}} \Big) \muR\big( E \big) \\
&+ C_1^\prime \big( 1+N^{10\delta_2} \big)  \exp\Big( C_1^\prime (1+N^{10\delta_2}) N^{10\delta_2} N^{-\frac{3\delta_1}{4}} \Big) \exp\Big( -c_1^\prime N^{-20\delta_2 \eta} N^{2\delta_2} \Big) \\ 
&\leq \exp\Big(  C_1^\prime N^{-\frac{5\delta_1}{8}}\Big) \muR\big( E \big) 
+ C_1^\prime \exp\big(  C_1^\prime \big) N^{20\delta_2} \exp\Big( -c_1^\prime N^{2(\delta_2-\eta)} \Big).
\end{aligned}
\end{equation}
From \eqref{main:eq-inv-3-p3}, we then directly obtain \eqref{main:eq-inv-3-p1}. Furthermore,  we also obtain that 
\begin{equation}\label{main:eq-inv-3-p4}
\begin{aligned}
&\muNR_{\zeta,t}\big( E \big)  - \muR\big( E \big) \\
\leq&\,  \Big( \exp\Big( C_1^\prime N^{-\frac{5\delta_1}{8}}\Big) -1 \Big) \muR\big( E \big) 
+  C_1^\prime \exp\big(  C_1^\prime \big) N^{20\delta_2} \exp\Big( -c_1^\prime N^{2(\delta_2-\eta)} \Big) 
\leq \tfrac{1}{2} C_1^{\prime\prime} N^{-\frac{5\delta_1}{8}}. 
\end{aligned}
\end{equation}
By applying \eqref{main:eq-inv-3-p4} to both $E$ and $\Omegas\backslash E$, it then follows that
\begin{equation}\label{main:eq-inv-3-p5}
\Big| \muNR_{\zeta,t}\big( E \big)  - \muR\big( E \big) \Big| \leq \tfrac{1}{2} C_1^{\prime\prime} N^{-\frac{5\delta_1}{8}}.
\end{equation}
Due to the definition of the total variation distance (Definition \ref{prelim:def-tv-distance}), this implies \eqref{main:eq-inv-3-p2}.\\

\emph{Fourth step: Scaling argument.} Under the condition $\coup \leq N^{-4\delta_2}$, we previously obtained \eqref{main:eq-inv-3-p1} and \eqref{main:eq-inv-3-p2}, which are stronger than the desired estimates \eqref{main:eq-almost-invariance-TV} and \eqref{main:eq-almost-invariance-event}. Thus, in order to prove the statement in the proposition, it remains to treat the case 
\begin{equation}\label{main:eq-ai-3-p1}
N^{-4\delta_2} < \coup \leq N^{\delta_2}, \qquad T\leq N^{\delta_2}, \qquad \text{and} \qquad R\leq N^{\delta_2}. 
\end{equation}
To this end, we choose a scaling parameter $\kappa \in 2^{\Z}$ as the unique dyadic number satisfying
\begin{equation}\label{main:eq-ai-3-p2}
\frac{1}{2} N^{-5\delta_2} < \kappa \leq N^{-5\delta_2}
\end{equation}
and then define $(N^\prime,T^\prime,R^\prime,\coup^\prime)\in \Dyadiclarge \times [1,\infty) \times [1,\infty) \times (0,\infty)$ and $\zeta^\prime \in \Cuttilde$ as
\begin{equation}\label{main:eq-ai-3-p3}
N^\prime = \kappa N, \qquad  T^\prime := \frac{T}{\kappa},  \qquad R^\prime = \frac{R}{\kappa}, \qquad \coup^\prime = \kappa \coup, \qquad \text{and} \qquad
\zeta^\prime(\cdot):= \zeta(\kappa \, \cdot).
\end{equation}
We emphasize that, due to Definition \ref{prelim:def-cut-off}, $\zeta^\prime$ is still an element of $\Cuttilde$.
Due to \eqref{main:eq-ai-3-p1}, \eqref{main:eq-ai-3-p2}, and \eqref{main:eq-ai-3-p3}, it holds that
\begin{equation}\label{main:eq-ai-3-p4}
(N^\prime)^{4\delta_2} \coup^\prime = \kappa^{1+4\delta_2} N^{4\delta_2} \coup \leq \kappa N^{5\delta_2} \leq 1, 
\quad 
 T^\prime = \frac{R}{\kappa} \leq  (N^{\prime})^{10\delta_2}, 
 \qquad \text{and} \qquad
 R^\prime = \frac{R}{\kappa} \leq  (N^\prime)^{10\delta_2}.
\end{equation}
Using \eqref{main:eq-inv-2-p1}, \eqref{main:eq-inv-2-p2}, and \eqref{main:eq-ai-3-p4}, we then obtain for all $t^\prime \in [0,T^\prime]$ that
\begin{align}
\Big\| \mupNR_{\zeta^\prime,t^\prime} - \mupR \Big\|_{\TV} &\leq C_1^{\prime\prime} (N^\prime)^{-\frac{5}{8} \delta_1} \leq C_2 N^{-\frac{1}{2} \delta_1} \label{main:eq-ai-3-p5} 
\end{align}
and obtain for all $t^\prime \in [0,T^\prime]$ and events $E^\prime \subseteq \Omegasprime$ that 
\begin{align}
 \mupNR_{\zeta^\prime,t^\prime}(E^\prime) &\leq C_1^{\prime\prime} \Big( \mupR(E^\prime) + \exp\big( - c_1^{\prime\prime} (N^\prime)^{2(\delta_2-\eta)} \big) \Big) \leq C_2 \Big( \mupR(E^\prime)  + \exp \big( - c_2 N^{\delta_2}\big) \Big). 
\label{main:eq-ai-3-p6} 
\end{align}
Furthermore, due to Lemma \ref{main:lem-scaling}, it holds that
\begin{equation}\label{main:eq-ai-3-p7}
\muNR_{\zeta,t} = (\Scaling_{\kappa^{-1}})_\# \mupNR_{\zeta^\prime,\kappa^{-1}t} 
\qquad \text{and} \qquad \muR = (\Scaling_{\kappa^{-1}})_\# \mupR. 
\end{equation}
We can therefore combine the identity \eqref{main:eq-ai-3-p7} with the estimates  \eqref{main:eq-ai-3-p5} and \eqref{main:eq-ai-3-p6} for $t^\prime=\kappa^{-1}t$, and then obtain 
\eqref{main:eq-almost-invariance-TV} and \eqref{main:eq-almost-invariance-event}. 
\end{proof}

\begin{remark}\label{main:rem-scaling-argument}
\revision{In the proof of Proposition \ref{main:prop-almost-invariance}, we were able to reduce the case $0<\coup \leq N^{\delta_2}$ to the low-temperature case $0<\coup \leq N^{-4\delta_2}$ by using an iteration and scaling argument. The reader may wonder why this reduction is possible in the proof of Proposition \ref{main:prop-almost-invariance}, but not possible in the proof of Proposition~\ref{main:prop-lwp}. In order words, why the local well-posedness in Proposition \ref{main:prop-lwp}.\ref{main:item-lwp-local} cannot be obtained from the short-time well-posedness in Proposition \ref{main:prop-lwp}.\ref{main:item-lwp-short} using an iteration and scaling argument, which would allow us to dispense with the Bourgain-Bulut event. The reason for this is that the growth resulting from each iteration of \eqref{main:eq-inv-1-p1} is much smaller than the growth resulting from each iteration of the short-time well-posedness. To be more precise, we saw in \eqref{main:eq-inv-2-p1} and \eqref{main:eq-inv-3-p1} that \eqref{main:eq-inv-1-p1} can be iterated up to times $T\leq N^{10\delta_2}$. Since the low-temperature condition only imposes $\coup \leq N^{-4\delta_2}$, it is therefore possible to make $\coup T$ large, which allows us to reach a large time even after using the scaling argument. \\
In contrast, our estimates of the nonlinearity in \eqref{main:eq-WM-null-N} have Lipschitz-constant $\sim \Dc$ (see Proposition~\ref{lipschitz:prop-remainder-no-jcb}). If one attempts to obtain Proposition \ref{main:prop-lwp}.\ref{main:item-lwp-local} by iterating Proposition \ref{main:prop-lwp}.\ref{main:item-lwp-short}, each step of the iteration therefore may increase potential errors by a factor of $1+C\Dc$, where $C>0$ is a constant\footnote{\revision{In fact, this would be slightly better than the difference-estimate in Definition \ref{main:def-lwp}.\ref{main:item-lwp-C2}.}}. In order to reach a large time after using a scaling argument, one needs to iterate roughly $\coup^{-1}$ times, which leads to a factor of $(1+C\Dc)^{1/\coup}=(1+C\hcoup \Ac \Bc)^{1/\coup}$. For any $0<\coup \leq N^{-4\delta_2}$, this is prohibitively large.}
\end{remark}

\begin{remark}[Comparison with results on quasi-invariant measures]
In \cite{DT21,FS22,FT22}, similar formulas as in Proposition \ref{structure:prop-Gibbs} have been used to prove the quasi-invariance of Gaussian measures under the flow of nonlinear dispersive equations. As in this article, the proofs in \cite{DT21,FS22,FT22} rely on estimates of space-time integrals such as in \eqref{main:eq-inv-1-p2}. However, there is also an alternative approach to \revision{quasi-invariance in~\cite{T15}}, which only relies on estimates at the time $t=0$. In the setting of this article, this method would require bounds on 
\begin{equation}\label{main:eq-quasi-e1}
\frac{\beta}{2}\frac{\mathrm{d}}{\mathrm{d}t} \int_{\bT_R} \Big( \big\| A^{(\Nscript,\Rscript,\coup,\zetascript)}(t,x) \big\|_{\frkg}^2 + \big\| B^{(\Nscript,\Rscript,\coup,\zetascript)}(t,x) \big\|_{\frkg}^2 \Big) \dx \bigg|_{t=0},
\end{equation}
where $(A^{(\Nscript,\Rscript,\coup,\zetascript)},B^{(\Nscript,\Rscript,\coup,\zetascript)})$ is the solution of \eqref{main:eq-WM-N}. From a direct calculation, \revision{it follows that}
\begin{equation}\label{main:eq-quasi-e2}
\revision{\eqref{main:eq-quasi-e1}=}\frac{1}{4} \zeta(0)^2 \int_{\bT_R} \big\langle \Renorm[N] W_0^{(\Rscript,\coup)},  W_1^{(\Rscript,\coup)} \rangle_\frkg \,  \dx. 
\end{equation}
Since $W_0^{(\Rscript,\coup)}$ and $W_1^{(\Rscript,\coup)}$ are independent, $\frkg$-valued white noises, it follows from Lemma \ref{ansatz:lem-renormalization} that 
\begin{equation*}
\E \bigg[ \Big| \int_{\bT_R} \big\langle \Renorm[N] W_0^{(\Rscript,\coup)},  W_1^{(\Rscript,\coup)} \rangle_\frkg \, \dx \Big|^2 \bigg] 
= \E \Big\| \Renorm[N] W_0^{(\Rscript,\coup)} \Big\|_{L_x^2(\mathbb{T}_R\rightarrow \frkg)}^2 \sim \coup RN. 
\end{equation*}
Thus, \eqref{main:eq-quasi-e1} diverges as $N$ tends to infinity, and can therefore not be used for a proof of Proposition \ref{main:prop-almost-invariance}. \revision{Recent articles inspired by \cite{T15}, such as \cite{GOTW22,OT20Q}, not only involve estimates at time $t=0$, but also rely on modified energies. However, the author was unable to find suitable modified energies for \eqref{main:eq-WM-N}.}
\end{remark}

As an application of Proposition \ref{main:prop-almost-invariance}, we now bound the probability of the Bourgain-Bulut event from Definition \ref{main:def-Bourgain-Bulut}.

\begin{corollary}[Probability of the Bourgain-Bulut event]\label{main:cor-bourgain-bulut-probability}
Let $N\in \Dyadiclarge$, let $1\leq R \leq N^{\delta_2}$, and let $0<\coup \leq N^{\delta_2}$. Furthermore, let $\chi=\chifixed$ be as in Definition \ref{prelim:def-cut-off}. Then, it holds that 
\begin{equation}\label{main:eq-bourgain-bulut-probability}
\muR \big( \Omega \backslash \BBAchi \big) \leq  c_3^{-1}  \exp\Big( - c_3 N^{\delta_2} \Big).
\end{equation}
\end{corollary}

\begin{proof} 
In order to shorten several expressions, we capture the supremum from \mbox{Definition \ref{main:def-Bourgain-Bulut}} using the following notation. For any $A,B\colon \bT_R \rightarrow \frkg$, we write 
\begin{align*}
&\CHHLNsup(A,B)  \\
:=&\,   
\sup_{y\in \R} \, \langle Ny \rangle^{-4} \max \bigg( \Big\| \CHHLN_y(A,B) \Big\|_{L_x^\infty(\bT_R)},
\Big\| \CHHLN_y(B,A) \Big\|_{L_x^\infty(\bT_R)}, \\
&\hspace{5ex}\Big\| \CHHLN_y\big( A ,A \big) - 8\coup \CNbd(y) \Kil  \Big\|_{L_x^\infty(\bT_R)}, 
 \Big\| \CHHLN_y(B,B) - 8\coup \CNbd(y) \Kil  \Big\|_{L_x^\infty(\bT_R)}  \bigg), 
\end{align*}
where $\CHHLN_y$ is as in Definition \ref{jacobi:def-chhl}. To further simplify the notation, we denote the solution of \eqref{main:eq-WM-N} with initial data $(W^{(\Rscript,\coup)}_0, W^{(\Rscript,\coup)}_1)\colon \bT_R \rightarrow \frkg^2$ and cut-off function $\chi$ by $(\An,\Bn)$ and denote the corresponding flow by $\Phin(t)=\Phin(t;0)$.  
Due to our definition of $\CHHLNsup$, it then suffices to show that the estimate
\begin{equation}\label{main:eq-bb-p1}
\begin{aligned}
\sup_{t\in [-4,4]} \CHHLNsup \big( \An(t) , \Bn(t) \big) \leq N^{-\frac{1}{2}+\delta} \coup 
\end{aligned}
\end{equation}
holds on an event satisfying \eqref{main:eq-bourgain-bulut-probability}.
The remainder of the argument is split into two steps. \\

\emph{First step: Meshing argument.} Let $c=c(\delta_\ast)$ be a sufficiently small constant.
In order to perform a meshing argument, we first let $(t_j)_{j=0}^J$ be a grid on $[-4,4]$ with spacing $|t_j-t_{j-1}|\lesssim c N^{-10}$ and less than $\sim c^{-1}N^{10}$ grid points. We also introduce the event
\begin{align*}
 \BBAchi_0 
:=\Big\{ &\big( W_0^{(\Rscript,\coup)},W_1^{(\Rscript,\coup)}\big) \colon 
\big\|  W_0^{(\Rscript,\coup)} \big\|_{\C_x^{s-1}(\bT_R)} \leq \hcoup N, \,
\big\|  W_1^{(\Rscript,\coup)} \big\|_{\C_x^{s-1}(\bT_R)} \leq \hcoup N, \\
 & \, \,\CHHLNsup \big(  W_0^{(\Rscript,\coup)},  W_1^{(\Rscript,\coup)} \big) \leq \frac{1}{2} N^{-\frac{1}{2}+\delta} \coup \Big\}. 
\end{align*}
Equipped with the above definitions and notations, we now claim that 
\begin{equation}\label{main:eq-bb-p2}
\bigcap_{j=0}^J \Phin(t_j)^{-1} \big( \BBAchi_0 \big) \subseteq \BBAchi. 
\end{equation}
In order to see this, we let $t\in [-4,4]$ be arbitrary. Then, we choose an index $0\leq j\leq J$ such that $|t-t_j|\lesssim c N^{-10}$. By definition of $\BBAchi_0$, we obtain on the event $\Phin(t_j)^{-1}(\BBAchi_0)$ that 
\begin{equation}\label{main:eq-bb-p3}
\big\| \An(t_j) \big\|_{\C_x^{s-1}}, \big\| \Bn(t_j) \big\|_{\C_x^{s-1}} \leq \hcoup N. 
\end{equation}
Using the algebra property of $L_x^\infty$, the corresponding well-posedness theory of \eqref{main:eq-WM-N}, and $|t-t_j|\lesssim c N^{-10}$, it is then straightforward to prove that 
\begin{equation}\label{main:eq-bb-p4}
\big\| P_{\leq N}^x\big( \An(t) - \An(t_j) \big) \big\|_{L_x^\infty}, 
\big\| P_{\leq N}^x\big( \Bn(t) - \Bn(t_j) \big) \big\|_{L_x^\infty} \lesssim c N^{-10} (\hcoup N \times N^{1-s})^2 \lesssim  c \coup N^{-6}.
\end{equation}
Using Lemma \ref{prelim:lem-integral}, \eqref{main:eq-bb-p3}, and \eqref{main:eq-bb-p4}, it then follows that
\begin{equation*}
\Big| \CHHLNsup \big( \An(t), \Bn(t) \big)  - \CHHLNsup \big( \An(t_j), \Bn(t_j) \big)  \Big| 
\lesssim  c \coup N^{-6} \times \hcoup N N^{1-s} \lesssim c \coup N^{-4} \leq \tfrac{1}{2} \coup N^{-4}.
\end{equation*}
On the event $\Phin(t_j)^{-1}(\BBAchi_0)$, we then obtain that 
\begin{equation}\label{main:eq-bb-p5}
\begin{aligned}
\CHHLNsup \big( \An(t), \Bn(t) \big) &\leq \tfrac{1}{2} \coup \big( N^{-\frac{1}{2}+\delta} +  N^{-4} \big).
\end{aligned}
\end{equation}
Since $t\in [-4,4]$ was arbitrary, this implies the desired inclusion \eqref{main:eq-bb-p2}. \\

\emph{Second step: Probability estimate.} 
Equipped with \eqref{main:eq-bb-p2}, it remains to prove that 
\begin{equation}\label{main:eq-bb-p6}
\muR \Big( \bigcup_{j=0}^J \Phin(t_j)^{-1} \big( \Omegas \backslash \BBAchi_0 \big) \Big) \leq c_3^{-1} \exp\Big( - c_3 N^{\delta_2} \Big).
\end{equation}
Using a union bound, Proposition \ref{main:prop-almost-invariance}, and $J\sim c^{-1} N^{10}$, we can estimate
\begin{align*}
\muR  \Big( \bigcup_{j=0}^J \Phin(t_j)^{-1} \big( \Omegas \backslash \BBAchi_0 \big) \Big) 
\leq&\, \sum_{j=0}^J \muR \Big(  \Phin(t_j)^{-1} \big( \Omegas \backslash \BBAchi_0 \big) \Big) \\ 
=&\,  \sum_{j=0}^J \muNR_{t_j} \Big( \Omegas \backslash \BBAchi_0 \Big) \\
\lesssim&\, c^{-1} C_2 N^{10} \Big( \muR \big( \Omegas \backslash \BBAchi_0 \big) + \exp\big( - c_2 N^{\delta_2} \big) \Big).
\end{align*}
Thus, the desired estimate \eqref{main:eq-bb-p6} can be reduced to 
\begin{equation*}
\muR \big( \Omegas \backslash \BBAchi_0 \big) \leq c_2^{-1} \exp\big( -c_2 N^{\delta_2} \big), 
\end{equation*}
which follows directly from Lemma \ref{killing:lem-bourgain-bulut} and $R\leq N^{\delta_2}$.
\end{proof}

\subsection{Proof of Theorem \ref{intro:thm-rigorous-A-B}}\label{section:main-gwp} 
In this subsection, we present a proof of Theorem \ref{intro:thm-rigorous-A-B}, i.e., the global well-posedness and invariance of the Gibbs measures for the $\frkg$-valued derivatives $A,B\colon \R^{1+1} \rightarrow \frkg$. In fact, Proposition \ref{main:prop-refined-gwp} contains a refined version of Theorem \ref{intro:thm-rigorous-A-B}, which will later be useful to 
lift our result to the $\frkG$-valued map $\phi\colon \R^{1+1} \rightarrow \frkG$ (see Subsection \ref{section:lifting-theorem}). 
Since we no longer need cut-off functions in time, we focus on the wave maps equation 
\begin{equation}\label{main:eq-WM-nocut}
\partial_t A = \partial_x B \qquad \text{and} \qquad \partial_t B =\partial_x A -  \big[ A, B \big]
\end{equation}
and its finite-dimensional approximation 
\begin{equation}\label{main:eq-WM-N-nocut}
\begin{aligned}
\partial_t A^{(\Nscript,\coup)} &= \partial_x B^{(\Nscript,\coup)} 
\qquad \text{and} \qquad 
\partial_t B^{(\Nscript,\coup)} &= \partial_x A^{(\Nscript,\coup)} 
- \big[ A^{(\Nscript,\coup)}, B^{(\Nscript,\coup)} \big]_{\leq N} + 2 \coup \Renorm[N] A^{(\Nscript,\coup)}.
\end{aligned}
\end{equation}
The rest of this subsection is split into two parts. In Definition \ref{main:def-gwp-small} and Lemma \ref{main:lem-gwp-small}, 
we consider the finite-dimensional approximation \eqref{main:eq-WM-N-nocut} in the low-temperature case $0<\lambda \ll 1$. By combining our local well-posedness result from Proposition \ref{main:prop-lwp} and Bourgain's globalization argument \cite{B94,B96}, we obtain the almost global well-posedness of \eqref{main:eq-WM-N-nocut}. 

In Proposition \ref{main:prop-refined-gwp}, we then prove the global well-posedness of \eqref{main:eq-WM-nocut}. After using Lemma \ref{main:lem-gwp-small} as a starting point, we need to eliminate the frequency-truncation in \eqref{main:eq-WM-N-nocut} and remove the low-temperature condition $0<\lambda\ll 1$. This is done using classical estimates (Proposition \ref{classical:prop-main}) and scaling (Lemma \ref{ansatz:lem-scaling-symmetry}), respectively. \\

In order to state our almost global well-posedness result for \eqref{main:eq-WM-N-nocut}, we first make the following definition.

\begin{definition}[Almost global well-posedness]\label{main:def-gwp-small}
Let $N\in \Dyadiclarge$, let $T,R \geq 1$, and let $\coup>0$. Furthermore, let $W^{(\Rscript,\coup)}_0,W^{(\Rscript,\coup)}_1\colon \bT_R \rightarrow \frkg$. For all $\Nd\in \Dyadiclarge$, let 
$Z^{(\Ndscript,\Rscript,\coup)}_{0},Z^{(\Ndscript,\Rscript,\coup)}_{1}\colon \bT_R\rightarrow \frkg$ be smooth
and define
\begin{equation*}
\Big( \widetilde{W}^{(\Ndscript,\Rscript,\coup)}_{0}  ,\widetilde{W}^{(\Ndscript,\Rscript,\coup)}_{1} \Big) := \Big( P_{\leq \Nd} W\Rlz + Z^{(\Ndscript,\Rscript,\coup)}_{0}, 
P_{\leq \Nd} W\Rlo +Z^{(\Ndscript,\Rscript,\coup)}_{1} \Big). 
\end{equation*}
Finally, let $(A^{(\Nscript,\Rscript,\coup)},B^{(\Nscript,\Rscript,\coup)})\colon \R \times \bT_R \rightarrow \frkg^2$ 
and $(\widetilde{A}^{(\Nscript,\Ndscript,\Rscript,\coup)},\widetilde{B}^{(\Nscript,\Ndscript,\Rscript,\coup)})\colon \R \times \bT_R \rightarrow \frkg^2$
be the unique global solutions of \eqref{main:eq-WM-N-nocut} with initial data $(W^{(\Rscript,\coup)}_0,W^{(\Rscript,\coup)}_1)$
and $( \widetilde{W}^{(\Ndscript,\Rscript,\coup)}_{0}  ,\widetilde{W}^{(\Ndscript,\Rscript,\coup)}_{1})$, respectively. 
For each $\Ld \in \Dyadiclarge$, we then define the event $\GWPN_{\Ld}(T)\subseteq \Omega$ by 
\begin{align*}
&\GWPN_{\Ld}(T) \\
:= \bigg\{ &\, 
\sup_{\substack{\Nd\in\Dyadiclarge\colon \\ \Ld \leq \Nd \leq N^{1-\delta}}} 
\Big\| \big( \widetilde{A}^{(\Nscript,\Ndscript,\Rscript,\coup)},\widetilde{B}^{(\Nscript,\Ndscript,\Rscript,\coup)} \big);
\big( 0, 0 \big) \Big\|_{C_t^0 \Xi_x^{s-1}([0,T]\times \bT_R\rightarrow \frkg^2)} \leq C_4^{4+RT},  \\ 
&\, \sup_{\substack{\Nd\in\Dyadiclarge\colon \\  \Ld \leq  \Nd \leq N^{1-\delta}}} \Nd^{\varsigma} 
\Big\| \big( \widetilde{A}^{(\Nscript,\Ndscript,\Rscript,\coup)},\widetilde{B}^{(\Nscript,\Ndscript,\Rscript,\coup)} \big);
\big( A^{(\Nscript,\Rscript,\coup)},B^{(\Nscript,\Rscript,\coup)}\big) \Big\|_{C_t^0 \Xi_x^{s-1}([0,T]\times \bT_R\rightarrow \frkg^2)}
\leq  C_4^{4+RT}, \\ 
&\, \sup_{\substack{\Nd\in\Dyadiclarge\colon \\  \Ld \leq  \Nd \leq N^{1-\delta}}} \Nd^{\varsigma} 
\Big\| \ac\big[ \widetilde{A}^{(\Nscript,\Ndscript,\Rscript,\coup)}(\cdot,0)\big];
\ac\big[ A^{(\Nscript,\Rscript,\coup)}(\cdot,0)\big]
\Big\|_{C_t^0([0,T]\rightarrow \frkG)} \leq
 C_4^{4+RT} \bigg\}. 
\end{align*}
\end{definition}

Equipped with Definition \ref{main:def-gwp-small}, we can now state and prove the following lemma.

\begin{lemma}[Almost global well-posedness and almost invariance at low temperature]\label{main:lem-gwp-small}
Let $T,R \geq 1$ and let $0<\coup \leq \coup_0$, where $\coup_0=\coup_0(c_1,c_2,c_3,\delta_\ast)$ is sufficiently small. Let $W\Rlz,W\Rlo\colon \bT_R \rightarrow \frkg$ be $2\pi R$-periodic, $\frkg$-valued white noises at temperature $8\coup$. For each $\Nd\in \Dyadiclarge$, let 
$Z^{(\Ndscript,\Rscript,\coup)}_{0},Z^{(\Ndscript,\Rscript,\coup)}_{1}\colon \bT_R\rightarrow \frkg$ be smooth, $\frkg$-valued maps and define
\begin{equation*}
p^Z_{\Nd} := \mathbb{P} \Bigg( \bigcup\displaylimits_{\substack{\Kd \in \Dyadiclarge\colon \\  \Kd \geq \Nd}} \Big\{ \max \Big( \big\| Z^{(\Kdscript,\Rscript,\coup)}_{0} \big\|_{\C_x^{r-1}(\bT_R)}, 
\big\| Z^{(\Kdscript,\Rscript,\coup)}_{1} \big\|_{\C_x^{r-1}(\bT_R)} \Big) \geq \Kd^{-\varsigma} \Big\} \Bigg).
\end{equation*}
Then, we have the following estimates:
\begin{enumerate}[label=(\roman*)]
\item\label{main:item-gwp-difference} (Pointwise estimate) For each $\Ld \in \Dyadiclarge$ satisfying $\Ld \geq C_4^{(2+RT)/\varsigma}$, it holds that 
\begin{equation}\label{main:eq-gwp-pointwise}
\begin{aligned}
 \limsup_{N\rightarrow \infty} \mathbb{P} \Big( \Omega \backslash \GWPN_{\Ld}(T)\Big) 
\leq C_4 T \exp\big(  - c_4R^{-2\eta} \coup^{-1}  \big) + p^Z_{\Ld}.    
\end{aligned}
\end{equation}
\item\label{main:item-gwp-wasserstein} (Wasserstein estimate) For each $N,\Nd\in \Dyadiclarge$, 
let $(\widetilde{A}^{(\Nscript,\Ndscript,\Rscript,\coup)},\widetilde{B}^{(\Nscript,\Ndscript,\Rscript,\coup)})$
be as in Definition \ref{main:def-gwp-small}. If $\Nd \geq C_4^{(2+RT)/\varsigma}$, then it holds for all $t\in [0,T]$ that
\begin{equation}\label{main:eq-gwp-wasserstein}
\begin{aligned}
&\limsup_{N\rightarrow \infty} \WassersteinR \Big( \operatorname{Law} \big( \widetilde{A}^{(\Nscript,\Ndscript,\Rscript,\coup)}(t),\widetilde{B}^{(\Nscript,\Ndscript,\Rscript,\coup)}(t)\big), \muR \Big)  \\
\leq&\,  C_4^{4+RT} \Nd^{-\varsigma} + C_4 T \exp \big( - c_4R^{-2\eta} \coup^{-1}  \big) + p^Z_{\Nd},
\end{aligned}
\end{equation}
where $\WassersteinR$ is the Wasserstein distance from Definition \ref{prelim:def-Wasserstein}.
\end{enumerate}
\end{lemma}

\begin{remark}
We make the following remarks regarding Lemma \ref{main:lem-gwp-small}.
\begin{enumerate}[label=(\alph*)]
    \item The error terms $Z^{(\Ndscript,\Rscript,\coup)}_{0},Z^{(\Ndscript,\Rscript,\coup)}_{1}\colon \bT_R\rightarrow \frkg$ are used in the proof of Proposition \ref{main:prop-refined-gwp} to switch between data which is localized on the spatial scale $\sim R$ and data which is $2\pi R$-periodic. 
    \item Due to Lemma \ref{ansatz:lem-scaling-symmetry}, $\coup$ can be interpreted as a local timescale. Thus, the term  $T\exp(-cR^{-2\eta} \coup^{-1} )$ is the natural loss from Bourgain's globalization argument, see e.g. \cite[(145)]{B96}. In \cite{B96}, this loss is eliminated by taking the local timescale $\coup$ to zero, which will be replicated in this paper via scaling (see the proof of Proposition \ref{main:prop-refined-gwp}).
    \item The use of the Wasserstein metric in Lemma \ref{main:lem-gwp-small} is not essential and \eqref{main:eq-gwp-wasserstein} could likely be replaced by several other estimates. However, due to the properties listed in Lemma \ref{prelim:lem-Wasserstein-properties}, using the Wasserstein metric is convenient.
\end{enumerate}
\end{remark}

The following proof of Lemma \ref{main:lem-gwp-small} is a variant of Bourgain's globalization argument \cite{B94,B96}. 

\begin{proof}[Proof of Lemma \ref{main:lem-gwp-small}:]
During this proof, we use the notation from Definition \ref{main:def-gwp-small}. 
For expository purposes, we separate the proof into four steps. \\ 

\emph{First step: Introduction of parameters.}
We choose the parameters $\Ac$ and $\Bc$ as 
\begin{equation}\label{main:eq-gwp-small-q1}
\Ac := \tfrac{c_1}{C_1}  \, \coup^{-\frac{1}{2}} \qquad \text{and} \qquad \Bc:=C_1. 
\end{equation}
Provided that $\coup_0$ is sufficiently small depending on $c_1$ and $C_1$, it holds that $\Ac \geq C_1$. Furthermore, due to \eqref{main:eq-gwp-small-q1}, it also holds that
\begin{equation}\label{main:eq-gwp-small-q2}
\Dc = \hcoup \Ac \Bc = \hcoup \times \tfrac{c_1}{C_1}  \coup^{-\frac{1}{2}} \times C_1 =c_1. 
\end{equation}

\emph{Second step: Global control of $A^{(\Nscript,\Rscript,\coup)}$ and $B^{(\Nscript,\Rscript,\coup)}$.} During this step,  we may always assume that 
\begin{equation}\label{main:eq-gwp-small-qq4}
N \geq C_4^{(2+RT)/\varsigma}.
\end{equation}
To simplify the notation, we assume that $T=J\in \mathbb{N}$, which can always be satisfied by increasing $T\geq 1$. We then let $\LWPD$ be as in \eqref{main:eq-lwpd-one} from Definition \ref{main:def-lwp} and define the event\footnote{The notation ``$\textup{BG}$" stands for Bourgain's globalization argument.}
\begin{equation}\label{main:eq-gwp-small-q4}
\BG(T) := \bigcap_{j=0}^{J-1} \Big\{ \big( A^{(\Nscript,\Rscript,\coup)}(j),B^{(\Nscript,\Rscript,\coup)}(j) \big) \in \LWPD \Big\}. 
\end{equation}
We now claim that 
\begin{equation}\label{main:eq-gwp-small-q5}
\mathbb{P} \big( \Omega \backslash \BG(T) \big) \leq C_4 T \Big( \exp\big( -c_4 R^{-2\eta} \coup^{-1} \big) + \exp\big( -c_4 N^{\delta_2} \big) \Big).
\end{equation}
Indeed, using \eqref{main:eq-gwp-small-q1}, Proposition \ref{main:prop-lwp}.\ref{main:item-lwp-local}, and Corollary \ref{main:cor-bourgain-bulut-probability}, we first obtain for the cut-off function $\chi=\chifixed$ from Definition \ref{prelim:def-cut-off} that
\begin{equation}\label{main:eq-gwp-small-q6}
\begin{aligned}
\mathbb{P} \Big( \big( W\Rlz,W\Rlo \big) \not \in \LWPD \Big) 
&\leq \muR \Big( \Omegas \backslash \BBAchi \Big) + c_1^{-1} \exp \Big( - c_1 R^{-2\eta} \Ac^2 \Big) \\ 
&\leq c_3^{-1} \Big( \exp\big(-c_3 N^{\delta_2}\big)  + \exp\big( -c_3R^{-2\eta} \coup^{-1} \big) \Big).    
\end{aligned}
\end{equation}
Then, using a union bound, Proposition \ref{main:prop-almost-invariance}, \eqref{main:eq-gwp-small-qq4}, and \eqref{main:eq-gwp-small-q6}, we then obtain that
\begin{align*}
\mathbb{P} \big( \Omega \backslash \BG(T) \big) 
&\leq \sum_{j=0}^{J-1} \mathbb{P} \Big( \big( A^{(\Nscript,\Rscript,\coup)}(j),B^{(\Nscript,\Rscript,\coup)}(j) \big) \not \in \LWPD \Big) \\
&\leq C_2 \sum_{j=0}^{J-1} \Big( \mathbb{P} \Big( \big( A^{(\Nscript,\Rscript,\coup)}(0),B^{(\Nscript,\Rscript,\coup)}(0) \big) \not \in \LWPD \Big) + \exp\big( - c_2 N^{\delta_2}\big) \Big) \\ 
&\leq 2 C_2 c_3^{-1} T  \Big( \exp\big( -c_3 R^{-2\eta} \coup^{-1} \big) + \exp\big(-c_3 N^{\delta_2}\big) \Big). 
\end{align*}
Thus, we obtain the desired estimate \eqref{main:eq-gwp-small-q5}. \\

\emph{Third Step: Proof of the pointwise estimate.} We now turn towards the proof of the pointwise estimate, i.e., the proof of \eqref{main:eq-gwp-pointwise}. To this end, we first fix $\Ld\in \Dyadiclarge$ and then define $\ZGd$ as the event on which all $Z^{(\Ndscript,\Rscript,\coup)}_{0}$ and $Z^{(\Ndscript,\Rscript,\coup)}_{1}$-terms with $\Nd \geq \Ld$ are good, i.e., we define
\begin{equation}\label{main:eq-zgd}
\ZGd := \bigcap_{\Nd \geq \Ld} \Big\{ \max \Big( \big\| Z^{(\Ndscript,\Rscript,\coup)}_{0} \big\|_{\C_x^{r-1}(\bT_R)}, 
\big\| Z^{(\Ndscript,\Rscript,\coup)}_{1} \big\|_{\C_x^{r-1}(\bT_R)} \Big) < \Nd^{-\varsigma} \Big\}.
\end{equation}
Using the probability $p^Z_{\Ld}$ defined in the statement of this lemma, we then obtain that 
\begin{equation}\label{main:eq-gwp-small-q8}
\mathbb{P} \big( \Omega \backslash \ZGd \big) \leq p^Z_{\Ld}.
\end{equation}
We now claim that the difference between $(  \widetilde{A}^{(\Nscript,\Ndscript,\Rscript,\coup)},\widetilde{B}^{(\Nscript,\Ndscript,\Rscript,\coup)})$ and 
$(A^{(\Nscript,\Rscript,\coup)},B^{(\Nscript,\Rscript,\coup)})$ is controlled on the event $\BG(T) \medcap \ZGd$, i.e., we claim that 
\begin{equation}\label{main:eq-gwp-small-q9}
\begin{aligned}
\BG(T) \medcap \ZGd 
\subseteq  \, \GWPN_{\Ld}(T). 
\end{aligned}
\end{equation}
To see \eqref{main:eq-gwp-small-q9}, we henceforth restrict to the event $\BG(T) \medcap \ZGd$. We further let $\Nd\in \Dyadiclarge$ satisfy $\Ld \leq \Nd \leq N^{1-\delta}$ and, for all $0\leq j \leq J$, define the parameter
\begin{equation}\label{main:eq-small-q10m}
\epsilon_j := C_4^{1+Rj} \Nd^{-\varsigma}.
\end{equation}
In order to use Proposition \ref{main:prop-lwp} repeatedly, we need to verify the three conditions
\begin{align}
\max \Big( \big\| Z^{(\Ndscript,\Rscript,\coup)}_{0} \big\|_{\C_x^{r-1}(\bT_R)}, 
\big\| Z^{(\Ndscript,\Rscript,\coup)}_{1} \big\|_{\C_x^{r-1}(\bT_R)} \Big)  
&\leq \hcoup \Ac \epsilon_0, 
\label{main:eq-epsilon-e1}\\
\epsilon_j &\leq C_1^{-1}, 
\label{main:eq-epsilon-e2}\\
C_1 \Bc R^{\eta} \big(\epsilon_j + \Bc \Nd^{-\varsigma} \big)
&\leq \epsilon_{j+1}. 
\label{main:eq-epsilon-e3}
\end{align}
The first condition \eqref{main:eq-epsilon-e1} can be obtained using the definitions of $\Ac$ and $\ZGd$ from \eqref{main:eq-gwp-small-q1} and \eqref{main:eq-zgd}, which imply that 
\begin{equation*}
\max \Big( \big\| Z^{(\Ndscript,\Rscript,\coup)}_{0} \big\|_{\C_x^{r-1}(\bT_R)}, 
\big\| Z^{(\Ndscript,\Rscript,\coup)}_{1} \big\|_{\C_x^{r-1}(\bT_R)} \Big)  
\leq \Nd^{-\varsigma} \leq \tfrac{c_1}{C_1} C_4 \Nd^{-\varsigma} =\hcoup \Ac \epsilon_0. 
\end{equation*}
The second condition \eqref{main:eq-epsilon-e2} can be obtained from our assumption $\Nd \geq \Ld \geq C_4^{(2+RT)/\varsigma}$, 
which yields
\begin{equation*}
\epsilon_j \leq C_4^{1+RT} \Ld^{-\varsigma} \leq C_4^{-1} \leq C_1^{-1}.
\end{equation*}
The third condition \eqref{main:eq-epsilon-e2} follows from the definition of $\Bc$ in \eqref{main:eq-gwp-small-q1}, which directly implies that
\begin{equation*}
C_1 \Bc R^{\eta} \big(\epsilon_j + \Bc \Nd^{-\varsigma} \big) \leq 2 C_1^2 R^\eta C_4^{1+Rj} \Nd^{-\varsigma} \leq C_4^{1+R(j+1)} \Nd^{-\varsigma} =\epsilon_{j+1}.
\end{equation*}
Thus, all three conditions in \eqref{main:eq-epsilon-e1}, \eqref{main:eq-epsilon-e2}, and \eqref{main:eq-epsilon-e3} are satisfied. 
As a result of \eqref{main:eq-epsilon-e1}, we now first obtain that 
\begin{equation*}
 \Big( P_{\leq \Nd} W\Rlz + Z^{(\Ndscript,\Rscript,\coup)}_{0}, 
P_{\leq \Nd} W\Rlo +Z^{(\Ndscript,\Rscript,\coup)}_{1} \Big) 
\in \Per \Big( W\Rlz, W\Rlo; \Ac, \epsilon_0\Big),
\end{equation*}
or, equivalently, 
\begin{equation}\label{main:eq-small-q12} 
\Big( \widetilde{A}^{(\Nscript,\Ndscript,\Rscript,\coup)}(0),\widetilde{B}^{(\Nscript,\Ndscript,\Rscript,\coup)}(0) \Big)  \in 
\Per \Big( A^{(\Nscript,\Rscript,\coup)}(0),B^{(\Nscript,\Rscript,\coup)}(0); \Ac,\epsilon_0 \Big).
\end{equation}
Due to \eqref{main:eq-gwp-small-q2}, \eqref{main:eq-gwp-small-q4}, \eqref{main:eq-epsilon-e2}, \eqref{main:eq-epsilon-e3}, and \eqref{main:eq-small-q12}, we can now iterate our local well-posedness result from Proposition \ref{main:prop-lwp}. Due to \ref{main:item-lwp-C4} from Definition \ref{main:def-lwp}, we obtain for all $0\leq j\leq J-1$ that
\begin{equation*}
\Big( \widetilde{A}^{(\Nscript,\Ndscript,\Rscript,\coup)}(j),\widetilde{B}^{(\Nscript,\Ndscript,\Rscript,\coup)}(j) \Big)  \in 
\Per \Big( A^{(\Nscript,\Rscript,\coup)}(j),B^{(\Nscript,\Rscript,\coup)}(j); \Dc,\epsilon_j \Big),
\end{equation*}
which is needed to perform the iteration. Furthermore, due to \ref{main:item-lwp-C1}, \ref{main:item-lwp-C2}, and \ref{main:item-lwp-C3} from Definition \ref{main:def-lwp} and the choice of $\Ac$ from \eqref{main:eq-gwp-small-q1}, we also obtain for all $0\leq j \leq J-1$ that 
\begin{align}
\Big\| \big( \widetilde{A}^{(\Nscript,\Ndscript,\Rscript,\coup)},\widetilde{B}^{(\Nscript,\Ndscript,\Rscript,\coup)} \big);
\big( 0, 0 \big) \Big\|_{C_t^0 \Xi_x^{s-1}([j,j+1]\times \bT_R)} 
&\leq C_1 \Dc,  \label{main:eq-small-q13-a} \\ 
\Big\| \big( \widetilde{A}^{(\Nscript,\Ndscript,\Rscript,\coup)},\widetilde{B}^{(\Nscript,\Ndscript,\Rscript,\coup)} \big);
\big( A^{(\Nscript,\Rscript,\coup)},B^{(\Nscript,\Rscript,\coup)}\big) \Big\|_{C_t^0 \Xi_x^{s-1}([j,j+1]\times \bT_R)}
&\leq \hcoup \Ac \epsilon_{j+1} \leq  \epsilon_{j+1} , \label{main:eq-small-q13-b}  \\ 
\Big\| \ac\big[ j, \widetilde{A}^{(\Nscript,\Ndscript,\Rscript,\coup)}(\cdot,0)\big];
\ac\big[j, A^{(\Nscript,\Rscript,\coup)}(\cdot,0)\big]
\Big\|_{C_t^0([j,j+1]\rightarrow \frkG)}  
&\leq \hcoup \Ac  \epsilon_{j+1} \leq \epsilon_{j+1}.\label{main:eq-small-q13-c} 
\end{align}
By taking the maximum of \eqref{main:eq-small-q13-a}, \eqref{main:eq-small-q13-b}, and \eqref{main:eq-small-q13-c} over all $0\leq j \leq J-1$ and all $\Nd\in \Dyadiclarge$ satisfying  $\Ld \leq \Nd \leq N^{1-\delta}$ and by using Lemma \ref{prelim:lem-gluing-estimate}, we obtain the estimates in the definition of $\GWPN$ and therefore obtain the inclusion \eqref{main:eq-gwp-small-q9}. Equipped with the inclusion \eqref{main:eq-gwp-small-q9}, we can now prove the estimate \eqref{main:eq-gwp-pointwise}. Indeed, using \eqref{main:eq-gwp-small-q5}, \eqref{main:eq-gwp-small-q8}, and \eqref{main:eq-gwp-small-q9}, we obtain that 
\begin{equation}\label{main:eq-small-q13-d}
\begin{aligned}
\mathbb{P} \big( \Omega \backslash \GWPN_{\Ld}(T) \big)  
&\leq  
\mathbb{P} \big( \Omega \backslash \BG_{\Ld}(T) \big) 
+ \mathbb{P} \big( \Omega \backslash \ZGd \big)  \\ 
&\leq   C_4 T  \exp\big( -c_4 R^{-2\eta} \coup^{-1} \big) + C_4 T \exp\big( -c_4 N^{\delta_2} \big)  + p^{Z}_{\Ld}. 
\end{aligned}
\end{equation}
By taking $N\rightarrow \infty$, we therefore obtain  \eqref{main:eq-gwp-pointwise}. \\ 

\emph{Fourth step: Wasserstein estimate.} We now turn towards the Wasserstein estimate, i.e., \eqref{main:eq-gwp-wasserstein}.
Let $\Nd \in \Dyadiclarge$ satisfy $\Nd \geq C_4^{(2+RT)/\varsigma}$ and let $t\in [0,T]$. Using the triangle inequality for the Wasserstein metric, we first estimate 
\begin{align}
&\, \limsup_{N\rightarrow \infty} \WassersteinR \Big( \operatorname{Law} \big( \widetilde{A}^{(\Nscript,\Ndscript,\Rscript,\coup)}(t),\widetilde{B}^{(\Nscript,\Ndscript,\Rscript,\coup)}(t)\big), \muR \Big)  \notag \\
\leq &\, \limsup_{N\rightarrow \infty}
\WassersteinR \Big( \operatorname{Law} \big( \widetilde{A}^{(\Nscript,\Ndscript,\Rscript,\coup)}(t),\widetilde{B}^{(\Nscript,\Ndscript,\Rscript,\coup)}(t)\big),  \operatorname{Law} \big( A^{(\Nscript,\Rscript,\coup)}(t),B^{(\Nscript,\Rscript,\coup)}(t) \big) \Big)
\label{main:eq-small-q15} \\
+ & \, \limsup_{N\rightarrow \infty}
\WassersteinR \Big( \operatorname{Law} \big( A^{(\Nscript,\Rscript,\coup)}(t),B^{(\Nscript,\Rscript,\coup)}(t) \big),  \muR \Big)
\label{main:eq-small-q16}.
\end{align} 
Using Lemma \ref{prelim:lem-Wasserstein-properties}.\ref{prelim:item-Wasserstein-coupling}, the definition of $\GWPN_{\Nd}(T)$, and \eqref{main:eq-small-q13-d}, it holds that 
\begin{align*}
\eqref{main:eq-small-q15} 
\leq \limsup_{N\rightarrow \infty} \Big( C_4^{4+RT} \Nd^{-\varsigma} + \mathbb{P}\big( \Omega \backslash \GWPN_{\Nd}(T) \big) \Big)
\leq C_4^{4+RT} \Nd^{-\varsigma} + C_4 T  \exp\big( -c_4 R^{-2\eta} \coup^{-1} \big) + p^{Z}_{\Nd}. 
\end{align*}
Due to Lemma \ref{prelim:lem-Wasserstein-properties}.\ref{prelim:item-Wasserstein-TV} and Proposition \ref{main:prop-almost-invariance}, it holds that 
\begin{align*}
\eqref{main:eq-small-q16} 
&= \limsup_{N\rightarrow \infty}  \WassersteinR \big( \muNR_t , \muR \big) \leq 
\limsup_{N\rightarrow \infty}  \big\| \muNR_t - \muR \big\|_{\TV} =0.
\end{align*}
By combining the estimates of \eqref{main:eq-small-q15} and \eqref{main:eq-small-q16}, we obtain the desired estimate \eqref{main:eq-gwp-wasserstein}.
\end{proof}

Equipped with Lemma \ref{main:lem-gwp-small}, we can now prove the following refined version of Theorem \ref{intro:thm-rigorous-A-B}. 

\begin{proposition}[Refined global well-posedness and invariance]\label{main:prop-refined-gwp}
Let $\coup >0$ and let $W^{(\coup)}_0,W^{(\coup)}_1\colon \R \rightarrow \frkg$ be independent, $\frkg$-valued white noises at temperature $8\coup$. For each $\Nd\in \Dyadiclarge$, let $(A^{(\Ndscript,\coup)},B^{(\Ndscript,\coup)})\colon \R^{1+1}\rightarrow \frkg^2$ be the unique smooth global solution of \eqref{main:eq-WM-nocut} with initial data $(P_{\leq \Nd}^x W^{(\coup)}_0,P_{\leq \Nd}^xW^{(\coup)}_1)$. Then, the limits
\begin{equation*}
\big( A^{(\coup)}, B^{(\coup)} \big) = \lim_{\Nd \rightarrow \infty} \big( A^{(\Ndscript,\coup)},B^{(\Ndscript,\coup)} \big) 
\qquad \text{and} \qquad 
\accoup = \lim_{\Nd\rightarrow \infty} \ac\Big[ A^{(\Ndscript,\coup)}(\cdot,0)\Big]
\end{equation*}
almost surely exist in $C_t^0 \Xi^{s-1}_x([0,T]\times [-R,R]\rightarrow \frkg^2)$ and in $C_t^0([0,T]\rightarrow \frkG)$ for all $T,R\geq 1$.
Furthermore, for each $t\in \R$, the law of $\big( A^{(\coup)}(t), B^{(\coup)}(t) \big)$ is given by $\mu^{(\coup)}$.
\end{proposition}

\begin{remark} It is natural to expect that $\accoup = \ac \big[ A^{(\coup)}(\cdot,0)\big]$. In order to define $\ac \big[ A^{(\coup)}(\cdot,0)\big]$, however, we would need to prove that $A^{(\coup)}(t,0)\in \Xi_t^{s-1}([0,T]\rightarrow \frkg)$ for all $T\geq 1$. While this is also expected, a proof poses the technical difficulties previously discussed in Remark \ref{main:rem-lwp-definition}. Since the identity $\accoup = \ac \big[ A^{(\coup)}(\cdot,0)\big]$ will not be needed in our proof of Theorem \ref{intro:thm-rigorous-A-B} or Theorem \ref{intro:thm-rigorous-phi}, we therefore do not pursue it here. 
\end{remark}

\begin{proof}[Proof of Proposition \ref{main:prop-refined-gwp}]
Since any constants from this proof will not appear in the statement of the proposition, and therefore not be referenced in the rest of this article, we now let $C\geq 1$ be sufficiently large and $0<c\leq 1$ be sufficiently small constants depending on $(C_j)_{j=0}^4$, $(c_j)_{j=0}^4$, and $(\delta_j)_{j=0}^5$. We also let $\varphi,\widetilde{\varphi}\colon \R \rightarrow [0,1]$ be smooth cutoff functions satisfying
\begin{alignat*}{3}
\varphi\big|_{[-1,1]} &= 1, \hspace{10ex}& \varphi\big|_{\R \backslash [-\frac{9}{8},\frac{9}{8}]} &=0, \\
\widetilde{\varphi}\big|_{[-\frac{5}{4},\frac{5}{4}]} &= 1, \qquad & \widetilde{\varphi}\big|_{\R \backslash [-\frac{11}{8},\frac{11}{8}]} &=0,
\end{alignat*}
and define $\varphi_R(x)=\varphi(x/R)$ and $\widetilde{\varphi}_R(x)=\widetilde{\varphi}(x/R)$. Throughout this proof, we may always assume that 
\begin{equation}\label{main:eq-refined-q-z1}
T \leq \frac{1}{8} R, 
\end{equation}
since otherwise we can simply increase $R$. For expository purposes, we split the argument into five steps. \\

\emph{First step: Periodization.} Since the statement of this proposition only depends on the law of $W^{(\coup)}_0$ and $W^{(\coup)}_1$, 
we can assume that $W^{(\coup)}_0$ and $W^{(\coup)}_1$ are given by the representation from 
Lemma \ref{prelim:lem-white-noise-representation}.\ref{prelim:item-representation-1}.
For each $R\geq 1$, we can then define $W_0^{(\Rscript,\coup)}$ and $W_1^{(\Rscript,\coup)}$ using the representations from Lemma \ref{prelim:lem-white-noise-representation}.\ref{prelim:item-representation-2}. 
For each $\Nd \in \Dyadiclarge$, we now define  $Z^{(\Ndscript,\Rscript,\coup)}_{0},Z^{(\Ndscript,\Rscript,\coup)}_{1}\colon \bT_R\rightarrow \frkg$ as the $2\pi R$-periodic extensions of 
\begin{equation}\label{main:eq-refined-q-a1}
\Big( \widetilde{\varphi}_R P_{\leq \Nd} \big( W^{(\coup)}_0 - W^{(R,\coup)}_0 \big), 
\widetilde{\varphi}_R P_{\leq \Nd} \big( W^{(\coup)}_1 - W^{(R,\coup)}_1 \big) \Big).
\end{equation}
Due to Lemma \ref{prelim:lem-white-noise-difference}, it then holds for all $\Nd \geq C^{RT}$ that 
\begin{equation}\label{main:eq-refined-q-a2}
\mathbb{P} \Big( \big\| Z^{(\Ndscript,\Rscript,\coup)}_{0} \big\|_{\C_x^{r-1}} \geq \Nd^{-\varsigma}  
\quad \textup{or} \quad  \big\| Z^{(\Ndscript,\Rscript,\coup)}_{1} \big\|_{\C_x^{r-1}} \geq \Nd^{-\varsigma} \Big) \lesssim \Nd^{-100} \leq \Nd^{-10}.
\end{equation}
For each $\Nd\in \Dyadiclarge$, we further define
\begin{equation}\label{main:eq-refined-q-a3}
\Big( \widetilde{W}^{(\Ndscript,\Rscript,\coup)}_{0}, \widetilde{W}^{(\Ndscript,\Rscript,\coup)}_{1} \Big) 
:= \Big( P_{\leq \Nd} W^{(\Rscript,\coup)}_{0} + Z^{(\Rscript,\coup)}_{0}, 
P_{\leq \Nd} W^{(\Rscript,\coup)}_{1} + Z^{(\Rscript,\coup)}_{1} \Big).
\end{equation}
In the following, we often think of \eqref{main:eq-refined-q-a3} as an approximation of 
$(P_{\leq \Nd} W^{(\Rscript,\coup)}_{0},P_{\leq \Nd} W^{(\Rscript,\coup)}_{1})$. 
We note that, due to the definition of $Z^{(\Ndscript,\Rscript,\coup)}_{0}$ and $Z^{(\Ndscript,\Rscript,\coup)}_{1}$, it then holds that 
\begin{equation}\label{main:eq-refined-q-a4}
\Big( \widetilde{W}^{(\Ndscript,\Rscript,\coup)}_{0}, \widetilde{W}^{(\Ndscript,\Rscript,\coup)}_{1} \Big)(x) 
= \Big( P_{\leq \Nd} W^{(\coup)}_{0}, 
P_{\leq \Nd} W^{(\coup)}_{1}  \Big)(x)  \qquad \text{for all } x\in \big[ - \tfrac{5}{4} R, \tfrac{5}{4} R\big].
\end{equation}
Finally, we let $(\widetilde{A}^{(\Ndscript,R,\coup)},\widetilde{B}^{(\Ndscript,R,\coup)})\colon \R \times \bT_R\rightarrow \frkg^2$ be the unique global solution of \eqref{main:eq-WM-nocut} with the initial data from \eqref{main:eq-refined-q-a3}. 
Due to finite speed of propagation, our assumption $T\leq \frac{1}{8} R$ from \eqref{main:eq-refined-q-z1}, and \eqref{main:eq-refined-q-a4}, it then follows 
for all $t\in [ - T , T ]$ that 
\begin{equation}\label{main:eq-refined-q-a5}
 \Big(\varphi_R \widetilde{A}^{(\Ndscript,\Rscript,\coup)}, \varphi_R \widetilde{B}^{(\Ndscript,\Rscript,\coup)} \Big)(t,x) 
= \Big( \varphi_R A^{(\Ndscript,\coup)}, \varphi_R B^{(\Ndscript,\coup)} \Big)(t,x). 
\end{equation}

\emph{Second step: Finite-dimensional approximation \eqref{main:eq-WM-N-nocut}.}
We now compare $(\widetilde{A}^{(\Ndscript,R,\coup)},\widetilde{B}^{(\Ndscript,R,\coup)})$ 
with solutions of the finite-dimensional approximation of the wave maps equation in \eqref{main:eq-WM-N-nocut}. 
For any $N\in \Dyadiclarge$, we therefore let $(\widetilde{A}^{(\Nscript,\Ndscript,\Rscript,\coup)},\widetilde{B}^{(\Nscript,\Ndscript,\Rscript,\coup)})\colon \R \times \bT_R \rightarrow \frkg^2$ be the unique global solution of \eqref{main:eq-WM-N-nocut} with the initial data 
$(\widetilde{W}^{(\Ndscript,\Rscript,\coup)}_{0}, \widetilde{W}^{(\Ndscript,\Rscript,\coup)}_{1})$
from \eqref{main:eq-refined-q-a3}. We emphasize that $(\widetilde{W}^{(\Ndscript,\Rscript,\coup)}_{0}, \widetilde{W}^{(\Ndscript,\Rscript,\coup)}_{1})$ does not depend on $N\in \Dyadiclarge$ and that, for any fixed $\Nd\in \Dyadiclarge$,
$(\widetilde{W}^{(\Ndscript,\Rscript,\coup)}_{0}, \widetilde{W}^{(\Ndscript,\Rscript,\coup)}_{1})$ is smooth. Using classical well-posedness results (Proposition \ref{classical:prop-main}), we obtain almost surely that
\begin{equation}\label{main:eq-refined-q-b1}
\lim_{N\rightarrow \infty} \Big\| \big( \widetilde{A}^{(\Nscript,\Ndscript,\Rscript,\coup)},\widetilde{B}^{(\Nscript,\Ndscript,\Rscript,\coup)}\big) -
\big( \widetilde{A}^{(\Ndscript,\Rscript,\coup)},\widetilde{B}^{(\Ndscript,\Rscript,\coup)}\big) \Big\|_{C_t^0 C_x^0([0,T]\times \bT_R)} = 0.
\end{equation}
Since the $C_t^0 C_x^0$-metric is much stronger than the $C_t^0 \Xi_x^{s-1}$ and $C_x^0 \Xi_t^{s-1}$-metrics, we obtain almost surely that 
\begin{align}
\lim_{N\rightarrow \infty} \Big\| \big( \widetilde{A}^{(\Nscript,\Ndscript,\Rscript,\coup)},\widetilde{B}^{(\Nscript,\Ndscript,\Rscript,\coup)}\big); 
\big( \widetilde{A}^{(\Ndscript,\Rscript,\coup)},\widetilde{B}^{(\Ndscript,\Rscript,\coup)}\big) \Big\|_{C_t^0 \Xi_x^{s-1}([0,T]\times \bT_R)} &= 0, \label{main:eq-refined-q-b2} \\ 
\lim_{N\rightarrow \infty} 
\Big\| \ac\big[ \widetilde{A}^{(\Nscript,\Ndscript,\Rscript,\coup)}(\cdot,0)\big] ; 
\ac\big[ \widetilde{A}^{(\Ndscript,\Rscript,\coup)}(\cdot,0)\big] \Big\|_{C_t^0([0,T]\rightarrow \frkG)}
&= 0. \label{main:eq-refined-q-b3}
\end{align}

\emph{Third step: Low temperature.} In this step, we only consider the low-temperature case
\begin{equation}\label{main:eq-refined-q-c1}
0 < \coup \leq \coup_0, 
\end{equation}
where $\coup_0$ is as in Lemma \ref{main:lem-gwp-small}. For all $\Ld\in \Dyadiclarge$ satisfying $\Ld \geq C^{RT}$, we now prove that 
\begin{equation}\label{main:eq-refined-q-c2}
\begin{aligned}
&\mathbb{P} \bigg( \sup_{\substack{\Nd \in \Dyadiclarge\colon \\  \Ld \leq \Nd}} 
 \Big\| \varphi_R \big( A^{(\Ndscript,\coup)}, B^{(\Ndscript,\coup)} \big); 
\big( 0, 0 \big) \Big\|_{C_t^0 \Xi^{s-1}_x([0,T]\times \R \rightarrow \frkg^2)}
> C^{RT}, \\ 
&\hspace{3ex}
\sup_{\substack{\Md,\Nd \in \Dyadiclarge\colon \\ \Ld \leq \Md \leq \Nd}} 
\Md^\varsigma \Big\| \varphi_R \big( A^{(\Mdscript,\coup)}, B^{(\Mdscript,\coup)} \big); 
\varphi_R  \big( A^{(\Ndscript,\coup)}, B^{(\Ndscript,\coup)} \big) \Big\|_{C_t^0 \Xi^{s-1}_x([0,T]\times \R \rightarrow \frkg^2)}
> C^{RT}, \\ 
&\hspace{3ex}\textup{or} \quad \sup_{\substack{\Md,\Nd\in\Dyadiclarge\colon \\  \Ld \leq \Md \leq \Nd}} \Md^{\varsigma}
\Big\| \ac\big[ A^{(\Mdscript,\coup)}(\cdot,0)\big] ; 
\ac\big[ A^{(\Ndscript,\coup)}(\cdot,0)\big] \Big\|_{C_t^0([0,T]\rightarrow \frkG)} > C^{RT} \bigg) \\
\leq&\, C T \exp \big( - c R^{-2\eta} \coup^{-1} \big) + C \Ld^{-10}.
\end{aligned}
\end{equation}
Furthermore, we also prove for all $\Nd\in \Dyadiclarge$ satisfying $\Nd \geq C^{RT}$ and all $t\in [0,T]$ that 
\begin{equation}\label{main:eq-refined-q-c3}
\begin{aligned}
&\, \WassersteinR \Big( \Law \big( \varphi_R A^{(\Ndscript,\coup)}(t), \varphi_R B^{(\Ndscript,\coup)}(t) \big), 
(\varphi_R \otimes \varphi_R)_\# \mu^{(\coup)} \Big)  \\ 
\leq&\,  C T \exp \big( - c R^{-2\eta} \coup^{-1} \big) +C^{RT}  \Nd^{-\varsigma}.
\end{aligned}
\end{equation}
In order to prove \eqref{main:eq-refined-q-c2}, we introduce an additional frequency-scale $\Hd \in \Dyadiclarge$, which will 
temporarily serve as an upper bound for the frequency-scales in \eqref{main:eq-refined-q-c2}.
By first restricting to frequency-scales $\leq \Hd$, then using \eqref{main:eq-refined-q-a5}, and then using \eqref{main:eq-refined-q-b2}-\eqref{main:eq-refined-q-b3}, we obtain that
\begin{align}
&\textup{LHS of } \eqref{main:eq-refined-q-c2}  \notag \\
\leq&\, \limsup_{\Hd \rightarrow \infty} \limsup_{N\rightarrow \infty} \mathbb{P} \bigg( \sup_{\substack{\Nd \in \Dyadiclarge\colon \\  \Ld \leq \Nd\leq \Hd}} 
 \Big\| \varphi_R \big( \widetilde{A}^{(\Nscript,\Ndscript,\Rscript,\coup)}, \widetilde{B}^{(\Nscript,\Ndscript,\Rscript,\coup)} \big); 
\big( 0, 0 \big) \Big\|_{C_t^0 \Xi^{s-1}_x([0,T]\times \R \rightarrow \frkg^2)}
> C^{RT}, 
\notag \\ 
&\hspace{0ex} 
\sup_{\substack{\Md,\Nd \in \Dyadiclarge\colon \\ \Ld \leq \Md \leq \Nd \leq \Hd}}  \hspace{-2ex}
\Md^\varsigma \Big\| \varphi_R \big( \widetilde{A}^{(\Nscript,\Mdscript,\Rscript,\coup)}, \widetilde{B}^{(\Nscript,\Mdscript,\Rscript,\coup)} \big); 
\varphi_R  \big( \widetilde{A}^{(\Nscript,\Ndscript,\Rscript,\coup)}, \widetilde{B}^{(\Nscript,\Ndscript,\Rscript,\coup)} \big) \Big\|_{C_t^0 \Xi^{s-1}_x([0,T]\times \R \rightarrow \frkg^2)}
> C^{RT}, \notag \\ 
&\hspace{3ex}\textup{or} \quad \sup_{\substack{\Md,\Nd\in\Dyadiclarge\colon \\  \Ld \leq \Md \leq \Nd \leq \Hd}} \Md^{\varsigma}
\Big\| \ac\big[ \widetilde{A}^{(\Nscript,\Mdscript,\Rscript,\coup)}(\cdot,0)\big] ; 
\ac\big[ \widetilde{A}^{(\Nscript,\Ndscript,\Rscript,\coup)}(\cdot,0)\big] \Big\|_{C_t^0([0,T]\rightarrow \frkG)} > C^{RT} \bigg).
\label{main:eq-refined-q-c4}
\end{align}
Using the triangle inequality, the distances between   
$( \widetilde{A}^{(\Nscript,\Mdscript,\Rscript,\coup)}, \widetilde{B}^{(\Nscript,\Mdscript,\Rscript,\coup)})$
and $( \widetilde{A}^{(\Nscript,\Ndscript,\Rscript,\coup)}, \widetilde{B}^{(\Nscript,\Ndscript,\Rscript,\coup)})$
in \eqref{main:eq-refined-q-c4} can be controlled using their distances to
$(A^{(\Nscript,\Rscript,\coup)}, B^{(\Nscript,\Rscript,\coup)})$, which is defined as the 
unique global solution \eqref{main:eq-WM-N-nocut} with initial data $(W_0^{(\Rscript,\coup)},W_1^{(\Rscript,\coup)})$. 
Provided that $C$ is sufficiently large depending on $C_4$, it then follows from  
Definition \ref{main:def-gwp-small}, Lemma \ref{main:lem-gwp-small}, and  \eqref{main:eq-refined-q-a2} that 
\begin{equation}\label{main:eq-refined-q-c5}
\eqref{main:eq-refined-q-c4} 
\leq \limsup_{N\rightarrow \infty} \mathbb{P} \big(  \revision{\Omega \backslash \GWPN_{\Ld}(T)} \big)
\leq  C T \exp \big( - c R^{-2\eta} \coup^{-1} \big) + C \Ld^{-10}.
\end{equation}
By combining \eqref{main:eq-refined-q-c4} and \eqref{main:eq-refined-q-c5}, we then obtain the desired estimate \eqref{main:eq-refined-q-c2}.
It remains to prove the Wasserstein estimate \eqref{main:eq-refined-q-c3}. To this end, we first note that 
\begin{equation*}
(\varphi_R \otimes \varphi_R)_\# \mu^{(\coup)} 
= (\varphi_R \otimes \varphi_R)_\# \mu^{(\Rscript,\coup)},
\end{equation*}
which follows directly from Definition \ref{prelim:def-white-noise-real-valued} and Definition \ref{prelim:def-white-noise-g-valued}.
By first using \eqref{main:eq-refined-q-a5} and then using \eqref{main:eq-refined-q-b1}, we then obtain that 
\begin{align*}
&\, \WassersteinR \Big( \Law \big( \varphi_R A^{(\Ndscript,\coup)}(t), \varphi_R B^{(\Ndscript,\coup)}(t) \big), 
(\varphi_R \otimes \varphi_R)_\# \mu^{(\coup)} \Big) \\ 
=&\,\WassersteinR \Big( \Law \big( \varphi_R \widetilde{A}^{(\Ndscript,\Rscript,\coup)}(t), \varphi_R \widetilde{B}^{(\Ndscript,\Rscript,\coup)}(t) \big), 
(\varphi_R \otimes \varphi_R)_\# \mu^{(\Rscript,\coup)} \Big) \\
\leq&\, \limsup_{N\rightarrow \infty}\WassersteinR \Big( \Law \big( \varphi_R \widetilde{A}^{(\Nscript,\Ndscript,\Rscript,\coup)}(t), \varphi_R \widetilde{B}^{(\Nscript,\Ndscript,\Rscript,\coup)}(t) \big), 
(\varphi_R \otimes \varphi_R)_\# \mu^{(\Rscript,\coup)} \Big). 
\end{align*}
After using Lemma \ref{prelim:lem-Wasserstein-properties}, Lemma \ref{main:lem-gwp-small}, and \eqref{main:eq-refined-q-a3}, we then obtain the desired estimate \eqref{main:eq-refined-q-c3}.\\

\emph{Fourth step: General temperature.} 
In the fourth step, we use the scaling symmetry (Lemma \ref{ansatz:lem-scaling-symmetry}) to remove the condition $\coup \leq \coup_0$. To this end, let $\kappa \in \inversedyadic$ satisfy
$0<\kappa \leq \kappa_0$, where $\kappa_0=\kappa_0(\coup,T,R)$ is sufficiently small, and let $\Kd=\Kd(\kappa,\coup,T,R)$ be sufficiently large.
We now claim for all $\Ld \geq \Kd$ that
\begin{equation}\label{main:eq-refined-q-d1}
\begin{aligned}
&\mathbb{P} \bigg( \sup_{\substack{\Nd \in \Dyadiclarge\colon \\ \Nd \geq \Ld}} 
\Big\| \varphi_R \big( A^{(\Ndscript,\coup)}, B^{(\Ndscript,\coup)} \big); 
\big( 0, 0 \big) \Big\|_{C_t^0 \Xi^{s-1}_x([0,T]\times \R \rightarrow \frkg^2)}
> C \kappa^{-10} C^{\frac{2RT}{\kappa^2}}, \\ 
&\hspace{3ex} \sup_{\substack{\Md,\Nd \in \Dyadiclarge\colon \\ \Ld \leq \Md \leq \Nd}} 
\Md^\varsigma \Big\| \varphi_R \big( A^{(\Mdscript,\coup)}, B^{(\Mdscript,\coup)} \big); 
\varphi_R  \big( A^{(\Ndscript,\coup)}, B^{(\Ndscript,\coup)} \big) \Big\|_{C_t^0 \Xi^{s-1}_x([0,T]\times \R \rightarrow \frkg^2)}
> C \kappa^{-10} C^{\frac{2RT}{\kappa^2}},  \\ 
&\hspace{3ex}\textup{or} \quad \sup_{\substack{\Md,\Nd\in\Dyadiclarge\colon \\  \Ld \leq \Md \leq \Nd}} \Md^{\varsigma}
\Big\| \ac\big[ A^{(\Mdscript,\coup)}(\cdot,0)\big] ; 
\ac\big[ A^{(\Ndscript,\coup)}(\cdot,0)\big] \Big\|_{C_t^0([0,T]\rightarrow \frkG)} > C \kappa^{-10} C^{\frac{2RT}{\kappa^2}} \bigg) \\
\leq&\, C T \kappa^{-1} \exp \Big( - c R^{-2\eta} \coup^{-1} \kappa^{-(1-2\eta)} \Big) + C \kappa^{-10} \Ld^{-10}
\end{aligned}
\end{equation}
and that, for all $\Nd\geq \Kd$, it holds that
\begin{equation}\label{main:eq-refined-q-d2}
\begin{aligned}
&\, \WassersteinR \Big( \Law \big( \varphi_R A^{(\Ndscript,\coup)}(t), \varphi_R B^{(\Ndscript,\coup)}(t) \big), 
(\varphi_R \otimes \varphi_R)_\# \mu^{(\coup)} \Big) \\
\leq&\,  C \kappa^{-10} \Big( C^{\frac{2RT}{\kappa^2}} \kappa^{-\varsigma} \Nd^{-\varsigma} + C T \kappa^{-1} \exp \Big( - c R^{-2\eta} \coup^{-1} \kappa^{-(1-2\eta)} \Big) \Big).
\end{aligned}
\end{equation}
In order to obtain \eqref{main:eq-refined-q-d1} and \eqref{main:eq-refined-q-d2}, we recall from Lemma \ref{ansatz:lem-scaling-symmetry} that  $(\Scaling_\kappa A^{(\Ndscript,\coup)}, \Scaling_\kappa B^{(\Ndscript,\coup)})$ solves \eqref{main:eq-WM-nocut} with $(\Nd,\coup)$ replaced by $(\kappa \Nd,\kappa \coup)$ and with $(W_0^{(\coup)},W_1^{(\coup)})$ replaced by 
\begin{equation*}
\big( \Scaling_\kappa  W_0^{(\coup)}, \Scaling_\kappa W_1^{(\coup)} \big),
\end{equation*}
which is a pair of independent, $\frkg$-valued white noises at temperature $8 \kappa \coup$. Provided  that $\kappa\leq \kappa_0$, we also have that  $\kappa \coup \leq \coup_0$. 
We can therefore apply our estimates \eqref{main:eq-refined-q-c2} and \eqref{main:eq-refined-q-c3} from the third step to 
$(\Scaling_\kappa A^{(\Ndscript,\coup)}, \Scaling_\kappa B^{(\Ndscript,\coup)})$. Then, we write 
\begin{align}\label{main:eq-refined-q-d3} 
\Big( A^{(\Ndscript,\coup)},  B^{(\Ndscript,\coup)} \Big) (t,x) 
&= \Big( (\Scaling_\kappa)^{-1} \Scaling_\kappa A^{(\Ndscript,\coup)}, (\Scaling_\kappa)^{-1} \Scaling_\kappa B^{(\Ndscript,\coup)}\Big)(t,x), \\  
\ac\big[A^{(\Ndscript,\coup)}(\cdot,0)\big](t) 
&= \ac\big[\Scaling_\kappa A^{(\Ndscript,\coup)}(\cdot,0)\big](\kappa^{-1} t).
\label{main:eq-refined-q-d4} 
\end{align}
Due to Lemma \ref{prelim:lem-Wasserstein-properties}.\ref{prelim:item-Wasserstein-Lipschitz}, \eqref{main:eq-refined-q-d3}, \eqref{main:eq-refined-q-d4}, and Lemma \ref{prelim:lem-scaling-estimate}, we can then transfer our estimates from $(\Scaling_\kappa A^{(\Ndscript,\coup)}, \Scaling_\kappa B^{(\Ndscript,\coup)})$ to $(A^{(\Ndscript,\coup)}, B^{(\Ndscript,\coup)})$, which yields  \eqref{main:eq-refined-q-d1} and \eqref{main:eq-refined-q-d2}.\\

\emph{Fifth step: Conclusion.} Using \eqref{main:eq-refined-q-d1} and \eqref{main:eq-refined-q-d2}, we now complete the proof of the proposition. We first prove the convergence of $(A^{(\Ndscript,\coup)}, B^{(\Ndscript,\coup)})$. Using \eqref{main:eq-refined-q-d1}, it follows that
\begin{align*}
&\, \mathbb{P} \bigg( \limsup_{\Md,\Nd\rightarrow \infty}  \Big\| \varphi_R \big( A^{(\Mdscript,\coup)}, B^{(\Mdscript,\coup)} \big); 
\varphi_R  \big( A^{(\Ndscript,\coup)}, B^{(\Ndscript,\coup)} \big) \Big\|_{C_t^0 \Xi^{s-1}_x([0,T]\times \R \rightarrow \frkg^2)} >0 \bigg) \\ 
\leq & \, \liminf_{\Ld \rightarrow \infty} \mathbb{P} \bigg( \sup_{\substack{\Md,\Nd \in \Dyadiclarge\colon \\ \Ld \leq \Md \leq \Nd}} 
\Md^\varsigma \Big\| \varphi_R \big( A^{(\Mdscript,\coup)}, B^{(\Mdscript,\coup)} \big); 
\varphi_R  \big( A^{(\Ndscript,\coup)}, B^{(\Ndscript,\coup)} \big) \Big\|_{C_t^0 \Xi^{s-1}_x([0,T]\times \R \rightarrow \frkg^2)}
> C \kappa^{-10} C^{\frac{2RT}{\kappa^2}} \bigg) \\
\leq&\,\liminf_{\Ld \rightarrow \infty} \Big( C T \kappa^{-1} \exp \Big( - c R^{-2\eta} \coup^{-1} \kappa^{-(1-2\eta)} \Big) + C \kappa^{-10} \Ld^{-10} \Big) 
\leq \, C T \kappa^{-1} \exp \Big( - cR^{-2\eta} \coup^{-1} \kappa^{-(1-2\eta)}\Big).
\end{align*}
By letting $\kappa\rightarrow 0$, it then follows that 
\begin{equation*}
\mathbb{P} \bigg( \limsup_{\Md,\Nd\rightarrow \infty}  \Big\| \varphi_R \big( A^{(\Mdscript,\coup)}, B^{(\Mdscript,\coup)} \big); 
\varphi_R  \big( A^{(\Ndscript,\coup)}, B^{(\Ndscript,\coup)} \big) \Big\|_{C_t^0 \Xi^{s-1}_x([0,T]\times \R \rightarrow \frkg^2)} >0 \bigg) 
=0, 
\end{equation*}
which implies the existence of the limit $(A^{(\coup)},B^{(\coup)})$. The convergence of $\ac\big[ A^{(\Ndscript,\coup)}(\cdot,0)\big]$ can be obtained from \eqref{main:eq-refined-q-d1} in a similar fashion. Furthermore, for any $T,R\geq 1$ and $t\in [0,T]$, it follows that 
\begin{align*}
&\, \WassersteinR \Big( \Law \big( \varphi_R A^{(\coup)}(t), \varphi_R B^{(\coup)}(t)\big), 
(\varphi_R \otimes \varphi_R)_\# \mu^{(\coup)} \Big)  \\ 
\leq &\, \liminf_{\Nd\rightarrow \infty}
\WassersteinR \Big( \Law \big( \varphi_R A^{(\Ndscript,\coup)}(t), \varphi_R B^{(\Ndscript,\coup)}(t)\big), 
(\varphi_R \otimes \varphi_R)_\# \mu^{(\coup)} \Big) \\ 
\leq&\,  C \kappa^{-10} \liminf_{\Nd \rightarrow \infty} 
 \Big( C^{\frac{2RT}{\kappa^2}} \kappa^{-\varsigma} \Nd^{-\varsigma} + C T \kappa^{-1} \exp \Big( - c R^{-2\eta} \coup^{-1} \kappa^{-(1-2\eta)} \Big) \Big) \\ 
 \leq&\,  C^{2} T \kappa^{-11} \exp \Big( - c R^{-2\eta} \coup^{-1} \kappa^{-(1-2\eta)} \Big).
\end{align*}
By letting $\kappa\rightarrow 0$, it then follows that 
\begin{equation*}
 \WassersteinR \Big( \Law \big( \varphi_R A^{(\coup)}(t), \varphi_R B^{(\coup)}(t)\big), 
(\varphi_R \otimes \varphi_R)_\# \mu^{(\coup)} \Big) =0.
\end{equation*}
Since $R\geq 1$ is arbitrary, this implies that the law of $\big(A^{(\coup)}(t),B^{(\coup)}(t)\big)$ is given by the Gibbs measure $\mu^{(\coup)}$. 
\end{proof}

Equipped with Proposition \ref{main:prop-refined-gwp}, we now readily obtain Theorem \ref{intro:thm-rigorous-A-B}. 

\begin{proof}[Proof of Theorem \ref{intro:thm-rigorous-A-B}:]  In Proposition \ref{main:prop-refined-gwp}, we considered the solutions 
$ \big( A^{(\Ndscript,\coup)},B^{(\Ndscript,\coup)} \big)$ and the initial data $W^{(\coup)}_0,W^{(\coup)}_1\colon \R \rightarrow \frkg$, which are written in the notation of Sections \ref{section:ansatz}-\ref{section:Lipschitz}. After relabelling $\Nd$ as $N$ and recalling from \eqref{intro:eq-coup} that $\coup=(8\beta)^{-1}$, the solutions and initial data in Proposition \ref{main:prop-refined-gwp} correspond to the solution $(A_{\leq N},B_{\leq N})$ and initial data $(\beta^{-\frac{1}{2}} W_0,\beta^{-\frac{1}{2}} W_1)$ from the statement of Theorem \ref{intro:thm-rigorous-A-B}. Since the $C_t^0 \Xi_x^{s-1}$-distance from Proposition \ref{main:prop-refined-gwp} controls the $C_t^0 \C_x^{s-1}$-norm from the statement of Theorem \ref{intro:thm-rigorous-A-B}, Proposition \ref{main:prop-refined-gwp} therefore implies Theorem \ref{intro:thm-rigorous-A-B}.
\end{proof}

\begin{remark}[Flow property]\label{main:rem-flow}
As previously mentioned in Remark \ref{intro:rem-properties}, a modification of our arguments can likely be used to prove that the limit from
Theorem \ref{intro:thm-rigorous-A-B} satisfies the flow property. We now briefly discuss the flow property and necessary modifications of our argument, but we do not provide a full proof. 
We first let $\big( \mathbb{A},\mathbb{B}\big)(t)$ be the flow corresponding to smooth solutions of \eqref{main:eq-WM-nocut}. Then, we define a set 
\begin{equation}\label{main:eq-flow-e1}
\begin{aligned}
\Sigma:=\Big\{ &\big( A_0, B_0 \big) \colon \textup{ For all $k\geq 1$ and $T,R\geq 1$, the limit of} \\ 
&
\big( \mathbb{A}, \mathbb{B} \big)(t_1) P^x_{\leq \Nd^{(1)}}  \big( \mathbb{A}, \mathbb{B} \big)(t_2) P^x_{\leq \Nd^{(2)}} \hdots  
\big( \mathbb{A}, \mathbb{B} \big)(t_k) P^x_{\leq \Nd^{(k)}} \big( A_0, B_0 \big)  \\ 
&\textup{exists in } C_{t_1,\hdots,t_k}^0 \C_x^{s-1} \big( [-T,T]^k \times [-R,R] \rightarrow \frkg^2\big) \textup{ as } 
\Nd^{(1)}, \Nd^{(2)},\hdots, \Nd^{(k)} \rightarrow \infty 
\Big\}. 
\end{aligned}
\end{equation}
From the definition of $\Sigma$, it directly follows that the flow $\big(\mathbb{A},\mathbb{B}\big)(t)$ can be extended to the set $\Sigma$ via $\big(\mathbb{A},\mathbb{B}\big)(t) := \lim_{\Nd \rightarrow \infty}  \big(\mathbb{A},\mathbb{B}\big)(t)P_{\leq \Nd}^x$. Furthermore, it directly follows that the flow satisfies the self-mapping property $\big(\mathbb{A},\mathbb{B}\big)(t) \Sigma \subseteq \Sigma$ and 
the flow property 
\begin{equation*}
\big(\mathbb{A},\mathbb{B}\big)(t_1) \big(\mathbb{A},\mathbb{B}\big)(t_2) = \big(\mathbb{A},\mathbb{B}\big)(t_1+t_2). 
\end{equation*}
In order to show that the limiting solution in Theorem \ref{intro:thm-rigorous-A-B} is induced by a flow, it therefore remains to show that the $\frkg$-valued white noise $\big(W_0^{(\coup)},W_1^{(\coup)}\big)$ from the proof of Proposition \ref{main:prop-refined-gwp} is an element of $\Sigma$. 
While the statement of Proposition \ref{main:prop-refined-gwp} only yields the existence of the limit in \eqref{main:eq-flow-e1} for $k=1$, the general case $k\geq 1$ can likely be treated using a modified version of the proof of Proposition \ref{main:prop-refined-gwp}. In this modified version, it is important to approximate 
\begin{equation}\label{main:eq-flow-e2}
\big( \mathbb{A}, \mathbb{B} \big)(t_1) P^x_{\leq \Nd^{(1)}}  \big( \mathbb{A}, \mathbb{B} \big)(t_2) P^x_{\leq \Nd^{(2)}} \hdots  
\big( \mathbb{A}, \mathbb{B} \big)(t_k) P^x_{\leq \Nd^{(k)}} \big( W_0^{(\coup)},W_1^{(\coup)} \big) 
\end{equation}
by 
\begin{equation}\label{main:eq-flow-e3}
\big( \mathbb{A}^{(\Nscript,\coup)}, \mathbb{B}^{(\Nscript,\coup)} \big)(t_1) 
\big( \mathbb{A}^{(\Nscript,\coup)}, \mathbb{B}^{(\Nscript,\coup)} \big)(t_2) 
\hdots 
\big( \mathbb{A}^{(\Nscript,\coup)}, \mathbb{B}^{(\Nscript,\coup)} \big)(t_k) 
\big( W_0^{(\Rscript,\coup)},W_1^{(\Rscript,\coup)} \big),
\end{equation}
where $N$ is much bigger than $\Nd^{(1)},\hdots, \Nd^{(k)}$ and $\big( \mathbb{A}^{(\Nscript,\coup)}, \mathbb{B}^{(\Nscript,\coup)} \big)(t)$ 
is the flow of \eqref{main:eq-WM-N-nocut}. Since \eqref{main:eq-flow-e3} equals $\big( \mathbb{A}^{(\Nscript,\coup)}, \mathbb{B}^{(\Nscript,\coup)} \big)(t_1+t_2+\hdots+t_k) \big( W_0^{(\Rscript,\coup)},W_1^{(\Rscript,\coup)} \big)$, 
it is then clear that \eqref{main:eq-flow-e3} still leaves the Gibbs measure almost invariant and still satisfies the estimates from the Bourgain-Bulut argument. 

\revision{We also mention that the flow-property is not obtained in the standard Bourgain-Bulut argument \cite{BB14}. This is due to an important difference between the argument in this article and in \cite{BB14}, which was previously mentioned in Remark \ref{intro:rem-Bourgain-Bulut}. In our argument, we obtain the convergence of the solutions with frequency-truncated initial data, which appear in \eqref{main:eq-flow-e1}. In \cite{BB14}, however, the authors only obtain the convergence of the finite-dimensional approximations from~\cite[(1.2)]{BB14}. Of course, one can try to introduce a version of \eqref{main:eq-flow-e1} involving the flows of the finite-dimensional approximations, but this may cause problems. As far as the author can tell, the arguments in \cite{BB14} do not yield control on the push-forward of the frequency-truncated Gibbs measures from \cite[(1.3)]{BB14} under flows with different frequency-truncation parameters. That is, in the notation\footnote{\revision{To make the underlying function spaces in \eqref{main:eq-flow-e4} consistent, the frequency-projection $P_N$ in \cite[(1.2)]{BB14} should appear in the nonlinearity, and not in the initial data. That is, the high frequency initial data should be evolved linearly. Similarly, the high frequencies in the Gibbs measure from \cite[(1.3)]{BB14} should be drawn from the limiting free measure $\mu_F$.}} of \cite[Proposition 3.4]{BB14}, the argument in \cite{BB14} does not control push-forwards of the form
\begin{equation}\label{main:eq-flow-e4}
V_{N_1}(t_1)_\# V_{N_2}(t_2)_{\#} \hdots V_{N_k}(t_k)_{\#} \mu^{(N)}_G, 
\end{equation}
where $N,N_1,N_2,\hdots, N_k \in \dyadic$ are distinct.}
\end{remark}

\subsection{Proof of Theorem \ref{intro:thm-rigorous-phi}}\label{section:lifting-theorem}

In the last subsection, we obtained the global well-posedness and invariance of the Gibbs measure at the level of derivatives (Proposition \ref{main:prop-refined-gwp}). In Subsection \ref{section:lifting-lifts} and \mbox{Subsection \ref{section:lifting-Gibbs}} above, we also studied the lifts of $\frkg$-valued maps $A$ and $B$ to $\frkG$-valued maps $\phi$. By combining both ingredients, we can now prove Theorem \ref{intro:thm-rigorous-phi}.

\begin{proof}[Proof of Theorem \ref{intro:thm-rigorous-phi}:] 
In order to better match the notation from Section \ref{section:ansatz}-\ref{section:Lipschitz}, we denote the frequency-truncation parameter from the statement of Theorem \ref{intro:thm-rigorous-phi} by $\Nd\in \Dyadiclarge$ instead of $N\in \Dyadiclarge$. We also recall that the temperature $\coup>0$ used in this section and the inverse temperature $\beta>0$ from the statement of Theorem \ref{intro:thm-rigorous-phi} are related via the identity $\coup=(8\beta)^{-1}$ from \eqref{intro:eq-coup}. 

Let $g\in \frkG$ and $\phi_{\leq \Nd}\colon \R^{1+1}\rightarrow \frkG$ be as in the statement of the theorem and let $A_{\leq \Nd},B_{\leq \Nd}\colon \R^{1+1}\rightarrow \frkG$ be defined as
\begin{equation}\label{lifting:eq-theorem-1}
A_{\leq \Nd} = \phi_{\leq \Nd}^{-1} \partial_t \phi_{\leq \Nd} 
\qquad \text{and} \qquad
B_{\leq \Nd} = \phi_{\leq \Nd}^{-1} \partial_x \phi_{\leq \Nd}.
\end{equation}
Using Lemma \ref{lifting:lem-equivalence}, it follows that $A_{\leq \Nd}$ and $B_{\leq \Nd}$ solve the derivative formulation of the wave maps equation, i.e., $A_{\leq \Nd}$ and $B_{\leq \Nd}$ solve the initial value problem
\begin{equation}\label{lifting:eq-theorem-2}
\begin{cases}
\begin{aligned}
&\partial_t A_{\leq \Nd} = \partial_x B_{\leq \Nd}, \\ 
&\partial_t B_{\leq \Nd} = \partial_x A_{\leq \Nd} - \big[ A_{\leq \Nd}, B_{\leq \Nd} \big], \\
&A_{\leq \Nd}(0,\cdot) = \beta^{-\frac{1}{2}} P_{\leq \Nd} W_0, \, \, B_{\leq \Nd}(0,\cdot) = \beta^{-\frac{1}{2}} P_{\leq \Nd} W_1. 
\end{aligned}
\end{cases}
\end{equation}
We note that, in the notation of Proposition \ref{main:prop-refined-gwp}, $A_{\leq \Nd}$ and $B_{\leq \Nd}$ 
are written as $A^{(\Ndscript,\coup)}$ and $B^{(\Ndscript,\coup)}$. 
Using Proposition \ref{main:prop-refined-gwp}, it follows that the limit 
\begin{equation}\label{lifting:eq-theorem-3-a}
(A,B)(t,x) = \lim_{\Nd\rightarrow \infty} (A_{\leq \Nd},B_{\leq \Nd})(t,x) 
\end{equation}
exists in $C_t^0 \Xi_x^{s-1}([0,T]\times [-R,R]\rightarrow \frkg^2)$ for all $T,R\geq 1$ and that the limit
\begin{equation}\label{lifting:eq-theorem-3-b}
\ac(t) = \lim_{\Nd \rightarrow \infty} \ac \big[ A_{\leq \Nd}(\cdot,0)](t)
\end{equation}
exists in $C_t^0([0,T]\rightarrow \frkG)$ for all $T\geq 1$.
Furthermore, it also follows from Proposition \ref{main:prop-refined-gwp} that 
\begin{equation}\label{lifting:eq-theorem-5}
\operatorname{Law} \Big( \big(A(t),B(t)\big)\Big) = \mu_\beta
\end{equation}
for all $t\in \R$.
Using Lemma \ref{lifting:lem-representation}, we can write the wave map $\phi_{\leq \Nd}$ and its derivative $\partial_t \phi_{\leq \Nd}$ as  
\begin{align}
\phi_{\leq \Nd}(t,x) &= g \cdot \ac\big[ A_{\leq \Nd}(\cdot,0)\big](t) \cdot \bc\big[ B_{\leq \Nd}(t,\cdot)\big](x), \label{lifting:eq-theorem-6} \\
\partial_t \phi_{\leq \Nd}(t,x) &= g \cdot \ac\big[ A_{\leq \Nd}(\cdot,0)\big](t) \cdot \bc\big[ B_{\leq \Nd}(t,\cdot)\big](x) \cdot A_{\leq \Nd}(t,x).  \label{lifting:eq-theorem-7}
\end{align}
Due to \eqref{lifting:eq-theorem-3-a}, \eqref{lifting:eq-theorem-3-b}, \eqref{lifting:eq-theorem-6}, and \eqref{lifting:eq-theorem-7}, we can apply Lemma \ref{lifting:lem-lift-properties} to $\phi_{\leq \Nd}$ and $\partial_t \phi_{\leq \Nd}$, which yields the desired claim regarding the convergence of $(\phi_{\leq \Nd},\partial_t \phi_{\leq \Nd})$ to a limit $(\phi,\partial_t \phi)$ and also yields the identities 
\begin{align}
\phi(t,x) &= g \cdot \ac(t) \cdot \bc\big[ B(t,\cdot)\big](x), \label{lifting:eq-theorem-8} \\
\partial_t \phi(t,x) &= g \cdot \ac(t) \cdot \bc\big[ B(t,\cdot)\big](x) \cdot A(t,x).  \label{lifting:eq-theorem-9}
\end{align}
The claim regarding the convergence of $\phi_{\leq \Nd}^{-1} \partial_t \phi_{\leq \Nd}$ follows directly from \eqref{lifting:eq-theorem-1} and the convergence of $A_{\leq \Nd}$. It therefore only remains to prove for all $t\in \R$ that the limit satisfies
\begin{equation}\label{lifting:eq-theorem-10}
\operatorname{Law}\Big( \big( \phi(t),\partial_t \phi(t) \big) \Big) = \muGb.
\end{equation}
However, this follows directly from the invariance from \eqref{lifting:eq-theorem-5}, the representations from \eqref{lifting:eq-theorem-8}-\eqref{lifting:eq-theorem-9}, and Lemma \ref{lifting:lem-Gibbs-representation}. 
\end{proof}

\begin{appendix}\label{section:appendix}

\section{Auxiliary analytical estimates}\label{section:appendix-auxiliary}

In this section, we state and prove several analytical estimates which are either more technical or used less frequently in this article than the analytical estimates from Subsection \ref{section:function-spaces}.

\subsection{Additional para-product, integral, and commutator estimates}
The first lemma allows us to ignore the precise form in which we restrict to high$\times$high$\rightarrow$low-interactions. 

\begin{lemma}[\protect{Insertion of $\Para[v][sim]$}]\label{prelim:lem-insertion-parasim}
Let $K,L\in \dyadic$ satisfy $K\sim L$ and let $\alpha_1,\alpha_2,\beta_1,\beta_2,\gamma_1,\gamma_2\in (-1,1)\backslash \{0\}$ be regularity parameters satisfying 
\begin{equation*}
\gamma_1 \leq \min (\alpha_1, \beta_1) \qquad \text{and} \qquad \alpha_1 + \beta_1 >0. 
\end{equation*}
Then, it holds that
\begin{equation}\label{prelim:eq-insertion-parasim}
\big\| P_K^v f \, P_L^v g - P_K^v f \Para[v][sim] P_L^v g \big\|_{\Cprod{\gamma_1}{\gamma_2}} 
\lesssim K^{\gamma_2-\alpha_2-\beta_2} \big\|  f \big\|_{\Cprod{\alpha_1}{\alpha_2}}
\big\|  g   \big\|_{\Cprod{\beta_1}{\beta_2}}. 
\end{equation}
\end{lemma}

\begin{proof}
Due to the definition of the para-product operators, it holds that 
\begin{equation}\label{prelim:eq-insertion-p1}
P_K^v f \, P_L^v g - P_K^v f \Para[v][sim] P_L^v g 
=  P_K^v f \Para[v][ll] P_L^v g + P_K^v f \Para[v][gg] P_L^v g.
\end{equation}
We only estimate the first summand in \eqref{prelim:eq-insertion-p1}, since the second summand can be controlled using a similar estimate. For any $0<\epsilon\ll 1$, the low$\times$high-estimate (Lemma \ref{prelim:lem-paraproduct}) implies 
\begin{align*}
\big\|  P_K^v f \Para[v][ll] P_L^v g  \big\|_{\Cprod{\gamma_1}{\gamma_2}}
\lesssim \big\| P_K^v f \big\|_{\Cprod{\alpha_1}{-\epsilon}}
\big\| P_L^v g \big\|_{\Cprod{\beta_1}{\gamma_2+\epsilon}} 
\lesssim K^{-\epsilon-\alpha_2} L^{\gamma_2+\epsilon-\beta_2} 
\big\|  f \big\|_{\Cprod{\alpha_1}{\alpha_2}}
\big\|  g   \big\|_{\Cprod{\beta_1}{\beta_2}} . 
\end{align*}
Since $K\sim L$, this yields the desired estimate \eqref{prelim:eq-insertion-parasim}. 
\end{proof}

\begin{lemma}[\protect{Commutator estimate involving $\Para[v][sim]$}]\label{prelim:lem-chi-commutator}
Let $f,g,h\in C^\infty_b(\R^{1+1})$. Then, it holds that 
\begin{align*}
\Big\| \big( fg \big) \Para[v][sim] h - f \Big( g \Para[v][sim] h \Big) \Big\|_{\Cprod{r-1}{r-1}}
&\lesssim \big\| f \big\|_{\Cprod{s}{1-\eta}} \big\| g \big\|_{\Cprod{r-1}{r-1}} \big\| h \big\|_{\Cprod{s}{r-1}}, \\ 
\Big\| \big( fg \big) \Para[v][sim] h - f \Big( g \Para[v][sim] h \Big) \Big\|_{\Cprod{s}{\scrr-1}}
&\lesssim \big\| f \big\|_{\Cprod{s}{s}} \big\| g \big\|_{\Cprod{s}{s}} \big\| h \big\|_{\Cprod{s}{s-1}}. 
\end{align*}
\end{lemma}

Similar to Lemma \ref{prelim:lem-para-product-trilinear}, this follows from a modification of \cite[Lemma 2.4]{GIP15} and we omit the details.

\begin{lemma}[Integration by parts in high$\times$high$\rightarrow$low-interactions]\label{prelim:lem-integration-by-parts}
Let $\alpha_1,\alpha_2,\beta_1,\beta_2,\gamma_1,\gamma_2\in (-1,1)\backslash \{0\}$ be regularity parameters satisfying
\begin{equation}\label{prelim:eq-integration-by-parts-condition}
\begin{aligned}
\gamma_1\leq \min(\alpha_1,\beta_1),  \qquad \alpha_1 +\beta_1 >0, \qquad \text{and} \qquad 
\gamma_2\leq \alpha_2 + \beta_2 +1. 
\end{aligned}
\end{equation}
Then, it holds that 
\begin{equation}
\big\| \Int^{v}_{u\rightarrow v} f \Para[v][sim] g + f \Para[v][sim] \Int^{v}_{u\rightarrow v} g \big\|_{\Cprod{\gamma_1}{\gamma_2}}
\lesssim \big\| \Int^{v}_{u\rightarrow v} f \big\|_{\Cprod{\alpha_1}{\alpha_2+1}}
\big\| \Int^{v}_{u\rightarrow v} g   \big\|_{\Cprod{\beta_1}{\beta_2+1}}. 
\end{equation}
\end{lemma}

\begin{proof}
Using the product rule, it holds that 
\begin{align*}
  \Int^{v}_{u\rightarrow v} f \Para[v][sim] g + f \Para[v][sim] \Int^{v}_{u\rightarrow v} g
  &=  \Int^{v}_{u\rightarrow v} f \Para[v][sim]   \partial_v \Int^{v}_{u\rightarrow v}  g +  \partial_v \Int^{v}_{u\rightarrow v}   f \Para[v][sim] \Int^{v}_{u\rightarrow v} g \\
  &= \partial_v \big( \Int^{v}_{u\rightarrow v} f \Para[v][sim] \Int^{v}_{u\rightarrow v} g  \big). 
\end{align*}
Due to the conditions in \eqref{prelim:eq-integration-by-parts-condition}, the high$\times$high-estimate from Lemma \ref{prelim:lem-paraproduct} implies that 
\begin{align*}
\big\|  \partial_v \big( \Int^{v}_{u\rightarrow v} f \Para[v][sim] \Int^{v}_{u\rightarrow v} g  \big) \big\|_{\Cprod{\gamma_1}{\gamma_2}} 
\lesssim \big\| \Int^{v}_{u\rightarrow v} f \Para[v][sim] \Int^{v}_{u\rightarrow v} g   \big\|_{\Cprod{\gamma_1}{\gamma_2+1}} 
\lesssim \big\| \Int^{v}_{u\rightarrow v} f \big\|_{\Cprod{\alpha_1}{\alpha_2+1}}
\big\| \Int^{v}_{u\rightarrow v} g   \big\|_{\Cprod{\beta_1}{\beta_2+1}}. 
\end{align*}
This completes the proof. 
\end{proof}

\begin{lemma}[Integral commutator estimate]\label{prelim:lem-commutator-integral}
For all $f,g \in C^\infty_b(\R^{1+1})$, it holds that 
\begin{equation*}
\Big\| \Int^v_{u\rightarrow v} \big( f \Para[v][gg] g \big) -  \Int^v_{u\rightarrow v} \big( f \big) \Para[v][gg] g \Big\|_{\Cprod{s}{\scrr}}
\lesssim \big\| f \big\|_{\Cprod{s}{s-1}} \big\| g \big\|_{\Cprod{s}{s}}. 
\end{equation*}
Furthermore, it holds for all $\varphi\in C^\infty_b(\R)$ that 
\begin{equation*}
\Big\| \Int^v_{u\rightarrow v} \big( \varphi(v-u) f \big) - \varphi(v-u) \Int^v_{u\rightarrow v} f \Big\|_{\Cprod{s}{\scrr}} 
\lesssim \big\| \varphi \big\|_{\Schwartz} \big\| f \big\|_{\Cprod{s}{s-1}},
\end{equation*}
where $\| \cdot \|_{\Schwartz}$ is as in \eqref{prelim:eq-Schwartz}.
\end{lemma}

This follows essentially from the Duhamel integral estimate (Lemma \ref{prelim:lem-Duhamel-integral}) and commutator estimates (as in Lemma \ref{prelim:lem-commutator-basic}), and we omit the proof. We also refer the reader to \cite[Lemma B.2]{GIP15} for a similar estimate. 

We now record an estimate for the Cartesian integral operator. In contrast to Lemma \ref{prelim:lem-integral}, however, we measure regularity with respect to the null coordinates $u$ and $v$, i.e., with respect to coordinates that differ from the integration variable. 

\begin{lemma}[The Cartesian integral operator in null coordinates]\label{prelim:lem-Cartesian-integral}
Let $K,L\in \dyadic$ satisfy $K\ll L$ and let $\alpha_1,\alpha_2,\beta_1,\beta_2\in (-1,1)\backslash \{0\}$. If 
\begin{equation}\label{prelim:eq-Cartesian-integral-a1} 
\alpha_1+\alpha_2 \leq 1 + \beta_1 + \beta_2 
\qquad \text{and} \qquad 
\alpha_2 \leq 1+\beta_2, 
\end{equation}
then it holds that 
\begin{equation}\label{prelim:eq-Cartesian-integral-e1} 
\big\| P^x_{>1} \Int^x_{0\rightarrow x} P_K^u P_L^v f \big\|_{\Cprod{\alpha_1}{\alpha_2}} \lesssim \big\| f \big\|_{\Cprod{\beta_1}{\beta_2}}. 
\end{equation}
Furthermore, if 
\begin{equation}\label{prelim:eq-Cartesian-integral-a2}
\alpha_1 + \alpha_2 \leq 1+\beta_1+\beta_2, \qquad 
\alpha_1 \leq \min(1+\beta_1+\beta_2,1+\beta_2), \qquad \text{and} \qquad \alpha_2 \leq 2+\beta_2, 
\end{equation}
then it holds that 
\begin{equation}\label{prelim:eq-Cartesian-integral-e2} 
\big\| \big( P^x_{>1}  \Int^x_{0\rightarrow x} - \Int^v_{u\rightarrow v} \big) P_K^u P_L^v f \big\|_{\Cprod{\alpha_1}{\alpha_2}} \lesssim \big\| f \big\|_{\Cprod{\beta_1}{\beta_2}}. 
\end{equation}
\end{lemma}

\begin{remark}
We briefly describe the motivation behind \eqref{prelim:eq-Cartesian-integral-e2}. To this end, we recall from \eqref{ansatz:eq-null-coordinates} 
that $u=x-t$ and $v=x+t$ and therefore $P^u_K P^v_L f(u,v)= P^u_K P^v_L f(x-t,x+t)$. Since $K\ll L$, integration in the $x$-variable should therefore be close to integration in the second argument, i.e., in the $v$-variable. 
\end{remark}

\begin{proof}
We denote the frequency variables corresponding to the $u$ and $v$-variables by $\xi$ and $\zeta$, respectively. We now turn towards the proof of the first estimate \eqref{prelim:eq-Cartesian-integral-e2}. To this end, we first define $\Int_{K,L}$ as the Fourier integral operator with symbol 
\begin{equation*}
\frac{1}{i(\xi+\zeta)}  \widecheck{\rho}_{>1}(\xi+\zeta) \widecheck{\rho}_{K}(\xi) \widecheck{\rho}_{L}(\zeta)
= \frac{1}{i(\xi+\zeta)} \widecheck{\rho}_{K}(\xi) \widecheck{\rho}_{L}(\zeta).
\end{equation*}
Using our assumption $K\ll L$ and direct estimates, it easily follows for all $\gamma_1,\gamma_2\in (-1,1)\backslash \{0\}$ that 
\begin{equation}\label{prelim:eq-Cartesian-integral-p1}
\Big\| \Int_{K,L} f \Big\|_{\Cprod{\gamma_1}{\gamma_2}} \lesssim L^{-1} \big\| f \big\|_{\Cprod{\gamma_1}{\gamma_2}}. 
\end{equation} 
 For any $\xi,\zeta\in \R$ satisfying $\xi+\zeta\neq 0$, a direct calculation yields that
\begin{equation}\label{prelim:eq-Cartesian-integral-p2} 
\begin{aligned}
&\, P^x_{>1} \Int^x_{0\rightarrow x}  \big( e^{i\xi u} e^{i\zeta v} \big)
=P^x_{>1}  \int_{0}^x \dy \, e^{i\xi (y-t)} e^{i\zeta (y+t)} \\ 
=&\,  \frac{1}{i(\xi+\zeta)}  P^x_{>1} \Big( e^{i\xi (x-t)} e^{i\zeta (x+t)} - e^{i \xi (-t)} e^{i\zeta t} \Big)
= \frac{1}{i(\xi+\zeta)} \widecheck{\rho}_{>1}(\xi+\zeta) e^{i\xi u} e^{i\zeta v}.
\end{aligned}
\end{equation}
As a result of \eqref{prelim:eq-Cartesian-integral-p2}, we obtain the identity 
\begin{equation}\label{prelim:eq-Cartesian-integral-q2} 
\big( P^x_{>1} \Int^x_{0\rightarrow x} P_K^u P_L^v f \big)(u,v)
= \big( \Int_{K,L} f \big)(u,v). 
\end{equation}
By combining our estimate \eqref{prelim:eq-Cartesian-integral-p1} and the identity from \eqref{prelim:eq-Cartesian-integral-q2}, it then follows that 
\begin{align*}
\big\| P^x_{>1} \Int^x_{0\rightarrow x} P_K^u P_L^v f \big\|_{\Cprod{\alpha_1}{\alpha_2}}  
\lesssim K^{\alpha_1-\beta_1} L^{\alpha_2-\beta_2-1} \big\| f \big\|_{\Cprod{\beta_1}{\beta_2}}. 
\end{align*}
Due to our assumption \eqref{prelim:eq-Cartesian-integral-a1}, this implies the first estimate \eqref{prelim:eq-Cartesian-integral-e1},
and it therefore only remains to prove the second estimate \eqref{prelim:eq-Cartesian-integral-e2}. Arguing similarly as in \eqref{prelim:eq-Cartesian-integral-p2}, we obtain for all $\xi,\zeta\in \R$ satisfying $\zeta,\zeta+\xi\neq 0$ that
\begin{equation}\label{prelim:eq-Cartesian-integral-p3} 
\big( P^x_{>1}  \Int^x_{0\rightarrow x} - \Int^v_{u\rightarrow v} \big) \big( e^{i\xi u} e^{i\zeta v} \big) = 
\frac{1}{i} \Big( \frac{\widecheck{\rho}_{>1}(\xi+\zeta)}{\xi+\zeta} - \frac{1}{\zeta} \Big)  e^{i\xi u} e^{i\zeta v} 
+ \frac{1}{i\zeta} e^{i(\xi+\zeta)u}. 
\end{equation}
As a result of  \eqref{prelim:eq-Cartesian-integral-p3}, the difference of the integral operators can be decomposed into two operators 
\begin{equation}\label{prelim:eq-Cartesian-integral-q3}
\big(  P^x_{>1} \Int^x_{0\rightarrow x} - \Int^v_{u\rightarrow v} \big) P_K^u P_L^v f
= \Int_{K,L}^{(1)} f  + \Int_{K,L}^{(2)} f
\end{equation}
as follows: $\Int_{K,L}^{(1)}$ is the Fourier multiplier with symbol 
\begin{equation*}
\frac{1}{i} \Big( \frac{\widecheck{\rho}_{>1}(\xi+\zeta)}{\xi+\zeta} - \frac{1}{\zeta} \Big)  \widecheck{\rho}_{K}(\xi) \widecheck{\rho}_{L}(\zeta) = \frac{i\xi}{(\xi+\zeta)\zeta}  \widecheck{\rho}_{K}(\xi) \widecheck{\rho}_{L}(\zeta),
\end{equation*}
where we used that $\widecheck{\rho}_{>1}(\xi+\zeta)=1$ on the support of $\widecheck{\rho}_{K}(\xi) \widecheck{\rho}_{L}(\zeta)$. Furthermore, 
  $\big( \Int_{K,L}^{(2)}f \big) (u,v) = \big( \widetilde{\Int}_{K,L}^{(2)}f \big) (u,u)$, where $\widetilde{\Int}_{K,L}^{(2)} $ is the Fourier multiplier with symbol 
\begin{equation*}
 \frac{1}{i\zeta}  \widecheck{\rho}_{K}(\xi) \widecheck{\rho}_{L}(\zeta). 
\end{equation*}
From direct estimates, it follows that 
\begin{align*}
\big\| \Int_{K,L}^{(1)} f \big\|_{\Cprod{\alpha_1}{\alpha_2}} 
&\lesssim K^{\alpha_1-\beta_1+1} L^{\alpha_2-\beta_2-2} 
\big\| f \big\|_{\Cprod{\beta_1}{\beta_2}}, \\  
\big\| \Int_{K,L}^{(2)} f \big\|_{\Cprod{\alpha_1}{\alpha_2}} 
&\lesssim K^{-\beta_1} L^{\alpha_1-\beta_2-1}
\big\| f \big\|_{\Cprod{\beta_1}{\beta_2}}. 
\end{align*}
Due to our assumption \eqref{prelim:eq-Cartesian-integral-a2}, this implies \eqref{prelim:eq-Cartesian-integral-e2}.
\end{proof}

\subsection{Scaling and gluing estimate}

In the following, we state two auxiliary estimates which are used in Subsection \ref{section:main-gwp}.

\begin{lemma}[Scaling estimate]\label{prelim:lem-scaling-estimate}
Let $\kappa \in 2^{-\mathbb{N}_0}$ and let $\Scaling_\kappa$ be as in \eqref{ansatz:eq-scaling-transform}. For all $\gamma \in (-1,0)$ and all $A,\widetilde{A},B,\widetilde{B}\colon \R \rightarrow \frkg$, it then holds that 
\begin{align}
\big\| \Scaling_\kappa^{-1} A \big\|_{\C_x^{\gamma}(\R)} &\lesssim \kappa^{-10} \big\| A \big\|_{\C_x^{\gamma}(\R)}, 
\label{prelim:eq-scaling-estimate-e1} \\ 
\big\| \Scaling_\kappa^{-1}A; \Scaling_\kappa^{-1}\widetilde{A} \big\|_{\Xi_x^{s-1}(\R\rightarrow \frkg)} 
&\lesssim \kappa^{-10} \Big( 1+ \| A \|_{\C_x^{s-1}} +  \| \widetilde{A} \|_{\C_x^{s-1}} \Big)  \big\| A; \widetilde{A} \big\|_{\Xi_x^{s-1}(\R\rightarrow \frkg)}, \label{prelim:eq-scaling-estimate-e2} 
\end{align}
and 
\begin{align}
&\big\| \big(\Scaling_\kappa^{-1}A,\Scaling_\kappa^{-1}B\big); \big(\Scaling_\kappa^{-1}\widetilde{A},\Scaling_\kappa^{-1}\widetilde{B}\big) \big\|_{\Xi_x^{s-1}(\R\rightarrow \frkg)} \label{prelim:eq-scaling-estimate-e3}  \\ 
\lesssim\, & \kappa^{-10} \Big( 1+ \| A \|_{\C_x^{s-1}} +  \| B \|_{\C_x^{s-1}}+  \| \widetilde{A} \|_{\C_x^{s-1}}
\| \widetilde{B} \|_{\C_x^{s-1}}\Big)  \big\| \big( A, B \big); \big( \widetilde{A},\widetilde{B}\big) \big\|_{\Xi_x^{s-1}(\R\rightarrow \frkg)}.
\notag 
\end{align}
\end{lemma}

\begin{proof}
We observe that 
\begin{equation}\label{prelim:eq-scaling-estimate-e4}
(\Scaling_\kappa)^{-1} = \Scaling_{\kappa^{-1}}, \quad 
\Int^x_{0\rightarrow x} \Scaling_{\kappa^{-1}}  = \kappa \Scaling_{\kappa^{-1}} \Int^x_{0\rightarrow x}, 
\quad \text{and} \quad P_L \Scaling_{\kappa^{-1}} =  \Scaling_{\kappa^{-1}} P_{\kappa L}
\end{equation}
for all $L>\kappa^{-1}$. 
By using \eqref{prelim:eq-scaling-estimate-e4} for all frequency scales  $>\kappa^{-1}$ and 
using crude estimates for all frequency scales $\leq \kappa^{-1}$, one then readily obtains \eqref{prelim:eq-scaling-estimate-e1}, \eqref{prelim:eq-scaling-estimate-e2}, and \eqref{prelim:eq-scaling-estimate-e3}, and we omit the details. 
\end{proof}

\begin{lemma}[Gluing estimate]\label{prelim:lem-gluing-estimate}
Let $t_0,t_1,t_2\in \R$ satisfy $t_0 < t_1 < t_2$, let $\ac_0,\widetilde{\ac}_0 \colon [t_0,t_1] \rightarrow \frkG$, 
let $\ac_1,\widetilde{\ac}_1\colon [t_1,t_2] \rightarrow \frkG$, and assume that $\ac_1(t_1)=\widetilde{\ac}_1(t_1)=\ec$, where $\ec$ is the unit element in $\frkG$. Furthermore, let $\ac,\widetilde{\ac}\colon [t_0,t_2]\rightarrow \frkG$ be defined as 
\begin{equation*}
\ac(t) := 
\begin{cases}
\begin{aligned}
\ac_0(t) &\qquad\text{if } t\in [t_0,t_1], \\ 
\ac_0(t_1) \ac_1(t) &\qquad\text{if } t \in (t_1,t_2]
\end{aligned}
\end{cases}
\hspace{5ex}
\text{and} \hspace{5ex}
\widetilde{\ac}(t) := 
\begin{cases}
\begin{aligned}
\widetilde{\ac}_0(t) &\qquad\text{if } t\in [t_0,t_1], \\ 
\widetilde{\ac}_0(t_1) \widetilde{\ac}_1(t) &\qquad\text{if } t \in (t_1,t_2].
\end{aligned}
\end{cases}
\end{equation*}
Then, it holds that
\begin{equation*}
\big\| \ac(t); \widetilde{\ac}(t) \big\|_{C_t^0([t_0,t_2]\rightarrow \frkG)}
\lesssim \big\| \ac_0(t); \widetilde{\ac}_0(t) \big\|_{C_t^0([t_0,t_1]\rightarrow \frkG)}
+  \big\| \ac_1(t); \widetilde{\ac}_1(t) \big\|_{C_t^0([t_1,t_2]\rightarrow \frkG)}.
\end{equation*}
\end{lemma}

\begin{proof}
This follows directly from the compactness of $\frkG$, which implies uniform bounds on $\ac_0,\ac_1,\widetilde{\ac}_0$, and $\widetilde{\ac}_1$, as well as the triangle inequality. 
\end{proof}

\subsection{Weighted H\"{o}lder spaces and lattice partition}

In this subsection, we first introduce weighted H\"{o}lder norms. Then, we record estimates for weighted H\"{o}lder norms and for sums over the lattice $\Lambda$ from Definition \ref{prelim:def-lattice-partition}.

\begin{definition}[Weighted H\"{o}lder spaces]\label{prelim:def-weighted-hoelder}
Let $\varphi\colon \R \rightarrow [0,1]$ be a smooth, compactly supported function such that $\{ \varphi(\cdot-t_0)\}_{t_0\in \Z}$ is a partition of unity. Then, for all $\gamma_1,\gamma_2\in (-1,1)\backslash \{0\}$ and all $f\in C^\infty_b(\R^{1+1})$, we define the weighted H\"{o}lder norms as 
\begin{equation*}
\big\| f \big\|_{\WCprod{\gamma_1}{\gamma_2}}
:= \sum_{t_0 \in \Z} \langle t_0 \rangle^{-2} \Big\| \varphi \big( \tfrac{v-u}{2} -t_0 \big) f \, \Big\|_{\Cprod{\gamma_1}{\gamma_2}}.
\end{equation*}
\end{definition}

The next lemma states that the weighted H\"{o}lder norms behave as expected under multiplication by cut-off functions and integration. 

\begin{lemma}[Properties of weighted H\"{o}lder norms]\label{prelim:lem-weighted-hoelder-properties}
Let $\gamma_1,\gamma_2 \in (-1,1)\backslash \{0\}$, let $\zeta \in C^\infty_b(\R)$, and let $f\in C^\infty_b(\R^{1+1})$. Then, it holds that 
\begin{equation*}
\Big\| \zeta \big( \tfrac{v-u}{2} \big) f \Big\|_{\WCprod{\gamma_1}{\gamma_2}}
\lesssim \big\| \zeta\big\|_{\Schwartz} \big\| f \big\|_{\Cprod{\gamma_1}{\gamma_2}},
\end{equation*}
where $\| \cdot\|_{\Schwartz}$ is as in \eqref{prelim:eq-Schwartz}. Furthermore, if $\gamma_1\leq \gamma_2 +1$, then it holds that
\begin{equation*}
\Big\| \Int^v_{u\rightarrow v} f \Big\|_{\Cprod{\gamma_1}{\gamma_2+1}} \lesssim \big\| f \big\|_{\WCprod{\gamma_1}{\gamma_2}}.
\end{equation*}
\end{lemma}
The estimates in Lemma \ref{prelim:lem-weighted-hoelder-properties} follow from minor modifications of Lemma \ref{prelim:lem-paraproduct} and 
Lemma \ref{prelim:lem-Duhamel-integral}, and we therefore omit the proof. In the next lemma, we state estimates for sums over $\Lambda$.

\begin{lemma}\label{prelim:lem-psi-sum}
Let $\gamma_1,\gamma_2\in (-1,1)\backslash \{ 0\}$ and, for each $u_0\in \Lambda$, let $f_{u_0}\in \Cprod{\gamma_1}{\gamma_2}$. Furthermore, let $\zeta\in C^\infty_b(\R)$. Then, it holds that 
\begin{equation*}
\Big\| \sum_{u_0 \in \Lambda} \zeta\big( u-u_0 \big) f_{u_0}(u,v) \Big\|_{\Cprod{\gamma_1}{\gamma_2}}
\lesssim \| \zeta \|_{\Schwartz} \sup_{u_0 \in \Lambda} \big\| f_{u_0} \big\|_{\Cprod{\gamma_1}{\gamma_2}}
\end{equation*}
and 
\begin{equation*}
\Big\| P_{\leq N}^x \sum_{u_0 \in \Lambda} \zeta\big( u-u_0 \big) f_{u_0}(u,v) 
- \sum_{u_0 \in \Lambda} \zeta\big( u-u_0 \big) \big( P_{\leq N}^x f_{u_0}\big)(u,v) 
\Big\|_{\Cprod{\gamma_1}{\gamma_2}}
\lesssim N^{-1} \| \zeta \|_{\Schwartz} \sup_{u_0 \in \Lambda} \big\| f_{u_0} \big\|_{\Cprod{\gamma_1}{\gamma_2}}.
\end{equation*}
Similar estimates also hold if $\Lambda$ is replaced with $\LambdaR$ and $\zeta$ is replaced with its $2\pi R$-periodization 
\begin{equation*}
\zeta^{(\Rscript)}(u)=\sum_{\ell \in \Z} \zeta(u+2\pi R\ell). 
\end{equation*}
\end{lemma}
The estimate follows essentially from the proof of Lemma \ref{prelim:lem-paraproduct} and Lemma \ref{prelim:commutator-PNx}, and we therefore omit the argument.

\subsection{Estimates of sharp frequency-projections}

We next state and prove two estimates involving the sharp frequency cut-offs $(\Psharp_{R;K})_{K\in \dyadic}$ from \eqref{prelim:eq-Psharp}. 

\begin{lemma}[Almost-boundedness of $\Psharp_{R;K}$]\label{prelim:lem-PKsharp}
Let $\alpha,\beta\in \R\backslash \{0\}$ and $\epsilon>0$. Then, it holds for all $R\geq 1$, $f\colon \bT_R\rightarrow \bC$, and $K\in \dyadic$ that
\begin{equation}
\big\| \Psharp_{R;K} f \big\|_{\C_x^\alpha(\bT_R)} \lesssim (RK)^\epsilon K^{\alpha-\beta} \big\| f \big\|_{\C_x^\beta(\bT_R)}. 
\end{equation}
\end{lemma}

\begin{proof}
Since $ \Psharp_{R;K} f$ is supported on frequencies $\sim K$, it suffices to prove that
\begin{equation}\label{prelim:eq-PKsharp-p1}
\big\| \Psharp_{R;K} f \big\|_{L_x^\infty(\bT_R)} \lesssim (RK)^{\epsilon} \big\| f \big\|_{L_x^\infty(\bT_R)}.
\end{equation}
By writing $\Psharp_{R;K}=\Psharp_{R;\leq K}-\Psharp_{R;\leq K/2}$, it further suffices to prove \eqref{prelim:eq-PKsharp-p1} with $\Psharp_{R;K}$ replaced by $\Psharp_{R;\leq K}$. For each $n\in \mathbb{N}$, we let $D_n(x):= \sum_{k=-n}^n e^{ikx}$
be the $n$-th Dirichlet kernel.
From \eqref{prelim:eq-Psharp}, it directly follows that
\begin{equation}
\Psharp_{R;\leq K} f (x) 
= \frac{1}{2\pi R} \int_{\bT_R} \dy \,  D_{KR}\Big( \frac{y}{R} \Big) \Theta^x_y f.
\end{equation}
Thus, the desired estimate \eqref{prelim:eq-PKsharp-p1} follows from the standard $L^1$-bound
\begin{equation*}
\frac{1}{2\pi R} \Big\| D_{KR}\Big( \frac{y}{R} \Big) \Big\|_{L_y^1(\bT_R)}
= \frac{1}{2\pi} \Big\| D_{KR} \big( z \big) \Big\|_{L_z^1(\bT)} \lesssim \log(KR). \qedhere
\end{equation*}
\end{proof}

\begin{lemma}[Sharp Fourier-cutoffs and white noise]\label{prelim:lem-commutator-sharp}
Let $R\geq 1$  and let $W^{(\Rscript),+}:=W^{(\Rscript,1),+}$  be as in \eqref{ansatz:eq-Wpm}. Then, the following estimates hold: 
\begin{enumerate}[label=(\roman*)]
\item\label{prelim:item-sharp-bound} (Bounds on $\Psharp_{R;K}$) For all $K\in\dyadic$, $\gamma \in (-1,1)\backslash \{0\}$, $\epsilon>0$, and $p \geq 1$, 
it holds that
\begin{equation*}
\E \Big[ \Big\|  \Psharp_{R;K} W^{(\Rscript),+} \Big\|_{\C_x^\gamma}^p \Big]^{\frac{1}{p}}
\lesssim  \sqrt{p}  R^\epsilon K^{\gamma+\frac{1}{2}+\epsilon}.
\end{equation*}
\item\label{prelim:item-sharp-action} (Action on representation) For all $K,L\in \dyadic$, $\gamma \in (-1,1)\backslash \{0\}$, $\epsilon>0$, and $p \geq 1$, it holds that 
\begin{equation*}
\E \Big[ \Big\|  P_L \Psharp_{R;K} W^{(\Rscript),+} 
-  \sum_{x_0\in \LambdaRR} \sum_{k\in \Z_K} \rho_L(k) \psiRxK  G_{x_0,k}^+ e^{ikx} \Big\|_{\C_x^\gamma}^p \Big]^{\frac{1}{p}}
\lesssim \sqrt{p} R^\epsilon K^{\gamma+\frac{1}{2}+\frac{\vartheta}{2}+\epsilon} \max\big( K,L \big)^{-\frac{1}{2}}.
\end{equation*}
\item\label{prelim:item-sharp-commutator} (Commutator) 
For all $K,L,M\in \dyadic$ satisfying $K\sim L \gg M$, $\alpha,\gamma \in (-1,1)\backslash \{0\}$, $\epsilon>0$, and $p \geq 1$, it holds that 
\begin{equation*}
\E \bigg[ \sup_{\| f \|_{\C_x^{\alpha} \leq 1}} \Big\| \Com \big(  \Psharp_{R;L}, (P_M f) \big) \Psharp_{R;K} W^{(\Rscript),+} \Big\|_{\C_x^{\gamma}}^p \bigg]^{\frac{1}{p}} 
\lesssim \sqrt{p} R^\epsilon K^{\gamma+\frac{1}{2}+\frac{\vartheta}{2}+\epsilon} \max\big( K,L \big)^{-\frac{1}{2}} M^{\frac{1}{2}-\alpha}. 
\end{equation*}
\end{enumerate}
\end{lemma}

\begin{remark}
While sharp frequency-truncations often do not satisfy useful commutator estimates, the estimates in \ref{prelim:item-sharp-action} and \ref{prelim:item-sharp-commutator} are exceptions. The reason is that white noise is evenly distributed in frequency space, and thus cannot be concentrated on frequencies near the boundary of the region of truncation.
\end{remark}

\begin{proof}
Since \ref{prelim:item-sharp-bound} follows directly from Lemma \ref{killing:lem-regularity-modulated-linear} and \ref{prelim:item-sharp-action}, it suffices to prove \ref{prelim:item-sharp-action} and \ref{prelim:item-sharp-commutator}. In order to prove \ref{prelim:item-sharp-action}, it suffices to treat the case $K\sim L$, since otherwise the argument in \ref{prelim:item-sharp-action} is identically zero. We now recall from 
\eqref{ansatz:eq-Wpm} that
\begin{equation}\label{prelim:eq-sharp-p1}
W^{(\Rscript),+} = \sum_{x_0 \in \LambdaRR} \psiRx W_{x_0}, 
\qquad \text{where} \qquad W_{x_0} =   \sum_{k\in \Z} G_{x_0,k}^+ e^{ikx}. 
\end{equation}
Using the definitions of $P_L$ and $\Psharp_{R;K}$, we may then write 
\begin{align}
 &\, P_L \Psharp_{R;K} W^{(\Rscript),+} 
-  \sum_{x_0\in \LambdaRR} \sum_{k\in \Z_K} \rho_L(k) \psiRxK  G_{x_0,k}^+ e^{ikx} \notag \\
=&\, P_L \Psharp_{R;K}   \sum_{x_0\in \LambdaRR}   \psiRx W_{x_0}  - \sum_{x_0\in \LambdaRR}   \psiRxK P_L  \Psharp_{R;K}W_{x_0}  \notag \\ 
=&\, P_L \Psharp_{R;K} \sum_{x_0\in \LambdaRR} \big( \psiRx - \psiRxK \big) W_{x_0} 
\label{prelim:eq-sharp-p2} \\ 
+&\, P_L  \sum_{x_0\in \LambdaRR} \Big(  \Psharp_{R;K}   \big( \psiRxK W_{x_0} \big) - \psiRxK \big(  \Psharp_{R;K}W_{x_0} \big) \Big)
\label{prelim:eq-sharp-p3} \\
+&\,  \sum_{x_0\in \LambdaRR} \Big( P_L \big( \psiRxK  \Psharp_{R;K}W_{x_0} \big) 
-  \psiRxK P_L \big(   \Psharp_{R;K}W_{x_0} \big) \Big). 
\label{prelim:eq-sharp-p4}
\end{align}
The first and third term \eqref{prelim:eq-sharp-p2} and \eqref{prelim:eq-sharp-p4} can easily be treated using chaos estimates, Lemma \ref{prelim:lem-psi-sum}, and Lemma \ref{prelim:lem-PKsharp}, and we omit the details. In order to treat \eqref{prelim:eq-sharp-p3}, which is the hardest term, we introduce the operator $\Psharpbd_{R;K}$. The operator $\Psharpbd_{R;K}$ is a sharp frequency-projection to frequencies which are at a distance of at most $\leq 4K^{\vartheta}$ from the boundary of $[-K,-\frac{1}{2} K]\medcup [\frac{1}{2}K,K]$. Since $\psiRxK$ is supported on frequencies $\leq K^{\vartheta}$, it then follows that 
\begin{equation}\label{prelim:eq-sharp-p5}
\begin{aligned}
&\, \Psharp_{R;K}   \big( \psiRxK W_{x_0} \big) - \psiRxK \big(  \Psharp_{R;K}W_{x_0} \big) \\ 
=&\, \Psharp_{R;K}   \big( \psiRxK \Psharpbd_{R;K} W_{x_0} \big) - \psiRxK \big(  \Psharp_{R;K} \Psharpbd_{R;K}W_{x_0} \big).
\end{aligned}
\end{equation}
By first using \eqref{prelim:eq-sharp-p5} and then using Lemma \ref{prelim:lem-psi-sum} and Lemma \ref{prelim:lem-PKsharp}, we obtain that
\begin{align*}
\big\| \eqref{prelim:eq-sharp-p3}\big\|_{\C_x^\gamma} 
&\lesssim 
\Big\| \Psharp_{R;K} \sum_{x_0\in \LambdaRR}  \psiRxK \Psharpbd_{R;K} W_{x_0} \Big\|_{\C_x^\gamma}
+ \Big\| \sum_{x_0\in \LambdaRR} \psiRxK \big(  \Psharp_{R;K} \Psharpbd_{R;K}W_{x_0} \big) \Big\|_{\C_x^\gamma} \\
&\lesssim (KR)^{\frac{\epsilon}{2}} \max_{x_0\in\LambdaRR} \big\| \Psharpbd_{R;K}W_{x_0}\big\|_{\C_x^\gamma}.
\end{align*}
We now recall that $\Psharpbd_{R;K}W_{x_0}$ is supported on $\sim K^{\vartheta}$-frequencies at a distance $\sim K$ from the origin. Using chaos estimates and Lemma \ref{prelim:lem-maxima}, we then obtain that 
\begin{align*}
\E \Big[ \max_{x_0\in\LambdaRR} \big\| \Psharpbd_{R;K}W_{x_0}\big\|_{\C_x^\gamma}^p \Big]^{\frac{1}{p}}
\lesssim \sqrt{p} (KR)^{\frac{\epsilon}{2}} K^{\gamma} \big( K^\vartheta \big)^{\frac{1}{2}}, 
\end{align*}
which completes the proof of \ref{prelim:item-sharp-action}.
It remains to prove \ref{prelim:item-sharp-commutator}, which can be obtained using a similar argument as in the proof of \ref{prelim:item-sharp-action}. The only difference is that, instead of localizing to frequencies at a distance $\leq 4 K^\vartheta$ near the boundary of $[-K,-\frac{1}{2} K]\medcup [\frac{1}{2}K,K]$, the distance is now given by $4 \max(K^\vartheta,M)$.
\end{proof}

\section{Classical well-posedness theory}\label{section:classical}

In this section, we recall classical estimates for the wave maps equation. 
The classical estimates play an important role in the proof of Proposition \ref{main:prop-refined-gwp} where, for a fixed $\Nd\in \Dyadiclarge$, they are used to eliminate the frequency-truncation parameter in the nonlinearity. 
We now let $N,\Nd \in \Dyadiclarge$ and let $\coup >0$. We then consider the wave maps equation
\begin{equation}\label{classical:eq-wm}
\begin{cases}
\begin{aligned}
\partial_t A &= \partial_x B, \\ 
 \partial_t B&= \partial_x A - \big[ A , B \big]
\end{aligned}
\end{cases}    
\end{equation}
and its finite-dimensional approximation
\begin{equation}\label{classical:eq-wm-N}
\begin{cases}
\begin{aligned}
\partial_t A^{(\Nscript)} &= \partial_x B^{(\Nscript)} , \\ 
 \partial_t B^{(\Nscript)} &= \partial_x A^{(\Nscript)}  - \big[ A^{(\Nscript)}  , B^{(\Nscript)}  \big]_{\leq N} + 2\coup \Renorm[\Nscript,\Ndscript] A^{(\Nscript)}.
\end{aligned}
\end{cases}    
\end{equation}
We note that, while the solution $(A^{(\Nscript)},B^{(\Nscript)})$ depends on $\Nd$ and $\coup$, this dependence is not reflected in our notation. For all $k\in \mathbb{N}$ and $f\colon \R \rightarrow \bC$, we define local Sobolev-norms by 
\begin{equation}\label{classical:eq-Sobolev}
\big\| f \big\|_{H^k_{\loc}(\R)}^2 := \sup_{x_0\in \R} \sum_{j=0}^k  \big\| \partial_x^j f \big\|_{L^2([x_0-1,x_0+1])}^2.
\end{equation}
Equipped with \eqref{classical:eq-wm}, \eqref{classical:eq-wm-N}, and \eqref{classical:eq-Sobolev}, we can now state the main result of this section. 

\begin{proposition}[Difference estimate]\label{classical:prop-main}
Let $N,\Nd\in \dyadiclarge$, let $T,R\geq 1$, let $\coup>0$, and let $A_0,B_0 \colon \bT_R \rightarrow \frkg$ be smooth. Furthermore, 
let $(A,B)\colon \R \times \bT_R\rightarrow \frkg^2$ and $(A^{(\Nscript)},B^{(\Nscript)})\colon \R \times \bT_R\rightarrow \frkg^2$
be the unique global solutions of \eqref{classical:eq-wm} and \eqref{classical:eq-wm-N} with initial data $(A_0,B_0)$. Then, there exist parameters
\begin{equation*}
N_0=N_0\big(T,R,\coup,\|A_0\|_{H^2_{\loc}(\R)},\|B_0\|_{H^2_{\loc}(\R)}\big)
\qquad \text{and} \qquad 
C=C\big(T,R,\coup,\|A_0\|_{H^2_{\loc}(\R)},\|B_0\|_{H^2_{\loc}(\R)}\big)
\end{equation*}
such that, for all $N\geq N_0$, it holds that 
\begin{equation}\label{classical:eq-main}
\big\|  \big( A, B \big) - \big( A^{(\Nscript)},B^{(\Nscript)}\big) \big\| _{C_t^0 C_x^0([-T,T]\times \R)} 
\leq C N^{-1}.
\end{equation}
\end{proposition}

Since Proposition \ref{classical:prop-main} assumes that the initial data lies in $H^2_{\loc}$, 
its proof is well-within the scope of the classical methods in \cite{G80,GV82,KT98,LS81,MNT10,P76}. Due to the frequency-truncated Lie-bracket and Killing-renormalization in \eqref{classical:eq-wm-N}, however, it cannot be directly cited from the literature. For this reason, we decided to sketch the argument below, but omit standard details. \\ 

For all $f,g\colon \R^{1+1} \rightarrow \bC$ and $1 \leq q,p\leq \infty$, we define the local space-time norms
\begin{equation}\label{classical:eq-local-space-time}
\begin{aligned}
\big\| f \big\|_{(L_{u}^q L_v^p)_{\loc}}=\big\| f \big\|_{(L_{u}^q L_v^p)_{\loc}(\R^{1+1})}
&:= \sup_{u_0,v_0\in \R} \big\| f \big\|_{L_u^q L_v^p([u_0-1,u_0+1]\times [v_0-1,v_0+1])}, \\
\big\| g \big\|_{(L_{v}^q L_u^p)_{\loc}}=\big\| g \big\|_{(L_{v}^q L_u^p)_{\loc}(\R^{1+1})}
&:= \sup_{u_0,v_0\in \R} \big\| g \big\|_{L_v^q L_u^p([v_0-1,v_0+1]\times [u_0-1,u_0+1])}. 
\end{aligned}
\end{equation}
Equipped with \eqref{classical:eq-local-space-time}, we can now state the following nonlinear estimates.

\begin{lemma}[Nonlinear estimates]\label{classical:lem-nonlinear}
Let $f,g\colon \R^{1+1} \rightarrow \bC$, let $\chi\in C^\infty_c(\R)$, and let $\tau \in (0,1]$. Then, we have the following estimates:
\begin{align}
\big\| f g \big\|_{(L_u^2 L_v^2)_{\loc}} &\lesssim \big\| f \big\|_{(L^2_u L^\infty_v)_{\loc}} \big\| g \big\|_{(L^2_v L^\infty_u)_{\loc}}, 
\label{classical:eq-nonlinear-e1} \\ 
\big\| \chi\big( \tfrac{v-u}{\tau}\big) \,  \Int^v_{u\rightarrow v} f \big\|_{(L_u^2 L_v^\infty)_{\loc}}
&\lesssim \big\| f \big\|_{(L_u^2 L_v^1)_{\loc}}, 
\label{classical:eq-nonlinear-e2} \\ 
\big\| \chi \big( \tfrac{v-u}{\tau} \big) f \big\|_{(L_u^2 L_v^1)_{\loc}} 
&\lesssim \tau^{\frac{1}{2}} \big\| f \big\|_{(L_u^2 L_v^2)_{\loc}}, 
\label{classical:eq-nonlinear-e3} \\ 
\sup_{t\in \R} \big\| f\big(x-t,x+t\big) \big\|_{L^2_{\loc}}
&\lesssim \big\| f \big\|_{(L_u^2 L_v^\infty)_{\loc}}. 
\label{classical:eq-nonlinear-e4} 
\end{align}
Furthermore, for all $1\leq q,p \leq \infty$, it also holds that 
\begin{align}
\big\| P_{\leq N}^x  f \big\|_{(L_u^q L_v^p)_{\loc}}, 
\big\| \Renorm[\Nscript,\Ndscript] f \big\|_{(L_u^q L_v^p)_{\loc}} 
\lesssim \big\|  f \big\|_{(L_u^q L_v^p)_{\loc}}. 
\label{classical:eq-nonlinear-e5}  
\end{align}
\end{lemma}

\begin{proof}
The first estimate \eqref{classical:eq-nonlinear-e1} follows by first using H\"{o}lder's inequality and then Minkowski's integral inequality.
The second estimate \eqref{classical:eq-nonlinear-e2} follows directly from the definition of $\Int^v_{u\rightarrow v}$. The third estimate \eqref{classical:eq-nonlinear-e3} follows from H\"{o}lder's inequality. To obtain the fourth estimate \eqref{classical:eq-nonlinear-e4}, we first fix $t_0,x_0\in \R$. For $u_0:= x_0-t_0$ and $v_0:=x_0+t_0$, it then follows that
\begin{align*}
 &\, \big\| f\big(x-t_0,x+t_0\big) \big\|_{L^2_x([x_0-1,x_0+1])} \\
 \leq&\,   \Big\| \mathbf{1}\big\{ |x-t_0-u_0|\leq 1, \, |x+t_0-v_0|\leq 1\big\} f\big(x-t_0,x+t_0\big)\Big\|_{L^2_x(\R)} \\
 \leq&\, \Big\| \sup_{v\in \R} \mathbf{1}\big\{ |x-t_0-u_0|\leq 1, \, |v-v_0|\leq 1\big\}\big|f\big(x-t_0,v\big)\big|\Big\|_{L^2_x(\R)}.
\end{align*}
After using the change of variables $u:=x-t_0$, this implies \eqref{classical:eq-nonlinear-e4}. Finally, due to the translation-invariance of $(L_u^q L_v^p)_{\loc}$, it holds that
\begin{equation*}
    \sup_{y\in \R} \big\| \Theta^x_y f \big\|_{(L_u^q L_v^p)_{\loc}}
\leq \big\|  f \big\|_{(L_u^q L_v^p)_{\loc}}.
\end{equation*}
Together with Lemma \ref{ansatz:lem-renormalization}, this readily implies \eqref{classical:eq-nonlinear-e5}.
\end{proof}

In the following lemma, we state local estimates for solutions of \eqref{classical:eq-wm} and \eqref{classical:eq-wm-N}.

\begin{lemma}[Local theory]\label{classical:lem-local}
Let $N,\Nd\in \dyadiclarge$, let $R\geq 1$, and let $\coup>0$. Let $0<\tau\leq 1$, let $D\geq 1$, and
let $A_0,B_0,A_0^{(\Nscript)},B_0^{(\Nscript)}\colon \bT_R \rightarrow \frkg$. Furthermore, assume that 
\begin{equation}
\big\|  (A_0,B_0) \big\| _{L^2_{\loc}(\R)}, 
\big\|  (A_0^{(\Nscript)},B_0^{(\Nscript)})\big\| _{L^2_{\loc}(\R)}  \leq D 
\qquad \text{and} \qquad \tau \leq c D^{-2},
\end{equation}
where $c\in (0,1)$ is a sufficiently small universal constant. 
Finally, let $(A,B)$ be global solution of \eqref{classical:eq-wm} with initial data $(A_0,B_0)$
and let $(A^{(\Nscript)},B^{(\Nscript)})$ be the global solution of \eqref{classical:eq-wm-N} with initial data 
$(A_0^{(\Nscript)},B_0^{(\Nscript)})$. Then, we have the following estimates: 
\begin{enumerate}[label=(\roman*)]
\item\label{classical:item-Hk} ($H^k_{\loc}$-bounds) For all $k=0,1,2$, it holds that 
\begin{align}
\big\|  (A(t),B(t)) \big\| _{C_t^0 H_{\loc}^k([-\tau,\tau]\times \R)}
&\leq 2 \big\|  (A_0,B_0) \big\| _{H^k_{\loc}(\R)}, \label{classical:eq-hk} \\
\big\|  (A^{(\Nscript)}(t),B^{(\Nscript)}(t)) \big\| _{C_t^0 H_{\loc}^k([-\tau,\tau]\times \R)}
&\leq 2 \big\|  (A^{(\Nscript)}_0,B^{(\Nscript)}_0) \big\| _{H^k_{\loc}(\R)}. \label{classical:eq-hk-N}
\end{align}
\item\label{classical:item-difference} (Difference) It holds that 
\begin{equation}
\begin{aligned}
&\, \big\|   (A(t),B(t)) - (A^{(\Nscript)}(t),B^{(\Nscript)}(t)) \big\| _{C_t^0 H_{\loc}^1([-\tau,\tau]\times \R)} \\
\leq&\,  C_0 \big\|  (A_0,B_0) - (A^{(\Nscript)}_0,B^{(\Nscript)}_0) \big\| _{H^1_{\loc}(\R)} 
+ C_0 \big( \coup+ \big\|  (A^{(\Nscript)}_0,B^{(\Nscript)}_0) \big\| _{H^2_{\loc}(\R)} \big)^2 N^{-1}, 
\end{aligned}
\end{equation}
where $C_0$ is a universal constant. 
\end{enumerate}
\end{lemma}

\begin{proof}
We first sketch the proof of \ref{classical:item-Hk}. To this end, let $\chi=\widebar{\chi}$ be as in Definition \ref{prelim:def-cut-off}. 
After localizing to $t\in [-\tau,\tau]$ and switching to the null-coordinates from \eqref{ansatz:eq-null-coordinates} 
and \eqref{ansatz:eq-null-unknowns}, we can write \eqref{classical:eq-wm} and \eqref{classical:eq-wm-N} as 
\begin{equation}\label{classical:eq-local-q1}
\begin{cases}
\begin{aligned}
U(u,v) &= \tfrac{1}{4} \big(A_0-B_0\big)(u) 
+ \chi\big( \tfrac{v-u}{\tau}\big) \Int^v_{u\rightarrow v} \Big(  \big[ U,V \big] \Big), \\ 
V(u,v) &= \tfrac{1}{4} \big(A_0+B_0\big)(v) 
+ \chi\big( \tfrac{v-u}{\tau}\big) \Int^u_{v\rightarrow u} \Big(   \big[ U,V \big] \Big)
\end{aligned}
\end{cases}
\end{equation}
and
\begin{equation}\label{classical:eq-local-q2}
\hspace{-0.25ex}
\begin{cases}
\begin{aligned}
U^{(\Nscript)}(u,v) &= \tfrac{1}{4} \big(A^{(\Nscript)}_0-B^{(\Nscript)}_0\big)(u) 
+ \chi\big( \tfrac{v-u}{\tau}\big) \Int^v_{u\rightarrow v} \Big(  \big[ U^{(\Nscript)},V^{(\Nscript)} \big]_{\leq N}  
- \coup \Renorm[\Nscript,\Ndscript] (U^{(\Nscript)}+V^{(\Nscript)}) \Big), \\ 
V^{(\Nscript)}(u,v) &= \tfrac{1}{4} \big(A^{(\Nscript)}_0+B^{(\Nscript)}_0\big)(v) 
+ \chi\big( \tfrac{v-u}{\tau}\big) \Int^u_{v\rightarrow u} \Big(   \big[ U^{(\Nscript)},V^{(\Nscript)} \big]_{\leq N} 
- \coup \Renorm[\Nscript,\Ndscript] (U^{(\Nscript)}+V^{(\Nscript)})\Big),
\end{aligned}
\end{cases}
\end{equation}
respectively. Using Lemma \ref{classical:lem-nonlinear}, both \eqref{classical:eq-local-q1} and \eqref{classical:eq-local-q2} can be solved in $(L^2_u L^\infty_v)_{\loc}\times (L^2_v L^\infty_u)_{\loc}$  using a contraction-mapping argument. Together with \eqref{classical:eq-nonlinear-e4}, this implies the estimates in \ref{classical:item-Hk} for $k=0$. The estimates for $k=1,2$ can then be obtained using a standard persistence\footnote{The important aspect is that, even for $k=1,2$, the condition on $\tau$ only depends on the $L^2_{\loc}$-norm of the initial data.} of regularity argument. 

In order to obtain \ref{classical:item-difference}, we treat the nonlinearity in \eqref{classical:eq-wm-N} as a perturbation of the nonlinearity in \eqref{classical:eq-wm}. Using the estimate $\| P_{>N}^x \|_{H_{\loc}^2\rightarrow H_{\loc}^1}\lesssim N^{-1}$ and the algebra properties of $H^1_{\loc}$ and $H^2_{\loc}$, we easily obtain that 
\begin{equation}\label{classical:eq-local-p1}
\big\| \big[ A^{(\Nscript)}  , B^{(\Nscript)}  \big]_{\leq N} - \big[ A^{(\Nscript)}  , B^{(\Nscript)}  \big] \big\|_{L_t^2 H^1_{\loc}([-\tau,\tau]\times \R)}
\lesssim N^{-1} \big\| A^{(\Nscript)}\big\|_{C_t^0 H^2_{\loc}([-\tau,\tau]\times \R)} \big\| B^{(\Nscript)}\big\|_{C_t^0 H^2_{\loc}([-\tau,\tau]\times \R)}.
\end{equation}
Using a similar argument as in the proof of \eqref{ansatz:eq-renormalization-3} from Lemma \ref{ansatz:lem-renormalization},  we also obtain that
\begin{equation}\label{classical:eq-local-p2}
\coup \big\| \Renorm[\Nscript,\Ndscript] A^{(\Nscript)} \big\|_{L_t^2 H^1_{\loc}([-\tau,\tau]\times \R)} 
\lesssim \coup N^{-1} \big\| A^{(\Nscript)} \big\|_{C_t^0 H^2_{\loc}([-\tau,\tau]\times \R)}.
\end{equation}
Together with the $H^2_{\loc}$-bounds from \eqref{classical:eq-hk-N}, the equivalence of the local $L_t^2 L_x^2$ and $L_u^2 L_v^2$-norms, and the contraction-mapping argument described above, 
this readily yields \ref{classical:item-difference}.   
\end{proof}

While Lemma \ref{classical:lem-local} only yields local bounds, the following lemma yields global $L^2$-bounds.

\begin{lemma}[\protect{Global $L^2$-bounds}]\label{classical:lem-global}
Let $N,\Nd\in \dyadiclarge$, let $R\geq 1$, let $\coup>0$, and let $A_0,B_0 \colon \bT_R \rightarrow \frkg$ be smooth. Furthermore, 
let $(A,B)\colon \R \times \bT_R\rightarrow \frkg^2$ and $(A^{(\Nscript)},B^{(\Nscript)})\colon \R \times \bT_R\rightarrow \frkg^2$
be the unique global solutions of \eqref{classical:eq-wm} and \eqref{classical:eq-wm-N} with initial data $(A_0,B_0)$. Then, it holds 
for all $t\in \R$ that 
\begin{equation}\label{classical:eq-global}
\begin{aligned}
\big\|  \big( A(t) , B(t) \big) \big\| _{L^2(\bT_R)}
&= \big\|  \big( A_0 , B_0 \big) \big\| _{L^2(\bT_R)}, \\ 
 \big\|  \big( A^{(\Nscript)}(t) , B^{(\Nscript)}(t) \big) \big\| _{L^2(\bT_R)}
 &\leq \exp\big( C_0 \coup |t| \big) \big\|  \big( A_0 , B_0 \big) \big\| _{L^2(\bT_R)}, 
\end{aligned}
\end{equation}
where $C_0$ is a universal constant.
\end{lemma}

\begin{proof}
From a direct calculation (see Lemma \ref{structure:lem-global-flow} and Subsection \ref{section:conservative-wave-maps}), it follows that
\begin{align}
\frac{1}{2} \frac{\mathrm{d}}{\mathrm{d}t} \int_{\bT_R} \dx \Big( 
\| A(t,x) \|_{\frkg}^2 + \| B(t,x) \|_{\frkg}^2 \Big)  &=0, \label{classical:eq-global-p1} \\  
\frac{1}{2} \frac{\mathrm{d}}{\mathrm{d}t} \int_{\bT_R} \dx\,  \Big( 
\| A^{(\Nscript)}(t,x) \|_{\frkg}^2 + \| B^{(\Nscript)}(t,x) \|_{\frkg}^2 \Big) &=
2\coup \int_{\bT_R} \dx \, \big\langle B^{(\Nscript)}(t,x), \Renorm[\Nscript,\Ndscript]A^{(\Nscript)}(t,x) \big \rangle_{\frkg}.\label{classical:eq-global-p2} 
\end{align}
From Definition \ref{ansatz:def-Killing} and the estimate \eqref{ansatz:eq-renormalization-1} from Lemma \ref{ansatz:lem-renormalization}, it also follows that
\begin{equation}
\Big\|  \Renorm[\Nscript,\Ndscript] \Big\|_{L^2(\bT_R)\rightarrow L^2(\bT_R)}
\leq \sum_{M\lesssim N} \int_{\R} \dy \, \big| \big( \widecheck{\rho}_{\leq N} \ast \widecheck{\rho}_{\leq N}\big)(y)\big| 
\big| \Cf^{(\Nscript,\Ndscript)}_M(y) \big| \lesssim \sum_{M\lesssim N} M N^{-1} \lesssim 1.\label{classical:eq-global-p3}
\end{equation}
The claims in \eqref{classical:eq-global} now follow from \eqref{classical:eq-global-p1}-\eqref{classical:eq-global-p3} and Gronwall's inequality. 
\end{proof}

Equipped with Lemma \ref{classical:lem-local} and Lemma \ref{classical:lem-global}, we can now prove Proposition \ref{classical:prop-main}. 

\begin{proof}[Proof of Proposition \ref{classical:prop-main}]
Using Lemma \ref{classical:lem-global} and the trivial estimate $\|\cdot\|_{L^2_{\loc}(\R)}\leq \| \cdot\|_{L^2(\bT_R)} \lesssim R^{\frac{1}{2}} \| \cdot\|_{L^2_{\loc}(\R)}$, we obtain that
\begin{equation*}
\big\| \big( A(t),B(t) \big) \big\|_{C_t^0 L^2_{\loc}([-T,T]\times \R)}, 
\big\| \big( A^{(\Nscript)}(t),B^{(\Nscript)}(t) \big) \big\|_{C_t^0 L^2_{\loc}([-T,T]\times \R)}
\leq C\big( T,R,\coup, \big\| A_0 \big\|_{L^2_{\loc}(\R)}, \big\| B_0 \big\|_{L^2_{\loc}(\R)}\big). 
\end{equation*}
By iterating Lemma \ref{classical:lem-local}.\ref{classical:item-Hk}, it follows that 
\begin{equation*}
\big\| \big( A(t),B(t) \big) \big\|_{C_t^0 H^2_{\loc}([-T,T]\times \R)}, 
\big\| \big( A^{(\Nscript)}(t),B^{(\Nscript)}(t) \big) \big\|_{C_t^0 H^2_{\loc}([-T,T]\times \R)}
\leq C\big( T,R,\coup, \big\| A_0 \big\|_{H^2_{\loc}(\R)}, \big\| B_0 \big\|_{H^2_{\loc}(\R)}\big). 
\end{equation*}
By iterating Lemma \ref{classical:lem-local}.\ref{classical:item-difference} and using $\|\cdot\|_{C^0(\R)}\lesssim \|\cdot\|_{H^1_{\loc}(\R)}$,  we obtain the desired \mbox{estimate \eqref{classical:eq-main}}.
\end{proof}

\section{Symbolic index}\label{section:appendix-symbolic}

{\small 
\begin{longtblr}[
entry=none
]{
  colspec = {!{\vrule width 2pt}p{0.2\textwidth}p{0.6\textwidth}p{0.08\textwidth}!{\vrule width 2pt}},
  rowhead = 1,
  hlines,
  row{even} = {gray9},
  row{1} = {olive9},
} 
\textbf{Symbol} & \textbf{Description} & \textbf{Page} \\[0ex] 
$\dyadic$       & Dyadic integers 
&\begin{NoHyper}\pageref{prelim:eq-dyadiclarge}\end{NoHyper} \\[0ex] 
$\dyadiclarge$       & Large dyadic integers 
&\begin{NoHyper}\pageref{prelim:eq-dyadiclarge}\end{NoHyper} \\[0ex] 
$\ac$           & Lift
&\begin{NoHyper}\pageref{lifting:def-lifts}\end{NoHyper} \\[0ex] 
$\Ac$           & Parameter from probabilistic hypothesis
&\begin{NoHyper}\pageref{hypothesis:probabilistic}\end{NoHyper} \\[0ex] 
$\A[N]$         & Solution in Cartesian coordinates
&\begin{NoHyper}\pageref{intro:eq-system-A-B-Kil}\end{NoHyper} \\[0ex]   
$A_{\leq N}$         & Solution in Cartesian coordinates with truncated data
&\begin{NoHyper}\pageref{intro:eq-system-A-B-N}\end{NoHyper} \\[0ex] 
$\ad$           & adjoint map 
&\begin{NoHyper}\pageref{prelim:eq-adjoint-map}\end{NoHyper} \\[0ex] 
$\bc$           & Lift
&\begin{NoHyper}\pageref{lifting:def-lifts}\end{NoHyper} \\[0ex] 
$\Bc$           & Parameter from pre-modulation hypothesis
&\begin{NoHyper}\pageref{hypothesis:pre}\end{NoHyper} \\[0ex] 
$\Bcin$           & Size of initial modulation operators
&\begin{NoHyper}\pageref{modulation:def-bin}\end{NoHyper} \\[0ex] 
$\B[N]$         & Solution in Cartesian coordinates
&\begin{NoHyper}\pageref{intro:eq-system-A-B-Kil}\end{NoHyper} \\[0ex]   
$B_{\leq N}$         & Solution in Cartesian coordinates with truncated data
&\begin{NoHyper}\pageref{intro:eq-system-A-B-N}\end{NoHyper} \\[0ex] 
$\BBA$          & Event on which Bourgain-Bulut estimate is satisfied
&\begin{NoHyper}\pageref{main:def-Bourgain-Bulut}\end{NoHyper} \\[0ex]
$\beta$         & Inverse temperature 
&\begin{NoHyper}\pageref{intro:eq-Gibbs}\end{NoHyper} \\[0ex]  
$\Cprod{\gamma_1}{\gamma_2}$    & H\"{o}lder space 
&\begin{NoHyper}\pageref{prelim:def-hoelder}\end{NoHyper} \\[0ex] 
$\Cf^{(\Nscript)},\Cf^{(\Ncs)}$           & Covariance functions 
&\begin{NoHyper}\pageref{ansatz:def-Killing}\end{NoHyper} \\[0ex]
$\CNbd, \CNcsbd$           & Covariance functions corresponding to frequency-boundary
&\begin{NoHyper}\pageref{jacobi:def-cf} \end{NoHyper} \\[0ex]
$\Cas$          & Casimir 
&\begin{NoHyper}\pageref{prelim:def-casimir}\end{NoHyper} \\[0ex] 
$\CHHLN$          & Cartesian high$\times$high$\rightarrow$low-matrices  
&\begin{NoHyper}\pageref{jacobi:def-chhl}\end{NoHyper} \\[0ex]
$\chi$ & Cut-off function 
&\begin{NoHyper}\pageref{prelim:def-cut-off}\end{NoHyper} \\[0ex] 
$\chi_K, \chi_M, \chi_{K,M}$ & Frequency-truncated cut-off functions
&\begin{NoHyper}\pageref{ansatz:def-cutoff-frequency-truncated}\end{NoHyper} \\[0ex] 
$\Cut,\Cuttilde$ & Sets of cut-off functions 
&\begin{NoHyper}\pageref{prelim:def-cut-off}\end{NoHyper} \\[0ex] 
$\Dc$           & Parameter from pre-modulation hypothesis 
&\begin{NoHyper}\pageref{prelim:eq-parameter-regularities}\end{NoHyper} \\[0ex] 
$\Der$          & Discrete derivative   
&\begin{NoHyper}\pageref{chaos:def-derivative}\end{NoHyper} \\[0ex] 
$\delta,\delta_0,\delta_1,\delta_2,\delta_3,\delta_4,\delta_5$   & Small parameters
&\begin{NoHyper}\pageref{prelim:eq-parameter-sizes}\end{NoHyper}\\[0ex] 
$\eta$          & Regularity parameter 
&\begin{NoHyper}\pageref{prelim:eq-parameter-regularities}\end{NoHyper} \\[0ex] 
$\Wfuv[\gamma_1][\gamma_2]$          & Function space for chaos estimates with frequency constraints 
&\begin{NoHyper}\pageref{killing:def-FW}\end{NoHyper} \\[0ex] 
$\frkG$         & Lie group         
&\begin{NoHyper}\pageref{section:introduction-result}\end{NoHyper}\\[0ex] 
$\frkg$         & Lie algebra       
&\begin{NoHyper}\pageref{section:introduction-result}\end{NoHyper}\\[0ex] 
$\Gain$         & Gain 
&\begin{NoHyper}\pageref{ansatz:def-Gains}\end{NoHyper} \\[0ex]
$\Haar$         & Haar measure  
&\begin{NoHyper}\pageref{lifting:def-Haar}\end{NoHyper} \\[0ex]
$\HHLN$         & High$\times$high$\rightarrow$low-matrices 
&\begin{NoHyper}\pageref{ansatz:def-hhl}\end{NoHyper}\\[0ex] 
$\HHLNErr$         & High$\times$high$\rightarrow$low-errors
&\begin{NoHyper}\pageref{ansatz:def-hhlerr}\end{NoHyper}\\[0ex]
$\EIN$          & Energy increment 
&\begin{NoHyper}\pageref{increment:def-energy-increment}\end{NoHyper}\\[0ex] 
$\JcbNErr$         & Jacobi errors
&\begin{NoHyper}\pageref{ansatz:def-jacobi-errors}\end{NoHyper}\\[0ex]
$\JcbNErr[][\dagger]$         & Modified Jacobi errors
&\begin{NoHyper}\pageref{jacobi:def-modified-jacobi}\end{NoHyper}\\[0ex]
$\Kil$          & Killing map 
&\begin{NoHyper}\pageref{prelim:eq-Killing}\end{NoHyper}\\[0ex] 
$\coup$         & Temperature parameter
&\begin{NoHyper}\pageref{intro:eq-coup}\end{NoHyper} \\[0ex] 
$\Lambda,\LambdaR$         & Lattices
&\begin{NoHyper}\pageref{prelim:def-lattice-partition}\end{NoHyper} \\[0ex] 
$\LON$         & Low-frequency terms
&\begin{NoHyper}\pageref{ansatz:def-lo}\end{NoHyper}\\[0ex] 
$\LWPD$         & Event on which local well-posedness holds 
&\begin{NoHyper}\pageref{main:def-lwp}\end{NoHyper}\\[0ex] 
$\mu^{(\coup)}$ & Gibbs measure at temperature $8\coup$
&\begin{NoHyper}\pageref{prelim:def-Gibbs}\end{NoHyper}\\[0ex]
$\muR$ & Periodic Gibbs measure at temperature $8\coup$
&\begin{NoHyper}\pageref{prelim:def-Gibbs}\end{NoHyper}\\[0ex]
$\muNR_t$ & Time-evolved, periodic Gibbs measure at temperature $8\coup$
&\begin{NoHyper}\pageref{structure:def-Gibbs-time-evolved}\end{NoHyper}\\[0ex]  
$\mu_\beta$     & Gibbs measure at inverse temperature $\beta$
&\begin{NoHyper}\pageref{intro:eq-Gibbs-AB-rigorous}\end{NoHyper} \\[0ex]
$\muGb$     & Gibbs measure for $\frkG$-valued maps at inverse temperature $\beta$
&\begin{NoHyper}\pageref{lifting:def-Gibbs}\end{NoHyper}\\[0ex]
$\Nd$       & Frequency-truncation parameter for initial data 
&\begin{NoHyper}\pageref{ansatz:eq-WRNd}\end{NoHyper}\\[0ex]
$\Nlarge$   & Large dyadic integer 
&\begin{NoHyper}\pageref{prelim:eq-dyadiclarge}\end{NoHyper}\\[0ex]
$P_{\leq N}$    & Littlewood-Paley projection to frequencies $\leq N$
&\begin{NoHyper}\pageref{prelim:eq-P-N}\end{NoHyper}\\[0ex] 
$\Pbd$          & Littlewood-Paley projection to frequencies $N^{1-2\delta_1}\leq \cdot \leq N$
&\begin{NoHyper}\pageref{prelim:eq-pbd}\end{NoHyper}\\[0ex]
$\Per$             & Perturbations of initial data
&\begin{NoHyper}\pageref{main:def-perturbations}\end{NoHyper}\\[0ex]
$\Perpm$             & Perturbations of initial modulation operators and nonlinear remainders
&\begin{NoHyper}\pageref{main:def-perturbations-null}\end{NoHyper}\\[0ex]
$\PIN$          & Perturbative interactions
&\begin{NoHyper}\pageref{ansatz:def-perturbative}\end{NoHyper}\\[0ex] 
$ \PHA$         & Event on which probabilistic hypothesis is satisfied
&\begin{NoHyper}\pageref{killing:def-probabilistic-hypothesis}\end{NoHyper}\\[0ex]
$ \psix,\psiRx$         & Lattice partition
&\begin{NoHyper}\pageref{prelim:def-lattice-partition}\end{NoHyper}\\[0ex]
$\Renorm[\Nscript],\Renorm[\Ncs]$   & Killing-renormalization 
&\begin{NoHyper}\pageref{ansatz:def-Killing}\end{NoHyper}\\[0ex]   
$\Renormbd[\Nscript],\Renormbd[\Ncs]$   & Killing-renormalization corresponding to frequency-boundary
&\begin{NoHyper}\pageref{jacobi:def-cf}\end{NoHyper}\\[0ex]
$\RenNErr$   & Renormalization error
&\begin{NoHyper}\pageref{ansatz:def-renormalization-error}\end{NoHyper} \\[0ex]   
$r$             & Regularity parameter 
&\begin{NoHyper}\pageref{prelim:eq-parameter-regularities}\end{NoHyper} \\[0ex] 
$\scrr$         & Regularity parameter 
&\begin{NoHyper}\pageref{prelim:eq-scrr}\end{NoHyper} \\[0ex] 
$\rho_{\leq N}$ & Symbol of Littlewood-Paley projection 
&\begin{NoHyper}\pageref{prelim:eq-rho-leqN}\end{NoHyper}\\[0ex]
$\scalebox{0.8}{$\sum\limits^{\hspace{0.15ex}\scalebox{0.6}{$\bullet$}}$}$        & Sum over large dyadic integers  
&\begin{NoHyper}\pageref{prelim:eq-sumlarge}\end{NoHyper} \\[0ex] 
$\Schwartz$        &Schwartz-type norm  
&\begin{NoHyper}\pageref{prelim:eq-Schwartz}\end{NoHyper} \\[0ex] 
$\Scaling_\kappa$   & Scaling transformation 
&\begin{NoHyper}\pageref{ansatz:eq-scaling-transform}\end{NoHyper}    \\[0ex] 
$\SN[][\pm][]$         & Modulation operators
&\begin{NoHyper}\pageref{ansatz:def-pure}\end{NoHyper}\\[0ex]  
$\pSN[][\pm][]$         & Pure modulation operators
&\begin{NoHyper}\pageref{ansatz:def-pure}\end{NoHyper}\\[0ex]  
$\SEN[][u][],\SEN[][v][]$         & Structural errors 
&\begin{NoHyper}\pageref{ansatz:def-structural-error}\end{NoHyper}\\[0ex]  
$\SHHLN[][]$         & Simplified high$\times$high$\rightarrow$low-interactions 
&\begin{NoHyper}\pageref{ansatz:def-shhl}\end{NoHyper}\\[0ex]  
$s$         & Regularity parameter 
&\begin{NoHyper}\pageref{prelim:eq-parameter-regularities}\end{NoHyper} \\[0ex] 
$\varsigma$         & Regularity parameter 
&\begin{NoHyper}\pageref{prelim:eq-parameter-regularities-less}\end{NoHyper} \\[0ex] 
$U_{\leq N}$         & Solution in null-coordinates with truncated data
&\begin{NoHyper}\pageref{intro:eq-system-U-V-N-data}\end{NoHyper} \\[0ex]
$\UN[][]$       &   Solution in null-coordinates          
&\begin{NoHyper}\pageref{ansatz:eq-wave-maps}\end{NoHyper}    \\[0ex] 
$\UN[][+]$      &      Modulated linear wave            
&\begin{NoHyper}\pageref{ansatz:def-modulated-linear}\end{NoHyper}                   \\[0ex]
$\IUN[][+]$      &      Integrated modulated linear wave            
&\begin{NoHyper}\pageref{ansatz:def-integrated-modulated-waves}\end{NoHyper}                   \\[0ex]
$\UN[][+-]$       &    Modulated bilinear wave 
&   \begin{NoHyper}\pageref{ansatz:def-modulated-bilinear}\end{NoHyper}    \\[0ex]
$\UN[][-]$       &     Reversed modulated linear wave              
&  \begin{NoHyper}\pageref{ansatz:def-modulated-linear-reversed}\end{NoHyper}     \\[0ex] 
$\UN[][+\fs]$       &    Mixed modulated object               
&\begin{NoHyper}\pageref{ansatz:def-mixed}\end{NoHyper}     \\[0ex]
$\UN[][\fs-]$       &    Mixed modulated object                 
&\begin{NoHyper}\pageref{ansatz:def-mixed}\end{NoHyper}     \\[0ex]
$\UN[][\fs ]$       &      Nonlinear remainder            
&\begin{NoHyper}\pageref{ansatz:eq-UN-decomposition}\end{NoHyper}     \\[0ex]
$ \UN[][\fsc]$       &      Combined smooth remainder           
&\begin{NoHyper}\pageref{ansatz:def-combined}\end{NoHyper}     \\[0ex]
$V_{\leq N}$         & Solution in null-coordinates with truncated data
&\begin{NoHyper}\pageref{intro:eq-system-U-V-N-data}\end{NoHyper} \\[0ex] 
$\VN[][]$       &    Solution in null-coordinates             
&\begin{NoHyper}\pageref{ansatz:eq-wave-maps}\end{NoHyper}    \\[0ex]        
$\VN[][-]$      &     Modulated linear wave             
&   \begin{NoHyper}\pageref{ansatz:def-modulated-linear}\end{NoHyper}                   \\
$\IVN[][-]$      &      Integrated modulated linear wave            
&\begin{NoHyper}\pageref{ansatz:def-integrated-modulated-waves}\end{NoHyper}                   \\[0ex]
$\VN[][+-]$       &   Modulated bilinear wave 
&\begin{NoHyper}\pageref{ansatz:def-modulated-bilinear}\end{NoHyper}      \\[0ex] 
$\VN[][+]$       &     Reversed modulated linear wave        
&\begin{NoHyper}\pageref{ansatz:def-modulated-linear-reversed}\end{NoHyper}     \\[0ex] 
$\VN[][\fs-]$       &    Mixed modulated object             
&\begin{NoHyper}\pageref{ansatz:def-mixed}\end{NoHyper}     \\[0ex]
$\VN[][+\fs]$       &     Mixed modulated object                
&\begin{NoHyper}\pageref{ansatz:def-mixed}\end{NoHyper}     \\[0ex]
$\VN[][\fs ]$       &      Nonlinear remainder            
&\begin{NoHyper}\pageref{ansatz:eq-VN-decomposition}\end{NoHyper}     \\[0ex] 
$ \VN[][\fcs]$       &      Combined smooth remainder      
&\begin{NoHyper}\pageref{ansatz:def-combined}\end{NoHyper}     \\[0ex]
$\vartheta$         & Regularity parameter 
&\begin{NoHyper}\pageref{prelim:eq-parameter-regularities-less}\end{NoHyper}     \\[0ex]
$W_0^{(\coup)},W_1^{(\coup)}$       &      $\frkg$-valued white noise at temperature $8\coup$ 
&\begin{NoHyper}\pageref{prelim:def-white-noise-g-valued}\end{NoHyper}  \\[0ex] 
$W^{(\coup),\pm}$       &      $\frkg$-valued white noise at temperature $\coup$        
&\begin{NoHyper}\pageref{prelim:def-white-noise-g-valued}\end{NoHyper} \\[0ex] 
$\Wuv[\gamma_1][\gamma_2][]$        & Function space for chaos estimates 
&\begin{NoHyper}\pageref{chaos:def-wc}\end{NoHyper}\\[0ex] 
$\WassersteinR$     & Wasserstein metric 
&\begin{NoHyper}\pageref{prelim:def-Wasserstein}\end{NoHyper}
\end{longtblr}
}

\end{appendix}

\bibliography{Wave_Maps_Library}
\bibliographystyle{myalpha}
\nopagebreak

\end{document}